\documentclass[reqno]{amsart}
\address{Department of Mathematics,  Stanford University,
Stanford, California 94305, USA}
\email{abouzaid@stanford.edu}
\address{Yau Mathematical Sciences Center, Tsinghua University 
and Beijing Institute of Mathematical Sciences and Applications,  Beijing, China}
\email{fukayakenji@mail.tsinghua.edu.cn}
\address{Center for Geometry and Physics, Institute for Basic Sciences (IBS), Pohang, Korea
\\
Department of Mathematics,
POSTECH, Pohang, Korea}
\email{yongoh1@postech.ac.kr}
\address{Graduate School of Mathematics,
Nagoya University, Nagoya, Japan}
\email{ohta@math.nagoya-u.ac.jp}
\address{Research Institute for Mathematical Sciences, 
Kyoto University, Kyoto, Japan}
\email{ono@kurims.kyoto-u.ac.jp}

\usepackage[all]{xy}
\usepackage{euscript}
\usepackage{comment}
\usepackage{yhmath}
\usepackage{mathrsfs}
\usepackage{bbding}
\usepackage[dvipdfmx]{graphicx}
\usepackage{amssymb}
\usepackage{amscd}
\usepackage{amsbsy} %used for making boldface greek letters
\usepackage{color}
\usepackage[normalem]{ulem}% used for \uwave{K}

\newcommand{\blue}{\textcolor{black}}

\newcommand{\newred}{\textcolor{black}}

\DeclareMathOperator*{\tensor}{\otimes} %Used to make subscript
\usepackage[original]{imakeidx}
\makeindex % will get the same name as the main file
\makeindex[name=syindex, title=List of notations]

\begin{document}

\def\R{\ifmmode{\mathbb R}\else{$\mathbb R$}\fi} %real numbers
\def\Q{\ifmmode{\mathbb Q}\else{$\mathbb Q$}\fi} %rational numbers
\def\C{\ifmmode{\mathbb C}\else{$\mathbb C$}\fi} %complex numbers
\def\Z{\ifmmode{\mathbb Z}\else{$\mathbb Z$}\fi} %integers
\def\F{\ifmmode{\mathbb F}\else{$\mathbb F$}\fi} %field
\def\N{\ifmmode{\mathbb N}\else{$\mathbb N$}\fi} %natural numbers

\theoremstyle{theorem}
\newtheorem{thm}{Theorem}[section]
\newtheorem{cor}[thm]{Corollary}
\newtheorem{lem}[thm]{Lemma}
\newtheorem{claim}[thm]{Claim}
\newtheorem{sublem}[thm]{Sublemma}
\newtheorem{prop}[thm]{Proposition}
\newtheorem{ax}[thm]{Axiom}
\newtheorem{conj}[thm]{Conjecture}

\theoremstyle{definition}
\newtheorem{defn}[thm]{Definition}
\newtheorem{rem}[thm]{Remark}
\newtheorem{ques}[thm]{Question}
\newtheorem{cons}[thm]{\rm\bfseries{Construction}}
\newtheorem{asump}[thm]{\rm\bfseries{Assumption}}
\newtheorem{exm}[thm]{Example}
\newtheorem{conds}[thm]{Condition}

\newtheorem{notation}[thm]{\rm\bfseries{Notation}}

\newtheorem*{thm*}{Theorem}

\numberwithin{equation}{section}

\def\osc{{\hbox{\rm osc}}}
\def\Crit{{\hbox{Crit}}}
\def\tr{{\hbox{\bf tr}}}
\def\Tr{{\hbox{\bf Tr}}}
\def\DTr{{\hbox{\bf DTr}}}

% Organising principle for notations:
% "standard" 
% mathbf for collection of objects (Lagrangians, ...)
% \mathfrak for "Floer theoretic data" (PO, bounding cochains)
% \EuScript for categories (also for bi-modules?)
% mathcal for modules, bi-modules, functors. Also for moduli spaces
% Commands added by M.A (14 April 2011)
\newcommand{\Fuk}{\EuScript F}
\newcommand{\Cat}{{\EuScript C}}
\newcommand{\Dat}{{\EuScript D}}
\newcommand{\Bat}{{\EuScript B}}
\newcommand{\Aat}{{\EuScript A}}
\newcommand{\modu}{\operatorname{mod}}
\newcommand{\Catbi}{\Cat\mbox{-}\! \modu \!\mbox{-} \Cat}
\newcommand{\Ob}{\operatorname{Ob}}
\newcommand{\Hom}{\operatorname{Hom}}
\newcommand{\Spin}{\operatorname{Spin}}
\newcommand{\cP}{{\mathcal P}}
\newcommand{\cQ}{{\mathcal Q}}
\newcommand{\cY}{{\mathcal  Y}}
\newcommand{\cZ}{{\mathcal  Z}}
\newcommand{\field}{\mathbf k}
\newcommand{\fm}{\mathfrak m}
\def\co{\colon\thinspace}
% Commands added by M.A (18 April 2011)
\newcommand{\cL}{{\mathcal  L}}
\newcommand{\fb}{\mathfrak b}
\newcommand{\PO}{{\mathfrak{PO}}_{\fb}}
\newcommand{\wq}{\widehat{\mathfrak q} }
\newcommand{\id}{\operatorname{id}}
% Commands added by M.A (26 July 2011)
\newcommand{\cX}{{\mathcal  X}}
\newcommand{\Tot}{{\mathcal  T}}
\newcommand{\Mat}{{\mathrm{Mat}}}
\newcommand{\Tw}{{\mathrm{Tw}}}
\newcommand{\unit}{{\mathbf e}}
\newcommand{\Po}{{\mathfrak{PO}}}
%Commands added by M.A (2 Sep 2011)
%Collections of Lagrangians
\newcommand{\bL}{{\mathbf  L}}
\newcommand{\bU}{{\mathbf  U}}
% Categories with these objects
\newcommand{\sL}{{\EuScript L}}
\newcommand{\sU}{{\EuScript  U}}

%Commands added by M.A (29 Sep 2011)
\newcommand{\liminv}{\varprojlim}
\newcommand{\Iota}{\iota}
\newcommand{\bfx}{{\mathbf x}}
\newcommand{\bfy}{{\mathbf y}}
\newcommand{\hCat}{\underleftarrow{\Cat}}
\newcommand{\cS}{{\mathcal  S}}
\newcommand{\superscript}[1]{\ensuremath{^{\textrm{#1}}} }
\newcommand{\subscript}[1]{\ensuremath{_{\textrm{#1}}} }
\renewcommand{\th}[0]{\superscript{th}}
\newcommand{\bfa}{{\mathbf a}}
\newcommand{\bfb}{{\mathbf b}}
\newcommand{\bfc}{{\mathbf c}}
\newcommand{\bfdelta}{{\boldsymbol \delta}}
%Commands added by M.A (14 Oct 2011)
\newcommand{\bfd}{{\mathbf d}}

%Commands added by M.A (2 Dec 2011)
\newcommand{\hBat}{\underleftarrow{\Bat}}
\newcommand{\Nat}{N}
\newcommand{\Bhat}{\hat{B}}
%Commands added by M.A (22 Dec 2011)
\newcommand{\Cone}{\operatorname{Cone}}
%Commands added by M.A (08 Dec 2016)
\newcommand{\coev}{\operatorname{Coev}}
\newcommand{\fq}{\mathfrak{q}}
\newcommand{\fp}{\mathfrak{p}}
\newcommand{\fqhat}{\hat{\fq}}
\newcommand{\fphat}{\hat{\fp}}

\makeatletter
\providecommand*{\cupdot}{%
  \mathbin{%
    \mathpalette\@cupdot{}%
  }%
}
\newcommand*{\@cupdot}[2]{%
  \ooalign{%
    $\m@th#1\cup$\cr
    \hidewidth$\m@th#1\cdot$\hidewidth
  }%
}
\makeatother

\makeatletter
\providecommand*{\capdot}{%
  \mathbin{%
    \mathpalette\@capdot{}%
  }%
}
\newcommand*{\@capdot}[2]{%
  \ooalign{%
    $\m@th#1\cap$\cr
    \hidewidth$\m@th#1\cdot$\hidewidth
  }%
}
\makeatother

\title[Quantum cohomology
and split generation]
{Quantum cohomology
and split generation
in Lagrangian Floer theory}

\author[M. Abouzaid, K. Fukaya, Y.-G. Oh, H. Ohta, K.
Ono]{Mohammed Abouzaid, Kenji Fukaya, Yong-Geun Oh, Hiroshi
Ohta, Kaoru Ono}
\thanks{MA was partially supported by a Clay Research Fellowship, by NSF grants DMS-1308179, DMS-1564172, and DMS-1609148, and by the Simons Foundation grant 385571, KF is supported
partially by JSPS Grant-in-Aid for Scientific Research
Nos. 30165261, 23224002, Global COE Program G08, Simons collaboration on homological Mirror symmetry, NSF Grant No. 1406423. 
and 2024 National Recruitment Program for High-level Overseas Talents, the People's Republic of China, YO by
IBS project IBS-R003-D1
HO by JSPS Grant-in-Aid
for Scientific Research Nos.19340017, 23340015, 15H02054, 21K18576, 21H00983, 23K20796, 23K20798, 26K00605
and KO by JSPS Grant-in-Aid
for
Scientific Research, Nos. 17654009, 18340014,21244002, 26247006, 19H00636, 24H00182.
KF thanks to Simons Center for Geometry and Physics where many of his research for this paper is performed.
KO is also grateful for National Center for Theoretical Sciences, Taiwan 
and Czech Academy of Sciences. 
The authors would like to thank various institutions, 
such as Simons Center for Geometry and Physics, IBS center for Geometry and Physics, and Research Institute for Mathematical Sciences, Kyoto University, 
which provided them the opportunities of discussions on this work.}

\begin{abstract} Given a finite collection of Lagrangian
submanifolds $\mathscr L$ in a compact symplectic manifold $X$, 
we construct a cyclic, filtered, strictly unital curved $A_{\infty}$ category $\mathcal L$
and develop Floer theory of closed-open maps and open-closed maps.
Using them, we prove that, whenever the map from the quantum cohomology of $X$
to the Hochschild cohomology of the Fukaya category $\mathcal L$ with objects
$\mathscr L$ is injective, the following consequences follow: (1) any
other Lagrangian submanifold equipped with a weak bounding cochain lies in the category
split-generated by $\mathscr L$, and (2) the Hochschild homology and
cohomology of the Fukaya category are isomorphic to quantum
cohomology.
In the exact case a similar result was obtained in
\cite{abouzaid:IHES}.
We also provide some applications.
\end{abstract}

\date{Revision June 10 2026}


\maketitle
\setcounter{tocdepth}{1}
\newpage
\tableofcontents

\newpage
\section{Introduction}
\label{sec:intro}
In \cite{fooo09,fooo092}, the last four authors studied the Floer theory of embedded Lagrangians $L$ in a general closed symplectic manifold $(X,\omega)$.  One of their main results takes the following form in a simplified setting: if $L$ is equipped with a $\mathrm{Spin}$ structure then %each choice of $\omega$-compatible almost complex structure $J$ on $M$ determines a 
Floer cohomology group $ HF^*(L,b)$ is defined for each element
 %\marginpar{Related to Mahammed's comment.  Floer homology 
%is independent of $J$ and perturbation. So I modified the sentence. KF. 2025 July}
\begin{equation} \label{eq:bounding_cochain}
b \in   \widehat{\mathcal M}(L;\Lambda_0) \subset H^{\rm odd}(L; \Lambda_0)
\end{equation}
of the set of \emph{bounding cochains}\footnote{The Floer theory associated to $b$ only depends on the \emph{gauge equivalence class} of this element in the Maurer-Cartan moduli space associated to $L$.
The quotient space of  $\widehat{\mathcal M}(L;\Lambda_0)$ by gauge equivalence 
is independent of $J$ and perturbations.}\footnote{We include the case when $b$ is a weak bounding cochain, 
that is, the case Maurer-Cartan equation is satisfied modulo unit.} associated to $L$, an $\omega$-compatible almost 
complex structure $J$ and perturbations, where $\Lambda_0$ denotes the Novikov ring over a field of characteristic $0$ (see the discussion of notation in 
Section \ref{Sec3Not}). The basic idea for understanding the space of bounding cochains is that the moduli space of $J$-holomorphic disks in $M$ with boundary on $L$, equipped  with one boundary marked point, defines an even dimensional cochain on $L$ which obstructs the equation $d^2 = 0$ when considering the Floer differential. The space $ H^{\rm odd}(L; \Lambda_0)$ parametrises ways of deforming the Floer differential, and the space of bounding cochains is the subset of this parameter space consisting of choices for which the resulting obstruction vanishes, yielding a well-defined cohomology group.

The Floer cohomology group was shown to be the underlying graded vector space of an $A_\infty$ algebra, and an $A_\infty$ bi-module  $CF^*((L_0,b_0),(L_1,b_1))$ 
with homology 
\begin{equation} \label{eq:Floer_cohomology-pair}
  HF^*((L_0,b_0),(L_1,b_1))
\end{equation}
was associated to each pair of Lagrangians equipped with bounding cochains.\footnote{Here $\Lambda_0$ is a Novikov ring over ground field $\F$ which is 
either $\C$ or $\R$.}
 %\marginpar{I comment out 
%`this bi-module was denoted $C(L_0,L_1; \Lambda_0)$, but we change the notation to the one used in this paper'. FOOO wrote it as 
%$CF^*((L_0,b_0),(L_1,b_1))$.  Since there is still nontrivial boundary I do not want to use $HF$.  KF 2025 Aug.}

The main goal in \cite{fooo09,fooo092} was the study of applications of Floer cohomology to symplectic topology and the mirror symmetry, so the focus was to formulate and prove its invariance properties, as well as to develop computational methods. Among aspects of the theory which were left for subsequent work was the construction of a category, usually referred to as the Fukaya category, whose objects are Lagrangians equipped with (weak) bounding cochains, and whose morphism spaces are given by Equation \eqref{eq:Floer_cohomology-pair}. The construction of an $A_{\infty}$ algebra associated to an embedded Lagrangian was generalized by Akaho and Joyce \cite{AkahoJoyce} to the case of immersed Lagrangians, relying the foundations of \cite{fooo09,fooo092}. In case one considers a finite collection of transversely intersecting Lagrangian submanifolds, the construction of an $A_{\infty}$ category is a consequence of 
this generalization, and one can also extract it from the work of the second author \cite{fukaya:functor}.

Starting from \cite{fooo:toric1}, the last four authors have rebuilt the foundations of Lagrangian Floer theory, replacing the use of singular cochains in \cite{fooo09,fooo092} by differential forms. While largely motivated by the use of averaging arguments with respect to Hamiltonian torus actions in \cite{fooo:toric1}, the outcome is a simplified construction of Floer homology, ultimately because the chain-level functoriality of pushforward and pullback maps, which is used to construct operations in Floer theory, is easier to establish in the de Rham than in the singular theory.%\marginpar{de-Rham is changed to de Rham in many places.
%\bbblue{2026 May 29}}

This paper constructs the $A_\infty$ category of Lagrangians using the de Rham model, and establishes a fundamental \emph{generation criterion} which is used to compute it, and which has its origin in the first author's work \cite{abouzaid:IHES} studying Floer theory in the technically much simpler setting of exact Lagrangians. Later in this introduction, we shall give a precise statement of the generation criterion, but it is useful at this stage to comment about what it means to compute a category: naively, this would entail a list of all (quasi)-isomorphism classes of objects, of the cohomology group of morphism spaces between them, and of all compositions. Unfortunately, such a task is in general hopeless even for problems arising in linear algebra \cite{Drodz}.

Instead, the sense in which we compute the $A_\infty$ categories arising from Lagrangian Floer theory is that we establish conditions under which we can identify a fully faithful embedding in a category of modules over some $A_\infty$ algebra, i.e. one can assign to each Lagrangian equipped with a bounding cochain a module over this algebra, and an isomorphism from the Lagrangian Floer cohomology for pairs to the space of morphisms between the corresponding modules, in such a way that compositions agree. In fact, the result we prove establishes both that the algebra is the $A_\infty$ algebra associated to a finite collection of Lagrangians, and that the modules we obtain are perfect, i.e. finitely built via extensions and retractions from free ones. Once an embedding in a category of modules has been established, one can translate (some) questions about Lagrangians to questions about modules, for which tools of representation theory or algebraic geometry can be used. For example, our result is applicable to all compact toric manifolds (see Theorem \ref{maintheorem4}), and we expect that this will be used to prove the appropriate version of Homological Mirror Symmetry in this context. We also discuss some applications to quantitative symplectic topology in Section \ref{sec: quatitative}.

In the remainder of this introduction, we give a precise statement of the main results of the paper, compare our results with related ones in the literature, and finally give a summary of the results established in each section.

% Since the first appearance of \cite{fooo092}, an extensive literature has constructed Fukaya categories in restricted contexts, starting with Seidel's work for closed exact Lagrangians in \cite{Sei06}, but also in different contexts, such as in the first author's work with Seidel \cite{abouzaidSeidel} on Lagrangians with Legendrian boundary.

% The starting point of this paper is the construction of the long-awaited Fukaya $A_\infty$-category associated to any closed symplectic manifold, with objects supported on embedded Lagrangians. 

% We construct this category in this paper.

% In this paper, we complete the formal part

% In this paper, we shall denote such a Floer cohomology group $HF^*(L;b) $ 

% The key construction of \cite{fooo092} is of a Floer cohomology group $HF^*(L;b)$ for each Lagrangian equipped with a bounding cochain.  The cohomology  that the space $M_{\text{weak}}(L;\Lambda_0) $ 

% $\fb \in H^{even}(X;\Lambda_0)$ a bulk class,
% ${\rm st} \in H^2(X;\Z_2)$  
% a background class,
% {
% $\sL$  a finite set of ${\rm st}$-relatively spin Lagrangian submanifolds}, and $\bL$ a finite collection whose elements are pairs
% $(L,b)$, with $L \in \sL$, and $b$ a class in $ H^{odd}(L, \Lambda_0)$  for which 
% \begin{equation}
% [b] \in \mathcal M_{\text{weak}}(L;\fb;\Lambda_0)
% \end{equation}
% is a weak bounding cochain (the right hand side is defined in Remark \ref{MCoddelementrem}).
% We fix \rcng{an $\omega$-compatible} almost complex structure $J$ and use it 
% throughout the paper.

\subsection*{Summary of the results}

As we explained earlier, the most straightforward constructions in Lagrangian Floer theory in general yield obstructed complexes. It  has been proved to be beneficial to consider deformations and twists of Floer complexes parametrised by geometric data on the ambient symplectic manifold, because this makes it more likely that one can find an unobstructed Floer complex. To this end, we fix for the remainder of this paper a pair of cohomology classes\index[syindex]{sct@${\rm st}$}
\begin{equation}
(  \fb, {\rm st})  \in H^{\rm even}(X;\Lambda_0) \times  H^2(X;\Z_2),
\end{equation}
called the \emph{bulk}\index{bulk class} and \emph{background} classes\index{background class}, which can be thought of as respectively parametrising ways of changing the area and sign contributions of moduli spaces of holomorphic curves in Floer theory. 
More precisely we fix a differential form representing $\fb$.\footnote{If we change $\fb$ by an element of $2\pi\sqrt{-1} H^{2}(X;\Z)$ the resulting theory does not change. See Remark \ref{rem611}.}
 %\marginpar{footnote is added. KF 2025 July}
Here $\Lambda_0$ is the Novikov ring (See (\ref{formula114}).) over the ground ring $\F$ (which is $\R$ or $\C$ in this paper.)

We fix an $\omega$-compatible almost complex structure $J$, that we use.
(The filtered $A_{\infty}$ category we obtain is indepenent of $J$ as will be proved in    Theorem  \ref{thm111}.)
 %\marginpar{Changed according to MA's comment.  KF 2025 July.}
We fix as well as a finite set $\mathscr L$ of immersed
 %\marginpar{immersed is added.  KF 2025 July} 
Lagrangian submanifolds of $X$  such that the summand $\fb_2$ of $\fb$ lying in $  H^{2}(X;\F)$ vanishes upon restriction to each element $L \in \mathscr L$ (more precisely, we choose a primitive for this restriction),\index[syindex]{Lscr@$\mathscr L$}
 and each of which is equipped with an ${\rm st}$-relatively $\Spin$ structure (see Definition \ref{def:relative_Spin}); this datum determines an orientation of the moduli spaces of holomorphic curves with these boundary conditions, and implies in particular that the restriction of ${\rm st}$ to each such Lagrangians agrees with the second Stiefel-Whitney class. Associated to a Lagrangian $L$ equipped with this data is the set 
\begin{equation}
\widehat{\mathcal M}_{\rm weak}(L;\fb;\Lambda_0) \subset H^{\rm odd}(L; \Lambda_0)
\end{equation}
of \emph{weak bounding cochains}\index{bounding cochain}\index{weak bounding cochain}, which is defined in Definition \ref{boundingcochain}.  See also Remark \ref{MCoddelementrem}. This implicitly depends on the background class ${\rm st}$, the choice of relative spin structure, the choice of primitive for $\fb_2$,\index[syindex]{MweakLbhat@$\widehat{\mathcal M}_{\rm weak}(L;\fb;\Lambda_0) $}
the $\omega$-compatible almost complex structure $J$ and the perturbations. 
(However after taking gauge equivalence classes it is independent to all of them.)
 %\marginpar{modified following MA's comment. KF 2025 July} 
As discussed after Equation \eqref{eq:bounding_cochain}, a {\it bounding cochain} is one for which the cohomology class $\fm_{0}^{(L,b)}(1)$ associated
to the deformed count of holomorphic disks with one boundary marked point vanishes. The adjective \emph{weak} refers, in this context, to the condition that this cohomology class is a scalar multiple of the unit in $H^0(L; \Lambda_0)$. This scalar is called the \emph{potential value}, denoted $ \PO(L,b)$\index{potential value}; writing $ \unit_{(L,b)}$ for the unit, the above definitions imply that
\begin{equation}
  \fm_{0}^{(L,b)}(1) = \PO(L,b)  \cdot \unit_{(L,b)}.
\end{equation}

{\it Hereafter we sometimes omit the word `weak' from `weak bounding cochain'. So in this paper bounding cochain means weak bounding cochain unless 
otherwise mentioned explicitely.}
 %\marginpar{A sentence in italic is added.  
%the word weak is removed in many places.  KF 2025 March.  }
\index[syindex]{Lcalcan@$\cL$}

Let $\bL$ be a \index[syindex]{Lbd@$\bL$} finite collection consisting of the pairs $(L,b)$,\footnote{We call 
the triple $\theta_L$, relative spin structure of $L$ and a weak bounding cochain $b$ to be  a 
{\it brane datum}.\index{brane datum}} with $L \in \mathscr L$, and $b$ a  bounding cochain. 
We denote by ${\bf b}$ the pair $(\frak b,\{b\})$ where $\frak b$ is the bulk class and $\{b\}$ is the set of weak bounding cochains of various $L \in \mathscr L$. As discussed earlier in the introduction, previous work has indicated that there is an $A_\infty$ category with objects elements of $\bL$, with morphisms between pairs of objects given by the Floer cochain complexes. In Sections \ref{sec:cyclicfil}-\ref{sec:CRperturb} we construct this category, which we denote $\cL^{\bf b}$, and establish some of its basic properties. 
More precisely, the outcome of Section \ref{canonical} is summarised as follows:\index[syindex]{Lcalb@$\sL^{\bf b}$}
\begin{thm}\label{thm11}
The category $\cL^{\bf b}$ is a cyclic, filtered, strictly unital curved $A_{\infty}$ category over the Novikov ring $\Lambda_0$, 
 %\marginpar{Changed coefficient to $\Lambda_0$ KF 2024 Dec} 
in which every object is weakly curvature free, and such that the underlying graded module of morphisms between objects is given by: 
\begin{enumerate}
\item[(i)] the free $\Lambda_0$ modules generated by the intersection point if the underlying Lagrangians are distinct,  
\item[(ii)] the de Rham cohomology of the $L$ with coefficient in $\Lambda_0$ if both objects are supported on $L$.
(See Remark $\ref{flatbundlerem}$.)
\end{enumerate}
\end{thm}
The independence of $\cL^{\bf b}$ of the choices involved (including an almost complex structure) is stated as Theorem  \ref{thm111}.

A few words
 %\marginpar{We can prove uniqueness up to {\bf cyclic} pseudo-isotopy, except $\Pi$ (harmonic projection) dependence. 
%I will add the statement about uniqueness if one can do that remaining part.  KF 2024 Dec.
%Now we can do it. KF 2025 Jan.} 
of clarification are in order: the fact that every object is weakly curvature free is simply a restatement of the fact that all objects that we consider are equipped with weak bounding cochains. The statement that the category is filtered is a basic consequence of the fact that holomorphic curves have positive area, and this is a feature of the category which is largely independent of the choice of model used to construct it. On the other hand, the statement that we construct a category which is simultaneously cyclic and strictly unital depends very strongly on our use de Rham theory: the cyclicity statement amounts to a strict implementation of the compatibility between the operations which we construct with Poincar\'e duality, and the unitality statement asserts that all higher $A_\infty$ operations identically vanish when one of the inputs is the unit. Even though we expect that there will be alternative approaches to our other main results which do not depend on these two properties, we rely very heavily on them in later constructions, and it seems essentially impossible to construct the operation 
which is cyclic and strictly unital directly   via geometry on the singular chain complex, that was used in \cite{fooo09,fooo092}. 
\begin{rem}\label{flatbundlerem} %\marginpar{Remark added.  KF 2024 Dec.}
$  $ \par
\begin{enumerate}
\item
For the purpose of studying bounding cochains which contain an element of $H^1(L;\Lambda_0)$ with 
non-vanishing leading order term $\in H^1(L;\Lambda_{0}/\Lambda_{+}) =H^1(L;\F)$, we will include the case where there is a non-trivial flat $\F$ bundle with a connection form $\theta$ on $L$. 
In other words the holonomy representation of the flat bundle is given by 
$\gamma \in \pi_1(L) \mapsto \exp(\gamma \cap \theta) \in \F^*$. (See Definition \ref{boundingcochain}.) When we consider two objects 
$(L,\theta)$, $(L,\theta')$ (where the Lagrangian $L$  is common) the morphism space is the cohomology 
with local coefficeint $H^*(L;\theta'-\theta)$.  (So the above item (ii) is slightly modified in such a case.) 
Here the closed one form $\theta'-\theta$ is regarded as a connection of an $\F$ bundle associated to the flat 
$\F^*$  bundle\footnote{This is an idea going back to Cho \cite{Cho2}.} and $H^*(L;\theta'-\theta)$ is the cohomology with coefficients 
in this flat bundle.
 %\marginpar{Rewritted a bit 
%following MA's comment.  KF 2025 July.} 
See Subsection \ref{localcoefficient}. 
\item
We include the case an element $L$ of $\mathscr L$ is immersed.  We assume $L$ is self-transversal in the following sense.
We identify $L$ to a pair $(\tilde L,i_L)$ of a smooth manifold $\tilde L$ and an immersion $i_L : \tilde L \to X$.
We require $i_L$ to be {\it self-transversal},\index{self-transversal} that is, the intersection of $(i_L,i_L) : \tilde L \times \tilde L \to X \times X$ with the diagonal 
$\cong X \subset X \times X$, 
is transversal outside the image of the diagonal embedding $\tilde L \subset \tilde L \times \tilde L$.  In case $L$ is immersed the morphism space from $L$ to itself is $H(\tilde L \times_X \tilde L;\Lambda_0)$,
that is, the direct sum of $H(\tilde L;\Lambda_0)$  and the 
free $\Lambda_0$ module with two generators for each self-intersection.
See \cite{AkahoJoyce,fukaya:functor}.
 %\marginpar{(2) is added.  KF 2025 July.}
\item
Even if we start from a finite set of relatively spin Lagrangians which are not necessarily assumed to be weakly unobstructed (i.e., without weak bounding cochains $b$), we can still construct a cyclic filtered strictly unital {\it curved} $A_{\infty}$ category over $\Lambda_0$. The same construction in the current paper works up to Section \ref{sec:unit}.
Then   we skip the process in Section \ref{sec:elicurv} and can move to Section \ref{canonical} to obtain 
curved structure. Here we incorporate {\it tad pole} to construct the canonical model for the curved $A_{\infty}$ category.  See \cite[Subsection 5.4.4]{fooo09}.
\end{enumerate}
\end{rem}
\begin{rem}
We expect that it is possible to equip the singular chain model of the category $\cL$ with the structure of a \emph{proper Calabi-Yau} category, as introduced by Kontsevich and Soibelman \cite{KS-Aoo}, which is a homotopically invariant version of the notion of a cyclic $A_\infty$ category. 
It is proved in  \cite[Theorem 10.2.2]{KS-Aoo} that 
one can apply formal algebraic constructions to replace a proper Calabi-Yau category with 
homology unit over a field of characteristic $0$ by an equivalent strictly unital cyclic category in the case curvature is $0$. To the best knowledge of ours, the same kind of statement is not considered in the literature 
for the curved case and/or in the case of a filtered $A_{\infty}$ category over $\Lambda_0$ coefficients.
We do not try to prove it\footnote{
%It might be difficult to work out the proof of Theorem \ref{maintheorem3.5} 
%in a similar way as this paper if we take 
%that way.
There might be an issue in our situation where we use weak bounding cochain to eliminate the curvature. 
In fact the 
definition of weak bounding cochain involves strict unit.} but use more geometry to prove Theorem \ref{thm11}.
Namely we construct a curved homotopically unital and cyclic filtered $A_{\infty}$ category 
and use algebraic constructions to replace it by an equivalent curved strictly unital and cyclic filtered $A_{\infty}$ category.
\end{rem}
\begin{rem}
Using de Rham model we can construct:
 %\marginpar{Remark 1.3 is rewritten and Remark 1.4 is added.  KF 2025 July.}
\begin{enumerate}
\item a curved unital and cyclic filtered $A_{\infty}$ algebra associated to a single embedded Lagrangian submanifold,
\item a curved cyclic filtered $A_{\infty}$ category associated to a finite set of immersed Lagrangian submanifolds,
\item a curved strictly unital filtered $A_{\infty}$ category associated to a finite set of immersed Lagrangian submanifolds.
\end{enumerate}
We do not know how to perform the construction as (1), (2), (3) above in the chain model other than de Rham model,
(such as singular homology or Morse homology and etc.)  
We can not realize (2) and (3) simultaneously on de Rham chain complex directly by a geometric 
construction.  However we can realize (2) and (3) simultaneously by combining it with algebraic constructions.
\end{rem}
As a final clarification about Theorem \ref{thm11}, observe that we may extract a subcategory $\cL^{\bf b}_{\lambda}$ consisting of those objects $\bL_{\lambda}$ 
 %\marginpar{ $\cL$ is changed to  $\cL^{\bf b}$ in several places.  
%$\sL_{\lambda}$ is changed to $\bL^{\bf b}$. KF 2024 Dec}  
with a prescribed potential value $\lambda \in \Lambda_+$; essentially by definition, $\cL^{\bf b}_{\lambda}$ is an $A_\infty$ category in which all contributions of the curvature $ \fm_{0}^{\bf b}$ to the relations among $A_\infty$ operations vanish. This means that the $A_\infty$ relations still hold if we set the value $\fm_{0}^{\bf b}$ to identically vanish; i.e. we may consider $\cL_{\lambda}^{\bf b}$ as an ordinary 
(uncurved) $A_\infty$ category.

We shall go back and forth between considering $\cL$ as a curved $A_\infty$ category, and considering the direct sum of uncurved categories\index[syindex]{Lcalcanlambda@$\cL_{\lambda}$}
\begin{equation}
  \bigoplus_{\lambda \in \Lambda} \cL^{\bf b}_{\lambda} \subset \cL^{\bf b}.
\end{equation}
The former point of view is of course more natural for Floer theory, but the latter simplifies the homological algebra constructions which we perform. In each context where it is possible to switch between the two points of view, without any loss of information, the underlying principle is that the contributions of terms of mixed curvature vanish; we refer the reader to Section \ref{sec:matr-fact-dimens}, where the most elementary version of a vanishing theorem is proved in a purely algebraic context.

Having limited the set of objects in $\cL$ to a finite collection, it is natural to ask about the interaction of this category with other Lagrangians equipped with bounding cochains. In general, little of substance can be said, as $\mathscr L$ may be the empty collection, which leads us to first consider the relationship of $\cL$ with the symplectic topology of the ambient symplectic manifold. 
%Denoting by {\bf b} the data of the bulk class $\frak b$ and the bounding cochains $b$ that enter in the definition of $\cL^{\rm can}$, 
We are thus led, in Section \ref{sec:open-closed map}, to construct the \emph{open-closed map}
 %\marginpar{section 14 is changed to 15 since 
%open closed map is introduced in Section 15 in the present version. It is possible to change 15 to 14 and open close to close open at the same time.
%KF 2025 Feb.}
\index{open-closed map}
\begin{equation}\label{mapphat(1.3)}
\widehat{\frak p}^{\text{\bf b}} : HH_*(\cL^{\bf b}, \cL^{\bf b}) \to QH^*_{\fb}(X;\Lambda_0)
\end{equation}
from the Hochschild homology of the category $\cL^{\bf b}$ to the bulk deformed quantum cohomology group of $X$ (see for example 
Formula (\ref{closedgwdef})). 
Here we put $\text{\bf b}= (\frak b,b,\theta)$.  The definition of Hochschild homology is recalled in {Subsection} \ref{sec:hochschild-homology}, and its construction uses 
the decomposition of $\cL^{\bf b}$ into curved categories, in the sense that we define it as a direct sum
\begin{equation} \label{eq:HH_defined_by_decomposition}
  HH_*(\cL^{\bf b}, \cL^{\bf b}) \equiv \bigoplus_{\lambda \in \Lambda}  HH_*(\cL^{\bf b}_\lambda, \cL^{\bf b}_\lambda). 
\end{equation}
Intuitively, this group is obtained from $\cL^{\bf b}$ by considering a chain complex whose generators correspond to arranging composable morphisms around the boundary of a disk, equipped with the differential given by applying the $A_\infty$ operations to any cyclically consecutive subsequence. The map $\widehat{\frak p}^{\text{\bf b}}$, which extends the analogously denoted map constructed in  \cite{fooo09,fooo092}, is then obtained by considering holomorphic disks with Lagrangian boundary conditions specified by the composable sequence, equipped with an interior marked point, which defines the evaluation map to $QH^*_{\fb}(X;\Lambda_0)$ (there may be additional marked points to account for the bulk deformation $\frak b$).
\begin{rem} \label{rem:open-closed-terminology}
The motivation for the terminology we use is that $\widehat{\frak p}^{\text{\bf b}}$ is a map from an invariant associated to the open sector controlled by holomorphic curves with Lagrangian boundary conditions, to the closed sector, in which only curves without boundary are considered. We shall later use the term \emph{closed-open map}\index{closed-open map} to refer to a map in the other direction.
\end{rem}
\begin{rem}
The Hochschild homology group $HH_*(\cL^{\bf b},\cL^{\bf b})$ (over Novikov field $\Lambda$) can be defined without reference to the decomposition into potential values, in which case Equation \eqref{eq:HH_defined_by_decomposition} becomes the consequence of a computation. Since we do not need such a result, we do not prove it in this paper.
\end{rem}
\par

We now give the first indication that much of the Floer theory of Lagrangians in $X$ can be extracted from a finite collection of Lagrangians. For the statement, let $U \subset X$ be a Lagrangian submanifold of $X$ (not necessarily lying in $\mathscr L$) equipped with a  bounding cochain $b_{U} \in \mathcal M(U,\fb;\Lambda_0)$. We recall that the pair $(U, b_U)$ is said to be Floer theoretically essential whenever the Floer cohomology group $HF((U;(\frak b,b_{U})),(U;(\frak b,b_{U}));\Lambda)$ does not vanish.

\begin{thm}\label{maintheorem1}
If the multiplicative unit ${\bf e}_{X} \in QH^0_{\fb}(X;\Lambda)$  lies in the image of $(\ref{mapphat(1.3)})$
after tensoring with $\Lambda$, then the potential value of any Floer theoretically essential pair $(U,b_U)$ agrees with the potential value of some element of $\bL$, i.e. there exists  $(L,b) \in \bL$  such that
\begin{equation}
  \label{eq:equal_PO}
\PO(L,b)  = \PO(U,b_U) .
\end{equation}
\end{thm}

Let $\lambda$ denote the common potential value in Equation (\ref{eq:equal_PO}). The construction of the category $\cL_\lambda^{\bf b}$ took an arbitrary finite collection of weakly unobstructed Lagrangians as input, and can thus be implemented with the object set  $\bL_{\lambda} \cup \{(U;b_U)\}$. We denote by  $\Fuk_{\lambda}^{\bf b}$ the corresponding cyclic 
unital filtered $A_{\infty}$ category, and observe in Subsection \ref{relative} that the natural inclusion on sets of objects extends to a filtered $A_{\infty}$ functor
\begin{equation}
I_{\lambda}:  \cL_{\lambda}^{\bf b} \to \Fuk_{\lambda}^{\bf b}.
\end{equation}

\par

Our main theorem asserts that $I_{\lambda}$ becomes an equivalence upon passing to an appropriate algebraic closure of the source category.  To be more precise, recall that the \emph{triangulated closure}\index{triangulated closure} $D(\Cat)$ of any $A_{\infty}$ category $\Cat$ is the $A_{\infty}$ category obtained from $\Cat$ by iteratively adding:
\begin{enumerate}
\item shifts of each object 
\item cones of each morphism.
\end{enumerate}
An explicit construction of $D(\Cat)$ using twisted complexes was first introduced in the setting of differential graded categories by Bondal and Kapranov \cite{BondalKapranov}, and is reviewed in Section \ref{sec:twisted-complexes}.    The \emph{idempotent closure}\index{idempotent closure} $D^{\pi}(\Cat)$ is then obtained by adding a summand representing each idempotent of an object in $D(\Cat)$. 
(See \cite[Section 4]{Sei06} and Section \ref{sec:idempotent-closure} of present paper.)
We are now ready to state our main result:\index[syindex]{DpiLcan@$D^{\pi}(\cL_{\lambda})$}

\begin{thm}\label{maintheorem2}
In the situation of Theorem $\ref{maintheorem1}$, 
the functor
$I_{\lambda}$ induces an equivalence of categories
\begin{equation}
D^{\pi}{I_\lambda} : D^{\pi}(\cL_{\lambda}^{\bf b}) \to D^{\pi}(\Fuk_{\lambda}^{\bf b}).   
\end{equation}
\end{thm}
\begin{rem} \label{rem:split-generation-summand}
As we explain in Section \ref{sec:idempotent-closure}, the above theorem is equivalent to the statement that $(U,b_{U})$, as an object of $D(\Fuk_{\lambda}^{\bf b})$, is a summand of an object obtained by iteratively taking cones of morphisms among objects in the original collection $\bL_{\lambda}$.  
\end{rem}

Theorem \ref{maintheorem2} allows us to interpret the hypothesis that ${\bf e}_{X}$ lies in the image of $\hat{\frak p}^{\text{\bf b}}$ as a split-generation criterion for a collection of Lagrangians in a symplectic manifold.  Since we are assuming that $X$ is compact, we can state a dual criterion in terms of Hochschild cohomology.  More precisely, we construct a homomorphism
\begin{equation}\label{mapqhat}
\wq^{\text{\bf b}}: QH^*_{\fb}(X;\Lambda_0)
\to  HH^*(\cL^{\bf b}, \cL^{\bf b})
\end{equation}
of algebras over the Novikov ring, where the right hand side is Hochschild cohomology, which we define in Section \ref{sec:hochsch-cohom} as the product of the Hochschild cohomology of the categories $ \cL_{\lambda}$ for various values of the potential function.  The map $\wq^{\text{\bf b}}$ is the closed-open map alluded to in Remark \ref{rem:open-closed-terminology}.
% \footnote{
% \rcng{Note $HH^*(\cL, \cL)$ 
% is the direct sum of $HH^*(\cL_{\lambda}, \cL_{\lambda})$ for various $\lambda$.}}

In Section \ref{sec:duality_bi-modules}, we show that the cyclic invariance of the $A_{\infty}$ structure on $ \cL^{\bf b}$ gives an isomorphism from $ HH^*(\cL^{\bf b}, \cL^{\bf b}) $ to the linear dual of $HH_{n-*}(\cL^{\bf b}, \cL^{\bf b})$.  Identifying $QH^*_{\fb}(X;\Lambda_0)$ with its dual using Poincar\'e duality, we prove the following result in Section \ref{dualpq}:
\begin{thm}\label{pqduealthe}
The map $\hat{\frak p}^{\text{\bf b}}$ is the linear dual to $\hat{\frak q}^{\text{\bf b}}$.
\end{thm}
We immediately conclude:
\begin{cor} \label{cor:generation_injective_statement}
  If $\hat{\frak q}^{\text{\bf b}}\otimes_{\Lambda_0} \Lambda$ is injective, then the hypothesis of Theorem $\ref{maintheorem1}$ is satisfied.  
 \par
 In the special case when $c_{1}(X) = 0$, $\fb\in H^2(X;\Lambda_0)$ 
  and $b \in H^1(L;\Lambda_0)$ for each $(L,b) \in \bL$,  and all elements of $\mathscr L$ are graded with respect to a chosen complex volume form,  it suffices to assume that the restriction of $\hat{\frak q}^{\text{\bf b}}$ to
  \begin{equation}
    QH^{2n}_{\fb}(X;\Lambda)
  \end{equation}
does not vanish, where $2n$ is the real dimension of $X$.
\end{cor}
In fact injectivity of $\frak q^{\bf b}:  QH^{2n}_{\fb}(X;\Lambda) \to HH^{2n}(L;\Lambda_0)$
is equivalent to the surjectivity of $\frak p^{\bf b} : HH_{-n}(L;\Lambda_0) \to QH^{0}_{\fb}(X;\Lambda)$.

\begin{rem}
In the situation of the second half of Corollary \ref{cor:generation_injective_statement},
 %\marginpar{I think $2n$ is the correct degree 
%for Corollary \ref{cor:generation_injective_statement}.  I add remark to clarify this point.  
%Actually I am still confused with $\Z$ graded case.
%KF 2025 Aug.}
$\frak q$ sends $QH_{\fb}^{2n}(X;\Lambda)$ to $HH^{2n}(\mathcal L)$.
(See (\ref{form460new}), (\ref{form464new}) etc. for the degree of Hochschild (co)chomology.)
When $\mathscr L$ consists of a single embedded Lagrangian submanifold,
$HH^{2n}(\mathcal L)$ is necessary $0$.\footnote{This is because only unit is an element of degree $0$.}  However in the case of single immersed 
Lagrangian, 
$HH^{2n}(\mathcal L)$ can be nonzero. For example if $p_i$ ($i=1,\dots, n$) is a self-intersection point of Maslov index $0$ and $p_0$ is a 
self-intersection point of Maslov index $n$ then 
$p_1 \otimes \dots \otimes p_{n} \mapsto p_0$ has 
shifted degree $2n$.\footnote{A similar situation appears in \cite{AS}.}
%A Hochishld chain which has
%non-trivial inner product with this Hochishld cohain 
%is $\tau(p_0) \otimes p_1 \otimes \dots \otimes p_{n+1}$.
%Here $\tau(x,y) = (y,x)$.  This Hochishld chain has shifted degree $-n-2$.
%So it is sent by $\frak p$ to an element of degree $-2$. 
%Taking into account degree shift of the ambient homology class, 
%this is the correct degree of $[1_{X}]$.
\end{rem}

We also prove (Theorem \ref{thm:duality}, Proposition \ref{nontrivialq}):
\begin{thm}\label{ringhomothm}
The closed-open  map $\hat{\frak q}^{\bf b}$ is a uinital ring homomorphism.
\end{thm}

We next prove the following result:
\begin{thm}\label{maintheorem3}
In the situation of Theorem $\ref{maintheorem1}$, the open-closed map  $(\ref{mapphat(1.3)})$ and the closed-open map $(\ref{mapqhat})$ are both isomorphisms.
\end{thm}

\begin{rem}
We mention that Floer cohomology with $\Lambda_0$ coefficient is used in Theorems \ref{pqduealthe} and \ref{ringhomothm} and 
while the $\Lambda$ coefficient
is used in Corollary \ref{cor:generation_injective_statement} and Theorems \ref{maintheorem1}, \ref{maintheorem3}.
\end{rem}

  As proved by Ganatra in the exact setting \cite{Ga}, the assumptions of  Theorem $\ref{maintheorem1}$ in addition imply that the category $\cL_{\lambda}$ is smooth in the sense of \cite{KS-Aoo}, i.e. that the diagonal bi-module is quasi-isomorphic to a summand of a perfect bi-module. In Floer theory, and using pseudo-holomorphic quilts \cite{mww},\cite{fukaya:functor} to relate Lagrangians in the product symplectic manifold $(X \times X, - \omega \oplus \omega)$ with bi-modules. 
In fact a filtered $A_{\infty}$ bi-functor from
$\sL \times \sL$ to a filtered $A_{\infty}$ category 
whose object set is ${\bf L} \times {\bf L}$ is constructed in \cite{fukaya:functor}. 
We will use it to prove the following.

\begin{thm}\label{maintheorem3.3}
In the situation of Theorem $\ref{maintheorem1}$ the filtered $A_{\infty}$ category 
$\cL^{\bf b}\otimes \Lambda$ is cohomologically smooth in the sense of \cite{KS-Aoo}.\footnote{See Definition \ref{defn261}.}
\end{thm} %\marginpar{Added. KF 2025 Aug 27}

We remark that the filtered $A_{\infty}$ category 
$\cL^{\bf b}\otimes \Lambda$ is always proper, that is, the morphism space is finite dimensional.
We will prove Theorem \ref{maintheorem3.3} in Section \ref{sec:smooth}.

Another application of operator $\frak p$ and Theorem \ref{maintheorem3} is 
the next result which is related to a conjecture by  Kontsevitch-Soibelman \cite[Conjecture 9.2]{KS}.

\begin{thm}\label{maintheorem3.5}
In the situation of Theorem $\ref{maintheorem1}$ there exists a
canonical isomorphism
$$
HC_*(\cL^{\bf b}, \cL^{\bf b})/(u^n) \cong H(X;\Lambda)\otimes \Lambda[\![u]\!]/(u^n).
$$
\end{thm}
Here $HC_*(\cL^{\bf b}, \cL^{\bf b})$ is the negative cyclic homology over $\Lambda$, which is 
a $\Lambda[\![u]\!]$ module.
Similar results for positive and periodic cyclic homologies also hold.
 %\marginpar{Added. KF 2025 July}
\index[syindex]{HC@$HC_*(\cdot,\cdot)$}

\begin{rem}
Kaledin \cite{kal} proved  the conclusion of Theorem \ref{maintheorem3.5}
under the assumption that the DG category is cohomologically smooth and proper, is defined over a field of characteristic $0$, and is $\Z$ graded.
So in case $\cL^{\bf b}$ is $\Z$ graded, Theorem  \ref{maintheorem3.5} may follow from Theorem \ref{maintheorem3.3}.\footnote{
A possible issue is  whether a proper $A_{\infty}$ category is quasi equivalent to a {\it proper} DG category or not.}
 %\marginpar{Footnote 11 is added. KF 2025 Sep 23.}
We do not know a similar result in the literature in the case when $A_{\infty}$ category is only $\Z_2$ graded.
\end{rem}

\begin{rem}\label{rem119}
If $\frak p \otimes \Lambda$ is injective (but not necessary surjective) we can still prove:
\begin{equation}
\aligned
&HC_*(\cL^{\bf b}, \cL^{\bf b})/(u^n)  \\
&\cong {\rm Im}\,(\frak p: 
HH_*(\cL^{\bf b}, \cL^{\bf b};\Lambda) \to H_*(X;\Lambda)) \otimes \Lambda[\![u]\!]/(u^n),
\endaligned
\end{equation}
in the same way.
\end{rem}

The above results complete the general Floer theoretic package which this paper constructs. In order to illustrate some applications of this work, we consider the case when $X$ is a compact toric manifold, possibly equipped with a bulk deformation. In \cite{fooo:toric1,fooo:bulk,toric3}, the last four authors studied the Floer theory of Lagrangian tori which arise as fibres of the moment map, and extracted from the  moduli spaces of holomorphic disks of Maslov index $2$, with boundary on such fibers, a potential function on the space of pairs $(L,b)$ consisting of a fibre, equipped with a first cohomology class $b$. The relationship with Floer theory is that the critical points of this function are exactly the pairs $(L,b)$ for which the Floer cohomology $HF^*((L,b),(L,b);\Lambda)$ is nontrivial, and the results of \cite{fooo:bulk} imply that the set of critical values is finite; this yields a finite collection $\bf L$ of such pairs $(L,b)$, which define a non-trivial category $\cL^{\bf b}$.

In \cite{toric3}, the last four authors proved that the quantum cohomology of a toric manifold can be recovered from the potential function on the space of fibers equipped with bounding cochains; in fact, only the formal expansion of this function near the critical points is required. Given the relationship between the potential function and Floer cohomology, this suggests that quantum cohomology can be recovered from the category $\cL^{\bf b}$; this is what we prove in Section \ref{sec:toric}:
 
\begin{thm}\label{maintheorem4}
Let $(X,\omega)$ be a compact toric manifold equipped with the trivial background class and an arbitrary bulk class $\frak b \in H^{\rm even}(X;\Lambda^{\C}_0)$.
Let $\bf L$ be the set of all pairs $(L,b)$ where $L$ is a torus orbit of $X$ and $b \in H^1(L;\Lambda^{\C}_0)/2\pi \sqrt{-1}H^1(L;\Z)$, such that the Floer cohomology 
$HF((L,(\frak b,b)),(L,(\frak b,b));\Lambda^{\C})$ does not vanish. % \cong H^*(L;\Lambda)$.
Then the map
$$
\widehat{\frak p}^{\text{\bf b}} : HH_*(\cL^{\bf b}, \cL^{\bf b};\Lambda^{\C}) \to QH^*_{\fb}(X;\Lambda^{\C})
$$
is surjective.
\end{thm}

\begin{rem}
In this paper we take and fix bulk class $\frak b$ and bounding cochains $\{b\}$.  So hereafter we omit the symbol ${\bf b}$ from the notation 
$\cL^{\bf b}$, $\hat{\frak p}^{\bf b}$, $\hat{\frak q}^{\bf b}$ sometimes.\index[syindex]{Lcal@$\sL$}
\end{rem}

\subsection*{Related results}
\label{sec:related-results}

Many of the results of this paper were announced in 2011--12; we apologize for the long delay in the appearance of a publicly available text. In the meantime, several related results have appeared:
\begin{itemize}
\item In \cite{sheridan-fano}, Sheridan proved a version of Theorem \ref{maintheorem2} for monotone Lagrangians in closed monotone symplectic manifolds (i.e. those for which the first Chern class is positively proportional to the symplectic form). A similar result for monotone symplectic manifolds with contact-type boundary was proved by Ritter and Smith in \cite{RS}.
\item  In the setting of relative Fukaya categories (i.e. categories associated to a symplectic manifold with a choice of symplectic hypersurface Poincar\'e dual to the symplectic form), Perutz and Sheridan proved in \cite{PS} that an equivalence of a subcategory of the derived Fukaya category with the category of coherent sheaves on a maximally degenerate mirror family implies that the generation criterion holds.
\item Ganatra proved in \cite{Ga2}, that whenever a subcategory of the Fukaya category is \emph{cohomologically smooth} in the sense of Kontsevich (i.e. admits a split-resolution of the diagonal bi-module), the open-closed map is necessarily injective. In particular, the hypothesis of Theorem \ref{maintheorem2}  can be verified by comparing ranks. We postpone the discussion of how Ganatra's results apply in our setting to Remark \ref{rem:ganatra-result}.
 %\marginpar{We need to quote more. KF 2025 July.
%I added last item.  Need to check whether this is OK.  KF 2025 Aug.}
\item See \cite{AT,Ga3,GTS,OS,PS2,Sd,she2,she3} and etc. for recent references related to open-closed Gromov-Witten theory.
The papers \cite{Ga3,GTS,PS2} etc. discuss exact Lagrangians in an exact symplectic manifold, and then study Floer theory of a compact 
symplectic manifold $X$ with 
a divisor $D$, which is Poincar\'e dual to the symplectic form, as a deformation of the Floer theory on $X \setminus D$.  
Such Floer theory is used in \cite{Se4,sheridan-fano} etc. to prove homological Mirror symmetry in various cases.
We have no doubt 
that the $A_{\infty}$ category (together with open-closed and closed-open maps) obtained in their method  will be embedded to 
the one we construct directly on $X$ in this paper, over Novikov field, and they are quasi isomorphic in various cases. 
(We do not try to prove it in this paper.)
\end{itemize}

\subsection*{Outline of the paper} %\marginpar{Outline is rewritten.  KF. 2024 Dec.}

The results contained in  Part \ref{part1} of this paper concern constructions in homological algebra which are for the most part well-established in the setting of ordinary categories, but for which no adequate reference seems to exist in the presence of curvature.   The main purpose of Section \ref{sec:curved-cat} is to recall the notion of a curved cyclic $A_{\infty}$ category, and of a weakly curvature free object in such a category.   In Sections \ref{sec:Hochschild} and \ref{sec:a_infty-categ-over} we define the Hochschild homology and cohomology of such categories with coefficients in an $A_{\infty}$ bi-module, and prove that the cyclic structure defines a duality between these groups.  In Section \ref{sec:triang-clos-split} we define the category of twisted complexes, which will be our model for the bounded derived category of an $A_{\infty}$ category, and discuss the notion of split-generation of an $A_{\infty}$ category.  Finally, Section \ref{sec:a_infty-categ-over} studies the problem of constructing a cyclic $A_{\infty}$ category from a homotopy inverse system: we solve the problem under the further assumption that the maps in the inverse system arise from pseudo-isotopies (see Subsection \ref{sec:pseudo-isotopy}).

In Part \ref{part2}, we study Lagrangian Floer theory in a compact symplectic manifold.\footnote{Some of the results of Section \ref{sec:cyclicfil}-\ref{sec:ring} are generalisations, to the case of multiple Lagrangians, of results which have been established, in \cite{fooo09,fooo:bulk,toric3} either in the case of a single Lagrangian, or in the special situation of toric manifolds.  There are several other results 
which prove these results under additional assumption such as exactness or monotonicity. 
In the generality we work in this paper, those results are new.} Section \ref{sec:cyclicfil} 
(together with Sections \ref{sec:unit}-\ref{sec:CRperturb}) provides an essentially complete construction of the $A_{\infty}$ structure on the Lagrangian Floer cohomology of a finite collection of Lagrangian submanifolds; this includes in particular the construction of pseudo-isotopies which is implemented in Section \ref{sec:homotopyequiv}.  
In Section \ref{sec:unit} we discuss the construction of homotopy unit in our de Rham setting.
(The reason we need to work with homotopy unit (but not with strict unit) is explained in Section \ref{sec:two perturbation}.)
In Section \ref{sec:elicurv} we discuss how to use bounding cochain to eliminate the curvature.
In Section \ref{canonical} we reduce the $A_{\infty}$ structure on a chain model (the de Rham complex in our case)
to its homology.  The new point  in this paper is to perform this construction in the situation when the 
inner product (with cyclic symmetry) and the homotopy unit are given.

The results of Sections \ref{sec:cyclicfil}-\ref{canonical}   rely on the constructions of Kuranishi structures on the moduli spaces 
of pseudo-holomorphic discs with Lagrangian boundary conditions, and of CF-perturbations.
These constructions are given
in Sections \ref{sec:Kuraconst} and \ref{sec:CRperturb},
after reviewing the necessary notions. 
Among the new results which we introduce in these sections is a method of obtaining 
Kuranishi structures and CF-perturbations  in a way compatible with forgetful maps.
They are essential for the $A_{\infty}$ structure we obtain to be unital.

In Part \ref{part3} we study  the relationship between the Lagrangian Floer theory to quantum cohomology. 
In Section \ref{sec:frakq}, we define the map $\hat{\frak q}^{\text{\bf b}}$ from the quantum cohomology of the ambient space to the Hochschild cohomology of the $A_{\infty}$ category constructed in the previous part.  In Section \ref{sec:open-closed map} we define a map $\hat{\frak p}^{\text{\bf b}}$ from 
Hochschild homology\footnote{Since we use de Rham model we identify $(\dim X - k)$-th homology with $k$-th cohomology.} to the homology of the ambient space.
In Section \ref{dualpq} we prove that the map $\hat{\frak p}^{\text{\bf b}}$ is dual to $\hat{\frak q}^{\text{\bf b}} $.
In Section \ref{nontrivialsec} we show that those maps are non-trivial under the assumption that the Floer homology 
is non-trivial.
In Section \ref{sec:ring}, we prove that $\hat{\frak q}^{\text{\bf b}}  $ is a ring homomorphism with respect to the natural structures on either side.   
In Section \ref{sec:annuli} we begin implementing the homological interpretation of the annulus argument (see \cite{abouzaid:IHES}) to the setting of this paper. 
Section \ref{sec:ori} is devoted to the discussion on orientation and sign.
The novel point in this section is the part related to the Cardy relation.

In Part \ref{part4} of the paper, we provide applications of the algebraic structures introduced in Part \ref{part1} to the study of Lagrangian Floer theory.  Sections \ref{sec:Theorem1}, \ref{sec:Theorem2}, \ref{sec:Theorem3}, \ref{sec:smooth}, \ref{sec:Hodge-to-de Rham} and \ref{sec:toric}  prove Theorems \ref{maintheorem1},    \ref{maintheorem2}  \ref{maintheorem3}, \ref{maintheorem3.3}, \ref{maintheorem3.5} and \ref{maintheorem4}, respectively; these are the main results of the paper.  They show that, under the assumption of Theorem \ref{maintheorem1}, the $A_{\infty}$ category and the quantum cohomology decompose as a direct sum of the potential values appearing non-trivially in the decomposition $\bL = \cup_{\lambda} \bL_{\lambda}$.  The next three sections concern this decomposition in the absence of this strong assumption:  In Section \ref{sec:Theorem2+1} we show that the images of the Hochschild homologies of the categories $\cL_{\lambda}$ under $ \hat{\frak p}^{\text{\bf b}} $ are orthogonal with respect to quantum multiplication.  In Section \ref{sec: decomp}, we compare the decomposition of quantum cohomology into indecomposable factors with the decomposition of the category $\cL$ by potential values, while in Section \ref{sec: value} we prove, under a Calabi-Yau assumption, a conjecture of Kontsevich about the relation of these decompositions with the  eigenvalues of the quantum product with the first Chern class.  The paper ends with Section \ref{sec: quatitative}, in which we state a quantitative version of the main results, that allows us to provide bounds from below on the minimal areas of holomorphic disks with boundary on Lagrangian submanifolds.
\par
The authors thank to L. Amorim who pointed out (in 2019 February) an important issue related to the unitality and cyclic symmetry  
in the case of immersed Lagrangian or several mutually transversal Lagrangians. (Such an issue does not exist 
when studying disks bounding a single embedded Lagrangian.\footnote{It however occurs when one studies 
higher genus bordered curves bounding a single embedded Lagrangian.}) (See Section \ref{sec:two perturbation}.) 
They also  thank Kyoji Saito and Atsushi Takahashi for their discussions regarding the topics studied
in Sections \ref{sec:smooth}, \ref{sec:Hodge-to-de Rham}.  
K. Fukaya wants to thank Aliakbar Daemi  for discussions especially those related to Subsection \ref{diskforget}.
\par\bigskip
\section{Notations  and Conventions}\label{Sec3Not} 
\subsection{Labels of Lagrangians}
We fix a class $\frak b_2 \in H^2(X;\F)$.\footnote{$\frak b$ is a bulk deformation and 
$\frak b_2$ is a part of it. See (\ref{decomposefrakb}).} %\marginpar{modified a bit.  KF. 2025 Jan}
We  fix a finite set $\mathscr L$ of pairs $(L,\theta)$ of relatively spin immersed Lagrangian submanifolds $L$ 
equipped with $1$ forms $\theta$ with $d\theta = \frak b_2$,\footnote{More precisely if 
$L = (\tilde L,i_L)$ is immersed then we require $i_L^*(\frak b^2) = d\theta$ } which we consider in this paper
and write:\index[syindex]{Lscr@$\mathscr L$}
\begin{equation*}
\mathscr L = \{(L_{\kappa},\theta_{\kappa}) \mid \kappa = 1,\dots,\#\mathscr L\}.
\end{equation*}

We require that $L_i$ are self-transversal. For a natural number $K$ we put\index[syindex]{Kunderine@$\uwave{K}$}
$$
\uwave{K} = \{0,\dots,K\}.
$$

\subsection{Immersed Lagrangians}
An immersed Lagrangian submanifold $L$ of $X^{2n}$ is identified with 
a pair $(\tilde L,i_L)$ of an  $n$ dimensional manifold $\tilde L$ and an immersion\index[syindex]{iL@$i_L$}\index[syindex]{Ltilde@$\tilde L$}
$i_L : \tilde L \to X^{2n}$ such that $i_L^*\omega_X = 0$.  

\subsection{De Rham theory}
In this paper, {we shall consider a field} $\F$ which is either $\R$ or $\C$; since we work with de Rham theory, we do not work over the rational numbers. {The choice of field is often suppressed from the notation, so that} we denote the de Rham complex of smooth forms on a manifold $X$ by 
$\Omega(X)$, {when considering a context in which both real and complex coefficients can be used.}\index[syindex]{Lambda0@$\Lambda_0$}

\subsection{Novikov rings}
The universal Novikov ring\index{universal Novikov ring}\index{Novikov ring}  $\Lambda_0$ is defined by 
\begin{equation}\label{formula114}
\Lambda_0
= \left\{ \sum_{i=1}^{\infty} a_i T^{\lambda_i} \mid a_i \in \F, \lambda_i \in \R_{\ge 0},
\lim_{i\to\infty} \lambda_i = \infty \right\}.
\end{equation}
Here $\F =\C$ or $\R$\index[syindex]{F@$\F$} in this paper, which we call the {\em ground field}\index{ground field}.
The Novikov field $\Lambda$ \index{Novikov field}\index[syindex]{Lambda@$\Lambda$} is its field of fractions, and $\Lambda_+$
\index[syindex]{Lambdaplus@$\Lambda_+$} its maximal ideal.

\subsection{Shift of graded vector spaces}
Let $C$ be a graded vector space. {We shall use the standard notation in homological algbera} $C[k]$ {to indicate} the degree shift of a $\Z$ graded
$\F$-vector space $C$ by {an integer} $k$. This is defined by\index[syindex]{CF[1]@$C[1]$}
\begin{equation}
(C[k])^d = C^{k+d}.  
\end{equation}

\subsection{Bar construction for algebras}
{At various points, we shall consider coalgebras denoted $BC$ and $EC$, which are associated to each graded vector space $C$, and are defined as follows.}   We
put
$
B_kC = \underbrace{C\otimes \cdots \otimes C}_{\text{$k$ times}}.
$
The symmetric group $\mathfrak S_k$ of order $k!$ acts on $B_kC$ by
$$
\sigma \cdot (x_1 \otimes \cdots \otimes x_k)
= (-1)^* x_{\sigma(1)} \otimes \cdots \otimes x_{\sigma(k)},
$$
where
$
* = \sum_{i<j; \sigma(i)>\sigma(j)} \deg x_i \deg x_j$.
We define $E_kC$ to be the quotient of $B_kC$ by the subspace generated by
$(x_1 \otimes \cdots \otimes x_k) 
-  \sigma\cdot (x_1 \otimes \cdots \otimes x_k)$.
\par
Let $B_kC$ be as above and define
$
BC = \bigoplus_{k=0}^{\infty} B_kC.
$
(We remark $B_0C = \F$.) This direct sum has the structure of a coassociative
coalgebra with its coproduct $\Delta: BC \to BC \otimes BC$ defined\index[syindex]{Delta@$\Delta$}
by\index[syindex]{BC@$BC$}
$$
\Delta(x_1 \otimes \cdots \otimes x_k)
= \sum_{i=0}^k (x_1 \otimes \cdots \otimes x_i)
\otimes (x_{i+1} \otimes \cdots\otimes x_k).
$$
There exists a similar coproduct $\Delta: EC \to EC \otimes EC$ on
$EC$ (using shuffles, see \cite[(1.23), (1.25)]{spectre}) which becomes a coassociative and graded\index[syindex]{EC@$EC$}
cocommutative. 
\par
We shall denote the iterated coproduct by $\Delta^{n-1}: BC \to (BC)^{\otimes n}$ or $EC \to
(EC)^{\otimes n}$; this is defined by the formula
$$
\Delta^{n-1} = (\Delta \otimes  \underbrace{{\rm id} \otimes \cdots
\otimes {\rm id}}_{n-2}) \circ (\Delta \otimes  \underbrace{{\rm id} \otimes
\cdots \otimes {\rm id}}_{n-3}) \circ \cdots \circ \Delta.
$$
Using Sweedler's notation,\index{Sweedler's notation} an element $\text{\bf x} \in B_{k}C$ can be expressed as\index[syindex]{xn:1c@$\text{\bf x}^{n;1}_c$}
\begin{equation}\label{deltawritebyc}
\Delta^{n-1}(\text{\bf x}) = \sum_c \text{\bf x}^{n;1}_c \otimes
\cdots \otimes \text{\bf x}^{n;n}_c
\end{equation}
where $c$ runs over an  index set depending on $\text{\bf x}$, which consists of the set of 
ordered partitions of $n$ if $\text{\bf x} $ is an indecomposable element of $B_nC$.   
%\blue{(When ${\mathbf x}$ is a ``monomial'', then the index set is 
%the one consisting partitions of $n$.)}  
 
\subsection{Shifted degrees}
For an element  %\marginpar{MA190515: Before, this was only defined for $C = H(L;R)$}
$
\text{\bf x} = x_1 \otimes \cdots \otimes x_k \in B_k(C[1])
$
we denote the shifted degree $\deg'x_i = \deg x_i - 1$ and\index[syindex]{deg'@$\deg'$}
$
\deg' \text{\bf x} = \sum \deg'x_i = \deg \text{\bf x} - k.
$
(Here $\deg x_i$ is the degree of $x_i$ before introducing the shift by $1$.)

\subsection{Bar construction for categories}
We use the notation $BC$ in the following slightly more general situation: assume that we are given, for each pair  $\kappa_1,\kappa_2 \in \uwave{K}$, {a graded vector space} $C_{\kappa_1,\kappa_2}$.
Then for a sequence $\vec{\kappa} = (\kappa_0,\dots,\kappa_k)$ we put:
$$
B_{\vec{\kappa}}(C[1]) =  C_{\kappa_0,\kappa_1}[1] \otimes 
\cdots \otimes C_{\kappa_{k-1},\kappa_k}[1].
$$
We then define
\begin{align*}
B_{k}(C[1]) & = \bigoplus_{\vec\kappa \in \uwave{K}^k}
B_{\vec{\kappa}}(C[1]) \\
B(C[1]) & = \bigoplus_{k}B_{k}(C[1]) .
\end{align*}
There is a coalgebra structure $\Delta \co B(C[1]) \to B(C[1]) \otimes B(C[1])$ defined as before. We also use the same notation (\ref{deltawritebyc}) in this case.
\par
Note 
$
B_{0}(C[1]) \cong \Lambda_0^{\#\mathscr L}.
$
Each of the generator\index[syindex]{B0C1@$B_{0}(C[1])$} corresponds to an element of $\mathscr L$.
Keeping in mind that the set $\uwave{K}$ labels objects of the category, assuming that $i$ labels an object $c = (L,\theta,b)$ of the category, we shall write $1_c$\index[syindex]{1c@$1_c$} the generator associated to the $i$-th component of $B_{0}(C[1]) = \Lambda_0^{\#\mathscr L}$.

\subsection{Kuranishi structures}
Let $\mathcal M$ be a space with Kuranishi structure and let $\text{\rm ev}_s:
\mathcal M \to M_s$, $\text{\rm ev}_t: \mathcal M \to M_t$ be strongly
continuous smooth maps. (See  \cite[Definition 6.6]{FO}.)
 %\marginpar{MA190515: If this is also defined in FOOO9 or FOOO10, we should give an additional reference} 
 (Here $s$ and $t$ stand for source and target,
respectively.) We assume that our smooth manifolds $M_s, M_t$ are compact
and oriented without boundary. We also assume that $\mathcal M$ has a
tangent bundle and is oriented in the sense of Kuranishi structure.
(See  \cite[Definition A1.14]{fooo09}.)
We include the case when $\mathcal M$ has a boundary or corners.
We assume that $\text{\rm ev}_t$ is weakly submersive. (See  \cite[A1.13]{fooo09}
and the description below.)
{Under the above assumptions,} we use the work of \cite[Section 12]{fooo:bulk} and \cite{tech2} to define a
{\it smooth correspondence map}\index{smooth correspondence map}\index[syindex]{Corr(evsevt@$\text{\rm Corr}(\mathcal M;\text{\rm ev}_s,\text{\rm ev}_t)$}
\begin{equation}\label{smoothcorrespondence}
\text{\rm Corr}(\mathcal M;\text{\rm ev}_s,\text{\rm ev}_t): \Omega^kM_s \to \Omega^{k+\dim M_t - \dim \mathcal M}M_t
\end{equation}
associated to  $(\mathcal M;\text{\rm ev}_s,\text{\rm ev}_t)$. This map is defined as a composition
$$
\text{\rm Corr}(\mathcal M;\text{\rm ev}_s,\text{\rm ev}_t)(g)
=
(\text{\rm ev}_t)_! (\text{\rm ev}_s)^*(g)
$$
where $ \text{\rm ev}_s^*(g)$ is the 
pull back of the form $g \in \Omega(M_s)$ to a form on 
$\mathcal M$. { More precisely, we use the notion of  CF-perturbations, introduced in \cite{tech2} and which is a variant of the notion of continuous family of perturbations from \cite{fooo:bulk}, and define $\text{\rm ev}_s^*(g)$ as a form on the zero set of a chosen \blue{CF-perturbation}.} The map
$(\text{\rm ev}_t)_! $ is the integration along the fiber.

Note that (\ref{smoothcorrespondence}) actually depends on the choice of the 
perturbation (\newred{CF-perturbation}).
If we take the \newred{CF-perturbations} in a compatible way, then  
(\ref{smoothcorrespondence}) satisfies 
Stokes' formula \cite[Lemma 12.13]{fooo:bulk}, \cite[Theorem 8.11]{springer}. It is also compatible with 
composition of correspondences \cite[Lemma 12.15]{fooo:bulk}, \cite[Theorem 10.21]{springer}.

\subsection{Shuffles}
For the description of the boundaries of various moduli spaces of Riemann surfaces, we introduce the notion of an {\it $\ell$-shuffle} $(\Bbb L_1,\Bbb L_2)$,\index{shuffle} for each $\ell \in \Z_{+}$, which  is 
a decomposition of $\underline{\ell}$ to a disjoint union 
$\Bbb L_1 \cup \Bbb L_2$.
Let $\text{\rm Shuff}(\ell)$ be the set of all shuffles of $\underline{\ell}$\index[syindex]{lxellunderline@$\underline{\ell}$}.
We enumerate  $\Bbb L_1$, $\Bbb L_2$ by its order and regard them 
as an ordered sets.
A {\it triple shuffle} $(\Bbb L_1,\Bbb L_2,\Bbb L_3)$ is defined in the same way.

\subsection{Categories of Lagrangians}
\label{subsec:CatLag}
 %\marginpar{Rewritten.  KF. 2024 Dec. More KF. 2025 Jan. and Feb.}

Recall that $\mathscr L$ is the fixed set of relatively Spin immersed Lagrangians (equipped with 
$\theta$) that we consider. We write 
$$
\mathscr L = \{ (L_{\kappa},\theta_{\kappa})\}.
$$
Note it may happen $L_{\kappa} = L_{\kappa'}$ for $\kappa' \ne \kappa$. (In such case 
$\theta_{\kappa} \ne \theta_{\kappa'}$.)
We denote by ${\bf L}$ a finite set of triples 
$(L,\theta,b)$ where $(L,\theta) \in \mathscr L$ and $b_+$ is a  bounding cochain of 
$(L,\theta)$ such that degree $1$ part of $b_+$ is $0 \mod \Lambda_+$ and 
$b = \theta + b_+$.
(See Definition \ref{boundingcochain}.) 
We write $\cL$ for the (unital and cyclic filtered) $A_{\infty}$ category whose objects are the elements of 
${\bf L}$, and whose morphism spaces are de Rham cohomology groups. We write ${\bf L}_{\lambda}$ for the subset of ${\bf L}$ consisting of triples $(L,\theta,b)$ such that 
$\frak{PO}_{\frak b}(b) = \lambda$ (See Subsection \ref{cvnoncomon}). As discussed earlier in the introduction, we write
$\cL_{\lambda}$ for the $A_{\infty}$ subcategory of $\cL$ whose objects are the elements of 
${\bf L_{\lambda}}$ (Condition \ref{conds614}).

Note that it may happen $(L,\theta,b), (L',\theta',b') \in {\bf L}$ and $(L,\theta)= (L',\theta')$ but $b \ne b'$.
\par\medskip
For a given finite set $\mathscr L$ we study three slightly different $A_{\infty}$ categories. (In Parts \ref{part2} and \ref{part3}.)
\par
$\mathcal L^{\rm form}_{\rm curve}$ is a cyclic curved filtered $A_{\infty}$ category whose modules of morphisms are de Rham complex.
(Definition \ref{defn847}.)  When the bulk deformation by $\frak b \in H^{\rm even}(X;\Lambda_0)$ is included we write,
 $\cL^{\rm form}_{\rm curve}(\frak b;\Lambda_0) := \mathcal L^{\rm form}_{\rm curve}$. 
\par
We write $\cL^{\rm form}_{\rm c.u.}$ for the curved 
homotopically unital and cyclic filtered $A_{\infty}$ category whose objects are elements of $\mathscr L$, and whose morphisms are constructed using differential forms.  
For the morphism module of $\mathcal L^{\rm form}_{\rm c.u.}$ we add extra generators ${\bf e}^+_c$ and ${\bf f}_c$ to the 
morphism modules $\mathcal L^{\rm form}_{\rm curve}(L,L;\F)$, 
 and define 
 $$
CF(\mathcal L^{\rm form};L,L;\F)^+\widehat{\otimes}\Lambda_0  = CF(\mathcal L^{\rm form};L,L;\F) \widehat{\otimes} \Lambda_0 \oplus \Lambda_0 [{\bf e}^+_{L}] \oplus \Lambda_0 [{\bf f}_{L}].
 $$
 (Begining of Subsection \ref{subsec:opeartor}.)
\par
We use weak bounding cochains to eliminate the curvature and obtain a homotopically unital and cyclic filtered $A_{\infty}$ category 
$\mathcal L^{\rm form}_{\rm uni}$, whose objects are all weakly curvature free.  (Definition \ref{defn105}.) 
We note (See Theorem \ref{cAinfconstuni}) 
$$
CF(\mathcal L^{\rm form};L,L;\F)^+\widehat{\otimes}\Lambda_0 =
\mathcal L^{\rm form}_{\rm uni}((L,b),(L,b)) = \mathcal L^{\rm form}_{\rm c.u.}((L,b),(L,b))^+.
$$
 \par
The morphism modules of the above $A_{\infty}$ categories are all de Rham {\it complex}.
We then replace de Rham {\it complex} by de Rham {\it cohomology}
to obtain $\mathcal L$ in Section \ref{canonical}.  It is strictly unital cyclic and its objects are all weakly curvature free.

\subsection{$\frak q$ and $\widehat{\frak q}$}
\label{labelqqhat}

The closed-open map ${\frak q}$ appears 
in two slightly different ways as follows.
Let $H(X)$ be the cohomology of the ambient space, 
$CF(L,L')$  the morphism space of the $A_{\infty}$ category 
and $BCF(L,L')$ its bar complex.
Then 
$\hat{\frak q}$
is a map $H(X)^{\otimes \ell} \to 
Hom(BCF(L,L'),CF(L,L'))$ 
and 
${\frak q}$ is a map
$H(X)^{\otimes \ell} \otimes  BCF(L,L')
\to CF(L,L')$.
There is an obvious identification between these 
two types of maps.
We will use symbols $
\frak q
$ and $\widehat{\frak q}$ (and their variants)  in this way
sometimes without mentioning it.

\subsection{Hochschild complex of the diagonal bi-module}
\label{hochchild}

We will define Hochschild chain and cochain complexes 
$CH_*(\mathcal C,\mathcal Q)$, $CH^*(\mathcal C,\mathcal Q)$
of a $\mathcal C$ bi-module $\mathcal Q$, where $\mathcal C$
is an $A_{\infty}$ category.  (See Definition \ref{defn412}.)
In case $\mathcal Q$ is the diagonal bi-module $\mathcal C$, (See the begining of Subsection \ref{sec:examples-bi-modules}.)
we sometimes write $CH_*(\mathcal C)$, $CH^*(\mathcal C)$, instead.\index[syindex]{CH*C2@$CH_*(\mathcal C)$}
The same applies to Hochschild (co)homology
$HH_*(\mathcal C)$, $HH^*(\mathcal C)$\index[syindex]{HHupper1@$HH^{*}(\mathcal C)$}
and cyclic (co)homology $HC_*(\mathcal C)$, $HC^*(\mathcal C)$.

\newpage
\part{Homological algebra of $A_{\infty}$ categories.}
\label{part1}

\section{Curved $A_{\infty}$ categories.}
\label{sec:curved-cat}
This is the first of three sections studying the homological algebra of curved $A_{\infty}$ categories.  While most of the results we describe are well-known in the absence of curvature, there seems to be no account in the literature that proves precisely what we want.

We shall consider a strictly unital and weakly curvature free $A_{\infty}$ category $\Cat$ defined over  a commutative ring $R$ (see Definition \ref{defn:a_infty-categories}).  In our applications, this ring will be either the Novikov field $\Lambda$,  the Novikov ring $\Lambda_{0}$, or a quotient of the Novikov ring by a closed ideal.   Note that, to construct categories over $\Lambda_0$, we shall need the more specialized notion of a \emph{gapped filtered category,} which we discuss in Section \ref{sec:a_infty-categ-over}.

\subsection{$A_{\infty}$ categories}

\label{sec:a-oo-categories}
Since the categories appearing in symplectic topology differ slightly from the standard setting, we start by describing our framework.  We consider a collection of objects denoted by $\Ob(\Cat)$,\index[syindex]{ObC@$\Ob(\Cat)$} and, for each pair, a morphism space $\Cat(X,Y)$ which is a $\Z_{2}$-graded module over $R$.\index[syindex]{BC1@$B(\Cat[1])$}  To define the structure maps of the $A_{\infty}$ category $\Cat$, we start by introducing the {\it bar complex}\index{bar complex}
\begin{align}
 B(\Cat[1]) & = \bigoplus_{X,Y \in \Ob(\Cat)} \bigoplus_{d = 0}^{\infty}   B_{d}(\Cat[1])(X,Y) \\
  B_{d}(\Cat[1])(X,Y) & = \bigoplus_{ X_{i} \in  \Ob(\Cat) } \Cat(X,X_1)[1] \otimes \Cat(X_1,X_2)[1] \otimes \cdots \otimes \Cat(X_{d-1},Y)[1],
\end{align}
with the convention that
\begin{equation}
B_{0}(\Cat[1])(X,Y)= \begin{cases}  R   & \textrm{ if } X=Y \\
 0 & \textrm{otherwise.} \end{cases}
\end{equation}
\begin{defn} \label{def:length_filtration_bar_complex}
 The \emph{length filtration}\footnote{This filtration is called 
 a number filtration in \cite{fooo09} and subsequent work.}\index{length filtration} on $B(\Cat[1]) $ is the filtration by submodules
\begin{equation}
 B^{(\leq N)}(\Cat[1])  = \bigoplus_{X,Y \in \Ob(\Cat)} \bigoplus_{d = 0}^{N}   B_{d}(\Cat[1])(X,Y)
\end{equation}
\end{defn}

In terms of the bar complex, an $A_{\infty}$ category $\Cat$ consists of a collection of maps
\begin{equation}
  \fm^{\Cat}_{d} \co B_{d}(\Cat[1])(X,Y) \to \Cat(X,Y)[1]
\end{equation}
for $ 0 \leq d$ of degree $1$, satisfying the $A_{\infty}$ equation.\index{$A_{\infty}$ equation} Writing\index[syindex]{mxCk@$\fm^{\Cat}$}
\begin{equation}
   \fm^{\Cat} \co B(\Cat[1])(X,Y) \to \Cat(X,Y)[1]
\end{equation}
for the sum of all structure maps $ \fm^{\Cat}_{d} $,  and using the notation of Equation \eqref{deltawritebyc}, we express this equation concisely as the condition that, for each $\bfa \in B(\Cat[1])$
  \begin{equation}
    \label{eq:a_oo-category-equation}
    \sum_{\alpha} (-1)^{\deg'(\bfa^{3;1}_{\alpha(3)})} \fm^{\Cat}(\bfa^{3;1}_{\alpha(3)}, \fm^{\Cat}( \bfa^{3;2}_{\alpha(3)}),\bfa^{3;3}_{\alpha(3)} ) = 0.
  \end{equation}
Whenever $\bfa$ is a word of length $d$, $\alpha(3)$ is an index running over all ordered 3-partitions of the integer $d$ into $3$ non-zero integers.
\begin{rem}
Note we do not assume $\frak m^{\Cat}_0 = 0$.
We use the name \emph{curved $A_{\infty}$ category}\index{curved $A_{\infty}$ category} sometimes to highlight this fact.
\end{rem}

\begin{defn} \label{defn:a_infty-categories}
The $A_{\infty}$ category $\Cat$  is \emph{strictly unital}\footnote{
We say `unital' in place of `strictly unital' sometimes in this paper.}\index{strictly unital}\index{unital}\index[syindex]{eBFX@${\bf e}_X$}
 if there exists an element ${\bf e}_{X} \in \Cat(X,X)$ satisfying the following conditions %\marginpar{$\id_{X}$ is ${\bf e}_{X}$ in Part 3.}
\begin{align}
   \fm^{\Cat}_{1}({\bf e}_{X}) & = 0 \\
\fm^{\Cat}_{2}({\bf e}_{X}, a) = (-1)^{\deg(a)} \fm^{\Cat}_{2}(a , {\bf e}_{X})  & =a \\
\fm^{\Cat}_{d}(a_1, \ldots,a_{k}, {\bf e}_{X} , a_{k+2}, \ldots, a_d) & = 0 \textrm{ if } 3 \leq d.
\end{align}

An object $X$  is \emph{weakly curvature free}\index{weakly curvature free} if there exists a scalar $\Po(X) \in R$ (called the \emph{potential value} of $X$) such that\index[syindex]{PO@$\Po$}
\begin{equation}
  \fm^{\Cat}_{0}(1) = \Po(X) \cdot {\bf e}_{X}.
\end{equation}
\end{defn}

When studying the homological algebra of curved categories, we only consider categories in which all objects are weakly curvature free.  Note that the $A_{\infty}$ equation \eqref{eq:a_oo-category-equation} when $d=2$ reads
  \begin{equation}
 \fm_{1}^{\Cat}\left( \fm_{1}^{\Cat}(a)\right) + \fm_{2}^{\Cat}\left( \fm_{0}^{\Cat}(1), a\right) + (-1)^{\deg'(a)}  \fm_{2}^{\Cat}\left(a,  \fm_{0}^{\Cat}(1)\right) = 0,
  \end{equation}
so that the endomorphism $   \fm_{1}^{\Cat}$ can only be guaranteed to define a differential on the  morphism spaces $\Cat(X,Y)$ if $X$ and $Y$ are weakly curvature free and $\Po(X) = \Po(Y)$; if these scalars do not agree, we obtain a matrix factorisation (see Section \ref{sec:matr-fact-dimens}). %\marginpar{A sentence removed.}

\begin{defn} \label{def:quasi-isomorphic}
A  weakly curvature free object  $Y$ is a \emph{summand}\index{summand} of a weakly curvature free object $X$ if there are closed morphisms
\begin{align}
f \in & \Cat(X,Y)   \\
g \in & \Cat(Y,X)
\end{align}
such that
\begin{equation}
  [\fm_{2}^{\Cat}(g,f)] = [{\bf e}_{Y}] \in H^{*}( \Cat(Y,Y),  \fm_{1}^{\Cat}).
\end{equation}
If, in addition,
\begin{equation}
  [\fm_{2}^{\Cat}(f,g)] = [{\bf e}_{X}] \in H^{*}( \Cat(X,X),  \fm_{1}^{\Cat})
\end{equation}
we say that $X$ and $Y$ are \emph{quasi-isomorphic}\index{quasi-isomorphic}.
\end{defn}

 %\marginpar{A paragraph removed.  KF. Feb 2025.}

\subsection{Matrix factorisations in dimension $0$} \label{sec:matr-fact-dimens}
In this section, we consider the category of matrix factorisations, which plays the same r\^ole, in the study of weakly curvature free $A_{\infty}$ categories, as the category of complexes does in the study of ordinary categories.

Recall that  a {\it matrix factorisation}\index{matrix factorisation} $V$ over a ring  $R$ consists of a $\Z_{2}$ graded free $R$-modules $V$, and a degree $1$  endomorphism $\partial_{V}$ which squares to a multiple  of the identity by a scalar $\Po(V)$.  We write  $V$ as the direct sum
\begin{equation}
  V^{0} \oplus V^{1}[1]
\end{equation}
and $\partial_{V}^{i}$ for the restriction of $\partial_{V}$ to $V^{i}[i]$.
\begin{rem}
In the algebraic literature, the condition of vanishing Euler characteristic is sometimes imposed on the category of matrix factorisations. For non-trivial values of the potential,  the condition that the differential square to a non-trivial multiple of the identity implies that the ranks of $V^0$ and $V^1$ agree, hence the Euler characteristic vanishes. For the trivial value, we shall not impose this, because  Floer cochain complexes $CF((L_1, b_1), (L_0, b_0))$ with $\mathfrak{PO}(L_0,b_0)=\mathfrak{PO}(L_1,b_1)$ need not satisfy this property.
%may  have not vanishing Euler characteristic.  
 \end{rem}

The collection of matrix factorisations forms a category \cite{eisenbud,Orlov},
in fact a curved differential graded category in which every object has curvature $\Po(V) \cdot \id_{V} $, and with vanishing higher multiplications.   The morphism spaces are the graded vector spaces $\Mat(V,W) =  \Hom_{R}(V,W)$,\index[syindex]{M1at@$\Mat(V,W)$} which are equipped with an endomorphism\index[syindex]{d1Mat@$d^{\Mat}$}
\begin{equation}\label{Matsign}
  d^{\Mat} \phi = \partial_{W} \circ \phi - (-1)^{\deg(\phi)} \phi \circ \partial_{V}
\end{equation}
and the composition is given by the usual composition of linear maps.  Note that we are using the standard conventions from  differential graded algebra, so that the integer $\deg(\phi)$ appearing in the above equation is the ordinary (not shifted) grading. %every morphism $\phi$ is not equipped with a shifted grading.

The tensor product of matrix factorisations equips $\Mat$ with a monoidal structure:  $V \otimes W$ is defined to be
\begin{equation}
   \left( (V^{0} \otimes W^{0}) \oplus (V^{1} \otimes W^{1}) \right) \oplus \left(   (V^{1} \otimes W^{0}) \oplus (V^{0} \otimes W^{1}) \right)[1]
\end{equation}
with differential
\begin{equation}
  \partial_{V \otimes W}(v \otimes w) = \partial_{V}(v) \otimes w   + (-1)^{\deg(v)}  v \otimes \partial_{W}(w)
\end{equation}
which has curvature $ \Po(V) +  \Po(W)   $.   In particular, we obtain an honest differential whenever the sum vanishes:
\begin{lem} \label{lem:tensor_product_matrix_vanishes}
  If $\Po(V) = r$ and $\Po(W) = -r$, then
  \begin{equation}
  r \cdot  H^{*}\left( V \otimes W, \partial_{V \otimes W} \right) = 0.
  %\partial_{V} \otimes \id \pm  \id \otimes \partial_{W} 
  \end{equation}
\end{lem}
\begin{proof}
Let $\partial^{0}_{V}$ denote the component of $\partial_{V}$ mapping $V^{0}$ to $V^{1}$.  The reader may directly compute that the endomorphism of $ V \otimes W $ which vanishes on $V^{1} \otimes W$ and is given by $\partial_{V}^{0} \otimes \id_{W} $  on $V^{0} \otimes W$ is a null-homotopy for $ r \cdot \id_{V \otimes W} $.
\end{proof}

 For each scalar $r \in R$, we obtain a curvature-free subcategory  $\Mat_{r}$.\index[syindex]{M1atr@$\Mat_{r}$}
\begin{lem} \label{lem:hom_factorisations_trivial}
The cohomology groups of all morphism spaces $\Mat_{r} $ are annihilated by $r$, i.e.:
\begin{equation}
  r \cdot H^{*}\left(\Hom_{R}(W,V), d^{\Mat} \right) = 0
\end{equation}
whenever $\Po(W) = \Po(V) = r$.
\end{lem}
\begin{proof}
It suffices to prove that the identity in $ H^{*}\left(\Hom_{R}(V,V), d^{\Mat} \right)  $ is annihilated by $r$.  Consider the map
\begin{equation}
\partial_{V}^{0}  \in \Hom_{R}(V^{0},V^{1}[1])    \subset \Hom_{R}(V,V).
\end{equation}
The equation $\partial_{V} \circ \partial_{V}  = r \cdot \id_{V} $ can be reinterpreted as $d^{\Mat}(\partial_{V}^{0}) = r \cdot \id_{V}$, i.e. that the class represented by $ \id_{V}  $  in $ H^{*}\left(\Hom_{R}(V,V), d^{\Mat} \right) $    vanishes when multiplied by $r$.  Since $  H^{*}\left(\Hom_{R}(V,V), d^{\Mat} \right)  $ is a unital associative algebra over $R$, we conclude the desired result.
\end{proof}

We conclude:
\begin{lem} \label{lem:factorisations_trivial}
If $R$ is a field, then the category $ \Mat_{r} $ is equivalent to the trivial category if and only if $r\ne0$. \qed
\end{lem}
%\blue{
%}\marginpar{Is this way of writing okay?  KO}

\subsection{$A_{\infty}$ functors}
\label{sec:a_infty-homom}

Let $\Cat$ and $\Dat$ be curved $A_{\infty}$ categories.  A \emph{curved $A_{\infty}$ functor}  $\Phi$\index{$A_{\infty}$ functor}
\index{curved $A_{\infty}$ functor} consists of a map $\Phi \co \Ob(\Cat) \to \Ob(\Dat)$, together with a map
\begin{equation}
  \Phi \co B(\Cat[1])(X,Y) \to \Dat [1] (\Phi(X),\Phi(Y))
\end{equation}
 satisfying the relation
\begin{equation}
\aligned
&\sum_{\substack{ 0 \leq k \leq \infty \\ \alpha(k)}} \fm_{k}^{\Dat}\left(\Phi(\bfa^{k;1}_{\alpha(k)}), \ldots,  \Phi(\bfa^{k;k}_{\alpha(k)}) \right)
\\
&=
\sum_{\alpha} (-1)^{\deg'(\bfa^{3;1})} \Phi\left(\bfa^{3;1}_{\alpha}, \fm^{\Cat}( \bfa^{3;2}_{\alpha}),\bfa^{3;3}_{\alpha}  \right).
\endaligned
\end{equation}
Whenever $\bfa$ has length $d$, $\alpha(k)  $ is an index running over all compositions of $d$ consisting of $k$ terms.  Writing $\Phi^{d}$ for the component of $\Phi$ defined on words of length $d$, the first above relation, for $d=0$, asserts that for each pair $(X,Y)$ of objects of $\Cat$, we have
\begin{equation}
\fm_{0}^{\Dat} (1)+  \sum_{k \geq 1} \fm_{k}^{\Dat}\left( \Phi^{0}(1), \ldots, \Phi^{0}(1) \right) = \Phi^{1}(\fm_{0}^{\Cat}(1)).
\end{equation}

\begin{rem}
The left hand side is an infinite sum so we need some assumption so that the sum is
well defined. In case we work over the Novikov ring $\Lambda_0$
we assume $\Phi^{0}(1) \in \Lambda_+$ and then it converges in the $T$-adic topology.
 %\marginpar{I put this remark Kenji.}
\end{rem}

In this paper, we shall only consider those functors for which the curvature $\Phi^{0}$ vanishes; we simply call such maps  \emph{$A_{\infty}$ functors}.  The map $\Phi^{1}$ is called the \emph{linear term}\index{linear term} of $\Phi$.

\begin{defn}
\begin{enumerate}
\item An \emph{$A_{\infty}$ quasi-embedding}\index{quasi-embedding} is an $A_{\infty}$ functor such that the map
  \begin{equation}
    \Phi^{1} \co \Cat[1](X,Y) \to \Dat [1] (\Phi(X),\Phi(Y))
  \end{equation}
is a quasi-isomorphism of chain complexes whenever the potential values of $X$ and $Y$ agree. 
 %\marginpar{In the filtered case 
%one might need to assume that $\Phi^{1}$ reduces to a chain homotopy equivalence over $R$ the ground ring, after reduction.
%Assume it for all $X,Y$.  KF. 2025 Dec.}
%for all pairs $(X,Y)$ of objects \blue{such that $\frak{PO}(X)=\frak{PO}(Y)$}.
 %\marginpar{Is this change okay?  See also Definitions 2.10, 2.11. KO}\marginpar
\item An \emph{$A_{\infty}$ quasi-isomorphism}\index{quasi-isomorphism} is an $A_{\infty}$ quasi-embedding such that the map $\Phi$ on objects is a bijection.
\item An \emph{$A_{\infty}$ quasi-equivalence}\index{quasi-equivalence}  is an $A_{\infty}$ quasi-embedding such that every object of the target is quasi-isomorphic to an object which is in the image of $\Phi$.
\end{enumerate}
\end{defn}
\begin{rem}
  In \cite{fooo09}, the word \emph{homotopy equivalence} was used instead of quasi-equivalence.
\end{rem}

An $A_{\infty}$ functor is  \emph{linear} if all higher order terms $\{ \Phi^{k} \}_{k \geq 2}$ vanish.
\begin{defn}
\begin{enumerate}
\item An \emph{$A_{\infty}$ embedding}\index{$A_{\infty}$ embedding} is a linear $A_{\infty}$ functor such $\Phi^{1}$ is an isomorphism of chain complexes for all pairs $(X,Y)$ of objects  with the same potential value.
\item An \emph{$A_{\infty}$ isomorphism}\index{$A_{\infty}$ isomorphism} is a linear $A_{\infty}$ quasi-isomorphism.
\end{enumerate}
\end{defn}

\subsection{Cyclic structure}
\label{sec:cycl-a_infty-algebra} 

The $A_{\infty}$-categories we consider have an additional structure which comes from Poincar\'e duality on the cohomology of a compact manifold:
\begin{defn} \label{defn:cyclic_structure}
  A \emph{\newred{cyclic} pairing of dimension}\index{cyclic pairing} $n$ on $A_{\infty}$ category $\Cat$ is a collection of linear maps
  \begin{equation} \label{eq:pairing_definition_CY}
    \langle \cdot\,,\, \cdot \rangle \co \Cat(X,Y) \otimes \Cat(Y,X) \to R
  \end{equation}
of degree $-n$, which are graded anti-symmetric with respect to the shifted grading in the sense that
\begin{equation} \label{eq:switch_factors_pairing}
  \langle x, y \rangle = (-1)^{1 + \deg'(x) \deg'(y)} \langle y, x \rangle
\end{equation}
and satisfying the cyclic invariance condition
\begin{equation} \label{eq:sign_cyclic_invariance}
  \langle \fm_{d}^{\Cat}(a_1, \ldots, a_{d}), a_0 \rangle = (-1)^{n \deg'(a_0)}  \langle \fm_{d}^{\Cat}(a_0, \ldots, a_{d-1}), a_d \rangle.
\end{equation}
We say that the pairing is \emph{\newred{non-degenerate}}\index{non-degenerate} if it descends to a perfect pairing on cohomology with respect to $\fm_1^{\Cat}$  for any pair $X,Y$ such that $\frak{PO}(X)=\frak{PO}(Y)$.
\par
An $A_{\infty}$ category $\Cat$ is called {\it cyclic}\index{cyclic} if it is equipped with a non-degenerate cyclic 
pairing.
\end{defn} 
\begin{rem}   %\marginpar{160127: MA clarified - no substantial change}
To clarify the statement on perfectness of the pairing on cohomology, we equip $ \Cat(X,Y)$ with shifted gradings and differential $\fm_1^{\Cat}$.   % When we regard $x \in \Cat (X,Y)$ as an element in $\Cat [1](X,Y)$, we denote it by $S^{1} x$.
 %  \marginpar{We use ${\mathfrak m}_k$ for operation on the map $B \Cat [1]$. Thus I changed the phrase here. KO} 
%; i.e. if we write $S^1 x$ for the element of $\Cat(X,Y)[1]$ 
%corresponding to $x \in \Cat(X,Y)$, then the image of $S^1x$ under the differential is $S^1  \fm_1^{\Cat}(x)$.   
Equations \eqref{eq:switch_factors_pairing} and \eqref{eq:sign_cyclic_invariance} and the Koszul conventions for the differential on a tensor product then imply that the map
\begin{align}\nonumber
  \Cat(X,Y) \otimes \Cat(Y,X) & \to R \\
x \otimes  y & \mapsto \langle x , y \rangle
\nonumber\end{align}
is a co-chain map, and hence induces a map on cohomology.
\end{rem}
\begin{rem} 
It is more common to find the sign
\begin{equation}\nonumber
  \deg'(a_0) \left( \sum_{i=1}^{d} \deg'(a_i) \right)
\end{equation}
in Equation \eqref{eq:sign_cyclic_invariance}.  This is equivalent to our sign, 
because the pairing is non-vanishing only when
\begin{equation}\nonumber
n = \deg(a_0) + \deg( \fm^{\Cat}_{d}(a_1, \ldots, a_{d})).
\end{equation}
Applying a similar consideration to Equation \eqref{eq:switch_factors_pairing} yields
\begin{equation}
  \label{eq:231bew}
  \langle x, y \rangle = (-1)^{1 + (n+1) \deg'(x) } \langle y, x \rangle
\end{equation}
\end{rem}

We shall find the following Lemma useful when associating a pairing on Hochschild homology to a cyclic structure:
\begin{lem}\label{lem316May26}
Given a pair $c$ and $d$ of morphisms, and sequences $\bfa$ and $\bfb$, we have
\begin{equation}\nonumber
  \label{eq:sign_cyclic_multiple_rotations}
  \langle \fm^{\Cat}_{r+s+1}(\bfb,  c, \bfa), d \rangle = (-1)^{1+ \deg'(c) +  n \deg'(\bfb)}  \langle  c ,  \fm^{\Cat}_{r+s+1}(\bfa, d, \bfb) \rangle .
\end{equation}
\end{lem}
\begin{proof}
 Applying Equation \eqref{eq:sign_cyclic_invariance} for $d$ and $\bf a$, we find that
 \begin{equation}\nonumber
       \langle \fm^{\Cat}_{r+s+1}(\bfb,  c, \bfa), d \rangle = (-1)^{n(\deg'(d) + \deg'(\bfa))} \langle  \fm^{\Cat}_{r+s+1}(\bfa, d, \bfb) , c \rangle.
 \end{equation}
Applying Equation \eqref{eq:231bew} gives
\begin{equation}\nonumber
  \langle  \fm^{\Cat}_{r+s+1}(\bfa, d, \bfb) , c \rangle = (-1)^{1 + (n+1) \deg'(c)}  \langle  c ,  \fm^{\Cat}_{r+s+1}(\bfa, d, \bfb) \rangle.
\end{equation}
Using the relation
\begin{equation}\nonumber
  n -3 = \deg'(c) + \deg'(d) + \deg'(\bfa) + \deg'(\bfb),    
\end{equation}
we deduce the desired result.
\end{proof}

Kajiura \cite{Ka} defined an appropriate notion of morphisms which respect cyclic structures:
\begin{defn} \label{def:cyclic-morphism}
If $\Bat$ and $\Cat$ are both cyclic $A_{\infty}$ categories, an $A_{\infty}$ functor $\Phi \co \Bat \to \Cat$ is \emph{cyclic}\index{cyclic} if
\begin{equation} \label{eq:cyclic_morphism}
\sum_c  \langle \Phi(\bfx^{2:1}_c), \Phi(\bfx^{2:2}_c) \rangle = \begin{cases}  \langle x, y \rangle  & \textrm{if } \bfx = x\otimes y \\
0 & \textrm{otherwise.}
\end{cases}
\end{equation}
Here we assume ${\bf x}$ is of the form ${\bf x} = x_1\otimes \dots \otimes x_n$ and put $\Delta{\bf x} = \sum_c \bfx^{2:1}_c \otimes \bfx^{2:2}_c$.
 %\marginpar{Formula corrected.  See [Fu2] p556. KF 2024 Dec.}
\end{defn}

\section{Hochschild (co)-homology.}
\label{sec:Hochschild}
 Most of our focus will be on Hochschild homology and cohomology as invariants of $A_{\infty}$ categories which depend only on their quasi-isomorphism type, and which are respectively obtained by deriving the tensor product with the diagonal bi-module, and its space of morphisms.

Since the setting is not exactly standard outside of symplectic topology, we shall give a complete definition of the Hochschild homology and cohomology for a curved $A_{\infty}$ category $\Cat$ in which all objects are weakly unobstructed.  We then prove the following two results, which are well known in slightly different contexts:
\begin{prop}
The  Hochschild cohomology of $\Cat$ is an algebra with respect to the cup product, over which Hochschild homology admits a natural module structure via the cap product.
\end{prop}
\begin{proof}
The algebra structure is established in Lemma \ref{lem:cup_product_HH}, and the module structure in Lemma \ref{lem:cap_product_HH}.
\end{proof}

Next, we assume that $\Cat$ admits a  pairing of dimension $n$ (see Definition \ref{defn:cyclic_structure}), and prove the following result (see Lemma \ref{lem:pairing_HH} for a precise formula).
\begin{prop} \label{prop:pairing_HH_cyclic}
The  pairing on $\Cat$ induces a pairing
\begin{equation}
  HH^{*}(\Cat,\Cat) \otimes HH_{*-n}(\Cat,\Cat) \to \Lambda.
\end{equation}
If $\Cat$ is cyclic, this pairing is perfect.
\qed
\end{prop}

\subsection{$A_{\infty}$ bi-modules} \label{sec:a_infty-bi-modules}
The proper setting for defining Hochschild (co)-homology are categories of bi-modules.  Given  curved $A_{\infty}$ categories $\Aat$ and $\Bat$, an object of the category of $\Aat$-$\Bat$-bi-modules\index{bi-module} consists of $R$ modules $\cP(X,Y)$ for each pair of objects $X \in \Ob(\Aat)$ and $Y \in \Ob(\Bat)$, together with structure maps $\fm_{r|1|s}^{\cP}$ defined whenever  $r$ and $s$ are non-negative integers, which generalize the notion of left and right multiplication.  Our conventions for signs agree with those in \cite{seidel2007}.
\begin{rem}
The reader may notice that all bi-modules we shall consider will have vanishing curvature (this is implicit in the definition as the curvature would correspond to $\fm_{0|0|0}$).  There is a notion of \emph{curved bi-module}, which has a r\^{o}le to play in symplectic topology, but with whose investigation we shall not be concerned.
\end{rem}

In order to define the structure maps, it is convenient to introduce the double sided bar complex of bi-modules
\begin{align} \label{eq:bar-cplex_bi-module}
  B(\Aat[1],\cP,\Bat[1]) (X, Y) & = \bigoplus_{r=0}^{\infty} \bigoplus_{s=0}^{\infty} B_{r,1,s}(\Aat[1],\cP,\Bat[1])(X, Y) \\
B_{r,1,s}(\Aat[1],\cP,\Bat[1]) (X, Y) & =   \bigoplus_{X_{0}, Y_0 }  B_{r}(\Aat[1])(X,X_0) \otimes \cP(X_0, Y_0) \otimes  B_{s}(\Bat[1])( Y_{0}, Y).
\end{align}
With this notation\index[syindex]{B1r1s@$B_{r,1,s}(\Aat[1],\cP,\Bat[1])$} at hand, the structure\index[syindex]{mxmr1sPB@$\fm_{r\vert 1\vert s}^{\cP}$} maps are linear maps of degree $1$:
\begin{equation}
  \label{eq:bi-module_structure_maps}
  \fm_{r|1|s}^{\cP} \co  B_{r,1,s}(\Aat[1],\cP,\Bat[1]) (X, Y)  \to  \cP(X, Y)
\end{equation}
whose direct sum will be denoted by $\fm^{\cP}  $ when no confusion can arise.   The cases $(r,s)=(1,0)$ and $(r,s)=(0,1)$ correspond respectively to left and right multiplication in the setting of ordinary algebras, and the condition that the multiplications commute is replaced by the sequence of equations
\begin{equation}
  \label{eq:bi-module_A_oo_equation}
  \aligned
&0 = \sum (-1)^{\deg'(\bfa^{2;1}_{\alpha(2)})}   \fm^{\cP}(\bfa^{2;1}_{\alpha(2)}\otimes \fm^{\cP}(\bfa^{2;2}_{\alpha(2)} \otimes p \otimes \bfb^{2;1}_{\beta(2)} ) \otimes \bfb^{2;2}_{\beta(2)}  ) \\
&+ \sum  (-1)^{\deg'(\bfa) + \frak{deg}(p) + \deg(\bfb^{3;1}_{\beta(3)}) }   \fm^{\cP}\left(\bfa \otimes  p \otimes  \bfb^{3;1}_{\beta(3)} \otimes  \fm^{\Bat}(  \bfb^{3;2}_{\beta(3)} ) \otimes   \bfb^{3;3}_{\beta(3)} \right) \\
&+ \sum (-1)^{ \deg'( \bfa^{3;1}_{\alpha(3)})  }  \fm^{\cP}(\bfa^{3;1}_{\alpha} \otimes \fm^{\Aat}( \bfa^{3;2}_{\alpha}) \otimes \bfa^{3;3}_{\alpha} \otimes  p \otimes \bfb ),
\endaligned
\end{equation}
where the signs are given by a Koszul sign corresponding to permuting the second structure map past the inputs which precede it.
Here and hereafter we use the symbol $\frak{deg}(p)$\index[syindex]{degfrak@$\frak{deg}$}  for the degree of an element $p$ of bi-module. %\marginpar{Notation 
%$\frak{deg}$ is introduced. KF 2025 Feb}
\begin{rem}\label{difference of Lag PO}
The first instance of Equation \eqref{eq:bi-module_A_oo_equation} implies that $  \cP(X, Y) $ is a matrix factorisation with potential value
\begin{equation}
  \Po(Y) - \Po(X).
\end{equation}
\end{rem}
We next consider the shift functor on bi-modules: 
 %\marginpar{MA: I would like not to use this at all} 
recall that for each integer $n$, one defines the $n$-fold shift of a graded module $V^{\bullet}$ via the formula
\begin{equation}
  V[n]^{i} = V^{i+n}.
\end{equation}
We then define the $n$-fold shift\index{nfoldshift@$n$-fold shift} $\cS^{n} \cP$\index[syindex]{Scaln@$\cS^{n} \cP$} of a bi-module $\cP$ to have underlying $R$-modules $ \cP(X,Y)[n] $, and operations
\begin{equation} \label{eq:shift_bi-module_operations}
   \fm^{\cS^{n} \cP}(\bfa \otimes S^{n }p \otimes \bfb  ) = (-1)^{n (1+ \deg'(\bfa))}  S^{n} \fm^{\cP}(\bfa \otimes  p \otimes \bfb ).
\end{equation}

\subsection{Morphisms and tensor products of bi-modules} \label{sec:morph-tens-prod}
The set of bi-modules forms a category, where morphisms\index{morphism of bi-module} are computed using a $2$-sided bar complex.  Given a pair of bi-modules $\cP$ and $\cQ$, we define the space of morphisms\index{morphism of bi-module} to be
\begin{equation}
\Aat-\Bat(\cP,\cQ) = \prod_{\substack{X \in \Ob(\Aat) \\  Y \in \Ob(\Bat)}} \Hom_{R}\left( B(\Aat[1],\cP,\Bat[1]) (X, Y) , \cQ(X, Y) \right).
 \end{equation}
We write an element of this morphism space as a sequence $\{ T^{r|1|s}\}_{r,s=0}^{+\infty}$, with $T^{r|1|s}$ defined on the direct sum of $B_{r,1,s}(\Aat[1],\cP,\Bat[1])(X, Y) $ over all pairs $( X,   Y) $, and use the structure maps of $\cP$ and $\cQ$ to define the differential
% The signs come from the fact that the signs in the $A_{\infty}$ structure come from assuming that all inputs are shifted by $1$.  If you shift $p$ by one, you have to permute that shift past all the elements to the left of it.
\begin{equation} \label{eq:differential_map_bi-modules}
\aligned
&(\partial T)( \bfa,  p, \bfb) = \\
&\sum (-1)^{\deg(T)  \deg'(\bfa^{2;1}_{\alpha(2)} )}   \fm^{\cQ}(\bfa^{2;1}_{\alpha(2)} \otimes T(\bfa^{2;2}_{\alpha(2)} \otimes p \otimes\bfb^{2;1}_{\beta(2)}) \otimes \bfb^{2;2}_{\beta(2)}   ) \\
&+ \sum (-1)^{\deg'(\bfa^{2;1}_{\alpha(2)} )+\deg(T)+1}  T(\bfa^{2;1}_{\alpha(2)} \otimes \fm^{\cP}(\bfa^{2;2}_{\alpha(2)} \otimes p \otimes \bfb^{2;1}_{\beta(2)}) \otimes \bfb^{2;2}_{\beta(2)} ) \\
&+ \sum  (-1)^{\deg'( \bfa) + \frak{deg}(p) + \deg'(\bfb^{3;1}_{\beta(3)} ) +\deg(T)+1} T(\bfa \otimes  p  \otimes \bfb^{3;1}_{\beta(3)} \otimes \fm^{\Bat}(\bfb^{3;2}_{\beta(3)} ) \otimes \bfb^{3;2}_{\beta(3)}  )  \\
&+ \sum (-1)^{\deg'(\bfa^{3;1}_{\alpha(3)} )+\deg(T)+1} T(\bfa^{3;1}_{\alpha(3)} \otimes \fm^{\Aat}( \bfa^{3;2}_{\alpha(3)}) \otimes \bfa^{3;3}_{\alpha(3)} \otimes  p \otimes \bfb ).
\endaligned
\end{equation}
\begin{defn}
  A \emph{quasi-isomorphism} of bi-modules of degree $n$ is a morphism $T$ of degree $n$ such that $T^{1}$ is a quasi-isomorphism
  \index{quasi-isomorphism of bi-module} of matrix factorisations in the sense of Definition \ref{def:quasi-isomorphic}, for every pair of objects $X$ and $Y$.
%in the sense of Definition \blue{i.e., an isomorphism in the homotopy category of matrix factorisations}, 
   Such a quasi-isomorphism is \emph{naive} if all higher order terms vanish.
\end{defn}
\begin{rem}
Note that a naive quasi-isomorphism $T$ of degree $n$ satisfies:
\begin{equation} \label{eq:sign_for_naive_quasi-iso-degree-n}
 (-1)^{\deg(T)(1+  \deg'(\bfa))}   \fm^{\cQ}(\bfa \otimes T( p ) \otimes \bfb  ) = T(\fm^{\cP}(\bfa \otimes p \otimes \bfb)).
\end{equation}
\end{rem}

We shall also use the tensor product of bi-modules:  Given $\cP$ and $\cQ$, which are respectively $\Aat-\Bat$ and $\Bat-\Cat$ bi-modules, we define an $\Aat-\Cat$ bi-module $\cP \tensor_{\Bat} \cQ$.  The $R$-modules underlying $\cP \tensor_{\Bat} \cQ$ 
are\index[syindex]{PtensorBQ@$\cP \tensor_{\Bat} \cQ$}
\begin{equation}\label{tensorprodbi-modules}
  \cP \tensor_{\Bat} \cQ (X,Y) = \bigoplus_{\substack{d \geq 0 \\ Z_0, Z_d  \in \Ob(\Bat)}}  \cP(X,Z_0) \otimes B_{d}(\Bat[1])(Z_0,Z_d) \otimes \cQ (Z_d,Y).
\end{equation}

The differential is given by
\begin{equation}\label{eq:differential_derived_tensor_product}
\aligned
  &\fm^{  \cP \tensor_{\Bat} \cQ }( p \otimes \bfb \otimes q) \\
=   &\sum \fm^{\cP}\left( p \otimes \bfb^{2;1}_{\beta(2)}\right) \otimes \bfb^{2;2}_{\beta(2)} \otimes q  \\
&+ \sum (-1)^{ \deg(p) + \deg'(\bfb^{3;1}_{\beta(3)}) } p \otimes \bfb^{3;1}_{\beta(3)} \otimes \fm^{ \Bat }\left(  \bfb^{3;2}_{\beta(3)} \right)   \otimes  \bfb^{3;3}_{\beta(3)} \otimes q \\
 &+  \sum (-1)^{ \deg(p) + \deg'(\bfb^{2;1}_{\beta(2)}) }  p \otimes \bfb^{2;1}_{\beta(2)} \otimes \fm^{ \cQ }\left(  \bfb^{2;2}_{\beta(2)} \otimes q \right).
\endaligned
\end{equation}

The structure maps incorporate the internal differential in $B_{d}(\Bat[1])(Z_0,Z_d)   $ together with the higher products $ \fm_{0|1|r}^{\cP}$ and $ \fm_{s|1|0}^{\cQ}$. (We remark $ \fm_{s|1|r}^{\cQ} = 0$ if $rs \ne 0$.)
\begin{align} \label{eq:structure_maps_bi-module}
  \fm^{  \cP \tensor_{\Bat} \cQ }( \bfa, p, \bfb, q, \bfc) & = 0 \textrm{ if } \bfa \neq 0 \textrm{ and }
\bfc \neq 0 \\ \label{eq:structure_maps_bi-module2}
 \fm^{  \cP \tensor_{\Bat} \cQ }( \bfa, p, \bfb, q) & =  \sum \fm^{  \cP }\left( \bfa \otimes p \otimes \bfb^{2;1}_{\beta(2)}\right) \otimes \bfb^{2;2}_{\beta(2)} \otimes q \\
\label{eq:structure_maps_bi-module3}
\fm^{  \cP \tensor_{\Bat} \cQ }( p, \bfb, q, \bfc) & = \sum (-1)^{ \frak{deg}(p) + \deg'(\bfb^{2;1}_{\beta(2)}) }  p \otimes \bfb^{2;1}_{\beta(2)} \otimes \fm^{ \cQ }\left(  \bfb^{2;2}_{\beta(2)} \otimes q \otimes \bfc \right).
\end{align}

%STATE (Prove?) DERIVED INVARIANCE OF TENSOR PRODUCT?

\subsection{Examples of bi-modules}
\label{sec:examples-bi-modules}

Given a category $\Cat$, we first consider the \emph{shifted diagonal bi-module}\index{diagonal bi-module} $\Cat_\Delta$, denoted by $\Cat$ if no confusion arises,\index[syindex]{CzDelta@$\Cat_\Delta$}
which assigns to every pair of objects the graded module $ \Cat (X,Y)$ equipped with the shifted gradings, for which the structure maps are given by
\begin{equation} \label{eq:product_diagonal_bi-module}
  \fm_{r|1|s}^{\Cat_\Delta} (  \bfa \otimes  p \otimes \bfb  ) =  \fm_{r+1+s}^{\Cat} ( \bfa \otimes  p \otimes  \bfb ).
\end{equation}
In other words $\frak{deg}(p) = \deg'p$, when a morphism $p$ of $\Cat$ is regarded as element of diagonal module.
 %\marginpar{
%The last two lines added. KF. 2025 Feb.}
\begin{rem}
While it would be abstractly more natural to discuss the diagonal bi-module first, our sign conventions are arranged so that this would introduce a Koszul sign of $(-1)^{\deg' \bfa}$ in the above equation. 
%In order to understand the above sign, it is helpful to keep in mind that the sign conventions we use are the Koszul signs associated to thinking of the $A_{\infty}$ operations as differentials on a complex built from copies of the morphism spaces \emph{shifted by $1$.}  When constructing a bi-module from morphism spaces, we take a morphism in $p$ with its natural grading, which accounts for the difference in signs.  The sign in Equation \eqref{eq:product_diagonal_bi-module} differs from the one appearing in \cite[Definition 3.7.73]{fooo09}.  The difference arises because of the different shift in \cite[Proposition 3.7.71]{fooo09}.
\end{rem}

Note that, whenever $\Aat$ and $\Bat$ are subcategories of $\Cat$, we may restrict this bi-module to an $\Aat-\Bat$-bi-module which we denote ${_\Aat}\Cat_{\Bat}$ when confusion might arise.

We next consider the tensor product of the left Yoneda module\index{Yoneda module} \newred{${_\Cat}\cZ$}
and right Yoneda module \newred{$\cZ_{\Cat}$} associated to an object $Z$ of $\Cat$: 
\begin{equation}
  \left( {_\Cat}\cZ \otimes   \cZ_{\Cat} \right) ( X,Y)  = \Cat(X,Z) \tensor_{R}  \Cat(Z, Y),
\end{equation}
where each factor in the right hand side is equipped with shifted gradings.  The differential is given by
\begin{equation} \label{eq:differential_tensor_yoneda}
   \fm_{0|1|0}^{ {_\Cat}\cZ \otimes   \cZ_{\Cat}} (p_{l} \otimes p_{r}) =  \fm_{1}^{\Cat}(p_{l}) \otimes p_{r}   + (-1)^{\deg'(p_{l})}p_{l} \otimes  \fm_{1}^{\Cat}(p_{r})
\end{equation}
while the higher structure maps vanish whenever both $r$ and $s$ are non-zero, and are otherwise given by
\begin{align} \label{eq:higher_product_tensor_yoneda}
  \fm_{r|1|0}^{ {_\Cat}\cZ \otimes   \cZ_{\Cat}} (\bfa \otimes  p_{l} \otimes p_{r}) & = %(-1)^{1+\deg'(\bfa)}
 \fm^{\Cat}(\bfa \otimes p_{l}) \otimes p_{r} \\
  \fm_{0|1|s}^{ {_\Cat}\cZ \otimes   \cZ_{\Cat}} ( p_{l} \otimes p_{r} \otimes \bfb  ) & = (-1)^{\deg'(p_{l})} p_{l} \otimes  \fm^{\Cat} (p_{r},\bfb).
\end{align}
The verification that the operations $   \fm_{r|1|s}^{ {_\Cat}\cZ \otimes   \cZ_{\Cat}} $ satisfy Equation \eqref{eq:bi-module_A_oo_equation} is a direct computation using the $A_{\infty}$ relations on $\Cat$.  The key fact is that the differentials on $ \Cat(X,Z) $  and $ \Cat(Z,X) $ square to opposite multiples of the identity.

Note that this example can be generalised to the following setting:  Given a subcategory $\Bat \subset \Cat$, we obtain a $\Cat$-bi-module $\Bat \tensor_{R} \Bat$ by taking the direct sum of the above bi-modules for all objects in $\Bat$.

We shall also consider a pullback construction for bi-modules:  Whenever $\Phi \co \Bat \to \Cat$ is an $A_{\infty}$ homomorphism, and $\cP$ a $\Cat$-bi-module, we define a $\Bat$ bi-module $\Phi^{*}(\cP)$ to have underlying modules
\begin{equation}\nonumber
  \Phi^{*}(\cP)(X,Y) = \cP(\Phi(X), \Phi(Y)),
\end{equation}
and operations
\begin{equation}\nonumber
  \fm^{ \Phi^{*}(\cP)  }(\bfa \otimes p \otimes \bfb) = \sum_{\alpha, \beta} \fm^{\cP}\left(\Phi(\bfa^{k;1}_{\alpha}) \otimes \cdots \otimes \Phi(\bfa^{k;k}_{\alpha}) \otimes p \otimes   \Phi(\bfb^{l;1}_{\beta}) \otimes  \cdots  \Phi(\bfb^{l;l}_{\beta}) \right)
\end{equation}

Note that, in the special case $\cP$ is the diagonal bi-module, we have a natural map of $\Bat$-bi-modules
\begin{align} \label{eq:natural_map_to_pullback_bi-module}
\Nat_{\Phi} \co  \Bat_\Delta & \to  \Phi^{*}(\Cat_\Delta)  \\
\Nat_{\Phi}(\bfa \otimes x \otimes \bfb) & \mapsto % (-1)^\text{{\rm deg}' ({\bf a})} 
\Phi( \bfa \otimes x \otimes \bfb )
\end{align}
 %\marginpar{sign must appears when we convert an $A_{\infty}$-morphism to a morphism of $A_{\infty}$-bi-modules. KO}

\subsection{Duality for bi-modules}
\label{sec:duality-bi-modules}

If $\cP$ is an $\Aat-\Bat$-bi-module, we also define the ($R$-linear) dual $\Bat-\Aat$-bi-module $\cP^{\vee}$ to have underlying $R$-modules
\begin{equation}
  \cP^{\vee}(Y,X) = \Hom_{R}( \cP(X,Y), R)
\end{equation}
and operations\index[syindex]{Pvee@$ \cP^{\vee}$}
\begin{equation} \label{eq:structure_maps_dual_bi-module}
  \fm^{\cP^{\vee}}(\bfb \otimes  \psi \otimes \bfa)(q) = (-1)^{\frak{deg}(\psi) +1}  \psi\left( \fm^{\cP}(\bfa \otimes  q \otimes \bfb)\right).
\end{equation}
\begin{rem} \label{rem:use_pairing_simplify_sign}
The sign in Equation \eqref{eq:structure_maps_dual_bi-module} can be derived from the Koszul conventions by requiring that the evaluation map
\begin{equation}\nonumber
    \cP^{\vee}(Y,X) \otimes   \cP(X,Y) \to R  
\end{equation}
extend to the pairing on Hochschild (co)-homology groups which we will later introduce in Equation \eqref{eq:pairing_cyclic_chains_cochains}.  The reader interested in going through the derivation will use the fact that the right hand side in Equation \eqref{eq:structure_maps_dual_bi-module} vanishes unless
\begin{equation}\nonumber
    \frak{deg}(\psi) = 1 + \deg'(\bfa) +\frak{deg}(q) + \deg'(\bfb).
\end{equation}
\end{rem}

We now relate the pairing appearing in the cyclic structure to duality of bi-modules:

\begin{lem} \label{cor:cyclic_structure_dual_bi-modules}
A cyclic structure on $\Cat$ induces a degree $-n$ quasi-isomorphism of bi-modules
\begin{align}
  \label{eq:3}
   \psi \co \Cat_\Delta & \to  \Cat_\Delta^{\vee} \\
p & \mapsto \psi(p)=  \langle p, \_ \rangle.
\end{align}
\end{lem}
\begin{proof}
According to Equation \eqref{eq:sign_for_naive_quasi-iso-degree-n}, we must prove that
\begin{align} \label{dual3.31}
\fm^{\Cat_\Delta^{\vee}}(\bfa \otimes \psi( p ) \otimes \bfb  )(q) & =   (-1)^{n(1+  \deg'(\bfa))}    \psi(\fm^{\Cat}(\bfa \otimes p \otimes \bfb))(q) \\
& =  (-1)^{n(1+  \deg'(\bfa))}   \langle \fm^{\Cat}(\bfa \otimes p \otimes \bfb),  q \rangle \label{dual3.312}
\end{align}
for each morphism $q$. On the left hand side, we use the definition of the dual bi-module to compute
\begin{align}\nonumber
\fm^{\Cat_{\Delta}^{\vee}} ({\bf a} \otimes \psi (p) \otimes {\bf b}) (q) & = (-1)^{1 + n + \deg' p} \langle p,  \fm^{\Cat}({\bf b} \otimes q \otimes {\bf a}) \rangle.
\end{align}
Putting this sign together with that appearing in Equation \eqref{dual3.312}, we get
\begin{equation}\nonumber
   1 + \deg' p + n \deg'(\bfa)
\end{equation}
which is exactly the sign appearing in Lemma \ref{lem316May26}. \end{proof}
 
If $\cP^{\vee}$ and $\cQ^{\vee}$ are $\Bat-\Aat$-bi-modules which are respectively dual to $\cP$ and $\cQ$, we can define the dual of a bi-module morphism
\begin{equation}\nonumber
  T \in  \Aat-\Bat(\cP,\cQ)
\end{equation}
via the formula 
\begin{equation} \label{eq:dual_map_bi-modules}
T^{\vee}(  \bfb \otimes \phi \otimes  \bfa) (p) = (-1)^{\deg'(\bfb)+ \deg(T) (\deg'(\bfb) + \frak{deg}(\phi))} \phi\left(  T(\bfa \otimes p \otimes \bfb)  \right).
\end{equation}
\begin{rem}
The same principle as in Remark \ref{rem:use_pairing_simplify_sign} determines the sign; namely, we require that the pairings of Hochschild (co)-cochain complexes be natural.
\end{rem}

\subsection{Hochschild cohomology}
\label{sec:hochsch-cohom}
Given a bi-module $\cQ$ over a category $\Cat$, the Hochschild cohomology of $\Cat$ with coefficients in $\cQ$ can be defined as the homology of the space of bi-module maps from $\Cat$ to $\cQ$.  However,  in order to maintain consistency with previous accounts on the r\^ole of this invariant in Floer theory, we use, instead, the  Hochschild complex.  
Given a category $\Cat$, let us denote by $\Cat_{r}$ the full subcategory whose objects have potential value $r$.
  
\begin{equation}\nonumber
CH^{*}(\Cat_{r},\cQ) = \prod_{\substack{X \in \Ob(\Cat_{r}) \\  Y \in \Ob(\Cat_{r})}} \Hom_{R}\left( B(\Cat_{r}[1]) (X, Y) , \cQ(X, Y) \right)
\end{equation}
and\index[syindex]{CH*C@$CH^{*}(\Cat, \cQ)$}
\begin{equation}\nonumber
CH^{*}(\Cat, \cQ) = \prod_{r \in R}  CH^{*}(\Cat_{r},\cQ).
\end{equation}

\begin{rem}
In Subsection \ref{frakqondeRham} we also use  a different version 
$$
CH^{\prime *}(\Cat,\cQ)
= \prod_{\substack{X \in \Ob(\Cat) \\  Y \in \Ob(\Cat)}} \Hom_{R}\left( B(\Cat[1]) (X, Y) , \cQ(X, Y) \right)
$$
of Hochschild cochain complex. There is an obvious map 
$$
CH^{*}(\Cat, \cQ) \to CH^{\prime *}(\Cat,\cQ)
$$
which we can easily prove to be a chain map. Based on the vanishing results from Section \ref{sec:matr-fact-dimens}, we expect that this map induces 
isomorphism in cohomology. However we do not try to prove it in this paper since we do not use it.
\end{rem}

We write an element of this morphism space as a sequence $\{ T^d \}_{d=0}^{\infty}$, with $T^{d}$ defined on the direct sum of $B_{d}(\Cat[1])(X, Y) $ over all pairs $( X,   Y) $, and use the structure maps of $\Cat$ and $\cQ$ to define the differential
\begin{multline} \label{eq:diff_CC^*}
  (\partial T)(\bfa)  = \sum (-1)^{\deg(T)\deg'( \bfa^{3;1}_{\alpha(3)}) }   \fm^{\cQ}(\bfa^{3;1}_{\alpha(3)} \otimes T( \bfa^{3;2}_{\alpha(3)}) \otimes \bfa^{3;3}_{\alpha(3)} )  \\
+ \sum  (-1)^{\deg'( \bfa^{3;1}_{\alpha(3)})  +\deg(T)+1}   T(\bfa^{3;1}_{\alpha(3)} \otimes  \fm^{\Cat}( \bfa^{3;2}_{\alpha(3)}) \otimes \bfa^{3;3}_{\alpha(3)} ).
\end{multline}
% The extra +1 in the last two signs comes from the fact that a chain map satisfies df-fd=0
% Note that the difference with Equation 1.8 in Seidel is exactly accounted for by the use of  \fm_{k|1|d-\ell} instead of  \fm_{k+1+d-\ell}
\begin{defn}\label{defn412}
The \emph{Hochschild cohomology of $\Cat$ with coefficients in a bi-module}\index{Hochschild cohomology} $\cQ$ is $\cQ$\index[syindex]{HHupper@$HH^{*}(\cdot,\cdot)$}
defined by:
 \begin{equation}
   \label{eq:HH_co-definition}
   HH^{*}(\Cat,\cQ) = H^*(  CH^{*}(\Cat,\cQ), \partial).
 \end{equation}
\end{defn}

We write the Hochschild cochain complex (resp. Hochschild cohomology) with coefficients in the diagonal bi-module as 
$CH^*(\Cat, \Cat)$ (resp. $HH^*(\Cat, \Cat)$) instead of $CH^*(\Cat, \Cat_\Delta)$ (resp. $HH^*(\Cat, \Cat_\Delta)$).

A map $\Psi$ of $\Cat$-bi-modules induces a map on cyclic cochains:
\begin{equation} \label{eq:induced_map_CH^*}
  CH^{*}(\Psi)(T)( \bfb ) = \sum (-1)^{\deg(T) \deg'(\bfb^{3;1}_{\beta})  }\Psi\left(  \bfb^{3;1}_{\beta} \otimes T(\bfb^{3;2}_{\beta}) \otimes \bfb^{3;3}_{\beta} \right).
\end{equation}

Hochschild cohomology is an invariant of the quasi-equivalence type of the coefficient bi-module:
\begin{lem} \label{lem:naive_equivalence_iso}
Equivalences of bi-modules induce isomorphisms of Hochschild cohomology groups.
\end{lem}
\begin{proof}
We prove the statement for  each $r \in R$.
  The length filtration on the bar complex gives rise to a decreasing filtration on $CH^{*}(\Cat_{r},\cQ)$, whose $N$\th step consisting of those sequences  $\{ T^d \}_{d=0}^{\infty}$ for which $T^{0} = \ldots = T^{N} = 0$.  The $E^{2}$ page of the associated  spectral sequence is
\begin{equation}
 \prod  H^{*}  \Hom_{R}\left( B_{N}(\Cat_{r}[1]) (X, Y)) , \cQ(X, Y)) \right).
\end{equation}
The assumption  implies that the spectral sequences on this page is isomorphic.
\end{proof}

It is well known that an $A_{\infty}$ homomorphism does not induce maps of Hochschild cohomology groups with coefficients in the diagonal bi-module.  However, given a $\Cat$ bi-module $\cQ$, such a homomorphism $\Phi \co \Bat \to \Cat$ induces a map from the Hochschild cohomology of $\Cat$ with $\cQ$ coefficients to the cohomology of $\Bat$ with coefficients in the pullback bi-module.  At the chain level, the map is given by
\begin{equation}
\aligned
CH^{*}_{\cQ}( \Phi) \co  CH^{*}(\Cat, \cQ) & \to CH^{*}(\Bat, \Phi^{*}( \cQ))  \\
CH^{*}_{\cQ}( \Phi)(T)(\bfb) & = \sum_{\beta} T \left( \Phi( \bfb^{k;1}_{\beta} ) \otimes \cdots \otimes \Phi (\bfb^{k;k}_{\beta} ) \right).  \label{eq:pulback_bi-module_map_cohomology}
\endaligned
\end{equation}
Using the length filtration as in the proof of Lemma \ref{lem:naive_equivalence_iso}, we conclude:
\begin{lem} \label{lem:quasi-iso-categories-cohomology-iso}
  If $\Phi$ is a quasi-isomorphism, then $HH^{*}_{\cQ}( \Phi)$ is an isomorphism for every bi-module $\cQ$.
\end{lem}

\begin{comment}

\begin{equation}
  CH^{*}(\Cat,\Cat) \to  CH^{*}(\Cat_{r},\Cat_{r}).
\end{equation}

\begin{lem} \label{lem:HH-splits-product}
  If $R$ is a field, the induced map
  \begin{equation}
    HH^{*}(  \Cat,\Cat) \to  \prod_{r \in R}  HH^{*}(\Cat_{r},\Cat_{r})
  \end{equation}
is an isomorphism.
\end{lem}
\begin{proof}
Consider the map on the spectral sequence associated to the length filtration.  The $E_{2}$ page of the spectral sequence for the Hochschild cohomology of $\Cat$ is
\begin{equation}
 \prod  H^{*}  \Hom_{R}\left( B_{N}(\Cat[1]) (X, Y)) , \Cat(X, Y)) \right).
\end{equation}
The result therefore follows if we can prove that (i) the above factors vanish unless $\Po(X) = \Po(Y)$, and (ii) if $  \Po(X) = \Po(Y) = r $, the inclusion
\begin{equation}
 B_{N}(\Cat_{r}[1]) (X, Y)) \to B_{N}(\Cat[1]) (X, Y))
\end{equation}
induces a quasi-isomorphism of chain complexes.  In the first case, $\Cat(X, Y)$   and $ B^{(N)}(\Cat[1]) (X, Y)) $ are a matrix factorisations with non-vanishing potential, so that (i) follows from Lemma \ref{lem:factorisations_trivial}.  To prove the second part, we note that the quotient of $B_{N}(\Cat[1]) (X, Y))  $  by the image of $ B_{N}(\Cat_{r}[1]) (X, Y))  $ is a direct sum of tensor products of matrix factorisations with non-zero potential values.  Lemma \ref{lem:tensor_product_matrix_vanishes} implies that each summand has trivial homology, hence that the inclusion is a quasi-isomorphism.
\end{proof}
\end{comment}

\subsection{Hochschild homology}
\label{sec:hochschild-homology}

To each pair of bi-modules $\cP$ and $\cQ$ one may assign a derived tensor product $\cP \tensor_{\Catbi} \cQ$.  Whenever the first bi-module is the diagonal, one can define this using the Hochschild chain complex\footnote{Hochschild chain complex here 
is called cyclic bar complex in various literatures.}\index{Hochschild chain complex}
 %\marginpar{The name cyclic bar complex 
%and Hochschild chain complex should be unified. In Part 3 etc. Hochschild boundary is not $b$ but is $\delta_H$. KF 2025.
%$B_{k}^{\text{\rm cyc}}$ is $CH_*$ in certain other places.  KF 2025 Aug  These notations are unified.  KF 2025 Aug}
\index[syindex]{CH*CQ@$CH_{*}(\Cat, \cQ)$}\index[syindex]{CHk@$CH_{k}(\Cat_{r}[1], \cQ)$}
\begin{align} \label{eq:cyclic_bar-cplex}
CH_{*}(\Cat, \cQ) & = \bigoplus_{r \in R} \bigoplus _{k=0} CH_k(\Cat_r[1], \cQ) \\
CH_k(\Cat_{r}[1], \cQ)  & =  \bigoplus_{X , Y}  \cQ(Y,X) \otimes B_{k}(\Cat_{r}[1])(X,Y)
\end{align}
with differential\footnote{The Hochschild differential $\delta_H$ is written as $b$ in various literatures.}
\begin{equation}\nonumber
\aligned
&\delta_H(q \otimes \bfa) \\
=&  \sum (-1)^{\frak{deg}(q) + \deg'(\bfa^{3;1}_{\alpha(3)} ) }  q \otimes  \bfa^{3;1}_{\alpha(3)} \otimes \fm^{\Cat}( \bfa^{3;2}_{\alpha(3)} ) \otimes \bfa^{3;3}_{\alpha(3)} \\
&+   \sum  (-1)^{\deg'( \bfa^{3;3}_{\alpha(3)}  ) \left( \frak{deg}(q) + \deg'(\bfa)  +1 \right) } \fm^{\cQ}(  \bfa^{3;3}_{\alpha(3)} \otimes q  \otimes  \bfa^{3;1}_{\alpha(3)}) \otimes  \bfa^{3;2}_{\alpha(3)}.
\endaligned
\end{equation}
\begin{defn}
 The \it{Hochschild homology}\index{Hochschild homolog} of $\Cat$ with coefficients in $\cQ$ is the cohomology of the cyclic bar complex with differential $\delta_H$:\index[syindex]{HHlower@$HH_{*}(\cdot,\cdot)$}
 \begin{equation}
   \label{eq:HH_homology}
   HH_{*}(\Cat, \cQ) = H^{*}(  CH_{*}(\Cat, \cQ), \delta_H).
 \end{equation}
\end{defn}

Note that the cyclic bar complex has an increasing filtration by subcomplexes
\begin{equation} \label{eq:length_filtraiton_cyclic}
  CH_{*}^{(N)}(\Cat, \cQ)  = \bigoplus _{k=1}^{N-1} CH_{k}(\Cat[1], \cQ).
\end{equation}
Writing $HH_{*}^{(N)}(\Cat, \cQ)  $ for the homology of this subcomplex, we conclude:
\begin{lem}
The Hochschild homology\index{Hochschild homology of bi-module} of a bi-module is  a direct limit:
  \begin{equation}
     HH_{*}(\Cat, \cQ) = \varinjlim_{N}   HH_{*}^{(N)}(\Cat, \cQ).
\end{equation} \qed
\end{lem}

We write the Hochschild chain complex (resp. Hochschild homology) with coefficents in the diagonal bi-module as 
$CH_*(\Cat, \Cat)$ (resp. $HH_*(\Cat, \Cat)$) instead of $CH_*(\Cat, \Cat_\Delta)$ (resp. $HH_*(\Cat, \Cat_\Delta)$).

Next, we assume that we are given a map $T \co \cP \to \cQ$ of bi-modules as defined in subsection \ref{sec:a_infty-bi-modules} which is closed with respect to the differential of Equation \eqref{eq:differential_map_bi-modules}, and consider the induced map
\begin{align}
  \label{eq:CH_map_bi-modules}
  CH_{*}(T) \co CH_{*}(\Cat, \cP) & \to CH_{*}(\Cat, \cQ) \\ \notag
 p \otimes \bfa & \mapsto \sum (-1)^{\deg'( \bfa^{3;3}_{\alpha(3)}  ) \left( \frak{deg}(p) + \deg'(\bfa)  +1 \right) }  T(  \bfa^{3;3}_{\alpha(3)} \otimes p  \otimes  \bfa^{3;1}_{\alpha(3)}) \otimes  \bfa^{3;2}_{\alpha(3)}.
\end{align}
\begin{lem}
The  map $ CH_{*}(T) $ commutes with the differential, and hence descends to a map
\begin{equation}
  \label{eq:HH_map_bi-modules}
HH_{*}(T) \co  HH_{*}(\Cat, \cP)  \to HH_{*}(\Cat, \cQ) .
\end{equation} \qed
\end{lem}

Moreover, the same argument as in Lemma \ref{lem:naive_equivalence_iso} implies
\begin{lem} \label{lem:equivalence_iso_HH}
  If $T$ is an equivalence of bi-modules then $  HH_{*}(T) $ is an isomorphism.
\end{lem}

An $A_{\infty}$ homomorphism induces a map on Hochschild homology with coefficients in the diagonal bi-module, but it will be convenient for our purpose to decompose this as maps
\begin{equation}
HH_{*}(\Bat,\Bat) \overset{HH_{*}( \Nat_{\Phi} )}\longrightarrow HH_{*}(\Bat,   \Phi^{*}(\Cat)) \overset{HH_{*, \Cat}( \Phi) }\longrightarrow  HH_{*}(\Cat,\Cat).
\end{equation}

The first map is the induced on Hochschild homology by Equation \eqref{eq:natural_map_to_pullback_bi-module}.  The second is the Hochschild homology analogue of Equation (\ref{eq:pulback_bi-module_map_cohomology}).  For any $\Cat$ bi-module $\cQ$, we have a chain level map:
\begin{align} \label{eq:pulback_bi-module_map_homology}
CH_{*,\cQ}( \Phi) \co  CH_{*}(\Bat, \Phi^{*}(\cQ)) & \to CH_{*}(\Cat, \cQ)  \nonumber\\
CH_{*,\cQ}( \Phi)(p \otimes \bfb) & = \sum_{\beta} p \otimes  \Phi( \bfb^{k;1}_{\beta} ) \otimes  \cdots \otimes \Phi(  \bfb^{k;k}_{\beta}).
\end{align}

Again, we have the analogue of Lemma \ref{lem:quasi-iso-categories-cohomology-iso}:
\begin{lem} \label{lem:quasi-iso-categories-cohomology-iso2}
  If $\Phi$ is a quasi-isomorphism, then $HH_{*,\cQ}( \Phi)$ is an isomorphism for every bi-module $\cQ$. \qed
\end{lem}

\begin{comment}

\begin{lem} \label{lem:HH_*-splits}
If $R$ is a field, then the natural map
\begin{equation}
 \bigoplus_{r \in R}  HH_{*}(\Cat_{r}, \Cat_r) \to HH_{*}(\Cat,\Cat)
\end{equation}
is an isomorphism. \qed
\end{lem}
\end{comment}

\subsection{Pairing on Hochschild (co)-homology} 
\label{sec:duality_bi-modules}
Recall that if $\cP$ is a $\Cat$-bi-module, then so is the linear dual $\cP^{\vee}$ defined in Equation \eqref{eq:structure_maps_dual_bi-module}. We obtain a pairing, $\langle *,*\rangle_{\rm HH}$  \index[syindex]{<*>HH@$\langle *,*\rangle_{\rm HH}$},
\begin{align} \label{eq:pairing_cyclic_chains_cochains}
CH^{-*}(\Cat,\cP^{\vee}) \otimes CH_{*}(\Cat,\cP) & \to R \nonumber\\
T \otimes (q \otimes \bfa) & \mapsto \langle T,q \otimes \bfa \rangle_{\rm HH}: = (-1)^{\frak{deg}(q) \deg'(\bfa)} T( \bfa )(q).
\end{align}

Since we constructed the bi-module structure on $\cP^{\vee} $ to be the dual of that on $\cP$, this map induces an isomorphism of chain complexes
\begin{equation}
  CH^{-*}(\Cat,\cP^{\vee}) \cong \Hom(CH_{*}(\Cat,\cP), R),
\end{equation}
so that this pairing is perfect.

If $\Cat$ is a cyclic category, then Lemma \ref{cor:cyclic_structure_dual_bi-modules} provides a degree $n$ equivalence of bi-modules $\Cat_\Delta  \cong \Cat_\Delta^{\vee}$. Combined with the previous result, we conclude:
\begin{lem} \label{lem:pairing_HH} 
 If $\Cat$ is equipped with a pairing of degree $-n$, then the degree $2-n$ map 
 \begin{align} 
   CH^{n-2-*}(\Cat,  \Cat) \otimes CH_{*}(\Cat, \Cat) &\to R  \nonumber\\
\langle T, q  \otimes \bfa \rangle & 
  = (-1)^{\deg'(q) \deg'(\bfa) }  \langle T( \bfa ), q \rangle  \\  \label{eq:pairing_HH_q_first}
& = (-1)^{ 1 + \deg'(q)(1 + n + \deg'(\bfa) )}  \langle q, T( \bfa ) \rangle
 \end{align}
descends to a pairing between Hochschild homology and cohomology.  If $\Cat$ is cyclic, this pairing is non-degenerate. \qed
\end{lem}

Let $\Cat$ and $\Bat$ be two cyclic $A_{\infty}$ categories, and consider an $A_{\infty}$ functor
\begin{equation}
  \Phi \co \Bat \to \Cat
\end{equation}
which is cyclic in the sense of Definition \ref{def:cyclic-morphism}.  
The map $\Nat_{\Phi}$ of $\Bat$ bi-modules  (see  Equation \eqref{eq:natural_map_to_pullback_bi-module}) induces a map on both Hochschild homology and cohomology.  
Together with the maps defined in Equations \eqref{eq:pulback_bi-module_map_cohomology} and \eqref{eq:pulback_bi-module_map_homology}, we obtain a diagram 
  \begin{equation} \label{eq:commuting_pairings}
    \xymatrix@R=2ex{ HH^{n-2-*}(\Bat, \Bat) \otimes   HH_{*}(\Bat, \Bat) \ar[dr] \ar[rr] & &  HH^{n-2-*}(\Bat, \Phi^{*}( \Cat)) \otimes   HH_{*}(\Bat,\Phi^{*}( \Cat) )\ar[ld] \\
& R & \\
HH^{n-2-*}(\Cat, \Cat) \otimes   HH_{*}(\Cat, \Cat) \ar[ur] & &  HH^{n-2-*}(\Cat, \Cat) \otimes   HH_{*}(\Bat, \Phi^{*}(\Cat)) \ar[ll] \ar[uu]. }
  \end{equation}
\par\medskip
\noindent
The following result gives a comparison between the Hochschild homology and cohomology pairings for $\Bat$ and $\Cat$:
\begin{prop}
 Diagram \eqref{eq:commuting_pairings} commutes.
\end{prop}
\begin{proof}
  The key point is to consider the degree $-n$ pairing
  \begin{align}
    CH^{-*}(\Bat, \Phi^{*} \Cat ) \otimes   CH_{*}(\Bat, \Phi^{*}( \Cat)) & \to R \\ \notag
T \otimes q \otimes \bfb & \mapsto(-1)^{\deg'(q) \deg'(\bfb) } \langle  \Phi\left(  T\left( \bfb\right)  \right), q \rangle
  \end{align}
which defines a map from the top right corner of Diagram (\ref{eq:commuting_pairings}) to $R$.  The cyclicity condition on the morphism $\Phi$ implies that the diagram   

\begin{equation}
     \xymatrix@R=2ex{ HH^{n-2-*}(\Bat, \Bat) \otimes   HH_{*}(\Bat, \Bat) \ar[dr] \ar[rr] & &  HH^{n-2-*}(\Bat, \Phi^{*}(\Cat)) \otimes   HH_{*}(\Bat, \Phi^{*}(\Cat)) \ar[dl]   \\
& R &  }
\end{equation}
commutes.  On the other hand, using  Equations \eqref{eq:pulback_bi-module_map_cohomology} and \eqref{eq:pulback_bi-module_map_homology}, one may readily check that we have a commutative square: 
\begin{equation} \label{eq:commuting_pairings2}
    \xymatrix@R=2ex{ R &  HH^{n-2-*}(\Bat, \Phi^{*}( \Cat)) \otimes   HH_{*}(\Bat,\Phi^{*}( \Cat) ) \ar[l]  \\
HH^{n-2-*}(\Cat, \Cat) \otimes   HH_{*}(\Cat, \Cat) \ar[u] &  HH^{n-2-*}(\Cat, \Cat) \otimes   HH_{*}(\Bat, \Phi^{*}(\Cat)) \ar[l] \ar[u]. }
  \end{equation}
\end{proof}

\subsection{Cup and cap products}
\label{sec:cup-cap-products}

We have a well-defined cup product\index{cup product} on the Hochschild cohomology of a category with coefficients in the diagonal bi-module:
\begin{equation}
  HH^{*}(\Cat,\Cat) \otimes   HH^{*}(\Cat,\Cat) \to  HH^{*}(\Cat,\Cat).
\end{equation}
Because of the conventions of using shifted gradings on the diagonal bi-module,
 %\marginpar{Small change on this subsection 
%because of a slight issue of degree and sign change.  KF 2025 Aug.} 
this product has degree $1$. At the chain level, one possible formula for this product is
\begin{align} \label{eq:cup_product}
  (S\cupdot T)^{d}( \bfa) &  = \sum_{ \alpha } (-1)^{\eqref{eq:sign_cup_product_HH}} \fm^{\Cat}(\bfa^{5;1}_{\alpha} \otimes  S( \bfa^{5;2}_{\alpha}) \otimes \bfa^{5;3}_{\alpha} \otimes T( \bfa^{5;4}_{\alpha} )  \otimes \bfa^{5;5}_{\alpha}   )  \\
\label{eq:sign_cup_product_HH} & \hspace{-0.8in} \deg(T) \deg' \left( \bfa^{5;1}_{\alpha} \otimes   \bfa^{5;2}_{\alpha} \otimes \bfa^{5;3}_{\alpha}  \right) +  \deg(S) \deg'( \bfa^{5;1}_{\alpha} ) 
\end{align}
We will explain the reason we use $\cupdot$ soon.\index[syindex]{cupdot@$\cupdot$}

A tedious but straightforward computation  using only the $A_{\infty}$ relations on $ \fm^{\Cat}$ shows that this product, together with the differential from Equation \eqref{eq:diff_CC^*} satisfy the relation
\begin{equation}\label{form457}
  \partial (S \cupdot T) + (\partial S) \cupdot T + (-1)^{\frak{deg} S+1} S \cupdot \partial T = 0,
\end{equation}
which is the first non-trivial $A_\infty$ relation.\footnote{Here $\frak{deg}$ of the element of 
Hochschild cohomology is the degree as an element of 
$\Hom(B_*\Cat[1],\Cat[1])$.} We can extend this product to an $A_\infty$ structure, with the tertiary operation (which provides a homotopy between $R \cupdot (S \cupdot T)$ and $(R \cupdot S) \cupdot T$), obtained as the appropriately signed sum of
\begin{equation}\label{form458}
 \fm^{\Cat}\left( \bfa^{7;1}_{\alpha} \otimes  R( \bfa^{7;2}_{\alpha}) \otimes \bfa^{7;3}_{\alpha} \otimes S( \bfa^{7;4}_{\alpha} )  \otimes \bfa^{7;5}_{\alpha}  \otimes  T( \bfa^{7;6}_{\alpha}) \otimes \bfa^{7;7}_{\alpha} \right).
\end{equation}
More precisely we have
\begin{equation}\label{form459}
(R \cupdot S) \cupdot T + (-1)^{\frak{deg} R+1}R \cupdot (S \cupdot T) = \pm \partial(\ref{form458}). 
\end{equation}
We remark that (\ref{form459}) and associativity of $\cupdot$ up to homotopy \index[syindex]{cup@$\cup$}
is different by sign.  In other words, (\ref{form459}) is the same as a special case of $A_{\infty}$ formula.
\par
An $A_\infty$ structure induces, on the cohomology level, an associative algebra structure with non-shifted gradings
and sign, that is to say, 
if we define\index[syindex]{cvap@$\cap$}
\begin{equation}\label{form460new}
\deg S = \frak{deg} S +1, 
\qquad S \cup T = (-1)^{\deg S} S \cupdot T
\end{equation}
 then $\cup$ induces a structure of associative 
algebra. We call $\cup$  the {\it cup product}.  We thus conclude:\footnote{Then $\deg S$ is the degree of $S$ 
as an element of $\Hom(B_*\Cat[1],\Cat)$.  Note there is no degree shift for the codomain.}
\begin{lem} \label{lem:cup_product_HH}
The cup product descends to an associative product on Hochschild cohomology $  HH^{*}(\Cat,\Cat) $ equipped with non-shifted gradings.   \qed
\end{lem}

The tensor product of bi-modules is functorial in both coefficients, but our choice of Hochschild complex has broken the symmetry between $\Cat$ and $\cQ$. 
 %\marginpar{`we shall only consider the second terms' is commented out. KF 2025 July  }%, and we shall only consider the second terms.  \marginpar{}
First, we define a cap product\index{cap product} \index[syindex]{cvapdot@$\capdot$} (dot version, $\capdot$)
\begin{align}
  \label{eq:cap_product_chains}
  CH^{*}(\Cat,\Cat) \otimes CH_{*}(\Cat, \cQ) & \to CH_{*}(\Cat, \cQ) \\ \notag
T \capdot \left( q \otimes \bfa  \right) & = \sum_{\alpha} (-1)^{ \eqref{eq:sign_cap_product}} \fm^{\cQ}( \bfa^{5;3}_{\alpha} \otimes T(\bfa^{5;4}_{\alpha})  \otimes \bfa^{5;5}_{\alpha} \otimes q \otimes \bfa^{5;1}_{\alpha}) \otimes \bfa^{5;2}_{\alpha}
\end{align}
\begin{equation}
  \label{eq:sign_cap_product}
  \deg'\left(\bfa^{5;3}_{\alpha} \otimes  \bfa^{5;4}_{\alpha} \otimes \bfa^{5;5}_{\alpha}  \right) \cdot \left( \frak{deg}(q) +  \deg' \left( \bfa^{5;1}_{\alpha}  \otimes \bfa^{5;2}_{\alpha} \right) \right) + \deg(T) \deg'\left(\bfa^{5;3}_{\alpha} \right)
\end{equation}
At the chain level, the cap product does not commute with the ring structure on Hochschild cohomology, but there is a homotopy between $(S \cupdot T) \capdot (q \otimes \bfa)  $  and $S \capdot \left(T \capdot (q  \otimes \bfa \right)  $ given by a signed sum of the terms
\begin{equation}
   \fm^{\cQ}\left( \bfa^{7;3}_{\alpha} \otimes S( \bfa^{7;4}_{\alpha} )  \otimes \bfa^{7;5}_{\alpha}  \otimes  T( \bfa^{7;6}_{\alpha}) \otimes \bfa^{7;7}_{\alpha} \otimes q \otimes \bfa^{7;1}_{\alpha} \right) \otimes   \bfa^{7;2}_{\alpha}.
\end{equation}

Then in a similar way as (\ref{form460new}) we put
\begin{equation}\label{form464new}
\deg T = \frak{deg} T +1,
\quad T \cap (q  \otimes \bfa) = (-1)^{\deg T} T \capdot (q  \otimes \bfa).
\end{equation}
Here $\deg$ is the degree as an element of $\Cat(c,c') \otimes B_*\Cat[1](c',c)$ 
and $\frak{deg}$ is the degree as an element of $\Cat[1](c,c') \otimes B_*\Cat[1](c',c)$.
We conclude:

\begin{lem} \label{lem:cap_product_HH}
The chain level cap product \eqref{eq:cap_product_chains},\eqref{form464new} defines an $HH^{*}(\Cat,\Cat)$-module structure on Hochschild homology:
\begin{equation}
  \label{eq:cap_product_homology}
  HH^{*}(\Cat,\Cat) \otimes HH_{*}(\Cat, \cQ)  \to HH_{*}(\Cat, \cQ).
\end{equation} \qed
\end{lem}

In a similar way, when $\cQ$ is a  $\Cat$ bi-module, we define: %\marginpar{Another version is added.  It is used in Section 23.3.  KF 25 June}

\begin{align}
  \label{eq:cap_product_chains2}
CH_{*}(\Cat, \Cat)   \otimes CH^{*}(\Cat,\cQ) & \to CH_{*}(\Cat, \cQ) \\ \notag
\left( a_0 \otimes \bfa  \right) \capdot  T& = \sum_{\alpha} (-1)^{ \eqref{eq:sign_cap_product2}} \fm^{\cQ}( \bfa^{4;5}_{\alpha} \otimes a_0 \otimes \bfa^{4;1}_{\alpha}  \otimes T(\bfa^{4;2}_{\alpha})  
\otimes \bfa^{4;3}_{\alpha}) \otimes \bfa^{4;4}_{\alpha}
\end{align}
where
\begin{equation}
  \label{eq:sign_cap_product2}
  \aligned
  &\deg' (\bfa^{5;1}_{\alpha} ) \cdot \left( {\rm deg}'(a_0) +  \deg' \left( \bfa^{5;2}_{\alpha}  \otimes \bfa^{5;3}_{\alpha} \otimes \bfa^{5;4}_{\alpha}  \otimes \bfa^{5;5}_{\alpha} \right) \right) \\
  &+ \deg(T) \deg'\left(\bfa^{5;2}_{\alpha} \otimes  \bfa^{5;3}_{\alpha} \otimes \bfa^{5;4}_{\alpha}  \right) .
  \endaligned
\end{equation}
and define $\cap$ by (\ref{form464new}).
We also define
\begin{align}
  \label{eq:cap_product_chains3}
CH^{*}(\Cat,\cQ)   \otimes CH_{*}(\Cat, \Cat) & \to CH_{*}(\Cat, \cQ) \\ \notag
T \capdot\left( a_0 \otimes \bfa  \right) =   & \sum_{\alpha} (-1)^{ \eqref{eq:sign_cap_product3}}\fm^{\cQ}\left(  \bfa^{5;3}_{\alpha} \otimes T( \bfa^{5;4}_{\alpha})\otimes  \bfa^{5;5}_{\alpha} \otimes a_0 \otimes \bfa^{5;1}_{\alpha}  \right) \otimes \bfa^{5;2}_{\alpha}\end{align}
where
\begin{equation}
  \label{eq:sign_cap_product3}
    \aligned
  & \left( {\rm deg}'(a_0) + \deg'\left(\bfa^{5;1}_{\alpha} \otimes  \bfa^{5;2}_{\alpha}\otimes  \bfa^{5;3}_{\alpha}   \right) \right) \cdot \deg' \left( \bfa^{5;4}_{\alpha} \otimes \bfa^{5;5}_{\alpha} \right)\\
   &+ \deg(T)\deg'( \bfa^{5;3}_{\alpha})
  \endaligned
\end{equation}  
and define $\cap$ by (\ref{eq:cap_product_chains3}).
We remark that  (\ref{eq:cap_product_chains3}) coincides with (\ref{eq:cap_product_chains}) when $\cQ =\Cat$.
They induce
\begin{equation}
  \label{eq:cap_product_homology2}
  \aligned
 & HH_{*}(\Cat, \Cat) \otimes  HH^{*}(\Cat,\cQ) \to HH_{*}(\Cat, \cQ), \\
& HH^{*}(\Cat, \cQ) \otimes  HH_{*}(\Cat,\Cat) \to HH_{*}(\Cat, \cQ).
 \endaligned
\end{equation}

\section{Triangulated closure and split-generation.}
\label{sec:triang-clos-split}

Given a closed morphism $f \co X \to Y$ in an $A_{\infty}$ category $\Cat$, the notion of a {\it cone} $\Cone(f)$ 
\index{cone} \index[syindex]{Cone@$\Cone(f)$} is well defined up to quasi-isomorphism:  this is an object of $\Cat$, equipped with morphisms $\Cone(f) \to X$ and $Y \to \Cone(f)$, which induce an exact sequence on morphism spaces.

While the notion is well-defined, such an object may not necessarily exist; an $A_{\infty}$ category is \emph{triangulated}\index{triangulated} if all closed morphisms admit cones.  Moreover, every category admits a quasi-embedding into a triangulated category.  In a sense, the most canonical such choice is the Yoneda embedding into the category of left modules.

The main disadvantage of the Yoneda embedding is that it is only a quasi-embedding, making computations sometimes cumbersome; i.e. the morphism spaces between Yoneda modules do not in general have finite rank even if all morphism spaces in $\Cat$ have finite rank.  It is partly for this reason that we prefer to work with \emph{twisted complexes}, following \cite{Sei06}.

\subsection{Twisted complexes} \label{sec:twisted-complexes}

In this subsection, we define an enlargement of $\Cat$ which we call the category of \emph{curved twisted complexes}\index{twisted complex}
\index{curved twisted complex} $\Tw(\Cat)$\index[syindex]{TwC@$\Tw(\Cat)$} over $\Cat$, and which serves as a model for the derived category of $\Cat$.

The category of twisted complexes is enriched over the category of matrix factorisations, and an object $  X_{\bullet}\boxtimes V_{\bullet}$ of such a category consists of the following data:
\begin{enumerate}
\item For each integer $i$ an object $X_{i}$ of $\Cat$, and a matrix factorisation $V_{i}$.  We require that $V_{i} = 0$ for all but finitely many $i$.
\item Whenever $i < j$, a degree $1$ linear map\footnote{The degree is 
defined by $\deg(x \otimes \phi) = \deg x + \deg \phi$.}
  \begin{equation}
    \delta_{i,j} \in  \Cat(X_i,X_j) \otimes \Mat(V_{i}, V_{j})  .
  \end{equation}  
We write $\bfdelta = {\displaystyle{\sum_{\substack{ \blue{0} \leq k \\ i_0 < i_1 < \ldots < i_{k}}} }} \delta_{i_0,i_1} \otimes  \delta_{i_1,i_2} \otimes \cdots  \otimes \delta_{i_{k-1},i_k} $.
\end{enumerate}

%Here we do not insert a term involving $\delta_{i,j}$ in the case that $k=0$.
These data are supposed to satisfy the following conditions:
\begin{enumerate}
  \item There is a scalar $r_{X_{\bullet}\boxtimes V_{\bullet}  } $ such that
    \begin{equation} \label{eq:twisted_complex_potentials_constant}
\Po(V_i) + \Po(X_i) = r_{X_{\bullet}\boxtimes V_{\bullet}  }.
    \end{equation}
\item The self-composition of $\bfdelta$ vanishes
 %\marginpar{It seems that $\fm^{\Sigma \Cat}$ 
%is not defined at this stage.  It is I think defined by (\ref{eq:tensor_composition_A_infty}).  I feel the way this is written is a bit confusing. So I rewrite it a %bit. KF 2025 Aug.}
  \begin{equation} \label{eq:d^2=0}
0 =     \fm^{\Sigma \Cat}( \bfdelta).
  \end{equation}
\end{enumerate}

Here the operation $\fm^{ \Sigma \Cat } $ in Equation \eqref{eq:d^2=0} is a sum of operations $\fm^{ \Sigma \Cat }_{k}$ for each integer $k$ greater than or equal to $1$.\index[syindex]{mxfrakSigmaC@$\fm^{ \Sigma \Cat } $} The operations $\fm^{ \Sigma \Cat }_{k}$  are constructed from the structure maps on $\Mat$ and $\Cat$
as follows.   The case $k=1$ is the natural differential on a tensor product
 %\marginpar{Formula corrected.  (left hand side.)  KF 2025 Feb}
\begin{align} \label{eq:differential_hom_enriched_cat}
  \fm_{1}^{\Sigma \Cat} \co  \Cat(X_i,X_j) \otimes \Mat(V_{i}, V_{j}) & \to \Cat(X_i,X_j)\otimes \Mat(V_{i}, V_{j}) \\
x \otimes \phi & \mapsto   \fm_{1}^{\Cat}(x) \otimes \phi  + (-1)^{\deg'(x)} x \otimes d^{\Mat} \phi .
\end{align}
The cases $k>1$ are obtained by higher composition by the general formula for the tensor product of two $A_{\infty}$ algebras, one of which has vanishing higher order terms.
\begin{equation} \label{eq:tensor_composition_A_infty}
\fm_{k}^{\Sigma \Cat}\left(x_1 \otimes \phi_1,  \ldots  , x_k \otimes \phi_k  \right)  = (-1)^{ \eqref{eq:sign_higher_product_matrix}} \fm^{\Cat}(\bfx) \otimes  \left( \phi_k \circ \cdots \circ \phi_{1}  \right).
\end{equation}
The sign in the above formula is the Koszul sign associated to permuting the morphisms $x_i$ and the maps $\phi_j$
\begin{equation}
  \label{eq:sign_higher_product_matrix}
\sum_{1 \leq i < k} \deg(\phi_i)  \left( \sum_{i < j \leq k} \deg'(x_j) + \deg(\phi_j)  \right).
\end{equation}
Having defined the objects of $\Tw(\Cat) $, we now define the morphisms.  Given a pair $ X_{\bullet}^{0} \boxtimes V_{\bullet}^{0} $  and $ X_{\bullet}^{1} \boxtimes V_{\bullet}^{1}$, we define a graded vector space
 %\marginpar{
%I changed left hand side of (\ref{eq:morphism_twisted_complexes}) from $\Sigma{\Cat}$ 
%to $ \Tw(\Cat)$.  KF 2025 Aug.}
\begin{equation} \label{eq:morphism_twisted_complexes}
  \Tw(\Cat) (X_{\bullet}^{0} \boxtimes V_{\bullet}^{0},  X_{\bullet}^{1}\boxtimes V_{\bullet}^{1} ) = \bigoplus_{i,j} \Cat( X_{i}^{0},  X_{j}^{1}) \otimes \Mat\left(V_{i}^{0},  V_{j}^{1} \right),
\end{equation}
We denote
 %\marginpar{$\Sigma{\Cat}$ appears in the next page.  I think its definition is missing. So I put it here.
%KF 2025 Aug.}
\index[syindex]{SigmaC@$\Sigma{\Cat}$}
$$
\Sigma{\Cat}(X \boxtimes V,X' \boxtimes V') = \Cat( X,  X') \otimes \Mat(V,  V')
$$
and then
$$
\Tw(\Cat) (X_{\bullet}^{0} \boxtimes V_{\bullet}^{0},  X_{\bullet}^{1}\boxtimes V_{\bullet}^{1} ) = \bigoplus_{i,j} \Sigma{\Cat}(X_i\boxtimes 
V_i, X_j\boxtimes V_j).
$$
In order to construct the $A_{\infty}$ structure on $\Tw(\Cat) $, recall from Subsection \ref{sec:matr-fact-dimens} that $\Mat\left(V_{i}^{0},  V_{j}^{1} \right)$ is a matrix factorisation with potential value $\Po(V_{j}^{1}) - \Po(V_{i}^{0}   )  $, and that $ \Cat( X_{i}^{0},  X_{j}^{1})  $ is a matrix factorisation with potential value $\Po(X_{j}^{1}) - \Po(X_{i}^{0}   ) $.  Since the value of the potential function is additive under tensor product, Equation \eqref{eq:twisted_complex_potentials_constant} implies that the morphism space in Equation \eqref{eq:morphism_twisted_complexes} is again a matrix factorisation, with potential function $ r_{X_{\bullet}^{1} \boxtimes V_{\bullet}^{1}  } - r_{X_{\bullet}^{0} \boxtimes V_{\bullet}^{0}  }$.  This justifies defining the curvature of objects in $\Tw(\Cat)   $ to be:
\begin{equation}
  \fm_0^{X_{\bullet}\boxtimes V_{\bullet} }(1) = r_{X_{\bullet}\boxtimes V_{\bullet}  } \sum_{i}  \unit_{X_i}\otimes  \id_{V_{i} }  \in  \Tw(\Cat)(X_{\bullet}\boxtimes V_{\bullet},  X_{\bullet}\boxtimes V_{\bullet}).
\end{equation}

Next, we define the higher operations.  These will be defined using the maps introduced in Equation \eqref{eq:tensor_composition_A_infty} and the differentials $\delta_{i,j}$:\index[syindex]{mxfraktwistC@$\fm_{k}^{\Tw (\Cat)}$}
\begin{equation} \label{eq:higher_product_A_infty}
\aligned
&\fm_{k}^{\Tw (\Cat)}\left(x_1 \otimes \phi_1,  \ldots ,x_k \otimes \phi_k \right) \\
&=  \fm^{\Sigma \Cat}( \bfdelta,  x_1 \otimes \phi_1,  \bfdelta, \cdots, \bfdelta, x_k \otimes \phi_k,   \bfdelta).
\endaligned
\end{equation}

The fact that the operations $\fm_{k}^{\Tw{\Cat}}$ on $\Tw(\Cat) $ satisfy the $A_{\infty}$ relations is a direct consequence of Equation \eqref{eq:d^2=0} and the fact that $\Cat$ is an $A_{\infty}$ category, and our observation that the morphism spaces in Equation \eqref{eq:morphism_twisted_complexes} are themselves matrix factorisations.  We conclude:

\begin{lem}\label{embCattoTwCat}
The operations $ \fm_{k}^{\Tw(\Cat)} $ define a curved $A_{\infty}$ structure on $ \Tw(\Cat)  $ such that the assignment
\begin{equation}
 X \mapsto X  \boxtimes R
\end{equation}
extends to a fully faithful embedding of $\Cat$ as a subcategory of $  \Tw(\Cat) $.  \qed
\end{lem}
Here $X  \boxtimes R$ is the twisted complex such that $X_0 = X$, $V_0  = R$ and all other $X_i = 0$, $V_i =0$.
 %\marginpar{
%A line added. KF 2025 Feb.}

As a consequence of Lemma \ref{lem:factorisations_trivial}, we have 
\begin{lem} \label{lem:maps_twisted_complex_torsion}
 If $Y$ is an object of $\Cat$ with $\Po(Y) =\Po(X_{\bullet}\boxtimes V_{\bullet} )$ and  $\Po(V_i) = r$ for all $i$, then
  %\marginpar{I think  the relation between $\Tw (\Cat)$ and $\Sigma\Cat$ is that, $\Sigma\Cat$ is $X \boxtimes V$ for a single $X,V$ but 
 %$\Tw (\Cat)$ contains a complex of such objects. So I changed the left hand side of (\ref{form512newnew}) to $\Tw (\Cat)$. 
 %(This is also suggested by Yong-Geun. KF 2025 Aug.}
 \begin{equation}\label{form512newnew}
 r \cdot H^{*} \left( {\Tw}(\Cat)(Y, X_{\bullet}\boxtimes V_{\bullet}),   \fm_{1}^{\Tw(\Cat)}  \right) = 0.
 \end{equation}
\end{lem}
\begin{proof}
 By construction, $ X_{\bullet}\boxtimes V_{\bullet} $ is an iterated extension of the objects $  X_{i}  \boxtimes V_{i} $; so we have a spectral sequence
  \begin{equation} \label{eq:spectral_sequence_twisted_complex_endomorphism}
   \bigoplus_{i} H^{*}(\Sigma \Cat (Y, X_{i}  \boxtimes V_{i} ) ;\fm_{1}^{\Sigma \Cat})   \Rightarrow  H^{*} \left( \Tw(\Cat)(Y, X_{\bullet}\boxtimes V_{\bullet}),   \fm_{1}^{\Tw(\Cat)}  \right).
  \end{equation}
It suffices therefore to show that each cohomology group $  H^{*}( \Sigma \Cat(Y, X_{i}  \boxtimes V_{i} ) ;\fm_{1}^{\Sigma \Cat}) $ is annihilated by $r$.  
Note that $\Sigma \Cat (Y, X_i \boxtimes V_i)= \Cat (Y, X_i) \otimes {\rm Mat}(R, V_i)$.  
By the hypothesis of the Lemma, the potential value of $\Cat (Y, X_i)$ (resp.  ${\rm Mat}(R, V_i)$) is $-r$ (resp. $r$).
Taking Remark \ref{difference of Lag PO} and the definition \eqref{eq:differential_hom_enriched_cat} of $\fm_1^{\Sigma \Cat}$ into account, 
Lemma \ref{lem:tensor_product_matrix_vanishes} implies that
$$ r \cdot H^{*}( \Sigma \Cat(Y, X_{i}  \boxtimes V_{i} ) ;\fm_{1}^{\Sigma \Cat}) = 0.$$
\end{proof}

\subsection{A universal complex} \label{sec:universal-complex}
Let $\Bat \subset \Cat$ be a fully faithfully embedded subcategory and $D$ a positive integer. 
We denote by $\chi$ the set of pairs $(X_0,d)$ where $X_0$ is an object of $\Bat$ and $d \in \{0,\dots,D\}$.
We take a linear order $<$ on $\chi$ such that $d_1 > d_2$ implies $(X_0,d_1) < (X'_0,d_2)$.
We fix an order preserving embedding from $(\chi,<)$ to $(\Z,<)$ and regard  $\chi$ as a subset of $\Z$.
 %\marginpar{
%A bit modified so that it exactly fits to the definition of twisted complex given before.
%The mathematical contents do not change.  KF 2025Feb.}
Given an object $Y$ of $\Cat$ and $(X_0,d) \in \chi$, let us introduce the notation
 \begin{equation} \label{eq:definition_bar_complex_Bat-Y}
 \aligned
Y_{X_0,d} :&= 
X_{0}  \boxtimes  \bigoplus_{X_d \in \rm{Ob}(\Bat)}B_{d}(\Bat[1])(X_0,X_{d}) \otimes  \Cat[1](X_{d},Y)\\
&= X_{0}  \boxtimes   B_{d+1}(\Bat [1])(X_0, Y)_{\Cat}.
\endaligned
 \end{equation}
Here $B_{d+1}(\Bat[1])(X_0,Y))_{\Cat}$ (which is defined by Formula (\ref{eq:definition_bar_complex_Bat-Y})) is regarded as a matrix factorization by $d^{\Mat} = \widehat{\frak m}_1$.
 \par
Note that when $Y$ is an object in $\Bat$, 
$B_{d+1}(\Bat[1])(X_0,Y)_{\Cat}$ coincides with the earlier definition.  An element of the vector space 
$B_{d+1}(\Bat[1])(X_0,Y)_{\Cat}$ consists of a collection of $d+1$ composable morphisms in $\Cat$, starting at $X_0$ and ending at $Y$, with the first $d$ being morphisms in the subcategory $\Bat$.

We shall define a twisted complex $\Tot_{\Bat}^{D}(Y) $ whose underlying sequence of  objects
are
\begin{equation} \label{eq:underlying_vector_space_}
\{ Y_{X_0,d} \mid (X_0,d) \in \chi \subset \Z\}.
\end{equation}
\begin{rem}
The main justification for studying this collection of twisted complexes (other than the fact it appears naturally when studying the Floer theory of Lagrangian manifolds), is that $Y$ lies in the category split-generated by $\Bat$ if and only if it is a summand of  $\Tot_{\Bat}^{D}(Y)  $ for some integer $D$.  Since we shall not use such a result, we shall not provide a proof, but  the case of differential graded categories may be extracted from Drinfeld's description of the DG quotient
\cite{Drinfeld} together with Remark A.3.(ii) in \cite{Drinfeld}. This complex can be studied more generally for any subcategory $\Bat \subset \Cat$, but we shall only discuss it in the case all morphism spaces are finite dimensional.
\end{rem}

Before defining the differential on $  \Tot_{\Bat}^{D}(Y)$, it is useful to keep in mind the following representation
\begin{equation}
  \xymatrix{  \displaystyle{\bigoplus_{X_0}}   X_0  \boxtimes \Cat[1](X_0,Y) \\
\displaystyle{\bigoplus_{X_0,X_1}}X_0  \boxtimes \Bat[1](X_0,X_1)   \otimes  \Cat[1](X_{1},Y)   \ar[u] 
 \\  \displaystyle{\bigoplus_{X_0, X_1,X_2} }X_0  \boxtimes \Bat[1](X_0,X_1)   \otimes  \Bat[1](X_1,X_2)   \otimes  \Cat[1](X_{2},Y) \ar[u] \ar@/^10pc/[uu] \\
\cdots \ar[u] \ar@/^11pc/[uu]  \ar@/^15pc/[uuu] }
\end{equation}
Using the coderivation $\hat{\fm}_1$, we consider $B_{d+1}(\Bat [1])(X_0, Y)_{\Cat}$ as a matrix factorisation.   
By the definition of  a twisted complex, the structure maps  $ \delta_{ \Tot } $  are elements of
 \begin{equation} \label{eq:morphism_spaces_total_complex}
\bigoplus_{ X'_0} \Bat(X_0,X'_0)  \boxtimes \Hom\left( B_{d+1}(\Bat[1])(X_0,Y)_{\Cat}  ,  B_{k+1}(\Bat[1])(X'_0,Y)_{\Cat} \right)
\end{equation}
whenever $0 \leq k < d$.

We shall define the differential as a sum of two terms
\begin{equation}\nonumber
  \delta_{ \Tot } = \delta_{ \Tot }^{0} + \delta_{ \Tot }^{\unit}.
\end{equation}
The first term $ \delta_{ \Tot }^{0} $ occurs only when $k=d-1$.  Since $\Bat(X, X')$ is finite dimensional, we may express this term as a sum 
\begin{equation}\nonumber
\delta_{ \Tot }^{0} = \sum_{i} f^{i} \boxtimes \left( f^{i}_{*} \otimes \id \right),
\end{equation}
where  $\{ f^{i} \}$ is a basis of the morphism spaces in $\Bat(X,X')$,
 %\marginpar{I think 
%$f^{i}$ is a basis of $\Bat(X,X')$. So $\Cat$ is changed to $\Bat$. KF 2025 Feb} 
and $\{f^{i}_*\}$ is the dual basis of the dual vector space $\Bat(X, X')^{*}$. 
 %\marginpar{I changed $\vee$ to $*$. In part 4 $\vee$ is used for $\Cat(X,X') = \Cat(X',X)^{\vee}$
%Here we are simply taking dual vector space.  KF 2025 Feb}
More explicitly
$$
\left( f^{i}_{*} \otimes \id \right)(a_1 \otimes \dots \otimes a_d \otimes a_{d+1})
= f^{i}_{*}(a_1)  \,\, a_2 \otimes \dots \otimes a_d \otimes a_{d+1}
$$
 %\marginpar{Formula added. KF. 2025 Feb} 
A straightforward computation shows that the map is independent of the choice of basis.
\par
We observe that $f^{i}_*$ decreases the degree by $\deg'f^{i}$.  Therefore $\delta_{ \Tot }^{0}$ has degree $1$.
\par
%$C(X,X') = C(X',X)^{\vee}$.}
Letting ${\bf e} = {\rm id}$, denote the identity morphism, the second term is defined as follows.
It consists of $\delta_{ \Tot }^{\unit} = \sum \delta_{X_0,d,k}^{\unit}$ with
$$
\delta_{X_0,d,k}^{\unit} : Y_{X_0,d} \to Y_{X_0,k}.
$$
Here $X_0$ are common in the domain and the target  and $0\le k < d$. 
 %\marginpar{Rewritted. 
%according to the slightly modified definition of our 
%twisted complex.  KF  2025 Feb.}
We define:
\begin{equation} \label{eq:differential_complex_tot}
\delta_{X_0,d,k}^{\unit} =  - \unit \boxtimes \widehat{\frak m}_{d-k+1} 
\end{equation}
where 
$$
\aligned
&\widehat{\frak m}_{d-k+1}(a_1 \otimes \dots \otimes a_d \otimes a_{d+1}) \\
&= 
\sum_{i=1}^{k+1} (-1)^* a_1 \otimes \dots \otimes a_{i-1} \otimes \frak m_{d-k+1}(a_i,\dots,a_{d-k+i}) \otimes
\\
& \qquad\qquad\qquad\qquad\qquad\qquad \otimes a_{d-k+i+1} \otimes \dots \otimes a_d \otimes a_{d+1}.
\endaligned
$$
Note that  $\widehat{\frak m}_1 + \sum_{k=0}^{d-1}\widehat{\frak m}_{d-k+1}$ is the differential 
on the bar complex $\partial_{B(\Cat[1])}$.
By this reason we slightly abuse the notation to write 
$$
\delta_{ \Tot }^{\unit} = \sum_{k=0}^{d-1}\delta_{X_0,d,k}^{\unit} = \unit \boxtimes (\widehat{\frak m}_1 - \partial_{B(\Cat[1])}).
$$
%This is a slight abuse of notation since the codomain of $\delta_{X_0,d,k}^{\unit}$ is
%$Y_{X_0,k}$ and is $k$ dependent. So the sum in the left hand side has an issue. However this does not cause a problem 
%during the calculation we perform below.

\begin{lem}\label{lem5454}
The map $  \delta_{ \Tot } $ defines a twisted complex.
\end{lem}
\begin{proof}
The condition that $\Cat$ is strictly unital implies that Equation \eqref{eq:d^2=0} for a twisted complex is given by
\begin{equation} \label{eq:sum_all_terms_twisted_complex}
\aligned
 \fm_{1}^{\Sigma \Cat}(\delta_{ \Tot }^{\unit}) &+ \fm_{2}^{\Sigma \Cat}(\delta_{ \Tot }^{\unit}, \delta_{ \Tot }^{\unit}   ) +  \fm_{2}^{\Sigma \Cat}( \delta_{ \Tot }^{0}, \delta_{ \Tot }^{\unit}   )\\
  &+ \fm_{2}^{\Sigma \Cat}( \delta_{ \Tot }^{\unit}, \delta_{ \Tot }^{0}   ) +  \sum_{1 \leq k} \fm_{k}^{\Sigma \Cat}(\delta_{ \Tot }^{0}, \ldots, \delta_{ \Tot }^{0} )  = 0.
 \endaligned
\end{equation}
An easy computation shows that the first two terms above vanish.  
Accounting for the signs in Equation \eqref{eq:tensor_composition_A_infty} we compute that% \marginpar{YO requests less numbered equations in this page.
%\rcng{I removed some of the number of the formula which is not lavelled. KF}}
\begin{align}
   \fm_{2}^{\Sigma \Cat}( \delta_{ \Tot }^{0}, \delta_{ \Tot }^{\unit} )  & =  \sum \fm_{2}^{\Cat}( f_i, \unit ) \boxtimes  (\hat{\fm}_1 -   \partial_{B(\Cat[1])}) \circ (f^{i}_* \otimes \id )   \nonumber
   \\ \label{eq:first_terms_twisted_complex}
& =  (-1)^{\frak{deg}(f_i)}  f_{i}  \boxtimes  (\hat{\fm}_1 -  \partial_{B(\Cat[1])}) \circ (f^{i}_* \otimes \id )   \\
\fm_{2}^{\Sigma \Cat}(\delta_{ \Tot }^{\unit},  \delta_{ \Tot }^{0})  & =  \sum  \fm_{2}^{\Cat}( \unit ,  f_i) \boxtimes  (f^{i}_* \otimes \id ) \circ (\hat{\fm}_1 - \partial_{B(\Cat[1])}) \nonumber
\\ \label{eq:second_terms_twisted_complex}
& =  \sum  f_i \boxtimes   (f^{i}_* \otimes \id ) \circ (\hat{\fm}_1 - \partial_{B(\Cat[1])})
\end{align}
 %\marginpar{$\deg f_i$ is changed to $\frak{deg}(f_i)$.  KF 2025 Feb.}
To compute further, let us write $ \partial_{B(\Bat[1])}' $  for the map:
\begin{align}\nonumber
  B_{d+1}(\Bat[1])(X_0,Y)_{\Cat} & \to B_{d+1}(\Bat[1])(X_0,Y)_{\Cat} \\
a_{0} \otimes \bfa_{1} & \mapsto (-1)^{\deg'(a_0)} a_{0} \otimes (\hat{\fm}_1 - \partial_{B(\Cat[1])})(\bfa_{1}).
\end{align}
This map essentially consists of those terms of $\partial_{B(\Cat[1])}  $  which act by the identity on the first term.  With this notation, we see that
\begin{equation}\nonumber
(f^{i}_* \otimes \id) \circ \partial'_{B(\Cat[1])} = (-1)^{ \deg'(f_i) } (\hat{\fm}_1 -   \partial_{B(\Cat[1])}) \circ (f^{i}_* \otimes \id ).
\end{equation}
Returning to Equations \eqref{eq:first_terms_twisted_complex} and \eqref{eq:second_terms_twisted_complex} we find that:
\begin{equation}\nonumber
     \fm_{2}^{\Sigma \Cat}( \delta_{ \Tot }^{0}, \delta_{ \Tot }^{\unit} )  +  \fm_{2}^{\Sigma \Cat}(\delta_{ \Tot }^{\unit},  \delta_{ \Tot }^{0})  = 
     -  \sum_i  f_i \boxtimes  ( f^{i}_* \otimes \id ) \circ \left(  \sum_{k \geq 2}  \partial_{B(\Cat[1]); k} \right) ,
\end{equation}
where each map $    \partial_{B(\Cat[1]); k}$ is given by
\begin{equation}
  \partial_{B(\Cat[1]); k}(a_0 \otimes \cdots \otimes a_{d}) = \fm_{k}(a_0 \otimes \cdots \otimes a_{k-1}) \otimes a_{k} \otimes \cdots \otimes a_{d}.
\end{equation}
Returning to Equation \eqref{eq:sum_all_terms_twisted_complex}, note that the terms involving only $ \delta_{ \Tot }^{0}  $ are unaccounted for: using \eqref{Matsign}, we have that $\fm_1^{\Sigma \Cat}(\delta^0_{ \Tot })=0$, and the reader may easily check that 
for $k \geq 2$
\begin{equation}
\sum_i   f_i \boxtimes   (f^{i}_* \otimes \id ) \circ \partial_{B(\Bat[1]); k}= \fm_{k}^{\Sigma \Cat}(\delta_{ \Tot }^{0}, \ldots, \delta_{ \Tot }^{0} ),
\end{equation}
which completes the proof of the Lemma \ref{lem5454}.
\end{proof}

\subsection{Hochschild homology and morphisms of twisted complexes}
Our goal is to understand when $Y$, as an object of $\Tw(\Cat)$, is a summand of $  \Tot_{\Bat}^{D}(Y) $.  
We begin with defining the evaluation map
\begin{equation}\nonumber
  X_0  \boxtimes \Cat[1](X_0,Y)  \to Y
\end{equation}
by
$$
\sum f^i \boxtimes f^i_* \in \Cat(X_0,Y)  \boxtimes \Hom(\Cat(X_0,Y),R),
$$
where $\{f^i\}$ is a basis of $\Cat(X_0,Y)$and 
$\{f^i_*\}$ is its dual basis of the vector space $\Cat(X_0,Y)^* = \Hom(\Cat(X_0,Y),R)$. %\marginpar{Definition corrected.  KF 2025 Feb}
We extend it by $0$ and denote it by\index[syindex]{tau@$\tau$} 
\begin{equation}\nonumber
  \tau \co  \Tot_{\Bat}^{D}(Y)  \to Y.
\end{equation}
It is straightforward to check that $\tau$ is a closed morphism.

We now consider morphisms in the other direction. Returning to the definition of the category of twisted complexes (see Section \ref{sec:twisted-complexes}), and using Equation \eqref{eq:definition_bar_complex_Bat-Y}, we compute that the space of such morphisms is a direct sum
\begin{equation}
\Tw \Cat[1](Y,\Tot_{\Bat}^{D}(Y) ) \cong  \bigoplus_{\stackrel{0 \leq d \leq D}{X_0, X_{d}}} \Cat[1](Y,X_0) \otimes B_{d}(\Bat[1])(X_0,X_{d}) \otimes  \Cat[1](X_{d},Y).
\end{equation}

We shall interpret these morphisms as subcomplexes of a derived tensor product of modules:  Indeed, let $\cY^{r}$ and $ \cY^{l}$ respectively denote the (shifted) Yoneda right and left $\Bat$-modules associated to $Y$ (see \cite{Has03,fukaya:mirII}) and whose underlying vector spaces are respectively  %\marginpar{YO prefers $Y^l$ and $Y^r$}
\begin{align}\nonumber
  \cY^{r}(X) & =  \Cat[1](Y,X)  \\
 \cY^{l}(X) & =  \Cat[1](X,Y).
\end{align}
The derived tensor product\index{derived tensor product} $ \cY^{r} \otimes_{\Bat} \cY^{l}   $  over $\Bat$ is the direct sum
\begin{equation}\nonumber
\bigoplus_{\stackrel{0 \leq d}{X_0, X_{d}}} \Cat[1](Y,X_0) \otimes B_{d}(\Bat[1])(X_0,X_{d}) \otimes  \Cat[1](X_{d},Y).
\end{equation}
By considering only elements which are tensor products of less than $D+2$-factors, we obtain an ascending filtration by subcomplexes
\begin{equation}\label{eq:filtration_tensor_product}
 \left( \cY^{r} \otimes_{\Bat} \cY^{l}  \right)^{(D)} \subset \cY^{r} \otimes_{\Bat} \cY^{l}
\end{equation}
Comparing the differential on the complex defined in Equation \eqref{eq:differential_derived_tensor_product}, we shall conclude:
\begin{lem}\label{lem45}
  There is a natural isomorphism of complexes
  \begin{equation}
    \Tw \Cat[1](Y,\Tot_{\Bat}^{D}(Y) ) \cong \left( \cY^{r} \otimes_{\Bat} \cY^{l}  \right)^{(D)}.
  \end{equation}
\end{lem}
\begin{proof}
  The strict unitality of $\Cat$ implies that the differential on $  \Tw \Cat[1](Y,\Tot_{\Bat}^{D}(Y) ) $ (see Equation \eqref{eq:higher_product_A_infty}), is given by
  \begin{equation} \label{eq:differential-morphisms-to-universal}
    \phi \mapsto \fm^{\Sigma \Cat}_{2}(\phi, \delta_{ \Tot }^{\unit} ) + \sum_{d=1}   \fm^{\Sigma \Cat}_{d}(\phi, \delta_{ \Tot }^{0}, \ldots,  \delta_{ \Tot }^{0} ).
  \end{equation}
The first term is given by
\begin{equation}\nonumber
\aligned
p \otimes \bfb \otimes q  \mapsto
 &   \sum (-1)^{ \deg'(p) + \deg'(\bfb^{3;1}_{\beta(3)}) } p \otimes \bfb^{3;1}_{\beta(3)} \otimes \fm^{ \Bat }\left(  \bfb^{3;2}_{\beta(3)} \right)   \otimes  \bfb^{3;3}_{\beta(3)} \otimes q \\
 &+  \sum (-1)^{ \deg'(p) + \deg'(\bfb^{2;1}_{\beta(2)}) }  p \otimes \bfb^{2;1}_{\beta(2)} \otimes \fm^{ \Cat}\left(  \bfb^{2;2}_{\beta(2)} \otimes q \right)
 \endaligned
\end{equation}
which exactly corresponds to the last two terms in Equation \eqref{eq:differential_derived_tensor_product}, 
%\eqref{eq:higher_product_A_infty}
after accounting for the sign change due to the shift of the Yoneda module.  The remaining terms in Equations \eqref{eq:differential-morphisms-to-universal} and  \eqref{eq:differential_derived_tensor_product} cancel term by term, using the formula
 %\marginpar{A formula is added.  KF 2025 Feb}
$$
\aligned
&\sum_{i_1,\dots,i_k}f^*_{1,i_1}(b_1) \dots f^*_{k,i_k}(b_k)  \frak m_k(p, f_{1,i_1},\dots,f_{k,i_k})
\otimes b_{k+1} \otimes
\dots \otimes b_{d} \\
&=  
\frak m_k(p,b_1,\dots,b_k) \otimes b_{k+1} \otimes \dots \otimes b_{d}.
\endaligned
$$
\end{proof}

Note that the $A_{\infty}$ structure maps on $\Cat$ define a chain map
\begin{align}
  \fm^{\Cat} \co \cY^{r} \otimes_{\Bat} \cY^{l}    & \to \Cat[1](Y,Y) \\ \label{eq:evaluation_map_derived_tensor}
  p \boxtimes \bfa \otimes q & \mapsto  \fm^{ \Cat }(   p \otimes \bfa \otimes q ).
\end{align}
On $ \Tw \Cat[1](Y,\Tot_{\Bat}^{D}(Y) ) $, composition with the canonical map $\tau$ defines a map to the same target:
\begin{lem}
  The following diagram commutes
  \begin{equation}
    \xymatrix{\Tw \Cat[1](Y,\Tot_{\Bat}^{D}(Y) ) \ar[rr] \ar[dr]^{\fm_{2}^{\Tw(\Cat)}(\_, \tau)} & &  \cY^{r} \otimes_{\Bat} \cY^{l}   \ar[dl]^{ \fm^{ \Cat } }  \\
&  \Cat[1](Y,Y). &
}
  \end{equation}
\end{lem}

\begin{proof}
Using again the fact that $ \Cat $ is strictly unital, we find that $ \delta_{ \Tot }^{\unit} $ does not contribute to the composition with $\tau$, which is given by
\begin{equation}
 p \boxtimes \bfa \otimes q   \mapsto  \sum_{k}  \fm^{\Sigma \Cat}_{k}( p \boxtimes \bfa \otimes q, \delta_{ \Tot }^{0}, \cdots, \delta_{ \Tot }^{0}, \tau).
\end{equation}
Since $\tau$ was obtained by extending the evaluation map from $ \Cat[1](X,Y) \otimes X$ to $Y$ by $0$, the summand $d$ in the above expression can only be non-zero if $\bfa = a_0 \otimes \cdots \otimes a_{d}$.  From the definition of the product in the category of twisted complexes, we conclude that Equation \eqref{eq:evaluation_map_derived_tensor} is the expression for the composition with $\tau$.
\end{proof}
Using the fact that $\Tw \Cat[1](Y,\Tot_{\Bat}^{D}(Y) )  $ is an ascending filtration of $\cY^{r} \otimes_{\Bat} \cY^{l} $, we conclude
\begin{cor} \label{cor:criterion_summand_tensor_bi-modules}
$Y$ is a summand of $ \Tw \Cat[1](Y,\Tot_{\Bat}^{D}(Y) ) $ if and only if $[\unit] \in  H^{*}\Cat(Y,Y) $ lies in the image of the composition
\begin{equation}
H^{*}  \left( \cY^{r} \otimes_{\Bat} \cY^{l}  \right)^{(D)} \to H^* \left( \cY^{r} \otimes_{\Bat} \cY^{l}  \right)  \to H^{*}\Cat(Y,Y).
\end{equation} \qed
\end{cor}

\subsection{Idempotent closures}
\label{sec:idempotent-closure}

Let $\Cat$ be an $A_{\infty}$ category, and $X$ a weakly unobstructed object:
\begin{defn}
A \emph{cohomological idempotent}\index{cohomological idempotent} $\iota \in H^{*}(\Cat)(X,X)$ is a morphism satisfying\index[syindex]{iota@$\Iota$}
\begin{equation} \label{eq:idempotent}
  \iota \cdot \iota = \iota.
\end{equation}
\end{defn}

Whenever $Y$ is a summand of $X$, the composition $X \to Y \to X$ is an idempotent.  In an idempotent closed category, every idempotent arises from a summand:
\begin{defn}
  $\Cat$ is \emph{idempotent closed}\index{idempotent closed} if every idempotent is represented by an object, i.e. if for every morphism $\iota$ satisfying Equation \eqref{eq:idempotent}, there exists an object $Y$, and morphisms
  \begin{equation}
    [f] \otimes [g] \in H^{*}(\Cat)(X,Y) \otimes H^{*}(\Cat)(Y,X)
  \end{equation}
such that
\begin{align}
  [f] \cdot [g] & = \iota \\
[g] \cdot [f] & = {\bf e}_{Y}.
\end{align}
\end{defn}

We say that $\Dat$ is an \emph{idempotent closure}\index{idempotent closure} of $\Cat$ if (i) there is a quasi-embedding  $\Tw(\Cat) \to \Dat$, such that every object of $\Dat$ is quasi-isomorphic to a summand of an object in the image of $\Phi$ and (ii) $\Dat$ is idempotent closed.   The existence and uniqueness of idempotent closures up to quasi-isomorphism is a standard fact in the theory of triangulated categories; the $A_{\infty}$ analogue was proved by Seidel in   \cite[Lemma 4.7 and Remark 4.10]{Sei06} which we summarize:
\begin{lem}\label{lem:seidel_existence_idempotent_closure}
  Every $A_{\infty}$ category $\Cat$ has an idempotent closure $\Tw^{\pi}(\Cat)$ which admits an embedding $\Tw(\Cat) \to  \Tw^{\pi}(\Cat)$, and which is quasi-equivalent to any idempotent closure of $\Cat$. \qed
\end{lem}

Observe that, whenever $\Bat \subset \Cat$ is an $A_{\infty}$ embedding, we have a commuting diagram of $A_{\infty}$ embeddings:
\begin{equation}
  \xymatrix{ \Bat \ar[r] \ar[d] & \Cat \ar[d] \\
\Tw(\Bat) \ar[r] & \Tw(\Cat). }
\end{equation}

\begin{defn} \label{def:split-generate}
An object $Y$ of $\Cat$, \emph{lies in the category split-generated by} $\Bat$ if its image in $\Tw(\Cat)$, is a summand of an object in $\Tw(\Bat)$.  If this property is satisfied by every object, we say that $\Bat$ \emph{split-generates}\index{split-generate} $\Cat$.
\end{defn}
\begin{lem}
  If a subcategory $\Bat \subset \Cat$ split-generates, the inclusion $\Tw^{\pi}(\Bat) \subset \Tw^{\pi}(\Cat) $ is an $A_{\infty}$ equivalence.
\end{lem}
\begin{proof}
Consider the smallest subcategory of $ \Tw^{\pi}(\Cat)$ which contains $\Bat$, and is closed under quasi-isomorphisms, cones, and summands.   By assumption, this subcategory includes every object of $\Cat$, hence must agree with $ \Tw^{\pi}(\Cat)$.  The result then follows from   \cite[Corollary 4.9]{Sei06}.
\end{proof}

\section{Gapped filtered $A_{\infty}$ categories.}
\label{sec:a_infty-categ-over}
The purpose of this section is to discuss a special class of $A_{\infty}$ categories over $\Lambda_{0}$ which we call \emph{gapped filtered} following \cite{fooo09}, and which are required for the inductive constructions of the categories in Section \ref{sec:cyclicfil}.

In Section \ref{sec:cyclicfil} there is a situation where we need to study filtered $A_{\infty}$
categories which are {\it not} unital.  So in Section \ref{sec:a_infty-categ-over} and 
Parts \ref{part2} and \ref{part3} we do not assume a filtered $A_{\infty}$
category is unital unless we say it is unital explicitly.
In Sections \ref{sec:curved-cat}-\ref{sec:triang-clos-split} and Part \ref{part4}, 
all $A_{\infty}$ categories are assumed to be (strictly) unital 
unless otherwise mentioned explicitly.

\subsection{Discrete monoids and gapped morphisms}
A discrete submonoid $G$ of $\R_{\ge 0}$ is a subset of  $\R_{\ge 0}$
such that
\begin{enumerate}
\item
$0 \in G$,
\item
if $g_1,g_2 \in G$ then $g_1 + g_2 \in G$,
\item
$G$ is a discrete subset.
\end{enumerate}
Note that the last condition implies that
$
\{g \in G \mid g \le E\}
$
is a finite set for any $E$.
\par
We denote by $\Lambda_G$ the subring of $\Lambda_0$ consisting of the\index[syindex]{LambdaG@$\Lambda_G$} 
formal sums $\sum a_i T^{\lambda_i}$ with $\lambda_i \in G$.
\par
We say {\it discrete monoid}\index{discrete monoid} instead of
discrete submonoid of $\R_{\ge 0}$  for simplicity.
\par

Given an $\F$ vector space $\overline V$, we denote by $V$ the completion of the algebraic tensor product with $\Lambda_0$, i.e.
\begin{equation} \label{eq:complete_tensor_product}
  V  = \overline V \,\widehat{\otimes}_{\F}\, \Lambda_0.
\end{equation}
This is given by the usual tensor product if $\overline V$ is finite dimensional.  For the next definition, we consider an ideal $I$ in $\Lambda_0$, and a pair of $\F$ vector spaces $\overline V_1$ and $\overline V_2$:
\begin{defn}
A map of modules  $ f \co V_1/I \to V_2/I$ is $G$-\emph{gapped}\index{gapped map} if there exists $\F$-linear maps $f_g : \overline V_1 \to \overline V_2$  indexed by elements of $G$  such that
\begin{equation} \label{eq:gapped_sum}
 f = \sum_{g\in G} T^g f_g.
\end{equation}
\end{defn}
\begin{rem}
Assuming that the ideal $I$ is non-trivial, the module $  V_i /I  $ is isomorphic to the tensor product $ \overline V_i \otimes_{\F} \Lambda_0/I$, and the expression in Equation \eqref{eq:gapped_sum} consists of only finitely many non-trivial terms.
\end{rem}

Various other structures can also be said to be $G$-gapped, for example an element $x$ of $V_i /I$ is $G$-gapped if
$$
x = \sum_{g\in G} T^g x_g.
$$
with $x_g \in \overline V_i$.

\subsection{Gapped and filtered categories} \label{sec:gapp-filt-categ}
%In the definition below, and in the rest of this Section, we shall continue to assume that all categories are strictly unital. 

We do not assume any conditions on the curvature, other than those imposed by the gapped property.
%, and all objects are weakly 
%curvature free.

An $A_{\infty}$ category $\Cat$ over $\Lambda_0/I$ is \emph{filtered}\index{filtered $A_{\infty}$ category} if the morphism spaces are isomorphic to the quotient by $I$ of the completion of a free module over $\Lambda_0$ and,
\begin{equation}
  \fm^{\Cat}_{0} \equiv 0 \mod \Lambda_{+},
\end{equation}
for each object of $\Cat$.

% Using the condition that each object be weakly unobstructed, the above condition implies that the curvature of each object is a multiple of the identity by a scalar in $\Lambda_+$.

We define
\begin{equation}
 \Bhat(\Cat[1])(X,Y)
\end{equation}
to be the completion of  $B(\Cat[1])(X,Y)$\index[syindex]{BhatC@$ \Bhat(\Cat[1])$} with respect to the $T$-adic filtration. If $I$ is non-trivial, this completion has no effect.  Our inductive construction of $A_{\infty}$ categories on $E$ 
will require the following more restricted notion:

\begin{defn}
  A \emph{gapped} $A_{\infty}$ category\index{gapped $A_{\infty}$ category} is a filtered $A_{\infty}$ category $\Cat$ over $\Lambda_0 / I$, for which there exists a monoid $G$ such that the structure map
  \begin{equation}
    \fm^{\Cat} \co \Bhat(\Cat)[1](X,Y) \to \Cat(X,Y)[1]
  \end{equation}
is $G$-gapped  for each pair of objects $X$ and $Y$.

A gapped category is \emph{cyclic}\index{cyclic} if it is equipped with a \newred{cyclic} pairing in the sense of Definition  \ref{defn:cyclic_structure}  such that the reduction modulo the maximal ideal
\begin{equation}\nonumber
  \Cat/ \Lambda_{+}
\end{equation}
is a cyclic $A_{\infty}$ category over $\F$.
\newred{In other words we require the pairing to be non-degenerate over $\F$.}
\end{defn}
\begin{rem}
One can state the cyclic condition more explicitly by noting that, if we write
\begin{equation}
  \fm^{\Cat} = \sum_{g \in G} T^{g} \left(\sum_{d \geq 0}  \fm_{g; d}\right),
\end{equation}
then the $A_{\infty}$ equation for $\Cat$ implies that the maps $\fm_{0; d}$ define an $A_{\infty}$ structure over $\F$.
\end{rem}

It is straightforward to extend the notion of a gap to morphisms of categories as well: an $A_{\infty}$ morphism $\Phi \co \Cat \to \Dat $ is gapped if there exists a monoid $G$ such that, for all objects $X$ and $Y$ of $\Cat$, the map
\begin{equation}
  \Phi \co \Bhat(\Cat[1])(X,Y) \to  \Dat[1] (  \Phi(X), \Phi(Y))
\end{equation}
is $G$-gapped.  Such a morphism is cyclic if it satisfies Equation \eqref{eq:cyclic_morphism}.

\begin{rem}
  In \cite{fukaya:cyc}, a gapped morphism $\Phi$ was required, in addition, to satisfy
  \begin{equation} \label{eq:condition_gapped_morphism}
   \sum_{g_1+g_2 =g} \langle \Phi^{1}_{g_1}(x), \Phi^{1}_{g_2}(y) \rangle  = 0
  \end{equation}
 \newred{for $g \ne 0$},
where $\Phi^{1}$ is the linear term of $\Phi$ given by 
\begin{equation}
\Phi^1 = \sum_{g \in G} T^{g}\Phi^{1}_{g}
\end{equation}
is the decomposition of $\Phi^1$ as a gapped morphism.  This condition follows 
from Equation \eqref{eq:cyclic_morphism} which can be seen 
as follows:
since the morphism spaces in $\Cat$ are free $\Lambda_0$ modules, we write $\overline{\Cat}$ for the underlying category over $\F$.  We may then regard
$\Bhat(\Cat[1])(X,Y)$ as the completion of $B(\overline{\Cat}[1])(X,Y)\otimes_{\F} \Lambda_0$.
For pairs of elements $x \in B(\overline{\Cat}[1])(X,Y)$ and $y 
\in B(\overline{\Cat}[1])(Y,X)$, we have
  \begin{equation}
\langle x, y \rangle \in \F \subset \Lambda_0.
  \end{equation}
Equation \eqref{eq:cyclic_morphism} implies that all higher order terms in
  \begin{equation}
\langle \Phi(x), \Phi(y) \rangle
  \end{equation}
vanish, which is precisely the content of Equation \eqref{eq:condition_gapped_morphism}.
\end{rem}

\subsection{Pseudoisotopy of $A_{\infty}$ category}
\label{sec:pseudo-isotopy}
The notion of pseudo-isotopy of cyclic filtered $A_{\infty}$
algebra was introduced in \cite[Definition 8.7]{fukaya:cyc}.
We can generalize the notion to the cyclic filtered $A_{\infty}$
{\it category} as follows.
\par
Let $\Cat$ be a cyclic filtered $A_{\infty}$ category, which
we assume is $G$-gapped.\footnote{We consider both the unital and non-unital versions.}
We consider one of the following situations.
In an algebraic situation we assume that
for each pair of objects $X,Y$
the morphism space $\Cat(X,Y)$ is a
{\it finitely generated} free $\Lambda_0$ module.
In the geometric situation
we includethe case when some morphism spaces are given by the space of differential forms on a Lagrangian. More precisely, in
Section \ref{sec:cyclicfil}, we shall construct a category whose objects $X = (L,\theta_L)$ are Lagrangian submanifolds equipped with extra data. Given a pair $X = (L,\theta_L)$ and $Y = (L,\theta_L)$ whose underlying Lagrangians agree, we shall have
\begin{equation}\label{derhamcaseCat}
\Cat(X,Y) = \Omega(L) \widehat{\otimes} \Lambda_0.
\end{equation}
% Here we are studying the category constructed in
% Section \ref{sec:cyclicfil} and  $X = (L,b)$, $Y = (L,b')$.
% (Namely $X, Y$ both are Lagrangian submanifold $L$
% plus extra data.)
In this setting, the differential we consider agrees 
 with the de Rham differential module the maximal ideal, i.e.
\begin{equation}
\frak m_1(u) \equiv  du \mod \Lambda_+.
\end{equation}

\par
Given a pair $X = (L,\theta_L)$ and $Y = (L',\theta_{L'})$ so that the underlying Lagrangian submanifolds do not agree, 
we shall assume that $L$ is transversal to $L'$.  In this case, $\Cat(X,Y)$
is a finitely generated free $\Lambda_0$ module, as in the algebraic situation.
\par
We also consider the case when
$\Cat$ is a category over $\Lambda_0/T^{E}$
for some $E>0$.
\par
We now discuss the notion of pseudo-isotopy between $\Cat$, and another cyclic unital and filtered $A_{\infty}$ category $\Dat$ (over $\Lambda_0$ or $\Lambda_0/T^{E}$).
We suppose that the set of objects of $\Dat$ is the same as the set of
objects of $\Cat$, and that
%Moreover we assume
\begin{equation}
\Cat(X,Y) =  \Dat(X,Y)
\end{equation}
as $\Lambda_0$ module. We do {\it not} assume the operator $\frak m^{\Cat}_k$ coincides with $\frak m^{\Dat}_k$.
\par
We put
$$
\overline\Cat(X,Y) =  \overline\Dat(X,Y)  = \Cat(X,Y)/\Lambda_+.
$$
It is either a finite dimensional $\Z_2$-graded $\mathbb F$ vector space
or the de Rham complex $\Omega(L) \otimes_{\R} \mathbb F$ of a Lagrangian $L$.
\par
We denote by
$
C^{\infty}([0,1] \times \overline \Cat(X,Y)^k)
$ the  $\F$ vector space of formal sums
\begin{equation}\label{elementdt}
a(t) + dt \wedge b(t)
\end{equation}
where $a(t) \in C^{\infty}([0,1],\overline \Cat(X,Y)^k)$,
$b(t) \in C^{\infty}([0,1],\overline \Cat(X,Y)^{k-1})$. % and denote the set thereof
\par

 In order to construct a filtered $A_{\infty}$ category with morphisms given by the associated completed $\Lambda_0$-modules, we consider the data of operations\index[syindex]{mxCtkbeta@$\frak m^{\frak C, t}_{k,\beta}$}
\begin{equation}\label{param}
\frak m^{\frak C, t}_{k,\beta} : B_k(\overline \Cat(X,Y)[1]) \to \overline \Cat(X,Y)[1]
\end{equation}
for each $t\in [0,1]$ and $\beta\in G$, as well as maps\index[syindex]{c1Ctkbeta@$\frak c^{\frak C, t}_{k,\beta}$}
\begin{equation}\label{parac}
\frak c^{\frak C, t}_{k,\beta} :  B_k(\overline \Cat(X,Y)[1]) \to \overline \Cat(X,Y)[1],
\end{equation}
where
$$
B_k(\overline \Cat(X,Y)[1])
= \bigoplus_{X_0=X,X_1,\dots,X_{k-1},X_k=Y}
\overline \Cat(X_0,X_1)[1] \,\otimes_{\mathbb F}\,\dots
\,\otimes_{\mathbb F}\,\overline \Cat(X_{k-1},X_k)[1].
$$
\par
In the setting where the morphism spaces are finite-dimensional, we say that $\frak m^{\frak C, t}_{k,\beta}$    and $\frak c^{\frak C, t}_{k,\beta}$ are {\it smooth}\index{smooth} if for
each $x_1,\ldots,x_k$ the maps
$$
t\mapsto \frak m^{\frak C, t}_{k,\beta}(x_1,\ldots,x_k) \textrm{, and } t\mapsto \frak c^{\frak C, t}_{k,\beta}(x_1,\ldots,x_k)
$$
lie in $C^{\infty}([0,1],\overline \Cat(X,Y))$.
\par
In the case $\overline \Cat(X,Y)$ is the
de Rham complex $\Omega(L)$, the definition of smoothness is slightly different.
Namely we require the maps 
$$
t\mapsto \frak m^{\frak C, t}_{k,\beta}(x_1,\ldots,x_k) \textrm{, and } t\mapsto \frak c^{\frak C, t}_{k,\beta}(x_1,\ldots,x_k)
$$
to be smooth sections of the pullback of $\Omega L$ to the product $L \times [0,1]$, i.e.\ elements of
$$
\Omega^{\blue{\rm vert}}([0,1] \times L )
$$
without $dt$ component, where $t$ is the parameter of $[0,1]$.
\par

We define
\begin{equation}
 \overline{\frak C}(X,Y)  =
\begin{cases}
\Omega^{\rm vert}([0,1] \times L), & \text{if $\overline \Cat(X,Y) = \Omega(L)$}\\
C^{\infty}([0,1],\overline \Cat(X,Y)),& \text{otherwise}.
\end{cases}
\end{equation}
\begin{rem} \label{rem:finite_dim_case_is_diff_forms}
 The two definitions are compatible in the following way: if $L_{X,Y}$ is a finite set indexing a basis of ${\frak C}(X,Y) $ then we have
  \begin{equation}
    \Omega^{\blue{\rm vert}}([0,1] \times L_{X,Y}) = C^{\infty}([0,1], \overline{\frak C}(X,Y) ).
  \end{equation}
 This implies that we can avoid any case by case analysis.
\end{rem}

\begin{defn}\label{pisotopydef}  %\marginpar{MA190518: I reordered this definition}
We say $(\frak C,\langle\, ,\,\rangle,\{\frak m^{\frak C, t}_{k,\beta}\},\{\frak c^{\frak C, t}_{k,\beta}\})$
is a {\it cyclic pseudo-isotopy}\index{pseudo-isotopy}\index{cyclic pseudo-isotopy} between $G$-gapped cyclic filtered $A_{\infty}$ categories
$\Cat$, $\Dat$ if $\frak m^{\frak C, t}_{k,\beta}$ and $\frak c^{\frak C, t}_{k,\beta}$ are smooth, and the following conditions hold for any value of
  % each (but fixed)
  $t \in [0,1]$:
\smallskip
\begin{enumerate}
%\item  
\item
  % each (but fixed)
  The triple  $(\frak C,\langle \, ,\,\rangle,\{\frak m^{\frak C, t}_{k,\beta}\})$ defines a cyclic filtered $A_{\infty}$
category.  At $t=0$ and $t=1$, this  $A_{\infty}$
category is respectively given by $ \frak m_{k,\beta}^{\frak C, 0} =  \frak m_{k,\beta}^{\Cat}$ and
$ \frak m_{k,\beta}^{\frak C, 1} =  \frak m_{k,\beta}^{\Dat}$.
\item \label{item:t-derivative-of-a_oo-operations}  The $t$-derivative of $ \frak m_{k,\beta}^{\frak C, t}$ is determined by the operations $\frak c^{\frak C, t}_{k,\beta}$ and the relation
\begin{equation}\label{isotopymaineq}
\aligned
&\frac{d}{dt} \frak m_{k,\beta}^{\frak C, t}(x_1,\ldots,x_k) =  \\
& \sum_{k_1+k_2=k}\sum_{\beta_1+\beta_2=\beta}\sum_{i=1}^{k-k_2+1}
\frak m^{\frak C, t}_{k_1,\beta_1}(x_1,\ldots, \frak c_{k_2,\beta_2}^{\frak C, t}(x_i,\ldots),\ldots,x_k)  \\
&- \sum_{k_1+k_2=k}\sum_{\beta_1+\beta_2=\beta}\sum_{i=1}^{k-k_2+1}
(-1)^{\maltese}\frak c^{\frak C, t}_{k_1,\beta_1}(x_1,\ldots, \frak m_{k_2,\beta_2}^{\frak C, t}(x_i,\ldots),\ldots,x_k), %\\
%&=0,
\endaligned
\end{equation}
for each $x_i \in \overline \Cat[1]$,  where the parity of the sign is given by $\maltese = \deg' x_1 + \ldots + \deg'x_{i-1}$.
\item The operation $\frak c_{k,\beta_0}^{\frak C, t}$ corresponding to $\beta_0 = 0\in G$ identically vanishes (hence the operation
$\frak m_{k,\beta_0}^{\frak C, t}$  is independent of $t$ by Equation \eqref{isotopymaineq}).
\item For any $x_i \in \overline \Cat[1]$, we have
\begin{equation}\label{cycsymforc}
\langle
\frak c^{\frak C, t}_{k,\beta}(x_1,\ldots,x_k),x_0
\rangle
=
(-1)^{\maltese}\langle
\frak c^{\frak C, t}_{k,\beta}(x_0,x_1,\ldots,x_{k-1}),x_{k}
\rangle,
\end{equation}
where the sign parity is given by $\maltese = (\deg \blue{x_0}+1)(\deg \blue{x_1}+\ldots + \deg \blue{x_k} + k)$.
\end{enumerate}
We say that $\Cat$ is {\it cyclically pseudo-isotopic}\index{cyclically pseudo-isotopic} to $\Dat$ if there exists a
cyclic pseudo-isotopy between them.
\par
If all the $A_{\infty}$ categories involved are unital and if the 
next (5) holds as well we say that our cyclic pseudo-isotopy is unital:
\begin{enumerate}
\item[(5)]Letting $\text{\bf e}$ denote the unit of an object, we have
\begin{equation}\label{candeequ}
\frak c^{\frak C, t}_{k,\beta}(x_1,\ldots,\text{\bf e},\ldots,x_{k-1})=0.
\end{equation}
\end{enumerate}

\end{defn}
%{\color{blue} Since the notion of cyclic pseudo-isotopy is the only one considered in this paper, we shall omit the qualifier ``cyclic'' unless it is required for clarification.

 We can interpret the above definition in terms of a single $A_\infty$ algebra $\frak C$, over $\Lambda_0$, whose morphism spaces are the complete $\Lambda_0$ modules associated to $ \overline{\frak C}(X,Y) $. The operations on $\frak C$ are given as follows: we define
\begin{align}
\frak m^{\frak C}_{k,\beta} \co
B_k(\overline{\frak C}(X,Y)[1]) & \to \overline{\frak C}(X,Y)[1] \\
{\frak m}^{\frak C}_{k,\beta}(\text{\bf x}_1,\ldots,\text{\bf x}_k)
& = x(t) + dt \wedge y(t) 
\end{align}
where
\begin{subequations}\label{combineainf}
\begin{equation}
x(t) = {\frak m}^{\frak C, t}_{k,\beta}(x_1(t),\ldots,x_k(t))
\end{equation}
and
\begin{equation}\label{combineainfmain}
\aligned
y(t)
=
& \frak c^{\frak C, t}_{k,\beta}
(x_1(t),\ldots,x_k(t)) \\
&-\sum_{i=1}^k (-1)^{\maltese_i} \frak m^{\frak C, t}_{k,\beta}
(x_1(t),\ldots,x_{i-1}(t),y_i(t),x_{i+1}(t),\ldots,x_k(t))
\endaligned\end{equation}
if $(k,\beta) \ne (1,\beta_0)$ and
\begin{equation}
y(t) = \frac{d}{dt} x_1(t) + \frak m_{1,\beta_0}^{\frak C, t} (y_1(t))
\end{equation}
if $(k,\beta) = (1,\beta_0)$. Here $\maltese_i$ in (\ref{combineainfmain}) is
$\maltese_i = \deg' x_1 +\ldots+\deg'x_{i-1}$.
\end{subequations}
The conditions imposed in Definition \ref{pisotopydef} imply:
\begin{lem}
 The category  $\frak C$, with operations ${\frak m}^{\frak C}_{k,\beta}$, is a filtered $A_{\infty}$ category. 
  %\marginpar{ cyclic is removed in two places .  KF 2024 Dec.}
 The evaluation maps
 \begin{equation}
\begin{CD}
 \Dat @<<{\text{\rm Eval}_{t=1}}<   \frak C @>>{\text{\rm Eval}_{t=0}}> \Cat ,
\end{CD}
 \end{equation}
 given by
$$
x = x(t) + dt  \wedge y(t)  \mapsto \text{\rm Eval}_{t=0}(x) = x(0)  \qquad (\text{\rm resp.}      \,\,\,\text{\rm Eval}_{t=1}(x) =x(1)).
$$
are  $A_\infty$ functors.\index[syindex]{Eval@$\text{\rm Eval}$}
The unital version also holds.
\qed
\end{lem}
\par

If $R \to S$ is a ring homomorphism, and $\Cat$ an $A_{\infty}$ category over $R$ in which all morphism spaces are free modules, we obtain an $A_{\infty}$ category $\Cat \tensor_{R} S$ whose objects are the same as those of $\Cat$, and whose morphism spaces are
\begin{equation}\label{basechangetoS}
  \Cat(X,Y)  \tensor_{R} S.
\end{equation}
In case $R=\Lambda_0$ or $\Lambda_0/T^E$ and $S = \Lambda_0/T^{E'}$ (with $E' < E$),
we call (\ref{basechangetoS}) {\it the reduction of $\Cat$ modulo $T^{E'}$}.
We can define the notion of reduction modulo $T^{E'}$ of a pseudo-isotopy in the same way.
\par
We say that two $A_{\infty}$ categories are {\it strongly isomorphic}\index{strongly isomorphic} if there exists a linear $A_{\infty}$ functor between them
that induces a bijection on objects and isomorphisms on
morphism spaces.
We can define the notion of strong isomorphism of pseudo-isotopy in the same way.
\begin{rem}
In \cite{fooo09} the isomorphism between two $A_{\infty}$ algebras is
defined to be an $A_{\infty}$ homomorphism that has an inverse.
The above notion of strong isomorphism requires much stronger conditions than that.
In fact we assume that it is linear. Namely it is $\{\Phi_k\}$ such that $\Phi_k = 0$ for $k\ne 1$.
\end{rem}
We can prove that pseudo-isotopy is an equivalence relation in the same way as
\cite[Lemma 8.2]{fukaya:cyc}.

We shall also prove that every cyclic pseudo-isotopy between $\Cat$ and $\Dat$ induces a cyclic filtered $A_\infty$ functor from $\Cat$ to $\Dat$. Letting $\Cat^t$ denote the cyclic and unital filtered $A_{\infty}$ category whose structure operations are $\frak m^{t,\frak C}_k$, the functor we construct will be obtained by evaluating, at $t=1$, an $A_\infty$ functor $\Phi^t$ from $\Cat$ to $\Cat^t$ which is defined as follows: for each $\text{\bf x} = x_1 \otimes \cdots \otimes x_k \in B_k\overline{\Cat}(X,Y)$, with
$k\ge 1$ and $\beta \in G$, we define $\Phi^t_{k,\beta}(\text{\bf x})$ as the solution to the differential equation
\begin{equation}\label{Psieq}
\aligned
&\frac{d}{dt}\Phi^t_{k,\beta}(\text{\bf x})
\\
&=-
\sum_{m=1}^k\sum_{\beta'+\beta_1+\dots+\beta_m = \beta}\sum_{c}
\frak c^{\frak C,t}_{m,\beta'}
\left(
\Phi_{\beta_1}^t(\text{\bf x}^{(m;1)}_c),\dots,
\Phi_{\beta_m}^t(\text{\bf x}^{(m;m)}_c)
\right)
\endaligned
\end{equation}
with initial condition
\begin{equation}
\Phi^0_{k,\beta}(\text{\bf x})
=
\begin{cases}
x_1 & \text{if $k=1, \beta= \beta_0$}\\
0    & \text{otherwise.}
\end{cases}
\end{equation}
Note that $\beta_i$ in the right hand side of (\ref{Psieq}) satisfies $\beta_i < \beta$, so we can solve (\ref{Psieq}) by induction on $\beta$.
\par
Taking the sum of the above map over all elements $\beta$ of $G$, we define
$$
\Phi^t_{k}(\text{\bf x}) = \sum_{\beta}T^{\omega(\beta)}\Phi^t_{k,\beta}(\text{\bf x}).
$$
\begin{lem}\label{phisisafunctor} \label{cyclicityhoo}
For each $t$, the sequence of maps $\{\Phi^t_{k}\}$ defines a filtered cyclic $A_{\infty}$ functor
$\Phi^t : \Cat \to \Cat^t$. 
\end{lem}
\begin{proof}
We first prove that $\Phi^t$ defines a filtered functor. Since $\Cat = \Cat^0$ and $\{\Phi^0_{k}\}$ is identity, this statement is obvious for $t=0$.
{By the construction, the operations $\Phi^t_k$ are gapped.  Let $G=\{E_i\}$.      
We prove 
\begin{equation}\label{inductiveclaim}
\aligned
&\sum_m\sum_c  \frak m^{t,\frak C}_m
\left(
\Phi^t(\text{\bf x}^{(m;1)}_c), \dots,
\Phi^t(\text{\bf x}^{(m;m)}_c)
\right)\\
& \equiv 
\sum_m\sum_c  (-1)^{\deg' \text{\bf x}^{(3;1)}_c} \Phi^t
\left( \text{\bf x}^{(3;1)}_c,
\frak m^{\blue{t=0},\frak C}(\text{\bf x}^{(3;2)}_c)
\text{\bf x}^{(3;3)}_c
\right)   \ \ \mod T^{E_i},
\endaligned
\end{equation}
by induction on the energy level $E_i$, noting that the base case $E=0$ holds trivially. }

We calculate
\begin{equation}\label{dffphi1}
\aligned
&\frac{d}{dt} \sum_m\sum_c  \frak m^{t,\frak C}_m
\left(
\Phi^t(\text{\bf x}^{(m;1)}_c), \dots,
\Phi^t(\text{\bf x}^{(m;m)}_c)
\right)
\\
=&
\sum_m\sum_c\sum_{1\le i\le j \le m}
\frak m^{t,\frak C}_{m+i-j}
\left(
\Phi^t(\text{\bf x}^{(m;1)}_c), \dots,
\right.
\\
&\qquad\qquad\qquad\qquad\qquad\left.
\frak c^{t,\frak C}_{j-i+1}
\left(
\Phi^t(\text{\bf x}^{(m;i)}_c), \dots,
\Phi^t(\text{\bf x}^{(m;j)}_c)
\right),
\dots,
\Phi^t(\text{\bf x}^{(m;m)}_c)
\right) \\
& +
\sum_m\sum_c\sum_{1\le i\le j \le m}
(-1)^{\maltese}
\frak c^{t,\frak C}_{m+i-j}
\left(
\Phi^t(\text{\bf x}^{(m;1)}_c), \dots,
\right.
\\
&\qquad\qquad\qquad\qquad\qquad\left.
\frak m^{t,\frak C}_{j-i+1}
\left(
\Phi^t(\text{\bf x}^{(m;i)}_c), \dots,
\Phi^t(\text{\bf x}^{(m;j)}_c)
\right),
\dots,
\Phi^t(\text{\bf x}^{(m;m)}_c)
\right)
\\
&-\sum_m\sum_m\sum_{ i=1}^m
\frak m^{t,\frak C}_m
\left(
\Phi^t(\text{\bf x}^{(m;1)}_c), \dots,
\frac{d\Phi^t}{dt}
(\text{\bf x}^{(m;i)}_c),
\dots,
\Phi^t(\text{\bf x}^{(m;m)}_c)
\right),
\endaligned
\end{equation}
where $\maltese =\deg' \text{\bf x}^{(m;1)}_c + \dots + \deg' \text{\bf x}^{(m;i-1)}_c$.
\par
On the other hand, we have
\begin{equation}\label{dffphi2}
\aligned
&\frac{d}{dt} \sum_m\sum_c
(-1)^{\deg' \text{\bf x}^{(3;1)}_c}\Phi^t
\left( \text{\bf x}^{(3;1)}_c,
\frak m^{\blue{t=0},\frak C}(\text{\bf x}^{(3;2)}_c),
\text{\bf x}^{(3;3)}_c
\right) \\
&=
\sum_m\sum_c\sum_{1\le i \le m-2}
(-1)^{\maltese}
\frak c^{t,\frak C}_m
\left(
\Phi^t(\text{\bf x}^{(m;1)}_c), \dots,
\right.
\\
&\
\qquad\qquad\qquad\left.
\blue{  \Phi^t
\left(
\text{\bf x}^{(m;i)}_c, 
\mathfrak{m}^{t=0, \frak C}(\text{\bf x}^{(m;i+1)}_c), 
\text{\bf x}^{(m;i+2)}_c
\right), }
\dots,
\Phi^t(\text{\bf x}^{(m;m)}_c)
\right),   
\endaligned
\end{equation}
where $\maltese =\deg' \text{\bf x}^{(m;1)}_c + \dots + \deg' \text{\bf x}^{(m;i-1)}_c$.
\par
By part \ref{item:t-derivative-of-a_oo-operations}  of Definition \ref{pisotopydef} and \eqref{inductiveclaim} for modulo $T^{E_i}$, 
(\ref{dffphi1}) is equal to (\ref{dffphi2}) modulo $T^{E_{i+1}}$.  
Hence Equation \eqref{inductiveclaim} holds modulo $T^{E_{i+1}}$.  
The proof that $\Phi^t$ is a filtered $A_\infty$ functor is complete.
\par
We now turn to the proof of cyclicity. The case $t=0$ is obvious. We calculate
\begin{equation}\label{diffinner}
\aligned
&-\frac{d}{dt} \sum_c
\left\langle
\Phi^t(\text{\bf x}^{(2;1)}_c),
\Phi^t(\text{\bf x}^{(2;2)}_c)
\right\rangle \\
=&
\sum_c\sum_{m=1}^{\infty}
\left\langle
\frak c^{t,\frak C}\left(\blue{\Phi^t(\text{\bf x}^{(m+1;1)}_c)} ,\dots,
\blue{\Phi^t(\text{\bf x}^{(m+1;m)}_c)} \right),
\Phi^t(\text{\bf x}^{(m+1;m+1)}_c)
\right\rangle \\
&+
\sum_c\sum_{m=1}^{\infty}
\left\langle
\Phi^t(\text{\bf x}^{(m+1;1)}_c),
\frak c^{t,\frak C}\left(\blue{\Phi^t(\text{\bf x}^{(m+1;2)}_c)} ,\dots,
\blue{\Phi^t(\text{\bf x}^{(m+1;m+1)}_c)} \right)
\right\rangle.
\endaligned
\end{equation}
The sum (\ref{diffinner}) is zero because of the cyclic symmetry (that is, \eqref{cycsymforc}) 
and anti-symmetry of the pairing
$\langle \, ,\,\rangle$ (that is, (\ref{eq:switch_factors_pairing})).
The proof %of Lemma \ref{cyclicityhoo} is complete.
of Lemma \ref{phisisafunctor} is complete.
\end{proof}

Specialising the above result to the case $t=1$, we have:
\begin{prop}\label{pisoimplyhomo}
If $\Cat$ is cyclically pseudo-isotopic to $\Dat$ then there exists a
cyclic and filtered $A_{\infty}$ quasi-isomorphism  $\Psi^{\frak C} : \Cat \to \Dat$, which is the identity map on the space of objects. Moreover, the following diagram commutes up to homotopy:
\begin{equation}\label{isotophomodia}
\begin{CD}
\frak C @= \frak C\\
@V{\text{\rm Eval}_{t=0}}VV @VV{\text{\rm Eval}_{t=1}}V \\
\Cat  @>>{\Psi^{\frak C}}>
\Dat.
\end{CD}
\end{equation}
If $\frak C$ and $\frak C'$ are cyclic pseudo-isotopies whose reductions modulo $T^E$
are strongly isomorphic, then the following diagram {\rm strictly} commutes:
\begin{equation}\label{isotophomstrict}
\begin{CD}
\Cat\otimes \Lambda_0/T^E @>{\Psi^{\frak C}}>> \Dat\otimes \Lambda_0/T^E\\
@V{=}VV @VV{=}V \\
\Cat\otimes \Lambda_0/T^E @>{\Psi^{\frak C'}}>> \Dat\otimes \Lambda_0/T^E.
\end{CD}
\end{equation}
\end{prop}
\begin{proof}
We define $\Psi^{\frak C} = \Psi^{1} : \Cat \to \Dat$; this is a cyclic filtered $A_{\infty}$ functor by Lemma \ref{cyclicityhoo}. The strict commutativity of (\ref{isotophomstrict}) is immediate from the construction.
\par
To prove the homotopy commutativity of (\ref{isotophomodia}), we associate to the collection of maps $\Psi^t_k$ the differential form $\tilde{\Psi}_k$ given by
\begin{equation}
\tilde{\Psi}_k(t) =  \Psi^t_k(x_1,\dots,x_k),
\end{equation}
with trivial $dt$ component. This is either an element of $\Omega^{\rm vert}(L)$, or of the space of smooth maps from $[0,1]$ to morphism spaces in $\Cat$. Lemma \ref{phisisafunctor} implies that the collection $\tilde{\Psi} = \{ \tilde{\Psi}_k\}_{k=1}^{\infty} $ defines an $A_{\infty}$ functor
\begin{equation}
  \tilde{\Psi} : \Cat \to \frak C.
\end{equation}
% $\tilde{\Phi} : \Cat \to \frak C$ by
% $$
% \tilde{\Phi}_k(x_1,\dots,x_k)
% = \Phi^t_k(x_1,\dots,x_k).
% $$
% Here the right hand side is regarded as
% a morphism in $\frak C$, that is a differential form on
% $([0,1] \times L_{\kappa}\cap L_{\kappa'}$ without $dt$ component.
\par

By construction we have
$$
\text{\rm Eval}_{t=0} \circ \tilde{\Psi} = \text{identity},
\qquad
\text{\rm Eval}_{t=1} \circ \tilde{\Psi} =\Psi^{\frak C}.
$$
Since $\text{\rm Eval}_{t=0}$ is a homotopy equivalence,
the homotopy commutativity of (\ref{isotophomodia}) follows.
\end{proof}
\begin{rem}
For the cyclic filtered $A_{\infty}$ algebra, a similar result
was proved in \cite[Theorem 8.2]{fukaya:cyc}.
There the morphism $\Psi^{\frak C}$ is constructed explicitly by summation over trees.
The same proof as that of \cite[Theorem 8.2]{fukaya:cyc} can be also applied to our case,
which in fact will produce the same functor as the one we define here.
\end{rem}
\begin{rem}
Existence of the homotopy equivalence $\Cat \to \Dat$ for which Diagram {\rm (\ref{isotophomodia})}
commutes is an immediate consequence of the fact that
$\text{\rm Eval}_{t=0} :\frak C \to \Cat$ has a homotopy inverse.
(See \cite[Proposition 8.9]{fukaya:mirII}.)
However there does not seem to be an appropriate inner product on $\frak C$ to make it a
cyclic $A_{\infty}$ category. So to obtain a {\it cyclic} functor we need to work a bit more.
This is a reason why we introduced the notion of cyclic pseudo-isotopy. % of cyclic $A_{\infty}$ category.
\end{rem}

\begin{prop} \label{prop:pseudo-isotopy-to-quasi-iso}
Let $\Cat_{E'} $ and $\Dat_{E}$ be gapped cyclic  filtered $A_{\infty}$ categories, respectively defined over $\Lambda_0/T^{E'} $ and  $\Lambda_0/T^{E} $, with $E' < E$.
Let $\frak C_{E'}$ be a pseudo-isotopy over $\Lambda_0/T^{E'}$ between
$
\Dat_{E} \otimes_{\Lambda_0/T^{E} } \Lambda_0/T^{E'}
$
and $  \Cat_{E'}$.
\par
Then there exists a gapped cyclic $A_{\infty}$ category   $\tilde{\Cat}_{E} $ with the following properties
\begin{enumerate}
\item $\Cat_{E'}$ is strongly isomorphic to the reduction of $\tilde{\Cat}_{E} $ modulo $T^{E'}$.
\item There exists a cyclic pseudo-isotopy $\tilde{\frak C}_{E}$ between $\tilde{\Cat}_{E} $
and $\Dat_{E}$ such that $\frak C_{E'}$ is strongly isomorphic to the reduction of $\tilde{\frak C}_{E}$.
\end{enumerate}
The unital version also holds.
\end{prop}
\begin{proof}
The proof is the same as that of \cite[Theorem 8.1]{fukaya:cyc}, which we sketch below for completeness.
It suffices to consider the case when $[E',E) \cap G$ consists of a single element, which we denote $\beta$.
\par
We define
$$
\frak c^{\tilde{\frak C}_{E}, t}_{k,\beta} = 0
$$
and
$$
\frak c^{\tilde{\frak C}_{E}, t}_{k,\beta'}
=
\frak c^{{\frak C}_{E'}, t}_{k,\beta'}
\qquad \frak m^{\tilde{\frak C}_{E}, t}_{k,\beta'}
=
\frak m^{{\frak C}_{E'}, t}_{k,\beta'}
$$
for $\beta' \in G \cap [0,E']$.
We then define
$\frak m^{\tilde{\frak C}_{E}, t}_{k,\beta}$ by solving the
ordinary differential equation
(\ref{isotopymaineq}) with initial condition
$
\frak m^{\tilde{\frak C}_{E}, 1}_{k,\beta}
=
\frak m^{{\Dat}_{E}\blue{, 1}}_{k,\beta}.
$
We then define
$
\frak m^{\tilde{\Cat}_{E}}_{k,\beta}
=
\frak m^{\tilde{\frak C}_{E},0}_{k,\beta}.
$
We can check that they have the required properties by a
straightforward calculation.
(See \cite[Section 9]{fukaya:cyc}.)
\end{proof}

\subsection{Homotopy inverse limits}
\label{sec:inverse-limits-categ}

\begin{defn}
A \emph{homotopy inverse system}\index{homotopy inverse system} of $G$-gapped cyclic  filtered $A_{\infty}$ 
categories consists of the following data: 
\begin{enumerate}
\item an  increasing sequence of non-negative real numbers $E_{i}$ for which $\lim_i E_i = +\infty$,
\item a $G$-gapped cyclic  filtered $A_{\infty}$ category $\Cat_{E_i}$ defined over $\Lambda_0 / T^{E_i}$, and 
\item a $G$-gapped  quasi-isomorphism between $\Cat_{E_i} $ and the reduction of $\Cat_{E_{i+1}}$ 
modulo $T^{E_i}$, which is cyclic.
\end{enumerate}
\end{defn}

Note that we do not assume filtered $A_{\infty}$ category to be unital here.
So the notion of quasi {\it equivalence} is not defined.  %\marginpar{A paragraph added.}
In fact we need unit to define the notion that two objects are quasi equivalent.
On the other hand, the notion of quasi {\it isomorphism} is defined without using 
unit, as follows:
\begin{defn}
Let $\Cat, \Dat$ be non-unital curved filtered  $A_{\infty}$ categores. A filtered 
$A_{\infty}$ functor $\Phi : \Cat \to \Dat$ is said to be a {\it quasi-isomorphism}\index{quasi-isomorphism} if 
it is bijection on objects and if  
$\Phi : \Cat(X,Y)/\Lambda_+ \to \Dat(X,Y)/\Lambda_+$ induces a chain homotopy 
equivalence on the chain complex with boundary operator $\overline{\frak m}_1$. 
\end{defn}
In the same way as the case of filtered $A_{\infty}$ algebra (proved in \cite[Theorem 4.2.45]{fooo09}),
we can show that it implies that $\Phi$ has homotopy inverse.

We say that a gapped cyclic $A_{\infty}$ category $ \hCat $ defined over $\Lambda_0$ is a \emph{homotopy inverse limit} of a given  homotopy inverse system, if 
the reduction of $\hCat$ modulo $T^{E_i}$ is gapped quasi-isomorphic 
to $\Cat_{E_i}$, which is cyclic, for each $i$  
such that the following diagram strictly commutes as functors
filtered $A_{\infty}$ categories over $\Lambda_0/T^{E_i}$:
\begin{equation} \label{eq:commuting_diagram_inverse_limit}
  \xymatrix{  \hCat / T^{E_{i+1}} \ar[r] \ar[d]  & \Cat_{E_{i+1}} \ar[d] \\
\hCat / T^{E_i} \ar[r] &  \Cat_{E_i} }
\end{equation}

It is reasonable to expect that every inverse system has a homotopy inverse limit,
but in the context we describe (of gapped cyclic filtered $A_{\infty}$ categories), such a result has  not been proved.
(In case we do not consider cyclic symmetry, the existence of a  homotopy inverse limit is proved in
\cite[Subsection 7.2.8]{fooo09}). % A related problem is also studied in \cite{abouzaidSeidel}).
\par
However if the homotopy equivalences are induced by pseudo-isotopies,
the next result asserts the existence of homotopy inverse limit
(the case of cyclic and unital filtered $A_{\infty}$ {\it algebras}
was proved in \cite[Section 12]{fukaya:cyc}):
\begin{prop}\label{htpyinvlim}
An inverse system  $\Cat_{E_i}$ of
gapped cyclic filtered $A_{\infty}$ 
categories has a homotopy inverse limit, which is cyclic, if the quasi-isomorphism between $\Cat_{E_i}$ 
and the reduction of $\Cat_{E_{i+1}}$ modulo $T^{E_i}$ is induced by a cyclic pseudo-isotopy for each $i$,
 %\marginpar{cyclic is added. KF 2024 Dec.}
in the sense of Proposition {\rm \ref{pisoimplyhomo}}.
\par
The unital version also holds.
\end{prop}

\begin{proof}
This is a corollary of Proposition \ref{prop:pseudo-isotopy-to-quasi-iso}.
Let $\frak C_{E_{i+1}}$ be the pseudo-isotopy between
$\Cat_{E_i}$ and the reduction modulo $T^{E_i}$ of $\Cat_{E_{i+1}}$.
We will prove the following by induction on $j$.
\begin{lem}\label{lemC13}
For any $j' \le j$ and $i\le j'$, there exists
$\Cat_{E_{i},{j'}}$, $\frak C_{E_{i},{j'}}$ with the following properties:
\begin{enumerate}
\item
$\Cat_{E_{i},{j'}}$ is a cyclic and filtered $A_{\infty}$ category over $\Lambda_0/T^{E_{j'}}$.
\item
$\frak C_{E_i,j'}$ is a cyclic pseudo-isotopy modulo $T^{E_{j'}}$ between %\marginpar{cyclic is added. KF 2024 Dec.}
$\Cat_{E_{i-1},{j'}}$ and $\Cat_{E_{i},{j'}}$.
\item
$\Cat_{E_{i},{i}} = \Cat_{E_{i}}$ and
$\frak C_{E_i,i} = \frak C_{E_i}$.
\item
The reduction modulo $T^{E_{j'-1}}$ of $\Cat_{E_{i},{j'}}$
and of $\frak C_{E_i,j'}$  are strongly isomorphic to
$\Cat_{E_i,j'-1}$ and $\frak C_{E_i,j'-1}$  respectively.
\end{enumerate}
\end{lem}
\begin{proof}
Suppose the lemma is proved for $j-1$.
We then define $\Cat_{E_{i},{j}}$ and $\frak C_{E_i,j}$
by downward induction on $i$ starting with $\Cat_{E_{i},{i}} = \Cat_{E_{i}}$,
$\frak C_{E_i,i} = \frak C_{E_i}$ and using Proposition \ref{prop:pseudo-isotopy-to-quasi-iso}.
\end{proof}
Now we consider $\Cat_{E_{1},{j}}$.
For objects $X,Y$ we can identify
$\overline{\Cat}_{E_{1},{j}}(X,Y)$ with $\overline{\Cat}_{E_{1},{1}}(X,Y)$
by the strong isomorphism given by item (4) above.
We put
$$
\hCat(X,Y) =  \overline{\Cat}_{E_{1},{1}}(X,Y) \otimes_{\mathbb F} \Lambda_0.
$$
We may regard the operators $\frak m^{\overline{\Cat}_{E_{1},{j}}}_k$
as a homomorphism from its bar complex. Then
we have
$$
\frak m^{\overline{\Cat}_{E_{1},{j}}}_k \equiv
\frak m^{\overline{\Cat}_{E_{1},{j'}}}_k
\mod T^{E_{j'}}
$$
if $j' \le j$. Thus we can define
$$
\frak m^{\hCat}_k
=  \liminv_{j\to \infty}  \frak m^{\overline{\Cat}_{E_{1},{j}}}_k.
$$
Let us check that $\hCat$ has the required property by constructing the diagram (\ref{eq:commuting_diagram_inverse_limit}):
Let $j \in \Z_{>0}$.
By Lemma \ref{lemC13} (\blue{4}), $\hCat \otimes {\Lambda_0}/T^{E_j}$
is strongly isomorphic to $\Cat_{E_1,j}$.
We have a cyclic pseudo-isotopy  $\frak C_{E_i,j'}$ modulo $T^{E_{j'}}$ between
$\Cat_{E_{i-1},j'}$ and $\Cat_{E_i,j'}$.
By Proposition \ref{pisoimplyhomo} it induces a cyclic  filtered $A_{\infty}$
functor modulo $T^{E_{j'}}$ that is a quasi-isomorphism.
We denote it by $\Psi_{i,j'} : \Cat_{E_{i-1},j'} \to \Cat_{E_i,j'}$.
We now define $\Psi_i : \hCat(X,Y)\otimes \Lambda_0/{T^{E_i}} \to \Cat_{E_{i}}$ as the composition
$$
\hCat(X,Y)\otimes \Lambda_0/{T^{E_i}}
\cong
\Cat_{E_1,i} \overset{\Psi_{2,i}}\longrightarrow \Cat_{E_2,i}
\overset{\Psi_{3,i}}\longrightarrow
\cdots
\overset{\Psi_{i-1,i}}\longrightarrow \Cat_{E_i,i}  = \Cat_{E_{i}}.
$$
This is a quasi-isomorpism.
\par
We remark that Lemma \ref{lemC13} (2) implies that
$\frak C_{E_i,j'} \otimes \Lambda_0/T^{E_{j'-1}}$
is strongly isomorphic to $\frak C_{E_i,j'-1}$.
It implies that
$\Psi_{i,j'} = \Psi_{i,j'-1} \mod T^{E_{j'-1}}$.
The commutativity of Diagram  (\ref{eq:commuting_diagram_inverse_limit}) follows easily.
\end{proof}

Note we did not use the commutativity of (\ref{isotophomodia}) in the last step of the above proof, because we need to show that (\ref{eq:commuting_diagram_inverse_limit})
commutes strictly, while (\ref{isotophomodia}) commutes only up to homotopy.

It is possible to compose pseudo-isotopies after deforming them to ensure that they are locally constant near the endpoints (see \cite[Lemma 8.2]{fukaya:cyc}). However, it is sufficient for our purpose to consider the notion of \emph{formal composition}  of pseudo-isotopies from $\Cat$ to $\Dat$, i.e. a collection  $\{ \Cat_0\}_{j=0}^{n}$ of $A_\infty$ categories with $\Cat_0 = \Cat$ and $\Cat_n = \Dat$, as a collection  $\frak C = \{ \frak C_j \}_{j=1}^{n}$ with each $\frak C_j$ a pseudo-isotopy from a $\Cat_{j-1}$ to $\Cat_{j}$. The argument used to prove Proposition \ref{htpyinvlim} establishes the following:
  \begin{prop} \label{prop:inverse_system_pseudo-isotopy}
    Given a collection $\{\Cat_{E_i}\}$ of cyclic curved $A_\infty$ categories defined over  $\Lambda_0/T^{E_i}$, and cyclic pseudo-isotopies $\{ \frak C_{E_i,E_{i+1}} \} $  between $ \Cat_{E_i} $ and the reduction of $ \Cat_{E_{i+1}} $ modulo  $T^{E_i}$, there exists a cyclic curved $A_\infty$ category $ \hCat $ defined over $\Lambda_0$, 
    equipped with a formal composition $\frak C_{\infty,E_{i}} $ of cyclic pseudo-isotopies between the reduction of $\hCat$ modulo $T^{E_i}$ and $\Cat_{E_i}  $  such that the reduction of $ \frak C_{\infty,E_{i+1}} $ modulo $T^{E_i}$ is the formal composition of $ \frak C_{\infty,E_{i}}$ with $\frak C_{E_i,E_{i+1}} $.
\par
The unital version also holds.
  \end{prop}

\subsection{Transfer of bounding cochains}
\label{sec:transfer_bounding_cochains}

In the context of the previous section, we would like to describe a choice of bounding cochains on an object of the homotopy inverse limit of an inverse system in terms of bounding cochains on the underlying categories of the inverse system. The starting point is to be able to relate bounding cochains on pseudo-isotopic $A_\infty$ algebras.

Let $E$ and $E'$ be positive numbers with $E < E'$, and let $\Cat$ and $\Dat$ denote gapped filtered $A_\infty$ categories, respectively defined modulo $T^E$ and $T^{E'}$. Let $\frak C$ be a pseudo-isotopy modulo $T^E$ between $\Cat$ and the reduction modulo $T^E$ of $\Dat$.

Recall that the objects of $\Dat$ and $\Cat$ are identified by the pseudo-isotopy. Let $X$ be such an object, and let $b_{+}^{\Dat} \in \Dat(X,X)$ and $b_{+}^{\Cat} \in \Cat(X,X)$ be  bounding cochains modulo $T^{E'}$ and $T^E$. 
\begin{defn}
  A \emph{pseudo-isotopy between bounding cochains}\index{pseudo-isotopy between bounding cochains} $b_+^{\Dat}$ and $b_+^{\Cat}$ is a  bounding cochain $b_+^{\frak C} $ whose image under the evaluation maps at $t = \{0,1\}$ are  $b_+^{\Dat}$ and $b_+^{\Cat}$. Such a pseudo-isotopy has \emph{temporal gauge} if the $dt$ component of $b_+^{\frak C} $ is trivial.
\end{defn}

% In the above construction 
% we require that $b^{\frak C}_{\kappa, +} $ has no $dt$ component.
% (In other words we take temporal gauge.)
% Under this additional condition $b^{\Cat}_{\kappa} $ is uniquely determined from $b^{\Dat}_{\kappa} $.

\begin{lem}\label{lem:solve_ODE_to_get_bounding_cochain}
  Given a bounding cochain $b_+^{\Dat}$ on $\Dat$, there exists a uniquely determined  bounding cochain $b_+^{\Cat}$ on $\Cat$ so that the pair ($b_+^{\Dat}, b_+^{\Cat})$ are pseudo-isotopic in temporal gauge (the corresponding pseudo-isotopy $b_+^{\frak C}$ is itself uniquely determined).
\end{lem}
\begin{proof}
  The strategy is to first construct the pseudo-isotopy  $b_+^{\frak C}$ as a solution to an ODE, with boundary condition $b_+^{\Dat} $ at $t=1$, then evaluate at $t=0$ to obtain $b_+^{\Cat} $. To this end, we define  
$$
\frak m^{\frak C,t}_k = \sum_{g \in G}  \frak m^{\frak C,t}_{k,\beta}, 
\qquad
\frak c^{\frak C,t}_k = \sum_{g \in G}  \frak c^{\frak C,t}_{k,\beta},
$$
and recall that  $\frak m^{\frak C}_k$ is induced by  $\frak m^{\frak C,t}_k$ and $\frak c^{\frak C,t}_k$ by Formula (\ref{combineainf}).

We define the desired  bounding cochains $b^{\frak C}_{+} \in \frak C(X, X)$ on $\frak C$ by
solving the ordinary differential equation
\begin{equation}\label{ODEbounding}
\frac{d b^{t,\frak C}_{+} }{dt}
+ 
\sum_{k} \frak c^{\frak C,t}_{k}((b^{t,\frak C}_{+} )^{k}) = 0,
\end{equation}
with initial condition 
$$
b^{1,\frak C}_{+}  = b^{\Dat}_{+}.
$$

To see that $ b^{\frak C}_{+}$ is a  bounding cochain, observe that Equation (\ref{ODEbounding}) is a part of the equation 
$\sum \frak m_k^{\frak C}((b^{\frak C}_{+})^{\otimes k}) \equiv \frak{PO}_{\frak b}(b) \cdot \text{\bf e}_L \mod T^{E}$.
The other part is 
\begin{equation}\label{formnew640}
\sum_{k} \frak m^{\frak C,t}_{k}((b^{t,\frak C}_{+} )^{\otimes k}) 
\equiv \frak{PO}_{\frak b}(b_{\kappa} ) \cdot \text{\bf e}_{L_{\kappa}} \mod T^{E}
\end{equation}
that can be proved by differentiating the left hand side of (\ref{formnew640}).

The bounding cochain $b_+^{\Cat}$ is then obtained by evaluation. 
\end{proof}

\begin{rem}
In case we allow $b^{\frak C}_{+} $ to have a nontrivial $dt$ component 
the chain $b^{\Cat}_{+} $ is not uniquely determined. But its gauge equivalence 
class is determined uniquely.
\end{rem}

Note that, in the above construction, it is straightforward to replace the pseudo-isotopy $\frak C$ by a formal composition. We may thus return to the context of Proposition \ref{prop:inverse_system_pseudo-isotopy}: given a system  $\{\Cat_{E_i}\}$ of curved $A_\infty$ categories defined over  $\Lambda_0/T^{E_i}$, and pseudo-isotopies $\{ \frak C_{E_i,E_{i+1}} \} $  between $ \Cat_{E_i} $ and the reduction of $ \Cat_{E_{i+1}} $ modulo  $T^{E_i}$, we denote by $\Cat$ the inverse limit, which we equip with the formal composition $\frak C_{\infty,E_{i}} $ pseudo-isotopies modulo $T^{E_i}$ to  $ \Cat_{E_i} $.
\begin{lem} \label{lem:bounding_cochains_inverse_limit}
Given an object $X$ of $\Cat$, there is a bijective correspondence between the following data: {\rm (i)} a bounding cochain $b$ on $X$, and 
{\rm (ii)} a collection $\{b(E_i)\}_{i=0}^{\infty}$ of bounding cochains on the object corresponding to $X$ in $\Cat_{E_i}$, together with a pseudo-isotopy $b_{E_i, E_{i+1}}$ modulo $T^{E_i}$ in temporal gauge between $b(E_i)$ and $b(E_{i+1})$ (i.e. a bounding cochain of $X$ as an object of $ \frak C_{E_i,E_{i+1}}$, with trivial $dt$ component, and which restricts to $b_{E_i}$ at the ends.)
\end{lem}
\begin{proof}
  In one direction, the data of $b$, determines, via applying Lemma \ref{lem:solve_ODE_to_get_bounding_cochain} to the formal composition $ {\frak C}_{\infty,E_{i}} $ of pseudo-isotopies modulo $T^{E_i}$ between $\Cat$ and $\Cat_{E_i}$, a bounding cochain $b(E_i)$ on $\Cat_{E_i}$. The fact that $b(E_i)$ and $b(E_{i+1})$ are themselves pseudo-isotopic then follows from the fact that the reduction modulo $T^{E_i}$ of $ {\frak C}_{\infty,E_{i+1}} $ is the formal composition of $ {\frak C}_{\infty,E_{i}} $ and ${\frak C}_{E_{i}, E_{i+1}}$. The argument in the reverse direction is similar.
\end{proof}

\begin{rem}
In Part \ref{part2} we first   construct a (curved) filtered $A_{\infty}$ category inductively and then use bounding cochains to eliminate the %\marginpar{Remark added. KF 2024 Dec}
curvature.  So we do not use  Lemma \ref{lem:bounding_cochains_inverse_limit} there.
Lemma \ref{lem:bounding_cochains_inverse_limit} is used for example in the situation when we use 
the fact that the obstructions for the existence of bounding cochains (See \cite[Theorem C]{fooo09}) vanish to inductively construct 
a bounding cochain, such as the case 
when $H^{\rm even}(X;\Q) \to H^{\rm even}(L;\Q)$ is surjective. (See  \cite[Theorem H]{fooo09}.)
\end{rem}

\subsection{Reduced Hochschild homology and cohomology}

Let $\Cat$ be a (strictly) unital and $G$-gapped filtered $A_{\infty}$ category
or an $A_{\infty}$ category over a filed.
 %\marginpar{This subsection is new.  KF 2025 March}
In later sections we need to study reduced version of Hochschild homology and cohomology.
We will define it in this subsection and show some of its properties.

\begin{defn}\label{defn525}
The {\it reduced Hochschild cochain complex}\index{reduced Hochschild cochain complex} $CH^*_{\rm red}(\Cat)$ is a subcomplex of Hochschild cochain complex\index[syindex]{CHredC@$CH^*_{\rm red}(\Cat)$}
consisting of elements 
$\varphi = (\varphi_0,\varphi_1,\dots,\varphi_k,\dots)$, 
$\varphi_k : B_k(\Cat[1])(c,c') \to \Cat[1](c,c')$ such that
$$
\varphi_k(\dots, {\bf e}_c,\dots) = 0.
$$
\end{defn}
It is easy to see from the definition of the unit that 
this subset is preserved by the Hochschild differential.
We denote its cohomology by $HH^*_{\rm red}(\Cat)$.

The inclusion defines a homomorphism
$HH^*_{\rm red}(\Cat) \to HH^*(\Cat)$ between 
Hochschild cohomology and its reduced version.

\begin{prop}\label{prop526}
If $\Cat$ is a cyclic and unital $A_{\infty}$ category over a field 
and is weakly curvature free, then the homomorphism
$HH^*_{\rm red}(\Cat) \to HH^*(\Cat)$ is an isomorphism.
 %\marginpar{Maybe cyclicity is not necessary to prove this proposition.
%Since we use only that case and since I still have an issue to prove it without cyclicity I leave it as a part of assumptions.
%KF 2025 Aug.}
\par
The same is true for a cyclic, unital and weakly curvature free  filtered $A_{\infty}$ category over a Novikov ring 
with ground ring being a field.\end{prop}

In the case of differential graded algebra over the ground ring this is 
well-known (See \cite[Proposition 1.1.15]{Loday}.)  
Its $A_{\infty}$ category version (over the field) seems to be well known to 
expert. Since it seems difficult to find a reference which discusses the curved (but weakly curvature free) case, we will prove 
Proposition \ref{prop526} later in this subsection.
See \cite{Cho12} for 
some discussion on the Hochschild (co)homology for the obstructed case.

We next discuss Hochschild homology.
For $c \in {\rm{Ob}}(\Cat)$, we define\index[syindex]{Cxredoverline@$\Cat^{\rm red}$}
$$
\Cat^{\rm red}(c,c) = \Cat(c,c)/\Lambda_0 {\bf e}_c.
$$
We also define $\Cat^{\rm red}(c,c') = \Cat(c,c')$ if $c\ne c'$.
For $c, c' \in {\rm{Ob}}(\Cat)$ we define reduced 
Bar complex  by 
$$
B_k^{\rm red}\Cat[1](c,c') = \bigoplus \Cat^{\rm red}[1](c_0,c_1) \otimes \cdots \otimes \Cat^{\rm red}[1](c_{k-1},c_k)
$$

Here the direct sum is taken over all 
$c_0,\dots,c_k \in  {\rm{Ob}}(\Cat)$ with $c_0 = c$, $c_k = c'$.
In the case of filtered $A_{\infty}$ category we replace $\otimes$ by $\widehat{\otimes}_{\Lambda_0}$.

\begin{defn}
We define {\it reduced Hochschild chain complex}\index{reduced Hochschild chain complex} $CH^{\rm red}_*(\Cat,\Cat)$
by\index[syindex]{CH^redC@$CH^{\rm red}_*(\Cat,\Cat)$}
$$
CH^{\rm red}_{*}(\Cat,\Cat) = \bigoplus_{c,c' \in \rm{Ob}(\Cat)}\bigoplus_{k=0}^{\infty}
\Cat[1](c',c) \otimes B_k\Cat^{\rm red}[1](c,c').
$$
In the case of filtered $A_{\infty}$ category we replace $\otimes$ by $\widehat{\otimes}_{\Lambda_0}$.
\end{defn}
Note that there exists a surjective homomorphism 
$$
\pi: CH_*(\Cat,\Cat) \to CH^{\rm red}_*(\Cat,\Cat)
$$
By the definition of strict unit, it is easy to see that
there exists a boundary operator $\delta_H$ on $CH^{\rm red}_*(\Cat,\Cat)$ 
such that $\pi$ is a chain map.

\begin{prop}\label{prop528}
If $\Cat$ is weakly curvature free then $\pi$ induces an isomorphism 
$HH_*(\Cat,\Cat) \to HH^{\rm red}_*(\Cat,\Cat)$.\footnote{We do not need cyclicity for 
this proposition.}
\footnote{The proposition holds both for filtered category over Novikov ring or category 
over a field.}
\end{prop}

\begin{proof}
Let ${\bf e}_c$ be the unit of $\Cat$ and $k \in \Z_{\ge 0}$.
We consider 
$$
{\bf x} = x_0 \otimes (x_1 \otimes \dots \otimes x_k) \in CH_{k+1}(\Cat,\Cat)
$$
Let $\ell \in \Z_{\ge 0}$.
We define
$$
s_{\ell}({\bf x}) = (-1)^{\deg'x_0 + \dots + \deg' x_{\ell}}x_0 \otimes  \dots \otimes x_{\ell} \otimes {\bf e}_c \otimes x_{\ell+1} \otimes \dots \otimes x_k,
$$
if $\ell \le k$, $x_{\ell} \in \Cat(c',c)$ and 
$$
s_{\ell}({\bf x}) = 0
$$
if $k <\ell$.
Let $F_{m} CH_*(\Cat,\Cat)$ be the submodule of $CH_*(\Cat,\Cat)$ generated 
by the images of $s_{\ell}$ with $\ell \le m$.
It is easy to see that $\delta_H$ preserves $F_{m} CH_*(\Cat,\Cat)$.
\begin{lem}\label{lem529}
$$
(\delta_H \circ s_{\ell} + s_{\ell} \circ \delta_H) \circ s_{\ell} -  s_{\ell }  \equiv 0  \mod  F_{\ell-1}  CH_*(\Cat,\Cat).
$$
\end{lem}
The proof of the lemma is a straight forward calculation and is omitted. %\marginpar{Check the sign with Loday.  KF 2025 Aug 29}
(Compare \cite[Lemma 1.6.6]{Loday} which is the case of associative ring.)
\par
Lemma \ref{lem529} implies that the identity map 
$F_{\ell} CH_*(\Cat,\Cat) \to F_{\ell}  CH_*(\Cat,\Cat)$
is chain homotopic to a map to $F_{\ell-1} CH_*(\Cat,\Cat)$.
We next apply the lemma to $\ell-1$ and repeat.
Then the identity map $F_{\ell} CH_*(\Cat,\Cat) \to  F_{\ell} CH_*(\Cat,\Cat)$
is chain homotopic to $0$.
Therefore $F_{\ell} CH_*(\Cat,\Cat)$ is acyclic.
It follows that 
$
\bigcup_{\ell}F_{\ell} CH_*(\Cat,\Cat)
$ is acyclic.
Since the kernel of $
\pi: CH_*(\Cat,\Cat) \to CH^{\rm red}_*(\Cat,\Cat)
$
is $
\bigcup_{\ell}F_{\ell} CH_*(\Cat,\Cat)
$
the proposition follows.
\end{proof}

\begin{proof}[Proof of Proposition \ref{prop526}]
In the case of cyclic $A_{\infty}$ category over a field, 
Hochschild cohomology is dual to Hochschild homology.
The same holds for the reduced versions. 
Therefore the proposition follows from Proposition \ref{prop528}.
\par
We consider the case of a filtered $A_{\infty}$ category over the Novikov ring.
\par
We first prove the surjectivity.
Let $\varphi = (\varphi_0,\varphi_1,\dots,\varphi_k,\dots) \in CH^*(\Cat)$
be a Hochschild cocycle.
We put $G = \{\lambda_j \in \R_{\ge 0} \mid j \in \Z_{\ge 0} \,\text{with $\lambda_j < \lambda_{j+1}$, $\lambda_0=0$}\}$.
We prove the next statement by induction on $j$.
We put
$$
\varphi = \sum T^{\lambda_i} \varphi^i
$$
with $\varphi^i \in CH^*(\overline{\Cat})$.
We will define 
$$
\varphi' = \sum T^{\lambda_i} \varphi^{\prime i}, \quad \psi = \sum T^{\lambda_i} \psi^i
$$
with $\varphi^{\prime i} \in CH^*_{\rm red}(\overline{\Cat})$, 
$\psi^i \in CH^*(\overline{\Cat})$
such that the next lemma holds for
$$
\varphi^{(j)} = \sum_{i=0}^{j} T^{\lambda_i} \varphi^i, 
\quad 
\varphi^{\prime (j)} = \sum_{i=0}^{j} T^{\lambda_i} \varphi^{\prime i}, 
\quad 
\psi^{(j)} = \sum_{i=0}^{j} T^{\lambda_i} \psi^i.
$$
\begin{lem}
\begin{equation}\label{indhohred}
\varphi^{(j)} - \delta^H(\psi^{(j)}) - \varphi^{\prime (j)} \equiv 0 \mod  T^{\lambda_{j+1}}.
\end{equation}
\end{lem}
\begin{proof}
The proof is by induction on $j$.
We assume (\ref{indhohred}) for $j$ and will prove it in the case of $j+1$.
Let $\phi^j \in CH(\overline{\Cat})$ with
$$
\varphi^{(j+1)} - \delta^H(\psi^{(j)}) - \varphi^{\prime (j)} \equiv T^{\lambda_{j+1}}\phi^j \mod  T^{\lambda_{j+2}}
$$
Since $\delta^H(\varphi) = 0$ 
we have
$$
\overline{\delta}^H(\phi^j) \in CH^*_{\rm red}(\overline{\Cat}).
$$
Therefore there exists $\varphi^{\prime, j+1,1} \in CH^*_{\rm red}(\overline{\Cat})$ such that
$$
\overline{\delta}^H(\phi^j) = \overline{\delta}^H(\varphi^{\prime, j+1,1}).
$$
Since 
$\overline{\delta}^H(\phi^j - \varphi^{\prime, j+1,1}) = 0$
we can
use $HH^*_{\rm red}(\overline{\Cat}) \cong HH^*(\overline{\Cat})$
to modify $\varphi^{\prime, j,1}$ to 
$\varphi'^{j+1}\in CH^*_{\rm red}(\overline{\Cat})$ such that $[\phi^j - \varphi'^{j+1}]$ is $0$ in $HH^*(\overline{\Cat})$.
Namely there exists $\psi^{j+1} \in  CH^*(\overline{\Cat})$  such that
$$
\phi^j - \varphi'^{j+1} = \overline\delta^H(\psi^{j+1}).
$$
Thus
$$
\varphi^{(j+1)} - \delta^H(\psi^{(j)}+ T^{\lambda_{j+1}}\psi^{j+1}) - \varphi^{\prime (j)} - T^{\lambda_{j+1}} \varphi^{\prime j+1} \equiv 0 \mod  T^{\lambda_{j+2}}
$$
as required.
The proof is complete by induction.
\end{proof}
We put 
$$ 
\varphi^{\prime} = \sum_{i=0}^{\infty} T^{\lambda_i} \varphi^{\prime i}, 
\quad 
\psi = \sum_{i=0}^{\infty} T^{\lambda_i} \psi^i.
$$
Then 
$$
\varphi = \delta^H(\psi) + \varphi^{\prime}. 
$$
We have proved the surjectivity.
\par
We next prove the injectivity.
 %\marginpar{Proof of the next lemma is given so that 
%we have isomorphism in the filtered case also.  KF 2025 Aug.}
Let 
$\varphi \in CH^*_{\rm red}({\Cat})$, $\psi \in CH({\Cat})$ such that 
$
\varphi = \delta^H( \psi)
$.
\begin{lem}
There exists $\phi^{(k)} \in CH({\Cat})$ such that
\begin{eqnarray}
&\psi - \delta^H(\phi^{(k)}) \in CH^*_{\rm red}({\Cat}) + T^{\lambda_k}CH({\Cat})  \label{form641}\\
& \phi^{(k+1)} -  \phi^{(k)} \in T^{\lambda_k}CH({\Cat}). \label{form642}
\end{eqnarray}
\end{lem}
\begin{proof}
The proof is by induction. The case $k=1$ follows from 
$HH^*_{\rm red}(\overline {\Cat}) \cong HH^*(\overline {\Cat})$.
Suppose we have $\phi^{(k)}$. 
We put $\psi - \delta^H(\phi^{(k)}) = a_k + T^{\lambda_k}b_k$ with $a_k \in  CH^*_{\rm red}({\Cat})$.
Then
$$
\varphi - \delta^H(a_k) = T^{\lambda_k} \delta^H b_k  \in CH^*_{\rm red}({\Cat})
$$
Using  $HH^*_{\rm red}(\overline {\Cat}) \cong HH^*(\overline {\Cat})$  we find 
$c_k \in T^{\lambda_k}CH^*_{\rm red}({\Cat})$
such that $T^{\lambda_k} \delta^H b_k - \delta^H c_k \in T^{\lambda_{k+1}}CH^*({\Cat})$.
Then 
$$
\delta^H(T^{\lambda_k} b_k - c_k) \in T^{\lambda_{k+1}}CH_{\rm red}^*({\Cat}).
$$
Using $HH^*_{\rm red}(\overline {\Cat}) \cong HH^*(\overline {\Cat})$ again we find $d_k
\in  T^{\lambda_k}CH^*_{\rm red}({\Cat})$, $e_k \in T^{\lambda_k}CH({\Cat})$
such that
$$
(T^{\lambda_k} b_k - c_k) -(d_k + \delta^H e_k) \in T^{\lambda_{k+1}}CH^*({\Cat}).
$$
Now
$$
\psi - \delta^H(\phi^{(k)} + e_k)- (a_k + c_k + d_k) 
 \in T^{\lambda_{k+1}}CH^*({\Cat})
$$
and $a_k + c_k + d_k \in CH^*_{\rm red}({\Cat})$.
Therefore $\phi^{(k+1)}: = \phi^{(k)} + e_k$ has the required property.
\par
The proof is complete by induction.
\end{proof}
We take $\phi = \lim_{k\to \infty} \phi^{(k)}$. Here the right hand side converges by 
(\ref{form642}). 
Then
(\ref{form641}) implies $\psi - \delta^H(\phi) \in CH^*_{\rm red}({\Cat})$.
Therefore $\varphi = \delta^H(\psi - \delta^H(\phi))$ is 
zero in $HH^*_{\rm red}({\Cat})$ as required.
\end{proof}

\newpage
\part{$A_{\infty}$ category in Lagrangian Floer theory.}
\label{part2}

\section{Filtered $A_{\infty}$ algebra on de Rham complex associated to a 
single Lagrangian submanifold.}
\label{ainfalgasssingle}

This section is a review of the 
construction of \cite{fooo09}. We rewrite it using 
de Rham cohomology and also include some enhancements.

We first introduce moduli spaces of pseudo-holomorphic disks.

Let $L$ be a relatively spin {\it embedded} Lagrangian submanifold of $X$ 
and $\beta \in H_2(X,L;\Z)$,   and recall that we have fixed an $\omega$-compatible almost complex structure $J$:
\begin{defn}\label{diskmoduli1}
We denote by ${\mathcal M}_{\ell;k+1}(L;\beta)$\index[syindex]{Mellk+1Lbeta@${\mathcal M}_{\ell;k+1}(L;\beta)$}
the set of all $\sim$ equivalence classes of quartets 
$(\Sigma;u;z_1^{+},\dots,z_{\ell}^+;w_0,\dots,w_k)$  consisting of the following data:
\begin{enumerate}
\item
$\Sigma$ is a connected bordered curve of genus zero with one boundary component, which has only 
double points as singularity.
\item
$u : \Sigma \to X$ is a pseudo-holomorphic map 
such that $u(\partial D^2) \subset L$.
\item $z_1^{+},\dots,z_{\ell}^+$ are mutually distinct points in the interior of $\Sigma$, 
which are not nodal.
\item
$w_0,\dots,w_k$ are mutually distinct points on the boundary  
$\partial \Sigma$ of $\Sigma$, whose ordering respects the counterclockwise cyclic order on $\partial \Sigma$.\footnote{We use the symbol $w$ in place of $z$ which was used in the previous literatures, 
for the sake of consistency with the later sections.}
They are not nodal.
\item The relative homology class $u_*([\Sigma,\partial \Sigma])$ 
is $\beta \in H_2(X,L;\Z)$.
\item The automorphism group of $(\Sigma;u;z_1^{+},\dots,z_{\ell}^+;w_0,\dots,w_k)$ 
is a finite group.  Here an automorphism of $(\Sigma;u;z_1^{+},\dots,z_{\ell}^+;w_0,\dots,w_k)$ 
is a biholomorphic map $v : \Sigma \to \Sigma$ which preserves $z^+_i$, $w_j$ and $u \circ v = u$.
\end{enumerate}
We call such $(\Sigma;u;z_1^{+},\dots,z_{\ell}^+;w_0,\dots,w_k)$ a \emph{stable holomorphic disk}.\index{stable holomorphic disk}
\par\medskip
The equivalence relation\index{equivalent} is given by 
$$(\Sigma;u;z_1^{+},\dots,z_{\ell}^+;w_0,\dots,w_k)
\sim 
(\Sigma';u';z_1^{\prime +},\dots,z_{\ell}^{\prime +};w'_0,\dots,w'_k)$$
if there exists a biholomorphic map $v : \Sigma \to \Sigma'$
such that
$$
u' \circ v = u,
\quad
v(z_i^+) = z_i^{\prime +},
\quad
v(w_i) = w_i^{\prime}.
$$
\par
In case $(\Sigma';u';z_1^{\prime +},\dots,z_{\ell}^{\prime +};w'_0,\dots,w'_k) =
(\Sigma;u;z_1^{+},\dots,z_{\ell}^+;w_0,\dots,w_k)$
such $v$ is called an automorphism.\index{automorphism}
An  extended automorphism\index{extended automorphism} is defined in a similar way 
but relax the condition 
$v(z_i^+) = z_i^{\prime +}$ to 
$v(z_i^+) = z_{\sigma(i)}^{\prime +}$
where $\sigma$ is a permutation of $\{1,\dots,\ell\}$.
\par
We denote by $\overset{\circ}{\mathcal M}_{\ell;k+1}(L;\beta)$ the subset of  ${\mathcal M}_{\ell;k+1}(L;\beta)$
consisting of elements such that $\Sigma$ is a disk with or without sphere bubbles (but not disk bubbles).
\end{defn}
\begin{rem}
{In other accounts of Lagrangian Floer theory (e.g. \cite{Sei06}), one often specifies that the boundary marked point $w_0$ corresponds to the output of Floer theoretic operations, thus treating it differently from the other boundary marked points. However, as we shall later require that all our constructions be cyclically invariant, the reader will find that, in our account, all boundary marked points are ultimately treated in the same way (i.e. we may break the cyclic symmetry to write a formula which seems to distinguish $w_0$ from the other marked points, but we could alternatively use any other marked point, and obtain the same operations in the end).}
\end{rem}

The moduli space ${\mathcal M}_{\ell;k+1}(L;\beta) $ is equipped with natural evaluation maps\index{evaluation map}
\index[syindex]{evrmi@$\text{\rm ev}_i$}\index[syindex]{evrmiplus@$\text{\rm ev}^+_i$}
$$
(\text{\rm ev}^+,\text{\rm ev}) 
= (\text{\rm ev}^+_1,\dots,\text{\rm ev}^+_{\ell};
\text{\rm ev}_0,\dots,\text{\rm ev}_k) 
: {\mathcal M}_{\ell;k+1}(L;\beta)
\to X^{\ell} \times L^{k+1}
$$
{whose factors with value in $X$ and $L$ are respectively given by}
$$
\text{\rm ev}^+_i([u;z_1^{+},\dots,z_{\ell}^+;w_0,\dots,w_k])
= u(z_i^+), 
\quad
\text{\rm ev}_{j}([u;z_1^{+},\dots,z_{\ell}^+;w_0,\dots,w_k])
= u(w_{{}j}).
$$
\begin{prop}\label{diskkura}
$  $ \par
\begin{enumerate}
\item  We can topologize the moduli space of stable holomorphic disks ${\mathcal M}_{\ell;k+1}(L;\beta)$ 
so that it is compact, Hausdorff  and metrizable.
\item
The space ${\mathcal M}_{\ell;k+1}(L;\beta)$ has an orientable Kuranishi structure with corners. 
\item The evaluation maps at boundary marked points 
 %\marginpar{Compatibility with forgetful map with 
%interior marked points are removed. (In other places also.) See right above Proposition \ref{multiproductdisk}
%KF2024}
are underlying continuous maps of strongly smooth maps %\marginpar{Check name and reference.} 
with respect to this 
Kuranishi structure.
\item \label{diskkura:weak_submersive} ${\rm ev}_0$ is weakly submersive.\footnote{In fact, Item (10) of Proposition \ref{diskkura} implies that ${\rm ev}_i$ is weakly submersive for every boundary marked point $i$.}
\item \label{diskkura:boundary_strata}
The normalized\footnote{See \cite[Definition 8.4]{springer}.} boundary of ${\mathcal M}_{\ell;k+1}(L;\beta)$ in the 
sense of Kuranishi structures is described by the following disjoint union of fiber products 
over $L$:\footnote{The fiber product is well defined by Item (4).}
\begin{equation}\label{eq177}
\partial{\mathcal M}_{\ell;k+1}(L;\beta)
= \bigcup {\mathcal M}_{\# {\mathbb L}_1;k_1+1}(L;\beta_1)
{}_{\text{\rm ev}_{0}}\times_{\text{\rm ev}_i} {\mathcal M}_{\# {\mathbb L}_2;k_2+1}(L;\beta_2), 
\end{equation}
where the union is taken over all $({\mathbb L}_1,{\mathbb L}_2) \in \text{\rm Shuff}(\ell)$, $i = 1, \dots, k_2$,
$k_1,k_2\in \Z_{\ge 0}$ with $k_1 + k_2 = k$ and $\beta_1,\beta_2
\in H_2(X,L;\Z)$ with $\beta_1 +\beta_2 = \beta$.
\item
{The (virtual) dimension of the moduli space of stable holomorphic disks is}
\begin{equation}\label{dimension}
\dim {\mathcal M}_{\ell;k+1}(L;\beta) = n+\mu_L(\beta) 
+k -2 + 2\ell,
\end{equation}
where $\mu_L : H_2(X,L;\Z) \to  2\Z$ is the Maslov index.
\item
{A choice of relative $\mathrm spin$ structure on $L$ determines} orientations of ${\mathcal M}_{\ell;k+1}(L;\beta)$ {so that the induced orientation of the boundary of the moduli space in Equation \eqref{eq177} is
compatible in the sense of \cite[Definition 21.3 (IX)]{springer}, when restricted to each stratum of the right hand side, with the fibre product orientation.}
\item \label{diskkura:forget}
  The Kuranishi structure is compatible with forgetful maps\index{forgetful map}\index[syindex]{forgetfrak@$\frak{forget}$}
  (See the beginning of Subsection {\rm \ref{diskforget}} for the definition)
  \begin{equation*}
 \frak{forget}:   {\mathcal M}_{\ell;k+1}(L;\beta) \to {\mathcal M}_{\ell;k}(L;\beta)
  \end{equation*}
  of the 
boundary marked points. (See \cite[Definition 3.1]{fukaya:cyc}.)
\item
The Kuranishi structure is invariant under the permutation of interior marked 
points.
\item \label{diskkura:cyclic}
The Kuranishi structure is invariant under the cyclic permutation of the 
boundary  marked 
points.
\end{enumerate} 
\end{prop}
Parts (1)--(7),(9),(10) of Proposition \ref{diskkura} are proved
in \cite[Propositions 7.1.1,7.1.2,8.5.1]{fooo09}.
Some more detail is written  \cite{const1, const2}, including the case $\ell \ne 0$.
The Kuranishi structure satisfying the additional properties (8)
is constructed in  \cite[Corollary 3.1]{fukaya:cyc}.
We will give the construction following \cite{const1, const2,fukaya:cyc} 
in Subsections \ref{diskreview} and \ref{diskforget}.
(Note that in Subsections \ref{diskreview} and \ref{diskforget} we construct 
a system of Kuranishi structures on the outer collaring of the moduli spaces, 
by a certain technical reason.)

\begin{rem}\label{runningoutrem}
There is one technical point to take care of in the statements, 
such as Proposition \ref{diskkura}, appearing many times in this paper. 
Namely to prove Proposition \ref{diskkura} as stated above, we need to 
study infinitely many moduli spaces and construct 
Kuranishi strutures or CF-perturbations etc. on them which are all compatible.
The construction of such objects is done by induction on the energy 
 $\omega(\beta)$ of homology classes $\beta$. We then meet a trouble 
which is described in detail in \cite[Subsection 7.2.3]{fooo09}.
In \cite[Subsection 7.2]{fooo09} a method to go around this trouble is 
given. Namely we prove such a statement only up to an arbitrary but fixed energy level $E$
for each $E$ and then we use the technique of homological algebra to promote the 
structure to the higher energy level.
This method is applicable to the situation of this paper. 
\par
We actually do it in the following way:
\begin{enumerate}
\item
We construct Kuranishi structures and \blue{CF-perturbations}  for the 
moduli spaces of holomorphic maps with energy $\le E$ for 
each but fixed $E \ge 0$.
\item
We use these data to define operations or maps up to energy level $E$.
\item
Then the equalities hold only modulo $T^E$.
\item 
We take a sequence $E_i$ of positive numbers that diverges to infinity
and construct a filtered $A_{\infty}$ category $\Cat_{E_i}$ 
over $\Lambda_{0}/T^{E_i}$ for each $i$.
\item
For $E_i < E_j$ we prove that $\Cat_{E_i}$ is homotopy equivalent to 
$\Cat_{E_j}$ over $\Lambda_{0}/T^{E_i}$.
(We moreover prove that they are pseudo-isotopic in the sense of 
Definition \ref{pisotopydef}.)
\item
We use homological algebra to take the inverse limit of $\Cat_{E_i}$ 
and obtain a filtered $A_{\infty}$ category over $\Lambda_0$.
\end{enumerate}
%{We give detail of this argument in Appendices \ref{sec:a_infty-categ-over} and \ref{sec:homotopyequiv}.}
The mod $T^E$ version has some additional applications 
which we discuss in Section \ref{sec: quatitative}.
\par
We implement this strategy in detail in Subsections \ref{sec:constr-cycl-a_infty4} and 
\ref{sec:homotopyequiv} for the construction of the category 
in Section \ref{sec:cyclicfil}. Meanwhile for simplicity we state 
Proposition \ref{diskkura} as above, leaving implicit the fact that the following line should be added before the statement:
\par
``For each $E \ge 0$ the following holds for all $\beta$ with $\omega(\beta) \le E$.''
\par
Many of the equalities appearing in this paper are actually hold mod $T^E$ only 
during Steps (1)-(3) above.
\end{rem}
\begin{rem}\label{Rmk:7.5}
Actually the Kuranishi structures and CF-perturbations we use are {\it not}
compatible with the forgetful map of the interior marked points.
See Remark \ref{rem819} for its reason.  So we actually need 
a double induction on the energy $E$ and  the number $\ell$ of interior marked points. %\marginpar{Remark added. KF 2024Oct}
\end{rem}

The notion of CF-perturbations\footnote{In \cite{fooo:bulk}, a slightly different version, ``continuous families of perturbations'', is used to integrate 
 along the fibers.}  is introduced to make sense of the integration along fibers in the realm of 
Kuranishi structures (see \cite[Chapter 7]{springer}). 
We follow the way of \cite{springer}
to construct a system of CF-perturbations on the   outer collaring
${\mathcal M}_{\ell;k+1}(L;\beta)^{\boxplus 1}$ of  ${\mathcal M}_{\ell;k+1}(L;\beta)$.
(See Definition \ref{outcolar} for the definition of the outer collaring.)

\begin{prop}\label{existmkulti1}
There exists
a system of CF-perturbations
on the moduli spaces ${\mathcal M}_{\ell;k+1}(L;\beta)^{\boxplus 1}$
such that the following holds:
\begin{enumerate}
\item
It is transversal to zero.
\item
It is compatible with the
description of the boundary in Proposition $\ref{diskkura}$ $(\ref{diskkura:boundary_strata})$ above.
\item
It is compatible
with the forgetful map of the boundary marked points.\footnote{See
Subsection \ref{CFforget} a line below Proposition \ref{prop1315}.
See also  \cite[Definition 5.1]{fukaya:cyc}.}
\item
It is compatible
with the permutation of interior marked points.
\item
It is compatible with cyclic permutation of the boundary marked points.
\item
The evaluation map
$\text{\rm ev}_0$ is strongly submersive with respect to this CF-perturbation.\footnote{See \cite[Definition-Lemma 7.50]{springer}.  It means that the restriction of $\text{\rm ev}_0$ to the 
zero set of the CF-perturbation is a submersion.}
\end{enumerate}
\end{prop}
The \newred{CF-perturbations} satisfying (1), (2), (4), (6) above 
is constructed in \cite[Chapter 22]{springer} 
(and Subsection \ref{CFconst}), while the system of \newred{CF-perturbations} satisfying (3)(5) in addition 
is constructed in \cite[Corollary 5.1]{fukaya:cyc}.  
See Subsection \ref{CFforget}.

\begin{rem}
$  $ \par
\begin{enumerate}
\item
In \cite{springer},
to carry out the inductive construction of a system of CF-pertrubations 
in  \cite[Chapter 22]{springer} we used  the consistency of 
Kuranishi structures and CF-perturbations not only along the boundary 
but also along the corners. In other words we used the corner compatibility 
conditions. 
However we do not need the corner compatibility 
conditions in this paper by the reason explained in Remark \ref{remark1313}.
\item
There is an issue of smoothness of the forgetful map on the outer collar.
We will discuss it in Section \ref{sec:CRperturb}.
\end{enumerate} %\marginpar{\rcng{Remark added. KF}}
\end{rem}
\par
{For the remainder of this section, we consider a fixed Lagrangian $L$ and recall the construction of the curved 
$A_\infty$ structure on $\Omega(L)$. This will be the first step in the inductive construction of the category with multiple Lagrangians, which we return to in Subsection \ref{cycliccategorymoduli} below.}

Let $g_1,\dots,g_{\ell} \in \Omega(X)$ and
$h_1,\dots,h_k \in \Omega(L)$ and $\beta$ with
$(\beta,k,\ell) \ne (0,0,0)$. We define\index[syindex]{qbetakell@$\frak q^{\rm form}_{\beta;\ell,k}$}
\begin{equation}\label{defqformula}
\aligned
&\frak q^{\rm form}_{\beta;\ell,k}(g_1,\dots,g_{\ell}; h_1,\dots,h_k)\\
&= (-1)^{\epsilon(h_1,\dots,h_k) + \delta(g_1,\dots,g_{\ell})} 
\frac{1}{\ell!}
\text{\rm Corr}({\mathcal M}_{\ell;k+1}(L;\beta);(\text{\rm ev}^+_1
\times\dots \times\text{\rm ev}^+_{\ell}\times
\\
&\qquad\qquad\qquad\qquad \times \text{\rm ev}_1
\times\dots \times\text{\rm ev}_{k}),
\text{\rm ev}_{0})
(g_1\wedge \dots \wedge g_{\ell}\wedge h_1\wedge \dots \wedge h_k).
\endaligned
\end{equation}
Here we use the notation (\ref{smoothcorrespondence}) in (\ref{defqformula}) and
$$
\epsilon(h_1,\dots,h_k) = \sum_{i=1}^k i\deg'h_i +1, 
\qquad 
 \delta(g_1,\dots,g_{\ell}) = \sum_{i=1}^{\ell} \deg g_i.
$$
See \cite{ono} for the sign $\epsilon(h_1,\dots,h_k)$. %and the appendix ???\marginpar{appendix to be added.} for $\delta(g_1,\dots,g_{\ell})$.
For $\beta = \beta_0 = 0$, $\ell = 0$ we put
\begin{equation}\label{defqformula2}
\aligned
&\frak q^{\rm form}_{\beta_0;0,k}(h_1,\dots,h_k)
=
\begin{cases}
0 & k\ne 1,2 \\
dh_1 & k=1\\
(-1)^{\deg h_1} h_1 \wedge h_2
& k=2.
\end{cases}
\endaligned
\end{equation}
\begin{prop}\label{qproperties}
The operators $\frak
q^{\rm form}_{\beta;\ell,k}$ have the following properties:
\begin{enumerate}
\item
For each $\beta$ and $\text{\bf x} \in B_k(\Omega(L)[1])$,
$\text{\bf y} \in E_{\ell}(\Omega(X)[2])$, we have the following:
\begin{equation}\label{qmaineq}
0 =
\sum_{\beta_1+\beta_2=\beta}\sum_{c_1,c_2}
(-1)^\maltese
\frak q^{\rm form}_{\beta_1}(\text{\bf y}^{2;1}_{c_1};
\text{\bf x}^{3;1}_{c_2} \otimes
\frak q^{\rm form}_{\beta_2}(\text{\bf y}^{2;2}_{c_1};\text{\bf x}^{3;2}_{c_2})
\otimes \text{\bf x}^{3;3}_{c_2})
\end{equation}
where
$
\maltese = \deg'\text{\bf x}^{3;1}_{c_2} +
\deg'\text{\bf x}^{3;1}_{c_2} \deg \text{\bf y}^{2;2}_{c_1}
+\deg \text{\bf y}^{2;1}_{c_1}.
$
In $(\ref{qmaineq})$ and hereafter, we write $\frak q_{\beta}(\text{\bf y};\text{\bf x})$ in place of
$\frak q_{\beta;\ell,k}(\text{\bf y};\text{\bf x})$ if
$\text{\bf y} \in E_{\ell}(\Omega(X)[2])$, $\text{\bf x} \in B_{k}(\Omega(L)[1])$.
We use notation $(\ref{deltawritebyc})$ in $(\ref{qmaineq})$.
\item We have
\begin{equation}\label{qism}
\langle\frak q^{\rm form}_{\beta;\ell,k}(\text{\bf y};x_1,\dots,x_k),x_0\rangle_{\text{\rm cyc}} 
= (-1)^{\maltese}\langle\frak q^{\rm form}_{\beta;\ell,k}(\text{\bf y};x_0,x_1,\dots,x_{k-1}),x_k\rangle_
{\text{\rm cyc}},
\end{equation}
where $\maltese = \deg'x_0(\deg'x_1 + \dots + \deg' x_k)$.
Here\index[syindex]{<>_cyc@$\langle *,* \rangle_{\text{\rm cyc}}$}
\begin{equation}\label{pairing}
\langle a,b\rangle_{\text{\rm cyc}}  = (-1)^{\deg a}\int_L a \wedge b.
\end{equation}
\item Let $\text{\bf e}_L$ be the constant function on $L$ with value $1$, which we consider as a $0$-form. Let $\text{\bf x}_i \in B(\Omega(L)[1])$ and we put
$\text{\bf x} = \text{\bf x}_1 \otimes \text{\bf e}_L \otimes \text{\bf x}_2
\in B(\Omega(L)[1])$. Then
\begin{equation}\label{unital}
\frak q^{\rm form}_{\beta}(\text{\bf y};\text{\bf x}) = 0.
\end{equation}
except in the following cases:
\begin{equation}\label{unital2}
\frak q^{\rm form}_{\beta_0}(1_X;\text{\bf e}_L \otimes x) =
(-1)^{\deg x}\frak q^{\rm form}_{\beta_0}(1_X;x \otimes \text{\bf e}_L) = x,
\end{equation}
where $\beta_0 = 0 \in H_2(X,L;\Z)$ and $x \in \Omega(L)[1]
= B_1(\Omega(L)[1])$.
Note $1_X$ in $(\ref{unital2})$ is $1_X \in E_0(\Omega(X)[2])$. 
\begin{equation}\label{unital22}
\frak q^{\rm form}_{\beta_0}(y;\text{\bf e}_L) = (-1)^{\deg y}i_L^* (y).
\end{equation}
Here $y$ is a differential form on $X$ and $i_L^* (y)$ is its pull back to $L$. %\marginpar{(\ref{unital22}) is added. The sign to be checked.  KF 2025 July.}
\end{enumerate}
\end{prop}
In (\ref{unital}) we omit the suffix $k,\ell$ to $\frak q^{\rm form}_{\beta}$. It is automatically determined by the inputs 
$\text{\bf y}$, $\text{\bf x}$.  We will not repeat this kinds of remark.
\begin{proof}
Proposition \ref{qproperties} (1) is a consequence of
Proposition \ref{diskkura} (\ref{diskkura:boundary_strata}) and the compatibility of \blue{the CF-perturbations} 
with this boundary identification.
\par
Equation (\ref{qism}) follows from Proposition \ref{diskkura} (\ref{diskkura:cyclic})
and Proposition \ref{existmkulti1} (5).
\par
Proposition \ref{qproperties} (3) is a consequence of
Proposition \ref{diskkura} (\ref{diskkura:forget}) and the compatibility of the
CF-perturbation to this forgetful map.
See Section \ref{sec:CRperturb} and \cite[section 7]{fukaya:cyc}  for the detail of this point.
The proof of Proposition \ref{qproperties} is complete.
\end{proof}
We next explain the way how we use the map $\frak q$ to deform {the de Rham differential and the wedge product of differential forms to a filtered $A_{\infty}$ structure $\frak m$ on $\Omega(L)$.} In this section we use the {universal Novikov ring $\Lambda_{0}$, and its quotient by $T^{E}$ for some constant $E$. }
\par
We split $\frak b \in H^2(X;\Lambda_0)$ as a sum\index[syindex]{bxfra0@$ \frak b_0$}\index[syindex]{bxfra2@$ \frak b_2$}\index[syindex]{bxfra+@$ \frak b_+$}
\begin{align}\label{decomposefrakb}
  \frak b & = \frak b_0 + \frak b_2+\frak b_+ \\ \notag
  \frak b_0 & \in H^0(X;\Lambda_0) \\ \notag
  \frak b_2 & \in H^2(X;\F) \\ \notag
  \frak b_+ & \in H^2(X;\Lambda_+)\oplus \bigoplus_{k\ge 2} H^{2k}(X;\Lambda_0).
\end{align}
\par
For the rest of the paper, we fix representatives of the cohomology classes $\frak b_0$, $\frak b_+$, and $\frak b_2$ by
closed differential forms, which are denoted by the same letters. We assume that $\frak b_0$ is represented by a constant function.
%We fix the representatives throughout.
\begin{rem}
  We omit the proof (which follows the same argument as Chapter 27 in \cite{spectre}) that the structures of a filtered {curved} $A_{\infty}$ algebra 
on the de Rham complex of $L$ which we will define is independent  up to homotopy equivalence of the choice of the de Rham representative of $\frak b $. 
\end{rem}
\par
We assume $i_L^*[\frak b_2] = 0$ in $H^2(L;\F)$, 
where $i_L : L \to X$ is the inclusion map.  
Then there exists $\theta_L \in \Omega^1(L) \otimes \F$\index[syindex]{thetaL@$\theta_L$}
such that
\begin{equation}\label{fixthetaL}
d\theta_L+ i_L^*\frak b_2 = 0,
\end{equation}
and consider the data $(L,\theta_L)$ as fixed in this section.\index[syindex]{bxfrak2thetaLcap@$(\frak b_2,\theta_L) \cap \beta$}
\par
Let  $u : (\Sigma,\partial \Sigma) \to 
(X,L)$ be a map representing a homology class $\beta \in H_2(X,L;\Z)$. Then, by Stokes' theorem, it is easy to see that
\begin{equation}\label{capwiththetaL}
(\frak b_2,\theta_L) \cap \beta
:= \int_{\Sigma} u^*\frak b_2 + \int_{\partial\Sigma} u^*\theta_L
\end{equation}
depends only on the relative homology class $\beta$. 
\par
 {Recall that we have fixed an energy level $E$ for which our choice of $CF$-perturbation is given:}
\begin{defn}\label{deformedqdef}
For each $k \ne 0$, we define $\frak m_{k}^{{\rm form}, \frak b}$ by\index[syindex]{mxkbfra@$\frak m_{k}^{{\rm form},\frak b}$}
\begin{equation}\label{mkdefeq}
\aligned
&\frak m_{k}^{{\rm form},\frak b}(x_1,\ldots,x_k) \\
&= \sum_{\stackrel{\beta\in H_2(X,L:\Z)}{\omega(\beta) \leq E}}
\sum_{\ell=0}^{\infty} T^{\omega\cap \beta}
\frac{\exp((\frak b_2,\theta_L) \cap \beta) 
  )}{\ell!}
\frak
q^{\rm form}_{\beta;\ell,k}(\frak b_{+}^{\otimes\ell};
x_1,\ldots,x_k),
\endaligned
\end{equation}
where $x_i \in \Omega(L)$.
We extend the definition to {$\Omega(L)
\otimes \Lambda_0/T^E$}  by $\Lambda_0$-linearity. 
\par
For $k=0$, we define $\frak m_0^{\frak b}$ by
\begin{equation}\label{mkdefeq2}
\frak m_{0}^{{\rm form},\frak b}(1)
= {i_L^*}\frak b_0 + \sum_{\stackrel{\beta\in H_2(X,L:\Z)}{\omega(\beta) \leq E}}
\sum_{\ell=0}^{\infty} T^{\omega\cap \beta}
\frac{\exp((\frak b_2,\theta_L) \cap \beta)
}{\ell!} \frak
q^{\rm form}_{\beta;\ell,0}(\frak b_{+}^{\otimes\ell}; 1_L).
\end{equation}

\end{defn}
\begin{rem}\label{rem611}
We may relax the assumption $i^*_L([\frak b_2]) = 0$
to $i^*_L([\frak b_2]) \in H^2(L;2\pi \sqrt{-1}\Z)$.
In that case we take a complex line bundle on $L$ whose first Chern class 
is $i^*_L([\frak b_2])$ and take a $\C^*$-connection whose first
Chern form is $i^*_L(\frak b_2)/2\pi \sqrt{-1}$. This connection can then be added to the 1-form $\theta_L$.
We do not use this generalization here so do not discuss it further.
\end{rem}
\begin{prop}\label{Algebraonderhamcomplex}
The formal series $(\ref{mkdefeq})$ and $(\ref{mkdefeq2})$ converge in the $T$-adic topology. The operations $\{\frak m_{k}^{{\rm form},\frak b}\}_{k=0}^{\infty}$ define the 
structure of cyclic filtered $A_{\infty}$ algebra
on {$\Omega(L){\otimes} \Lambda_0$}, with strict unit
$\text{\bf e}_L \in \Omega^0(L)$.
\end{prop}
\begin{proof}
We can prove this proposition by a straightforward calculation 
using Proposition \ref{qproperties}.
\end{proof}
\begin{rem}\label{rem61313}
{If we did not decompose the forms $\frak b$ as in (\ref{decomposefrakb}), we would have to define the $A_\infty$ operations $\frak m_{k}^{{\rm form},\frak b}(x_1,\ldots,x_k)$ using the formula}
$$
\sum_{\beta\in H_2(X,L:\Z)}
\sum_{\ell=0}^{\infty}
\frac{T^{\omega\cap \beta}}{\ell!} \frak
q^{\rm form}_{\beta;\ell,k}(\frak b^{\otimes\ell};
x_1,\ldots,x_k).
$$
{However this expression may have infinitely many terms of energy less than $E$ even with the proviso that the energy of $\beta$  is bounded; so that it would not define a map over $\Lambda_0/T^E$. 
So we rewrite it to (\ref{mkdefeq}), which is a finite sum in this setting. } 
\end{rem}
The operations $\{\frak m_{k}^{{\rm form},\frak b}\}_{k=0}^{\infty}$ 
on {$\Omega(L)\widehat{\otimes} \Lambda_0$} together with 
homotopy  limit argument (explained in Subsection \ref{sec:inverse-limits-categ}) 
using pseudo-isotopy (explained in Subsection \ref{sec:pseudo-isotopy}),
we obtain a structure of cyclic filtered $A_{\infty}$ algebra
on {$\Omega(L) \widehat{\otimes} \Lambda_0$}, with strict unit
$\text{\bf e}_L \in \Omega^0(L)$.
We postpone the discussion about homotopy  limit argument 
until Subsections \ref{sec:constr-cycl-a_infty4} and \ref{sec:homotopyequiv}.
 
We call the Kuranishi structure of Proposition \ref{diskkura} and the CF-perturbation  of Proposition \ref{existmkulti1}
the ${\frak  q}$-Kuranishi structure etc.  
See Proposition \ref{Kuraeistspoly} for its category version.
Here $\frak  q$ is the 
symbol we use for the closed-open map.
We put superscript ${\frak  q}$ for the ${\frak  q}$-Kuranishi structure etc. 
in case we need to distinguish it from the other.
In Section \ref{sec:frakq} we introduce ${\frak  p}$-Kuranishi structure 
which is used to define the open-closed map ${\frak  p}$.

Using the operations $\frak m_{k}^{{\rm form},\frak b}$ we define the 
set of weak bounding cochains of $((L,\theta_L),\frak b)$ below.

\begin{defn}\label{boundingcochain}
Let $((L,\theta_L),\frak b)$ be given.\index[syindex]{brm1@$b_1$}\index[syindex]{brm+@$b_+$}
An element $b \in \Omega \widehat{\otimes} \Lambda_0$ is said to be a 
weak bounding cochain of  $((L,\theta_L),\frak b)$ 
if it satisfies the following:
\begin{enumerate}
\item
$b$ is a sum of $b = b_1 + b_+$ where 
$b_1$ is a closed one form that does not contain $T$ and
$b_+$ is in the direct sum
$$
b_+ \in ( \Omega^1  \widehat{\otimes}  \Lambda_+) \oplus \bigoplus_{k \ge 1}( \Omega^{2k+1}  \widehat{\otimes}  \Lambda_0)
$$
\item 
We consider a pair $(L,\theta_L + b_1)$. Note that $d(\theta_L + b_1) + i^*(\frak b_2) = 0$.
Let $\frak m_{k}^{{\rm form},b_1,\frak b}$\index[syindex]{mxkbfrab@$\frak m_{k}^{{\rm form},b_1,\frak b}$} is the $A_{\infty}$ operation (\ref{mkdefeq}) 
with $\theta_L$ replaced by $\theta_L + b_1$.
We then require:
\begin{equation}\label{S7MCeq}
\sum_{k=0}^{\infty} \frak m_{k}^{{\rm form},b_1,\frak b}(b_+,\dots,b_+) 
= c {\bf e}_L.
\end{equation}
Here ${\bf e}_L$ is the differential $0$-form $1$ on $L$ and
$c \in \Lambda_+$.
\end{enumerate}
\end{defn}

We denote by $\widehat{\mathcal M}_{\text{\rm weak}}(\Omega(L);(\frak b,\theta_L);\Lambda_0)$
\index[syindex]{MweakhatLbtheta@$\widehat{\mathcal M}_{\text{\rm weak}}(\Omega(L);(\frak b,\theta_L);\Lambda_0)$} the set of all weak bounding cochains.
\begin{lem}\label{lem:Tadic-convergence}
 The left hand side of \eqref{S7MCeq} converges in $T$-adic topology.
\end{lem}
\begin{proof}
We decompose $b_+=b^1_{+} + b^{\rm high}$ as above
where $b_+^1 \in  \Omega^1  \widehat{\otimes}  \Lambda_+$
and $b^{\rm high} \in \bigoplus_{k>1}( \Omega^{2k+1}  \widehat{\otimes}  \Lambda_0)$.
We remark that
$$
\lim_{\ell \to \infty}\inf \{ \beta \cap \omega \mid \mu(\beta) \ge \ell, \mathcal M(L;\beta) \ne \emptyset\} = \infty.
$$
We consider the decomposition of the left hand side of  (\ref{S7MCeq})
as
$$
\sum_{\beta}\sum_{k=0}^{\infty} T^{\beta\cap \omega}\frak m_{k,\beta}^{{\rm form},b_1,\frak b}(b_+,\dots,b_+).
$$
Then modulo $T^E$ for any $E$, all except finitely many homology classes $\beta$ appear.
We next consider the sum:
\begin{equation}\label{form:716}
\sum_{k=0}^{\infty} T^{\beta\cap \omega}\frak m_{k,\beta}^{{\rm form},b_1,\frak b}(b_+,\dots,b_+)
\end{equation}
and decompose it as 
$$
\sum_{k=0}^{\infty} \sum_{j_1,\dots,j_m \in \{0,\dots,k\}}
\frak m_{k,\beta}^{{\rm form},b_1,\frak b}(b_+^1,\dots,b^{\rm high},\dots,\dots,b^{\rm high}, \dots, b^1_+)
$$
where $b^{\rm high}$ are the $j_1,\dots,j_m$-th inputs and all other inputs are 
$b_+^1$. By the degree reason we have
$$
2m \le \mu(\beta) + 2.
$$
On the other hand there exists $\rho > 0$ such that 
$b_+^1 \equiv 0 \mod T^{\rho}$.
Therefore (\ref{form:716}) converges in $T$-adic topology.
Thus (\ref{S7MCeq}) also converges.
\end{proof}

\begin{defn}
Let $((L,\theta_L),\frak b)$ be given.
\begin{enumerate}
\item
A pair $(b(0), b(1))$ of  bounding cochains
are {\it gauge equivalent} if
there exists 
$$
b(t) + dt \wedge c(t)  \in \Omega^{\rm odd}(L\times \R) 
\widehat{\otimes} \Lambda_0
$$
such that %\marginpar{Degree of $b$ is changed to odd.}
$b(t) \in \widehat{\mathcal M}_{\text{\rm weak}}(\Omega(L);(\frak b,\theta_L);\Lambda_0)$ for each $t$ and that
$$
\frac{d b(t)}{dt} + \frak m^{{\rm form},\text{\bf b}(t)}_{1}(c(t)) = 0
$$
where $\text{\bf b}(t) = (\frak b,\theta_L,b(t))$.
(See  \cite[Definition 4.3.1 and Proposition 4.3.5]{fooo09}.)
We can show that this is an equivalence relation.
(See \cite[Lemma 4.3.4]{fooo09}.)
\item
Let ${\mathcal M}_{\text{\rm weak}}(\Omega(L);(\frak b,\theta_L);\Lambda_0)$
be the set of all gauge equivalence classes of 
 bounding cochains for 
$(L,\frak b,\theta_L)$.
\item
For $b \in \widehat{\mathcal M}_{\text{\rm weak}}(\Omega(L);(\frak b,\theta_L);\Lambda_0)$
we define $\frak{PO}^{\frak b}(b) \in \Lambda_+$ by the equation
$$
\frak m_{0}^{\rm form,{\text{\bf b}}}(1) = \frak{PO}^{\frak b}(b)\cdot\text{\bf e}_L.
$$
It defines a map
$$
\widehat{\mathcal M}_{\text{\rm weak}}(\Omega(L);(\frak b,\theta_L);\Lambda_0)
\to 
\Lambda_+.
$$
In the same way as \cite[Lemma 4.3.22]{fooo09}, we  show that it descends to a map
$$
{\mathcal M}_{\text{\rm weak}}(\Omega(L);(\frak b,\theta_L);\Lambda_0)
\to 
\Lambda_+,
$$
which we also write $\frak{PO}^{\frak b}$ and call the
{\it potential function}.
\end{enumerate}
\end{defn}
Hereafter we write $
{\mathcal M}_{\text{\rm weak}}(L;\frak b;\Lambda_0)$ \index[syindex]{MweakLb@${\mathcal M}_{\text{\rm weak}}(L;\frak b;\Lambda_0)$}
in place of 
$
{\mathcal M}_{\text{\rm weak}}(\Omega(L);(\frak b,\theta_L);\Lambda_0)$
for simplicity.

\section{Cyclic symmetry and unitality.}
\label{sec:two perturbation}
 %\marginpar{This section is new.}
We start with a bit technical discussion about the construction 
of (abstract) perturbation to obtain a cyclically symmetric, unital, and filtered $A_{\infty}$ category.
We use the moduli spaces ${\mathcal M}_{\ell;k+1}(L;\beta)$ of pseudo-holomorphic disks 
 with a boundary condition given 
by a single (embedded) Lagrangian submanifold $L$.
We also use the
moduli space of  pseudo-holomorphic polygons ${\mathcal M}_{\ell;\vec k}((\vec{\kappa},\vec p);B)$
(see Section \ref{sec:cyclicfil} for this notation.)
with the boundary 
condition given by a finite set of Lagrangian submanifolds.
The former gives an $A_{\infty}$ algebra and the latter 
an $A_{\infty}$ category.
We want to use some of the properties mentioned in Section \ref{sec:a_infty-categ-over}.
\begin{equation}\label{cyclicsymdef}
\aligned
\langle \frak m^{\rm form}_k(x_0,\dots,x_{k-1}),x_k\rangle = (-1)^{\maltese}
\langle \frak m^{\rm form}_k(x_k,x_0,\dots,x_{k-2}),x_{k-1}\rangle  \\
\frak m^{\rm form}_k(x_0,\dots,x_{i-1},{\bf e},x_{i+1}\dots,x_{k})
= 
\begin{cases}
0   &  k\ne 2 \\
(-1)^{\deg x_0}x_0   &  k=2, i=1 \\
x_1   &  k=2, i=0
\end{cases}
\endaligned
\end{equation}
where $\maltese = \deg'x_k (\deg'x_0+\dots+\deg'x_{k-1})$.  %\marginpar{sign added.  KF.}
${\bf e}$ is the unit which is the degree $0$-form $1$ on 
$L_{\kappa_i}$.  ($L_{\kappa_i} = L_{\kappa_{i+1}}$).
The first formula is the cyclic symmetry (\ref{eq:sign_cyclic_invariance}), the second is the unitality.
For the cyclic symmetry we need to make the perturbation we use to be symmetric under the cyclic 
permutation of the boundary marked points. 
For the unitality, we need the perturbation to be compatible 
with respect to the forgetful map of the boundary marked points.
In the case of a single embedded Lagrangian submanifold 
we can find such a perturbation as we  discussed in Section \ref{ainfalgasssingle}.
However in the case of several (mutually transversal) 
Lagrangian submanifolds (or of a single immersed Lagrangian 
submanifold) it is actually impossible to find 
a perturbation so that both of the properties are satisfied,
in case we use de Rham model.\footnote{For other model 
such as singular homology or Morse homology, it seems even harder.}

In fact the following explicit example illustrates such a situation.\footnote{The authors thank L. Amorim 
for explaining this point.}
Let $L_0$ and $L_1$ be two embedded Lagrangian submanifolds which intersect 
transversally at one point $p$.
We consider the moduli space of pseudo-holomorphic maps 
$
u: D^2 \to X
$
such that 
$$
u(\exp (2\pi\sqrt{-1}\theta)) 
\in 
\begin{cases} 
L_0      & \text{ for $-1/3 < \theta < 1/3$}   \\
L_1       & \text{ for  $1/3 < \theta <2/3$}.
\end{cases}
$$
We regard $z_0 = 1$, $z_1 = \exp (2\pi\sqrt{-1}/3)$, $z_2 = \exp (4\pi\sqrt{-1}/3)$
as marked points in $\partial D^2$.
\begin{figure}[h]
\centering
\includegraphics[scale=0.4]{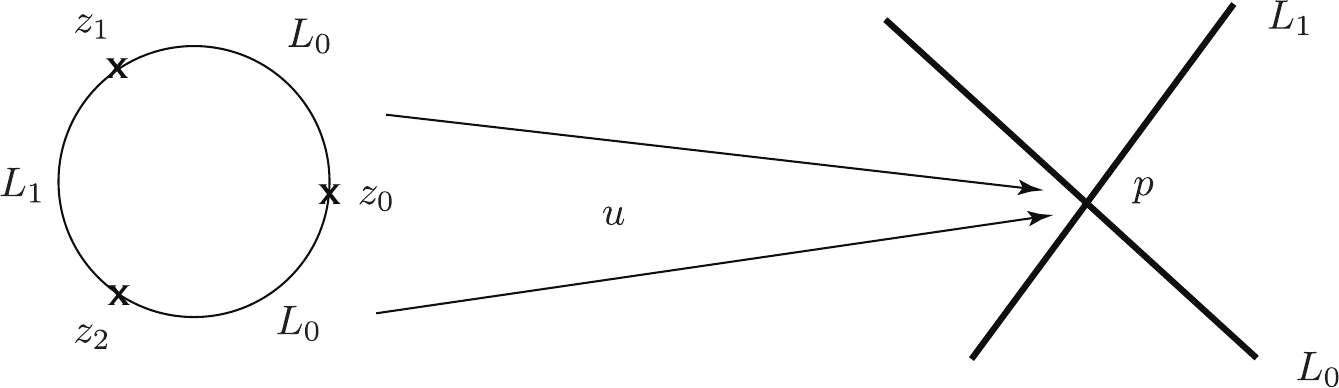}
\caption{Constant map}
\label{Figure7-1}
\end{figure}
We consider the case when the homology class of $u$ is zero.
So actually this moduli space consists of one point that is a constant map 
with value $p$. See Figure \ref{Figure7-1}.
This moduli space is transversal.
The issue is that in order to define operations in the de Rham model, the
evaluation map at the marked point corresponding to the output 
(that is, $z_0$) should be submersive.
This is because the push out or the integration along the fiber of a 
differential form becomes a smooth form only under such assumptions.
The way employed in \cite{springer} and etc. to resolve this issue is 
to use a continuous family of perturbations (CF-perturbation).
In this situation, we add an extra obstruction bundle  $E$
so that the thickened moduli space, $\mathcal M^+$ 
has positive dimension and that 
the evaluation map ${\rm ev}_0 :  \mathcal M^+ \to L_0$ 
is a submersion. (Here ${\rm ev}_0(u) = u(1)$.)
Now we take a parametrized family of sections $s_{\xi}$ of $E$
such that the map from $\{(u,\xi) \mid u \in s_{\xi}^{-1}(0)\}$
to $L_0$ defined by ${\rm ev}^+_0(u,\xi)= u(1) \in L_0$ is a submersion.
(We also require that $s_{0}^{-1}(0)$ is the original moduli space 
consisting of one point.)
We take a compactly supported smooth form $\chi$ of top degree with $\int \chi = 1$ on the parameter space.
Then the virtual fundamental chain associated to this CF-perturbation is
$$
\rho = ({\rm ev}^+_0)_! \chi.
$$
This is a smooth differential $\dim L_0$-form on $L_0$ supported on a small 
neighborhood of $p$ and $\int_{L_0}\rho = 1$.
In other words, it is a smooth approximation of the delta form supported at $p$.
Thus after the perturbation our operation is
$$
\frak m^{\rm form}_2(p_{01},p_{10}) = \rho
$$
Here $p_{01}$ is $p$ regarded as an intersection point $L_0 \cap L_1$ 
and $p_{10}$ is $p$ regarded as an intersection point $L_1 \cap L_0$.
Note that in the singular homology model we can take 
the right hand side to be the $0$ chain $p$ of $L_0$.
\par
Now we require the operation to be cyclically symmetric 
and unital.
It means that we may cyclically change the enumeration 
of the marked points to $z_1 = 1$, $z_2 = \exp (2\pi\sqrt{-1}/3)$, $z_0 = \exp (4\pi\sqrt{-1}/3)$.
Then $z_1 = 1$ corresponds to the input of the operation.
So for the sake of unitality we need the perturbation 
to be compatible with the forgetful map.
The CF-perturbation we used above is {\it not} compatible with forgetful map. %\marginpar{Figures}
\par
\begin{figure}[h]
\centering
\includegraphics[scale=0.4]{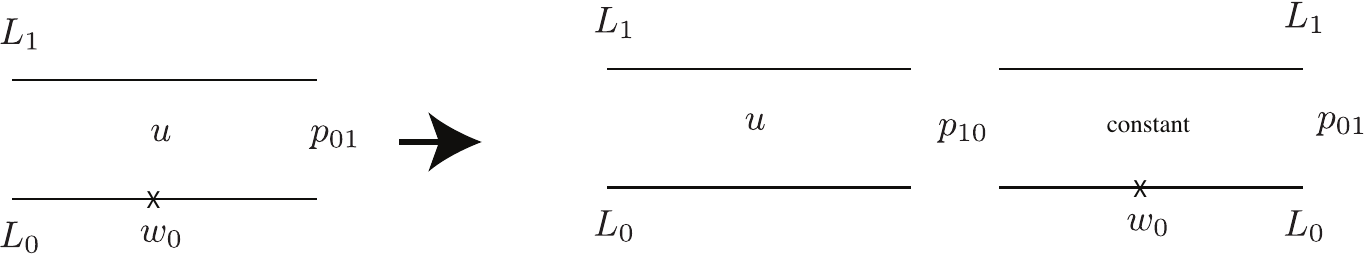}
\caption{Constant triangle bubbles}
\label{Figure7new2}
\end{figure}

In fact we consider the configuration of 
one strip $[0,1]\times \R \to X$ between $L_0$ and $L_1$ with the ends being asymptotic to $q$ 
and $p$, respectively.
We consider one marked points $w_0 \in \{0\} \times \R$.
We want to find a perturbation so that it is compatible with 
the forgetful map at $w_0$.
It can be induced from the translation invariant 
perturbation of the strip (without marked point).
We however consider the limit when $w_0$ goes to $+\infty$.
Then the limit in the stable map topology
consists of two components, one is 
a strip without marked point, 
the other is a constant map to $p$ with 3 marked points 
2 of which are mapped to $p_{01}$, $p_{10}$, and the third is mapped
to somewhere on $L_0$.
(See Figure \ref{Figure7new2}.)
This is exactly the configuration discussed above.
Therefore if we introduce non-trivial perturbation 
for this second component, then 
the perturbation of the moduli space of strip with one marked point 
can not be compatible with forgetful map 
of $w_0$.

A  way to resolve this issue is to replace 
various symmetries by homotopy (or infinite homotopy) symmetries.
In this paper we take the route of replacing the unitality by the homotopy unitality 
but keeping cyclic symmetry to be exact. 
(There is a homotopy version of cyclic symmetry, which is 
called Calabi-Yau property. We do not use it in this paper.) 
Strict unitality will be restored when we go to the canonical model in 
Section \ref{canonical}.

\begin{rem}
In \cite{fukaya:functor} a  strictly unital filtered $A_{\infty}$ category 
is constructed from a finite set of relatively spin Lagrangian submanifolds in a similar situation.
In \cite{fukaya:functor} a strictly unital filtered $A_{\infty}$ category 
without cyclic symmetry is constructed.
\end{rem} %\marginpar{Remark added.  KF. 2025 Feb}

\section{Cyclic filtered $A_{\infty}$ category.}
\label{sec:cyclicfil}
  %\marginpar{Beginning is a copy from previous Section 5}
We fix a bulk class $\mathfrak b \in H^*(X;\Lambda^{\C}_0)$ in this section. We start from a finite collection ${\mathscr L}$ as in the 
introduction and construct a curved cyclic filtered $A_{\infty}$
category $\cL^{\rm form}_{\rm curve}$, whose object set ${\mathscr L}$ is a finite collection  
${\mathscr L} = \{(L_{\kappa},\theta_{\kappa})\mid {\kappa}=1,\dots, \# {\mathscr L}\}$ such that:
 %\marginpar{
%Actually we may take several $\theta_{\kappa}$ and it seems possible to take infinitely many of them. KF 2024 June}
\begin{conds}\label{conds614}
$  $ \par
\begin{enumerate}
\item  $L_{\kappa}$ is an immersed Lagrangian submanifold equipped with 
a relative spin structure with  background class $\text{\rm st}
\in H^2(X;\Z_2)$ independent of $\kappa$.
We assume $L_{\kappa}$ is connected and self-transversal.
\item $\theta_{\kappa}$ is a $1$-form on $L_{\kappa}$ satisfying Equation \eqref{fixthetaL}.
\item
If $L_{\kappa}\ne L_{\kappa'}$, then the Lagrangian submanifolds $L_{\kappa}$ and $L_{\kappa'}$ intersect transversely.
The intersection is disjoint from self-intersection. 
\end{enumerate}
\end{conds}
\begin{rem}
It may happen that $L_{\kappa} = L_{\kappa'}$ but  $\theta_{\kappa} \ne \theta_{\kappa'}$.
We need to include such cases so that we consider two bounding cochains $b$ of $L_{\kappa}$ 
which have  different $H^2(L_{\kappa};\F)$ components.  (See Definition \ref{boundingcochain} (2).) %\marginpar{Remark added KF 2024 Dec.}
\end{rem}

In this section we prove:
\begin{thm}\label{cAinfconst}
There exists a  curved cyclic gapped filtered $A_{\infty}$ category $\cL^{\rm form}_{\rm curve}$ \index[syindex]{Lcateformcurve@{$\cL^{\rm form}_{\rm curve}$}}
over $\Lambda_0$ whose set of objects  
is ${\mathscr L}$, {with morphism spaces from $L_i$ to $L_j$ given by the completion of the space of differential forms on 
$\tilde L_i \times_X \tilde L_j$ with coefficients in $\Lambda_0$, and with operations induced by choices of virtual fundamental chains on the moduli spaces of holomorphic disks with boundary conditions on the elements of ${\mathscr L}$.}
 %\marginpar{unital is removed}
\end{thm}
Corollary \ref{existsderham} is a version of Theorem \ref{cAinfconst} and is proved in 
Subsection \ref{constcyclic}.
The proof of Theorem \ref{cAinfconst} is completed in Subsection \ref{sec:constr-cycl-a_infty4}.
\begin{rem}
We remark that the $A_{\infty}$ category in Theorem \ref{cAinfconst} is 
{\it not} unital.  The reason was explained in Section \ref{sec:two perturbation}.
We need unitality for  applications of our main theorems, since 
in many important cases such as the toric case, 
we need to study weakly unobstructed Lagrangian 
submanifolds rather than unobstructed ones. 
For this purpose we need unit.

We recall: %\marginpar{The lest of this subsection is rewritted.}
\begin{defn}\label{defn84}
A {\em homotopically unital curved filtered $A_{\infty}$ category}\index{homotopically unital}
is a pair of curved filtered $A_{\infty}$ category $\Cat$
and a strictly unital curved filtered $A_{\infty}$ category $\Cat^+$\index[syindex]{Czatplus@ $\Cat^+$}
together with a linear filtered $A_{\infty}$ functor $\iota : \Cat \to \Cat^+$ 
with the following properties.
\begin{enumerate}
\item $\iota$ is a bijection on objects.
\item  If $c \ne c'$, $c,c' \in \Ob(\Cat)$ then $\Cat^+(c,c') =  \Cat(c,c')$.
\item For any object $c \in \Ob(\Cat)$ there exist ${\bf e}_c \in \Cat(c,c)$
of degree $0$ and ${\bf f}_c \in \Cat^+(c,c)$ of degree $-1$, 
such that ${\bf e}_c$ together with ${\bf f}_c$ and the unit ${\bf e}^+_c \in \Cat^+(c,c)$
makes $\Cat(c,c)$ to be homotopically unital filtered $A_{\infty}$ algebra 
in the sense of \cite[Definition 3.3.2]{fooo09}.
In particular\index[syindex]{eBplusbf@${\bf e}^+_c$}\index[syindex]{fplusbf@${\bf f}_c$}
$$
\Cat^+(c,c) =  \Cat(c,c) \oplus \Lambda_0{\bf e}^+_c  \oplus \Lambda_0{\bf f}_c.
$$
\item ${\bf e}^+_c,{\bf f}_c \notin \Lambda_+\Cat^+(c,c)$.  
\item
The image of all the operations including ${\bf f}_c$ are in $ \Cat(c,c')$
except $\frak m^+_1({\bf f}_c)$, $\frak m^+_2({\bf e}^+_c,{\bf f}_c)$, 
$\frak m^+_2({\bf f}_c,{\bf e}^+_c)$.
\end{enumerate}
\end{defn}
Note that other than unitality of ${\bf e}^+_c$ item (3) means:
\begin{equation}\label{condforrf}
\frak m^+_1({\bf f}_c) - ({\bf e}^+_c - {\bf e}_c) \equiv 0 \mod \Lambda_+ \Cat(c,c)
\end{equation}
here $\frak m^+_*$ is the structure operation of $\Cat^+$.
In particular the left hand side does not contain 
${\bf f}_c$ or ${\bf e}^+_c$.
In the other cases, when input of $\frak m^+_*$  contains ${\bf f}_c$, the output is in 
$\Lambda_+ \Cat(c,c)$.\footnote{In \cite[Section 3.3]{fooo09} the definition was written 
in such a way that only  {\it existence} of $\Cat^+$ is required for $\Cat$ to be homotopically 
unital.  In the above definition $\iota$ and $\Cat^+$ are parts of the data consisting homotopically 
unital $A_{\infty}$ category.  This difference however is not important because of \cite[Lemma 4.2.55]{fooo09}.
The formulation of this paper (that is, including 
the operations involving ${\bf f}$ as a 
part of the data) is better suited 
for the definition of cyclicity of 
homotopically unital $A_{\infty}$ category. (Definition \ref{def98}.)
Nevertheless we sometimes say $\Cat$ is homotopically unital (instead of saying 
$(\Cat,\Cat^+,\iota)$ is homotopically unital) by an abuse of notation.}
\par
From now on we omit $^+$ in $\frak m^+_*$ and simply write it as $\frak m_*$.

In Section \ref{sec:unit} we prove:
\begin{thm}\label{homotopyunitality}
The curved cyclic filtered $A_{\infty}$ category $\cL^{\rm form}_{\rm curve}$
in Theorem {\rm \ref{cAinfconst}} is  homotopically unital.
\end{thm}
We write $\cL^{\rm form}_{\rm c.u.}$ \index[syindex]{Lcateformcu@{$\cL^{\rm form}_{\rm c.u.}$}} in place of $\cL^{\rm form}_{\rm curve}$
\index[syindex]{Lcateformcurve@{$\cL^{\rm form}_{\rm curve}$}} when we include 
the structure of homotopy unit.\footnote{Here ${\rm c.u.}$ stands for curved and homotopically unital. 
$\rm c$ does not stand for cyclic.  In this paper our moduli spaces with structures are always 
symmetric with respect to the cyclic permutation of boundary marked points.}

\begin{rem}
See Definition \ref{def98} for the definition of the  cyclicity of a homotopically unital  $A_{\infty}$ category.
\end{rem}

Note the $A_{\infty}$ category in Theorem \ref{homotopyunitality} is still curved.
In Section \ref{sec:elicurv}  we include bounding cochain and prove the next theorem.

\begin{defn}
A unital filtered $A_{\infty}$ category is said to be \emph{weakly 
curvature free}\index{weakly curvature free} if all objects are weakly 
curvature free in the sense of Definition \ref{defn:a_infty-categories}.
\end{defn}

\begin{thm}\label{cAinfconstuni}
There exists a unital  gapped filtered $A_{\infty}$ category $\cL^{\rm form}_{\rm uni}$\index[syindex]{Lformuni@$\cL^{\rm form}_{\rm uni}$} 
over $\Lambda_0$ such that:
\begin{enumerate}
\item
$\cL^{\rm form}_{\rm uni}$ is weakly 
curvature free.
\item
Its set of objects $\Ob(\cL^{\rm form}_{\rm uni})$ is 
the set of all pairs $(L_{\kappa},b)$ 
where $L_{\kappa}  \in \mathscr L$ and $b$ is a  bounding cochain, 
with respect to the filtered $A_{\infty}$ structure $\frak m_k^{\rm form,\frak b}$ 
obtained in Section $\ref{ainfalgasssingle}$.
\item 
Let  $c = (L_{\kappa},b),c' = (L_{\kappa'},b') \in \Ob(\cL^{\rm form}_{\rm uni})$.  If  $c \ne c'$
then the morphism space  $\cL^{\rm form}_{\rm uni}(c,c')$ is 
equal to $\cL^{\rm form}_{\rm c.u.}(L_{\kappa},L_{\kappa'})$.
\item 
If $c = c'$ then
$$
\cL^{\rm form}_{\rm uni}(c,c) = \cL^{\rm form}_{\rm c.u.}(L_{\kappa},L_{\kappa})^+.
$$
Here $^+$ is as in Definition {\rm \ref{defn84}}.
Namely the right hand side has two more generators ${\bf e}^+_c$, ${\bf f}_c$
compared to $\cL^{\rm form}_{\rm c.u.}(L_{\kappa},L_{\kappa})$.
\item
${\bf e}^+_c$ is the strict unit of $c$. The formula $(\ref{condforrf})$ holds.
\item
The structure operations of $\cL^{\rm form}_{\rm uni}$ among 
elements other than ${\bf e}^+_c$, ${\bf f}_c$ 
are obtained from those of $\cL^{\rm form}_{\rm c.u.}$
by deforming the latter via  bounding cochains.
\end{enumerate}
\end{thm}
The proof of Theorem \ref{cAinfconstuni} is completed in Section \ref{sec:elicurv}.

We remark that we do not claim $\cL^{\rm form}_{\rm uni}$ to be cyclically symmetric.
The reason is that the inner product $\langle *,* \rangle_{\rm cyc}$ is not defined in case 
${\bf e}^+_c$, ${\bf f}_c$ are involved.  However a certain cyclic symmetry holds.
See (\ref{qismcatprime}).  

In Section \ref{canonical}, we finally obtain a unital cyclic gapped filtered $A_{\infty}$ category $\cL$ 
so that the morphism space is the cohomology group $H(L_{\kappa})$
in the case of $c = (L_{\kappa},b),c' = (L_{\kappa'},b')$ with 
$L_{\kappa} = L_{\kappa'}$, $b \equiv b' \mod \Lambda_+$  (Corollary \ref{cor117}).

Because the choice of a virtual fundamental chain on the moduli space of holomorphic curves is not unique, the construction of the categories $\cL^{\rm form}_{\rm curve}$,
 $\cL^{\rm form}_{\rm uni}$ depends on various choices. 
 We will prove that this resulting category does not depend on these choices up to cyclic pseudo-isotopy 
 as Theorem \ref{thm111}.
 %Our techniques can be used to show that this resulting category does not depend on these choices up to cyclic pseudo-isotopy, though we do not include the proof in this paper. 
 In Subsection \ref{relative}, %we do however 
we discuss a relative version of the above statement, in which we show that we can extend the choice of data to define a category with objects $\mathscr L$ to include an additional Lagrangian.
\end{rem}\

\subsection{Construction of the cyclic $A_{\infty}$ category 1 : 
the moduli space}
\label{cycliccategorymoduli}
Let $\mathscr L$ be a finite collection of Lagrangians satisfying conditions laid out in
Condition \ref{conds614}.
\begin{defn}\label{def2110}
Let $\vec{\kappa} = ({\kappa_0}, \dots,\kappa_K)$, $\kappa_i \in \{ 1,\dots, \#\mathscr L\}$.
Here and hereafter we put $\kappa_{K+1} = \kappa_{\blue{0}}$ 
by convention 
\par
 We take also $\vec p =\{ p_{i}\mid i= {0}, \dots,K\}$ such that:
 \begin{enumerate}
 \item
 If  $L_{\kappa_{i-1}} \ne L_{\kappa_i}$ then
$p_{i} \in L_{\kappa_{i-1}} \cap L_{\kappa_{i}}$.
\item
 If  $L_{\kappa_{i-1}} = L_{\kappa_i}$ then $p_i = (p_i^1,p_i^2)$ 
 such that $p_i^1,p_i^2 \in \tilde L_{\kappa_i}$ with 
 $p_i^1  \ne p_i^2$ and $i_{L_{\kappa_i}}(p_i^1) = i_{L_{\kappa_i}}(p_i^2)$. In other words it is a self-intersection 
 point.\footnote{We remark that when $L_{\kappa}$ is immersed we identify it with a pair $(\tilde L_{\kappa},i_{L_{\kappa}})$
such that $i_{L_{\kappa}} : \tilde L_{\kappa} \to X$ is an immersion with image $L_{\kappa}$.}
\end{enumerate}
In the exceptional case $K=0$ we include the case $\vec p = \emptyset$. Namely we may or may not 
take $p_0$ as in (2).
\par
We denote by $\text{\rm Seq}_{K}$\index[syindex]{Seq@$\text{\rm Seq}_{K}$}
the collection of all such $(\vec{\kappa},\vec p)$.
 %\marginpar{Definitions \ref{def2110} and \ref{def211} are modified so that immersed cases is 
%included.  KF 2025 Aug}
\end{defn}
 
\begin{defn}\label{def211}
Let $(\vec{\kappa},\vec p) \in \text{\rm Seq}_{K}$.
We denote by $\Pi_2(\vec{\kappa},\vec p)$ the 
set of equivalence classes of pairs 
$(u,\vec z)$ consisting of a  map 
$u : D^2 \to X$ and a sequence $\vec z
= ({z_0}, \dots,z_{K})$
with the following properties:
\begin{enumerate}
\item {The points $z_i$ lie in 
$\partial D^2$, and the ordering 
$(\blue{z_0}, \dots,z_{K})$ respects the} counter-clockwise 
orientation of $\partial D^2$.
\item
Let $\{\overline{z_{{i}} z_{{i+1}}}\}_{i=0}^{K}$ be { the decomposition of  
$\partial D^2$ into intervals with endpoints
$\{z_{\blue{i}}, z_{\blue{i+1}}\}$ and disjoint interiors}. Then
$u(\overline{z_{\blue{i}} z_{\blue{i+1}})} \subset L_{\kappa_i}$.
Here $z_{K+1} = z_0$ by convention.
\item
If $L_{\kappa_{i-1}} \ne L_{\kappa_i}$ then $u(z_{i}) = p_{i}$.
\item 
If  $L_{\kappa_{i-1}} = L_{\kappa_i}$ the following {\it switching condition} is satisfied.  \index{switching condition}
Let $z_j^{\uparrow} \in \overline{z_{{i-1}} z_{{i}}}$ and $z_j^{\downarrow} \in \overline{z_{{i}} z_{{i+1}}}$ 
be sequences converging to $z_i$.  We have uniquely the limit $\tilde u(z_j^{\uparrow}) \in \tilde L_{\kappa_{i-1}}$, 
$\tilde u(z_j^{\downarrow})  \in \tilde L_{\kappa_{i}}$ 
such that 
$
i_{L_{\kappa_{i-1}}}(\tilde u(z_j^{\uparrow})) =  u(z_j^{\uparrow})$, 
$i_{L_{\kappa_{i}}}(\tilde u(z_j^{\downarrow})) =  u(z_j^{\downarrow}).
$
We then require:
\begin{equation}\label{formswitch}
\lim_{j\to \infty}(\tilde u(z_j^{\uparrow}),\tilde u(z_j^{\downarrow})) 
= (p_i^1,p_i^2).
\end{equation}
(See \cite[Definition 3.17 (5)]{fukaya:cyc}.)
\end{enumerate}
\par
We say two elements $(u_1,\vec z)$, 
$(u_2,\vec z\,')$ are equivalent each other if
there exists a continuous map $f: W \to X$ enjoying the following properties: 
\begin{enumerate}
\item  $W$ is a compact 3-manifold with boundary and corners. \\
\item {We have a decomposition of the boundary}
  \begin{equation}
    \partial W = \Sigma_1 \cup \Sigma_2 
\cup \bigcup_{i \in {\mathbb Z}/(k+1){\mathbb Z}} S_{i, i+1}
  \end{equation}
into compact 2-manifolds such that  
\begin{enumerate}
\item[(i)] $\Sigma_a \cap S_{i, i+1}$ is the part of $\Sigma_a$ between $z_i$ and $z_{i+1}$ in $\partial \Sigma_a$, $a=1,2$ and 
\item
[(ii)] $S_{i, i+1} \cap S_{j, j+1} = \emptyset$ if $i$ and $j$ are not adjacent, $S_{i, i+1} \cap S_{i+1, i+2} = I_{i+1} \cong [0,1]$.  
\end{enumerate}
\item The restriction of $f$ to $\Sigma_a$ coincides with $u_a$. 
\item The map $f$ sends $S_{i, i+1}$ to $L_{\kappa_i}$.  In particular, $f(I_i)=p_i$ if $L_{\kappa_{i-1}} \ne L_{\kappa_i}$.  
A similar switching condition as (\ref{formswitch}) holds if $L_{\kappa_{i-1}} = L_{\kappa_i}$.
\end{enumerate} 
{We denote by $\Pi_2(\vec{\kappa},\vec p)$\index[syindex]{Pi2kappap@$\Pi_2(\vec{\kappa},\vec p)$} 
the set of all equivalence classes of such pairs $(u,\vec z)$.}
\end{defn}

\begin{defn}\label{defn615615}
Let 
$(\vec{\kappa}',\vec p') \in \text{\rm Seq}_{K'}$, 
$(\vec{\kappa}'',\vec p'') \in \text{\rm Seq}_{K''}$.
Assume
$p'_{i+1}=p''_0$, $\kappa'_i=\kappa''_0$ and $\kappa'_{i+1}=\kappa''_{K''}$.  
Then we define
\begin{equation}\label{213formula}
(\vec{\kappa}'',\vec p'')\#_{i,i+1} (\vec{\kappa}',\vec p')
=: (\vec{\kappa},\vec p)
\end{equation}
where we arrange
$$
\aligned
\vec{\kappa} &= (\blue{\kappa'_0},\dots, \kappa'_{i} \blue{=\kappa''_0}, \kappa''_{\blue{1}},\dots,\kappa''_{K''}\blue{=\kappa'_{i+1}}, 
\kappa'_{i+2}, \dots,\kappa'_{K'}),
\\
\vec p
&= (\blue{p'_{0}},\dots,p'_{\blue{i}}, p''_{1},
\dots,p''_{\blue{K''}},p'_{i+\blue{2}},\dots,p'_{K'}).
\endaligned$$
\par
Then it induces the map
\begin{equation}
\#_{i, i+1} 
:
\Pi_2(\vec{\kappa}'',\vec p'') 
\times \Pi_2(\vec{\kappa}',\vec p') 
\to
\Pi_2((\vec{\kappa}'',\vec p'')\#_{i,i+1} (\vec{\kappa}',\vec p'))
\end{equation}
by the obvious concatenation at   %
$p'_{i+1} = p''_{0}$. 
\end{defn}
\begin{figure}[h]
\centering
\includegraphics[scale=0.7]{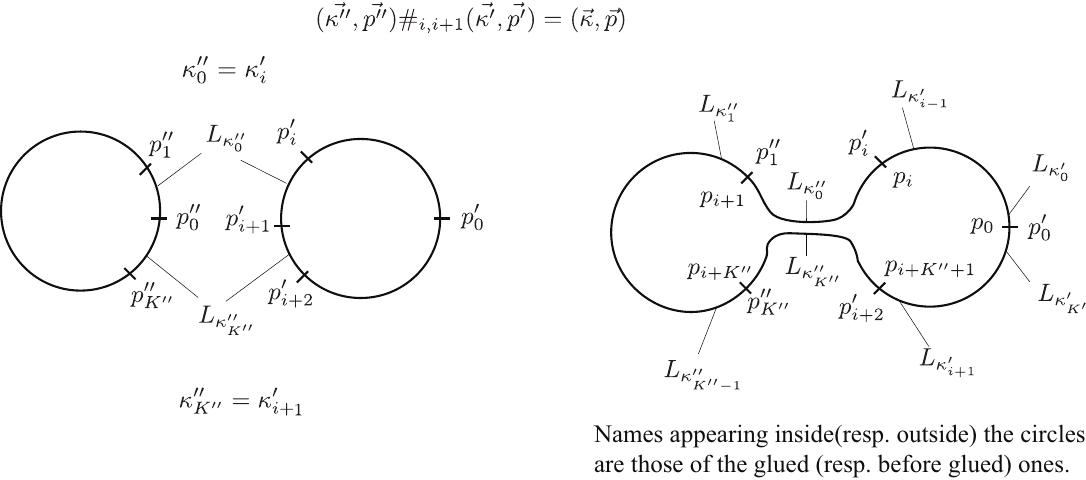}
\caption{Gluing in Definition \ref{defn615615}}
\label{Figure6-2-1}
\end{figure}
\begin{rem}
  {Note that the description of this operation on relative homology classes breaks the symmetry between the boundary marked points, since gluing always takes place at the marked point labelled by $0$ on the disk to the left of Figure \ref{Figure6-2-1}. One could restore symmetry by allowing gluing at arbitrary marked points for which the given boundary conditions agree, but the result would give a redundant list of boundary strata when we describe the compactification of the moduli space. }
\end{rem}

{Because we have excluded the possibility that successive Lagrangians agree and 
$p_i$ is in the diagonal, the procedure for connect sum described in the above definition does not account for all possible degenerations of moduli spaces of holomorphic disks. Additional degenerations will be explicitly described in Definition \ref{defn615616} below in terms of additional marked points along each interval  $\overline{z_{\blue{i}} z_{\blue{i+1}}}$, but we can describe the effect of the degenerations on relative homotopy classes without introducing these marked points, as follows: }
\begin{defn}\label{defn615616}
Suppose 
$(\vec{\kappa}',\vec p') \in \text{\rm Seq}_{K'}$, 
$(\vec{\kappa}'',\vec p'') \in \text{\rm Seq}_{K''}$.
We assume 
$\kappa'_i = \kappa''_{\blue{0}}$.
Then we define
\begin{equation}\label{form215}
(\vec{\kappa}'',\vec p'')\#_{i} (\vec{\kappa}',\vec p')
=(\vec{\kappa},\vec p)
\end{equation}
where 
$$
\aligned
\vec{\kappa} &= (\kappa'_{\blue{0}},\dots, \kappa'_{i}=\kappa''_0 ,\kappa''_{1},\dots,\kappa''_{K''},
\kappa''_0=\kappa'_{i},\kappa'_{i+1},\dots,\kappa'_{K'}),
\\
\vec p
&= (p'_{\blue{0}},\dots,p'_{\blue{i}},p''_{1},
\dots,p''_{K''}, \blue{p''_0}, p'_{\blue{i+1}},\dots,p'_{K'}).
\endaligned$$
\par
Then we define the map
\begin{equation}
\#_{i} 
:
\Pi_2(\vec{\kappa}'',\vec p'') 
\times \Pi_2(\vec{\kappa}',\vec p') 
\to
\Pi_2((\vec{\kappa}'',\vec p'')\#_{i} (\vec{\kappa}',\vec p'))
\end{equation}
by concatenating the maps of the left hand side by 
the band sum of 
$u'(\overline{z'_{\blue{i}} z'_{\blue{i+1}}}) \subset L_{\kappa'_i}$ 
and $u''(\overline{z''_{\blue{0}} z''_{\blue{1}}}) \subset L_{\kappa''_{\blue{0}}}$ 
in $L_{\kappa'_i} = L_{\kappa''_{\blue{0}}}$.
The equivalence class of the concatenation is 
independent of the choice of this arc.
\end{defn}
Let $\ell \in \Z_{\ge 0}$, $\vec k = (k_0,k_1,\dots,k_{K})$ where $k_i \in \Z_{\ge 0}$,
$(\vec{\kappa},\vec p) \in \text{\rm Seq}_K$ and
$B \in \Pi_2(\vec{\kappa},\vec p)$. 
\begin{figure}[h]
\centering
\includegraphics[scale=0.7]{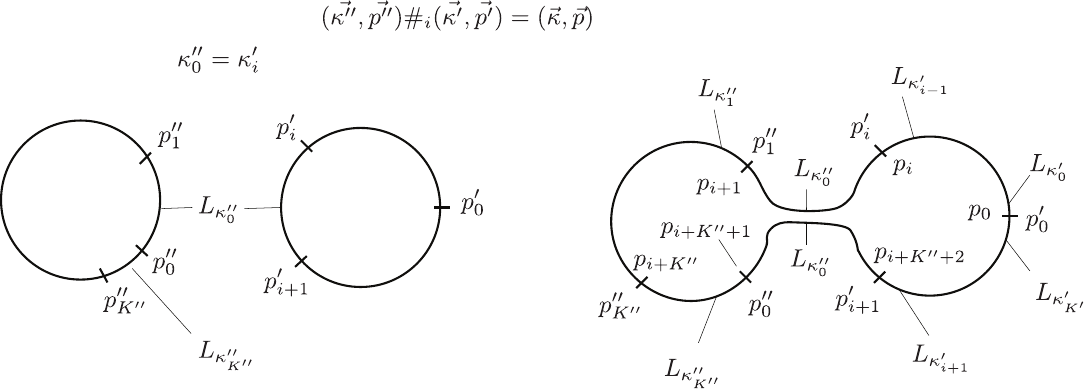}
\caption{Gluing in Definition \ref{defn615616}}
\label{Figure6-2-2}
\end{figure}

\begin{figure}[h]
\centering
\includegraphics[scale=0.25]{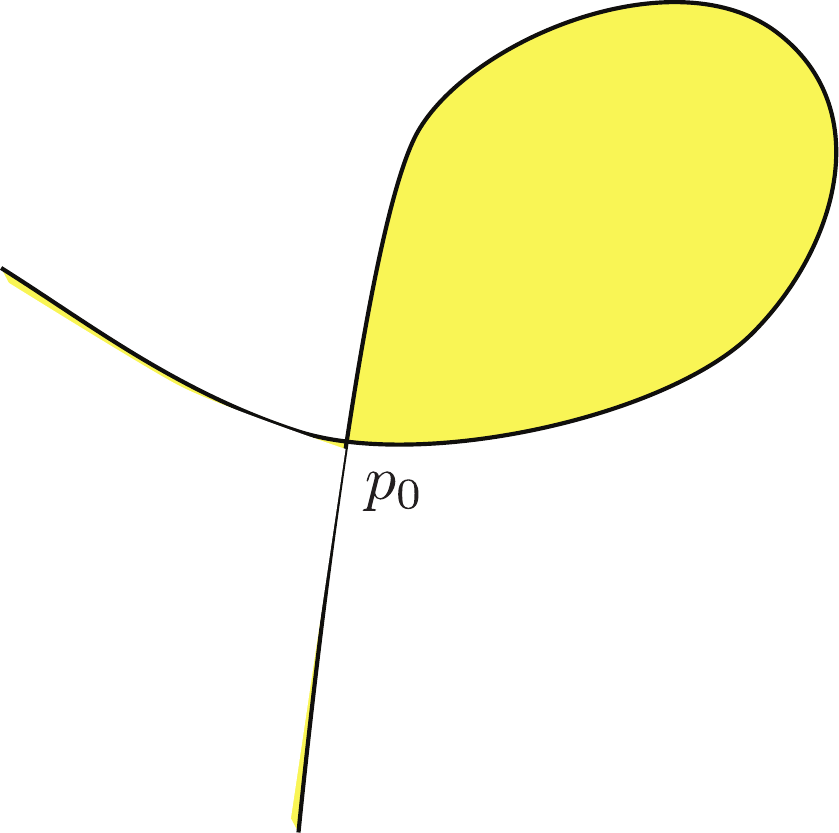}
\caption{Teardrops}
\label{eardrops}
\end{figure}

\begin{defn}\label{defMellveck}
We define $\overset{\circ}{\mathcal M}_{\ell;\vec k}((\vec{\kappa},\vec p);B)$\index[syindex]{MkellLLaE@${\mathcal M}_{\ell;\vec k}((\vec{\kappa},\vec p);B)$}
to be the set of all $\sim$ equivalence classes of the quartet $(u;\vec z^+,\vec z;\vec{\vec w})$ with the following properties:
\begin{enumerate}
\item $(u,\vec z)$ satisfies (1),(2),(3),(4) of Definition $\ref{def211}$;
its homotopy class is in the equivalence class $B$. 
\item $u : D^2 \to X$ is $J$-holomorphic.
\item $\vec{\vec w} = (\vec w_{\blue{0}},\dots,\vec w_K)$ where
$\vec w_i  =(w_{i,1},\dots,w_{i,k_i})$  $w_{i,j} \in \overline{z_{\blue{i}} z_{\blue{i+1}}}$.
They are consistent with  the counter-clockwise orientation of $\partial D^2$.
\item
$\vec z^+ = (z^+_1,\dots,z^+_{\ell})$ is a collection of  mutually disjoint marked points 
$z^+_i \in \text{\rm Int}D^2$.
\end{enumerate}
\par
We say that $(u;\vec z^+,\vec z;\vec{\vec w}) 
\sim (u';\vec z^{\prime +},\vec z\,';\vec{\vec w}')$
if there exists a biholomorphic map $v : D^2 \to D^2$ such that
$u' \circ v = u$, $v(w_{i,j}) = w'_{i,j}$, 
$v(z^+_i) = z^{\prime +}_i$, $v(z_i) = z^{\prime}_i$.
\end{defn}
\begin{defn}
We call $z_i$ a {\it switching boundary marked point}\index{switching boundary marked point} and 
$w_{i,j}$ a {\it diagonal boundary marked point}\index{diagonal boundary marked point}.
\end{defn}
Note in case $K=\blue{0}$, i.e., only one Lagrangian boundary, this moduli space coincides with the 
moduli space in Definition \ref{diskmoduli1} if $\vec p = \emptyset$.
 %\marginpar{Rewritten so that  teardrop is included. KF 2025 Aug.}
If $K= 0$ and $\vec p = p_0$ then it is the moduli space of 
teardrops.\index{teardrop} (Figure \ref{eardrops}.)
 Hereafter in this subsection we assume $\vec p \ne \emptyset$.
(This is because of the way how we organize the induction to 
define Kuranishi structures and \blue{CF-perturbations}.  
Namely they are fixed already in the case $K=0$, $\vec p = \emptyset$ in Section \ref{ainfalgasssingle}.)
\par
By evaluating $u$ at $z_i^{+}$ and $w_{i,j}$ we obtain 
an evaluation map
\begin{equation}\label{evevev}
\text{\rm ev} 
= (\text{\rm ev}^+, (\text{\rm ev}_i)_{i=\blue{0}}^K)
:
\overset{\circ}{\mathcal M}_{\ell;\vec k}((\vec{\kappa},\vec p);B)
\to
X^{\ell} \times \prod_{i=\blue{0}}^K \tilde L_{\kappa_i}^{k_i},
\end{equation}
where $\text{\rm ev}^+ = (\text{\rm ev}^+_1,\dots,\text{\rm ev}^+_{\ell})$,
$\text{\rm ev}_i = (\text{\rm ev}_{i,1},\dots,\text{\rm ev}_{i,k_i})$.

\begin{defn}\label{defn917}
Let $(\vec{\kappa},\vec p) \in  \text{\rm Seq}_{K}$ and $\vec k = (k_0,\dots,k_K) \in \Z_{\ge 0}^K$.
We consider $m \in \{0,\dots,k_0\}$.  %\marginpar{Definition is added.  KF 2025 Sep.}\index[syindex]{MkellLLaEm@${\mathcal M}_{\ell;\vec k}((\vec{\kappa},\vec p,m);B)$}
We put
$$
\overset{\circ}{\mathcal M}_{\ell;\vec k}((\vec{\kappa},\vec p,m);B)
= \overset{\circ}{\mathcal M}_{\ell;\vec k}((\vec{\kappa},\vec p);B).
$$
For an element $(u;\vec z^+,\vec z;\vec{\vec w}) \in \overset{\circ}{\mathcal M}_{\ell;\vec k}((\vec{\kappa},\vec p,m);B)$
we define its $0$-th marked point $z(0)$\index{0-th marked point}  as follows.\index[syindex]{zz(0)@ $z(0)$}\index[syindex]{zz(i)@ $z(i)$}
\begin{enumerate}
\item If $m=0$ then $z(0) = z_{0}$.
\item If $m>0$ then $z(0) = w_{0,m}$.
\end{enumerate}
Then we enumerate $\vec z \cup \vec{\vec w}$ by the counter clockwise ordering
and define $z(0),z(1),\dots,z(k+1)$, where $k = K + \sum k_i$.
\par
Note that an element of $\overset{\circ}{\mathcal M}_{\ell;\vec k}((\vec{\kappa},\vec p);B)$ 
does not contain  data which specify $0$-th marked point.  By adding $m \in \{0,\dots,k_0\}$
we specify it as above.
\end{defn}

We describe the boundary of the stable map compactification 
${\mathcal M}_{\ell;\vec k}((\vec{\kappa},\vec p);B)$ of 
$\overset{\circ}{\mathcal M}_{\ell;\vec k}((\vec{\kappa},\vec p);B)$ below.
There are three kinds of boundaries.
\par
\begin{figure}[h]
\includegraphics[scale=0.7]{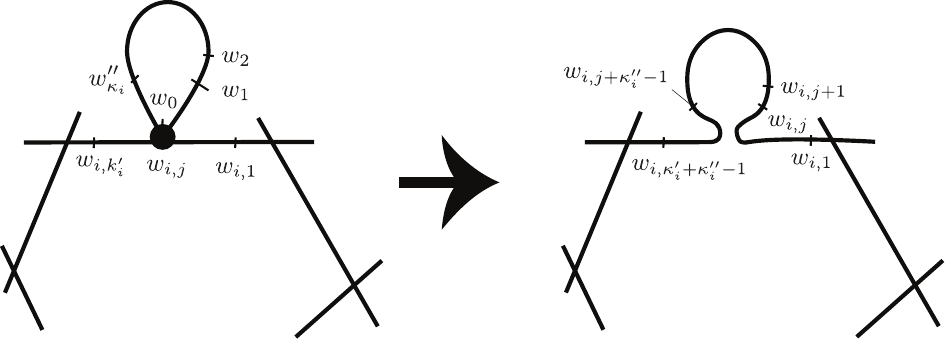}
\caption{Boundary of Type I}
\label{Figure6-2-6-2}
\end{figure}

\medskip
\noindent{\bf (Boundary of Type I)}:
Let $(\vec{\kappa},\vec p) \in \text{\rm Seq}_K$,
$B' \in \Pi_2(\vec{\blue{\kappa}} ,\vec{\blue{p}})$ and $\beta \in \pi_2(X;L_{\kappa_i})$.
We assume $B' \# \beta =B$, where $\#$ is the obvious  
concatenation.  
Let $\vec{k}=(k_0, \dots, k_K)$.  
Let $i$ be an integer such that $0 \leq i \leq K$, $k'_i$, $k''_i$ integers such that $k'_i + k''_i  = k_i + 1$ 
and $j \in \underline k'_i$.  
Write $\vec k' = (k_0,\dots,k'_i,\dots,k_K)$.  
We consider the fiber product
\begin{equation}\label{218}
{\mathcal M}_{\ell'';k''_i+1}(L_{\kappa_i};\beta)
\,{}_{\text{\rm ev}_0} \times_{\text{\rm ev}_{i,j}}
{\mathcal M}_{\ell';\vec k'}((\vec{\kappa},\vec p);B')
\end{equation}
Here the fiber product is taken over $\tilde L_{\kappa_i}$.
See Figure \ref{Figure6-2-6-2}.  %\marginpar{Figure is redrawn. One of the figures are removed KF}
\par
\begin{figure}[h]
\centering
\includegraphics[scale=0.7]{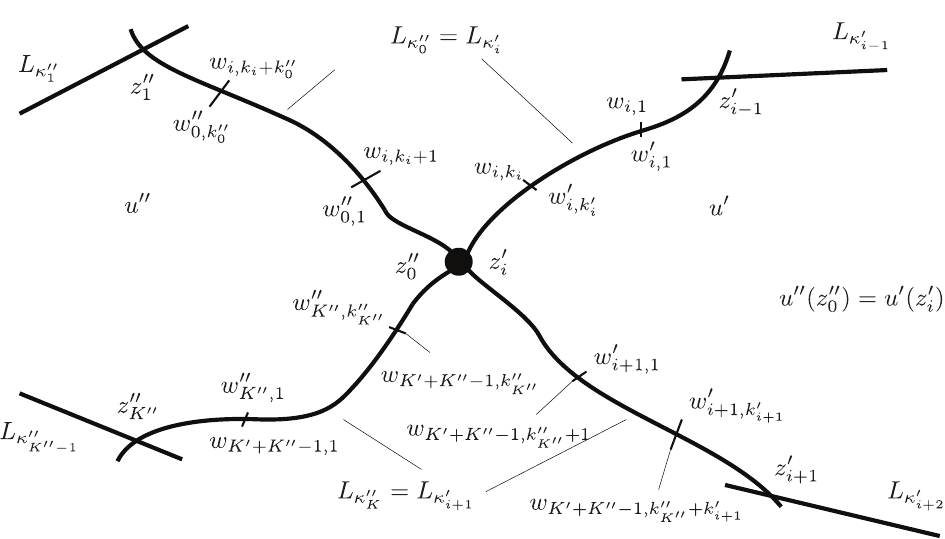}
\caption{Boundary of Type II}
\label{Figure6-2-7-2}
\end{figure}

\medskip
\noindent{\bf (Boundary of Type II)}:
We assume (\ref{213formula}).
Let $B' \in \Pi_2(\vec{\kappa}',\vec p')$, 
$B'' \in \Pi_2(\vec{\kappa}'',\vec p'')$ with
$B'' \#_{i,i+1} B' = B$.
Let $\vec k' = (k'_0,\dots,k'_{K'})$,
$\vec k'' = (k''_0,\dots,k''_{K''})$.
We define $\vec k = \vec k'' \#_{i,i+1} \vec k'$ by
$$
k_j = 
\begin{cases}
k'_j   & j=\blue{0},\dots,i-1, \\
k'_i  + k''_{\blue{0}}, &j=i, \\
k''_{\blue{j-i}}  & j=i+1,\dots, \blue{i+K'' -1}, \\
k'_{i+1} + k''_{K''}  &j=\blue{i+K''}, \\
k'_{\blue{j-K''+1}} & j=i+K'' +1, \dots, \blue{K'+K''-1}.
\end{cases}
$$
(See Figure \ref{Figure6-2-7-2}.) 
We then consider the direct product:
\begin{equation}\label{219}
{\mathcal M}_{\ell'';\vec k''}((\vec{\kappa}'',\vec p'');B'')
\times
{\mathcal M}_{\ell';\vec k'}((\vec{\kappa}',\vec p');B').
\end{equation} 
 %\marginpar{{ Figure are rewritten. KF}}

\par
\begin{figure}[h] 
\centering
\includegraphics[scale=0.7]{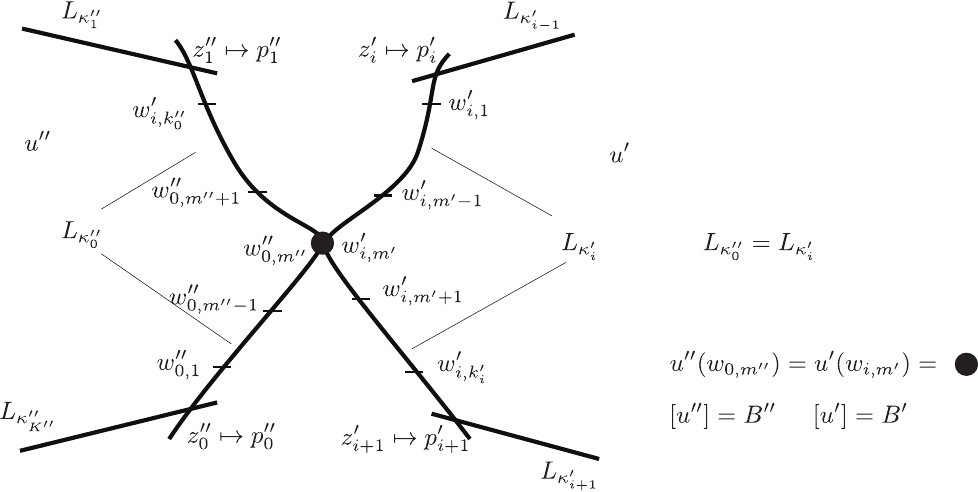}
\caption{Boundary of Type III}
\label{Figure6-2-8-2}
\end{figure}

\medskip
\noindent{\bf (Boundary of Type III)}:
We assume (\ref{form215}).
Let $B' \in \Pi_2(\vec{\kappa}',\vec p')$, 
$B'' \in \Pi_2(\vec{\kappa}'',\vec p'')$ with
$B'' \#_i B' = B$.
Let $\vec k' = (k'_{\blue{0}},\dots,k'_{K'})$,
 $\vec k'' = (k''_{\blue{0}},\dots,k''_{K''})$.
 Let $m' \in \{ 1,\dots, k'_i\}$,
$m'' \in \{1,\dots,k''_{0}\}$.
We define $\vec k = \vec k'' \#_{i,(m',m'')} \vec k'$ by
$$
k_j = 
\begin{cases}
k'_j   & j=\blue{0},\dots, i-1, \\
m' +k''_{\blue{0}} - m'' -\blue{1} &j=i, \\
k''_{\blue{j-i}}  & j=i+1,\dots, \blue{i+K''}, \\
m'' + k'_{i} - m' -\blue{1} & j= \blue{i+K'' +1}, \\
k'_{\blue{j-K''-1}} & j=i+K''+\blue{2},\dots, K'+K''\blue{+1}.
\end{cases}
$$
(See Figure  \ref{Figure6-2-8-2}.)
We now consider the fiber product
\begin{equation}\label{220}
{\mathcal M}_{\ell'';\vec k''}((\vec{\kappa}'',\vec p'');B'')
\,{}_{\text{\rm ev}_{0,m''}}\times_{\text{\rm ev}_{i,m'}}
{\mathcal M}_{\ell';\vec k'}((\vec{\kappa}',\vec p');B'),
\end{equation}
where the fiber product is over $\tilde L_{\kappa_i}
= \tilde L_{\kappa'_i}=\tilde L_{\kappa''_{\blue{0}}}$.
 %\marginpar{Figure redrawn. KF}

An element of the compactification ${\mathcal M}_{\ell;\vec k}((\vec{\kappa},\vec p);B)$ 
can be regarded also as an isomorphism classes of $(\Sigma;u;\vec z^+,\vec z;\vec{\vec w})$ 
where $\Sigma$ is a bordered genus zero curve with one boundary component, 
$u : \Sigma \to X$ is a pseudo-holomorphic map with appropriate boundary 
conditions, $\vec z^+$ are interior marked points, $\vec z;\vec{\vec w}$ are boundary 
marked points.
The outer collaring ${\mathcal M}_{\ell;\vec k}((\vec{\kappa},\vec p);B)^{\boxplus 1}$ 
is by definition the set of pairs $([(\Sigma;u;\vec z^+,\vec z;\vec{\vec w}) ],\vec s)$
where $[(\Sigma;u;\vec z^+,\vec z;\vec{\vec w})] \in {\mathcal M}_{\ell;\vec k}((\vec{\kappa},\vec p);B)$ 
and $\vec s$ assigns a number $\in [-1,0]$ to each boundary node of $\Sigma$. (See Definition \ref{outcolar}.)
(\ref{evevev}) is extended to the outer collaring in an obvious way.\footnote{We do not 
repeat this kinds of remark.}
\par
The above description of the boundary induces one for outer collaring in an obvious way.
Now we are ready to state the main result of this subsection.
\begin{prop}\label{Kuraeistspoly}
$  $ \par
\begin{enumerate}
\item
Each of the spaces $\overset{\circ}{\mathcal M}_{\ell;\vec k}((\vec{\kappa},\vec p);B)$
has a compactification 
${\mathcal M}_{\ell;\vec k}((\vec{\kappa},\vec p);B)$ that is Hausdorff
and metrizable.
\item
The space ${\mathcal M}_{\ell;\vec k}((\vec{\kappa},\vec p);B)^{\boxplus 1}$  has a
Kuranishi structure.
 %\marginpar{compatibility with forgetful map 
%is removed from the statement.}
\item 
The boundary of ${\mathcal M}_{\ell;\vec k}((\vec{\kappa},\vec p);B)^{\boxplus 1}$
is a union of the three types of fiber or direct products 
$(\ref{218})$, $(\ref{219})$, $(\ref{220})$. 
Here the union is taken over the data described above 
together with the shuffles $(\mathbb L'',\mathbb L')$ of $\underline\ell$
such that $\# \mathbb L'' = \ell''$, $\# \mathbb L' = \ell'$.
\par
The Kuranishi structure on the first factor of $(\ref{218})$ is the one in Proposition $\ref{diskkura}$.\par
The orientations are compatible in the same sense as Proposition $\ref{diskkura}$ $(3)$.
\item
The evaluation maps $(\ref{evevev})$ extend to the 
compactification and are compatible with 
the description of the boundary in $(3)$.
\item
For each $i \in \uwave{K}$ and $j \in \{1,\dots,k_i\}$, the evaluation map
$$
\text{\rm ev}_{i,j} : 
{\mathcal M}_{\ell;\vec k}((\vec{\kappa},\vec p);B)^{\boxplus 1}
\to \tilde L_{\kappa_i}
$$
is weakly submersive.
\item
The Kuranishi structure is invariant under  permutation of 
the interior marked points.
\item
The Kuranishi structure is invariant under cyclic permutation 
of the boundary marked points in the sense we describe below. 
\end{enumerate}
\end{prop}
Let us describe the statement (7) precisely.
Let $\vec\kappa = (\kappa_0,\kappa_1,\dots,\kappa_K)$, 
$\vec p = (p_{0},\dots,p_{K})$, 
$\vec k = (k_0,k_1,\dots,k_K)$.
We put 
\begin{align} \label{eq:cyclic_rotation_kappa}
c\vec\kappa & = (\kappa_K,\kappa_0,\kappa_1,\dots,\kappa_{K-1}) \\ 
c\vec p & = (p_{K},p_{0},p_{1},\dots,p_{K-1}) \\ \label{eq:cyclic_rotation_k}
c\vec k & = (k_K,k_1,\dots,k_{K-1}).  
\end{align}
\par
Then in an obvious way we can define a  homeomorphism
$$
\overset{\circ}{\mathcal M}_{\ell;\vec k}((\vec{\kappa},\vec p);B)
\to
\overset{\circ}{\mathcal M}_{\ell;c\vec k}((c\vec{\kappa},c\vec p);B)
$$
that can be enhanced to a homeomorphism 
on their compactifications and its outer collaring. 
Statement (7) claims that it is an underlying homemorphism of an isomorphism 
of Kuranishi structures.

We define Kuranishi structures on ${\mathcal M}_{\ell;\vec k}((\vec{\kappa},\vec p,m);B)$
by the equality
\begin{equation}\label{formdef915}
{\mathcal M}_{\ell;\vec k}((\vec{\kappa},\vec p,m);B) 
=
{\mathcal M}_{\ell;\vec k}((\vec{\kappa},\vec p);B)
\end{equation}
and Proposition \ref{Kuraeistspoly}.
In other words the Kuranishi structures on ${\mathcal M}_{\ell;\vec k}((\vec{\kappa},\vec p,m);B)$
are independent of $m$.
Together with Proposition \ref{Kuraeistspoly} (7) it implies that 
the Kuranishi structures on ${\mathcal M}_{\ell;\vec k}((\vec{\kappa},\vec p,m);B)$
are invariant under the cyclic permutation of all the boundary marked points of it. 
\begin{rem}
In Proposition \ref{Kuraeistspoly} we do not 
 %\marginpar{The rest of this subsection is rewritted.}
claim  the compatibility of Kuranishi structure 
with the forgetful map of the {\it boundary} marked points.
The compatibility with such forgetful maps will be discussed 
in Section \ref{sec:unit}.
\end{rem}
\par
\begin{proof}
The construction of the Kuranishi structures that have the required properties except (5) and (7)
is the same as the argument of \cite[Section 7.2]{fooo09} \newred{or of \cite{const1,const2}}.
Proposition \ref{Kuraeistspoly} is a special case of Proposition \ref{prop92} 
which we prove in Subsection \ref{polygonforget}.
\end{proof}
\begin{prop}\label{existmultipolu}
There exists a  system of CF-perturbations of the outer collaring\footnote{
Outer collaring is defined in Definition \ref{outcolar}.}
$\{{\mathcal M}_{\ell;\vec k}((\vec{\kappa},\vec p);B)^{\boxplus 1}\}$ 
of the moduli spaces in Proposition \ref{Kuraeistspoly} with 
the following properties:\index[syindex]{oplusupper@${}^{\boxplus 1}$}
\begin{enumerate}
\item
They are transversal to $0$.
\item 
They are compatible with the description of the 
boundary in Proposition $\ref{Kuraeistspoly}$ $(3)$.
The CF-perturbation on the first factor of $(\ref{218})$
is the one of Proposition $\ref{existmkulti1}$.
\item
Let $i \in \uwave{K}$, $j \in  \{1,\dots,k_i\}$.
Then the evaluation map
$$
\text{\rm ev}_{i,j} : 
{\mathcal M}_{\ell;\vec k}((\vec{\kappa},\vec p);B)^{\boxplus 1}
\to \tilde L_{\kappa_i}
$$
is strongly submersive with respect to our  \blue{CF-perturbation}.
\item
The \blue{CF-perturbation} is invariant under the permutation of 
interior marked points.
\item
The \blue{CF-perturbation} is invariant under the cyclic permutation 
of the data in the same sense as Proposition $\ref{Kuraeistspoly}$ $(7)$. 
\end{enumerate}
\end{prop}
\begin{proof}
In the previous literature,  Proposition \ref{existmultipolu}  (1)-(5) 
are proved 
by the same 
induction as Proposition $\ref{Kuraeistspoly}$.
\newred{See \cite[Proof of Proposition 19.1]{springer}}.
We prove Proposition $\ref{existmultipolu}$ in Section \ref{sec:CRperturb}.
\end{proof}
We define CF-pereturbation of $\{{\mathcal M}_{\ell;\vec k}((\vec{\kappa},\vec p,m);B)^{\boxplus 1}\}$
by using the fact that it is isomorphic to 
$\{{\mathcal M}_{\ell;\vec k}((\vec{\kappa},\vec p);B)^{\boxplus 1}\}$
and Proposition \ref{existmultipolu}.

\begin{defn}
We call the Kuranishi structure of Proposition \ref{Kuraeistspoly},  
the ${\frak q}$-Kuranishi structure\index{qKuranishistructure@${\frak q}$-Kuranishi structure} and 
CF-perturbation in Proposition \ref{existmultipolu}
the $\frak q$-CF-perturbation\index{qCFperturabtion@$\frak q$-CF-perturbation}.
We write ${\mathcal M}_{\ell;\vec k}((\vec{\kappa},\vec p);B)^{\frak q}$ etc.
when we need to distinguish it from $\frak p$-Kuranishi structure introduced in Section \ref{sec:frakq}.
\end{defn}

\begin{rem}\label{rem819}
We explain the reason why we do {\it not} claim our Kuranishi structures and CF-perturbations 
are  compatible with the forgetful map of the {\it interior} marked points.\footnote{See also 
\cite[Remark 3.1]{fukaya:cyc}.}
(The reason is different from the reason why CF-perturbations 
are not compatible with the forgetful map of the {\it boundary} marked points.
The latter is explained in Section \ref{sec:two perturbation}.) %\marginpar{Remark added}
\par
We consider an element ${\bf p} = (\Sigma;u;\vec z^+,\vec z;\vec{\vec w})$ of ${\mathcal M}_{\ell;\vec k}((\vec{\kappa},\vec p);B)$.
To define its Kuranishi neighborhood we consider 
$u'$ close to $u$ and study equation
\begin{equation}
\overline{\partial} u' \in  \mathcal E_{\bf p}({\bf x}).
\end{equation}
Here ${\bf x}$ is $u'$ together with $\Sigma', \vec z^{\prime +},\vec z^{\prime}, \vec{\vec w}^{\prime}$.
The linear space $\mathcal E_{\bf p}({\bf x})$ is an obstruction space 
that is a finite dimensional subspace of 
$\Gamma(\Sigma';(u')^*TX \otimes \Lambda^{01})$.
We require that the support of elements  $\mathcal E_{\bf p}({\bf x})$ are 
away from interior nodes.  
The reason why we do so is that, then, we can use the trivialization of the
 %\marginpar{Explanation of 
%$\delta$-thick part is added in the appendix since it is removed from Kuranishi structure section.  KF 2025 Aug.} 
universal family of Deligne-Mumford moduli space to identify $\delta$-thick parts $\Sigma'(\delta)$ and $\Sigma(\delta)$
of $\Sigma'$ and of $\Sigma$.\footnote{Here $\delta$-thick part\index{thick part} $\Sigma(\delta)$ \index[syindex]{Sigma(delta)@$\Sigma(\delta)$} is by definition a connected component of
$
\{x \in \Sigma \mid i_{\Sigma}(x) > \delta \}
$.
(See \cite{Hu}.)
The notation $i_{\Sigma}$ in the above formula is defined as follows.
We take the double $\widehat{\Sigma}$ of $\Sigma$. We remove 
 boundary marked points and the doubles of interior marked points from $\widehat{\Sigma}$ to obtain $\widehat{\Sigma}_0$.
When $\Sigma$ together with marked points are 
stable the irreducible components of $\widehat{\Sigma}_0$ have  unique metrics of constant negative 
curvature $-1$.  The notation $i_{\Sigma}$ stands for its injectivity radious.
We say a complement of ${\Sigma}(\delta)$ a $\delta$-thin part. \index{thin part}
In case $\Sigma$ together with marked points is unstable we add additional marked points 
using $u$.
}
To require this property we also require that the support of an element of 
$\mathcal E_{\bf p}({\bf x})$ is away from interior marked points.
In fact when we consider the case where two interior marked points $z_i^+$ and $z_j^+$
become the same point in the limit, there will be a new sphere bubble appearing at that point.
If the support of an element of 
$\mathcal E_{\bf p}({\bf x})$ contains $z_i^+$ and $z_j^+$ then the support of an element of
the obstruction space (= fiber of the obstruction bundle) of the limit necessarily contains the interior nodes.
\par
Now it is easy to observe that there is no reasonable way 
to require simultaneously that the support of an element of 
$\mathcal E_{\bf p}({\bf x})$ is away from interior marked points 
and that it is compatible with the forgetful map of the interior marked points.
\par
Note that we always require that the  the support of an element of 
$\mathcal E_{\bf p}({\bf x})$ is away from boundary.
So in particular it is away from the boundary nodes or boundary marked points.
So the above mentioned reason does not apply to the 
forgetful map of boundary marked points.
In fact in Proposition \ref{diskkura} we did claim that the 
Kuranishi structure is compatible with the forgetful map of the 
boundary marked points.
 \end{rem}

\subsection{Construction of the cyclic $A_{\infty}$ category 2 : Orientation}
\label{oricyclic}

In this subsection, we describe the data
which determine orientation of 
the moduli space $ {\mathcal M}_{\ell;\vec k}((\vec{\kappa},\vec p);B)$ introduced in the last subsection.
We also describe the compatibility of the orientations 
with the identification mentioned in Proposition $\ref{Kuraeistspoly}$ (3) of the boundary.
The discussion here follows
\cite{anchor}.
\par
We recall that we fix a background class $\text{\rm st} \in H^2(X;\Z_2)$ and assume our 
Lagrangian submanifolds are $\text{\rm st}$-relatively spin.
More precisely we need to fix a relatively spin structure.
Let us review its definition now.
We fix a smooth triangulation of $X$ which includes each Lagrangian submanifold 
we study as a simplicial subcomplex.
We take an oriented 
real vector bundle $V(\text{\rm st})$ on the  3-skeleton $X^{(3)}$ of $X$ whose 
second Stiefel-Whiteney class is $\text{\rm st}$.\index[syindex]{Vst@ $V(\text{\rm st})$}
\begin{defn} \label{def:relative_Spin}
An $\text{\rm st}$-{\it relative spin structure}\index{relative spin structure} on a Lagrangian submanifold 
$L$ is a trivialization of $TL \oplus V(\text{\rm st})$ on the  2-skeleton $L^{(2)}$ 
of $L$.
\end{defn}
If $L$ is $\text{\rm st}$-relative spin,  allowing the trivialization of $TL \oplus V(\text{\rm st})$ on the  2-skeleton $L^{(2)}$ to vary,  the set of all 
$\text{\rm st}$-relative structures forms a principal homogeneous space of 
the group $H^2(X,L;\Z_2)$ and this set is independent of the choice of 
triangulation. (See \cite[Proposition 8.1.6]{fooo09}.) 
\par
In this paper, we choose a real vector bundle $V(\text{\rm st})$ and fix it.   The set of spin structures on  $TL \oplus V(\text{\rm st})$ 
is a principal homogeneous space of $H^1(L;\Z_2)$.   The notion of stable conjugacy among them is defined in \cite[Definition 8.1.5]{fooo092} with 
$V_1=V_2=V(\text{\rm st})$.  
Then the set of stable conjugacy classes of spin structures 
on $TL \oplus V(\text{\rm st})$ is a principal homogeneous space of 
$H^1(L;\Z_2)/ (\text{\rm Im}(H^1(X;\Z_2) \to H^1(L;\Z_2)))$, which is a subgroup of $H^2(X,L;\Z_2)$.  
Note also that the orientation on the moduli spaces of pseudo-holomorphic disks only depends on the stable conjugacy 
class of relative spin structures, see \cite[Proposition 8.1.16]{fooo092}.
\par
In this paper we consider various finite sets $\mathscr L$ of $\text{\rm st}$-relative spin 
Lagrangian submanifolds.
We always choose and fix an $\text{\rm st}$-relative spin structure on each element $L$ of $\mathscr L$.
%Actually taking a different choice of $\text{\rm st}$-relative spin structure on $L$
%is equivalent to taking an appropriate homomorphism $\pi_1(L) \to \Z_2$.
%In other words it corresponds to an element $\Delta \in H^1(L;\Z_2)$.  
%\blue{If $\Delta$ has a lift in $H^1(L;\Z)$, 
%it is easy to see that this is canceled by changing $b \in \mathcal M_{\text{\rm weak}}(L;\Lambda_0)$
%to $b + \pi\sqrt{-1}  \Delta$. Thus the choice of $\text{\rm st}$-relative structure 
%does not change the filtered $A_{\infty}$ category $\cL$ etc.  
%In general, a different choice of st-relative structure corresponds to a twisting by a local system 
%associated with the homomorphism $\pi_1(L) \to \Z_2$.}
\par
We next take $(L,\theta), (L',\theta') \in \mathscr L$, which are transversal to each other.
We take a path $\lambda_p : [0,1] \to {\mathcal{LAG}}_+(T_{p}(X))$ 
for each $p\in L \cap L'$ such that
\begin{equation} \label{eq:path_Lagrangian}
\lambda_p(0) = T_{p}L, \quad
\lambda_p(1) = T_{p}L'.
\end{equation}
Here ${\mathcal{LAG}}_+(T_{p}(X))$ is the oriented Lagrangian Grassmannian, that is the set of oriented $n$ dimensional subspaces
of $T_pX$ on which the symplectic form $\omega$ vanishes.
\par
We make a similar choice for self-intersection points of immersed Lagrangians $L$.
\begin{rem}
{In case the symplectic manifold $X$ satisfies $c_1(X)=0$, and we have chosen a complex volume form on $X$, we may restrict attention to Lagrangians equipped with a trivialisation of the phase map to $S^1$ (i.e. graded Lagrangians). Assuming that all $(L,\theta_L) \in \mathscr L$ are graded Lagrangians, this data determines a unique homotopy class of paths $\lambda_p$, associated to each transverse intersection $p$.} 
\end{rem}
\begin{defn}
A relative spin structure determines trivializations of 
$$
T_pL \oplus V(\text{\rm st})_p, \quad T_pL' \oplus V(\text{\rm st})_p.
$$
We take a trivialization of
$\lambda_p(t) \oplus V(\text{\rm st})_p$ which depends smoothly on $t$ and 
coincides with  the given ones at $t=0,1$.
Pick and fix oriented frames of $T_pL$ and $T_pL'$.  
There are two homotopy classes of trivializations
relative to the frames at $t=0,1$.
We write $\tilde{\lambda_p}$
in place of $\lambda_p$, when we include the datum of a choice of
a homotopy class of such trivializations. 
%We write $\tilde \lambda_p$ in place of $\lambda_p$ when 
%we include the datum of homotopy class of this trivialization.
\par
Consider $(-\infty,0] \times [0,1] \subset \C$.
Let $L^{1,q}(T_pX;L,L';\lambda_p)$ be the set of 
maps $\xi : (-\infty,0] \times [0,1] \to T_pX$ of $L^{1,q}$ class 
(that is the set of maps whose first derivatives are of $L^q$ class, $q > 2$) satisfying the boundary conditions 
\begin{equation}\label{form815}
\xi(\tau,0) \in T_pL, \quad \xi(\tau,1) \in T_pL',
\quad
\xi(0,t) \in \lambda_p(t).
\end{equation}
 Let $L^{q}(T_pX\otimes \Lambda^{0,1})$ be the set of sections $(-\infty,0] \times [0,1] \to T_pX\otimes \Lambda^{0,1}$
of $L^q$ class of the bundle $T_p X \otimes \Lambda^{0,1}\C$ on $(-\infty,0] \times [0,1]$.
(Actually the bundle $\Lambda^{0,1}\C$ is trivial.)
The standard Cauchy-Riemann operator defines a Fredholm map
\begin{equation}\label{orientationcomplex}
\overline{\partial}_p : L^{1,q}(T_pX;L,L';\lambda_p) \to L^{q}(T_pX
\otimes \Lambda^{0,1}).
\end{equation}
Let $o_p$ be an orientation of the determinant bundle of the index of 
(\ref{orientationcomplex}).
\end{defn}
\begin{thm}\label{orithem}
The choice of $V(\text{\rm st})$-relative spin structures of each of $L_{\kappa}$ and 
$\tilde \lambda_p$, $o_p$ for each of $p \in L_{\kappa} \cap L_{\kappa}'$ determines 
a system of orientations on ${\mathcal M}_{\ell;\vec k}((\vec{\kappa},\vec p);B)$
that is compatible with Proposition $\ref{Kuraeistspoly}$ $(3)$.
\end{thm}
This is \cite[Theorem 6.7]{anchor}.
We can prove a similar statement on the orientation of the other moduli spaces 
which appear in this paper, in the same way.
\par
We use this orientation throughout the paper.
We remark that if we change $o_p$ to the opposite orientation,
then the orientation of all  the moduli spaces which contain $p$ are affected.
This is actually equivalent to replacing the generator $[p] \in CF(L,L')$ by 
$-[p] \in CF(L,L')$. Therefore, the filtered $A_{\infty}$ category does not change 
up to isomorphism. 
The same remark applies to the choice of $\tilde \lambda_p$.
So we do not write those choices  $\tilde \lambda_p$, $o_p$ in our 
notation for filtered $A_{\infty}$ category.

\blue{In this paper, we construct filtered $A_{\infty}$-structures on the de Rham complexes with $\Lambda_0$ coefficients as in \cite{fooo:toric1}, 
\cite{fooo:bulk}, \cite{toric3}, \newred{\cite{tech22}} while 
we did so on certain subcomplexes of singular chain complexes in \cite{fooo09}, \cite{anchor}.  
The signs in operations $\frak m$, $\frak q$, $\frak p$ arising from the correspondence are obtained in \cite{fooo09} 
as we explained in the last paragraph of \cite{fooo:bulk}.  
}

\subsection{Construction of the cyclic $A_{\infty}$ category 3 : The operations}
\label{constcyclic}

We now use the moduli space and its \newred{CF-perturbations} 
constructed in Subsection \ref{cycliccategorymoduli} 
to define the operators of the cyclic filtered $A_{\infty}$ category {$\cL^{\rm form}_{\rm curve}$}.
(More precisely its $A_{n,k}$ version.)
\par
\begin{figure}[h] 
\centering
\includegraphics[scale=0.35]{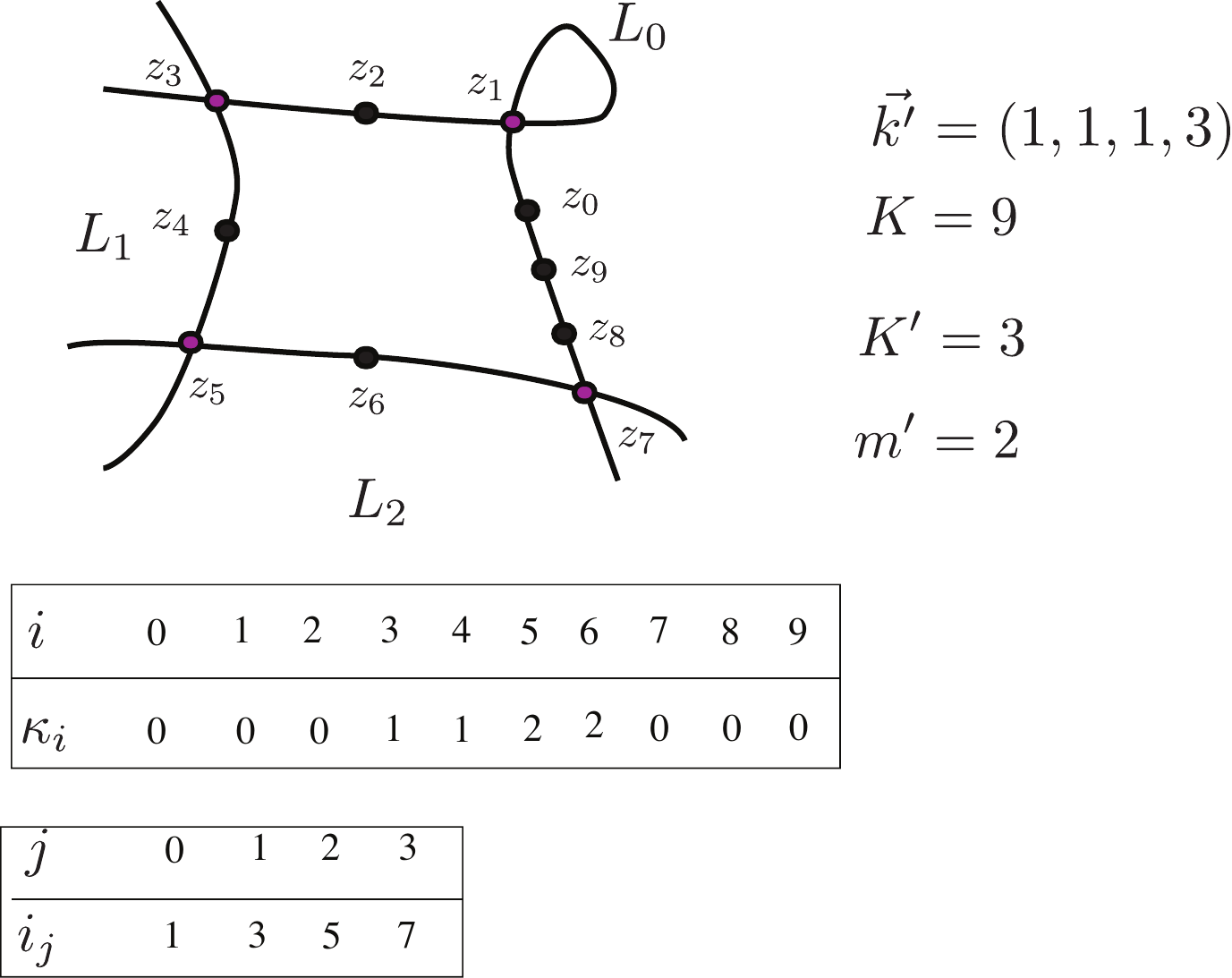}
\caption{Forgetful map Red}
\label{Figure(8.18)}
\end{figure}

For this purpose we change our notation a bit.
\begin{defn}\label{newdef926}
Let $\vec{\kappa} = (\kappa_0,\kappa_1,\dots,\kappa_K)$, $\kappa_0,\dots,\kappa_K \in \{1,\dots,\#\mathscr L\}$,
%This time we do {\it not} assume ${\kappa_i} \ne {\kappa_{i+1}}$.  
%{We also choose $\theta_{L_{\kappa}} \in \Omega^1(L_{\kappa})$ satisfying  
%(\ref{fixthetaL}) for each $L_{\kappa}$ appearing in this set of Lagrangians.
and let $\vec{p} = (p_{i_0},\dots,p_{i_{K'}})$
 be a sequence
 such that:
 \begin{enumerate}
 \item  $0\le i_0  <\dots < i_{K'} \le K$.  $K' \ge -1$.  In case $K' = -1$, $\vec{p} =\emptyset$.
 \item If $i \notin \{ i_0, \dots, i_{K'}\}$ then $L_{\kappa_{i-1}} = L_{\kappa_{i}}$.
 \item If $L_{\kappa_{i_j-1}} \ne L_{\kappa_{i_j}}$ then $p_{i_j} \in L_{\kappa_{i_j-1}} \cap L_{\kappa_{{i_j}}}$.
 \item  If  $L_{\kappa_{i_j-1}} = L_{\kappa_{i_j}}$ then $p_{i_j} = (p_{i_j}^1,p_{i_j}^2)$ 
 such that $p_{i_j}^1,p_{i_j}^2 \in \tilde L_{\kappa_{i_j}}$ with 
 $p_{i_j}^1  \ne p_{i_j}^2$ and $i_{L_{\kappa_i}}(p_i^1) = i_{L_{\kappa_i}}(p_i^2)$. In other words it is a self-intersection point.
 \end{enumerate}
We denote the set of such  pairs $(\vec{\kappa},\vec{p})$ by $\widetilde{\text{\rm Seq}}_K$. \index[syindex]{SeqtildeK@$\widetilde{\text{\rm Seq}}_K$}
\par
We define a forgetful map \index[syindex]{Red@${\rm Red}$}
  \begin{equation} \label{eq:reduced_sequence}
{\rm Red} :  \widetilde{\text{\rm Seq}}_K \to 
\bigcup_{K'} \left({\text{\rm Seq}}_{K'} \times   \Z^{1 + K'}_{\geq 0} \times \Z_{\geq 0}\right).
  \end{equation}
by $\text{\rm Red}(\vec{\kappa},\vec{p}) = (\vec{\kappa}', \vec{p}, \vec k', m')$, such that    
\begin{itemize}
    \item $\vec{\kappa}' = (\kappa_{i_0},\dots,\kappa_{i_{K'}})$.
    \item $k'_j = i_{j+1} - i_{j} -1  \in \Z_{\ge 0}$. (Here $i_{-1} = i_{K'}$ as convention.)
      \item $m' = K-i_{K'}$. 
    \end{itemize}
  If $K' = -1$ then ${\rm Red}(\vec{\kappa},\vec{p}) = ((\kappa,\emptyset),K,m')$.
 %When $K = K'$, we call the pair $(\vec{\kappa}, \vec{p})$ a reduced pair.  
 \end{defn}
 Here $\vec k' = (k'_0,\dots,k'_{K'})$ determine the numbers $k'_j  \in \Z_{\ge 0}$ of the 
 diagonal marked point on $L_{\kappa_{i_j}}$, between two switching marked points (that is, $z_{i_{j}}$ and $z_{i_{j+1}}$).  The number $m'$ determines which marked point on $L_{\kappa_0}$
 is the $0$-th marked point.  If $m' = 0$ then the $0$-th marked point is a switching marked point.
 If $m' > 0$ then  the $0$-th marked point is a diagonal marked point.
 See Figure \ref{Figure(8.18)} and compare (\ref{eq:change_notation_moduli_space}). 

\par
We also put\index[syindex]{BCFformcurve@$
BCF(\cL^{\rm form};\vec{\kappa})$}
$$
BCF(\cL^{\rm form};\vec{\kappa})
= \bigotimes_{i=1}^K CF(L_{\kappa_{i-1}},L_{\kappa_{i}};\F) \blue{[1]}
$$
where
 %\marginpar{rewritten so that it include immersed case. KF 2025 Aug.}\index[syindex]{CF(F@$
%CF(L_{\kappa_{i-1}},L_{\kappa_{i}};\F)$}
\begin{equation}\label{eq920}
\aligned
&CF(L_{\kappa_{i-1}},L_{\kappa_{i}};\F)\\
&= 
\begin{cases}
\displaystyle
\bigoplus_{p \in L_{\kappa_{i-1}} \cap L_{\kappa_{i}}} \F[p]
& \text{if $L_{\kappa_{i-1}} \ne L_{\kappa_{i}}$}, \\
\displaystyle(\Omega(\tilde L_{\kappa_{i}}) \otimes \F) \oplus \bigoplus_{p \in (\tilde L_{\kappa_i} \times_{X} \tilde L_{\kappa_i}) \setminus \tilde L_{\kappa_i}} \F[p]
& \text{ if $L_{\kappa_{i-1}} = L_{\kappa_{i}}$}.
\end{cases}
\endaligned
\end{equation}
We recall that the immersed Lagrangian $L_{\kappa_i}$ is identified with $(\tilde L_{\kappa_i},i_{L_{\kappa_i}})$, 
where $\tilde L_{\kappa_i}$ is an $n = \frac{1}{2}\dim_{\R} X$ dimensional manifold and $i_{L_{\kappa_i}}: \tilde L_{\kappa_i} \to X$ 
is a Lagrangian immersion.
Note the set $ (\tilde L_{\kappa_i} \times_{X} \tilde L_{\kappa_i}) \setminus \tilde L_{\kappa_i}$ is a finite set 
consisting of self-intersections of $L_{\kappa_i}$ (each self-intersection corresponds two points in this set.)
In the case of  self-intersection point we also use Maslov-Viterbo index to define its degree.
\begin{rem}\label{rem927}
We can also write 
\begin{equation}\label{eq921}
CF(L_{\kappa_{i-1}},L_{\kappa_{i}};\F)
= \Omega(\tilde L_{\kappa_{i-1}} \times_{X} \tilde L_{\kappa_i}) \otimes \F.
\end{equation}
We use both (\ref{eq920}) and (\ref{eq921}) according to the situation. 
\end{rem}
We define the degree of an element of 
$CF(L_{\kappa_{i-1}},L_{\kappa_{i}};\F)$ as follows. 
In case $L_{\kappa_{i-1}} \ne L_{\kappa_{i}}$
the degree of $p  \in L_{\kappa_{i-1}}  \cap L_{\kappa_{i}}$
is its Maslov-Viterbo index. \footnote{See \cite[Definition 2.2.12]{fooo09}.
In  \cite[Definition 2.2.12]{fooo09} an extra data $w$ (that is a path joining constant path to $p$ 
to a given base point of the space of arcs) is used. 
Here we write $[p]$ in place of $[p,w]$ since we consider $\Z_2$ grading only
and $\Z_2$ grading is independent of $w$.}
 %\marginpar{MA190527: Is $w$ the same as the choice of $\lambda_p$ in Equation \eqref{eq:path_Lagrangian}? 
%Maybe we can add a small remark about the CY case.
%I moved it to footnote and add short explanation. A remark is added. KF}

In case $L_{\kappa_{i-1}} = L_{\kappa_{i}}$ the degree of an
element of $CF(L_{\kappa_{i-1}},L_{\kappa_{i}};\F)$ 
is its degree as a differential form.
\par
We put 
$$
\aligned
\deg' x &= \deg x -1, \\
\deg (x_1\otimes \dots \otimes x_K)
& =\deg x_1+ \dots + \deg x_K,
\\
\deg' (x_1\otimes \dots \otimes x_K)
&=
\deg (x_1\otimes \dots \otimes x_K) - K.
\endaligned
$$
\begin{rem}
In this paper we mainly use the $\Z_2$ grading.
However in certain special cases we use the $\Z$ grading.
Suppose $c_1(TX) = 0$ and Maslov index homomorphism 
$\pi_2(X,L) \to \Z$ is trivial.  
$c_1(TX) = 0$  implies that the Lagrangian Grassmannian 
bundle $\mathcal{LAG}(X) \to X$ over $X$ (whose fiber at $p$ is the Lagrangian Grassmannian of 
the symplectic vector space $T_p(X)$)  has a $\Z$ fold cover $\widetilde{\mathcal{LAG}}(X)$
which restricts to the universal cover for each fiber of $p \in X$.
$p \mapsto T_pL$ defines a section of the pull back of $\mathcal{LAG}(X)$ to $L$.
The vanishing of Maslov index implies that there is a lift of this section to 
the pull back of $\widetilde{\mathcal{LAG}}(X)$ to $L$.
Lagrangian submanifold $L$ equipped  with a lift is called graded Lagrangian submanifold 
(\cite{Se}).  If $L_0,L_1$ are graded Lagrangian submanifolds intersecting transversally,
then for each $p \in L_0 \cap L_1$ there is a canonical choice of 
(a homotopy class of) a path from $T_p(L_0)$ to $T_p(L_1)$ in 
the Lagrangian Grassmannian of $T_p(X)$.  Using this choice of path 
we can obtain a well defined $\Z$ degree of $p$.
See \cite{Se} for detail.
\par
In the case when $c_1(TX)$ and the half of the Maslov index are divisible by a natural number $m$
we obtain a $\Z/2m$ grading in a similar way.
\end{rem}
\par
\begin{defn}\label{defn927777}
For $(\vec{\kappa},\vec{p}) \in \widetilde{\text{\rm Seq}}_K$,
 %\marginpar{This definition is modified
%so that immersed case is included.  KF 2025 Aug.}
we define $L^{\vec{\kappa},\vec p}_{\text{\rm source}}$\index[syindex]{Lwkappavecpsource@$L^{\vec{\kappa},\vec p}_{\text{\rm source}}$}
to be the open closed subset of 
$\prod_{i=1}^K (\tilde L_{\kappa_{i-1}} \times_{X} \tilde L_{\kappa_i})
$ consisting $(x_i,y_i)_{i=1}^{K}$ such that
\begin{enumerate}
\item If $L_{\kappa_{i-1}} \ne L_{\kappa_i}$ then $x_i = y_i \in 
L_{\kappa_{i-1}}\cap L_{\kappa_i}$.
\item If $L_{\kappa_{i-1}} = L_{\kappa_i}$ but $i=i_j$ (See Definition \ref{newdef926} (1)) 
then $(x_i,y_i) = p_{i_j}$.  In other words $(x_i,y_i)$ is a self-intersection point of $L_{\kappa_i}$
\item Otherwise $x_i = y_i \in \tilde L_{\kappa_i}$. In other words $(x_i,y_i)$ corresponds to a diagonal marked point.
\end{enumerate}
By Remark \ref{rem927}, we may regard 
$CF(\cL^{\rm form};\vec{\kappa})$ as the 
set of $\F$-valued differential forms on $\prod_{i=1}^K (\tilde L_{\kappa_{i-1}} \times_{X} \tilde L_{\kappa_i})
$ and so it determines an  $\F$-valued differential forms on $L^{\vec{\kappa},\vec p}_{\text{\rm source}}$.
\end{defn}
We put
$$
\Pi_2(\vec{\kappa},\vec{p}) = \Pi_2(\pi_{\vec{\kappa}, \vec{p}} \circ \text{\rm Red}(\vec{\kappa},\vec{p})),
$$
where $\pi_{\vec{\kappa}, \vec{p}}$ is the projection to $\text{\rm Seq}_{K'}$ and the right hand side is defined in 
Definition \ref{def211}.
Let $(\vec{\kappa},\vec{p}) \in \widetilde{\text{\rm Seq}}_K$ and $B \in \Pi_2(\vec{\kappa},\vec p)$. We  define\index[syindex]{mellkappaB@$\mathcal M_{\ell}(\vec \kappa,\vec p;B)$}
\begin{equation}
  \label{eq:change_notation_moduli_space}
\mathcal M_{\ell}(\vec \kappa,\vec p;B)
=
\mathcal M_{\ell,\vec k'}((\vec{\kappa}',\vec p,m');B)  
\end{equation}
where
 %\marginpar{(\ref{eq:change_notation_moduli_space})
%and the discussion below is corrected.  KF 2025 Aug. 31}
$\text{\rm Red}(\vec{\kappa},\vec p) = (\vec{\kappa}', \vec{p}, \vec k', m')$.
We define the Kuranishi structure and CF-perturbation on 
$\mathcal M_{\ell}(\vec \kappa,\vec p;B)$
by the equality (\ref{eq:change_notation_moduli_space})
and Propositions \ref{Kuraeistspoly} and \ref{existmultipolu}.
\par
We will use $\mathcal M_{\ell}(\vec{\kappa},\vec p;B)$
to define a map\index[syindex]{qformellvec@$\frak q^{\text{\rm form}}_{\ell;\vec{\kappa},\vec p;B}$}
$$
\frak q^{\text{\rm form}}_{\ell;\vec{\kappa},\vec p;B}
: 
\Omega(X)\blue{[2]}^{\otimes \ell}
\otimes
BCF(\cL^{\rm form};\vec{\kappa})
\to CF(L_{\kappa_{0}},L_{\kappa_K};\F)[1],
$$
which will in turn induce the de Rham {\it complex} version 
$\frak m_K^{\text{\rm form},\frak b}$
of the operator 
$\frak m_K$.
(Namely the operators $\frak m_K$ that 
will be defined in Section \ref{canonical}, is its 
de Rham {\it cohomology} version.)  

Here we apply Koszul's rule with respect to shifted gradings to arrange 
the order of marked points $\{z_{i'}\} \cup \{w_{i',j} \}$ to the one starting with {
  $w_{0,m'}$} and respecting the counter-clockwise order.  {We write this new sequence $\{\tilde{z}_i\}_{i=0}^{K}$, and note that the boundary condition along $\overline{\tilde{z}_{i} \tilde{z}_{i+1}} $ is the Lagrangian $L_{\kappa_i}$.}\index[syindex]{evpartial@$\text{\rm ev}^{\partial}$} 

\par
We define an evaluation map
\begin{equation}\label{evmap2}
\text{\rm ev} = (\text{\rm ev}^+;\text{\rm ev}^{\partial}): %;\text{\rm ev}^{\partial,2}) 
\mathcal M_{\ell}(\vec{\kappa},\vec p;B) \to X^{\ell} \times L^{\vec{\kappa},\vec p}_{\text{\rm source}}.
\end{equation}
Here ${\rm ev}^{\partial}_i$ is the evaluation at the $i$-th marked point $z(i)$ of
$\mathcal M_{\ell}(\vec \kappa,\vec p;B)
=
\mathcal M_{\ell,\vec k'}((\vec{\kappa}',\vec p,m');B)$.
$m'$ is used to define $z(i)$ as in Definition \ref{defn917}.
In other words, the $\tilde L_{\kappa_{i-1}} \times_{X} \tilde L_{\kappa_i}$ component 
of $\text{\rm ev}^{\partial}$  assigns $p_{i}$  in case  Definition \ref{defn927777} (1),(2).
Otherwise, this component is given by the boundary evaluation map at ${z}(i)$ to $\tilde L_{\kappa_{i}}$, which agrees with the evaluation map at the corresponding the marked points $w_{i',j}$ under the relabeling of elements.
The map
$\text{\rm ev}^+ : \mathcal M_{\ell}(\vec{\kappa},\vec p;B) 
\to X^{\ell}$ is induced by the evaluation map at $z^+_i$  on 
$\mathcal M_{\ell, \vec k'}((\vec{\kappa}',\vec p,m');B)$.
\begin{rem}
Note that the value of the $(\tilde L_{\kappa_{i-1}} \times_{X} \tilde L_{\kappa_i})$-component of $\text{\rm ev}$  with $L_{\kappa_{i-1}}  \ne L_{\kappa_{i}} $ 
is determined by the data $B$, so is just a constant. 
Later on we will sum up the contributions of these evaluation maps  for various $B$.
\end{rem}
\par
We decompose 
$\text{\rm ev}^{\partial} : 
\mathcal M_{\ell}(\vec{\kappa},\vec p;B) \to L^{\vec{\kappa},\vec p}_{\text{\rm source}}$
as $\text{\rm ev}_{0} \times \text{\rm ev}^{\partial}_{>0}$
where the first factor is the $0$-th factor of $\text{\rm ev}^{\partial}$.
(Its target is a subset of $\tilde L_{\kappa_{0}} \times_{X} \tilde L_{\kappa_{K}}$, that is either a finite set or $\tilde L_{\kappa_{0}}$.)
\par
Let $\text{\bf g} = g_1\otimes \dots\otimes g_{\ell} \in \Omega(X)^{\otimes \ell}$ and 
$\text{\bf h} = h_1 \otimes \dots \otimes h_{K} \in BCF(\cL^{\rm form};\vec{\kappa})$.%  and 
\begin{defn}\label{def827}
$$
\aligned
\frak q^{\text{\rm form}}_{\ell;\vec{\kappa},\vec p;B}(\text{\bf g},\text{\bf h} )
%,\text{\bf h}')
&=
(-1)^{\maltese}(\text{\rm ev}_{0})_!
\left((\text{\rm ev}^+)^*\text{\bf g}
\wedge
(\text{\rm ev}^{\partial}_{>0})^*\text{\bf h}
% \wedge
% (\text{\rm ev}^{\partial,2})^*\text{\bf h}'
\right) \\
&
= (-1)^{\maltese}\text{\rm Corr}
(\mathcal M_{\ell}(\vec{\kappa},\vec p;B);(\text{\rm ev}_{0};
\text{\rm ev}^+\times \text{\rm ev}^{\partial}_{>0}% \times 
% \text{\rm ev}^{\partial,2}
))
(\text{\bf g}\wedge \text{\bf h} % \wedge \text{\bf h}'
)\\
&
\in CF(L_{\kappa_{0}},L_{\kappa_{K}};\F),
\endaligned$$
where the notation $\text{\rm Corr}$ in the
second line is defined in (\ref{smoothcorrespondence}).
The sign ${\maltese}$ is $\epsilon(h_1,\dots,h_k) + \delta(g_1,\dots,g_{\ell})$ as in (\ref{form21ten1}), (\ref{form21ten2}).
Here and hereafter we omit the symbol ${}^{\boxplus 1}$ in the notation of 
integration along the fiber.
\par
In case all Lagrangians agree with $L_{\kappa_0}$,  %$\vec k = k$, $\vec{\kappa} = i$,
we define $\frak q^{\text{\rm form}}_{\ell,K;\beta}$\index[syindex]{qformellvecbeta@$ \frak q^{\text{\rm form}}_{\ell,K;\beta}$} by:
$$ \frak q^{\text{\rm form}}_{\ell,K;\beta}
= \sum_{[B] \in \beta}\sum_{\vec p}\frak q^{\text{\rm form}}_{\ell;\vec{\kappa},\vec p;B}.
$$
Here $B$ determines $[B]
\in H_2(X;L_{\kappa};\Z)/\sim$ in an obvious way 
in case $K=\blue{0}$. 
($B \sim B'$ is defined by
$B(\omega) = B(\omega')$ and 
$\mu(B) = \mu(B')$, where $\mu$ is an appropriate 
polygonal Maslov index. See \cite[Definition 5.10]{anchor}.)
\par
In case $\ell = 0$, $\vec{\kappa} = (\kappa_1,\kappa_2)$, $B=0$, $L_{\kappa_1} = L_{\kappa_2}$,
$\theta_{\kappa_1} \ne \theta_{\kappa_2}$,
we define
\begin{equation}\label{locacoeffdif}
\frak q^{\text{\rm form}}_{0;(\kappa_1,\kappa_2);0}(1,h )
= dh +  (\theta_{\kappa_2} - \theta_{\kappa_1})\wedge h.
\end{equation}
Note that this definition coincides with (\ref{defqformula2}) in case 
$\theta_{\kappa_1} = \theta_{\kappa_2}$.
 %\marginpar{Formula (\ref{locacoeffdif}) is added.
%KF. 2024 Dec.}
\end{defn}
\par
These operators have the properties stated in Lemma 
\ref{qpropertiescat} below.

\begin{defn}
We define the pairing \index[syindex]{<>_cyc@$\langle *,* \rangle_{\text{\rm cyc}}$}
$$
\langle*,*\rangle_{\rm cyc} :
CF(L;L';\F) \otimes_{\F} CF(L';L;\F) \to \F.
$$

In case $L\ne L'$ (and so $L$ is transversal to $L'$)
$CF(L;L';\F)$ and $CF(L';L;\F)$ both are vector space 
whose basis is identified with $L\cap L'$. We put
\begin{equation}\label{form823}
\langle [p_1], [p_2] \rangle_{\rm cyc} 
= 
\begin{cases}
\pm 1   & \text{$p_1$, $p_2$  correspond to the same points.} \\
0            & \text{otherwise}.
\end{cases}
\end{equation}
The sign $\pm$ is defined in Subsection \ref{sec:signpair}.
In case $L=L'$ 
$$
CF(L;L;\F) = \Omega(\tilde L) \oplus  \bigoplus_{p \in (\tilde L \times_X \tilde L) \setminus \tilde L} \F [p].
$$
The paring among the elements of $\Omega(\tilde L)$ is defined by (\ref{pairing}).
For the self-intersection point $p = (p^1,p^2) \in (\tilde L \times_X \tilde L) \setminus \tilde L$
we put
\begin{equation}\label{form8232}
\langle [(p^1,p^2)], [(q^1,q^2)] \rangle_{\rm cyc} 
= 
\begin{cases}
\pm 1   & \text{$p_1 = q_2$, $p_2 = q_1$.} \\
0            & \text{otherwise}.
\end{cases}
\end{equation}
The sign $\pm$ is defined in Subection \ref{sec:signpair}.
\end{defn}
We  make the following choice of the 
orientation.
We recall that
we fixed an orientation of the index bundle of the operator (\ref{orientationcomplex}).
\begin{equation}\label{operatorpLL'}
\overline{\partial}_p : L^{1,q}(T_pX;L,L';\lambda_p) \to L^{q}(T_pX\otimes \Lambda^{0,1}).
\end{equation}
See (\ref{form815}), (\ref{orientationcomplex}) for the definition of the domain.

We regard the same $p \in L \cap L'$ as an intersection of $L'$ and $L$.
In this case we write $p^o$ to distinguish it from $p$.
We took a path $\lambda_p : [0,1] \to {\mathcal{LAG}}_+(T_pX)$ such that
$\lambda_p(0) = T_{p}L$, 
$\lambda_p(1) = T_{p}L'$.
We define 
$\lambda_{p^o}$ by $\lambda_{p^o}(t) = \lambda_p(1-t)$.
Then we have a Fredholm operator:

\begin{equation}\label{operatorpL'L}
\overline{\partial}_{p^o} : L^{1,q}(T_{p^o}X;L',L;\lambda_{p^o}) \to L^{q}(T_{p^o}X
\otimes \Lambda^{0,1}).
\end{equation}
We fixed $o_p$ an orientation of the index bundle of (\ref{operatorpLL'}).
It induces an orientation $o_{p^o}$ of the index bundle of (\ref{operatorpL'L}) as follows.
\par
We define
$$
I : [0,\infty) \times [0,1] \to (-\infty,0] \times [0,1]
$$
by $I(\tau,t) = (-\tau,1-t]$. This is a biholomorphic map.
Using it we can identify (\ref{operatorpL'L}) with an operator
on $[0,\infty) \times [0,1]$, that we write
$I^*(\overline{\partial}_{p^o} )$.
We glue $I^*(\overline{\partial}_{p^o} )$ and (\ref{operatorpLL'}) at 
$\infty \times [0,1]$ and $-\infty \times [0,1]$.
We then obtain a Fredholm operator, which we describe below.
\par
We consider the rectangle $[-T,T] \times [0,1]$ for large $T$. We consider the 
maps $\xi : [-T,T] \times [0,1] \to T_pX$ with the following boundary conditions.
\begin{equation}\label{gluedboundary}
\xi(\tau,0) \in T_pL, \quad \xi(\tau,1) \in T_pL',
\quad \xi(-T,t) = \lambda_p(1-t), \quad
\xi(T,t) = \lambda_p(t).
\end{equation}
Then the Cauchy-Riemann operator defines a Fredholm operator
\begin{equation}\label{operatorpglued}
\overline{\partial}_{p^o\#p} : L^{1,q}([-T,T]\times [0,1];T_{p^o}X;(\ref{gluedboundary})) \to L^{q}([-T,T]\times [0,1];T_pX \otimes \Lambda^{0,1}).
\end{equation}
It is now standard to find the canonical isomorphism:
\begin{equation}\label{indexsum}
\text{\rm Index}(I^*(\overline{\partial}_{p^o} )) \oplus \text{\rm Index}(\overline{\partial}_p)
\cong 
\text{\rm Index}(\overline{\partial}_{p^o\#p}).
\end{equation}
On the other hand the path $\partial([-T,T]\times [0,1]) \to {\mathcal{LAG}}_+(T_pX)$ defining 
the boundary condition (\ref{gluedboundary}) is homotopic to zero.
In case when this path becomes trivial the index of the corresponding operator 
becomes $T_pL$. Therefore the  orientation of $\text{\rm Index}(\overline{\partial}_{p^o\#p})$ 
canonically identified with the orientation of $T_pL$. Since we assumed that 
$L$ is oriented,  $\text{\rm Index}(\overline{\partial}_{p^o\#p})$ is oriented.
Therefore by (\ref{indexsum}) the orientation $o_p$ canonically 
induces $o_{p^o}$.

\begin{lem}\label{qpropertiescat}
The operators $\frak q^{\text{\rm form}}_{\ell;\vec{\kappa};B}$ have the following properties:
\begin{enumerate}
\item
The following equality holds.
\begin{equation}\label{qmaineqcat}
\aligned
0= &\sum
(-1)^{\maltese_1}
\frak q^{\text{\rm form}}_{\ell_1;\vec{\kappa}_1,\vec p_1;B_1}(\text{\bf g}^{2;1}_{c_1};
\text{\bf h}^{3;1}_{c_2} \otimes \frak q^{\text{\rm form}}_{\ell_2;\vec{\kappa}_2,\vec p_2;B_2}(\text{\bf g}^{2;2}_{c_1};\text{\bf h}^{3;2}_{c_2})
\otimes \text{\bf h}^{3;3}_{c_2})
\\
&+
 \frak q^{\text{\rm form}}_{\ell;\vec{\kappa},\vec p;B}(d\text{\bf g},\text{\bf h})
% \\
% & + 
% \sum (-1)^{\maltese_2}
% \frak q^{\text{\rm form}}_{\ell';\vec{\kappa};B'}(\text{\bf g}^{2;1}_{c_1}; \text{\bf h};\widehat{\frak q}^{\blue{\text{\rm form}}}_{\beta}(\text{\bf g}^{2;2}_{c_2};\text{\bf h}'))
% \\
% &= 0
,
\endaligned
\end{equation}
where
$
\maltese_1 = \deg'  \text{\bf h}^{3;1}_{c_2} + \deg \text{\bf g}^{2;1}_{c_1}
+ \deg\text{\bf g}^{2;2}_{c_1} \deg'  \text{\bf h}^{3;1}_{c_2}
%+  \deg'  \text{\bf h}^{2;1}_{c_3} \deg'  \text{\bf h}^{3;3}_{c_2}
$.
% and
% $
% \maltese_2 = \deg'  \text{\bf h} \deg \text{\bf g}^{2;1}_{c_1}
% +  \deg'  \text{\bf h} + \deg \text{\bf g}^{2;1}_{c_1}
% $.
The operator $\text{\bf g} \mapsto d\text{\bf g}$ is defined from the de Rham differential $d$ by 
$$
d(g_1 \otimes \cdots \otimes g_{\ell})
= 
\sum_{i=1}^{\ell} (-1)^\maltese g_1 \otimes \cdots \otimes
dg_i \otimes \cdots \otimes g_{\ell}
$$
where $\maltese = \deg g_1 + \dots + \deg g_{i-1}$.
\par
The summation in the first line of $(\ref{qmaineqcat})$ is taken over all $c_1,c_2$ and $B_1, B_2$
with $B = B_2\#_j B_1$, where $j-1$ is the length filtration of 
$\text{\bf h}^{3;1}_{c_2}$.
The integers $\ell_1, \ell_2$ are determined by $\text{\bf g}^{2;1}_{c_1}$ 
and $\text{\bf g}^{2;2}_{c_1}$,%  while
% $\vec k_1, \vec k_2$ are determined by 
% $(\text{\bf h}')^{2;1}_{c_3}, (\text{\bf h}')^{2;2}_{c_3}$,
and
$\vec{\kappa}_1, \vec{\kappa}_2$ are determined by 
$\text{\bf h}^{3;1}_{c_2}, \text{\bf h}^{3;2}_{c_2},
\text{\bf h}^{3;3}_{c_2}$.
\par
% The summation in the second line of $(\ref{qmaineqcat})$ is taken over all $c_1, c_2$ and 
% $\beta, B'$ with $B' \# \beta = B$, and
% $\ell'$ is the length filtration of 
% $\text{\bf g}^{2;1}_{c_1}$.
% Precisely speaking,
% $\vec k'$ depends on the term of 
% $\widehat{\frak q}_{\beta}(\text{\bf g}^{2;2}_{c_2};\text{\bf h}')$
% and determined by it.

\item The following equality holds.
\begin{equation}\label{qismcat}
\aligned
&\langle\frak q^{\text{\rm form}}_{\ell;\vec{\kappa},\vec p;B}
(\text{\bf g};h_0,h_1,\dots,h_{K-1}% ;\text{\bf h}'
),h_K\rangle \\
&= (-1)^\maltese
\langle\frak q^{\text{\rm form}}_{\ell;c\vec{\kappa},c\vec p;cB}
(\text{\bf g};h_1,\dots,h_K% ;\text{\rm cyc}(\text{\bf h}'
)),h_0\rangle,
\endaligned
\end{equation}
where 
$\maltese = (\deg' h_0)(\deg' h_1 + \dots + \deg' h_K)$. 
{We recall the definition of $c \vec \kappa$ from Equation % \eqref{eq:cyclic_rotation_k} and
\eqref{eq:cyclic_rotation_kappa}, and note that this map induces an isomorphism} % Note $c\vec k = (k_1,\dots,k_K,k_0)$ 
% where $\vec k = (k_0,\dots,k_K)$. $c\vec{\kappa}$ is defined in the same way. It induces
$c : \Pi_2(\vec{\kappa},\vec p) \to \Pi_2(c\vec{\kappa},c\vec p)$.
 %\marginpar{(3) (the unitality) is removed and 
%is moved to the next section. KF}
% \end{equation}
% 
\end{enumerate}
\end{lem}
\begin{proof}
This a consequence of Propositions \ref{Kuraeistspoly}, \ref{existmultipolu}
and Theorem \ref{orithem}.
Formula (\ref{qmaineqcat}) is a consequence of the 
description of the boundary 
Proposition \ref{Kuraeistspoly} (3), compatibility of 
family of CF-perturbations 
(Proposition \ref{existmultipolu} (2)), Stokes' theorem 
\cite[Theorem 8.11]{springer} 
and the Composition formula \cite[Theorem 10.21]{springer}.
The formula (\ref{qismcat}) follows from cyclic symmetry of the moduli space 
(Proposition \ref{Kuraeistspoly} (7)) and 
its \blue{CF-perturbation} (Proposition \ref{existmultipolu} (5)).
\end{proof}
We remark that the unitality is not claimed in Lemma \ref{qpropertiescat}.
We will discuss it later.

Below we explain a case where the perturbation of the constant map is involved.
 %\marginpar{An example is added.}
\begin{exm}\label{exam;constmap1}
Let us consider the case when we have two Lagrangian submanifolds $L_0$ and $L_1$.
Let $p \in L_0 \cap L_1$. We write $p_{01}$ (resp. $p_{10}$) when 
we regard it as an intersection point $L_0 \cap L_1$ (resp. $L_1 \cap L_0$).
We consider 
$
\frak q^{\text{\rm form}}_{\ell;\vec{\kappa},\vec p;B}
$
with $\ell = 0$, $\vec{\kappa} = (0,1)$, $B = 0$
and study
\begin{equation}\label{eq823}
\frak q^{\text{\rm form}}_{0;(0,1),(p_{01},p_{10});0}([p_{01}],[p_{10}]).
\end{equation}
It will be a differential form on $L_0$ and easy 
degree counting shows that it is a differential $n$-form where $n= \dim L_0$.
The moduli space we use to define this operation consists of 
pseudo-holomorphic triangle (of area $0$), 
so that the three edges (in counter clockwise order) are mapped to $L_0$, $L_1$, $L_0$.
The 0-th marked point is a vertex between the last and the first edges. 
We require first and second marked points are mapped to $p_{01}$, $p_{10}$
respectively.  Such a moduli space consists of one point which is a 
constant map to $p = p_{01} = p_{10}$.
This space itself is transversal but does {\it not} satisfies 
Proposition \ref{Kuraeistspoly} (5).  Namely the evaluation map at the $0$-th marked 
point is not a submersion to $L_0$.  
If we use this (unperturbed) moduli space 
then (\ref{eq823}) would be a delta form supported at $p \in L_0$.
Since we are constructing structure operations on the space of 
smooth forms this is not the correct one.
We need to include an obstruction bundle and an appropriate 
CF-perturbation.  Then (\ref{eq823}) becomes an $n$-form 
which is a smoothing of the delta form.
The $n$-form itself depends on the CF-perturbation.
\par
We remark that here we need to perturb the moduli space of constant maps.
The compatibility with forgetful map then becomes an issue.
This point is further discussed in Sections \ref{sec:unit} and \ref{sec:Kuraconst}.
\par
If we regard the diagonal marked point (which is mapped to $L_0$) as 
an input then we obtain a formula
$$
\overline{\frak m}^{\rm form}_2(h,[p_{01}]) = \int_{L_0}h \rho  [p_{01}]
$$
Here $h$ is a differential $0$ form (function) on $L_0$ and $\rho$ is an appropriate smoothing of delta form at $p$.
This is the formula after we perturb.
Before perturbation the formula is
$$
\frak m^{\rm form}_2(h,[p_{01}]) = h(p) [p_{01}]
$$
They coincide if $h$ is harmonic (that is, constant).
This calculation appears in the study of $\frak m_1\circ \frak m_1$.
\end{exm}

We next use the bulk class %and weak bounding cochains 
to define a curved cyclic filtered $A_{\infty}$ category {$\cL^{\rm form}_{\rm curve}$.} 
We mention that the endomorphism space of any object is given by the de Rham {\it complex} of the underlying Lagrangian.
\par
For the next statement, we recall the splitting (\ref{decomposefrakb}) of the ambient cohomology 
class $\frak b$ as a sum of elements $\frak b_0$, $\frak b_2$, $\frak b_+$, and recall that we fixed a representative of each of these classes by closed forms (with coefficients in $\Lambda_0$, $\F$, or $\Lambda_+$). We also recall that each Lagrangian $L_\kappa$ is equipped with a primitive $\theta_{L_\kappa}$ for the restriction of $\frak b_2$. %  We then in Subsection \ref{ainfalgasssingle} obtain the structure  
% (This involves various choices. But the resulting 
%  filtered $A_{\infty}$ algebra is independent of the 
%  choice up to homotopy equivalence.)

%  {\color{blue} As in Section \ref{ainfalgasssingle}, we consider the additional data of odd, closed, differential forms $b_{\kappa}$ supported on Lagrangians  $L_{\kappa} \in \sL$.} \marginpar{MA190528: I believe that we should be picking all such odd forms $b_\kappa$ at this stage, rather than enumerate them. I implicitly assume that this is done below.} We denote by $\frak b$ the totality of the data $\frak b$, $\theta_{L_{\kappa}}$, 
% $b_{\kappa}$. 
% Recall that we decomposed  
% $
% b_{\kappa} = b_{\kappa,0} + b_+
% $
% as in (\ref{bdecomp}).
% and fixed the closed forms representing $b_{\kappa,0}$, $b_+$.
\begin{defn}\label{def936}
We define \index[syindex]{rhobtheta@$\rho_{\frak b,\theta}$}
$
\rho_{\frak b,\theta} : \Pi_2(\vec{\kappa},\vec p) \to \Lambda_0
$
by
\begin{equation}
\rho_{\frak b,\theta}(B) 
= \exp\left(\int_{D^2} u^* \frak b_2 
+\sum \int_{z_{i}}^{z_{i+1}} u^* \theta_{{\kappa_i}}\right).
\end{equation}
Here $(u;\vec z)$ is a representative of the equivalence class of $B$ in the sense of Definition \ref{def211}.
The equality (\ref{fixthetaL}) implies that the right 
hand side is independent of the representative.
\par
% We next define 
% $
% \rho_{b} : \Pi_2(\vec{\kappa},\vec k) \to \F
% $
% by
% \begin{equation}
% \rho_{b}(B) 
% =  \exp\left( \blue{ \sum_{i=0}^K }\int_{z_i}^{z_{i+1}} u^* b_{\kappa_i,0}\right).
% \end{equation}
% Since $b_{\kappa,0}$ are closed one forms this is invariant of the 
% representative $(u;\vec z)$.
\par
We also define\index[syindex]{omegaB@$\omega(B)$}
$$
\omega(B) = \int_{D^2} u^*\omega.
$$
This is independent of the representative since $L_{\kappa}$ are 
Lagrangian submanifolds.
\end{defn}

{We now introduce the maps which yield the curved filtered $A_\infty$ structure:}
\begin{defn}\label{defn626}
For given $\vec{\kappa}$ we define\index[syindex]{qformfrabellkappa@$\frak q^{\text{\rm form},\frak b}_{\ell,\vec{\kappa}}$} 
$$
\frak q^{\text{\rm form},\frak b}_{\ell,\vec{\kappa}}
: 
\Omega(X)^{\otimes \ell}
\otimes 
BCF(\cL^{\rm form};\vec{\kappa})
\to CF(\newred{L_{\kappa_{\blue{0}}},L_{\kappa_{\blue{K}}}};\F)
\widehat{\otimes} \Lambda_0,
$$
by
\begin{equation}\label{qcat}
\frak q^{\text{\rm form},\frak b}_{\ell,\vec{\kappa}}(\text{\bf g};\text{\bf h})
=
\sum_{\vec p,B,\ell'}T^{\omega(B) }
\frac{\rho_{\frak b,\theta}(B)% \rho_{b}(B)
}{\ell'!}
\frak q^{\text{\rm form}}_{\ell+\ell';\vec{\kappa},\vec p;B}
(\text{\bf g}\otimes \frak b_+^{\ell'};
\text{\bf h}) %;b_+^{\vec k})
\end{equation}
{where the sum is taken over all classes $B$ and $\ell'$with $\omega(B) + \ell' e  \leq E$.}
Here $e>0$ is defined such that $\frak b_+ \equiv 0 \mod T^e$.
\begin{rem}
Note that in case $K=0$ the operation (\ref{qcat}) coincides with  $\frak q$, 
in Proposition \ref{qproperties}.  
 %\marginpar{remark is added.}
The one in Proposition \ref{qproperties} is unital but 
(\ref{qcat}) is not unital in the case of $K>0$.
\end{rem}
We also define \index[syindex]{mxformfrabellkappa@$\frak m^{\text{\rm form},\frak b}_{\vec{\kappa}}$}
$$
\frak m^{\text{\rm form},\frak b}_{\vec{\kappa}}
: 
BCF(\cL^{\rm form};\vec{\kappa})
\to CF(L_{\kappa_0},L_{\kappa_{K}};\F),
$$
by
\begin{equation} \label{eq:a-oo-operation-trivial-gauge}
\frak m^{\text{\rm form},\frak b}_{\vec{\kappa}}(\text{\bf h})
= 
\frak q^{\text{\rm form},\frak b}_{0,\vec{\kappa}}(1;\text{\bf h}).
\end{equation}
Note that in case $K=0$ this becomes $\frak m^{\text{\rm form},\frak b}_k$ in 
Definition \ref{deformedqdef}.
\end{defn}

Now Lemma \ref{qpropertiescat} and  the equality
$$
\rho_{\frak b,\theta}(B_1\#_jB_2) 
= \rho_{\frak b,\theta}(B_1)\rho_{\frak b,\theta}(B_2), 
% \quad
% \rho_{b}(B_1\#_jB_2) 
% = \rho_{b}(B_1)\rho_{b}(B_2)
$$
imply the following:
\begin{prop}\label{prop627}
\blue{When $\text{\bf g}$ consists of closed differential forms, 
the} operators $\frak q^{\text{\rm form},\frak b}_{\ell,\vec{\kappa}}$ have the following properties:
\begin{enumerate}
\item
The  equality 
\begin{equation}\label{qmaineqcat2}
\aligned
\sum
(-1)^{\maltese}
\frak q^{\text{\rm form},\frak b}_{\ell_1;\vec{\kappa}_1}(\text{\bf g}^{2;1}_{c_1};
\text{\bf h}^{3;1}_{c_2} \otimes
\frak q^{\text{\rm form},\frak b}_{\ell_2,\vec{\kappa}_2}(\text{\bf g}^{2;2}_{c_1};\text{\bf h}^{3;2}_{c_2})
\otimes \text{\bf h}^{3;3}_{c_2})
= 0
\endaligned
\end{equation}
holds, where
$
\maltese = (\deg'\text{\bf g}^{2;2}_{c_1} +1)\deg' \text{\bf h}^{3;1}_{c_2}.
$
\par
The summation is taken over all $c_1,c_2$, where 
the integers $\ell_1, \ell_2$ are determined by $\text{\bf g}^{2;1}_{c_1}$, 
$\text{\bf g}^{2;2}_{c_2}$ and 
$\vec{\kappa}_1, \vec{\kappa}_2$ are determined by 
$\text{\bf h}^{3;1}_{c_2}, \text{\bf h}^{3;2}_{c_2},
\text{\bf h}^{3;3}_{c_2}$.
\item The following equality holds. %\marginpar{Unitality is removed 
%from the statement.}
\begin{equation}\label{qismcat2}
\aligned
&\langle\frak q^{\text{\rm form},\frak b}_{\ell,\vec{\kappa}}
(\text{\bf g};h_0,h_1,\dots,h_{K-1}),h_K\rangle_{\rm cyc} \\
&= (-1)^\maltese
\langle\frak q^{\text{\rm form},\frak b}_{\ell;c\vec{\kappa}}
(\text{\bf g};h_1,\dots,h_K)),h_0\rangle_{\rm cyc} ,
\endaligned
\end{equation}
where $c\vec{\kappa}$ is defined as before and 
$\maltese = \deg'h_0 (\deg'h_1+\dots+\deg'h_K)$.
\end{enumerate}
\end{prop}
\begin{rem}
Since Kuranishi structures and CF perturbations are {\it not}
compatible with forgetful map of the interior and boundary marked points, we
first construct the operator $\frak q^{\text{\rm form},\frak b}_{\ell;\vec{\kappa}}$
for
$\ell \le \ell_0$ and
$\vert\vec{\kappa}\vert  \le K$ with given $\ell_0$ and $K$ only
and then use inductive argument.
See Remark \ref{Rmk:7.5} for the case of a single Lagrangian submanifold.
\end{rem}
\begin{cor}\label{existsderham}
The operators
$\{\frak m^{\text{\rm form},\frak b}_{\vec{\kappa}}\}$ with $\vert\vec{\kappa}\vert  \le K$
define a structure of cyclic filtered $A_{K}$
category {modulo $T^{E}$}. %\marginpar{Unitality is removed 
%from the statement.}
\end{cor}
%{We denote the category constructed above by $\cL^{\rm form}_{\rm curv}(\Lambda_0/T^E)$.}
%\par
Note we have not yet discussed the unitality at this stage.

\subsection{Construction of  cyclic $A_{\infty}$ category 4: Homotopy equivalences of $A_{n,k}$ structures} \label{sec:constr-cycl-a_infty4}
We begin by adapting Kajiura's definition  \cite{Ka} of cyclic $A_\infty$ functors to the  filtered setting. %\marginpar{This subsection 
%is written using double projective limit.}
\begin{defn}
We say a filtered $A_{\infty}$ functor $\{\frak f_{\vec{\kappa}}\}$ 
between cyclic filtered $A_{\infty}$ categories to be {\it cyclic}\index{cyclic filtered $A_{\infty}$ functor}\index{cyclic}
if the following holds.
\begin{equation}\label{836form}
\sum_{i=0}^{k}
\langle
\frak f_i(x_1,\dots,x_{i}), \frak f_{k-i}(x_{i+1},\dots,x_k)
\rangle
=
\begin{cases}
0  & k>2 \\
\langle x_1,x_2 \rangle & k = 2.
\end{cases}
\end{equation}
\end{defn}
Since our construction of Kuranishi structure is not compatible with the 
forgetful map of boundary marked point we need to make
double induction over the energy and the number of inputs, to construct
the relevant operations. For this purpose,
we apply the method of \cite[Subsection 7.2.6]{fooo092}.

We fix a discrete submonoid $G$ of $\Z$ and our $A_{\infty}$ structures ($A_{\infty}$ category, 
$A_{\infty}$ functor and etc.) are all $G$-gapped.
For $E \in G$ we put
\begin{equation}
\Vert E \Vert = \sup \{ n \mid \exists E_i \in G, \sum_{i=1}^n E_i = E \}
\end{equation} 
\begin{defn}\label{defn837837}
We define a partial order $<$ on $(G \times \Z_{\ge 0} \times \Z_{\ge 0}) \setminus \{(0,0)\}\times\Z_{\ge 0}$ 
by defining $(E_1,k_1,\ell_1) > (E_2,k_2,\ell_2)$  to be
\begin{enumerate}
\item $2\Vert E_1\Vert + k_1 + \ell_1 > 2\Vert E_2\Vert + k_2+ \ell_2$, $\ell_1 \ge \ell_2$ or
\item  $2\Vert E_1\Vert + k_1 + \ell_1 = 2\Vert E_2\Vert + k_2 + \ell_2$, $\Vert E_1\Vert > \Vert E_2\Vert$, $\ell_1 \ge \ell_2$ or
\item  $2\Vert E_1\Vert + k_1 + \ell_1 = 2\Vert E_2\Vert + k_2 + \ell_2$, $\Vert E_1\Vert = \Vert E_2\Vert$, $\ell_1 > \ell_2$
\end{enumerate}
We write $(n_1,k_1,\ell_1) > (n_2,k_2,\ell_2)$ if there exist $E_1,E_2$ 
with $\Vert E_i \Vert = n_i$ and $(E_1,k_1,\ell_1) > (E_2,k_2,\ell_2)$.
\par
We define the notation $(n_1,k_1,\ell_1) > (E_2,k_2,\ell_2)$ etc.  in a similar way.
For $B \in \Pi_2(\vec{\kappa},\vec p)$ or $\beta \in H_2(X,L)$
we define the notations $(n_1,k_1,\ell_1) > (B,k_2,\ell_2)$, $(\beta,k_1,\ell_1) > (E_2,k_2,\ell_2)$
etc. by $(n_1,k_1,\ell_1) > (B\cap \omega,k_2,\ell_2)$, $(\beta\cap \omega,k_1,\ell_1) > (E_2,k_2,\ell_2)$ etc.
\end{defn}
\begin{lem}
If $\mathcal M_{\ell'}(\vec{\kappa}',\vec p';B')$
appears as a fiber product factor of  the boundary of 
$\mathcal M_{\ell}(\vec {\kappa},\vec p;B)$ then
$(B',\vert \vec{\kappa}'\vert,\ell') < (B,\vert \vec{\kappa}\vert,\ell)$, where 
$k = \vert \vec\kappa\vert$ and $k' = \vert \vec\kappa'\vert$. %\marginpar{The statement corrected. KF 2025 Aug. 27}
\end{lem}
The proof is easy and is omitted.\footnote{Note $k' > k$ can occur.
This is the case when there is a disk bubble (without boundary marked point).}
\begin{defn}\label{defn838}
An {\it $A_{n,k_0,\ell_0}$ structure}\index{$A_{n,k_0,\ell_0}$ structure} with bulk is a series of operators:\index[syindex]{qformfrabellkappaB@$\frak q^{\text{\rm form},\frak b}_{\ell;\vec{\kappa};B}$}
$$
\frak q^{\text{\rm form},\frak b}_{\ell;\vec{\kappa};B}:
\Omega(X)^{\otimes \ell} \otimes  BCF(\cL^{\rm form};\vec{\kappa})
\to CF(L_{\kappa_{0}},L_{\kappa_{K}};\F)
$$
for $(B,k,\ell) \le (n,k_0,\ell_0)$, 
such that (\ref{qmaineqcat}) holds.
\end{defn}

In the case $\ell = 0$ an $A_{n,k,\ell}$ structure is called an $A_{n,k}$ structure.\index{$A_{n,k}$ structure@$A_{n,k}$ structure}
We consider $\frak b = \frak b_0 + \frak b_+$ and 
let $\frak b_+ = \sum T^{\lambda_i} \frak b_{+,i}$.
We replace $G$ by the monoid  $G_+$ generated by 
$G$ and $\{\lambda_i \mid i=1, 2,\dots\}$. Then
for a {\it $G$-gapped filtered $A_{n,k_0,\ell_0}$ structure} with bulk $\frak q^{\text{\rm form}}_{\ell;\vec{\kappa};B}$, 
the operations\index[syindex]{mxkbfraB@$\frak m^{{\rm form},\frak b}_{B,k}$}
$$
\frak m^{{\rm form},\frak b}_{B,k}
(\text{\bf h})=
\sum_{\ell,\vert\vec k\vert \le k}T^{\omega(B) }
\frac{\rho_{\frak b,\theta}(B)% \rho_{b}(B)
}{\ell!}
\frak q^{\text{\rm form}}_{\ell;\vec{\kappa};B}
(\frak b_+^{\ell};
\text{\bf h})
$$
define a $G_+$-gapped filtered $A_{n,k}$ structure.
Note that we take $\ell$ sufficiently large determined by $n$ (or $E_n$).
We remark that the upper bound of this $\ell$ has no relation to $\ell_0$ appearing 
in the $A_{n,k,\ell_0}$ structures above.
(The symbol $\ell_0$ there is the number of inputs, which are  differential forms of $\Omega$ that are not $\frak b$.
The symbol $\ell$ here is the number of inputs $\frak b$.)
We denote it by $\cL^{\rm form}_{\frak b;n,k}$.\index[syindex]{Lformbnk@$\cL^{\rm form}_{\frak b;n,k}$}
\begin{defn}
We consider  $\cL^{\rm form}_{\frak b;n,k}$ with 
two $G_+$-gapped filtered $A_{n,k}$ structures 
$\frak m^{{\rm form},\frak b}_{B,k}$ and $\frak m^{\prime \, {\rm form},\frak b}_{B,k}$.
A $G_+$-gapped {\it filtered $A_{n,k}$ functor}t\index{filtered $A_{n,k}$ functor} between 
$G_+$-gapped filtered $A_{n,k}$ categories is defined by
$$
\Phi_{B,k} : BCF(\cL^{\rm form};\vec{\kappa})
\to CF(L_{\kappa_{0}},L_{\kappa_{K}};\F)
$$
such that 
$$
\aligned
&\sum \frak m^{\prime\, {\rm form},\frak b}_{B_0,k_0}
\left(
\Phi_{B_1,k^c_{k_0;1}}(\text{\bf h}_{c}^{k_0;1}),\dots,\Phi_{B_{k_0},k^c_{k_0;k_0}}(\text{\bf h}_{c}^{k_0;k_0})
\right) \\
&= 
\sum (-1)^{\deg' \text{\bf h}_c^{3;1}}\Phi_{B'_c,k_c^{3;1}}(\text{\bf h}_{c}^{3;1},
\frak m^{{\rm form},\frak b}_{B''_c,k^c_{3;2}}(\text{\bf h}_{c}^{3;2}),\text{\bf h}_{c}^{3;3}).
\endaligned
$$
Here the sum of the left hand side is taken over all $B_i$, $k_0$, $c$ with
$B_0+B_1 + \dots+B_{k_0} = B$,
$$
\underbrace{((\Delta \otimes \dots \otimes {\rm id})\circ \dots\circ \Delta)}_{k_0-1} ({\bf h}) = \sum_c \text{\bf h}_{c}^{k_0;1} \otimes \dots \otimes \text{\bf h}_{c}^{k_0;k_0}.
$$ 
and $k_{k_0,j}^c$ is the number of tensor factors of $\text{\bf h}_{c}^{k_0:j}$.
\par
The sum of the right hand side is taken over all 
$B',B''$, $k',k''$, $c$ with $B'+B'' = B$, $(\vert B\vert,\vert \vec{\kappa}\vert) \le (n,k)$.
$$
((\Delta \otimes {\rm id})\circ \Delta) ({\bf h}) = \sum_c \text{\bf h}_{c}^{3;1} \otimes \text{\bf h}_{c}^{3;2} \otimes \text{\bf h}_{c}^{3;3}.
$$
Here $k''$ is the number of tensor factors of $ \text{\bf h}_{c}^{3:2}$.
$k'$ is the sum of the number of tensor factors of $ \text{\bf h}_{c}^{3:1}$,
the number of tensor factors of $ \text{\bf h}_{c}^{3:3}$,  and $1$.
\par
It is called cyclic if (\ref{836form}) holds (after obvious modification).
\end{defn}

\begin{defn}
We define a {\it cyclic pseudo-isotopy}\index{cyclic pseudo-isotopy} between $G$-gapped filtered $A_{n,k}$ categories in the same way as 
Definition \ref{pisotopydef} but requiring $\frak C$ to have operations $\frak m^{\frak C, t}_{k',\beta}$ and $\frak c^{\frak C, t}_{k',\beta}$ 
for $(\beta\cap \omega,k') \le (n,k)$ only.
\par
Non-cyclic version is defined in the same way.
\end{defn}
We will define homotopically unital version later in Definition \ref{defn913}.

\begin{lem}\label{homotopyequiv}
For $n,k$, $n',k'$ with $(n',k')<(n,k)$, there exists a cyclic pseudo-isotopy between $G_+$-gapped cyclic filtered $A_{n',k'}$ categories
 %\marginpar{unital
%is removed from statement.  homotopy equivalence is changed to pseudo-isotopy.}  
$\cL^{\rm form}_{\frak b;n,k}$ and $\cL^{\rm form}_{\frak b;n',k'}$.
\end{lem}
We will prove Lemma \ref{homotopyequiv} in Subsection \ref{sec:homotopyequiv}.

\begin{rem} In Lemma \ref{homotopyequiv} we regard $\cL^{\rm form}_{\frak b;n,k}$ (which is an $A_{n,k}$ category) as an $A_{n',k'}$ category by 
forgetting some of the operations.
\end{rem}
\begin{defn}
Let $(n',k') < (n,k)$.  We say a $G_+$-gapped cyclic filtered $A_{n,k}$ category $\tilde\cL^{\rm form}_{\frak b;n',k'}$
is a {\it promotion}\index{promotion} of a $G_+$-gapped cyclic filtered $A_{n',k'}$ category $\cL^{\rm form}_{\frak b;n',k'}$ 
if $\tilde\cL^{\rm form}_{\frak b;n',k'}$ regarded as a cyclic filtered $A_{n',k'}$ category is strongly isomorphic 
to $\cL^{\rm form}_{\frak b;n',k'}$.
\end{defn}

We next
define a curved filtered $A_{\infty}$ category $\cL^{\rm form}_{\rm curve}(\Lambda_0;\frak b)$ with objects $\mathscr L$ with coefficients in $\Lambda_0$ as an `inverse limit' over the categories $\cL^{\rm form}_{\frak b;n,k}$.  For this purpose, we construct this limit by replacing these categories by ones, for which we have a pseudo-isotopy. 
The proof of the next result is similar to the argument in \cite[Theorem 7.2.72]{fooo09} and \blue{similar to Proposition \ref{htpyinvlim}}.

\begin{prop} \label{prop:commutative_diagram_energy_induction}
In the situation of Lemma $\ref{homotopyequiv}$, there exists a 
$G_+$-gapped cyclic and filtered $A_{n,k}$ category $\tilde{\cL}^{\rm form}_{\frak b;n,k}$ such that 
it is a promotion of $\tilde{\cL}^{\rm form}_{\frak b;n;,k;}$ and  we have a cyclic pseudo-isotopy 
of filtered $A_{n,k}$ categories between  $\tilde{\cL}^{\rm form}_{\frak b;n,k}$ and $\cL^{\rm form}_{\frak b;n,k}$.
Moreover
the restriction of a cyclic pseudo-isotopy (as $A_{n,k}$ categories) is one in Lemma $\ref{homotopyequiv}$ 
when we regard it as a cyclic pseudo-isotopy of $A_{n',k'}$ categories.
\end{prop}
Proposition \ref{prop:commutative_diagram_energy_induction} claims that 
we can add operations $\frak m'_{B'',k''}(\text{\bf h})$ for $(B'',k'') \le (n,k)$, but not $(B'',k'') \le (n',k')$ 
such that together with $\frak m'_{B'',k''}(\text{\bf h})$ which was previously defined for $(B'',k'')\le (n',k')$
to obtain a $G_+$-gapped cyclic $A_{n,k}$ category $\tilde{\cL}^{\rm form}_{\frak b;n,k}$.

We can now define the central objects of our study:
\begin{defn}\label{defn847}
  The category $\cL^{\rm form}_{\rm curve}(\frak b;\Lambda_0)$
  \index[syindex]{Lcatformcurv@$\cL^{\rm form}_{\rm curve}(\frak b;\Lambda_0)$} is the `inverse limit'
  \begin{equation} \label{eq:inverse_limit_categories_energy}
    \liminv_{(n,k)} \tilde{\cL}^{\rm form}_{\frak b;n,k}.
  \end{equation}
The category $ \cL^{\rm form}_{\rm curve}(\frak b;\blue{\Lambda})$ \index[syindex]{Lcatformcurvfield@$\cL^{\rm form}_{\rm curve}(\frak b;\Lambda)$} is the localisation of $\cL^{\rm form}_{\rm curve}(\frak b;\Lambda_0)$ at the maximal ideal:
  \begin{equation}
    \cL^{\rm form}_{\rm curve}(\frak b;\Lambda_0) \otimes_{\Lambda_{0}}  \Lambda.
  \end{equation}
\end{defn}
This construction is mostly the same as Subsection \ref{sec:inverse-limits-categ}, except we are taking double induction here.
\begin{rem}
Since our situation is slightly more complicated than in Subsection \ref{sec:inverse-limits-categ}, 
we explain the meaning of the `inverse limit' in more details.  %\marginpar{Remark added.}

We take a sequence $(n_i,k_i)$ such that $(n_i,k_i) < (n_{i+1},k_{i+1})$.   
We fix $i_0, k_0$ and will construct a structure of $G_+$-gapped cyclic and filtered $A_{n_j,k_j}$ category
on ${\cL}^{{\rm form},(n_j,k_j)}_{\frak b;n_{i_0},k_{i_0}}$  such that
${\cL}^{{\rm form},(n_{j+1},k_{j+1})}_{\frak b;n_{i_0},k_{i_0}}$ is a promotion of 
${\cL}^{{\rm form},(n_j,k_j)}_{\frak b;n_{i_0},k_{i_0}}$.
We can construct such a sequence by induction on $j$ using Lemma \ref{homotopyequiv}.
\par
Now we consider the set of objects and morphism modules of ${\cL}^{{\rm form}}_{\frak b;n_{i_0},k_{i_0}}$.
We define $A_{n,k}$ operations on it for arbitrary $n,k$ as follows. We take $j$ with
$(n,k) < (n_{j},k_{j})$.  Then the $A_{n_{j},k_{j}}$ operations of ${\cL}^{{\rm form},(n_j,k_j)}_{\frak b;n_{i_0},k_{i_0}}$
restrict to  $A_{n,k}$ operations on it.  By construction such $A_{n,k}$ operations are 
independent of the choice of such $j$.  Now the totality of such $A_{n,k}$ operations 
obviously define an $A_{\infty}$ structure on (the morphisms modules of) ${\cL}^{{\rm form}}_{\frak b;n_{i_0},k_{i_0}}$.
\par
We can use Lemma \ref{homotopyequiv} in a similar way that this $G_+$-gapped cyclic and filtered $A_{\infty}$ category
is independent of various choices involved (which we need to take `inverse limit') up to  cyclic pseudo-isotopy of 
 filtered $A_{\infty}$ categories.
\end{rem}

The proof of Theorem \ref{cAinfconst} is now complete admitting Lemma \ref{homotopyequiv}.
\qed

\begin{rem}\label{Remark849}
The proof of Theorem \ref{cAinfconst} shows the following.
Let us fix $(n',k')$ and  
construct a structure of a filtered $A_{n',k'}$ category with object set $\mathscr L$. 
For this purpose, we first need to directly carry out a geometric construction 
by employing a system of Kuranishi structures  and CF-perturbations on the 
relevant moduli spaces.
\par
Then we promote the filtered $A_{n',k'}$ category structure
to a filtered  $A_{\infty}$ category $\cL^{\rm form}_{\rm curve}$
by taking its $A_\infty$ operations $\frak m_k$ in a filtered $A_{n',k'}$ category that we obtain directly from geometry for a sufficient large pair $(n', k')$ relative to the given $k$.
\par
In other words we can take $\cL^{\rm form}_{\rm curve}$ so that its 
operations as a filtered $A_{n',k'}$ category is one we obtain directly from geometry. %\marginpar{Remark added.}
\end{rem}

 %\marginpar{ The rest of this subsection is moved to the next section:}

\subsection{Construction of Pseudo-isotopy}
\label{sec:homotopyequiv}

%input{homotopyequiv.tex}
% \subsection{Statement}
% \label{subsec:isotopystatement}

In this subsection we prove Lemma \ref{homotopyequiv}.
More precisely we will prove  
Proposition \ref{existpisomain} below, from which Lemma \ref{homotopyequiv} will follow
by adapting the proof of Proposition \ref{pisoimplyhomo}.
To state Proposition \ref{existpisomain} we recap the situation 
we work with.
We consider a finite set of immersed
Lagrangian submanifolds $\mathscr L = \{(L_{\kappa},\theta_{\kappa})\}$.
Assume that $L_{\kappa}$ is transversal to $L_{\kappa'}$
if $L_{\kappa'}$ is different from $L_{\kappa}$.
\par
We first fix a discrete submonoid $G$ we work with.
Given a compatible almost complex structure $J$, we consider the set of $E \in \R_{\ge 0}$
which satisfies one of the following:
\begin{enumerate}
\item
There exist $\kappa$ and a $J$-holomophic map $u : (D^2,\partial D^2) \to
(X,L_{\kappa})$ such that 
\begin{equation}\label{areaisbeta}
E  = \int_{D^2} u^* \omega.
\end{equation}
\item
There exist $\vec{\kappa} = (\kappa_0,\dots,\kappa_K)$
and $\vec p$ as in Subsection  \ref{cycliccategorymoduli}
such that there exists a $J$-holomorphic map 
$u : D^2 \to X$ satisfying Items (1)-(3) of Definition \ref{def211}
and that (\ref{areaisbeta}) holds.
\item
There exists a $J$-holomorphic map 
$u : S^2 \to X$ with symplectic area $E$.
\end{enumerate}
Let $G$ be the discrete monoid that is generated by the set of 
$E$ satisfying (1)-(3) above.
The symplectic areas of all the $J$-holomorphic maps 
appearing in this section are elements of $G$.
\par
We next take $\frak b = \frak b_0 + \frak b_2 + \frak b_+$ and define $G_+$ as 
in right after Definition \ref{defn838}.
Therefore the structures and morphisms (after bulk deformation)  we will obtain are all $G_+$-gapped.
We put $G_+ = \{E_0,E_1,\dots\}$, $E_0 =0$, $E_i < E_{i+1}$.
\par
Let $E,E' \in G_+$ and $k,k'$ with $(E',k') < (E,k)$. We define $n,n'$  by $E_n = E$, $E_{n'} = E'$.
\par
Now by the construction of  Subsections \ref{ainfalgasssingle}-\ref{constcyclic}, 
we obtain $\cL^{\rm form}_{\Dat}$\index[syindex]{LormeDat@$\cL^{\rm form}_{\Dat}$} that is a cyclic and curved   $G_+$-gapped filtered $A_{n,k}$ category  
whose object set is $\mathscr L$. 
We denote this structure by $\{\frak m^{\Dat}_k\}$.
Performing the same construction but 
using possibly different perturbation, we  obtain a cyclic curved  $G_+$-gapped filtered $A_{n',k'}$ category  $\cL^{\rm form}_{\Cat} $, with the same\index[syindex]{LormeCat@$\cL^{\rm form}_{\Cat}$} set $\mathscr L$ of objects. We denote this structure by $\{\frak m^{\Cat}_k\}$.
(In other words, $\cL^{\rm form}_{\Cat} = \mathcal L^{\rm form}_{\frak b;n',k'}$,
$\cL^{\rm form}_{\Dat} = \mathcal L^{\rm form}_{\frak b;n,k}$.)
Now the main result of this subsection is stated as follows:
\begin{prop}\label{existpisomain}
The category $\cL^{\rm form}_{\Dat}$ is cyclically pseudo-isotopic to $\cL^{\rm form}_{\Cat} $
as $G_+$-gapped filtered $A_{n',k'}$ category.
\end{prop}
This is the analogue of \cite[Theorem 11.1]{fukaya:cyc} for multiple Lagrangians, and the proof, which we now recall, is entirely parallel:  %\marginpar{MA190629: I removed a discussion of bounding cochains from here}
%\par
% \subsection{Proof}
% \label{subsec:isotopyporoof}

\begin{proof}[Proof of Proposition \ref{existpisomain}]
The proof is a 1-parametric version of the construction of 
Subsections \ref{cycliccategorymoduli}-\ref{constcyclic}
 and proceeds as follows.
 \par
As in Subsection \ref{cycliccategorymoduli}, we start with $(\vec{\kappa},\vec p)$ as in  Definition \ref{def2110}. 
We consider the moduli space ${\mathcal M}_{\ell}((\vec{\kappa},\vec p);B)$ 
obtained in Proposition \ref{Kuraeistspoly} (1).
For each $(B,\vec{\kappa})$ with $(B,\vert\vec{\kappa}\vert) \le (E,K)$, we write ${\mathcal M}_{\ell}((\vec{\kappa},\vec p);B)^{\Dat}$ when we equip this moduli space with the Kuranishi structure on it as in Proposition \ref{Kuraeistspoly}, for the purpose of constructing 
$\cL^{\rm form}_{\Dat}$. Similarly, for each $(B,\vec{\kappa})$ with $(B,\vert\vec{\kappa}\vert) \le (E',K')$, we write ${\mathcal M}_{\ell}((\vec{\kappa},\vec p);B)^{\Cat}$ when using the Kuranishi structure that gives rise to the category $\cL^{\rm form}_{\Cat}$. 
\par
We also took the CF-perturbation on them to construct $\cL^{\rm form}_{\Dat}$, $\cL^{\rm form}_{\Cat}$
that satisfy Proposition \ref{existmultipolu}.
We write them as $\widehat{\frak S}^{\Dat}$, $\widehat{\frak S}^{\Cat}$ respectively.
\par
We consider $[0,1]_t\times {\mathcal M}_{\ell}((\vec{\kappa},\vec p);B)$.
Besides the evaluation maps given in $(\ref{evevev})$ we also consider\index[syindex]{01Mellveck@$[0,1]_t\times {\mathcal M}_{\ell}((\vec{\kappa},\vec p);B)$}\index[syindex]{evt@$\text{\rm ev}_{t}$}
\begin{equation}\label{evaluationmaptime}
\text{\rm ev}_{t} : [0,1]_t\times {\mathcal M}_{\ell;\vec k}((\vec{\kappa},\vec p);B)
\to [0,1]_t
\end{equation}
that is the projection to the first factor.
We denote by $\vert(\vec{\kappa},\vec k)\vert$ the number of boundary marked points of 
an element of ${\mathcal M}_{\ell;\vec k}((\vec{\kappa},\vec p);B)$.

\begin{lem}\label{Kuraonepara}
For each $(n_0,k_0,\ell_0)$, consider the product $[0,1]_t \times {\mathcal M}_{\ell;\vec k}((\vec{\kappa},\vec p);B)$
for a relative homology class $B$ and $n,\vec{\kappa}$ with $(B,\vert(\vec{\kappa},\vec k)\vert,\ell) \le (n_0,k_0,\ell_0)$.
 %\marginpar{
%Check whether change is OK.  KF2025. Aug 27}
\begin{enumerate}
\item
The space $[0,1]_t \times {\mathcal M}_{\ell;\vec k}((\vec{\kappa},\vec p);B)^{\boxplus 1}$  has a
Kuranishi structure.
\item 
The boundary of $[0,1]_t \times {\mathcal M}_{\ell;\vec k}((\vec{\kappa},\vec p);B)^{\boxplus 1}$
is a union of 
$$
\{0,1\} \times {\mathcal M}_{\ell;\vec k}((\vec{\kappa},\vec p);B)^{\boxplus 1}
$$
and the three types of fiber or direct products $(\text{\rm \ref{218para}})$, 
$(\text{\rm \ref{219para}})$, $(\text{\rm \ref{220para}})$.
(Below, we use the notations of $(\ref{218})$,  $(\ref{219})$, $(\ref{220})$, respectively.)
\begin{equation}\label{218para}
\aligned
&([0,1]_t \times{\mathcal M}_{\ell'';k''_i+1}(L_{\kappa_i};\beta)^{\boxplus 1}) \\
&\,{}_{(\text{\rm ev}_t,\text{\rm ev}_0)} \times_{(\text{\rm ev}_t,\text{\rm ev}_{i,j})}
([0,1]_t\times {\mathcal M}_{\ell';\vec k'}((\vec{\kappa},\vec p);B')^{\boxplus 1}),
\endaligned
\end{equation}
\begin{equation}\label{219para}
([0,1]_t\times {\mathcal M}_{\ell'';\vec k''}((\vec{\kappa}'',\vec p'');B'')^{\boxplus 1})
{}_{\text{\rm ev}_t}\times_{\text{\rm ev}_t}
([0,1]_t\times {\mathcal M}_{\ell';\vec k'}((\vec{\kappa}',\vec p');B')^{\boxplus 1}),
\end{equation}
\begin{equation}\label{220para}
\aligned
&([0,1]_t\times {\mathcal M}_{\ell'';\vec k''}((\vec{\kappa}'',\vec p'');B'')^{\boxplus 1})\\
&\,{}_{(\text{\rm ev}_t,\text{\rm ev}_{1,m''})}\times_{(\text{\rm ev}_t,\text{\rm ev}_{i,m'})}
([0,1]_t\times {\mathcal M}_{\ell';\vec k'}((\vec{\kappa}',\vec p');B')^{\boxplus 1}).
\endaligned
\end{equation}
Here the union is taken over the data appearing in $(\ref{218})$,  $(\ref{219})$, $(\ref{220})$, 
respectively
together with the shuffles $(\mathbb L'',\mathbb L')$ of $\underline\ell$
such that $\# \mathbb L'' = \ell''$, $\# \mathbb L' = \ell'$.
\item
The evaluation maps $(\ref{evevev})$, $\text{\rm(\ref{evaluationmaptime})}$ extend to the 
compactification and are compatible with 
the description of the boundary in $(2)$.
\item
Let $i \in \uwave{K}$, $j \in  \{1,\dots,k_i\}$.
Then the evaluation map
$$
(\text{\rm ev}_t,\text{\rm ev}_{i,j}) : 
[0,1]_t\times {\mathcal M}_{\ell;\vec k}((\vec{\kappa},\vec p);B)^{\boxplus 1}
\to [0,1]_t\times L_{\kappa_i}
$$
is weakly submersive.
\item
The restriction of the Kuranishi structure to $
\{0\} \times {\mathcal M}_{\ell;\vec k}((\vec{\kappa},\vec p);B)^{\boxplus 1}
$ and 
$
\{1\} \times {\mathcal M}_{\ell;\vec k}((\vec{\kappa},\vec p);B)^{\boxplus 1}
$
coincides with ${\mathcal M}_{\ell;\vec k}((\vec{\kappa},\vec p);B)^{\Cat,\boxplus 1}$  and
${\mathcal M}_{\ell;\vec k}((\vec{\kappa},\vec p);B)^{\Dat,\boxplus 1}$,  respectively. 
\item
The Kuranishi structure is invariant under the permutation of 
interior marked points.
\item
The Kuranishi structure is invariant under the cyclic permutation 
of the data in the same sense as Proposition $\ref{Kuraeistspoly}$ $(7)$. 
\item
The Kuranishi structures are oriented so that it is compatible with item $(2)$.
\end{enumerate}
\end{lem}
The proof is the same as the proof of Proposition $\ref{Kuraeistspoly}$ and so is omitted.
\par
For $m \in \{0,\dots,k_0\}$ we define $[0,1]_t \times {\mathcal M}_{\ell;\vec k}((\vec{\kappa},\vec p,m);B)
= [0,1]_t \times {\mathcal M}_{\ell;\vec k}((\vec{\kappa},\vec p,m);B)$.\index[syindex]{01Mellveckm@$[0,1]_t\times {\mathcal M}_{\ell}((\vec{\kappa},\vec p,m);B)$}
We use $m$ to specify the $0$-th (boundary) marked point of an object in 
$[0,1]_t \times {\mathcal M}_{\ell;\vec k}((\vec{\kappa},\vec p,m);B)$ in the 
same way as Definition $\ref{defn917}$.
\par
The compactification of $[0,1]_t \times {\mathcal M}_{\ell;\vec k}((\vec{\kappa},\vec p,m);B)$ 
is defined by identifying it with $[0,1]_t \times {\mathcal M}_{\ell;\vec k}((\vec{\kappa},\vec p);B)$
and using Proposition \ref{Kuraonepara}.
\begin{rem}
In the special case $\vec{\kappa} = (\kappa,\dots,\kappa)$, 
the moduli space ${\mathcal M}_{\ell;\vec k}((\vec{\kappa},\vec p);B)$ is 
the one used to define the filtered $A_{\infty}$ algebra for $L_{\kappa}$.
In that case Lemma \ref{Kuraonepara} is used to prove 
Proposition \ref{existpisomain} in the case of a single Lagrangian.
More precisely, Proposition \ref{existpisomain} is proved by induction on the number of distinct Lagrangians $L_{\kappa}$, with the inductive step being the pseudo-isotopy extracted from Lemma \ref{Kuraonepara} in a way 
so that the Kuranishi structures obtained are consistent 
with the earlier choice.
The same remark applies to the choice of \blue{CF-perturbations} in Lemma \ref{existmultipolupara}.
\end{rem}
\begin{rem}
In Lemma \ref{Kuraonepara} and in similar situations we do {\it not} take 
outer collaring of the $[0,1]_t$ factor. Actually it is easy to take Kuranishi structure, CF-perturbation
etc. to be constant in $t$ direction in a neighborhood of $\partial[0,1]_t$.
\end{rem}
\begin{lem}\label{existmultipolupara}
There exists a system of {CF-perturbations} 
%multisections 
of $[0,1] \times {\mathcal M}_{\ell;\vec k}((\vec{\kappa},\vec p);B)^{\boxplus 1}$ 
for  $B,n,\vec{\kappa}$ with $(B,\vert(\vec{\kappa},\vec k)\vert,\ell) \le (n_0,k_0,\ell_0)$,  that has
the following properties.
\begin{enumerate}
\item
They are transversal to $0$.
\item 
They are compatible with the description of the 
boundary in Lemma {\rm \ref{Kuraonepara}} $(2)$.
\item
At $\{0\} \times {\mathcal M}_{\ell;\vec k}((\vec{\kappa},\vec p);B)^{\boxplus 1}$ 
it coincides with $\widehat{\frak S}^{\Cat}$.
At $\{1\} \times {\mathcal M}_{\ell;\vec k}((\vec{\kappa},\vec p);B)^{\boxplus 1}$ 
it coincides with $\widehat{\frak S}^{\Dat}$.
 %\marginpar{compatibility with forgetful map is removed.}
\item
Let $i \in \uwave{K}$, $j \in  \{1,\dots,k_i\}$.
Then the evaluation map
$$
(\text{\rm ev}_t,\text{\rm ev}_{i,j}) : {\mathcal M}_{\ell;\vec k}((\vec{\kappa},\vec p);B)^{\boxplus 1}
\to [0,1] \times L_{\kappa_i}
$$
is strongly submersive with respect to our \blue{CF-perturbations}.
 %\marginpar{ $: 
%[0,1] \times \frak{forget}$ is removed from formula in (4).}
\item
The CF-perturbation is invariant under the permutation of 
interior marked points.
\item
The CF-perturbation is invariant under the cyclic permutation 
of the data in the same sense as Lemma {\rm \ref{Kuraonepara}} $(7)$. 
\end{enumerate}
The same holds for $[0,1]_t\times {\mathcal M}_{\ell}((\vec{\kappa},\vec p,m);B)$.
\end{lem}
The proof is the same as Proposition \ref{existmultipolu} and so is omitted.
\par
Now we use the moduli spaces and its  \blue{CF-perturbations} %multisections 
constructed in 
Lemmas \ref{Kuraonepara} and \ref{existmultipolupara} 
to construct a pseudo-isotopy  $\frak C$  between $\cL^{\rm form}_{\Cat}$ and $\cL^{\rm form}_{\Dat}$ and $A_{n,k}$\index[syindex]{Cxfrak@$\frak C$}
structure. We consider $n_0$, $\vec{\kappa}$, $\ell_0$ 
such that $(n_0,\vert\vec{\kappa}\vert) \le (n,k)$ and $\ell_0 \le n$.
The object set of $\frak C$ is $\mathcal L$, and the morphism set is 
\index[syindex]{CxfrakLkappa@${\frak C}(L_{\kappa}, L_{\kappa'};\F)$}
$$
{\frak C}(L_{\kappa}, L_{\kappa'};\F) = \Omega(\blue{[0,1] \times (\tilde L_{\kappa} \times_X \tilde L_{\kappa'})}) \otimes \F.
$$
We shall use the notation\index[syindex]{BCxfrakLkappa@$B\frak C(\cL^{\rm form};\vec{\kappa})$}
$$
B\frak C(\cL^{\rm form};\vec{\kappa} )
= \bigotimes_{i=1}^K {\frak C}(L_{\kappa_{i-1}},L_{\kappa_{i}};\F),
$$
where we change the notation as in Section \ref{constcyclic} and we adopt the notation from Equation \eqref{eq:change_notation_moduli_space}, denoting by $\mathcal M_{\ell}(\vec \kappa,\vec p;B)$ the moduli space of holomorphic disks with these boundary conditions.
We modify the evaluation map (\ref{evmap2}) to\index[syindex]{evtilde@$\widetilde{\text{\rm ev}}$}
\begin{equation}\label{evmap2para}
\widetilde{\text{\rm ev}} = (\text{\rm ev}^+;\widetilde{\text{\rm ev}}): % ;\widetilde{\text{\rm ev}}^{\partial,2}):
[0,1] \times \mathcal M_{\ell}(\vec{\kappa},\vec p;B)^{\boxplus 1}
\to X^{\ell} \times [0,1] \times L^{\vec{\kappa},\vec p}_{\text{\rm source}}
\end{equation}
(See Definition \ref{defn927777} for the definition of the right hand side.)
in an obvious way.
(Namely we map the first factor $[0,1]$ diagonally to the product of $[0,1]$ factors of the target.)
\par
Let $\text{\bf g} = g_1\otimes \dots\otimes g_{\ell} \in \Omega(X)^{\otimes \ell}$, 
$\text{\bf h} = h_1 \otimes \dots \otimes h_{K} \in B\frak C(\cL^{\rm form};\vec{\kappa})$. \index[syindex]{qcfrakellkappaB@$\frak q^{\frak C}_{\ell;\vec{\kappa};B}$}
\begin{defn}
$$
\aligned
\frak q^{\frak C}_{\ell;\vec{\kappa};B}(\text{\bf g},\text{\bf h}
)
&=
(-1)^{\maltese}(\text{\rm ev}^{\partial}_{0})_!
\left((\text{\rm ev}^+)^*\text{\bf g}
\wedge
(\widetilde{\text{\rm ev}} _{>0})^*\text{\bf h}
\right) \\
&
= (-1)^\maltese{}\text{\rm Corr}
([0,1] \times \mathcal M_{\ell}(\vec{\kappa},\vec p;B);(\widetilde{\text{\rm ev}}_{0};
\text{\rm ev}^+\times \widetilde{\text{\rm ev}}_{>0}))% \times 
(\text{\bf g}\wedge \text{\bf h} % \wedge \text{\bf h}'
)\\
&
\in {\frak C}(L_{\kappa_{0}},L_{\kappa_{K}};\F).
\endaligned$$
\end{defn}
Here the sign $\maltese$ is the same as (\ref{defqformula}).
Now we define:
\begin{defn}\label{defn78}
For given $\vec{\kappa}$ we define\index[syindex]{qcfrakellkappa@$\frak q^{\frak C}_{\ell;\vec{\kappa}}$} 
$$
\frak q^{\frak C,\frak b}_{\ell,\vec{\kappa}}
: 
\Omega(X)^{\otimes \ell}
\otimes 
B\frak CF(\cL^{\rm form};\vec{\kappa})
\to {\frak C}(L_{\kappa_{0}},L_{\kappa_{K}};\F) \widehat{\otimes} \Lambda_0,
$$
by
\begin{equation}\label{qcat2}
\frak q^{\frak C,\frak b}_{\ell}(\text{\bf g};\text{\bf h})
=
\sum_{B,\ell',\vec k}T^{\omega(B) }
\frac{\rho_{\frak b,\theta}(B)}{\ell'!} %\rho_{b}(B)
\frak q^{\frak C}_{\ell+\ell',\vec k;\vec{\kappa};B}
(\text{\bf g}\otimes \frak b_+^{\ell'};
\text{\bf h}).
\end{equation}
% where $b_+^{\frak C, \vec k} = \bigotimes_{i=0}^K (b_{\kappa_i,+}^{\frak C})^{\otimes k_i}$
% and $b_+^{\frak C, \vec k}$ is the solution  $b^{t,\frak C}_{\kappa, +} \in \Omega^{odd}(L) \otimes \Lambda_+$  of (\ref{ODEbounding}) 
% that we regard as an element of $\Omega([0,1] \times L) \widehat{\otimes} \Lambda_+$.
Other notations are the same as (\ref{qcat2}).
\par
We also define\index[syindex]{mxcfrakkappaB@$\frak m^{\frak C,\frak b}_{\vec{\kappa}}$}

$$
\frak m^{\frak C,\frak b}_{\vec{\kappa}}
: 
B\frak C(\cL^{\rm form};\vec{\kappa})
\to {\frak C}(L_{\kappa_{0}},L_{\kappa_{K}};\F) \widehat{\otimes} \Lambda_0,
$$
by
$$
\frak m^{\frak C,\frak b}_{\vec{\kappa}}(\text{\bf h})
= 
\frak q^{\frak C,\frak b}_{0}(1;\text{\bf h}).
$$
\end{defn}
Hereafter we omit $\frak b$ and just write $\frak m^{\frak C}_{\vec{\kappa}}$.
We take $\text{\bf h} = h_1 \otimes \cdots \otimes h_k$ so that $h_i$ is 
induced by the form on $CF(L_{\kappa_{i-1}},L_{\kappa_{i}};\F)$.
(Namely $h_i$ does not contain a $dt$ component and 
the coefficients are independent of $t$ component.)
We have\index[syindex]{mxcfrakkappat@$\frak m^{\frak C,t}_{\vec{\kappa}}$}\index[syindex]{c1Ctkbeta@$\frak c^{\frak C,\frak b,t}_{\vec{\kappa},\beta}$}
$$
\frak m^{\frak C,\frak b}_{\vec{\kappa}}(\text{\bf h})
= \frak m^{\frak C,\frak b,t}_{\vec{\kappa}}(\text{\bf h}) + 
 dt  \wedge \frak c^{\frak C,\frak b,t}_{\vec{\kappa}}(\text{\bf h}).
$$
By our choice they define $G_+$-gapped filtered $A_{n,k}$ structure.
Using Lemmas \ref{Kuraonepara} and \ref{existmultipolupara} it is straightforward to check 
that they define the required pseudo-isotopy.
The proof of Proposition \ref{existpisomain} is complete.
\end{proof}
\begin{rem}
$\frak q^{\frak C,\frak b}_{\ell;\vec{\kappa};B}$ is invariant under the action of the symmetric group 
on $\Omega(X)^{\otimes \ell}$. Therefore it induces a map
$$
E_{\ell}(\Omega(X))
\otimes 
B\frak C(\cL^{\rm form};\vec{\kappa})
\to {\frak C}(L_{\kappa_{0}},L_{\kappa_{K}};\F) \widehat{\otimes} \Lambda_0.
$$
\end{rem}

 %\marginpar{ The rest of this section is moved to latter sections.}

\section{The forgetful map and unitality.}
\label{sec:unit}
 %\marginpar{This section is new.}

As we mentioned several times, the filtered $A_{\infty}$
structure constructed in Section \ref{sec:cyclicfil} is not 
unital.
The purpose of this section is 
to describe the way to obtain a filtered homotopically
unital $A_{\infty}$
structure.

\subsection{Cyclicity of homotopically unital filtered $A_{\infty}$ category.}
The notion of homotopically unital curved filtered $A_{\infty}$ category
is defined in Definition \ref{defn84}.
We define cyclicity of such category below.

Let $\Cat$ be a filtered $A_{\infty}$ category and $c$  its object.
We define
$CH_k\Cat(c,c)$\index[syindex]{CHCkcc@$CH_k\Cat(c,c)$} as the set of finite sums of elements of the form
$
{\bf h} = h_0 \otimes \dots \otimes h_k.
$ 
Here  $c_0,\dots,c_k$ are objects with $c_0 = c$  and $h_i \in \Cat(c_i,c_{i+1})$, where $c_{k+1} = c_0$ by convention.
Let $\Cat^+(c,c) = \Cat(c,c) \oplus \Lambda_0[{\bf e}^+_c] \oplus \Lambda_0[{\bf f}_c]$
as in Definition \ref{defn84}.

\begin{defn}\label{def98}
A homotopically unital filtered $A_{\infty}$ category $(\Cat,\Cat^+,\iota)$ is said to be\index[syindex]{iota@$\Iota$} 
{\it cyclic}\index{cyclic homotopically unital filtered $A_{\infty}$ category}\index{cyclic} if the following holds.

\begin{enumerate}
\item
There is an inner product $\langle *,*\rangle : \Cat(c,c') \otimes \Cat(c',c)$ as in Definition \ref{defn:cyclic_structure}.\footnote{We remark that 
inner product is not defined for ${\bf e}^+_c$, ${\bf f}_c$.}
\item
There exists an operator\index[syindex]{mxzk+@${}^+\frak m_k$}
$$
{}^+\frak m_k : CH_k\Cat^+(c,c) \to \Lambda_0.
$$
such that:
\begin{enumerate}
\item 
\begin{equation}\label{qismcatprimeaa}
{}^{+}\frak m_{k}
(h_0,h_1,\dots,h_{k-1},h_k) = (-1)^\maltese
{}^{+}\frak m_{k}
(h_1,\dots,h_k,h_0),
\end{equation}
where 
$\maltese = (\deg' h_0)(\deg' h_1 + \dots + \deg' h_k$). 
\item
In case $h_0$ does not contain ${\bf e}^+_{c}$ or ${\bf f}_{c}$,
we have:
\begin{equation}\label{form9533}
{}^{+}{\frak m}_{k}(h_1,\dots,h_k,h_0)
= 
\langle {\frak m}_{k}(h_1,\dots,h_k),h_0  \rangle.
\end{equation}
Here ${\frak m}_{k}$ is the structure operation of $\Cat$.
There are exceptions when ${\frak m}_{k}(h_1,\dots,h_k)$ 
contais ${\bf f}_c$ or ${\bf e}_c^+$.
(In such cases the right hand side of (\ref{form9533}) is not defined.)
In such cases  we require:\footnote{Note 
${}^{+}{\frak m}_{1}({\bf e}^+_c,h_0) = 0$ follows from (\ref{form9533}).}
\begin{equation}\label{fmnew103}
\aligned
{}^{+}{\frak m}_{1}({\bf f}_c,h_0)
&= 
\langle {\frak m}_{1}({\bf f}_c) - {\bf e}_c^+ ,h_0  \rangle,
\\
{}^{+}{\frak m}_{2}({\bf e}_c^+,{\bf e}_c^+,h_0) &= {}^{+}{\frak m}_{2}({\bf e}_c^+,{\bf f}_c,h_0) ={}^{+}{\frak m}_{2}({\bf f}_c,{\bf e}_c^+,h_0) = 0.
\endaligned
\end{equation}
\item
${}^{+}{\frak m}_{k}(\dots,{\bf e}^+_{c},\dots) = 0$
when all the inputs are either ${\bf e}^+_{c}$ or ${\bf f}_{c}$.
\item For $x_0 \otimes {\bf x} \in CH\Cat^+(c,c)$   
\begin{equation}\label{form9533+1}
\sum_c (-1)^{\deg' x_0 + \deg'  {\bf x}_c^{3;1}}\,\, {}^{+}{\frak m}_*(x_0 \otimes {\bf x}_c^{3;1} \otimes \frak m_*({\bf x}_c^{3;2}), {\bf x}_c^{3;3})
= 0
\end{equation}
in case $x_0$ does not include ${\bf f}_c$.  In that case (\ref{form9533+1}) becomes 
\begin{equation}\label{form9533+2}
\aligned
\sum_c (-1)^{1 + \deg'  {\bf x}_c^{3;1}}\,\, {}^{+}{\frak m}_*({\bf f}_c \otimes {\bf x}_c^{3;1} \otimes &\frak m_*({\bf x}_c^{3;2}), {\bf x}_c^{3;3})\\
& 
+  {}^{+}{\frak m}_*({\bf e}_c^+,{\bf x}) =  0.
\endaligned
\end{equation}
(Compare (\ref{fmnew103}).)\footnote{The extra term in (\ref{form9533+2}) is mostly zero but
$
{}^{+}{\frak m}_2({\bf e}_c^+,x_1,x_2) = \langle \frak m_2({\bf e}_c^+,x_1),x_2 \rangle
$
may be nonzero.}
\end{enumerate}

\end{enumerate}
\end{defn}
In case at least one input is not ${\bf e}_c^+$ or ${\bf f}_c$, Item (d) follows from 
(\ref{qismcatprimeaa}), (\ref{form9533}) and the $A_{\infty}$ relation of $\frak m$.

\subsection{Forgetable and unforgetable marked points.}
\label{subsec:forgetunforget}

We use the notation of Subsection \ref{cycliccategorymoduli} and
consider the moduli space 
${\mathcal M}_{\ell;\vec k}((\vec{\kappa},\vec p);B)$
in Proposition \ref{Kuraeistspoly}.

We recall that for a category to be homotopically unital we add two additional generators ${\bf e}^+_c$ and 
${\bf f}_c$ for each object $c$. 
 The former will be the (strict) unit and the latter will be a `homotopy' between 
${\bf e}^+_c$ and the original (homotopy) unit ${\bf e}_c$. We introduce two auxiliary data  ${\frak f}_{\rm{\bf e}}$ and 
$\frak f_{\rm{\bf f}}$ for this purpose.
Detailed definitions thereof are now in order.

We start with  ${\frak f}_{\rm{\bf e}}$. It is a map
which assigns a subset $\frak f_{{\rm{\bf e}},i}$ of $\{1,\dots,k_i\}$ to each $i=0,\dots,K$.
(Recall that $\{1,\dots,k_i\}$ is the index set of the diagonal marked points $\vec w_i = 
(w_{i,1},\dots,w_{i,k_i})$.)
We use marked points labeled by ${\frak f}_{\rm{\bf e}}$
to define operations including ${\bf e}^+_c$. \index[syindex]{ffrake@${\frak f}_{\rm{\bf e}}$}

To include the degree $-1$ element ${\bf f}_c$ in our discussion, 
we use another forgetability data ${\frak f}_{\bf f}$.\index[syindex]{ffrakf@${\frak f}_{\rm{\bf f}}$}
We will write
\begin{equation}\label{not(9.1)}
{\frak f}_{\bf f} = \{(i,j)  \in \Z_{\geq 0} \times \Z_{>0} \mid j \in {\frak f}_{{\bf f},i}\}.
\end{equation}
We have already used a similar notation for ${\frak f}_{\bf e} = \{\frak f_{{\bf e},i}\} $
above for which we put
\begin{equation}\label{not(9.12)}
{\frak f}_{\bf e} =  \{(i,k) \in \Z_{\geq 0} \times \Z_{> 0}\mid k \in \frak f_{{\bf e},i}\},
\end{equation}
and $\vec{\frak f} = ({\frak f}_{\bf e},{\frak f}_{\bf f})$. 
\begin{defn}
We define %\marginpar{Two definitions are added. KF 2025 Sep.}
$$
{\mathcal M}_{\ell;\vec k}((\vec{\kappa},\vec p);B;{\frak f}_{\bf e}) = 
{\mathcal M}_{\ell;\vec k}((\vec{\kappa},\vec p);B;\vec{\frak f}) =
{\mathcal M}_{\ell;\vec k}((\vec{\kappa},\vec p);B).
$$
\end{defn}

The moduli spaces \index[syindex]{MellkkappapBfrakfe-@${\mathcal M}_{\ell;\vec k}((\vec{\kappa},\vec p);B;{\frak f}_{\rm{\bf e}})$}
\index[syindex]{MellkkappapBfrakfeff-@${\mathcal M}_{\ell;\vec k}((\vec{\kappa},\vec p);B;\vec{\frak f})$}
${\mathcal M}_{\ell;\vec k}((\vec{\kappa},\vec p);B;{\frak f}_{\bf e})$, 
${\mathcal M}_{\ell;\vec k}((\vec{\kappa},\vec p);B;\vec{\frak f})$
are actually the same as ${\mathcal M}_{\ell;\vec k}((\vec{\kappa},\vec p);B)$.
We however use different Kuranishi structures and CF-perturbations on those moduli spaces.
\begin{defn}
$$
{\mathcal M}_{\ell;\vec k}((\vec{\kappa},\vec p,m);B;{\frak f}_{\bf e}) = 
{\mathcal M}_{\ell;\vec k}((\vec{\kappa},\vec p,m);B;\vec{\frak f}) =
{\mathcal M}_{\ell;\vec k}((\vec{\kappa},\vec p,m);B).
$$
When we write ${\mathcal M}_{\ell;\vec k}((\vec{\kappa},\vec p,m);B;{\frak f}_{\bf e})$ 
or ${\mathcal M}_{\ell;\vec k}((\vec{\kappa},\vec p,m);B;\vec{\frak f})$, 
we require the $0$-th boundary marked point $z(0)$ determined from $m$ in the same way as 
Definition \ref{defn917}
is not a marked point labeled by ${\frak f}_{\bf e}$ or by ${\frak f}_{\bf f}$ unless all the 
boundary marked points are labeled by ${\frak f}_{\bf e}$ or  ${\frak f}_{\bf f}$. 
\par
We then put\index[syindex]{MellkkappapBfrakfenove@${\mathcal M}_{\ell}(\vec{\kappa},\vec p;B;{\frak f}_{\rm{\bf e}})$}
\index[syindex]{MellkkappapBfrakfeffnove@${\mathcal M}_{\ell}(\vec{\kappa},\vec p;B;\vec{\frak f})$}
$$
{\mathcal M}_{\ell}(\vec{\kappa},\vec p;B;{\frak f}_{\bf e}) = 
{\mathcal M}_{\ell}(\vec{\kappa},\vec p;B;\vec{\frak f}) =
{\mathcal M}_{\ell}(\vec{\kappa},\vec p;B).
$$
\end{defn}

\begin{defn}
A diagonal marked point $w_{i,j}$ of an element of 
${\mathcal M}_{\ell;\vec k}((\vec{\kappa},\vec p);B;{\frak f}_{\bf e})$
is said to be {\it forgetable}\index{forgetable}   if $j \in  \frak f_{{\rm{\bf e}},i}$.
Otherwise it is said to be {\it unforgetable.}\index{unforgetable}
Switching marked points $z_i$ 
are, by definition, unforgetable.
We call $\frak f_{\rm{\bf e}}$, {\it forgetability data}.\index{forgetability data}
\end{defn}

We will pull back the $0$-form $1$ by the evaluation maps at the marked points 
labeled by $\frak f_{\rm{\bf e}}$. 
We also require various data and perturbations to be compatible with the 
forgetful maps of the marked points labeled by $\frak f_{\rm{\bf e}}$. 
Then we will have the strict unitality.
We will carry out this process precisely below.

\begin{center}
\begin{figure}[h]
 \begin{tabular}{cc}
 \begin{minipage}[t]{0.5\hsize}
\centering
\includegraphics[scale=0.3]
{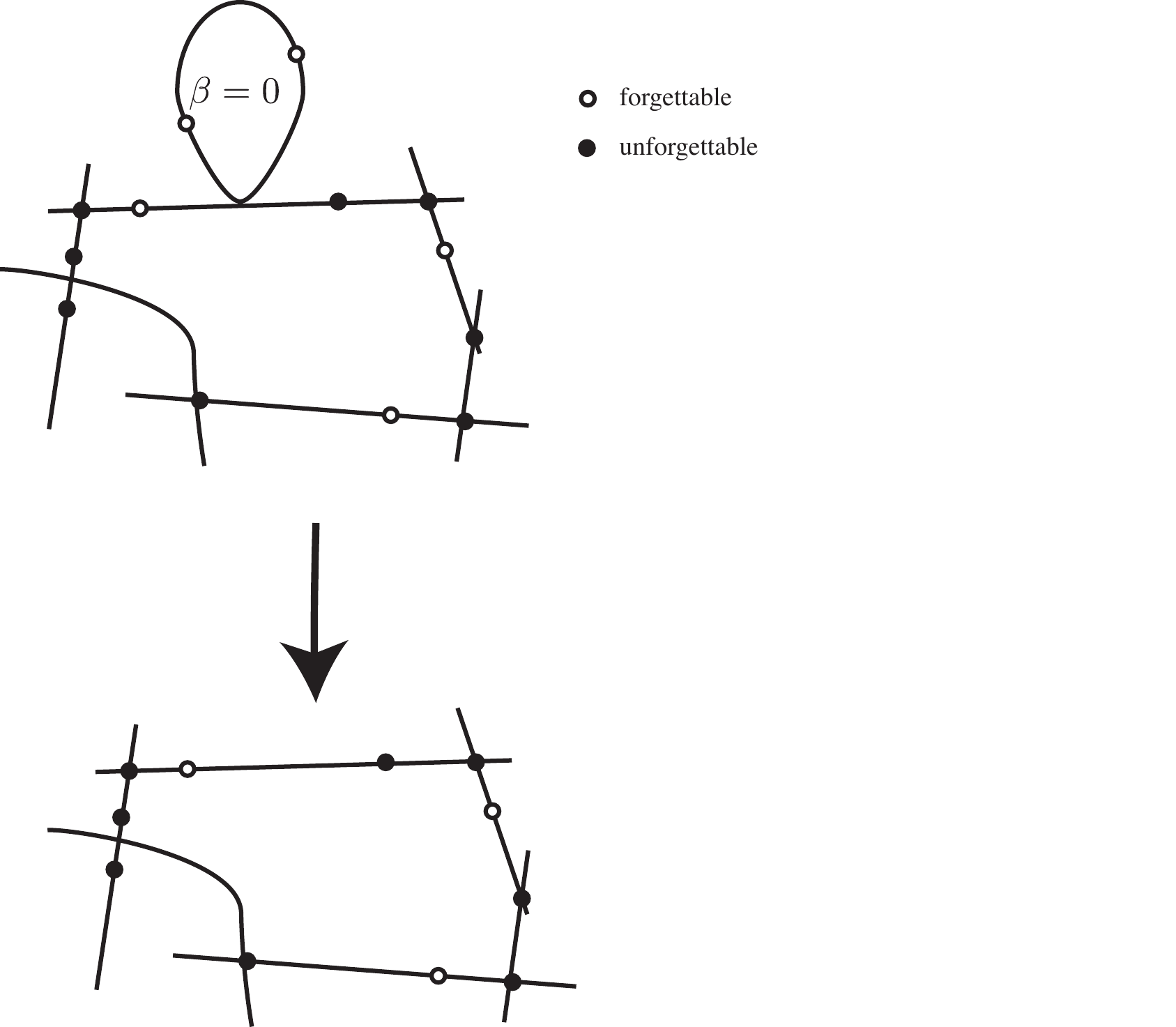}
\caption{(tricom1)}
\label{Figuretrocom1}
\end{minipage} &
 \begin{minipage}[t]{0.5\hsize}
\centering
\includegraphics[scale=0.35]
{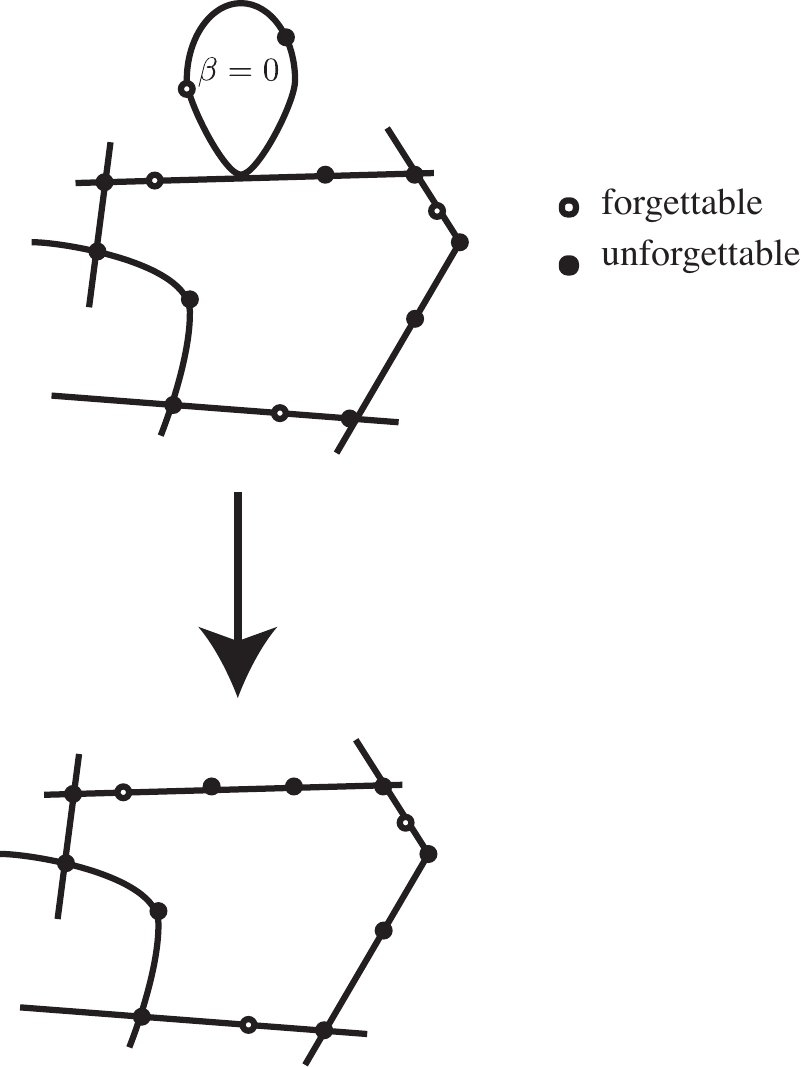}
\caption{(tricom2)}
\label{Figuretrocom2}
\end{minipage}
\end{tabular}
\end{figure}
\end{center}

We consider three types of the boundary 
$(\ref{218})$, $(\ref{219})$, $(\ref{220})$ 
of ${\mathcal M}_{\ell;\vec k}((\vec{\kappa},\vec p);B)$.
We will define the forgetability data $\frak f^1_{\rm{\bf e}}$, $\frak f^2_{\rm{\bf e}}$
of the factors of them when we regard the boundary as that  
of ${\mathcal M}_{\ell;\vec k}((\vec{\kappa},\vec p);B;\frak f_{\rm{\bf e}})$.
We remark that in the factors appearing in the formulas
$(\ref{218})$, $(\ref{219})$, $(\ref{220})$, 
there are two extra boundary marked (or one nodal) 
points which are not  marked points
as an element of ${\mathcal M}_{\ell;\vec k}((\vec{\kappa},\vec p);B;\frak f_{\rm{\bf e}})$. 
For the marked points of the factors of $(\ref{218})$, $(\ref{219})$, $(\ref{220})$ 
other than those two marked points
we define their forgetability data induced from
$\frak f_{\rm{\bf e}}$.
For the two new marked points 
we regard them as  unforgetable unless we are in case $(\ref{218})$ and the following 
situation happens.\footnote{In cases $(\ref{219})$ and $(\ref{220})$, two new marked points are by definition 
unforgetable.} See Figures \ref{Figuretrocom1}, \ref{Figuretrocom2}.
\begin{enumerate}
\item[(tricom1)]
All the marked points on the first factor are forgetable 
and the homology class $\beta$ is $0$.
\item[(tricom2)]
All the marked points on the first factor are forgetable 
except one marked point
and the homology class $\beta$ is $0$.
\end{enumerate}
We say the irreducible component of the first factor is a {\it trivial component}\index{trivial component} if (tricom1) or (tricom2) happens.
In the situation (tricom1) we declare the new marked point on the second factor to be 
forgetable. 
(Note that the object of the first factor is the one 
associated to a {\it single embedded}
Lagrangian appearing in Section \ref{ainfalgasssingle}. 
So we do not need to distinguish forgetable and unforgetable marked points 
for the marked points on this factor.)
In the situation (tricom2) we declare the new marked point on the second factor to be 
unforgetable.
We will discuss the way to handle this {\it trivial component} of the first factor in Proposition \ref{prop92} (3) below.
\par
We have thus defined the forgetability data $\frak f^1_{\rm{\bf e}}$ of the first factor and $\frak f
^2_{\rm{\bf e}}$ of the second factor of $(\ref{218})$, $(\ref{219})$, $(\ref{220})$.
We call them the {\it induced forgetability data}.\index{induced forgetability data}

\begin{prop}\label{prop92}
For given $E_0,k_0,\ell_0$, the following holds for ${\mathcal M}_{\ell;\vec k}((\vec{\kappa},\vec p);B;{\frak f}_{\rm{\bf e}})^{\boxplus 1}$
with $(B,\vert\vec{\kappa}\vert,\ell) \le (E_0,k_0,\ell_0)$.\index[syindex]{MellkkappapBfrakfe@${\mathcal M}_{\ell;\vec k}((\vec{\kappa},\vec p);B;{\frak f}_{\rm{\bf e}})^{\boxplus 1}$}
\begin{enumerate}
\item The same as Proposition $\ref{Kuraeistspoly}$ $(1)$.
\item The same as Proposition $\ref{Kuraeistspoly}$ $(2)$.
\item The boundary of ${\mathcal M}_{\ell;\vec k}((\vec{\kappa},\vec p);B;\frak f_{\bf e})^{\boxplus 1}$
is a union of the three types of fiber or direct products 
$(\ref{218})$, $(\ref{219})$, $(\ref{220})$ with the Kuranishi structure 
via induced forgetability data $\frak f^1_{\bf e}$ and $\frak f^2_{\bf e}$. 
Here the union is taken over the data described above 
together with the shuffles $(\mathbb L'',\mathbb L')$ of $\underline\ell$
such that $\# \mathbb L'' = \ell''$, $\# \mathbb L' = \ell'$.
In case of the trivial component of the first factor of $(\ref{218})$, we take 
the Kuranishi structure whose obstruction bundle is $0$.
\item The same as Proposition $\ref{Kuraeistspoly}$ $(4)$.
\item If $w_{i,j}$ is unforgetable then the  evaluation map
$$
\text{\rm ev}_{i,j} : 
{\mathcal M}_{\ell;\vec k}((\vec{\kappa},\vec p);B;{\frak f}_{\rm{\bf e}})^{\boxplus 1}
\to L_{\kappa_i}
$$
is weakly submersive.
\item The same as Proposition $\ref{Kuraeistspoly}$ $(6)$. (In the case of this proposition, cyclic symmetry changes the 
forgetability data ${\frak f}_{\rm{\bf e}}$.) 
\item The same as Proposition $\ref{Kuraeistspoly}$ $(7)$.
\item
The Kuranishi structure is compatible with the forgetful map\index{forgetful map}\index[syindex]{forgetfrak@$\frak{forget}$}
\begin{equation}\label{form106new}
\frak{forget} : {\mathcal M}_{\ell;\vec k }((\vec{\kappa},\vec p);B;{\frak f}_{\rm{\bf e}})^{\boxplus 1}
\to {\mathcal M}_{\ell;\vec k\setminus {\frak f}_{\bf e}}((\vec{\kappa},\vec p);B;\emptyset)^{\boxplus 1}
\end{equation}
which forgets the marked points $w_{i,j}$ labeled by ${\frak f}_{\bf e}$.\footnote{See Definition \ref{def1221} for the definition of compatibility.}
%\item
%The Kuranishi structure is compatible the permutation  
%of the interior marked points.
%\item
%The Kuranishi structure is compatible the cyclic permutation   
%of the boundary marked points.
\item
In case $K=0$, which is the case when there is only one Lagrangian submanifold,
the Kuranishi structure coincides with the one in Proposition $\ref{diskkura}$.
\item If ${\frak f}_{\rm{\bf e}}$ is empty then the Kuranishi structure coincides with the one 
of Proposition $\ref{Kuraeistspoly}$.
\item The description of the boundary is compatible with orientation in the sense of \cite[(21.15)]{springer}.
\end{enumerate}
The compactification ${\mathcal M}_{\ell,\vec k}((\vec{\kappa},\vec p,m);B;{\frak f}_{\rm{\bf e}})^{\boxplus 1}$ of 
$\overset{\circ}{\mathcal M}_{\ell,\vec k}((\vec{\kappa},\vec p,m);B;{\frak f}_{\rm{\bf e}})^{\boxplus 1}$
is similar.
We require $0$-th marked point is not labeled by ${\frak f}_{\rm{\bf e}}$.
The Kuranishi structures are independent of the choice of $0$-th marked point, as far as it is {\it not} labeled by  ${\frak f}_{\rm{\bf e}}$.
\end{prop}

\begin{rem}\label{rem106}
We elaborate on the cyclic symmetry of the Kuranishi structure of the 
moduli space ${\mathcal M}_{\ell,\vec k}((\vec{\kappa},\vec p,m);B;{\frak f}_{\rm{\bf e}})^{\boxplus 1}$.
It is the cyclic symmetry of the permutation of marked points 
which are unforgetable, that is, not labeled by ${\frak f}_{\rm{\bf e}}$.  This is the requirement in case 
there exists at least one unforgetable marked points.
In case all the marked points are forgetable, the space
of ${\mathcal M}_{\ell,\vec k}((\vec{\kappa},\vec p,m);B;{\frak f}_{\rm{\bf e}})^{\boxplus 1}$ is not 
defined.
In this case the cyclic symmetry  
 is a consequence of compatibility of the Kuranishi structure with the
forgetful map of  ${\mathcal M}_{\ell,\vec k}((\vec{\kappa},\vec p);B;{\frak f}_{\rm{\bf e}})^{\boxplus 1}$
(Proposition \ref{prop92} (6)).
\end{rem}

Item (8) implies unitality with respect to the input at the forgetable marked points.
The notion of homotopical unitality is a suitable algebraic framework to use 
forgetable and unforgetable marked points.
To work it out we need to introduce  additional 
Kuranishi structures which interpolate among them.
We here follow the  method  used in \cite[Section 7.3]{fooo092}.\footnote{
In \cite{fooo092} the singular chain model was used. Here we use the 
de Rham model.}

We consider ${\frak f}_{\rm{\bf f}}$ which is disjoint from $\frak f_{\rm{\bf e}}$.
We put $\vert{\frak f}_{\rm{\bf f}}\vert = \sum \vert{\frak f}_{{\rm{\bf f}},i}\vert \in \Z_{\ge 0}$
and consider 
$
[0,1]^{\vert{\frak f}_{\rm{\bf f}}\vert}
$.
If $\vec s \in [0,1]^{{\frak f}_{\rm{\bf f}}\vert}$ is given and $w_{i,j}$
is a marked point  in ${\frak f}_{\rm{\bf f}}$, we denote by 
$s_{i,j}$ the corresponding factor of $\vec s$.

We consider the situation of $(\ref{218})$, $(\ref{219})$, $(\ref{220})$.
We will define induced forgetability data $\frak f^1_{\rm{\bf e}}$ , $\frak f^2_{\rm{\bf e}}$ 
and $\frak f^1_{\rm{\bf f}}$ , $\frak f^2_{\rm{\bf f}}$  of the 
factors.
 $\frak f^1_{\rm{\bf f}}$ , $\frak f^2_{\rm{\bf f}}$ are defined in a similar way as the case of Proposition \ref{prop92}
 as follows.
 (tricom1) and (tricom2) are as before. 
 (In those cases we require the forgetable marked points 
 are not in $\frak f_{\rm{\bf f}}$.)
 We have one more case below:
 \begin{enumerate}
\item[(tricom3)] 
We are in Case $(\ref{218})$. All the marked points on the first factor are forgetable 
and in ${\frak f}_{\rm{\bf e}}$
except one marked point which is in  $\frak f_{\rm{\bf f}}$.
The homology class $\beta$ is $0$.
 \end{enumerate}
 \begin{figure}[h] 
\centering \hskip2.5cm
\includegraphics[scale=0.4]{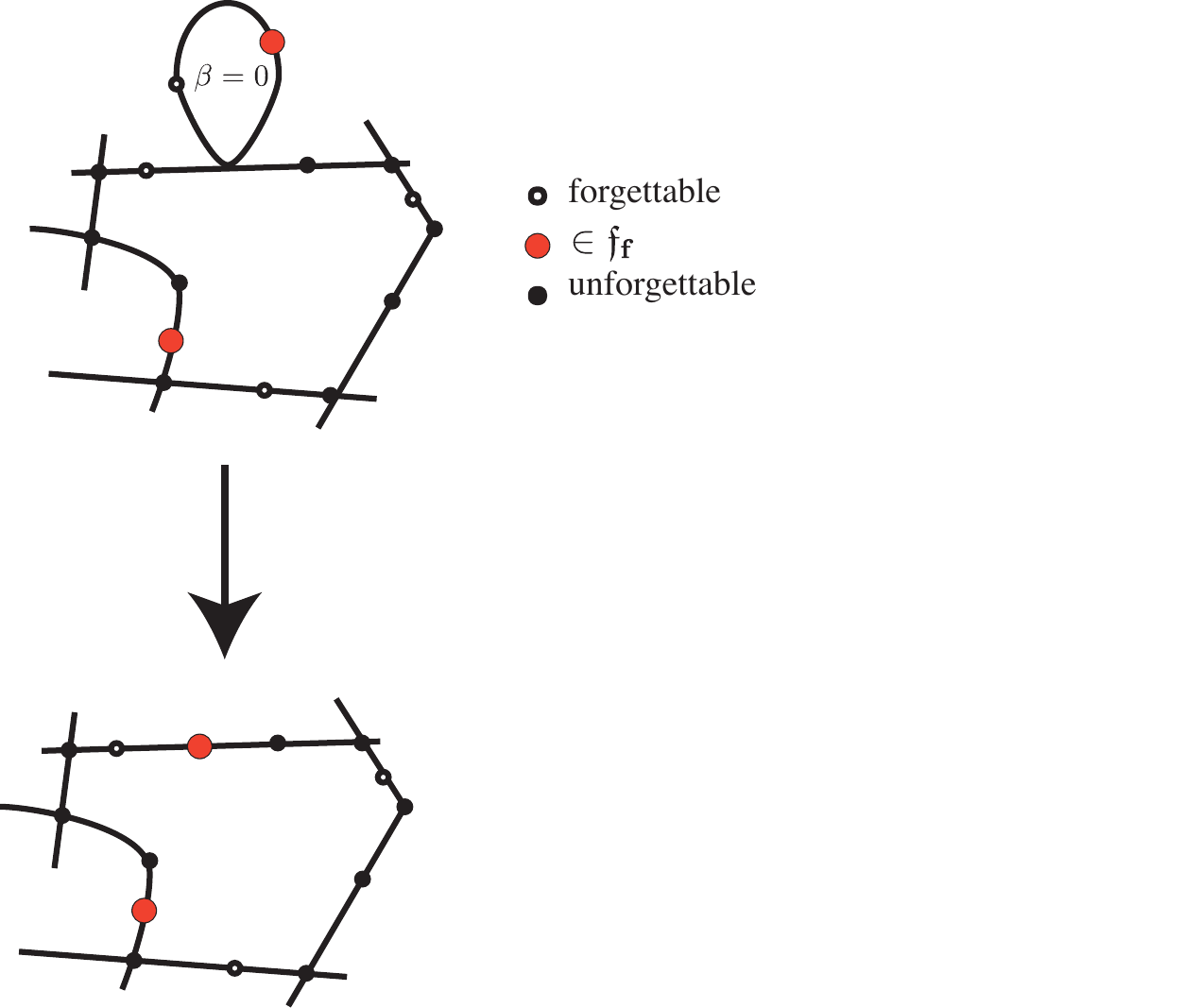}
\caption{(tricom3)}
\label{Figure(tricom3).1}
\end{figure}
In the case (tricom1) (resp. (tricom2)) we regard the new 
marked point on the second factor as a forgetable marked point in ${\frak f}_{\rm{\bf e}}^2$
(resp. unforgetable marked point).
(Recall that we require that the marked points in the first factor 
are all in ${\frak f}_{\rm{\bf e}}^1$ for (tricom1) and 
is in ${\frak f}_{\rm{\bf e}}^1$ except one for (tricom2).)
In the case (tricom3) we regard the new marked point 
of the second factor to be an element of $\frak f^2_{\rm{\bf f}}$.
\par
Except the above marked points we set $\frak f^1_{\rm{\bf f}}$ , $\frak f^2_{\rm{\bf f}}$ to be the marked points of $\frak f_{\rm{\bf f}}$ in 
each of the components.

\begin{prop}\label{prop94}
We put $\vec{\frak f} = ({\frak f}_{\rm{\bf e}},{\frak f}_{\rm{\bf f}})$.\index[syindex]{ffrakvec@$\vec{\frak f}$}
The following holds for $[0,1]^{\vert {\frak f}_{\rm{\bf f}}\vert} \times {\mathcal M}_{\ell;\vec k}((\vec{\kappa},\vec p);B;\vec{\frak f})^{\boxplus 1}
$.\index[syindex]{01MellkkappapBfrakff@$[0,1]^{\vert {\frak f}_{\rm{\bf f}}\vert} \times {\mathcal M}_{\ell;\vec k}((\vec{\kappa},\vec p);B;\vec{\frak f})^{\boxplus 1}
$}

\begin{enumerate}
\item This space is compact, Hausdorff and metrizable.
\item It admits a Kuranishi structure.
\item 
Its normalized boundary is decomposed as 
$$
\aligned
\partial([0,1]^{\vert{\frak f}_{\rm{\bf f}}\vert} \times  {\mathcal M}_{\ell;\vec k}((\vec{\kappa},\vec p);B;\vec{\frak f})^{\boxplus 1})
=&
[0,1]^{\vert{\frak f}_{\rm{\bf f}}\vert} \times \partial({\mathcal M}_{\ell;\vec k}((\vec{\kappa},\vec p);B;\vec{\frak f})^{\boxplus 1})  \\
&\cup 
\partial([0,1]^{\vert{\frak f}_{\rm{\bf f}}\vert})  \times  {\mathcal M}_{\ell;\vec k}((\vec{\kappa},\vec p);B;\vec{\frak f})^{\boxplus 1}.
\endaligned
$$ 
The following compatibility condition holds for their Kuranishi structures.
\begin{enumerate}
\item
The second factor of the first summand is decomposed as 
 a union of the three types of fiber or direct products 
$(\ref{218})$, $(\ref{219})$, $(\ref{220})$.
We consider the induced forgetability data $\frak f^1_{\rm{\bf e}}$ , $\frak f^2_{\rm{\bf e}}$ and 
$\frak f^1_{\rm{\bf f}}$ , $\frak f^2_{\rm{\bf f}}$ as above.
Put $\vec{\frak f}^1 = (\frak f^1_{\rm{\bf e}},\frak f^1_{\rm{\bf f}})$, $\vec{\frak f}^2= (\frak f^2_{\rm{\bf e}},\frak f^2_{\rm{\bf f}})$.
Then in the case of $(\ref{218})$ this summand is
\begin{equation}\label{218+I}
\aligned
&( [0,1]^{\vert \frak f^1_{\rm{\bf f}} \vert} \times  {\mathcal M}_{\ell'';k''_i+1}(L_{\kappa_i};\beta; \vec{\frak f}^1)^{\boxplus 1})\\
&\,{}_{\text{\rm ev}_0} \times_{\text{\rm ev}_{i,j}}
([0,1]^{\vert \frak f^2_{\rm{\bf f}} \vert}  \times {\mathcal M}_{\ell';\vec k'}((\vec{\kappa},\vec p);B'; \vec{\frak f}^2)^{\boxplus 1}).
\endaligned
\end{equation}
The fiber product is well defined by Item $(5)$.  Note that by definition the $0$-th marked point is always unforgetable.

We require that the restriction of the Kuranishi structure 
to this part of the boundary coincides with the 
fiber product Kuranishi structure.
In the cases {\rm (tricom1), (tricom2), (tricom3)}, 
the Kuranishi structure of the first factor is one whose obstruction bundle is 
trivial.
\item
The  other two   cases of the first summand are similar.
\item
The case of the second summand is described in items $(9)$ and $(10)$.
\end{enumerate}
\item The same as Proposition $\ref{Kuraeistspoly}$  $(4)$.
\item If $w_{i,j}$ is unforgetable then the  evaluation map
$$
\text{\rm ev}_{i,j} : 
[0,1]^{\vert \frak f_{\rm{\bf f}} \vert} \times {\mathcal M}_{\ell;\vec k}((\vec{\kappa},\vec p);B;\vec{\frak f})^{\boxplus 1}
\to L_{\kappa_i}
$$
is weakly submersive.
\item The same as Proposition $\ref{Kuraeistspoly}$ $(6)$. 
\item The same as Proposition $\ref{Kuraeistspoly}$ $(7)$. (In the case of this proposition, cyclic symmetry changes the 
forgetability data ${\frak f}_{\rm{\bf e}}$, ${\frak f}_{\rm{\bf f}}$.) 
\item
The Kuranishi structure is compatible with the forgetful map
$$
\frak{forget} : [0,1]^{\vert\frak f_{\bf f}\vert} \times {\mathcal M}_{\ell;\vec k}((\vec{\kappa},\vec p);B;\vec{\frak f})^{\boxplus 1} 
\to {\mathcal M}_{\ell;\vec k\setminus {\frak f}_{\bf e}}((\vec{\kappa},\vec p);B;(\emptyset, {\frak f}_{\bf f}))^{\boxplus 1}
$$
which forgets the marked points $w_{i,j}$ labeled by $\frak f_{\bf e}$.\footnote{See Definition \ref{def1221}  for the definition of compatibility.}
\item 
Let $(i,j) \in{\frak f}_{\rm{\bf f}}$. (We use the notation $(\ref{not(9.1)})$ here.) We consider the subset of 
$[0,1]^{\vert{\frak f}_{\rm{\bf f}}\vert}\times {\mathcal M}_{\ell;\vec k}((\vec{\kappa},\vec p);B;\vec{\frak f})^{\boxplus 1}$ where $s_{i,j} = 1$.
We consider ${\frak f}'_{\rm{\bf e}} = {\frak f}_{\rm{\bf e}} \cup (i,j)$, ${\frak f}'_{\rm{\bf f}} = {\frak f}_{\rm{\bf f}} \setminus (i,j)$,
$\vec{\frak f}' = ({\frak f}'_{\rm{\bf e}},{\frak f}'_{\rm{\bf f}})$
and identify this subset with 
$[0,1]^{\vert{\frak f}_{\rm{\bf f}}'\vert} \times {\mathcal M}_{\ell;\vec k}((\vec{\kappa},\vec p);B;\vec{\frak f}')^{\boxplus 1}
$.
Then the Kuranishi structure 
obtained by restriction coincides with 
one of $[0,1]^{\vert{\frak f}_{\rm{\bf f}}'\vert} \times  {\mathcal M}_{\ell;\vec k}((\vec{\kappa},\vec p);B;\vec{\frak f}')^{\boxplus 1}
$.\footnote{In other words 
at $s_{i,j} = 1$ the Kuranishi structure coincides with one 
where we regard $w_{i,j}$ as a forgetable marked point.}
\item
Let $(i,j) \in {\frak f}_{\rm{\bf f}}$. We consider the subset of 
$[0,1]^{\vert{\frak f}_{\rm{\bf f}}\vert} \times {\mathcal M}_{\ell;\vec k}((\vec{\kappa},\vec p);B;\vec{\frak f})^{\boxplus 1}$ where $s_{i,j} = 0$.
We consider $\frak f'_{{\rm{\bf f}}} =  {\frak f}_{\rm{\bf f}} \setminus (i,j)$ and 
$\vec{\frak f}'' = (\frak f_{\bf e},\frak f'_{{\rm{\bf f}}})$,
and identify this subset with 
$[0,1]^{\vert {\frak f}_{\rm{\bf f}}'\vert} \times {\mathcal M}_{\ell;\vec k}((\vec{\kappa},\vec p);B;\vec{\frak f}'')^{\boxplus 1}
$.
Then the Kuranishi structure 
obtained by restriction coincides with $[0,1]^{\vert{\frak f}_{\rm{\bf f}}'\vert}  \times {\mathcal M}_{\ell;\vec k}((\vec{\kappa},\vec p);B;\vec{\frak f}'')^{\boxplus 1}  
$.\footnote{In other words 
at $s_{i,j} = 0$ the Kuranishi structure coincides with one 
where we regard $w_{i,j}$ as an unforgetable marked point.}
\item In the case when $K=0$, that is, there is only one 
Lagrangian submanifold involved, 
the Kuranishi structure is the direct product 
of the one in Proposition $\ref{diskkura}$ and $[0,1]^{\vert{\frak f}_{\bf f}\vert}$ with trivial Kuranishi structure.
\item If $\frak f_{\bf f}$ is empty the Kuranishi structure coincides with the one in Proposition $\ref{prop92}$.
\end{enumerate}
The compactification of 
${\mathcal M}_{\ell;\vec k}((\vec{\kappa},\vec p,m);B;{\frak f}_{\rm{\bf e}},{\frak f}_{\rm{\bf f}})^{\boxplus 1}$
is similar.
\end{prop}
\begin{rem}
This is a continuation of Remark \ref{rem106}.
In our situation we need to consider 
${\mathcal M}_{\ell;\vec k}((\vec{\kappa},\vec p);B;\vec{\frak f})^{\boxplus 1}$
such that all the marked points are labeled by ${\frak f}_{\rm{\bf f}}$, for the proof 
of Definition \ref{def98} (d). 
In that case we require that the Kuranishi structure is invariant under the cyclic 
permutation of such marked points.
\end{rem}

The proof is similar to the proof of Proposition \ref{prop92} which is proved in Subection \ref{polygonforget}.

Now we can find a system of CF-perturbations of those moduli spaces 
with Kuranishi structures which satisfy the appropriate compatibility 
with various additional data as follows.
\begin{prop}\label{prop95}
Suppose we are in the situation of Proposition $\ref{prop92}$.
Then for any $E_0,k_0,\ell_0$ there exists a system of 
CF-perturbations on 
${\mathcal M}_{\ell;\vec k}((\vec{\kappa},\vec p);B;{\frak f}_{\bf e})^{\boxplus 1}$
with $(B,\vert\vec{\kappa}\vert,\ell) \le (B_0,k_0,\ell_0)$, such that the following holds.
%\footnote{
%%The second and fourth equalities 
%%mean that the number of unforgetable marked points is bounded.
%\bf We might simply remove this footnote.}
\begin{enumerate}
\item
The CF-perturbations are transversal to $0$.
\item
They are compatible with the description of the boundary in
Proposition $\ref{prop92}$ $(3)$.
\item
If $w_{i,j}$ is unforgetable then 
the evaluation map in Proposition $\ref{prop92}$ $(5)$ is 
strongly submersive with respect to our CF-perturbations.
\item The symmetry in Proposition $\ref{prop92}$ $(6)$, $(7)$ preserves 
the CF-perturbations.
\item In the situation of Proposition $\ref{prop92}$
$(8)$ the CF-perturbations are compatible with the 
forgetful maps.
\item
In the case of Item $(9)$ of  Proposition $\ref{prop92}$ 
the CF-perturbations coincide with the ones in Proposition $\ref{existmkulti1}$.
\item
If ${\frak f}_{\bf e}$ is empty then the CF-perturbations coincide with the ones
of Proposition $\ref{existmultipolu}$.
\end{enumerate}
The ${\mathcal M}_{\ell,\vec k}((\vec{\kappa},\vec p,m);B;{\frak f}_{\rm{\bf e}})^{\boxplus 1}$ version also holds.
\end{prop}
Once the compatible choice of Kuranishi structures is given as in Proposition \ref{prop94},
the proof of the construction of the CF-perturbation of Proposition \ref{prop95} is similar to those in the literature such as \cite[Section 22.1]{springer}.
We explain the point concerning the forgetful map on the outer collaring in Section \ref{sec:CRperturb}.
\begin{prop}\label{exstCFUNI}
Suppose we are in the situation of Proposition $\ref{prop94}$.
Then for any $E_0,k_0,\ell_0$ there exists a system of 
CF-perturbations on 
$ [0,1]^{\vert {\frak f}_{\rm{\bf f}}\vert} \times {\mathcal M}_{\ell;\vec k}((\vec{\kappa},\vec p);B;\vec{\frak f})^{\boxplus 1}$
with $(B,\vert\vec{\kappa}\vert,\ell) \le (E_0,k_0,\ell_0)$, that has
the following properties.\footnote{We apply outer collaring only to ${\mathcal M}_{\ell;\vec k}((\vec{\kappa},\vec p);B;\vec{\frak f})$
factors.  Actually for $[0,1]^{\vert {\frak f}_{\rm{\bf f}}\vert}$ factor we can easily make the Kuranishi structure and CF-perturbation to be 
collared.}
\begin{enumerate}
\item
The CF-perturbations are transversal to $0$.
\item 
They are compatible with the description of the boundary in
Proposition $\ref{prop94}$ $(3)$.
\item
If $w_{i,j}$ is unforgetable then 
the evaluation map in Proposition $\ref{prop94}$ $(5)$ is 
strongly submersive with respect to our CF-perturbations.
\item The symmetry in Proposition $\ref{prop94}$ $(6)$, $(7)$ preserves 
the CF-perturbations.
\item
The compatibility of Kuranishi structures
in Proposition $\ref{prop94}$ $(8)$, $(9)$, $(10)$  is enhanced to the compatibility of
CF-perturbations. 
\item 
In the case of Proposition $\ref{prop94}$ $(11)$  the CF-perturbations 
are the ones in Proposition $\ref{existmkulti1}$.  They are constant in the $[0,1]^{\vert {\frak f}_{\rm{\bf f}}\vert}$
direction.
\item 
If ${\frak f}_{\bf f}$ is empty then the CF-perturbations coincide with the ones in Proposition $\ref{prop95}$.
\end{enumerate}
The ${\mathcal M}_{\ell;\vec k}((\vec{\kappa},\vec p,m);B;\vec{\frak f})^{\boxplus 1}$ version 
also holds.
\end{prop}
The proof is the same as the proof of Proposition \ref{prop95}.

We omit the discussion on the orientation since it is 
the same as Subsection \ref{oricyclic}.

\subsection{Construction of operations.}
\label{subsec:opeartor}

Suppose we are in the situation of 
Theorem \ref{cAinfconst}.
We will construct $A_{\infty}$ operations 
in Theorem  \ref{cAinfconst} in this subsection.
We use the notation described at the beginning of  Subsection \ref{constcyclic}
in this subsection.
We put\index[syindex]{BCFformcurve+@$BCF(\cL^{\rm form};\vec{\kappa})^+$}
\begin{equation}\label{form93}
BCF(\cL^{\rm form};\vec{\kappa})^+
= \bigotimes_{i=1}^K CF(L_{\kappa_{i-1}},L_{\kappa_{i}};\F)^+ [1]
\end{equation}
where
$$
\aligned
&CF(L_{\kappa_{i-1}},L_{\kappa_{i}};\F) \\
&= 
\begin{cases}
\displaystyle
\bigoplus_{p \in L_{\kappa_{i-1}} \cap L_{\kappa_{i}}} \F[p]
& \text{if $L_{\kappa_{i-1}} \ne L_{\kappa_{i}}$}, \\
\displaystyle (\Omega(\tilde L_{\kappa_{i}}) \otimes\F) \oplus \bigoplus_{p \in (\tilde L_{\kappa_i} \times_{X} \tilde L_{\kappa_i}) \setminus \tilde L_{\kappa_i}} \F[p]
& \displaystyle\text{ if $L_{\kappa_{i-1}} = L_{\kappa_{i}}$}
\end{cases}
\endaligned
$$
and\index[syindex]{CF(Fplus@$
CF(L_{\kappa_{i-1}},L_{\kappa_{i}};\F)^+$}
$$
\aligned
&CF(L_{\kappa_{i-1}},L_{\kappa_{i}};\F)^+
&= \begin{cases} 
\displaystyle
 CF(L_{\kappa_{i-1}},L_{\kappa_{i}};\F)     &\text{if $\kappa_{i-1} \ne \kappa_i$}, \\
  CF(L_{\kappa_{i}},L_{\kappa_{i}};\F) \oplus \F [{\bf e}^+_{\kappa_i}]  
  \oplus \F [{\bf f}_{\kappa_i}]    &\text{if $\kappa_{i-1} = \kappa_i$}.
 \end{cases}
 \endaligned
$$
Note that $\kappa_{i-1} = \kappa_i$ is equivalent to 
$(L_{\kappa_{i-1}},\theta_{\kappa_{i-1}}) = (L_{\kappa_{i}},\theta_{\kappa_{i}})$.

 %\marginpar{$CF(L_{\kappa_{i-1}},L_{\kappa_{i}};\Lambda_0)^+$ is changed 
%to $CF(L_{\kappa_{i-1}},L_{\kappa_{i}};\Lambda_0)$ in the second formula.  KF 2024 Dec.
%Coefficient is changed to $\F$. 2025 June KF}

For $B \in \Pi_2(\vec{\kappa},\vec p)$, we  define
\begin{equation}
  \label{eq:change_notation_moduli_space2}
\aligned
\mathcal M_{\ell}(\vec \kappa,\vec p;B;\vec{\frak f})
&=
\mathcal M_{\ell;\vec k'}((\vec{\kappa}',\vec p,m');B;\vec{\frak f})  \\
 [0,1]^{\vert {\frak f}_{\rm{\bf f}}\vert} \times\mathcal M_{\ell}(\vec \kappa,\vec p;B;\vec{\frak f}) 
&= 
 [0,1]^{\vert {\frak f}_{\rm{\bf f}}\vert} \times \mathcal M_{\ell;\vec k'}((\vec{\kappa}',\vec p,m');B;\vec{\frak f}) 
\endaligned
\end{equation}
in the same way as (\ref{eq:change_notation_moduli_space}).
Here ${\rm Red}(\vec{\kappa},\vec p) = (\vec{\kappa}',\vec k,\vec p,m')$.
\begin{rem}
When $L_{\kappa_{i-1}} = L_{\kappa_{i}}$ but $\theta_{\kappa_{i-1}} \ne \theta_{\kappa_{i}}$,
the boundary marked point $z_i$ is regarded to be unforgetable, by definition.
 %\marginpar{Remark added. KF. 2024 Dec.}
\end{rem}
We use the moduli space to define the $A_{\infty}$ operations:\footnote{The notation f.c.u. in the next formula stands for 
`form, curved, (homotopically) unital'.}
$$
\frak q^{\text{\rm f.c.u.}}_{\ell;\vec{\kappa};B}
: 
\Omega(X)\blue{[2]}^{\otimes \ell}
\otimes
BCF(\cL^{\rm form};\vec{\kappa})^+
\to CF(L_{\kappa_{0}},L_{\kappa_K};\F)^+[1].
$$
Let $\text{\bf g} = g_1\otimes \dots\otimes g_{\ell} \in \Omega(X)^{\otimes \ell}$ and 
$\text{\bf h} = h_1 \otimes \dots \otimes h_{K} \in BCF(\cL^{\rm form};\vec{\kappa})^+$.
We consider the moduli space with $K$ boundary marked points.
We define forgetability data $\vec{\frak f}_{\rm{\bf e}}$ and $\vec{\frak f}_{\rm{\bf f}}$ as follows:
\begin{enumerate}
\item
If $h_i$ is not one of ${\bf e}_{\kappa}^+$ or 
${\bf f}_{\kappa}$ then we regard the $i$-th marked point 
to be unforgetable, that is, $i \notin  {\frak f}_{\rm{\bf e}} \cup  {\frak f}_{\rm{\bf f}}$.
\item
If $h_i = {\bf e}_{\kappa}^+$ then $i \in  {\frak f}_{\rm{\bf e}}$.
\item
If $h_i = {\bf f}_{\kappa}$ then $i \in  {\frak f}_{\rm{\bf f}}$.
\end{enumerate}
We put $\text{\bf h}' = h'_1 \otimes \dots \otimes h'_{K}$ where
$h'_i = h_i$ in case (1) above and $h'_i = {\bf e}_{\kappa} = 1$ (the constant function $1$ on $L_{\kappa}$)
in case (2) or (3).
Now we define:\index[syindex]{qfcuellkappaB@$\frak q^{\text{\rm f.c.u.}}_{\ell;\vec{\kappa};B}$}
\begin{defn}\label{def96}
We put $\text{\rm ev}_{>0} = (\text{\rm ev}_1,\dots,\text{\rm ev}_k)$ and define:\index[syindex]{evl1>0@$\text{\rm ev}_{>0}$}
$$
\aligned
\frak q^{\text{\rm f.c.u.}}_{\ell;\vec{\kappa};B}(\text{\bf g},\text{\bf h} )
&=
(-1)^{\maltese}(\text{\rm ev}_{0})_!
\left((\text{\rm ev}^+)^*\text{\bf g}
\wedge
(\text{\rm ev}_{>0})^*\text{\bf h}'
\right) \\
&
= (-1)^{\maltese}\text{\rm Corr}
([0,1]^{\vert {\frak f}_{\rm{\bf f}}\vert}  \times \mathcal M_{\ell}(\vec \kappa,\vec p;B;\vec{\frak f}),(\text{\rm ev}_{0};
\text{\rm ev}^+\times \text{\rm ev}_{>0}
))
(\text{\bf g}\wedge \text{\bf h}' 
)\\
&
\in CF(L_{\kappa_{0}},L_{\kappa_{K}};\F).
\endaligned$$
other than the following cases. 
\begin{equation}\label{excepfff}
\aligned
\frak q^{\text{\rm f.c.u.}}_{0;(\kappa);0}(1,\text{\bf f}_{\kappa})
&= \text{\bf e}^+_{\kappa} - \text{\bf e}_{\kappa}
 \\
 \frak q^{\text{\rm f.c.u.}}_{0;(\kappa);0}(1,\text{\bf f}_{\kappa}\otimes \text{\bf f}_{\kappa})
&= 0
 \\
 \frak q^{\text{\rm f.c.u.}}_{0;(\kappa,\kappa);0}(1,\text{\bf f}_{\kappa}\otimes\text{\bf e}^+_{\kappa})
&= -\text{\bf f}_{\kappa}\\
 \frak q^{\text{\rm f.c.u.}}_{0;(\kappa,\kappa);0}(1,\text{\bf e}^+_{\kappa}\otimes\text{\bf f}_{\kappa})
&= \text{\bf f}_{\kappa}\\
 \frak q^{\text{\rm f.c.u.}}_{0;(\kappa,\kappa);0}(1,\text{\bf e}^+_{\kappa}\otimes \text{\bf e}^+_{\kappa})
&= \text{\bf e}^+_{\kappa}
\endaligned
\end{equation}
Here $\text{\bf e}_{\kappa}$ is the function $1$ on $L_{\kappa}$.
\end{defn}
\begin{rem}
The sign ${\maltese}$ is obtained by modifying the sign in (\ref{defqformula})
by the sign to move each parameter in $[0,1]^{\vert {\frak f}_{\rm{\bf f}}\vert}$ to the 
position of the corresponding ${\bf f}_{\kappa}$.
\end{rem}

\begin{lem}\label{qpropertiescatuni}
The operators $\frak q^{\text{\rm f.c.u.}}_{\ell;\vec{\kappa};B}$ have the following properties:
\begin{enumerate}
\item
The equality $(\ref{qmaineqcat})$ in Lemma $\ref{qpropertiescat}$ holds.
\item The equality $(\ref{qismcat})$ in Lemma $\ref{qpropertiescat}$ holds.
In fact in the case when $h_0$ or $h_1$  is either ${\bf e}^+_{\kappa}$ or ${\bf f}_{\kappa}$
the formula is not defined since the inner product for such elements 
are not defined.  We will explain this point during the proof.
\item If the inputs do not contain ${\bf e}^+_{\kappa}$  or ${\bf f}_{\kappa}$ then  the operators 
coincide with one defined in Definition $\ref{defn626}$.
\end{enumerate}
\end{lem}
\begin{proof}
The proof is rather immediate from 
Propositions \ref{exstCFUNI} together with Stokes' theorem (\cite[Theorem 9.28]{springer}) and 
the composition formula (\cite[Theorem 10.21]{springer}).
\par
When none of the factor of the input ${\bf h}$ is ${\bf f}_{\kappa}$,
the $A_{\infty}$ relation follows from the description of the boundary 
in Proposition \ref{prop92} in the same way as the proof of Lemma \ref{qpropertiescat}.
\par
If one of the inputs ${\bf h}$ is ${\bf f}_{\kappa}$ then we have an extra boundary 
corresponding to the case when $s_{i,j}$ in the factor $ [0,1]^{\vert \frak f_{\bf f}\vert}$ 
goes to $0$ or $1$.
When $s_{i,j}$ goes to $1$, Proposition \ref{prop94} (9) implies that the output 
will be the case when ${\bf f}_{\kappa}$ is replaced by ${\bf e}^+_{\kappa}$.
When $s_{i,j}$ goes to $0$ Proposition \ref{prop94} (10) implies that the output 
will be the case when ${\bf f}_{\kappa}$ is replaced by ${\bf e}_{\kappa}$.
Thus taking the formula (\ref{excepfff}) into account, the $A_{\infty}$
relation holds including the case when some inputs are ${\bf f}_{\kappa}$
using the formula (\ref{condforrf}).
\par
Proposition \ref{prop92} (11) and Proposition \ref{prop95} (6) imply Item (3) of Lemma \ref{qpropertiescatuni}.
\par
We finally discuss cyclic symmetry, Formula  (\ref{qismcat}).
We first consider Hochschild chain complex\index[syindex]{CHCFLformcu@$CHCF_*(\cL^{\rm form};\vec{\kappa})^+$}
\begin{equation}\label{CHCFform}
CHCF_*(\cL^{\rm form};\vec{\kappa})^+
= \bigotimes_{i=1}^{K+1} CF(L_{\kappa_{i-1}},L_{\kappa_{i}};\F)^+ \blue{[1]}
\end{equation}
where $\kappa_{K+1} = \kappa_{0}$ by definition.
This is a particular case of (\ref{eq:cyclic_bar-cplex}).

We will define\index[syindex]{qfvcuplus@${}^{+}{\frak q}^{\text{\rm f.c.u.}}_{\ell;\vec{\kappa};B}$}
$$
{}^{+}{\frak q}^{\text{\rm f.c.u.}}_{\ell;\vec{\kappa};B}
:
\Omega(X)\blue{[2]}^{\otimes \ell}
\otimes
CHCF(\cL^{\rm form};\vec{\kappa})^+ \to \F.
$$

Let $\text{\bf g} = g_1\otimes \dots\otimes g_{\ell} \in \Omega(X)[2]^{\otimes \ell}$, 
$\text{\bf h} =  h_1 \otimes \dots \otimes h_{K} \otimes h_0 \in CHCF_*(\cL^{\rm form};\vec{\kappa})^+$. 
(Note that $h_0$ is included.)
We put

\begin{enumerate}
\item
If $h_i$ is not one of ${\bf e}_{\kappa}^+$ or 
${\bf f}_{\kappa}$ then we regard the $i$-th marked point 
to be unforgetable, that is, $i \notin {\frak f}_{\rm{\bf e}} \cup {\frak f}_{\rm{\bf f}}$.
\item
If $h_i = {\bf e}_{\kappa}^+$ then $i \in  {\frak f}_{\rm{\bf e}}$.
\item
If $h_i = {\bf f}_{\kappa}$ then $i \in  {\frak f}_{\rm{\bf f}}$.
\end{enumerate}
We put $\text{\bf h}' = h'_0 \otimes h'_1 \otimes \dots \otimes h'_{K}$ where
$h'_i = h_i$ in case (1) above and $h'_i = {\bf e}_{\kappa} = 1$ (the constant function $1$ on $L_{\kappa}$)
in case (2) or (3).
\begin{defn}\label{defn98}
We  define
\begin{equation}\label{form950}
{}^{+}{\frak q}^{\text{\rm f.c.u.}}_{\ell;\vec{\kappa};B}(\text{\bf g},\text{\bf h}
)
= (-1)^{\maltese}\int_{[0,1]^{\vert\frak f_{\bf f}\vert} \times \mathcal M_{\ell}(\vec{\kappa},\vec p;B;\vec{\frak f})}
(
\text{\rm ev}^+\times \widetilde{\text{\rm ev}} )^*% \times 
(\text{\bf g}\wedge \text{\bf h}').
\end{equation}
Here $\widetilde{\text{\rm ev}} = {\text{\rm ev}} \circ \pi$ and 
$\pi : [0,1]^{\vert\frak f_{\bf f}\vert} \times \mathcal M_{\ell}(\vec{\kappa},\vec p;B) \to \mathcal M_{\ell}(\vec{\kappa},\vec p;B)$ 
is the projection. We can define the sign $\maltese$ so that (\ref{form95}) holds with sign.
\par
There are a few exceptions.
\par
In case $B=0$, ${\bf g} = 1$ and all $h_i$'s are either ${\bf e}_{\kappa}^+$ or $\bf f_{\kappa}$ 
we put ${}^{+}{\frak q}^{\text{\rm f.c.u.}}_{0;\vec{\kappa};0}(1,\text{\bf h}) = 0$. 
\par
In case $B=0$, ${\bf g} = 1$ and ${\bf h} =  ({\bf f}_{\kappa},h_2)$  with $h_2 \ne {\bf f}_{\kappa},{\bf e}_{\kappa}^+$  we put
\begin{equation}\label{form950exc}
{}^{+}{\frak q}^{\text{\rm f.c.u.}}_{0;(\kappa,\kappa);0}(1,{\bf f}_{\kappa},h_2) 
=
- \langle  {\bf e}_{\kappa},h_2 \rangle_{\rm cyc}.
\end{equation}
${}^{+}{\frak q}^{\text{\rm f.c.u.}}_{0;(\kappa,\kappa);0}(1,h_2,{\bf f}_{\kappa})$ is then determined by 
cyclic symmetry. \footnote{The  definition (\ref{form950exc}) is related to the equality 
$\frak q_{0;(\kappa);0}(1;{\bf f}_{\kappa}) = {\bf e}_{\kappa}^+ - {\bf e}_{\kappa}$ and (\ref{fmnew103}).}

\end{defn}
Comparing Definitions \ref{def96} with \ref{defn98}   and using the fact
that integration along the fiber is 
adjoint to the pull back (\cite[Lemma 7.84 (2)]{springer}) 
we find the next formula in case $h_0$ is neither ${\bf e}^+_{L_{\kappa}}$ nor ${\bf f}_{L_{\kappa}}$.
\begin{equation}\label{form95}
{}^{+}{\frak q}^{\text{\rm f.c.u.}}_{\ell;\vec{\kappa};B}(\text{\bf g},\text{\bf h})
= 
\langle {\frak q}^{\text{\rm f.c.u.}}_{\ell;\vec{\kappa};B}(\text{\bf g},\text{\bf h}_-),h_0  \rangle_{\rm cyc}
\end{equation}
where $\text{\bf h}_- = (h_1,\dots,h_K)$.

Moreover  
Proposition \ref{prop92} (8), 
Proposition \ref{prop94} (7), Proposition  \ref{prop95} (4) and 
Proposition  \ref{exstCFUNI} (4)
imply
\begin{equation}\label{qismcatprime1}
{}^{+}\frak q^{\text{\rm f.c.u.}}_{\ell;\vec{\kappa};B}
(\text{\bf g};h_0,h_1,\dots,h_{K-1},h_K) = (-1)^\maltese
{}^{+}\frak q^{\text{\rm f.c.u.}}_{\ell;c\vec{\kappa};cB}
(\text{\bf g};h_1,\dots,h_K,h_0),
\end{equation}
where 
$\maltese = (\deg' h_0)(\deg' h_1 + \dots + \deg' h_K)$. 
(\ref{qismcatprime1}) is actually Lemma \ref{qpropertiescatuni} (2).
\end{proof}
We also define:\index[syindex]{mxzfcukappabplus@${}^{+}{\frak m}^{\text{\rm f.c.u.}}_{\vec{\kappa};B}$}
\begin{equation}\label{form952}
{}^{+}{\frak m}^{\text{\rm f.c.u.}}_{\vec{\kappa};B}(\text{\bf h})
= 
{}^{+}{\frak q}^{\text{\rm f.c.u.}}_{\ell;\vec{\kappa};B}(1,\text{\bf h}).
\end{equation}
We have:
\begin{equation}\label{qismcatprime}
{}^{+}\frak m^{\text{\rm f.c.u.}}_{\vec{\kappa};B}
(h_0,h_1,\dots,h_{K-1},h_K) = (-1)^\maltese
{}^{+}\frak m^{\text{\rm f.c.u.}}_{c\vec{\kappa};cB}
(h_1,\dots,h_K,h_0),
\end{equation}
where 
$\maltese = (\deg' h_0)(\deg' h_1 + \dots + \deg' h_K)$.

Now in the same way as Definition \ref{defn626}
(\ref{form95}) we define:

\begin{defn}\label{defn626uni}
For given $\vec{\kappa}$ we define\index[syindex]{qfcuellkappafrakb@$\frak q^{\text{\rm f.c.u.}\frak b}_{\ell,\vec{\kappa}}$}
$$
\frak q^{\text{\rm f.c.u.}\frak b}_{\ell,\vec{\kappa}}
: 
\Omega(X)^{\otimes \ell}
\otimes 
BCF(\cL^{\rm form};\vec{\kappa})^+ \widehat{\otimes} \Lambda_0
\to CF(\newred{L_{\kappa_{\blue{0}}},L_{\kappa_{\blue{K}}}};\F)^+
\widehat{\otimes} \Lambda_0,
$$
by
\begin{equation}\label{qcatuni}
\frak q^{\text{\rm f.c.u.}\frak b}_{\ell,\vec{\kappa}}(\text{\bf g};\text{\bf h})
=
\sum_{B,\ell'}T^{\omega(B) }
\frac{\rho_{\frak b,\theta}(B)% \rho_{b}(B)
}{\ell'!}
\frak q^{\text{\rm f.c.u.}}_{\ell+\ell';\vec{\kappa};B}
(\text{\bf g}\otimes \frak b_+^{\ell'};
\text{\bf h}) %;b_+^{\vec k})
\end{equation}
% $b_+^{\vec k} = \bigotimes_{i=0}^K b_{\kappa_i,+}^{\otimes k_i}$,
where the sum is taken over all classes $B$ and parts $T^{\lambda_i} \frak b_{+,i}$
of $\frak b_+$ such that
$\omega(B) + \sum \lambda_i \le E$.
\par
We define
and 
$$
{}^+\frak q^{\text{\rm f.c.u.}\frak b}_{\ell,\vec{\kappa}}(\text{\bf g};\text{\bf h}):
\Omega(X)\blue{[2]}^{\otimes \ell}
\otimes
CHCF(\cL^{\rm form};\vec{\kappa})^+\to \Lambda_0^{\F},
$$
by taking a weighted sum of ${}^{+}\frak q^{\text{\rm f.c.u.}}_{\ell;\vec{\kappa};B}$ in the same way as (\ref{qcatuni}).
\par
We also define\index[syindex]{mxfcubfrak@$\frak m^{\text{\rm f.c.u.}\frak b}_{\vec{\kappa}}$} 
$$
\aligned
\frak m^{\text{\rm f.c.u.}\frak b}_{\vec{\kappa}}
&: 
BCF(\cL^{\rm form};\vec{\kappa})^+\widehat{\otimes} \Lambda_0
\to CF(\newred{L_{\kappa_{\blue{0}}},L_{\kappa_{\blue{K}}}})^+ \widehat{\otimes} \Lambda_0,
\\
{}^+\frak m^{\text{\rm f.c.u.}\frak b}_{\vec{\kappa}}
&: 
CHCF(\cL^{\rm form};\vec{\kappa})^+\widehat{\otimes} \Lambda_0
\to \Lambda_0^{\F},
\endaligned
$$
by
\begin{equation} \label{eq:a-oo-operation-trivial-gauge}
\frak m^{\text{\rm f.c.u.}\frak b}_{\vec{\kappa}}(\text{\bf h})
= 
\frak q^{\text{\rm f.c.u.}\frak b}_{0}(1;\text{\bf h}).
\end{equation}
and  
$$
{}^+\frak m^{\text{\rm f.c.u.}\frak b}_{\vec{\kappa}}({\bf h} \otimes h_0)
={}^+\frak q^{\text{\rm f.c.u.}\frak b}_{\ell,\vec{\kappa}}(1;\text{\bf h} \otimes h_0).
$$
\end{defn}

\begin{figure}[h] 
\centering
\includegraphics[scale=0.3]{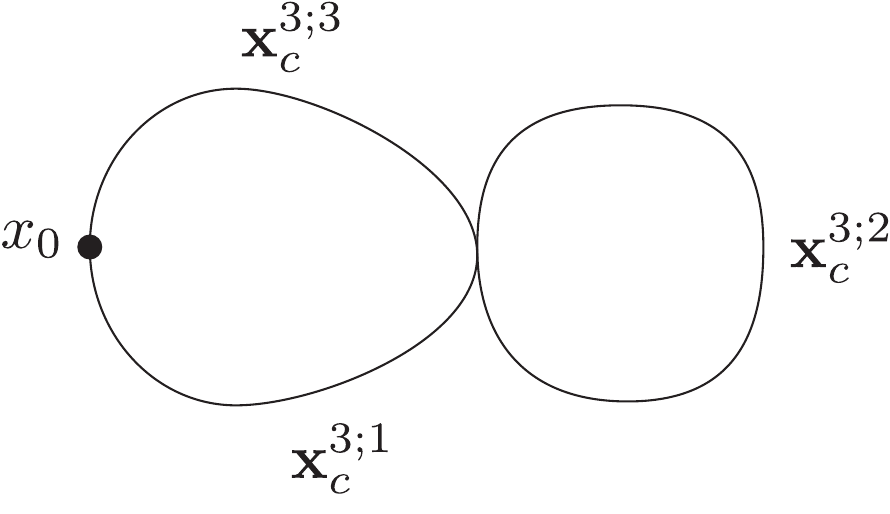}
\caption{Definition \ref{def98} (d).}
\label{Figurecyclic(d)}
\end{figure}

The next statement is immediate. (Definition \ref{def98} (d) also follows from Stokes' formula. See Figure \ref{Figurecyclic(d)}):
 %\marginpar{Figure is added.  KF 2025 Aug.}
\begin{prop}\label{prop910}
The operations $\frak m^{\text{\rm f.c.u.}\frak b}_{\vec{\kappa}}(\text{\bf h})$
define a structure of homotopically unital and cyclic $G_+$-gapped filtered $A_{n,k}$ category\footnote{whose morphism 
module is $CF(L_{\kappa},L_{\kappa'};\F) \otimes_{\F}\Lambda^{\F}_0.$},
for given $n,k$.\footnote{The cyclicity of $A_{n,k}$ category is defined by modifying 
the $A_{\infty}$ case in an obvious way.}
\end{prop}
We then can work out the homotopy  limit process in the same way 
as Subsections \ref{sec:constr-cycl-a_infty4} and \ref{sec:homotopyequiv}  to obtain a structure over $\Lambda_0$.
We denote it by $\mathcal L^{\rm form}_{\rm c.u.}$. \index[syindex]{Lcateformcu@{$\cL^{\rm form}_{\rm c.u.}$}}
In fact the notion of cyclic and homotopically unital pseudo-isotopy is defined as follows.

\begin{defn}\label{defn913} %\marginpar{Definition added.  KF 2025 Jan.}
Let $(\frak C,\langle\, ,\,\rangle,\{\frak m^{\frak C, t}_{k,\beta}\},\{\frak c^{\frak C, t}_{k,\beta}\})$
be a cyclic pseudo-isotopy between $G$-gapped filtered $A_{\infty}$ categories
$\Cat$, $\Dat$. Namely it satisfies Definition \ref{pisotopydef} (1)-(4).
Suppose $\Cat$, $\Dat$ are homotopically unital and is cyclic as 
in Definition \ref{def98}.
\par
Now we say  $(\frak C,\langle\, ,\,\rangle,\{\frak m^{\frak C, t}_{k,\beta}\},\{\frak c^{\frak C, t}_{k,\beta}\})$
can be extended to a cyclic and homotopically unital pseudo-isotopy between $G$-gapped filtered $A_{\infty}$ categories
if the following holds.
\begin{enumerate}
\item $\frak m^{\frak C, t}_{k,\beta}$ and $\frak c^{\frak C, t}_{k,\beta}$ 
extend to the operators 
$\frak m^{\frak C,+, t}_{k,\beta}$ and $\frak c^{\frak C,+, t}_{k,\beta}$
which take ${\bf e}^+_c$ and  ${\bf f}_c$ also as inputs.
\item 
For each $t$ the operators $\frak m^{\frak C,+, t}_{k,\beta}$ define a 
homotopically unital filtered $A_{\infty}$ structure.
Moreover ${}^+\frak m^{\frak C,+, t}_{k,\beta}$ is given as in Definition \ref{def98}
such that for each $t$ the above homotopically unital filtered $A_{\infty}$ structure is 
cyclic.
\item (\ref{isotopymaineq}) holds with $\frak m^{\frak C, t}_{k,\beta}$, 
$\frak c^{\frak C, t}_{k,\beta}$ replaced by
$\frak m^{\frak C,+, t}_{k,\beta}$ and $\frak c^{\frak C,+, t}_{k,\beta}$, respectively.
\item ${}^+\frak c^{\frak C,+, t}_{k,\beta}$ is given and plays the same role 
for $\frak c^{\frak C,+, t}_{k,\beta}$
as ${}^+\frak m^{\frak C,+, t}_{k,\beta}$ does for $\frak m^{\frak C,+, t}_{k,\beta}$.
\item (\ref{candeequ}) with $\frak c^{\frak C, t}_{k,\beta}$ replaced by
$\frak c^{\frak C,+, t}_{k,\beta}$, ${}^+\frak c^{\frak C,+, t}_{k,\beta}$
holds. 
\end{enumerate}
We can define the $A_{n,k}$ version in a similar way.
\end{defn}
Then Proposition \ref{prop:commutative_diagram_energy_induction} can be 
generalized to homotopically unital and cyclic case in the same way.
We use the proposition to work out the  homotopy  limit process.

Presuming Propositions  \ref{prop92}, \ref{prop94}, \ref{prop95} and \ref{exstCFUNI}  we have finished the proof of Theorem  \ref{homotopyunitality}.
The proofs of Propositions  \ref{prop92} and \ref{prop94} are postponed until Subsection \ref{polygonforget}
and the proofs of Propositions \ref{prop95} and \ref{exstCFUNI} will complete in Subsection \ref{CFforget}.
\qed
\begin{rem}
Note that the cyclicity in Proposition \ref{prop910} is not so satisfactory, since the inner product is not defined for elements ${\bf e}^+_{{\kappa}}$ and ${\bf f}_{{\kappa}}$. 
After taking the canonical model where we take the de Rham cohomology instead of de Rham chain complex 
(or its $^+$ version), this point is resolved and we obtain a genuine  unital and cyclic filtered $A_{\infty}$ category.
See Section \ref{canonical}.
\end{rem}

\section{Elimination of the curvature.}
\label{sec:elicurv}
 %\marginpar{several things are moved to here.}

In this section, we explain how to pass from the above curved category to the category whose objects
are pairs 
of a Lagrangian submanifold and its  bounding cochain. 
In Definition \ref{boundingcochain} we defined the notion of 
weak bounding cochain of an $A_{\infty}$ algebra $\{\frak m_k\}$.
Because of Lemma \ref{lem104} we state below, this structure is mostly the same 
as the restriction of the one we defined in the last section.
A slight difference is that the structure in the last section contains two 
extra generators ${\bf f}_{L_{\kappa}}$ and ${\bf e}^+_{L_{\kappa}}$.
We first give a few simple algebraic lemmas 
to clarify this point.
\begin{lem}\label{uniishuni}
Let $(\mathcal C,\{\frak m_k\})$ be a curved and gapped 
filtered $A_{\infty}$ algebra which has a strict unit ${\bf e}_c$.
Then it is  homotopically unital in a canonical way.
\end{lem}
\begin{proof}
We consider 
$\mathcal C^+ = \mathcal C \oplus \Lambda_0 {\bf e}_c^+ \oplus \Lambda_0 {\bf f}_c$
and define an $A_{\infty}$ structure on $\mathcal C^+$ as follows:
When all the inputs are in  $\mathcal C$ the $A_{\infty}$ operation is the same 
as one of $\mathcal C$.
If the input contains ${\bf e}^+_c$ then $A_{\infty}$ operation is determined 
from the fact ${\bf e}^+_c$ is the unit.
If the input contains ${\bf f}_c$ then $A_{\infty}$ operation is zero 
except the case it is determined 
from the fact ${\bf e}^+_c$ is the unit and the next formula.
$$
\frak m_1({\bf f}_c) = {\bf e}_c^+ - {\bf e}_c.
$$
\end{proof}
\begin{lem}\label{CandC+hom}
In the situation of Lemma $\ref{uniishuni}$, 
$\mathcal C$ is unitally, filtered and gapped homotopy equivalent to $\mathcal C^+$.\footnote{See
\cite[Subsection 4.2.3]{fooo09} for the definition.}
\end{lem}
\begin{proof}
We define a unital, filtered and gapped $A_{\infty}$ homomorphism 
$\frak F : \mathcal C^+ \to \mathcal C$ as follows.
When restricted to $ \mathcal C \subset \mathcal C^+$ the map
$\frak F$ is the identity.  (In particular it is linear there.)
If the input contains ${\bf e}^+$ the map $\frak F$ is determined 
by the fact that it is unital.
If the input contains ${\bf f}$ then $\frak F$ is zero.
It is easy to check that $\frak F$   is a unitally, filtered and gapped $A_{\infty}$
homomorphism.  It is also easy to see that it induces an isomorphism on 
$\overline{\frak m}_1$ cohomology.  Therefore it is a 
homotopy equivalence by \cite[Theorem 4.2.45a]{fooo09}.
\end{proof}
\begin{lem}
Suppose we are in the situation of Lemma $\ref{uniishuni}$.
If $b$ is a  bounding cochain of  $\mathcal C$ with potential value $c$ then
$b + c{\bf f}$ is a  bounding cochain of  $\mathcal C^+$ with potential value $c$.
\par
The  bounding cochain $b$ is gauge equivalent to $b'$  if and only if $b + c{\bf f}$ is gauge equivalent to $b' + c{\bf f}$.
\end{lem}
The proof is an easy calculation. (See \cite[Remark 3.6.31]{fooo09}.)
Lemma \ref{CandC+hom} implies that an arbitrary   bounding cochain 
of $\mathcal C^+$ is gauge equivalent to the one of the form $b + c{\bf f}$.

We observe:
\begin{lem}\label{lem104}
The unital $A_{\infty}$ algebra structure on $\Omega(L_{\kappa}) \hat\otimes \Lambda_0$
in Proposition $\ref{Algebraonderhamcomplex}$ induces one on  $(\Omega(L_{\kappa}) \hat\otimes \Lambda_0)^+$
in Subsection $\ref{subsec:opeartor}$ by Lemma $\ref{uniishuni}$.
This unital $A_{\infty}$ algebra in Section $\ref{ainfalgasssingle}$ 
coincides with one obtained by restricting the $A_{\infty}$ structure 
in Subsection $\ref{subsec:opeartor}$.
\end{lem}
This is a consequence of Proposition \ref{prop92} (9) and Proposition  \ref{prop94} (11).

Thus if we have a  bounding cochain $b$ of $\Omega(L_{\kappa}) \hat\otimes \Lambda_0$
with potential value $c$ then 
$b + c{\bf f}$ is a  bounding cochain of $(\Omega(L_{\kappa}) \hat\otimes \Lambda_0)^+$.
Sometimes we omit ${\bf f}$ and write $b$ in stead of $b + c{\bf f}$ .
\begin{defn}\label{defn105}
We define $\bL$ to be a collection of 
  $\{(L_{\kappa},\theta_{\kappa}, b_{\kappa})\mid \kappa =1,\dots,K\}$ such that $b_{\kappa}$ is a  bounding cochain. These are the objects of an $A_\infty$ category which we denote by $\cL^{\rm form}_{\rm uni}$, with $A_\infty$ operations given by the formula\index[syindex]{Lformuni@$\cL^{\rm form}_{\rm uni}$}\index[syindex]{mxfubbdk@$\frak m^{\text{\rm f.u.}\text{\bf b}}_k$}
  \begin{equation} \label{eq:mk-operatioin-forms-bounding-cochain}
    \frak m^{\text{\rm f.u.}\text{\bf b}}_k( x_1 \otimes \cdots \otimes x_k ) \equiv \sum    \frak m^{\text{\rm f.c.u.}\frak b}(b^{m_0}_{\kappa_0,+} \otimes x_1 \otimes b^{m_1}_{\kappa_1,+} \otimes \cdots \otimes x_k \otimes  b^{m_k}_{\kappa_k,+}), 
  \end{equation}
  whenever $x_i$ is a morphism from the object labelled by $\kappa_{i-1}$ to the object labelled by $\kappa_{i}$.
  This  filtered $A_{\infty}$ category is homotopically unital.
  \par
  Note in (\ref{eq:mk-operatioin-forms-bounding-cochain}), $b_{\kappa} = \theta_{\kappa} + b_{\kappa,+}$
  and $b_{\kappa,+} \in H^1(L_{\kappa},\Lambda_+) \oplus \bigoplus_{2j+1 \ge 3}H^{2j+1}(L_{\kappa},\Lambda_0)$.
  (See Definition \ref{boundingcochain}.)
  $\theta_{\kappa}$ is used in the definition of $\frak m^{\text{\rm f.c.u.}\frak b}$.  (See (\ref{qcat}).)
 The definition of ${}^+\frak m_{k+1}^{\text{\rm f.u.}\text{\bf b}}$  is similar from  ${}^+\frak m^{\text{\rm f.c.u.}\frak b}_{\vec{\kappa}}$.
 \end{defn} 

 All the objects of $\cL^{\rm form}_{\rm uni}$  are weakly curvature free.  
 We have thus proved Theorem \ref{cAinfconstuni}.
 \qed
 \begin{rem}
 Note that the  bounding cochain $b+c{\bf f}$ contains an element ${\bf f}$ of negative degree.
In case the bounding cochain contains an element of negative degree, 
the proof of Lemma \ref{lem:Tadic-convergence} which claims that (\ref{S7MCeq}) in Definition \ref{boundingcochain} converges given there does not work.
However using the fact that the coefficient of ${\bf f}$, which is $c = \frak{PO}(b)$ is in $\Lambda_+$ we 
can easily modify the argument there to show  (\ref{S7MCeq}) and (\ref{eq:mk-operatioin-forms-bounding-cochain}) still converge.
 \end{rem}
 \par\medskip
 We denote by $\cL_{\lambda}^{\rm form}$ the subcategory\index[syindex]{Llambdaform@$\cL_{\lambda}^{\rm form}$} of $\cL^{\rm form}_{\rm uni}$ whose set of objects $\bL_\lambda$ consists of those with potential value $\lambda$.  In particular, Equation (\ref{qmaineqcat2}) implies:

\begin{lem}\label{lambdaandlambdaprime}
For a pair of objects $(L_i,\theta_i,b_i)$ and $(L_{i'},\theta_{i'},b_{i'})$ of $\cL_{\lambda_i}^{\rm form} $ and  $\cL_{\lambda_{i'}}^{\rm form}$, we have
\begin{equation}\label{ddnonzero}
\frak m^{\rm f.u.\text{\bf b}}_{1}
\circ \frak m^{\rm f.u.\text{\bf b}}_{1}
= (\lambda_i -\lambda_{i'}) {\rm id}.
\end{equation} \qed
\end{lem}
  
As a special case, if both objects have the same curvature, we have 
\begin{equation}\label{ddzero}
\frak m^{\rm f.u.\text{\bf b}}_{1}
\circ \frak m^{\rm f.u.\text{\bf b}}_{1}
= 0.
\end{equation}
See \cite[Proposition 3.7.17]{fooo09}.

\begin{rem}\label{MCoddelementrem}
As we prove in Section \ref{canonical} below (and 
also as in \cite[Theorem 5.4.2]{fooo09}) 
there exists a structure of filtered $A_{\infty}$ algebra 
$\{\frak m^{\text{\bf b}}_k\}$ on 
the cohomology group $H(L;\Lambda_0)$
that is homotopy equivalent to 
$(\Omega(L),\{\frak m^{\text{\rm f.u.\bf b}}_k\})$.
By the homotopy invariance of Maurer-Cartan moduli spaces 
(see \cite[Theorem 4.3.22]{fooo09}) we have an isomorphism\index[syindex]{MweakomegaL@$\mathcal M_{\text{\rm weak}}(\Omega(L),\{\frak m^{\text{\bf b}}_k\})$}\index[syindex]{MweakHL@$\mathcal M_{\text{\rm weak}}(H(L;\Lambda_0),\{\frak m^{\text{\bf b}}_k\})$}
$$
\mathcal M_{\text{\rm weak}}(\Omega(L),\{\frak m^{\text{\bf b}}_k\})
=
\mathcal M_{\text{\rm weak}}(H(L;\Lambda_0),\{\frak m^{\text{\bf b}}_k\}).
$$
The right hand side is the set of gauge equivalence classes  of elements of  
$\widehat{\mathcal M}_{\text{\rm weak}}(H(L;\Lambda_0),\{\frak m^{\text{\bf b}}_k\})
$, 
and 
$\widehat{\mathcal M}_{\text{\rm weak}}(H(L;\Lambda_0),\{\frak m^{\text{\bf b}}_k\})
\subset
H^{\rm odd}(L;\Lambda_0).
$
Therefore an element of 
$\mathcal M_{\text{\rm weak}}(\Omega(L),\{\frak m^{\rm f.u.\text{\bf b}}_k\})$
is represented by a cohomology class of odd degree.
\end{rem}

\section{Canonical model.}
\label{canonical}

In the  filtered $A_{\infty}$ category that we defined in Section
\ref{sec:elicurv},
the morphism space 
from $(L_{\kappa},b_{\kappa})$ to $(L_{\kappa'},b_{\kappa'})$  is the de Rham complex in case
$L_{\kappa} = L_{\kappa'}$ (or its $+$ version when $b_{\kappa} = b_{\kappa'}$). 
In certain situations, we prefer to work with a  filtered $A_{\infty}$ category
whose morphism spaces are finite dimensional.
We use the canonical model, which we construct in this subsection, for this purpose. 
Another reason to go to the canonical model is an issue that 
the inner product are not defined for ${\bf e}^+_{L_{\kappa}}$ or  ${\bf f}_{L_{\kappa}}$.  %\marginpar{A sentence added. KF}
In this subsection we follow the strategy explained in \cite[Section 5.4]{fooo09}
and \cite[Section 10]{fukaya:cyc}.
See also \cite{KS}. We restrict our construction to the category of weakly unobstructed objects, even though it would be possible to work at the level of curved categories.

\subsection{Reduction to the canonical model.}
\label{subsec:red}

%In particular, we shall write $\bf b$ for the choice of the bulk class $\frak b$, together with all 
%weak bounding cochains in the category we consider.

\begin{defn}\label{decribbon}
We define the following objects $\Gamma = (\mathcal T,\mathcal K(\cdot),p(\cdot),B(\cdot))$
 %\marginpar{The definition is modified. Since related one in previous section is modified.  KF 2025 Aug.}
to be a {\it decorated 
ribbon tree}\index{decorated ribbon tree}.
\begin{enumerate}
\item $\mathcal T$ is a ribbon tree. 
\item
Let $C_0(\mathcal T)$ be the set of its vertices. 
We assume that we are given its decomposition 
$C_0(\mathcal T) = C^{\text{\rm int}}_0(\mathcal T) \sqcup C^{\text{\rm ext}}_0(\mathcal T)$ 
into the sets of {\it interior vertices}\index{interior vertex} and of {\it exterior vertices}\index{exterior vertex}.
We enumerate exterior vertices as $v^{\text{\rm ext}}_0,v^{\text{\rm ext}}_1,\dots,v^{\text{\rm ext}}_N$, 
which respect the counter-clockwise cyclic order 
induced by the ribbon structure of $\mathcal T$.
Denote by $e^{\text{\rm ext}}_i$ the unique edge which contains $\blue{v_i^{\text{\rm ext}}}$.
Those edges are called {\it exterior edges}\index{exterior edge}.
All other edges are called {\it interior}\index{interior edge}.
We call $v^{\text{\rm ext}}_0$ the {\it root} of $\mathcal T$.
\item
Let $C_1(\mathcal T)$ be the set of all edges of $\mathcal T$.
We orient the edges so that 
for each vertex $v$ there exists a unique path 
from $v$ to $v^{\text{\rm ext}}_0$ respecting the 
orientation of the edges.
\item
Note that $\bL$  is a finite set. We put
$\bL = \{(L_{\kappa},\theta_{\kappa}, b_{\kappa}) \mid \kappa \in \mathscr K\}$.
We take a map
$\mathcal K = (\mathcal K_1,\mathcal K_2) : C_1(\mathcal T) \to \mathscr K^2$ that satisfies the property
\begin{equation}\label{compatibility}
\mathcal K_2(e_i(v)) = \mathcal K_1(e_{i+1}(v)), \qquad i = 0,1,\dots,m(v),
\end{equation}
at each vertex $v$ 
where $\{e_0(v),e_1(v),\dots,e_{m(v)}(v)\}$ is the set of all the edges 
containing $v$ which are enumerated with respect to the counter-clockwise cyclic order given by the ribbon structure,  so that 
$e_0(v)$ is the unique edge which starts at $v$ with respect to the 
given orientation of the edges. In addition, we set $e_{m(v)+1}(v) = e_0(v)$.
\item
We assume a decomposition of $C_1(\mathcal T) =  C^s_1(\mathcal T)
\sqcup C_1^d(\mathcal T)$ 
into the set of edges of {\it s-type} and of  {\it d-type}\footnote{s stands for 
switching and d stands for diagonal.} is given.
If $e$ is d-type we require $L_{\mathcal K_1(e)} = L_{\mathcal K_2(e)}$.\footnote{It may happen $\mathcal K_1(e)
\ne \mathcal K_2(e)$ for d-type edge $e$.}
\item
For each $e \in C^s_1(\mathcal T)$ we are given $p(e)$
such that:
\begin{enumerate}
\item
If $L_{\mathcal K_1(e)} \ne L_{\mathcal K_2(e)}$ then $p(e) \in L_{\mathcal K_1(e)} \cap L_{\mathcal K_2(e)}$.
\item If
$L_{\mathcal K_1(e)} = L_{\mathcal K_2(e)}$ then $p(e) \in (\tilde L_{\mathcal K_1(e)} \times_X \tilde L_{\mathcal K_2(e)})
\setminus  \tilde L_{\mathcal K_1(e)}$.  Namely $p(e)$ is a self-intersection point.
\end{enumerate}
\item We define
\begin{equation}\label{veckapadef}
\vec{\kappa}(v) = (\mathcal K_1(e_{0}(v)),\dots,\mathcal K_{1}(e_{m(v)}(v))).
\end{equation}
We define $i_0,\dots,i_{K'(v)} \in \mathscr K$ 
as follows.  $K'(v)+1$ is the number of edges of $s$-type containing $v$.
We enumerate the edges containing $v$ as above.
Then $e_{i_j}(v)$ is an  edge of $s$-type for $j=0,\dots,K'(v)$.
\par
We define $\vec p(v)$ by $p_{i_j}(v): = p(e_{i_j}(v))$. 
In case all the edges containing $v$ is of $d$-type we put $\vec p(v)  = \emptyset$.  ($K'(v)=-1$.)
\item
$(\vec{\kappa}(v),\vec p(v))$ determine an element of $\widetilde{\text{\rm Seq}}_{m(v)}$
as in Definition \ref{newdef926}.\footnote{In Definition \ref{newdef926} Lagrangians are labeled by 
$\{1,\dots,\#\mathscr L\}$ here it is by $\mathscr K$.}
From $(\vec{\kappa}(v),\vec p(v))$, we obtain a reduced sequence $\vec{\kappa}'(v)$, a sequence of non-negative integers $\vec k'(v)$, and a 
non-negative integer $m'(v)$ as in Equation (\ref{eq:reduced_sequence}): 
 %\marginpar{MA190629: I don't understand the notation. What does $(0, \ldots, 0)$ stand for?
%$,(0,\dots,0)$ is removed  Some explation is added right after the next definition. KF.}
$$
{\rm Red}(\vec{\kappa}(v)) = (\vec{\kappa}'(v),\vec k'(v),m'(v)).
$$
\item
For each $v$ an element 
$B(v) \in \Pi_2(\vec{\kappa}(v),\vec p(v))$ is given.
\item
We also require  the following condition: 
if $e$ is an interior edge and $\mathcal K(e) = (\kappa,\kappa')$ then 
$L_{\kappa} = L_{\kappa'}$.
\end{enumerate}
\end{defn}

We explain a geometric background of this definition.
A decorated ribbon tree as in Definition \ref{decribbon} can be used 
to describe a combinatorial type of a bordered stable disk.
Suppose $(\Sigma,(z_0,z_1,\dots,z_{N}))$ is a tree like union of disks 
with double points together with boundary marked points enumerated in counter-clock-wise 
order. Let $u : \Sigma \to X$ be a pseudo-holomorphic map with a certain 
boundary condition, described below.
Its combinatorial type is described by 
$\Gamma = (\mathcal T,\mathcal K(\cdot),p(\cdot),B(\cdot))$.
An interor vertex $v$ of the $\mathcal T$ is identified with an 
irreducible component $\Sigma_v$ of $\Sigma$.  An exterior vertex is 
identified with one of the marked points.
The exterior edge joins an exterior vertex corresponding to a
marked point $z_i$ 
to the vertex $v$ such that $z_i \in \Sigma_v$.
An interior edge $e$ corresponds to a double point $z_{e}$
and joins two vertices $v,v'$ such that $z_{e} \in \Sigma_v \cap \Sigma_{v'}$.
$B(v)$ is the homology class of the map $u$ restricted 
to $\Sigma_v$.
The edge $e$ is of type $s$ if $u$ sends $z_{e}$ to an intersection point 
of two Lagrangian submanifolds or is a self-intersection point of an immersed Lagrangian submanifold, 
that is, the point $p(e)$.
The edge $e$ is of type $d$ if $u$ sends $z_{e}$ to one of the Lagrangian submanifolds.
Let $v$ be one of the interior vertices  and $\Sigma_v$ be the 
corresponding irreducible component (that is, a disk).
$\Sigma_v$ contains $m(v)+1$ double or marked points $z_{e_{i}(v)}$ ($i=0,\dots,m(v)$)
that corresponds to the edges $e_{i}(v)$.
We denote them by $z_i \in \partial \Sigma_v = S^1$.  
They are 
on the boundary in the counter-clock-wise orientation.
The part of $\partial \Sigma_v$ between $z_i,z_{i+1}$ is denoted by $\overline{z_iz_{i+1}}$.
We require $u$ sends $\overline{z_iz_{i+1}}$ to
$L_{\mathcal K_{2}(e_i(v))}$. 
\par\medskip
\begin{figure}[h] 
\centering
\includegraphics[scale=0.5]{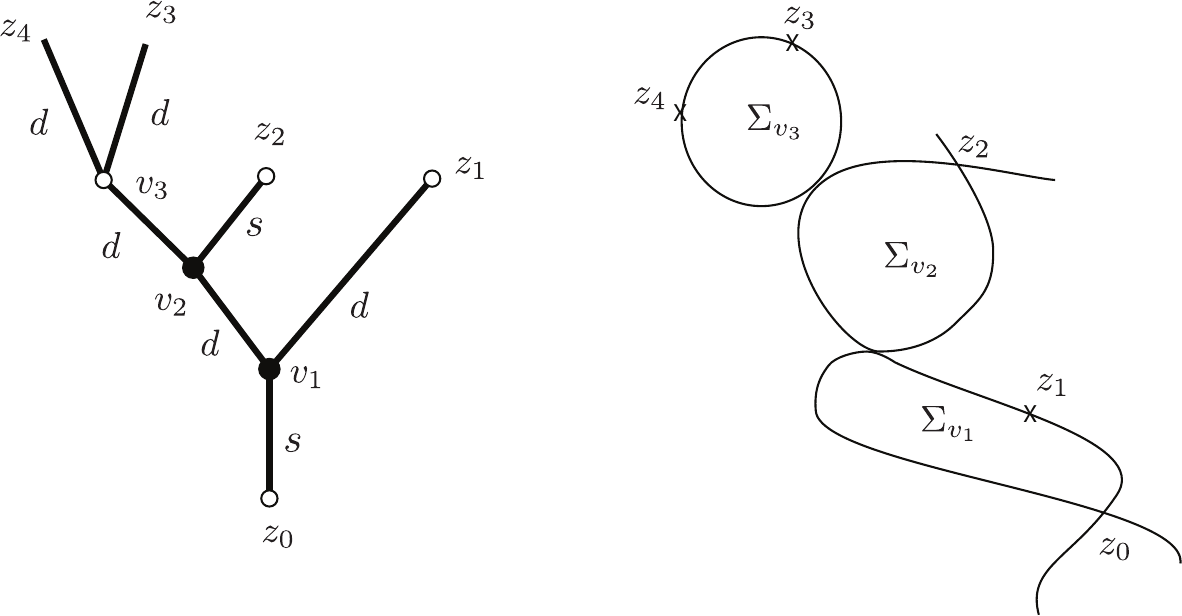}
\caption{Ribbon tree and stable disk}
\label{Figuredef11.1}
\end{figure}
We first consider the part which does not contain ${\bf e}^+_{L_{\kappa}}$ or ${\bf f}_{L_{\kappa}}$,
that is the cyclic $A_{\infty}$ category in Section \ref{sec:cyclicfil}.
\par
For a while we also assume that $L_{\kappa} \ne L_{\kappa'}$ for $\kappa \ne \kappa'$.
(Namely we exclude the case $L_{\kappa} = L_{\kappa'}$, $\theta_{\kappa} \ne \theta_{\kappa'}$
for a while.)
We remove this assumption in Subsection \ref{localcoefficient}.
\par
For each $L_{\kappa}$ 
we fix a linear space $H^k(L_{\kappa};\F)$ of closed differential forms 
on $L_{\kappa}$ so that elements of them uniquely represent the de Rham cohomology
classes. (For example we may fix a Riemannian metric and 
take the set of harmonic forms.)
We also  choose maps 
$\Pi_{L_{\kappa}} : \Omega^k(\tilde L_{\kappa}) \to \Omega^k(\tilde L_{\kappa})$
and $G_{L_{\kappa}} : \Omega^k(\tilde L_{\kappa}) \to \Omega^{k-1}(\tilde L_{\kappa})$ ($\Pi$ and $G$ respectively stand for projector and Green's kernel), 
satisfying the following conditions:\index[syindex]{PiLkappa@$\Pi_{L_{\kappa}}$}
\index[syindex]{GLkappa@$G_{L_{\kappa}}$}
\begin{conds}\label{cond112}
$  $ \par
\begin{enumerate}
\item $\Pi_{L_{\kappa}} \circ \Pi_{L_{\kappa}} = \Pi_{L_{\kappa}}$.
The image of $\Pi_{L_{\kappa}}$ is  $\bigoplus_k H^k(\tilde L_{\kappa};\F)$.
\item 
$\text{\rm identity} - \Pi_{L_{\kappa}} = -(\frak m_{1,0}\circ G_{L_{\kappa}} + G_{L_{\kappa}} 
\circ \frak m_{1,0})$.
\item $G_{L_{\kappa}}\circ G_{L_{\kappa}} = 0$.
\item $\langle \Pi_{L_{\kappa}} x, y \rangle_{\rm cyc} = \langle x,\Pi_{L_{\kappa}} y\rangle_{\rm cyc}$.
\item
$\langle y,G_{L_{\kappa}} x\rangle_{\rm cyc} = (-1)^{\deg x\deg y}\langle x,G_{L_{\kappa}} y\rangle_{\rm cyc}$.
\end{enumerate}
Here $\frak m_{1,0}$ is the part of $\frak m^{\text{\rm f.u.}\text{\bf b}}_1$ which does not contain $T$.
It coincides with de Rham differential $d$.
\end{conds}
The existence of such maps $G_{L_{\kappa}}, \Pi_{L_{\kappa}}$ is classical and is proved, for example, in 
\cite[Lemma 10.1]{fukaya:cyc}.\footnote{The definition of $\langle *,*\rangle_{\rm cyc}$ 
and of $\frak m_{1,0}$ in \cite{fukaya:cyc} is different from ours by sign.  However the 
proof of \cite[Lemma 10.1]{fukaya:cyc} can be adapted to our situation with minor change.}
 
\par
Let $\Gamma = (\mathcal T,\mathcal K(\cdot),p(\cdot),B(\cdot))$ be a  decorated 
ribbon tree and $e^{\text{\rm ext}}_0,\dots,e^{\text{\rm ext}}_N$
be its exterior edges. 
By (\ref{compatibility}) we have
\begin{equation}\label{compatibility2}
\mathcal K_2(e^{\text{\rm ext}}_i) = \mathcal K_1(e^{\text{\rm ext}}_{i+1}), 
\qquad i = 0,1,\dots, N,
\end{equation}
where $e^{\text{\rm ext}}_{N+1} = e^{\text{\rm ext}}_{0}$
by convention.
We define

\begin{equation}\label{vecGammakappa}
\vec{\kappa}(\Gamma) = (\mathcal K_1(e^{\text{\rm ext}}_{0}),
\dots,\mathcal K_{1}(e^{\text{\rm ext}}_{N}))
\in \uwave{K}^{N+1}.
\end{equation}
We also define
$$
\rho_{\frak b,\theta}(\Gamma) =
\prod_{v\in C_0^{\text{\rm int}}({\mathcal T})}\rho_{\frak b,\theta}(B(v)),
\quad
\rho_{b}(\Gamma) =
\prod_{v\in C_0^{\text{\rm int}}({\mathcal T})}\rho_{b}(B(v)).
$$
Moreover we put
$$
\omega(\Gamma) 
= \prod_{v\in C_0^{\text{\rm int}}({\mathcal T})}\omega(B(v)).
$$
We put\index[syindex]{BCFanLkappa@$BCF(\cL;\vec{\kappa})$}
$$
BCF(\cL;\vec{\kappa})
= \bigotimes_{i=1}^N CF^{\text{\rm can}}(L_{\kappa_{i-1}},L_{\kappa_{i}};\F) \blue{[1]}
$$
where\index[syindex]{CFanLkappa@$CF^{\text{\rm can}}(L_{\kappa_{i-1}},L_{\kappa_{i}};\F)$}
$$
CF^{\text{\rm can}}(L_{\kappa_{i-1}},L_{\kappa_{i}};\F)
= 
\begin{cases}
\displaystyle
\bigoplus_{p \in L_{\kappa_{i-1}} \cap L_{\kappa_{i}}} \F[p]
& \text{if $L_{\kappa_{i-1}} \ne L_{\kappa_{i}}$}, \\
\displaystyle H(L_{\kappa_i};\F)\oplus \bigoplus_{p \in (\tilde L_{\kappa_i} \times_{X} \tilde L_{\kappa_i}) \setminus \tilde L_{\kappa_i}} \F[p]
& \text{if $L_{\kappa_{i-1}} = L_{\kappa_{i}}$}.
\end{cases}
$$
We will define\index[syindex]{fGamma@$\frak f^{\bf b}_{\Gamma}$} 
\begin{equation}\label{form115}
\frak f^{\bf b}_{\Gamma} : 
BCF(\cL;\vec{\kappa}(\Gamma))
\to CF(L_{\kappa(\Gamma)_{0}},L_{\kappa(\Gamma)_{1}};\F) \blue{[1]}
\otimes \Lambda_0,
\end{equation}
and \index[syindex]{mxGamma@$\frak m^{\bf b}_{\Gamma}$} 
\begin{equation}\label{form116}
\frak m^{\bf b}_{\Gamma} : 
BCF(\cL;\vec{\kappa}(\Gamma))
\to CF^{\text{\rm can}}(L_{\kappa(\Gamma)_{0}},L_{\kappa(\Gamma)_{1}};\F) \blue{[1]}
\otimes \Lambda_0.
\end{equation}
(Here $\vec{\kappa}(\Gamma)$ is defined in (\ref{vecGammakappa}).)
We then put
\begin{equation}\label{fandmdefcan}
\aligned
\frak f^{\text{\bf b}}_{\vec{\kappa}} = \sum_{\Gamma; \vec{\kappa} = \vec{\kappa}(\Gamma)} 
T^{\omega(\Gamma)}\rho_{b}(\Gamma) \rho_{\frak b,\theta}(\Gamma)
\frak f^{\bf b}_{\Gamma},
\\
\frak m^{\text{\bf b}}_{\vec{\kappa}} = \sum_{\Gamma; \vec{\kappa} = \vec{\kappa}(\Gamma)} 
T^{\omega(\Gamma)}\rho_{b}(\Gamma) \rho_{\frak b,\theta}(\Gamma)
\frak m^{\bf b}_{\Gamma}.
\endaligned
\end{equation}

The maps $\frak f^{\bf b}_{\Gamma}$ 
and $\frak m^{\bf b}_{\Gamma}$ are  defined by induction over the 
number of edges of $\Gamma$ as follows:
\par
As a special case we include $\Gamma=\Gamma_{0,\mathcal K(\cdot),p(\cdot)}$
which has only two (exterior) vertices and one 
exterior edge $e$. It comes with $\mathcal K(e) = (\kappa,\kappa') \in \underline{\kappa}^2$.
In case $L_{\kappa} \ne L_{\kappa'}$, the point $p(e) \in L_{\kappa} \cap L_{\kappa'}$ is included 
in the data defining $\Gamma$.  %\marginpar{I added a figure KF 2024Dec}
$B(v) = 0$ in 
this case.
\begin{figure}[h] 
\centering
\includegraphics[scale=0.3]{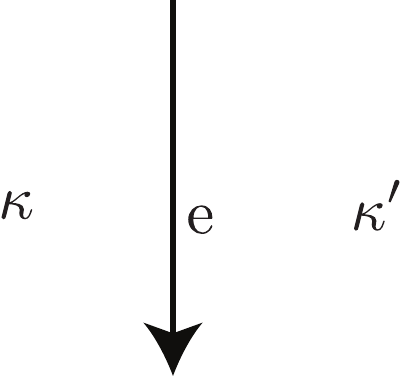}
\caption{$\Gamma_{0,\mathcal K(\cdot),p(\cdot)}$}
\label{Figure(11.6)}
\end{figure}
Now we put
\begin{equation}\label{0thcase}
\frak f^{\bf b}_{\Gamma_{0,\mathcal K(\cdot),p(\cdot)}}  =
\begin{cases}
\text{\rm inclusion} : H(\tilde L_{\kappa};\F) \to  \Omega(\tilde L_{\kappa})
& \text{if $L_{\kappa} = L_{\kappa'}$,}
\\
p \mapsto
\begin{cases}
p  & p=p(e) \\
0 & p\ne p(e)
\end{cases}
& \text{if $L_{\kappa} \ne L_{\kappa'}$,}
\end{cases}
\end{equation}
and
\begin{equation}\label{0thcase2}
\frak m^{\bf b}_{\Gamma_{0,\mathcal K(\cdot),p(\cdot)}}  =
0.
\end{equation}

\par
The next step of our inductive construction is to consider a tree
$\Gamma$ with exactly one interior vertex $v$.
It has $N+1$ exterior edges and $\vec{\kappa}(\Gamma) 
= \vec{\kappa}(v)  \in \uwave{K}^{N+1}$.
We put
\begin{equation}\label{257+1}
\frak f^{\bf b}_{\Gamma} = G_{L_{\kappa}} \circ \frak m_{\vec{\kappa}(v),B(v)}^{{\rm f.u.}\text{\bf b}},
\qquad
\frak m^{\bf b}_{\Gamma} = \Pi_{L_{\kappa}} \circ \frak m_{\vec{\kappa}(v),B(v)}^{{\rm f.u.}\text{\bf b}}.
\end{equation}
\begin{figure}[h] 
\centering
\includegraphics[scale=0.3]{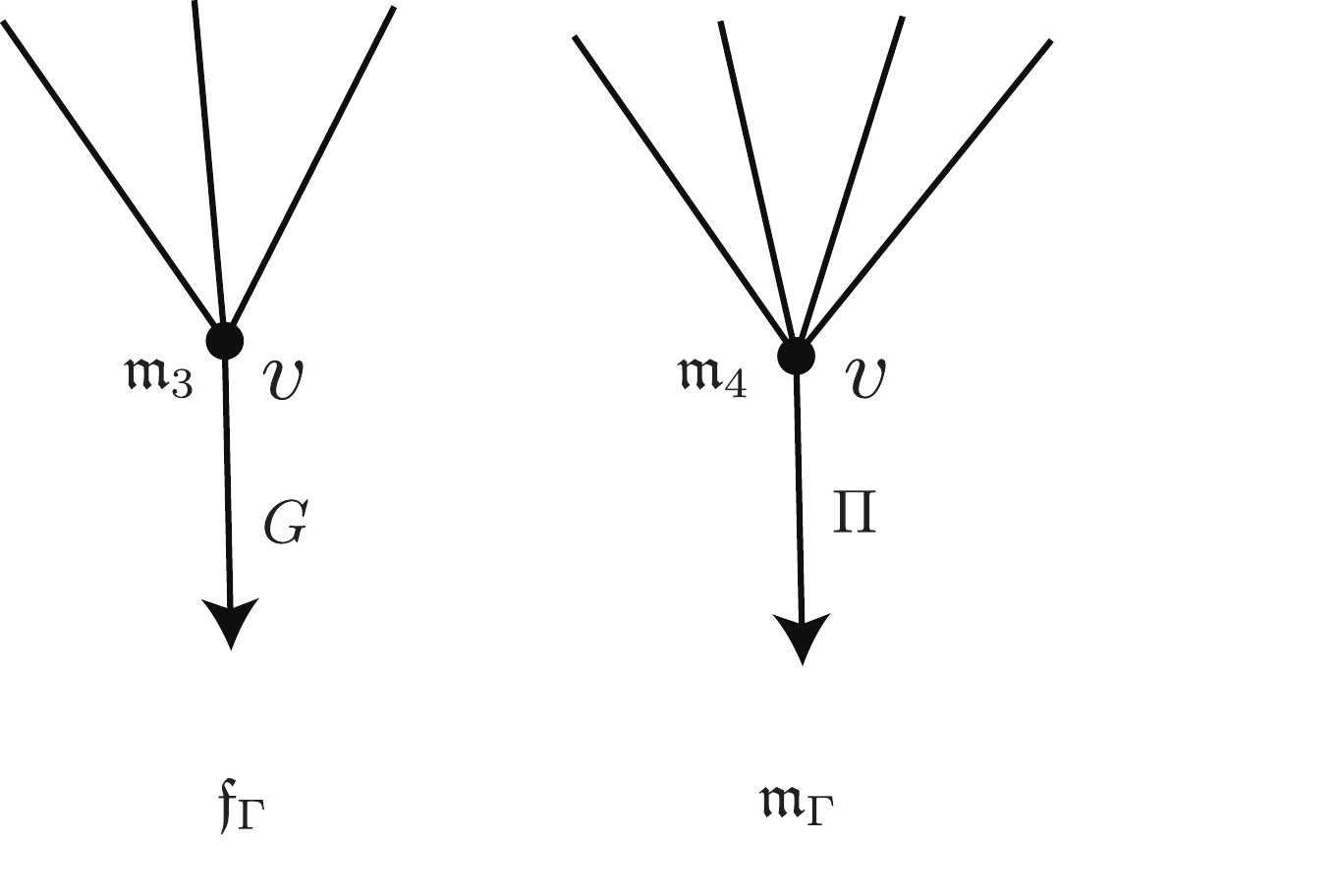}
\caption{$\frak f_{\Gamma}$ and $\frak m_{\Gamma}$}
\label{Figure(11.7)}
\end{figure}
Here we write $\frak m_{\vec{\kappa},B}^{{\rm f.u.}\text{\bf b}}$ for the contribution of curves in homology class $B$ to the $A_\infty$ operation $\frak m^{{\rm f.u.}\text{\bf b}} $ defined in Equation \eqref{eq:mk-operatioin-forms-bounding-cochain}.
 %\marginpar{Figure added. KF 2024 Dec.} 
Explicitly,  each class $B \in \Pi_2(\vec{\kappa}^-(v),\vec p(v))$ contributes a term to the operation  $ \frak q^{{\rm f.u.}\frak b}_{\ell}$ defined in Equation \eqref{qcat}. Using Equations \eqref{eq:a-oo-operation-trivial-gauge} and \eqref{eq:mk-operatioin-forms-bounding-cochain}, we find that the $A_\infty$ operation $\frak m^{{\rm f.u.}\text{\bf b}}$ is also decomposed as a sum indexed by such classes.
\par
Now suppose we have defined $\frak f_{\Gamma'}$ 
and $\frak m_{\Gamma'}$ when $\omega(\Gamma') = \sum_{v\in C_0^{\text{\rm int}}(\Gamma')} \omega(B(\Gamma'))$
is strictly smaller than $\omega(\Gamma)$
or $\omega(\Gamma') = \omega(\Gamma)$ and
the number of edges of $\Gamma'$ is strictly smaller than 
that of $\Gamma$.
We will define $\frak f_{\Gamma}$ 
and $\frak m_{\Gamma}$.
We take the $0$-th exterior vertex $v^{\text{\rm ext}}_0$ and 
the edge $e^{\text{\rm ext}}_0$.
Let $v'$ be the other vertex contained in $e^{\text{\rm ext}}_0$.
We remove $e^{\text{\rm ext}}_0$ from $\Gamma$.
The connected components of the complement 
consists of $k$ ribbon trees $\Gamma_i$.
($i=1,\dots,k$.) Here $k+1$ is the number of edges containing 
$v'$. We regard $v'$ as the root
of $\Gamma_i$. With the other data induced from $\Gamma$ 
in an obvious way, $\Gamma_i$ becomes a 
decorated ribbon tree.  %\marginpar{Include a figure}
We now define
\begin{equation}\label{form119}
\aligned
\frak f^{\bf b}_{\Gamma} &= G_{ L_{\kappa}} \circ \frak m_{\vec{\kappa}(v'),B(v')}^{{\rm f.u.}\text{\bf b}}
\circ \bigotimes_{i=1}^k \frak f^{\bf b}_{\Gamma_i},
\\
\frak m^{\bf b}_{\Gamma} &= \Pi_{ L_{\kappa}} \circ \frak m_{\vec{\kappa}(v'),B(v')}^{{\rm f.u.}\text{\bf b}}
\circ  \bigotimes_{i=1}^k \frak f^{\text{\bf b}}_{\Gamma_i},
\endaligned
\end{equation}
in case $L_{\kappa(\Gamma)_0} = L_{\kappa(\Gamma)_1} = L_{\kappa}$.
Otherwise we put $\frak f^{\bf b}_{\Gamma} = 0$, 
$$
\frak m^{\bf b}_{\Gamma} = \frak m_{\vec{\kappa}(v'),B(v')}^{{\rm f.u.}\text{\bf b}}
\circ  \bigotimes_{i=1}^k \frak f^{\bf b}_{\Gamma_i}.
$$
\par
We thus defined $\frak f^{\bf b}_{\Gamma} $ and 
$\frak m^{\bf b}_{\Gamma}$, and hence $\frak f^{\text{\bf b}}_{\vec{\kappa}}$ and 
$\frak m^{\text{\bf b}}_{\vec{\kappa}}$ 
by (\ref{fandmdefcan}).

\begin{lem}\label{prop230}
$\frak m^{\text{\bf b}}_{\vec{\kappa}}$ defines a structure of 
filtered $A_{\infty}$ category, which we denote by 
$\cL$.\index[syindex]{Lcal@$\sL$}
\par
$\frak f^{\text{\bf b}}_{\vec{\kappa}}$ 
defines a  filtered $A_{\infty}$ functor
$\cL \to \cL^{\rm form}_{\rm{c.u.}}$,
that is a quasi-equivalence.
\end{lem}
\begin{proof}
The proof is the same as that of \cite[Theorem 5.4.1]{fooo09}
and is omitted.
\end{proof}
\begin{lem}\label{lem114}
The (curved) filtered $A_{\infty}$ category {$\cL$}  
is cyclic.
\par
The filtered $A_{\infty}$ functor $\frak f^{\text{\bf b}}_{\vec{\kappa}}$
is cyclic. 
\end{lem}
\begin{proof}
The  proof is the same as the case of filtered $A_{\infty}$
algebra that is given in \cite[Lemmas 10.3 and 10.4]{fukaya:cyc}.
\end{proof}
\begin{rem}\label{lem1115}
If we fix $\lambda$ then 
the filtered $A_{\infty}$ structure in Lemmas \ref{prop230} and \ref{lem114} are curvature free. 
See Subsection \ref{cvnoncomon}.  Then 
$\frak f^{\text{\bf b}}_{\vec{\kappa}}$ 
defines a  filtered $A_{\infty}$ functor
$\cL \to \cL^{\rm form}_{\rm{uni}}$.
\end{rem}
\subsection{Homotopy unit and canonical model.}
\label{hunitcansubsec}
Now we include ${\bf e}^+_{{\kappa}}$ and ${\bf f}_{{\kappa}}$ and  %\marginpar{This subsection is new. KF}
construct a unital and cyclic filtered $A_{\infty}$ category.
In fact we prove the following:
\begin{prop}\label{prop115}
Suppose a gapped filtered $A_{\infty}$ category $\Cat$ is homotopically unital and cyclic, 
then its canonical model is strictly unital and cyclic.
\end{prop}
\begin{proof}
We first modify $\Pi_{L_{\kappa}}$ and $G_{L_{\kappa}}$.
The homotopy unit ${\bf e}_{\kappa}$ is the zero form $1$ on $L_{\kappa}$.
 %\marginpar{More detail 
%of the proof is added. Sign should be checked.  KF 2025 Aug.
%Sign Checked  KF 2025 Aug}
We may take $\Pi_{L_{\kappa}}({\bf e}_{\kappa}) = {\bf e}_{\kappa}$.
Moreover we may take the degree zero part and degree $n=\dim L_{\kappa}$ part of $\Pi_{L_{\kappa}}$ as 
\begin{equation}
\Pi_{L_{\kappa}}(h) =\langle h ,{\rm vol}_{L_{\kappa}}\rangle_{\rm cyc}   {\bf e}_{\kappa}, \qquad \Pi_{L_{\kappa}}(h) =(-1)^n\langle h ,{\bf e}_{\kappa}\rangle_{\rm cyc}   {\rm vol}_{L_{\kappa}}.
\end{equation}
See (\ref{pairing}).
Here we take a Riemannian metric on $L_{\kappa}$ and take $\Pi_{L_{\kappa}}$ to be the 
harmonic projection. The top form  ${\rm vol}_{L_{\kappa}}$ is the volume form.
We change degree $0$ part of the harmonic projection to\index[syindex]{PiprimeLkappa@$\Pi'_{L_{\kappa}}$}
\begin{equation}
\Pi'_{L_{\kappa}}(h) =\langle h ,{\rm vol}_{L_{\kappa}}\rangle_{\rm cyc}   {\bf e}^+_{\kappa}
= (-1)^n \langle {\rm vol}_{L_{\kappa}},h\rangle_{\rm cyc}  {\bf e}^+_{\kappa}.
\end{equation}
We put
$$
\Pi'_{L_{\kappa}}({\bf f}_{\kappa}) = 0.
$$ 
Note that an element of $\alpha = \sum \alpha_i \otimes \alpha'_i \in C \otimes C$  for a 
module $C$ with inner product defines a map 
$\alpha_*: C \to C$ by 
$$
\alpha_*(h) = \sum \langle \alpha'_i, h\rangle  \alpha_i .
$$
Using this notation, the degree $0$ part of $\Pi'_{L_{\kappa}}$ is:
$$
\Pi'_{L_{\kappa}} =  (-1)^n({\bf e}^+_{\kappa} \otimes {\rm vol}_{L_{\kappa}})_*.
$$
We `change'  the degree $n$ part  of $\Pi'_{L_{\kappa}}$ to
$$
\Pi'_{L_{\kappa}} =  ({\rm vol}_{L_{\kappa}} \otimes {\bf e}^+_{\kappa})_*.
$$
Actually the notation $({\rm vol}_{L_{\kappa}} \otimes {\bf e}^+_{\kappa})_*$ is not defined 
since the inner product with ${\bf e}^+_{\kappa}$, which appears, is not defined. This is the reason we put `change' in the quote.
This point will soon be explained.
When we put 
$$
\Delta \Pi_{L_{\kappa}} = (-1)^n({\bf e}_{\kappa}^+ - {\bf e}_{\kappa}) \otimes {\rm vol}_{L_{\kappa}} +{\rm vol}_{L_{\kappa}} 
\otimes ({\bf e}_{\kappa}^+ - {\bf e}_{\kappa})
$$
we might write
$$
\Pi'_{L_{\kappa}} = \Pi_{L_{\kappa}} + (\Delta \Pi_{L_{\kappa}})_*.
$$
We define
$$
\Delta G_{L_{\kappa}} =(-1)^n({ \bf f}_{\kappa} \otimes {\rm vol}_{L_{\kappa}}) - ({\rm vol}_{L_{\kappa}} \otimes {\bf f}_{\kappa}).
$$
Since $d{\bf f}_{L_{\kappa}} = {\bf e}_{L_{\kappa}}^+ - {\bf e}_{L_{\kappa}}$ (where $d$ is the part of the boundary operator which 
does not contain $T$), we have
$$
 ({\bf f}_{\kappa} \otimes {\rm vol}_{L_{\kappa}})_* \circ d +  d \circ ({\bf f}_{\kappa} \otimes {\rm vol}_{L_{\kappa}})_* = 
(({\bf e}_{\kappa}^+ - {\bf e}_{\kappa}) \otimes {\rm vol}_{L_{\kappa}})_*.
$$
In fact

Namely $(-1)^n({\bf f}_{\kappa} \otimes {\rm vol}_{L_{\kappa}})_*$  is a chain homotopy between $\Pi'_{L_{\kappa}}$ and  $\Pi_{L_{\kappa}}$ 
on degree $0$ part.
\par
Moreover using the fact that 
$$
\langle x,(a \otimes b)_*y\rangle_{\rm cyc} = \langle b,y\rangle \langle x,a\rangle_{\rm cyc} 
= (-1)^{\deg' a\deg' x + \deg' b\deg' y}\langle y,(b \otimes a)_*x\rangle_{\rm cyc}
$$
we `have' \footnote{We also remark $$
(d(a \otimes b))_*
=  d\circ (a\otimes b)_*  + (-1)^{\deg a + \deg b}(a \otimes b)_* \circ d. 
$$.} 
\begin{equation}\label{neweq124}
\aligned
&\langle (\Delta \Pi_{L_{\kappa}})_* x,y\rangle_{\rm cyc} =  \langle x,(\Delta \Pi_{L_{\kappa}})_*y\rangle_{\rm cyc} \\
&\langle x,(\Delta G_{L_{\kappa}})_*y\rangle_{\rm cyc} = (-1)^{\deg x\deg y} \langle y,(\Delta G_{L_{\kappa}})_*x\rangle_{\rm cyc}.
\endaligned
\end{equation}
\par
Now we `change' 
$G_{L_{\kappa}}$ to\index[syindex]{GprimeLkappa@$G'_{L_{\kappa}}$} 
$$
G'_{L_{\kappa}} = G_{L_{\kappa}} + (\Delta G_{L_{\kappa}})_*.
$$
We set $G'_{L_{\kappa}}({\bf e}^+_{\kappa})= G'_{L_{\kappa}}({\bf f}_{\kappa}) = 0$.
(Note that $G'_{L_{\kappa}}({\bf e}_{\kappa}) = {\bf f}_{\kappa}$ and is non-zero.)
\par
Actually the notation $(\Delta G_{L_{\kappa}})_*$ is {\it not} defined 
since the inner product with ${\bf f}_{\kappa}$, which appears 
in the second term, is not defined.
This is the reason we put `have' and `change' in the quote.
This point will soon be explained.

We will define\index[syindex]{fGamma+@$\frak f^+_{\Gamma}$} \index[syindex]{mxGamma+@$\frak m^+_{\Gamma}$} 
\begin{equation}\label{eq1214}
\frak f^+_{\Gamma} : 
BCF(\cL;\vec{\kappa}(\Gamma))^+
\to CF(L_{\kappa(\Gamma)_{0}},L_{\kappa(\Gamma)_{1}};\F) \blue{[1]}^+
\otimes \Lambda_0,
\end{equation}
and 
\begin{equation}\label{eq1215}
\frak m_{\Gamma}^+ : 
BCF(\cL;\vec{\kappa}(\Gamma))^+
\to CF^{\text{\rm can}}(L_{\kappa(\Gamma)_{0}},L_{\kappa(\Gamma)_{1}};\F)[1]
\otimes \Lambda_0.
\end{equation}
Note  that there is $^{+}$ in the domain and the codomain, that is to say, we will include ${\bf e}^+_{\kappa}$, ${\bf f}_{\kappa}$
and then take the canonical model.
The definition is by the same induction scheme as before.
First the formula (\ref{0thcase}) and (\ref{0thcase2}) are the same.
For the formula (\ref{257+1}) we change $\Pi_{L_{\kappa}}$ to $\Pi'_{L_{\kappa}}$ and 
$G_{L_{\kappa}}$ to $G'_{L_{\kappa}}$. 
\par
Now we come to the main point of our proof of Proposition \ref{prop115}.
As we mentioned already the operators $(\Delta G_{L_{\kappa}})_*$ and 
$(\Delta\Pi_{L_{\kappa}})_*$ are actually not defined since they contain
an inner product with ${\bf f}_c$ or ${\bf e}^+_c$.
The key idea to resolve this issue is the observation 
that in the definition of $\frak f^+_{\Gamma}$ and $\frak m_{\Gamma}^+$  the expressions 
$(\Delta G_{L_{\kappa}})_*$ and 
$(\Delta \Pi'_{L_{\kappa}})_*$ appear always in the form after they are composed 
with $\frak m_{\vec{\kappa}(v),B(v)}^{\rm f.u.\text{\bf b}}$.
Namely  if we formally apply the definition, the inner product with ${\bf f}_{\kappa}$ appears in the form:
\begin{equation}\label{equation121666-}
( {\rm vol}_{L_{\kappa}} \otimes {\bf f}_{\kappa})_*\circ \frak m_{\vec{\kappa}(v),B(v)}^{{\rm f.u.}\text{\bf b}}({\bf h})
= \langle {\bf f}_{\kappa}, \frak m_{\vec{\kappa}(v),B(v)}^{{\rm f.u.}\text{\bf b}}(\bf h)\rangle_{\rm cyc} {\rm vol}_{\rm L_{\kappa}}
\end{equation}
We remark that in this situation $\frak m_{\vec{\kappa}(v),B(v)}^{{\rm f.u.}\text{\bf b}}({\bf h})$ 
has degree $n+1$ and so is zero.
Moreover the inner product with ${\bf f}_{\kappa}$ is {\it not} defined.
Nevertheless we can and will `define'  (\ref{equation121666-}) as follows. (It can be non-zero).
\par
We consider the `equality'
\begin{equation}\label{equation121666}
\langle  {\bf f}_{\kappa}, \frak m_{\vec{\kappa}(v),B(v)}^{{\rm f.u.}\text{\bf b}}({\bf h})\rangle_{\rm cyc} =
- {}^{+}{\frak m}^{{\rm f.u.}\text{\bf b}}_{\vec{\kappa};B}({\bf f}_{\kappa}\otimes\bf h).
\end{equation}
This would be a consequence of (\ref{form9533}) if ${\bf f}_{\kappa}$ is replaced by 
an element of $CF(L_{\kappa},L_{\kappa};\F)$ (without $^+$) of odd degree.
Therefore we use the right hand side of (\ref{equation121666}) 
as the {\it definition} of the left hand side.
Namely when (\ref{equation121666-}) appears 
in the definition of $\frak f^+_{\Gamma}$ and $\frak m_{\Gamma}^+$
we {\it define}
$$
( {\rm vol}_{L_{\kappa}} \otimes {\bf f}_{\kappa})_*\circ \frak m_{\vec{\kappa}(v),B(v)}^{{\rm f.u.}\text{\bf b}}({\bf h})
= - {}^{+}{\frak m}^{{\rm f.u.}\text{\bf b}}_{\vec{\kappa};B}({\bf f}_{\kappa}\otimes\bf h){\rm vol}_{\rm L_{\kappa}}.
$$
\par
In the case of the second formula in (\ref{257+1}), the composition
$
\Pi'_{L_{\kappa}}\circ \frak m_{\vec{\kappa}(v),B(v)}^{{\rm f.u.}\text{\bf b}}({\bf h})
$
contains a term
$$
({\rm vol}_{L_{\kappa}} \otimes {\bf e}^+_{\kappa})_*\frak m_{\vec{\kappa}(v),B(v)}^{{\rm f.u.}\text{\bf b}}({\bf h})
=
\langle {\bf e}^+_{\kappa},\frak m_{\vec{\kappa}(v),B(v)}^{{\rm f.u.}\text{\bf b}}({\bf h})\rangle_{\rm cyc} {\rm vol}_{L_{\kappa}}.
$$
We {\it define} it as
$$
{}^{+}{\frak m}^{{\rm f.u.}\text{\bf b}}_{\vec{\kappa};B}({\bf e}^+_{\kappa}\otimes{\bf h}){\rm vol}_{L_{\kappa}}.
$$
Thus we can modify (\ref{257+1}) using 
$\Pi'_{L_{\kappa}}$ and  $G'_{L_{\kappa}}$ in place of  $\Pi_{L_{\kappa}}$ and 
$G_{L_{\kappa}}$. 

The inductive step can be carried out in the same way.
Namely we replace $\Pi_{L_{\kappa}}$ and $G_{L_{\kappa}}$ in  (\ref{form119})
by $\Pi'_{L_{\kappa}}$ and $G'_{L_{\kappa}}$. 
The `composition' $( {\rm vol}_{L_{\kappa}} \otimes {\bf f}_{\kappa})_* \circ \frak q^{{\rm f.u.}\text{\bf b}}_{**}$ makes sense by the same reason.
\par
We can check that $(\Delta G_{L_{\kappa}})_*$ is the chain homotopy between 
$\Pi'_{L_{\kappa}}$ and  $\Pi_{L_{\kappa}}$ 
on the intersection of degree $n$ part and the image of $\frak m$ by the 
following calculation:
 %\marginpar{Sign should be put to this calculation.  KF 2025 Aug.
%Sign added KF 2025 Aug.}
$$
\aligned
& (({\rm vol}_{L_{\kappa}} \otimes {\bf f}_{\kappa})_* \circ d)(\frak m^{\rm f.u.\text{\bf b}}(\bf x)) \\
 &= {\rm vol}_{L_{\kappa}} \langle  {\bf f}_{\kappa}, (d\circ \frak m^{\rm f.u.\text{\bf b}}_*(\bf x))\rangle_{\rm cyc} \\
 & = {\rm vol}_{L_{\kappa}} \sum_{c,*' \ne (\kappa,0)} (-1)^{\deg' {\bf x}^{3;1}_c +1}\langle  {\bf f}_{\kappa}, \frak m^{{\rm f.u.}\text{\bf b}}_{*'}({\bf x}^{3;1}_c, \frak m^{{\rm f.u.}\text{\bf b}}_{*}({\bf x}^{3;2}_c),{\bf x}^{3;3}_c)\rangle_{\rm cyc}\\
 & =  {\rm vol}_{L_{\kappa}} \sum_{c,*' \ne (\kappa,0)} (-1)^{\deg' {\bf x}^{3;1}_c} {}^+\frak m^{{\rm f.u.}\text{\bf b}}_{*'}({\bf f}_{\kappa},{\bf x}^{3;1}_c, \frak m^{{\rm f.u.}\text{\bf b}}_*({\bf x}^{3;2}_c),{\bf x}^{3;3}_c)
 \\   & =  {\rm vol}_{L_{\kappa}} ( {}^+\frak m^{{\rm f.u.}{\bf b}}_{(\kappa,0)}({\bf f}_{\kappa},\frak m^{{\rm f.u.}\text{\bf b}}_*({\bf x}))   
 +  {}^+\frak m^{{\rm f.u.}\text{\bf b}}_*({\bf e}_{\kappa}^+,{\bf x}))
  \\   & =  {\rm vol}_{L_{\kappa}} \langle ({\bf e}_{\kappa}^+ - {\bf e}_{\kappa}),\frak m^{{\rm f.u.}\text{\bf b}}_*({\bf x})\rangle_{\rm cyc}.
  \endaligned
 $$
 Here  the second equality is $A_{\infty}$ relation of $\frak m$, the third equality is the definition of $\langle {\bf f}, \cdot \rangle_{\rm cyc}$, 
 the fourth equality is Definition \ref{def98} (d).
In the last line, we put
 $
  \langle {\bf e}_{\kappa}^+,\frak m^{{\rm f.u.}\text{\bf b}}_*({\bf x})\rangle_{\rm cyc}
  = 
  {}^+\frak m^{{\rm f.u.}\text{\bf b}}_*({\bf e}_{\kappa}^+,{\bf x})
 $
 as definition.

Now it is straightforward to check that the resulting maps 
have the required properties.
(The proof of cyclic symmetry is similar to the proof of 
\cite[Proposition 10.1]{fukaya:cyc}, using the fact that (\ref{neweq124}) 
hold in case $x,y$ are in the image of $\frak m$.)
\end{proof}

\begin{cor}\label{cor117}
The gapped cyclic filtered $A_{\infty}$ category $\cL$ in Lemma $\ref{lem114}$
is strictly unital.
\end{cor}

\subsection{The case when $L_{\kappa} = L_{\kappa'}$ but $\theta_{\kappa} \ne \theta_{\kappa'}$.}
\label{localcoefficient}

We consider the case  $L_{\kappa} = L_{\kappa'}$ but $\theta_{\kappa} \ne \theta_{\kappa'}$
then we need to replace the de Rham differential $d$ by  %\marginpar{Subsection added. 2025 Jan. KF}
$
d + (\theta_{\kappa'}-\theta_{\kappa}) \wedge
$. (See (\ref{locacoeffdif}).) We put $d' = d + (\theta_{\kappa'}-\theta_{\kappa}) \wedge$.
The $L^2$ adjoint of $d'$ is
$
\delta' = (-)^{n+np+1} * d' *$ on $p$ forms. 
Hence we replace the usual Laplacian by $\Delta' = - d'\delta'  - \delta' d'$ 
and use Harmonic analysis in the same way.
Thus we can obtain $\Pi_{\kappa,\kappa'}$, $G_{\kappa,\kappa'}$
satisfying a similar condition as 
Condition \ref{cond112}.
We use them to modify the argument of Subsection \ref{subsec:red} 
in an obvious way\footnote{Namely we use $G_{\kappa,\kappa'}$ in place of $G$ 
for the edge $e$ such that 
$\mathcal K(e) =(\kappa,\kappa')$.} to include such cases.
\par
The image of $\Pi_{\kappa,\kappa'}$ is 
the cohomology $H^*(\tilde L_{\kappa};d')$, that is, the 
de Rham cohomology with local coefficient.
Therefore $H^*(\tilde L_{\kappa};d') \otimes_{\F} \Lambda_0$ is (the diagonal part of) the 
underlying module of $\mathcal L((L_{\kappa},\theta_{\kappa}),(L_{\kappa'},\theta_{\kappa'}))$.\footnote{
In the case when $L_{\kappa}$ is immersed, there appears another generator corresponding to the 
switching point.  The operators on those generators are defined in the same way as before.}
The proof of Theorem
\ref{thm11} is now complete.
\qed
\par\medskip
We remark that we assumed that $i_L^*[\frak b_2] = 0 \in H_2(L;\F)$.
(This is necessary for the connection $\theta_L$ of the trivial bundle 
to exist.)
On the other hand, the cohomology class $i_L^*[\frak b_+] \in H_2(L;\Lambda_+)$
is not assumed to be zero.  In fact sometimes we use non-trivial $i_L^*[\frak b_+] \in H_2(L;\Lambda_+)$
to cancel the obstruction class $\frak m_0(1)$ (for the existence of the bounding cochains).
(This is the way to show that if $i_L^* : H^{\rm even}(X;\F) \to H^{\rm even}(L;\F)$ 
is surjective then $L$ is unobstructed after suitable bulk deformation.
See \cite[Theorem C]{fooo09}.) 

\subsection{Note on divisor axioms.}
\label{divisor}

In (\ref{decomposefrakb}) we decompose the bulk class $\frak b$  %\marginpar{Subsection added. KF 2024 Dec}
to the sum $\frak b_0 + \frak b_2 + \frak b_+$ such that $\frak b_2 \in H^2(X;\F)$ 
and $\frak b_0  \in H^0(X;\Lambda_0)$.
The reason for decomposing the bulk class this way
is related to the  convergence issue as explained in Remark \ref{rem61313}.

Also in Definition \ref{boundingcochain}, 
the bounding cochain $b$ is 
decomposed as $b = b_1 + b_+$ where 
$b_1$ is a closed one form that does not contain $T$.
Then we replace $\theta_L$ by $\theta_L + b_1$.
So $b_1$ is used in the different way from $b_+$.
See (\ref{form:716}).  
We remark replacing $\theta_L$ by $\theta_L + b_1$ 
changes $\beta\cap (\frak b,\omega_L)$ 
by $\partial\beta\cap b_1$.

The {\it divisor axiom}\index{divisor axiom} for an ambient cycle is written as 

\begin{equation}\label{divbulk}
\sum_{k=0}^{\infty} \frac{1}{k!}\frak q^{\rm form}_{\beta}(\frak a^{k};{\bf x}) = 
\exp(\beta \cap \frak a)\frak q^{\rm form}_{\beta}(1;{\bf x})
\end{equation}
for  $\frak a \in H^2(X;\Lambda_+)$ for which the left hand side converges in $T$-adic topology.
Therefore if we split $\frak b_+$ into 
$\frak b_2^+ + \frak b_{\rm high}$
where $\frak b_2^+$ is degree $2$ and $\frak b_{\rm high}$ is in higher degree,
then (\ref{divbulk}) would imply that 
\begin{equation}\label{613rev}
\aligned
\text{RHS of }(\ref{mkdefeq}) = \sum_{\stackrel{\beta\in H_2(X,L:\Z)}{\omega(\beta) \leq E}}
\sum_{\ell=0}^{\infty} T^{\omega\cap \beta} 
&\frac{\exp((\frak b_2+\frak b_2^+,\theta_L) \cap \beta) 
  )}{\ell!}\\
&\frak
q^{\rm form}_{\ell,k;\beta}(\frak b_{\rm high}^{\otimes\ell};
x_1,\ldots,x_k).
\endaligned
\end{equation}
(Here for simplicity we assume $\frak b_2^+ = 0$ on $L$.)
(\ref{613rev}) may look more natural since 
put the degree 2 parts $\frak b_2$ and $\frak b_2^+$ together here.
\par
The divisor axiom for a cycle in $L$ is written as:
\begin{equation}\label{divbundary}
\aligned
&\sum_{n_0,\dots,n_k=0}^{\infty} \frak q^{\rm form}_{\beta}({\frak a};b_{1,+}^{n_0}, x_1, 
b_{1,+}^{n_1},\dots,b_{1,+}^{n_{k-1}},x_k,b_{1,+}^{n_{k}})\\
&=
\exp(b_{1,+}\cap \partial \beta)\frak q^{\rm form}_{\beta}({\frak a};x_1, 
,\dots,x_k).
\endaligned
\end{equation}
Here $b = b_1 + b_{1,+} + b_{\rm high}$, $b_{1,+} \in H^1(L;\Lambda_+)$,
$b_{\rm high} \in \bigoplus_{k\ge 1} H^{2k+1}(L;\Lambda_0)$.

We put
$$
\theta'_L = \theta_L + b_1 + b_{1,+}
$$
and regard it as a connection of $\Lambda_0^*$ valued (rank one) bundle 
on $L$. 
Then we replace (\ref{capwiththetaL}) by
\begin{equation}
(\frak b_2,\theta'_L) \cap \beta
:= \int_{\Sigma} u^*\frak b_2 + \int_{\partial\Sigma} u^*\theta'_L \in \Lambda_0
\end{equation}
and (\ref{mkdefeq}) by
\begin{equation}\label{mkdefeqprime}
\aligned
&\frak m_{k}^{\prime \frak b,\theta'_L 
}(x_1,\ldots,x_k) \\
& := \sum_{\stackrel{\beta\in H_2(X,L:\Z)}{\omega(\beta) \leq E}}
\sum_{\ell=0}^{\infty} T^{\omega\cap \beta}
\frac{\exp((\frak b_2,\theta'_L) \cap \beta) 
  )}{\ell!}
\frak
q^{\rm form}_{\ell,k;\beta}(\frak b_{+}^{\otimes\ell};
x_1,\ldots,x_k).
\endaligned
\end{equation}
Then (\ref{divbundary}) would make two definitions coincide in that
$$
\frak m_{k}^{\prime \frak b,\theta'_L 
}(x_1,\ldots,x_k)
=
\frak m_{k}^{\frak b,b_1+ b_{1,+}
}(x_1,\dots,x_k).
$$
Here the right hand side is the $A_{\infty}$ operation deformed 
by the bounding cochain $b_1+ b_{1,+}$, that is,
$$
\sum_{j_0,\dots,j_k=0}^{\infty} \frak m_{k+j_0+\dots + j_k}^{\frak b,\theta_L+b_1}(b_{1,+}^{j_0},x_1,b_{1,+}^{j_1},\dots,
b_{1,+}^{j_{k-1}},x_k,b_{1,+}^{j_k}).
$$
In (\ref{mkdefeqprime}) we did not decompose $\theta'_{L}$ to 
$\theta_L + b_1$ and $b_{1,+}$. So it could be more natural.

However it turns out that we cannot make both divisor aximons (\ref{divbulk}) and 
(\ref{divbundary})  hold simultaneously in our situation.
Formula (\ref{divbulk}) holds if Kuranishi structure and CF-perturbation 
can be chosen to be compatible with forgetful map of interior marked points. On the other hand, formula
(\ref{divbundary})  holds  if Kuranishi structure and CF-perturbation 
can be chosen to be compatible with forgetful map of boundary marked points.
Both compatibilites do not hold in general in our situation.\footnote{Divisor 
axiom (\ref{divbundary}) holds when we study only finitely many 
mutually disjoint set of embedded Lagrangian submanifolds. This is the case of \cite{toric3}.
See \cite[Section 4.1]{toric3}.}

In this regard, the bulk deformation defined by (\ref{613rev}) is 
not the same as   those defined in (\ref{mkdefeq}).
The deformation by  weak bounding cochains defined by  
(\ref{mkdefeqprime}) is not the same as 
$\frak m_{k}^{\frak b,b_1+ b_{1,+}}$.
Nevertheless, even we use (\ref{613rev}) and/or (\ref{mkdefeqprime}) 
we still obtain the structure which satisfies the same properties as one we use in this paper.

However some properties are lost when we use (\ref{613rev}) and/or (\ref{mkdefeqprime}).

Let us first discuss the case of bulk deformation.
Let $\frak b(s): = \frak b + s \Delta \frak b$ be an $s \in (-\epsilon,+\epsilon)$ parametrized 
family of bulk deformations. We assume $\frak b \in H^{\rm even}(X;\Lambda_0)$
and $\Delta\frak b \in H^{\rm 2}(X;\Lambda_+) \oplus H^{\rm even > 2}(X;\Lambda_0)$.
Then it is easy to see from (\ref{mkdefeq}) that
\begin{equation}\label{1777}
\frac{d}{ds}\frak m_{k}^{\frak b(s)}(x_1,\ldots,x_k)\vert_{s=0}
= \frak q_{1;k}^{\rm f.c.u.\frak b}(\Delta\frak b;x_1,\ldots,x_k).
\end{equation}
See (\ref{qqmcat}), (\ref{eq145+}) for the right hand side. Actually since the left hand side is homology version 
we need to take the canonical model version of (\ref{qqmcat}), (\ref{eq145+})  for the right hand side.
Formula (\ref{1777}) does not hold if we use (\ref{613rev}).
If we take $\Delta \frak b \in H^{2}(X;\Lambda_0)$
instead of $\Delta\frak b \in H^{2}(X;\Lambda_+)$ 
then (\ref{1777})  does not hold either.
See Subsection \ref{perturbcomp}.

We next consider the case of deformation by bounding cochain.
We put $b(s) = b + s\Delta b$ with 
$\Delta b \in H^{\rm 1}(L;\Lambda_+) \oplus H^{\rm odd > 1}(L;\Lambda_0)$.
Then it is easy to see from the definition that:
\begin{equation}\label{form1118}
\frac{d}{ds} \frak m_0^{\frak b,b(s)}(1)\vert_{s=0} =  \frak m_1^{\frak b,b}(\Delta b).
\end{equation}
However if we use  (\ref{mkdefeqprime}) to define $\frak m_k^{\frak b,b(s)}(1)$,
Formula (\ref{form1118}) does not hold.
However we can choose our system of perturbations so that (\ref{form1118}) holds mod $T^E$ 
for any fixed $E$.

In \cite{toric3}, two formulas (\ref{1777}) and (\ref{form1118}) played very important roles.
(\ref{1777}) together with the fact that $\frak q$ is a ring homomorphism 
is used to prove that Kodaira-Spencer map from the ambient quantum cohomology 
to the Jacobian ring is a ring homomorphism.
(\ref{form1118}) is used to show a relation between criticality of the potential function and 
non-vanishing of Lagrangian Floer cohomology.
We will go back to the discussion of this point at Subsection \ref{perturbcomp}.

\subsection{The case when the critical values are  common.}
\label{cvnoncomon}

In our discussion so far we did not assume that  %\marginpar{The name of subsection is changed.}
$\frak{PO}_{\frak b}(b)$ is independent of  $(L,\theta,b) \in {\bf L}$ 
in our collection of Lagrangian submanifolds equipped with a  bounding cochain $b$.
($\theta$ include the leading order term of $b$.)

Then as we mentioned before we have
\begin{equation}\label{mmne0}
\frak m^{\bf b}_1 \circ \frak m^{\bf b}_1 = (\frak{PO}_{\frak b}(b) - \frak{PO}_{\frak b}(b'))
\,\, \text{\rm id}
\end{equation}
on $CF((L,\theta,b),(L',\theta',b');\Lambda_0)$.
%So precisely speaking we have a curved filtered $A_{\infty}$ category.\marginpar{
%I feel the description around here is already written before in this paper.}  
%\newred{See \cite[Definition 2.2]{fukaya:functor} for its definition, for example.}
 %\marginpar{include the definition?}  

\par
Note in a curved  filtered $A_{\infty}$ category the operator $\frak m_0$ may not be $0$.
In fact
$ \frak m^{\bf b}_0 = 
\frak m^{(L,\theta,b)}_0 : \Lambda_0 \to CF((L,\theta,b),(L,\theta,b);\Lambda_0)
$
is given by 
$$
\frak m^{\bf b}_0(1) = \frak{PO}_{\frak b}(b)
\cdot \text{\bf e}_L \in H^0(L;\Lambda_0).
$$
The identity (\ref{mmne0}) is a consequence of the $A_{\infty}$ relations and this formula.
\par
Let us divide the collection $\bL$ into the disjoint union\index[syindex]{Lbflambda@$\bL_{\lambda}$} 
$$
\bL = \bigcup_{\lambda} \bL_{\lambda},
$$
where
$$
\bL \subseteq \{(L, \theta, b) \mid (L,\theta) \in \mathscr L, \,\, b \in \mathcal M_{\rm weak}(\Omega(L),\frak b, \theta;\Lambda_0)\}
$$
and
$$
 \bL_{\lambda} = \{(L_\kappa, \theta_\kappa, b_\kappa) \in \bL \mid \frak{PO}_{\frak b}(b_\kappa) = \lambda\}.
$$
Let $\cL_{\lambda}$ be the full subcategory of $\cL$ whose object is an element of 
$\bL_{\lambda}$. Then if $(L,\theta, b),(L',\theta',b')$ are objects of $\bL_{\lambda}$
we have $\frak m^{\bf b}_1 \circ \frak m^{\bf b}_1 =0$ on $CF((L,\theta,b),(L',\theta',b');\Lambda_0)$.
When we consider the filtered $A_{\infty}$ category $\cL_{\lambda}$
we {\it re-define} $\frak m^{\bf b}_0$ so that $\frak m^{\bf b}_0 = 0$.
Then the facts that $\frak{PO}_{\frak b}(b)$ is independent of \blue{$(L,\theta, b) \in \bL_{\lambda}$} and that
$\frak m^{\bf b}_0(1)$  is 
proportional to the unit imply that the $A_{\infty}$ formula holds with this 
re-defined $\frak m^{\bf b}_0$.
Thus $\cL_{\lambda}$ becomes a {\it curvature-free}  filtered $A_{\infty}$ category.

\begin{rem}
Note that there are two slightly different ways to split a curved $A_{\infty}$ categories 
according to the value of $\frak m^{\bf b}_0$. %\marginpar{Remark added KF.2024Nov}
\begin{enumerate}
\item
One way is to first go to the canonical model (taking all elements of $\bL$ as objects), and then 
split $\bL$ into a sum of $\bL_{\lambda}$. We then define the
$A_{\infty}$ category generated by the subset consisting of the elements of $\bL_{\lambda}$
by redefining $\frak m^{\bf b}_0$ to be $0$.
\item The other way is to   split $\bL$ into a sum of $\bL_{\lambda}$. 
We then define  the $A_{\infty}$ category whose objects are  elements of $\bL_{\lambda}$ by
redefining $\frak m^{{\rm f.u.} \bf b}_0$ to be $0$. Finally
we take the canonical models of them for various $\lambda$ individually.
\end{enumerate}
Actually those two slightly different ways give exactly the same $A_{\infty}$ category $\cL_{\lambda}$.
We now give an explanation on why this is the case.

A priori, the two results may be different from $\cL_{\lambda}$ because of the following reason:
While we take the canonical model with $\frak m^{{\rm f.u.} \bf b}_0 = 0$  in Case (2), 
we take $\frak m^{{\rm f.u.} \bf b}_0 \ne 0$ in Case (1).

However this difference does not change the construction.
Let us elaborate on this point.
The curvature value  $\frak m^{{\rm f.u.} \bf b}_0 = \lambda {\bf e}^+$ affect the definition of 
$\frak f^+_{\Gamma}$ etc. so that the graph $\Gamma$ may contain a tad pole.
(An interior vertex which has one edge.) See Figures \ref{Figure117(i)} and \ref{Figure117(ii)}.
An output of the tad pole (the edge) contain other vertex $v'$.
There are two cases. 
\begin{enumerate}
\item[(i)]
The operation associated to $v'$ is $\frak m^{{\rm f.u.} \bf b}_{2,0}$, 
which is the wedge product (up to sign).
\item[(ii)]
Other cases.
\end{enumerate}
See Figures \ref{Figure117(i)} and \ref{Figure117(ii)}. 
\begin{figure}[h]
\centering
\includegraphics[scale=0.2]{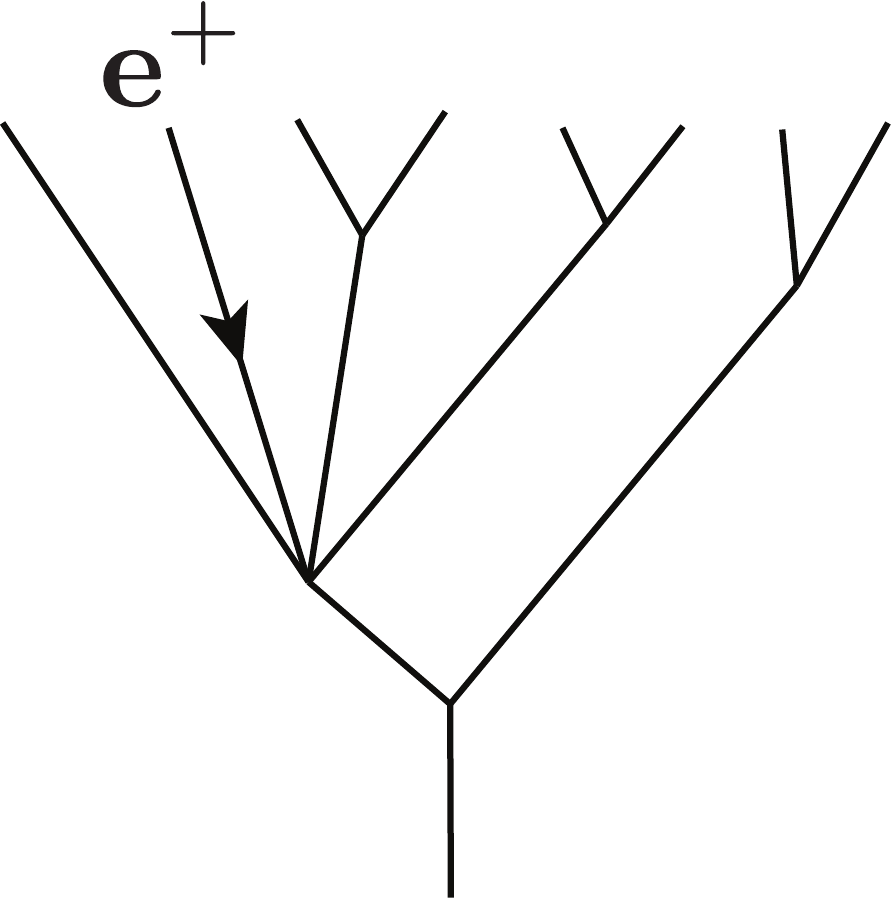}
\caption{Case (ii)}
\label{Figure117(i)}
\end{figure}
\begin{figure}[h]
\centering
\includegraphics[scale=0.2]{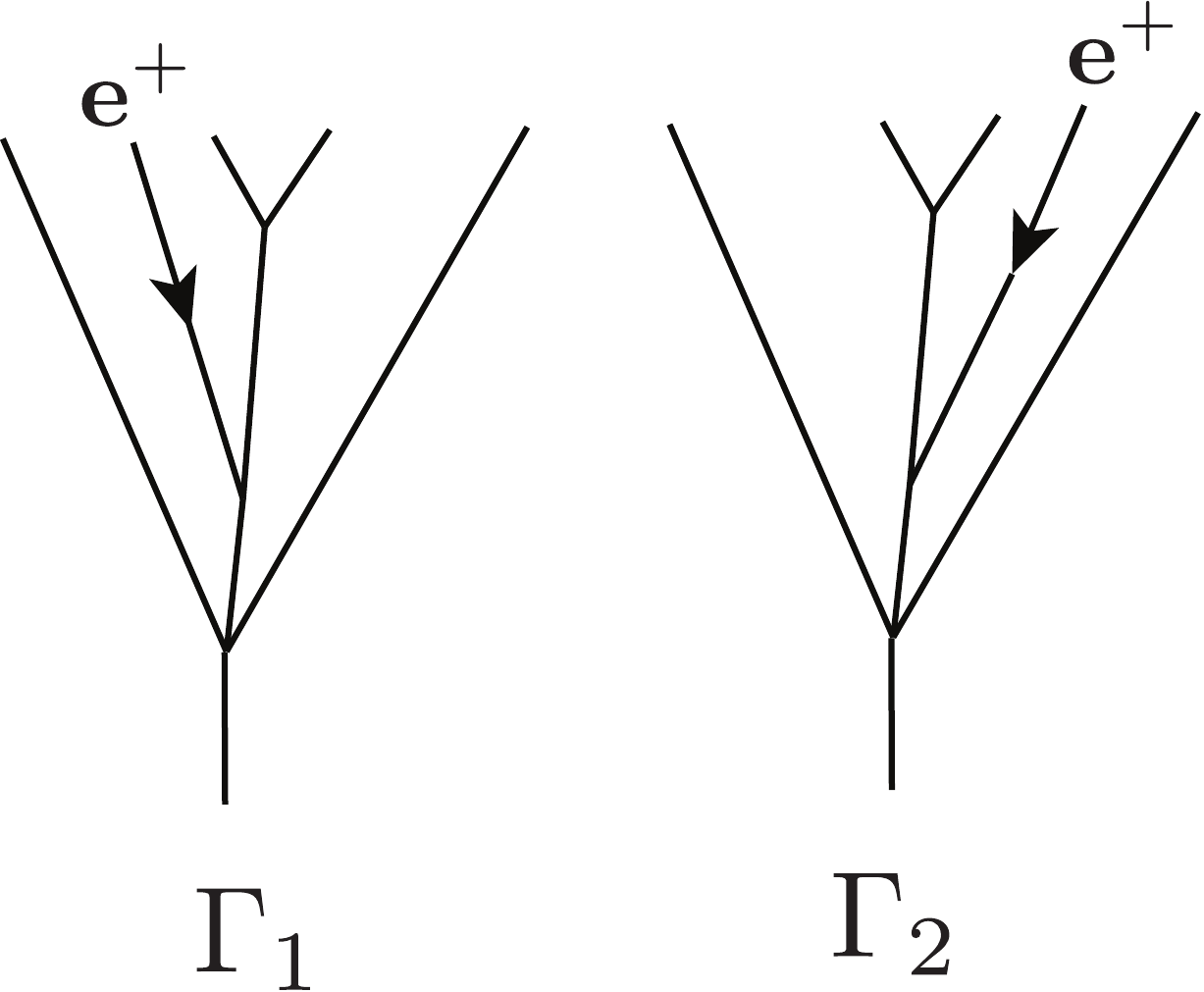}
\caption{Case (i)}
\label{Figure117(ii)}
\end{figure}

In case (ii) the map $\frak f^+_{\Gamma}$ is zero because of the 
definition of the (exact) unit ${\bf e}^+$.
In case (i) there are two diagrams $\Gamma_1$, $\Gamma_2$ as in Figure \ref{Figure117(ii)},
such that $\frak f^+_{\Gamma_1} + \frak f^+_{\Gamma_2} = 0$.
Therefore the tad pole does not affect the construction.
\end{rem}

\subsection{Relative version.}
\label{relative}

We fix a bulk cohomology class $\frak b$ and a background class $\text{\rm st}$.
Let  $\bL$ be a finite collection $\{(L_{\kappa}, \theta_{\kappa}, b_{\kappa}) \mid \kappa=1,2,\dots\}$. 

We take another collection $\bU$ of $\{(U_{\upsilon},\theta_{\upsilon}, b_{\upsilon})\mid 
\upsilon =1,2,\dots\}$.

\begin{prop}\label{sec2relative}
Suppose we are given choices to define the structure of 
unital cyclic filtered $A_{\infty}$ category $\cL$ 
with $\bL$ as the set of objects and 
$\sU$ 
with $\bU$ as the set of objects.

\par
Then we can define the structure of unital cyclic filtered $A_{\infty}$ category
with  $\bL \cup \bU$ as the set of objects,
extending the given structures.
\end{prop}

The proof is immediate from the way of constructing 
unital cyclic filtered $A_{\infty}$ structures.

\subsection{Uniqueness up to cyclic and unital pseudo-isotopy.}
\label{unique}

In this subsection we prove the uniqueness of cyclic and unital filtered 
$A_{\infty}$ category of Theorem \ref{thm11}, 
 %\marginpar{This subsection is new.
%KF. 2025 Dec}
up to cyclic and unital pseudo-isotopy.

We begin with:

\begin{prop}\label{formunique}
The cyclic and homotopically unital filtered $A_{\infty}$ category $\mathcal L^{\rm form}_{\rm c.u.}$
in Proposition $\ref{prop910}$ is independent of the choices (including almost complex structures) 
up to cyclic and homotopically unital pseudo-isotopy.
\end{prop}

\begin{proof}
Uniqueness as cyclic and homotopically unital filtered $A_{n,k}$  category
up to cyclic and homotopically unital pseudo-isotopy of filtered $A_{n,k}$ category 
is a consequence of Proposition \ref{existpisomain} and its analogue 
including homotopical unitality.  The latter can be proved in the same way.
\par
We can slightly modify the proof of Proposition \ref{existpisomain}
to include independence of the compatible almost complex structure as follows.
Let $J,J'$ be two compatible almost compatible structures. 
We can find a one parameter family of almost complex structures 
$\{J_t\}$ joining them.
Then in place of the moduli space $ [0,1]_t\times {\mathcal M}_{\ell;\vec k}((\vec{\kappa},\vec p);B)$
appearing in (\ref{evaluationmaptime}) we take 
\begin{equation}\label{form11242}
\bigcup_{t \in [0,1]} \{t\} \times  {\mathcal M}_{\ell;\vec k}((\vec{\kappa},\vec p);B;J_t)
\end{equation}
where ${\mathcal M}_{\ell;\vec k}((\vec{\kappa},\vec p);B;J_t)$ is the moduli space 
${\mathcal M}_{\ell;\vec k}((\vec{\kappa},\vec p);B)$ obtained 
by using the almost complex structure $J_t$.
\par
We use (\ref{form11242}) in the same way as the proof of  Proposition \ref{existpisomain}
to show independence under the change of almost complex structure 
of cyclic and homotopically unital filtered $A_{n,k}$ category 
up to cyclic and homotopically unital pseudo-isotopy. %\marginpar{Independence of $J$ is discussed here.  KF 2025 July}
\par
To prove uniqueness in the $A_{\infty}$ sense we need to apply 
homotopy limit argument to the pseudo-isotopy.
For this purpose we use pseudo-isotopy of pseudo-isotopies. 
Such an argument is in \cite[Chapter 22]{springer} and \cite[Section 11]{fukaya:cyc} 
in the case of single embedded Lagrangian submanifold.
In the same way as we discussed several times already in this paper, we can 
modify it
to prove Proposition \ref{formunique}.
\end{proof}

We next use Proposition \ref{formunique} to prove the uniqueness of the 
canonical model $\mathcal L$.
Such a uniqueness is discussed in  \cite[Section 12]{fukaya:cyc}.
Actually one point is not discussed there.
It is independence of the harmonic projection $\Pi$ and the propagator 
(Green kernel) $G$.  (See \cite[Remark 12.1]{fukaya:cyc}.)
In \cite[Section 7]{HY}, Hang Yuan proved such uniqueness up to pseudo-isotopy.
He did not discuss cyclic symmetry there. 
However, actually we can easily adapt the argument of  \cite{HY}   to complete the 
proof of the next uniqueness result:

\begin{thm}\label{thm111}
The cyclic and unital filtered $A_{\infty}$ category $\mathcal L$
is independent of various choice involved  (including almost complex structures) up to 
cyclic and unital pseudo-isotopy.
\end{thm}

\begin{proof}
We use Riemannian metric on $L_{\kappa}$ to obtain $\Pi_{\kappa}$ and $G_{\kappa}$.
We will discuss dependence of them on Riemannian metric.
We take two Riemannian metrics $g_0,g_1$ on $L_{\kappa}$ and join them by a
one parameter family $g_t$ $(t \in [0,1])$ of  metrics.
We then obtain  $\Pi_{\kappa,t}$ and $G_{\kappa,t}$.
We decompose\index[syindex]{pzismkappat@$\pi_{\kappa,t}$}\index[syindex]{pikappat@$\Pi_{\kappa,t} $}  
\index[syindex]{ikappat@$i_{\kappa,t}$}
$\Pi_{\kappa,t} = i_{\kappa,t} \circ \pi_{\kappa,t}$
where
$$
\pi_{\kappa,t} : \Omega(L_{\kappa}) \to H(L_{\kappa},\R),
\qquad
i_{\kappa,t} : H(L_{\kappa},\R) \to \Omega(L_{\kappa}).
$$
($i_{\kappa,t}$ is the canonical map to the subspace of $g_t$ harmonic forms.)
\begin{rem}
In Subsection \ref{subsec:red}, we identified $H(L_{\kappa};\R)$ with the set of harmonic forms 
and embed $H(L_{\kappa};\R)$ to $\Omega(L_{\kappa})$ by this identification.
Here we use a family of metrics and so the embedding  of 
$H(L_{\kappa};\R)$ to $\Omega(L_{\kappa})$ depends on $t$.
This is the reason why we decompose $\Pi_{\kappa,t} = i_{\kappa,t} \circ \pi_{\kappa,t}$.
\end{rem}
We use the next Lemma \ref{lem113} (\cite[Lemmas 7.4 and 7.5]{HY}).\index[syindex]{h1kappat@$h_{\kappa,t}$}
\index[syindex]{kkappat@$k_{\kappa,t}$}\index[syindex]{szigmakappat@$\sigma_{\kappa,t}$}
\begin{lem}\label{lem113}
There exists smooth family of operators:
$h_{\kappa,t} : H(L_{\kappa},\R) \to \Omega(L_{\kappa})$,
$k_{\kappa,t} : \Omega(L_{\kappa}) \to H(L_{\kappa},\R)$ and
$\sigma_{\kappa,t} : \Omega(L_{\kappa}) \to\Omega(L_{\kappa})$,
such that $h_{\kappa,t}$, $k_{\kappa,t}$ are of degree $-1$ and 
$\sigma_{\kappa,t}$ is of degree $-2$. They have the following properties.
\begin{equation}\label{form1120}
\aligned
&\frac{d i_{\kappa,t}}{dt} = d \circ h_{\kappa,t}, 
&\frac{d \pi_{\kappa,t}}{dt} = k_{\kappa,t} \circ d, \\
& \pi_{\kappa,t} \circ h_{\kappa,t} = 0, 
 k_{\kappa,t}  \circ i_{\kappa,t} = 0, \\
&  G_{\kappa,t} \circ h_{\kappa,t} = 0, 
\qquad  &k_{\kappa,t}  \circ G_{\kappa,t} = 0, 
\\
& \frac{d G_{\kappa,t}}{dt} - i_{\kappa,t} \circ k_{\kappa,t} 
- h_{\kappa,t} \circ \pi_{\kappa,t} = d \circ \sigma_{\kappa,t} - \sigma_{\kappa,t} \circ d,
\\
&\sigma_{\kappa,t} \circ G_{\kappa,t} = G_{\kappa,t} \circ \sigma_{\kappa,t}  = 0, 
\\
& \pi_{\kappa,t} \circ \sigma_{\kappa,t} = 0,  &\sigma_{\kappa,t} \circ  i_{\kappa,t} =0.
\endaligned
\end{equation}
\end{lem}
Moreover we have:
\begin{lem}\label{lem1113}
For $x,x' \in \Omega(L_{\kappa})$ and $y \in H(L_{\kappa};\R)$ we have
\begin{equation}\label{form1121}
\aligned
\langle k_{\kappa,t}(x),y\rangle_{\rm cyc} &= (-1)^{\deg' x} \langle x, i_{\kappa,t}(y)\rangle_{\rm cyc} \\
\langle \sigma_{\kappa,t}(x),x'\rangle_{\rm cyc} &= \langle x, \sigma_{\kappa,t}(x')\rangle_{\rm cyc}.
\endaligned
\end{equation}  
\end{lem}
We can prove Lemma \ref{lem1113} by examining the construction of those operators given 
in \cite{HY}, which is based on the harmonic analysis. 

Let $H^d(L_{\kappa})_{[0,1]_t}$ be the set of 
elements $a = a_0(t) + dt \wedge a_1(t)$ where 
$a_0(t)$ is a smooth map $[0,1]_t \to H^d(L_{\kappa};\R)$ 
and $a_1(t)$ is a smooth map $[0,1]_t \to H^{d-1}(L_{\kappa};\R)$.
We also consider $\Omega([0,1]_t \times L_{\kappa})$.
Its element is written as 
$a = a_0(t) + dt \wedge a_1(t)$ where 
$a_0(t)$, $a_1(t)$ are smooth sections of  the pull back 
bundle $\pi_2^*\Lambda(L_{\kappa})$  of differential forms on $L_{\kappa}$
to $[0,1]_t \times L_{\kappa}$.
\par
Now we define the `$[0,1]$ parametrized version' of the operators
$i$, $\pi$, $G$ as follows.
\begin{equation}\label{form1122}
\aligned
&i^{\kappa}_{\rm para} : H^d(L_{\kappa})_{[0,1]} \to \Omega^d([0,1] \times L_{\kappa}),
\\
&i^{\kappa}_{\rm para}(a_0(t) + dt \wedge a_1(t))
= i_{\kappa,t}(a_0(t)) + dt \wedge (h_{\kappa,t}(a_0(t)) + i_{\kappa,t}(a_1(t))).
\endaligned
\end{equation}
\begin{equation}\label{form1123}
\aligned
&\pi^{\kappa}_{\rm para} :  \Omega^d([0,1] \times L_{\kappa}) \to H^d(L_{\kappa})_{[0,1]}
\\
&\pi^{\kappa}_{\rm para}(a_0(t) + dt \wedge a_1(t))
=
\pi_{\kappa,t}(a_0(t)) + dt \wedge (k_{\kappa,t}(a_0(t)) + \pi_{\kappa,t}(a_1(t))).
\endaligned
\end{equation}
\begin{equation}\label{form1124}
\aligned
&G^{\kappa}_{\rm para} :  \Omega^d([0,1] \times L_{\kappa}) \to  \Omega^{d-1}([0,1] \times L_{\kappa})
\\
&G^{\kappa}_{\rm para}(a_0(t) + dt \wedge a_1(t))
=
G_{\kappa,t}(a_0(t)) + dt \wedge (\sigma_{\kappa,t}(a_0(t)) - G_{\kappa,t}(a_1(t))).
\endaligned
\end{equation}
This definition is taken from \cite[(129)]{HY}.
\par
Suppose we have a cyclic pseudo-isotopy
$(\frak C^F,\langle\, ,\,\rangle,\{\frak m^{F, t}_{k,\beta}\},\{\frak c^{F, t}_{k,\beta}\})$
as in Definition \ref{evmap2para}.
(This is a pseudo-isotopy between two cyclic filtered $A_{\infty}$ structures 
obtained from two different choices of Kuranishi structures and CF-perturbations.\footnote{
Here they are $A_{\infty}$ structures not $A_{n,k}$ structures.
Namely we already have taken homotopy limit.
To construct a cyclic pseudo-isotopy in $A_{\infty}$ sense, we use pseudo-isotopy of pseudo-isotopies
as in \cite[Subsection 22.2]{springer} or \cite[Section 14]{fukaya:cyc}.})
The morphism module $\frak C^F(L_{\kappa},L_{\kappa'})$
is given by $\Omega([0,1] \times (\tilde L_{\kappa} \times_X \tilde L_{\kappa'}))$.
\par
We use the operators (\ref{form1122}), (\ref{form1123})
and (\ref{form1124})  in place of $i$, $\pi$ and $G$.
We also use  
$$
\frak M^F_{k,\beta} = \frak m^F_{k,\beta} + dt \wedge  \frak c^F_{k,\beta}
$$
in place of $\frak m_{k,\beta}$.
We then proceed in the same way as 
Subsection \ref{subsec:red} to obtain:
$$
\frak M^H_k : B_k{\frak C}^H(L_{\kappa},L_{\kappa'})
\to {\frak C}(L_{\kappa},L_{\kappa'}).
$$
Here $B_k{\frak C}^H(L_{\kappa},L_{\kappa'})$ is the Bar complex of a filtered $A_{\infty}$ category 
$\frak C^H$. The morphism module $\frak C^H(L_{\kappa},L_{\kappa})$ is $H(L_{\kappa})_{[0,1]}$
and if  $L_{\kappa}\ne L_{\kappa'}$ the morphism module $\frak C^H(L_{\kappa},L_{\kappa'})$ is the set of elements 
$a(t)+ dt \wedge b(t)$, where $a(t)$, $b(t)$ are smooth functions on $[0,1]$ to the set $(L_{\kappa} \cap L_{\kappa'}) \times \R$.
(We also use the energy zero term that is obtained from de Rham differential of $L_{\kappa}$, exterior derivative on $[0,1]$
and wedge product.)
\par
We decompose 
$$
\frak M^H_k(h_1,\dots,h_k) = \frak m^{H,t}_k(h_1,\dots,h_k) + dt \wedge  \frak c^{H,t}_k(h_1,\dots,h_k)
$$
where $h_i$ does not contain $dt$ and is independent of $t$.
\par
We can prove the $A_{\infty}$ relation for $\frak M^H_k$ and therefore $\frak m^{H,t}$, $\frak c^{H,t}$
satisfy Definition \ref{pisotopydef} (1)-(3).
Thus they define a pseudo-isotopy.
Moreover we can use Lemma \ref{lem1113} and the fact that
$(\frak C^F,\langle\, ,\,\rangle,\{\frak m^{F, t}_{k,\beta}\},\{\frak c^{F, t}_{k,\beta}\})$ is cyclic
to prove  Definition \ref{pisotopydef} (4), in a way similar to 
the proof of Lemma \ref{lem114} (that is similar to the proof of \cite[Lemmas 10.3 and 10.4]{fukaya:cyc}).
Thus $(\frak C^H,\langle\, ,\,\rangle,\{\frak m^{H, t}_{k,\beta}\},\{\frak c^{H, t}_{k,\beta}\})$
is a cyclic pseudo-isotopy.
\par
We finally discuss unitality of our cyclic pseudo-isotopy.
We first recall that $(\frak C^F,\langle\, ,\,\rangle,\{\frak m^{F, t}_{k,\beta}\},\{\frak c^{F, t}_{k,\beta}\})$ 
is extended to a homotopically unital and cyclic pseudo-isotopy. (Proposition \ref{formunique}.)
We next proceed in the same way as Subsection \ref{hunitcansubsec}.
For this purpose we modify operators (\ref{form1122}), (\ref{form1123})
and (\ref{form1124}).
\par
We used the volume element ${\rm vol}_{L_{\kappa}}$ in Subsection \ref{hunitcansubsec}.
The volume element depends on the metric. Actually 
$$
{\rm vol}^t_{L_{\kappa}} = i_{{\kappa},t}([L_{\kappa}])
$$
is the one parameter family of volume elements we use.
Here $[L_{\kappa}] \in H^n(L_{\kappa};\F)$ is the fundamental cohomology class. 
The (homotopy) unit ${\bf e}_{L_{\kappa}}$ is independent of the metric,
since it is the differential $0$ form $1$ on $L_{\kappa}$. 
We also remark that the inner product $\langle \cdot,\cdot\rangle_{\rm cyc}$ is 
{\it independent} of the Riemannian metric.
(It is the Poincar\'e duality and is not the $L^2$ inner product etc.)
\par
The degree $0$ part of $\pi_{{\kappa},t}$ is
$$
\pi_{L_{{\kappa},t}}(h) = \langle {\rm vol}^t_{L_{\kappa}},h\rangle_{\rm cyc} [{\bf e}^+_{\kappa}].
$$
Here $[{\bf e}^+_{\kappa}]$ is the unit of $H(L_{\kappa};\R)$ which is 
the same as $[{\bf e}_{\kappa}]$.
We need to modify the degree $0$ part of $i_{{\kappa},t}$ to
$$
i_{{\kappa},t}^{\prime}([1]) = {\bf e}^+_{\kappa}.
$$
(Note that 
$
i_{\kappa,t}([1]) = {\bf e}_{\kappa}
$.)
We have $\frac{d}{dt}i'_{\kappa,t}([1]) = 0$.
So we do not need to change 
$k_{\kappa,t}$ and $h_{\kappa,t}$, that is, $k'_{\kappa,t} = k_{\kappa,t}$ and $h'_{\kappa,t} = h_{\kappa,t}$.
\par
We next  `change' $G_{L_{\kappa}}$ to $G'_{L_{\kappa}} = G_{L_{\kappa}} + (\Delta G^t_{L_{\kappa}})_*$
where
$$
\Delta G^t_{L_{\kappa}} = {\bf f}_{\kappa} \otimes {\rm vol}^t_{L_{\kappa}} + (-1)^n {\rm vol}^t_{L_{\kappa}} \otimes {\bf f}_{\kappa}.
$$
Note that $\Delta G^t_{L_{\kappa}}$ is $t$-dependent.
We have
$$
\frac{d}{dt} \Delta G^t_{L_{\kappa}} 
= d \left( 
 {\bf f}_{\kappa} \otimes h_{\kappa,t}([L_{\kappa}]) - h_{\kappa,t}([L_{\kappa}]) \otimes {\bf f}_{\kappa}
\right).
$$
In fact 
$$
\frac{d}{dt} {\rm vol}^t_{L_{\kappa}} = 
\frac{d}{dt} i_{{\kappa},t}([L_{\kappa}]) = (d\circ h_{t,\kappa})([L_{\kappa}]).
$$
We put
$$
\Delta\sigma_{\kappa,t} = 
 {\bf f}_{\kappa} \otimes h_{\kappa,t}([L_{\kappa}]) - (-1)^n h_{\kappa,t}([L_{\kappa}]) \otimes {\bf f}_{\kappa}
$$
and `put'
$$
\sigma'_{\kappa,t} = \sigma_{\kappa,t} + (\Delta\sigma{\kappa,t})_*.
$$
Those operators satisfy (\ref{form1120}), (\ref{form1121}).
\par
Actually the map $(h_{\kappa,t}([L_{\kappa}]) \otimes {\bf f}_{L_{\kappa}})_*$ is not defined 
but its composition with $\frak m^{F}$ is defined in the same way as in Subsection \ref{hunitcansubsec}.

Now we use a homotopically unital and cyclic pseudo-isotopy 
(which extends the cyclic pseudo-isotopy $(\frak C^F,\langle\, ,\,\rangle_{\rm cyc},\{\frak m^{F, t}_{k,\beta}\},\{\frak c^{F, t}_{k,\beta}\})$)
and modified versions of (\ref{form1122}) - (\ref{form1124})
in the same way as in Subsection \ref{hunitcansubsec}.
We then get a strictly unital and cyclic pseudo-isotopy between canonical models.
The proof of Theorem \ref{thm111} is complete.\footnote{
In case $L_{\kappa} = L_{\kappa'}$ but $\theta_{\kappa} \ne \theta_{\kappa'}$
we replace $d$ by $d_{\theta' - \theta} = d + (\theta' - \theta)\wedge$ and 
modify (\ref{form1120}), (\ref{form1121}) appropriately.  Then 
in the same way as Subsection \ref{localcoefficient} we can include such cases.} 
\end{proof}
\begin{rem}
In Subsection \ref{hunitcansubsec}, we prove that there exists a cyclic filtered $A_{\infty}$ functor 
from $\mathcal L$ to $\mathcal L^{\rm form}_{\rm c.u.}$ 
which is a quasi-isomorphism.  Unfortunately 
it does not so immediately imply the independence of  $\mathcal L$ 
up to cyclic homotopy equivalence.
In fact at the time of writing of this paper the authors do not know 
whether existence of homotopy inverse of a quasi-isomorphism
(a version of Whitehead theorem for $A_{\infty}$ category) holds in the cyclic category.
\end{rem}

\section{Construction of Kuranishi structures.}
\label{sec:Kuraconst}

 %\marginpar{This section is new.}
 %\marginpar{Sections 13 and 14 are rewritten. 
%The method to realize compatibility at boundary is changed from one in \cite{const2} to 
%one in \cite{linear}.
%Two sections become shorter.
%I believe it becomes easier to read.  
%Also it works for immersed Lagrangian.
% KF 2025 August}
In this and the next sections, we construct (systems of) Kuranishi 
structures and CF-perturbations thereon 
which are used in the previous sections of Part \ref{part2}. Similar arguments are 
used to construct (systems of) Kuranishi 
structures and CF-perturbations appearing in Part \ref{part3}.
The arguments are rather technical and make heavy use 
of the tool kits developed by 2nd-5th authors.
Readers who are not interested in such detail 
can skip this section  and the next.

In this section we prove Propositions \ref{diskkura}, \ref{Kuraeistspoly}, \ref{prop92} and \ref{prop94}
which claim existence of Kuranishi structures 
on various moduli spaces involved.
We employ the strategy established in \cite{const1,const2,linear}.
In fact, except Items (8), Proposition \ref{diskkura}
was proved in \cite{const1,const2}, in mostly the same way as this section.
(The method taken in this paper to achieve the compatibility at the boundary is different 
from \cite{const2} but is similar to \cite{linear}.)

\subsection{Construction of Kuranishi structures of the 
moduli spaces of disks: review.}
\label{diskreview}

In this subsection, we review the construction in 
\cite{const1,const2} and in the next subsection 
we explain how we modify it  so that 
Items (8), (9), (11) are satisfied.

We consider the situation of Section \ref{ainfalgasssingle}, that is, we have a single embedded Lagrangian submanifold $L$.
The   set
${\mathcal X}_{\ell;k+1}(L;\beta)$\index[syindex]{Xellk+1L@${\mathcal X}_{\ell;k+1}(L;\beta)$},
which we call an
{\it ambient {\bf set}}\index{ambient set},  is
defined in the same way
as ${\mathcal M}_{\ell;k+1}(L;\beta)$ is defined 
in Definition \ref{diskmoduli1}, 
except we do not require $u$ to be pseudo-holomorphic.
We do not define topology etc. on it and regard it as a set.
The inclusion ${\mathcal M}_{\ell;k+1}(L;\beta)
\subset {\mathcal X}_{\ell;k+1}(L;\beta)$ is obvious.

Suppose  a pair $(\mathcal X,\mathcal M)$
of a metrizable space $\mathcal M$ and a set $\mathcal X$ 
which contains it is given.
\begin{defn}{\rm (\cite[Definition 4.1]{const1}])}\label{defnpartialtopo}
A {\it partial topology}\index{partial topology} of such a pair, by definition, \index[syindex]{BepsilonXp@$B_{\epsilon}(\mathcal X,{\bf p})$}
assigns $B_{\epsilon}(\mathcal X,{\bf p}) \subset \mathcal X$
for each $\epsilon > 0$ and ${\bf p} \in \mathcal M$,
such that:
\begin{enumerate}
\item ${\bf p}$ is an element of $B_{\epsilon}(\mathcal X,{\bf p})$
and 
$
\{ B_{\epsilon}(\mathcal X,{\bf p}) \cap \mathcal M \mid {\bf p}, \epsilon\}
$
is a basis of the topology of $\mathcal M$.
\item 
For each $\epsilon, {\bf p}$ and ${\bf q} \in B_{\epsilon}(\mathcal X,{\bf p}) \cap \mathcal M$,
there exists $\delta >0$ such that
$
B_{\delta}(\mathcal X,{\bf q}) \subseteq B_{\epsilon}(\mathcal X,{\bf p}).
$
\item
If $\epsilon_1 < \epsilon_2$ then
$B_{\epsilon_1}(\mathcal X,{\bf p}) \subseteq B_{\epsilon_2}(\mathcal X,{\bf p})$.
Moreover
$$
\bigcap_{\epsilon} B_{\epsilon}(\mathcal X,{\bf p}) = \{{\bf p}\}.
$$
\end{enumerate}
Given a partial topology, a subset of $\mathcal X$ is said to be a 
neighborhood of ${\bf p} \in \mathcal M$ if 
it contains $B_{\epsilon}(\mathcal X,{\bf p})$ for some $\epsilon >0$.
Two partial topologies are said to be equivalent 
if the notion of neighborhood coincides.
\end{defn}

It is proved in \cite[Proposition 4.3]{const1}
that the pair $({\mathcal X}_{\ell;k+1}(L;\beta),{\mathcal M}_{\ell;k+1}(L;\beta))$
has a partial topology.  A brief outline of the proof is now in order.
From now on for an element ${\bf p} \in {\mathcal M}_{\ell;k+1}(L;\beta)$ or 
${\bf x} \in {\mathcal X}_{\ell;k+1}(L;\beta)$ we write its representative as:\index[syindex]{Sigmasubp@$\Sigma_{\bf p}$}
\index[syindex]{usubp@$u_{\bf p}$}
$$
\aligned 
{\bf p} &= (\Sigma_{\bf p};u_{\bf p};\vec z^+_{\bf p};\vec w_{\bf p})\\
{\bf x} &= (\Sigma_{\bf x};u_{\bf x};\vec z^+_{\bf x};\vec w_{\bf x}).
\endaligned
$$
Let ${\bf p} \in {\mathcal M}_{\ell;k+1}(L;\beta)$, 
${\bf x} \in {\mathcal X}_{\ell;k+1}(L;\beta)$.
We take additional interior marked points 
$\vec w^+_{\bf p}$ on 
$\Sigma_{\bf p}$ such that $(\Sigma_{\bf p};\vec z^+_{\bf p}
\cup \vec w^+_{\bf p};\vec w_{\bf p})$ (without map) is 
stable. We also require the map $u$ to be an immersion at the 
points $\vec w^+_{\bf p}$.
We say ${\bf x}$ is $\epsilon$-close to $\bf p$
(that is, ${\bf x} \in B_{\epsilon}(\mathcal X,{\bf p})$) 
if there exists $\vec w^+_{\bf x}$ on 
$\Sigma_{\bf x}$
such that:
\begin{enumerate}
\item
$(\Sigma_{\bf x};\vec z^+_{\bf x}
\cup \vec w^+_{\bf x};\vec w_{\bf x})$
is 
in the $\epsilon$-neighborhood of  $(\Sigma_{\bf p};\vec z^+_{\bf p}
\cup \vec w^+_{\bf p};\vec w_{\bf p})$
in the moduli space of marked disks.
\item
The $C^2$ distance between $u_{\bf x}$ and $u_{\bf p}$
is smaller than $\epsilon$ in the thick part.\footnote{See footnote 41.}
\item
The diameter of the image by $u_{\bf x}$ of each connected component of the thin part\footnote{See footnote 41.} of $(\Sigma_{\bf x};\vec z^+_{\bf x}
\cup \vec w^+_{\bf x};\vec w_{\bf x})$
is smaller than $\epsilon$.
\end{enumerate}
Let us elaborate on these conditions.
Since $(\Sigma_{\bf p};\vec z^+_{\bf p}
\cup \vec w^+_{\bf p};\vec w_{\bf p})$ 
and $(\Sigma_{\bf x};\vec z^+_{\bf x}
\cup \vec w^+_{\bf x};\vec w_{\bf x})$
are stable, they have thick-thin decomposition.\footnote{A reference which 
is easy to access for sympletic geometers is \cite[Chapter IV]{Hu}.}
%See also Lemma \ref{lem1221}.
Thus (1), (3) make sense.
We fix various data (such as the analytic family of coordinates at 
the nodes  (See \cite[Definitions 3.1, 3.2, 3.3]{const1}.) and local trivialization of the 
universal family of marked disks (See \cite[Definition 3.6]{const1}.) 
then we obtain a diffeomorphism 
between the thick part of 
$(\Sigma_{\bf p};\vec z^+_{\bf p}
\cup \vec w^+_{\bf p};\vec w_{\bf p})$ 
and that of 
$(\Sigma_{\bf x};\vec z^+_{\bf x}
\cup \vec w^+_{\bf x};\vec w_{\bf x})$. (See \cite[Definition 3.9]{const1}, 
where the map which gives this identification is written as $\hat{\Phi}$.)
Thus (2) makes sense.

We refer \cite{const1} for the precise definition of this 
partial topology.

Now we define the notion of 
obstruction bundle data as follows.
Let 
${\bf x} = (\Sigma_{\bf x};\vec z^+_{\bf x};\vec w_{\bf x}) \in {\mathcal X}_{\ell;k+1}(L;\beta)$.
We consider the direct sum of 
the space of $C^2$ sections of $u_{\bf x}^*TX \otimes \Lambda^{01}_{\Sigma_{\bf x}}$
of irreducible components of $\Sigma_{\bf x}$ and denote it by 
$C^2({\bf x};u_{\bf x}^*TX \otimes \Lambda^{01})$.

\begin{defn}\label{defn51}
{\it Obstruction bundle data}\index{Obstruction bundle data} of the
moduli space ${\mathcal M}_{\ell,k+1}(X,L;\beta)$
assign to each 
${\bf p} \in {\mathcal M}_{\ell,k+1}(X,L;\beta)$
a neighborhood $\mathscr U_{\bf p}$ of ${\bf p}$ in ${\mathcal X}_{k+1,\ell}(X,L;\beta)$
and an object
$E_{\bf p}(\bf x)$ to each ${\bf x} \in \mathscr U_{\bf p}$.\index[syindex]{Uscriptp@$\mathscr U_{\bf p}$} We require that they have the following properties.
\begin{enumerate}
\item
$E_{\bf p}({\bf x})$  is a finite dimensional linear subspace of\index[syindex]{Epx@$E_{\bf p}({\bf x})$} 
$C^2(\Sigma_{\bf x};u_{\bf x}^*TX \otimes \Lambda^{01})$, whose support is away from nodal or marked points 
and from the boundary. %\marginpar{and from the boundary is added. KF 2024 Dec.}
\item (Smoothness)
$E_{\bf p}({\bf x})$ depends smoothly on ${\bf x}$.
\item (Transversality)
$E_{\bf p}({\bf x})$ satisfies the transversality condition.
\item (Semi-continuity)
$E_{\bf p}({\bf x})$ is semi-continuous on  ${\bf p}$.
\item (Invariance under extended automorphisms)
$E_{\bf p}({\bf x})$ is invariant under the extended automorphism group of 
${\bf x}$.\footnote{Extended automorphism is defined in Definition \ref{diskmoduli1}.}
\item (Effectivity)   The action of the automorphism group ${\rm Aut}({\bf p})$
on the linear space $(D_{u_{\bf p}}\overline{\partial})^{-1}E_{\bf p}({\bf p})/{\rm aut}(\Sigma_{\bf p},\vec z_{\bf p}^+,\vec w_{\bf p})$
is effective.
\end{enumerate}
In case $Z \subseteq {\mathcal M}_{\ell,k+1}(X,L;\beta)$ an obstruction bundle data on $Z$ is 
by definition $\{E_{\bf p}({\bf x})\}$, such that $E_{\bf p}({\bf x})$ is defined only for ${\bf p} \in Z$, and it satisfies Conditions (1)-(6) above.
\end{defn}
Let us elaborate on those conditions.
(See \cite{const1} for detail.)
The meanings of Conditions (1), (5) are clear.
Condition (4) says that if ${\bf q} \in B_{\epsilon}(\mathcal X,{\bf p})
\cap  {\mathcal M}_{\ell,k+1}(X,L;\beta)$
and ${\bf x}$ is in a sufficiently small neighborhood of ${\bf q}$ 
then 
$E_{\bf q}({\bf x}) \subseteq E_{\bf p}({\bf x})$.
\par
Condition (3) says that the sum of  $E_{\bf p}({\bf x})$ 
and the image of the linearization operator of pseudo-holomorphic 
curve equation generates $C^2({\bf x};u_{\bf x}^*TX \otimes \Lambda^{01})$.
(There is a certain additional assumption so that
evaluation maps become weakly submersive.)

Condition (2) is a certain condition on the behavior 
of the obstruction spaces $E_{\bf p}({\bf x})$
when we fix ${\bf p}$ and move ${\bf x}$.
Since the support of $E_{\bf p}({\bf x})$ is away from 
nodal points, we can use appropriate 
identification of the thick part of $\Sigma_{\bf x}$
and of $\Sigma_{\bf p}$ and can regard
$E_{\bf p}({\bf x})$ as a family of 
subspaces of the sections of $u_{\bf p}^*TX \otimes \Lambda^{01}$.
(We use parallel transport also.)
\par
Condition (6) is assumed so that the Kuranishi neighborhood we obtain 
will be an {\it effective} orbifold.  See  \cite[Section 3]{const2}.

Thus we can formulate these conditions.
See \cite[Section 8]{const1} for the precise definition.

\begin{thm}\label{constthm}{\rm (\cite[Theorem 7.1]{const1})}
To an obstruction bundle data of the 
moduli space ${\mathcal M}_{\ell,k+1}(X,L,J;\beta)$
we can associate a Kuranishi structure on ${\mathcal M}_{\ell,k+1}(X,L,J;\beta)$.
\end{thm}
Moreover the Kuranishi structure obtained is canonical in the sense of 
germs of Kuranishi structures.\footnote{See \cite[Definition 6.6]{const1}
for the definition of the notion of a germ of Kuranishi strucutre.}
By choosing appropriate obstruction bundle data, 
the evaluation maps become weakly submersive with respect to
the Kuranishi structure obtained. (See Definition \ref{defn1223new}.)

Construction of the Kuranishi structure in Theorem \ref{constthm}
roughly goes as follows.

Let ${\bf p} \in{\mathcal M}_{\ell,k+1}(X,L,J;\beta)$.
Its Kuranishi neighborhood (an orbifold) $U_{\bf p}$
is the set of all ${\bf x} \in  \mathscr U_{\bf p}$
such that $\overline{\partial}u_{\bf x} \in E_{\bf p}({\bf x})$.
Condition (3) implies that we can apply the implicit function 
theorem to this equation. Then together with Condition (2) 
we find that $U_{\bf p}$ has a structure of orbifold.
The obstruction bundle $\mathcal E_{\bf p}$\index{obstruction bundle}
(an orbibundle on $U_{\bf p}$) 
is defined so that the fiber at ${\bf x}$ is  $E_{\bf p}({\bf x})$.
The Kuranishi map $s_{\bf p}$ (a section of $E_{\bf p}$)
is defined by $s_{\bf p}({\bf x}) = \overline{\partial}u_{\bf x}$.
Condition (2) implies that it is a smooth section.
If $s_{\bf p}({\bf x}) = 0$ then  $u_{\bf x}$ is pseudo-holomorphic.
Namely ${\bf x} \in {\mathcal M}_{\ell,k+1}(X,L,J;\beta)$.
Thus we obtain a parametrization map
$\psi_{\bf p} : s_{\bf p}^{-1}(0) \to  {\mathcal M}_{\ell,k+1}(X,L,J;\beta)$.
The quadruplet $(U_{\bf p},\mathcal E_{\bf p},s_{\bf p},\psi_{\bf p})$\index[syindex]{Excalp@$\mathcal E_{\bf p}$}
becomes a Kuranishi chart\index{Kuranishi chart} in the sense of 
\cite[Definition 3.1]{springer}.

We define coordinate change\index{coordinate change} as follows.
Let ${\bf q} \in {\mathcal M}_{\ell,k+1}(X,L,J;\beta)$
be in a small neighborhood  $\mathscr U_{\bf p}$ of
${\bf p} \in {\mathcal M}_{\ell,k+1}(X,L,J;\beta)$.
We can shrink $\mathscr U_{\bf q}$
a bit so that we can use Condition (4) 
to show $U_{\bf q} \subseteq U_{\bf p}$
as subsets of  ${\mathcal X}_{\ell,k+1}(X,L,J;\beta)$.
We can use Condition (2) to show the inclusion:
$U_{\bf q} \to U_{\bf p}$ is a smooth embedding.\footnote{See 
\cite{Corrigendum} for the definition of smooth structures of $U_{\bf q}$, $U_{\bf p}$.}
The bundle map  $\mathcal E_{\bf q} \to \mathcal E_{\bf p}$ which covers 
$U_{\bf q} \to U_{\bf p}$ is easily obtained by using Condition (4).
The compatibility with Kuranishi maps and parametrization maps 
are immediate from the definition.
We have thus defined a coordinate change in the sense of
\cite[Definition 3.2]{springer}.
\par
The cocycle condition of the coordinate change (\cite[(3.3)]{springer})
is actually  obvious in our situation since 
coordinate change map is the inclusion map as subsets of 
the ambient set ${\mathcal X}_{\ell,k+1}(X,L,J;\beta)$.
\par
This is an outline of the proof of Theorem \ref{constthm}.
See \cite{const1} for detail.
\par
To obtain a Kuranishi structure as in Propositions \ref{diskkura},
we need to show an existence of appropriate 
obstruction bundle data.
An important point is that,
to prove consistency at the boundary Propositions \ref{diskkura} (5)
and consistency with forgetful map Propositions \ref{diskkura} (8),
we need to choose our obstruction bundle data so that
it has corresponding properties.
We review the construction of obstruction bundle data
so that Propositions \ref{diskkura} (5)
(and others except (8), (9), (11)) is satisfied.
This was worked out in \cite{const2}.
\par
An additional argument needed to show Propositions \ref{diskkura} (8)
is explained in the next subsection.
\par
To show Propositions \ref{diskkura} (5)
we need our obstruction bundle data
to be disk-component-wise.
We define this notion now.

Consider the situation of (\ref{eq177}).
We consider one of the summands
\begin{equation}\label{form121}
{\mathcal M}_{\#{\mathbb L}_1;k_1+1}(L;\beta_1)
{}_{\text{\rm ev}_{0}}\times_{\text{\rm ev}_i} {\mathcal M}_{\#{\mathbb L}_2,k_2+1}(L;\beta_2)
\end{equation}
of the right hand side.
This is contained in the fiber product\footnote{This fiber product is taken 
in the category of sets.} of the 
ambient set:
\begin{equation}\label{form122}
{\mathcal X}_{\#{\mathbb L}_1,k_1+1}(L;\beta_1)
{}_{\text{\rm ev}_{0}}\times_{\text{\rm ev}_i} {\mathcal X}_{\#{\mathbb L}_2,k_2+1}(L;\beta_2).
\end{equation}
Let $({\bf p}_1,{\bf p}_2)$ be in (\ref{form121}).
We take a neighborhood  $\mathscr U_{{\bf p}_i}$ of ${\bf p}_i$ 
in its ambient set.
Then it is easy to see that the fiber product 
$\mathscr U_{{\bf p}_1} {}_{\text{\rm ev}_{0}}\times_{\text{\rm ev}_i} \mathscr U_{{\bf p}_2} $
is a neighborhood of $({\bf p}_1,{\bf p}_2)$ in (\ref{form122}).
(See \cite[Lemma 4.10]{const2}.)
\begin{defn}\label{diskcomponentwiseness}{\rm (\cite[Definition 5.1]{const2})}
A system of obstruction bundle data $\{E_{\bf p}({\bf x})\}$
for ${\mathcal M}_{k+1}(L;\beta)$ is said to be 
{\it disk-component-wise}\index{disk-component-wise} if
$$
E_{({\bf p}_1,{\bf p}_2)}({\bf x}_1,{\bf x}_2)
=
E_{{\bf p}_1}({\bf x}_1) \oplus E_{{\bf p}_2}({\bf x}_2)
$$
holds in the above situation and ${\bf x}_i \in \mathscr U_{{\bf p}_i}$.
\end{defn}
\begin{rem}
In \cite[Definition 5.1]{const2}, disk-component-wise-ness is 
defined in a slightly different way including the case of multiple fiber product.
However the definition there is equivalent to Definition \ref{diskcomponentwiseness}.
\end{rem}
Now the next lemma is intuitively clear.
\begin{lem}\label{obsgtoKrua}{\rm (\cite[Theorem 5.3]{const2})}
The system of Kuranishi structures on ${\mathcal M}_{k+1}(L;\beta)$
induced from a disk-component-wise system of obstruction bundle data 
via Theorem $\ref{constthm}$, satisfies Proposition $\ref{diskkura}$ $(5)$.
\end{lem}
\begin{proof}
The proof is mostly obvious. 
In fact, if 
${\bf x} = ({\bf x}_1,{\bf x}_2) \in \mathscr U_{{\bf p}_1} {}_{\text{\rm ev}_{0}}\times_{\text{\rm ev}_i} \mathscr U_{{\bf p}_2} $
then 
$$
\overline{\partial}u_{\bf x} \in E_{({\bf p}_1,{\bf p}_2)}({\bf x}_1,{\bf x}_2)
$$
if and only if
$$
\overline{\partial}u_{{\bf x}_1} \in E_{{\bf p}_1}({\bf x}_1) \quad \text{and} \quad
\overline{\partial}u_{{\bf x}_2} \in E_{{\bf p}_2}({\bf x}_2).
$$
\end{proof}
Now to complete the proof of Proposition \ref{diskkura} except 
Item (8), it suffices to find a disk-component-wise system of obstruction bundle data.
 (Items (9), (10) follow from corresponding symmetry of the obstruction bundle data.)
The construction of such a system was carried out in \cite[Sections 7,8]{const2}.
In \cite[Sections 7,8]{const2} we used the notion of quasi-component 
of an object of ${\mathcal M}_{\ell;k+1}(L;\beta)$.
Here we provide a different method, which is similar to one 
used in \cite{linear} (and \cite[Remark 4.3.89]{toric3}).\footnote{One reason we provide a different method is 
that it works better when we make it compatible with the forgetful map.}
In this method we construct a system of 
Kuranishi structures not on ${\mathcal M}_{\ell,k+1}(L;\beta)$ itself 
but on its outer collaring  ${\mathcal M}_{\ell,k+1}(L;\beta)^{\boxplus 1}$, 
introduced in \cite{springer}. (See also \cite[Subsection 8.6]{Farber}.)\index{outer collaring}
We first review the definition of the space  ${\mathcal M}_{\ell,k+1}(L;\beta)^{\boxplus 1}$.

\begin{defn}\label{outcolar}
An element of $\mathcal M_{\ell,k+1}(X,L;\beta)^{\boxplus 1}$ 
is a pair $({\bf p},\vec s)$ of an element $
{\bf p} \in \mathcal M_{\ell,k+1}(X,L;\beta)$
and $\vec s = (s_\frak z)_{\frak z}$ where $\frak z$ is a 
boundary node of ${\bf p}$ and $s_{\frak z} \in [-1,0]$. 
\par
We define a map 
$\mathcal M_{\ell,k+1}(X,L;\beta) \to \mathcal M_{\ell,k+1}(X,L;\beta)^{\boxplus 1}$
by sending ${\bf p}$ to $({\bf p},(0,\dots,0))$.
\par
We define the ambient set $\mathcal X_{\ell,k+1}(X,L;\beta)^{\boxplus 1}$
in the same way.
\end{defn}
We can define an analogue of stable map topology on $\mathcal M_{\ell,k+1}(X,L;\beta)^{\boxplus 1}$
in an obvious way
so that the map $\mathcal M_{\ell,k+1}(X,L;\beta) \to \mathcal M_{\ell,k+1}(X,L;\beta)^{\boxplus 1}$ 
becomes a topological embedding.  (See \cite[Definition 5.1]{linear}.)
We define a partial topology of the pair
$(\mathcal X_{\ell,k+1}(X,L;\beta)^{\boxplus 1}, \mathcal M_{\ell,k+1}(X,L;\beta)^{\boxplus 1})$
as follows.\index[syindex]{oplusupper@${}^{\boxplus 1}$}\index[syindex]{Mellk+1Lbetaoplus@$\mathcal M_{\ell,k+1}(X,L;\beta)^{\boxplus 1}$}
\index[syindex]{Xellk+1Lbetaoplus@$\mathcal X_{\ell,k+1}(X,L;\beta)^{\boxplus 1}$}
\begin{defn}\label{def137}
Let $\widehat{\bf p} = ({\bf p},\vec s) \in  \mathcal M_{\ell,k+1}(X,L;\beta)^{\boxplus 1}$ 
and  $\widehat{\bf x} = ({\bf x},\vec s') \in \mathcal X_{\ell,k+1}(X,L;\beta)^{\boxplus 1}$.
We put $\min s =-\max\{s_{\frak z} \mid s_{\frak z} < 0\}$.
We define $B_{\epsilon}(\widehat{\bf p})$ when $\epsilon < \min s$.
We cut ${\bf p}$ at the boundary nodes $\frak z$ for which $s_{\frak z} < 0$.
Let $\{{\bf p}_a\}$ ($a \in \mathcal A$) be the set of all the maximal connected unions of extended 
disk components of ${\bf p}$\index{extended disk components} obtained by this cut.
Here an extended disk component is by definition an irreducible disk component together with 
tree of spheres attached to it.
We define  $\widehat{\bf x} \in B_{\epsilon}(\widehat{\bf p})$
if the following holds.
\begin{enumerate}
\item
$\widehat{\bf x}$ is decomposed into ${\bf x}_a$ with 
$
{\bf x}_a \in B_{\epsilon}(\widehat{\bf p}_a)
$. 
The way ${\bf x}_a$ are glued to give ${\bf x}$ is the same as the way 
${\bf p}_a$ are glued to give ${\bf p}$. 
\item
Let $\frak z$ be a boundary node of ${\bf p}$ such that $s_{\frak z} < 0$.  Then, 
by Item (1), there exists a corresponding boundary node of ${\bf x}$, which we denote 
by $\frak z$. We require:
$$
\vert s_{\frak z} - s'_{\frak z} \vert < \epsilon.
$$
\item
Let $\frak z$ be a boundary node of ${\bf p}$ such that $s_{\frak z} = 0$.
Then there may or may not exist a corresponding boundary node of ${\bf x}$.
When it exists we write it $\frak z$. We then require
$
s'_{\frak z}  > -\epsilon.
$
\end{enumerate}
\end{defn}

It is easy to see that Definition \ref{def137} defines a partial topology.
Now we  modify Definition \ref{defn51} as follows
\begin{defn}\label{defn51col}
{\it Obstruction bundle data}\index{Obstruction bundle data} of the
moduli space ${\mathcal M}_{\ell,k+1}(X,L;\beta)^{\boxplus 1}$
assign to each 
$\widehat{\bf p} = ({\bf p},\vec s) \in {\mathcal M}_{\ell,k+1}(X,L;\beta)^{\boxplus 1}$
a neighborhood $\mathscr U_{\hat{\bf p}}$ of $\hat{\bf p}$ in ${\mathcal X}_{\ell,k+1}(X,L;\beta)^{\boxplus 1}$
and an object
$E_{\hat{\bf p}}(\hat{\bf x})$ to each $\widehat{\bf x} = ({\bf x},\vec s') \in \mathscr U_{\hat{\bf p}}$ . We require that they have the following properties.
\begin{enumerate}
\item
$E_{\hat{\bf p}}(\hat{\bf x})$  is a finite dimensional linear subspace of 
$C^2({\bf x};u_{\bf x}^*TX \otimes \Lambda^{01})$, whose support is away from nodal or marked points 
and from the boundary.
\item (Smoothness)
$E_{\hat{\bf p}}(\hat{\bf x})$ depends smoothly on $\hat{\bf x}$.
\item (Transversality)
$E_{\hat{\bf p}}(\hat{\bf x})$ satisfies the transversality condition.
\item (Semi-continuity)
$E_{\hat{\bf p}}(\hat{\bf x})$ is semi-continuous on  $\hat{\bf p}$.
\item (Invariance under extended automorphisms)
$E_{\hat{\bf p}}(\hat{\bf x})$ is invariant under the extended automorphism group of 
${\bf x}$.
\item (Effectivity)   The action of the automorphism group ${\rm Aut}({\bf p})$
on the linear space $(D_{u_{\bf p}}\overline{\partial})^{-1}E_{\hat{\bf p}}(\hat{\bf p})/{\rm aut}(\Sigma_{\bf p},\vec z_{\bf p}^+,\vec w_{\bf p})$
is effective.
\end{enumerate}
\end{defn}

The next proposition can be proved in the exactly the same way as Theorem \ref{constthm}.

\begin{prop}\label{constthmout}
To an obstruction bundle data of the 
moduli space ${\mathcal M}_{\ell,k+1}(X,L,J;\beta)^{\boxplus 1}$
we can associate a Kuranishi structure on ${\mathcal M}_{\ell,k+1}(X,L,J;\beta)^{\boxplus 1}$.
\end{prop}

\begin{figure}[h]
\centering
\includegraphics[scale=0.3]{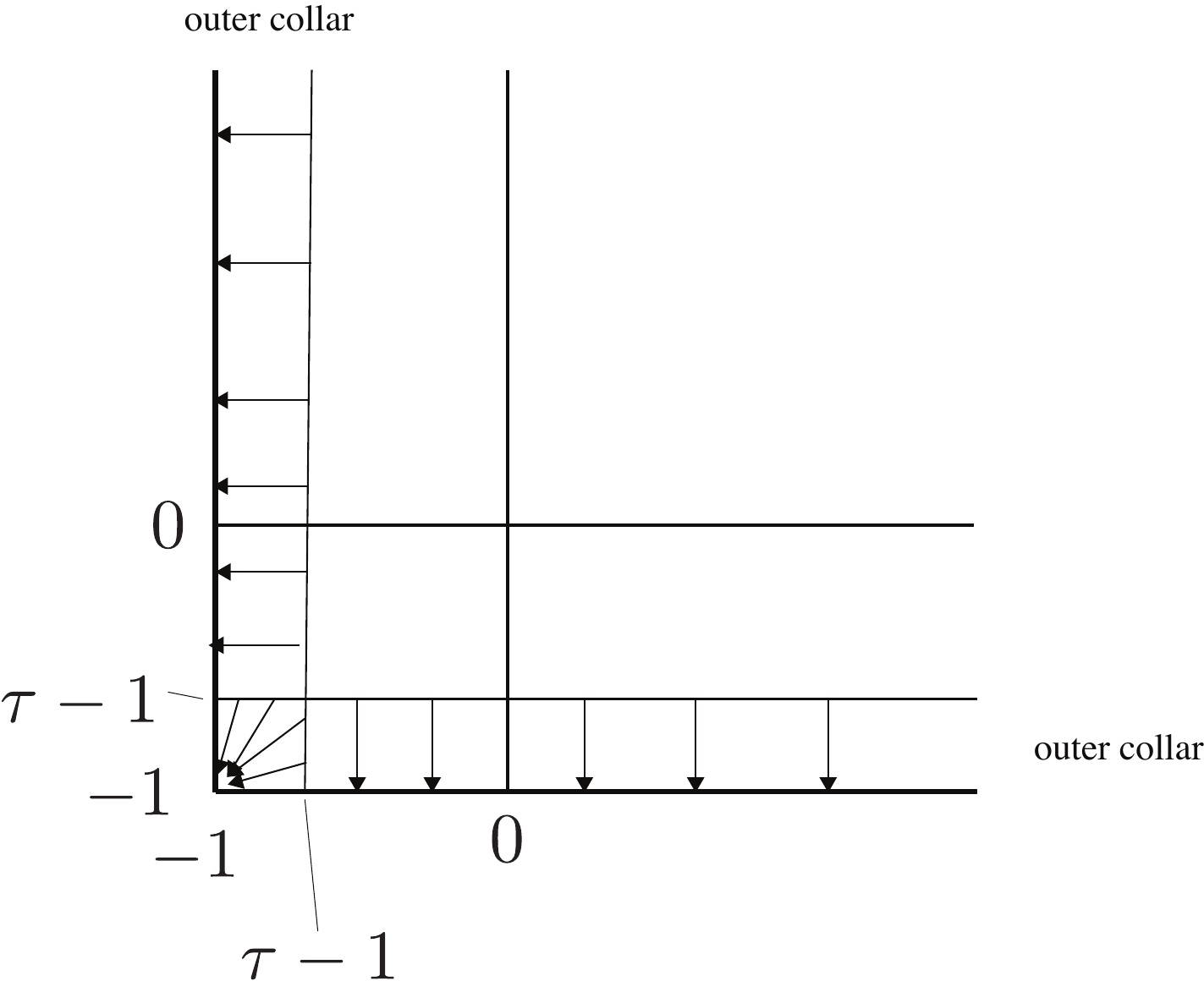}
\caption{Projection $\Pi_{\tau}$}
\label{zu10}
\end{figure}

We next define $\tau$-collared-ness of obstruction bundle data.
For $\tau > 0$ we define a map\index[syindex]{Pitau@$\Pi_{\tau}$} 
$$
\Pi_{\tau}: 
\mathcal M_{k+1,\ell}(X,L;\beta)^{\boxplus 1} \to \mathcal M_{k+1,\ell}(X,L;\beta)^{\boxplus 1}
$$ 
as follows. Let $({\bf p},\vec s) \in \mathcal M_{k+1,\ell}(X,L;\beta)^{\boxplus 1}$ .
We define $\vec s\,^{\prime}$ by
\begin{enumerate}
\item $s'_{\frak z} = s_{\frak z}$ if $s_{\frak z} > \tau-1$.
\item
$s'_{\frak z} = -1$ if $s_{\frak z} \le \tau -1 $.
\end{enumerate}  
This map is discontinuous.  See Figure \ref{zu10}.     %\marginpar{Figure added} 
 In the figure, we partition the 
square $[-1,0)^2$ into the union
\begin{eqnarray*}
[-1,0)^2 & = & [-1,-1 + \tau)^2 \\
& {} & \quad \cup \left(([-1+\tau, 0) \times [-1,-1 +\tau)) \cup ([-1,-1+\tau) \times [-1,\tau,0))\right)\\
& {} & \quad \cup (-1 + \tau, 0)^2.
\end{eqnarray*}
Note that this decomposition is consistent with the dimension stratification of $[0,1)^2$ in that each
component is a semi-open tubular neighborhood of the associated stratum chosen
inductively starting from the zero dimensional strata. In this regard,
the map $\Pi_\tau$ is the union of the deformation retractions of the relevant tubular neighborhood to
the associated stratum.
We define 

$$
\Pi_{\tau}: 
\mathcal X_{\ell,k+1}(X,L;\beta)^{\boxplus 1} \to \mathcal X_{\ell,k+1}(X,L;\beta)^{\boxplus 1}
$$
in the same way.
\par
It is simpler defining $\tau$-collared-ness of obstruction bundle data than 
defining $\tau$-collared-ness of 
Kuranishi structure.\footnote{See \cite[Section 17.5]{springer} for the latter.}

Let $\partial \mathcal M_{\ell,k+1}(X,L;\beta)^{\boxplus 1}$
be the boundary, which is the set of $({\bf p},\vec s_{\bf p})$ such that 
at least one of $s_{{\bf p},i}$ is $-1$.
We say $\Pi_{\tau}^{-1}(\partial \mathcal M_{\ell,k+1}(X,L;\beta)^{\boxplus 1})$
the $\tau$-collar of the boundary.
In other words it is the set of $({\bf p},\vec s_{\bf p})$ such that 
$s_{{\bf p},i} \le \tau-1$ holds for at least one $i$.

\begin{defn}
Obstruction bundle data $E_{\hat{\bf p}}(\widehat{\bf x})$ is said to be $\tau$-{\it collared}\index{collared obstruction bundle data} 
if 
\begin{equation}\label{formnewnew1233}
E_{\Pi_{\tau}(\hat{\bf p})}(\Pi_{\tau}(\hat{\bf x}))  =   E_{\hat{\bf p}}(\widehat{\bf x}).
\end{equation}
\end{defn}

A Kuranishi structure is said to be $\tau$-{\it collared}\index{collared Kuranishi structure} if it is obtained from 
 $\tau$-collared obstruction bundle data.

When we replace ${\mathcal M}_{\ell;k+1}(L;\beta)$  by
${\mathcal M}_{\ell;k+1}(L;\beta)^{\boxplus 1}$  
its normalized boundary $\partial{\mathcal M}_{\ell;k+1}(L;\beta)^{\boxplus 1}$  
is 
\begin{equation}\label{form135rev}
\bigcup {\mathcal M}_{\# {\mathbb L}_1;k_1+1}(L;\beta_1)^{\boxplus 1}
{}_{\text{\rm ev}_{0}}\times_{\text{\rm ev}_i} {\mathcal M}_{\# {\mathbb L}_2;k_2+1}(L;\beta_2)^{\boxplus 1}.
\end{equation}

We define disk-component-wise-ness of a system of obstruction bundle data on $\mathcal M_{k+1,\ell}(X,L;\beta)^{\boxplus 1}$
in the same way as Definition \ref{diskcomponentwiseness}.

Lemma \ref{obsgtoKrua} still holds for the outer collared version.
Now we prove:

\begin{lem}\label{obsgtoKruaout}
There exists a  disk-component-wise system of obstruction bundle data  ${\mathcal M}_{\ell,k+1}(L;\beta)^{\boxplus 1}$.
Those obstruction bundle data are chosen to be collared.
\end{lem}

\begin{proof}
The proof is by induction on the partial order $<$ on $\{(\omega(\beta),k,\ell)\}$ as in Definition \ref{defn837837}.
We remark that
$(\omega(\beta_i),k_i,\# {\mathbb L}_i) < (\omega(\beta),k,\# {\mathbb L})$ in the situation of (\ref{form135rev}).
\par
Therefore the first step of the induction follows  from Lemma \ref{lem12222} below.
\par
We next discuss an inductive step.
Suppose we have required obstruction bundle data $\{E_{\hat{\bf p}}(\hat{\bf x})\}$ on ${\mathcal M}_{\ell',k'+1}(L;\beta')^{\boxplus 1}$
with $(\omega(\beta'),k',\ell') < (\omega(\beta),k,\ell)$.
Suppose that those are $\tau'$-collared for $\tau' > 0$.
We will construct one on ${\mathcal M}_{\ell,k+1}(L;\beta)^{\boxplus 1}$, which is $\tau$-collared 
for a certain $0< \tau < \tau'$.
\par
Suppose $\hat{\bf p} \in \partial{\mathcal M}_{\ell,k+1}(L;\beta)^{\boxplus 1}$.
Then $\hat{\bf p} = [\hat{\bf p}_1,\hat{\bf p}_2]$ is the equivalence class in (\ref{form135rev})
represented by the tuple $({\bf p}_1,\hat{\bf p}_2)$. 
We can take $U_{\hat{\bf p}}$ so that it is contained in 
$
U_{\hat{\bf p}_1} {}_{\text{\rm ev}_{0}}\times_{\text{\rm ev}_i} U_{\hat{\bf p}_2}.
$
Then for $\hat{\bf x} = (\hat{\bf x}_1,\hat{\bf x}_2) \in U_{\hat{\bf p}}$ we define:
\begin{equation}\label{newform124}
E_{\hat{\bf p}}(\hat{\bf x}) = E_{\hat{\bf p}_1}(\hat{\bf x}_1) \oplus E_{\hat{\bf p}_2}(\hat{\bf x}_2).
\end{equation}
We remark that the union (\ref{form135rev}) is not the disjoint union.
(Here we consider the union as subsets of the ambient set.  If we consider the normalized boundary 
as a space with Kuranishi structure then the union is the disjoint union.)
Using the induction hypothesis, that is,  $\{E_{\hat{\bf p}}(\hat{\bf x})\}$ on ${\mathcal M}_{\ell',k'+1}(L;\beta')^{\boxplus 1}$
is disk-component-wise, it is easy to see that the right hand side of (\ref{newform124}) is independent of the representative
$(\hat{\bf p}_1,\hat{\bf p}_2)$ depending only on
the element $\hat{\bf p} = [\hat{\bf p}_1,\hat{\bf p}_2]$ itself.
\par
We have thus defined an obstruction bundle data on the boundary 
$\partial{\mathcal M}_{k+1}(L;\beta)^{\boxplus 1}$.
We can extend it to a neighborhood of $\partial{\mathcal M}_{k+1}(L;\beta)^{\boxplus 1}$ by requiring (\ref{formnewnew1233}).
Now the inductive step follows from the next lemma.

\begin{lem}\label{lem12222}
Let $K$ be a compact subset of ${\mathcal M}_{\ell,k+1}(L;\beta)^{\boxplus 1}$ and $U$  its neighborhood.
Let $Z \supset U$ be a compact set.
Suppose we are given an obstruction bundle data $\{E_{\bf}({\bf x})\}$ on $U$.\footnote{Namely
$E_{\bf p}({\bf x})$ is defined for ${\bf p} \in U$.} Then there exists an obstruction bundle data $\{E'_{\bf}({\bf x})\}$ on $
Z \subseteq {\mathcal M}_{k+1}(L;\beta)^{\boxplus 1}$ 
such that $E'_{\bf p}({\bf x}) = E_{\bf p}({\bf x})$ for ${\bf p} \in K$.
\end{lem}
\begin{proof}
If we replace ${\mathcal M}_{k+1}(L;\beta)^{\boxplus 1}$ by ${\mathcal M}_{k+1}(L;\beta)$ and take $K = \emptyset$ then 
Lemma \ref{lem12222} is 
\cite[Theorem 11.1]{const1}. We can modify the proof easily to prove Lemma \ref{lem12222}.
For completeness sake we briefly review the proof.
\par
For ${\bf p} \in {\mathcal M}_{k+1}(L;\beta)^{\boxplus 1}$ we obtain 
$E^0_{{\bf p}}({\bf x})$ (for ${\bf x}$ in a sufficiently small neighborhood of ${\bf p}$ in the ambient set)
such that Definition \ref{defn51} (1)(2)(4)(6) is satisfied.
Then there exists a neighborhood $V({\bf p})$  of ${\bf p}$ in ${\mathcal M}_{k+1}(L;\beta)^{\boxplus 1}$
such that if ${\bf q} \in V({\bf p})$ then the choice $E^0_{{\bf q}}({\bf x}): = E^0_{{\bf p}}({\bf x})$ satisfies Definition \ref{defn51} (1)(2)(4)(6)
for any ${\bf q} \in V({\bf p})$.
\par
We may assume that $V({\bf p})$ is a closed subset of ${\mathcal M}_{k+1}(L;\beta)^{\boxplus 1}$.
We take  a closed subset $U_0$ of $U$ which is a neighborhood of $K$.
We may take a finite subset $\{{\bf p}_a\}$ of ${\mathcal M}_{k+1}(L;\beta)^{\boxplus 1}$
such that
\begin{enumerate}
\item $U \cup \bigcup_a{\rm Int}V({\bf p}_a) \supset Z$.
\item $V({\bf p}_a) \cap U_0 = \emptyset$.
\end{enumerate}
(We may shrink $U$ and $V({\bf p}_a)$ if necessary.)
We put
$$
E'_{\bf{p}}({\bf x}) =
\begin{cases} 
\displaystyle \bigoplus_{a; {\bf p} \in V({\bf p}_a)}   E^0_{{\bf p}_a}({\bf x})    &\text{if ${\bf p} \notin U_0$},
\\
\displaystyle E_{\bf p}({\bf x}) \oplus \bigoplus_{a; {\bf p} \in V({\bf p}_a)}   E^0_{{\bf p}_a}({\bf x}) &\text{if ${\bf p} \in U_0$}.
\end{cases}  
$$
By perturbing $E^0_{{\bf p}_a}({\bf x})$ a bit we may assume that the sum appearing in the right hand side is the 
direct sum. (See \cite[Lemma 11.7]{const1}.)
\par
This satisfies Definition \ref{defn51} (1)(2)(4)(6) since $E_*(*)$ and $E^0$ have the corresponding properties.
Definition \ref{defn51} (3) then follows from the same property of $E_*(*)$ and the closed-ness of 
$V_{{\bf p}_a}$, $U_0$.
\end{proof}
The proof of Lemma \ref{obsgtoKruaout} is complete.\end{proof}
We have thus proved the outer-collared-version of Propositions \ref{diskkura}, except Item (8). 
\par
We can use  this outer-collared-version  in the same way for our application to the 
construction of $A_{\infty}$ structures.

\subsection{Construction of Kuranishi structure of the 
moduli space of disks: compatibility with forgetful map.}
\label{diskforget}

We now discuss the way to obtain a system of Kuranishi structures  
that is compatible with the forgetful maps.
The strategy is as follows.
We formulate the notion that a system of obstruction bundle data 
is compatible with the forgetful maps. Then we show that 
if a system of obstruction bundle data is compatible with forgetful maps  then the resulting system of 
Kuranishi structures  is compatible with forgetful maps.
\par
We begin with the first step.
We consider ${\mathcal M}_{\ell;k+1}(L;\beta)$
and a forgetful map\index{forgetful map}
$$
\frak{forget}: {\mathcal M}_{\ell;k+1}(L;\beta) \to {\mathcal M}_{\ell,0}(L;\beta)
$$
in the case $\beta \ne 0$.  (In the case when $\beta = 0$ the 
obstruction bundle is defined to be $0$.)
\par
Here we review the definition of the forgetful map $\frak{forget}$.\index[syindex]{forgetfrak@$\frak{forget}$}
 %\marginpar{Definition of forgetful map is reviewed here. KF 2025 Aug.}

Let 
${\bf p} = (\Sigma;u;z_1^{+},\dots,z_{\ell}^+;w_0,\dots,w_k)$
be an element of ${\mathcal M}_{\ell;k+1}(L;\beta)$.
We consider 
$(\Sigma;u;z_1^{+},\dots,z_{\ell}^+;w_0,\dots,w_{k-1})$.
If this object is stable we define
$$
\frak{forget}({\bf p}): = (\Sigma;u;z_1^{+},\dots,z_{\ell}^+;w_0,\dots,w_{k-1}).
$$
If $(\Sigma;u;z_1^{+},\dots,z_{\ell}^+;w_0,\dots,w_{k-1})$ is not stable, 
there is an irreducible component $\Sigma_a$ of $\Sigma$ 
such that there exists a set of automorphisms of positive dimension 
supported on $\Sigma_a$.  We shrink $\Sigma_a$ to a point and
obtain $\Sigma'$.  If one of the marked points 
$z_1^{+},\dots,z_{\ell}^+;w_0,\dots,w_{k-1}$ is on $\Sigma_a$ we take a corresponding marked 
point of $\Sigma'$ (See \cite[page 419]{fooo092}.). 
The map $u$ is constant on $\Sigma_a$. Therefore we obtain a map $u'$ on $\Sigma'$.  We thus obtain 
$(\Sigma';u';z^{\prime +}_1,\dots,z^{\prime +}_{\ell};w'_0,\dots,w'_{k-1})$.
If it is stable we put
$$
\frak{forget}({\bf p}): = (\Sigma';u';z^{\prime +}_1,\dots,z^{\prime +}_{\ell};w'_0,\dots,w'_{k-1}).
$$
Otherwise we repeat the process.
\par
We can also define
$$
\frak{forget}: {\mathcal X}_{\ell;k+1}(L;\beta) \to {\mathcal X}_{\ell,0}(L;\beta)
$$
in a similar way.
An issue is that the forgetful map is {\it not} continuous with respect to the
partial topology.
Namely for $\overline{\bf p} = \frak{forget}(\bf p)$ the following does {\it not}
hold.
\begin{enumerate}
\item[(*)]
For each neighborhood $\mathscr U_{\overline{\bf p}}$ of  $\overline{\bf p}$ in the ambient set
there exists a neighborhood  $\mathscr U_{{\bf p}}$ of ${\bf p}$ in the ambient set 
such that
$$
\frak{forget}(\mathscr U_{{\bf p}}) \subset \mathscr U_{\overline{\bf p}}.
$$
\end{enumerate}
\begin{exm}\label{exm1314}
A counter example is the following. 
Let ${\bf p} = (\Sigma_{\bf p},\vec z_{\bf p},u_{\bf p}) \in {\mathcal M}_{0,3}(L;\beta)$
consist of two disk components $\Sigma_1,\Sigma_2$.
$z_{{\bf p},0} \in \Sigma_1$ and $z_{{\bf p},1},z_{{\bf p},2} \in \partial\Sigma_2$.
Suppose $u_{\bf p}$ is non-constant on $\Sigma_1$ and is constant on 
$\Sigma_2$.  When we forget the marked points on the component 
$\Sigma_2$ and shrink the components. We thus obtain a map from $\Sigma_1$  only.
See Figure \ref{FigureSigma12}.
\begin{figure}[h] 
\centering
\includegraphics[scale=0.5]{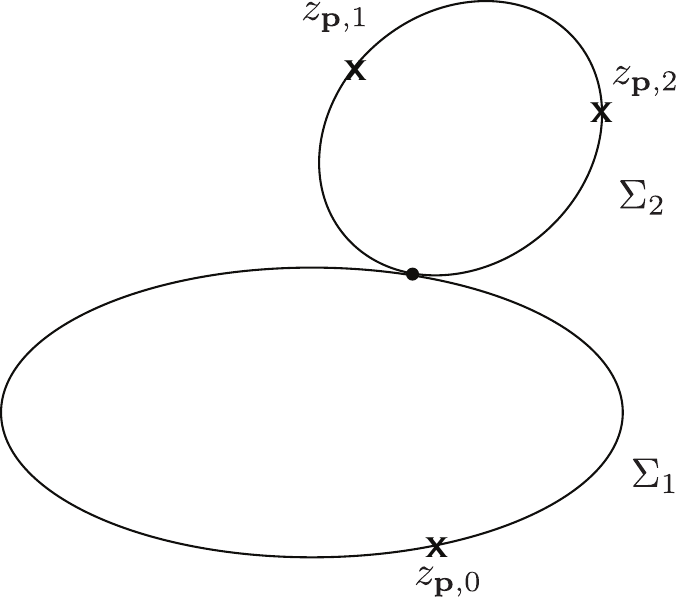}
\caption{Forgetful map collapses a disk bubble}
\label{FigureSigma12}
\end{figure}
\par
We consider ${\bf x}$ in an arbitrary small neighborhood of ${\bf p}$
such that the domain is the same and the map on $\Sigma_2$
is close to a constant map but is not a constant map.
Then the source curve $\frak{forget}({\bf x})$ still has two components.
(We shrink only an unstable component to define  $\frak{forget}$. 
The automorphism of a disk with non-constant map is finite  
even though the map is close to a constant map.)
From the definition of partial topology such an object  is not in a neighborhood of $\overline{\bf p}$.
\end{exm}
\par
The idea to resolve this issue is as follows.
We take an obstruction bundle $E_{{\bf p}}({\bf x})$ so that it is $0$ on the 
component where the map  is constant.
Then if ${\bf x}$ satisfies the equation 
$\overline{\partial}u_{{\bf x}} \in E_{{\bf p}}({\bf x})$
this phenomenon does not happen, if we take $E_{{\bf p}}({\bf x})$
in an appropriate way.
\par
We formulate a few conditions so that the above argument works.
Let ${\bf p} \in {\mathcal M}_{\ell;k+1}(L;\beta)$
and ${\bf x}$ be in a small neighborhood $\mathscr U_{{\bf p}}$ 
in the ambient set.
For each such pair $({\bf p}, {\bf x})$,
we have an open embedding
$$
\hat{\Phi}_{{\bf p},{\bf x}} : \Sigma_{\bf x}({\rm thick}) \to  \Sigma_{\bf p}.
$$
Here $ \Sigma_{\bf x}({\rm thick})$ is a thick part.\index{thick part}  (See footnote 41.)
This map is constructed for example in \cite[Lemma 3.9]{const1} (and is denoted by $\widehat{\Phi}_{\frak x;\vec \epsilon}$ there).
The construction goes roughly as follows.
In case $(\Sigma_{\bf p}, \vec z_{\bf p})$  is stable, 
then we can use a local trivialization of the universal family of the Deligne-Mumford moduli space 
to obtain $\hat{\Phi}_{{\bf p},{\bf x}}$.  (More precisely we use a $C^{\infty}$ class   local trivialization of the universal family over the
 Deligne-Mumford moduli space  away from nodal points of the singular fiber.)
 In the unstable case, we add marked points $
 \vec z_{\bf p}^{\,\rm add}$ to the unstable components of $\Sigma_{\bf p}$ where $u_{\bf p}$ is an
 immersion. We then take codimension 2 submanifolds of $X$ which intersect transversally to $u_{\bf p}$ at each point of $\vec z_{\bf p}^{\,\rm add}$.
 Then for each nearby ${\bf x}$ we obtain additional marked points $\vec z_{\bf x}^{\,\rm add}$ on $\Sigma_{\bf x}$ such that $u_{\bf x}(\vec z_{\bf x}^{\,\rm add})$ are on 
 those codimension 2 submanifolds.  Then since $[\Sigma_{\bf x}, \vec z_{\bf x}\cup  \vec z_{\bf x}^{\,\rm add}]$ 
 is close to $[\Sigma_{\bf p}, \vec z_{\bf p}\cup  \vec z_{\bf p}^{\,\rm add}]$  we can use  a local trivialization of the universal family of the 
 Deligne-Mumford moduli space 
to obtain $\hat{\Phi}_{{\bf p},{\bf x}}$.

Let $\{ E_{\bf p}(\bf x)\}$ be an obstruction bundle data.
We obtain a complex linear embedding  
$$
I_{{\bf p},{\bf x}} : E_{\bf p}({\bf x}) \to \Gamma(\Sigma_{\bf p},u_{\bf p}^*TX\otimes \Lambda^{01})
$$
by the parallel transport of the pull-back $\hat{\Phi}_{{\bf p},{\bf x}}^*(E_{\bf p}({\bf x}))$ along the short geodesics
from the target map of ${\bf p}$ to that of ${\bf x}$.
We compose it with the projection 
$$
 \Gamma(\Sigma_{\bf p},u_{\bf p}^*TX\otimes \Lambda^{01}) \to E_{\bf p}({\bf p})
$$
($L^2$ projection on each irreducible components) to obtain 
$$
I'_{{\bf p},{\bf x}} : E_{\bf p}({\bf x}) \to E_{\bf p}({\bf p}).
$$

\begin{defn}\label{defn1313}
We say an obstruction bundle data $\{ E_{\bf p}(\bf x)\}$
to enjoy the {\it support preserving property}\index{support preserving property} if for any $\epsilon>0$
we can choose $\mathscr U_{{\bf p}}$  so that 
for any $s \in  E_{\bf p}(\bf x)$ and ${\bf x} \in \mathscr U_{{\bf p}}$
we have
$$
B_{\epsilon} ({\rm Supp}(I'_{{\bf p},{\bf x}}(s)))
\supset \hat{\Phi}_{{\bf p},{\bf x}}({\rm Supp}(s)).
$$
Here $B_{\epsilon}$ denotes the $\epsilon$ neighborhood and 
${\rm Supp}$ denotes the support.
\end{defn}
It is obvious from definition that the obstruction bundle data we obtained in the 
last subsection enjoy the support preserving property.

\begin{defn}\label{defn1314}
We say that an  obstruction bundle data $\{ E_{\bf p}({\bf x})\}$ on ${\mathcal M}_{\ell;k+1}(L;\beta)$
{\it preserves triviality}\index{preserves triviality} if the following holds. If ${\bf p} \in {\mathcal M}_{\ell;k+1}(L;\beta)$
and if the map $u_{\bf p}$ on an extended disk component $\Sigma_v$ of $\Sigma_{\bf p}$
is constant,  all the sections in $E_{\bf p}({\bf p})$ are zero on $\Sigma_v$.
\end{defn}

Note that $\hat{\Phi}_{{\bf p},{\bf x}}$ depends on various choices such as the local trivialization
of the universal family over the Deligne-Mumford space 
or the additional marked points to stabilize the domain.  However two different choices give maps which are $C^0$ closed to each other.
Therefore Definitions \ref{defn1313}, \ref{defn1314} are independent of such choices.

Now the following lemma is  obvious.

\begin{lem}\label{lem1218}
Let $\{ E_{\bf p}(\bf x)\}$ be an obstruction bundle data on ${\mathcal M}_{\ell;k+1}(L;\beta)$
which satisfies the support preserving property and preserves triviality.
Then for each neighborhood $\mathscr U_{\overline{\bf p}}$ of  $\overline{\bf p}$ in the ambient set
there exists a neighborhood  $\mathscr U_{{\bf p}}$ of ${\bf p}$ in the ambient set 
such that
if ${\bf x} \in \mathscr U_{{\bf p}}$ and $\overline{\partial}u_{\bf x} \in E_{\bf p}({\bf x})$
then 
$$
\frak{forget}({\bf x}) \in \mathscr U_{\overline{\bf p}}.
$$
\end{lem}
\begin{exm}
Let us consider the situation of Example \ref{exm1314}, which is a counter example to (*).
Suppose ${\bf x}$ is in an arbitrary small neighborhood of $u_{\bf p}$
such that the domain is the same, namely it has two irreducible 
components $\Sigma_1$, $\Sigma_2$.
Let $\{ E_{\bf p}(\bf x)\}$ be an obstruction data 
which satisfies the support preserving property and preserves triviality.
Then it follows from the definition that the support of an element of $E_{\bf p}(\bf x)$
does not intersect with $\Sigma_2$.
Suppose $\overline{\partial}u_{\bf x} \in E_{\bf p}(\bf x)$.
Then $u_{\bf x}$ is constant on $\Sigma_2$.
Therefore, after forgetting the marked points the domain of 
$\overline{\bf x} = \frak{forget}({\bf x})$ is $\Sigma_1$.
(Namely we shrink $\Sigma_2$.)
Therefore $\overline{\bf x}$ is in a small neighborhood of 
$\overline{\bf p}$.
\end{exm}

We next formulate the condition for obstruction bundles to be 
compatible with a forgetful map.

\begin{defn}\label{forgetcompobs}
Let $\{ E_{\bf p}({\bf x})\}$ be an obstruction bundle data on ${\mathcal M}_{\ell;k+1}(L;\beta)$
and 
$\{  E_{\overline{\bf p}}(\overline{\bf x})\}$ be an obstruction bundle data on ${\mathcal M}_{0;0}(L;\beta)$.
We say they are {\it compatible with respect to the forgetful map}\index{compatible with respect to the forgetful map}
if the following holds for sufficiently small neighborhoods $\mathscr U_{\overline{\bf p}}$ of  $\overline{\bf p}$
and   $\mathscr U_{{\bf p}}$ of ${\bf p}$.
\par
If ${\bf x} \in \mathscr U_{{\bf p}}$ and $\frak{forget}({\bf x}) \in \mathscr U_{\overline{\bf p}}$
then we have
\begin{equation}\label{12777}
E_{{\bf p}}({\bf x}) = E_{\overline{\bf p}}(\frak{forget}({\bf x})).
\end{equation}
\end{defn}
Let us elaborate on the equality (\ref{12777}).
We remark that the source curve $\Sigma_{\bf x}$ and $\Sigma_{\overline{\bf x}}$
are related in such a way that $\Sigma_{\overline{\bf x}}$ is obtained by shrinking 
some of the components of $\Sigma_{\bf x}$.
(The components we shrink are those where the maps are constant.)
Therefore we can regard  an element of the right hand side (whose support 
is away from nodal points) as a section of $u_{\bf{x}}^*TX \otimes \Lambda^{01}$.
Thus the equality (\ref{12777}) makes sense.
\par
Note that in particular it implies that elements of $E_{\bf p}(\bf x)$ 
are zero on the components of $\Sigma_{\bf x}$ which we shrink 
to obtain $\Sigma_{\overline{\bf x}}$.
\par
Now suppose we have an obstruction bundle data
 $\{ E_{\bf p}({\bf x})\}$   on ${\mathcal M}_{\ell;k+1}(L;\beta)$
and 
$\{ \overline E_{\overline{\bf p}}(\overline{\bf x})\}$  on ${\mathcal M}_{0;0}(L;\beta)$
which are compatible with the forgetful map.
We also assume that  $\{ E_{\bf p}({\bf x})\}$ preserves triviality
and enjoy the support preserving property.
We now define maps between Kuranishi neighborhoods $U_{\bf p}$ of 
${\bf p}$ and $U_{\overline{\bf p}}$ of $\overline{\bf p}$.
\begin{defn}
We recall
$$
\aligned
U_{\bf p} &= \{ {\bf x} \in \mathscr U_{{\bf p}} \mid  \overline{\partial}u_{{\bf x}} \in E_{{\bf p}}({\bf x})\},
\\
U_{\overline{\bf p}} &= \{ \overline{\bf x} \in \mathscr U_{\overline{\bf p}} \mid  \overline{\partial}u_{\overline{\bf x}} \in E_{\overline{\bf p}}(\overline{\bf x})\}.
\endaligned
$$
for an appropriate $\mathscr U_{{\bf p}}$.
By Lemma \ref{lem1218} we can shrink $\mathscr U_{{\bf p}}$ so that 
$\frak{foget}({\bf x}) \in \mathscr U_{\overline{\bf p}}$ holds for all 
${\bf x} \in U_{\bf p}$.
Then (\ref{12777}) implies $\frak{foget}({\bf x}) \in U_{\overline{\bf p}}$.
Thus we obtain a map
$$
 F_{\bf p}: U_{\bf p} \to U_{\overline{\bf p}}.
$$
\end{defn}
We can show that this map is strata-wise smooth in the same way as the proof of 
smoothness of coordinate change map (See \cite[Chapter 8]{foooexp}.).
There is an issue of smoothness of $F_{\bf p}$ across the strata. We will discuss it later at the beginning of 
Subsection \ref{CFforget}.

Again by (\ref{12777}) and its definition, 
%we have
%$$
%E_{{\bf p}}({\bf x}) \cong E_{\overline{\bf p}}(\frak{foget}({\bf x})).
%$$ 
%Thus 
the map $ F_{\bf p}$ is covered by a bundle map
$E_{\bf p} \to E_{\overline{\bf p}}$, which is a
fiberwise  isomorphism.
Its compatibilities with the Kuranishi map and with the parametrization map 
are obvious from construction.
Compatibility of $F_{\bf p}$ with coordinate change map 
is also immediate from the fact that 
coordinate change map is the inclusion map in our situation.
\par
We can summarize these properties as Proposition \ref{prop1220} below.
We define:
\begin{defn}\label{def1221}
Let $\mathcal M_i$ ($i=1,2$) be spaces 
with Kuranishi structures and 
${\frak f} : \mathcal M_1 \to \mathcal M_2$ a continuous map.
We say Kuranishi structures on $\mathcal M_i$ are {\it compatible} with $\frak f$\index{compatible with a continuous map} 
if the following holds.
\par
Let ${\bf p}_i \in \mathcal M_i$ and ${\frak f}({\bf p}_1) = {\bf p}_2$.
Let $(U^i_{{\bf p}_i},\mathcal E^i_{{\bf p}_i},\psi^i_{{\bf p}_i},s^i_{{\bf p}_i})$
be the Kuranishi chart of ${\bf p}_i$.
Then after shrinking $U^1_{{\bf p}_1}$ we have the following.
\begin{enumerate}
\item There exists a strata-wise smooth map $F_{{\bf p}_1} : U^1_{{\bf p}_1}
\to U^2_{{\bf p}_2}$ which is covered by a bundle map
$\tilde F_{{\bf p}_1} : \mathcal E^1_{{\bf p}_1}
\to \mathcal E^2_{{\bf p}_2}$.
\item $\tilde F_{{\bf p}_1}$ is a fiberwise isomorphism.
\item
$\tilde F_{{\bf p}_1}\circ s^1_{{\bf p}_1} =  s^2_{{\bf p}_2} \circ F_{{\bf p}_1}$.
\item 
On $(s^1_{{\bf p}_1})^{-1}(0)$ we have 
$\psi^2_{{\bf p}_2} \circ F_{{\bf p}_1} = \frak f \circ \psi^1_{{\bf p}_1}$.
\item
Let $(U_{{\bf p}'_i{\bf p}_i},\varphi_{{\bf p}'_i{\bf p}_i},\hat\varphi_{{\bf p}'_i{\bf p}_i})$
be a coordinate change for $i=1,2$.  (Here $U_{{\bf p}'_i{\bf p}_i} \subset 
U_{{\bf p}_i}$ is open, $\varphi_{{\bf p}'_i{\bf p}_i} : U_{{\bf p}'_i{\bf p}_i} \to 
U_{{\bf p}'_i}$ is a smooth map and $\hat\varphi_{{\bf p}'_i{\bf p}_i} : 
\mathcal E_{{\bf p}_i}\vert_{U_{{\bf p}'_i{\bf p}_i}} \to 
\mathcal E_{{\bf p}'_i}$ is a bundle map which covers it. 
(See \cite[Definition 3.6]{springer}.))
We require that it is compatible with $F_{{\bf p}_1}$, $\tilde F_{{\bf p}_1}$, that is,
$$
\varphi_{{\bf p}'_2{\bf p}_2} \circ F_{{\bf p}_1}  = 
F_{{\bf p}'_1} \circ \varphi_{{\bf p}'_1{\bf p}_1}, \qquad
\hat\varphi_{{\bf p}'_2{\bf p}_2} \circ \tilde F_{{\bf p}_1}  = 
\tilde F_{{\bf p}'_1} \circ \hat\varphi_{{\bf p}'_1{\bf p}_1}
$$
when both sides are defined.
\end{enumerate}
\end{defn}

Thus we have proved:
\begin{prop}\label{prop1220}
In the situation of Definition $\ref{forgetcompobs}$,  assume that 
$\{ E_{\bf p}({\bf x})\}$ preserves triviality
and has the support preserving property.
Then the Kuranishi structures we obtain from 
the obstruction bundle data are 
compatible with the forgetful map.
\end{prop}
To complete the proof of Proposition \ref{diskkura} 
it suffices to find a system of obstruction bundle data 
which satisfies the assumption of Proposition \ref{prop1220}
and is disk-component-wise.
For this purpose,
we will follow the strategy we explained in the last subsection and 
use outer collaring.
More precisely we find the next lemma whose proof is given at the end of this subsection.
\begin{lem}\label{lem12220}
There exists a system of collared obstruction bundle data on 
$\mathcal M_{\ell,k}(L;\beta)^{\boxplus 1}$ for $\ell,k \in \Z_{\ge 0}$ and $\beta \in \Pi_2(X;L)$
such that: 
\begin{enumerate}
\item It is compatible with the forgetful map of boundary marked points.
\item It is invariant under the cyclic permutation of boundary marked points.
\item It is  invariant under the permutation of interior marked points.
\item It has the evaluation transversality property  laid out in Definition $\ref{defn1223new}$
below.
\end{enumerate}
\end{lem}
\begin{rem}
We use only $\mathcal M_{\ell,k}(L;\beta)$ with $k\ge 1$ for applications in this paper.
However we need to include the case $k=0$ so that 
induction works to prove (1)(2)(4).  See \cite[Theorem 3.1]{fukaya:cyc}.
\end{rem}
The rest of this subsection is occupied to the proof of Lemma \ref{lem12220}.

Before launching in the proof of Lemma \ref{lem12220}, we first provide the precise description of
 the forgetful map compatibility and the evaluation transversality property mentioned therein.
 
Let $\hat{\bf p} = ({\bf p},s) \in \mathcal M_{\ell,k}(L;\beta)^{\boxplus 1}$.
We consider $u_{\bf p} : \Sigma_{\bf p} \to X$ which is a part of the data consisting ${\bf p}$.
We decompose $\Sigma_{\bf p}$ into irreducible
components as
$$
\Sigma_{\bf p}
= \bigcup_{a \in \mathcal A_{\bf p}^{\rm s}} \Sigma_{\bf p}(a)
\cup \bigcup_{a \in \mathcal A_{\bf p}^{\rm d}} \Sigma_{\bf p}(a).
$$
Here $\Sigma_{\bf p}(a)$ is a disk (resp. sphere) if $a \in \mathcal A_{\bf p}^{\rm s}$
(resp. $a \in \mathcal A_{\bf p}^{\rm d}$).
Let $u_{{\bf p},a}$ be the restriction of $u_{\bf p}$ to $\Sigma_{\bf p}(a)$.
The linearization of the non-linear Cauchy-Riemann equation 
defines a linear elliptic operator
\begin{equation}\label{form51}
\aligned
D_{u_{{\bf p},a}}\overline\partial : &L^2_{m+1}(\Sigma_{\bf p}(a),\partial \Sigma_{\bf p}(a);
u_{{\bf p},a}^*TX,u_{{\bf p},a}^*TL)\\
&\to L^2_{m}(\Sigma_{\bf p}(a);u_{{\bf p},a}^*TX\otimes \Lambda^{01})
\endaligned
\end{equation}
for $a \in \mathcal A_{\bf p}^{\rm d}$  and
\begin{equation}\label{form52}
D_{u_{{\bf p},a}}\overline\partial : L^2_{m+1}(\Sigma_{\bf p}(a);
u_{{\bf p},a}^*TX) \to L^2_{m}(\Sigma_{\bf p}(a);u_{{\bf p},a}^*TX\otimes \Lambda^{01})
\end{equation}
for $a \in \mathcal A_{\bf p}^{\rm s}$.
Here $L^2_{m+1}(\Sigma_{{\bf p}}(a),\partial \Sigma_{{\bf p},a}(a);
u_{{\bf p},a}^*TX,u_{{\bf p},a}^*TL)$ is the space of all sections of the bundle 
$u_{{\bf p},a}^*TX$ of $L^2_{m+1}$-class  whose boundary values lie in $u_{{\bf p},a}^*TL$.
Other spaces are appropriate Sobolev spaces of the sections. Take $m$ sufficiently 
large.
We take a direct sum
\begin{equation}\label{form5333}
\aligned
&\bigoplus_{a \in \mathcal A_{\bf p}^{\rm s}}L^2_{m+1}(\Sigma_{\bf p}(a);
u_{{\bf p},a}^*TX\otimes \Lambda^{01}) \\
&\oplus
\bigoplus_{a \in \mathcal A_{\bf p}^{\rm d}}
L^2_{m+1}(\Sigma_{\bf p}(a),\partial \Sigma_{\bf p}(a);
u_{{\bf p},a}^*TX,u_{{\bf p},a}^*TL).
\endaligned
\end{equation}
We also consider
\begin{equation}\label{form5444}
\bigoplus_{a \in \mathcal A_{\bf p}^{\rm s} \cup \mathcal A_{\bf p}^{\rm d}}L^2_{m}(\Sigma_{\bf p}(a);
u_{{\bf p},a}^*TX\otimes \Lambda^{01}).
\end{equation}
\begin{defn}
We define 
$L^2_{m}(\Sigma_{\bf p};u_{{\bf p}}^*TX\otimes \Lambda^{01})$
\index[syindex]{L1^2_{m}(\Sigma_{\bf p};u_{{\bf p}}^*TX\otimes \Lambda^{01})@$L^2_{m}(\Sigma_{\bf p};u_{{\bf p}}^*TX\otimes \Lambda^{01})$}
to be the Hilbert space (\ref{form5444}).
\par
We define a Hilbert space $W^2_{m+1}(\Sigma_{\bf p},\partial\Sigma_{\bf p};u_{{\bf p}}^*TX,u_{{\bf p}}^*TL)$
\index[syindex]{W^2_{m+1}(\Sigma_{\bf p},\partial\Sigma_{\bf p};u_{{\bf p}}^*TX,u_{{\bf p}}^*TL)@$W^2_{m+1}(\Sigma_{\bf p},\partial\Sigma_{\bf p};u_{{\bf p}}^*TX,u_{{\bf p}}^*TL)$}
as the subspace of the Hilbert space (\ref{form5333}) consisting of  elements 
$
\sum_{a \in \mathcal A_{\bf p}^{\rm s} \cup 
\mathcal A_{\bf p}^{\rm d}}V_a
$
(where $V_a$ is a section on $\Sigma_{\bf p}(a)$)
with the following properties.
Let $p \in \Sigma_{\bf p}$ be a nodal point. We take $a_1(p)$, $a_2(p)$
such that
$\{p\} = \Sigma_{\bf p}(a_1(p)) \cap \Sigma_{\bf p}(a_2(p))$. We require
$$
V_{a_1(p)}(p) = V_{a_2(p)}(p).
$$
We require this condition at all the nodal points $p$.
\end{defn}
The operators (\ref{form51}), (\ref{form52}) induce a Fredholm  operator
\begin{equation}\label{form523}
D_{u_{{\bf p}}}\overline\partial : W^2_{m+1}(\Sigma_{\bf p},\partial\Sigma_{\bf p};u_{{\bf p}}^*TX,u_{{\bf p}}^*TL) \to L^2_{m}(\Sigma_{\bf p};u_{{\bf p}}^*TX\otimes \Lambda^{01}).
\end{equation}
\begin{defn}{\rm (Evaluation transversality)}\label{defn1223new}
We say that $E_{{\bf p}}({\bf p}) \subset L^2_{m}(\Sigma_{\bf p};u_{{\bf p}}^*TX\otimes \Lambda^{01})$
satisfies the {\it evaluation transversality}\index{evaluation transversality} if 
$$
{\rm EV}_0 : (D_{u_{{\bf p}}}\overline\partial)^{-1}(E_{{\bf p}}({\bf p})) \to T_{u_{{\bf p}}(w_0)}L
$$
defined by $(V_a) \mapsto V_{a_0}(w_0)$ is surjective.
Here $\Sigma_{{\bf p}}(a_0)$ is the irreducible component containing the $0$-th boundary marked point $w_0$.
\end{defn}
The next lemma is easy to show.
\begin{lem}
If $E_{{\bf p}}({\bf p})$ satisfies the evaluation transversality then ${\rm ev}_0 : U_{\bf p} \to L$ becomes a submersion, by shrinking $U_{\bf p}$ if necessary.
\end{lem}

We  now need to make precise  the compatibility with the  forgetful map of boundary marked points
for a system of the Kuranishi structures defined on the outer collaring $\mathcal M_{\ell,k}(L;\beta)^{\boxplus 1}$.
For this purpose, we first extend the forgetful map 
to the outer collar.
We consider the forgetful map
\begin{equation}\label{formula138}
\frak{forget} : {\mathcal M}_{\ell;k+1}(L;\beta) \to {\mathcal M}_{\ell;0}(L;\beta).
\end{equation}
Let ${\bf p} \in {\mathcal M}_{\ell;k+1}(L;\beta)$
and ${\frak I}_{\bf p}$ be the set of boundary nodes of ${\bf p}$.
Let $\frak{forget}({\bf p}) = \overline{\bf p}$
Since there exists a map 
$\Sigma_{\bf p} \to \Sigma_{\overline{\bf p}}$ between source curves which shrinks some of the extended disk components
(that becomes unstable after forgetting marked points), 
there exists a surjective map
$\mathcal S_{\bf p} : {\frak I}_{\bf p} \to {\frak I}_{\overline{\bf p}}$.
Note that ${\mathcal M}_{\ell;k+1}(L;\beta)^{\boxplus 1}$ 
is the set of pairs $({\bf p},\vec s)$ where
$\vec s$ may be regarded as a map
${\frak I}_{\bf p} \to [-1,0]^{\vert \vec s\vert}$.
We want to extend (\ref{formula138}) to a map:
\begin{equation}\label{formula138+}
\frak{forget}^{\boxplus 1} : {\mathcal M}_{\ell;k+1}(L;\beta)^{\boxplus 1} \to {\mathcal M}_{\ell;0}(L;\beta)^{\boxplus 1}
\end{equation}
by putting
$$
\frak{forget}^{\boxplus 1}({\bf p},\vec s) = (\overline{\bf p},\vec s\,^{\prime}).
$$
A naive choice is:
$$
\vec s\,^{\prime}(i) = \min \{ \vec s(j) \mid \mathcal S_{\bf p}(j) = i \}.
$$
It is easy to see that this choice defines a 
continuous map $\frak{forget}^{\boxplus 1}$.
However since $\min$ is not differentiable
this choice does not give a smooth map 
$\frak{forget}^{\boxplus 1}$.
It seems impossible to define $\frak{forget}^{\boxplus 1}$
appropriately so that it becomes smooth.
We next explain the way to define $\frak{forget}^{\boxplus 1}$
and define an appropriate class of obstruction bundle data so that we can 
pull back obstruction bundle data in such a class (on the outer collaring) 
by a forgetful map.
\par
\begin{figure}[h]
\centering
\includegraphics[scale=0.3]{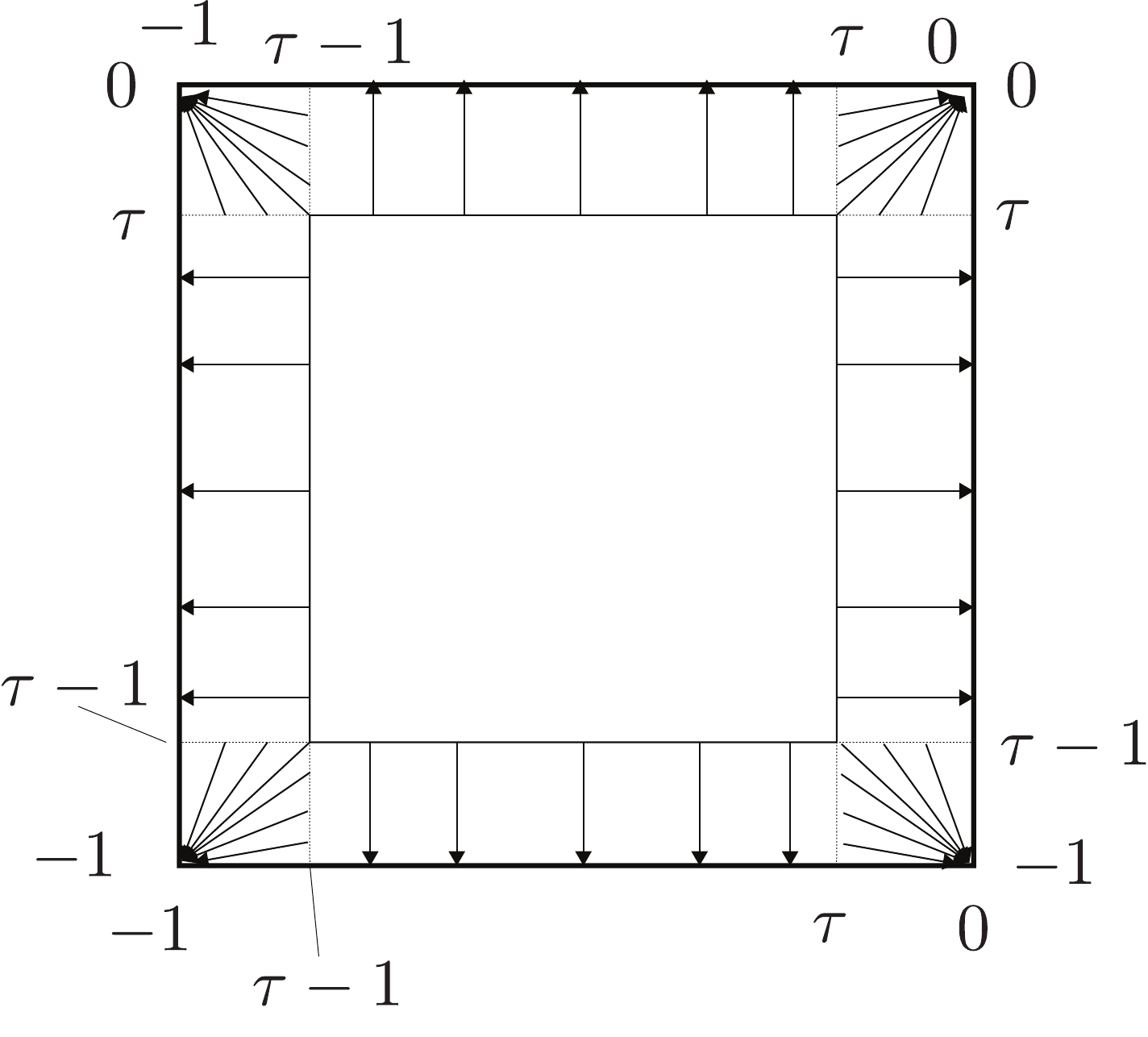}
\caption{Projection to bi-collar}
\label{zu11}
\end{figure}
\begin{figure}[h]
\includegraphics[scale=0.4]{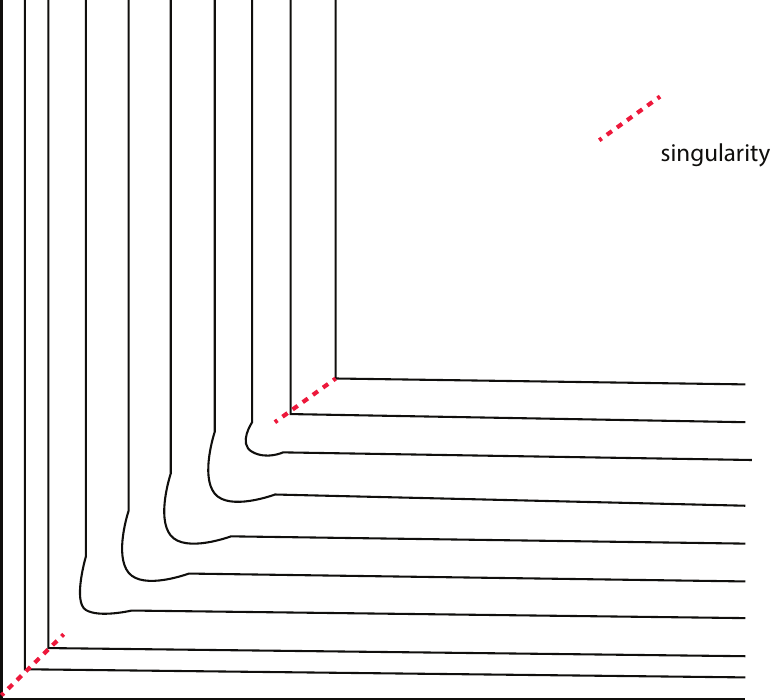}
\caption{The level set of $\mathscr M_2$ }
\label{FigureM2}
\end{figure}

We begin with defining $\Pi'_{\tau}$ which 
is similar to but is different from $\Pi_{\tau}$ we used in 
Subsection \ref{diskreview}.\index[syindex]{Pitauprime@$\Pi'_{\tau}$}
\begin{defn}
Let $\tau \in (0,1/2)$.
We define 
$$
\Pi'_{\tau}: 
{\mathcal M}_{\ell;k+1}(L;\beta)^{\boxplus 1} \to {\mathcal M}_{\ell;k+1}(L;\beta)^{\boxplus 1}
$$
by $\Pi'_{\tau}({\bf p},\vec s) = ({\bf p},\vec s\,^{\prime})$ such that:
\begin{enumerate}
\item $s'(i) = s(i)$ if $-\tau > s(i) > \tau-1$.
\item
$s'(i) = -1$ if $s(i) \le \tau -1$.
\item
$s'(i) = 0$ if $s(i) \ge  -\tau$. %\marginpar{Figure putted}
\end{enumerate}  
Here $\vec s, \vec s^{\,\prime}  \in [-1,0]^m$ ($m$ is the number of boundary nodes.)
See Figure \ref{zu11}
\end{defn}
\begin{lem}\label{lem1314}
For the above given $\tau \in (0,1/2)$, we
fix a constant $\epsilon > 0$ with $\epsilon < \tau$.
There exist functions 
$\mathscr M_m : [-1,0]^m \to [-1,0]$ with the following properties.\index[syindex]{Myscrk@$\mathscr M_k$}
\begin{enumerate}
\item
If there exists $s_i =-1$ then $\mathscr M_m(\vec s) = -1$.
\item
$\Pi'_{\tau} \circ \mathscr M_m = \mathscr M_m \circ \Pi'_{\tau}$.
\item
$\mathscr M_m$ is invariant under the permutation of the factors 
of the domain.
\item
$\mathscr M_m(s_1,\dots,s_{m-1},0) = \mathscr M_{m-1}(s_1,\dots,s_{m-1})$.
\item 
$\mathscr M_m$ is continuous.
\item
$\mathscr M_m$ is smooth except on the following region.
\begin{enumerate}
\item
There exists $i$ such that 
$s_i  < \epsilon -1$.
\item
For all $i$ 
$s_i >  -\epsilon$.
\end{enumerate}
\end{enumerate}
\end{lem}
We can prove the lemma by induction on $k$ 
and smoothing $\min$ appropriately.
The level set of $\mathscr M_2$ is drawn in Figure \ref{FigureM2}.
\par
Now we define\index[syindex]{forgetfrakboxplus1@$\frak{forget}^{\boxplus 1}$} 
\begin{equation}\label{formula138+2}
\frak{forget}^{\boxplus 1} : {\mathcal M}_{\ell;k+1}(L;\beta)^{\boxplus 1} \to {\mathcal M}_{\ell;0}(L;\beta)^{\boxplus 1}
\end{equation}
by putting
$$
\frak{forget}^{\boxplus 1}({\bf p},\vec s) = (\overline{\bf p},\vec s\,^{\prime})
$$
where 
$$
\vec s\,^{\prime}(i) = \mathscr M_{m} ((\vec s(j); {\mathcal S_{\bf p}(j) = i})).
$$
Here $m$ is the cardinarity of the set $\{j \mid \mathcal S_{\bf p}(j) = i\}$.
\par
We define 
\begin{equation}\label{formula138+2+1}
\frak{forget}^{\boxplus 1} : {\mathcal X}_{\ell;k+1}(L;\beta)^{\boxplus 1} \to {\mathcal X}_{\ell;0}(L;\beta)^{\boxplus 1}
\end{equation}
in the same way.
\par
We remark that the notions of `preserving triviality'
and `support preserving property' can be defined for the 
obstruction bundle data of ${\mathcal M}_{\ell;k+1}(L;\beta)^{\boxplus 1}$
in the same way.
Then Definition \ref{forgetcompobs} (compatibility of the obstruction bundle data with 
forgetful map) is adapted to our situation in an obvious way.
\par
We next introduce a class of obstruction bundle data 
which can be pulled back by the map $\frak{forget}^{\boxplus 1}$.
It is called a bi-collared\index{bi-collared} obstruction bundle data.

\begin{defn}\label{1210newnew}
An obstruction bundle data $\{E_{\hat{\bf p}}(\widehat{\bf x})\}$ is said to be $\tau$-bi-collared\index{bi-collared obstruction bundle data} 
if 
\begin{equation}
E_{\Pi'_{\tau}(\hat{\bf p})}(\Pi'_{\tau}(\hat{\bf x}))  =   E_{\hat{\bf p}}(\widehat{\bf x}).
\end{equation}
\end{defn}

Now we are in the position to complete the proof of Lemma \ref{lem12220}.
\begin{proof}[Proof of Lemma \ref{lem12220}]
The proof will be given by induction on the partial order $<$ on $(\omega(\beta),\ell,0)$.
\par
We first claim the following. Suppose an obstruction bundle data $E_{\hat{\bf p}}(\widehat{\bf x})$ is given on 
${\mathcal M}_{\ell;0}(L;\beta)^{\boxplus 1}$ which is bi-collared. Then 
we can define one on ${\mathcal M}_{\ell;k+1}(L;\beta)^{\boxplus 1}$ by 
\begin{equation}\label{form1210}
E_{{\bf p}}({\bf x}) = E_{\overline{\bf p}}(\frak{forget}({\bf x})).
\end{equation}
It is  bi-collared, preserves triviality
and has the support preserving property. 
\par
The claim is mostly obvious. The only point to be discussed is 
smoothness (Definition \ref{defn51} (2)).
This is because the forgetful map on the outer collar is not smooth.
However the singular set of
 $\mathscr M_m$
is already  described in Lemma \ref{lem1314} (6), (a) and (b).
We observe that in such a case the $\tau$-bi-collared-ness
implies that various objects appearing in $E_{\bf p}({\bf x})$ are constant 
in the $s_i$ direction.  
(We use $\epsilon <\tau$ here.)
Therefore, the non-smooth-ness of  $\mathscr M_m$
does not affect the smooth-ness of (\ref{form1210}).
\par
Note this matter of smoothness concerns only the case $\beta \ne 0$. For we put 
obstruction bundle data to be $0$ for the case $\beta = 0$.
\par
Now we will 
construct an obstruction bundle data on ${\mathcal M}_{\ell;k}(L;\beta)^{\boxplus 1}$ 
by the induction on the partial order $<$ on $(\omega(\beta),\ell,0)$.
\par
The first step of the induction is a consequence of 
Lemma \ref{lem12222}. Namely we use it for ${\mathcal M}_{\ell;0}(L;\beta)^{\boxplus 1}$
with  $(\omega(\beta),\ell,0)$ minimal and apply (\ref{form1210}).
\par
Let us discuss the inductive step.
Suppose that for all $(\beta',\ell',k')$ with 
$(\omega(\beta'),\ell',0) < (\omega(\beta),\ell,0)$
we constructed an obstruction bundle data on
${\mathcal M}_{\ell';k'}(L;\beta')^{\boxplus 1}$.
We will construct one on ${\mathcal M}_{\ell;k}(L;\beta)^{\boxplus 1}$.
Using (\ref{form1210}) it suffices to study the case $k=0$.
\par
We stratify ${\mathcal M}_{\ell;0}(L;\beta)^{\boxplus 1}$ as follows.
We define an element $\hat{{\bf p}} = ({\bf p},\vec s)$ to be   
in $G_m{\mathcal M}_{\ell;0}(L;\beta)^{\boxplus 1}$
if and only if:
$$
\# \{ i \mid s_i \ne 0 \} \le m.
$$
We also put
$$
U(G_m) = 
\{ \hat{{\bf p}} \in {\mathcal M}_{\ell;0}(L;\beta)^{\boxplus 1} 
\mid \Pi'_{\tau}(\hat{{\bf p}}) \in
G_m{\mathcal M}_{\ell;0}(L;\beta)^{\boxplus 1} \cup 
\partial {\mathcal M}_{\ell;0}(L;\beta)^{\boxplus 1} \}.
$$
See Figures \ref{FigureG0andUG0}, \ref{FigureG1andUG1}.
We will define $E_{\hat{{\bf p}}}(\hat{{\bf x}})$ for 
$\hat{{\bf p}} \in U(G_m)$ by an upward induction on $m$.
\par
\begin{figure}[h]
\includegraphics[scale=0.4]{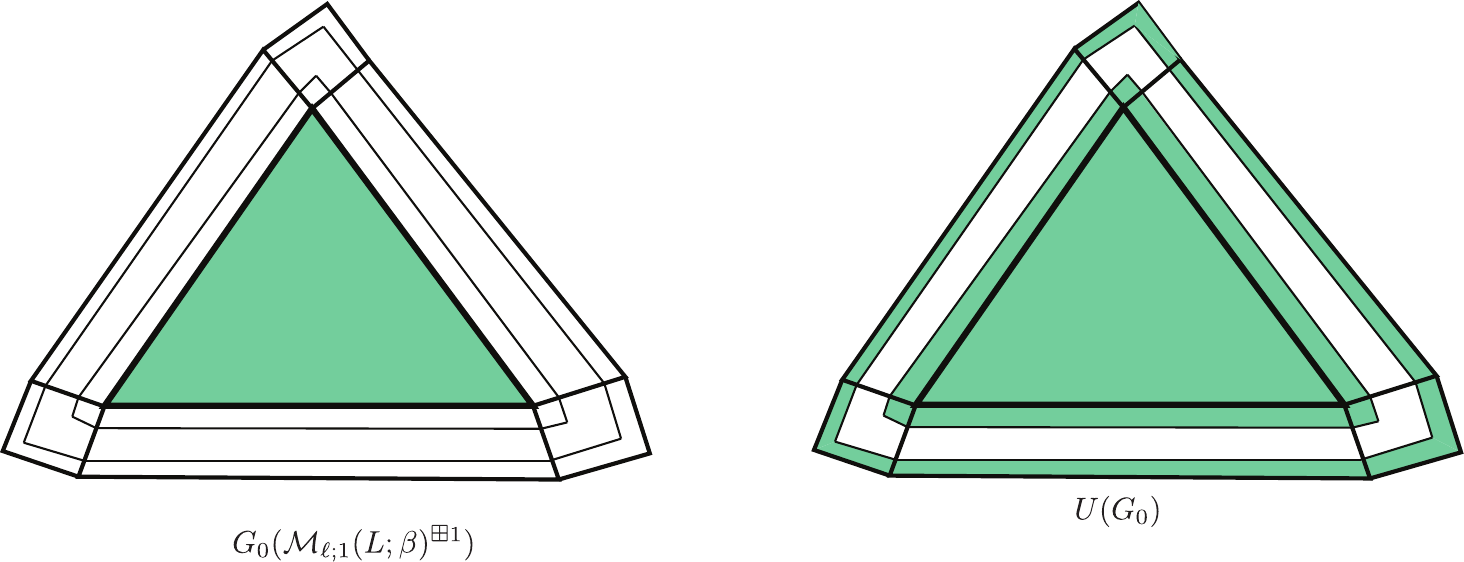}
\caption{$G_0{\mathcal M}_{\ell;0}(L;\beta)^{\boxplus 1}$ and $U(G_0)$}
\label{FigureG0andUG0}
\end{figure}
\begin{figure}[h]
\includegraphics[scale=0.4]{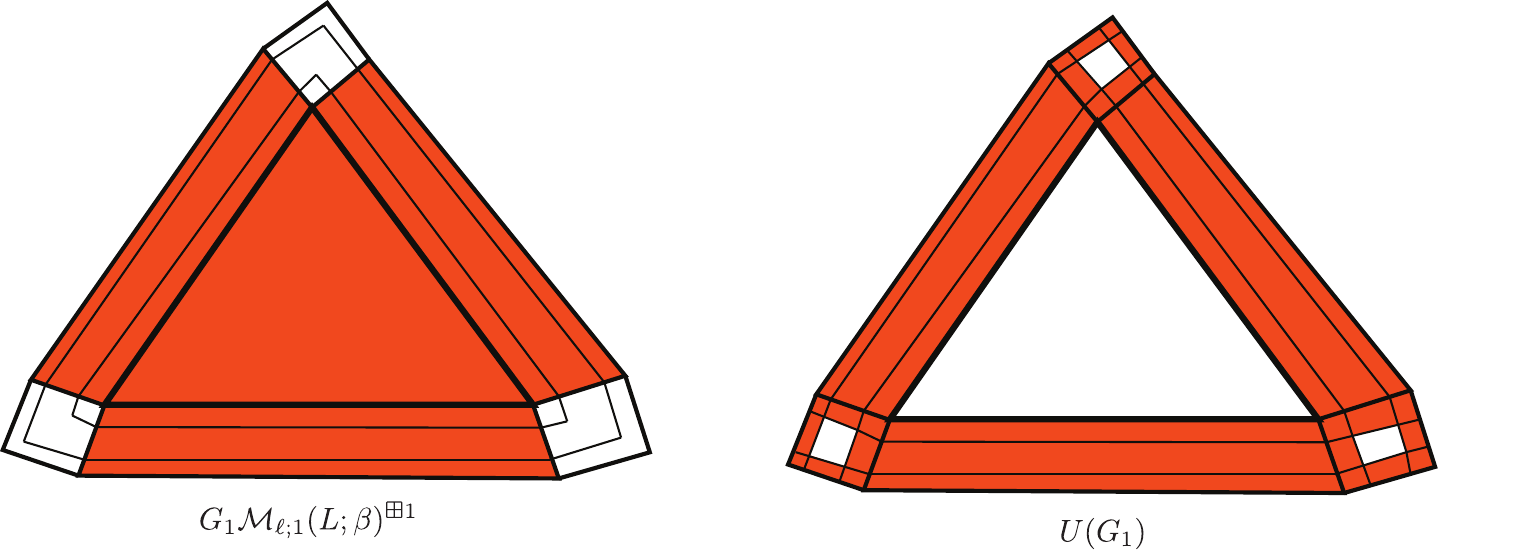}
\caption{$G_1{\mathcal M}_{\ell;0}(L;\beta)^{\boxplus 1}$ and $U(G_1)$}
\label{FigureG1andUG1}
\end{figure}
\par
We first consider the case when $m=0$. 
Note that $G_0{\mathcal M}_{\ell;0}(L;\beta)^{\boxplus 1}
= {\mathcal M}_{\ell;0}(L;\beta)$ which is disjoint from 
$\partial{\mathcal M}_{\ell;0}(L;\beta)^{\boxplus 1}$.
So we can construct required $E_{\hat{{\bf p}}}(\hat{{\bf x}})$ 
on $U(G_0)$ as follows. 
We first define it for $\hat{{\bf p}} \in \partial{\mathcal M}_{\ell;0}(L;\beta)^{\boxplus 1}$ by
(\ref{newform124}).  We construct it for $\hat{{\bf p}} \in  G_0{\mathcal M}_{\ell;0}(L;\beta)^{\boxplus 1}$
by Lemma \ref{lem12222}. Then we extend it 
to  $U(G_0)$ by (\ref{form1210}).
\par
Suppose we defined $E_{\hat{{\bf p}}}(\hat{{\bf x}})$  for $\hat{{\bf p}} \in U(G_{m-1})$.
We use Lemma \ref{lem12222} to extend it to ${\mathcal M}_{\ell;0}(L;\beta)^{\boxplus 1}$.
We restrict it to $\hat{{\bf p}} \in G_m{\mathcal M}_{\ell;0}(L;\beta)^{\boxplus 1}$.
We then use (\ref{form1210}) to extend it to $\hat{{\bf p}} \in U(G_{m})$.
Note that by taking $\tau$ smaller if necessary it is still bi-collared.
\par
Lemma \ref{lem12220} (2)   (invariance under  the cyclic permutation of boundary marked points)
is immediate from construction. Lemma \ref{lem12220} (3)
(invarianceby the permutation of interior marked points) can be achieved by 
requiring it inductively during construction.
\par
To prove Lemma \ref{lem12220} (4) we will construct obstruction bundle data 
on ${\mathcal M}_{\ell;0}(L;\beta)^{\boxplus 1}$ 
inductively so that the next (*) is satisfied.
\begin{enumerate}
\item[(*)] 
We consider $\frak{forget} : {\mathcal M}_{\ell;1}(L;\beta)^{\boxplus 1} 
\to {\mathcal M}_{\ell;0}(L;\beta)^{\boxplus 1}$. Then the obstruction 
bundle data on ${\mathcal M}_{\ell;1}(L;\beta)^{\boxplus 1}$ pulled back by $\frak{forget}$
from one on ${\mathcal M}_{\ell;0}(L;\beta)^{\boxplus 1}$
has the evaluation transversality property (Definition \ref{defn1223new}).
\end{enumerate}
To show (*) we go back to the proof of Lemma \ref{lem12222}.
We can choose $E^0_{{\bf p}}({\bf p})$ so that (*) is satisfied.
This is \cite[Lemma 3.1]{fukaya:cyc}. Then by construction (during the proof of Lemma \ref{lem12222})
the Kuranishi structure we obtain by applying Lemma \ref{lem12222} enjoy (*) as well.
\par
The proof of Lemma \ref{lem12220} is now complete.
\end{proof}
The proof  of Proposition \ref{diskkura} 
is complete.
\qed

\subsection{Construction of Kuranishi structure of the 
moduli space of pseudo-holomorphic polygons.}
\label{polygonforget}

In this subsection, we prove Proposition \ref{prop92}.
Proposition \ref{Kuraeistspoly} is proved at the same time.
The proof of Proposition \ref{prop94} is similar.
The proof of Proposition \ref{prop92} is similar to the proof of last two subsections but we need to modify it in a few points.
\par
We consider the moduli space
${\mathcal M}_{\ell;\vec k}((\vec{\kappa},\vec p);B;\vec{\frak f}_{\bf e})$
and define its ambient set 
${\mathcal X}_{\ell;\vec k}((\vec{\kappa},\vec p);B;\vec{\frak f}_{\bf e})$
in the same way.

The partial topology of the pair
$({\mathcal X}_{\ell;\vec k}((\vec{\kappa},\vec p);B;\vec{\frak f}_{\bf e}),
{\mathcal M}_{\ell;\vec k}((\vec{\kappa},\vec p);B;\vec{\frak f}_{\bf e}))$
is also defined in the same way. (See Definition \ref{defnpartialtopo}.)
The notions of obstruction bundle data 
and its disk-component-wise-ness are 
defined in the same way.
Note that we use the obstruction bundle data
obtained in Subsection \ref{diskforget} for the 
first factor of $(\ref{218})$.
\par
\begin{figure}[h] 
\centering
\includegraphics[scale=0.5]{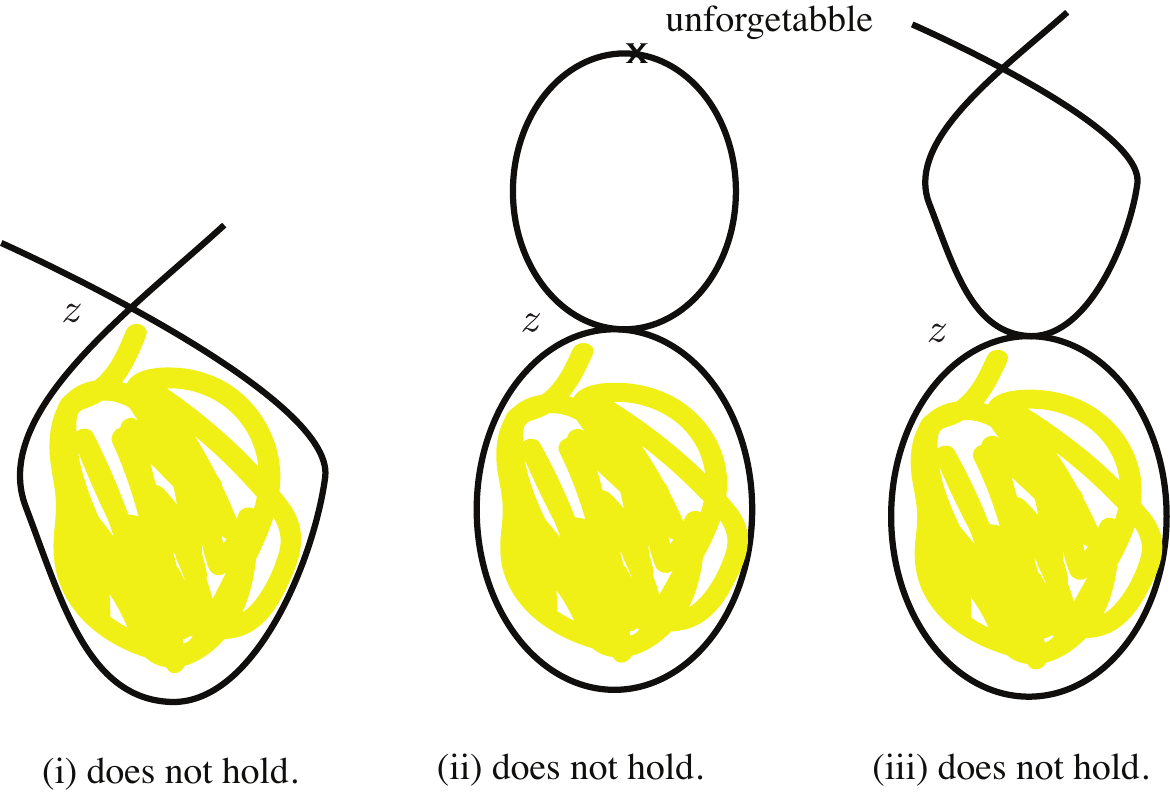}
\caption{$z$ is unforgetabble}
\label{Figureunforget}
\end{figure}

We define $\vec k \setminus \vec{\frak f}_{\bf e}$ by removing all the 
forgetable marked points.
We then have a forgetful map
\begin{equation}\label{1210form}
\frak{forget}:
{\mathcal M}_{\ell;\vec k}((\vec{\kappa},\vec p);B;\vec{\frak f}_{\bf e})
\to 
{\mathcal M}_{\ell;\vec k \setminus \vec{\frak f}_{\bf e}}((\vec{\kappa},\vec p);B;\emptyset)
\end{equation}
and its analogue for the ambient set.
We can define the notion of the support preserving property
in the same way.
We need to  modify the property `preserving triviality' \index{preserves triviality}
a bit as follows.
Let ${\bf p} \in {\mathcal M}_{\ell;\vec k}((\vec{\kappa},\vec p);B;\vec{\frak f}_{\bf e})$
and $\Sigma_{{\bf p}_0}$ be a connected union of extended disk components
of the source curve $\Sigma_{\bf p}$.
We regard each boundary node $z$ on $\Sigma_{{\bf p}_0}$ which intersects 
with the closure of $\Sigma_{\bf p} \setminus \Sigma_{{\bf p}_0}$
as a forgetable marked point of $\Sigma_{{\bf p}_0}$,
if the following 3 conditions hold. %\marginpar{I put a figure.}
\begin{enumerate}
\item[(i)]
$z$ is a diagonal marked point.
\item[(ii)]
Let $\Sigma_{{\bf p}'}$ be the closure of the connected component of   
$\Sigma_{\bf p} \setminus \Sigma_{{\bf p}_0}$ which contains $z$.
Then $\Sigma_{{\bf p}'}$ contains no unforgetable marked points.
\item[(iii)]
None of the marked points of the 
extended disk components of $\Sigma_{{\bf p}'}$ is a switching marked point.   
\end{enumerate}
See Figure \ref{Figureunforget}.
We say $z$ is forgetable if it is not unforgetable.
\par
We consider ${\bf p}_0$ together with unforgetable marked points, 
that is, either an unforgetable marked point on $\Sigma_{{\bf p}_0}$ 
or a nodal point which is regarded as unforgetable 
marked point by the above convention.
We say ${\bf p}_0$ is a trivial component if 
the map $u_{{\bf p}_0}$ is constant and 
the number of unforgetable marked points is at most two.
\par
We say an obstruction bundle data preserves 
triviality if the restriction of any element of $E_{{\bf p}}({\bf p})$ 
is $0$ on trivial components.
\par
Using these definitions, we can easily generalize Lemma \ref{lem1218} 
to our situation.
\par
We can formulate the compatibility of obstruction bundle
data with the forgetful map (\ref{1210form}) in the same way as in Subsection \ref{diskforget}.
Then we can show if 
obstruction bundle
data enjoy the support preserving property,  preserves 
triviality property, and are compatible 
with the forgetful map (\ref{1210form})
and  disk-component-wise, then 
the system of Kuranishi structures we obtain is compatible 
with the forgetful map (\ref{1210form}) and enjoys other properties in 
Proposition \ref{prop92}.
We can construct such a system of obstruction bundle
data in the same way as Subsection \ref{diskforget} (Lemma \ref{lem12220}).
We have thus proved Proposition \ref{prop92}.
The proofs of Propositions  \ref{Kuraeistspoly}, \ref{prop92} and \ref{prop94} are similar.

\section{Construction of CF perturbations.}
\label{sec:CRperturb}

 %\marginpar{\color{blue} This section is new.}

\subsection{A brief review of a CF-perturbation and the integration along the fiber.}
\label{CFreview}

This subsection and the next are reviews of \cite{springer} and explain how to construct 
a system of CF-perturbations and use it to define the integration along the fiber.
We first briefly review the definition of a CF-perturbation.
Let $\Gamma $ be a finite group and $\mathcal E \to V$  a 
$\Gamma$ equivariant vector bundle on a manifold $V$ 
with effective $\Gamma$ action.
It is, by definition, a vector bundle on $U = V/\Gamma$.
Let $s$ be a section of  ${\mathcal E} \to V$.  We require that $s$ is smooth 
and $\Gamma$ equivariant.
We consider the situation where $U = V/\Gamma$ is a Kuranishi neighborhood, 
${\mathcal E}$ is an obstruction bundle and $s$ is a Kuranishi map.
\par
We also consider the situation where we have a map 
$f: U \to N$ to a manifold $N$ that is $\Gamma$ equivariant. 
(The $\Gamma$ action on $N$ is trivial.)
We will apply the situation to the case when $f$ is an evaluation map.
We also assume that $f$ is a submersion.
(This condition will correspond to the weak submersivity in the case of Kuranishi structure.)
A CF-perturbation\index{CF-perturbation} on a Kuranishi chart is a triple $(W,\{\frak s_{\epsilon}\},\chi)$.
Here $W$ is a finite dimensional $\Gamma$-vector space, 
$\chi$ is a smooth $\Gamma$ invariant differential form on $W$ of degree $\dim W$ of compact support on $W$ 
whose support is a neighborhood of $0$ and $\int_W \chi = 1$,
and $\frak s_{\epsilon}$ is an $\epsilon \in [0,1)$ 
parametrized family of ($\Gamma$-equivariant) 
sections of ${\mathcal E} \times W \to  V \times W$, 
such that 
$\frak s_0(x,\xi) = s(x)$.
We say $f$ is strongly submersive with respect to 
$(W,\{\frak s_{\epsilon}\},\chi)$ 
if $(f \circ \pi)\vert_{\frak s_{\epsilon}^{-1}(0)}: \frak s_{\epsilon}^{-1}(0) \to N$ is a submersion, 
on the support of $\chi$ when $\epsilon \ne 0$ is sufficiently small.  Here $\pi : V \times W \to V$ is the 
projection.
\par
When $f$ is strongly submersive with respect to 
$\frak S = (W,\{\frak s_{\epsilon}\},\chi)$
and $h$ is a $\Gamma$ equivariant differential form on $V$ with compact support,
we define the integration along the fiber\index{integration along the fiber} by the formula:
\begin{equation}\label{form13-1}
f_{!}(h;\frak S) = \frac{1}{\#\Gamma} ((f \circ \pi)\vert_{\frak s_{\epsilon}^{-1}(0)})_{!}(\pi^* h \wedge \chi).
\end{equation}
Here $((f \circ \pi)\vert_{\frak s_{\epsilon}^{-1}(0)})_{!}$ 
is the usual integration along the fiber, which is well defined 
since  $(f \circ \pi)\vert_{\frak s_{\epsilon}^{-1}(0)}$ is a submersion and $\pi^* h \wedge \chi$ 
has a compact support. 
\begin{rem}
The value (\ref{form13-1}) is defined for each sufficiently small positive $\epsilon$.
It depends on $\epsilon$ and defines a smooth function on an unspecified open interval
$(0,\epsilon_0)$ for some $\epsilon_0 > 0$.
\end{rem}
\par
We consider an appropriate equivalence relation between 
$(W,\{\frak s_{\epsilon}\},\chi)$, $(W',\{\frak s'_{\epsilon}\},\chi')$
an equivalence class of which is called a CF-perturbation.
For example if $W' = W \times W''$,
$\chi' = \chi \times \chi''$
and $\frak s'_{\epsilon}$ is the pull back of $\frak s_{\epsilon}$
on the support of $\chi'$ then 
$(W,\{\frak s_{\epsilon}\},\chi)$ is equivalent to $(W',\{\frak s'_{\epsilon}\},\chi')$.
See \cite[Definition 7.6]{springer} for the definition of the equivalence relation.
Integration along the fiber does not change if we change the 
representative of an equivalence class (\cite[Lemma 7.13]{springer}).
\par
If a moduli space $\mathcal M$ has a Kuranishi structure, 
we can define the notion of a CF-perturbation\index{CF-perturbation} globally,
roughly, as follows.
We first review the notion of a coordinate change of 
Kuranishi charts.
Let $p \in \mathcal M$.
Its Kuranishi chart is $\mathcal U_p = (U_p,{\mathcal E}_p,s_p,\psi_p)$
where $U_p = V_p/\Gamma_p$, ${\mathcal E}_p$ and $s_p$ 
are as above. 
(Namely $(V,{\mathcal E},s,\Gamma) = (V_p,{\mathcal E}_p,s_p,\Gamma_p)$.)
Here $\Gamma_p$ is the group of automorphisms of the object $p$.
$\psi_p$ is a map from $s_p^{-1}(0)$ to $\mathcal M$ 
which is a homeomorpism onto a neighborhood of $p$.
We call ${\mathcal E}_p$ an obstruction bundle, $s_p$ a Kuranishi map 
and $\psi_p$ a parametrization.
\par
Let $q \in \mathcal M$ be in the image of $\psi_p$.
Suppose $\mathcal U_q = (U_q,{\mathcal E}_q,s_q,\psi_q)$ is a 
Kuranishi chart of $q$, where $U_q = V_q/\Gamma_q$.
A coordinate change from $\mathcal U_q$ to $\mathcal U_p$ 
is $(\varphi_{pq},\widehat{\varphi}_{pq})$.
Here $\varphi_{pq} : U_q \to U_p$ is an embedding of 
orbifolds, that is to say, a pair of a 
smooth embedding $\tilde{\varphi}_{pq}  : 
V_q \to V_p$ of manifolds and an injective 
homomorphism $\psi_{pq} : \Gamma_q \to \Gamma_p$
of groups, such that for each $x \in V_q$, 
$\psi_{pq}$ induces an isomorphism between
$\{\gamma \in \Gamma_q \mid \gamma(x) = x\}$
and 
$\{\gamma \in \Gamma_p \mid \gamma(\tilde{\varphi}_{pq}(x)) = \tilde{\varphi}_{pq}(x)\}$.
\par
$\widehat{\varphi}_{pq}$ is an injective bundle map ${\mathcal E}_q \to {\mathcal E}_p$
which covers $\varphi_{pq}$.
\par
We require an obvious compatibility condition between $(\varphi_{pq},\widehat{\varphi}_{pq})$
and Kuranishi maps, parametrizations.
(See \cite[Definition 3.3]{springer}.)
Moreover we require (tangent bundle condition) that the differential $D_{\varphi_{pq}(x)}s_p$ of the Kuranishi map induces an isomorphism
$$
\frac{T_{\varphi_{pq}(x)}U_p}{(D_{x}\varphi_{pq})(T_xU_q)}
\cong
\frac{({\mathcal E}_p)_{\varphi_{pq}(x)}}{\widehat{\varphi}_{pq}(({\mathcal E}_q)_{x})}
$$
if $s_{p}(x) = 0$. See \cite[Definition 3.3 (5)]{springer}.
\par
We require that such a coordinate change exists after shrinking $U_q$, if necessary.
\par
Kuranishi structure $\widehat{\mathcal U}$ on $\mathcal M$
assigns a Kuranishi chart $\mathcal U_p = (U_p,{\mathcal E}_p,s_p,\psi_p)$
for each $p \in \mathcal M$
and coordinate change $(\varphi_{pq},\widehat{\varphi}_{pq})$ for $q 
\in \mathcal M$ sufficiently 
close to $p$.
An appropriate cocycle condition is required for the coordinate changes.
See \cite[Definition 3.9]{springer}.
\par
We next discuss a {\it CF-perturbation}\index{CF-perturbation} of a 
Kuranishi structure.
It assigns $\frak S_p: = (W_p,\{\frak s_{p,\epsilon}\},\chi_p)$ as above to each $p \in \mathcal M$.  More precisely 
it assigns an equivalence class of such objects.
If $q \in \mathcal M$ is close to $p \in \mathcal M$, 
we need compatibility of 
$(W_p,\{\frak s_{p,\epsilon}\},\chi_p)$,
$(W_q,\{\frak s_{q,\epsilon}\},\chi_q)$ with coordinate change 
$(\varphi_{pq},\widehat{\varphi}_{pq})$,
which can be described as follows.
\begin{enumerate}
\item[(i)]
We can take representatives so that $W_p = W_q$ and  $\chi_p = \chi_q$.
\item[(ii)]
We have
$$
\frak s_{p,\epsilon}(\varphi_{pq}(x),\xi) = (\widehat{\varphi}_{pq},{\rm id})(\frak s_{q,\epsilon}(x,\xi)).
$$
Here $x \in U_q$ and $\xi \in W_p = W_q$.
(Note that $(\widehat{\varphi}_{pq},{\rm id})$ is a map from $\mathcal E_q \times W_q$ to 
$\mathcal E_p \times W_p$.)
\end{enumerate}
See \cite[Definition 7.49]{springer}.  The definition there looks somewhat different from 
but is equivalent to the one we described above in our situation.
We say such $\{\frak S_p\}$ is a CF-perturbation and use the symbol $\widehat{\frak S}$.
\par
Suppose we have a Kuranishi structure  $\widehat{\mathcal U}$ on $\mathcal M$ 
and a continuous map $f: \mathcal M \to N$ to a manifold $N$.
(We will use it in the case of an evaluation map.)
We say $f$ is strongly smooth\index{strongly smooth} if there exists a
smooth map $f_p : U_p \to N$ for each $p\in \mathcal M$
such that $f \circ \psi_p = f_p$ on $\frak s_p^{-1}(0)$ and
$f_p \circ \varphi_{pq} = f_q$ for a coordinate change map $\varphi_{pq}$.
\par
We say a strongly smooth map $\{f_p\}$ to be 
weakly submersive\index{weakly submersive} if all $f_p$ are submersions.
\par
Suppose we are given a CF-perturbation $\widehat{\frak S} =  \{\frak S_p\}$ ($\frak S_p = [W_p,\{\frak s_{p,\epsilon}\},\chi_p]$)
of $\widehat{\mathcal U}$. 
Then we say $\{f_p\}$ is strongly submersive\index{strongly submersive} with respect to $\widehat{\frak S}$ if 
$f_p \circ \pi\vert_{\frak s_{p,\epsilon}^{-1}(0)} :  \frak s_{p,\epsilon}^{-1}(0) \to N$ is a submersion 
for each $p$.
\par
A differential form on $\widehat{\mathcal U}$ is $\hat h = \{h_p\}$
where $h_p$ is a differential form on $U_p$ such that
$\varphi_{pq}^*h_p = h_q$.
\par
Given a CF-perturbation $\widehat{\frak S}$ on $\widehat{\mathcal U}$,
a map $\hat f = \{f_p\}$ to $N$ which is strongly submersive 
with respect to $\widehat{\frak S}$, and a smooth differential form 
$\hat h = \{h_p\}$ on $\widehat{\mathcal U}$, 
we can define the integration along the fiber\index{integration along the fiber} 
\begin{equation}\label{intfiber}
(\hat f )_! (\hat h;\widehat{\frak S}_{\epsilon})
\end{equation}
which is a smooth differential form on $N$ as follows.
We first define an appropriate notion of partition of unity 
and then use it to reduce to the case of a single Kuranishi chart.
Then we use Formula (\ref{form13-1}).
Actually we can not do it directly with a Kuranishi 
structure but first construct a 
good coordinate system. See \cite[Definition 3.15]{springer}
for the definition of a good coordinate system and
 \cite[Theorem 3.35]{springer} for the existence of
 a good coordinate system compatible with 
 a Kuranishi structure.
 In the case of a manifold, to define the integration of a differential form
 we take a certain locally finite cover by coordinate 
 charts. A good coordinate system consists of finitely many 
 Kuranishi charts (in our situation $\mathcal M$ is compact), 
 so plays the role of such a  locally finite cover.
 \begin{rem} The integration along the fiber
(\ref{intfiber}) is defined for each sufficiently small positive $\epsilon$.
It depends on $\epsilon$ and $\widehat{\frak S}$. However it is 
independent of other choices such as, the choice of partition of unity
and of good coordinate system compatible with $\widehat{\mathcal U}$.
See \cite[Theorem 9.14]{springer}.
\end{rem}

\subsection{Construction of CF-perturbations using outer collaring: review}
\label{CFconst}

In this subsection we review the proof of Proposition \ref{existmkulti1}
except Item (3), the forgetful map compactibility.
The proof was given in \cite[Chapter 22]{springer} using various results 
in the earlier chapters of \cite{springer} especially its Chapter 17.
In this subsection we review it.

Proposition \ref{existmkulti1} is an existence theorem of a system of
CF-perturbations. We can prove it by using two 
basic  existence theorems, which we will explain below in order.
\par
The first existence theorem is an analogy of the standard transversality 
theorem for sections of vector bundles on manifolds.
Let $E \to M$ be a smooth vector bundle on a smooth manifold.
Let $K \subset M$ be a compact subset.
Suppose there exists a section $s_0$ of $E$ on a neighborhood of $K$ 
 which is transversal to $0$.
Then there exists a smooth section $s$ of $E$ on $M$
which is transversal to $0$ and $s=s_0$ on $K$.  
\par
We want to obtain a version of this classical result 
for a CF-perturbation of a Kuranishi structure.
To prove such a variant, we need to employ a good coordinate system 
rather than a Kuranishi structure.
So we first prove the result by employing a good coordinate system
compatible to the given Kuranishi structure and 
obtain a good coordinate system version of CF-perturbation.
We then go back to the Kuranishi structure to obtain a (Kuranishi structure version of)
CF-perturbation.
A technical point to deal with in the process is that the Kuranishi structure 
which we extract from a good coordinate system compatible to the given Kuranishi structure
is not the original Kuranishi structure but is its thickening in general.
We explain the notion of thickening below.
\par
We take a short cut using the fact that our Kuranishi structure on the moduli space 
$\mathcal M$
is obtained from an obstruction bundle data $\{E_{{\bf p}}({\bf x})\}$
as we explained in Subsection \ref{diskreview}.
We denote by $\widehat{\mathcal U}^{E_*(*)}$
the  Kuranishi structure
obtained from an obstruction bundle data $\{E_{{\bf p}}({\bf x})\}$.
We then have the next lemma:
\begin{lem}\label{thickeningdef}
Let $\{E_{{\bf p}}({\bf x})\}$ and $\{E'_{{\bf p}}({\bf x})\}$ be 
two obstruction bundle data
which induce Kuranishi structures 
$\widehat{\mathcal U}^{E_*(*)}$ and
$\widehat{\mathcal U}^{E'_*(*)}$
on our moduli space $\mathcal M$
respectively.
We assume: 
\begin{enumerate}
\item[(*)]
For each ${\bf p}\in \mathcal M$ there exists a neighborhood $\mathfrak U_{\bf p}$
of ${\bf p}$ in $\mathcal M$ such that the following holds.
If ${\bf q} \in \mathfrak U_{\bf p}$ and ${\bf x} \in \mathcal X$ is in the ambient set and is in a 
sufficiently small neighborhood of ${\bf q}$ then
\begin{equation}\label{form13333}
E_{{\bf p}}({\bf x}) \subseteq E'_{{\bf q}}({\bf x}).
\end{equation}
\end{enumerate}
If $(*)$ is satisfied then $\widehat{\mathcal U}^{E'_*(*)}$
is a thickening\index{thickening} of $\widehat{\mathcal U}^{E_*(*)}$
in the sense of \cite[Definition 5.3]{springer}.
\end{lem}

Once \cite[Definition 5.3]{springer} is understood,  
Lemma \ref{thickeningdef} is easy to prove.
However, for the purpose of this paper, 
we can (and will) take a short cut
and regard (*) as the definition of thickening.
\begin{rem}
$  $ \par
\begin{enumerate}
\item
Semi-continuity of obstruction bundle data (Definition \ref{defn51} (4)) 
implies $E_{{\bf p}}({\bf x}) \supseteq E_{{\bf q}}({\bf x})$.
The direction of inclusion is opposite from (\ref{form13333}).
\item
Semi-continuity implies that $\{{\bf p} \in \mathcal M \mid {\rm rank}\,E_{{\bf p}}({\bf x}) \ge d
\}$ is closed. So typically this set is not open. So  (\ref{form13333}) is not satisfied for $E'_*(*) = E_*(*)$ 
typically.
\end{enumerate}
\end{rem}
Now the aforementioned relative existence theorem of 
CF-perturbations reads as follows.
\begin{prop}\label{relativeexisCF}
Suppose $\mathcal M$ has a Kuranishi structure 
$\widehat{\mathcal U}^{E_*(*)}$
and $f : \mathcal M \to N$ is an underlying 
continuous map of a strongly smooth map $\hat f$.
We also assume that $\hat f$ is weakly submersive.
\par
Let $K$ be a compact subset of $\mathcal M$.
We assume that there exists a $CF$-perturbation 
$\widehat{\frak S}$ such that $\hat f$ is 
strongly submersive with respect to $\widehat{\frak S}$
in a neighborhood of $K$.
\par
Then there exists a thickening $\widehat{\mathcal U}^{E'_*(*)}$ of $\widehat{\mathcal U}^{E_*(*)}$
and a $CF$-perturbation 
$\widehat{\frak S}'$ of $\widehat{\mathcal U}^{E'_*(*)}$
with the following properties.
\begin{enumerate}
\item
$\hat f$ is extended to  $\widehat{\mathcal U}^{E'_*(*)}$ and is
strongly submersive with respect to $\widehat{\frak S}'$
everywhere.
\item
On a neighborhood of $K$, CF-perturbations $\widehat{\frak S}$ and $\widehat{\frak S}'$
are compatible in the sense of {\rm (i)}, {\rm (ii)} in $Subsection\,  \ref{CFreview}$.
\end{enumerate}
\end{prop}
Proposition \ref{relativeexisCF} is a direct consequence 
of \cite[Proposition 7.59]{springer} using  \cite[Lemmas 9.10, 9.11]{springer}.\footnote{\cite[Proposition 7.59]{springer} is a similar result for a good coordinate 
system. We use \cite[Lemmas 9.10, 9.11]{springer} to 
relate a CF-perturbation of a Kuranishi structure with one 
of a good coordinate system.}
We omit the proof in this review.
\par\medskip
To prove Proposition \ref{existmkulti1} we need another kind of 
extension result of CF-perturbations, which we now explain.
The construction of a system of CF-perturbations on $\mathcal M_{\ell,k+1}(X,L;\beta)$ is by induction
on the numbers $\omega(\beta)$, $k$, $\ell$.  
The moduli space $\mathcal M_{\ell,k+1}(X,L;\beta)$ has a boundary and corners (in the sense of 
Kuranishi structure).  The boundary is described by (\ref{eq177}).
The corner is described by multiple fiber product. See \cite[(21.1) and Definition 21.3]{springer}.
To prove that the obtained 
CF-perturabtions are compatible along the boundary and corners, 
we need to show the following.
The restriction of the CF-perturbation 
of $\mathcal M_{\ell,k+1}(X,L;\beta)$ to the boundary coincides with the fiber product CF-perturbations 
of the boundary and corners.  See \cite[Section 10.2]{springer} for the definition 
of the fiber product  CF-perturbation. In our situation ${\rm ev}_0$ is strongly 
submersive. Therefore the fiber product  CF-perturbation is well defined
by \cite[Lemma-Definition 10.12]{springer}.
To prove the existence of such CF-perturbations by induction, 
we need to extend a CF-perturbation defined at the boundary 
of $\mathcal M_{\ell,k+1}(X,L;\beta)$ to its neighborhood.
In fact the assumption of Proposition \ref{relativeexisCF} 
requires that the CF-perturbation is given on {\it a neighborhood of} $K$.
The proof of this extendability is somewhat cumber some.
\par
\begin{figure}[h]
\centering
\includegraphics[scale=0.3]{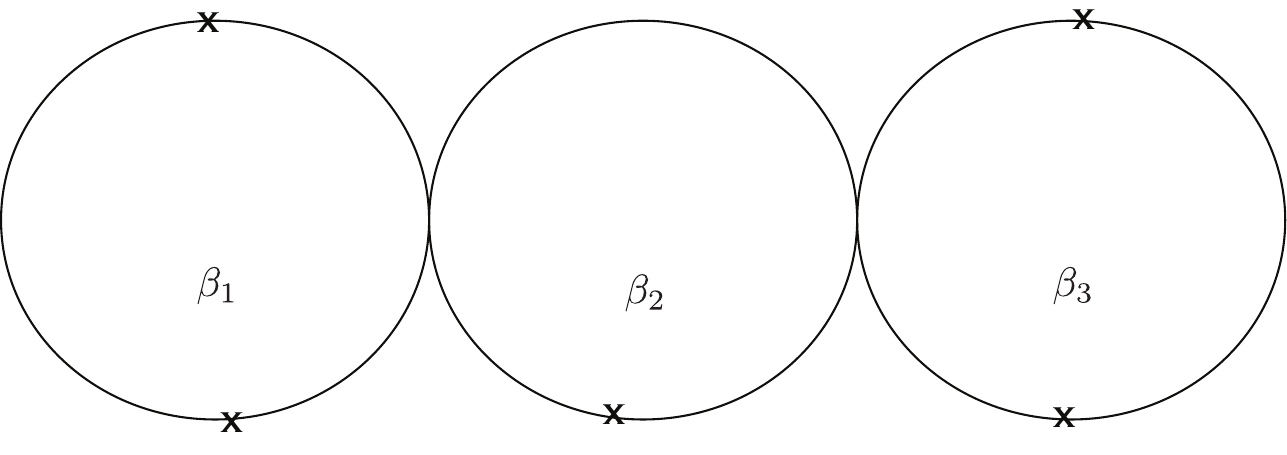}
\caption{$p \in \mathcal M_5(X,L;\beta)$}
\label{zu9}
\end{figure}

Let us elaborate on  the issue that arises while we give the proof of this extendability.
We remark that the CF-perturbation is an {\it equivalence class 
of} $(W,\{\frak s_{\epsilon}\},\chi)$.
Let us consider the case when $p \in \mathcal M_{5}(X,L;\beta)$ be a point of 
codimension 2 corner.
It is in the intersection of two fiber products
\begin{equation}\label{formula133}
\aligned
\mathcal M_{3}(X,L;\beta_1) &\times_{L} \mathcal M_{4}(X,L;\beta_2+\beta_3), \\
\mathcal M_{4}(X,L;\beta_1+\beta_2) &\times_{L} \mathcal M_{3}(X,L;\beta_3).
\endaligned
\end{equation}
Here $\beta = \beta_1 + \beta_2 +\beta_3$.
See Figure \ref{zu9}.
By assumption the CF-perturbation $(W_p,\{ \frak s_{p,\epsilon}\},\chi_p)$
extends to the CF-perturbations of the two fiber products 
in (\ref{formula133}).
We can use this fact to extend $(W_p,\{\frak s_{p,\epsilon}\},\chi_p)$ 
to a neighborhood of $p$ in $\mathcal M_{5}(X,L;\beta)$.
(We use the fact that a smooth function defined on the boundary of 
a manifold with corners can be extended smoothly to its 
neighborhood.  (See \cite[Lemma 7.2.121]{fooo092}.))
Then we can use a partition of unity to globally extend the 
CF-perturbation on the boundary to its neighborhood.
\par
However we need ensure the following in order for this argument to work:
\begin{itemize}
\item
We can choose a representative $(W_p,\{{\frak s}_{p,\epsilon}\},\chi_p)$
that extends to a neighborhood of $p$ each of the two fiber products 
in (\ref{formula133}). 
\end{itemize}
It is important to take a {\it single} representative which 
 extends to {\it both} fiber products simultaneously.
Since the equivalence relation between CF-perturbations 
is not so simple, finding such a representative 
is a nontrivial issue.\footnote{In the case of multisection, such a 
choice is easier.}
\par
A way to go around this trouble is to utilize the collar of the 
boundary and corners.  
Let us consider again the element $p$ above.
It is contained in the codimension 2 corner which is 
identified with the fiber product:
\begin{equation}\label{form134}
\mathcal M_{3}(X,L;\beta_1) \times_{L} \mathcal M_{3}(X,L;\beta_2) 
\times_{L} \mathcal M_{3}(X,L;\beta_3).
\end{equation}
A Kuranishi neighborhood of $p$ in 
$\mathcal M_{5}(X,L;\beta)$ is written as 
$U_p \times [0,1) \times [0,1)$
where $U_p$ is a Kuranishi neighborhood of $p$ 
in (\ref{form134}).
We may take an obstruction bundle on $U_p \times [0,1) \times [0,1)$
by taking a pull back of those on $U_p$.
\par
The Kuranishi neighborhood of $p$ in the two fiber products 
in (\ref{formula133}) can be identified with 
$U_p \times \{0\} \times [0,1)$ and
$U_p \times [0,1) \times \{0\}$, respectively.
Suppose the CF perturbation on them 
are both collared. 
It means that there exists $\tau > 0$ such that 
$\frak s_{p,\epsilon}((x,t,0),\xi) = \frak s_{p,\epsilon}((x,0,0),\xi)$ 
and $\frak s_{p,\epsilon}((x,0,t),\xi) = \frak s_{p,\epsilon}((x,0,0),\xi)$
for $t < \tau$.  Here $x \in U_p$, $\xi \in W$.
We can then extend it to a neighborhood of 
$U_p \times \{(0,0)\}$ in 
$U_p \times [0,1) \times [0,1)$ 
by defining
$$
\frak s_{p,\epsilon}((x,t_1,t_2),\xi) = \frak s_{p,\epsilon}((x,0,0),\xi)
$$
for $t_1,t_2 \in [0,\tau)$. 
This way of extension is simpler and works better with the 
equivalence relation of CF-perturbations.
\par
The problem is the way to find a global collar to a Kuranishi structure.
The existence of a collar of a given manifold with corners
is classical.\footnote{It is not easy to find a reference however.}
In the case of Kuranishi structure an issue occurs since 
collar of a manifold with corners is  not enough canonical.
(See \cite[Section 18.4]{springer} for a discussion on a related problem.)
In the case of a Kuranishi structure, we need  to equip each
Kuranishi chart with a collar simultaneously in the collection of Kuranishi charts associated to
the given Kuranishi structure so that 
the coordinate changes respect the collar structure.
This is cumbersome work to carry out.
\par
The way taken in \cite[Chapter 17]{springer} is to use outer collaring
which we explained in Subsection \ref{diskforget}.
In the case of the previous example we consider 
$(U_p\times [0,1) \times [0,1))^{\boxplus \tau}: = U_p \times [-\tau,1) \times [-\tau,1)$,
where $\tau > 0$.

\begin{defn}
Let  $\widehat{\mathcal U}$ be a $\tau$-collared Kuranishi structure.
Its CF-perturbation $\widehat{\frak S}$ is said to be $\tau$-collared\index{collared CF-perturbation}
if for each  $\widehat{\bf p}_1, \widehat{\bf p}_2 \in \mathcal M_{\ell,k+1}(X,L;\beta)^{\boxplus 1}$ 
with $\Pi_{\tau}(\widehat{\bf p}_1) = \Pi_{\tau}(\widehat{\bf p}_2)$ we can take a representative 
$(W_{\hat{\bf p}_i},\{\frak s_{\hat{\bf p}_i,\epsilon}\},\chi_{\hat{\bf p}_i})$ ($i=1,2$)
 %\marginpar{Definition corrected.
%KF 2025 July}
such that $W_{\hat{\bf p}_1} = W_{\hat{\bf p}_2}$ and that the following holds.
\par
Suppose $\widehat{\bf x}_1, \widehat{\bf x}_2$ are in a sufficiently small neighborhood 
$\mathscr U_{{\bf p}_i}$ of  $\hat{\bf p}_i$ for $i=1,2$.
in the ambient set.  We assume $U_{{\bf p}_i}$ is a Kuranishi neighborhood of ${\bf p}_i$
contained in $\mathscr U_{{\bf p}_i}$.
If  $\Pi_{\tau}(\widehat{\bf x}_1) = \Pi_{\tau}(\widehat{\bf x}_2)$ then
$$
\frak s_{{\bf p}_1,\epsilon}(\widehat{\bf x}_1,\xi) = \frak s_{{\bf p}_2,\epsilon}(\widehat{\bf x}_2,\xi).
$$
\end{defn}
We next recall the notion of normalized boundary.
For a manifold with corners $M$ its normalized boundary is $\partial M$ 
which is also 
a manifold with corners.  
There is a map $\pi : \partial M \to M$. If $p$ is interior of $M$ then 
$\pi^{-1}(p)$ is the empty set. 
If $p$ is in the boundary of $M$ but is not in the corner, then 
$\pi^{-1}(p)$ consists of one point. 
If $p$ is in the codimension $k$ corner but not in 
in the codimension $k+1$ corner, 
then $\pi^{-1}(p)$ consists of $k$ points. 
See  \cite[Definition 8.4 and Lemma 8.2]{springer}.
By taking normalized boundary of each Kuranishi chart, 
a space with Kuranishi structure (with a boundary and corners) 
has a normalized boundary which is again a 
Kuranishi structure with a boundary and corners.
\par
In the case of 
$\partial{\mathcal M}_{\ell;k+1}(L;\beta)$ its normalized boundary is a 
{\it disjoint union}
\begin{equation}\label{form135}
\bigcup {\mathcal M}_{\# {\mathbb L}_1;k_1+1}(L;\beta_1)
{}_{\text{\rm ev}_{0}}\times_{\text{\rm ev}_i} {\mathcal M}_{\# {\mathbb L}_2;k_2+1}(L;\beta_2).
\end{equation}
The disjoint union is taken over all $({\mathbb L}_1,{\mathbb L}_2) \in \text{\rm Shuff}(\ell)$, $i = 1, \dots, k_2$,
$k_1,k_2\in \Z_{\ge 0}$ with $k_1 + k_2 = k$ and $\beta_1,\beta_2
\in H_2(X,L;\Z)$ with $\beta_1 +\beta_2 = \beta$. (Compare (\ref{eq177}).)
\par
Note that the boundary of (\ref{form135}) is 
\begin{equation}\label{form136}
\bigcup {\mathcal M}_{\# {\mathbb L}_1;k_1+1}(L;\beta_1)
{}_{\text{\rm ev}_{0}}\times_{\text{\rm ev}_i} {\mathcal M}_{\# {\mathbb L}_2;k_2+1}(L;\beta_2)
{}_{\text{\rm ev}_{0}}\times_{\text{\rm ev}_j} {\mathcal M}_{\# {\mathbb L}_3;k_3+1}(L;\beta_3).
\end{equation}
More precisely each element of (\ref{form136}) appears 
twice in the boundary of (\ref{form135}).
We say a Kuranishi structure, a CF-perturbation etc. 
of (\ref{form135}) is compatible at the corner 
if they coincide at the two points which corresponds to a point in (\ref{form136}).

\begin{rem}\label{remark1313}
Here we say {\it coincide}.  In the situation of this paper we can safely say so.
In fact a Kuranishi neighborhood is, for example, a subset of the 
ambient set.  So whether it coincides or not is a well defined notion.
The same is true for obstruction bundles, CF-perturbations and etc.
In \cite{springer} the boundary compatibility 
(See \cite[Condition 16.17 (X)]{springer}.) requires that two Kuranishi structures are {\it isomorphic}, 
instead of requiring them to coincide.  In the story of a system of abstract Kuranishi structures 
which is developed in  \cite{springer}, we can not define the notion that two Kuranishi structures  coincide.
By this reason, we need to require further consistency of isomorphisms. 
(See \cite[Remark 15.6]{springer}.)
Such a further consistency condition is the corner compatibility condition. (See \cite[Condition 16.17 (XI), (XII)]{springer}.)
\par
Thus the story here is simplified compared to \cite[Chapter 17]{springer}.
\end{rem}
Now we can state the second extension theorem as follows.
Let $\rho_0 > 0$ be the smallest symplectic  energy of non-constant pseudo-holomorphic sphere or  
disk which bounds $L$.

\begin{prop}\label{prop13.12}
We fix $(\beta,k,\ell)$.
We assume that for each $(\beta',k',\ell')$ with 
$(\beta',k',\ell') \le (\beta,k,\ell)$, 
we have a CF-perturbation $\widehat{\mathfrak S}$
on ${\mathcal M}_{\ell';k'+1}(L;\beta')^{\boxplus 1}$
with the following properties.
\begin{enumerate}
\item
Proposition $\ref{existmkulti1}$ $(1),(2),(4),(5),(6)$ hold.
\item
It is $\tau$-collard.
\end{enumerate}
Then there exists a CF-perturbation on the $\tau$-collar of the boundary of $\partial{\mathcal M}_{\ell;k+1}(L;\beta)^{\boxplus 1}$
that is $\tau$-collard 
and satisfies Proposition $\ref{existmkulti1}$ $(1),(2),(4),(5),(6)$.
\end{prop}
\begin{proof}
This follows from \cite[Proposition 17.73]{springer}.
As we mentioned in Remark \ref{remark1313}
the proof in our situation is easier and proceed as follows.
\par
Note that by assumption on the moduli space appearing 
in the boundary of ${\mathcal M}_{\ell;k+1}(L;\beta)^{\boxplus 1}$
CF-perturbation is given.
Therefore each of summand of the boundary has a CF-perturbation, 
that is, the fiber product CF-perturbation.
By induction hypothesis Proposition \ref{existmkulti1} (2)
and associativity of fiber product, 
at the codimension 2 boundary (which is an outer-collared-version of (\ref{form136})), the restrictions of the 
two summands (the boundary components) coincide.
\par
Therefore we can pull back CF-perturbation 
by $\Pi_{\tau}$ in an obvious way to obtain 
the CF-perturbation with claimed properties.
\end{proof}
\begin{proof}[Proof of Proposition  \ref{existmkulti1} except (3)]
Using Propositions \ref{relativeexisCF} and \ref{prop13.12} the proof is easy.
The proof is by induction on the partial order  as in Definition \ref{defn837837}.
\par
The first step of induction (that is, the case of the smallest $(\beta,k,\ell)$)  is a case when
 ${\mathcal M}_{\ell;k+1}(L;\beta)^{\boxplus 1} =  {\mathcal M}_{\ell;k+1}(L;\beta)$ has no boundary.
 So we can apply Proposition \ref{relativeexisCF} with $K = \emptyset$.
 \par
 For the inductive step, we apply Proposition \ref{prop13.12} to obain a CF-perturbation 
 on the $\tau$-collar of the boundary of $\partial{\mathcal M}_{\ell;k+1}(L;\beta)^{\boxplus 1}$.
 We then apply Proposition \ref{relativeexisCF}  taking $K$ to be the $\tau'$-collar where $\tau'$ is slightly smaller 
than $\tau$ to obtain the required CF-perturbation.
\par
Note that Proposition \ref{relativeexisCF} does not claim explicitly that if $\widehat{\frak S}$ is
$\tau'$-collared then $\widehat{\frak S}'$ is
$\tau'$-collared as well. However this is obvious from the proof.
\end{proof}

\subsection{Smooth-ness of the pull back of CF-perturbation by 
the forgetful map.}
\label{CFforget}

In this subsection, we prove Proposition  \ref{existmkulti1}  (3).
For this purpose, we need to pull back 
a CF-perturbation by a forgetful map.
The main purpose of this subsection is to explain how to pull back a CF-perturbation 
defined on the outer collar.
\par
Before going there we review an issue of smoothness of the forgetful map 
which was studied in \cite[pages 777-778]{fooo092}.
We consider the case of Figure \ref{Figure13-14} below.
\begin{figure}[h]
\includegraphics[scale=0.4]{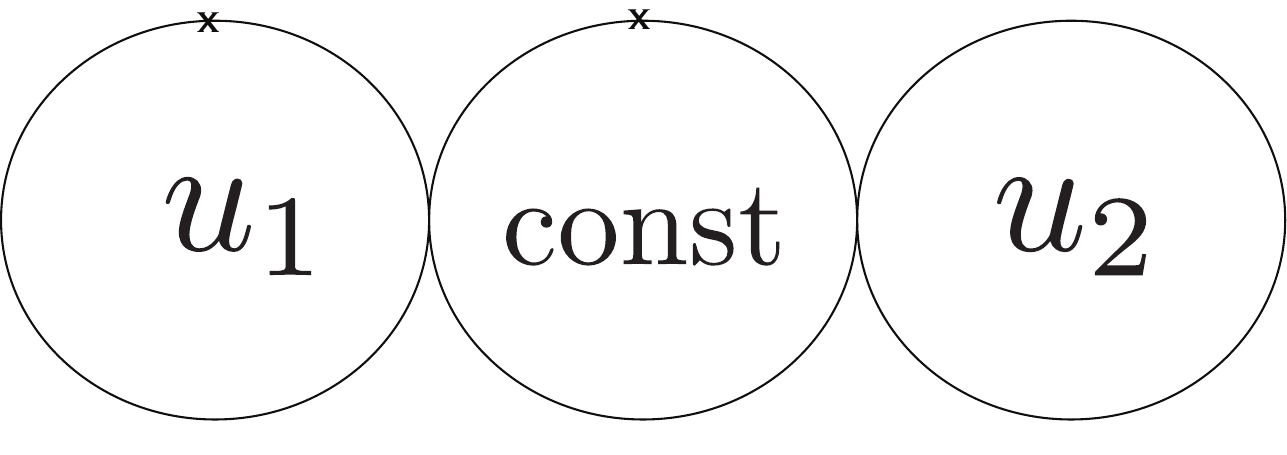}
\caption{An element of a codimension two stratum}
\label{Figure13-14}
\end{figure}
It is an element $\bf p$ of $\mathcal M_{0,2}(L,\beta)$.
If we forget the two marked points then 
the component in the middle (on which the map
is a constant map) becomes unstable and so we shrink it.
There are two boundary nodes in $\bf p$
and so ${\bf p}$ lies in the codimension 2 corner.
We assume $u_1$, $u_2$ are smooth points of and are isolated
in the moduli spaces, 
for simplicity. We also assume that they are non-constant.
The parameter to deform two boundary nodes are 
small real numbers $r_1$, $r_2$.
We fix coordinates of the two irreducible 
components containing the nodes.
We denote them by $z_1$, $w_1$ and 
$z_2$, $w_2$.
To deform the node we put the relation 
$z_1w_1 = r_1$, $z_2w_2 = r_2$
and obtain a smooth domain for $r_1,r_2 > 0$.
However these gluing parameters $r_1$, $r_2$ are not the
one in our Kuranishi structure.
We consider the map
$$
u_{r_1,r_2} : D^2 \to X
$$
which is obtained by gluing. (By the assumption such 
pseudo-holomorphic map with boundary condition $L$ 
exists uniquely.
We fix domains $D_i \subset D^2$ which can be 
identified to a part where $u_i$ was defined.
Suppose $D_i$ is away from the node.
We restrict $u_{r_1,r_2}$ to $D_i$.
Then the following is proved in \cite[Theorem 6.4]{foooexp},
where $T_i = - \log r_i$
\begin{equation}\label{form138}
\left\vert \frac{\partial}{\partial T_1^{k_1}}\frac{\partial}{\partial T_2^{k_2}} u_{r_1,r_2} - u_i \right\vert_{C^k(D_i)} < C_k \exp (-\delta_k
(k_1 T_1 + k_2 T_2 + \inf \{T_1,T_2\})).
\end{equation}
This estimate does {\it not} imply that 
the family $u_{r_1,r_2}$ (restricted to $D_i$) is smooth 
with respect to the coordinate $r_i$ at $(r_1,r_2) = (0,0)$.
This could cause a problem to show, for example, smoothness of the
coordinate change map or the Kuranishi map.
\par
A way to resolve this issue is to change the coordinate to
$t_i =1/ \log T_i.$
Then (\ref{form138}) implies that 
the restriction of $u_{r_1,r_2}$  to $D_i$ is a smooth family 
with respect to the $t_i$ coordinate at $t_i = 0$.
\par
We next consider $\overline{\bf p} = \frak{forget}({\bf p})$.
Its domain has one double point.
So in a similar way we have one gluing parameter $r$.
The coordinate we use then is
$t = 1/\log T$,  $T = -\log r$.
\par
It is easy to see that the forgetful map in the coordinates $r_i,r$ is written as 
$(r_1,r_2) \mapsto r_1r_2 = r$.
Therefore by an easy calculation 
it becomes
$$
(t_1,t_2) \mapsto \frac{1}{\log(e^{1/t_1}+e^{1/t_2})}.
$$
This map is {\it not} smooth at $(0,0)$.
So it could cause a problem 
to pull back a CF-perturbation etc. 
by the forgetful map.
\par
The way to resolve this problem is as follows.
(See \cite[pages 777-778]{fooo092}.)
For a family of sections ${\frak s}_{\epsilon}(x,\xi)$
which is a part of a CF-perturbation, 
we consider only those with the properties
$$
\vert \frak s_{\epsilon}((T_1,T_2,{\overline x}),\xi) - \frak s_{\epsilon}((\infty,\infty,{\overline x}),\xi)\vert_{C^k} 
< C_k \exp (-\delta_k T_i)
$$
Here we take derivatives with respect to $T_i$ to define the $C^k$ norm.
(The parameter $x$ of the Kuranishi neighborhood contains those 
corresponding to the gluing parameters. So we write $x = (T_1,T_2,\overline x)$.)
Such a function is called {\it admissible}\index{admissible} in \cite[Chapter 25]{springer}.
(See \cite[Definition 25.3]{springer} for its precise definition.)
We can define orbifold, vector bundle, sections and etc. in the category 
of admissible functions and can prove various results on usual smooth functions 
such as implicit function theorem etc. in the same way in the admissible 
category.  So the construction of a Kuranishi structure, 
CF-perturbations etc. can be carried out in the same way in the admissible category.
(See \cite[Chapter 25]{springer}.)
If a CF-perturbation etc. is admissible 
then it is smooth with respect to the coordinates $t_i$ at $(t_1,t_2) = (0,0)$.
Moreover since $T = T_1 + T_2$ is the forgetful map,
the pull back of an admissible CF-perturbation is 
again admissible.
\par
Another reason to work in the category of admissible 
functions concerns the outer collaring.
At the boundary or corners, which is given by $T_i = \infty$,
that is, $t_i = 0$, we attach outer collar 
that is $t_i \in [-1,0]$. The coordinate change maps, 
CF-perturbations and etc. are extended to the outer collar so 
that they are constant in $t_i$ direction. (for $t_i \in [0,\epsilon)$.)
For this extension to become smooth at $t_i = 0$, 
the coordinate change maps, 
CF-perturbations and etc. must have vanishing $t_i$ derivatives 
of all order at $t_i = 0$.
The coordinate $t_i$ does not have this 
property.
The coordinate $1/T_i$ has this properties 
for CF-perturbations because of the estimates of the form 
(\ref{form138}).
For coordinate change maps, we need to take the logarithum 
once more and use $t_i = 1/\log T_i$. 
Then using $t_i$ cordinate, the coordinate change maps, 
CF-perturbations etc. have vanishing $t_i$ derivatives 
of all order at $t_i = 0$.
Therefore using this coordinate 
outer collar has a smooth structure. 
See \cite[Chapters 25 and 17]{springer}, \cite{Corrigendum} for detail.
\par\medskip
Now we discuss the way to pull back the CF-perturbation defined 
on the outer collar by the forgetful map $\frak{forget}^{\boxplus 1}$.
\par
We  introduce a class of CF-perturbations 
which can be pulled back by the map $\frak{forget}^{\boxplus 1}$.
It is called a bi-collared\index{bi-collared} CF-perturbation.
\begin{defn}
Let $\{E_{\hat{\bf p}}(\widehat{\bf x})\}$ be a $\tau$-bi-collared obstruction bundle data on
the moduli space
${\mathcal M}_{\ell;k+1}(L;\beta)^{\boxplus 1}$.
It induces a $\tau$-bi-collared Kuranishi structure $\widehat{\mathcal U}^{E_*(*)}$.
A CF-perturbation 
$\widehat{\frak S}$ of $\widehat{\mathcal U}^{E_*(*)}$ is said to be $\tau$-bi-collared\index{bi-collared CF-perturbation}
if for each  $\widehat{\bf p} \in \mathcal M_{\ell,k+1}(X,L;\beta)^{\boxplus 1}$ we can take a representative 
$(W_{\hat{\bf p}},\{\frak s_{\hat{\bf p},\epsilon}\},\chi_{\hat{\bf p}})$
such that the following holds.
\par
Suppose $\widehat{\bf x}_1, \widehat{\bf x}_2$ are in a neighborhood of  $\hat{\bf p}$
in the ambient set.  We assume they are in $U_{\bf p}$ the Kuranishi neighborhood of ${\bf p}$.
If  $\Pi'_{\tau}(\widehat{\bf x}_1) = \Pi'_{\tau}(\widehat{\bf x}_2)$ then
$$
\frak s_{{\bf p},\epsilon}(\widehat{\bf x}_1,\xi) = \frak s_{{\bf p},\epsilon}(\widehat{\bf x}_2,\xi).
$$
\end{defn}

\begin{prop}\label{prop1315}
Suppose a system of $\tau$-bi-collared obstruction bundle data
on ${\mathcal M}_{\ell;k+1}(L;\beta)^{\boxplus 1}$ 
and on ${\mathcal M}_{\ell;0}(L;\beta)^{\boxplus 1}$
are given such that they are compatible with the forgetful map.
\par
In this situation a $\tau$-bi-collared admissible CF-perturbation 
on ${\mathcal M}_{\ell;0}(L;\beta)^{\boxplus 1}$ can 
be pulled back by the forgetful map 
to a (smooth and) $\tau$-bi-collared admissible  CF-perturbation
on ${\mathcal M}_{\ell;k+1}(L;\beta)^{\boxplus 1}$.
\end{prop}
We say that two CF-perturbations are compatible if 
one is obtained from the other via the pull-back.

The notion of pull-back of a CF-perturbation by a forgetful map is the obvious one.
Using the notation of Definition \ref{def1221}, 
the pull back of a CF-perturbation $(W,\{\frak s_{\epsilon}\},\chi)$ on 
 $(U^2_{{\bf p}_2},\mathcal E^2_{{\bf p}_2},\psi^2_{{\bf p}_2},s^2_{{\bf p}_2})$
to  $(U^1_{{\bf p}_1},\mathcal E^1_{{\bf p}_1},\psi^1_{{\bf p}_1},s^1_{{\bf p}_1})$
is $(W,\{\tilde F_{{\bf p}_1}^*\frak s_{\epsilon}\},\chi)$.
\begin{proof}
The only issue to take care of is the non-smoothness 
of $\mathscr M_m$.
\par
The  singular locus of  $\mathscr M_m$ 
is described in Lemma \ref{lem1314} (6), (a) and (b).
We observe that in such a case the $\tau$-bi-collared-ness
implies that the CF-perturbation is constant 
in the $s_i$ direction.  
(We use $\epsilon <\tau$ here.)
Therefore, the non-smooth-ness of  $\mathscr M_m$
does not affect the smooth-ness of the pull back.
\end{proof}

Now we are in the position to complete the proof of an outer-collared-version of
Proposition  \ref{existmkulti1}  (3).
We observe that the inductive construction 
of the system of CF-perturbations, which we described in 
Subsection \ref{CFconst}, can be performed 
in the category of bi-collared CF-perturbation.
Therefore 
using Proposition \ref{prop1315} 
we can perform the inductive construction 
of CF-perturbations compatible with forgetful map
by the induction scheme 
of the construction of Kuranishi structures in Subsection \ref{diskforget},
to obtain a system of CF-perturbations
as in Proposition  \ref{existmkulti1} including Item (3).
(We use the inductive scheme of the proof of Lemma \ref{lem12220}
so that the CF-perturbations become bi-collared.)

The proofs of Propositions \ref{prop95}, \ref{exstCFUNI} etc. 
are similar.  To prove Proposition \ref{exstCFUNI} we use 
Proposition \ref{relativeexisCF}.
\qed
\newpage

\part{Floer theory of closed and open strings.}
\label{part3}

In Section \ref{sec:frakq} we define the {\it closed-open map}\index{closed-open map} (\ref{mapqhat}).
In Section \ref{sec:open-closed map}, we define the {\it open-closed map}\index{open-closed map}
\begin{equation}\label{mapphat}
\widehat{\frak p}^{\text{\bf b}} : HH_*(\cL)
\to 
H(X;\Lambda_0).
\end{equation}
We prove that $\widehat{\frak p}^{\text{\bf b}}$ is dual to $\widehat{\frak q}^{\text{\bf b}}$ 
in Section \ref{dualpq}.
We also prove the non-triviality of $\widehat{\frak p}$ under the assumption 
of nontriviality of Floer cohomology in Section \ref{nontrivialsec}.
We will show that $\widehat{\frak q}^{\text{\bf b}}$ is a ring homomorphism in Section \ref{sec:ring}.
In Sections \ref{sec:annuli}, \ref{sec:ori} we will prove the key property of $\widehat{\frak p}^{\text{\bf b}}$ which we use 
to prove our main theorem.

\section{The closed-open map $\widehat{\frak q}$.}
\label{sec:frakq}
%input{section3}

\subsection{Construction of  $\widehat{\frak q}$.}
\label{frakqondeRham}

In this subsection we construct the closed-open map $\widehat{\frak q}^{\rm f.u.{\bf b}}$ in the de Rham setting. Recall that we are interested in the Hochschild invariants of the category $\cL ^{\rm form}_{\rm c.u.}$ from Section \ref{sec:unit} (right after Proposition \ref{prop910}), 
with objects $\mathscr L$ consisting of a collection $\{(L_{\kappa},\theta_{\kappa})\mid \kappa =0,1,\dots,K\}$  of Lagrangians with flat connections.
We then construct $\widehat{\frak q}^{\rm f.c.u.\frak b}$.
Next we eliminate the curvature by bounding cochains and obtain $\widehat{\frak q}^{\rm f.u.{\bf b}}$.

Note that the morphism module $\cL ^{\rm form}_{\rm uni}((L_{\kappa},\theta_{\kappa}),(L_{\kappa},\theta_{\kappa}))$
between the same element $c = (L_{\kappa},\theta_{\kappa})$ contains ${\bf e}^+_c$ and ${\bf f}_c$.

The map $\widehat{\frak q}^{\rm f.c.u.\frak b}$ will be defined using the operations $\frak q^{{\rm f.c.u.}\frak b}_{1,\vec{\kappa}}$ from Definitions \ref{defn626} and \ref{def96}. $\frak q^{{\rm f.c.u.}\frak b}_{\ell,\vec{\kappa}}$ is defined as a sum (\ref{qcat}) using 
$\frak q^{\text{\rm form}}_{\ell,\vec k;\vec{\kappa};B}$ (and including unit as in Section \ref{sec:unit}).  
Since the operators $\frak q^{\text{\rm form}}_{\ell,\vec k;\vec{\kappa};B}$ are defined inductively on $(B,\vert\vec{\kappa}\vert,\ell)$, 
the construction of $\widehat{\frak q}^{\rm f.c.u.\frak b}$ is also by induction. 
For the inductive construction we need to work out the following point as usual.
For given $(E,k,\ell_0)$ the operations $\frak q^{\text{\rm form}}_{\ell,\vec k;\vec{\kappa};B}$ 
are defined for $B,\vec{\kappa},\ell$ with $(B,\vert\vec{\kappa}\vert,\ell) \le (E,k,\ell_0)$ only.
We  must first relate the choices of various data defining them to work out the induction. 
\begin{rem}
Note that to define the structure operation $\frak m^{{\rm f.c.u.}\frak b}_{\vec{\kappa}}
= \frak q^{{\rm f.c.u.}\frak b}_{0,\vec{\kappa}}$ we used 
$\frak q^{\text{\rm form}}_{\ell;\vec{\kappa};B}$ for various $\ell$.
We plug in $\frak b$ to the variables of $\frak q^{\text{\rm form}}_{\ell;\vec{\kappa};B}$
corresponding to the interior marked points.
To define $\widehat{\frak q}^{\rm f.c.u.\frak b}$ we also use $\frak q^{\text{\rm form}}_{\ell;\vec{\kappa};B}$.
This time we plug in $\frak b$ to the variables corresponding to the interior marked points except for 
one of them, to which we plug in an element of $\Omega(X)$ that is the domain of 
$\widehat{\frak q}^{\rm f.c.u.\frak b}$.
\end{rem}
The result we obtain in this section can be summarized as follows.

\begin{prop}\label{prop142}
There exists a series of maps %\marginpar{Put the name.  KF.2025 Aug.}
\index[syindex]{qformelqlvechat@$\widehat{\frak q}^{\rm f.c.u.\frak b}_{\ell}$}
\begin{equation}\label{formLinfmap}
\widehat{\frak q}^{\rm f.c.u.\frak b}_{\ell}: (\Omega(X) \hat\otimes \Lambda_0)^{\otimes \ell}  \to CH^{*}_{\rm red}(\cL ^{\rm form}_{\rm uni},\cL ^{\rm form}_{\rm uni})
\end{equation}
which induce a filtered $L_{\infty}$ homomorphism.
\end{prop}
To obtain a filtered $L_{\infty}$ homomorphism from the maps $\widehat{\frak q}^{\rm f.c.u.\frak b}_{\ell}$  we need to reverse the order of input and  put appropriate sign.
(See (\ref{formnew88}) and \cite[Proposition 7.3.34]{fooo092}.)
The $L_{\infty}$ structure on $\Omega(X) \hat\otimes \Lambda_0$ is trivial except the boundary operator (which is de Rham differential).
The $L_{\infty}$ structure on $CH^{*}_{\rm red}(\cL ^{\rm form}_{\rm uni},\cL ^{\rm form}_{\rm uni})$ is a differential graded Lie algebra 
structure, where differential is the Hochschild differential and bracket is the Gerstenhaber bracket. 

To prove Proposition \ref{prop142}, we first introduce `$A_{n,k}$ version' of the Hochschild complex.
We fix a discrete monoid $G \subset \R_{\ge 0}$ and put $G = \{E_i \mid i=0,1,\dots,\}$ where $E_0 = 0$ and 
$E_{i} < E_{i+1}$.  We assume all the structure operations 
$\frak m^{{\rm f.c.u.}\frak b}_{\vec{\kappa}} 
= \frak q^{{\rm f.c.u.}\frak b}_{0,\vec{\kappa}}$ of our filtered $A_{\infty}$ category 
is $G$-gapped.  (Note that the bulk deformation is included. The bounding cochain is not yet
included at this stage.)
We consider elements
\begin{equation}\label{form142142}
\frak T_{E_i,\vec\kappa} \in
\Hom(BCF(\cL^{\rm form};\vec{\kappa})^+,{CF}(L_{\kappa_{0}},L_{\kappa_{K}};\F)^+)
\end{equation}
for $E_i,\vec{\kappa}$ with $(E_i,\vert\vec{\kappa}\vert) \le (E_n,k)$.\footnote{
Here and hereafter we write $(E_i,\vert\vec{\kappa}\vert) \le (E_n,k)$ in place of $(E_i,\vert\vec{\kappa}\vert,0) \le (E_n,k,0)$.}
\begin{defn}\label{defn142142}
We denote by
$CH^*_{E,k}(\cL^{\rm form}_{\rm uni})$ the set of all
$\frak T = (\frak T_{E_i,\vec{\kappa}})_{(E_i,\vert\vec{\kappa}\vert) \le (E,k)}$
where $\frak T_{E_i,\vec\kappa}$ is as in (\ref{form142142}).\index[syindex]{CHEkLformcu@$CH^*_{E,k}(\cL^{\rm form}_{\rm c.u.})$}
\end{defn}
It is obvious that $CH^*_{E,k}(\cL^{\rm form}_{\rm uni})$ is a graded $\F$ module.
Moreover it is filtered.
\begin{defn}
Let $G \subset \R_{\ge 0}$ be a discrete monoid.
We denote by $\Lambda_G$ the subring of $\Lambda_0$ consisting 
of elements $\sum a_i T^{\lambda_i}\in \Lambda_0$ such that
$\lambda_i \in G$ for all $i$.
\end{defn}
Then $CH^*_{E,k}(\cL^{\rm form}_{\rm uni})$  is a $\Lambda_G$-module.
In fact if $E' \in G$ then 
$$
T^{E'}(\frak T_{E_i,\vec{\kappa}}) = (\frak T'_{E_i,\vec\kappa})
$$
where 
$$
\frak T'_{E_i+E',\vec{\kappa}} = \frak T_{E_i,\vec\kappa}
$$
when $(E_i+E',\vert\vec{\kappa}\vert) \le (E,k)$.
We remark that $\frak T'_{E'',\vec{\kappa}} = 0$ unless $E''$ is of the form $E_i+E'$.
\par
We next define Hochschild differential on $CH^*_{E,k}(\cL^{\rm form}_{\rm uni})$.
We put
\begin{equation}\label{mmcat}
\frak m^{{\rm f.c.u.}\frak b}({\bf h})
=
\sum_{B,\ell,\vec k}T^{\omega(B) }
\frac{\rho_{\frak b,\theta}(B)% \rho_{b}(B)
}{\ell!}
\frak q^{\text{\rm form}}_{\ell;\vec{\kappa};B}
(\frak b_+^{\ell};
\text{\bf h}). %;b_+^{\vec k})
\end{equation}
See Definition \ref{def827}.
For ${\bf h} = h_1\otimes \dots \otimes h_k$ which does not contain $T$, 
we rewrite it as\index[syindex]{mxfcubfrakk@$\frak m^{\text{\rm f.c.u.}\frak b}_{k}$} 
$$
\frak m_k^{{\rm f.c.u.}\frak b}(\text{\bf h}) 
= 
\sum_{E \in G} T^E
\frak m^{{\rm f.c.u.}\frak b}_{k,E}(\text{\bf h})
$$
to define $\frak m^{{\rm f.c.u.}\frak b}_{k,E}$.\index[syindex]{mxfcubfrakkE@$\frak m^{\text{\rm f.c.u.}\frak b}_{k,E}$} 
Here $\frak m^{{\rm f.c.u.}\frak b}_{k,E}$ is $\F$ linear.
Now in the same way as (\ref{eq:diff_CC^*}), we define
\begin{equation}\label{form144144}
\aligned
&(\delta^H \frak T)_{E,\vec{\kappa}}(\bf a) \\
=&\sum_{E_1 + E_2 = E} (-1)^{\deg(\frak T)\deg'( \bfa^{3;1}_{\alpha(3)}) }   \frak m^{{\rm f.c.u.}\frak b}_{E_1}(\bfa^{3;1}_{\alpha(3)} \otimes \frak T_{E_2}( \bfa^{3;2}_{\alpha(3)}) \otimes \bfa^{3;3}_{\alpha(3)} )  \\
&+ \sum_{E_1+E_2=E}  (-1)^{\deg'( \bfa^{3;1}_{\alpha(3)})  +\deg(\frak T)+1}   \frak T_{E_1}(\bfa^{3;1}_{\alpha(3)} \otimes   
\frak m^{{\rm f.c.u.}\frak b}_{E_2}( \bfa^{3;2}_{\alpha(3)}) \otimes \bfa^{3;3}_{\alpha(3)} ).
\endaligned
\end{equation}
Here we omit the subscript $k$, $\vec{\kappa}$ etc. of 
operators $\frak T$, $\frak m^{{\rm f.c.u.}\frak b}$ since it is automatically determined 
by the input.
It is easy to see that $\delta^H\circ \delta^H = 0$. We thus obtain the
$A_{n,k}$ analogue of Hochschild complex $(CH^*_{E,k}(\cL^{\rm form}_{\rm uni}),\delta^H)$.
We can define a structure of Lie algebra on it in the same way as \cite[Definition 7.4.26]{fooo092} so that it becomes a differential graded Lie algebra.
\begin{rem}
Actually we are in the situation that 
$\frak q^{\text{\rm form}}_{\ell;\vec{\kappa};B}$ 
is defined only partially.
It is obvious that if we have $\frak q^{\text{\rm form}}_{\ell;\vec{\kappa};B}$  
for $(B,\vert\vec{\kappa}\vert,\ell) \le (E,k,\ell_0)$ 
for certain $(E,k,\ell_0)$ sufficiently large compared to $(E_n,k)$ then 
(\ref{mmcat}) is defined .
\end{rem}

We now define maps\index[syindex]{qformellvechat@$\widehat{\frak q}^{\rm f.c.u.\frak b}_{\ell;\le k,E}$}
\begin{equation}\label{map83}
\widehat{\frak q}^{\rm f.c.u.\frak b}_{\ell;k,\le E} \co \Omega(X)^{\otimes \ell} \to  (CH^*_{E,k}(\cL^{\rm form}_{\rm c.u.}),\delta^H), 
\end{equation}
as follows.  We fix $\ell_0$ and consider $\ell$ with $\ell \le \ell_0$ only.

Let ${\bf g} \in \Omega(X)^{\otimes \ell}$.
We put
\begin{equation}\label{qqmcat}
\frak q_{\ell,k}^{{\rm f.c.u.}\frak b}({\bf g},\text{\bf h})
=
\sum_{B,\ell',\vec{\kappa}}T^{\omega(B) }
\frac{\rho_{\frak b,\theta}(B)% \rho_{b}(B)
}{\ell'!}
\frak q^{\text{\rm form}}_{\ell+\ell';\vec{\kappa};B}
({\bf g},\frak b_+^{\ell'};
\text{\bf h}). %;b_+^{\vec k})
\end{equation}
For ${\bf h} = h_1\otimes \dots \otimes h_k$ which does not contain $T$, \index[syindex]{qformellkecu@${\frak q}_{\ell,k}^{{\rm f.c.u.}\frak b}$}
\index[syindex]{qformellkecuE@${\frak q}_{\ell;k,E}^{{\rm f.c.u.}\frak b}$}
we rewrite it as 
\begin{equation}\label{eq145+}
{\frak q}_{\ell,k}^{{\rm f.c.u.}\frak b}({\bf g},\text{\bf h}) 
= 
\sum_{E \in G} T^E
\frak q^{{\rm f.c.u.}\frak b}_{\ell;k,E}({\bf g};\text{\bf h}).
\end{equation}
Here $\frak q^{{\rm f.c.u.}\frak b}_{\ell;k,E}$ is $\F$ linear.
Now we define
$$
\widehat{\frak q}^{{\rm f.c.u.}\frak b}_{\ell;\le k,E}({\bf g}) = (\frak T_{\ell;E_i,\vec{\kappa}}({\bf g}))_{(E_i,\vert\vec{\kappa}\vert) \le (E,k)}
$$
where
\begin{equation}\label{form147147}
\frak T_{\ell;E_i,\vec{\kappa}}({\bf g})(\text{\bf h})= \frak q^{{\rm f.c.u.}\frak b}_{\ell;k,E_i,}({\bf g};\text{\bf h}).
\end{equation}
\begin{lem}
$\wq^{\rm f.c.u.\frak b}_{1;k,\le E}$ is a chain map.
\end{lem}
\begin{proof}
This is a consequence of Propositions \ref{Kuraeistspoly}, \ref{existmultipolu}
and Lemma \ref{qpropertiescat} and its variants containing ${\bf e}^+_c$ and ${\bf f}_c$ (See Section \ref{sec:unit}).
\end{proof}
\par
As we explained several times, at this stage we obtain $\wq^{\rm f.c.u.\frak b}_{\ell;k,E}$ for each fixed $(E,k)$.
In other words its target is $CH^*_{E,k}(\cL^{\rm form}_{\rm uni})$, the $A_{n,k}$ version of Hochschild complex, 
and is not the full  Hochschild complex.
A similar issue appeared during the construction of the $A_{\infty}$ operations 
$\frak m_k$ in Section \ref{sec:cyclicfil}, which we resolved by the use of homotopy limit argument. 
We use the same strategy to construct a map \begin{equation}\label{mapqhatde2}
\wq^{{\rm f.c.u.}\frak b}_1:\Omega(X)
\to 
CH^*(\cL ^{\rm form}_{\rm uni}, 
\cL ^{\rm form}_{\rm uni})
\end{equation} 
from a
sequence of maps (\ref{form147147}).
To work out the homotopy  limit we need to show 
the independence of the map $\wq^{\rm f.c.u.\frak b}_{\ell,k,E}$ 
of the choices up to homotopy.
\par
The main problem to implement this inductive scheme is the fact that 
Hochschild cohomology is neither covariant nor contravarant.
We use a variant of pseudo-isotopy of Hochschild cochain complex to 
go around it.
We introduce this variant below.
\begin{defn}
We use the notation in Definitions \ref{defn626} and \ref{defn78}.
\par
For given $\vec{\kappa} = (\kappa_0,\dots,\kappa_K)$ we define\index[syindex]{HomBCFainfty@$\Hom^{\infty}(BCF(\cL^{\rm form};\vec{\kappa})^+,{CF}(L_{\kappa_{0}},L_{\kappa_{K}};\F)^+)$} 
$$
\Hom^{\infty}(BCF(\cL^{\rm form};\vec{\kappa})^+,{CF}(L_{\kappa_{0}},L_{\kappa_{K}};\F)^+)
$$
to be the set of all $\F$-linear maps
$$
\frak{a} : 
BCF(\cL^{\rm form};\vec{\kappa})^+ \to {CF}(L_{\kappa_{0}},L_{\kappa_{K}};\F)^+
$$
such that there exists a distributional differential form 
$\mathcal K_{\frak a}$ on $\prod_{i=1}^K (\tilde L_{\kappa_{i-1}} \times_X \tilde L_{\kappa_{i}}) \times  (\tilde L_{\kappa_{0}} \times_X \tilde L_{\kappa_{K}})$ 
with Condition \ref{cod149} below and such that
 %\marginpar{Codition for $\mathcal K_{\frak a}$ is modified. 
%Actually it may not be smooth since evaluation maps to $L_{\kappa_{i-1}} \cap L_{\kappa_{i}}$ may not be weakly submersive 
%except the last facor.}
$$
\aligned
&\int_{\tilde L_{\kappa_{0}} \times_X \tilde L_{\kappa_{K}}}\frak{a}(h_1,\dots,h_K) \wedge h_0 \\
&= 
\int_{\prod_{i=1}^K (\tilde L_{\kappa_{i-1}} \times_X \tilde L_{\kappa_{i}}) \times  (\tilde L_{\kappa_{0}} \times_X \tilde L_{\kappa_{K}})}
(h_1 \times \dots \times h_K \times h_0) \wedge \mathcal K_{\frak a}
\endaligned
$$
for any smooth differential forms $h_i$ on $\tilde L_{\kappa_{i-1}} \times_X \tilde L_{\kappa_{i}}$ and 
$h_0$ on $\tilde L_{\kappa_{0}} \times_X \tilde L_{\kappa_{K}}$.
We call $\mathcal K_{\frak a}$ the {\it Schwartz kernel}\index{Schwartz kernel} of $\frak a$.\footnote{We only 
require existence of $\mathcal K_{\frak a}$. So we do not need to specify the sign in the formula.}
\par
\begin{conds}\label{cod149}
We consider the wave front set\index{wave front set} $WF(\mathcal K_{\frak a})$\index[syindex]{WF@$WF$} of $\mathcal K_{\frak a}$ in the sense of 
\cite{Ho}.  For each $x \in L_{\kappa_{0}} \cap L_{\kappa_{K}}$ 
let $F_x$ be the submanifold $\prod_{i=1}^K (\tilde L_{\kappa_{i-1}} \times_X \tilde L_{\kappa_{i}})  \times \{x\}$ of ${\prod_{i=1}^K (\tilde L_{\kappa_{i-1}} \times_X \tilde L_{\kappa_{i}}) \times  (\tilde L_{\kappa_{0}} \times_X \tilde L_{\kappa_{K}})}$.

We consider the co-normal bundle 
$
N^*_{F_x}\left({\prod_{i=1}^K (\tilde L_{\kappa_{i-1}} \times_X \tilde L_{\kappa_{i}}) \times  (\tilde L_{\kappa_{0}} \times_X \tilde L_{\kappa_{K}})}\right).
$
We then require
\begin{equation}
WF(\mathcal K_{\frak a}) \cap N^*_{F_x}\left({\prod_{i=1}^K (\tilde L_{\kappa_{i-1}} \times_X \tilde L_{\kappa_{i}}) 
\times  (\tilde L_{\kappa_{0}} \times_X \tilde L_{\kappa_{K}})}\right)
= \{0\}
\end{equation}
for all $x$. We also assume similar conditions  for the part where ${\bf e}^+_c$, ${\bf f}_c$ are involved.
\end{conds}
Note that the existence of smooth differential form $\frak{a}(h_1,\dots,h_K)$ is then a consequence of
\cite[Theorem 2.5.12]{Ho}.

We next define\index[syindex]{HomBCFainftyfrac@$\Hom^{\infty}(BCF(\cL^{\rm form};\vec{\kappa})^+,{\frak C}(\newred{L}_{\kappa_{0}},\newred{L}_{\kappa_{K}};\F)^+)$}  
$$
\Hom^{\infty}(BCF(\cL^{\rm form};\vec{\kappa})^+,{\frak C}(\newred{L}_{\kappa_{0}},\newred{L}_{\kappa_{K}};\F)^+)
$$
to be the set of all $t \in [0,1]$ parametrized family of elements 
$$
\frak f_t \in \Hom^{\infty}(BCF(\cL^{\rm form};\vec{\kappa})^+,
{CF}(\newred{L}_{\kappa_{0}},\newred{L}_{\kappa_{K}};\F)^+)
$$
which depends smoothly on $t$, that is, its Schwartz kernel depends smoothly on $t$.
More precisely we require the Schwartz kernel $\mathcal K(t;*,*)$ 
satisfies the next condition. %\marginpar{Condition added.  KF 2025 Sep.}
\begin{conds}\label{cond159}
Let $\mathcal K$ be a differential form on $[0,1]_t \times \prod_{i=1}^K (\tilde L_{\kappa_{i-1}} \times_X \tilde L_{\kappa_{i}}) \times  (\tilde L_{\kappa_{0}} \times_X \tilde L_{\kappa_{K}})$.
We consider its wave front set $WF(\mathcal K)$. For each $x \in L_{\kappa_{0}} \cap L_{\kappa_{K}}$ and $t \in [0,1]$ 
let $F_{t,x}$ be the submanifold $\{t\} \times \prod_{i=1}^K (\tilde L_{\kappa_{i-1}} \times_X \tilde L_{\kappa_{i}})  \times \{x\}$ of 
$[0,1]_t \times {\prod_{i=1}^K (\tilde L_{\kappa_{i-1}} \times_X \tilde L_{\kappa_{i}}) \times  (\tilde L_{\kappa_{0}} \times_X \tilde L_{\kappa_{K}})}$.
\par
Consider the co-normal bundle 
$
N^*_{F_{t,x}}\left([0,1]_t \times {\prod_{i=1}^K (\tilde L_{\kappa_{i-1}} \times_X \tilde L_{\kappa_{i}}) \times  (\tilde L_{\kappa_{0}} \times_X \tilde L_{\kappa_{K}})}\right).
$
We then require
\begin{equation}
WF(\mathcal K) \cap N^*_{F_{t,x}}\left([0,1]_t\times {\prod_{i=1}^K (\tilde L_{\kappa_{i-1}} \times_X \tilde L_{\kappa_{i}}) 
\times  (\tilde L_{\kappa_{0}} \times_X \tilde L_{\kappa_{K}})}\right)
= \{0\}
\end{equation}
for all $t,x$.
\end{conds}
\par
Note that 
$$
\aligned
&{\frak C}(\newred{L}_{\kappa_{0}},\newred{L}_{\kappa_{K}};\F)^+((L,\theta_L,b),(L,\theta_L,b))\\
&= 
{\frak C}(\newred{L}_{\kappa_{0}},\newred{L}_{\kappa_{K}};\F)((L,\theta_L,b),(L,\theta_L,b))
\oplus
\Lambda^{\F}_0[{\bf e}^+]
\oplus 
\Lambda^{\F}_0[{\bf f}].
\endaligned
$$ 
%contains an element such as $f(t) {\bf e}^+_c$ (or those with $ {\bf e}^+_c$ replaced by $ {\bf f}_c$.)
%\par
 
We define\index[syindex]{CHinftyG@$CH_{\infty,G}^*(\cL_{\rm c.u.} ^{\rm form},\cL ^{\rm form}_{\rm c.u.})$} 
$$
CH_{\infty;G}^*(\cL_{\rm c.u.} ^{\rm form},\cL ^{\rm form}_{\rm c.u.})^+
$$
to be the set of all elements $\sum T^{\lambda_i}\alpha_i$
where 
$$
\alpha_i \in 
\Hom^{\infty}(BCF(\cL^{\rm form};\vec{\kappa})^+,{CF}(L_{\kappa_{0}},L_{\kappa_{K}};\F)^+)
$$
$\lambda_i \in G$.  
We also define \index[syindex]{CHinftyE@$CH_{\infty,E,k}^*(\cL_{\rm c.u.} ^{\rm form},\cL ^{\rm form}_{\rm c.u.})$} 
$$
CH_{\infty;E,k}^*(\cL_{\rm c.u.} ^{\rm form},\cL ^{\rm form}_{\rm c.u.})^+
$$
to be the set of such $\sum T^{\lambda_i}\alpha_i$ with $(\lambda_i,\vert \vec{\kappa}\vert) \le (E,k)$.
\end{defn}

We remark that $CH_{\infty;G}^*(\cL_{\rm c.u.} ^{\rm form},\cL ^{\rm form}_{\rm c.u.})^+$
and $CH_{\infty;E,k}^*(\cL_{\rm c.u.} ^{\rm form},\cL ^{\rm form}_{\rm c.u.})^+$
are the same as 
$CH_{\infty;G}^*(\cL_{\rm uni} ^{\rm form},\cL ^{\rm form}_{\rm uni})$
and $CH_{\infty;E,k}^*(\cL_{\rm uni} ^{\rm form},\cL ^{\rm form}_{\rm uni})$.

\begin{defn}\label{defn148148}
Let $G \subset \R_{\ge 0}$ be a discrete monoid.\index[syindex]{CH*infinityG@$CH_{\infty;G}^*(\cL_{\rm c.u.} ^{\rm form},\cL ^{\rm form}_{[0,1]})$}
We define\index[syindex]{CHinftyLcu[0,1]@$CH_{\infty;G}^*(\cL_{\rm c.u.} ^{\rm form},\cL ^{\rm form}_{[0,1]})$}  
$$
CH_{\infty;G}^*(\cL_{\rm c.u.} ^{\rm form},\cL ^{\rm form}_{[0,1]})^+
$$
to be the set of all collections $(\widehat{\widehat{\frak{a}}}_{\vec{\kappa}})_{\vec{\kappa}}$
where $\widehat{\widehat{\frak{a}}}_{\vec{\kappa}} = {\widehat{\frak{a}}}^0_{\vec{\kappa}} + dt \wedge {\widehat{\frak{a}}}^1_{\vec{\kappa}}$
and 
$$
\widehat{\frak{a}}^0_{\vec{\kappa}}
= \sum_{E \in G} T^{E}{\frak{a}}^0_{\vec{\kappa};E}, \qquad
\widehat{\frak{a}}_{\vec{\kappa}}^1
= \sum_{E \in G} T^{E}{\frak{a}}^1_{\vec{\kappa};E}
$$
with
$$
{\frak{a}}^0_{\vec{\kappa};E}, {\frak{a}}^1_{\vec{\kappa};E}
\in \Hom^{\infty}(BCF(\cL^{\rm form};\vec{\kappa})^+,
{\frak C}(\newred{L}_{\kappa_{0}},\newred{L}_{\kappa_{K}};\F)^+).
$$
For $E \in G$ and $k$ we define 
$$
CH_{\infty;E,k}^*(\cL_{\rm c.u.} ^{\rm form},\cL ^{\rm form}_{[0,1]})^+
$$
in a similar way. This is a subset of $CH^*_{E,k}(\cL^{\rm form}_{\rm c.u.})^+$.
\end{defn}
We define evaluation maps:\index[syindex]{Euvt=i@${\rm Ev}_{t=i}$}
$$
{\rm Ev}_{t=i} : 
CH_{\infty;G}^*(\cL_{\rm c.u.} ^{\rm form},\cL ^{\rm form}_{[0,1]})^+
\to
CH_{\infty;G}^*(\cL ^{\rm form}_{\rm c.u.},\cL _{\rm c.u.}^{\rm form})^+
$$
by
\begin{equation}\label{form1512}
{\rm Ev}_{t=i}(\widehat{\frak{a}}^0_{\vec{\kappa}} + dt \wedge \widehat{\frak{a}}^1_{\vec{\kappa}})
=
\widehat{\frak{a}}^0_{\vec{\kappa}}\vert_{t=i},
\end{equation}
for $i=0,1$.
Condition \ref{cond159} implies the restriction in the right hand side of (\ref{form1512}) exists. 
Here we use the obvious restriction map 
$$
\aligned
&\Hom^{\infty}(BCF(\cL^{\rm form};\vec{\kappa})^+,{\frak C}(\newred{L}_{\kappa_{0}},\newred{L}_{\kappa_{K}};\F)^+) \\
&\qquad\to
\Hom^{\infty}(BCF(\cL^{\rm form};\vec{\kappa})^+,C(L_{\kappa_{0}},L_{\kappa_{K}};\F)^+)
\endaligned
$$
in the right hand side.
The evaluation map on $
CH_{\infty;E,k}^*(\cL_{\rm c.u.} ^{\rm form},\cL ^{\rm form}_{[0,1]})^+
$ is defined in a similar way.
\begin{defn}
The differential 
 %\marginpar{Definition added. KF 2024 Nov.  position of $dt$
%changed so maybe sign affected. KF 2025 Jan}  
$\delta^H$ on $CH_{\infty;G}^*(\cL ^{\rm form}_{\rm c.u.},\cL ^{\rm form}_{[0,1]})^+$ is defined as 
follows. Let $\widehat{\widehat{\frak{a}}}_{\vec{\kappa}} = {\widehat{\frak{a}}}^0_{\vec{\kappa}} + dt \wedge{\widehat{\frak{a}}}^1_{\vec{\kappa}}$
be as in Definition \ref{defn148148}. 
Let $\{\frak m^{\frak C, t}_{k,\beta}\}$, $\{\frak c^{\frak C, t}_{k,\beta}\}$ be the 
structure operations of the pseudo-isotopy $\frak C = \cL ^{\rm form}_{[0,1]}$ as in Definition \ref{pisotopydef}.  
(This is the definition of $\cL ^{\rm form}_{[0,1]}$.  Note that $\frak C$ is defined in Subsection \ref{sec:homotopyequiv}.)
We put
$$
\delta^H(\widehat{\widehat{\frak{a}}}_{\vec{\kappa}}) = (\widehat{\widehat{\frak{a}}}'_{\vec{\kappa}})
$$
where $\widehat{\widehat{\frak{a}}}'_{\vec{\kappa}} = {\widehat{\frak{a}}}^{0 \prime}_{\vec{\kappa}} + dt
\wedge {\widehat{\frak{a}}}^{1 \prime}_{\vec{\kappa}}$,
with 
$
\widehat{\frak{a}}^{0 \prime}_{\vec{\kappa}}
= \sum_{E \in G} T^{E}{\frak{a}}^{0 \prime}_{\vec{\kappa};E}$, 
$\widehat{\frak{a}}^{1 \prime}_{\vec{\kappa}}
= \sum_{E \in G} T^{E}{\frak{a}}^{1 \prime}_{\vec{\kappa};E}
$
and 
$$
\aligned
{{\frak{a}}}^{0 \prime}_{\vec{\kappa};E}(\bf a) 
=&\sum_{E_1 + E_2 = E} (-1)^{\deg(\frak{a}_{E_2}^{0})\deg'( \bfa^{3;1}_{\alpha(3)}) }   \frak m^{\frak C, t}_{E_1}(\bfa^{3;1}_{\alpha(3)} \otimes 
{{\frak{a}}}^{0}_{E_2}( \bfa^{3;2}_{\alpha(3)}) \otimes \bfa^{3;3}_{\alpha(3)} )  \\
&+ \sum_{E_1 + E_2 = E}  (-1)^{\deg'( \bfa^{3;1}_{\alpha(3)})  +\deg(\frak{a}^{0}_{E_1})+1}   {{\frak{a}}}^{0}_{E_1}(\bfa^{3;1}_{\alpha(3)} \otimes   
\frak m^{\frak C, t}_{E_2}( \bfa^{3;2}_{\alpha(3)}) \otimes \bfa^{3;3}_{\alpha(3)} ) .  
 \endaligned
$$
(Compare (\ref{form144144}).) Moreover
$$
\aligned
{{\frak{a}}}^{1 \prime}_{\vec{\kappa};E}(\bf a) 
=&\sum_{E_1 + E_2 = E} (-1)^{\maltese_1}   \frak c^{\frak C, t}_{E_1}(\bfa^{3;1}_{\alpha(3)} \otimes 
{{\frak{a}}}^{0}_{E_2}( \bfa^{3;2}_{\alpha(3)}) \otimes \bfa^{3;3}_{\alpha(3)} )  \\
&+ \sum_{E_1 + E_2 = E}  (-1)^{\maltese_2}   {{\frak{a}}}^{0}_{E_1}(\bfa^{3;1}_{\alpha(3)} \otimes   
\frak c^{\frak C, t}_{E_2}( \bfa^{3;2}_{\alpha(3)}) \otimes \bfa^{3;3}_{\alpha(3)} )\\
&+ \sum_{E_1 + E_2 = E} (-1)^{\maltese_3}   \frak m^{\frak C, t}_{E_1}(\bfa^{3;1}_{\alpha(3)} \otimes 
{{\frak{a}}}^{1}_{E_2}( \bfa^{3;2}_{\alpha(3)}) \otimes \bfa^{3;3}_{\alpha(3)} )\\
&+ \sum_{E_1 + E_2 = E}  (-1)^{\maltese_4}   {{\frak{a}}}^{1}_{E_1}(\bfa^{3;1}_{\alpha(3)} \otimes   
\frak m^{\frak C, t}_{E_2}( \bfa^{3;2}_{\alpha(3)}) \otimes \bfa^{3;3}_{\alpha(3)} )
\\
& + \frac{d \widehat{\frak{a}}^0_{\vec{\kappa},E}}{dt}  (\bf a).
\endaligned
$$
The signs are by Koszul rule,\footnote{Koszul sign here includes the sign caused by 
exchanging $dt$ with operators or variables. $\deg'\bf a$ 
appearing in the formula is caused by it.} that is,
$$
\aligned
\maltese_1 &= \deg(\frak{a}_{E_2}^{0})\deg'( \bfa^{3;1}_{\alpha(3)})\\
\maltese_2 &= \deg'({\bf a}^{3;1}_{\alpha(3)}) + \deg'(\frak{a}_{E_1}^{0})\\
\maltese_3 &= (\deg(\frak{a}_{E_2}^{1})+1)\deg'( \bfa^{3;1}_{\alpha(3)}) \\
\maltese_4 &=  \deg'( \bfa^{3;1}_{\alpha(3)}) + \deg'(\frak{a}_{E_1}^{1}) .
\endaligned
$$
Note that $\delta^H(\widehat{\widehat{\frak{a}}}_{\vec{\kappa}})$ is the (super)
 %\marginpar{Sign modified. Still to be checked.  KF. 2025 Jan}
commutator of $\widehat{\widehat{\frak{a}}}$ and $\frak m^{\frak C, t}+ dt \wedge \frak c^{\frak C, t} $, that is,
$$
(\frak m^{\frak C, t}+ dt \wedge \frak c^{\frak C, t}) \circ \widehat{\widehat{\frak{a}}}
+ (-1)^{\deg' \widehat{\widehat{\frak{a}}}
}\,\,\, \widehat{\widehat{\frak{a}}}
\circ (\frak m^{\frak C, t}+ dt \wedge \frak c^{\frak C, t}).
$$
(In the case when  $\frak m^{\frak C, t}$ is $\frak m^{\frak C, t}_1$ the term containing 
the $t$-derivative of $\widehat{\widehat{\frak{a}}}$ should also be included in an obvious way.)
\par
The equality $\delta^H \circ \delta^H = 0$ follows from this fact.
\par
It is easy to see that 
$\delta^H$ induces a differential on $
CH_{\infty;E,k}^*(\cL_{\rm c.u.} ^{\rm form},\cL ^{\rm form}_{[0,1]})^+
$.
\end{defn}

We work in the situation of Subsection \ref{sec:homotopyequiv}.
Consider  a finite set $\{ (L_{\kappa}, \theta_\kappa) \}$ of Lagrangian submanifolds of $X$ equipped with 
relative spin structures and primitives $\theta_\kappa$ for the degree $2$ part of the leading order terms of the bulk.
Suppose $E ,E' \in G$, $E'=E_{n'} \le E=E_{n}$ with $(E',k') \le (E,k)$.
We have a curved filtered $A_{n',k'}$ category 
and a curved filtered $A_{n,k}$ category, the set of 
whose objects are identified with the set $\{( L_{\kappa}, \theta_\kappa) \}$,
which we denote by $\cL^{\rm form}_{\mathcal C}$ and $\cL^{\rm form}_{\mathcal D}$.
In other words, $\cL^{\rm form}_{\Cat} = \cL^{\rm form}_{\frak b;n',k'}$, 
$\cL^{\rm form}_{\Dat} = \cL^{\rm form}_{\frak b;n,k}$.
(See Definition \ref{defn838} for the right hand side.)

There is also a pseudo-isotopy of $A_{n',k'}$ categories between $\cL^{\rm form}_{\mathcal C}$ and 
$\cL^{\rm form}_{\mathcal D}$.
This pseudo-isotopy is also a curved filtered $A_{n',k'}$ category, 
which we write $\cL^{\rm form}_{[0,1]}$.
We consider their ($A_{n',k'}$ or $A_{n,k}$  versions of) Hochschild complexes and denote them by
$$
CH_{n',k'}^*(\cL^{\rm form}_{\mathcal C},\cL^{\rm form}_{\mathcal C})^+,
\quad
CH_{n,k}^*(\cL^{\rm form}_{\mathcal D},\cL^{\rm form}_{\mathcal D})^+,
\quad
CH_{n',k'}^*(\cL_{\rm c.u.}^{\rm form},\cL^{\rm form}_{[0,1]})^+,
$$
respectively.  (Here $CH_{n,k}^* = CH_{E_n,k}^*$ and etc.  $+$ means that it includes ${\bf e}^+$ and ${\bf f}$.)
\par
We first observe there exists a chain map
$$
CH_{n,k}^*(\cL^{\rm form}_{\mathcal D},\cL^{\rm form}_{\mathcal D})^+
\to CH_{n',k'}^*(\cL^{\rm form}_{\mathcal D},\cL^{\rm form}_{\mathcal D})^+
$$
which is defined by forgetting $\frak T_{E'',\kappa''}$ with 
$(E'',\kappa'') \nleqq (E',k')$, $(E'',\kappa'') \leqq (E,k)$.
\par
We also observe that there exist obvious $\Lambda_G$ linear maps:
\begin{equation}\label{form86}
\aligned
&CH_{\infty;n',k'}^*(\cL_{\mathcal C}^{\rm form},\cL_{\mathcal C}^{\rm form})^+
\to CH^*_{n',k'}(\cL^{\rm form}_{\mathcal C},\cL^{\rm form}_{\mathcal C})^+,
\\
&CH_{\infty;n,k}^*(\cL_{\mathcal D}^{\rm form},\cL_{\mathcal D}^{\rm form})^+
\to CH^*_{n,k}(\cL^{\rm form}_{\mathcal D},\cL^{\rm form}_{\mathcal D})^+,
\\
&CH_{\infty;n',k'}^*(\cL_{\rm c.u.}^{\rm form},\cL ^{\rm form}_{[0,1]})^+
\to CH^*_{n,k}(\cL^{\rm form}_{\rm c.u.},\cL^{\rm form}_{[0,1]})^+,
\endaligned
\end{equation}
which are injective.
By construction the Hochschild differential defined on the right hand side 
preserves the image and hence define a differential 
on the left hand side. (Its square is zero.)
We observe that the {\it Gerstenhaber bracket}\index{Gerstenhaber bracket} on the right hand side 
preserves the image and so the left hand sides are differential graded 
Lie algebras.
\par
In the case of $CH^*_{n,k}(\cL_{\rm c.u.}^{\rm form},\cL^{\rm form}_{[0,1]})^+$ 
the Lie algebra structure is defined by
 %\marginpar{Definiiton added KF 2024 Nov.
%Modified 2025 Jan. KF}
$$
[ {\widehat{\frak{a}}}^0 + dt \wedge{\widehat{\frak{a}}}^1,
 {\widehat{\frak{a}}}^{0 \prime} + dt \wedge{\widehat{\frak{a}}}^{1 \prime}]
 =
 [{\widehat{\frak{a}}}^{0},{\widehat{\frak{a}}}^{0 \prime}]
 +(-1)^{\deg' \widehat{\frak{a}}^{0}} dt \wedge
 [{\widehat{\frak{a}}}^{0},{\widehat{\frak{a}}}^{1 \prime}]
- dt \wedge
 [{\widehat{\frak{a}}}^{1},{\widehat{\frak{a}}}^{0 \prime}].
$$
Here
$[{\widehat{\frak{a}}}^{0},{\widehat{\frak{a}}}^{0 \prime}]
= {\widehat{\frak{a}}}^{0} \circ {\widehat{\frak{a}}}^{0 \prime} \pm 
 {\widehat{\frak{a}}}^{0 \prime} \circ {\widehat{\frak{a}}}^{0}$.
 We can prove that the Schwartz kernel  %\marginpar{An argument added.  KF 2025 June}
 $\mathcal K_{{\widehat{\frak{a}}}^{0} \circ {\widehat{\frak{a}}}^{0 \prime}}$
 of ${\widehat{\frak{a}}}^{0} \circ {\widehat{\frak{a}}}^{0 \prime}$ satisfies 
 Condition \ref{cod149} as follows.
 Let  $\mathcal K_{{\widehat{\frak{a}}}^{0}}(t;x_1,\dots,x_k,x_0)$
 and $\mathcal K_{t;{\widehat{\frak{a}}}^{0 \prime}}(y_1,\dots y_{k'},y_0)$
 be the Schwartz kernels of  ${\widehat{\frak{a}}}^{0}$ 
 and ${\widehat{\frak{a}}}^{0 \prime}$, respectively.
 (Here $x_i,y_j \in L$ and $x_0$, $y_0$ correspond to the $0$-th marked 
 points.)
 We consider the product 
$\mathcal K_{{\widehat{\frak{a}}}^{0}} \mathcal K_{{\widehat{\frak{a}}}^{0 \prime}}$,
 which is a distributional form on $L^{k+k'+2}$.  We can restrict it 
to a submanifold defined by $y_0 = x_i$. (Here $i=1,\dots,k$.)
Using  Condition \ref{cod149} of ${\widehat{\frak{a}}}^{0 \prime}$ 
this is a consequence of \cite[Theorem 2.5.11]{Ho}.
\par
The Schwartz kernel $\mathcal K_{{\widehat{\frak{a}}}^{0} \circ {\widehat{\frak{a}}}^{0 \prime}}$
 of ${\widehat{\frak{a}}}^{0} \circ {\widehat{\frak{a}}}^{0 \prime}$ is a signed sum of 
\begin{equation}\label{form1413}
 \int_{y_0}  \mathcal K_{\widehat{\frak{a}^{0}}}(t;x_1,\dots,x_{i-1},y_0,x_{i+1},\dots,x_0) \mathcal K_{\widehat{\frak{a}^{0 \prime}}}(t;y_1,\dots,y_{k'},y_0).
\end{equation}

Then by Condition \ref{cond159} of ${\widehat{\frak{a}}}^{0}$
we can use \cite[Theorems 2.5.11' and 2.5.12]{Ho} to estimate the 
wave front set of  (\ref{form1413}).

Therefore ${\widehat{\frak{a}}}^{0} \circ {\widehat{\frak{a}}}^{0 \prime}$  
satisfies Condition \ref{cond159}.

The next lemma is also obvious from construction: %\marginpar{Where ${\rm Ev}_{t=i}$ is defined. KF 2025 Aug.}

\begin{lem}
The maps ${\rm Ev}_{t=i}$ (defined by $(\ref{form1512})$) are (linear) differential graded Lie algebra homomorphisms. 
\end{lem}

\begin{rem}
The filtered $A_{\infty}$ functor $\cL^{\rm form}_{[0,1]} 
\to \cL^{\rm form}_{\mathcal C}$ does {\it not} induce a
differential Lie algebra homomorphism
$CH(\cL^{\rm form}_{[0,1]},\cL^{\rm form}_{[0,1]})^+
\to CH(\cL^{\rm form}_{\mathcal C},\cL^{\rm form}_{\mathcal C})^+$.
\end{rem}

We now claim that the map $\frak q^{{\rm f.c.u.}\frak b}_{\ell}$ from Equation \eqref{qqmcat} yields an $L_{\ell_0}$ homomorphism from the cochains of $X$ to the above model for the Hochschild cochains. 
(Note we fixed $\ell_0$ and considered only $\ell$ with $\ell \le \ell_0$.) 
More precisely, we need to modify the map a bit in the same way as 
in \cite[Subsection 7.4.3]{fooo092}. Namely for $\text{\bf g} 
= g_1 \otimes \dots \otimes g_{\ell}$ we put\index[syindex]{qformellvechato@$\widehat{\frak q}^{\rm f.c.u.\frak b}_{\ell;\vec{\kappa};o}$}
\begin{equation}\label{formnew88}
\widehat{\frak q}^{\rm f.c.u.\frak b}_{\ell;\vec{\kappa};o}(\text{\bf g})
= 
\widehat{\frak q}^{\rm f.c.u.\frak b}_{\ell;\vec{\kappa}}(\text{\bf g}^{\rm op})
\end{equation}
with
$
\text{\bf g}^{\rm op} = (-1)^{\sum_{i<j}\deg' g_i \deg' g_j}g_{\ell} \otimes \dots \otimes g_{1}
$.  %\marginpar{Sign put.  KF 2025 Aug.}
Here $\widehat{\frak q}^{\rm f.c.u.\frak b}_{\ell;\vec{\kappa}}$ is induced from 
${\frak q}^{\rm f.c.u.\frak b}_{\ell;\vec{\kappa}}$ in an obvious way.\footnote{${\frak q}^{\rm f.c.u.\frak b}_{\ell;\vec{\kappa}}$
is a map $\Omega(X)^{\otimes \ell}
\otimes 
BCF(\cL^{\rm form};\vec{\kappa})
\to CF(\newred{L_{\kappa_{\blue{0}}},L_{\kappa_{\blue{K}}}};\F)
\widehat{\otimes} \Lambda_0$.  $\widehat{\frak q}^{\rm f.c.u.\frak b}_{\ell;\vec{\kappa}}$ is a map
$\Omega(X)^{\otimes \ell}
\to  \Hom(
BCF(\cL^{\rm form};\vec{\kappa}),
CF(\newred{L_{\kappa_{\blue{0}}},L_{\kappa_{\blue{K}}}};\F)
\widehat{\otimes} \Lambda_0)$.  As we mentioned in Subsection \ref{labelqqhat}, this is the way $\frak q$ and $\hat{\frak q}$ are related in several other places.}
The next result is then a consequence of  Proposition \ref{prop627} (and a similar results for operations 
in Definition \ref{defn626uni}), in the same way as in \cite[Corollary 7.4.10]{fooo092}.

We also remark that the Schwartz kernels of $T^{\lambda}$ coefficients of $\widehat{\frak q}^{\rm f.c.u.\frak b}_{\ell;\vec{\kappa}}(\text{\bf g}^{\rm op})$ 
satisfy Condition \ref{cod149}.  This is a consequence of the facts that 
it is locally a push out of a smooth form by the evaluation map at the marked points and evaluation map 
at the $0$-th marked point is a submersion.\footnote{
In fact for a proper smooth map $F : M \to N$ and a current $\eta$ on $M$ we have:
$$
WF(F_*\eta) \subseteq \{(x,v) \in T^*N \mid \exists (y,w) \, y\in M, F(y) = x, w \in WF(\eta) \cap T^*_yM, (D_yF)^*v = w\}.
$$
%In the case when $F$ is an embedding this formula is a cosequence of \cite[Theorem 2.5.10]{Ho}.  
%In the case when $F$ is a submersion this formula is a consequence of \cite[Theorem 2.5.12]{Ho}. 
%The general case then follows since any $F$ is a composition of an embedding and a submersion.
See \cite[Proposition 1.3.4]{duister}.}  (The evaluation map at other marked points may not be a 
submersion and so the Schwartz kernel may not be smooth.) %\marginpar{A sentence added.  KF. 2015 June.}

In our situation we have two different systems of perturbations, corresponding to the choices made up to energies $E'$ and $E$. We denote the corresponding categories $\cL^{\rm form}_{\Cat}$ and $\cL^{\rm form}_{\Dat}$, respectively.
For each of them,  we obtain a map as in Lemma \ref{lem8787}.
We respectively denote them $\widehat{\frak q}^{\rm f.c.u.\frak b}_{o;\Cat}$
(Here $(E',k') = (E_{n'},k')$, which we omit from notation.)
and $\widehat{\frak q}^{\rm f.c.u.\frak b}_{o;\Dat}$ (Here $(E,k) = (E_{n},k)$, which we omit from notation.)
 to distinguish them. We also have a system of perturbations, $\frak C$ associated to the pseudo-isotopy between them constructed in Subsections \ref{sec:homotopyequiv} and \ref{subsec:opeartor}. 
\begin{lem}\label{lem8787}
The maps\index[syindex]{qfcubohatD@$\widehat{\frak q}^{\rm f.c.u.\frak b}_{o;\Dat}$} \index[syindex]{qfcubohatC@$\widehat{\frak q}^{\rm f.c.u.\frak b}_{o;\Cat}$} 
$$
\aligned
&\widehat{\frak q}^{\rm f.c.u.\frak b}_{o;\Dat} : 
\Omega(X)^{\otimes \ell} 
\to CH_{\infty;n,k}^*(\cL_{\Dat}^{\rm form},\cL_{\Dat}^{\rm form})^+ \\
&
\widehat{\frak q}^{\rm f.c.u.\frak b}_{o;\Cat} : 
\Omega(X)^{\otimes \ell} 
\to CH_{\infty;n',k'}^*(\cL_{\mathcal C}^{\rm form},\cL_{\Cat}^{\rm form})^+
\endaligned
$$
for $\ell \le \ell_0$
define  $L_{\ell_0}$ homomorphisms from 
$ \Omega(X)$ equipped  
with the trivial $L_{\ell_0}$ structure (except for the differential 
which is the de Rham 
differential). 
\end{lem} 

See \cite[Corollary 7.4.40]{fooo092}.
 
 The next result is the pseudo-isotopy analogue of Lemma \ref{lem8787}  the proof of which is omitted:
\begin{lem}\label{lem88}
The maps 
$\frak q^{\frak C,{\frak b}}_{\ell,\vec{\kappa}}$
in Definition $\ref{defn78}$ 
induce a $\Lambda_0$ homomorphism\index[syindex]{qfcubohatcfra@$\widehat{\frak q}^{\frak C,\frak b}_{\ell;o;n',k'}$}
$$
\widehat{\frak q}^{\frak C,\frak b}_{\ell;o;n',k'} :
\Omega(X)^{\otimes \ell} 
\to 
CH_{\infty;n',k'}^*(\cL_{\rm c.u.} ^{\rm form},\cL ^{\rm form}_{[0,1]})^+.
$$
The collection of such maps over all $\ell \le \ell_0$ defines an $L_{\ell_0}$ homomorphism 
$\widehat{\frak q}^{\frak C,\frak b}_{o;n',k'}$.
\end{lem}
%The proof is mostly the same as Lemma \ref{lem8787} and so is omitted.

\begin{lem}\label{lem8ten9}
For $n',k',\ell'_0$, the next diagram (of $L_{\ell'_0}$ homomorphisms) commutes.
\begin{equation}\label{diag88}
\xymatrix{ 
&&&&
CH_{\infty;n',k'}^*(\cL_{\Dat}^{\rm form},\cL_{\Dat}^{\rm form})^+ \\
\Omega(X)   
\ar[urrrr]^{\widehat{\frak q}^{\rm f.c.u.\frak b}_{o;\Dat}} 
\ar[rrrr]_{\widehat{\frak q}^{\frak C,\frak b}_{o;n',k'}}
\ar[drrrr]_{\widehat{\frak q}^{\rm f.c.u.\frak b}_{o;\Cat}} 
&&&&CH_{\infty;n',k'}^*(\cL ^{\rm form}_{\rm c.u.},\cL ^{\rm form}_{[0,1]})^+ 
\ar[u]_{{\rm Ev}_{t=1}}
\ar[d]^{{\rm Ev}_{t=0}}
\\
&&&&
CH_{\infty;n',k'}^*(\cL^{\rm form}_{\Cat},\cL ^{\rm form}_{\Cat})^+  
}
\end{equation}
\par\smallskip
\end{lem}
This is immediate from our construction.
\par
Note that in our situation the $L_{\ell_0}$ homomorphism 
$\widehat{\frak q}^{\rm f.c.u.\frak b}_{o;\Dat}$ 
is actually defined as a map to 
$CH_{\infty;n,k}^*(\cL_{\Dat}^{\rm form},\cL_{\Dat}^{\rm form})^+$
(where $(n',k')< (n,k)$)
and the one appearing in Diagram (\ref{diag88}) is its reduction.
Furthermore during the construction of the filtered $A_{\infty}$ structure 
we promoted $\cL ^{\rm form}_{\Cat}$  to a 
filtered $A_{n,k}$ category.
We now prove the next:
Let $(n',k') < (n,k)$ and $\ell'_0 \le \ell_0$.

\begin{lem}\label{lem810}
We can promote Diagram $(\ref{diag88})$ to the next Diagram (of $L_{\ell_0}$ homomorphisms) :
\begin{equation}\label{diag89}
\xymatrix{ 
&&&&
CH_{\infty;n,k}^*(\cL_{\Dat}^{\rm form},\cL_{\Dat}^{\rm form})^+\\
\Omega(X) 
\ar[urrrr]^{\widehat{\frak q}^{\rm f.c.u.\frak b}_{o;\Dat}} 
\ar[rrrr]_{\widehat{\frak q}^{\frak C,\frak b}_{o;n,k}}
\ar[drrrr]_{\widehat{\frak q}^{\rm f.c.u.\frak b}_{o;\Cat}} 
&&&&CH_{\infty;n,k}^*(\cL _{\rm c.u.}^{\rm form},\cL ^{\rm form}_{[0,1]})^+
\ar[u]_{{\rm Ev}_{t=1}}
\ar[d]^{{\rm Ev}_{t=0}}
\\
&&&&
CH_{\infty;n,k}^*(\cL ^{\rm form}_{\Cat},\cL ^{\rm form}_{\Cat})^+
}
\end{equation}
Here:
\begin{enumerate}
\item
The differential graded Lie algebra structures of the three modules in the right hand side are 
defined by the above mentioned promoted filtered $A_{n,k}$ structures.
\item
The map $\widehat{\frak q}^{{\rm f.c.u.}\frak b}_{o;\Dat}$ is the one in Lemma $\ref{lem8787}$.
\item
The maps $\widehat{\frak q}^{\frak C,\frak b}_{o;n,k}$ and $\widehat{\frak q}^{{\rm f.c.u.}\frak b}_{o;\Cat}$ 
become the ones in Diagram $(\ref{diag88})$ after reduction.
\end{enumerate}
\end{lem}
\begin{proof}
The proof is similar to the proof of Proposition \ref{prop:pseudo-isotopy-to-quasi-iso}.
We first consider the case when $\ell_0 = \ell_0'$.
We take $(E(i),k(i))$ ($i=1,\dots,I$) such that
$(E(1),k(1)) = (E_{n'},k')$, 
$(E(I),k(I)) = (E_{n},k)$,
$(E(i),k(i)) < (E(i+1),k(i+1))$,
and there is no $(E,k)$ with 
$(E(i),k(i)) < (E,k) < (E(i+1),k(i+1))$.
We lift maps to 
$CH^*_{\infty;E(i),k(i)}$ inductively on $i$.
Suppose lifts to $CH^*_{\infty;E(i),k(i)}$ are defined.
We need to define a lift to  $CH^*_{\infty;E(i+1),k(i+1)}$.
\par
The map $\frak q^{\frak C,\frak b}_{\ell;o;\le E(i+1),k(i+1)}$ 
sends an element of $\Omega(X)^{\otimes \ell}$
to 
$$
\Hom_{\infty}(BCF(\cL^{\rm form};\vec{\kappa})^+,{\frak C}(\newred{L}_{\kappa_{0}},\newred{L}_{\kappa_{K}};\F)^+).
$$
We define it from $\frak q^{\frak C,\frak b}_{\ell;o}$  in a similar way as (\ref{eq145+}).

The maps $\frak q^{\frak C,\frak b}_{\ell;o;\le E(i),k(i)}$  which\index[syindex]{qCbelloEEelll@$\frak q^{\frak C,\frak b}_{\ell;o;\le E,k}$}
we have by the inductive hypothesis can be written as\index[syindex]{qzQCbelloEEle@$\frak Q^{j,\frak C,\frak b}_{\ell;o;\le E,k}$} 
$$
\frak q^{\frak C,\frak b}_{\ell;o;\le E(i),k(i)}
= 
\frak Q^{0;\frak C,\frak b}_{\ell;o;\le E(i),k(i)} 
+dt \wedge 
\frak Q^{1;\frak C,\frak b}_{\ell;o;\le E(i),k(i)}
$$
where 
$$
\frak Q^{j;\frak C,\frak b}_{\ell;o;\le E(i),k(i)} (\text{\bf g})
\in CH_{\infty;E(i),k(i)}^*(\cL _{\rm c.u.}^{\rm form},\cL_{\rm c.u.} ^{\rm form})^+
$$
for $j=0,1$.
\par
$\frak Q^{j;\frak C,\frak b}_{\ell;o;\le E(i),k(i)}$ consists of 
$\frak Q^{j;\frak C,\frak b}_{\ell;o;E,k}$ with $(E,k) \le (E(i),k(i))$.\index[syindex]{qzQCbelloEE@$\frak Q^{j,\frak C,\frak b}_{\ell;o;E,k}$} 
\par

The map we look for is also written as
$$
\frak q^{\frak C,\frak b}_{\ell;o;E(i+1),k(i+1)}
= 
\frak Q^{0;\frak C,\frak b}_{\ell;o;E(i+1),k(i+1)} 
+dt \wedge 
\frak Q^{1;\frak C,\frak b}_{\ell;o;E(i+1),k(i+1)}.
$$
We put 
$$
\frak Q^{1;\frak C,\frak b}_{\ell;o;E(i+1),k(i+1)} = 0.
$$
Then the condition that $\frak q^{\frak C,\frak b}_{\ell;o;E(i+1),k(i+1)}$
together with  $\widehat{\frak q}^{\frak C,\frak b}_{\ell;o;\le E(i),k(i)}$  
is an $L_{\ell}$ homomorphism  
can be written as an ordinary differential equation:
\begin{equation}\label{newform1518}
\frac{d}{dt}\frak Q^{0;\frak C,\bf b}_{\ell;o;E(i+1),k(i+1)}
+ \frak P(\frak Q^{j;\frak C,\bf b}_{\ell';o;\le E(i),k(i)}) = 0
\end{equation}
where the second term is a quadratic expression 
which contains $\frak Q^{j;\frak C,\bf b}_{\ell';o;E,k}$ with $(E,k)
< (E(i+1),k(i+1))$, $\ell' \le \ell$, $j=0,1$ only.
(Such $\frak Q^{j;\frak C,\bf b}_{\ell';o;E,k}$ is a part of $\widehat{\frak q}^{\frak C,\frak b}_{\ell';o;n(i),k(i)}$ with $E_{n(i)} = E(i)$.)

Under the  initial condition, that is,
$$
\frak Q^{0;\frak C,\frak b}_{\ell;o;E(i+1),k(i+1)}\vert_{t=1} = 
\frak q^{\rm f.c.u.\frak b}_{o;E(i+1),k(i+1)},
$$
$\frak Q^{0;\frak C,\frak b}_{\ell;o;E(i+1),k(i+1)}$ is uniquely determined.
We thus obtain a lift  
$\frak q^{\frak C,\frak b}_{\ell;o;E(i+1),k(i+1)}$.
The lift $\frak q^{\rm f.c.u.\frak b}_{o;E(i+1),k(i+1)}$  
is its restriction at $t=0$.
\par
We can prove Condition \ref{cod149} is satisfied using the fact that the process of solving (\ref{newform1518}) 
corresponds to solving an ODE for the Schwartz kernels of the operators and 
so it preserves Condition \ref{cod149}. 
\par
Thus we have proved the lemma in the case $\ell_0 = \ell'_0$.
In case $\ell'_0 < \ell_0$ we first lift 
with $\ell_0 = \ell'_0$. Then we obtain the lift of $L_{\ell'_0}$ homomorphism 
to the $L_{\ell_0}$ homomorphism in a similar way.
\end{proof}
Let $\Cat$ be a filtered $A_{n',k'}$ category obtained by 
taking a system of Kuranishi structures and CF perturbations 
(of finitely many moduli spaces needed to define $A_{n',k'}$ structure.)
Then as we explained in 
Remark \ref{Remark849} we promote it to a filtered $A_{\infty}$ structure.
For a fixed $\ell_0$ we can also obtain an $L_{\ell_0}$ 
homomorphism $\widehat{\frak q}^{\rm f.c.u.\frak b}_{o;\Cat;E_{n'},k'}$
(using system of Kuranishi structures and CF perturbations 
on finitely many moduli spaces).

We then use Lemma \ref{lem810} inductively 
we can promote $\widehat{\frak q}^{\rm f.c.u.\frak b}_{o;\Cat;E_{n'},k'}$ 
to homomorphisms\index[syindex]{qfcuboEkell@$\frak q^{\rm f.c.u.\frak b}_{o;\Cat;E,k,\ell}$}
$$
\frak q^{\rm f.c.u.\frak b}_{o;\Cat;E,k,\ell} :
\Omega(X)^{\otimes\ell}
\to 
CH_{\infty;G}^*(\cL_{\rm c.u.} ^{\rm form},\cL ^{\rm form}_{[0,1]})^+
$$
for each $(E,k,\ell)$ which gives an $L_{\infty}$ homomorphism 
to the Hochschild cohomology of  the promoted filtered $A_{\infty}$ category.
\par
Composed with (\ref{form86}) it induces maps:\index[syindex]{qformelqlvechat@$\widehat{\frak q}^{\rm f.c.u.\frak b}_{\ell}$}
\begin{equation}\label{form14151415}
\wq^{\rm f.c.u.\frak b}_{\ell} : \Omega(X)^{\otimes \ell} 
\to 
CH^*(\cL^{\rm form}_{\Cat},\cL^{\rm form}_{\Cat})^+.
\end{equation}
It is easy to see that they induce the required map (\ref{mapqhatde2}).

We finally include the bounding cochains.
We study the Hochschild invariants of the category $\cL ^{\rm form}$ from Section \ref{sec:elicurv}, with objects $\bL$ consisting of a collection $\{(L_{\kappa},\theta_{\kappa}, b_{\kappa})\mid \kappa =1,\dots,K\}$  of Lagrangians with  bounding cochains.
Note that the morphism complex $\cL ^{\rm form}_{\rm c.u.}((L_{\kappa},\theta_{\kappa}, b_{\kappa}),(L_{\kappa'},\theta_{\kappa'}, b_{\kappa'}))$ is the same as 
$\cL ^{\rm form}((L_{\kappa},\theta_{\kappa}),(L_{\kappa'},\theta_{\kappa'}))$
except the case $(L_{\kappa},\theta_{\kappa}) = (L_{\kappa'},\theta_{\kappa'})$
but $b_{\kappa} \ne b_{\kappa'}$.
In such a case 
the latter  contains ${\bf e}^+$ and 
${\bf f}$ but the former not.
\par
The map\index[syindex]{qformfubdfell@$\wq^{\rm f.u.\bf b}_{\ell}$}
\begin{equation}\label{map8322}
\wq^{\rm f.u.\bf b}_{\ell} \co \Omega(X)^{\otimes \ell} \to  CH^*(\cL^{\rm form}_{\rm c.u.},\cL^{\rm form}_{\rm c.u.})^+, 
\end{equation}
is obtained from (\ref{form14151415}) by\index[syindex]{qformfubdfellnohat@$\frak q^{\rm f.u.\bf b}_{\ell}$} putting
\begin{equation}  \label{qcat-bounding-cochain-deformed}
\frak q^{\rm f.u.\bf b}_{\ell}(\text{\bf g};\text{\bf h})
=
\sum_{0 \leq m_i} \frak q^{\rm f.u.\frak b}_{\ell}(\text{\bf g};b^{m_0}_{\kappa_0} \otimes h_1 \otimes b^{m_1}_{\kappa_1} \otimes \cdots \otimes h_k \otimes  b^{m_k}_{\kappa_k})
\end{equation}
Here we denote by $\bf b$ the total data of the bulk deformation $\frak b$ and the bounding cochains $\{b_{\kappa}(E)\}$.
The case $\ell = 1$ is $\wq$.\footnote{Here $\cL^{\rm form}_{\rm uni}$ is a homotopically unital $A_{\infty}$
category whose objects are all weakly curvature free.}

We  have the following.  %\marginpar{The rest of this subsection is added in 2025 March KF}
Note that $CH^*(\cL _{\rm uni}^{\rm form},\cL _{\rm uni}^{\rm form})
=CH^*(\cL _{\rm c.u.}^{\rm form},\cL _{\rm c.u.}^{\rm form})^+$.

\begin{prop}\label{prop1416}
If ${\bf h} = {\bf h}_1 \otimes {\bf e}^+_c \otimes {\bf h}_2$
then
\begin{equation}\label{eq1417}
\wq^{\rm f.u.\bf b}_{\ell}(g;{\bf h}) = 0.
\end{equation}
In other words $\wq^{\rm f.u.\bf b}_{\ell}$ factors as
$$
\Omega(X)^{\otimes \ell} \to CH^*_{\rm red}(\cL _{\rm uni}^{\rm form},\cL _{\rm uni}^{\rm form})
\to CH^*(\cL _{\rm uni}^{\rm form},\cL _{\rm uni}^{\rm form}).
$$
Here $ CH^*_{\rm red}(\cL _{\rm uni}^{\rm form},\cL _{\rm uni}^{\rm form})$ is the reduced 
version of the Hochschild  complex. (Definition $\ref{defn525}$.)
\end{prop}
\begin{proof}
Formula (\ref{eq1417}) is a consequence of Proposition \ref{prop92} (8) and Proposition \ref{prop94} (8).\footnote{
While we perform induction argument via pseudo-isotopy we can and will do it so that unitality is 
preserved in an obvious sense.}
The second half then is a consequence of the definition.
\end{proof}
The proof of Proposition \ref{prop142} is complete. 
\qed

\subsection{Construction of  $\widehat{\frak q}$.}
\label{frakqoncoho}
The map (\ref{mapqhat}) is induced from $\wq^{\rm f.u.\bf b}_{1}$ (the $\ell=1$ case of (\ref{qcat-bounding-cochain-deformed})) 
and the invariance of Hochschild (co)-homology
under homotopy equivalence. 
It is however useful to define it directly.
(We will use this explicit construction  in Sections \ref{dualpq}, \ref{sec:ring} and \ref{sec: value}.)
\par
A {\it decorated ribbon tree with one interior marked point} 
is, by definition, a pair $(\Gamma,v)$ consisting of a decorated 
ribbon tree $\Gamma$ in the sense of Definition \ref{decribbon} and 
$v \in C^0_{\text{int}}(\Gamma)$.
Below we will define\index[syindex]{qcangammav@${\frak q}_{(\Gamma,v)}$} 
\begin{equation}
\aligned
{\frak q}_{(\Gamma,v)} :
\Omega(X) \otimes 
BCF^{\text{\rm can}}(\cL;\vec{\kappa}(\Gamma)) 
\to CF^{\text{\rm can}}(L_{\kappa(\Gamma)_0},L_{\kappa(\Gamma)_N};\F)
\otimes \Lambda_0.
\endaligned
\end{equation}
We will define\index[syindex]{gcangammav@${\frak g}_{(\Gamma,v)}$} 
\begin{equation}
\aligned
{\frak g}_{(\Gamma,v)} &:
\Omega(X) \otimes 
BCF^{\text{\rm can}}(\cL;\vec{\kappa}(\Gamma))
\to CF(L_{\kappa(\Gamma)_0},L_{\kappa(\Gamma)_N};\F)
\otimes \Lambda_0
\endaligned
\end{equation}
at the same time. %\marginpar{Removed $+$ in the above two formula, rhs.  KF 2025 June.}
(Here $\vec\kappa(\Gamma)$ is defined in (\ref{veckapadef}).
$N+1$ is the number of exterior edges of $\Gamma$.)
\begin{rem} The maps
${\frak q}_{(\Gamma,v)}$ and 
${\frak g}_{(\Gamma,v)}$ 
are the analogs of 
${\frak m}_{\Gamma}$ and 
${\frak f}_{\Gamma}$ 
in (\ref{257+1}), respectively, 
where a single differential form on the ambient space is 
included.
There is a version involving arbitrarily many 
differential forms on the ambient space. See the proof of \cite[Theorem 5.4.1]{fooo09}.
\end{rem}
Let $v'$ be the vertex that is contained in $e^{\text{\rm ext}}_0$ 
and is different from the root.
We remove $e^{\text{\rm ext}}_0$ from $\Gamma$.
The closures of the connected components of the complement 
consist of $k$ ribbon trees $\Gamma_i$.
($i=1,\dots,k$.)  Here $k+1$ is the number of edges containing 
$v'$. We regard $v'$ as the root
of $\Gamma_i$. With the other data induced from $\Gamma$ 
in an obvious way, $\Gamma_i$ becomes a 
decorated ribbon tree.
We consider the following two cases separately.
Let $g \in \Omega(X)$.
\par\medskip
\noindent(Case 1): 
$v' = v$. In this case we put : 
\begin{equation}\label{form14181418}
\aligned
{\frak q}_{(\Gamma,v)}(g)
&= \Pi'_{L_{\kappa}} \circ {\frak q}_{\vec{\kappa}(v'),B(v')}^
{\rm f.u.\bf b}(g)
\circ \bigotimes_{i=1}^k \frak f_{\Gamma_i}, \\
{\frak g}_{(\Gamma,v)}(g)
&= G'_{L_{\kappa}} \circ {\frak q}_{\vec{\kappa}(v'),B(v')}^{\rm f.u.\bf b}(g)
\circ \bigotimes_{i=1}^k \frak f_{\Gamma_i},
\endaligned
\end{equation}
if $L_{\kappa(\Gamma)_{\blue{K}}} = L_{\kappa(\Gamma)_{\blue{0}}} = L_{\kappa}$,
and 
\begin{equation}\label{form14191419}
{\frak q}_{(\Gamma,v)}(g)
={\frak q}_{\vec{\kappa}(v'),B(v')}^{\rm f.u.\bf b}(g) 
\circ \bigotimes_{i=1}^k \frak f_{\Gamma_i},
\qquad 
\widehat{\frak g}_{(\Gamma,v)}(g) = 0,
\end{equation}
if $L_{\kappa(\Gamma)_{\blue{K}}} \ne L_{\kappa(\Gamma)_{\blue{0}}}$.
\par
Here  $\Pi'_{L_{\kappa}}$ is defined during the proof of Proposition \ref{prop115}.
$G'_{L_{\kappa}}$ is defined also there. 
More precisely $G'_{L_{\kappa}}$ does not make sense sometimes 
but the composition 
$G'_{L_{\kappa}} \circ {\frak q}_{\vec{\kappa}(v'),B(v')}^{\rm f.u.\bf b}$ 
is defined as explained there.
\par
The operators ${\frak q}_{\vec{\kappa},B}^{\rm f.u.\bf b}$\index[syindex]{qformfubdfellnohatkappaB@$\frak q^{\rm f.u.\bf b}_{\vec{\kappa},B}$} are obtained by decomposing ${\frak q}_{1}^{\rm f.u.\bf b}$ 
in (\ref{map8322}) in an obvious way. %\marginpar{A sentence added.  KF 2025 Aug}
\footnote{We use $B \in \Pi_2(\vec{\kappa},\vec p)$
rather than $E \in G$.  This is because a decorated ribbon tree associates an element $B$ rather than $E$ to each vertices.}
\par\medskip
\noindent(Case 2): 
$v' \ne v$. 
In this case, 
there exists unique $j$ such that $v$ is a vertex of $\Gamma_j$.
We put
\begin{equation}\label{1421new}
\aligned
{\frak q}_{(\Gamma,v)}(g)
&= \Pi'_{L_{\kappa}} \circ \frak m_{\vec{\kappa}(v'),B(v')}^{\rm f.u.\bf b} 
\circ 
\left(\bigotimes_{i=1}^{j-1} \frak f_{\Gamma_i} \otimes {\frak g}
_{(\Gamma_j,v)} (g) 
\otimes \bigotimes_{i=j+1}^{k} \frak f_{\Gamma_i}\right) , \\
{\frak g}_{(\Gamma,v)}(g)
&= G'_{L_{\kappa}} \circ \frak m_{\vec{\kappa}(v'),B(v')}^{\rm f.u.\bf b} 
\circ \left(\bigotimes_{i=1}^{j-1} \frak f_{\Gamma_i} \otimes {\frak g}_{(\Gamma_j,v)} (g) 
\otimes \bigotimes_{i=j+1}^{k} \frak f_{\Gamma_i}\right),
\endaligned
\end{equation} %\marginpar{Sign put.  KF. 2024 Dec. It is actually always $+$}

if $L_{\kappa(\Gamma)_{\blue{K}}} = L_{\kappa(\Gamma)_{\blue{0}}} = L_{\kappa}$,
and 
\begin{equation}\label{1422new}
\aligned
{\frak q}_{(\Gamma,v)}(g)
&= \frak m_{\vec{\kappa}(v'),B(v')}^{\rm f.u.\bf b} \circ\left(\bigotimes_{i=1}^{j-1} \frak f_{\Gamma_i} \otimes {\frak g}_{(\Gamma_j,v)} (g)
\otimes \bigotimes_{i=j+1}^{k} \frak f_{\Gamma_i}\right), \\
{\frak g}_{(\Gamma,v)}(g) &= 0,
\endaligned
\end{equation}
if $L_{\kappa(\Gamma)_{\blue{K}}} \ne L_{\kappa(\Gamma)_{\blue{0}}}$.
\par
We remark that the tensor products of homomorphisms here are `super tensor product'.
Namely
$$
(\psi \otimes \phi)(x \otimes y) = (-1)^{\deg \phi \deg' x} \psi(x) \otimes \phi(y).
$$
The operators $\frak m_{\vec{\kappa},B}^{\rm f.u.\bf b}$\index[syindex]{mxfubbdkkappaB@$\frak m^{\text{\rm f.u.}\text{\bf b}}_{\vec{\kappa},B}$}
 is obtained by decomposing $\frak m_{k}^{\rm f.u.\bf b}$
in Definition \ref{defn105} in an obvious way. %\marginpar{A sentence added.  KF 2025 Aug.}
\par\smallskip
We now put:
\begin{equation}
{\frak q}^{\text{\bf b}}_{\vec{\kappa}} = 
\sum_{(\Gamma,v); \vec{\kappa} = \vec{\kappa}(\Gamma)} 
T^{\omega(\Gamma)}\rho_{b}(\Gamma) \rho_{\frak b,\theta}(\Gamma)
{\frak q}^{\text{\bf b}}_{(\Gamma,v)}.
\end{equation}
\begin{prop}\label{lem1417} The map
$\frak q^{\text{\bf b}}_{\vec{\kappa}}$ induces a cochain map\index[syindex]{qcanbbd@${\wq^{\text{\bf b}}}$}
$$
\wq^{\text{\bf b}} : \Omega(X) \to CH^*(\cL, \cL).
$$
\end{prop}
The proof is a straightforward calculation which is 
similar to the proof of \cite[Lemma 5.4.38]{fooo09}, and so is omitted.
 %\marginpar{Do we need to write the proof ?}

\begin{rem}
In the case of the filtered $A_{\infty}$ category $\cL$
(which we obtained by taking the canonical model), the space of morphisms is 
always a {\it finitely generated} free module over $\Lambda_0$.
In particular we do not need to take a completion with respect to the 
$T$ adic topology for the tensor products among them. 
\end{rem}

\begin{rem}
In the situation of Proposition \ref{sec2relative},
we can prove a similar relative version of the construction of $\widehat{\frak q}$.
Namely when $\widehat{\frak q}$ is defined with its value in the Hochschild cohomologies 
of $\cL$ and  $\mathcal U$ respectively, we can extend it to the Hochschild cohomology of the union.
The same remark applies for $\widehat{\frak p}$ which we will construct in the next section.
\end{rem}

For the purpose of using it later in Subsection \ref{subseqprodcan}, we rewrite the above discussion a bit.
 %\marginpar{The rest of this subsection is added. KF. 2024 Nov}
For $g \in \Omega(X)$ we have
${\frak q}^{{\rm f.u.} {\bf b}}_{\vec{\kappa},B}(g)$ for various $\vec{\kappa}$, $B$
which gives an element of $ CH^*(\cL^{\rm form}_{\rm uni}, \cL^{\rm form}_{\rm uni})$.
We denote it by $\widehat{\frak q}^{{\rm f.u.} {\bf b}}(g)$.

\begin{figure}[h]
\includegraphics[scale=0.3]{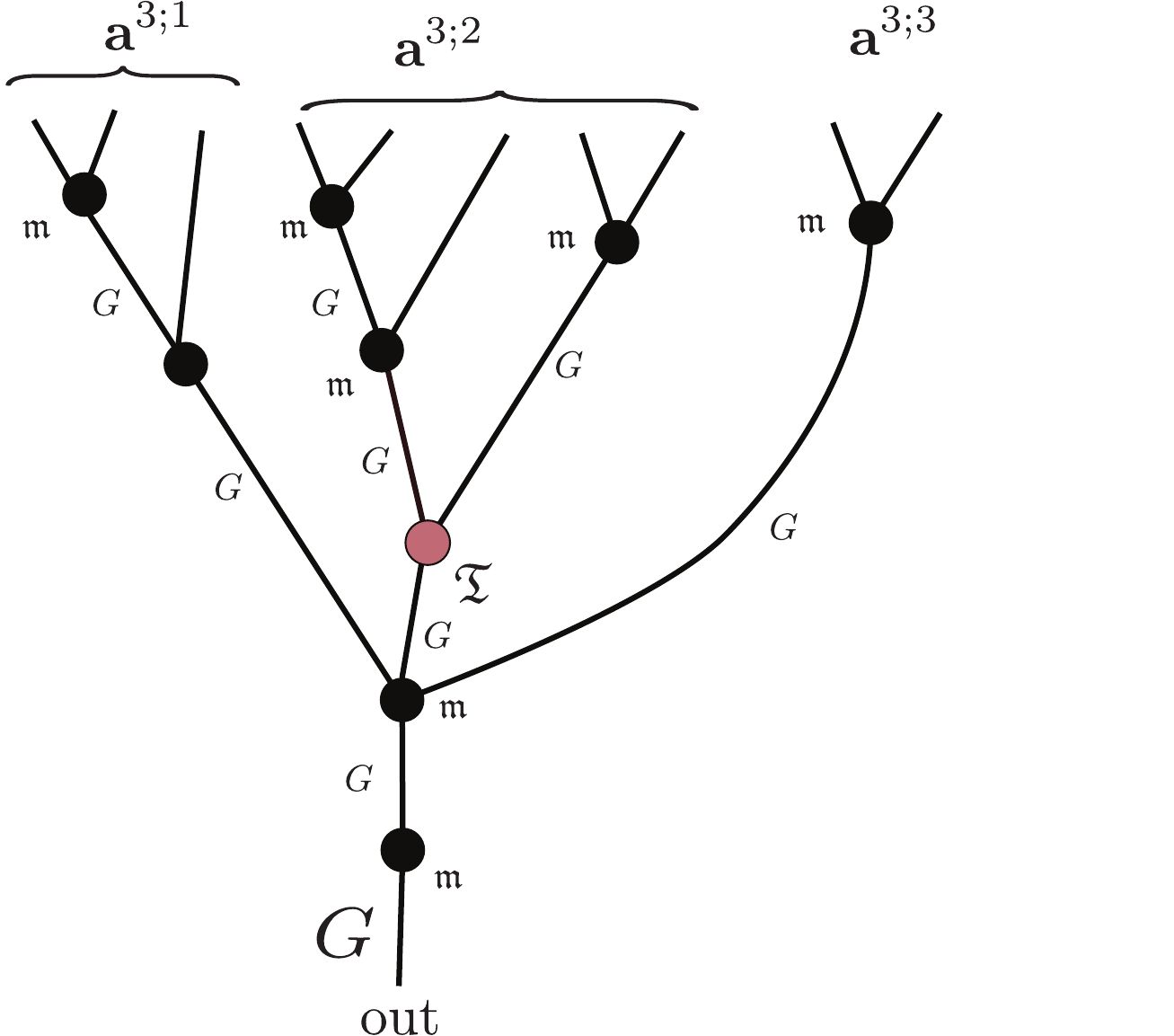}
\caption{$\frak H_{(\Gamma,v)}(\frak T)$}
\label{Figure14-17}
\end{figure}

\begin{lem}\label{lem142020}
There exists a chain map\index[syindex]{Jfrak@$\frak J$}
$$
\frak J: CH^*(\cL^{\rm form}_{\rm uni}, \cL^{\rm form}_{\rm uni}) \to  CH^*(\cL, \cL)
$$
such that $\frak J(\widehat{\frak q}^{{\rm f.u.} {\bf b}}(g)) = \wq^{\text{\bf b}}(g)$.
\end{lem}
\begin{proof}
Let $\frak T \in CH^*(\cL^{\rm form}_{\rm uni}, \cL^{\rm form}_{\rm uni})$. By replacing $G$, if necessary, we may assume 
that $\frak T$ is $G$-gapped. We write $\frak T$ as totality of $\frak T_{\vec{\kappa},E}$.
\par
In Definition \ref{decribbon} we assigned $B(v)$ to each vertex of a metric ribbon tree.
We change the definition slightly and assign $E(v)$ (an element of our discrete monoid $G$ and is a non-negative real 
number).  Then the discussions there and in the previous parts of Subsection \ref{frakqoncoho} work without change.
 %\marginpar{A sentence added.  KF. 2025 Aug.}
\par
We now replace ${\frak q}^{{\rm f.u.} {\bf b}}_{\vec{\kappa},E}(g): = \sum_{B, \omega(B) = E} {\frak q}^{{\rm f.u.} {\bf b}}_{\vec{\kappa},B}(g)$ 
by $\frak T_{\vec{\kappa},E}$ in Formulas (\ref{form14181418}), (\ref{form14191419}), (\ref{1421new}) (\ref{1422new}) 
(their versions with $B(v)$ replaced by $E(v)$).
Then we obtain $\frak J$.
\end{proof}
Note that we obtained ${\frak g}_{(\Gamma,v)}(g)$ also.
Replacing 
${\frak q}^{{\rm f.u.} {\bf b}}_{\vec{\kappa},E}(g)$ 
by $\frak T_{\vec{\kappa},E}$ 
we define 
$
\frak H_{(\Gamma,v)}(\frak T).
$
Its weighted sum over $(\Gamma,v)$ is, by definition, 
$\frak H(\frak T)$.
See Figure \ref{Figure14-17}.
We denote by $\widehat{\frak f} : BCF(\cL)
\to BCF(\cL^{\text{\rm form}}_{\rm uni})$ the coalgebra homomorphism induced by ${\frak f}$.
(See \cite[(3.2.6)]{fooo09}.)

\begin{lem}\label{lem1421}
We have
 %\marginpar{need a sign. Kozule KF 2024 Nov  I put some. need to be checked. KF 2024 Dec.
%Error in the formula is corrected.  KF 2025 Jan}
$$
\aligned
\frak H(\delta^H(\frak T))({\bf a}) 
= & \frak J(\frak T)({\bf a}) \\
&+ \sum_c (-1)^{\deg \frak T +\deg'{\bf a}_{c}^{3;1}+1}\frak H(\frak T)({\bf a}_{c}^{3;1},
\frak m^{\rm can}({\bf a}_{c}^{3;2}),{\bf a}_{c}^{3;3}) \\
&+
\sum_{c} (-1)^{\deg \frak T\deg'{\bf a}_{c}^{3;1}}\widehat{\frak f}({\bf a}_{c}^{3;1},
\frak J(\frak T)({\bf a}_{c}^{3;2}),{\bf a}_{c}^{3;3}) \\
&+ \sum_{c}(-1)^{(\deg \frak T +1)(1+\deg'{\bf a}_{c}^{3;1})}
\frak m^{{\rm f.u.} {\bf b}}
(\widehat{\frak f}({\bf a}_{c}^{3;1}),
\frak H(\frak T)({\bf a}_{c}^{3;2}),
\widehat{\frak f}({\bf a}_{c}^{3;3})) \\
& -  \frak T
(\widehat{\frak f}({\bf a})).
\endaligned
$$
\end{lem} %\marginpar{The sign in the third term of RHS is corrected.  KF 2025Feb}
Here $((\Delta \otimes 1) \circ \Delta)({\bf a}) = \sum_c
{\bf a}_{c}^{3;1} \otimes {\bf a}_{c}^{3;2} \otimes {\bf a}_{c}^{3;3}$.
\begin{proof}
The proof of Lemma \ref{lem1421} is also similar to 
Lemma \ref{lem1417}.  %\marginpar{Proof written.  KF 2025 Jan.}
See Figures \ref{Figure14-18}-\ref{Figure14-19}. %\marginpar{Figures added.}

\begin{figure}[h]
\includegraphics[scale=0.3]{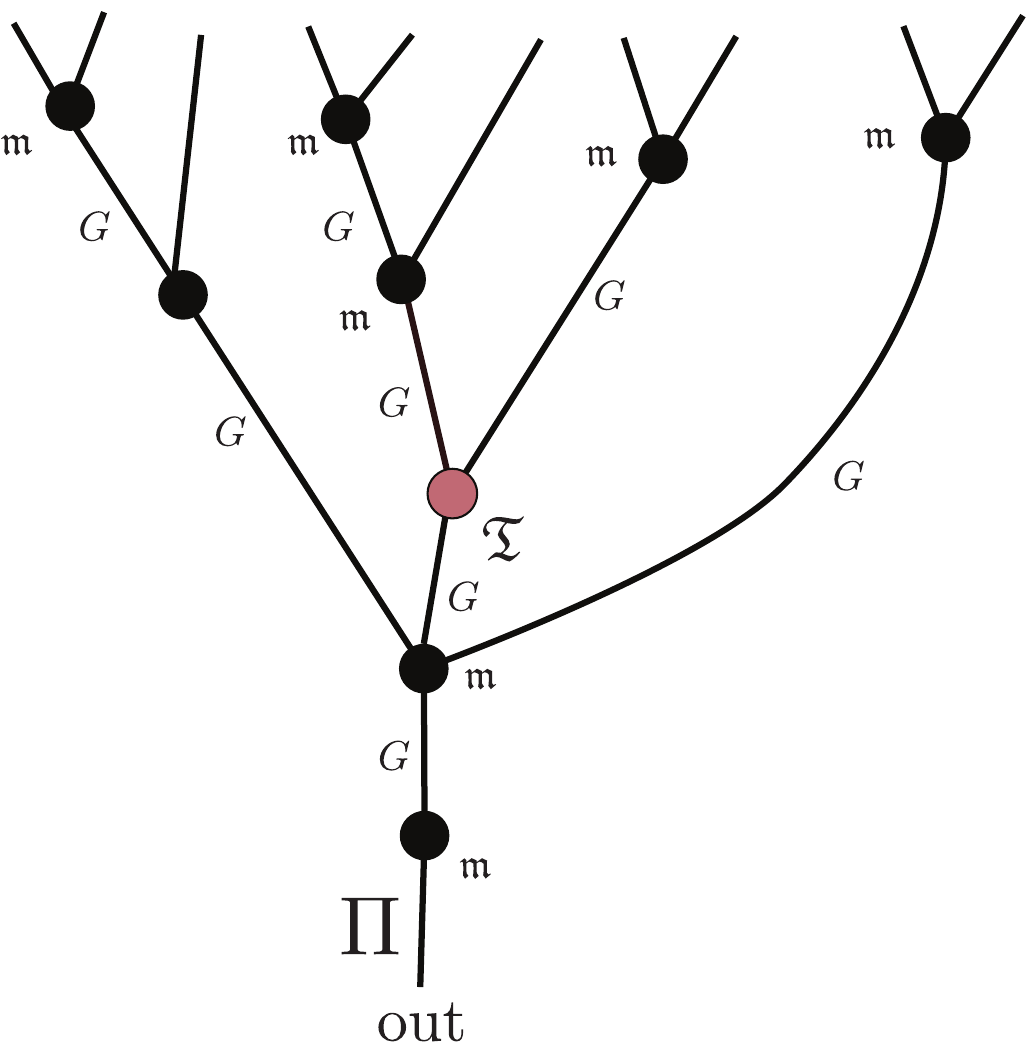}
\caption{$\frak J_{(\Gamma,v)}(\frak T)$}
\label{Figure14-18}
\end{figure}

\begin{figure}[h]
\includegraphics[scale=0.25]{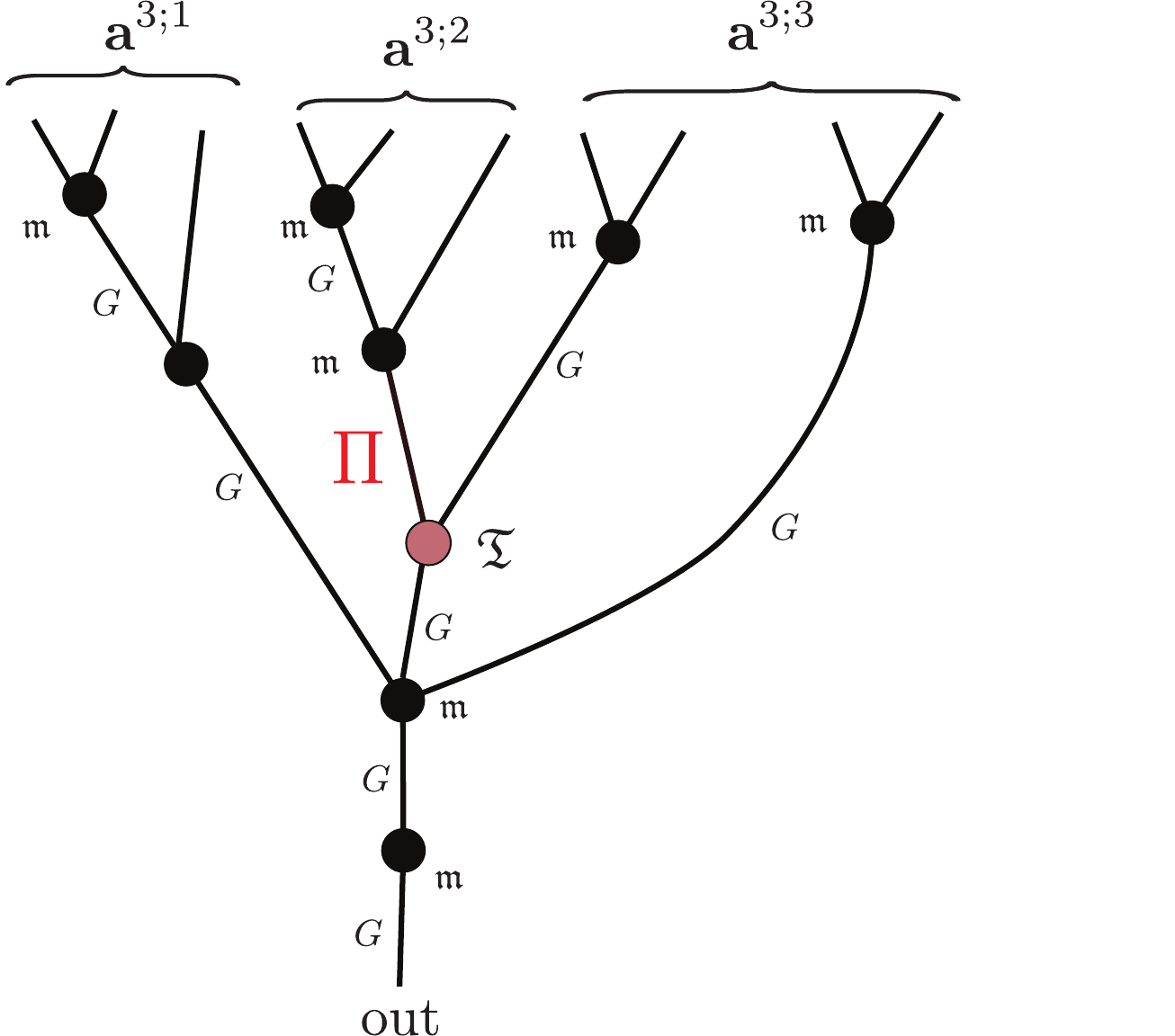}
\caption{2nd term}
\label{Figure14-20}
\end{figure}

\begin{figure}[h]
\includegraphics[scale=0.45]{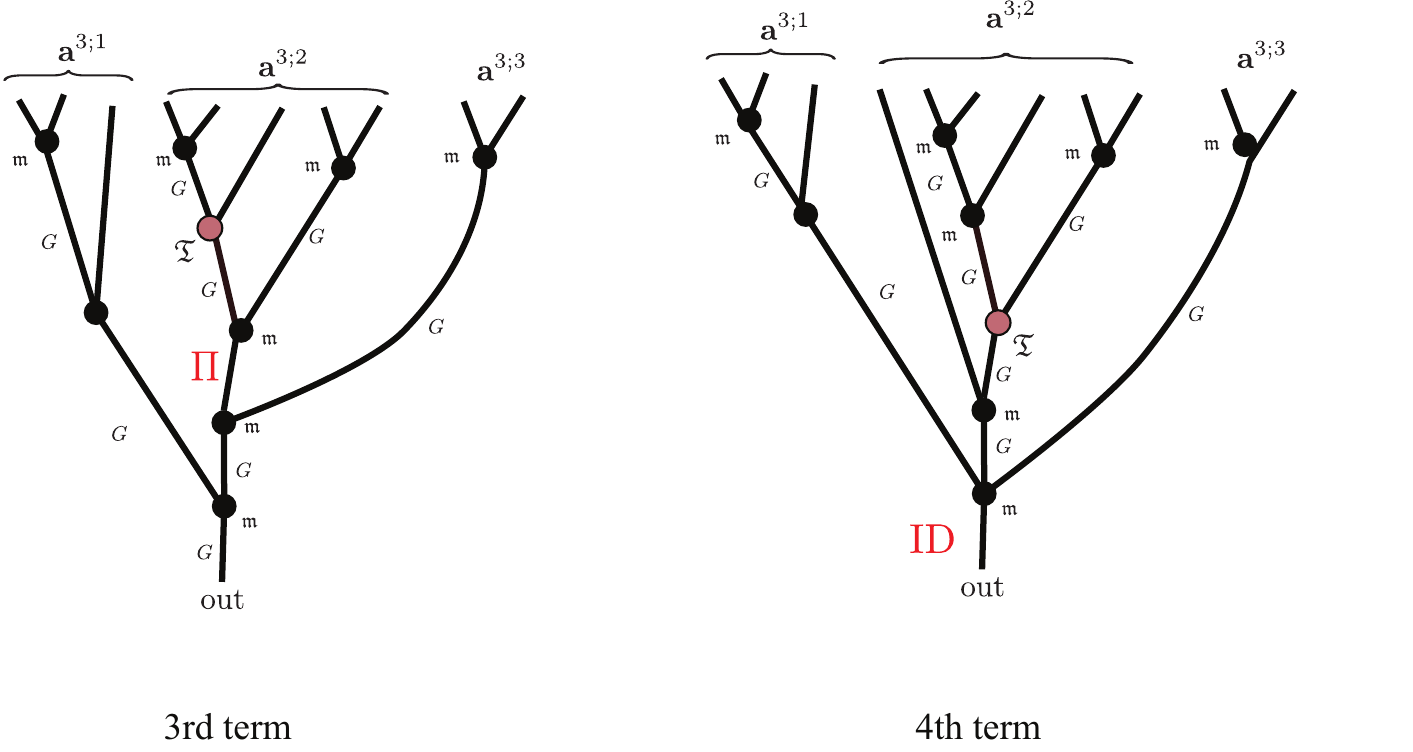}
\caption{3rd and 4th terms}
\label{Figure14-19}
\end{figure}

The idea of the proof is to study the commutator:
$$
[d,\frak H_{(\Gamma,v)}(\frak T)]
$$
of $\frak H_{(\Gamma,v)}(\frak T)$ with the de Rham differential $d$.  Recalling that  $\frak H_{(\Gamma,v)}(\frak T)$
is a composition of operators assigned to the vertices and the edges of $\Gamma$,
this commutator is expressed as a sum of the terms arising from the replacement of various operators 
assigned to the vertices and edges of $\Gamma$ for $\frak H(\frak T)$
by their commutators with $d$.
Taking commutators with those assigned to the vertices gives 
$$
[d,\frak m]  \quad {\rm or} \quad [d,\frak T].
$$
and those assigned to the edges gives 
$$
[d,G] = {\rm id} - \Pi.
$$
respectively. Note $[d,\frak m] = \frak m \circ \widehat{\frak m}$ by the $A_{\infty}$ relation.
Therefore the terms replacing $\frak m$ by $[d,\frak m]$ cancel the terms arising from the replacement of $G$ by ${\rm id}$.
\par
On the other hand,\footnote{Here we denote 
$\widehat{\frak T}(x_1,\dots,x_n) = \sum_i (-1)^{\deg\frak T(\deg'x_1 + \dots + \deg'x_{i-1})}x_1 \otimes \dots  \otimes\frak T(x_{i},\dots)
\otimes \dots \otimes x_n$.}
$$ 
[d,\frak T] + \frak m \circ \widehat{\frak T} + \frak T \circ \widehat{\frak m} = \delta_H(\frak T).
$$
The second and the third terms of the left hand side, cancel with terms replacing $G$ by ${\rm id}$ 
for edges adjacent to the vertex assigned with $\frak T$.
\par
The terms arising from the replacement of $G$ by $\Pi$ for edges become 
the second and third terms respectively of the right hand side of the formula to be proven. 
\par
There is an exceptional case, that is, the case when we replace $G$ with $\Pi$ or  ${\rm id}$
for the edge adjacent to the root.
They will give $\frak J_{(\Gamma,v)}(\frak T)$ (the first term) and the fourth term 
of the right hand side, respectively.
\par
Thus we have proved the required formula modulo the sign.
The fact the equality holds with sign is a consequence of the fact 
that all the signs here are by the  Koszul rule.
\end{proof}
\begin{rem}\label{rem425}
We can also define 
$$
\frak J_{\ell}: CH^*(\cL^{\rm form}_{\rm uni}, \cL^{\rm form}_{\rm uni})^{\otimes \ell} \to  CH^*(\cL, \cL)
$$
for $\ell \ge 2$
so that, together with $\frak J= \frak J_{1}$, the operations $\{\frak J_{\ell}\}_{\ell =1}^\infty$  define an $L_{\infty}$ homomorphism (between differential 
graded Lie algebras). See \cite[Theorem 7.4.41]{fooo092}.
\end{rem}
 %\marginpar{The rest of this subsection 
%is added in 2025 March.}

Using Remark \ref{rem425} and Proposition \ref{prop142} we have:
\begin{cor}\label{Cor1526}
We have a filtered $L_{\infty}$ homomorphism:
$$
\Omega(X)^{\otimes \ell} \to CH^*(\cL, \cL).
$$
\end{cor}
We can include the unit ${\bf e}^+$ and show the next lemma easily.
\begin{lem}
There exists
$$
\frak J_{\rm red}: CH_{\rm red}^*(\cL^{\rm form}_{\rm uni}, \cL^{\rm form}_{\rm uni}) \to  CH_{\rm red}^*(\cL, \cL)
$$
such that the next diagram commutes:
$$
\begin{CD}
CH_{\rm red}^*(\cL^{\rm form}_{\rm uni}, \cL^{\rm form}_{\rm uni}) @>\frak J_{\rm red} >>  CH_{\rm red}^*(\cL, \cL)
\\ 
@V VV   @ V VV   \\
CH^*(\cL^{\rm form}_{\rm uni}, \cL^{\rm form}_{\rm uni})  @>\frak J>> CH^*(\cL, \cL)
\end{CD}
$$
\end{lem}
Therefore Proposition \ref{prop1416} implies:
\begin{prop}\label{lem1425}
The map $
\wq^{\text{\bf b}} : \Omega(X) \to CH^*(\cL, \cL)
$  factors into
$$
\Omega(X) \to CH_{\rm red}^*(\cL, \cL) \to CH^*(\cL, \cL).
$$
\end{prop}

\begin{rem}
After we decompose $\cL$ into the sum of $\cL_{\lambda}$ with  critical values 
$\lambda$ we redefine $\frak m_0(1)$ to be zero.
Let us denote (in the rest of this remark)  $\cL_{\lambda,1}$ (resp.  $\cL_{\lambda,2}$)
the one with $\frak m_0(1) = \lambda {\bf e}^+$ (resp. $\frak m_0(1) = 0$).
The Hochschild coboundary operator of $CH^*(\cL_{\lambda,1}, \cL_{\lambda,1})$ 
does {\it not} coincide with the one of $CH^*(\cL_{\lambda,2}, \cL_{\lambda,2})$.
On the other hand it is easy to see that Hochschild coboundary operator of $CH_{\rm red}^*(\cL_{\lambda,1}, \cL_{\lambda,1})$ 
does coincide with that of $CH_{\rm red}^*(\cL_{\lambda,2}, \cL_{\lambda,2})$.
Therefore Proposition \ref{prop528} implies the canonical isomorphism:
\begin{equation}\label{eq1425}
HH^*(\cL_{\lambda,1} , \cL_{\lambda,1})\otimes_{\Lambda_0} \Lambda
\cong
HH^*(\cL_{\lambda,2}, \cL_{\lambda,2}) \otimes_{\Lambda_0} \Lambda.
\end{equation}
Proposition \ref{lem1425} implies that the isomorphism (\ref{eq1425}) is 
compatible with $\wq^{\text{\bf b}}$.
\end{rem}

\section{The open-closed map $\widehat{\frak p}$.}
\label{sec:open-closed map}

\subsection{Basic idea of the construction}
\label{frakponcoho}

In this section we define the {\it open-closed map}\index{open-closed map} $\widehat{\frak p}$.
 %\marginpar{This subsection is much modified in 2025 March. KF}
We use the same moduli space as the one we used in the definition of $\widehat{\frak q}$.
We however regard one of the {\it interior} marked points 
as the output. So we need to 
assume that the evaluation map at this marked point 
is a \blue{weak} submersion. By this reason, we need to use a slightly 
different perturbation from those used for the definition of  $\widehat{\frak q}$.
\par
We summarize the properties of the map $\widehat{\frak p}$ we obtain in the next theorem.
\begin{thm}\label{themp}
Consider the category $\cL$ associated to the choices of $\mathscr L$ and $\bf b$.
Then there exists a chain map\index[syindex]{pcanbbf@$\widehat{\frak p}^{{\bf b}}$} 
$$
\widehat{\frak p}^{{\bf b}} : CH_*(\cL, \cL) \to \Omega(X) \hat\otimes_{\F} {\Lambda_0}
$$
with the following properties.
\begin{enumerate}
\item 
We denote by
${\frak p}^{{\bf b}}_1 : \cL((L_{\kappa},b_{\kappa}),(L_{\kappa},b_{\kappa})) \to \Omega(X) \hat\otimes_{\F} {\Lambda_0}$
the restriction of $\widehat{\frak p}^{{\bf b}}$. It then satisfies
\begin{equation}\label{formnew161}
[{\frak p}^{{\bf b}}_1(x)] \equiv (-1)^{\deg' x}i_!(x)  \mod \Lambda_+.
\end{equation}
Here $x \in H(\tilde L_{\kappa}) \subseteq \cL((L_{\kappa},b_{\kappa}),(L_{\kappa},b_{\kappa}))$
and $i_!(x)$ is a differential form representing  the class $[i_!(x)] \in H^*(X)$ that is 
the class to which the Gysin homomorphism sends $x$.
\item 
$\widehat{\frak p}^{{\bf b}}(x_0\otimes {\bf x}) = 0$
if ${\bf x} = {\bf x}_1 \otimes {\bf e} \otimes {\bf x}_2$.
\item
$$
\widehat{\frak p}^{{\bf b}}(x_0,\dots,x_k) = (-1)^{\maltese} \widehat{\frak p}^{{\bf b}}(x_k,x_0,\dots,x_{k-1})
$$
with $\maltese = \deg'x_k (\deg'x_0 + \dots + \deg' x_{k-1})$.
\end{enumerate}
\end{thm}
\begin{rem}
$  $ \par
\begin{enumerate}
\item 
Theorem \ref{themp} (1)  is a de Rham version of \cite[Theorem 3.8.9 (3.8.10.1)]{fooo092}.
In \cite{fooo092} we used the singular homology. So the inclusion $i : L_{\kappa} \to X$ 
sends a singular chain to a singular chain.
Here we need a smoothing process to define the Gysin map in the chain level.
\par
The operator $i_!$ in (\ref{formnew161}) is defined as a special case of the integration along the fiber 
(of spaces with Kuranishi structures and CF-perturbations).  Its sign is defined as a special case of
one defined in \cite{springer}. So it is different from the sign appearing in \cite[Theorem 3.8.9 (3.8.10.1)]{fooo092}. 
\item 
Theorem \ref{themp} (2) corresponds to \cite[Theorem 3.8.9 (3.8.10.4)]{fooo092}.
In the current situation it implies that there exists a map
$$
\widehat{\frak p}^{{\bf b}}_{\rm red} : CH^{\rm red}_*(\cL, \cL) \to \Omega(X) \hat\otimes_{\F} {\Lambda_0}
$$
such that $\widehat{\frak p}^{{\bf b}}$ factors into 
$$
CH_*(\cL, \cL) \longrightarrow CH^{\rm red}_*(\cL, \cL) \overset{\widehat{\frak p}^{{\bf b}}_{\rm red}}\longrightarrow \Omega(X) \hat\otimes_{\F} {\Lambda_0}.
$$
\item The formula \cite[Theorem 3.8.9 (3.8.10.6)]{fooo092} is an error as mentioned in \cite[Remark 3.1.7]{toric3}.
Namely the right hand side of \cite[Theorem 3.8.9 (3.8.10.6)]{fooo092} should be $0$. 
This correction is incorporated into Theorem \ref{themp} (2).
\item
In \cite[Theorem 3.8.9]{fooo092} the map $\frak p$ is defined on the cyclic chain complex. 
So the statement corresponding to Theorem \ref{themp} (3) is included.
\item
See Remark \ref{rem1527} for a variant of  \cite[Theorem 3.8.9 (3.8.10.3)]{fooo092}.
\item
The degree of $\frak p$ is given by
$
\deg \frak p(S) = \frak{deg}(S) + n +1,
$
or $\deg \frak p(S) = \deg (S) + n$. See (\ref{form464new}) and footnotes 20 and 21  %\marginpar{Check the number of footnote. Checked 2026 Feb}
for $\frak{deg}$ and $\deg$. %\marginpar{(6) is added.  KF 2025 Aug.}
\end{enumerate}
\end{rem}
\par
Consider the 
moduli space $\mathcal M_{\ell+1;\vec k}((\vec{\kappa},\vec p);B)$.
Among the $\ell+1$ interior marked points 
we treat the $0$-th one $z_0^+$ in a different way 
from the other  $\ell$ interior marked points $z_i^+$, $i=1,\dots,\ell$.
We need to modify Proposition \ref{Kuraeistspoly} appropriately.
To include the generators ${\bf e}^+$ and ${\bf f}$ we proceed in the same way as Section \ref{sec:unit} 
as follows.  We include extra label $\vec{\frak f} = (\frak f_{\bf e},\frak f_{\bf f})$ as in Section \ref{sec:unit}.
We then modify Proposition \ref{prop94} also.

Another important modification we need to make in the argument of Section 9 is the following.
We need to consider 
$
\frak p({\bf e}^+)
$, which corresponds to the case when we have one boundary marked point that is 
forgetable. We then need to study the moduli space which has no boundary marked point.
We denote it by $\mathcal M_{\ell+1;\emptyset}(\kappa;B)$.
Here we require that the boundary of the domain will be mapped to $L_{\kappa}$.
An important case we need to study for such a moduli space
is the moduli space of constant maps with one interior marked point.
The group of automorphisms of an element of such a moduli space is $S^1$ which has positive dimension.
We however need to perturb it since otherwise the evaluation map 
at the 0-th  interior marked point would not be (weakly) submersive.
We will explain this case in detail later.

\begin{prop}\label{Kuraeistspoly2}
$  $ \par
\begin{enumerate}
\item
The space $[0,1]^{\vert \frak f_{\bf f}\vert} \times {\mathcal M}_{\ell+1;\vec k}((\vec{\kappa},\vec p);B;\vec{\frak f})^{\boxplus 1}$  has an
oriented Kuranishi structure, $\frak p$-Kuranishi structure.\index{pKuranishistructure@$\frak p$-Kuranishi structure}
(We write $([0,1]^{\vert \frak f_{\bf f}\vert} \times{\mathcal M}_{\ell+1;\vec k}((\vec{\kappa},\vec p);B)^{\boxplus 1})^{\frak p}$ when we equip $\frak p$-Kuranishi structure to this space.)
\item 
The boundary of $([0,1]^{\vert \frak f_{\bf f}\vert} \times{\mathcal M}_{\ell+1;\vec k}((\vec{\kappa},\vec p);B)^{\boxplus 1})^{\frak p}$
is a union of the following five types of fiber or direct products:
\par
{\bf (Boundary of type 1)}
We use the notation of $(\ref{218})$ and $(\ref{218+I})$.
\begin{equation}\label{2182}
\aligned
&([0,1]^{\vert \frak f_{\bf f}^2\vert} \times{\mathcal M}_{\ell''+1; k''_i+1}(L;\beta;\vec{\frak f}^2)^{\boxplus 1})^{\frak p}\\
&\,{}_{\text{\rm ev}_0} \times_{\text{\rm ev}_{i,j}}
(([0,1]^{\vert \frak f_{\bf f}^1\vert} \times{\mathcal M}_{\ell';\vec k'}((\vec{\kappa},\vec p);B';\vec{\frak f}^1)^{\boxplus 1})^{\frak q}
\endaligned
\end{equation}
(We put $\frak p$-Kuranishi structure to the first factor and 
$\frak q$-Kuranishi structure to the second factor.)
$\vec{\frak f}^1$ and $\vec{\frak f}^2$ are defined in the same way as Section $\ref{sec:unit}$ using {\rm (tricom1),(tricom2),
(tricom3)}.
\begin{equation}\label{2183}
\aligned
&([0,1]^{\vert \frak f_{\bf f}^2\vert} \times{\mathcal M}_{\ell''; k''_i+1}(L;\beta;\vec{\frak f}^2)^{\boxplus 1})^{\frak q} \\
&\,{}_{\text{\rm ev}_0} \times_{\text{\rm ev}_{i,j}}
([0,1]^{\vert \frak f_{\bf f}^1\vert} \times{\mathcal M}_{\ell'+1;\vec k'}((\vec{\kappa},\vec p);B';\vec{\frak f}^1)^{\boxplus 1})^{\frak p}
\endaligned
\end{equation}
In $(\ref{2182})$, $(\ref{2183})$, the union is taken over the data described in $(\ref{218})$  and
$(\ref{218+I})$
together with the shuffles $(\mathbb L'',\mathbb L')$ of $\underline\ell$
such that $\# \mathbb L'' = \ell''$, $\# \mathbb L' = \ell'$.
\blue{In \eqref{2182} (resp. \eqref{2183}), $z_0^+$ is in the first (resp. second) factor.}  
\par
We put the $\frak p$-Kuranishi structure of Proposition $\ref{Kuraeistspoly2}$
on the first factor of $(\ref{2182})$ and on the second factor of $(\ref{2183})$.
We put the $\frak q$-Kuranishi structure of Propositions $\ref{Kuraeistspoly}$ and $\ref{prop94}$
on the second factor of $(\ref{2182})$ and on the first factor of $(\ref{2183})$.
\par
{\bf (Boundary of type 2)} 
We use the notation of $(\ref{219})$.
\begin{equation}\label{2191}
\aligned
&(([0,1]^{\vert \frak f_{\bf f}^2\vert} \times{\mathcal M}_{\ell''+1;\vec k''}((\vec{\kappa}'',\vec p'');B'';\vec{\frak f}^2)^{\boxplus 1})^{\frak p} \\
&\times
([0,1]^{\vert \frak f_{\bf f}^1\vert} \times{\mathcal M}_{\ell';\vec k'}((\vec{\kappa}',\vec p');B';\vec{\frak f}^1)^{\boxplus 1})^{\frak q}.
\endaligned
\end{equation}
\begin{equation}\label{2192}
\aligned
&([0,1]^{\vert \frak f_{\bf f}^2\vert} \times{\mathcal M}_{\ell'';\vec k''}((\vec{\kappa}'',\vec p'');B'';\vec{\frak f}^2)^{\boxplus 1})^{\frak q}\\
&\times
([0,1]^{\vert \frak f_{\bf f}^1\vert} \times{\mathcal M}_{\ell'+1;\vec k'}((\vec{\kappa}',\vec p');B';\vec{\frak f}^1)^{\boxplus 1})^{\frak p}.
\endaligned
\end{equation}
In $(\ref{2191})$, $(\ref{2192})$, the union is taken over the data described in $(\ref{219})$ 
together with the shuffles $(\mathbb L'',\mathbb L')$ of $\underline\ell$
such that $\# \mathbb L'' = \ell''$, $\# \mathbb L' = \ell'$.  
\blue{In \eqref{2191} (resp. \eqref{2192})), $z_0^+$ is in the first (resp. second) factor.}  
\par
We put the $\frak p$-Kuranishi structure of Proposition $\ref{Kuraeistspoly2}$
on the first factor of $(\ref{2191})$ and on the second factor of $(\ref{2192})$.
We put the $\frak q$-Kuranishi structure of Propositions $\ref{Kuraeistspoly}$ and $\ref{prop94}$ 
on the second factor of $(\ref{2191})$ and on the first factor of $(\ref{2192})$.
\par
{\bf (Boundary of type 3)} 
We use the notation of $(\ref{220})$.
\begin{equation}\label{2201}
\aligned
&([0,1]^{\vert \frak f_{\bf f}^2\vert} \times{\mathcal M}_{\ell''+1;\vec k''}((\vec{\kappa}'',\vec p'');B'');\vec{\frak f}^2)^{\boxplus 1})^{\frak p}
\\
&{}_{\text{\rm ev}_{1,m''}}\times_{\text{\rm ev}_{i,m'}}
([0,1]^{\vert \frak f_{\bf f}^1\vert} \times{\mathcal M}_{\ell';\vec k'}((\vec{\kappa}',\vec p');B';\vec{\frak f}^1)^{\boxplus 1})^{\frak q},
\endaligned
\end{equation}
\begin{equation}\label{22021}
\aligned
&([0,1]^{\vert \frak f_{\bf f}^2\vert} 
{\mathcal M}_{\ell'';\vec k''}((\vec{\kappa}'',\vec p'');B'';\vec{\frak f}^2)^{\boxplus 1})^{\frak q}
\\
&{}_{\text{\rm ev}_{1,m''}}\times_{\text{\rm ev}_{i,m'}}
([0,1]^{\vert \frak f_{\bf f}^1\vert} \times{\mathcal M}_{\ell'+1;\vec k'}((\vec{\kappa}',\vec p');B';\vec{\frak f}^1)^{\boxplus 1})^{\frak p},
\endaligned
\end{equation}
In $(\ref{2201})$, $(\ref{22021})$, the union is taken over the data described in $(\ref{220})$ 
together with the shuffles $(\mathbb L'',\mathbb L')$ of $\underline\ell$
such that $\# \mathbb L'' = \ell''$, $\# \mathbb L' = \ell'$.  
\blue{In \eqref{2201} (resp. \eqref{22021}), $z_0^+$ is in the first (resp. second) factor.}  
\par
We put the Kuranishi structure of Proposition $\ref{Kuraeistspoly2}$
to the first factor of $(\ref{2201})$ and the second factor of $(\ref{22021})$.
We put the Kuranishi structure of Proposition $\ref{Kuraeistspoly}$
to the second factor of $(\ref{2201})$ and the first factor of $(\ref{22021})$.
\par
{\bf (Boundary of type 4)} 
$\partial([0,1]^{\vert \frak f_{\bf f}\vert}) \times{\mathcal M}_{\ell+1;\vec k}((\vec{\kappa},\vec p);B)^{\boxplus 1}$.
\par
{\bf (Boundary of type 5)} 
In case $B\vert_{\partial} = 0 \in H_1(L_{\kappa})$
The space ${\mathcal M}_{\ell+1;\emptyset}(\kappa;B)^{\boxplus 1}$
has an extra boundary which is identified with the fiber product
\begin{equation}\label{eq157}
L \times_{X} \mathcal M^{\rm cl}_{\ell+2:B'}
\end{equation}
here $B' \in H_{2}(X;\Z)$ is a homology class which goes to $B \in H_2(X,L_{\kappa})$.
$\mathcal M^{cl}_{\ell+2:B'}$ is the compactified moduli space of 
pseudo-holomorphic sphere with $\ell+2$ marked point and of homology 
class $B'$.  See Definition \ref{defn155555} and discussions thereafter.
\item
The orientation is compatible with the identification 
of the boundary in $(2)$.
\item
The evaluation map $(\ref{evevev})$ extends to the 
compactification and is compatible with 
the description of the boundary in $(2)$.
\item
We have a forgetful map 
$$
\aligned
\frak{forget} :&([0,1]^{\vert \frak f_{\bf f}\vert} \times {\mathcal M}_{\ell+1;\vec k}((\vec{\kappa},\vec p);B;\vec{\frak f})^{\boxplus 1})^{\frak p} \\
&\to ([0,1]^{\vert \frak f_{\bf f}\vert} \times{\mathcal M}_{\ell+1;\vec k \setminus \frak f_{\bf e}}((\vec{\kappa},\vec p);B;(\emptyset,{\frak f}_{\bf f}))^{\boxplus 1})^{\frak p}
\endaligned
$$
which forgets marked points labeled by $\frak f_{\bf e}$.  The $\frak p$-Kuranishi structure 
is compatible with this forgetful map. 
We will describe the case when $\vec k \setminus \frak f_{\bf e} = \emptyset =  \frak f_{\bf f}$ later at 
right above Proposition $\ref{existmultipolu2}$.
\item
At the point when the coordinate $s_{i,j}$ of the  $[0,1]^{\vert \frak f_{\bf f}\vert}$ factor is $0$ or $1$, a 
similar compatibility condition as  Proposition $\ref{prop94}$ $(9),(10)$ is satisfied. %\marginpar{This part is modified.  KF 2025 March}
\item
The evaluation map
$$
\text{\rm ev}_{0}^+ 
 : ([0,1]^{\vert \frak f_{\bf f}\vert} \times{\mathcal M}_{\ell+1;\vec k}((\vec{\kappa},\vec p);B;\vec{\frak f})^{\boxplus 1})^{\frak p}
\to X
$$
at the $0$-th interior marked point is weakly submersive.
\item
The Kuranishi structure is invariant under the permutation of $\ell$
interior marked points $z^+_1,\dots,z^+_{\ell}$.
\item
The Kuranishi structure is invariant under the cyclic permutation 
of the data. %in the sense we described in Proposition $\ref{Kuraeistspoly}$. 
\end{enumerate}
\end{prop}
\begin{rem}\label{remnew164}
Note that we do not assume any submersivity for the evaluation map 
at the boundary marked points.  In other words, there is nothing similar to 
Proposition \ref{Kuraeistspoly} (5).  %\marginpar{Remark added.  KF 2025 Sep.}
By this reason it is simpler to state cyclic symmetry Proposition \ref{Kuraeistspoly2} (9).
Namely we simply require the symmetry under the cyclic permutation of the 
boundary marked points.  In the situation of Proposition \ref{prop92}, for example,
we require the compatibility with the forgetful map of the marked points labeled by $\frak f_{\bf e}$
and submersivity at the evaluation map of the unforgetable marked points.
So we need to state cyclic symmetry a bit more carefully. In fact it is 
impossible to require the weak submersivity of the evaluation map and the compatibility with the forgetful map
simultaneously for the same marked point.
\end{rem}

\begin{figure}[h]
\centering
\includegraphics[scale=0.2]{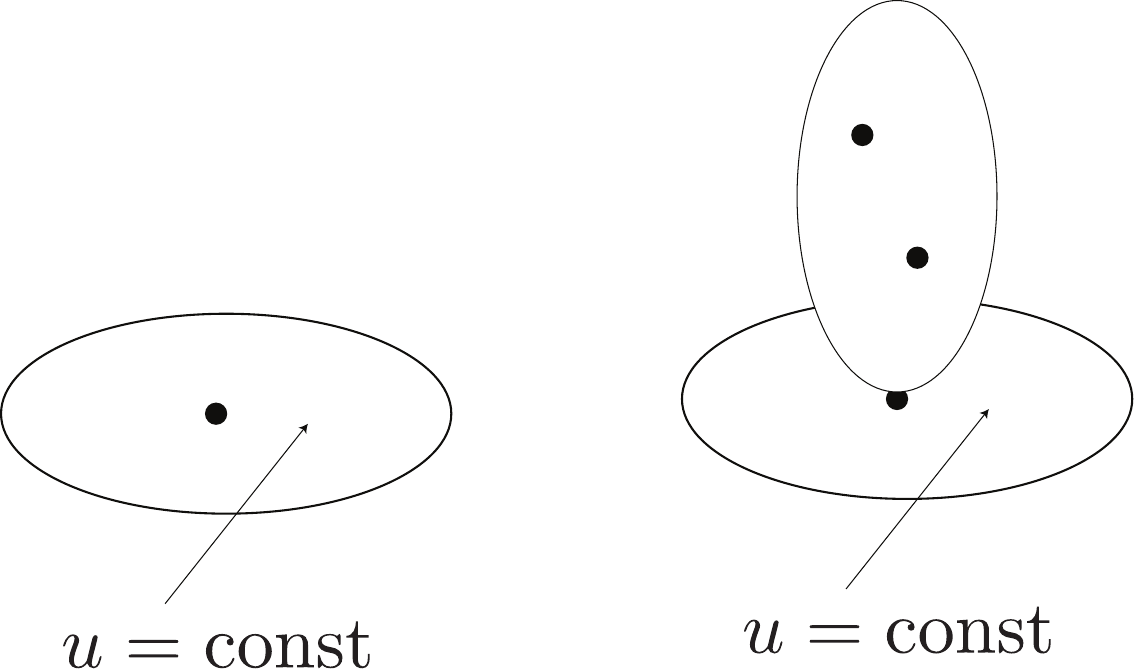}
\caption{An object with $S^1$ symmetry.}
\label{Figure15.6}
\end{figure}

We remark that Proposition \ref{Kuraeistspoly2} includes the case $K=\blue{0}$
(that is the case of one Lagrangian submanifold). 
Given the system of Kuranishi structures described in Proposition $\ref{Kuraeistspoly}$ 
the proof of Proposition $\ref{Kuraeistspoly2}$ 
proceeds by induction in the same way 
as the proof of Proposition $\ref{Kuraeistspoly}$,
except the point in which the moduli space $\mathcal M_{\ell+1;\emptyset}(\kappa;B)$ 
are involved.  We discuss this case in detail below.
\begin{rem}\label{rem36}
We remark that in Proposition \ref{Kuraeistspoly2} (2) 
the Kuranishi structure we put on the `bubble'\footnote{Here `bubble' means an 
extended disk component such that the $0$-th interior marked point is not on it.} is that 
of Propositions \ref{Kuraeistspoly}, \ref{prop94}, that is, the $\frak q$-Kuranishi structure.
\end{rem}
We now discuss the moduli space $\mathcal M_{\ell+1;\emptyset}(\kappa;B)$.
Its element is represented by $(\Sigma,(z^+_0,z^+_1,\dots,z^+_{\ell});u)$ 
where $\Sigma$ is a genus zero bordered curve with one boundary 
component and $z^+_0,z^+_1,\dots,z^+_{\ell}$ are interior marked points 
which are mutually distinct.  The map $u : \Sigma \to X$ is a pseudo-holomorphic map 
such that $u(\partial \Sigma) \subset L_{\kappa}$.
\begin{defn}\label{defn155555}
We say $(\Sigma,(z^+_0,z^+_1,\dots,z^+_{\ell});u)$  is {\it weakly stable}\index{weakly stable} if one of the 
following two conditions are satisfied.
\begin{enumerate}
\item The group of automorphisms of $(\Sigma,(z^+_0,z^+_1,\dots,z^+_{\ell});u)$ is a finite group.
\item  The connected component of the group of automorphisms of $(\Sigma,(z^+_0,z^+_1,\dots,z^+_{\ell});u)$ is
$S^1$.  An element of this $S^1$ acts non-trivially on $\partial\Sigma$.
\end{enumerate}
\end{defn}
We observe that case (2) occurs  (only) when the following 3 conditions are all satisfied.  
\begin{enumerate}
\item[(2.a)]  There exists a unique extended disk component of $\Sigma$.
\item[(2.b)]  On the unique disk component of $\Sigma$ there exists exactly 
one nodal or (interior) marked point.  In case it is a marked point,
$\ell = 0$, and the marked point on the disk component is the $0$-th one.
\item[(2.c)]  The map $u$ is a constant map on the unique disk component. 
\end{enumerate}
See Figure \ref{Figure15.6}.
The next lemma is easy to show.
\begin{lem}
The moduli space of objects satisfying {\rm (2.a), (2.b), (2.c)}  above 
coincides with the fiber product $(\ref{eq157})$.
\end{lem}
\begin{defn}
$\mathcal M_{\ell+1;\emptyset}(\kappa;B)$ is the set of 
isomorphism classes of objects $(\Sigma,(z^+_0,z^+_1,\dots,z^+_{\ell});u)$
which are weakly stable and $u_*([\Sigma]) = B$.
\end{defn}
Now we have:
\begin{lem}\label{lem158}
For each element  ${\bf p} = (\Sigma,(z^+_0,z^+_1,\dots,z^+_{\ell});u)$ satisfying {\rm (2.a),(2.b),(2.c)}
we can define an obstruction space $E_{{\bf p}}({\bf p})$ satisfying  Definition $\ref{defn51}$ $(1)(2)$.
Moreover the following holds:
\begin{enumerate}
\item[$(3)'$]
$$
{\rm EV} : D_u\overline{\partial}^{-1}(E_{{\bf p}}({\bf p})) \to T_{u(z^+_0)} X 
$$
obtained by the evaluation at $z^+_0$ is a submersion.  (This is replacement of Definition $\ref{defn51}$ $(3)$.)\index[syindex]{EV@${\rm EV}$}
\item[$(5)'$]
Moreover $E_{{\bf p}}({\bf p})$ is invariant under the action of the automorphism group ${\rm Aut}({\bf p})$ of ${\bf p}$.
(This is a replacement of Definition $\ref{defn51}$ $(5)$.)
\item[$(6)'$]
The action of $S^1 \subset {\rm Aut}({\bf p})$ on $E_{{\bf p}}({\bf p})$ 
is {\bf trivial}.
The action of ${\rm Aut}({\bf p})/S^1$ on $D_u\overline{\partial}^{-1}(E_{{\bf p}}({\bf p}))/{\rm aut}((\Sigma,(z^+_0,z^+_1,\dots,z^+_{\ell}))$
is effective.\footnote{We can choose $E_{{\bf p}}({\bf p})$ so that the $S^1$ action on $D_u\overline{\partial}^{-1}(E_{{\bf p}}({\bf p}))$ is 
not only well defined but also trivial.} (This is the replacement of Definition $\ref{defn51}$ $(6)$.)
\end{enumerate}
\end{lem}
By the compactness of  ${\rm Aut}({\bf p})$ with $S^1$  its connected component, we can easily find such  $E_{{\bf p}}({\bf p})$.
Then in the same way as the case of $\mathcal M_{\ell+1;\vec{k}}((\vec{\kappa},\vec p);B)$ with $\vec{\kappa} \ne \emptyset$
we can show:
\begin{lem}\label{lem159}
There exists an obstruction bundle data $\{E_{\bf p}({\bf x})\}$ for 
$\mathcal M_{\ell+1;\emptyset}(\kappa;B)^{\boxplus 1}$ that satisfies Definition $\ref{defn51}$  $(1)-(6)$.
At the point satisfying {\rm (2.a),(2.b),(2.c)} we replace $(5)$, $(6)$ by $(5)'$, $(6)'$ above
and $(3)$  by $(3)'$.
\par
The obstruction bundle data is symmetric under the permutation 
of the $1$st-$k$-th interior marked points, as well as the cyclic permutations of the 
boundary marked point.
%\par
%In the case ${\bf p} \in \mathcal M_{1;\emptyset}(\kappa;0)$ we choose $E_{\bf p}(\bf x)$ to be 2 dimensional.
\end{lem}
Using Lemma \ref{lem158} we can prove Lemma \ref{lem159} in the same way as 
\cite[Theorem 11.1]{const1}.
\par
Now we have:
\begin{lem}\label{lem1510}
Given the obstruction bundle data as in Lemma $\ref{lem159}$ we can associate 
Kuranishi structures on various $\mathcal M_{\ell+1;\vec{k}}((\vec{\kappa},\vec p);B)^{\boxplus 1}$'s so that they satisfy 
Proposition $\ref{Kuraeistspoly2}$ $(7)$, $(8)$, $(9)$.
\end{lem} 
\begin{rem}\label{rem1511}
Note that the moduli space $\mathcal M_{1;\emptyset}(\kappa;0)$ is not involved in our
definition of the open-closed 
map $\frak p$.  However, we need to study the  
obstruction bundle data of this moduli space to formulate 
the compatibility of the forgetful map of the unique boundary marked point of 
$\mathcal M_{1;(\kappa,\dots,\kappa)}(\kappa;0)$ while keeping the cyclic symmetry of the boundary marked points.
\end{rem}
\begin{proof}
The proof is mostly the same as that of \cite[Theorem 7.1]{const1}.
The Kuranishi neighborhood of ${\bf p} \in \mathcal M_{\ell+1;\vec{k}}((\vec{\kappa},\vec p);B)$
is given as 
$$
U_{\bf p} = \{ {\bf x} \mid \overline \partial u_{\bf x} \in E_{\bf p}({\bf x})\}.
$$
Here ${\bf x}$ is an object of the ambient set which is $\epsilon$-close to ${\bf p}$
in the sense of partial topology which is defined in \cite[Definition 4.1]{const1}.  
(See also (1)(2)(3) right above Definition \ref{defn51}.)
An object ${\bf x}$ consists of a 
marked bordered curve $(\Sigma_{\bf x},z_{\bf x},z^+_{\bf x})$ 
and a map $u_{\bf x}: \Sigma_{\bf x} \to X$. $U_{\bf p}$ consists of isomorphism classes of nearby objects ${\bf x}$ 
satisfying the equation $\overline \partial u_{\bf x} \in E_{\bf p}({\bf x})$.
Using the transversality, Definition \ref{defn51}  (3), and the smoothness of obstruction bundle data, 
Definition \ref{defn51}  (2),  we can show that 
$U_{\bf p}$ is a smooth orbifold.  
This is a consequence of \cite{foooexp}.
A new point to be discussed is the case when ${\bf p}$ satisfies (2.a),(2.b),(2.c).
We will explain this case later in this proof.
\par
The obstruction bundle $\mathcal E_{\bf p}$ is defined so that its fiber at ${\bf x}$ is 
$E_{\bf p}({\bf x})$. The Kuranishi map is $s_{\bf p}({\bf x}) = \overline \partial u_{\bf x}$.
Using the semi-continuity of the obstruction bundle data, Definition \ref{defn51} (4), we can define the coordinate change map.
(Actually it is the inclusion as  elements of the ambient set.)
\par
Proposition \ref{Kuraeistspoly2} (7) is a consequence of Lemma \ref{lem158} (3)$'$.
Proposition \ref{Kuraeistspoly2} (8)(9) are consequences of corresponding 
symmetry of the obstruction bundle data.
\par
Now we discuss the case when  ${\bf p}$ satisfies (2.a),(2.b),(2.c).
This case is similar to the discussion of \cite[Subsection 7.4.1]{fooo092}, which we explain now.
\par
We first consider the simplest case where ${\bf p} \in \mathcal M_{1;\emptyset}(\kappa;0)$.
Our ${\bf p}$ is $((D^2,0);i_p)$ where  $(D^2,0)$ is a disk together with the $0$-th marked point $0$.
$i_p$ is the constant map to $p \in L_{\kappa}$. 
We put a (non-trivial) obstruction bundle $E_{{\bf p}}({\bf x}) \subset C^{\infty}(u_{\bf x}^*(TX) \otimes \Lambda^{0,1})$. Here
${\bf x} = ((D^2,0);u_{\bf x})$ and $u_{\bf x} : (D^2,\partial D^2) \to (X,L_{\kappa})$ is close 
to a constant map but is not necessarily constant. We require  $\overline \partial u_{\bf x} \in  E_{\bf p}({\bf x})$.
Since constant maps are transversal, the set of such maps $u_{\bf x}$ is identified 
with the product of a neighborhood $U(L)$ of $p$ in $L_{\kappa}$ and a neighborhood of $0$ of $E_{\bf p}({\bf p})$. 
The group $S^1$ acts trivially on $E_{\bf p}({\bf p})$.  Therefore 
  $U_{\bf p} \cong  U(L) \times B_{\epsilon}(0;\R^m)$. If this neighborhood is small then the $S^1 = {\rm Aut}({\bf p})$ 
action is trivial on   $U_{\bf p}$. 
A neighborhood of $\bf p$ in $\mathcal M_{1;\emptyset}(\kappa;0)$ is $U(L)$.  We remark that the virtual dimension of $\mathcal M_{1;\emptyset}(\kappa;0)$ is $\dim L - 1$
because of $S^1$ symmetry.  On the other hand $\dim U_{\bf p}$ is $\dim L + m$ and the rank of $E_{\bf p}$ is $m$.
Hence $\dim U_{\bf p} - {\rm rank}(E_{\bf p}) = \dim L$ does not coincide with virtual 
dimension. This is because $S^1$ action is trivial on $U_{\bf p}$.
So in this case our `chart' is not Kuranishi chart in the usual sense.
\par
We next consider the case of $\mathcal M_{\ell+1;\emptyset}(\kappa;B)$ with $\ell > 0$ or $B \ne 0$.
Suppose ${\bf p} \in \mathcal M_{\ell+1;\emptyset}(\kappa;B)$ satisfies (2.a),(2.b),(2.c).
Then ${\bf p}$ is represented by $((D^2 \cup \Sigma^s;\vec z^+),u_{\bf p})$.
Here $\Sigma^s$ is a tree of sphere components attached to $0 \in D^2$ 
and $u_{\bf p} : D^2 \cup \Sigma^s \to X$ is holomorphic such that $u_{\bf p}(\partial D^2)
\in L_\kappa$.  $u_{\bf p}$ is constant on $D^2$.  Moreover 
$((\Sigma^s;\vec z^+),u_{\bf p}\vert_{\Sigma^s})$ is a stable map 
in the usual sense.
The set of elements with the same combinatorial type as ${\bf p}$ is described by the fiber product 
(\ref{eq157}). 
\par
To obtain a Kuranishi neighborhood we include a gluing parameter to smooth the double point 
$0 \in D^2$.  The set of gluing parameters is identified with a small neighborhood of $0$ in $\C$.
The automorphism group $S^1$ of ${\bf p}$ acts as a linear action of $U(1) = S^1$ on $\C$, which is the set of gluing parameters. 
In particular the action is free.  Therefore the actual 
gluing parameter (appearing in the Kuranishi neighborhood) is $B_{\epsilon}(\C)/S^1 \cong [0,\epsilon)$.
Thus $U_{\bf p}$ is an orbifold with boundary.  Since $S^1$ action on $E_{\bf p}(\bf p)$ is 
trivial $\mathcal E_{\bf p}$ is a well defined orbibundle on $U_{\bf p}$. Thus in this case 
we have a Kuranishi neighborhood in the usual sense.  See \cite[Subsection 7.4.1]{fooo092} for the detailed argument 
of this case.
\end{proof}
\begin{rem}
Note that we need to include non-trivial obstruction bundle 
$E_{\bf p}({\bf x})$ for ${\bf p} \in \mathcal M_{1;\emptyset}(\kappa;0)$.
This is because otherwise the evaluation map at the $0$-th {\it interior} marked 
point is not weakly submersive.
This point is a big difference from the case of $\frak q$-perturbation.
In that case we need the weak submersivity of the evaluation map 
at the $0$-th {\it boundary} marked point. The moduli space
$\mathcal M_{0;1}(\kappa;0)$ 
is transversal and is  identified with $\tilde L_{\kappa}$.
Therefore we do not need to put nontrivial obstruction bundle but have only to put zero there
so that the evaluation map $\mathcal M_{0;1}(\kappa;0) \to \tilde L_{\kappa}$  
becomes weakly submersive.
This point was very important
to make our system of $\frak q$-Kuranishi structures compatible 
with the forgetful map of boundary marked points.
\par
We remark that we put the $\frak q$-Kuranishi structure on the disk bubbles 
for the moduli spaces in Proposition \ref{Kuraeistspoly2} (which we use to define
the operator $\frak p$).
Therefore we equip the trivial obstruction bundle with the constant disk bubbles. On the other hand,
we equip the $\frak p$-Kuranishi obstruction bundle data with the extended disk component
which contains the $0$-th interior marked point. The latter obstruction bundle data is non-trivial
on this extended disk component where the map is constant.
\end{rem}

Now to prove Proposition \ref{Kuraeistspoly2} we need to take the obstruction bundle data 
(which induces ${\frak p}$-Kuranishi structure) so that it is compatible with the fiber product 
description of the boundary.  This step is similar to and easier than the case of ${\frak q}$ perturbation 
and so is omitted.
We however explain the  compatibility of Kuranishi structures with 
the forgetful map in the case $\vec k \setminus \frak f_{\bf e} = \emptyset =  \frak f_{\bf f}$,
since this case is new in our situation.
In this case we need to define a forgetful map
\begin{equation}\label{eq158}
\frak{forget} : {\mathcal M}_{\ell+1;\vec k}(((\kappa,\dots,\kappa),\emptyset);B;\vec{\frak f})
\to {\mathcal M}_{\ell+1;\emptyset}(\kappa;B).
\end{equation}
(Here $\emptyset$ in $((\kappa,\dots,\kappa),\emptyset)$ means that $\vec p$ is empty.)
We took the obstruction data of the target.
We define the obstruction data of the domain 
so that it is compatible with this forgetful map.
Then the induced Kuranishi structures are compatible with 
forgetful map.
\par
In case $\ell = 0$ and $B  = 0$, the `chart' we have on $ {\mathcal M}_{1;\emptyset}(\kappa;0)$ 
is not a Kuranishi chart.  However we do have obstruction bundle data and all the elements of the 
obstruction bundle are fixed by the $S^1$ action. (Lemma \ref{lem158} (6)$'$).
Therefore we can pull back obstruction bundle data.
\par
We next modify Propositions \ref{existmultipolu}, \ref{exstCFUNI} as follows.
\begin{prop}\label{existmultipolu2}
There exist \blue{CF-perturbations} of 
$([0,1]^{\vert \frak f_{\bf f}\vert} \times{\mathcal M}_{\ell+1;\vec k}((\vec{\kappa},\vec p);B;\vec{\frak f})^{\boxplus 1})^{\frak p}$ with 
the following properties.
\begin{enumerate}
\item
They are transversal to $0$.
\item 
They are compatible with the description of the 
boundary in Proposition $\ref{Kuraeistspoly2}$ $(2)$.
(Note that, on the factors equipped with the Kuranishi structures 
in Propositions $\ref{Kuraeistspoly}$ and $\ref{prop94}$,
we consider the \blue{CF-perturbations}  in 
Propositions $\ref{existmultipolu}$ and   $\ref{exstCFUNI}$
and on the factors equipped with the Kuranishi structures 
in Proposition $\ref{Kuraeistspoly2}$, 
we consider the \blue{CF-perturbations} in 
Proposition $\ref{existmultipolu2}$.)
\item 
They are compatible with the forgetful map 
appearing in Proposition $\ref{Kuraeistspoly2}$ $(5)$ and $(\ref{eq158})$.
\item
The compatibility in Proposition $\ref{Kuraeistspoly2}$ $(6)$
is enhanced to the compatibility of CF-perturbations.
\item
The evaluation map
$$
\text{\rm ev}_{0}^+ : 
 ([0,1]^{\vert \frak f_{\bf f}\vert} \times{\mathcal M}_{\ell+1;\vec k}((\vec{\kappa},\vec p);B;\vec{\frak f})^{\boxplus 1})^{\frak p}
\to X
$$
is strongly submersive with respect to our CF-perturbations.
\item
The \blue{CF-perturbation} is invariant under the permutation of 
$\ell$ interior marked points $z^+_1,\dots,z_{\ell}^+$.
\item
The \blue{CF-perturbation} is invariant under the cyclic permutation 
of the data. 
\end{enumerate}
\end{prop}
We call the perturbation in Proposition \ref{existmultipolu2} the {\it $\frak p$-perturbation}. \index{pperturbation@$\frak p$-perturbation} 
The proof is mostly the same as that of Proposition  \ref{existmultipolu} except the point 
we discuss below.
We remark that the Kuranishi chart we defined on $\mathcal M_{1;\emptyset}(\kappa;0)$
 is a `quotient' by trivial $S^1$ action. The Kuranishi neighborhood of this Kuranishi chart is $U(L_{\kappa}) \times \R^m$ with a trivial $S^1$ action.
Here $U(L_\kappa)$ is an open subset of our Lagrangian $L_\kappa$.  (The evaluation map ${\rm ev}_0$ is the projection to this factor.)
The $S^1$ action on the obstruction bundle is also trivial.
So we take an $S^1$ invariant CF-perturbation of the Kuranishi chart,\footnote{Actually any CF-perturbation is 
$S^1$ invariant.} which we can pull back.
This moduli space appears only once at the first step of the inductive 
construction of the CF-perturbations.
\par\smallskip

We now change our notation as in Section \ref{constcyclic}. 
\par
We first introduce ${\mathcal M}_{\ell+1;\vec k}((\vec{\kappa},\vec p,m);B;\vec{\frak f})$
(the version which include $m \in \{0,\dots,k_0\}$ that determines the choice of $0$-th marked point).  
We then define $\mathcal M_{\ell+1}(\vec \kappa,\vec p;B)$ using Equation \eqref{eq:change_notation_moduli_space} to be the moduli space associated to the reduction of the sequence $\vec \kappa$. This space is equipped with evaluation maps $\{{\rm ev}^+_j \}_{j=0}^{\ell}$ to $X$, and $\{ \text{\rm ev}^{\partial}_{i}\}_{i=0}^{K}$ to the fiber produce $\tilde L_{\kappa_{i-1}}\times_X \tilde L_{\kappa_{i}}$ (in the case $K=0$, we have $L_{\kappa_{K+1}} = L_{\kappa_{0}}$ by convention).

We put\index[syindex]{CHCFLformcunonplus@$CHCF_*(\cL^{\rm form};\vec{\kappa})$}
\begin{equation}
\aligned
CHCF(\cL^{\rm form};\vec{\kappa})
&=CF(L_{\kappa_{K}},L_{\kappa_{0}};\F) [1] \otimes \bigotimes_{i=1}^K CF(L_{\kappa_{i-1}},L_{\kappa_{i}};\F) [1]
 \\
&=
\bigotimes_{i=0}^K CF(L_{\kappa_{i}},L_{\kappa_{i+1}};\F) \blue{[1]},
\endaligned
\end{equation}
and its $+$ version (\ref{CHCFform}).
We then define\index[syindex]{pafcubfrak@${\frak p}^{\rm f.c.u. \frak b}_{\ell,\vec{\kappa}}$}
$$
{\frak p}^{\rm f.c.u. \frak b}_{\ell,\vec{\kappa}} 
:
\Omega(X) \blue{[2]}^{\otimes \ell} \otimes 
CHCF(\cL^{\rm form};\vec{\kappa} )^+
\to 
\Omega(X) \blue{[2]} \hat\otimes \Lambda_0
$$
as follows:
\begin{defn}
For $\text{\bf g} \in \Omega(X)^{\otimes \ell}$,
$\text{\bf h} \in CHCF(\cL^{\rm form};\vec{\kappa})$, 
%$\text{\bf h}' \in B_{k_i}(\Omega(L_{\kappa_{i}})[1])$
we put\index[syindex]{paformfrakellB@${\frak p}^{\rm form}_{\ell,\vec{\kappa},B}$} 
\begin{equation}\label{form1430}
\aligned
&{\frak p}^{\rm form}_{\ell,\vec{\kappa},B}
(\text{\bf g},\text{\bf h}) \\ %,\text{\bf h}')\\
&=
(-1)^{\maltese}\text{\rm Corr}({\mathcal M}_{\ell+1}(\vec{\kappa},\vec p;B)),
(\text{\rm ev}^{\partial}\times (\text{\rm ev}^+_1,\dots,\text{\rm ev}^+_{\ell}), \text{\rm ev}^+_0))(\text{\bf g},\text{\bf h}).%,\text{\bf h}').
\endaligned
\end{equation}
The sign $\maltese$ is the same as in (\ref{form21ten2}).
We then put 
\begin{equation}\label{form1611new}
{\frak p}^{\rm f.c.u. \frak b}_{\ell,\vec{\kappa}}
(\text{\bf g},\text{\bf h})
=\sum_{B,\ell',\vec k}
T^{\omega(B)}
\frac{\rho_{\frak b,\theta}(B)%\rho_{b}(B)
}{\ell'!}
{\frak p}^{\rm form}_{\ell+\ell';\vec{\kappa};B}
(\text{\bf g}\otimes\frak b_+^{\otimes\ell'},h_0, \ldots, h_K).
\end{equation}
%where $b_+^{\vec k} = \bigotimes_{i=0}^K b_{\kappa_i,1}^{\otimes k_i}$.
Using the marked points labeled by ${\frak f}_{\bf e}$ and ${\frak f}_{\bf f}$ we can 
extend it to\index[syindex]{pafcubfrakellB@${\frak p}^{\rm f.c.u. \frak b}_{\ell,\vec{\kappa},B}$} 
\begin{equation}
{\frak p}^{\rm f.c.u. \frak b}_{\ell,\vec{\kappa},B} %,\vec k}
:
\Omega(X) \blue{[2]}^{\otimes \ell} \otimes 
CHCF(\cL^{\rm form};\vec{\kappa})^+
\to 
\Omega(X) \blue{[2]}
\end{equation}
in the same way as Subsection \ref{subsec:opeartor}.
\end{defn}

In Subsection \ref{frakponcoho2} we will define maps
\begin{equation}\label{formmap826}
\widehat{\frak p}^{\rm f.c.u. \frak b}_{n,k}
: CH_*^{n,k}(\cL^{\rm form}_{\rm c.u.})^+ 
 \to \Omega(X) 
\otimes \Lambda_0/T^{E'_n}
\end{equation}
from the Hochschild chain complex to the de Rham complex,
using the maps
$\widehat{\frak p}^{\rm f.c.u. \frak b}_{0,\vec{\kappa}}$ defined  
for various $\vec{\kappa}$.
Here $E'_n$ is a  sequence  with $\lim_{n\to \infty} E'_n = \infty$.
Existence of such $E'_n$ follows from the Gromov compactness, 
in the same way as the proof of convergence of (\ref{form:716}).
 %\marginpar{$E_n$ is changed to $E'_n$. 
%KF 2025 Feb. Moved to here.}
 
We  define $CH_*^{n,k}(\cL^{\rm form}_{\rm c.u.})^+$
also in the next subsection. \subsection{Homotopy inverse limit for the construction of  $\widehat{\frak p}$.}
\label{frakponcoho2}
 
Because of the problem we explained in Remark \ref{runningoutrem}
the construction of Subsection \ref{frakponcoho} actually defines a map 
for each $(E,k) = (E_n,k)$ only. 

We can prove their compatibilities for various $(n,k)$ which enable us to
define the map\index[syindex]{paformcubbf@$\widehat{\frak p}^{\rm f.c.u. \frak b}$} 
\begin{equation}\label{pformdef}
\widehat{\frak p}^{\rm f.c.u. \frak b} 
: CH_*(\cL^{\rm form}_{\rm c.u.})^+ \to \Omega(X) 
\widehat\otimes \Lambda_0
\end{equation}
as their limit.
The details of this argument are similar to but are slightly different from the case of $\widehat{\frak q}$, as we now explain.
\par
We first define\index[syindex]{CH*upnkLform@$CH_*^{n,k}(\cL_{\rm c.u.}^{\rm form})$}
$
CH_*^{n,k}(\cL_{\rm c.u.}^{\rm form}) 
$
as the set of formal sums
\begin{equation}\label{form1434}
{\frak T} = \sum T^{\lambda_i} {\frak T}_i \in CH_*(\cL^{\rm form}_{\rm c.u.})^+
\end{equation}
where $\lambda_i \in G$, ${\frak T}_i \in CHCF(\cL^{\rm form};\vec{\kappa})$
such that $(\lambda_i,\vert \vec{\kappa}\vert) \le (E_n,k)$.
Note (\ref{form1434}) is necessary a finite sum.  (We recall that our collection $\mathscr L$ is 
assumed to be a finite set.)
\par
Note that Hochschild chain complex of the de Rham complex $CH_*(\cL_{\rm c.u.}^{\rm form})^+$
over Novikov ring $\Lambda_G$ is the set of all 
{\it infinite} sums of the form (\ref{form1434}) where $\lambda_i \in G$.
(We do not require a condition such as $(\lambda_i,\vert \vec{\kappa}\vert) \le (E_n,k)$.)
We consider an element $x$ of $CH_*^{n,k}(\cL_{\rm c.u.}^{\rm form})^+$ and take its 
Hochschild boundary as an element of $CH_*(\cL_{\rm c.u.}^{\rm form})^+$, 
which we write
$$
\delta_H(x) = \sum T^{\lambda_i} {\frak T}'_i.
$$
We consider only the terms $T^{\lambda_i} {\frak T}'_i$ where  $(\lambda_i,\vert \vec{\kappa}\vert) \le (E_n,k)$ 
and ${\frak T}'_i \in CHCF(\cL^{\rm form};\vec{\kappa})$
which we write
$$
\delta^{n,k}_H(x) = \sum_{i\in I}T^{\lambda_i} {\frak T}'_i.
$$
Thus we obtain 
$\delta^{n,k}_H: CH_*^{n,k}(\cL^{\rm form}_{\rm c.u.})^+ \to CH_*^{n,k}(\cL^{\rm form}_{\rm c.u.})^+$.
The next equality is a consequence of $\delta_H\circ \delta_H=0$ and 
the choice of the order $<$ (Definition \ref{defn837837}).
\begin{equation}
\delta^{n,k}_H \circ \delta^{n,k}_H = 0.
\end{equation}
For $
{\frak T} \in CH_*^{n,k}(\cL_{\rm c.u.}^{\rm form}) 
$
we define ${\frak p}^{\rm f.c.u. \frak b}_{\ell,\vec{\kappa}}
(\text{\bf g},\frak T)$
by (\ref{form1611new}).
Then ${\frak p}^{\rm f.c.u. \frak b}_{0,\vec{\kappa}}$ 
induces\index[syindex]{pankformcubbf@$\widehat{\frak p}_{n,k}^{\rm f.c.u. \frak b}$}  %\marginpar{error corrected.  KF 2025 Aug.}
\begin{equation}\label{pformdef2}
\widehat{\frak p}_{n,k}^{\rm f.c.u. \frak b} 
: CH_*^{n,k}(\cL^{\rm form}_{\rm c.u.})^+ \to \Omega(X) 
\otimes \Lambda_0/T^{E'_n}.
\end{equation}

\begin{prop}\label{prop142566}
$\widehat{\frak p}_{n,k}^{\rm f.c.u. \frak b}$ is a chain map. Namely
\begin{equation}
d \circ \widehat{\frak p}_{n,k}^{\rm f.c.u. \frak b} 
- \widehat{\frak p}_{n,k}^{\rm f.c.u. \frak b}  \circ \delta^{n,k}_H
\equiv 0  \mod T^{E'_n}.
\end{equation}
\end{prop}
\begin{proof}
This is a consequence of Propositions 
\ref{Kuraeistspoly2}, \ref{existmultipolu2} and 
Stokes' theorem.
(We use Remark \ref{rem36} in the proof.)
\end{proof}

Now we use homotopy inverse limit argument to construct $\widehat{\frak p}^{\rm f.c.u. \frak b}$.
Note that the domain of $\widehat{\frak p}^{\rm f.c.u. \frak b}$ is the
direct {\it sum} of various tensor products 
$BCF(\cL^{\rm form}_{\rm c.u.};\vec{\kappa}; \text{\rm whole})$. 
In the same way as in the case of $\widehat{\frak q}^{\rm form}$ we take 
its completion with respect to the Fr\'echet topology of the space of differential forms
as follows. 
Let $\vec\kappa = (\kappa_0,\dots,\kappa_K)$ be as above.
Recall
$$
BCF(\cL^{\rm form};\vec{\kappa})^+
= \bigotimes_{i=1}^{K+1} CF((L_{\kappa_{i-1}},b_{\kappa_{i-1}}),(L_{\kappa_{i}},b_{\kappa_{i}});\F)^+ \blue{[1]}
$$
where
$$
\aligned
&CF((L_{\kappa_{i-1}},b_{\kappa_{i-1}}),(L_{\kappa_{i}},b_{\kappa_{i}});\F)^+ \\
&= 
\begin{cases}
\displaystyle
\bigoplus_{p \in L_{\kappa_{i-1}} \cap L_{\kappa_{i}}} \F[p]
& \text{if $L_{\kappa_{i-1}} \ne L_{\kappa_{i}}$}, \\
(\Omega(L_{\kappa_{i}}) \otimes \F) \oplus  \displaystyle\bigoplus_{p \in (\tilde L_{\kappa_i} \times_{X} \tilde L_{\kappa_i}) \setminus \tilde L_{\kappa_i}} \F[p]
& \text{ if $L_{\kappa_{i-1}} = L_{\kappa_{i}}$ $b_{\kappa_{i-1}}\ne b_{\kappa_{i}}$},\\
\displaystyle(\Omega(L_{\kappa_{i}})\otimes \F) \oplus  \bigoplus_{p \in (\tilde L_{\kappa_i} \times_{X} \tilde L_{\kappa_i}) \setminus \tilde L_{\kappa_i}} \F[p]
\oplus \F {\bf e}^+_{\kappa_{i}}
\oplus
\F {\bf f}_{\kappa_{i}} &
\text{ if $L_{\kappa_{i-1}} = L_{\kappa_{i}}$ $b_{\kappa_{i-1}}= b_{\kappa_{i}}$}.
\end{cases}
\endaligned
$$
Let $BCF(\cL^{\rm form};\vec{\kappa})$ be the direct summand of 
$BCF(\cL^{\rm form}_{\rm c.u.};\vec{\kappa})^+$ which complements to the span of the part
containing ${\bf e}^+$ and ${\bf f}$.

We take\index[syindex]{Lveckappa@$L(\vec{\kappa})$}
$$
L(\vec{\kappa}) = \prod_{i=1}^{K+1} (\tilde L_{\kappa_{i-1}} \times_X \tilde L_{\kappa_{i}}).
$$
(Here we put $\kappa_{K+1} = \kappa_0$ as convention.)
 %\marginpar{Changed to $+$ version from 
%here KF. 2025 March.}
There is an obvious embedding
$$
\frak I : BCF(\cL^{\rm form};\vec{\kappa})^+ \to \Omega(L(\vec{\kappa})).
$$
We put\index[syindex]{BCFhatLcuform@${BCF}^{\infty}(\cL^{\rm form};\vec{\kappa})^+$} 
$$
{BCF}^{\infty}(\cL^{\rm form};\vec{\kappa})^+
= 
\Omega(L(\vec{\kappa})).
$$
We include ${\bf e}^+$ and ${\bf f}$ and obtain an $\F$ module
${BCF}^{\infty}(\cL^{\rm form};\vec{\kappa})^+$,
in which
${BCF}(\cL^{\rm form};\vec{\kappa})^+$ is dense 
with respect to the 
$C^{\infty}$ norm  of the differential forms.
\par
We remark that ${BCF}^{\infty}$ is a completion of ${BCF}$. We define
${CH}^{\infty}_*(\cL_{\rm c.u.}^{\rm form})$\index[syindex]{CH*upnkLformhat@${CH}^{\infty}_*(\cL_{\rm c.u.}^{\rm form})$} 
and ${CH}^{\infty}(\cL^{\rm form}_{\rm c.u.})^+$\index[syindex]{CH*upnkLformhat+@${CH}^{\infty}_*(\cL_{\rm c.u.}^{\rm form})^+$}  as completions of ${CH}_*(\cL_{\rm c.u.}^{\rm form})$ and
${CH}_*(\cL_{\rm c.u.}^{\rm form})^+$ in the same way. 
(The latter contains ${\bf e}^+$ and ${\bf f}$.)
Since the Hochschild differential is defined either by the de Rham differential 
or by a correspondence via distributional form satisfying Condition \ref{cod149}, it extends to the completions 
${CH}^{\infty}_*(\cL_{{\rm c.u.}}^{\rm form})$.
By the above mentioned denseness the extension is unique.
%We define $\widehat{CH}_*(\cL_{\rm c.u.}^{\rm form})$ to be the completion of 
%${CH}_*(\cL_{\rm c.u.}^{\rm form})$ in the same way.
We can further extend it to ${CH}^{\infty}(\cL^{\rm form}_{\rm c.u.})^+$ 
using the evaluation at the marked points labeled by $\frak f_{\bf e}$, $\frak f_{\bf f}$ 
in the same way as in Subsection  \ref{subsec:opeartor}, and obtain  ${CH}^{\infty}_*(\cL_{\rm c.u.}^{\rm form})^+$.
\par
For the sake of simplicity of notation we  denote by  
${CH}_{\infty}^*(\cL_{\rm c.u.}^{\rm form})$, 
${CH}_{n,k,\infty}^*(\cL_{\rm c.u.}^{\rm form})$ the Hochschild
complexes
$CH_{\infty,G}^*(\cL_{\rm c.u.} ^{\rm form},\cL ^{\rm form}_{\rm c.u.})
$, $CH_{\infty,G;n,k}^*(\cL_{\rm c.u.} ^{\rm form},\cL ^{\rm form}_{\rm c.u.})
$, defined in Definition \ref{defn148148}.  We use similar notations for $+$ versions.
\par
Since (\ref{formmap826}) is also defined by a smooth correspondence, 
it extends to a chain map\index[syindex]{pankformcubbfmmorehat@${\widehat{\frak p}}^{\rm f.c.u. \frak b}_{n,k,\infty}$} 
\begin{equation}\label{formmap82622}
{\widehat{\frak p}}^{\rm f.c.u. \frak b}_{n,k,\infty}
: {CH}^{n,k,\infty}_*(\cL_{\rm c.u.}^{\rm form})^+
\to \Omega(X) 
\otimes \Lambda_0/T^{E'_n},
\end{equation}
 as well.
Here ${CH}^{n,k,\infty}_*(\cL_{\rm c.u.}^{\rm form})^+$ is the completion of ${CH}^{n,k}_*(\cL_{\rm c.u.}^{\rm form})^+$
in ${CH}^{\infty}_*(\cL_{\rm c.u.}^{\rm form})^+$.
(For the part where ${\bf e}^+$ and ${\bf f}$ are involved we can define the completion in the same way as in 
Subsection \ref{subsec:opeartor}.\footnote{We do not repeat this remark below.})

We take a sequence $(E_j,k_j)$ such that $E_j \in G$, $k_j \in \Z_{+}$
and $(E_j,k_j) < (E_{j+1},k_{j+1})$.
(Note that then the discreteness of $G$ implies that for any $E \in G$ and $k \in \Z_+$, we have $(E,k) < (E_j,k_j)$ for sufficiently large $j$.)

For $j'< j$, by the inductive construction of the $A_{\infty}$ category, 
there exists a quasi-isomorphism
$$
\frak f_{j',j}:
(\cL^{\rm form}_{\rm c.u.},\frak m^{{\rm f.c.u.} \frak b;j})\vert_{(E_{j'},k_{j'})}
\to 
(\cL^{\rm form}_{\rm c.u.},\frak m^{{\rm f.c.u.} \frak b;j'})
$$
Here 
the domain  
is the $A_{E_{j},k_{j}}$ category  
that is obtained in Subsection \ref{constcyclic} 
and is regarded as an $A_{E_{j'},k_{j'}}$ category.\footnote{The notation 
$\vert_{(E_{j'},k_{j'})}$ means that we regard an $A_{E_{j},k_{j}}$ category   as an 
$A_{E_{j'},k_{j'}}$ category.}
The target is the $A_{E_{j'},k_{j'}}$ category .
The functor $\frak f_{j',j}$ is a quasi-isomorphism of $A_{E_{j'},k_{j'}}$ categories.
The $A_{E_{j'},k_{j'}}$ morphism $\frak f_{j',j}$ induces a chain homotopy equivalence:\footnote{We here use the fact 
that Hochschild chain complex is covariant. Note that the Hochschild cochain complex 
is neither covariant nor contravariant.}
\begin{equation}\label{form1620}
CH({\frak f}_{j',j}^{+}) : 
{CH}^{E_{j'},k_{j'}}_*(\cL_{\rm c.u.}^{\rm form},\delta_H^{E_{j},k_{j}})^+
\to {CH}^{E_{j'},k_{j'}}_*(\cL_{\rm c.u.}^{\rm form} ,\delta_H^{E_{j'},k_{j'}})^+.
\end{equation}
We can take its completion:
\begin{equation}\label{form1621new}
CH({\frak f}_{j',j}^{\infty,+}) : 
{CH}^{E_{j'},k_{j'},\infty}_*(\cL_{\rm c.u.}^{\rm form},\delta_H^{E_{j},k_{j}})^+
\to {CH}^{E_{j'},k_{j'},\infty}_*(\cL_{\rm c.u.}^{\rm form} ,\delta_H^{E_{j'},k_{j'}})^+.
\end{equation}
We can define this map on the completed Hochschild chain complex 
since $\frak f_{j',j}^+$ (that is, $\frak f_{j',j}$ with ${\bf f}$ and ${\bf e}^+$ included) is obtained by a  correspondence via distributional form satisfying Condition \ref{cod149}.
\begin{rem}
The Hochschild differential $\delta_H^{E_{j},k_{j}}$ 
(resp. $\delta_H^{E_{j'},k_{j'}}$) of the domain (resp. target) is induced from 
$\frak m^{{\rm f.c.u.} \frak b;j}$
(resp. $\frak m^{{\rm f.c.u.} \frak b;j'}$).
Note that in the domain we only take the terms `not greater' than $(E_{j'},k_{j'})$.
\end{rem}
\begin{prop}\label{prop818}
The chain map 
${\widehat{\frak p}}^{\rm f.c.u. \frak b}_{(E_{j},k_{j}),\infty}$\index[syindex]{phatfcubE(j)inf@${\widehat{\frak p}}^{\rm f.c.u. \frak b}_{(E_{j},k_{j}),\infty}$}
is 
chain homotopic to the chain map
${\widehat{\frak p}}^{{\rm f.c.u.} \frak b}_{(E_{j'},k_{j'}),\infty}\circ {\frak f}_{j',j}^{\infty, +}$
on ${CH}^{E_{j'},k_{j'},\infty}_*(\cL_{\rm c.u.}^{\rm form})^+$.\footnote{Here ${\widehat{\frak p}}^{\rm f.c.u. \frak b}_{(E_{j},k_{j}),\infty}$ 
is ${\widehat{\frak p}}^{\rm f.c.u. \frak b}_{n_j,k_{j},\infty}$ with $n_j$  determined by $E_j$.}
\end{prop}

We can use Proposition \ref{prop818} to complete the construction of 
$$
{\widehat{\frak p}}^{\rm f.c.u. \frak b} 
: CH_*(\cL_{\rm c.u.}^{\rm form})^+ \to \Omega(X) 
\widehat\otimes \Lambda_0
$$
as follows.
 %\marginpar{The construction of ${\widehat{\frak p}}^{\rm f.c.u. \frak b}$ is rewritten since 
%previous version has a problem.  KF 2025 Aug.}
\par
We recall from from our inductive construction that 
we have a  filtered $A_{\infty}$ structure 
$\frak m^{{\rm f.c.u.} \frak b;j;\#}$\index[syindex]{mxfrakfcub\#@$\frak m^{{\rm f.c.u.} \frak b;j;\#}$}
for each $j$ such that
the $A_{\infty}$ structure $\frak m^{{\rm f.c.u.} \frak b;j;\#}$
coincides (strongly isomorphic) with $A_{E_j,k_j}$ 
structure $\frak m^{{\rm f.c.u.} \frak b;j}$ as an 
$A_{E_j,k_j}$ 
structure.
(Note here these two $A_{E_j,k_j}$ structures coincide and not just weakly equivalent.)
Moreover we have a  quasi-isomorphism 
(of $A_{\infty}$ structures (filtered and homotopicaly unital))
 %\marginpar{change notation from $+$ to $\#$ since 
%$+$ is confusing with the process to add ${\bf f}$ and ${\bf e}^+$  KF 2025 Aug}
$$
{\frak f}_{j',j}^{\#} : 
(\cL_{\rm c.u.}^{\rm form},\frak m^{{\rm f.c.u.} \frak b;j,\#})
\to 
(\cL_{\rm c.u.}^{\rm form},\frak m^{{\rm f.c.u.} \frak b;j',\#})
$$
such that\index[syindex]{ffrakjj'*@${\frak f}_{j',j}^{\#}$}
$$
{\frak f}_{j',j}^{\#}\vert_{A_{E_{j'},k_{j'}}} = {\frak f}_{j',j}.
$$
Here the left hand side is ${\frak f}_{j',j}^{\#}$ regarded as 
an $A_{E_{j'},k_{j'}}$ functor.
It induces a chain homotopy equivalence
$$
CH({\frak f}_{j',j}^{\#,+}) : {CH}_*(\cL_{\rm c.u.}^{\rm form},\delta^j_H)^+
\to {CH}_*(\cL_{\rm c.u.}^{\rm form} ,\delta^{j'}_H)^+
$$
which extends (\ref{form1620}).  Here $\delta^j_H$ (resp. $\delta^{j'}_H$) 
is the Hochschild boundary of the $A_{\infty}$ structure $\frak m^{{\rm f.c.u.} \frak b;j,\#}$ 
(resp. $\frak m^{{\rm f.c.u.} \frak b;j',\#}$).

By construction this quasi-isomorphism is induced by smooth correspondences with appropriate 
Schwartz kernels. Therefore we can take a completion to obtain a chain homotopy equivalence
$$
CH({\frak f}_{j',j}^{\infty,\#,+}) : {CH}^{\infty}_{*}(\cL_{\rm c.u.}^{\rm form},\delta^j_H)^+
\to {CH}^{\infty}_*(\cL_{\rm c.u.}^{\rm form} ,\delta^{j'}_H)^+
$$
which extends (\ref{form1621new}).
\par
Then we can use Proposition \ref{prop818}
to obtain\index[syindex]{phatfcubE(j)E(j')inf@${\widehat{\frak p}}^{\rm f.c.u. \frak b}_{(E_j,k_j),(E_{j'},k_{j'}),\infty}$} 
$$
{\widehat{\frak p}}^{\rm f.c.u. \frak b}_{(E_j,k_j),(E_{j'},k_{j'}),\infty}
: CH^{\infty}_{(E_j,k_j)}(\cL_{\rm c.u.}^{\rm form}, \delta^{j'}_H)^+
\to \Omega(X) \widehat\otimes \Lambda_0
$$
such that 
$$
{\widehat{\frak p}}^{\rm f.c.u. \frak b}_{(E_j,k_j),(E_{j'},k_{j'}),\infty} = {\widehat{\frak p}}^{\rm f.c.u. \frak b}_{(E_{j'},k_{j'}),\infty}
$$
on $CH^{\infty}_{(E_{j'},k_{j'})}(\cL_{\rm c.u.}^{\rm form}, \delta^{j'}_H)^+$.
This is a simple algebraic fact. (See for example \cite[Proposition 19.33]{tech2}.)
\par
Now we take an increasing sequence $j(n)$ ($j(n) > j'$)  and use a similar induction as in the case or $\frak m$ or $\frak q$ 
to obtain 
$$
{\widehat{\frak p}}^{\rm f.c.u. \frak b}_{(E_{j(n)},k_{j(n)}),(E_{j'},k_{j'}),\infty} 
: CH^{\infty}_{(E_{j(n)},k_{(n)}j)}(\cL_{\rm c.u.}^{\rm form}, \delta^{j'}_H)^+
\to \Omega(X) \widehat\otimes \Lambda_0
$$
inductively such that for $n' < n$ 
$$
{\widehat{\frak p}}^{\rm f.c.u. \frak b}_{(E_{j(n')},k_{j(n')}),(E_{j'},k_{j'}),\infty} 
= 
{\widehat{\frak p}}^{\rm f.c.u. \frak b}_{(E_{j(n)},k_{j(n)}),(E_{j'},k_{j'}),\infty} 
$$
on $CH^{\infty}_{(E_{j(n)},k_{j(n)})}(\cL_{\rm c.u.}^{\rm form}, \delta^{j'}_H)^+$.
\par
We thus obtain\index[syindex]{phatfcubinf@${\widehat{\frak p}}^{\rm f.c.u. \frak b}_{\infty}$}
$$
{\widehat{\frak p}}^{\rm f.c.u. \frak b}_{\infty}
: CH^{\infty}(\cL_{\rm c.u.}^{\rm form}, \delta^{j'}_H)^+
\to \Omega(X) \widehat\otimes \Lambda_0
$$
by putting 
$$
{\widehat{\frak p}}^{\rm f.c.u. \frak b}_{\infty}
=
{\widehat{\frak p}}^{\rm f.c.u. \frak b}_{(E_{j(n)},k_{j(n)}),(E_{j'},k_{j'}),\infty}
$$
on $CH^{\infty}_{(E_{j(n)},k_{j(n)})}(\cL_{\rm c.u.}^{\rm form}, \delta^{j'}_H)^+$.
\par
We finally obtain ${\widehat{\frak p}}^{\rm f.c.u. \frak b}$ by restricting ${\widehat{\frak p}}^{\rm f.c.u. \frak b}_{\infty}$
to the (algebraic) tensor product.
\par
Thus to complete the construction of $\widehat{\frak p}^{\rm f.c.u. \frak b}$ 
it suffices to prove Proposition \ref{prop818}.
\begin{rem}
We remark that in the construction of closed-open map $\frak q$, we use a variant of the Hochschild
cochain complex, that is, those consisting of elements with  Schwartz kernels satisfying Condition \ref{cod149}.
This is a {\it subcomplex} of the algebraic Hochschild
cochain complex.
Here, during the construction of open-closed map $\frak p$, we use a variant of 
Hochschild chain complex, where we replace the algebraic tensor product $\Omega(L)\otimes\Omega(L')$
of $\Omega(L)$ and $\Omega(L')$ by its completion 
$\Omega(L \times L')$.
So our variant {\it contains} the algebraic Hochschild chain complex as a subcomplex. 
In other words, the direction of inclusion is opposite. This fact is consistent with the fact 
that Hochschild
cochain complex is the target of the closed-open map 
and Hochschild chain complex is the domain of the open-closed map.
\end{rem}

\subsection{Proof of Proposition \ref{prop818}}
We start with the following remark.
The Hochschild cochain complex $CH^*(\cL_{\rm c.u.}^{\rm form})^+$ is a 
differential graded Lie algebra 
whose differential is the Hochschild boundary and Lie algebra structure 
is the Gerstenhaber bracket.
The Hochschild chain complex 
$CH_*(\cL_{\rm c.u.}^{\rm form})^+$
is a differential graded module over the differential graded Lie algebra 
$CH^*(\cL_{\rm c.u.}^{\rm form})^+$.
In fact the operation is defined by
\begin{equation}\label{new828}
\aligned
\varphi \cdot {\bf x} 
= 
&\sum_{0< i < j \le k} (-1)^{\maltese_1} x_0 \otimes \dots \otimes \varphi_{j-i+1}(x_i, \dots, x_j) \otimes \dots
\otimes x_k \\
&+ \sum_{0\le i < j \le k} (-1)^{\maltese_2} \varphi_{k-j+i}(x_j, \dots, x_0,\dots,x_i)\otimes x_{i+1} \otimes \dots \otimes x_{j-1}.
\endaligned
\end{equation}
Here ${\bf x} = x_0 \otimes x_1 \otimes \dots \otimes x_k$,
$\varphi = (\varphi_i)_{i=1,2\dots}$ and
$$
\aligned
\maltese_1 &= \deg \varphi_{j-i+1}(\deg' x_0 + \dots + \deg' x_{i-1}), \\ 
\maltese_2 &= (\deg' x_j + \dots + \deg' x_k)(\deg' x_0 + \dots + \deg' x_{j-1})
\\ &\qquad + \deg \varphi_{k-j+i} (\deg' x_j + \dots + \deg'x_i).
\endaligned
$$
%(See \cite[Proposition 7.4.184]{fooo092}.)
\begin{rem}
This differential graded module structure is different from the one used in Lemma
\ref{lem:cap_product_HH}. In fact the latter one uses the structure operations of 
$A_{\infty}$ structure in the definition. On the other hand, (\ref{new828}) 
does not contain structure operations of the 
$A_{\infty}$ structure.
\end{rem}
In the same way, $CH^*_{n,k}(\cL_{\rm c.u.}^{\rm form})$ is a differential graded Lie algebra 
and $CH_*^{n,k}(\cL_{\rm c.u.}^{\rm form})$ is a differential graded module over $CH^*_{n,k}(\cL_{\rm c.u.}^{\rm form})$.
(The same holds for $+$ verions.)

We regard $\Omega(X) \widehat\otimes \Lambda_G$ as an $L_{\infty}$ algebra 
(where the differential is the de Rham differential and all the other 
structure operations are trivial)
and $\widehat{\frak q}$ induces an 
$L_{\infty}$ homomorphism\index[syindex]{qfcubellnk@$\widehat{\frak q}^{\rm f.c.u.\frak b}_{\ell,n,k}$}
$$
\widehat{\frak q}^{\rm f.c.u.\frak b}_{\ell,n,k}: \Omega(X)^{\otimes \ell} \widehat\otimes \Lambda_G \to CH^*(\cL_{\rm c.u.}^{\rm form})^+,
$$
by which we regard $CH_*(\cL_{\rm c.u.}^{\rm form})^+$ 
as an $L_{\infty}$ module over $\Omega(X) \widehat\otimes \Lambda_G$.
\par
In the same way $CH_*^{n,k}(\cL^{\rm form}_{\rm c.u.})^+$
is an  $L_{\ell}$ module over $\Omega(X) \widehat\otimes \Lambda_G$.
\begin{rem}
We say $CH_*^{n,k}(\cL_{\rm c.u.}^{\rm form})^+$ is an $L_{\ell}$ module rather than an $L_{\infty}$ module. In fact, during the inductive construction 
of structures, we consider moduli spaces of pseudo-holomorphic polygons, the number of whose interior marked 
points is not greater than $\ell$.
\end{rem}
\par
Then the map $\widehat{\frak p}_{n,k}^{\rm f.c.u.\frak b}$ is induced from the 
$L_{\ell}$ module homomorphism\index[syindex]{pellEnkfcub@${\frak p}_{\ell;E_n,k}^{\rm f.c.u.\frak b}$}
\begin{equation}\label{form1520}
{\frak p}_{\ell;E_n,k}^{\rm f.c.u.\frak b} : \Omega(X)^{\otimes \ell} \otimes  CH_*^{n,k}(\cL_{\rm c.u.}^{\rm form})^+
\to \Omega(X) \otimes \Lambda_0/T^{E'_n}.
\end{equation}
Here  $L_{\infty}$ module structure of $\Omega(X)$ over   $\Omega(X)$ is trivial.  (Note that $L_{\infty}$ algebra structure of $\Omega(X)$ 
was also trivial.)
The statement that (\ref{form1520}) is an $L_{\infty}$ module homomorphism is the next equality.
\begin{equation}\label{form1521}
\aligned
 d (\frak p({\bf w};x_0 \otimes {\bf x})) =&\sum_{c,c'} (-1)^{\maltese_1} \frak p({\bf w}_{c}^{2;1};x_0 \otimes {\bf x}_{c'}^{3;1} \otimes \frak q({\bf w}_{c}^{2;2}; {\bf x}_{c'}^{3;2})  \otimes {\bf x}_{c'}^{3;3} ) \\
&+ \sum_{c,c'} (-1)^{\maltese_2} \frak p({\bf w}_{c}^{2;1};\frak q({\bf w}_{c}^{2;2};  {\bf x}_{c'}^{3;3} \otimes x_0 \otimes {\bf x}_{c'}^{3;1} )\otimes {\bf x}_{c'}^{3;2}).
\endaligned
\end{equation}
Here the signs are by Koszul rule. We omit suffix such as ${\ell,E_n,k}$ in the formula for simplicity.
We define ${\bf w}_{c}^{2;1}$ by $\Delta({\bf w}) = \sum_c {\bf w}_{c}^{2;1} \otimes {\bf w}_{c}^{2;2}$,
where $\Delta$ is the {\it co-comutative} coalgebra structure on $\oplus_{\ell} \Omega(X)^{\otimes \ell}$.
The symbols ${\bf x}_{c'}^{3;1}$ etc. are defined by $((\Delta \otimes {\rm id}) \circ \Delta)({\bf x}) = \sum_{c'} {\bf x}_{c'}^{3;1} \otimes {\bf x}_{c'}^{3;2} \otimes {\bf x}_{c'}^{3;3}$,
where $\Delta$ is the non-cocommutative colagebra structure.
\par
Note that in case when ${\bf w}=1$ the right hand side becomes $\frak p(1;\delta^H(x_0 \otimes {\bf x}))$.
Therefore, Formula (\ref{form1520}) implies that $\frak p$ is a chain map.
%\begin{rem}
%This is similar to \cite[Theorem 7.4.192]{fooo092}.
%The cyclic chain complex (instead of the Hochschild chain complex) 
%was studied in \cite{fooo092}. The situation there is a bit more complicated than here 
%because of the phenomenon that boundary loop shrinks to a point.
%(See \cite[Subsubsection 3.8.3, Proposition 3.8.27]{fooo09}). In our situation of Hochschild chain complex, we 
%use moduli spaces of pseudo-holomorphic disks with at least one boundary 
%marked point. So this phenomenon occurs in codimension 2 and is not an issue.
%\end{rem}

We can take the completion (with respect to the Fr\'echet space structure of the space of smooth 
forms) to obtain\index[syindex]{qfcubnkinfhat@${\widehat{\frak q}}^{\rm f.c.u.\frak b}_{n,k,\infty}$}
$$
{\widehat{\frak q}}^{\rm f.c.u.\frak b}_{n,k,\infty} : \Omega(X) \to CH_{\infty;n,k}^*(\cL_{\rm c.u.}^{\rm form})^+,
$$
and\index[syindex]{pfcubellEnkinfhat@${\widehat{\frak p}}_{\ell;E_n,k,\infty}^{\rm f.c.u.\frak b}$}
$$
{\widehat{\frak p}}_{\ell;E_n,k,\infty}^{\rm f.c.u.\frak b} : \Omega(X)^{\otimes \ell} \otimes  {CH}_*^{n,k,\infty}(\cL_{\rm c.u.}^{\rm form})^+
\to \Omega(X) \otimes \Lambda_0/T^{E'_n}.
$$
Here $CH_{\infty;n,k}^*(\cL^{\rm form})^+$ is the Hochschild cochain complex consisting of linear maps 
with Schwartz kernels satisfying Condition \ref{cod149}  which we introduced in Subsection \ref{frakqondeRham}.
\par
Again using the fact that operations have Schwartz kernels satisfying appropriate condition on its wave front set we have an operation 
\begin{equation}\label{form1625new}
CH_{\infty;n,\vec{\kappa}}^*(\cL_{\rm c.u.}^{\rm form})^+ 
\otimes {CH}_{\vec{\kappa},\infty}(\cL_{\rm c.u.}^{\rm form})^+
\to \cL^{\rm form}(L_{\kappa_0},L_{\kappa_k};\F)\widehat{\otimes} \Lambda_0
\end{equation}
which defines the structure of a differential graded module. Here 
$CH_{\infty;n,\vec{\kappa}}^*(\cL_{\rm c.u.}^{\rm form})^+$ is a completion of 
$CHCF_*(\cL^{\rm form};\vec{\kappa})^+$. (See (\ref{CHCFform}).)
\par
The operation (\ref{form1625new}) is the completed version of differential graded module structure of $CH_*^{n,k}(\cL_{\rm c.u.}^{\rm form})^+$ over $CH^*_{n,k}(\cL_{\rm c.u.}^{\rm form})^+$.
\par
\begin{figure}[h]
\includegraphics[scale=0.5]{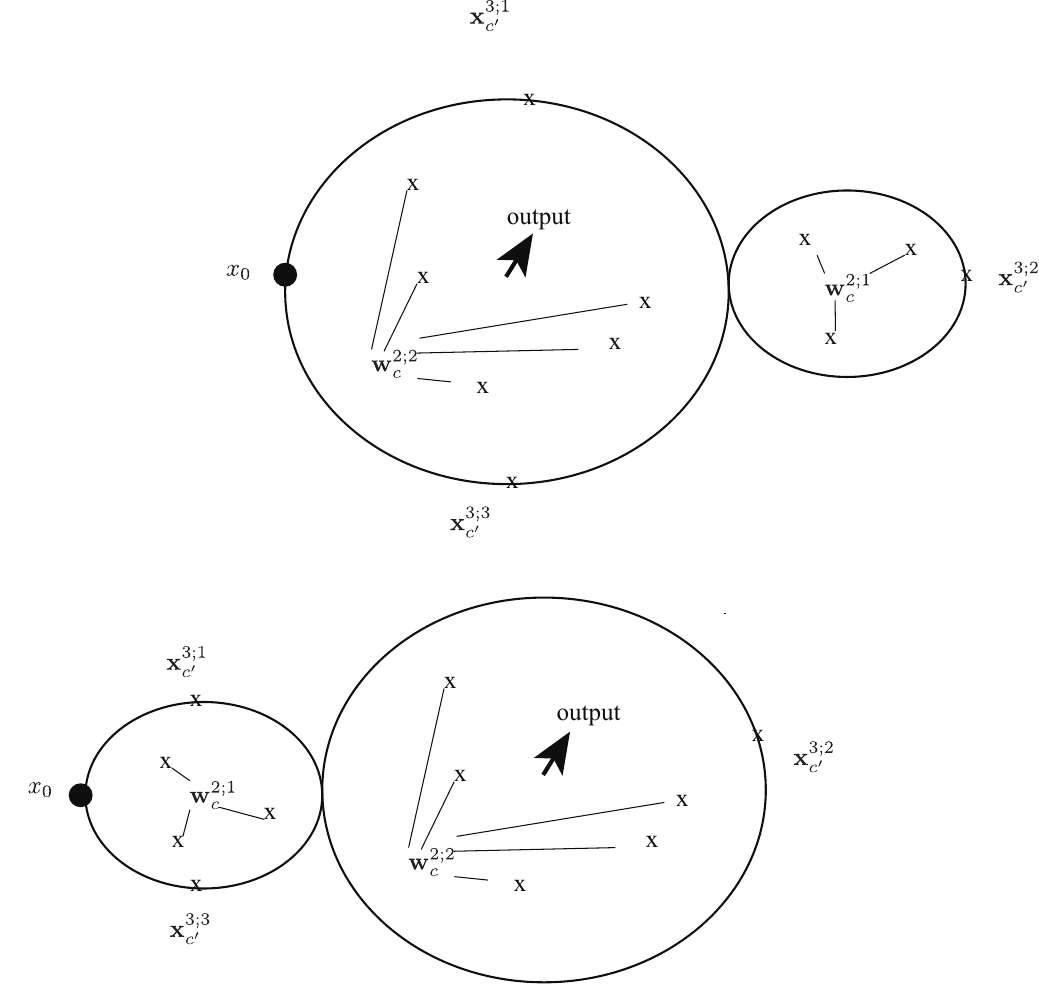}
\caption{$\frak p$ is an $L_{\infty}$ module homomorphism}
\label{FigureLinfp}
\end{figure}

For $(E_{n'},k') < (E_{n},k)$ we obtained a pseudo-isotopy and the following commutative diagram of $L_{\ell}$ homomorphism
in Lemma \ref{lem810}.  See Figure \ref{FigureLinfp}
\begin{equation}\label{diag828}
\xymatrix{ 
&&&&
CH^*_{\infty;n,k}(\cL^{\rm form}_{\frak b;n,k},\cL^{\rm form}_{\frak b;n,k})^+ \\
\Omega(X)    
\ar[urrrr]^{{\widehat{\frak q}}^{\rm f.c.u.\frak b}_{n,k}} 
\ar[rrrr]_{\widehat{\frak q}^{\frak C,\frak b}_{o;n',k'}}
\ar[drrrr]_{{\widehat{\frak q}}^{\rm f.c.u.\frak b}_{n',k'}} 
&&&&CH_{\infty;n,k}^*(\cL ^{\rm form}_{\frak b;n',k'},(\cL ^{\rm form}_{\frak b;n',k'})_{[0,1]})^+
\ar[u]_{{\rm Ev}_{t=1}}
\ar[d]^{{\rm Ev}_{t=0}}
\\
&&&&
CH_{\infty;n,k}^*(\cL ^{\rm form}_{\frak b;n',k'},\cL ^{\rm form}_{\frak b;n',k'}) ^+
}
\end{equation}
Note that here the maps $\widehat{\frak q}$ and the $A_{n,k}$ structures  with suffix $n', k'$ 
in Diagram (\ref{diag828}) denote those which are promoted to $A_{n,k}$ structures.

\begin{lem}\label{lem821}
There exists a commutative diagram
\begin{equation}\label{diag829}
\xymatrix{ 
{CH}^{n',k',\infty}_*(\cL^{\rm form}_{\frak b;n,k})^+
\ar[rrr]^{{\widehat{\frak p}}^{\rm f.c.u.\frak b}_{\ell,E_{n},k,\infty}}
&&&\Omega(X) \otimes \Lambda_0/T^{E'_{n'}}\\
{CH}^{n',k',\infty}_*((\cL ^{\rm form}_{\frak b;n',k'})_{[0,1]})^+
\ar[u]_{{\rm Ev}_{t=1}}
\ar[d]^{{\rm Ev}_{t=0}}
\ar[rrr]_{{\widehat{\frak p}}^{\frak C,\frak b}_{\ell,E_{n'},k',\infty}}
&&&\Omega([0,1] \times X) \otimes \Lambda_0/T^{E'_{n'}}\ar[u]_{t=1}\ar[d]^{t=0}
\\
{CH}_*^{n',k',\infty}(\cL ^{\rm form}_{\frak b;n',k'})^+
\ar[rrr]_{{\widehat{\frak p}}^{\rm f.c.u.\frak b}_{\ell,E_{n'},k',\infty}}
&&&\Omega(X) \otimes \Lambda_0/T^{E'_{n'}}
}
\end{equation}
with the following properties.

\begin{enumerate}
\item
${CH}^{n',k',\infty}_*((\cL ^{\rm form}_{\frak b;n',k'})_{[0,1]})^+$ 
and ${CH}^{n',k',\infty}_*((\cL ^{\rm form}_{\frak b;n,k})_{[0,1]})^+$
are differential graded modules over the Lie algebra
$CH_{\infty;n'k'}^*(\cL ^{\rm form}_{\frak b;n',k'},(\cL ^{\rm form}_{\frak b;n',k'})_{[0,1]})^+$.
\item
The three objects in the left hand side of Diagram $(\ref{diag829})$ are regarded as 
filtered $L_{\ell}$ modules over $\Omega(X)$ by the three 
morphisms which are horizontal arrows of Diagram $(\ref{diag828})$.
\item 
The horizontal arrow ${\widehat{\frak p}}^{\rm f.c.u.\frak b}_{\ell,E_n,k,\infty}$,
(resp. ${\widehat{\frak p}}^{\rm f.c.u.\frak b}_{E_{n'},k',\infty}$) is
a differential graded Lie module homomorphism over 
$CH^*_{\infty;n,k}(\cL^{\rm form}_{\frak b;n,k},\cL^{\rm form}_{\frak b;n,k})^+$ 
(resp. $CH^*_{\infty;n',k'}(\cL^{\rm form}_{\frak b;n',k'},\cL^{\rm form}_{\frak b;n',k'})^+$) as explained above.
We regard it as an $L_{\ell}$ module homomorphism over $\Omega(X)$ by the first (resp. third) horizontal 
arrow of Diagram $(\ref{diag828})$.
%\item
%The first horizontal arrow is one 
\item The middle horizontal arrow ${\widehat{\frak p}}^{\frak C,\frak b}_{\ell,E_{n'},k',\infty}$ \index[syindex]{phathatfrakCb@${\widehat{\frak p}}^{\frak C,\frak b}_{E_{n'},k',\infty}$} 
is a differential graded Lie module homomorphism over
$CH_{\infty;n',k'}^*(\cL ^{\rm form}_{\frak b;n',k'},(\cL ^{\rm form}_{\frak b;n',k'})_{[0,1]})^+$.
We regard it as an $L_{\ell}$ module homomorphism over $\Omega(X)$ by the second horizontal 
arrow of Diagram $(\ref{diag828})$.
\item
The Diagram $(\ref{diag829})$ commutes as a diagram of $L_{\ell}$ module homomorphisms over $\Omega(X)$.
\end{enumerate}
\end{lem}
\begin{proof}[Proof of Lemma \ref{lem821} $\Rightarrow$ Proposition \ref{prop818}]
We conider $t$-parameter family version of Equation (\ref{form1521}) and proceed in the same way as 
the proof of Lemma \ref{lem810}. Since the argument is completely parallel we do not repeat it.

\begin{proof}[Proof of Lemma \ref{lem821}]
Since the proof is similar to the proofs of other similar statements given before, our proof here will be somewhat sketchy.
\par
We first define
${CH}_*^{n',k',\infty}((\cL ^{\rm form}_{\frak b;n',k'})_{[0,1]})^+$.
For $\vec{\kappa} = (\kappa_0,\dots,\kappa_K)$
we consider 
$$
L(\vec{\kappa}) = \prod_{j=1}^{K+1} (\tilde L_{\kappa_{j-1}} \times_X \tilde L_{\kappa_{j}})
$$
and $\Omega([0,1] \times L(\vec{\kappa}))$.
Then we define ${CH}^{n',k',\infty}_*((\cL ^{\rm form}_{\frak b;n',k'})_{[0,1]})^+$ 
to be the set of sums\index[syindex]{CHnkinftyLform01@${CH}^{n',k',\infty}_*((\cL ^{\rm form}_{\frak b;n',k'})_{[0,1]})^+$}
$$
\sum_i T^{\lambda_i} \frak T_{i}
$$
where $\lambda_i \in G$, $\frak T_i \in  \Omega([0,1] \times L(\vec{\kappa}(i))$  and
$(\lambda_i,\vert\vec\kappa(i)\vert) \le (E_{n'},k')$.
(Here $\vec\kappa(i) = (\kappa_1(i),\dots, \kappa_{\vert\vec\kappa(i)\vert}(i))$.)
(We also include ${\bf e}^+$ and ${\bf f}$.)
It is a completion of the filtered module ${CH}^{n',k'}_*((\cL ^{\rm form}_{\frak b;n',k'})_{[0,1]})^+$.
Its Hocshchild boundary is also defined by  completion.

%\begin{equation}\label{form833}
%\widehat{CH}^{n,k}_*(\cL ^{\rm form} \times [0,1])= 
%\bigoplus_{\vec{\kappa}} \Omega([0,1] \times L(\vec{\kappa})) \widehat\otimes \Lambda_0.
%\end{equation}
%Here $\widehat\otimes$ is the $T$ adic completion of the algebraic tensor product.
%In other words an element of 
%$ \Omega([0,1] \times L(\vec{\kappa})) \widehat\otimes \Lambda_0$ is a sum
%$$
%\sum_i T^{\lambda_i} h_i
%$$
%where $\lambda_i \ge 0$, $\lim_{i\to \infty} \lambda_i = + \infty$ and 
%$h_i \in \Omega([0,1] \times L(\vec{\kappa}))$.
%The sum over $\vec{\kappa}$ in (\ref{form833}) is the direct sum.
\par
We next define ${\widehat{\frak p}}^{\frak C,\frak b}_{E_{n'},k',\infty}$.
We consider the direct product\index[syindex]{phathatfrakCb@${\widehat{\frak p}}^{\frak C,\frak b}_{E_{n'},k',\infty}$} 
$$
[0,1]_t \times  [0,1]^{\vert \frak f_{\bf f}\vert} \times{\mathcal M}_{\ell+1;\vec k}((\vec{\kappa},\vec p);B;\vec{\frak f})^{\boxplus 1}
$$
of the moduli space in 
Proposition \ref{Kuraeistspoly2} and the interval $[0,1]$ (with coordinate $t$).
By Propositions \ref{Kuraeistspoly2} and  \ref{existmultipolu2} 
we have constructed Kuranishi structures and CF-perturbations on the moduli spaces ${\mathcal M}_{\ell+1;\vec k}((\vec{\kappa},\vec p);B)^{\boxplus 1}$.
More precisely speaking, for $(E_{n'},k')$, we take  
Kuranishi structures and CF-perturbations on $ [0,1]^{\vert \frak f_{\bf f}\vert} \times{\mathcal M}_{\ell+1;\vec k}((\vec{\kappa},\vec p);B;\vec{\frak f})^{\boxplus 1}$
for $(\omega(B),\vert\vec\kappa\vert) \le (E_{n'},k')$ and $\ell \le \ell_0$.\footnote{Here $\ell_0$ is taken sufficiently larger than 
a number depending on 
$(E_{n'},k')$.}  We denote these Kuranishi structures and CF-perturbations by $\Xi_{(E_{n'},k')}$.
Let $i < j$ and $\omega(B) \le E_i$. 
In the same way as Propositions \ref{Kuraeistspoly2} and  \ref{existmultipolu2} 
we can construct a Kuranishi structure on 
$[0,1]_t \times  [0,1]^{\vert \frak f_{\bf f}\vert} \times{\mathcal M}_{\ell+1;\vec k}((\vec{\kappa},\vec p);B;\vec{\frak f})^{\boxplus 1}$\index[syindex]{01t01fMell+1kappa
@$[0,1]_t \times  [0,1]^{\vert \frak f_{\bf f}\vert} \times{\mathcal M}_{\ell+1;\vec k}((\vec{\kappa},\vec p);B;\vec{\frak f})^{\boxplus 1}$}
with the following properties. We call this Kuranishi structure the $\frak p$-para Kuranishi structure. 
\begin{enumerate}
\item
Its boundary is the union of the following types of moduli spaces: 
\begin{enumerate}
\item[(*)]The direct product of boundary of types 1,2,3,4,5 in Proposition \ref{Kuraeistspoly2}
and $[0,1]_t$. 
\item[(**)] $\{0,1\}  \times  [0,1]^{\vert \frak f_{\bf f}\vert} \times{\mathcal M}_{\ell+1;\vec k}((\vec{\kappa},\vec p);B;\vec{\frak f})^{\boxplus 1}$.
\end{enumerate}
\item
On $\{0\} \times [0,1]^{\vert \frak f_{\bf f}\vert} \times{\mathcal M}_{\ell+1;\vec k}((\vec{\kappa},\vec p);B;\vec{\frak f})^{\boxplus 1}$
it coincides with $\Xi_{i}$ and on 
$ \{1\} \times [0,1]^{\vert \frak f_{\bf f}\vert} \times{\mathcal M}_{\ell+1;\vec k}((\vec{\kappa},\vec p);B;\vec{\frak f})^{\boxplus 1}$
it coincides with $\Xi_{j}$.
\item
On the boundary given as the product of type 1 boundary in Proposition \ref{Kuraeistspoly2}
and $[0,1]_t$, the Kuranishi structure is given as the fiber products 
\begin{equation}\label{21822}
\aligned
&([0,1]_t \times [0,1]^{\vert \frak f_{\bf f}^2\vert} \times{\mathcal M}_{\ell''+1; k''_i+1}(L;\beta;\vec{\frak f}^2)^{\boxplus 1})^{\frak p} \\
&
\,{}_{(\text{\rm ev}_t,\text{\rm ev}_0)} \times_{(\text{\rm ev}_t,\text{\rm ev}_{i,j})}
([0,1]_t \times  [0,1]^{\vert \frak f_{\bf f}^1\vert}\times {\mathcal M}_{\ell';\vec k'}((\vec{\kappa},\vec p);B';\vec{\frak f}^1)^{\boxplus 1})^{\frak q},
\endaligned
\end{equation}
\begin{equation}\label{2202120}
\aligned
&([0,1]_t \times  [0,1]^{\vert \frak f_{\bf f}^2\vert} \times {\mathcal M}_{\ell'';\vec k''}(L;\beta;\vec{\frak f}^2))^{\boxplus 1})^{\frak q}\\
&\,{}_{(\text{\rm ev}_t,\text{\rm ev}_{1,m''})}\times_{(\text{\rm ev}_t,\text{\rm ev}_{i,m'})}
([0,1]_t \times [0,1]^{\vert \frak f_{\bf f}^1\vert}\times {\mathcal M}_{\ell'+1;\vec k'}((\vec{\kappa},\vec p);B';\vec{\frak f}^1))^{\boxplus 1})^{\frak p},
\endaligned
\end{equation}
over $[0,1]_t \times L_{\kappa}$ for some $\kappa$.
Here the Kuranishi structures on the second factor of (\ref{21822}) 
and the first factor of (\ref{2202120}) (and one parameter version of 
Proposition \ref{prop94})
 are the ones in Lemma \ref{Kuraonepara}
and the Kuranishi structures on the first factor of (\ref{21822}) 
and the second factor of (\ref{2202120}) are the $\frak p$-para Kuranishi structures.
\item
On the boundary given as the product of type 2 boundary in Proposition \ref{Kuraeistspoly2}
and $[0,1]_t$, the Kuranishi structure is given as the fiber products 
\begin{equation}\label{21912}
\aligned
&([0,1]_t\times [0,1]^{\vert \frak f_{\bf f}^2\vert} \times {\mathcal M}_{\ell''+1;\vec k''}((\vec{\kappa}'',\vec p'');B'';\vec{\frak f}^2
))^{\boxplus 1})^{\frak p}
 \\
&\times_{[0,1]_t}([0,1]_t\times [0,1]^{\vert \frak f_{\bf f}^1\vert} \times {\mathcal M}_{\ell';\vec k'}((\vec{\kappa}',\vec p');B';\vec{\frak f}^1))^{\boxplus 1})^{\frak q},
\endaligned
\end{equation}
\begin{equation}\label{219221}
\aligned
&([0,1]_t \times  [0,1]^{\vert \frak f_{\bf f}^2\vert} \times {\mathcal M}_{\ell'';\vec k''}((\vec{\kappa}'',\vec p'');B'';\vec{\frak f}^2))^{\boxplus 1})^{\frak q}
\\
&\times_{[0,1]_t}  ([0,1]_t\times [0,1]^{\vert \frak f_{\bf f}^1\vert} \times{\mathcal M}_{\ell'+1;\vec k'}((\vec{\kappa}',\vec p');B'
;\vec{\frak f}^1))^{\boxplus 1})^{\frak p},
\endaligned
\end{equation}
over $[0,1]_t$.
Here the Kuranishi structures on the second factor of (\ref{21912}) and the first factor of 
(\ref{219221})  are the ones in Lemma \ref{Kuraonepara} (and one parameter version of 
Proposition \ref{prop94})
and the Kuranishi structures on the first factor of (\ref{21912}) and 
the second factor of (\ref{219221}) are the $\frak p$-para Kuranishi structures.
\item
On the boundary given as the product of type 3 boundary in Proposition \ref{Kuraeistspoly2}
and $[0,1]_t$, the Kuranishi structure is given as the fiber products 
\begin{equation}\label{220121}
\aligned
&([0,1]_t\times  [0,1]^{\vert \frak f_{\bf f}^2\vert} \times {\mathcal M}_{\ell''+1;\vec k''}((\vec{\kappa}'',\vec p'');B'';
\vec{\frak f}^2))^{\boxplus 1})^{\frak p} \\
&\,{}_{(\text{\rm ev}_t,\text{\rm ev}_{1,m''})}\times_{(\text{\rm ev}_t,\text{\rm ev}_{i,m'})}
([0,1]_t\times [0,1]^{\vert \frak f_{\bf f}^1\vert} \times {\mathcal M}_{\ell';\vec k'}((\vec{\kappa}',\vec p');B';\vec{\frak f}^1))^{\boxplus 1})^{\frak q},
\endaligned
\end{equation}
\begin{equation}\label{220212}
\aligned
&([0,1]_t\times [0,1]^{\vert \frak f_{\bf f}^2\vert}\times {\mathcal M}_{\ell'';\vec k''}((\vec{\kappa}'',\vec p'');B'';\vec{\frak f}^2))^{\boxplus 1})^{\frak q}\\
&\,{}_{(\text{\rm ev}_t,\text{\rm ev}_{1,m''})}\times_{(\text{\rm ev}_t,\text{\rm ev}_{i,m'})}
([0,1]_t\times  [0,1]^{\vert \frak f_{\bf f}^1\vert}\times {\mathcal M}_{\ell'+1;\vec k'}((\vec{\kappa}',\vec p');B';\vec{\frak f}^1))^{\boxplus 1})^{\frak p},
\endaligned
\end{equation}
over $[0,1]_t \times L_{\kappa}$ for some $\kappa$.
Here the Kuranishi structures on the second factor of (\ref{220121}) and the first factor of 
(\ref{220212})  are the ones in Lemma \ref{Kuraonepara}
and the Kuranishi structures on the first factor of (\ref{220121}) and 
the second factor of (\ref{220212})   are the $\frak p$-para Kuranishi structures.

\item
The evaluation map 
$$
({\rm ev}_t,{\rm ev}^+_0): ([0,1]_t \times  [0,1]^{\vert \frak f_{\bf f}\vert} \times {\mathcal M}_{\ell+1;\vec k}((\vec{\kappa},\vec p);B;\vec{\frak f})^{\boxplus 1})^{\frak p}
\to [0,1]_t \times  X 
$$
is weakly submersive.
\end{enumerate}
This Kuranishi structure is orientable and have a similar property as 
Proposition \ref{Kuraeistspoly2} (3),(4),(5),(7),(8),(9).
\par
We then define a CF-perturbation, which we call the $\frak p$-para CF perturbations, on the $\frak p$-para Kuranishi structures.
We require that it is consistent with the fiber product description 
of its boundary. 
Namely on the boundary given as the product of types 1,2,3 boundary in Proposition \ref{Kuraeistspoly2}
and $[0,1]$ the restriction of the $\frak p$-para CF perturbation coincides with the fiber product 
CF-perturbation defined on the 
Kuranishi structures mentioned in (3), (4), (5) above.
We also require that on the boundary component $[0,1]_t \times [0,1]^{\vert \frak f_{\bf f}\vert} \times{\mathcal M}_{\ell+1;\vec k}((\vec{\kappa},\vec p);B;\vec{\frak f})^{\boxplus 1}$ labelled by $i \in \{0,1\}$, 
the restriction of the outer collaring of the $\frak p$-para CF perturbation coincides with $\Xi_{E_n,k}$.
\par
We also require that the evaluation map 
$$
({\rm ev}_t,{\rm ev}^+_0):  ([0,1]_t  \times  [0,1]^{\vert \frak f_{\bf f}\vert} \times {\mathcal M}_{\ell+1;\vec k}((\vec{\kappa},\vec p);B;\vec{\frak f})^{\boxplus 1})^{\frak p}
\to [0,1]  \times  X
$$
is strongly submersive with respect to the $\frak p$-para CF perturbation. 
\par
The construction of such a system of Kuranishi structures can be carried out 
in the same way as \cite{const2} and Sections \ref{sec:Kuraconst},\ref{sec:CRperturb} for example and the 
construction of such a system of  CF-perturbations 
can be carried out in exactly the same way as the construction in 
\cite[Proof of Proposition 19.1]{tech2} and Section \ref{sec:CRperturb} for example.
\par
Now we consider the following smooth correspondence:\footnote{For this purpose we need to 
add a data to specify the $0$-th marked point, and then to change the notation 
to one in Subsection \ref{constcyclic} in the same way we did several times.}
$$
\xymatrix{ 
&&  [0,1] \times  {\mathcal M}_{\ell+1}(\vec{\kappa},\vec p;B)^{\boxplus 1}
\ar[rd]_{({\rm ev}^+_0,{\rm ev}_t)}
\ar[ld]^{\,\,\,\,\,\,\,\,\,({\rm ev}_t,({\rm ev}^+_{1},\dots,{\rm ev}^+_{\ell}), {\rm ev}^{\partial})}
\\
& [0,1] \times 
X^{\ell} \times L(\vec{\kappa})  &&  [0,1] \times  X 
}
$$
It induces a map
$$
\Omega(X)^{\otimes\ell} \otimes \Omega([0,1] \times L(\vec{\kappa})))
\to \Omega([0,1]  \times  X)
$$
We use this map in the same way as for the definition of 
${\widehat{\frak p}}^{\rm f.c.u.\frak b}_{E_{n'},k',\infty}$ 
to obtain ${\widehat{\frak p}}^{\frak C,\frak b}_{E_{n'},k',\infty}$.
(We can include ${\bf e}^+$ and ${\bf f}$ by using the evaluation
maps labeled by ${\frak f}_{\bf e}$, ${\frak f}_{\bf f}$
in the same way as Subsection \ref{subsec:opeartor}.)
\par
The fact that ${\widehat{\frak p}}^{\frak C,\frak b}_{E_{n'},k',\infty}$
is an $L_{\ell}$ module homomorphism is a consequence 
of the compatibility at the boundary of the Kuranishi structures 
and the CF perturbations (at the boundaries 
described by items (3), (4), (5) above), together with Stokes' theorem and composition 
formula.
The commutativity of Diagram (\ref{diag829}) is a consequence 
of the compatibility at the boundary of Kuranishi structure 
and CF perturbation described by item (2) above.
We can perform induction on $\ell$ in the same way.
The proof of Lemma \ref{lem821} is now complete.
\end{proof}
The proof of Proposition \ref{prop818} is finally complete.
\end{proof}
We have thus obtained\index[syindex]{paformcubbf@$\widehat{\frak p}^{\rm f.c.u. \frak b}$} 
\begin{equation}
\widehat{\frak p}^{\rm f.c.u. \frak b} 
: CH_*(\cL^{\rm form}_{\rm c.u.})^+ \to \Omega(X) 
\widehat\otimes \Lambda_0.
\end{equation}
We include bounding cochain and define\index[syindex]{pfubfb@$\widehat{\frak p}^{\rm f.u. {\bf b}}$}
\begin{equation}
\widehat{\frak p}^{\rm f.u. {\bf b}} 
: CH_*(\cL^{\rm form}_{\rm uni}) \to \Omega(X) 
\widehat\otimes \Lambda_0.
\end{equation}
by
\begin{equation}
\widehat{\frak p}^{\rm f.u.\text{\bf b}}_{\vec{\kappa}}
(\text{\bf h})
=\sum_{\vec k}
\widehat{\frak p}^{\rm f.c.u. \frak b}_{\ell,\vec{\kappa} \sqcup \vec k}
(b_{\kappa_0,+}^{\otimes k_0}, h_0, \ldots, h_K ,b_{\kappa_K,+}^{\otimes k_K}).
\end{equation}
It is easy to see that $\widehat{\frak p}^{\rm f.u. {\bf b}}$ is a chain map.
 
We have thus proved a de Rham complex version of Theorem \ref{themp}. We have also proved the next proposition 
at the same time.

\begin{prop}\label{prop1522}
The maps $\widehat{\frak p}^{\rm f.c.u. \frak b}$ and  $\widehat{\frak p}^{\rm f.u. {\bf b}}$ 
are the linear parts of $L_{\infty}$ module homomorphisms.
\end{prop}

We used the homotopy inverse limit  to construct a chain map by utilizing 
the fact that the existence of the promotion of the chain complex depends only on the chain 
homotopy type.  To prove Proposition \ref{prop1522} we use appropriate obstruction theory 
to show that the existence of the promotion of the $L_{\ell}$ module structure depends only on 
homotopy type.  We omit the detail here and referring readers to \cite[Section 7.2.12]{fooo092}. %\marginpar{A sentence added.  KF 2025 Aug.}

\subsection{Going to the canonial model.}
In this subsection we  construct a similar map $\widehat{\frak p}^{\rm f.c.u. \frak b}$ from the canonical model.
We put
$$
CH_*(\cL) = CH_*(\cL,\cL)
, \quad
CH_*(\cL_{\rm uni}) = CH_*(\cL_{\rm uni},\cL_{\rm uni}),
$$
where the right hand side is defined in Subsection \ref{sec:hochschild-homology}.
(We use a similar notation for $+$ versions.)
\begin{defn}\label{def1523}
We define a map\index[syindex]{fhatH@$\widehat{\frak f}^H $}
$$
\widehat{\frak f}^H 
: CH_*(\cL)
\to CH_*(\cL^{\rm form}_{\rm uni})
$$
by the following formula:
\begin{equation}
\aligned
\widehat{\frak f}^H (x_0\otimes x_1\otimes\dots \otimes x_k)
=
\sum (-1)^{\maltese} 
&\frak f^+(x_{k_0+1},\dots,x_{k_1})
\otimes 
\frak f^+(x_{k_1+1},\dots,x_{k_2})
\\
&\otimes 
\cdots
\otimes 
\frak f^+(x_{k_{\ell}+1},\dots,x_{k_0}),
\endaligned
\end{equation}
where 
$$
\maltese = (\deg' x_{k_0+1}+ \dots + \deg' x_{k})
(\deg' x_0 + \dots + \deg' x_{k_0})
$$
and the sum is taken over all $\ell$ and $k_i$ with 
$$
0\le k_1 \le k_2 \le \cdots \le k_{\ell} \le k_0 \le k.
$$
\end{defn}
(We put $x_{k+1} = x_0$ as convention.)
Here $\frak f^+$ is defined in Section \ref{canonical}.
(See also  Remark \ref{lem1115}.)
Namely it is a weighted sum of $\frak f^+_{\Gamma}$ in Subsection \ref{hunitcansubsec}.
\begin{lem}
$\widehat{\frak f}^H$ is a chain homotopy equivalence.
\end{lem}
\begin{proof}
We can prove that $\widehat{\frak f}^H$ is a chain map by an easy calculation 
using Lemma \ref{prop230}.
\par
The map $\widehat{\frak f}^H$ preserves the length filtration of the Hochschild complex and induces an isomorphism on the cohomology of associated graded groups.
So $\widehat{\frak f}^H$  is a chain homotopy equivalence as required.
\end{proof}
By composing $\widehat{\frak f}^H$ and $\widehat{\frak p}^{\rm f.u.\text{\bf b}}$ we obtain $\widehat{\frak p}^{{\bf b}}$.\index[syindex]{pcanbbf@$\widehat{\frak p}^{{\bf b}}$}
\begin{proof}[Wrap up of the proof of Theorem \ref{themp}]
The properties claimed in Theorem \ref{themp} is a consequence of corresponding properties of $\widehat{\frak p}^{\rm f.u.\text{\bf b}}$.
 %\marginpar{The rest of this subsection is added in 2025 March. KF}
The $\widehat{\frak p}^{\rm f.u.\text{\bf b}}$ version of Theorem \ref{themp} (2) follows from the compatibility 
of the forgetful maps with the Kuranishi structures (Proposition \ref{Kuraeistspoly2} (5)) and of CF-perturbations 
(Proposition \ref{existmultipolu2} (3)). The cyclic symmetry (Theorem \ref{themp} (3))
is a consequence of Proposition \ref{Kuraeistspoly2} (9) and Proposition \ref{existmultipolu2} (7).
Theorem \ref{themp} (1) is a consequence of the description of local chart on $\mathcal M_{1;\emptyset}(\kappa;0)$
given during the proof of Lemma \ref{lem1510}.
The proof of Theorem \ref{themp} is complete.
\end{proof}
\begin{rem}
It seems likely that we can prove that the chain map $
\widehat{\frak f}^H $ in Definition \ref{def1523} is the linear part of an $L_{\infty}$ module homomorpism 
over the $L_{\infty}$ homomorphism mentioned in Remark \ref{rem425}. It will imply the de Rham cohomology 
version of Proposition \ref{prop1522}. We do not try to prove it here since we do not use it in this paper.
See \cite[Theorem 7.4.186]{fooo092} for a related result.
\end{rem}
\begin{rem}
By the construction of this section, it becomes possible to work with the 
filtered $A_{\infty}$ category where morphism spaces are finite dimensional 
$\Lambda$ vector spaces in Part \ref{part4} and etc. So taking a completion of the tensor products among them 
becomes unnecessary.
\end{rem}

\begin{rem}
Note that we can and will decompose the $A_{\infty}$ category $\cL$ 
into  $\cL_{\lambda}$ corresponding to various values $\lambda$ of the potential function $\frak{PO}$.
We then re-define $\frak m_0$ to be $0$.  (It was  $\frak m_0(1) = \lambda {\bf e}$ before 
re-definition.) These two versions of the choices of $\frak m_0$ give different Hochschild differentials 
on the Hochschild chain complex.  However Proposition \ref{prop528} implies that the
resulting Hochschild homologies are canonically 
isomorphic. Moreover Theorem \ref{themp} (2) implies that the open closed map is 
compatible with this (canonical) isomorphism.
\end{rem}

\begin{rem}\label{rem1527}
In \cite[Theorem 3.8.9]{fooo09} another property of the open-closed map, that is, \cite[(3.8.10.3)]{fooo09} was given.
In our de Rham version it will be
\begin{equation}\label{eq1533}
\widehat{\frak p}^{\bf b}_1(\frak m_0(1)) \pm d\, \widehat{\frak p}^{\bf b}_0(1) \pm GW^{\rm form}_{0,1}(X)(i_!(L)) = 0. 
\end{equation}
Here $GW^{\rm form}_{0,1}(X) : \Omega(X) \to \Omega(X)$ is a certain\index[syindex]{GW01X@$GW^{\rm form}_{0,1}(X)$}
(de Rham) chain level version of Gromov-Witten invariant and 
$i_!(L)$ is a certain differential form on $X$ which represents the Poincar\'e dual to $[L]$.
$\widehat{\frak p}^{\bf b}_0(1)$ is the differential form on $X$ obtained by the pushforward of the differential $0$ form $1$ on the moduli spaces 
$\mathcal M_{1;\emptyset}(\kappa;B)$ for various $B$.
(In case the bulk class $\frak b$ is non-zero we need to include the forms obtained 
from $\frak b^{\ell}$ by the smooth correspondence $\mathcal M_{\ell+1;\emptyset}(\kappa;B)$.)
We can prove (\ref{eq1533}) by studying the boundary of the moduli space ${\mathcal M}_{\ell+1;\vec k}((\vec{\kappa},\vec p);B)^{\frak p}$
appearing in Theorem \ref{Kuraeistspoly2} (2) (boundary of type 5).
We do not discuss the detail of the proof of (\ref{eq1533}) here because we do not use it in this paper,
also because
the proof thereof is mostly the same as that of \cite[Theorem 3.8.9 (3.8.10.3)]{fooo09}.
\end{rem}

\section{Duality between   $\widehat{\frak q}$ and 
$\widehat{\frak p}$.}
\label{dualpq}

In this section we prove the following:

\begin{thm}\label{thm:duality}
Let $g \in H(X;\Lambda_0)$ be a 
de Rham cohomology class and 
$\text{\bf h} \in HH_*(\cL;\Lambda_0)$
 a Hochschild homology class. Then we have
\begin{equation}
\langle
\widehat{\frak q}^{\text{\bf b}}(g),\text{\bf h}
\rangle_{\rm HH}
=
(-1)^{\deg g}\langle
g,\widehat{\frak p}^{\text{\bf b}}(\text{\bf h})
\rangle_{\text{\rm PD}_X}.
\end{equation}
Here the pairing in the left hand side is the one used in  %\marginpar{Sign to be checked.  KF 2025 Aug 30}
Lemma $\ref{lem:pairing_HH}$\index[syindex]{<*,*>PDX@$\langle*,*\rangle_{\text{\rm PD}_X}$}
%{\rm \ref{lem:dual_bi-modules_HH}} 
and the pairing in the right hand side is the
Poincar\'e duality pairing on $X$.
\end{thm}

In this section we prove Theorem \ref{thm:duality} except the sign.
The sign will be discussed in Subsection \ref{sec:signduality}.
To prove Theorem \ref{thm:duality}, we first prove a similar duality between the form versions.

\begin{prop}\label{dualityinformlevel}
{\rm (Compare \cite[Remark 7.4.197]{fooo09},
\cite[Theorem 3.3.8]{toric3})}
Let $g \in \Omega(X)$ with $dg =0$.
Then there exist\index[syindex]{Hxgfrak@$\frak H_g$} 
$$
\frak H_g \in 
\Hom(CH_*(\cL_{\rm c.u.}^{\rm form}), 
\Lambda_0)
$$
such that the following equality holds:
\begin{equation}\label{dualityform}
{}^+\widehat{\frak q}^{\rm f.u.\text{\bf b}}(g;\text{\bf h}\otimes h_0)
=
(-1)^{\deg g}\langle g,
\widehat{\frak p}^{\rm f.u.\text{\bf b}}(\text{\bf h}\otimes h_0)
\rangle_{\text{\rm PD}_X}
+ 
\frak H_g(\delta^H(\text{\bf h}\otimes h_0)).
\end{equation}
\end{prop}
Here ${}^+\widehat{\frak q}^{\rm f.u.\text{\bf b}}$
\index[syindex]{qfxormelqlvechat+@${}^+\widehat{\frak q}^{\rm f.u.\text{\bf b}}$}
is an operator which is related to $\widehat{\frak q}^{\rm f.u.\text{\bf b}}$
by the formula 
$$
{}^+\widehat{\frak q}^{\rm f.u.\text{\bf b}}(g;\text{\bf h}\otimes h_0)
= \langle \widehat{\frak q}^{\rm f.u.\text{\bf b}}(g;\text{\bf h}),h_0\rangle_{\rm cyc}
$$
when $h_0$ is not ${\bf f}$ or ${\bf e}^+$.
\begin{rem}
Since we constructed operators by induction, actually 
we prove the version where we replace 
$\widehat{\frak p}^{\rm f.u.\text{\bf b}}$ 
by $\widehat{\frak p}_{n,k}^{\rm f.u.\text{\bf b}}$ etc.
On the other hand, since  Theorem \ref{thm:duality} we are proving is an 
equality between two well defined numbers. So it suffices to prove the equality 
modulo $T^E$ for an arbitrary but fixed $E$.  Therefore it suffices to prove 
its $A_{n,k}$ version. So for the simplicity of exposition we state Proposition \ref{dualityinformlevel} 
in the above form without specifying these details.
\par
The same remark applies to various other places 
including Sections \ref{sec:ring}, \ref{sec:annuli}.
\end{rem}
\begin{proof}
Note that if we {\it could} use the same 
perturbation of
${\mathcal M}_{\ell;\vec k}((\vec{\kappa},\vec p);B)$ to define 
$\widehat{\frak q}^{\rm f.u.\text{\bf b}}$
and 
$\widehat{\frak p}^{\rm f.u.\text{\bf b}}$
then (\ref{dualityform}) would be rather immediate.
Namely, changing the notation so that $ \text{\rm ev}^{\partial} $ corresponds to evaluation at all boundary marked points, the left hand side would become
\begin{equation}\label{integralexpreqp}
\sum_{B,\ell',\vec k}
\pm\frac{\rho_{\frak b,\theta}(B)}{\ell!}%\rho_{b}(B)
\int_{{\mathcal M}_{\ell+1;\vec k}((\vec{\kappa},\vec p);B)}
(\text{\rm ev}^{\partial})^{*}(b_{\kappa_0}^{k_0} \otimes h_1 \otimes \cdots \otimes h_K \otimes b_{\kappa_K}^{k_K} \otimes h_0)
%(\text{\rm ev}^{\partial,2})^{*}(b_+^{\vec k})
(\text{\rm ev}^+)^*(g \times\frak b_+^{\otimes\ell}).
\nonumber
\end{equation}
This would be equal to 
$(-1)^{\deg g}  
\langle g,
\widehat{\frak p}^{\rm f.u.\text{\bf b}}(\text{\bf h}\otimes h_0)
\rangle_{\text{\rm PD}_X}.
$
Namely (\ref{dualityform}) would hold without the term 
$
\frak H_g(\delta^H(\text{\bf h}\otimes h_0))
$.\footnote{In the above formula we consider the case when 
${\bf h}$ does not contain ${\bf e}^+$ or ${\bf f}$.}
\par
However we remark that the above integration 
does depend on the choice of 
the perturbation
of ${\mathcal M}_{\ell+1;\vec k}((\vec{\kappa},\vec p);B)$.
For the left hand side of  (\ref{dualityform}) we use the 
\newred{CF-perturbation} of Proposition \ref{existmultipolu2} 
and for the right hand side we use one of
Proposition \ref{existmultipolu}. 
\par
So for the proof we need to interpolate the two CF-perturbations defined on the same moduli space.
The detail follows. 
\begin{lem}\label{Kuraeistspoly3}
There exists a system of Kuranishi structures on 
$[0,1]_{\tau} \times [0,1]^{\vert \frak f_{\bf f}\vert} \times{\mathcal M}_{\ell+1;\vec k}((\vec{\kappa},\vec p);B;\vec{\frak f})^{\boxplus 1}$
for various $\ell, \vec k$, $\vec{\kappa}$, $\vec p$, $B$, $\vec{\frak f}$ with the 
following properties.
\begin{enumerate}
\item
The restriction of this Kuranishi structure to 
$(\{0\} \times[0,1]^{\vert \frak f_{\bf f}\vert} \times {\mathcal M}_{\ell+1;\vec k}((\vec{\kappa},\vec p);B;\vec{\frak f})^{\boxplus 1})^{\frak p}$
coincides with the one in Propositions $\ref{Kuraeistspoly}$ and $\ref{prop94}$.
The restriction of this Kuranishi structure to 
$(\{1\} \times [0,1]^{\vert \frak f_{\bf f}\vert} \times {\mathcal M}_{\ell+1;\vec k}((\vec{\kappa},\vec p);B;\vec{\frak f})^{\boxplus 1})^{\frak p}$
coincides with the one in Proposition \ref{Kuraeistspoly2}.
\item 
The boundary of $([0,1]_{\tau} \times [0,1]^{\vert \frak f_{\bf f}\vert}\times {\mathcal M}_{\ell+1;\vec k}((\vec{\kappa},\vec p);B;\vec{\frak f})^{\boxplus 1})^{\frak p}$
is a union of the following five types of fiber or direct products:
\par
{\rm(\bf Boundary of type 0)}
$(\{0,1\}\times [0,1]^{\vert \frak f_{\bf f}\vert} \times {\mathcal M}_{\ell+1;\vec k}((\vec{\kappa},\vec p);B;\vec{\frak f})^{\boxplus 1})^{\frak p}$.
\par
{\rm (\bf Boundary of type 1)}
We use the notation of $(\ref{218})$ and Proposition $\ref{prop94}$.
\begin{equation}\label{21832}
\aligned
&([0,1]^{\vert \frak f^2_{\rm{\bf f}} \vert}  \times{\mathcal M}_{\ell'';k''_i+1};(L;\beta;\vec{\frak f}^2)^{\boxplus 1})^{\frak q}\\
&\,{}_{\text{\rm ev}_0} \times_{\text{\rm ev}_{i,j}}
([0,1]_{\tau}\times [0,1]^{\vert \frak f^1_{\rm{\bf f}} \vert}  \times{\mathcal M}_{\ell'+1;\vec k'}((\vec{\kappa},\vec p);B');\vec{\frak f}^1)^{\boxplus 1})^{\frak p}.
\endaligned
\end{equation}
In $(\ref{21832})$, the union is taken over the data described in $(\ref{218})$ and Proposition $\ref{prop94}$
together with the shuffles $(\mathbb L'',\mathbb L')$ of $\underline\ell$
such that $\# \mathbb L'' = \ell''$, $\# \mathbb L' = \ell'$.
\par
In this lemma, we put the Kuranishi structure in this lemma
%Lemma $\ref{Kuraeistspoly3}$
on the second factor of $(\ref{21832})$.
We put the Kuranishi structure of Proposition $\ref{Kuraeistspoly}$  and Proposition $\ref{prop94}$
on the first factor of $(\ref{21832})$.
\par
{\bf (Boundary of type 2)} 
We use the notation of $(\ref{219})$ and Proposition $\ref{prop94}$.
\begin{equation}\label{219220}
\aligned
&([0,1]_{\tau}\times [0,1]^{\vert \frak f^2_{\rm{\bf f}} \vert}  \times{\mathcal M}_{\ell''+1;\vec k''}((\vec{\kappa}'',\vec p'');B'';\vec{\frak f}^2)^{\boxplus 1})^{\frak p}\\
&\times
([0,1]^{\vert \frak f^1_{\rm{\bf f}} \vert} \times {\mathcal M}_{\ell';\vec k'}((\vec{\kappa}',\vec p');B';\vec{\frak f}^1)^{\boxplus 1})^{\frak q}.
\endaligned
\end{equation}
\begin{equation}\label{21922}
\aligned
& ([0,1]^{\vert \frak f^2_{\rm{\bf f}} \vert}  \times{\mathcal M}_{\ell'';\vec k''}((\vec{\kappa}'',\vec p'');B'';\vec{\frak f}^2)^{\boxplus 1})^{\frak q}\\
&\times
([0,1]_{\tau} \times[0,1]^{\vert \frak f^1_{\rm{\bf f}} \vert} \times {\mathcal M}_{\ell'+1;\vec k'}((\vec{\kappa}',\vec p');B';\vec{\frak f}^1)^{\boxplus 1})^{\frak p}.
\endaligned
\end{equation}
In $(\ref{219220})$, $(\ref{21922})$, the union is taken over the data described in $(\ref{219})$ 
and Proposition $\ref{prop94}$
together with the shuffles $(\mathbb L'',\mathbb L')$ of $\underline\ell$
such that $\# \mathbb L'' = \ell''$, $\# \mathbb L' = \ell'$.
\par
We put the Kuranishi structure \blue{in this lemma} and the first factor of $(\ref{219220})$ 
%Lemma $\ref{Kuraeistspoly3}$  
for the  second factor of $(\ref{21922})$.
We put the Kuranishi structure of Propositions $\ref{Kuraeistspoly}$ and $\ref{prop94}$
for  the   second factor of $(\ref{219220})$ and the first factor of $(\ref{21922})$.
%$$
%\text{\bf Figure 3.2}
%$$
\par
{\bf (Boundary of type 3)} 
We use the notation of $(\ref{220})$ and Proposition $\ref{prop94}$.
\begin{equation}\label{22012}
\aligned
&([0,1]_{\tau}\times [0,1]^{\vert \frak f^2_{\rm{\bf f}} \vert}  \times{\mathcal M}_{\ell''+1;\vec k''}((\vec{\kappa}'',\vec p'');B'';\vec{\frak f}^2)^{\boxplus 1})^{\frak p}\\
&\,{}_{\text{\rm ev}_{1,m''}}\times_{\text{\rm ev}_{i,m'}}
([0,1]^{\vert \frak f^1_{\rm{\bf f}} \vert} \times{\mathcal M}_{\ell';\vec k'}((\vec{\kappa}',\vec p');B';\vec{\frak f}^1)^{\boxplus 1})^{\frak q},
\endaligned
\end{equation}
\begin{equation}\label{22022a}
\aligned
& ([0,1]^{\vert \frak f^2_{\rm{\bf f}} \vert}  \times{\mathcal M}_{\ell'';\vec k''}((\vec{\kappa}'',\vec p'');B'';\vec{\frak f}^2)^{\boxplus 1})^{\frak q}\\
&\,{}_{\text{\rm ev}_{1,m''}}\times_{\text{\rm ev}_{i,m'}}
([0,1]_{\tau}\times [0,1]^{\vert \frak f^1_{\rm{\bf f}} \vert} \times{\mathcal M}_{\ell'+1;\vec k'}((\vec{\kappa}',\vec p');B';\vec{\frak f}^1)^{\boxplus 1})^{\frak p},
\endaligned
\end{equation}
In $(\ref{22012})$, $(\ref{22022a})$, the union is taken over the data described in $(\ref{220})$ 
together with the shuffles $(\mathbb L'',\mathbb L')$ of $\underline\ell$
such that $\# \mathbb L'' = \ell''$, $\# \mathbb L' = \ell'$.
\par
We put the Kuranishi structure in this lemma
for the first factor of $(\ref{22012})$ and the second factor of $(\ref{22022a})$.
We put the Kuranishi structure of Propositions $\ref{Kuraeistspoly}$  and $\ref{prop94}$
for the second factor of $(\ref{22012})$ and the first factor of $(\ref{22022a})$.
\par
{\bf (Boundary of type 4)} 
$([0,1]_{\tau} \times \partial([0,1]^{\vert \frak f_{\bf f}\vert}) \times{\mathcal M}_{\ell+1;\vec k}((\vec{\kappa},\vec p);B;\vec{\frak f})^{\boxplus 1})^{\frak p}$
%$$
%\text{\bf Figure 3.3}
%$$
\item
The Kuranishi structure can be oriented so that it is 
compatible with the identification $(2)$ of the boundaries.
\item
The evaluation maps $(\ref{evevev})$ extend to the 
compactification and are compatible with 
the description of the boundary in $(3)$.
\item 
A similar compatibility with the forgetful map as Proposition $\ref{prop94}$ $(8)$ holds.
\item When the coordinate $s_{i,j}$ of $ [0,1]^{\vert \frak f_{\bf f}\vert}$ factor 
becomes $0$ or $1$ a similar compatibility condition as Proposition $\ref{prop94}$ $(9)$, $(10)$ hold.
\item
The Kuranishi structure is invariant under the permutation of $\ell$
interior marked points $z^+_1,\dots,z^+_{\ell}$.
\item
The Kuranishi structure is invariant under the cyclic permutation 
of the data in the sense we described in Proposition $\ref{Kuraeistspoly}$. 
\end{enumerate}
\end{lem}
The proof is the same as the proof of Propositions \ref{Kuraeistspoly}, 
\ref{Kuraeistspoly2} and so is omitted.
\begin{lem}\label{existmultipolu3}
There exists a system of CF-perturbations on 
$[0,1]_{\tau}\times [0,1]^{\vert \frak f_{\bf f}\vert}\times {\mathcal M}_{\ell+1;\vec k}((\vec{\kappa},\vec p);B;\vec{\frak f})^{\boxplus 1}$ with 
the following properties.
\begin{enumerate}
\item
They are transversal to $0$.
\item 
They are compatible with the description of the 
boundary in Proposition $\ref{Kuraeistspoly3}$ $(1), (2)$ and Proposition $\ref{prop95}$, where we consider the \blue{CF-perturbation} in 
Proposition $\ref{existmultipolu}$, $\ref{prop95}$, $\ref{existmultipolu2}$,
$\ref{existmultipolu3}$ on the factors equipped with the Kuranishi structures 
in Propositions $\ref{Kuraeistspoly}$, $\ref{prop94}$, $\ref{Kuraeistspoly2}$,
$\ref{Kuraeistspoly3}$.
\item 
The compatibility with forgetful map in Lemma $\ref{Kuraeistspoly3}$ $(5)$ is  enhanced  to the compatibility of CF-perturbations.
\item
The compatibility in Lemma $\ref{Kuraeistspoly3}$ $(6)$ is enhanced to the compatibility of CF-perturbations.
\item
The CF-perturbation is invariant under the permutation of 
$\ell$ interior marked points $z^+_1,\dots,z_{\ell}^+$.
\item
The CF-perturbation is invariant under the cyclic permutation 
of the data in the same sense as Proposition $\ref{Kuraeistspoly}$ $(8)$. 
\end{enumerate}
\end{lem}
The proof is the same as the proof of Proposition  \ref{existmultipolu}, and is also omitted.
\par
Continuing with our proof of Proposition \ref{dualityinformlevel}, we now put
\begin{equation}\label{integralexpreH}
\aligned
\frak H_g(\text{\bf h})
= \sum_{B,\ell,\vec k}(-1)^{\maltese}
&\frac{\rho_{\frak b,\theta}(B)}{\ell!} 
\int_{\blue{[0,1] \times}  {\mathcal M}_{\ell+1}(\vec{\kappa},\vec p;B)}
\\
&
(\text{\rm ev}^+)^*(g \times\frak b_+^{\otimes\ell})
(\text{\rm ev}^{\partial})^{*}(b_{\kappa_0}^{k_0} \otimes h_1 \otimes \cdots \otimes b_{\kappa_K}^{k_K} \otimes h_0).
\endaligned
\end{equation}
Here  
$$
{\mathcal M}_{\ell+1}(\vec{\kappa},\vec p;B)
= 
{\mathcal M}_{\ell+1,\vec k'}((\vec{\kappa}',\vec p,m');B)
=
{\mathcal M}_{\ell+1,\vec k'}((\vec{\kappa}',\vec p);B),
$$
with $\text{\rm Red}(\vec{\kappa},\vec p) = (\vec{\kappa}', \vec{p}, \vec k', m')$ as in (\ref{eq:change_notation_moduli_space}).
The sign $\maltese$ is the same as (\ref{form21ten2}).
\par
The formula (\ref{integralexpreH}) is the case when $\bf h$ does not contain ${\bf e}^+$ or ${\bf f}$.
When it contains ${\bf e}^+$ and/or ${\bf f}$ we use the evaluation at the marked points 
labelled by $\frak f_{\bf e}$, $\frak f_{\bf f}$ in the same way as Subsection \ref{subsec:opeartor}.
\par
The equality (\ref{dualityform}) follows from Stokes' theorem 
and Lemmas \ref{Kuraeistspoly3}, \ref{existmultipolu3}.  
\end{proof}

We have thus proved the duality between $\widehat{\frak q}^{\rm f.u.\text{\bf b}}$
and  $\widehat{\frak p}^{\rm f.u.\text{\bf b}}$.
(The sign will be studied in Subsection \ref{sec:signduality}.) 
We next prove the same conclusion 
for the canonical model.
\par
Let $(\Gamma,v)$ be the decorated ribbon tree in the sense of 
Definition \ref{decribbon}.
Removing $v$ from $\Gamma$, the closure of the 
connected components of the complement are 
ribbon trees, which we write $\Gamma_0, \Gamma_1, \dots, \Gamma_m$.
Here $\Gamma_0$ contains the root of $\Gamma$ and enumeration 
respects the cyclic order of the edges containing $v$, that is 
induced by the ribbon structure.
We regard $v$ as the root of each of $\Gamma_i$ (including the case of $\Gamma_0$.).
Let $k+1$ be the number of exterior edges of $\Gamma$ we take elements of Floer complex   
$h_1,\dots,h_k$ such that
$$
h_1 \otimes \dots \otimes h_k \in BCF(\cL^{\rm form};\vec{\kappa};\text{\rm source})^+
$$
where $\vec{\kappa} = \vec\kappa(\Gamma)$ is as in (\ref{veckapadef}).
We also take 
$h_0 \in CF(L_{\kappa_0},L_{\kappa_K})$.
By construction each of the exterior vertices of $\Gamma_i$ other than $i$-th one 
is assigned one of the morphisms $h_j$'s.  We write $j_1,\dots, j_m,  j_0$ for the natural numbers such that $
h_{j_i+1},\dots,h_{j_{i+1}}
$
is assigned to the vertices of $\Gamma_i$ if $1 \leq i \leq m$, and differential forms

$
h_{j_{0}+1},\dots,h_{j_{k}}, h_{0}, h_1,\dots,h_{m_0}
$
are assigned to the vertices of $\Gamma_0$.
Here
$0\le j_1\le j_2 \le \dots \le j_m \le j_0 \le k$.
\par
We put
$$
\aligned
\text{\bf h}_1 &= h_{j_1+1} \otimes \dots \otimes h_{j_2}, \\
&\dots\\
\text{\bf h}_m &= h_{j_m+1} \otimes \dots \otimes h_{j_0},\\
\text{\bf h}_0 &= h_{j_0+1} \otimes \dots \otimes h_k \otimes h_0\otimes \dots \otimes h_{j_1}, \\
\text{\bf h} &= h_1 \otimes \dots \otimes h_k.
\endaligned
$$
\begin{lem}\label{lemmovevertex}
\begin{equation}\label{eqlem314}
{}^+\frak q^{\rm f.c.u.\frak b}_{(\Gamma,v)}(g,\text{\bf h},h_0)
=
(-1)^{\maltese}{}^+\frak q^{\rm f.c.u.\frak b}_{m,B(v)}(g,
\frak f_{\Gamma_1}(\text{\bf h}_1),\dots,\frak f_{\Gamma_m}(\text{\bf h}_m),
\frak f_{\Gamma_0}(\text{\bf h}_0))
\end{equation}
where $\maltese=
(\deg' h_1+\dots+\deg'h_{j_1})
(\deg' h_{j_1+1}+\dots+\deg' h_k + \deg' h_0 + \deg g)$.
\end{lem}
The proof is the same as that of 
\cite[Proof of Proposition 10.1]{fukaya:cyc}
and so is omitted.
\par
\begin{proof}[Proof of Theorem \ref{thm:duality}]
Note the sum  over $(\Gamma,v)$ of left hand side of  (\ref{eqlem314}) 
with weight 
$T^{\omega(\Gamma)}\rho_{b}(\Gamma) \rho_{\frak b,\theta}(\Gamma)$
is the left hand side of
Theorem \ref{thm:duality}.
The sum over $(\Gamma,v)$ of
$$
\langle
g,\frak p^{\rm f.c.u.\frak b}_{m+1,B(v)}(
\frak f_{\Gamma_1}(\text{\bf h}_1),\dots,\frak f_{\Gamma_m}(\text{\bf h}_m),
\frak f_{\Gamma_0}(\text{\bf h}_0))
\rangle_{\text{\rm PD}_X}
$$
with the same weight is 
the right hand side of 
Theorem  \ref{thm:duality}.
 
Therefore Theorem \ref{thm:duality} follows from 
Proposition \ref{dualityinformlevel} and 
Lemma \ref{lemmovevertex}.
The proof of Theorem  \ref{thm:duality} is complete.
\end{proof}

\section{Nontriviality of   $\widehat{\frak q}$ and 
$\widehat{\frak p}$.}
\label{nontrivialsec}

In this subsection we prove the following:

\begin{prop}\label{nontriialp}
If the Floer cohomology $HF((U,b_U),(U,b_U);\Lambda)$ 
with Novikov {\it field} $\Lambda$ coefficient is 
nonzero, then 
$$
\widehat{\frak p}^{(\frak b,b_U)}  : HH_*(HF((U,b_U),(U,b_U);\Lambda))
\to H(X;\Lambda)
$$
is nonzero.
\end{prop}
Here $U$ is a relatively spin Lagrangian submanifold in $X$ and 
$b_U$ is its  bounding cochain.
We consider the category with only one object $(U,b_U)$
and the Hochschild homology in the left hand side is that 
of this category.
\par
We remark that we use the Novikov field $\Lambda$ rather than $\Lambda_0$
as a coefficient in Proposition \ref{nontrivialq}. 
In fact, the non-triviality over  $\Lambda$ is stronger that 
the non-triviality over $\Lambda_0$.
\par
Proposition \ref{nontriialp} is a consequence of 
Theorem \ref{thm:duality} and the next result on $\widehat{\frak q}$.
We remark that the complex defining  Hochschild cohomology
$HH^*(HF((U,b_U),(U,b_U);\Lambda))$ is given by 
\begin{equation}\label{HHcomplex}
\bigoplus_{k=0}^{\infty}
\Hom(B_kH(U;\Lambda)[1] ,H(U;\Lambda)[1] ).
\end{equation}
Note that $k=0$ is included. So 
$\Hom(B_0H(U;\Lambda)[1] ,H(U;\Lambda)[1] )= H(U;\Lambda)[1]$
is a direct summand of (\ref{HHcomplex}).  The restriction of Hochschild coboundary to 
this summand is the usual boundary operator 
$\frak m_1^{(\frak b,b_U)}$.
So we have a natural homomorphism
\begin{equation}\label{HFtoHH}
\frak I_{(U,b_U)}:
HF((U,b_U),(U,b_U);\Lambda) \to HH^*(HF((U,b_U),(U,b_U);\Lambda)).
\end{equation}
This map is injective since 
the $k>0$ part of (\ref{HHcomplex}) 
does not go to $k=0$ part by the Hochschild coboundary map.
Let us take an element $\text{\bf e}_U$  of $ H^0(U;\Lambda)$ 
that gives a nonzero element of $HF((U,b_U),(U,b_U);\Lambda)$ 
as long as $HF((U,b_U),(U,b_U);\Lambda) \ne 0$.

The element $\text{\bf e}_U \in H^0(U;\Lambda) $ is represented by the constant function on $U$ with value $1$. Because the $A_\infty$ structure on $HF((U,b_U),(U,b_U);\Lambda)$ is strictly unital, this represents a closed element of $HH^*(HF((U,b_U),(U,b_U);\Lambda)) $.  We write it $\text{\bf e}_{(U,b_U)}$.

\begin{prop}\label{nontrivialq}
We have
\begin{equation}\label{qunital}
\widehat{\frak q}^{(\frak b,b_U)}(\text{\bf e}_X) 
= \text{\bf e}_{(U,b_U)},
\end{equation}
where $\text{\bf e}_X \in H^0(X;\Lambda)$ is the unit.
\end{prop}

Proposition \ref{nontriialp} follows from Proposition \ref{nontrivialq}, the 
injectivity of (\ref{HFtoHH}) and Theorem \ref {thm:duality}.
The rest of this subsection is devoted to the 
proof of Proposition \ref{nontrivialq}.
\begin{proof}
We first explain the idea of the proof.
The map 
$\widehat{\frak q}^{(\frak b,b_U)}$ 
is induced by the correspondence 
via the moduli spaces 
\begin{equation}\label{correspondencefor1}
\text{\rm Corr}(\mathcal M_{\ell+1;k+1}(U;\beta);
{\rm ev}^+\times ({\rm ev}_1\times \dots \times {\rm ev}_k),
{\rm ev}_0)
(\text{\bf  e}_X\otimes \frak b_+^{\otimes \ell}, \text{\bf h})
\end{equation}
with appropriate weight.
We consider the forgetful map
\begin{equation}\label{forget}
\frak{forget}_+^1 : \mathcal M_{\ell+1;k+1}(U;\beta)
\to \mathcal M_{\ell;k+1}(U;\beta).
\end{equation}
Then 
$$
({\rm ev}^+_1)^*(\text{\bf e}_X) = (\frak{forget}_+^1)^* 1,
$$
where $1$ is the constant function $1$ regarded as a zero form 
on $U$.
Moreover all the other evaluation maps factor 
through $\frak{forget}_+^1$.
Therefore (\ref{correspondencefor1}) is 
zero unless the dimension of the fiber of 
(\ref{forget}) is zero.
The dimension of the fiber of 
(\ref{forget}) is zero if and only if $k=\ell=0$ and 
$\beta =0$. Thus 
(\ref{correspondencefor1}) is nonzero only 
in case $k=\ell=0$ and 
$\beta =0$. It is easy to see that
$$
\text{\rm Corr}(\mathcal M_{1;k}(U;0);
{\rm ev}^+,
{\rm ev}^{\partial}_0))
(\text{\bf e}_X)
= \text{\bf e}_{(U,b_U)}.
$$
Thus (\ref{qunital}) holds in the chain level.
\par
Actually there is one gap in the above argument.
Namely the \blue{CF-perturbation} (and the Kuranishi structure) 
is not \blue{compatible with} the forgetful map (\ref{forget}), 
which is a map forgetting interior marked points.
(See Remark \ref{rem819} for the reason why it is difficult to do so.)
Therefore we need to take a homotopy between the two \blue{CF-perturbations}, one which is 
compatible with the forgetful map, the other which is used to define $\frak q$, and 
use the homotopy between them to prove  (\ref{qunital}) in the homology level.
The detail follows.
We remark that in this subsection we study a single Lagrangian submanifold $U$ so we do not need to involve 
${\bf e}^+$ or ${\bf f}$.

\begin{lem}\label{diskkurahomot}
$  $ \par
\begin{enumerate}
\item
The space $[0,1]_{\tau}\times {\mathcal M}_{\ell+1;k+1}(U;\beta)^{\boxplus 1}$ has an orientable Kuranishi structure with corners. 
\item
Its restriction to $\{0\} \times {\mathcal M}_{\ell+1;k+1}(U;\beta)^{\boxplus 1}$
coincides with the Kuranishi structure in Proposition $\ref{diskkura}$.
\item
Its restriction to $\{1\} \times {\mathcal M}_{\ell+1;k+1}(U;\beta)^{\boxplus 1}$
coincides with the Kuranishi structure induced by the forgetful map
$$
\frak{forget}_+^1 : \{1\} \times \mathcal M_{\ell+1;k+1}(U;\beta)^{\boxplus 1}
\to \mathcal M_{\ell;k+1}(U;\beta)^{\boxplus 1}
$$
from the Kuranishi structure on 
$\mathcal M_{\ell;k+1}(U;\beta)^{\boxplus 1}$ in Proposition $\ref{diskkura}$.
\item
The normalized boundary of $[0,1]_{\tau}\times {\mathcal M}_{\ell+1;k+1}(U;\beta)^{\boxplus 1}$
\index[syindex]{01MtauUbeta@$[0,1]_{\tau}\times {\mathcal M}_{\ell+1;k+1}(U;\beta)^{\boxplus 1}$} in the 
sense of Kuranishi structure is described by the following fiber product 
over $U$:
\begin{equation}\label{eq1772}
\aligned
&\partial([0,1]_{\tau}\times {\mathcal M}_{\ell+1;k+1}(U;\beta)^{\boxplus 1})\\
&= 
\partial[0,1]_{\tau}\times\mathcal M_{\ell+1;k+1}(U;\beta)^{\boxplus 1}\\
&\quad \cup \bigcup {\mathcal M}_{\# {\mathbb L}_1;k_1+1}(U;\beta_1)^{\boxplus 1}
{}_{\text{\rm ev}^{\partial}_{0}}\times_{\text{\rm ev}^{\partial}_i} 
([0,1]_{\tau} \times {\mathcal M}_{\# {\mathbb L}_2;k_2+1}(U;\beta_2)^{\boxplus 1}), 
\endaligned
\end{equation}
where the union is taken over all $({\mathbb L}_1,{\mathbb L}_2) \in \text{\rm Shuff}(\ell)$, $i = 1, \dots, k_2$,
$k_1,k_2\in \Z_{\ge 0}$ with $k_1 + k_2 = k$ and $\beta_1,\beta_2
\in H_2(X,U;\Z)$ with $\beta_1 +\beta_2 = \beta$.
\par
Here the Kuranishi structure of the first factor of the second term 
is the one in Proposition $\ref{diskkura}$ and the one of
the second factor of the second term 
is the Kuranishi structure provided by this lemma.
\item
We can define orientations of $[0,1]_{\tau}\times {\mathcal M}_{\ell+1;k+1}(U;\beta)^{\boxplus 1}$ so that 
$(4)$ above is compatible with this orientation.  
\item The evaluation maps are extended to the 
compactification so that it is compatible with $(\ref{eq177})$.
\item $ev_0^{\partial}$ is weakly submersive.
\item
The Kuranishi structure is compatible with the forgetful map of 
boundary marked points.
\item
The Kuranishi structure is invariant under the permutation of 
$1$st \dots $\ell$-th interior marked 
points. (It is not necessarilyy invariant under the permutation including the
$0$-th marked point.)
\item 
The Kuranishi structure is invariant under the cyclic permutation of the 
boundary  marked 
points. 
\end{enumerate} 
\end{lem}
The proof is the same as that of Proposition \ref{diskkura}.
\begin{lem}\label{existmkulti1homot}
There exists
a system of \blue{CF-perturbations} 
on the moduli spaces $[0,1]_{\tau}\times {\mathcal M}_{\ell+1;k+1}(U;\beta)^{\boxplus 1}$
such that the following holds:
\begin{enumerate}
\item
It is transversal to zero.
\item
It is compatible with the
description of the boundary in Proposition $\ref{diskkurahomot}$ $(2),(3),(4)$ above.
(Especially it is compatible with the forgetful map 
$\frak{forget}_+^1$ at $\{1\}\times {\mathcal M}_{\ell+1;k+1}(U;\beta)^{\boxplus 1}$.)
\item
It is compatible
with the forgetful map of the boundary marked points.
\item
It is compatible
with the permutation of 1st, \dots, $\ell$-th interior marked points.
\item
It is compatible with cyclic permutation of the boundary marked points.
\item
$\text{\rm ev}_0^{\partial}%\times \text{\rm ev}_{\text{\rm time}}
$ restricted to the
zero set of this \blue{CF-perturbation} is a submersion.
%Here $\text{\rm ev}_{\text{\rm time}}$ is the projection to the 
%$[0,1]$ direction.
\end{enumerate}
\end{lem}
The proof is the same as that of
Lemma \ref{existmkulti1}.
We now define
$$
\frak G_{k,\beta}(\text{\bf h})
=
\text{\rm Corr}([0,1]_{\tau}\times \mathcal M_{\ell+1;k+1}(U;\beta)^{\boxplus 1};
{\rm ev}^+\times ({\rm ev}^{\partial}_1\times  \cdots \times {\rm ev}^{\partial}_k),
{\rm ev}^{\partial}_0))
(\text{\rm id}_X\otimes \frak b_+^{\otimes \ell}, \text{\bf h}),
$$
where we use the Kuranishi structure and \blue{CF-perturbations} %multisection 
in Lemmas \ref{diskkurahomot} and \ref{existmkulti1homot}.
We then put
$$
\frak G_{k}(\text{\bf h})
=
\sum_{\beta} 
T^{\omega(\beta)}%\rho_{b}(\partial\beta)
\rho_{\frak b,\theta}(\partial\beta)
\frak G_{k,\beta}(\text{\bf h}).
$$
Now we have
\begin{equation}\label{unitalityinchain}
\widehat{\frak q}^{\rm f.u. {\bf b}}(\text{\bf e}_X) 
- \text{\bf e}_{(U,b_U)}
= \delta^H(\frak G).
\end{equation}
Here $\frak G$ is a map from the Hochschild complex to $\Lambda$ which restricts to 
$\frak G_{k}$ on $B_k(\Omega(X)[1])$. The identity
(\ref{unitalityinchain}) is a consequence of 
Lemmas \ref{diskkurahomot} and \ref{existmkulti1homot}, the 
argument we explained at the begining of the proof and 
the Stokes' theorem.
Using the unitality of $\frak f$ it is easy to see that (\ref{unitalityinchain}) 
implies Proposition \ref{nontrivialq}.
\end{proof}

\section{$\hat{\frak q}$ is a ring homomorphism.}
\label{sec:ring}
 
In this section we prove:

\begin{thm}\label{theoremD1}
$\hat{\frak q}^{\bf b}  : QH^*_{\frak q}(X;\Lambda_0) \to HH^*(\cL,\cL)$ 
is a unital ring homomorphism.
\end{thm}
The degree of $\hat{\frak q}^{{\bf b}}$ is given by:  %\marginpar{A paragraph added.  KF 2025 Aug}
$
{\frak{deg}} (\hat{\frak q}^{{\bf b}}(g)) = \deg g -1.
$\footnote{In \cite[Theorem 3.8.32]{fooo09} the degree is claimed to be $+1$. Note in 
\cite{fooo09} the degree is shifted by $2$ for $H^*(X;\Lambda_0)$ 
and $\frak{deg}$ is used for Hochschild cohomology.  Therefore it is consistent with the degree 
given here.}
Therefore
$$
\deg (\hat{\frak q}^{{\bf b}}(g)) = \deg g.
$$
In other words, $\hat{\frak q}^{{\bf b}}$ is degree preserving if we use $\deg$ for 
Hochschild cohomology.  (See (\ref{form460new}) 
and footnotes 20, 21 for $\frak{deg}$ and $\deg$.)
Since the cup product $\cup$ on Hochschild cohomology is associative and degree preserving using  $\deg$, this 
fact is consistent with Theorem \ref{theoremD1}.
\begin{rem}
We refer to \cite[Section 4.7]{toric3} for some of the history of this and related results.
 %\marginpar{Maybe we need to 
%quote more reference which appears after \cite{toric3}  KF 2025 July}
\end{rem}
\begin{rem}
Theorem \ref{theoremD1} and Corollary \ref{Cor1526} should be combined and 
be a part of the statement that closed-open map is a morphism over 
a little-disc-operad.  We do not to try to prove it in this paper.  Note that the existence of module structure over little-disc-operad on Hochschild cohomology
is called Deligne conjecture and is  proved by several people (\cite{KS2,MS,Ta,Vo}). 
 %\marginpar{Remak added.  
%Whether the name littile disc operad is OK or not should be checked. Reference to be
%checked.  KF 2025 Aug.}
\end{rem}
\begin{rem}
In order to handle infinitely many moduli spaces at the same time 
there occurs a  problem similar to the one mentioned in Remark \ref{runningoutrem}.
In the situation of this and the next sections however, we will need to study only finitely many of them.
In fact the statement we prove here is an equality between 
two operators that are independent of various choices involved. 
So if we prove that they coincide
modulo $T^E$ for each given $E$ then the equality holds. So 
we have  to study only the moduli spaces of pseudo-holomorphic maps 
of symplectic area $\le E_i$ for a sequence $E_i$ going to $+\infty$.
In this and the next sections we do not explicitly mention that we consider finitely 
many moduli spaces only in order to simplify the description.
\end{rem}

\subsection{The chain level formula}
We first prove the following chain level version of Theorem \ref{theoremD1}.
\begin{prop}\label{prop18141414}
The chain map $\hat{\frak q}^{\rm f.u. {\bf b}}$ induces a ring homomorphism $ QH^*_{\frak q}(X;\Lambda_0) \to HH^*(\cL^{\rm f.u. {\bf b}},\cL^{\rm f.u. {\bf b}})$.
\end{prop}
\begin{proof}
The proof is similar to the proof of 
\cite[Theorem 2.6.1]{toric3} which dealt with the toric case and proceeds as follows.
\par
Consider the moduli space 
${\mathcal M}_{\ell+2;\vec k}((\vec{\kappa},\vec p);B)$\index[syindex]{Mell+2kappap@${\mathcal M}_{\ell+2;\vec k}((\vec{\kappa},\vec p);B)$}
introduced in Definition \ref{defMellveck}.
Among the $\ell+2$ interior marked points we 
handle the first and the second ones 
in a different way from others. 
We denote them by $z^+_{\text{\rm I}}$ and $z^+_{\text{\rm II}}$.
Other interior marked points are denoted by 
$z^+_{1}, \dots, z^+_{\ell}$.
We add data to specify $0$-th boundary marked point and define
${\mathcal M}_{\ell+2;\vec k}((\vec{\kappa},\vec p,m);B)$.\index[syindex]{Mell+2kappapm@${\mathcal M}_{\ell+2;\vec k}((\vec{\kappa},\vec p,m);B)$}
\par
Let $\mathcal M_{1;2}$ be the compactified moduli space of genus zero\index[syindex]{Mcal12@$\mathcal M_{1;2}$}
bordered Riemann surface with connected boundary,
one boundary marked point and $2$ interior marked points.
The moduli space $\mathcal M_{1;2}$ can be identified with a 2-disk $D^2$.
We describe three of its special points. See Figure \ref{Figure9-1}.
\par
One point $[\Sigma_0]$ corresponds to the origin in $D^2$.
Here $\Sigma_0$ is obtained by gluing $D^2$ and $S^2$ at the origin of 
$D^2$. Two interior marked points are on $S^2$.
\par
The two other  distinguished points, denoted $[\Sigma']$ and  $[\Sigma'']$, lie on the boundary.
Both $\Sigma'$ and $\Sigma''$ have three disk components.
Let $D^2_0$, $D^2_1$, $D^2_2$ be those components.
The boundary marked point $z_0$ is on $D^2_0$.
The components
$D^2_1$, $D^2_2$ are both attached to $D^2_0$.
We enumerate them so that 
$z_0$, $D^2_0\cap D^2_1$, $D^2_0\cap D^2_2$ 
respect the counter-clockwise cyclic  order of $\partial D^2_0$.
\par
In the case of $\Sigma'$ we have $z^+_{\text{\rm I}} \in D^2_1$, 
 $z^+_{\text{\rm II}} \in D^2_2$.
\par
In the case of $\Sigma''$ we have $z^+_{\text{\rm I}} \in D^2_2$, 
 $z^+_{\text{\rm II}} \in D^2_1$.
\begin{figure}[h]
\centering
\includegraphics[scale=0.7]{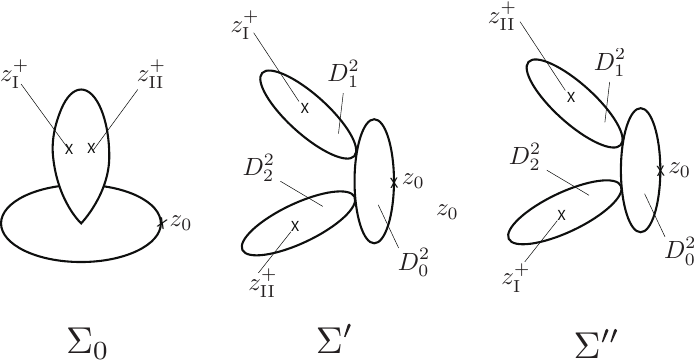}
\caption{Three points in $\mathcal M_{1,2}$}
\label{Figure9-1}
\end{figure}
\par
 We next consider the forgetful map
\begin{equation}
 \frak{forget} : 
 {\mathcal M}_{\ell+2;\vec k}((\vec{\kappa},\vec p,m);B)
 \to \mathcal M_{1;2}
\end{equation}
which forgets all but the first two interior marked points and the $0$\th boundary marked point. In more detail,  this map is obtained as follows:   let $((\Sigma,\vec z^+,\vec z),u) \in  {\mathcal M}_{\ell+2;\vec k}((\vec{\kappa},\vec p,m);B)$.
 Namely $(\Sigma,\vec z^+,\vec z)$ is a genus zero bordered Riemann surface with one boundary 
 component and marked points $\vec z^+,\vec z$, and $u : \Sigma \to X$ is a pseudo-holomorphic map.
 We assume  the boundary condition determined by $(\vec{\kappa},\vec p)$ and that $u$ is of 
 homology class $B$.
 We forget all the 3rd - ($\ell$+2)nd interior marked points except the first two and all the boundary marked points except
 the 0-th one. We also forget the map $u$. Then we obtain $(\Sigma,(z^+_{\text{\rm I}},z^+_{\text{\rm II}},z_0))$.
 It may  have unstable component. We shrink them in an obvious way and 
 obtain an element of $\mathcal M_{1;2}$.
 \par
$ {\mathcal M}_{\ell+2;\vec k}((\vec{\kappa},\vec p,m);B)$ has a Kuranishi structure 
and  the map $ \frak{forget}$ is smooth outside the fiber of three singular points 
$[\Sigma_0], [\Sigma'], [\Sigma'']$.
We need to put a Kuranishi structure on $ {\mathcal M}_{\ell+2;\vec k}((\vec{\kappa},\vec p,m);B)$
that satisfies appropriate compatibility properties at $[\Sigma_0], [\Sigma'], [\Sigma'']$.
We describe the compatibility condition below after recalling the description of the relevant spaces of stable maps.
\par
Let $ {\mathcal M}^{\text{\rm cl}}_{\ell+3}(\alpha)$ be the compactified moduli 
space of stable maps from closed Riemann surfaces of genus zero with $(\ell+3)$
marked points and of homology class $\alpha \in H_2(X;\Z)$.
Let 
\begin{equation}\label{evfromclosed}
\text{ev}^+ = (\text{ev}^+_{\text{\rm I}},\text{ev}^+_{\text{\rm II}},\text{ev}^+_0,\dots,\text{ev}^+_{\ell})
: {\mathcal M}^{\text{\rm cl}}_{\ell+3}(\alpha) \to X^{\ell+3}
\end{equation}
be the evaluation maps.
We handle the three marked points $z^+_{\text{\rm I}},  z^+_{\text{\rm II}}, z^+_0$
in a way different from other marked points.
\par
We also consider the case when only a subset (possibly empty)
%not both 
of $z^+_{\text{\rm I}}, z^+_{\text{\rm II}}$ %but only some (or none) 
%of them 
are included as marked points. (The marked point $z^+_0$
is always included.)
We write  ${\mathcal M}^{\text{\rm cl}}_{\ell+m}(\alpha)$, $m=1$ or $2$ 
in such case. We actually need to distinguish $z^+_{\text{\rm I}}$ from $z^+_{\text{\rm II}}$.
So there are two possibilities for 
${\mathcal M}^{\text{\rm cl}}_{\ell+2}(\alpha)$. 
Namely the case when $z^+_{\text{\rm I}}$  is included and 
the case when $z^+_{\text{\rm II}}$  is included.
We use the notations ${\mathcal M}^{\text{\rm cl}}_{\ell+1+{\text{\rm I}}}(\alpha)$
and ${\mathcal M}^{\text{\rm cl}}_{\ell+1+{\text{\rm II}}}(\alpha)$ to 
distinguish them.
\par
We define Kuranishi structures on them so that 
$\text{ev}_0^{\blue{+}}$ is weakly submersive.
More precisely, we have the following.\index[syindex]{Mclell@${\mathcal M}^{\text{\rm cl}}_{\ell+3}(\alpha)$}
\begin{prop}\label{Kuraproduct}
$  $ \par
\begin{enumerate}
\item
The spaces 
${\mathcal M}^{\text{\rm cl}}_{\ell+3}(\alpha)$,
${\mathcal M}^{\text{\rm cl}}_{\ell+1+{\text{\rm I}}}(\alpha)$,
${\mathcal M}^{\text{\rm cl}}_{\ell+1+{\text{\rm II}}}(\alpha)$,
${\mathcal M}^{\text{\rm cl}}_{\ell+1}(\alpha)$
have
Kuranishi structures without boundary.
\item 
The codimension two stratum 
of ${\mathcal M}^{\text{\rm cl}}_{\ell+3}(\alpha)$ 
is the union of the following fiber products over $X$:
\begin{equation}\label{formnew93}
\aligned
&{\mathcal M}^{\text{\rm cl}}_{\ell'+1}(\alpha_1)
{}_{\text{\rm ev}_0^+} \times_{\text{\rm ev}_{\ell''+1}^+}
{\mathcal M}^{\text{\rm cl}}_{(\ell''+1)+2}(\alpha_2)
\\
&
{\mathcal M}^{\text{\rm cl}}_{\ell'+1+{\text{\rm I}}}(\alpha_1)
{}_{\text{\rm ev}_0^+} \times_{\text{\rm ev}_{\ell''+1}^+}
{\mathcal M}^{\text{\rm cl}}_{(\ell''+1)+{\text{\rm II}}}(\alpha_2)
\\
&
{\mathcal M}^{\text{\rm cl}}_{\ell'+1+{\text{\rm II}}}(\alpha_1)
{}_{\text{\rm ev}_0^+} \times_{\text{\rm ev}_{\ell''+1}^+}
{\mathcal M}^{\text{\rm cl}}_{(\ell''+1)+{\text{\rm I}}}(\alpha_2)
\\
&
{\mathcal M}^{\text{\rm cl}}_{(\ell'+1)+2}(\alpha_1)
{}_{\text{\rm ev}_0^+} \times_{\text{\rm ev}_{\ell''+1}^+}
{\mathcal M}^{\text{\rm cl}}_{(\ell''+1)}(\alpha_2)
\endaligned
\end{equation}
Here the union is taken over the data described above 
together with the shuffles $(\mathbb L',\mathbb L'')$ of $\underline\ell$
such that $\# \mathbb L' = \ell'$, $\# \mathbb L'' = \ell''$, and $(\alpha_1,\alpha_2)$ with $\alpha_1 + \alpha_2 = \alpha$.
We assume that the domains are stable in the following sense: $\ell' \ge 2$, $\ell''\ge 0$ in the first case,
$\ell' \ge 1$ $\ell''\ge 1$ in the second and the third case, 
and
$\ell' \ge 0$ $\ell''\ge 2$ in the fourth case.
\par
We assume similar stability conditions for the 
codimension two stratum of 
${\mathcal M}^{\text{\rm cl}}_{\ell+1+{\text{\rm I}}}(\alpha)$, 
${\mathcal M}^{\text{\rm cl}}_{\ell+1+{\text{\rm II}}}(\alpha)$,
${\mathcal M}^{\text{\rm cl}}_{\ell+1}(\alpha)$.
\item
The Kuranishi structures have orientations that are compatible with 
$(2)$.
\item
The evaluation map $(\ref{evfromclosed})$ and
the evaluation maps of the moduli spaces appearing in $(\ref{formnew93})$  are compatible \blue{with respect to} 
the description of the boundary in $(2)$.
\item
The Kuranishi structure is invariant under permutations of 
interior marked points $z^+_1,\dots,z^+_{\ell}$.
\item
The evaluation map $\text{\rm ev}_0^{\blue{+}}$ is weakly submersive.
\end{enumerate}
\end{prop}
\begin{proof}
The construction of the Kuranishi structure is mostly the same as the one performed in 
\cite{FO}. The difference is as follows. Here, we assume the weak submersivity of 
the evaluation map $\text{\rm ev}_0^+$ only.
Since in the fiber product appearing in (2) we 
always use $\text{\rm ev}_0^+$ for the first factor 
this is enough for the inductive construction.
We also assume equivariance under the permutation of the marked points 
$z_1^+,\dots,z_k^+$ (the marked points  $z_{\text{\rm I}}^+, z_{\text{\rm II}}^+,z_0^+$ play different
roles).
% from them.
With the above remarks in mind,  the construction is the same as that of \cite{FO}.
\end{proof}
\begin{prop}\label{multiproduct} There exists a system of CF-perturbations on the moduli spaces
${\mathcal M}^{\text{\rm cl}}_{\ell+3}(\alpha)$,
${\mathcal M}^{\text{\rm cl}}_{\ell+1+{\text{\rm I}}}(\alpha)$,
${\mathcal M}^{\text{\rm cl}}_{\ell+1+{\text{\rm II}}}(\alpha)$, and 
${\mathcal M}^{\text{\rm cl}}_{\ell+1}(\alpha)$ satisfying the following properties:
\begin{enumerate}
\item
They are transversal to $0$.
\item
The evaluation map $\text{\rm ev}^+_0$ is strongly submersive with respect to  this \blue{CF-perturbation}.
\item
They are compatible at the codimension $2$ stratum 
described in Proposition {\rm \ref{Kuraproduct}}.
\item
The \blue{CF-perturbations} are invariant under the permutation of 
interior marked points $z^+_1,\dots,z^+_{\ell}$.
\end{enumerate}
\end{prop}
\begin{proof}
We can prove Proposition \ref{multiproduct} by the 
same induction process as the proof of 
Proposition \ref{Kuraproduct}.
To extend CF-perturbation of the codimension 2 strata to its neighborhood we use 
fattening described in \cite[Section 10]{linear}.
\end{proof}
\begin{rem}
In \cite{FO} we do not need to  carefully state  the compatibility at 
the codimension 2 stratum, since Gromov-Witten invariant is 
well-defined in the homology level and we can 
study each of the moduli spaces ${\mathcal M}^{\text{\rm cl}}_{\ell}(\alpha)$
separately.
In this section we also study the disk case  where the operation is 
well-defined only in the chain level. 
So it is necessary to start with stating the compatibility condition 
precisely for ${\mathcal M}^{\text{\rm cl}}_{\ell}(\alpha)$ also.
\end{rem}
Returning to Equation \eqref{decomposefrakb}, recall that we have $$
\frak b = \frak b_0 + \frak b_2 + \frak b_+
$$
where $\frak b_0\in H^0(X;\Lambda_0)$, $\frak b_2\in H^2(X;\F)$, and $\frak b_+$ vanishes in degree $0$, with degree $2$ part in $H^2(X;\Lambda_+)$, and higher degree part in $ \bigoplus_{k\ge 2} H^{2k}(X;\Lambda_0)$. We lift $\frak b_+$ to the chain level as follows:  letting $Q_i$ be closed differential forms 
 whose cohomology classes define a basis of the de Rham cohomology of $X$, we write $\frak b_+ = \sum_{i=1}^B w_i Q_i$ where  $w_i \in \Lambda_0$ 
if the degree of $Q_i$ is not 2 and $w_i \in \Lambda_+$ if the degree of $Q_i$ is 2.

Let $g_1,g_2$ be differential form on $X$. We put\index[syindex]{cupfrakb@$\cup^{\frak b}$}
%\marginpar{Need to put sign.  2026 April.}
\begin{equation}\label{closedgwdef}
\aligned
g_1 \cup^{\frak b} g_2
= 
&\sum_{i_1,\dots,i_B=0}^{\infty}\sum_{\alpha\in \pi_2(X)}
\frac{T^{\alpha\cap \omega}
\exp(\alpha\cap  \frak b_2)}{i_1!\cdots i_{B}!}w_1^{i_1} \dots w_B^{i_B} \\
&\text{\rm Corr}({\mathcal M}^{\text{\rm cl}}_{i_1+\dots + i_B+3}(\alpha);
((\text{ev}^+_{\text{\rm I}},\text{ev}^+_{\text{\rm II}},\text{ev}^+_1,
\dots,\text{ev}^+_{i_1+\dots + i_B}),\text{ev}^+_0)\\
&\hskip5.4cm (g_1,g_2,
Q_{1}^{i_1},\dots,Q_{B}^{i_B}).
\endaligned
\end{equation}
We use the system of  CF-perturbations provided in Proposition \ref{multiproduct} to define 
the correspondence in the right hand side.
\begin{lem}
The right hand side of {\rm(\ref{closedgwdef})} converges in $T$-adic topology.
\end{lem}
\begin{proof}
The proof is the same as the proof of \cite[Theorem 3.14]{fooo:bulk}, for example. 
\end{proof}
Since ${\mathcal M}^{\text{\rm cl}}_{i_1+\dots + i_B+3}(\alpha)$ 
has no boundary, it follows from Stokes' theorem \cite[Lemma 12.14]{fooo:bulk} 
\newred{and \cite[Proposition 9.26]{tech2}} that $\cup^{\frak b}$ induces a chain map
\begin{equation}\label{prod1chain}
\Omega(X) \otimes \Omega(X) \to \Omega(X)\widehat{\otimes} \Lambda_0.
\end{equation}
We can extend it to a product %structure 
on $\Omega(X)\widehat{\otimes} \Lambda_0$ 
by linearity and continuity.
Using the fact that the right hand side of (\ref{closedgwdef}) is well-defined 
at the cohomology level we can easily show that the map induced on cohomology from 
(\ref{prod1chain}) coincides with the 
$\frak b$-deformed (quantum) cup product.
\par

We now go back to the study of the moduli space of pseudo-holomorphic disks
${\mathcal M}_{\ell+2;\vec k}((\vec{\kappa},\vec p,m);B)$ and describe  the
compatibility condition mentioned above.
\begin{prop}\label{Kuraproductdisk}
There exists a Kuranishi structure on 
$[0,1]^{\vert {\frak f}_{\rm{\bf f}}\vert} \times{\mathcal M}_{\ell+2;\vec k}((\vec{\kappa},\vec p,m);B;\vec{\frak f})^{\boxplus 1}$ with 
the following properties.
For $\frak x \in \mathcal M_{1;2}$ we denote the 
fiber 
$$
\frak{forget}^{-1}(\frak x) \subset
[0,1]^{\vert {\frak f}_{\rm{\bf f}}\vert} \times{\mathcal M}_{\ell+2;\vec k}((\vec{\kappa},\vec p,m);B;\vec{\frak f})^{\boxplus 1}
$$
by
$[0,1]^{\vert {\frak f}_{\rm{\bf f}}\vert} \times{\mathcal M}_{\ell+2;\vec k}((\vec{\kappa},\vec p,m);B;\frak x;\vec{\frak f})^{\boxplus 1}$.
\begin{enumerate}
\item
If $\frak x \in \text{\rm Int}\,\,\mathcal M_{1;2} \setminus \{[\Sigma_0]\}$, 
then the Kuranishi structure on the space $[0,1]^{\vert {\frak f}_{\rm{\bf f}}\vert} \times{\mathcal M}_{\ell+2;\vec k}((\vec{\kappa},\vec p,m);B;\vec{\frak f})^{\boxplus 1}
$ induces one on $[0,1]^{\vert {\frak f}_{\rm{\bf f}}\vert} \times{\mathcal M}_{\ell+2;\vec k}
((\vec{\kappa},\vec p,m);B;\frak x;\vec{\frak f})^{\boxplus 1}$
such that its boundary is decomposed into a union of those
similar to $(\ref{218})$, $(\ref{219})$, $(\ref{220})$ (and to the corresponding statements in 
Proposition \ref{prop94}) and that 
the two distinguished interior marked points $z^+_{\text{\rm I}}$, 
$z^+_{\text{\rm II}}$ lie in the second factor. 
For example, one case among the analogue of  $(\ref{218})$ is 
\begin{equation}\label{D5formula}
{\mathcal M}_{\ell'';k''_i+1}(L_{\kappa_i};\beta)^{\boxplus 1}
\,{}_{\text{\rm ev}_0} \times_{\text{\rm ev}_{i,j}}
{\mathcal M}_{\ell'+2;\vec k'}((\vec{\kappa},\vec p,m);B';\frak x)^{\boxplus 1},
\end{equation}
when $i > 0$, for example.
Here $\beta\#B' = B$, and the union is taken over all $({\mathbb L}_1,{\mathbb L}_2) \in \text{\rm Shuff}(\ell)$, and $\#{\mathbb L}_1 = \ell'$, 
$\#{\mathbb L}_2 = \ell''$. 
\par
In Formula {\rm (\ref{D5formula})}, we put the Kuranishi structure of Proposition 
$\ref{Kuraeistspoly}$ on the first factor 
and one of this proposition on the second factor.
Then it is compatible with the identification.\footnote{
We remark that $z_{\rm I}^+$,  $z_{\rm II}^+$ must be in the second factor since $\frak x \in \text{\rm Int}\,\,\mathcal M_{1;2} \setminus \{[\Sigma_0]\}$.}
\begin{figure}[h]
\centering
\includegraphics[scale=0.7]{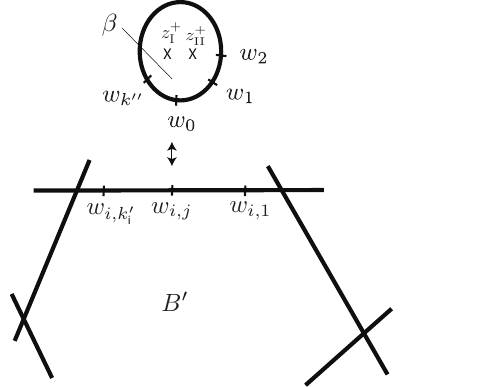}
\caption{An element of (\ref{D5formula})}
\label{Figure9-2}
\end{figure}
\item
The fiber $[0,1]^{\vert {\frak f}_{\rm{\bf f}}\vert} \times{\mathcal M}_{\ell+2;\vec k}((\vec{\kappa},\vec p,m);B;[\Sigma_0];\vec{\frak f})^{\boxplus 1}$ of the moduli space at $[\Sigma_0]$ 
is a union of fiber products
\begin{equation}\label{D6formula}
{\mathcal M}^{\text{\rm cl}}_{\ell''+2}(\alpha)
{}_{\text{\rm ev}_0^+}\times_{\text{\rm ev}_{\ell'+1}^+}
([0,1]^{\vert {\frak f}_{\rm{\bf f}}\vert} \times{\mathcal M}_{\ell'+1;\vec k}((\vec{\kappa},\vec p,m);B';\vec{\frak f}))^{\boxplus 1}
\end{equation}
where $\alpha\# B' = B$, and the union is taken over all $({\mathbb L}_1,{\mathbb L}_2) \in \text{\rm Shuff}(\ell)$, and $\#{\mathbb L}_1 = \ell'$, 
$\#{\mathbb L}_2 = \ell''$. 
The two distinguished interior marked points $z^+_{\text{\rm I}}$, 
$z^+_{\text{\rm II}}$ are in the first factor. 
\par
In Formula {\rm (\ref{D6formula})}, we put the Kuranishi structure of Proposition {\rm \ref{Kuraproduct}}
on the first factor and the Kuranishi structure of Propositions $\ref{Kuraeistspoly}$  and $\ref{prop94}$
on the second factor. Then it is compatible with this identification.
\begin{figure}[h]
\centering
\includegraphics[scale=0.7]{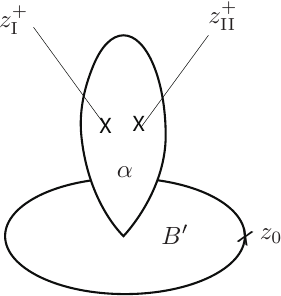}
\caption{An element of (\ref{D6formula})}
\label{Figure9-3}
\end{figure}
\item
The fiber $[0,1]^{\vert {\frak f}_{\rm{\bf f}}\vert} \times{\mathcal M}_{\ell+2;\vec k}((\vec{\kappa},\vec p,m);B;[\Sigma'];\vec{\frak f})^{\boxplus 1}$
is a union of several fiber products. There are two singular points in  the stable curve 
$[\Sigma']$. 
\blue{There are  $9$ cases.  Namely,  $z_{\rm I}^+$  (resp. $z_{\rm II}^+$) lies on one of the 
three types of components similar to the first factors of  the products $(\ref{218})$, $(\ref{219})$, $(\ref{220})$},
or those on Proposition $\ref{prop94}$.

In case when both are similar to  $(\ref{218})$ it is a fiber product 
\begin{equation}\label{case11quad}
\aligned
({\mathcal M}_{\ell^{(1)}+\text{\rm I};k^{(1)}+1}(L_{\kappa_{i(1)}};\beta_1)^{\boxplus 1}
\times 
{\mathcal M}_{\ell^{(2)}+\text{\rm II};k^{(2)}+1}&(L_{\kappa_{i(2)}};\beta_2))^{\boxplus 1}\\
\,{}_{(\text{\rm ev}_0,\text{\rm ev}_0)} \times_{(\text{\rm ev}_{i(1),j(1)}, 
\text{\rm ev}_{i(2),j(2)})}
&{\mathcal M}_{\ell'+2;\vec k'}((\vec{\kappa},\vec p,m);B')^{\boxplus 1},
\endaligned
\end{equation}
when $i(1) >  0$, for example.
\par
Here $\beta_1 \# \beta_2 \# B' = B$, and the union is taken over all 
triple shuffles $({\mathbb L}_1,{\mathbb L}_2,{\mathbb L}_3)$, and $\#{\mathbb L}_1 = \ell^{(1)}$, 
$\#{\mathbb L}_2 = \ell^{(2)}$, $\#{\mathbb L}_3 = \ell'$.
And $\vec k = (k_1,\dots,\vec k_K)$, 
$\vec k' = (k'_1,\dots,\vec k'_K)$ where
$
k'_j = k_j
$
for $j \ne i(1),i(2)$, 
$k'_{i(1)} = k_{i(1)} + k^{(1)}$, 
$k'_{i(2)} = k_{i(2)} + k^{(2)}$.
\par
In {\rm (\ref{case11quad})}  we take the Kuranishi structure of 
Proposition \ref{prop94} 
%$\ref{Kuraeistspoly}$ 
on the first and the second factors.
Then it is compatible with the identification.
\begin{figure}[h]
\centering
\includegraphics[scale=0.7]{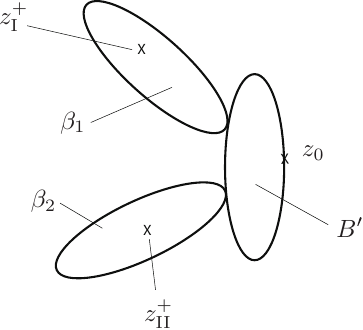}
\caption{An element of (\ref{case11quad})}
\label{Figure9-4}
\end{figure}
\par
The fiber ${\mathcal M}_{\ell+2;\vec k}((\vec{\kappa},\vec p,m);B;[\Sigma''])^{\boxplus 1}$
 is described in a similar way. For example in {\rm (\ref{case11quad})} we exchange 
$\text{\rm I}$ and $\text{\rm II}$.
\item
The Kuranishi structures are oriented compatibly with the decompositions in $(1)$, $(2)$ and $(3)$.
\item 
The Kuranishi structure is compatible with the forgetful map
of the boundary marked points labeled by $\frak f_{\bf e}$,
in a similar sense as Proposition $\ref{prop94}$ $(8)$.
\item When the coordinates $s_{i,j}$ of $[0,1]^{\vert {\frak f}_{\rm{\bf f}}\vert}$ factor 
becomes $0$ or $1$ a similar compatibility as Proposition $\ref{prop94}$ $(9), (10)$ holds.
 \item If $(i,j) \notin \vec{\frak f}$ the evaluation map
$$
\text{\rm ev}_{i,j} : 
[0,1]^{\vert \frak f_{\bf f}\vert} \times {\mathcal M}_{\ell+2;\vec k}((\vec{\kappa},\vec p,m);B;\vec{\frak f})^{\boxplus 1}
\to L_{\kappa_i}
$$
is weakly submersive.
\item
The Kuranishi structure is invariant under the permutation of 
interior marked points $z^+_1,\dots,z^+_{\ell}$.
\item
The Kuranishi structure is invariant under the cyclic permutation of the boundary marked points in the same sense 
as Proposition  $\ref{Kuraeistspoly}$.
\end{enumerate}
\end{prop}
The proof is a straightforward analogue of the proofs of 
Propositions  $\ref{Kuraeistspoly}$ and \ref{prop94} so is omitted.
We remark that 
in the case (2),(3) the union is not the disjoint union.
In case (2) the fiber ${\mathcal M}_{\ell+2;\vec k}((\vec{\kappa},\vec p,m);B;[\Sigma_0])^{\boxplus 1}$
is a `normal crossing divisor' as described in \cite[Lemma 2.6.13]{toric3}.
(The compatibility in Proposition \ref{Kuraproduct} is needed for the 
consistency at this normal crossing divisor.)
In case (3) the associated fibre is a union of several components each of which is a space with 
Kuranishi structure with corners, as defined in \cite[Lemma 2.6.22]{toric3}.
\par
We next equip them with relevant CF-perturbations as follows.
\begin{prop}\label{multiproductdisk}
For each $\frak x \in \mathcal M_{1;2}$ we can take a \blue{CF-perturbation} on $[0,1]^{\vert \frak f_{\bf f}\vert} \times{\mathcal M}_{\ell+2;\vec k}((\vec{\kappa},\vec p,m);B;\frak x;\vec{\frak f})^{\boxplus 1}$
with the following properties:
\begin{enumerate}
\item It is transversal to $0$.
\item 
It is compatible with the description of the fibers given in Proposition {\rm \ref{Kuraproductdisk}}
$(1),(2),(3)$.
\item 
It is compatible with the forgetful map 
in Proposition $\ref{Kuraproductdisk}$ $(5)$.
\item The compatibility in Proposition $\ref{Kuraproductdisk}$ $(6)$
is enhanced to the compatibility of CF-perturbations.
\item
If $(i,j) \notin \vec{\frak f}$ then the evaluation map
$$
\text{\rm ev}_{i,j} : 
[0,1]^{\vert \frak f_{\bf f}\vert} \times{\mathcal M}_{\ell+2;\vec k}((\vec{\kappa},\vec p,m);B;\frak x;\vec{\frak f})^{\boxplus 1}
\to L_{\kappa_i}
$$
is a strongly submersive with respect to our \blue{CF-perturbation}.
\item
The \blue{CF-perturbation} is invariant under the permutation of 
interior marked points $z^+_1,\dots,z^+_{\ell}$.
\item
The \blue{CF-perturbation} is invariant under the  cyclic permutation of the boundary marked points  
in the same sense as Proposition $\ref{Kuraeistspoly}$ $(7)$. 
\end{enumerate}
\end{prop}
The proof is again similar to that of Proposition \ref{existmultipolu} and is 
omitted.
See \cite[Lemma 2.6.16.2]{toric3} for the compatibility of \blue{CF-perturbations}  in case (2) and 
\cite[Lemma 2.6.16.4]{toric3}
for the compatibility of CF-perturbations in case (3).
We use `fattening' described in \cite[Proof of Proposition 10.9]{linear} to obtain a CF-perturbation satisfying Proposition \ref{multiproductdisk} (2).
\par
We next put
$$
{\mathcal M}_{\ell+2}(\vec{\kappa},\vec p;B) = {\mathcal M}_{\ell+2;\vec k'}((\vec{\kappa}',\vec p,m');B)
$$
for $\text{\rm Red}(\vec{\kappa},\vec p) = (\vec{\kappa}', \vec{p}, \vec k', m')$ as in (\ref{eq:change_notation_moduli_space}).
We define its $\frak f$ version in the same way.
\begin{defn}
For given $\frak x \in \mathcal M_{1;2}$, $\vec{\kappa}$,   we define  a map
$$
\frak Q^{\text{\bf b},\frak x}_{\vec{\kappa}}
: 
\Omega(X)\blue{[2]}^{\otimes 2}
\otimes 
BCF(\cL^{\rm form};\vec{\kappa})^+
\to CF(L_{\kappa_{0}},L_{\kappa_{1}};\C)^+ \blue{[1]},
$$
by %\marginpar{Sign to be put.  KF 2025 Aug 30}
\begin{equation}\label{Q2cat}
\aligned
\frak Q^{\text{\bf b},\frak x}_{\vec{\kappa}}&(g_1,g_2;\text{\bf h})
\\=
%\sum_{B,\ell',\vec k}
\sum_{i_1,\dots,i_B,\vec k}&T^{\omega(B) }
\frac{\rho_{\frak b,\theta}(B)}{i_1!\cdots i_{B}!}(-1)^{\maltese}
\text{\rm Corr}(
{\mathcal M}_{i_1+\dots + i_B+2}((\vec{\kappa},\vec p);B;\frak x);\\
&(\text{ev}^+_{\text{\rm I}},\text{ev}^+_{\text{\rm II}},\text{ev}^+_1,\dots,\text{ev}^+_{i_1+\dots + i_B},
\text{\rm ev}^{\partial}
, \text{ev}^{\partial}_0),\\
& (g_1,g_2,
Q_{1}^{i_1},\dots,Q_{B}^{i_B},b_{\kappa_0}^{k_0}, h_1, \ldots, h_K, b_{\kappa_K}^{k_K} ).
\endaligned
\end{equation}
(This is the case when $\text{\bf h}$ does not contain ${\bf e}^+$ or ${\bf f}$.
The case ${\bf e}^+$ and/or ${\bf f}$ is included we can use marked points 
labelled by $\frak f_{\bf e}$ and/or $\frak f_{\bf f}$ in the same way as 
Subsection \ref{subsec:opeartor} to define it.)
We define the sign $\maltese$ in the same way as  (\ref{form21ten1}). 
This collection of maps defines a cochain map\index[syindex]{Qfxrabx@$\frak Q^{\text{\bf b},\frak x}$}
$$
\frak Q^{\text{\bf b},\frak x}
:
\Omega(X) \blue{[2]}^{\otimes 2}
\to 
CH^*(\cL^{\rm form}_{\rm c.u.},\cL^{\rm form}_{\rm c.u.})^+.
$$
\end{defn}
\begin{lem}\label{lem1912}
If $\frak x = [\Sigma_0]$ then
$$
\frak Q^{\text{\bf b},\frak x}(g_1,g_2;\text{\bf x}) 
= 
\frak q^{\text{\bf b}}(g_1\cup^{\frak b} g_2;\text{\bf x}).
$$
If $\frak x = [\Sigma']$ then
$$
\frak Q^{\text{\bf b},\frak x}(g_1,g_2;\text{\bf x}) 
= 
\sum_{c}(-1)^{\maltese} \frak m^{\text{\bf b}}({\mathbf x}_c^{(5;1)}, \frak q^{\text{\bf b}}(g_1;\text{\bf x}_c^{(5;2)}),{\mathbf x}_c^{(5;3)}, \frak q^{\text{\bf b}}(g_2;\text{\bf x}_c^{(5;4)}), {\mathbf x}_c^{(5;5)}).
$$
with
$\maltese = \deg' \text{\bf x}_c^{(5;2)} +  \deg' \text{\bf x}_c^{(5;3)} + 
\deg'{\mathbf x}_c^{(5;1)}\deg g_1 + (\deg'{\mathbf x}_c^{(5;1)} +\deg'{\mathbf x}_c^{(5;2)} 
+ \deg'{\mathbf x}_c^{(5;3)})\deg g_2 +   \deg g_1$.
 %\marginpar{Sign added. KF 2025 Aug.}
\end{lem}
\begin{proof}
The first equality is a consequence of Proposition 
\ref{Kuraproductdisk} (2) and Proposition \ref{multiproductdisk} (2)
and the definition.
The second equality is a consequence of Proposition 
\ref{Kuraproductdisk} (3) and Proposition \ref{multiproductdisk} (2)
and the definition.
\end{proof}
\begin{rem}
We explain the sign $\maltese$ above.
The domain curve $\Sigma'$  has two boundary nodes.
The parameter to resolve them are $[0,\epsilon)^2$, which 
will be identified with the deformation parameter of the moduli space $\mathcal M_{1;2}$.
The moduli space $\mathcal M_{1;2}$ has a complex structure and so is canonically oriented.
When we identify $[0,\epsilon)^2$ with the tangent space of $\mathcal M_{1;2}$, 
we need to move those two parameters to the same place.
Those parameters are placed 
at the outputs of $ \frak q^{\text{\bf b}}(g_1;\text{\bf x}_c^{(5;2)})$ and of $\frak q^{\text{\bf b}}(g_2;\text{\bf x}_c^{(5;4)})$.
To move them to the same place we need to exchange one of the two factors of $[0,\epsilon)^2$
with the parameters of the chains appearing in the fiber product.
The sign we obtain is ${\deg' \text{\bf x}_c^{(5;2)} +  \deg' \text{\bf x}_c^{(5;3)} + \deg g_1}$.
To move $g_1$, $g_2$ to the beginning of the formula we have another sign, that is, 
$\deg'{\mathbf x}_c^{(5;1)}\deg g_1 + (\deg'{\mathbf x}_c^{(5;1)} +\deg'{\mathbf x}_c^{(5;2)} 
+ \deg'{\mathbf x}_c^{(5;3)})\deg g_2$.
The sum of these signs is $\maltese$.
\par
The Koszul 
sign appearing in the definition of $\cupdot$ (\ref{eq:cup_product})
is ${\deg' \text{\bf x}_c^{(5;2)} +  \deg' \text{\bf x}_c^{(5;3)}} + \deg'{\mathbf x}_c^{(5;1)}\deg g_1 + (\deg'{\mathbf x}_c^{(5;1)} +\deg'{\mathbf x}_c^{(5;2)}+\deg'{\mathbf x}_c^{(5;3)})\deg g_2$.
\par
$\maltese$  contains an extra term $\deg g_1$.
It coincides with the sign ${\deg S}$ %\marginpar{Remark added.  KF 2025 Aug.} 
in the formula (\ref{form457}). 
\par
Thus the right hand side coincides with 
$
\hat{\frak q}^{\rm f.u. {\bf b}}(g_1) \cup \hat{\frak q}^{\rm f.u. {\bf b}}(g_2).
$
Here we  obtain $\cup$ not $\cupdot$.
\end{rem}
\begin{prop}
Let $\frak x_i \in \text{\rm Int}\,\,\mathcal M_{1;2} \setminus \{[\Sigma_0]\}$.
\begin{enumerate}
\item
If $\lim_{i\to \infty}\frak x_i = [\Sigma_0]$ then
$$
\lim_{i\to \infty} 
\frak Q^{\text{\bf b},\frak x_i}_{\vec{\kappa}}(g_1,g_2;\text{\bf h})
=
\frak Q^{\text{\bf b},[\Sigma_0]}_{\vec{\kappa}}(g_1,g_2;\text{\bf h}).
$$
\item
If $\lim_{i\to \infty}\frak x_i = [\Sigma']$ then
$$
\lim_{i\to \infty} 
\frak Q^{\text{\bf b},\frak x_i}_{\vec{\kappa}}(g_1,g_2;\text{\bf h})
=
\frak Q^{\text{\bf b},[\Sigma']}_{\vec{\kappa}}(g_1,g_2;\text{\bf h}).
$$
\end{enumerate}
\end{prop}
\begin{proof}
Using the exponential 
decay estimate \cite[Theorem 6.4]{foooexp},  we can prove the proposition in the same way as \cite[Lemma 3.4.26 and 27]{toric3}.
\end{proof}
\begin{prop}
If $\frak x_i \in \text{\rm Int}\,\,\mathcal M_{1;2} \setminus \{[\Sigma_0]\}$, $i=1,2$,  then the maps
$
\frak Q^{\text{\bf b},\frak x_1}
:
\Omega(X)^{\otimes 2}
\to 
CH^*(\cL_{\rm c.u.}^{\rm form},\cL_{\rm c.u.}^{\rm form})
$
 and 
$
\frak Q^{\text{\bf b},\frak x_2}$ are chain homotopic.
\end{prop}
\begin{proof}
Choose a path $\gamma : [0,1] \to \mathcal M_{1;2} \setminus \{[\Sigma_0], 
[\Sigma'], [\Sigma'']\}$ such that 
$\gamma(0) = \frak x_1$, $\gamma(1) = \frak x_2$.
We define the parametrised moduli space 
$$
\aligned
&[0,1]^{\vert {\frak f}_{\rm{\bf f}}\vert} \times {\mathcal M}_{i_1+\dots + i_B+2;\vec k}((\vec{\kappa},\vec p);B;\gamma;\vec{\frak f})^{\boxplus 1}\\
&=
\bigcup_{t\in [0,1]} 
\{t\} \times[0,1]^{\vert {\frak f}_{\rm{\bf f}}\vert}\times {\mathcal M}_{i_1+\dots + i_B+2;\vec k}((\vec{\kappa},\vec p);B;\gamma(t);\vec{\frak f})^{\boxplus 1}.
\endaligned
$$
It has a Kuranishi structure that coincides with those defined on the space
$[0,1]^{\vert {\frak f}_{\rm{\bf f}}\vert} \times {\mathcal M}_{i_1+\dots + i_B+2;\vec k}((\vec{\kappa},\vec p);B;\frak x_i;\vec{\frak f})^{\boxplus 1}$
at $t=0,1$ (where $i=1,2$) and is also compatible with the description of the fibre in 
Proposition \ref{Kuraproductdisk} (2).
We take a \blue{CF-perturbation} with the same compatibility 
condition and the transversality condition for $\text{\rm ev}_{j,i}$.
We now put
\begin{equation}\label{H2cat}
\aligned
\frak H^{\text{\bf b},\frak x}_{\vec{\kappa}}&(g_1,g_2;\text{\bf h})
\\=
\sum_{i_1,\dots,i_B,\vec k}& (-1)^{\maltese}T^{\omega(B) }
\frac{\rho_{\frak b,\theta}(B)\rho_{b}(B)}{i_1!\cdots i_{B}!}
\text{\rm Corr}(
{\mathcal M}_{i_1+\dots + i_B+2;\vec k}((\vec{\kappa},\vec p);B;\gamma)\\
&(\text{ev}^+_{\text{\rm I}},\text{ev}^+_{\text{\rm II}},\text{ev}^+_1,\dots,\text{ev}^+_{i_1+\dots + i_B},
\text{\rm ev}^{\partial}_{1}% \times 
% \text{\rm ev}^{\partial,2}
),\text{ev}^{\partial}_0)\\
& (g_1,g_2,
Q_{1}^{i_1},\dots,Q_{B}^{i_B},b_{0}^{k_0}, h_1, \ldots, h_K, b_{K}^{k_K}).
\endaligned
\end{equation}
(This is the case when $\text{\bf h}$ does not contain ${\bf e}^+$ or ${\bf f}$.
The case ${\bf e}^+$ and/or ${\bf f}$ is included we can use marked points 
labelled by $\frak f_{\bf e}$ and/or $\frak f_{\bf f}$  in the same way as 
Subsection \ref{subsec:opeartor}.)
The sign $\maltese$ can be defined in a similar way as the beginning of Section \ref{sec:ori}.

It is straightforward to see that (\ref{H2cat}) induces the required chain homotopy.

The proof of Proposition \ref{prop18141414} is complete.\end{proof}
\subsection{Formula on canonical model.}
\label{subseqprodcan}

We next prove the same conclusion for the canonical model.  %\marginpar{This subsection is new.}
Let $g_1,g_2$ be de Rham cycles in $X$.
We define $g_1 \cup^{\frak b} g_2$ by (\ref{closedgwdef}).
Proposition \ref{prop18141414} implies that the difference
$$
\widehat{\frak q}^{{\rm f.u.} {\bf b}}(g_1 \cup^{\frak b} g_2)
-
\widehat{\frak q}^{{\rm f.u.} {\bf b}}(g_1)
 \cup
\widehat{\frak q}^{{\rm f.u.} {\bf b}}(g_2) 
$$ 
is a Hochschild coboundary.
Therefore to complete the proof of Theorem \ref{theoremD1} it suffices to prove the next lemma.
\begin{lem}
The map $\frak J$ in Lemma $\ref{lem142020}$ induces a ring homomorphism 
between Hochschild cohomologies.
\end{lem}
\begin{proof}
During this proof we compare the operators $\cupdot$ rather than comparing $\cup$.
 %\marginpar{The first line of the proof is added. KF 2025 Aug.}
\par
We use the notation $\frak M_2$ to denote the $\frak m_2$ operation in Hochschild complex.
Let $\frak T,\frak S \in CH^*(\cL^{\rm form}_{\rm c.u.}, \cL^{\rm form}_{\rm c.u.})^+$ be cycles.
By definition, the Hochschild cocycle
$
\frak U = \frak M_2(\frak T,\frak S)
$
is given by 
$$
\frak U({\bf a}) = \sum_c (-1)^{\maltese_1}\frak m^{{\rm f.u.} {\bf b}}({\bf a}_{c}^{5;1},\frak T({\bf a}_{c}^{5;2}),{\bf a}_{c}^{5;3},\frak S({\bf a}_{c}^{5;4}),{\bf a}_{c}^{5;5}).
$$
$\maltese_1= \deg \frak T\deg'{\bf a}_{c}^{5;1} + \deg \frak S \sum_{i=1}^3\deg'{\bf a}_{c}^{5;i}$.
See (\ref{eq:cup_product}).
Here $((\Delta \otimes 1\otimes 1\otimes 1) \circ (\Delta \otimes 1\otimes 1) \circ (\Delta \otimes 1) \circ \Delta)({\bf a}) = \sum_c
{\bf a}^{5;1}_{c} \otimes {\bf a}^{5;2}_{c} \otimes {\bf a}_c^{5;3}\otimes {\bf a}_c^{5;4}\otimes {\bf a}^{5;5}_{c}$.
Therefore, by definition %\marginpar{Sign put.  KF. 2024 Dec. It should be just Koszul. To be checked}
$$
\aligned
\frak J(\frak U)({\bf a}) 
=
\sum_c (-1)^{\maltese_2}{\frak m}^{\bf b}({\bf a}_{c}^{7;1},
(G'\circ\frak m^{\rm f.u.{\bf b}})(\widehat{\frak f}({\bf a}_{c}^{7;2}),
&\frak T(\widehat{\frak f}({\bf a}_{c}^{7;3})),
\widehat{\frak f}({\bf a}_{c}^{7;4}),\\
&\frak S(\widehat{\frak f}({\bf a}_{c}^{7;5})),\widehat{\frak f}({\bf a}_{c}^{7;6})),{\bf a}_{c}^{7;7}),
\endaligned
$$
with 
$$
\maltese_2 =  (\deg'{\bf a}_{c}^{7;1} + \deg'{\bf a}_{c}^{7;2})\deg \frak T + 
(\deg'{\bf a}_{c}^{7;1} + \dots + \deg'{\bf a}_{c}^{7;4})\deg \frak S.
$$ %\marginpar{The sign of this proof is corrected.  It should be OK now.  KF. 2025 Feb.}

See Figure \ref{Figure18-1}. (The definition of ${\bf a}_{c}^{7;i}$ is similar to ${\bf a}_{c}^{5;i}$.)
In case when ${\bf a}_{c}^{7;1} = {\bf a}_{c}^{7;7} = 1$ the right hand side should be replaced by 
\begin{equation}\label{newform1911}
\aligned
&(-1)^{\maltese_2}\left\{
(\Pi'\circ \frak m^{\rm f.u.{\bf b}})(\widehat{\frak f}({\bf a}_{c}^{7;2}),
\frak T(\widehat{\frak f}({\bf a}_{c}^{7;3})),
\widehat{\frak f}({\bf a}_{c}^{7;4}),\frak S(\widehat{\frak f}({\bf a}_{c}^{7;5})),\widehat{\frak f}({\bf a}_{c}^{7;6}))\right.
\\
&+
\left.{\frak m}^{{\bf b}}
\left((G'\circ \frak m^{\rm f.u.{\bf b}})(\widehat{\frak f}({\bf a}_{c}^{7;2}),
\frak T(\widehat{\frak f}({\bf a}_{c}^{7;3})),
\widehat{\frak f}({\bf a}_{c}^{7;4}),\frak S(\widehat{\frak f}({\bf a}_{c}^{7;5})),\widehat{\frak f}({\bf a}_{c}^{7;6}))\right)\right\}.
\endaligned
\end{equation}
\begin{figure}[h]
\centering
\includegraphics[scale=0.3]{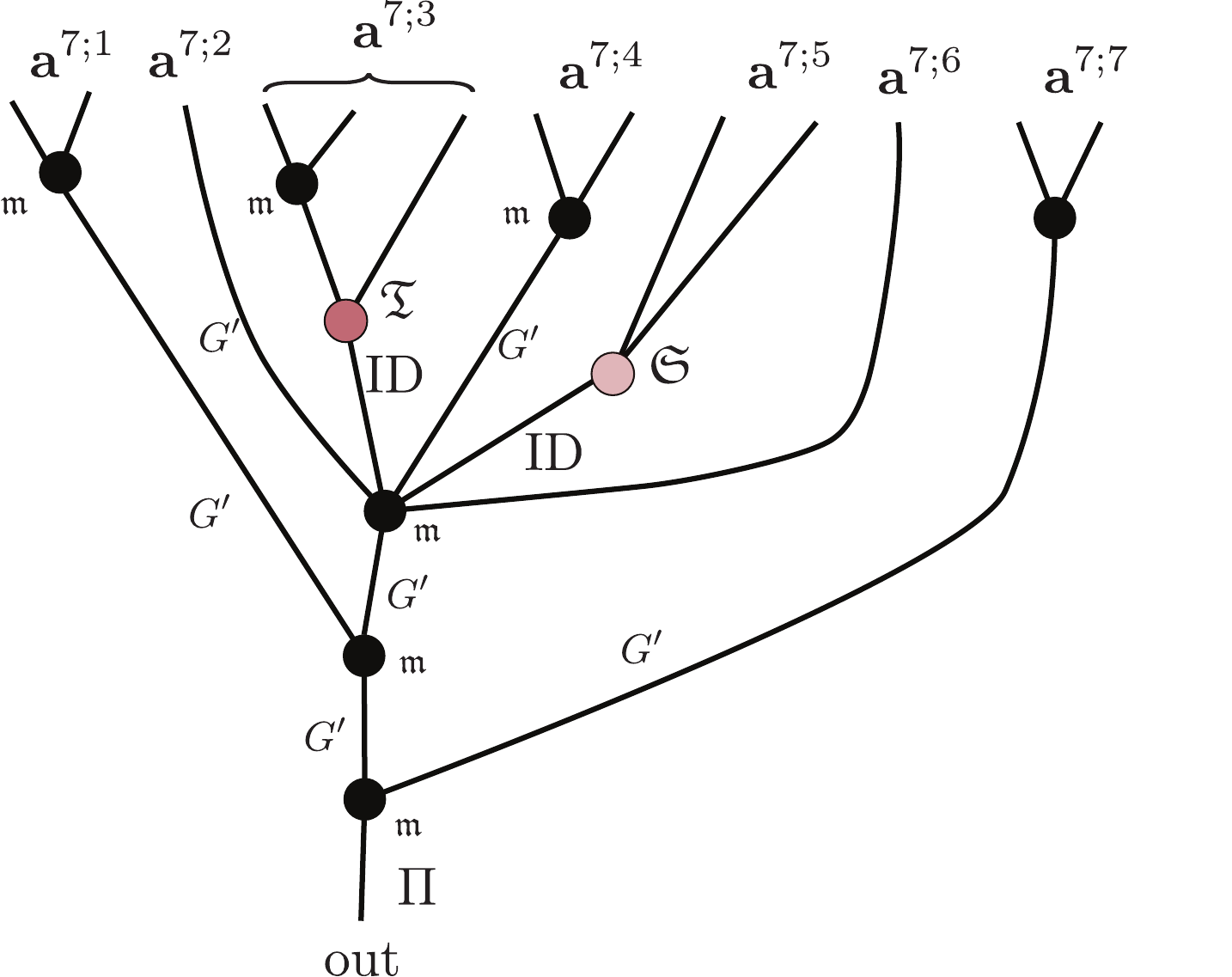}
\caption{$\frak J(\frak U)({\bf a})$}
\label{Figure18-1}
\end{figure}
Note that the operation ${\frak m}^{{\bf b}}$ is defined in Subsection \ref{subsec:red}.
Its inputs are supposed to be harmonic forms there. However we can 
plug in arbitrary forms without changing the formula.\footnote{We do not include $d$ the exterior derivative in it.} 
So ${\frak m}^{\bf b}$ in the right hand side makes sense.

We define $\mathscr K_1 \in  CH^*(\cL, \cL)$ by\index[syindex]{Kscr1@$\mathscr K_1$}
$$
\aligned
\mathscr K_1({\bf a}) 
= 
\sum_c (-1)^{\maltese_3}{\frak m}^{\bf b}({\bf a}_{c}^{7;1}, (G'\circ\frak m^{\rm f.u. {\bf b}})(\widehat{\frak f}({\bf a}_{c}^{7;2}),
&\frak H(\frak T)({\bf a}_{c}^{7;3}),
\widehat{\frak f}({\bf a}_{c}^{7;4}),\\
&\frak S(\widehat{\frak f}({\bf a}_{c}^{7;5})),\widehat{\frak f}({\bf a}_{c}^{7;6})),{\bf a}_{c}^{7;7})
\endaligned
$$
with 
$$
\maltese_3 = (\deg'{\bf a}_{c}^{7;1} + \deg'{\bf a}_{c}^{7;2})(\deg \frak T +1) + 
(\deg'{\bf a}_{c}^{7;1} + \dots + \deg'{\bf a}_{c}^{7;4})\deg \frak S.
$$
Here $\frak H$ is defined  right above Lemma \ref{lem1421}.
See Figure \ref{Figure18-2}.
In case when ${\bf a}_{c}^{7;1} = {\bf a}_{c}^{7;7} = 1$ we replace 
corresponding term in a similar way as (\ref{newform1911}). 
\par
Then using Lemma \ref{lem1421}\footnote{Or looking  Figure \ref{Figure18-2} and 
discussing in the same way as the proof of Lemma \ref{lem1421}} we can show:
 %\marginpar{Sign  put.  KF. 2025 Jan. It is just Koszul
%but to be checked.}
$$
\aligned
&\delta^H(\mathscr K_1)({\bf a})\\
=&
\sum_c (-1)^{\maltese_2}{\frak m}^{{\bf b}}({\bf a}_{c}^{7;1}, (G'\circ\frak m^{\rm f.u.{\bf b}})(\widehat{\frak f}({\bf a}_{c}^{7;2}),\frak J(\frak T)({\bf a}_{c}^{7;3}),
\widehat{\frak f}({\bf a}_{c}^{7;4}),\\
&\qquad\qquad\qquad\qquad\qquad\qquad\qquad\qquad\,\,\,
\frak S(\widehat{\frak f}({\bf a}_{c}^{7;5})),\widehat{\frak f}({\bf a}_{c}^{7;6})),{\bf a}_{c}^{7;7})\\
&+
\sum_c (-1)^{\maltese_2}{\frak m}^{{\bf b}}({\bf a}_{c}^{7;1}, (G'\circ\frak m^{\rm f.u.{\bf b}})(\widehat{\frak f}({\bf a}_{c}^{7;2}),\frak T(\widehat{\frak f}({\bf a}_{c}^{7;3})),
\widehat{\frak f}({\bf a}_{c}^{7;4}),\\
&\qquad\qquad\qquad\qquad\qquad\qquad\qquad\qquad\quad\,\,\,
\frak S(\widehat{\frak f}({\bf a}_{c}^{7;5})),\widehat{\frak f}({\bf a}_{c}^{7;6})),{\bf a}_{c}^{7;7}).
\endaligned
$$
Where ${\maltese_2}$ is as above.
\begin{figure}[h]
\centering
\includegraphics[scale=0.3]{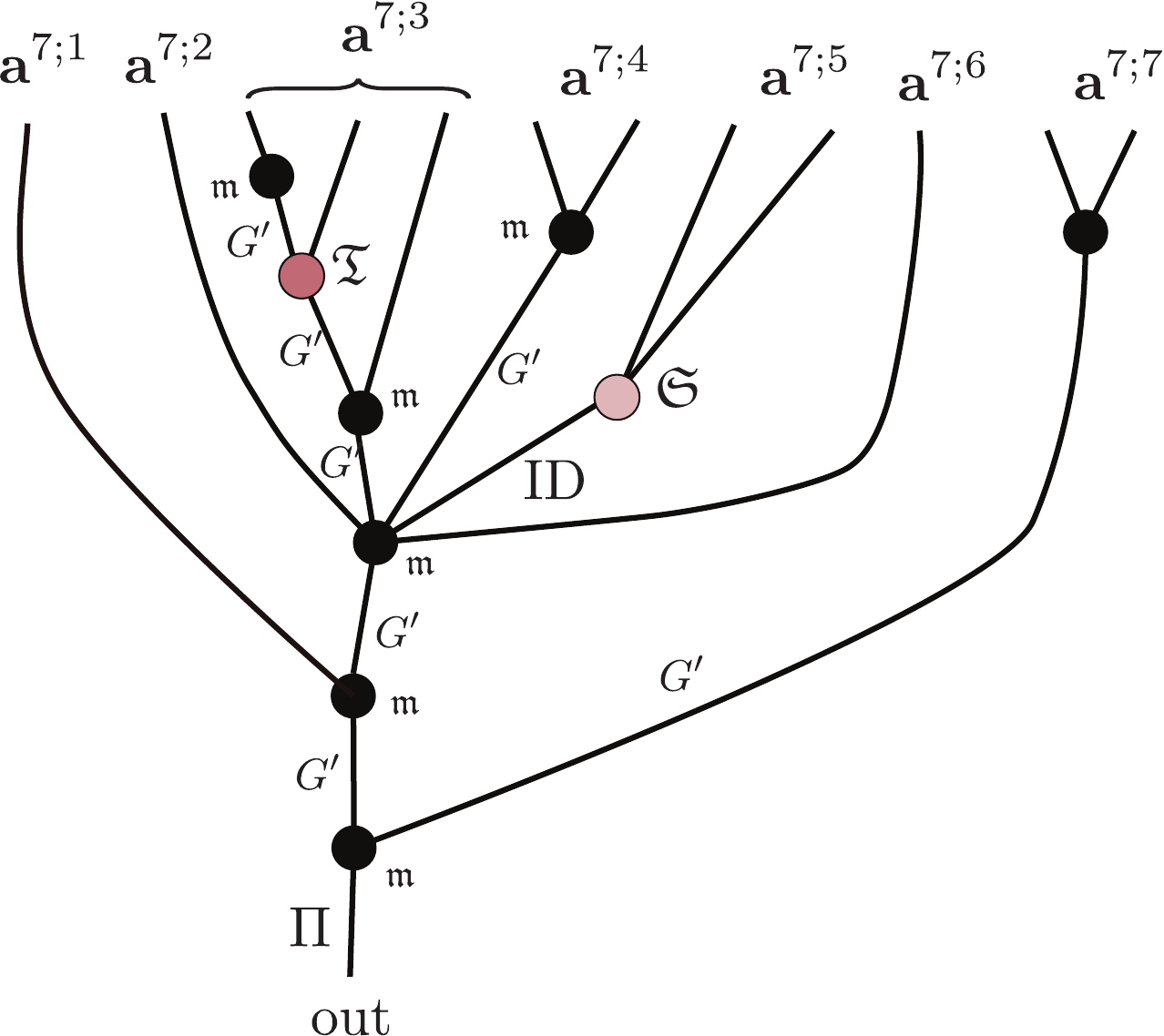}
\caption{$\mathscr K_1({\bf a})$}
\label{Figure18-2}
\end{figure}
In case when ${\bf a}_{c}^{7;1} = {\bf a}_{c}^{7;7} = 1$ we replace 
corresponding term in a similar way as (\ref{newform1911}). 

We next define $\mathscr K_2 \in  CH^*(\cL, \cL)$ by\index[syindex]{Kscr2@$\mathscr K_2$}
$$
\aligned
\mathscr K_2({\bf a}) 
= 
\sum_c (-1)^{\maltese_4}{\frak m}^{{\bf b}}(&{\bf a}_{c}^{9;1},(G'\circ\frak m^{\rm f.u.{\bf b}})(\widehat{\frak f}({\bf a}_{c}^{9;2}),
\frak f({\bf a}_{c}^{9;3},\frak J(\frak T)({\bf a}_{c}^{9;4}),{\bf a}_{c}^{9;5)}),
\\
&\widehat{\frak f}({\bf a}_{c}^{9;6}),\frak H(\frak S)({\bf a}_{c}^{9;7}),\widehat{\frak f}({\bf a}_{c}^{9;8})),
{\bf a}_{c}^{9;9})
\endaligned
$$
with 
$$
\maltese_4 =  (\deg'{\bf a}_{c}^{9;1} + \deg'{\bf a}_{c}^{9;2} + \deg'{\bf a}_{c}^{9;3})\deg \frak T + 
(\deg'{\bf a}_{c}^{9;1} + \dots + \deg'{\bf a}_{c}^{9;6})(\deg \frak S+1).
$$
Here $\frak H$ is as in right after Lemma \ref{lem142020}.
Note here we again abuse a notation a bit. 
In fact the input of $\frak m^{{\bf b}}$ was assumed to be harmonic forms.
The input 
$$(G'\circ\frak m^{\rm f.u.{\bf b}})(\widehat{\frak f}({\bf a}_{c}^{9;2}),
\frak f({\bf a}_{c}^{9;3},\frak J(\frak T)({\bf a}_{c}^{9;4}),{\bf a}_{c}^{9;5)}),
\widehat{\frak f}({\bf a}_{c}^{9;6}),\\
\frak H(\frak S)({\bf a}_{c}^{9;7}),\widehat{\frak f}({\bf a}_{c}^{9;8}))
$$
is actually not harmonic. However we can use exactly the same formula to define 
$\frak m^{{\bf b}}$ in this case. See Figure \ref{Figure18-3}.
(Note that the vertices to which $\frak T$ and $\frak S$ are assigned 
and the edge to which ${\rm ID}$ is assigned, both must be `above' the 
blue vertex,\footnote{The blue vertex is the (unique) `lowest' vertex 
such that if we remove it then $\frak T$ and $\frak S$ will be in  different connected
components.} which corresponds to the operation $\frak m^{\rm f.u.{\bf b}}$
appearing in the formula defining $\mathscr K_2$.)
In case when ${\bf a}_{c}^{9;1} = {\bf a}_{c}^{9;9} = 1$ we replace 
corresponding term in a similar way as (\ref{newform1911}).
We do not repeat this kinds of remark until the end of this subsection. 
\par
Using Lemma \ref{lem1421} (or discussing in the same way 
as the proof of Lemma \ref{lem1421}) again we can show:
$$
\aligned
&\delta^H(\mathscr K_2)({\bf a})\\
=&
\sum_c (-1)^{\maltese_2} {\frak m}^{{\bf b}}({\bf a}_{c}^{7;1}, (G'\circ\frak m^{\rm f.u.{\bf b}})(\widehat{\frak f}({\bf a}_{c}^{7;2}),\frak J(\frak T)({\bf a}_{c}^{7;3}),
\widehat{\frak f}({\bf a}_{c}^{7;4}),\\
&\qquad\qquad\qquad\qquad\qquad\qquad\qquad\qquad\,\,
\frak S(\widehat{\frak f}({\bf a}_{c}^{7;5})),\widehat{\frak f}({\bf a}_{c}^{7;6})),{\bf a}_{c}^{7;7})\\
&+
\sum_c (-1)^{\maltese_2}{\frak m}^{{\bf b}}({\bf a}_{c}^{7;1}, (G'\circ\frak m^{\rm f.u.{\bf b}})(\widehat{\frak f}({\bf a}_{c}^{7;2}),\frak J(\frak T)({\bf a}_{c}^{7;3}),
\widehat{\frak f}({\bf a}_{c}^{7;4}),\\
&\qquad\qquad\qquad\qquad\qquad\qquad\qquad\qquad\quad
\frak J(\frak S)({\bf a}_{c}^{7;5}),\widehat{\frak f}({\bf a}_{c}^{7;6})),{\bf a}_{c}^{7;7}).
\endaligned
$$
Here $\maltese_2$ is as above. 

\begin{figure}[h]
\centering
\includegraphics[scale=0.3]{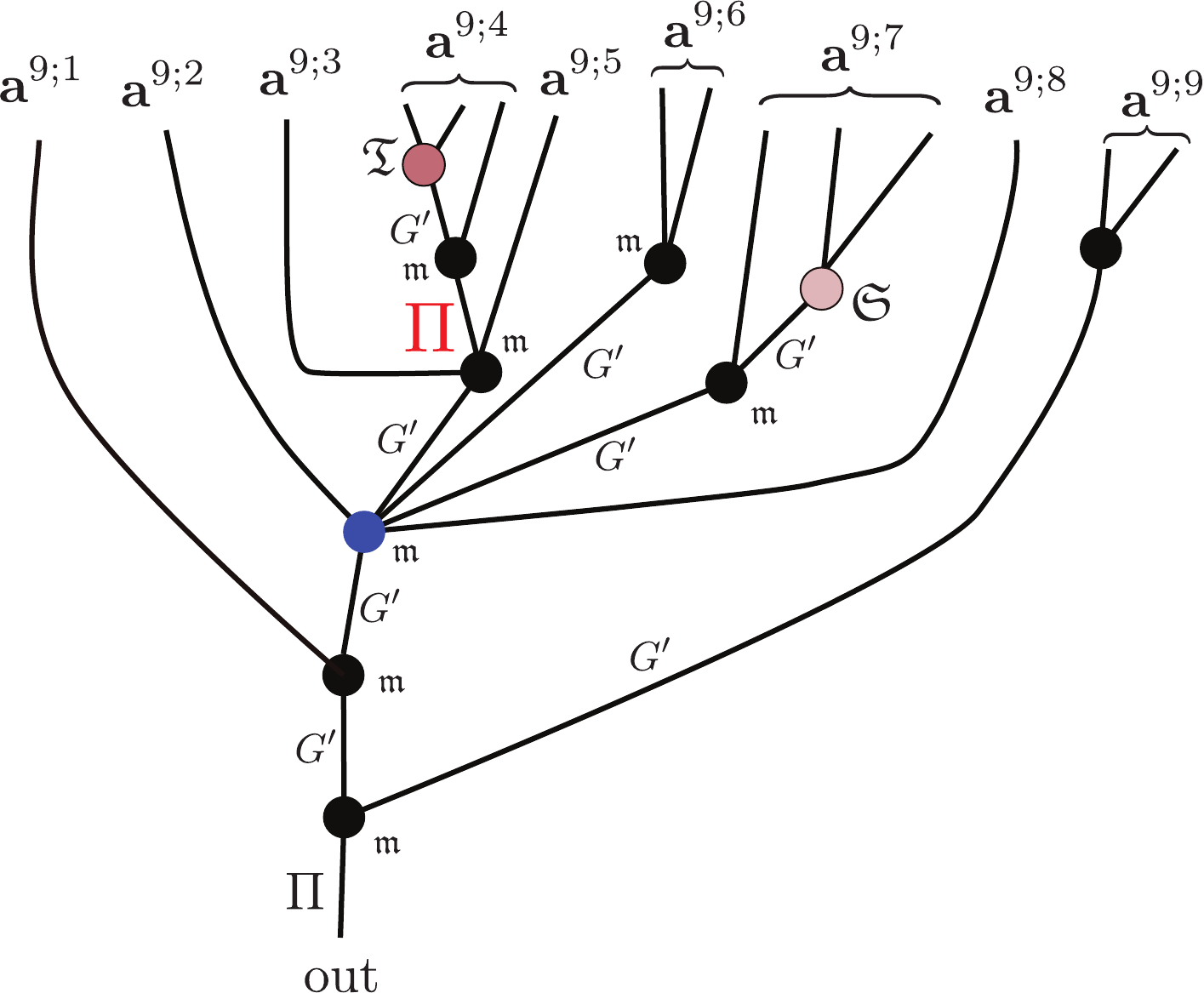}
\caption{$\mathscr K_2({\bf a})$}
\label{Figure18-3}
\end{figure}

Since the formula
$$
\aligned
&\sum_c (-1)^{\maltese_2}{\frak m}^{{\bf b}}({\bf a}_{c}^{7;1}, (G'\circ\frak m^{\rm f.u.{\bf b}})(\widehat{\frak f}({\bf a}_{c}^{7;2}),\frak J(\frak T)({\bf a}_{c}^{7;3}),
\widehat{\frak f}({\bf a}_{c}^{7;4}),\\
&\qquad\qquad\qquad\qquad\qquad\qquad\qquad\qquad
\frak J(\frak S)({\bf a}_{c}^{7;5}),\widehat{\frak f}({\bf a}_{c}^{7;6}),{\bf a}_{c}^{7;7}) \\
&= \frak M_2(\frak J(\frak T),\frak J(\frak S))({\bf a})
\endaligned
$$
is immediate from the definition of $A_{\infty}$ structure in $\cL$,\footnote{
In fact 
$$
\aligned
&\sum_c \pm{\frak m}^{{\bf b}}({\bf a}_{c}^{7;1}, (G'\circ\frak m^{\rm f.u.{\bf b}})(\widehat{\frak f}({\bf a}_{c}^{7;2}),\frak J(\frak T)({\bf a}_{c}^{7;3}),
\widehat{\frak f}({\bf a}_{c}^{7;4}),\frak J(\frak S)({\bf a}_{c}^{7;5}),\widehat{\frak f}({\bf a}_{c}^{7;6}),{\bf a}_{c}^{7;7})
\\
&=
\sum_c\pm \frak m^{\bf b}({\bf a}_{c}^{5;1},\frak J(\frak T)({\bf a}_{c}^{5;2}),{\bf a}_{c}^{5;3},\frak J(\frak S)({\bf a}_{c}^{5;4}),
{\bf a}_{c}^{5;5}).
\endaligned
$$}
the lemma follows.
\end{proof}
The proof of Theorem \ref{theoremD1} is complete.
\end{proof}

\section{The Cardy relation.}
\label{sec:annuli}

%input{section4}
In this section we define the moduli space of 
pseudo-holomorphic annuli, which is the   key ingredient in the proofs of 
Theorems \ref{maintheorem1} and  \ref{maintheorem2}
and discuss its basic properties.

\subsection{Inner product $Z$ among Hochschild homologies.}
\label{inproZ}

To state the conclusion of this section, that is, the Cardy relation, we first need to consider a pairing on Hochschild homology, which we denote by $Z = Z({\bf x}, {\bf y})$.\footnote{This pairing is called Shklyarov paring or Mukai paring 
in certain literatures.} Explanation of the construction of this pairing is in order. 
As before we fix a bulk class $\fb$ and a background class ${\rm st} \in H_2(X; \Z_2)$. 
Then let ${\bf  L} = \{(L_{\kappa}, \theta_{\kappa},b_{\kappa}) \mid 
{\kappa}_1, \dots, {\kappa}_K\}$ be a collection of ${\rm st}$-relatively spin pairwise transversal Lagrangian submanifolds enhanced by weak bounding cochains. We denote by $\cL$ the cyclic filtered $A_{\infty}$ category which they generate as in Section \ref{canonical}. We also consider another such collection 
$\bU = \{(U_{\upsilon},\theta_{\upsilon},b_{\upsilon}) \mid {\upsilon} =1,\dots, R\}$ and denote by $\sU$  the associated cyclic filtered $A_{\infty}$-category. 
\par
We also assume the following:
\begin{asump}\label{LUtrans}
For all $\kappa$ and ${\upsilon}$, $L_{\kappa}$ is transversal to 
$U_{{\upsilon}}$.
\end{asump}
We will define a pairing\index[syindex]{Z@$Z$} 
\begin{equation}\label{productZ}
Z : HH_*(\cL;\Lambda_0) \otimes
HH_*(\sU;\Lambda_0)
\to \Lambda_0,
\end{equation}
under this assumption.
We need to set up some notation to define it.
\par
Let
$$
\text{\bf x} = x_0 \otimes \cdots \otimes x_k 
\in CH_*(\cL)  
$$
be a Hochschild chain.\index[syindex]{DeltaH@$\Delta_H$} 
We define 
\begin{equation}\label{deltaH}
\Delta_H(\text{\bf x})
= \sum_{0\le i < j \le k}
(-1)^{\maltese}
(x_{j+1} \otimes \cdots \otimes  x_0 \otimes \cdots \otimes x_i)
\otimes
(x_{i+1} \otimes \cdots \otimes x_j),
\end{equation}
where
$
\maltese = (\deg' x_{j+1} + \dots + \deg' x_{k})
(\deg' x_{0} + \dots + \deg' x_{j})
$. Equation 
(\ref{deltaH}) defines a map %\marginpar{The image does not seem to be in the Hochschild.}
$$
\Delta_H : CH_*(\cL) \to BC_*(\cL)\otimes BC_*(\cL).
$$
%We remark that
%$$
%\delta_H(\text{\bf x}) 
%= (\frak m \widehat{\otimes} 1 + 1 \widehat{\otimes}\frak m)\circ \Delta_H.
%$$
We also put
$$
\Delta_H^{k-1} = 
\underbrace{
(1\otimes \cdots \otimes 1 \otimes \Delta)
\circ \cdots \circ (1 \otimes \Delta)}_{k-2} \circ  \Delta_H
: CH_*(\cL) \to BC_*(\cL)^{\otimes k}.
$$
Again, using Sweedler's notation, we write\index[syindex]{xxbdcHki@$\text{\bf x}_c^{(H;k;i)}$}
\begin{equation}\label{DeltaH}
\Delta_H^{k-1}(\text{\bf x})
=
\sum_c (-1)^{\maltese} \text{\bf x}_c^{(H;k;1)} \otimes \cdots \otimes 
\text{\bf x}_c^{(H;k;k)}.
\end{equation}
Here $\maltese$ is the Kozule sign appearing in (\ref{deltaH}). 
\begin{rem} 
We recall that the Hochschild differential is given by  
$$
\delta_H(\text{\bf x}) 
= \sum_{c_1} (-1)^{\maltese}\frak m(\text{\bf x}_{c_1}^{(H;2;1)}) \otimes
\text{\bf x}_{c_1}^{(H;2;2)}
+
\sum_{c_2} 
(-1)^{\deg' \text{\bf x}^{(3;1)}_{c_2}}
\text{\bf x}_{c_2}^{(3;1)} \otimes
\frak m(\text{\bf x}_{c_2}^{(3;2)}) \otimes \text{\bf x}_{c_2}^{(3;3)}.
$$
(Note $\text{\bf x}_{c_2}^{(3;1)}$ in the second term is not 
$\text{\bf x}_{c_2}^{(H;3;1)}$.) 
Here the second term is the sum with the case $\text{\bf x}_{c_2}^{(3;1)} = 1$ being excluded.
\end{rem}

In order to describe the sign in the pairing $Z$, we also need to recall a basic consequence of the Koszul sign conventions: assume that $V$ and $V^\vee$ are $n$-dual cochain complexes over $\Lambda$, in the sense that we have a cochain map of degree $-n$
\begin{equation}
\langle \_, \_ \rangle \co  V \otimes V^{\vee} \to \Lambda
\end{equation}
which induces a perfect pairing of vector spaces. Given bases $\{v_i\}$ and $\{v_i^\vee\}$, which are dual in the sense that $\langle v_i, v_j^\vee \rangle = \delta_{ij}$ we have:
\begin{lem} \label{lem:sign_coevaluation_pairing}
The degree $n$ cochain map\index[syindex]{coev@$\coev$}
\begin{align}
\coev \co  \Lambda & \to   V \otimes V^{\vee} \\
1 & \mapsto \sum_{i} (-1)^{\deg' v_i} v_i \otimes v_i^\vee
\end{align}
fits in a commutative diagram:
\begin{equation}
  \xymatrix{ V \ar[r]^{\cong} \ar[d]^{=} & V \otimes  \Lambda \ar[rrr]^{v \otimes 1 \mapsto (-1)^{n\deg'(v)} v \otimes \coev(1)} &  && V \otimes V \otimes V^{\vee} \ar[d]^{1 \otimes \tau} \\
V  & \Lambda \otimes V  \ar[l]^{\cong} & &&  V \otimes V^{\vee} \otimes V    \ar[lll]^{\langle \_, \_ \rangle \otimes \id }.}
\end{equation}
where $\tau(x\otimes y) = (-1)^{\deg' x\deg' y}y \otimes x$.
\end{lem}

\begin{defn}\label{defnZ}
For $\text{\bf x} \in CH_*(\cL)$ 
and $\text{\bf y} \in CH_*(\sU)$  we define\footnote{In (\ref{eq:formula_Z-pairing})
we omit the symbol ${\bf b}$ in $\frak m^{\bf b}$.  In this section we also 
omit ${\bf b}$ in $\widehat{\frak p}^{\bf b}$ etc.} 
 
\begin{equation} \label{eq:formula_Z-pairing}
\aligned
Z(\text{\bf x},\text{\bf y}) 
= \sum_{c_1,c_2}\sum_{f^1,f^2}
(-1)^{\maltese}
&\langle 
\frak m(\text{\bf x}_{c_1}^{(H;2;1)},f^1,
\text{\bf y}_{c_2}^{(H;2;2)}),f^{2\vee}
\rangle_{\rm cyc}\\
&\cdot\langle 
\frak m(\text{\bf x}_{c_1}^{(H;2;2)},f^2,\text{\bf y}_{c_2}^{(H;2;1)}
), f^{1\vee}
\rangle_{\rm cyc}.
\endaligned
\end{equation} %\marginpar{Formula changed so that it is closer to geometric argument in Section 21.  KF 2025 Nov}
See Figure \ref{Figure10-1new}.
Here $(-1)^\maltese$ is the product of: (i) the Koszul sign associated to the  permutation 
 
\begin{equation}\label{form199}
\aligned
&\langle \cdot \rangle_{(1)}, \langle \cdot \rangle_{(2)},\frak m_{(1)}, \frak m_{(2)},     \text{\bf x},  \text{\bf y},   f^1,  f^{1\vee},  f^2,  f^{2\vee}, \\
&\mapsto \langle \cdot \rangle_{(1)},\frak m_{(1)},  \text{\bf x}_{c_1}^{(H;2;1)}, f^1, \text{\bf y}_{c_2}^{(H;2;2)},   f^{2\vee},\langle \cdot \rangle_{(2)},
\frak m_{(2)},  \text{\bf x}_{c_1}^{(H;2;2)}, f^2, \text{\bf y}_{c_2}^{(H;2;1)}, f^{1\vee}
\endaligned
\end{equation}
(ii) with, the sign $\deg'f^1$ and $\deg'f^2$ from Lemma \ref{lem:sign_coevaluation_pairing}, (iii)
$(n+1)\deg' \text{\bf x}$.
 %\marginpar{This part are to be slightly modified.  KF 2025 Feb.}\marginpar{(\ref{form199}) etc is modified. KF Feb 2025
%More KF Jan 2026}
Here we count the degree of $ \frak m$ to be one and the degree of $\langle \cdot \rangle$ to be $n$.\footnote{To obtain the Koszul sign 
we use the shifted  degree of $ \text{\bf x}$, $ \text{\bf y}$, $f^1,f^2$ etc. but do not shift  degree of $\frak m$ and 
$\langle \cdot \rangle$. In other words, we use the degree $1$ and $n$ respectively.}
(Here we put suffix $(1)$ and $(2)$ to $\langle \rangle$ and $\frak m$ to distinguish them.)
\begin{rem}
The extra sign $(-1)^{(n+1) \deg' {\bf y}}=(-1)^{(n+1) \deg' {\bf x}}$  may be understood as the Koszul sign 
of the permutation
$$
\langle \rangle_{{\rm PD}_X}, \frak p, \text{\bf x}, \frak p, \text{\bf y}
\mapsto Z, \text{\bf x}, \text{\bf y}
$$
Note $\langle \rangle_{{\rm PD}_X}$ has degree $-2n$, $\frak p$ has degree $n+1$ 
(using  $\deg' {\bf x}$ as the domain degree)
and $Z$ has degree $0$. So we might identify $\langle \rangle_{{\rm PD}_X}$, $\frak p$, $\frak p$
with $Z$ to define the  Koszul sign.
\end{rem}
\begin{rem}
The Poincar\'e duality of de Rham cohomology of $X$ is anti symmetric with respect the grading.  
$\frak p$ has degree $n$ from Hochschild homology (with our degree convention) to 
the ambient homology.   So in case $\frak p$ is an isomorphism (\ref{form2015}) implies that 
$Z$ has symmetry
\begin{equation}
Z({\bf x},{\bf y}) = (-1)^{(n+ \deg' {\bf x} + 1)(n+ \deg' {\bf y} + 1)}Z({\bf y},{\bf x}).
\end{equation}
 %\marginpar{Remarks added but to be checked.  KF 2026 Jan.}
\par
In \cite{toric3} we studied the case when $ {\bf x} =  {\bf y}$ is a delta form supported at a point.
In this case since $n+ \deg' {\bf x} + 1 \equiv 0 \mod 2$, $Z$ is symmetric.
Since toric manifolds have only even degree cohomology the Poincar\'e duality is symmetric.
Note that in case for example the Potential function is Morse, 
the Hochschild homology of $A_{\infty}$ category splits into $1$ dimensional 
factors which are orthogonal to each other with respect to $Z$.  So $Z$ must be 
symmetric.
\end{rem}
\begin{figure}[h]
\centering
\includegraphics[scale=0.6]{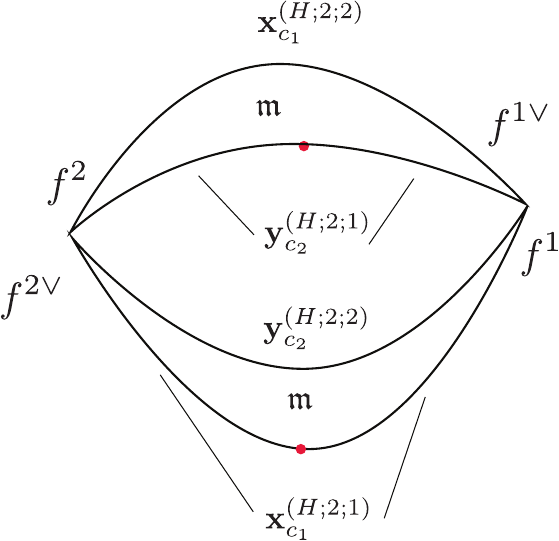}
\caption{Pairing $Z$}
\label{Figure10-1new}
\end{figure} %\marginpar{Figure changed a bit. KF 2025 Nov.}
\end{defn}
The sum over $f^1,f^2$ is taken as follows: we put
$$
\text{\bf x}_{c_1}^{(H;2;1)}
= x_1^{c_1} \otimes \cdots \otimes x_{a(c_1)}^{c_1}, 
\quad
\text{\bf y}_{c_2}^{(H;2;1)}
= y_1^{c_2} \otimes \cdots \otimes y_{b(c_2)}^{c_2}.
$$
We define $\kappa(c_1;1)$, $\kappa(c_1;2)$,
${\upsilon}(c_2;1)$, ${\upsilon}(c_2;2)$ by
$$
\aligned
x_1^{c_1} &\in CF(L_{\kappa'},L_{\kappa(c_1;1)}), 
\quad
x_{a(c_1)}^{c_1} \in CF(L_{\kappa''},L_{\kappa(c_1;2)}), \\
y_1^{c_2} &\in CF(U_{{\upsilon}(c_2;1)},U_{{\upsilon}'}), 
\quad
y_{b(c_2)}^{c_2} \in CF(U_{{\upsilon}(c_2;2)}, U_{{\upsilon}''}).
\endaligned
$$
Then $\sum_{f^1,f^2}$ is the sum over%Added by M.A. Aug 30
\begin{align}
f^1 & \textrm{ a generator of } CF^{*}(L_{\kappa(c_1;1)},U_{{\upsilon}(c_2;2)})  \\
f^2 & \textrm{ a generator of }CF^{*}(L_{\kappa(c_1;2)}, U_{{\upsilon}(c_2;1)}).
\end{align}
Note that $f^1$ corresponds to a point in $L_{\kappa(c_1;1)} \cap U_{{\upsilon}(c_2;2)}$.
We denote by $f^{1\vee}$ the same point regarded as 
a point in $U_{{\upsilon}(c_2;2)} \cap  L_{\kappa(c_1;1)}$, with sign correction so that $\langle f^1,f^{1\vee}\rangle_{\rm cyc} = 1$.
The meaning of the notation $f^2$, $f^{2\vee}$ is similar.
\par
We can rewrite \eqref{eq:formula_Z-pairing} to
\begin{equation} \label{eq:formula_Z-pairing2}
\aligned
&Z(\text{\bf x},\text{\bf y}) \\
&= \sum_{c_1,c_2}\sum_{f^1}
(-1)^{\maltese'}
\langle 
\frak m(\text{\bf x}_{c_1}^{(H;2;2)},\frak m(\text{\bf x}_{c_1}^{(H;2;1)},f^1,
\text{\bf y}_{c_2}^{(H;2;2)}),\text{\bf y}_{c_2}^{(H;2;1)}
), f^{1\vee}
\rangle.
\endaligned
\end{equation}
with sign
$$
\aligned
\maltese' =  &\deg f^1  + \deg \frak m\deg'  {\bf x}_{c_1}^{(H;2,2)} + \deg' {\bf x} \deg' f^1
\\
&+ \deg'  {\bf x}_{c_1}^{(H;2,2)}\deg'  {\bf x}_{c_1}^{(H;2,1)}
+ \deg'  {\bf y}_{c_2}^{(H;2,2)}\deg'  {\bf y}_{c_2}^{(H;2,2)}
\endaligned
$$
where 
$\maltese'$ (which depends on $c_1,c_2$)
is the sum of the Koszul sign of 
$$
\langle \cdot\rangle, \frak m, \frak m,  f^{1\vee},f^1,\frak m, {\bf x}, {\bf y}\,\,\,\,
\mapsto\,\,\,\,
\langle \cdot\rangle, \frak m, {\bf x}_{c_1}^{(H;2,2)}, \frak m, {\bf x}_{c_1}^{(H;2,1)}, f^1, {\bf y}_{c_2}^{(H;2,2)},{\bf y}_{c_2}^{(H;2,1)}, f^{1\vee},
$$
and $\deg f^1$ which appears to take super trace.
We will prove that (\ref{eq:formula_Z-pairing}) and (\ref{eq:formula_Z-pairing2}) are equivalent 
(together with sign)
at the last part of the proof of Theorem \ref{ZisPD} given at the end of Subsection \ref{sec:signCardy}.

\begin{lem}\label{lem205} %\marginpar{proof added KF. 2025 Nov.}
\begin{equation}\label{lem205form}
Z(\delta_H \text{\bf x},\text{\bf y}) 
+ (-1)^{\deg' \text{\bf x}}Z(\text{\bf x},\delta_H \text{\bf y})  = 0.
\end{equation}
\end{lem}
\begin{proof}
Schematically the left  hand side is described by Figure \ref{Figure29new}.
\begin{figure}[h]
\centering
\includegraphics[scale=0.5]{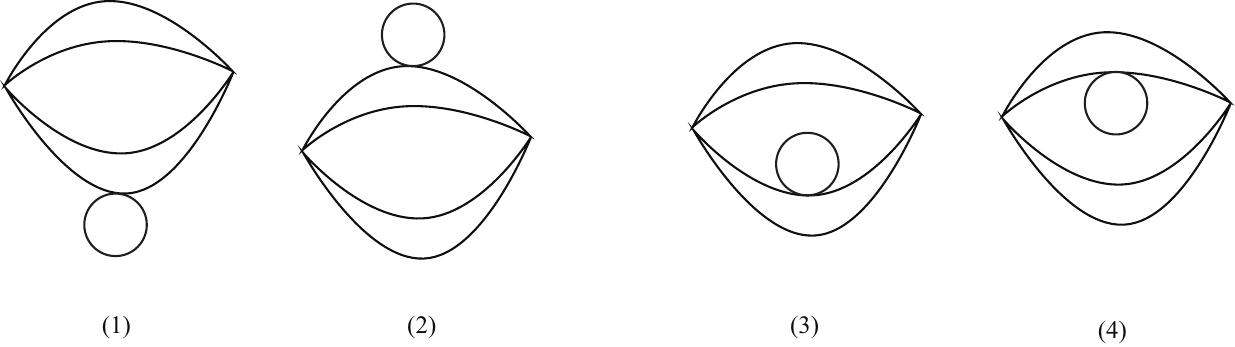}
\caption{Left hand side of (\ref{lem205form}).}
\label{Figure29new}
\end{figure}
(1) and (2) (resp. (3) and (4)) in the figure give the first term (resp. the second term) of  (\ref{lem205form}).
Here we do not distinguish whether $0$-th marked point is on the bubble or not.
\par
Applying the $A_{\infty}$ relation to the lower (resp. upper) part of Figure \ref{Figure10-1new}
gives Figure \ref{Figure30new} (resp. \ref{Figure31new}).
\begin{figure}[h]
\centering
\includegraphics[scale=0.5]{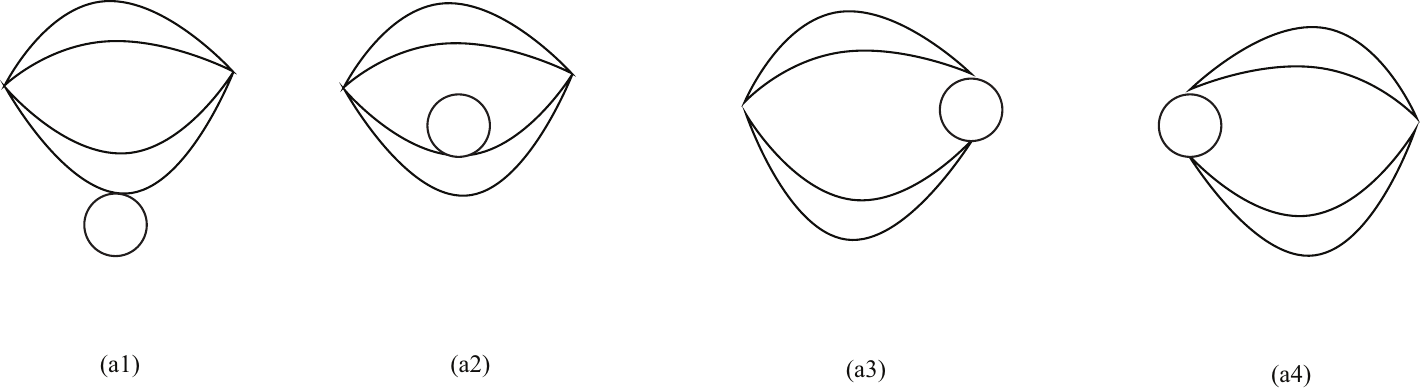}
\caption{$A_{\infty}$ relation I.}
\label{Figure30new}
\end{figure}
\begin{figure}[h]
\centering
\includegraphics[scale=0.5]{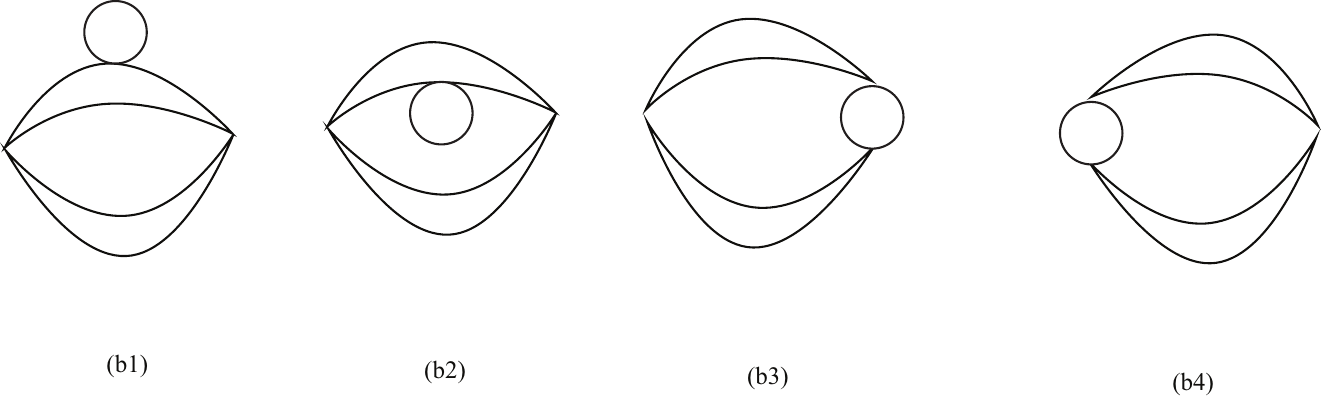}
\caption{$A_{\infty}$ relation II.}
\label{Figure31new}
\end{figure}
In other words the sum of the terms corresponding to Figures \ref{Figure30new} and \ref{Figure31new}
are zero.  The terms corresponding to (a3) and (a4) cancel with  terms corresponding to (b3) and (b4).
On the other hand, 
(a1), (a2), (b1), (b2) coincide with (1), (3), (2), (4), respectively.
Thus the equality (\ref{lem205form}) follows modulo sign.
Using the fact the definition of sign for $Z$ is by the Koszul rule, we can check the sign.
For example we can check the cancelation with sign of (a3) and (b3) as follows.
Both are super trace of
$$
f^2 \mapsto \frak m_*(y_* \dots \frak m_*(y_* \dots \frak m_*(y_* \dots y_*, f^2, x_* \dots x_*) \dots x_*) \dots x_*)
$$
modulo sign.  In (a3) it appears as a part of 
$$
\frak m_*(y_* \dots [d,\frak m_*](y_* \dots y_*, f^2, x_* \dots x_*) \dots x_*).
$$
In (b3) it appears as a part of 
$$
[d,\frak m_*](y_* \dots \frak m_*(y_* \dots y_*, f^2, x_* \dots x_*) \dots x_*).
$$
Therefore by the Koszul rule they appear with opposite sign and cancel.\footnote{To work out 
the detail of the sign part of the above proof it seems better to use 
`super trace formulation of the sign' appearing in (\ref{eq:formula_Z-pairing2}) etc.}
\end{proof}

\subsection{The main theorem of section \ref{sec:annuli}.}
\label{sec4maintheorem}

\begin{thm}\label{ZisPD}{\rm (Cardy relation)}\index{Cardy relation}
Under the Assumption $\ref{LUtrans}$ we have 
\begin{equation}\label{form2015}
 \langle
\widehat{\frak p}(\text{\bf x})
,
\widehat{\frak p}(\text{\bf y})
\rangle_{\text{\rm PD}_X} 
= (-1)^{n(n-1)/2}Z(\text{\bf x},\text{\bf y}), 
\end{equation}
for 
$[\text{\bf x}] \in HH_*(\cL;\Lambda_0)$,
$[\text{\bf y}] \in HH_*(\sU;\Lambda_0)$.
\end{thm}

\begin{rem}
$  $ \par
\begin{enumerate}
\item
The conclusion of Theorem \ref{ZisPD}
holds without Assumption $\ref{LUtrans}$.
We can prove it 
by perturbing $U_{\upsilon}$ so that Assumption \ref{LUtrans}
is satisfied.
\item
In case $\mathscr L = \mathscr U =\{L(u)\}$ and 
$\text{\bf x} = \text{\bf y} = PD([{\rm pt}]) \in HF((L(u),b),(L(u),b))$
Theorem \ref{ZisPD} coincides with 
\cite[Proposition 3.5.2]{toric3}.
In this case Definition \ref{defnZ} becomes 
\cite[Definition 1.3.22]{toric3}.
\end{enumerate}
\end{rem}

The proof of Theorem \ref{ZisPD} is based on the study of moduli space 
of pseudomorphic annuli and occupies the rest of this section.
The proof is based on the  
annulus argument  which appeared
previously in \cite{abouzaid:IHES,bircor,toric3} in the literature of 
symplectic topology.
The Cardy relation appeared earlier in the physics literature  such as in \cite{cardy}.

\subsection{Moduli space of annuli: compactification and Kuranishi structure}
\label{sec:moduli-space-annuli}
Let $\vec{\kappa} = (\kappa_0,\kappa_1,\dots,\kappa_K)$ and 
$\kappa_i \in \uwave{K}$ be as in Definition  \ref{def211}, so that we have $\kappa_i \neq \kappa_{i+1}$. 
Let $\vec k = (k_0,\dots,k_K)$ with $k_0,\dots,k_K \in \Z_{\ge 0}$.
We denote the set of such pairs $(\vec{\kappa},\vec k)$ by $\text{\rm seq}'_K(\cL)$.\index[syindex]{seqKcL@$\text{\rm seq}'_K(\cL)$}
We define a similar set $\text{\rm seq}'_{R}(\sU)$ of the pair $(\vec{{\upsilon}},\vec r)$
with
$\vec{{\upsilon}} = ({\upsilon}_0,{\upsilon}_1,\dots,{\upsilon}_{R})$, $\vec r = (r_0,\dots,r_R)$
in the same way.
\par
%In case that $L_{\kappa_{\blue{i-1}}} \ne L_{\kappa_{\blue{i}}}$ % the element $x_i$ 
% is identified with
We fix an intersection point $L_{\kappa_{\blue{i-1}}} \cap L_{\kappa_{\blue{i}}}$ (or self-intersection point of 
$L_{\kappa_{\blue{i-1}}} = L_{\kappa_{\blue{i}}}$),
which we write $p_i$.
Similarly %$y_i$ is identified with
we fix $q_i \in U_{\rho_{\blue{i-1}}} \cap U_{\rho_{\blue{i}}}$ (or self-intersection point). %whenever these Lagrangians are different.
We denote by $\vec p$ (resp. $\vec q$), the totality of $p$'s (resp. $q$'s).
\par
We put $\vert\vec k\vert  = \sum_{i=0}^K k_i$, 
$\vert\vec r\vert  = \sum_{i=0}^R r_i$.
\par
We now define the moduli space of the pseudo-holomorphic annuli which we use 
for the proof of Theorem \ref{ZisPD}.
\par
Consider the bordered semi-stable curve $\Sigma$ of genus {$0$}  with two
boundary components $\partial_1 \Sigma$,
$\partial_2 \Sigma$. We equip this curve with $K+\vert\vec k\vert +1$ and $R+\vert\vec r\vert+1$ boundary marked points, denoted
$\vec z^{\,1}$ and $\vec z^{\,2}$ on
the two components $\partial_1 \Sigma$ and
$\partial_2 \Sigma$ of the boundary, as well as $\ell$ interior marked points
$\vec z^{\,\,+} = (z^+_1,
\dots,z^+_{\ell})$.
We denote the totality of these data by $(\Sigma;\vec z^{\,1},\vec z^{\,2},\vec z^{\,+})$
or sometimes by $\Sigma$ for short.
When we regard $\Sigma$ as an annulus in $\C$, 
we take the counter-clockwise cyclic order for the outer circle and the
clockwise cyclic order for the inner circle (both agree with the natural orientation of the boundary of the annulus).
\par
Let $(\Sigma;\vec z^{\,1},\vec z^{\,2},\vec z^{\,+})$  be as above.
Let
\begin{equation}\label{formnew216}
\aligned
\vec z^{\,1} &= (z^1_0,w^1_{\blue{0},1},\dots,w^1_{\blue{0} ,k_1},
z^1_1,w^1_{\blue{1} ,1},\dots,w^1_{\blue{1} ,k_2},z^1_2,\dots,z^1_K,
w^1_{\blue{K} ,1},\dots,w^1_{\blue{K},k_K}), \\
\vec z^{\,2} &= (z^2_0,w^2_{\blue{0},1},\dots,w^2_{\blue{0},r_1},
z^2_1,w^2_{\blue{1},1},\dots,w^2_{\blue{1} ,r_2},z^2_2,\dots,z^2_R,
w^2_{\blue{R} ,1},\dots,w^2_{\blue{R},r_R}).
\endaligned
\end{equation}
\noindent
We
consider a holomorphic map
$u: \Sigma  \to X$ such that:
\begin{enumerate}
\item
$
u(\overline{z^1_{{i}} z^1_{{i+1}}}) \subset L_{\kappa_i}$, where we put $\kappa_{K+1}
= \kappa_0$.
\item
$
u(\overline{z^2_{{i}} z^2_{{i+1}}}) \subset U_{\upsilon_i}$, where we put $\upsilon_{R+1}
= \upsilon_0$.
\item
$u(z^1_i) = p_i$ and $u(z^2_i) = q_i$.
\item 
In case $p_i$, $q_i$ are self-intersection points we require a similar switching condition as  
Definition \ref{def211} (4).
\end{enumerate}
\begin{figure}[h]
\centering
\includegraphics[scale=0.7]{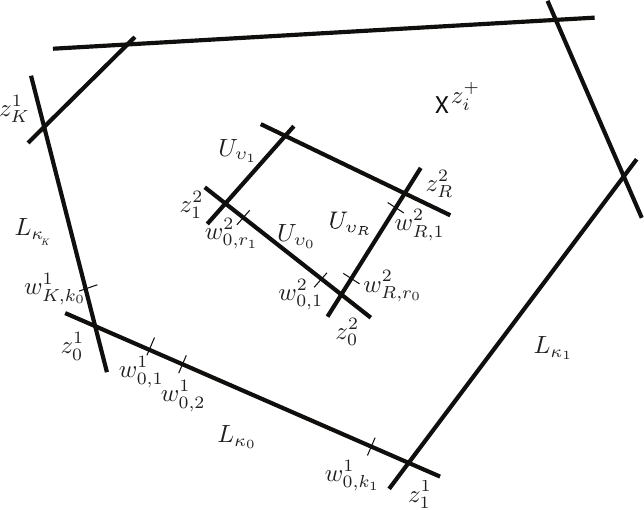}
\caption{$(\Sigma;\vec z^{\,1},\vec z^{\,2},\vec z^{\,+})$}
\label{Figure10-1}
\end{figure}
We say that $(u;\Sigma;\vec z_1,\vec z_2,\vec z^{\,\,+})$
is a {\it stable map} if its automorphism group is of finite order. 
We can define the notion of isomorphism between $(u;\Sigma;\vec z_1,\vec z_2,\vec z^{\,\,+})$ 
and $(u';\Sigma';\vec z'_1,\vec z'_2,(\vec z^{\,\,+})')$
in an obvious way.
\par
Let $\pi_2(X;\vec \kappa,\vec {\upsilon},\vec p,\vec q)$ be the set of 
homotopy classes of the maps satisfying the conditions 
(1), (2), (3) and (4) above.
Using the Maslov index (in other words the virtual dimension of 
the moduli space) and the energy we define an equivalence relation 
on $\pi_2(X;\vec \kappa,\vec {\upsilon},\vec p,\vec q)$,  and 
define $\Pi_2(X;\vec \kappa,\vec {\upsilon},\vec p,\vec q)$ to be the\index[syindex]{Pi2XkappaUpsilon@$\Pi_2(X;\vec \kappa,\vec {\upsilon},\vec p,\vec q)$} 
quotient by this equivalence relation, so that two homotopy classes are identified if they have the same energy and Maslov index.
We take forgetability data ${}_{\mathcal L}\vec{\frak f} = ({}_{\mathcal L}{\frak f}_{\bf e},{}_{\mathcal L}{\frak f}_{\bf f})$, 
${}_{\mathcal U}\vec{\frak f} = ({}_{\mathcal U}{\frak f}_{\bf e},{}_{\mathcal U}{\frak f}_{\bf f})$.

\begin{defn}
We put $B \in \Pi_2(X;\vec \kappa,\vec {\upsilon},\vec p,\vec q)$.
We denote by\index[syindex]{Mkappakann@$\mathcal M_{(\vec\kappa,\vec k),(\vec{\upsilon},\vec r);\ell}^{\text{ann}}
(\vec p,\vec q;B)$}
$$
\mathcal M_{(\vec\kappa,\vec k),(\vec{\upsilon},\vec r);\ell}^{\text{ann}}
(\vec p,\vec q;B)
$$
the set of all isomorphism classes of stable maps $(u;\Sigma;\vec z_1,
\vec z_2,\vec z^{\,\,+})$ in homology class $B$
that
satisfy Conditions $(1)$, $(2)$, and $(3)$ above.
\par
When we include ${}_{\mathcal L}\vec{\frak f}$ and ${}_{\mathcal U}\vec{\frak f}$ we write it as
$
\mathcal M_{(\vec\kappa,\vec k),(\vec{\upsilon},\vec r);\ell}^{\text{ann}}
(B;\vec p,\vec q;{}_{\mathcal L}\vec{\frak f},{}_{\mathcal U}\vec{\frak f}).
$
\end{defn}
We define the evaluation map
\begin{equation} \label{eq:evaluation_map_annuli}
\aligned
\text{\rm ev} : 
&\mathcal M_{(\vec\kappa,\vec k),(\vec{\upsilon},\vec r);\ell}^{\text{ann}}
(\vec p,\vec q;B)\\
&\to 
X^{\ell}
\times \prod_{i=0}^K ((\tilde L_{\kappa_{i-1}} \times_X \tilde L_{\kappa_{i}})\setminus {\rm Diagonal})
\times \prod_{i=0}^R ((\tilde U_{{\upsilon}_{i-1}}  \times_X \tilde U_{{\upsilon}_{i}})\setminus {\rm Diagonal})\\
& \qquad \times \prod_{i=0}^K \tilde L_{\kappa_i}^{k_i}
\times \prod_{i=0}^R \tilde U_{{\upsilon}_i}^{r_i}
\endaligned
\end{equation}
whose components are given as follows:
$$
\text{\rm ev} = (\text{\rm ev}^+;\text{\rm ev}^{\partial,1,1},\text{\rm ev}^{\partial,1,2};
\text{\rm ev}^{\partial,2,1},\text{\rm ev}^{\partial,2,2}).
$$
Here the map
$$
\text{\rm ev}^+ = (\text{\rm ev}^+_1,\dots,\text{\rm ev}^+_{\ell}) : \mathcal M_{(\vec\kappa,\vec k),(\vec{\upsilon},\vec r);\ell}^{\text{ann}}
(B;\vec p,\vec q) \to X^{\ell}
$$
is defined by
$$
\text{\rm ev}^+_i(u;\Sigma;\vec z_1,
\vec z_2,\vec z^{\,\,+}) = u(z^+_i),
$$
the map 
$$
\text{\rm ev}^{\partial,1,1} = (\text{\rm ev}^{\partial,1,1}_{0},\dots,\text{\rm ev}^{\partial,1,1}_{K}) : \mathcal M_{(\vec\kappa,\vec k),(\vec{\upsilon},\vec r);\ell}^{\text{ann}}
(\vec p,\vec q;B) \to \prod_{i=0}^K ((\tilde L_{\kappa_{i-1}} \times_X \tilde L_{\kappa_{i}})\setminus {\rm Diagonal})
$$
by %\marginpar{Error corrected.  KF 2025 Aug 27.}
$$
\text{\rm ev}^{\partial,1,1}_{i}(u;\Sigma;\vec z_1,
\vec z_2,\vec z^{\,\,+}) = u(z^1_i) = p_i,
$$
and the map 
$$
\text{\rm ev}^{\partial,2,1} = (\text{\rm ev}^{\partial,2,1}_{0},\dots,\text{\rm ev}^{\partial,2,1}_{R}) : \mathcal M_{(\vec\kappa,\vec k),(\vec{\upsilon},\vec r);\ell}^{\text{ann}}
(\vec p,\vec q;B) \to \prod_{i=0}^R ((\tilde U_{{\upsilon}_{i-1}} \times_X \tilde U_{{\upsilon}_{i}})\setminus {\rm Diagonal})
$$
by
$$
\text{\rm ev}^{\partial,2,1}_{i}(u;\Sigma;\vec z_1,
\vec z_2,\vec z^{\,\,+}) = u(z^2_i) = q_i.
$$
Namely $\text{\rm ev}^{\partial,1,1}$ and $\text{\rm ev}^{\partial,2,1}$ are determined by $\vec p$ and $\vec q$.

Moreover we define 
$$
\text{\rm ev}^{\partial,1,2} = (\vec{\text{\rm ev}}^{\partial,1,2}_{0},\dots,\vec{\text{\rm ev}}^{\partial,1,2}_{K}) : \mathcal M_{(\vec\kappa,\vec k),(\vec{\upsilon},\vec r);\ell}^{\text{ann}}
(\vec p,\vec q;B) \to \prod_{i=0}^K L_{{\kappa}_i}^{k_i}
$$
by
$$
\vec{\text{\rm ev}}^{\partial,1,2}(u;\Sigma;\vec z_1,
\vec z_2,\vec z^{\,\,+})
= (u(w^1_{i,1}),\dots,u(w^1_{i,k_i}))_{i=0}^K
$$
and 
$$
\vec{\text{\rm ev}}^{\partial,2,2} = (\vec{\text{\rm ev}}^{\partial,2,2}_{0},\dots,\vec{\text{\rm ev}}^{\partial,2,2}_{R}) : \mathcal M_{(\vec\kappa,\vec k),(\vec{\upsilon},\vec r);\ell}^{\text{ann}}
(\vec p,\vec q;B) \to \prod_{i=0}^R U_{{\upsilon}_i}^{r_i}
$$
by
$$
\vec{\text{\rm ev}}^{\partial,2,2}(u;\Sigma;\vec z_1,
\vec z_2,\vec z^{\,\,+})
= (u(w^2_{i,1}),\dots,u(w^2_{i,r_i}))_{i=0}^R.
$$

Now we change the notation as in Subsection \ref{constcyclic}.
 %\marginpar{\bf Corrected.
%Several places in this and the next pages. KF 2025 Aug. 27}
\par
We first introduce, for $m'_1 \in \{0,\dots,k'_0\}$,  $m'_2 \in \{0,\dots,r'_0\}$, 
the moduli spaces
$$
\mathcal M_{(\vec\kappa',\vec k',m'_1),(\vec{\upsilon}',\vec r',m'_2);\ell}^{\text{ann}}
(\vec p,\vec q;B) :=
\mathcal M_{(\vec\kappa',\vec k'),(\vec{\upsilon}',\vec r');\ell}^{\text{ann}}
(\vec p,\vec q;B),
$$
which specify $0$-th marked points $z_1(0) \in \partial_1\Sigma$,\index[syindex]{Mkappakannm1m2@$\mathcal M_{(\vec\kappa,\vec k,m_1),(\vec{\upsilon},\vec r,m_2);\ell}^{\text{ann}}
(\vec p,\vec q;B)$} 
$z_2(0) \in \partial_1\Sigma$ in the same way as Definition \ref{defn917}.\index[syindex]{zz01@$z_1(0)$}
\index[syindex]{zz02@$z_2(0)$}
It then determines the enumeration of the marked points $\vec z_{\,1}$ and of $\vec z_{\,2}$.
(It may be different from (\ref{formnew216}) since $z_1(0)$ may be $w^1_{0,i}$.)

Next, let $(\vec{\kappa},\vec p)\in \widetilde{\text{\rm Seq}}_K$, 
$(\vec{\upsilon},\vec q) \in \widetilde{\text{\rm Seq}}_R$. We put\index[syindex]{Mkappakannell@$\mathcal M_{\ell}^{\text{ann}}
(\vec\kappa,\vec{\upsilon};\vec p,\vec q;B)$}
\begin{equation}\label{newnew2018}
\mathcal M_{\ell}^{\text{ann}}
(\vec\kappa,\vec{\upsilon};\vec p,\vec q;B)
= 
\mathcal M_{(\vec\kappa',\vec k',m'_1),(\vec{\upsilon}',\vec r',m'_2);\ell}^{\text{ann}}
(\vec p,\vec q;B)
\end{equation}
where ${\rm Red}(\vec{\kappa},\vec p) = (\vec{\kappa}', \vec{p}, \vec k', m'_1)$
and ${\rm Red}(\vec{\upsilon},\vec q) = (\vec{\upsilon}', \vec{q}, \vec r', m'_2)$.
See (\ref{eq:change_notation_moduli_space}).

We next consider the forgetful map 
\begin{equation}\label{forgettoannul}
\aligned
&\frak{forget}_{\rm sc} : 
\mathcal M_{(\vec\kappa',\vec k',m'_1),(\vec{\upsilon}',\vec r',m'_2);\ell}^{\text{ann}}
(\vec p,\vec q;B)
\to 
\mathcal M_{(1,1);0},
\\
&\frak{forget}_{\rm sc} : 
\mathcal M_{\ell}^{\text{ann}}
(\vec\kappa,\vec{\upsilon};\vec p,\vec q;B)
\to 
\mathcal M_{(1,1);0}.
\endaligned
\end{equation}
Here $\mathcal M_{(1,1);0}$ is\index[syindex]{Ma110@$\mathcal M_{(1,1);0}$} 
the moduli space of bordered Riemann surfaces of genus 
zero with two boundary components and 
one marked point in each component of the boundary.
\par
The map (\ref{forgettoannul}) is defined by\index[syindex]{forgetsource@$\frak{forget}_{\rm sc}$}
$$
\frak{forget}_{\rm sc}(u;\Sigma;\vec z_1,
\vec z_2,\vec z^{\,\,+})
= (\overline{\Sigma},z_1(0),z_2(0)),
$$
when the space $\overline{\Sigma}$ is defined as follows.
\par
We consider the bordered Riemann surface $\Sigma$
together with two marked points $z_1(0),z_2(0)$ 
on the boundary. 
If $(\Sigma;z_1(0),z_2(0))$ is stable then 
$\overline{\Sigma} = \Sigma$.
If not we shrink the unstable components of 
$\Sigma$ to obtain $\overline{\Sigma}$.
We then get 
$(\overline{\Sigma},z_1(0),z_2(0)) 
\in \mathcal M_{(1,1);0}$.
\par
It is well-known (and is proved in \cite[Section C.4]{abouzaid:IHES} and 
\cite[Lemma 3.4.5]{toric3})
that $\mathcal M_{(1,1);0}$ is a 2 disk. It has a stratification 
by the combinatorial type of its elements $({\Sigma},z_1(0),z_2(0))$.
The $0$-dimensional stratum consists of  two particular points which we denote by $(\Sigma_1;z^1_1,z_2^1)$ and 
$(\Sigma_2;z_1^2,z^2_2)$, and which we now describe:\index[syindex]{Sigma1@$\Sigma_1$}
\par
\begin{figure}[h]
\centering
\includegraphics[scale=0.7]{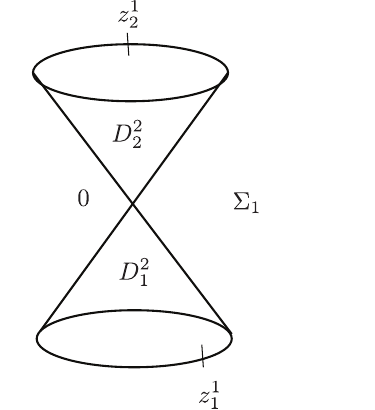}
\caption{$({\Sigma}_1,z^1_1,z^1_2)$}
\label{Figure10-2}
\end{figure}

({\bf Description of $\Sigma_1$}:) Let $D^2_1$ and $D^2_2$ be two copies of the disk $D^2 = \{z \in \C \mid
\vert z\vert \le 1\}$. We identify $0 \in D^2_1$ and $0 \in D^2_2$
to obtain $\Sigma_1$. We put
$\partial_1\Sigma_1 = \partial D^2_1$, $\partial_2\Sigma_1 = \partial D^2_2$,
$z^1_1 = 1 \in \partial D^2_1$, and $z_2^1 = 1 \in \partial D^2_2$.
Thus we have $(\Sigma_1;z^1_1,z^1_2) \in \mathcal M_{(1,1);0}^{\text{ann}}$.
See Figure \ref{Figure10-2}.
\par
\begin{figure}[h]
\centering
\includegraphics[scale=0.5]{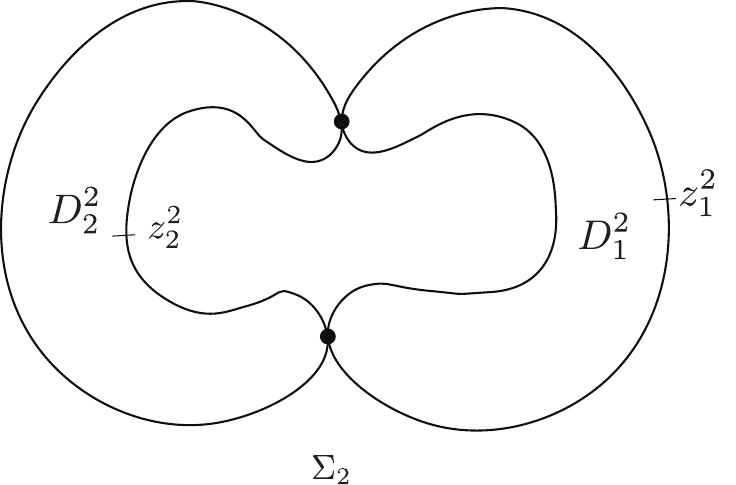}
\caption{$({\Sigma}_2,z^2_1,z^2_2)$}
\label{Figure10-3}
\end{figure}

({\bf Description of $\Sigma_2$}:)  We consider two copies\index[syindex]{Sigma2@$\Sigma_2$}
$D^2_1$, $D^2_2$ of the disk $D^2$.
We put $\frak z^i_k = e^{2k\pi\sqrt{-1}/3} \in \partial D^2_i$ for $i=1,2$, and
$k=0,1,2$. We identify $\frak z^1_1$ with $\frak z^2_1$ and
$\frak z^1_2$ with $\frak z^2_2$ respectively and obtain
$\Sigma_2$.
We put $\frak z^i_0 = \blue{z_i^2}$ and $i=1,2$. We call
$\partial_i \Sigma_2$ the component which contains ${z_i^2}$.
We have thus defined
$(\Sigma_2;z_1^2,z^2_2) \in \mathcal M_{(1,1);0}^{\text{\rm ann}}$.
See Figure \ref{Figure10-3}.
\par It is easy to see from the construction that
the forgetful map $(\ref{forgettoannul})$ is continuous and is a smooth submersion
on each stratum.
We remark that (\ref{forgettoannul}) may not be smooth across strata. See \cite[pp 777-778]{fooo09}.
\par
Hereafter we write
\begin{equation}
\frak{forget}_{\rm sc}^{-1}(X) \cap
\mathcal M_{\ell}^{\text{ann}}
(\vec\kappa,\vec{\upsilon};\vec p,\vec q;B)
\end{equation}
for the inverse image of a subset $X \subset \mathcal M_{(1,1);0}^{\text{ann}}$
by the map (\ref{forgettoannul}).
\par
We next describe Kuranishi structures and their
CF perturbations of $\mathcal M_{(\vec\kappa,\vec k),(\vec{\upsilon},\vec r);\ell}^{\text{ann}}
(B;\vec p,\vec q)^{\boxplus 1}$, $\mathcal M_{(\vec\kappa,\vec k,m_1),(\vec{\upsilon},\vec r,m_2);\ell}^{\text{ann}}
(B;\vec p,\vec q)^{\boxplus 1}$ and
$\mathcal M_{\ell}^{\text{ann}}
(\vec\kappa,\vec{\upsilon};B)^{\boxplus 1}$.
\begin{lem}\label{bdryMannu}
$  $ \par
\begin{enumerate}
\item
The space $\mathcal M_{(\vec\kappa,\vec k),(\vec{\upsilon},\vec r);\ell}^{\text{\rm ann}}
(\vec p,\vec q;B;{}_{\mathcal L}\vec{\frak f},{}_{\mathcal U}\vec{\frak f})$ is compact and Hausdorff.
\item
The space $[0,1]^{\vert {}_{\mathcal L}\frak f_{\bf f}\vert + \vert {}_{\mathcal U}\frak f_{\bf f}\vert} \times \mathcal M_{(\vec\kappa,\vec k),(\vec{\upsilon},\vec r);\ell}^{\text{\rm ann}}
(\vec p,\vec q;B;{}_{\mathcal L}\vec{\frak f},{}_{\mathcal U}\vec{\frak f})^{\boxplus 1}$ has an oriented Kuranishi structure with corners. 
\item
The boundary of $[0,1]^{\vert {}_{\mathcal L}\frak f_{\bf f}\vert + \vert {}_{\mathcal U}\frak f_{\bf f}\vert} \times \mathcal M_{(\vec\kappa,\vec k),(\vec{\upsilon},\vec r);\ell}^{\text{\rm ann}}
(\vec p,\vec q;B;{}_{\mathcal L}\vec{\frak f},{}_{\mathcal U}\vec{\frak f})^{\boxplus 1}$
is the union of the following $8$ types of fiber products:
\smallskip
\begin{enumerate}
\item
$
\frak{forget}_{\rm sc}^{-1}(\partial \mathcal M_{(1,1);0}^{\text{\rm ann}})
\cap ([0,1]^{\vert {}_{\mathcal L}\frak f_{\bf f}\vert + \vert {}_{\mathcal U}\frak f_{\bf f}\vert} \times \mathcal M_{(\vec\kappa,\vec k),(\vec{\upsilon},\vec r);\ell}^{\text{\rm ann}}
(\vec p,\vec q;B;{}_{\mathcal L}\vec{\frak f},{}_{\mathcal U}\vec{\frak f})^{\boxplus 1})$
\item
\begin{equation}\label{31822}
\aligned
&([0,1]^{\vert {}_{\mathcal L}\frak f^2_{\bf f}\vert} \times {\mathcal M}_{\ell'';k''_i+1}(L_{\kappa_i};\beta';{}_{\mathcal L}\vec{\frak f}^2)^{\boxplus 1})
\\
&{}_{\text{\rm ev}_0} \times_{\text{\rm ev}^{\partial,\blue{1,2}}_{i,j}}
([0,1]^{\vert {}_{\mathcal L}\frak f^1_{\bf f}\vert + \vert {}_{\mathcal U}\frak f_{\bf f}\vert} \times\mathcal M_{(\vec\kappa,\vec k'),(\vec{\upsilon},\vec r);\ell'}^{\text{\rm ann}}
(\vec p,\vec q;B';{}_{\mathcal L}\vec{\frak f}^1,{}_{\mathcal U}\vec{\frak f})^{\boxplus 1})
\endaligned
\end{equation}
The notations are as follows: 
{\rm (i)} the sequence $\vec k'$ is the same as $\vec k$ except for the $i$-th component which is 
$k'_i$ such that $k'_i+k''_i=k_i \blue{+1}$, {\rm (ii)} the integer $j$ lines in
$\underline k'_i$, and {\rm (iii)} the equivalence class
$\beta'\in \Pi_2(X;L_{\kappa_i})$ is such that $\beta$ and $B'$ satisfies
$\beta'\# B' = B$, where $\#$ is the obvious concatenation.
{\rm (iv)}  $ {}_{\mathcal L}\vec{\frak f}^1$ and $ {}_{\mathcal L}\vec{\frak f}^2$ are obtained from $ {}_{\mathcal L}\vec{\frak f}$ in the same way as Subsection  $\ref{subsec:forgetunforget}$.
\par
In $(\ref{31822})$, the union is taken over the data $i$, $k'_i$, $k''_i$, $j$, $\beta'$, and $B$
together with the shuffles $(\mathbb L'',\mathbb L')$ of $\underline\ell$
such that $\# \mathbb L'' = \ell''$, $\# \mathbb L' = \ell'$.
\par
We put the Kuranishi structure of Proposition $\ref{Kuraeistspoly2}$
on the first factor of $(\ref{31822})$ .
We put the Kuranishi structure in $(2)$
on the second factor of $(\ref{31822})$.
See Figures $\ref{Figure4-4-1}$ and $\ref{Figure4-4-2}$.
\begin{center}
\begin{figure}[h]
 \begin{tabular}{cc}
 \begin{minipage}[t]{0.45\hsize}
\centering
\includegraphics[scale=0.7]
{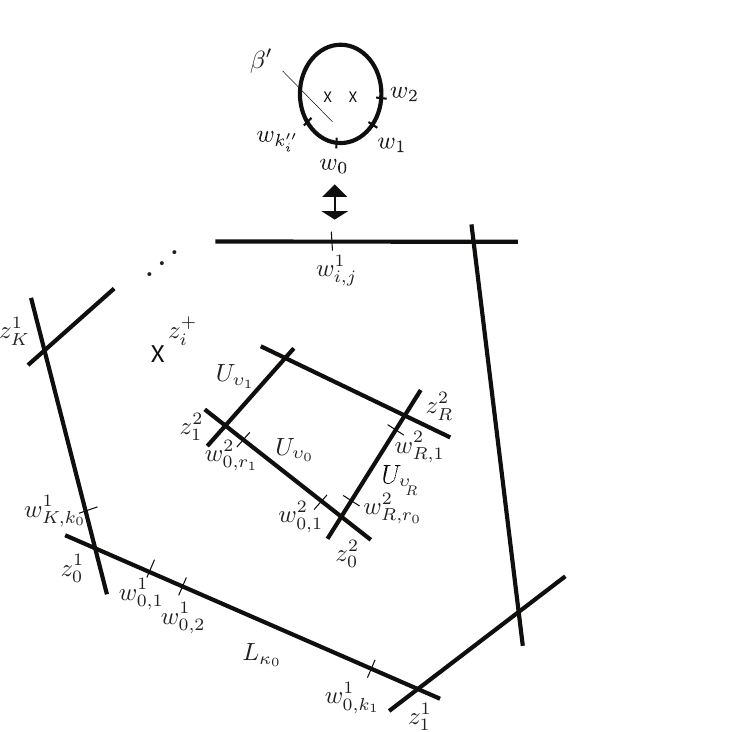}
\caption{An element of (\ref{31822})}
\label{Figure4-4-1}
\end{minipage} &
 \begin{minipage}[t]{0.45\hsize}
\centering
\includegraphics[scale=0.7]
{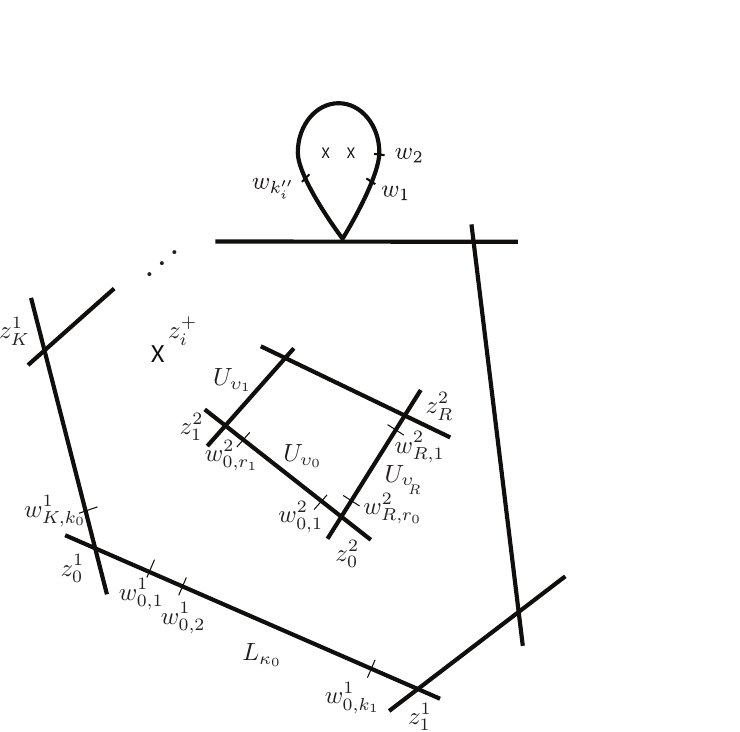}
\caption{An element of (\ref{31822})}
\label{Figure4-4-2}
\end{minipage}
\end{tabular}
\end{figure}
\end{center}
\item
\begin{equation}\label{3182}
\aligned
&([0,1]^{\vert {}_{\mathcal U}\frak f^2_{\bf f}\vert} \times {\mathcal M}_{\ell'';r''_i+1}(U_{\upsilon_i};\beta'';{}_{\mathcal U}\vec{\frak f}^2)^{\boxplus 1}) \\
&\,{}_{\text{\rm ev}_0} \times_{\text{\rm ev}^{\partial,\blue{2,2}}_{i,j}}
([0,1]^{\vert {}_{\mathcal L}\frak f_{\bf f}\vert +\vert {}_{\mathcal U}\frak f^1_{\bf f}\vert} \times\mathcal M_{(\vec\kappa,\vec k),(\vec{\upsilon},\vec r^{\,\prime});\ell'}^{\text{\rm ann}}
(\vec p,\vec q;B';{}_{\mathcal L}\vec{\frak f},{}_{\mathcal U}\vec{\frak f}^1)^{\boxplus 1}).
\endaligned
\end{equation}
The notations are as follows. The sequence
$\vec r'$ is the same as $\vec r$ except for the $i$-th component which is 
$r'_i$ such that $r'_i+r''_i=r_i \blue{+1}$, the integer
$j$ lies in $\underline r'_i$, and the class
%\par
$\beta''\in \Pi_2(X;U_{\upsilon_i})$is  such that $\beta''$ and $B'$ satisfies
$\beta''\# B' = B$, where $\#$ is the obvious concatenation.
$ {}_{\mathcal U}\vec{\frak f}^1$ and $ {}_{\mathcal U}\vec{\frak f}^2$ are obtained from $ {}_{\mathcal U}\vec{\frak f}$ in the same way as Subsection  $\ref{subsec:forgetunforget}$.
\par
In $(\ref{3182})$, the union is taken over the data $i$, $r'_i$, $r''_i$, $j$, $\beta''$, and $B$
together with the shuffles $(\mathbb L'',\mathbb L')$ of $\underline\ell$
such that $\# \mathbb L'' = \ell''$, $\# \mathbb L' = \ell'$.
\par
We put the Kuranishi structure of Proposition $\ref{Kuraeistspoly2}$
on the first factor of $(\ref{3182})$.
We put the Kuranishi structure in $(2)$
on the second factor of $(\ref{3182})$.
See Figure $\ref{Figure4-4-3}$.
\begin{figure}[h]
\centering
\includegraphics[scale=0.7]{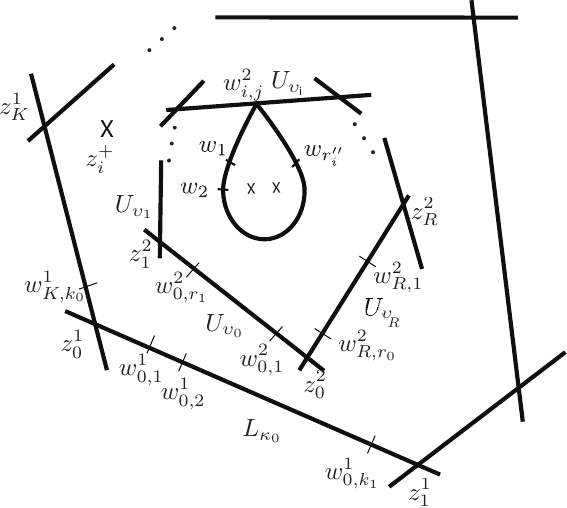}
\caption{An element of (\ref{3182})}
\label{Figure4-4-3}
\end{figure}
\item
Let $\vec{\kappa}^{\,\prime}$, $\vec{\kappa}^{\,\prime\prime}$, 
$\vec p^{\,\prime}$, and $\vec p^{\,\prime\prime}$,  satisfy 
$(\vec{\kappa}^{\,\prime\prime},\vec p^{\,\prime\prime}) \#_{i,i+1}   (\vec{\kappa}^{\,\prime},\vec p^{\,\prime})
=(\vec{\kappa},\vec p)$ in the sense of $(\ref{213formula})$.
Let $B'' \in \Pi_2(\vec{\kappa}'',\vec k'')$ and
$B' \in \Pi_2(X;\vec \kappa',\vec {\upsilon},\vec p^{\,\prime},\vec q)$ satisfy
$B'' \#_{i,i+1} B' = B$.
Let $\vec k' = (k'_1,\dots,k'_{K'})$ and
$\vec k'' = (k''_1,\dots,k''_{K''})$.
We set $\vec k = \vec k'' \#_{i,i+1} \vec k'$. Namely 
\blue{
$$
k_j = 
\begin{cases}
k'_j   & j=0 ,\dots,i-1, \\
k'_i  + k''_{0}, &j=i \\
k''_{j-i}  & j=i+1,\dots, i+K'' - 1, \\
k'_{i+1} + k''_{K''}  &j=\blue{i+K''}, \\
k'_{j-K''+1} & j=i+K'' +1, \dots, K'+K''-1.
\end{cases}
$$
}
%$$
%k_j = 
%\begin{cases}
%k'_j   & j=1,\dots,i, \\
%k''_{j-i}  & j=i+1,\dots, i+K'', \\
%k'_{j-K''} & j=i+K'',\dots, K'+K''.
%\end{cases}
%$$

We then consider the direct product:
\begin{equation}\label{3219}
\aligned
&([0,1]^{\vert {}_{\mathcal L}\frak f^2_{\bf f}\vert} \times {\mathcal M}_{\ell'';\vec k''}((\vec{\kappa}^{\,\prime\prime},\vec p^{\,\prime\prime});B'';{}_{\mathcal L}\vec{\frak f}^2)^{\boxplus 1}) \\
&\times
([0,1]^{\vert {}_{\mathcal L}\frak f^1_{\bf f}\vert +\vert {}_{\mathcal U}\frak f_{\bf f}\vert} \times\mathcal M_{(\vec\kappa^{\,\prime},\vec k),(\vec{\upsilon},\vec r);\ell'}^{\text{\rm ann}}
(\vec p^{\,\prime},\vec q;B';{}_{\mathcal L}\vec{\frak f}^1,{}_{\mathcal U}\vec{\frak f})^{\boxplus 1})
\endaligned
\end{equation}
We take the sum over all the above data 
together with the shuffles $(\mathbb L'',\mathbb L')$ of $\underline\ell$
such that $\# \mathbb L'' = \ell''$, $\# \mathbb L' = \ell'$.
\par
We put the Kuranishi structure 
of Proposition $\ref{Kuraeistspoly}$
on the first factor of $(\ref{3219})$.
We put the Kuranishi structure \blue{in $(2)$} 
on the second factor of $(\ref{3219})$.
See Figure $\ref{Figure4-4-6}$.
\begin{figure}[h]
\centering
\includegraphics[scale=0.7]{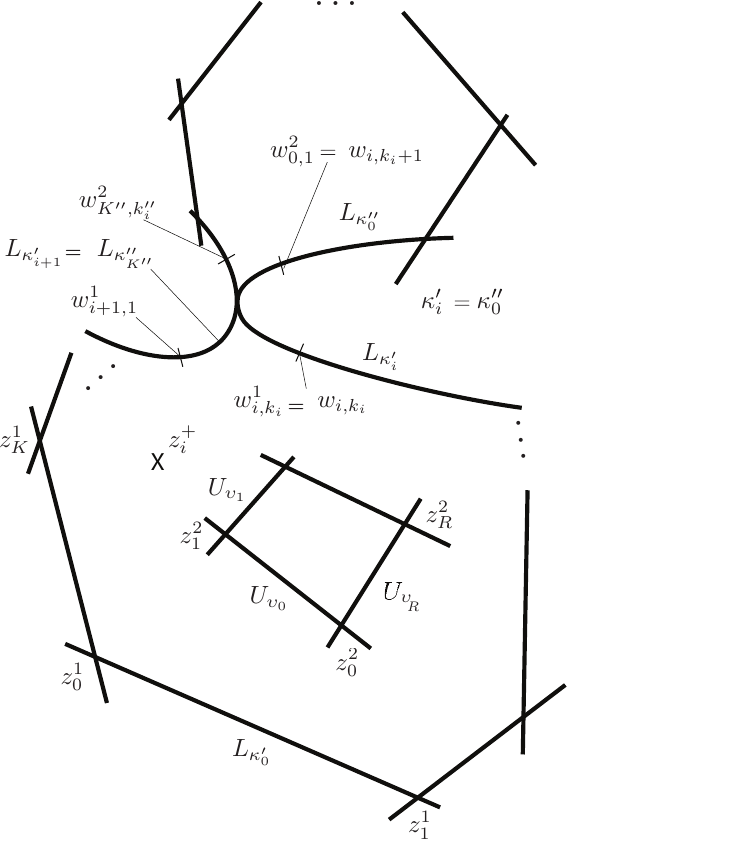}
\caption{An element of (\ref{3219})}
\label{Figure4-4-6}
\end{figure}
\item
The direct product:
\begin{equation}\label{32192}
\aligned
&([0,1]^{\vert {}_{\mathcal U}\frak f^2_{\bf f}\vert} \times 
{\mathcal M}_{\ell'';\vec r^{\,\prime\prime}}((\vec{\upsilon}'',\vec q^{\,\prime\prime});B'';{}_{\mathcal U}\vec{\frak f}^2)^{\boxplus 1})\\
&\times
([0,1]^{\vert {}_{\mathcal L}\frak f_{\bf f}\vert +\vert {}_{\mathcal U}\frak f^1_{\bf f}\vert} \times\mathcal M_{(\vec\kappa,\vec k),(\vec{\upsilon}',\vec r');\ell'}^{\text{\rm ann}}
(\vec p,\vec q^{\,\prime};B';{}_{\mathcal L}\vec{\frak f},{}_{\mathcal U}\vec{\frak f}^1)^{\boxplus 1}).
\endaligned
\end{equation}
Here the notation etc. are similar to $(\ref{3219})$ but 
we replace $L_{\kappa_i}$ etc. by $U_{\upsilon_i}$ etc.
See Figure $\ref{Figure4-8}$. 
\begin{figure}[h]
\centering
\includegraphics[scale=0.7]{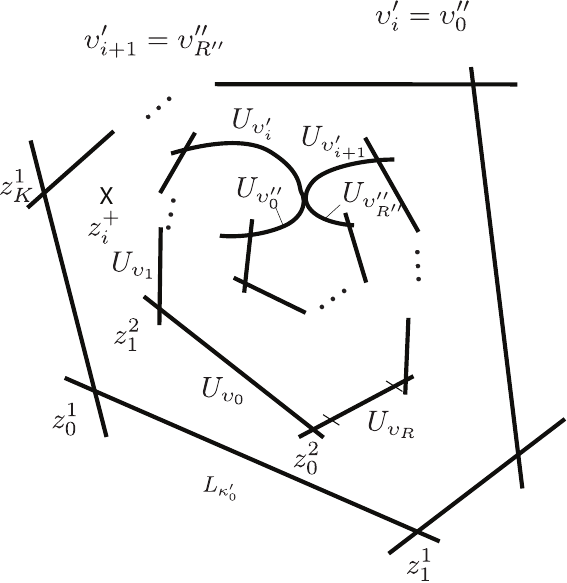}
\caption{An element of (\ref{32192})}
\label{Figure4-8}
\end{figure}

\item
We write $(\vec{\kappa},\vec p)$ for $(\vec{\kappa}'',\vec p'')\#_{i} (\vec{\kappa}',\vec p') $ as in $(\ref{form215})$.
Consider $B' \in \Pi_2(X;\vec \kappa',\vec {\upsilon},\vec p^{\,\prime},\vec q)$ and
$B'' \in \Pi_2(\vec{\kappa}'',\vec k'')$ such that 
$B'' \#_i B' = B$
(Here $ \#_i $ is the obvious concatenation).
Let $\vec k' = (k'_1,\dots,k'_{K'})$,  $\vec k'' = (k''_1,\dots,k''_{K''})$, and consider $m' \in \{ 1,\dots, k'_j\}$, and
$m'' \in \{1,\dots,k''_{1}\}$.
We set $\vec k = \vec k'' \#_{i,(m',m'')} \vec k'$ by 
\blue{
$$
k_j = 
\begin{cases}
k'_j   & j=0,\dots, i-1, \\
m' +k''_{0} - m'' -1 &j=i \\
k''_{j-i}  & j=i+1,\dots, i+K'', \\
m'' + k'_{i} - m' -\blue{1} & j= i+K'' +1, \\
k'_{j-K''-1} & j=i+K''+2,\dots, K'+K''+1.
\end{cases}
$$
}

We now consider the fiber product  over $L_{\kappa_i}$
\begin{equation}\label{2202}
\aligned
&([0,1]^{\vert {}_{\mathcal L}\frak f^2_{\bf f}\vert} \times {\mathcal M}_{\ell'';\vec k^{\,\prime\prime}}((\vec{\kappa}^{\,\prime\prime},\vec p^{\,\prime\prime});B''
;{}_{\mathcal L}\vec{\frak f}^2)^{\boxplus 1}) \\
&\,{}_{\text{\rm ev}_{1,m''}}\times_{\text{\rm ev}^{\partial,2,1}_{i,m'}}
([0,1]^{\vert {}_{\mathcal L}\frak f^1_{\bf f}\vert + \vert {}_{\mathcal U}\frak f_{\bf f}\vert} \times\mathcal M_{(\vec\kappa^{\,\prime},\vec k^{\,\prime}),(\vec{\upsilon},\vec r);\ell'}^{\text{\rm ann}}
(\vec p^{\,\prime},\vec q;B';{}_{\mathcal L}\vec{\frak f}^1,{}_{\mathcal U}\vec{\frak f})^{\boxplus 1}).
\endaligned
\end{equation}
%where the fiber product is.
\par
We put the Kuranishi structure 
of Proposition $\ref{Kuraeistspoly}$ 
%$\ref{Kuraeistspoly2}$
on the first factor of $(\ref{2202})$.
We put the Kuranishi structure 
Proposition $\ref{bdryMannu}$ $(2)$
on the second factor of $(\ref{2202})$.
See Figure $\ref{Figure4-4-5}$.
\begin{figure}[h]
\centering
\includegraphics[scale=0.7]{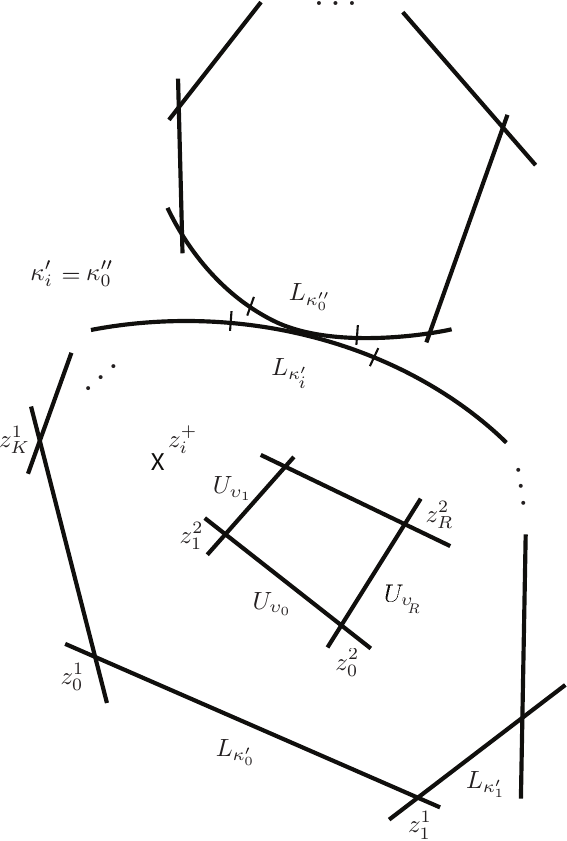}
\caption{An element of (\ref{2202})}
\label{Figure4-4-5}
\end{figure}

\item
We consider the fiber product
\begin{equation}\label{22022}
\aligned
&([0,1]^{\vert {}_{\mathcal U}\frak f^1_{\bf f}\vert} \times {\mathcal M}_{\ell'';\vec r^{\,\prime\prime}}((\vec{\upsilon}^{\,\prime\prime},\vec q^{\,\prime\prime});B'';{}_{\mathcal U}\vec{\frak f}^2)^{\boxplus 1})\\
&\,{}_{\text{\rm ev}_{1,m''}}\times_{\text{\rm ev}^{\partial,2,1}_{i,m'}}
\mathcal M_{(\vec\kappa,\vec k),(\vec{\upsilon}^{\,\prime},\vec r^{\,\prime});\ell'}^{\text{\rm ann}}
(\vec p,\vec q^{\,\prime};B';{}_{\mathcal L}\vec{\frak f},{}_{\mathcal U}\vec{\frak f}^1)^{\boxplus 1}),
\endaligned
\end{equation}
over $U_{\upsilon_i}$. 
Here the notation etc. are similar to $(\ref{3219})$ but 
we replace $L_{\kappa_i}$ etc. by $U_{\upsilon_i}$ etc.
See Figure $\ref{Figure4-9}$.
\item
$\partial([0,1]^{\vert {}_{\mathcal L}\frak f_{\bf f}\vert + \vert {}_{\mathcal U}\frak f_{\bf f}\vert}) \times \mathcal M_{(\vec\kappa,\vec k),(\vec{\upsilon},\vec r);\ell}^{\text{\rm ann}}
(\vec p,\vec q;B;{}_{\mathcal L}\vec{\frak f},{}_{\mathcal U}\vec{\frak f})^{\boxplus 1}$.
\end{enumerate}
\begin{figure}[h]
\centering
\includegraphics[scale=0.7]{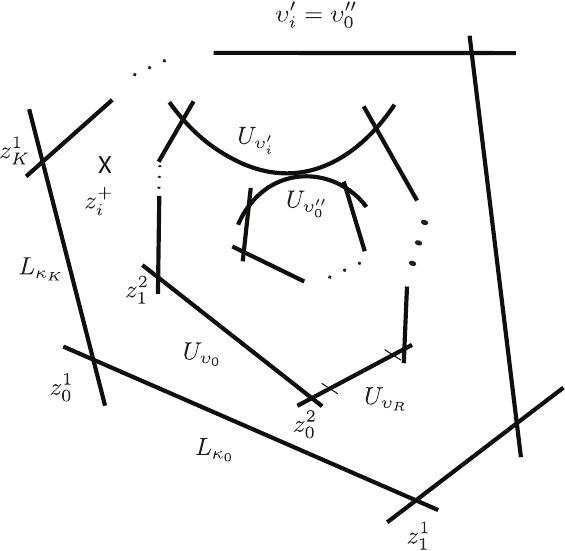}
\caption{An element of (\ref{22022})}
\label{Figure4-9}
\end{figure}

\item
We can define orientations of our Kuranishi structures 
so that they are compatible with the identification 
$(3)$ of the boundaries.
\item
The evaluation maps are compatible 
with the description of the boundaries in $(3)$.
\item
The Kuranishi structure is compatible with the map 
forgetting the marked points labelled by  ${}_{\mathcal L}\vec{\frak f}_{\bf e}$,
${}_{\mathcal U}\vec{\frak f}_{\bf e}$.
\item 
At the point where a coordinate of $[0,1]^{\vert {}_{\mathcal L}\frak f_{\bf f}\vert + \vert {}_{\mathcal U}\frak f_{\bf f}\vert}$ 
becomes $0$ or $1$, we require a similar compatibility condition as 
Proposition $\ref{prop94}$ $(9),(10)$.
\item
The Kuranishi structure is invariant under permutation of $\ell$ interior 
marked points.
\item
The Kuranishi structure is invariant under the cyclic permutations 
 of the datas associated to each of the boundary components.
 \end{enumerate}
\end{lem}
\begin{rem}
By the same reason as we mentioned in Remark \ref{remnew164}, it is simpler to state 
cyclic symmetry in our situation than, for example, in Proposition \ref{Kuraeistspoly} (7). %\marginpar{Remark added.  2025 Sep KF.}
\end{rem}
The proof is similar to  those discussed already in various places and 
so omitted.
Those Kuranishi structures induce one on $\mathcal M_{\ell}^{\text{ann}}
(\vec\kappa,\vec{\upsilon};\vec p,\vec q;B)^{\boxplus 1}$ 
\par
The same hold for versions including ${}_{\mathcal L}\vec{\frak f},{}_{\mathcal U}\vec{\frak f}$
and (a product of) $[0,1]$ factors. %\marginpar{A sentence added.  KF 2025 Aug 28}
\par
We need to take our system of Kuranishi structures on the space
$[0,1]^{\vert {}_{\mathcal L}\frak f_{\bf f}\vert + \vert {}_{\mathcal U}\frak f_{\bf f}\vert} \times \mathcal M_{\ell}^{\text{ann}}
(\vec\kappa,\vec{\upsilon};\vec p,\vec q;B;{}_{\mathcal L}\vec{\frak f},{}_{\mathcal U}\vec{\frak f})$ so that it is compatible 
with those that have been already been constructed 
at the fibers of $\Sigma_1, \Sigma_2 \in \mathcal M_{(1,1);0}^{\text{ann}}$.
We state this compatibility below.\footnote{Forgetful map to $\mathcal M_{(1,1);0}^{\text{ann}}$ 
is not defined for $\mathcal M_{(\vec\kappa,\vec k),(\vec{\upsilon},\vec r);\ell}^{\text{\rm ann}}
(B;\vec p,\vec q)$ etc.  This is because it does not contain an information to specify 
which marked point is the $0$-th one.}
\par
We first consider the fiber of  $\Sigma_1$.
We consider the moduli space $[0,1]^{\vert {}_{\mathcal L}{\frak f}_{\rm{\bf f}}\vert} \times{\mathcal M}_{\ell_1+1}((\vec{\kappa},\vec p,m_1);B_1;{}_{\mathcal L}\vec{\frak f})$
introduced in Definition \ref{defMellveck}.
Among the $\ell_1 +1$ interior marked points of elements of $[0,1]^{\vert {}_{\mathcal L}{\frak f}_{\rm{\bf f}}\vert} \times{\mathcal M}_{\ell_1+1;\vec k}((\vec{\kappa},\vec p,m_1);B_1;{}_{\mathcal L}\vec{\frak f})$
we handle one of them in a different way from the others as 
we did in Subsection \ref{frakponcoho}. We call this marked point the zero-th marked point and take 
an evaluation map
$$
\text{\rm ev}^+_0 : [0,1]^{\vert {}_{\mathcal L}{\frak f}_{\rm{\bf f}}\vert} \times{\mathcal M}_{\ell_1+1;\vec k}((\vec{\kappa},\vec p,m_1);B_1;{}_{\mathcal L}\vec{\frak f})^{\boxplus 1} \to X,
$$
at that marked point. 
 We define a group $\Pi_2(\vec{\upsilon},\vec q)$  of equivalences classes of homotopy classes of disks in a similar way as the group 
$\Pi_2(\vec{\kappa})$ that is defined in Definition \ref{def211} using $U_{\upsilon_i}$ in place of $L_{\kappa_i}$, and
define $[0,1]^{\vert {}_{\mathcal U}{\frak f}_{\rm{\bf f}}\vert} \times {\mathcal M}_{\ell_2+1;\vec r}((\vec{\upsilon},\vec q,m_2);B_2;{}_{\mathcal U}\vec{\frak f})$ by performing the same substitution in Definition \ref{defMellveck}.
We define 
$$
\text{\rm ev}^+_0 : [0,1]^{\vert {}_{\mathcal U}{\frak f}_{\rm{\bf f}}\vert} \times{\mathcal M}_{\ell_2+1;\vec r}((\vec{\upsilon},\vec q,m_2);B_2;{}_{\mathcal U}\vec{\frak f})^{\boxplus 1} \to X.
$$
in a similar way.
We put Kuranishi structures on $[0,1]^{\vert {}_{\mathcal L}{\frak f}_{\rm{\bf f}}\vert} \times{\mathcal M}_{\ell_1+1;\vec k}((\vec{\kappa},\vec p,m_1);B_1;{}_{\mathcal L}\vec{\frak f})^{\boxplus 1}$ 
and on $\text{\rm ev}^+_0 : [0,1]^{\vert {}_{\mathcal U}{\frak f}_{\rm{\bf f}}\vert} \times{\mathcal M}_{\ell_2+1;\vec r}((\vec{\upsilon},\vec q,m_2);B_2;{}_{\mathcal U}\vec{\frak f})^{\boxplus 1} \to X
$ which are 
defined in Proposition \ref{Kuraeistspoly2}.
We remark that it is {\it different} from the one in Proposition \ref{Kuraeistspoly}.
The evaluation map $\text{\rm ev}^+_0$ is weakly submersive with respect to this 
Kuranishi structure by Proposition \ref{Kuraeistspoly2} (5). Therefore the fiber product
$$
\aligned
 &([0,1]^{\vert {}_{\mathcal U}{\frak f}_{\rm{\bf f}}\vert} \times{\mathcal M}_{\ell_2+1;\vec r}((\vec{\upsilon},\vec q,m_2);B_2;{}_{\mathcal U}\vec{\frak f})^{\boxplus 1})\\
& {}_{\text{\rm ev}^+_0}
\times_{\text{\rm ev}^+_0} 
 ([0,1]^{\vert {}_{\mathcal U}{\frak f}_{\rm{\bf f}}\vert} \times{\mathcal M}_{\ell_2+1;\vec r}((\vec{\upsilon},\vec q,m_1);B_2;{}_{\mathcal U}\vec{\frak f})^{\boxplus 1})
 \endaligned
$$
over $X$ has a Kuranishi structure.
By definition we have the following identity
 %\marginpar{Since this forget is not defined. 
%We use version of Subsection 9.3.  KF 2025 Aug 28.}
\begin{equation}\label{sigma1comp}
\aligned
&\frak{forget}_{\rm sc}^{-1}(\Sigma_1) 
\cap
[0,1]^{\vert {}_{\mathcal L}\frak f_{\bf f}\vert + \vert {}_{\mathcal U}\frak f_{\bf f}\vert} \times \mathcal M_{(\vec\kappa,\vec k,m_1),
(\vec{\upsilon},\vec r,m_2);\ell}^{\text{ann}}
(\vec p,\vec q; B;{}_{\mathcal L}\vec{\frak f},{}_{\mathcal U}\vec{\frak f})^{\boxplus 1}\\
= 
\bigcup &
([0,1]^{\vert {}_{\sU}{\frak f}_{\rm{\bf f}}\vert} \times{\mathcal M}_{\ell_2+1;\vec r}((\vec{\upsilon},\vec q,m_2);B_2;{}_{\sU}\vec{\frak f})^{\boxplus 1})\\
 &\, {}_{\text{\rm ev}^+_0}
\times_{\text{\rm ev}^+_0} 
 ([0,1]^{\vert {}_{\sL}{\frak f}_{\rm{\bf f}}\vert} \times{\mathcal M}_{\ell_1+1;\vec k}((\vec{\kappa},\vec p,m_1);B_1;{}_{\sL}\vec{\frak f})^{\boxplus 1})
\endaligned
\end{equation}
as sets. Here the union is taken over all $B_1$, $B_2$, $\mathbb L_1$, $\mathbb L_2$ 
such that $B_1 \#B_2 = B$ with respect to the obvious concatenation $\#$ and 
$(\mathbb L_1,\mathbb L_2)$ is a shuffle of $\underline\ell$.

\begin{lem}\label{lemsigma1comp}
We may choose the Kuranishi structure of the space
$[0,1]^{\vert {}_{\mathcal L}\frak f_{\bf f}\vert + \vert {}_{\mathcal U}\frak f_{\bf f}\vert} \times \mathcal M_{(\vec\kappa,\vec k,m_1),(\vec{\upsilon},\vec r,m_2);\ell}^{\text{\rm ann}}
(B;\vec p,\vec q;{}_{\mathcal L}\vec{\frak f},{}_{\mathcal U}\vec{\frak f})$
so that it is compatible with the equality $(\ref{sigma1comp})$.
Here we use the Kuranishi structure of Proposition $\ref{Kuraeistspoly2}$ 
in the right hand side.
\end{lem}
We remark that the union in the right hand side of  (\ref{sigma1comp}) is not a disjoint union.
The intersections of different members of the union are of codimension 
2. More precisely the left hand side is a `normal crossing divisor' and 
right hand side (regarded as a disjoint union) is its `normalization'.
We refer the reader to \cite[Lemma 2.6.13 and Lemma 3.4.14]{toric3} for the detailed
explanation of this point.
The proof of Lemma \ref{lemsigma1comp} is similar to the proof of 
\cite[Lemma 2.6.13]{toric3} and so is omitted.
\par
We next study the fiber of $\Sigma_2$.
Let 
\begin{equation}\label{deff1f2}
f^1 \in L_{\kappa_{i_1}}\cap U_{{\upsilon}_{i_2}}, 
\quad 
f^2 \in U_{{\upsilon}_{j_2}}\cap  L_{\kappa_{j_1}}.
\end{equation}
with
$1\le i_1 \le j_1 \le K$, $1\le i_2 \le j_2 \le R$.
We put
\begin{equation}\label{defvarrho}
\aligned
\vec\varrho_{1;(i_1,i_2,j_1,j_2)}
&= (\kappa_{\blue{0}}, \dots,\kappa_{i_1},\upsilon_{i_2},\dots,\upsilon_{j_2},\kappa_{j_1},\dots,\kappa_{\blue{K}}),
\\
\vec\varrho_{2;(i_1,i_2,j_1,j_2)}
&= (\upsilon_{\blue{0}}, \dots,\upsilon_{i_2},\kappa_{i_1},\dots,\kappa_{j_1},\upsilon_{j_2},\dots,\upsilon_{\blue{R}}).
\endaligned
\end{equation}
We also define
\begin{equation}
\aligned
\vec o_{1;(f^1,f^2)} = (p_{\blue{0}}, \dots,p_{\blue{i_1}}, f^1 ,q_{\blue{i_2 +1}}, \dots,q_{\blue{j_2}}, f^2, p_{\blue{j_1 +1}}, \dots,p_{\blue{K}}), \\
\vec o_{2;(f^1,f^2)} = (q_{\blue{0}}, \dots,q_{\blue{i_2}} ,f^1, p_{\blue{i_1 +1}}, \dots,p_{\blue{j_1}}, f^2, q_{\blue{j_2+1}}, \dots,q_{\blue{R}}).
\endaligned
\end{equation}
\begin{figure}[h]
\centering
\includegraphics[scale=0.5]{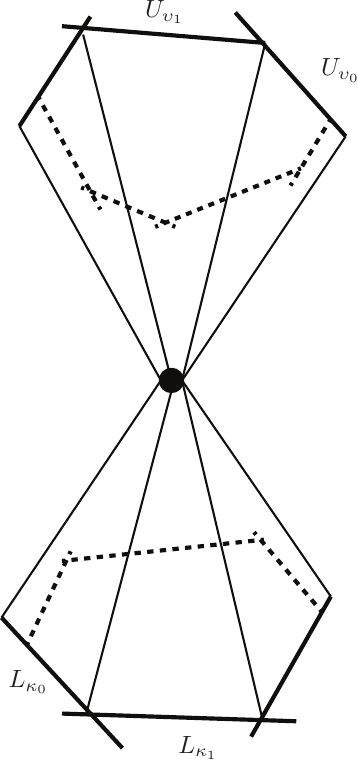}
\caption{An element of (\ref{sigma1comp})}
\label{Figure4-10}
\end{figure}
Again by definition we have the following identification as sets:
\begin{equation}\label{sigma2comp}
\aligned
&\frak{forget}_{\rm sc}^{-1}(\Sigma_2) 
\cap
([0,1]^{\vert {}_{\mathcal L}\frak f_{\bf f}\vert + \vert {}_{\mathcal U}\frak f_{\bf f}\vert} \times \mathcal M_{(\vec\kappa,\vec k,m_1),(\vec{\upsilon},\vec r,m_2);\ell}^{\text{ann}}
(\vec p,\vec q;B;{}_{\mathcal L}\vec{\frak f},{}_{\mathcal U}\vec{\frak f})^{\boxplus 1})
\\
&= 
\bigcup 
([0,1]^{\vert {}_{1}\frak f_{\bf f}\vert} \times {\mathcal M}_{\vert\mathbb L_1\vert;\vec s_1}((\vec\varrho_{1;(i_1,i_2,j_1,j_2)},\vec o_{1;(f^1,f^2)},m'_1);B_1;{}_{1}\vec{\frak f})^{\boxplus 1})\\
&\qquad \times ([0,1]^{\vert {}_{2}\frak f_{\bf f}\vert} \times
{\mathcal M}_{\vert\mathbb L_2\vert;\vec s_2}((\vec\varrho_{2;(i_1,i_2,j_1,j_2)},\vec o_{2;(f^1,f^2)},m'_2);B_2;{}_{2}\vec{\frak f})^{\boxplus 1}).
\endaligned
\end{equation}
${}_{1}\vec{\frak f}$ and $ {}_{2}\vec{\frak f}$ are induced from $ {}_{\mathcal L}\vec{\frak f}$ and $ {}_{\mathcal U}\vec{\frak f}$
in an obviousway.
 %\marginpar{Since this forget is not defined. 
%We use version of Subsection 9.3.  KF 2025 Aug 28.}
$m'_1$ and $m'_2$ are chosen so that the choice of $0$-th marked points of the  left hand side coincide with one  of the right hand side. 
(It is obvious that such a choice of $m'_1$ and $m'_2$ uniquely exists.)
Here the union is taken over all data $( B_1, B_2, \mathbb L_1, \mathbb L_2, 
i_1,i_2,j_1,j_2, f^1,f^2, \vec s_1, \vec s_2)$,
where $B_1 \#B_2 = B$ with respect to the obvious concatenation $\#$ and 
$(\mathbb L_1,\mathbb L_2)$ is a shuffle of $\underline\ell$.
The data $( i_1,i_2,j_1,j_2, f_1,f_2)$ are as above, and 
$(\vec s_1, \vec s_2)$ are defined by 
$$
\aligned
\vec s_1 &= (k_{\blue{0}}, \dots,k_{i_1-1},k'_{i_1},r''_{i_2},r_{i_2+1},\dots,r_{j_2-1},r'_{j_2},k''_{j_1},k_{j_1+1},\dots,k_{\blue{K}}), \\
\vec s_2 &= (r_{\blue{0}}, \dots,r_{j_2-1},r'_{\blue{i_2}},k''_{i_1},k_{i_1+1},\dots,k_{j_1-1},k'_{j_1},r''_{j_2},r_{j_2+1},\dots,r_{\blue{R}}),
\endaligned
$$
with $k'_{i_1} + k''_{i_1} = k_{i_1}$, $k'_{j_1} + k''_{j_1} = k_{j_1}$, $r'_{i_2} + r''_{i_2} = r_{i_2}$, and 
$r'_{j_2} + r''_{j_2} = r_{j_2}$.
\par
We use the relative version Proposition \ref{sec2relative} to define the Kuranishi structure on the 
right hand side of (\ref{sigma2comp}). 
\par
We remark that the product in the right hand side of  (\ref{sigma2comp}) is the 
direct product by Assumption $\ref{LUtrans}$.
Therefore the Kuranishi structure of the factors induce one on the 
product.
We also remark that to define the map $\frak{forget}_{\rm sc}$ we need to specify the $0$-th boundary 
marked points on $\partial_i\Sigma$ ($i=1,2$).
However the subset $\frak{forget}_{\rm sc}^{-1}(\Sigma_2)$ is independent of such a choice.
\begin{lem}\label{lemsigma2comp}
We may choose the Kuranishi structure of the space
$[0,1]^{\vert {}_{\mathcal L}\frak f_{\bf f}\vert + \vert {}_{\mathcal U}\frak f_{\bf f}\vert} \times \mathcal M_{(\vec\kappa,\vec k,m_1),(\vec{\upsilon},\vec r,m_2);\ell}^{\text{\rm ann}}
(B;\vec p,\vec q;{}_{\mathcal L}\vec{\frak f},{}_{\mathcal U}\vec{\frak f})^{\boxplus 1}$
so that it is compatible with the equality $(\ref{sigma2comp})$.
\end{lem}
We remark that the union of the right hand side of (\ref{sigma2comp}) is not a disjoint union.
In this case each member of the union has a Kuranishi structure with corners 
and they are glued along the boundaries and corners. 
The precise description of it is the same as that of \cite[Lemmas 2.6.22 and  3.4.14]{toric3}.
We refer the readers to the proof of \cite[Lemma 2.6.22]{toric3} for the 
proof of Lemma \ref{lemsigma2comp}.
\par
We next discuss CF-perturbations.
\begin{lem}\label{multisectionannul}
We can define a system of \blue{CF-perturbations} on 
$[0,1]^{\vert {}_{\mathcal L}\frak f_{\bf f}\vert + \vert {}_{\mathcal U}\frak f_{\bf f}\vert} \times \mathcal M_{(\vec\kappa,\vec k,m_1),(\vec{\upsilon},\vec r,m_2);\ell}^{\text{\rm ann}}
(B;\vec p,\vec q;{}_{\mathcal L}\vec{\frak f},{}_{\mathcal U}\vec{\frak f})^{\boxplus 1}$
so that the following holds:
\begin{enumerate}
\item
They are transversal to $0$.
\item 
They are compatible with the description of the 
boundary in Lemma $\ref{bdryMannu}$ $(3)$.
\item 
They are compatible with the forgetful maps 
forgetting the marked points labelled by  ${}_{\mathcal L}{\frak f}_{\bf e}$,
${}_{\mathcal U}{\frak f}_{\bf e}$.
\item
The compatibility of Proposition $\ref{bdryMannu}$ $(7)$ is enhanced to the compatibility of CF-perturbations.
\item
The \blue{CF-perturbations} are 
compatible with the identifications in 
Lemmas $\ref{lemsigma1comp}$ and $\ref{lemsigma2comp}$.
\item
The \blue{CF-perturbations} are invariant under the permutation of 
interior marked points.
\item
The \blue{CF-perturbations} are invariant under the cyclic permutation 
of the data in the same sense as Lemma $\ref{bdryMannu}$ $(7)$. 
\end{enumerate}
\end{lem}
We later use this \blue{CF-perturbation} to integrate the forms 
pulled back by evaluation maps.

\subsection{Proof of Theorem \ref{ZisPD}.}
\label{sec4maintheoremproof}
%In order to prove Theorem \ref{ZisPD}, we begin by changing the notation as in Section \ref{constcyclic}: we drop the assumption that $\kappa_i 
%\neq \kappa_{i+1}$, and ${\upsilon}_j \neq {\upsilon}_{j+1}$. We can associate to each such sequences $\vec{\kappa}$ and $\vec{\upsilon}$
 %\marginpar{Need to correct this subsection much. KF 2025 Sep.} reduced sequences $ \vec{\kappa}'$ and $ \vec{\upsilon}'$, by omitting repeated 
%elements, as well as sequence $ \vec k'$ and $ \vec r'$ recording the number of omitted elements. For each sequences $\vec p$ and $\vec q$ of 
%intersections among the Lagrangians in the reduced sequences $ \vec{\kappa}'$ and $ \vec{\upsilon}'$, and for each relative homotopy class $B %
%\in \Pi_2(X;\vec \kappa',\vec {\upsilon}',\vec p,\vec q)$, we then define
In this subsection we use the notation as in Section \ref{constcyclic}.
Namely we put:
$$
\mathcal M_{\ell}^{\text{ann}}
(\vec\kappa,\vec{\upsilon};\vec p,\vec q;B)
= 
\mathcal M_{(\vec\kappa',\vec k',m'_1),(\vec{\upsilon}',\vec r',m'_2);\ell}^{\text{ann}}
(\vec p,\vec q;B)
$$
as in (\ref{newnew2018}).

For each such $B$, we define  $\rho(B) \in \Lambda_0\setminus \Lambda_+$ via the formula
\begin{equation}
\aligned
\rho(B)
= 
T^{\int_{\Sigma}u^*\omega}
&\exp
\left(
\int_{\Sigma} u^*\frak b_2
+
\sum_{i=0}^K \int_{z^1_i}^{z^1_{i+1}}u^*\theta_{L_{\kappa_i}}
+
\sum_{i=0}^R \int_{z^2_i}^{z^2_{i+1}}u^*\theta_{U_{\upsilon_i}}
\right)% \\
% &\times\exp
% \left(
% \sum_{i=0}^K \int_{z^1_i}^{z^1_{i+1}}u^*b_{\kappa_i,0}
% +
% \sum_{i=0}^R \int_{z^2_i}^{z^2_{i+1}}u^*b_{\upsilon_i,0}
% \right)
,
\endaligned
\end{equation}
where we use the splitting (\ref{decomposefrakb}), the 1-forms $\theta_{L_{\kappa}}$ and $\theta_{U_{\upsilon}}$
as in (\ref{fixthetaL}) for each $L_{\kappa}$, $U_{\upsilon}$, and the bounding chains $b_{\kappa,+}$ and $b_{\upsilon,+}$
of $L_{\kappa}$ and $U_{\upsilon}$, respectively, which have strictly positive energy. % We decompose them into 

We consider the map
\begin{equation}
\aligned
  \text{\rm ev} = (\text{\rm ev}^+, \text{\rm ev}^{\partial,1},
\text{\rm ev}^{\partial,2}) &: 
\mathcal M_{\ell}^{\text{ann}}
(\vec\kappa,\vec{\upsilon};\vec p,\vec q;B)\\
&\to 
X^{\ell}
\times \prod_{i=0}^K (\tilde L_{\kappa_{i-1}} \times_X \tilde L_{\kappa_{i}})
\times \prod_{i=0}^R (\tilde U_{{\upsilon}_{i-1}} \times_X \tilde U_{{\upsilon}_{i}})
\endaligned
\end{equation}
which is given by the evaluation maps \eqref{eq:evaluation_map_annuli} at the  interior or boundary  
marked points.
We use the enumeration of the boundary marked points we have for 
an object of $\mathcal M_{\ell}^{\text{ann}}
(\vec\kappa,\vec{\upsilon};\vec p,\vec q;B)$.
(In other words on $\mathcal M_{(\vec\kappa',\vec k',m'_1),(\vec{\upsilon}',\vec r',m'_2);\ell}^{\text{ann}}
(\vec p,\vec q;B)$ such enumerations are obtained in a way similar to Definition \ref{defn917}.)

\begin{proof}[Proof of of Theorem \ref{ZisPD}]
At the level of differential forms, we can define a map 
$
Z^{\rm form}(\text{\bf x},\text{\bf y})
$,
%the form version, 
in the same way as in Subsection \ref{inproZ}.\footnote{
We use the fact that the Lagrangian submanifolds $U_{\upsilon}$ is transversal 
to the elements $L_{\kappa}$ of $\cL$.
In fact, then the basis $f^1, f^2$ appearing in (\ref{eq:formula_Z-pairing}) are of finite number.}
Here
$$
\text{\bf x} \in CH_*(\cL_{\rm c.u.}^{\rm form}) \textrm{ and }
\text{\bf y} \in CH_*(\sU_{\rm c.u.}^{\rm form})
$$
are Hochschild chains of the categories, defined using differential forms, corresponding to  $\cL$ and $\sU$.
In view of the definition it suffices to prove
\begin{equation}\label{form45}
\langle
\widehat{\frak p}^{\rm f.u.c.\bf b}(\text{\bf x})
,
\widehat{\frak p}^{\rm f.u.c.\bf b_{\sU}}(\text{\bf y})
\rangle_{\text{\rm PD}_X} 
=
(-1)^{n(n-1)/2} Z^{\rm form}(\text{\bf x},\text{\bf y}) 
\end{equation}
under the assumption 
\begin{equation}\label{delHvanish}
\delta_H(\text{\bf x}) = \delta_H(\text{\bf y}) = 0.
\end{equation}
\par
Let
$$
\text{\bf x} = x_0 \otimes x_1 \otimes \cdots \otimes x_K
\in CH_*(\cL^{\rm form}_{\rm c.u.}, \cL^{\rm form}_{\rm c.u.}),
\quad
\text{\bf y} = y_0 \otimes y_1 \otimes \cdots \otimes y_R
\in CH_*(\sU^{\rm form}_{\rm c.u.}, \sU^{\rm form}_{\rm c.u.})
$$
and 
$$
x_i \in CF(L_{\kappa_i},L_{\kappa_{i+1}};\F),
\quad
y_i \in CF(U_{\rho_i},U_{\rho_{i+1}};\F).
$$
%We remark that we do {\it not} assume (\ref{delHvanish}) here.

\par
%Let $\gamma : [0,1] \to \mathcal M_{(1,1);0}^{\text{\rm ann}}$ 
 %\marginpar{MA190710: I cannot see where this path is used.
%KF. Removed it is $P$ in the proof of Lemma \ref{lem414}.}
%be a smooth path 
%whose image is contained in a single stratum of 
%$\mathcal M_{(1,1);0}^{\text{\rm ann}}$.
%\par
For each %\marginpar{We need sign. 2026 April} 
curve $\Sigma$ in $\mathcal M_{(1,1);0}^{\text{\rm ann}}$, we define\index[syindex]{Ifrabsigma@$ \frak I^{\bf b}_{\Sigma}$}
\begin{equation}\label{form1634}
 \aligned
 &
 \frak I^{\bf b}_{\Sigma}( \text{\bf x}, \text{\bf y}) \equiv  \sum_{\vec k, \vec r, \ell} \sum_{B} \rho(B) \int_{\frak{forget}_{\rm sc}^{-1}(\Sigma)\cap\mathcal M_{\ell}^{\text{ann}}
(\vec\kappa,\vec{\upsilon};\vec p,\vec q;B)} \\
&\qquad\qquad\qquad\qquad\qquad (\vec{\text{\rm ev}}^{\partial,1})^*(b_{\kappa_0}^{k_0} \times x_0 \times \cdots \times x_K \times b_{\kappa_K}^{k_K})  \\
&\qquad\qquad\qquad\qquad\qquad
\wedge (\vec{\text{\rm ev}}^{+})^* \frak b_{+}^{\ell}  \wedge (\vec{\text{\rm ev}}^{\partial,2})^*(b_{\upsilon_0}^{r_0} \times y_0 \times \cdots \times y_R \times b_{\upsilon_R}^{r_R}).
\endaligned
\end{equation}

%In the case $L_{\kappa_i} \ne L_{\kappa_{i+1}}$, the element
%$x_i$ corresponds to an intersection point $\in L_{\kappa_i} \cap L_{\kappa_{i+1}}$; we identify it with the constant function $1  \in \Omega 
%( L_{\kappa_i} \cap L_{\kappa_{i+1}})$.
Here $x_i$, $y_i$ are either intersection point  $\in L_{\kappa_i} \cap L_{\kappa_{i+1}}$
self-intersection point of $L_{\kappa_i}$ or differential form on $L_{\kappa_i}$.
Actually we need to put signs as in Subsection \ref{inproZ} and Section \ref{sec:ori}.
We omit the sign in this subsection.
\par
(\ref{form1634}) is the case when ${\bf x}$, ${\bf y}$ do not contain ${\bf e}^+$ and/or ${\bf f}$.
To include ${\bf e}^+$ and/or ${\bf f}$ we use the evaluation maps at the marked points labelled by 
${}_{\mathcal L}\vec{\frak f}$, ${}_{\mathcal U}\vec{\frak f}$ in the same way as Subsection \ref{subsec:opeartor}.
\begin{lem}\label{lem411}
We have the equality
$$
\langle
\widehat{\frak p}^{\rm f.u.c.{\bf b}}(\text{\bf x})
,
\widehat{\frak p}^{\rm f.u.c.{\bf b}_{\sU}}(\text{\bf y})
\rangle_{\text{\rm PD}_X} =   \frak I^{\bf b}_{\Sigma_1}( \text{\bf x}, \text{\bf y}).
$$
Here ${\bf b}_{\sU} = (\frak b,b_{\sU})$.
\end{lem}
\begin{proof}
This is a consequence of Lemmas \ref{bdryMannu}, \ref{lemsigma1comp} and \ref{multisectionannul} and the definition.
\end{proof}
\begin{lem}\label{lem412}
We have the equality
$$
%\aligned
%&
Z^{\rm form}(\text{\bf x},\text{\bf y} ) =  \frak I^{\bf b}_{\Sigma_2}( \text{\bf x}, \text{\bf y}).
%\\
% &= \sum_{\vec k, \vec r, \ell} \sum_{B} \rho(B) \int_{\frak{forget}^{-1}(\Sigma_1)\cap\mathcal M_{\vec\kappa,\vec{\upsilon};\ell}^{\text{\rm ann}}
% (B;\vec p,\vec q)} 
% (\vec{\text{\rm ev}}^{\partial,1})^*(b_{\kappa_0}^{k_0} \times x_0 \times \cdots \times x_K \times b_{\kappa_K}^{k_K})  \\
% &\qquad\qquad\qquad\qquad\qquad
% {\color{blue} \wedge (\vec{\text{\rm ev}}^{+})^* \frak b_{+}^{\ell}}  \wedge (\vec{\text{\rm ev}}^{\partial,2})^*(b_{\upsilon_0}^{r_0} \times y_0 \times \cdots \times y_R \times b_{\upsilon_R}^{r_R}) 
% \endaligned
$$
\end{lem}
\begin{proof}
This is a consequence of Lemmas \ref{bdryMannu}, \ref{lemsigma2comp} and \ref{multisectionannul} and the definition.
\end{proof}

We remark that $\text{\rm Int}\mathcal M_{(1,1);0}^{\text{ann}} \setminus \{\Sigma_1\}$
consists of a single stratum of dimension 2.
We choose sequences $\Sigma_{1,m}$ and $\Sigma_{2,m}$ on this stratum 
so that 
$$
\lim_{m\to \infty} \Sigma_{1,m} = \Sigma_{1}, \quad \lim_{m\to \infty} \Sigma_{2,m} = \Sigma_{2}
$$
in the moduli space $\mathcal M_{(1,1);0}$. %stable surface $ \Sigma_{(1,1);0}$.

We may also assume that $\Sigma_{1,m}$ and $\Sigma_{2,m}$ are regular values of 
the restriction of $\frak{forget}$ on the zero set of the 
family of CF-perturbations defined in  Lemma \ref{multisectionannul}.
\begin{lem}\label{annlimit}
  We have the equality
  \begin{equation}
    \lim_{m\to\infty} \frak I^{\bf b}_{\Sigma_{c,m}}( \text{\bf x}, \text{\bf y}) = \frak I^{\bf b}_{\Sigma_c}( \text{\bf x}, \text{\bf y})
  \end{equation}
% $$
% \aligned
% &\lim_{m\to\infty}
% \int_{\frak{forget}^{-1}(\Sigma_{c,m})\cap\mathcal M_{(\vec\kappa,\vec k),(\vec{\upsilon},\vec r);\ell}^{\text{\rm ann}}
% (B;\vec p,\vec q)} 
% (\vec{\text{\rm ev}}^{\partial,2,1})^*(\vec b_{\kappa,+}) \wedge (\vec{\text{\rm ev}}^{\partial,2,2})^*(\vec b_{\upsilon,+})\\
% &\qquad\qquad\qquad\qquad\qquad\qquad\qquad\qquad\qquad \wedge
% (\vec{\text{\rm ev}}^{\partial,1,1})^*(\frak x) \wedge (\vec{\text{\rm ev}}^{\partial,1,2})^*(\frak y) \\
% &=\int_{\frak{forget}^{-1}(\Sigma_c)\cap\mathcal M_{(\vec\kappa,\vec k),(\vec{\upsilon},\vec r);\ell}^{\text{\rm ann}}
% (B;\vec p,\vec q)} 
% (\vec{\text{\rm ev}}^{\partial,2,1})^*(\vec b_{\kappa,+}) \wedge (\vec{\text{\rm ev}}^{\partial,2,2})^*(\vec b_{\upsilon,+})\\
% &\qquad\qquad\qquad\qquad\qquad\qquad\qquad\qquad\qquad \wedge
% (\vec{\text{\rm ev}}^{\partial,1,1})^*(\frak x) \wedge (\vec{\text{\rm ev}}^{\partial,1,2})^*(\frak y) 
% \endaligned
% $$
for $c=1,2$.
\end{lem}
\begin{proof}
This equality is slightly nontrivial only because 
$\frak{forget}$ is not smooth at the fiber of $\Sigma_1$, $\Sigma_2$.
We can however prove Lemma \ref{annlimit} using the exponential decay estimate of \cite[Lemma 1.58]{fooo09}.
The argument is the same as the proof of \cite[Lemmas 3.4.26 and 27]{toric3}.
\end{proof}
\begin{lem}\label{lem414}
There exists a sequence of operators $\mathcal K_m$ such that\index[syindex]{Kcalm@$\mathcal K_m$}
\begin{equation}\label{homotopyannul}
\mathcal K_m(\delta^H\text{\bf x},\text{\bf y}) + (-1)^{\deg'{\bf x}}\mathcal K_m(\text{\bf x},\delta^H\text{\bf y})
=  \frak I^{\bf b}_{\Sigma_{1,m}}( \text{\bf x}, \text{\bf y}) - \frak I^{\bf b}_{\Sigma_{2,m}}( \text{\bf x}, \text{\bf y}).
\end{equation}
\end{lem}
\begin{proof}
We take a path $P \subset \text{\rm Int}\mathcal M_{(1,1);0}^{\text{ann}} \setminus \{\Sigma_1\}$ 
such that $\partial P = \{\Sigma_{1,m},\Sigma_{2,m}\}$.
We may also assume that $P$ is  transverse to
  %regular value of 
the restriction of $\frak{forget}$ on the zero set of the \blue{CF-perturbation} defined in  Lemma \ref{multisectionannul}.
We now put %\marginpar{We need sign  2026 April.}
$$
\aligned
&
\mathcal K_m(\text{\bf x},\text{\bf y}) \\
 &= \sum_{\vec k, \vec r, \ell} \sum_{B} \rho(B) \int_{\frak{forget}_{\rm sc}^{-1}(P)\cap\mathcal M_{\ell}^{\text{ann}}
(\vec\kappa,\vec{\upsilon};\vec p,\vec q;B)} 
(\vec{\text{\rm ev}}^{\partial,1})^*(b_{\kappa_0}^{k_0} \times x_0 \times \cdots \times x_K \times b_{\kappa_K}^{k_K})  \\
&\qquad\qquad\qquad\qquad\qquad
\wedge (\vec{\text{\rm ev}}^{+})^* \frak b_{+}^{\ell}  \wedge (\vec{\text{\rm ev}}^{\partial,2})^*(b_{\upsilon_0}^{r_0} \times y_0 \times \cdots \times y_R \times b_{\upsilon_R}^{r_R}).
\endaligned
$$
This is the case when ${\bf x}$, ${\bf y}$ do not contain ${\bf e}^+$ and/or ${\bf f}$.
To include ${\bf e}^+$ and/or ${\bf f}$ we use the evaluation maps at the marked points labelled by 
${}_{\mathcal L}\vec{\frak f}$, ${}_{\mathcal U}\vec{\frak f}$ in the same way as Subsection \ref{subsec:opeartor}.
We omit the sign in this subsection. (See Proposition \ref{pairingdelta}.)

Equality (\ref{homotopyannul})
is then a consequence of Lemmas \ref{bdryMannu}, \ref{multisectionannul} and Stokes' theorem. 
In fact the first term of the left hand side is induced by Lemmas \ref{bdryMannu} (3) (b)(d)(f) 
and the second term of the left hand side is induced by Lemmas \ref{bdryMannu} (3) (c)(e)(g).
\end{proof}
The formula (\ref{form45}) under the assumption (\ref{delHvanish}) 
follows from Lemmas \ref{lem411}, \ref{lem412}, \ref{annlimit}, \ref{lem414}.
The proof of Theorem \ref{ZisPD} is now complete except the orientation and sign issue 
which will be discussed in the next section.
 %\marginpar{maybe we need to say something about the sign.}
\end{proof}

\section{Orientation and sign.}\label{sec:ori}

We define the sign of the operations $\frak q$ and $\frak p$ in Definition \ref{def827} and Theorem
\ref{themp} as follows, respectively. %\marginpar{This section is new KO 2025 Sep.}
\begin{eqnarray}
\frak q^{\text{\rm form}}_{\ell;\vec{\kappa},\vec p;B}(\text{\bf g},\text{\bf h} )
&=&(-1)^{\epsilon(h_1, \dots, h_k) + \delta(g_1, \dots, g_\ell)} 
(\text{\rm ev}_{0})_!
\left((\text{\rm ev}^+)^*\text{\bf g}
\wedge
(\text{\rm ev}_{>0})^*\text{\bf h}
\right) \label{form21ten1}
\\
{\frak p}^{\rm form}_{\ell;\vec{\kappa},\vec p;B}(\text{\bf g},\text{\bf h}) 
&=&(-1)^{\epsilon(h_1, \dots, h_k) + \delta(g_1, \dots, g_\ell)+ \mu(B) + 1} \nonumber \\
&&(\text{\rm ev}^+_{0})_!
\left((\text{\rm ev}^+)^*\text{\bf g}
\wedge
(\text{\rm ev}_{>0})^*\text{\bf h}
\right) \label{form21ten2}
\end{eqnarray}
Here\index[syindex]{epsilonh_1dots@$\epsilon(h_1,\dots,h_k)$}\index[syindex]{deltah_1dots@$\delta(g_1,\dots,g_\ell)$}
\begin{eqnarray} 
\epsilon(h_1, \dots, h_k) & = & \sum_{i=1}^{k} (i + \sum_{j=1}^{i-1} \mu_j) (\text{\rm d-}\deg h_i -1) + 1,\label{212plus1}  \\
\delta(g_1, \dots, g_\ell) & =  & \sum_{i=1}^\ell \text{\rm d-}\deg g_i. \nonumber
\end{eqnarray}
We define $\mu_j$ as follows.\index[syindex]{myui@$\mu_j$}  Let $z(j)$ be as in Definition \ref{defn917}.
We put:
$$
\mu_j =
\begin{cases} 
0      & \text{if $z(j)$ is a diagonal marked point},   \\
\mu(\lambda_{p_j})    & \text{if $z(j)$ is a switching marked point}.
\end{cases}
$$  
Here $\mu(\lambda_{p_j})$ is the index of the operator (\ref{orientationcomplex}).\footnote{$\mu_i$ is the one denoted by $\mu(R_{\alpha_i})$ in \cite[(3.1)]{ono}.} 
We then define 
\begin{equation}\label{defnform213}
 \mu(B)=\sum_{j=0}^{k} \mu_j.
\end{equation}
 \index[syindex]{myuB@$\mu(B)$}
\par
$\text{\rm d-}\deg g_i$ \index[syindex]{ddeg@$\text{\rm d-}\deg g$}is the degree of $g_i$ as a differential form on $X$ and 
$\text{\rm d-}\deg h_i$\footnote{$\deg h_i = {\rm d}\text{-}\deg h_i + \mu_i$.} 
the degree of $h_i$ as a differential form on $\tilde L_{\kappa_{i-1}} \times_X \tilde L_{\kappa_i}$.  
When $p_i$ is a self-intersection point, $h_i$ is a constant function on the self-intersection point $p_i=(p^1_i, p^2_i)$.
When $L_{\kappa_{i-1}} \ne L_{\kappa_i}$ $h_i$ is a constant function on the intersection point.
So $\text{\rm d-}\deg = 0$ in those cases.
\par
We remark that when we consider a single embedded Lagrangian submanifold $L$, the numbers
$\mu_j$ are $0$. Therefore (\ref{form21ten1}) coincides with (\ref{defqformula}).
\par
We use $\frak q$-perturbation for (\ref{form21ten1}) and $\frak p$-perturbation for (\ref{form21ten2}), but we abbreviate the superscript in this section.  
In \cite{ono}, we fix the sign and orientation convention and verified the filtered $A_{\infty}$-relations 
for $\mathfrak m$, namely, the case that $\ell = 0$.  We can also verify  formulas 
for $\mathfrak q$ and $\mathfrak p$ using
our convention above in 
a similar way.  

\begin{rem}
In the case $\ell = 0$ and $\vec{\kappa} = \kappa$ (a single embedded Lagrangian) $B =0$,
we have
$$
\frak p^{\rm form}_{0;\kappa;0}(\emptyset,h) = (-1)^{\epsilon(h)+1} ({i_L})_{!}(h).
$$ 
Since $\epsilon(h) = \deg h$ we obtain Theorem \ref{themp} (1).
\end{rem}

\subsection{Sign for the Cardy relation} 
\label{sec:signCardy}

In this subsection, we verify the sign in the Cardy relation (Theorem \ref{ZisPD}) based on \cite{ono}.  
Namely we prove Theorem \ref{ZisPD} together with sign.
We consider the moduli space $\mathcal M_{\ell}(\vec{\kappa},\vec p;B)$ in 
(\ref{eq:change_notation_moduli_space}).

We consider the case $\ell = 1$. Let $k+1$ be the number of boundary marked points. 
Note that
$k+1$ is determined by $\vec{\kappa}$.  Since $B$ determines $\vec{\kappa},\vec p$ 
we write $\mathcal M_{1;k+1}(B)$\index[syindex]{M1k+1B@$\mathcal M_{1;k+1}(B)$} in place of $\mathcal M_{1}(\vec{\kappa},\vec p;B)$
for simplicity in this section.
We remark that for an element of $\mathcal M_{1;k+1}(B)$ enumeration of its boundary marked points is given. 
\begin{rem}\label{rem211}
We recall that the orientation of $\mathcal M_{1;k+1}(B)$ is defined as follows. %\marginpar{Remark added 2025 Dec.}
(See \cite[Section 8.3]{fooo092}.)
Let $(u,(\vec z,z^+))$ be an element of $\overset{\circ}{\mathcal M}_{1;k+1}(B)$.
Here $\vec z = (z_0,\dots,z_k)$ are the boundary marked points and $z^+$ the interior marked point.
$u : D^2 \to X$ is a pseudo-holomorphic map satisfying the boundary condition 
specified by $B$. For each switching marked point $z_i$ (such that $u(z_i) = p_i$) we consider 
the operator $\overline{\partial}_{p_{i}}$ as in (\ref{orientationcomplex}).  We denote its index by $\mathcal {V}_{p_i}$.\footnote{We may choose 
the path $\lambda_{p_i}$ so that the index of  $\overline{\partial}_{p_i}$ is positive.   We may also assume that  $\overline{\partial}_{p_{i}}$ is surjective. 
So the index is a vector space.}  
Denote by $\Theta_{p_i}$ the determinant line bundle of $\mathcal {V}_{p_i}$.
Then the relative spin structure determines the orientation of 
$$
{\rm Index}D_u\overline{\partial} \oplus \bigoplus \mathcal{V}_{p_i}.
$$
We fix an orientation of $\mathcal{V}_{p_i}$. The orientation of ${\rm Index}D_u\overline{\partial}$ is then determined.
The tangent bundle $T_{((\vec z,z^+);u)}\mathcal M_{1;k+1}(B)$ 
is identified with 
$$
\left({\rm Index}D_u\overline{\partial} \oplus \C \oplus \R^{k+1} \right) / sl(2;\R).
$$
In other words
$$
T_{((\vec z,z^+);u)}\mathcal M_{1;k+1}(B) \oplus {sl}(2;\R)
\cong {\rm Index}D_u\overline{\partial} \oplus \C \oplus \R^{k+1}.
$$
We orient $T_{((\vec z,z^+);u)}\mathcal M_{1;k+1}(B)$ so that this 
isomorphism respects the orientation.  See \cite[(8.2.1.2)]{fooo092}.
\par
Alternatively we may proceed as follows.
Any element of $\overset{\circ}{\mathcal M}_{1;k+1}(B)$ have   unique representative 
such that $z_0 = 1$, $z^+ = 0$.  Therefore  we may identify
$$
T_{((\vec z,z^+);u)}\mathcal M_{1;k+1}(B)
= {\rm Index}D_u\overline{\partial} \oplus \R^{k}.
$$
Here $\R^k$ is the parameter space to move $z_1,\dots,z_k$.
We may orient $T_{((\vec z,z^+);u)}\mathcal M_{1;k+1}(B)$ 
so that this isomorphism respects the orientation.
We put the symbol $*$ to the moduli space and write $\mathcal M^*_{1;k+1}(B)$
etc. when we use this second convention of the orientation.\index[syindex]{MhatStar1k+1@$\mathcal M^*_{1;k+1}(B)$} 
\par
We claim that the difference of the two orientations is $(-1)^{k}$.
This is a consequence of the fact that we consider  the {\it right} action\footnote{See \cite[(8.2.1.2)]{fooo092}.}
of $\text{PSL}(2;\R)$ and the parameter to move $z_0$ is put on the top of $\R^{k+1}$
in the first method.
In other words we have 
$$
\mathcal M^*_{1;k+1}(B) = (-1)^{k}\mathcal M_{1;k+1}(B).
$$
\par
During the proof of Proposition \ref{pp} we will use the second convention 
of the orientation. Since we follow  \cite[Chapter 8]{fooo092} to take the first convention 
there is a term $(-1)^{k^{(1)} + k^{(2)}}$ in (\ref{form214}).
\end{rem}
\par
Let ${\rm ev}_{0,B_1}^+:{\mathcal M}_{1; k^{(1)}+1}(B_1) \to X$ and  ${\rm ev}_{0,B_2}^+ : {\mathcal M}_{1; k^{(2)}+1}(B_2) \to X$ be the evaluation maps 
at the 0-th interior marked points, respectively.\index[syindex]{ev0B+@${\rm ev}_{0,B}^+$}
We take their fiber product ${\mathcal M}_{1; k^{(1)}+1}(B_1) \times_X {\mathcal M}_{1; k^{(2)}+1}(B_2)$.  
Denote by $\pi_{B_1}$ (resp. $\pi_{B_2}$) the projection to the first (resp. second) factor.  
In the definition of $\frak p$ above, $(h_1, \dots, h_k)$ are the inputs.  Since we denote the inputs by $(h_0, \dots, h_k)$ 
for Hochschild chains, we have 
 $$
 \epsilon(h^{(s)}_0, \dots, h^{(s)}_k)= \sum_{i=0}^{k^{(s)}}((i+1)+\sum_{j=0}^{i-1} \mu^{(s)}_j)(\text{\rm d-}\deg h^{(s)}_i -1) + 1$$
for $s=1,2$.  
$
 \mu(B_s)
$
is defined in the same way as (\ref{defnform213}).
 
 \begin{prop}\label{pp}
 \begin{equation}\label{form214}
 \aligned & 
 (-1)^{\mu(B_1) + \mu(B_2) + k^{(1)} + k^{(2)} + (n+1) \deg' {\bf h}^{(2)}} \int_X {\mathfrak p}_{B_1}  ({\bf h}^{(1)}) \wedge {\mathfrak p}_{B_2} ({\bf h}^{(2)}) \\
 = & (-1)^{\delta_1} \int_{{\mathcal M}^*_{1; k^{(1)}+1}(B_1)^{\frak p} \times_X {\mathcal M}^*_{1; k^{(2)}+1}(B_2)^{\frak p}} 
 \pi_{B_1}^* {\rm ev}_{B_1}^* {\bf h}^{(1)} \wedge \pi_{B_2}^* {\rm ev}_{B_2}^* {\bf h}^{(2)}.  
 \endaligned
 \end{equation}
 Here 
 \begin{equation}
 \aligned \delta_1 = &  \epsilon({\bf h}^{(1)}) + \epsilon({\bf h}^{(2)}) + (\mu(B_1) + k^{(1)} + 1) \sum_{j=0}^{k^{(2)}} (\text{\rm d-} \deg h_j^{(2)} -1)  \\
 &  + (n+1) (\mu(B_1) + k^{(1)}) + \mu(B_1) \mu(B_2) + (k^{(1)} + 1) \mu(B_2).
 \endaligned
 \end{equation}
 \end{prop} 
 %$*$ is related to Remark \ref{rem212} and is discussed at the end of the proof.
 \begin{proof}
 Write $\delta_0= {\epsilon({\bf h}^{(1)}) + \epsilon({\bf h}^{(2)})}$.  
 \begin{eqnarray}
& &  (-1)^{\mu(B_1) + \mu(B_2) + k^{(1)} + k^{(2)}}
\int_X {\mathfrak p}_{B_1} ({\bf h}^{(1)}) \wedge {\mathfrak p}_{B_2}({\bf h}^{(2)}) \nonumber \\
& = & (-1)^{\delta_0}
\int_X ({\rm ev}_{0,B_1}^+)_! {\rm{ev}}_{B_1}^* {\bf h}^{(1)} \wedge ({\rm ev}_{0,B_2}^+)_! {\rm{ev}}_{B_2}^* {\bf h}^{(2)} \nonumber \\
 & = &  (-1)^{\delta_0}  \int_X ({\rm ev}_{0,B_2}^+)_!({\rm ev}_{0,B_2}^* ({\rm ev}_{0,B_1}^+)_!   ({\rm{ev}}_{B_1}^*{\bf h}^{(1)}) \wedge {\rm{ev}}_{B_2}^*({\bf h}^{(2)})) \nonumber \\
 & = &  (-1)^{\delta_0}  \int_X ({\rm ev}_{0,B_2}^+)_! (((\pi_{B_2})_! \pi_{B_1}^* {{\rm ev}}_{B_1}^*{\bf h}^{(1)}) \wedge {{\rm ev}}_{B_2}^*{\bf h}^{(2)}) 
 \nonumber \\
 & = & (-1)^{\delta_0} \int_{{\mathcal M}^*_{1; k^{(2)}+1}(B_2)^{\frak p}} ((\pi_{B_2})_! \pi_{B_1}^*{{\rm ev}}_{B_1}^*{\bf h}^{(1)}) \wedge 
 {\rm ev}_{B_2}^*{\bf h}^{(2)} \nonumber \\
 & = & (-1)^{\delta_0+ d_1} \int_{{\mathcal M}^*_{1; k^{(2)}+1}(B_2)^{\frak p}}  (\pi_{B_2})_! (\pi_{B_2}^*{{\rm ev}}_{B_2}^*{\bf h}^{(2)} \wedge \pi_{B_1}^* {{\rm ev}}_{B_1}^* {\bf h}^{(1)} )\nonumber \\
 & = & (-1)^{\delta_0+ d_1} \int_{{\mathcal M}^*_{1; k^{(2)}+1}(B_2)^{\frak p} \times_X {\mathcal M}^*_{1; k^{(1)}+1}(B_1)^{\frak p}} 
 \pi_{B_2}^* {{\rm ev}}_{B_2}^*{\bf h}^{(2)} \wedge \pi_{B_1}^* {{\rm ev}}_{B_1}^*  {\bf h}^{(1)} \nonumber \\
  & = & (-1)^{\delta_0+ d_1+d_2} \int_{{\mathcal M}^*_{1; k^{(2)}+1}(B_2)^{\frak p} \times_X {\mathcal M}^*_{1; k^{(1)}+1}(B_1)^{\frak p}} 
 \pi_{B_1}^* {{\rm ev}}_{B_1}^*{\bf h}^{(1)} \wedge \pi_{B_2}^* {{\rm ev}}_{B_2}^*  {\bf h}^{(2)} \nonumber \\
 & = & (-1)^{\delta_0+ d_1+ d_2 +d_3} \int_{{\mathcal M}^*_{1; k^{(1)}+1}(B_1)^{\frak p} \times_X {\mathcal M}^*_{1; k^{(2)}+1}(B_2)^{\frak p}} 
 \pi_{B_1}^* {{\rm ev}}_{B_1}^*{\bf h}^{(1)} \wedge \pi_{B_2}^* {{\rm ev}}_{B_2}^*  {\bf h}^{(2)}, \nonumber 
 \end{eqnarray}
 where
 \begin{eqnarray} 
 d_1 & = & \text{\rm d-}\deg (\pi_{B_2})_! \pi_{B_1}^* {{\rm ev}}_{B_1}^* {\bf h}^{(1)} \times \text{\rm d-}\deg {{\rm ev}}_{B_2}^* {\bf h}^{(2)} \nonumber \\
 d_2 & = & \text{\rm d-}\deg  \pi_{B_1}^* {{\rm ev}}_{B_1}^* {\bf h}^{(1)} \times \text{\rm d-}\deg {{\rm ev}}_{B_2}^* {\bf h}^{(2)} \nonumber \\
 d_3 & = & (\mu(B_1) + k^{(1)} - n) (\mu(B_2) + k^{(2)} -n). \nonumber 
 \end{eqnarray} 
The thrid inequality follows from \cite[Proposition 10.24]{springer}, \cite[Proposition 2.3 and the subsequent paragraphs]{ono}.   
 \par
 Note that $d_1+d_2 = (\text{rel}\dim \pi_{B_2}) \times \text{\rm d-}\deg {\bf h}^{(2)} = (\mu(B_1) + k^{(1)} - n) \text{\rm d-}\deg {\bf h}^{(2)}.$ 
 Then a simple calculation shows that 
 $\delta_0 + d_1 + d_2 + d_3 \equiv \delta_1 + (n+1) \deg' {\bf h}^{(2)}   \mod 2$.  
 Hence the proof is  complete.
 \end{proof}

 Recall that we have the forgetful map \eqref{forgettoannul} from the moduli space of holomorphic annuli to the moduli of domain curves.  
\begin{equation}
\frak{forget}_{\rm sc} : 
\mathcal M_{1}^{\text{ann}}
(\vec\kappa^{(1)},\vec\kappa^{(2)};\vec p,\vec q;B)
\to 
\mathcal M_{(1,1);0}. \nonumber
\end{equation}

Forgetting the orientation issue, $\frak{forget}_{\rm sc}^{-1}([\Sigma_1])$ is the union of fiber products 
${\mathcal M}_{1; k^{(1)}+1}(B_1) \times_X {\mathcal M}_{1; k^{(2)}+1}(B_2)$ for various 
$B_1,B_2$.  (See Figure \ref{Figure10-2}.)

We use the following convention.  
\begin{equation}\label{oriholann} T_{[u, \Sigma, \vec{z}_1 \vec{z}_2]} \mathcal M_{1}^{\text{ann}}(\vec\kappa^{(1)},\vec\kappa^{(2)};\vec p,\vec q;B) 
= T_u {\rm Hol}(B) \oplus T_{[\Sigma, \vec{z}_1, \vec{z}_2]} {\mathcal M}_{(k^{(1)} +1, k^{(2)}+1), 0}.
\end{equation} 
For the moduli space of the domain curves, i.e., annlui with boundary marked points, 
\begin{equation}\label{oriannuli}
T_{[\Sigma, \vec{z}_1, \vec{z}_2]} {\mathcal M}_{(k^{(1)} +1, k^{(2)}+1), 0}  
=  T_{[\Sigma, z_1(0), z_2(0)]}  {\mathcal M}_{(1,1), 0} \oplus {\R}^{k^{(1)}} \oplus {\R}^{k^{(2)}}, 
\end{equation}
where $\Sigma$ is an annlus, ${\R}^{k^{(1)}}$, (resp. ${\R}^{k^{(2)}}$) corresponds to 
$\bigoplus_{i=1}^{k^{(1)}} T_{z_1(i)} \partial_1 \Sigma$, (resp. $\bigoplus_{i=1}^{k^{(2)}} T_{z_2(i)} \partial_2 \Sigma$).

The orientations are related as follows.   

\begin{prop}\label{pp, ann}
\begin{equation}
(-1)^{\delta_2} {\mathcal M}^*_{1; k^{(1)}+1}(B_1) \times_X {\mathcal M}^*_{1; k^{(2)}+1}(B_2) \subset \frak{forget}_{\rm sc}^{-1}([\Sigma_1]).
\end{equation}
Here 
$$\delta_2=k^{(1)}(\mu(B_2)-n) + \mu(B_1)(\mu(B_2)-n) + k^{(2)}.$$
\end{prop}

\begin{proof} 
As for as  orientation issues is concerned,\footnote{We discuss orientations in the sense of Kuranishi structure and  the issue is on the liniearization level.} it is enough to discuss the  spaces  $\overset{\circ}{\mathcal M^*}_{1; k^{(s)}+1}(B_s)$, $s=1,2$.  
We may assume that the interior marked point $z^+_0 = 0 $, the origin,  and the 0-th boundary marked point $z_s(0) =1$.  
(See Remark \ref{rem211}.)
Since it is enough to study the orientation on the tangent space level, we may regard 
 the moduli space $\overset{\circ}{\mathcal M^*}_{1, k^{(s)}}(B_s)$ as  an open subset of 
$\text{Hol}(B_s) \times (\partial D)^{k^{(s)}}$, where $\text{Hol}(B_s)$\index[syindex]{HolB@$\text{Hol}(B)$} is the space of holomorphic polygons such that the parameter encoding the positions of switching marked points and their images are fixed.
The parameters encoding the positions of the diagonal marked points are also fixed.
Then we have 
$$\aligned
&\overset{\circ}{\mathcal M^*}_{1, k^{(1)}+1}(B_1) \times_X \overset{\circ}{\mathcal M^*}_{1, k^{(2)}+1}(B_2) \\
&= (-1)^{d_4} (\text{Hol}(B_1) \times_X \text{Hol}(B_2)) \times (\partial D)^{k^{(1)}} \times (\partial D)^{k^{(2)}}),
\endaligned
$$
where $d_4 = k^{(1)} (\mu(B_2) - n)$.  

Inserting $\widetilde{\lambda}^{(s)}_i$ at switching marked points, the boundary condition is equipped with a spin structure, which stably trivializes 
the boundary condition as a fixed Lagrangian space (unique up to homotopy).  
For $v \in \text{Hol}(B_s)$, the spin structure determines the orientation for 
$T_v \text{Hol}(B_s)  \oplus \mathcal{V}^{(s)}_0 \oplus \dots \oplus \mathcal{V}^{(s)}_{k^{(s)}}$.  
Here $\mathcal{V}^{(s)}_j$ is the index of the operator 
(\ref{orientationcomplex}) at the $j$-th switching marked point. (See \cite[begining of Section 3]{ono}.)

For the space of holomorphic annuli, the data $\widetilde{\lambda}^{(s)}_i$, $ i=0, \dots, k^{(s)}$ and $s=1,2$, determine the orientation for 
$T_u \text{Hol}^{\rm ann}(B) \oplus  \mathcal{V}^{(1)}_0 \oplus \dots \oplus \mathcal{V}^{(1)}_{k^{(1)}} \oplus  \mathcal{V}^{(2)}_0 \oplus \dots \oplus \mathcal{V}^{(2)}_{k^{(2)}}$.  
Here $\text{Hol}^{\rm ann}(B)$ is defined in a similar way as $\text{Hol}(B)$.

Now we consider the case of $\frak{forget}_{\rm sc}^{-1}([\Sigma_1])$.  
The (interior) node corresponds to the coincidence condition for two holomorphic disks at interior marked points.  
The orientation of 
$T_u \text{Hol}^{\rm ann}(B) \oplus  \mathcal{V}^{(1)}_0 \oplus \dots \oplus \mathcal{V}^{(1)}_{k^{(1)}} \oplus  \mathcal{V}^{(2)}_0 \oplus \dots \oplus \mathcal{V}^{(2)}_{k^{(2)}}$ 
and the fiber product orientation of 
$$(T_{u^{(1)}} \text{Hol}(B_1) \oplus \mathcal{V}^{(1)}_0 \oplus \dots \oplus \mathcal{V}^{(1)}_{k^{(1)}}) 
{}_{{\rm ev}^+_0}\times_{{\rm ev}^+_0} (T_{u^{(2)}} \text{Hol}(B_2) \oplus  \mathcal{V}^{(2)}_0 \oplus \dots \oplus \mathcal{V}^{(2)}_{k^{(2)}})$$ are the same.\footnote{
We remind the reader that we are using the convention to require $z_0 = 1$ and $z^+=0$ here.
See Remark \ref{rem211}.}
The latter differs from the fiber product orientation of 
$$(T_{u^{(1)}}\text{Hol}(B_1) {}_{{\rm ev}^+_0}\times_{{\rm ev}^+_0}
T_{u^{(2)}} \text{Hol}(B_2)) \oplus  \mathcal{V}^{(1)}_0 \oplus \dots \oplus \mathcal{V}^{(1)}_{k^{(1)}} \oplus  \mathcal{V}^{(2)}_0 \oplus \dots \oplus \mathcal{V}^{(2)}_{k^{(2)}}$$
 by $(-1)^{d_5}$, where $d_5=\mu(B_1) \times (\mu(B_2) - n).$  
 The evaluation maps ${\rm ev}^+_0$ is defined by ${\rm ev}^+_0(u^{(1)}) =u^{(1)}(0)$,  ${\rm ev}^+_0(u^{(2)}) =u^{(2)}(0)$.
 The fiber product is taken over $X$.
 
 For marked points, the forgetful map sends $[(z_1(0), \dots z_1(k^{(1)}), z_2(0), \dots, z_2(k^{(2)}))]$ to $[(z_1(0), z_2(0))]$.  
 In addition to permutation of $z_1(1), \dots, z_1({k^{(1)}})$ with $z_2(0)$, the automorphism group $S^1$ also jump over 
 $z_1(1), \dots, z_1(k^{(1)}), z_2(1), \dots, z_2(k^{(2)})$ in this process.   (See (\ref{oriannuli}).)
 Thus the fiber orientation is given by 
 $$(-1)^{d_6} (z_1(1), \dots, z_1(k^{(1)}), z_2(1), \dots, z_2(k^{(2)})),$$
 where $d_6=k^{(1)} + (k^{(1)}+k^{(2)}) \equiv k^{(2)} \mod 2$.\footnote{In other words
 our convention of the orientation of $\mathcal M_{k^{(1)}+1,k^{(2)}+1}^{\rm ann}$ 
 is $\R^{k^{(1)}+k^{(2)}+2}/\R$, where $\R^{k^{(1)}+k^{(2)}+2}$ encodes the positions of 
 $z_1(0),z_1(1),\dots,z_1(k^{(1)}),z_2(0),z_2(1),\dots,z_2(k^{(2)})$.}
 Set $\delta_2=d_4+d_5+d_6$.  Then we obtain Proposition \ref{pp, ann}.
 \end{proof}

Note that
$$\epsilon({\bf h}^{(1)}, {\bf h}^{(2)}) = \epsilon({\bf h}^{(1)}) + \epsilon({\bf h}^{(2)}) + (\mu(B_1)+k^{(1)}+1) 
\sum_{j=0}^{k^{(2)}} (\text{\rm d-}\deg h^{(2)}_j -1).$$
Using this equality, Propositions \ref{pp} and \ref{pp, ann} imply the following: 

\begin{prop}\label{ppSigma1}
\begin{equation}
\aligned 
& 
(-1)^{(n+1)\deg' {\bf h}^{(2)}} \sum_{B_1+B_2=B} \int_X {\mathfrak p}_{B_1} ({\bf h}^{(1)}) \wedge {\mathfrak p}_{B_2} ({\bf h}^{(2)}) \\
= & \sum_{B_1+B_2=B} (-1)^{\epsilon({\bf h}^{(1)}, {\bf h}^{(2)})} \int_{\frak{forget}_{\rm sc}^{-1}([\Sigma_1])} {\rm ev}^*({\bf h}^{(1)}, {\bf h}^{(2)}).
\endaligned
\end{equation}
\end{prop}

Let $\gamma:[0,1] \to {\mathcal M}_{(1,1);0}$ be a path with $\gamma(0) = [\Sigma']$ and $\gamma(1)=[\Sigma'']$.  
We define 

$$
\aligned
{\mathcal Z}^{\gamma}_B ({\bf h}^{(1)}, {\bf h}^{(2)}) = &\sum_{B_1+B_2=B} (-1)^{\epsilon({\bf h}^{(1)}, {\bf h}^{(2)})}   \\
&\int_{\gamma \times_{\mathcal M_{(1,1);0}}{\mathcal M_{1}^{\text{ann}}
(\vec\kappa^{(1)},\vec\kappa^{(2)};\vec p,\vec q;B)}} {\rm ev}^*({\bf h}^{(1)}, {\bf h}^{(2)}),
\endaligned
$$
$$
{\mathcal Z}^{[\Sigma]}_B  ({\bf h}^{(1)}, {\bf h}^{(2)}) = \sum_{B_1+B_2=B} (-1)^{\epsilon({\bf h}^{(1)}, {\bf h}^{(2)})} 
\int_{\frak{forget}_{\rm sc}^{-1}[\Sigma]} {\rm ev}^*({\bf h}^{(1)}, {\bf h}^{(2)}).
$$
Then we have the following:

\begin{prop}\label{pairingdelta}
\begin{equation}\label{Zdelta}
\aligned
&{\mathcal Z}^{[\Sigma'']}_B ({\bf h}^{(1)}, {\bf h}^{(2)}) - {\mathcal Z}^{[\Sigma']}_B({\bf h}^{(1)}, {\bf h}^{(2)}) \\
&= {\mathcal Z}^{\gamma}_B (\delta{\bf h}^{(1)}, {\bf h}^{(2)}) + (-1)^{\deg' {\bf h}^{(1)}} {\mathcal Z}^{\gamma}_B ({\bf h}^{(1)}, \delta {\bf h}^{(2)}).  
\endaligned
\end{equation}
\end{prop}

\begin{proof}
We have the following identity 
\begin{equation}\label{relp} 
d {\mathfrak p} ({\bf h}) + {\mathfrak p} \circ \widehat{\mathfrak m} ({\bf h}) = 0.  
\end{equation}
For the pairing ${\mathcal Z}^{[\Sigma_1]}$, by the Stokes formula, Proposition \ref{pp} and \eqref{relp}, we have 
we have the following equality 
$$
{\mathcal Z}^{[\Sigma_1]}_B (\delta {\bf h}^{(1)}, {\bf h}^{(2)}) + (-1)^{\deg' {\bf h}^{(1)} } {\mathcal Z}^{[\Sigma_1]} _B({\bf h}^{(1)}, \delta {\bf h}^{(2)}) = 0.
$$
This is \eqref{Zdelta} in the case where $\gamma$ is the constant path at $[\Sigma_1]$.  
Our concern in this subsection is the issue on signs.  It verifies \eqref{Zdelta}.
\end{proof}

\begin{figure}[h]
\centering
\includegraphics[scale=0.25]{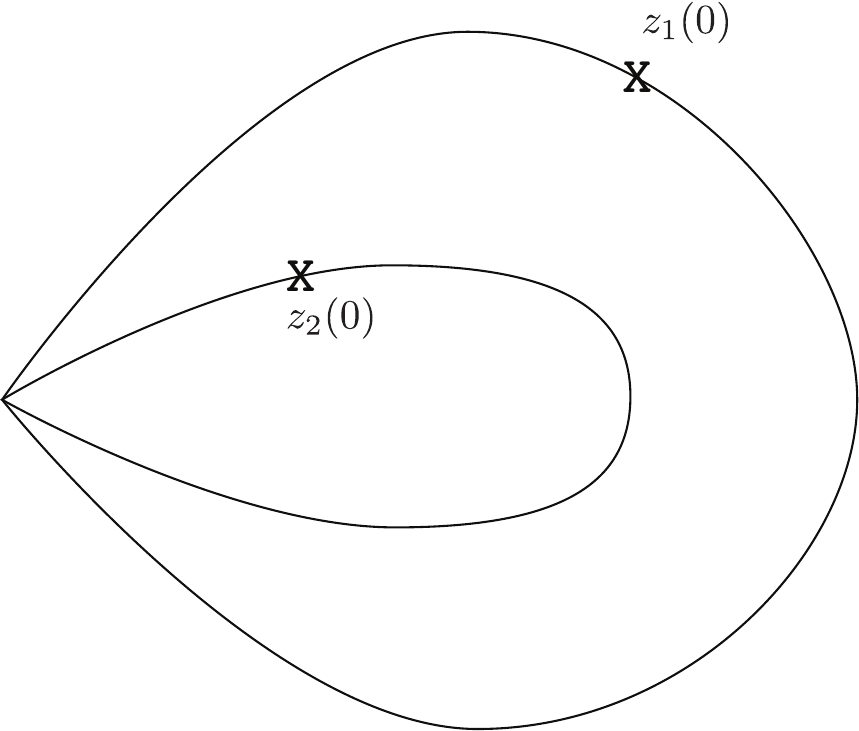}
\caption{$\Sigma_3$}
\label{Figureaaa}
\end{figure}

Consider a degenerate annulus $\Sigma_3$ with a boundary node and boundary marked points $z_1(0)$, $z_2(0)$, see Figure \ref{Figureaaa}.  
Here we fix the boundary marked points $z_1(0)$ and $z_2(0)$ and the boundary node up to automorphisms.  

\begin{figure}[h]
\centering
\includegraphics[scale=0.25]{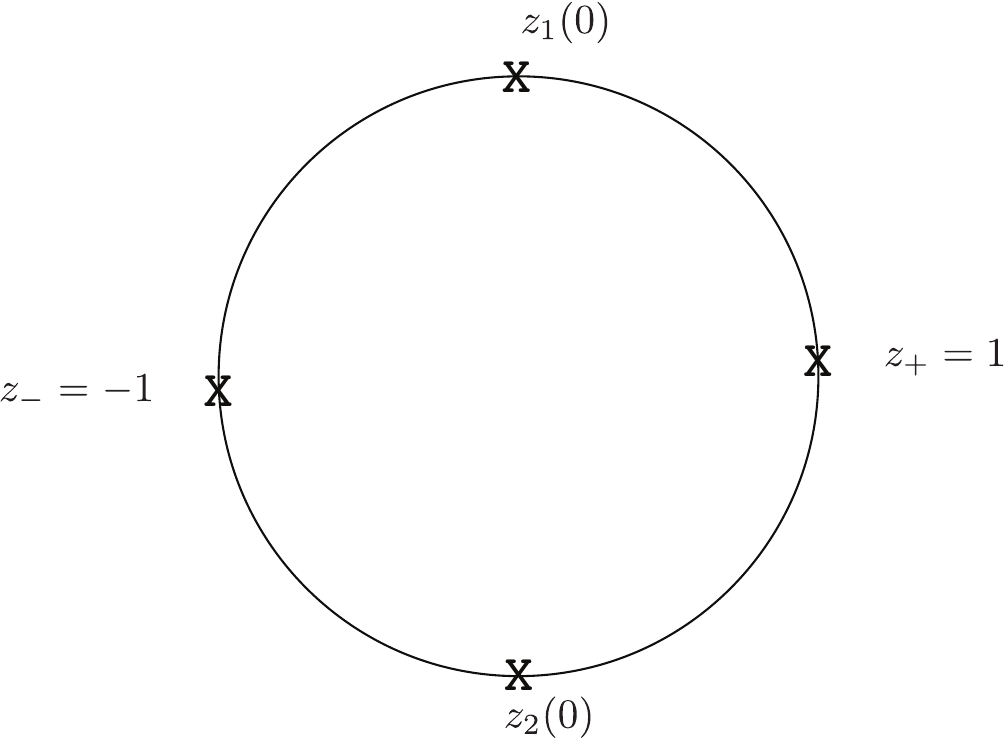}
\caption{$\tilde \Sigma_3$}
\label{Figureiii}
\end{figure}

In other words, 
after taking a normalization $\widetilde{\Sigma}_3$ of the boundary node, we obtain a unit disk with four boundary marked points, $z_+, z_1(0), z_-, z_2(0)$, where 
$z_+$ and $z_-$ are the inverse image of the boundary nodes via the normalization.  The cross ratio of  $z_+, z_1(0), z_-, z_2(0)$ is a fixed number depending 
only on $\Sigma_3$.  We may fix $z_{\pm}=\pm 1$.  See Figure \ref{Figureiii}.

\begin{figure}[h]
\centering
\includegraphics[scale=0.25]{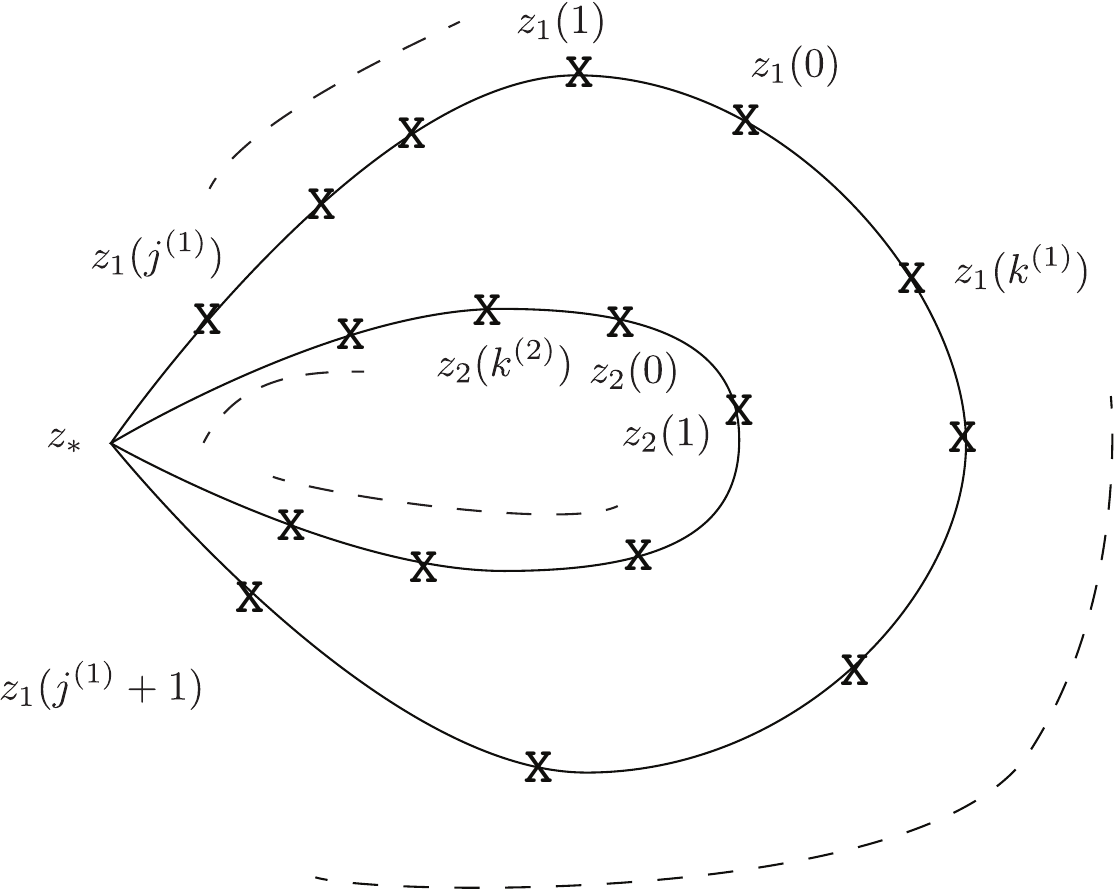}
\caption{$\mathcal{C}({j^{(1)}, j^{(2)}})$.}
\label{Figureuuu}
\end{figure}

Then $\frak{forget}_{\rm sc}^{-1}([\Sigma_3])$ is decomposed into the components $\mathcal{C}(j^{(1)}, j^{(2)})$ consisting of elements of  $\frak{forget}_{\rm sc}^{-1}([\Sigma_3])$, 
whose domain curves have boundary nodes $z_*$ as the intersection of  $\overline{z_1(j^{(1)}) z_1(j^{(1)}+1)}$ and $\overline{z_2(j^{(2)}) z_2(j^{(2)}+1)}$, 
see Figure \ref{Figureuuu}.  
Let $p \in L_{\kappa^{(1)}_{j^{(1)}}} \cap U_{\kappa^{(2)}_{j^{(2)}}}$, which is the image $u(z_*)$ of the boundary node under the holomorphic map $u$.  

\begin{rem}\label{rem2177}
We define the orientation of $\mathcal{C}({j^{(1)}, j^{(2)}})$ as follows.
 %\marginpar{Remark added.  KF 2025 Dec.
%Modified KF 2026 Feb.}
We take the index of the operator $D_u\overline{\partial}$
where $u$ is the map part of $\mathcal{C}({j^{(1)}, j^{(2)}})$. Using the relative spin structure and 
the orientation of the index ${\mathcal V}_{p_i}$ associated to the switching marked points (including $p$), 
the index of $D_u\overline{\partial}$ is oriented as follows. 
We consider the sum $D_u\overline{\partial} \oplus \bigoplus \mathcal {V}_{p_i}$.
Here we order $\bigoplus \mathcal {V}_{p_i}$  according to the way they appear in the configuration, as in \eqref{orientationpolygon}.  
This sum is identified with the index of an appropriate Cauchy-Rieman operator with Lagrangian boundary condition. 
The relative spin structures of $L_{\kappa}$, $U_{\kappa'}$ induces trivialization of the boundary condition and so the sum 
is oriented. Therefore together with the orientation of $\mathcal {V}_{p_i}$ we have obtained the orientation of the index of $D_u\overline{\partial}$.
(Note that when we regard $u$ as a holomorphic strip with switching points on the boundaries, we give an orientation of $\text{Index~} D_u\overline{\partial} = \text{Hol}^{\text{strip}}(B)$ 
as in \eqref{2111form}, which is used in Proposition \ref{comparisonwithconnorb}.   This is only an intermediate step.  
Then we adopt the orientation convention as in \eqref{orientationpolygon} for a later argument.  The difference between them is taken care of in Propostion  \ref{aspolygon}.)

We consider 
\begin{equation}\label{newnew2110}
\text{Index~} D_u\overline{\partial} \oplus \left(\bigoplus \R^{k^{(1)} + k^{(2)}+2}\right)/\R.
\end{equation}
Here $\R^{k^{(1)} + k^{(2)}+2}$ is the parameter encoding the positions of the points 
$$z_1(0), z_1(1),\dots,z_1(k^{(1)}),z_2(0),z_2(1),\dots,z_2(k^{(2)})$$ and $\R$ is the translation of them by the $\R$ action.  
There is a map
\begin{equation}\label{newnew2113}
\left(\bigoplus \R^{k^{(1)} + k^{(2)}+2}\right)/\R \to \R
\end{equation}
which associates $[z_1(0),z_2(0)]$ to $[z_1(0), z_1(1),\dots,z_1(k^{(1)}),z_2(0),z_2(1),\dots,z_2(k^{(2)})]$.

The tangent space of $\mathcal{C}({j^{(1)}, j^{(2)}})$ is identified to the direct sum of
$\text{Index~} D_u\overline{\partial}$ and the fiber of 
(\ref{newnew2113}) and so is oriented.
(We remark that the orientation of the fiber of (\ref{newnew2113}) is $(-1)^{k^{(2)}}$ times the standard orientation 
of $\R^{k^{(1)} + k^{(2)}}$.
This is the orientation which is consistent with the proof of Proposition \ref{pp, ann}.)
\end{rem}

We can regard an element of $\mathcal{C}({j^{(1)}, j^{(2)}})$ as a holomorphic map $u:{\mathbb R} \times [0,1] \to X$ of homology class $B$ with boundary marked points $z_1(i) \in {\mathbb R} \times \{1\}$, 
$z_2(i) \in {\mathbb R} \times \{0 \}$ with ${\rm Re} (z_1(0)) - {\rm Re}(z_2(0))$ is fixed such that $\lim_{\tau \to \pm \infty} u(\tau, t) = p$ and 
$u(z_1(i)) \in L_{\kappa^{(1)}_{i-1}} \cap L_{\kappa^{(1)}_i}$ for $i = 0, \dots, k^{(1)}$, $u(z_2(i)) \in U_{\kappa^{(2)}_{i-1}} \cap L_{\kappa^{(2)}_i}$ for 
$i = 0, \dots, k^{(2)}$\footnote{By the abuse of notation,  $L_{\kappa^{(1)}_{-1}}=L_{\kappa^{(1)}_{k^{(1)}}}$ and $U_{\kappa^{(2)}_{-1}}=U_{\kappa^{(2)}_{k^{(2)}}}$.} 
with the following boundary condition:
$$
u(\tau, 1) \in 
\begin{cases}
L_{\kappa^{(1)}_{j^{(1)}}} & \text{ for } \tau > {\rm Re}(z_1(j^{(1)}+1)) \\
L_{\kappa^{(1)}_i} & \text{ for  } {\rm Re}(z_1(i)) < \tau < {\rm Re}(z_1(i-1)), i=j^{(1)}+2, \dots, k^{(1)}, 0, \dots, j^{(1)}\\
L_{\kappa^{(1)}_{j^{(1)}}} & \text{ for } \tau <{\rm Re}(z_{(1)}(j^{(1)}))
\end{cases}
$$
and
$$
u(\tau, 0) \in 
\begin{cases}
U_{\kappa^{(2)}_{j^{(2)}}} & \text{ for } \tau < {\rm Re} (z_2(j^{(2)}+1)) \\
U_{\kappa^{(2)}_i} & \text{ for  } {\rm Re}(z_2(i-1)) < \tau < {\rm Re}(z_2(i)), i=j^{(2)}+2, \dots, k^{(2)}, 0, \dots, j^{(2)}\\
U_{\kappa^{(2)}_{j^{(2)}}} & \text{ for } \tau > {\rm Re}(z_{(2)}(j^{(2)})), 
\end{cases}
$$
see Figure \ref{Figureeee}.

Denote by $\overset{\circ}{\mathcal M}_{k^{(1)}+1, k^{(2)}+1}^{}(p,p;\{L_{\kappa^{(1)}_i}\}, \{U_{\kappa^{(2)}_i}\};B;[\widetilde{\Sigma}_3], j^{(1)}, j^{(2)})$ the moduli space of 
such holomorphic curves. 
Its stable map compactification is denoted by the symbol ${\mathcal M}_{k^{(1)}+1, k^{(2)}+1}(p,p;\{L_{\kappa^{(1)}_i}\}, \{U_{\kappa^{(2)}_i}\};B;[\widetilde{\Sigma}_3], j^{(1)}, j^{(2)})$.   
%We use the orientation by \eqref{2111form} below and Remark \ref{rem2177}.  
%Namely, we regard an element in ${\mathcal M}_{k^{(1)}+1, k^{(2)}+1}(p,p;\{L_{\kappa^{(1)}_i}\}, \{U_{\kappa^{(2)}_i}\};B;[\widetilde{\Sigma}_3], j^{(1)}, j^{(2)})$ as 
%a stable strip with switching points on the boundary.  
The orientation convention \eqref{2111form} is used in the intermediate step of the computation in Proposition \ref{comparisonwithconnorb}.  See also Remark \ref{rem2177}. 
We use the orientation convention  in \eqref{orientationpolygon} from Proposition \ref{aspolygon} on.  
\par
We put also
\begin{equation}
\mu(s;j^{(s)},k^{(s)}) =  \left(\sum_{j=j^{(s)}+1}^{k^{(1)}} \mu^{(s)}_j\right) \left(\sum_{j=0}^{j^{(s)}} \mu^{(s)}_j\right) 
\end{equation}
for $j^{(s)} \le k^{(s)}$, $s=1,2$.

\begin{prop}\label{comparisonwithconnorb}
The orientation of $\mathcal{C}({j^{(1)}, j^{(2)}}) \subset \frak{forget}_{\rm sc}^{-1}([\Sigma_3])$ and the orientation of 
${\mathcal M}_{k^{(1)}+1, k^{(2)}+1}(p,p;\{L_{\kappa^{(1)}_i}\}, \{U_{\kappa^{(2)}_i}\};B;[\widetilde{\Sigma}_3], j^{(1)}, j^{(2)})$ differs by 
\begin{equation}\label{exchangingpoints}
\aligned
& (j^{(1)}+1)(k^{(1)}-j^{(1)}) + (j^{(2)} + 1) (k^{(2)} -j^{(2)}) + \mu(p) 
+ n(n-1)/2 \\ &+ \mu(1;j^{(1)},k^{(1)}) + \mu(2;j^{(2)},k^{(2)}). 
\endaligned
\end{equation}
\end{prop} 

\begin{proof} 
As in the proof of Proposition \ref{pp, ann}, the boundary conditions along $\partial^{(s)} D$ with inserting $\widetilde{\lambda}^{(s)}_i$, $i=0, \dots, k^{(2)}$ and $s=1,2$,  
is trivialized by spin structures.  
We may assume that they are constant Lagrangian subspaces $T_p L_{\kappa^{(1)}_{j^{(1)}}}$ and $T_p U_{\kappa^{(2)}_{j^{(2)}}}$, respectively.  
Then the sign of the intersection number of $T_p L_{\kappa^{(1)}_{j^{(1)}}}$ and $T_p U_{\kappa^{(2)}_{j^{(2)}}}$ in $T_p X$ is $(-1)^{\mu(p)+n(n-1)/2}$.
Therefore, when we regard $u$ as a holomorphic annulus, we find that the orientation of 
$$
T_u \text{Hol}^{\rm ann}(B) \oplus \mathcal{V}^{(1)}_0 \oplus \dots \oplus \mathcal{V}^{(1)}_{k^{(1)}} \oplus \mathcal{V}^{(2)}_0 \oplus \dots \oplus \mathcal{V}^{(2)}_{k^{(2)}}
$$ 
is 
$(-1)^{\mu(p)+n(n-1)/2}$ times complex orientation.  

Now we regard the map $u$ (omitting the source curve from notation) as an element of the moduli space ${\mathcal M}_{k^{(1)}+1, k^{(2)}+1}^{}(p,p;\{L_{\kappa^{(1)}_i}\}, \{U_{\kappa^{(2)}_i}\};B; [\widetilde{\Sigma}_3], j^{(1)}, j^{(2)})$.  
Pick a Lagrangian path $\widetilde{\lambda}(p)$ 
from $T_p L_{\kappa^{(1)}_{j^{(1)}}}$ to  $T_p U_{\kappa^{(2)}_{j^{(2)}}}$ in $T_p X$ with a spin structure relative to the boundary.  
Denote by $\mathcal{V}(p)$ the index of $\overline{\partial}$ with boundary condition $\widetilde{\lambda(p)}$ 
and by $\Theta(p)$ its determinant line bundle.  

Since the trivialization of the boundary inserting $\widetilde{\lambda}(p)$ is homotopic to $\widetilde{\lambda}(p)$ relative to $T_{u(+\infty, 1)} L_{\kappa^{(1)}_{j^{(1)}}}$ 
and $T_{u(+\infty,0)} U_{\kappa^{(2)}_{j^{(2)}}}$ and the trivialization given by $\widetilde{\lambda}(p)$ are homotopic, 
we find that the orientation of 
\begin{eqnarray}\label{2111form}
T_u \text{Hol}^{\rm strip} (B) & \oplus & \mathcal{V}^{(1)}_{j^{(1)}+1} \oplus \dots \oplus \mathcal{V}^{(1)}_{k^{(1)}} \oplus \mathcal{V}^{(1)}_0 \oplus \dots \oplus \mathcal{V}^{(1)}_{j^{(1)}} \nonumber \\
& \oplus & \mathcal{V}^{(2)}_{j^{(2)}+1} \oplus \dots \oplus \mathcal{V}^{(2)}_{k^{(2)}} \oplus \mathcal{V}^{(2)}_0 \oplus \dots \oplus \mathcal{V}^{(2)}_{j^{(2)}} \oplus \mathcal{V}(p)  
\end{eqnarray}
is canonically identified with the one of $\mathcal{V}(p)$.   
Here $\text{Hol}^{\rm strip} (B)$ is defined in a similar way as $\text{Hol}^{\rm ann} (B)$.

Comparing these, the difference of orientation of $T_u \text{Hol}^{\rm ann}(B)$ and $T_u \text{Hol}^{\rm strip} (B)$ is $(-1)^{d_7}$, where 
$$d_7 = \mu(p) +  n(n-1)/2 + \mu(1;j^{(1)},k^{(1)}) + \mu(2;j^{(2)},k^{(2)}).
$$

Now we compare the orientations on the configurations of marked points.    
On the annulus, we have the moduli space of boundary marked points on a fixed annulus.  
$$\{(z_1(0), \dots, z_1(k^{(1)}), z_2(0), \dots, z_2(k^{(2)})\}/S^1.$$  Here $S^1$ acts on the annulus by rotations, See Figure \ref{Figureooo}.  

\begin{figure}[h]
\centering
\includegraphics[scale=0.35]{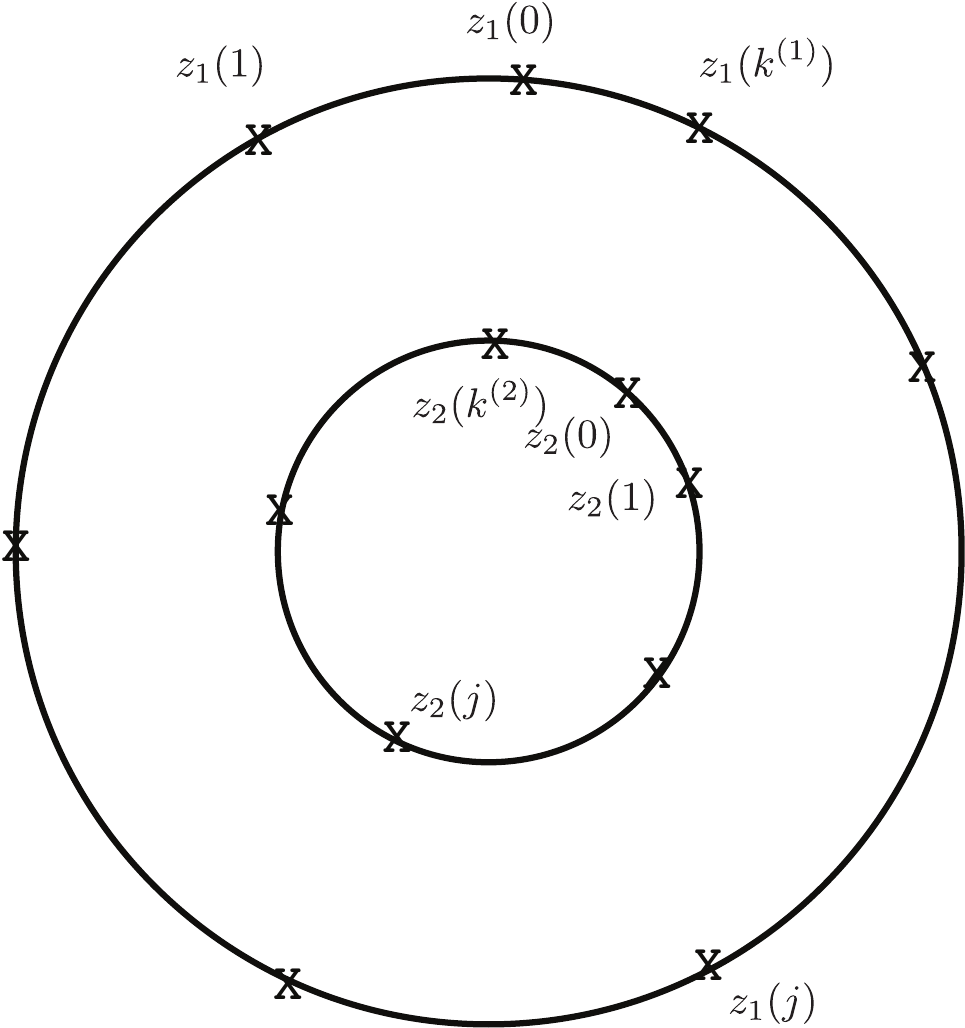}
\caption{Annulus with marked point}
\label{Figureooo}
\end{figure}

On the strip, we have the moduli space of boundary marked points 
\begin{equation}
\aligned
\{& z_1(j^{(1)}+1), \dots, z_1(k^{(1)}), z_1(0), \dots, z_1(j^{(1)}),  \\
& z_2(j^{(2)}+1), \dots, z_2(k^{(2)}), z_2(0), \dots,z_2(j^{(2)}))\}/{\mathbb R}. \nonumber
\endaligned
\end{equation}
Here ${\mathbb R}$ acts on the strip by translations, see Figure  \ref{FigureKakaka}.  
\begin{figure}[h]
\centering
\includegraphics[scale=0.35]{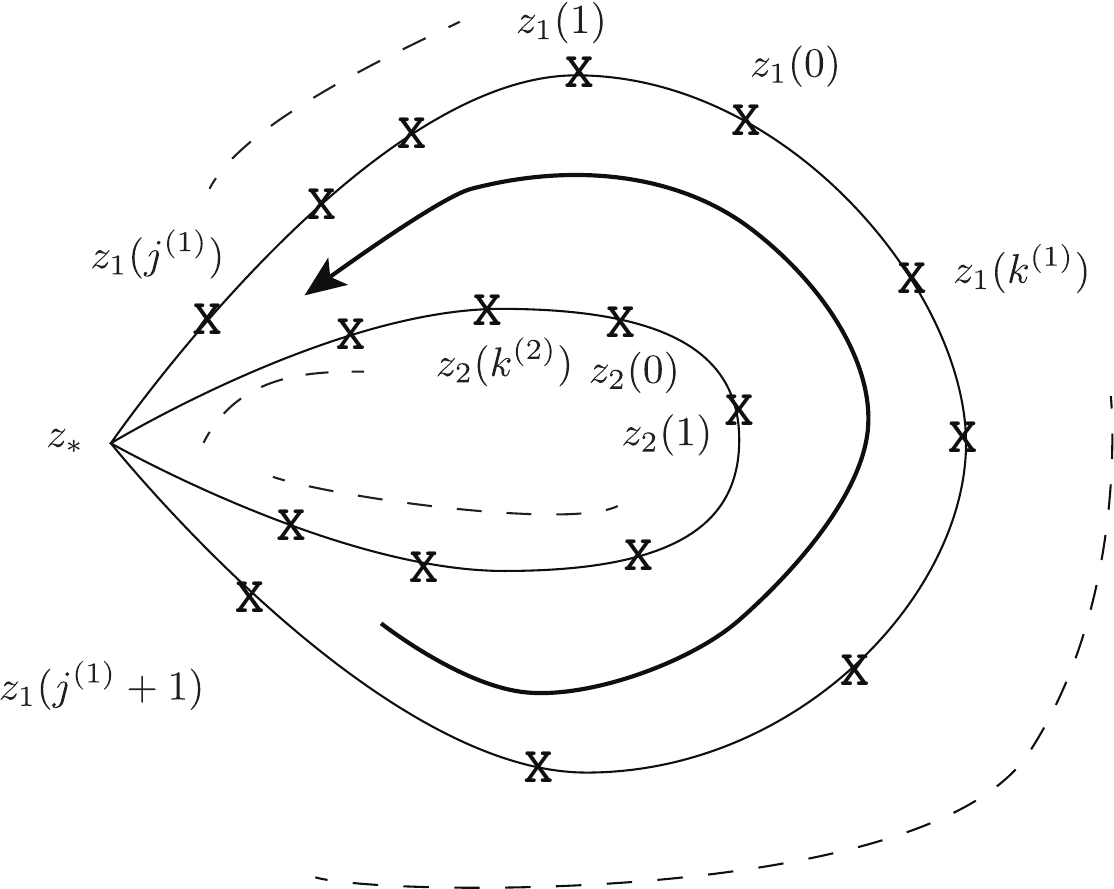}
\caption{$\R$ action on $\mathcal{C}({j^{(1)}, j^{(2)}})$}
\label{FigureKakaka}
\end{figure}
The difference of orientations is $(-1)^{d_8}$, 
with $d_8= (j^{(1)} +1)(k^{(1)}-j^{(1)}) + (j^{(2)}+1) (k^{(2)} -j^{(2)})$. 
Hence $(-1)^{d_7 + d_8}$ is also the difference of orientations of $\mathcal{C}({j^{(1)}, j^{(2)}}) \subset \frak{forget}^{-1}_{\rm sc}([\Sigma_3])$ and 
${\mathcal M}_{k^{(1)}+1, k^{(2)}+1}^{}(p,p;\{L_{\kappa^{(1)}_i}\}, \{U_{\kappa^{(2)}_i}\};B;[\widetilde{\Sigma}_3], j^{(1)}, j^{(2)})$.  
Hence the proof.  
\end{proof} 

The moduli  space ${\mathcal M}_{k^{(1)}+1, k^{(2)}+1}^{}(p,p;\{L_{\kappa^{(1)}_i}\}, \{U_{\kappa^{(2)}_i}\};B;[\widetilde{\Sigma}_3], j^{(1)}, j^{(2)})$ 
can be regarded 
as the moduli  space of holomorphic polygons (see Figure $\ref{Figureeee}$) and also as the moduli space of connecting orbits.

\begin{prop}\label{aspolygon}
The orientation as the moduli  space of holomorphic polygons differs from the orientation as connecting orbits by 
\begin{equation}\label{connorbpoly}
\mu(p) \cdot \mu(B_2) + k^{(1)} + 1.
\end{equation}
\end{prop}
\begin{figure}[h]
\centering
\includegraphics[scale=0.35]{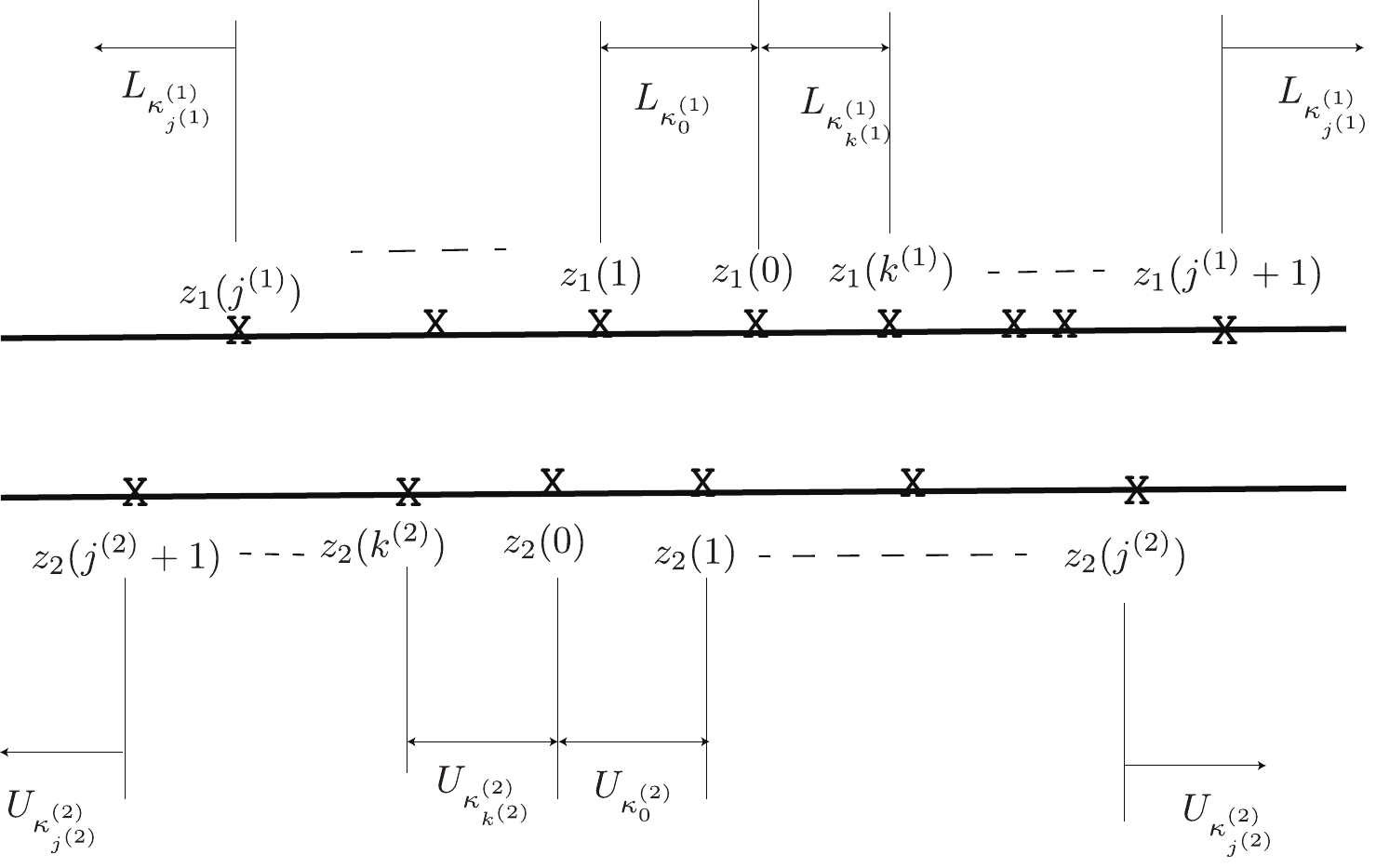}
\caption{${\mathcal M}_{k^{(1)}+1, k^{(2)}+1}^{}(p,p;\{L_{\kappa^{(1)}_i}\}, \{U_{\kappa^{(2)}_i}\};B;[\widetilde{\Sigma}_3], j_1, j_2)$}
\label{Figureeee}
\end{figure}

\begin{proof}
There are two sources for the sign difference.  
One comes from the orientation of the index of the linearization operaotor $D_u\overline{\partial}$ and the other from 
the order of the boundary marked points.  
As for the former, the orientation of  the index of $D_u\overline{\partial}$, we insert $\widetilde{\lambda}^{(1)}_{j^{(1)}+1}, 
\dots, \widetilde{\lambda}^{(1)}_{k^{(1)}}, \widetilde{\lambda}^{(1)}_0 \dots, \widetilde{\lambda}^{(1)}_{j^{(1)}}, \widetilde{\lambda}(p), 
\widetilde{\lambda}^{(2)}_{j^{(1)}+1}, \dots, \widetilde{\lambda}^{(2)}_{k^{(2)}}, \widetilde{\lambda}^{(2)}_0 \dots, \widetilde{\lambda}^{(2)}_{j^{(2)}}$ 
at marked points at boundary marked points.  The relative spin structures and $\widetilde{\lambda}^{(s)}_i$, $s=1,2$, $j=0, \dots, k^{(s)}$ determines the orientation of 
\begin{eqnarray}\label{orientationpolygon}
T_u \text{Hol}^{\rm strip} (B) & \oplus & \mathcal{V}^{(1)}_{j^{(1)}+1} \oplus \dots \oplus \mathcal{V}^{(1)}_{k^{(1)}} \oplus \mathcal{V}^{(1)}_0 \oplus \dots \oplus \mathcal{V}^{(1)}_{j^{(1)}} \oplus \mathcal{V}(p) \nonumber \\
& \oplus & \mathcal{V}^{(2)}_{j^{(2)}+1} \oplus \dots \oplus \mathcal{V}^{(2)}_{k^{(2)}} \oplus \mathcal{V}^{(2)}_0 \oplus \dots \oplus \mathcal{V}^{(2)}_{j^{(2)}}.   
\end{eqnarray}
Comparing with \eqref{2111form}, we find that 
the difference of orientation is $(-1)^{d_9}$, where $d_9 = \mu(p) (\sum_{j=0}^{k^{(2)}} \mu^{(2)}_j)$.  

As for the latter, we note that the strip is identified with the unit disk with $\{\pm 1\}$ removed.  
Then the action in Proposition \ref{comparisonwithconnorb} fixes $\pm 1$, which is $0$-th marked point and $k^{(1)}+2$ nd marked points 
in the cyclic order.  In  \cite[Convention 8.3.1]{fooo092}, the ${\mathbb R}$-action fixes $0$-th and $1$-st marked points.  The orientations 
of the quotient spaces by these ${\mathbb R}$-actions differes by $(-1)^{d_{10}}$, where $d_{10}=k^{(1)}+1$.  
Recall that $\mu(B_2)=(\sum_{i=0}^{k^{(2)}} \mu^{(2)}_i)$.  
Combining these effects, we obtain Proposition \ref{aspolygon}.  
\end{proof}

We write 
\begin{equation}
\vec \kappa^{(s)}_{\langle j \rangle} = (\kappa^{(s)}_{j+1}, \dots, \kappa^{(s)}_{k^{(s)}}, \kappa^{(s)}_0, \dots, \kappa^{(s)}_j),
\end{equation}
Here $\vec\kappa^{(s)}$ is determined by $B_s$.
Note in the formulation here $k^{(s)} + 1$ is the total number of boundary marked point, that is, 
both switching and boundary marked points. 
We put 
\begin{equation}
{\bf h}^{(s)}_{\langle j \rangle} = (h^{(s)}_{j+1}, \dots, h^{(s)}_{k^{(s)}}, h^{(s)}_0, \dots, h^{(s)}_j).
\end{equation}
To shorten the formula we put
\begin{equation}
\aligned
\deg'(s;j^{(s)}, k^{(s)}) & = \left(\sum_{j=j^{(s)}+1}^{k^{(s)}}  \deg'  h^{(s)}_j \right) \left(\sum_{j=0}^{j^{(s)}}  \deg'  h^{(s)}_j \right), \\
\text{\rm d-}\deg(s;j^{(s)},k^{(s)}) &=\left(\sum_{j=j^{(s)}+1}^{k^{(s)}} \text{\rm d-}\deg h^{(s)}_j\right) \left(\sum_{j=0}^{j^{(s)}} \text{\rm d-}\deg h^{(s)}_j\right). \
\endaligned
\end{equation}

Combining Propositions \ref{ppSigma1}, \ref{pairingdelta}, \ref{comparisonwithconnorb} and \ref{aspolygon}, we have:
\begin{prop}\label{21.9}  
For Hochschild cycles ${\bf h}^{(1)}$, ${\bf h}^{(2)}$, we have 
\begin{equation}\label{form21111}
\aligned 
& 
(-1)^{(n+1) \deg' {\bf h}^{(2)} } \sum_{B_1 + B_2 = B} \int_X {\mathfrak p}_{B_1} ({\bf h}^{(1)}) \wedge {\mathfrak p}_{B_2} ({\bf h}^{(2)}) \\
= & \sum_{j^{(1)}, j^{(2)}, p}  (-1)^{\delta_3} \int {\rm ev}^*({\bf h}^{(1)}_{\langle j^{(1)} \rangle}, 1_p, {\bf h}^{(2)}_{\langle  j^{(2)} \rangle}), 
\endaligned
\end{equation}
where the  integration is over 
${\mathcal M}_{k^{(1)}+1, k^{(2)}+1}^{}(p,p;\{L_{\kappa^{(1)}_i}\}, \{U_{\kappa^{(2)}_i}\};B;[\widetilde{\Sigma}_3], j^{(1)}, j^{(2)})$
and
\begin{eqnarray}\label{pppoly}
\delta_3 & = & \epsilon ({\bf h}^{(1)}, 1_p, {\bf h}^{(2)}) + \deg' {\bf h}^{(1)} \deg'  1_p  \nonumber \\
& & +  (j^{(1)}+1)(k^{(1)}-j^{(1)}) + (j^{(2)} + 1) (k^{(2)} -j^{(2)})\nonumber \\
& & + \text{\rm d-}\deg(1;j^{(1)},k^{(1)}) + \text{\rm d-}\deg(2;j^{(2)},k^{(2)}) \nonumber \\
&& +  \mu(1;j^{(1)},k^{(1)}) + \mu(2;j^{(2)},k^{(2)}) + \mu(B) \nonumber \\
& & + n(n-1)/2 + \mu(p) -1. 
\end{eqnarray} 
Here $1_p$ is the constant function $1$ at the intersection point $p$.  The summation is taken over $0 \leq j^{(1)} \leq k^{(1)}, 0\leq j^{(2)} \leq k^{(2)}$, intersection point $p$ of 
$L_{\kappa^{(1)}_{j^{(1)}}}$ and $U_{\kappa^{(2)}_{j^{(2)}}}$.
%and all equivalence classes $B$ of polygonal maps.  
\end{prop}

\begin{proof}
For Hochschild cycles, Proposition \ref{pairingdelta} implies that the integration does not depend on $[\Sigma] \in {\mathcal M}_{1,1;0}$.  

Note that 
$${\rm ev}^*({\bf h}^{(1)}_{\langle j^{(1)} \rangle} \times 1_p \times  {\bf h}^{(2)}_{\langle j^{(2)} \rangle}) 
= (-1)^{d_{11}} {\rm ev}^* ({\bf h}^{(1)} \times 1_p \times {\bf h}^{(2)}),$$ 
where 
$
d_{11}= \text{\rm d-}\deg(1;j^{(1)},k^{(1)}) + \text{\rm d-}\deg(2;j^{(2)},k^{(2)}).
$
Here ${\rm ev}$  on the left hand side is the evaluation map at marked points 
$$z_1(j^{(1)}+1), \dots, z_1(k^{(1)}), z_1(0), \dots, z_1(j^{(1)}), z_-, 
z_2(j^{(2)}+1), \dots, z_2(k^{(2)}), z_2(0), \dots,z_2(j^{(2)})$$ 
and ${\rm ev}$ on the right hand side is 
the evaluation map at marked points 
$$z_1(0), \dots, z_1(k^{(1)}), z_-, z_2(0), \dots, z_2(k^{(2)}).$$ 
Then the equality \ref{pppoly} holds with $\delta_3 = \epsilon({\bf h}^{(1)}, {\bf h}^{(2)}) + \eqref{exchangingpoints} + \eqref{connorbpoly} +d_{11}$.  

When the integral is non-zero, we have 
$$\text{\rm d-}\deg {\bf h}^{(1)} + \text{\rm d-}\deg {\bf h}^{(2)} = \mu(B_1)+\mu(B_2) + k^{(1)}+k^{(2)}.$$
Then a simple calculation tells that 
\begin{eqnarray}
& & \epsilon ({\bf h}^{(1)}, {\bf h}^{(2)}) + \mu(p) \cdot \mu(B_2) + k^{(1)} + 1 \nonumber \\
& = & \epsilon({\bf h}^{(1)}, 1_p, {\bf h}^{(2)}) + \deg' {\bf h}^{(1)}  \deg' 1_p  - (\mu(B_1) + \mu(B_2))  -1. \nonumber
\end{eqnarray}
Hence we have \eqref{pppoly}.  
\end{proof}
We set
\begin{eqnarray}
{\mathrm M}^{}({\bf h}^{(1)}_{\langle j^{(1)}\rangle}, 1_p, {\bf h}^{(2)}_{\langle j^{(2)} \rangle};B;[\widetilde{\Sigma}_3]) 
 =(-1)^{\epsilon({\bf h}^{(1)}_{\langle j^{(1)} \rangle}, 1_p, {\bf h}^{(2)}_{\langle j^{(2)} \rangle})} \int{\rm ev}^*({\bf h}^{(1)}_{\langle j^{(1)} \rangle} \times 1_p \times {\bf h}^{(2)}_{\langle j^{(2)} \rangle}),  \nonumber
\end{eqnarray}
where integration is taken over ${\mathcal M}_{k^{(1)}+1, k^{(2)}+1}^{}(p,p;\{L_{\kappa^{(1)}_i}\}, \{U_{\kappa^{(2)}_i}\};B;[\widetilde{\Sigma}_3], j^{(1)}, j^{(2)})$.

\begin{prop}\label{Moperation}
\begin{equation}
\aligned
&(-1)^{(n+1) \deg'  {\bf h}^{(2)} } \sum_{B_1 + B_2 = B}  \int_X {\mathfrak p}_{B_1} ({\bf h}^{(1)}) \wedge {\mathfrak p}_{B_2} ({\bf h}^{(2)}) \\
&=  \sum_{p} (-1)^{\delta_4} {\mathrm M}^{}({\bf h}^{(1)}_{\langle j^{(1)}\rangle}, 1_p, {\bf h}^{(2)}_{\langle j^{(2)} \rangle};B;[\widetilde{\Sigma}_3]),
\endaligned
\end{equation}
where 
$$
\aligned
\delta_4= \deg' {\bf h}^{(1)}  \deg' 1_p  + \deg' (1;j^{(1)},k^{(1)}) + \deg' (2;j^{(2)}, k^{(2)}) + \mu(B) + \mu(p) + n(n-1)/2.
\endaligned$$  
\end{prop}

For the proof, we prepare the following lemma. 
Denote the transposition of ${\bf h}^{(s)}$ in $j$-th and $j+1$-st entries by
$${\bf h}^{(s)}_{(j,j+1)} = (h^{(s)}_0, \dots, h^{(s)}_{j+1}, h^{(s)}_j, \dots, h^{(s)}_{k^{(s)}}).$$

\begin{lem}\label{transposition}
$$
\aligned
&\epsilon({\bf h}^{(1)}_{(j,j+1)}, 1_p, {\bf h}^{(2)}) - \epsilon({\bf h}^{(1)}, 1_p, {\bf h}^{(2)}) + \deg' h^{(1)}_j \deg' h^{(1)}_{j+1}  \\
&\equiv \mu^{(1)}_j  \mu^{(1)}_{j+1} + \text{\rm d-}\deg h^{(1)}_j \cdot \text{\rm d-}\deg h^{(1)}_{j+1} + 1.
\endaligned$$
$$
\aligned
&\epsilon({\bf h}^{(1)}, 1_p, {\bf h}^{(2)}_{(j,j+1))}) - \epsilon({\bf h}^{(1)}, 1_p, {\bf h}^{(2)}) + \deg' h^{(2)}_j \deg' h^{(2)}_{j+1}  \\
&\equiv \mu^{(2)}_j  \mu^{(2)}_{j+1} + \text{\rm d-}\deg h^{(2)}_j \cdot \text{\rm d-}\deg h^{(2)}_{j+1} + 1.
\endaligned$$
\end{lem}

The proof is a direct computation using the definitions of $\epsilon$-factor and the shifted degree of $h^{(s)}_i$.  

\begin{proof}[Proof of Proposition \ref{Moperation}]
From Lemma \ref{transposition}, we find the next formula  (\ref{effecttransposition}).
Now we have:
\begin{eqnarray}\label{effecttransposition} 
&  & \epsilon({\bf h}^{(1)}, 1_p, {\bf h}^{(2)}) + (j^{(1)}+1)(k^{(1)} - j^{(1)}) + (j^{(2)}+1)(k^{(2)} - j^{(2)})) \nonumber \\ 
& & + \mu(1;j^{(1)},k^{(1)}) + \mu(2;j^{(2)},k^{(2)}) \nonumber \\
& & + \text{\rm d-}\deg(1;j^{(1)},k^{(1)}) + \text{\rm d-}\deg(2;j^{(2)},k^{(2)}) \nonumber \\ 
& \equiv & \epsilon({\bf h}^{(1)}_{\langle j^{(1)} \rangle}, 1_p, {\bf h}^{(2)}_{\langle j^{(2)} \rangle}) + 
\deg'(1;j^{(1)}, k^{(1)}) + \deg'(2;j^{(2)}, k^{(2)})
\end{eqnarray}

\begin{comment}
Note that 
$${\rm ev}^*({\bf h}^{(1)}_{\langle j^{(1)} \rangle} \times 1_p \times  {\bf h}^{(2)}_{\langle j^{(2)} \rangle}) 
= (-1)^{d_{11}} {\rm ev}^* ({\bf h}^{(1)} \times 1_p \times {\bf h}^{(2)}),$$ 
where 
$
d_{11}= \text{\rm d-}\deg(1;j^{(1)},k^{(1)}) + \text{\rm d-}\deg(2;j^{(2)},k^{(2)}).
$
Here the evaluation map on the left hand side are taken at 
$$z_1(j^{(1)}+1), \dots, z_1(k^{(1)}), z_1(0), \dots, z_1(j^{(1)}), z_-, z_2(j^{(2)}+1), \dots, z_2(k^{(2)}), z_2(0), \dots, z_2(j^{(2)}).$$ 
The evaluation map on the right hand side is taken at 
$$z_1(0), \dots, z_1(k^{(1)}), z_-, z_2(0), \dots, z_2(k^{(2)}).$$  
\end{comment}
Using \eqref{pppoly}, \eqref{effecttransposition}, we obtain Proposition  \ref{Moperation}.  
\end{proof}

\begin{figure}[h]
\centering
\includegraphics[scale=0.35]{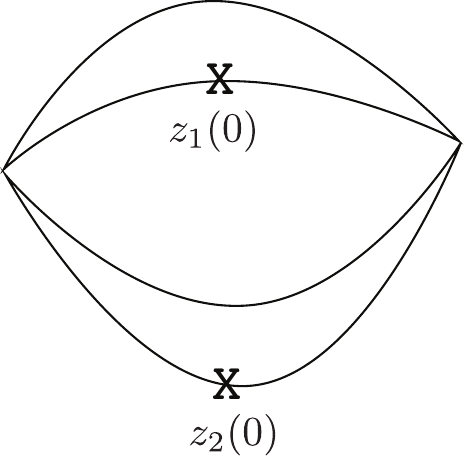}
\caption{$\Sigma_4$}
\label{FigureKikiki}
\end{figure}

\begin{figure}[h]
\centering
\includegraphics[scale=0.35]{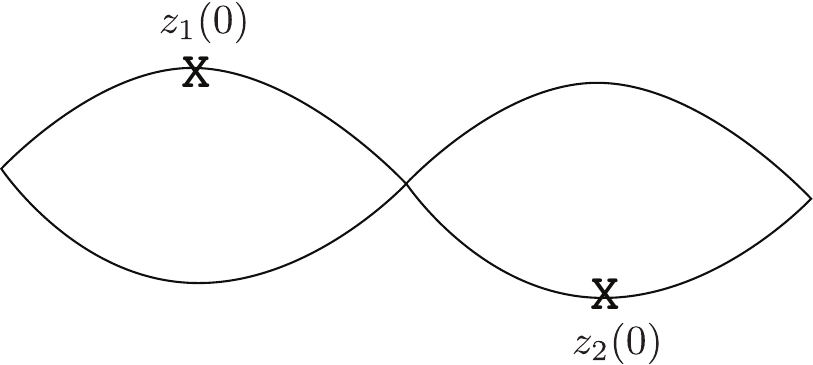}
\caption{$\tilde\Sigma_4$}
\label{FigureKukuku}
\end{figure}

\begin{figure}[h]
\centering
\includegraphics[scale=0.6]{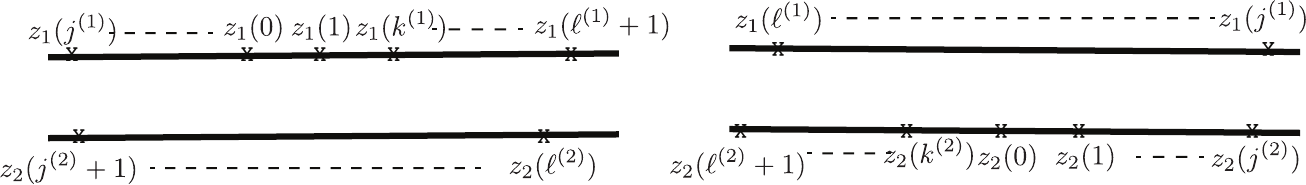}
\caption{$\tilde\Sigma_4$ with marked points.}
\label{FigureKekekeono}
\end{figure}

Now, we consider $[\Sigma_4] \in {\mathcal M}_{1,1; 0}$ with two boundary nodes, which is a degeneration of $\Sigma_3$, see Figure \ref{FigureKikiki}.  
The corresponding degeneration for $\widetilde{\Sigma}_3$ is denoted by $\widetilde{\Sigma}_4$, see Figures \ref{FigureKukuku}
and \ref{FigureKekekeono}.   


We set 
$${\bf h}^{(s)}_{j, \ell} = h^{(s)}_j \times \dots \times h^{(s)}_{\ell}.$$

For a generator $p \in CF(L_{\kappa^{(1)}_{j^{(1)}}}, U_{\kappa^{(2)}_{j^{(2)}}})$,  which is represented by an intersection point 
$p \in L_{\kappa^{(1)}_{j^{(1)}}} \cap U_{\kappa^{(2)} _{j^{(2)}}}$ equipped with $\widetilde{\lambda}_p$ and $o(\lambda)$, we denote by 
$p^{\vee} \in CF(U_{\kappa^{(2)}_{j^{(2)}}}, L_{\kappa^{(1)}_{j^{(1)}}})$ such that 
$$\langle p, p^{\vee} \rangle_{\rm cyc} = 1.$$

We define the number
${\mathrm M}^{}(\frak b^{\otimes \ell};{\bf h}^{(1)}_{\langle j^{(1)}\rangle}, 1_p, {\bf h}^{(2)}_{\langle j^{(2)} \rangle};B;[\widetilde{\Sigma}_3])$
by including $\ell$ bulk classes $\frak b$ to 
${\mathrm M}({\bf h}^{(1)}_{\langle j^{(1)}\rangle}, 1_p, {\bf h}^{(2)}_{\langle j^{(2)} \rangle};B;[\widetilde{\Sigma}_3])$
in an obvious way.  (Including bulk class, which is of even degree, does not affect the sign issue that is our concern in this 
section.)

We next prove Lemma \ref{lem2114} below.

The moduli space 
${\mathcal M}_{k^{(1)}+1, k^{(2)}+1}^{}(p,p;\{L_{\kappa^{(1)}_i}\}, \{U_{\kappa^{(2)}_i}\};B;[\widetilde{\Sigma}_4], j^{(1)}, j^{(2)})$  is defined in the same way as 
${\mathcal M}_{k^{(1)}+1, k^{(2)}+1}^{}(p,p;\{L_{\kappa^{(1)}_i}\}, \{U_{\kappa^{(2)}_i}\};B;[\widetilde{\Sigma}_3], j^{(1)}, j^{(2)})$ but replacing $\widetilde{\Sigma}_3$ by $\widetilde{\Sigma}_4$.
We consider two orientations on it.
The first orientation is one as the boundary of the moduli space 
of pseudo-holomorphic disks of homology class $B$.
We call it the boundary orientation.
The other orientation is one we have been using so far. See Remark \ref{rem2177}.

\begin{lem}\label{lem2114}
The difference between the two orientations is $(-1)^{\mu(B)}$.   %\marginpar{The lemma added.  KF 2025 Dec.,  revised 2026 Feb KO}
\end{lem}
\begin{proof}
Recall our orientation convention for $\mathfrak{forget}_{sc}([\Sigma])$ given at \eqref{oriholann}, \eqref{oriannuli}.  
In Propositions \ref{aspolygon}, \ref{21.9} and \ref{Moperation}, we give orientations on the moduli spaces in the same way.  
Namely, the moduli space ${\mathcal M}_4$, which is the stable curve compactification of $\{(D^2; z_+, z_1(0), z_-, z_2(0))\}/{\text{PSL}}(2,{\R})$ 
is put between $T_u \text{Hol}(B)$ and $\R^{k^{(1)}} \oplus \R^{k^{(2)}}$.  

Hence the boundary orientation of the modul space of pseudo-holomorphic polygons  of holomogy class $B$ and the orientation on which we use so far
differs by the effect of switching $T_u{\rm Hol}(B)$, which is $\mu(B)$-dimensional and the moduli space ${\mathcal M}_4$, 
which is one-dimensional.  Hence we obtain the lemma.    
\end{proof}

\begin{comment}
We orient $T_u \text{Hol}^{\rm strip} (B)$ using (\ref{2111form}).
Here $u$ is a map from $\Sigma_3$.
The boundary orientation is given as the limit of the orientation of the fiber of the map
when $\Sigma_3$ goes to $\Sigma_4$.
\begin{eqnarray}\label{2441form}
&&T_u \text{Hol}^{\rm strip} (B) \oplus  \mathbb R^{k^{(1)}+ k^{(2)}} \oplus \mathbb R^4/\text{psl}(2;\mathbb R) \nonumber
\\
& \oplus & \mathcal{V}^{(1)}_{j^{(1)}+1} \oplus \dots \oplus \mathcal{V}^{(1)}_{k^{(1)}} \oplus \mathcal{V}^{(1)}_0 \oplus \dots \oplus \mathcal{V}^{(1)}_{j^{(1)}} 
\nonumber \\
& \oplus & \mathcal{V}^{(2)}_{j^{(2)}+1} \oplus \dots \oplus \mathcal{V}^{(2)}_{k^{(2)}} \oplus \mathcal{V}^{(2)}_0 \oplus \dots \oplus \mathcal{V}^{(2)}_{j^{(2)}} \oplus \mathcal{V}(p)  \nonumber
\\
&& \to \mathbb R.
\end{eqnarray}
Here $\mathbb R^4$ is the parameter moving $z_1(0), z_2(0)$ and $(\pm \infty,0)$.
The automorphism   $\text{PSL}(2;\R)$ acts on the moduli space of marked disks from the {\it right}.
The quotient $\mathbb R^4/\text{psl}(2;\mathbb R)$ is $\R$ which is identified with the parameter 
to move $\Sigma_3$ to $\Sigma_4$. It then becomes the outer normal 
to ${\mathcal M}_{k^{(1)}+1, k^{(2)}+1}^{}(p,p;\{L_{\kappa^{(1)}_i}\}, \{U_{\kappa^{(2)}_i}\};B;[\widetilde{\Sigma}_4], j^{(1)}, j^{(2)})$
as the boundary of the moduli space 
of pseudo holomorphic disks of homology class $B$.
This parameter $\R$ exchanges with 
$T_u \text{Hol}^{\rm strip} (B) \oplus  \mathbb R^{k^{(1)}+ k^{(2)}}$ and so 
we have sign correction $(-1)^{k^{(1)}+k^{(2)}+\mu(B)}$.
\end{comment}

Using 
\begin{eqnarray}
& & {\mathfrak m}({\bf h}^{(1)}_{\ell^{(1)}+1, k^{(1)}}, {\bf h}^{(1)}_{0,j^{(1)}}, 1_p, {\bf h}^{(2)}_{j^{(2)}+1, \ell^{(2)}}) \nonumber \\
& = & \sum_q \langle {\mathfrak m}({\bf h}^{(1)}_{\ell^{(1)}+1, k^{(1)}}, {\bf h}^{(1)}_{0,j^{(1)}}, 1_p, 
{\bf h}^{(2)}_{j^{(2)}+1, \ell^{(2)}}), 1_{q^{\vee}} \rangle_{\text{\rm cyc}} \cdot 1_q,  \nonumber 
\end{eqnarray}
we have 
\begin{equation}
\aligned
& \langle {\mathfrak m}({\bf h}^{(1)}_{j^{(1)}+1, \ell^{(1)}}, {\mathfrak m}({\bf h}^{(1)}_{\ell^{(1)}+1, k^{(1)}}, {\bf h}^{(1)}_{0,j^{(1)}}, 1_p, 
{\bf h}^{(2)}_{j^{(2)}+1, \ell^{(2)}}), {\bf h}^{(2)}_{\ell^{(2)}+1, k^{(2)}}, {\bf h}^{(2)}_{0, j^{(2)}}), 1_{p^{\vee}} \rangle_{\text{\rm cyc}},  \\
& = \sum_{q} \langle {\mathfrak m}({\bf h}^{(1)}_{\ell^{(1)}+1, k^{(1)}}, {\bf h}^{(1)}_{0,j^{(1)}}, 1_p, 
{\bf h}^{(2)}_{j^{(2)}+1, \ell^{(2)}}), 1_{q^{\vee}} \rangle_{\text{\rm cyc}}   \\
& \qquad\times\langle {\mathfrak m}({\bf h}^{(1)}_{j^{(1)}+1, \ell^{(1)}}, 1_q, 
{\bf h}^{(2)}_{\ell^{(2)}+1, k^{(2)}}, {\bf h}^{(2)}_{0, j^{(2)}}), 1_{p^{\vee}} \rangle_{\text{\rm cyc}}.   
\endaligned 
\nonumber
\end{equation}

Now using these equalities, we have 

\begin{prop}
\begin{eqnarray}\label{form2110}
& & 
\sum_{B,\ell} T^{\omega(B)} \frac{\exp((\frak b_2,\theta_L) \cap B)}{\ell!} 
(-1)^{\delta_4}{\mathrm M}^{} 
(\frak b^{\otimes \ell};{\bf h}^{(1)}_ {\langle j^{(1)} \rangle}, 1_p, {\bf h}^{(2)}_{\langle j^{(2)} \rangle};B;[\widetilde{\Sigma}_3]) \nonumber 
  \\ 
& = &  \sum_{\ell^{(1)}, \ell^{(2)}, p, q} (-1)^{\delta_4 + \delta_5}   \langle {\mathfrak m}({\bf h}^{(1)}_{\ell^{(1)}+1, k^{(1)}}, {\bf h}^{(1)}_{0,j^{(1)}}, 1_p, 
{\bf h}^{(2)}_{j^{(2)}+1, \ell^{(2)}}), 1_{q^{\vee}} \rangle_{\text{\rm cyc}}  \nonumber \\
& & \qquad \qquad \times \langle {\mathfrak m}({\bf h}^{(1)}_{j^{(1)}+1, \ell^{(1)}}, 1_q,  
{\bf h}^{(2)}_{\ell^{(2)}+1, k^{(2)}}, {\bf h}^{(2)}_{0, j^{(2)}}), 1_{p^{\vee}} \rangle_{\text{\rm cyc}},  
\end{eqnarray}
where $\delta_5 = \sum_{j=j^{(1)}+1}^{\ell^{(1)}} \deg' h^{(1)}_j  + \mu(B)$.  
\end{prop}

 %\marginpar{(In Definition of $Z$, Definition \ref{defnZ}, the sign is as in the beginning of the proof of Theorem \ref{ZisPD} below. )
%KO 2025 Sep. }

\begin{proof}[Proof of Theorem \ref{ZisPD}]
\begin{eqnarray}
& & (-1)^{(n+1)\deg' {\bf h}^{(2)}}  \langle {\mathfrak p} ({\bf h}^{(1)}), {\mathfrak p} ({\bf h}^{(2)}) \rangle_{\text{PD}_X} \nonumber \\
& = & (-1)^{(n+1)\deg' {\bf h}^{(2)}} \sum_B T^{\omega(B)} \sum_{B_1 + B_2 = B} \int_X {\mathfrak p}_{B_1} ({\bf h}^{(1)}) \wedge {\mathfrak p}_{B_2} ({\bf h}^{(2)}) \nonumber \\
& = &  \sum_{\ell^{(1)}, \ell^{(2)}, p, q} (-1)^{\delta_4 + \delta_5}   \langle {\mathfrak m}({\bf h}^{(1)}_{\ell^{(1)}+1, k^{(1)}}, {\bf h}^{(1)}_{0,j^{(1)}}, 1_p, 
{\bf h}^{(2)}_{j^{(2)}+1, \ell^{(2)}}), 1_{q^{\vee}} \rangle_{\text{\rm cyc}}  \nonumber \\
& & \qquad \qquad \times\langle {\mathfrak m}({\bf h}^{(1)}_{j^{(1)}+1, \ell^{(1)}}, 1_q,  
{\bf h}^{(2)}_{\ell^{(2)}+1, k^{(2)}}, {\bf h}^{(2)}_{0, j^{(2)}}), 1_{p^{\vee}} \rangle_{\text{\rm cyc}},  
\nonumber
\end{eqnarray}
where 
$$
\aligned
\delta =  \delta_4 + \delta_5 = \mu(p) &+ \sum_{j=j^{(1)}+1}^{\ell^{(1)}} \deg' h^{(1)}_j  + \deg' {\bf h}^{(1)}   \deg' 1_p    
\\
&+ \deg'(1;j^{(1)}, k^{(1)}) + \deg'(2;j^{(2)}, k^{(2)}) + n(n-1)/2. 
\endaligned
$$
Recall that 
\begin{eqnarray}\label{Zsign}
& & (-1)^{(n+1) \deg' {\bf h}^{(2)} }Z({\bf h}^{(1)}, {\bf h}^{(2)})  \nonumber \\
& = & 
\sum_{\ell^{(1)}, \ell^{(2)}, p, q} (-1)^{\delta_4 + \delta_5}   \langle {\mathfrak m}({\bf h}^{(1)}_{\ell^{(1)}+1, k^{(1)}}, {\bf h}^{(1)}_{0,j^{(1)}}, 1_p, 
{\bf h}^{(2)}_{j^{(2)}+1, \ell^{(2)}}), 1_{q^{\vee}} \rangle  \nonumber \\
& & \qquad \qquad \times \langle {\mathfrak m}({\bf h}^{(1)}_{j^{(1)}+1, \ell^{(1)}}, 1_q,  
{\bf h}^{(2)}_{\ell^{(2)}+1, k^{(2)}}, {\bf h}^{(2)}_{0, j^{(2)}}), 1_{p^{\vee}} \rangle_{\text{\rm cyc}},  
\end{eqnarray}
where 
$$
\aligned
\delta' = \mu(p) &+ \sum_{j=j^{(1)}+1}^{\ell^{(1)}} \deg' h^{(1)}_j  + \deg' {\bf h}^{(1)}   \deg' 1_p    
\\
&+ \deg'(1;j^{(1)}, k^{(1)}) + \deg'(2;j^{(2)}, k^{(2)}). 
\endaligned
$$
Namely, $\delta' = \delta - n(n-1)/2$.  
Therefore we have 
$$\langle {\mathfrak p} ({\bf h}^{(1)}), {\mathfrak p} ({\bf h}^{(2)}) \rangle_{\text{PD}_X} = (-1)^{n(n-1)/2} Z({\bf h}^{(1)}, {\bf h}^{(2)}).$$ 

We will give an interpretation of the sign $(-1)^{\delta'}$ in terms of the Koszul rule.  
In fact, using the notation in Subsection \ref{inproZ} and putting ${\bf x} = {\bf h}^{(1)}$, ${\bf y} = {\bf h}^{(2)}$,
 %\marginpar{The rest 
%of this subsection is added.  KF 2026 Jan.}
$f^1= 1_p$,
we have
$$
\aligned
&{\bf h}^{(1)}_{\ell^{(1)}+1, k^{(1)}} \otimes {\bf h}^{(1)}_{0,j^{(1)}}  =
{\bf x}_{c_1}^{(H;2,1)},
\quad 
{\bf h}^{(1)}_{j^{(1)}+1, \ell^{(1)}} = {\bf x}_{c_1}^{(H;2,2)},
\\
&{\bf h}^{(2)}_{\ell^{(2)}+1, k^{(2)}} \otimes {\bf h}^{(2)}_{0,j^{(2)}}  =
{\bf y}_{c_2}^{(H;2,1)},
\quad 
{\bf h}^{(2)}_{j^{(2)}+1, \ell^{(2)}} = {\bf y}_{c_2}^{(H;2,2)}.
\endaligned
$$
Here the indices $c_1$ (appearing in Sweedler's notation)  (resp. $c_2$) is $(j^{(1)},k^{(1)},\ell^{(1)})$
(resp. $(j^{(2)},k^{(2)},\ell^{(2)})$)
and 
\begin{equation}\label{form2124pre}
\aligned
\delta'=  &\deg f^1  + \deg \frak m\deg'  {\bf x}_{c_1}^{(H;2,2)} + \deg' {\bf x} \deg' f^1
\\
&+ \deg'(1;j^{(1)}, k^{(1)}) + \deg'(2;j^{(2)}, k^{(2)}).
\endaligned
\end{equation}
Here $\deg'(1;j^{(1)}, k^{(1)})$ and $\deg'(2;j^{(2)}, k^{(2)})$ are Koszul signs in shifted degrees for the permutations 
${\bf x} \to {\bf x}_{c_1}^{(H;2,2)}, {\bf x}_{c_1}^{(H;2,1)}$ and ${\bf y} \to {\bf y}_{c_2}^{(H;2,2)}, {\bf y}_{c_2}^{(H;2,2)}$, 
respectively.  

The right hand side of (\ref{Zsign}) is equal to 
\begin{equation}\label{form2124}
\sum_{c_1,c_2, f_1} (-1)^{\delta'} \langle {\mathfrak m}({\bf x}_{c_1}^{(H;2,2)},
 {\mathfrak m}({\bf x}_{c_1}^{(H;2,1)}, f^1, 
{\bf y}_{c_2}^{(H;2,2)}),{\bf y}_{c_2}^{(H;2,1)}), f^{1\vee}\rangle_{\text{cyc}}.
\end{equation}
This is Koszul sign by the permutation
\begin{equation}\label{form212422}
\aligned
&\langle \rangle,\mathfrak m,\mathfrak m,f^1,f^{1\vee},{\bf x}, {\bf y}  \\
&\mapsto \langle \rangle,\mathfrak m, {\bf x}_{c_1}^{(H;2,2)},\mathfrak m,{\bf x}_{c_1}^{(H;2,1)},f^1,{\bf y}_{c_2}^{(H;2,2)},{\bf y}_{c_2}^{(H;2,1)},f^{1\vee}
\endaligned
\end{equation}
plus $\deg f^1$. Here we use the fact that the total shifted degree of the object
${\bf x}_{c_1}^{(H;2,1)},{\bf x}_{c_1}^{(H;2,2)},{\bf y}_{c_2}^{(H;2,1)},{\bf y}_{c_2}^{(H;2,2)}$ is even when 
(\ref{form2124}) is nonzero.

We further modify the formula so that it becomes the form we discussed in Subsection \ref{inproZ}.

For a generator $f^2 \in CF(L_{\kappa^{(1)}_{\ell^{(1)}}}, U_{\kappa^{(2)}_{\ell^{(2)}}})$,  which is represented by an intersection point 
$p \in L_{\kappa^{(1)}_{\ell^{(1)}}} \cap U_{\kappa^{(2)} _{\ell^{(2)}}}$, we denote by 
$f^{2\vee} \in CF(U_{\kappa^{(2)}_{\ell^{(2)}}}, L_{\kappa^{(1)}_{\ell^{(1)}}})$ such that 
$$\langle f^2, f^{2\vee} \rangle_{\rm cyc} = 1.$$

Then 
\begin{equation}\label{2126form}
\aligned
&\langle {\mathfrak m}({\bf x}_{c_1}^{(H;2,2)}, {\mathfrak m}({\bf x}_{c_1}^{(H;2,1)}, f^1 ,{\bf y}_{c_2}^{(H;2,2)}), {\bf y}_{c_2}^{(H;2,1)}), f^{1\vee} \rangle_{\text{cyc}}  \\
& =  \langle   {\mathfrak m}({\bf x}_{c_1}^{(H;2,1)}, f^1 ,{\bf y}_{c_2}^{(H;2,2)}), f^{2\vee} \rangle_{\text{cyc}} \langle {\mathfrak m}({\bf x}_{c_1}^{(H;2,2)}, f^2, {\bf y}_{c_2}^{(H;2,1)}), f^{1\vee} \rangle_{\text{cyc}} 
. \endaligned
\end{equation}
Therefore
$$
(\ref{form2124}) 
= \sum_{c_1,c_2, f^1,f^2} (-1)^{\delta'}
\langle   {\mathfrak m}({\bf x}_{c_1}^{(H;2,1)}, f^1 ,{\bf y}_{c_2}^{(H;2,2)}), f^{2\vee} \rangle_{\text{cyc}}\langle {\mathfrak m}({\bf x}_{c_1}^{(H;2,2)}, 
f^2,{\bf y}_{c_2}^{(H;2,1)}), f^{1\vee} \rangle_{\text{cyc}}.
$$

\begin{claim}
$\delta'$ in the sign above (that is, $(\ref{form2124pre})$) can also be interpreted  as the Koszul sign 
of the permutation
\begin{equation}\label{2126form+++}
\aligned
&\langle \rangle_{(1)},\langle \rangle_{(2)},\mathfrak m_{(1)},\mathfrak m_{(2)},f^1,f^{1\vee},f^2,f^{2\vee},{\bf x},{\bf y}  \\
&\mapsto \langle \rangle_{(1)},\mathfrak m_{(1)},{\bf x}_{c_1}^{(H;2,1)},f^1, {\bf y}_{c_2}^{(H;2,2)},f^{2\vee},\langle \rangle_{(2)},\mathfrak m_{(2)}, {\bf x}_{c_1}^{(H;2,2)}, f^2, 
 {\bf y}_{c_2}^{(H;2,1)},f^{1\vee} 
\endaligned
\end{equation}
plus $\deg f^1 + \deg f^2$.
(Here we put suffix $(1)$ and $(2)$ to $\langle \rangle$ and $\frak m$ to distinguish them.)
\end{claim}

To prove the claim we put 
$$
a =  {\mathfrak m}_{(1)}({\bf x}_{c_1}^{(H;2,1)}, f^1 ,{\bf y}_{c_2}^{(H;2,2)}).
$$
Note that $\deg' a = \deg' f^2$ in case (\ref{2126form}) is non-zero.
(\ref{form2124})  is
$$
\sum_{c_1,c_2, f^1,f^2} (-1)^{\delta'}
\langle {\mathfrak m}_{(2)}({\bf x}_{c_1}^{(H;2,2)}, \langle  a, f^{2\vee} \rangle_{(1)}f^2,{\bf y}_{c_2}^{(H;2,1)}), f^{1\vee} \rangle_{(2)}.
$$
(We omit ${\rm cyc}$ in $\langle \rangle_{\rm cyc}$ for the rest of this subsection.)
We remark that
\begin{equation}\label{form2128}
\langle \rangle_{(1)}, f^2, f^{2\vee}, a \mapsto \langle \rangle_{(1)}, a,  f^{2\vee}, f^2
\end{equation}
induces the Koszul sign\footnote{This calculation is basically the proof of Lemma \ref{lem:sign_coevaluation_pairing}.}
$$
n \deg' a + \deg' f^2 (n-\deg'f^2) \equiv \deg' f^2 \mod 2.
$$
This implies the claim.
In fact we first permute:
$$
\aligned
&\langle \rangle_{(1)},\langle \rangle_{(2)},\mathfrak m_{(1)},\mathfrak m_{(2)},f^1,f^{1\vee},f^2,f^{2\vee},{\bf x},{\bf y}  \\
&\mapsto \langle \rangle_{(1)},f^2,f^{2\vee},\langle \rangle_{(2)},\mathfrak m_{(1)},\mathfrak m_{(2)},f^1,f^{1\vee},{\bf x},{\bf y}\endaligned
$$
without sign. Then we exchange $\frak m_{(1)}$ and $\frak m_{(2)}$  and  permute in the same way as (\ref{form212422}) to
$$
\mapsto \langle \rangle_{(1)},f^2,f^{2\vee},\langle \rangle_{(2)},\mathfrak m_{(2)}, {\bf x}_{c_1}^{(H;2,2)},\mathfrak m_{(1)},{\bf x}_{c_1}^{(H;2,1)},f^1,{\bf y}_{c_2}^{(H;2,2)},{\bf y}_{c_2}^{(H;2,1)},f^{1\vee}
$$
with the Koszul sign $1+\delta' + \deg f^1$. The right side is identified with:
$$
\langle \rangle_{(1)},f^2,f^{2\vee},\langle \rangle_{(2)},\mathfrak m_{(2)}, {\bf x}_{c_1}^{(H;2,2)},a,{\bf y}_{c_2}^{(H;2,1)},f^{1\vee}
$$
We next permute it  to
$$
\mapsto \langle \rangle_{(2)},\mathfrak m_{(2)}, {\bf x}_{c_1}^{(H;2,2)},\langle \rangle_{(1)},f^2,f^{2\vee}, a, {\bf y}_{c_2}^{(H;2,1)},f^{1\vee}
$$
without sign.
We then apply permutation (\ref{form2128}) to obtain
$$
\mapsto \langle \rangle_{(2)},\mathfrak m_{(2)}, {\bf x}_{c_1}^{(H;2,2)},\langle \rangle_{(1)}, a,  f^{2\vee}, f^2,{\bf y}_{c_2}^{(H;2,1)},f^{1\vee}
$$
with the Koszul sign $\deg' f^2$.  It is equivalent to 
$$
\langle \rangle_{(2)},\mathfrak m_{(2)}, {\bf x}_{c_1}^{(H;2,2)},\langle \rangle_{(1)}, \mathfrak m_{(1)},{\bf x}_{c_1}^{(H;2,1)},f^1,{\bf y}_{c_2}^{(H;2,2)},  f^{2\vee}, f_2,{\bf y}_{c_2}^{(H;2,1)},f^{1\vee}
$$
We finally permute and move $\langle \rangle_{(1)}, \mathfrak m_{(1)},{\bf x}_{c_1}^{(H;2,1)},f^1,{\bf y}_{c_2}^{(H;2,2)},  f^{2\vee}$ to the 
top.  (This object has  degree $0$. So it does not cause sign.)  We thus get
$$
\langle \rangle_{(1)},\mathfrak m_{(1)},{\bf x}_{c_1}^{(H;2,1)},f^1, {\bf y}_{c_2}^{(H;2,2)},f^{2\vee},\langle \rangle_{(2)},\mathfrak m_{(2)}, {\bf x}_{c_1}^{(H;2,2)},f^2,
 {\bf y}_{c_2}^{(H;2,1)},f^{1\vee}
 $$
that is the right hand side of (\ref{2126form+++}).
\par
We thus have proved the claim.
The claim immediately implies (\ref{ZisPD}).
\end{proof}

\subsection{Sign for the duality between $\frak p$ and $\frak q$} 
\label{sec:signduality}

In this subsection, we check the sign in Theorem \ref{thm:duality}.
\par
\begin{figure}[h]
\centering
\includegraphics[scale=0.35]{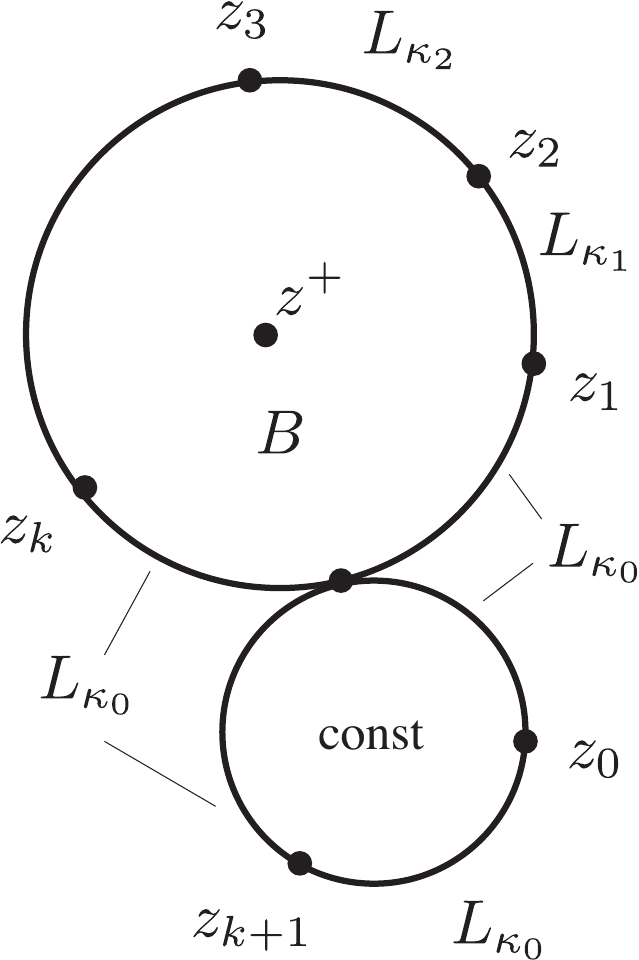}
\caption{Fiber product we use.}
\label{dualfig1}
\end{figure}

The proof is based on \cite[Lemma 4.3 and Proposition 4.1]{ono}
in the case $k'=2$, $j=1$, $B'=0$, $B''=B$, $R_{\alpha_0} = L_{\kappa_k} \cap L_{\kappa_k}$.
Note that we discuss the diagonal component here.
$B'=0$ means that we consider the constant map. 
Moreover $\mu(R_{\alpha_0}) = 0$ and $\dim R_{\alpha_0} = n$.
See Figure \ref{dualfig1}.
\par
In \cite{ono} the case $\frak m \circ \hat{\frak m}$ was studied.
Here we use its version for $\frak m \circ \hat{\frak q}$, that is,\footnote{We use the symbol $\xi$ instead of $h$, for the sake of consistency 
with \cite{ono}.}
$$
\frak m^{\rm form}_{2,0}(\frak q^{\rm form}_{k;1,B}(g;\xi_1,\dots,\xi_k),\xi_{k+1}).
$$
We can easily adapt the argument of \cite{ono} to our case as follows.
\begin{enumerate}
\item[(*)] In the case of  \cite[Lemma 4.3]{ono}
the numbers $\eta_1$, $\eta_2$ appearing in its proof are 
replaced by
$$
\eta_1= \bigl( ( (\text{\rm d-}\deg g + \sum_{i=1}^{k} \text{\rm d-}{\rm deg~}\xi_i) + (\mu_{R_{\alpha}} - \sum_{i=1}^{k} \mu(R_{\alpha_i}) + k -2) \bigr) (\text{\rm d-}{\rm deg~}\xi_{k+1})$$ 
$$\eta_2 =  (\text{\rm d-}\deg g + \sum_{i=1}^{k} \text{\rm d-}{\rm deg~}\xi_i)  (\text{\rm d-}{\rm deg~}\xi_{k+1}).$$ 
(Note $k$ in \cite{ono} is $k+1$ here.)
Namely the extra term $\text{\rm d-}\deg g$ appears both in $\eta_1$ and $\eta_2$.   However the sum $\eta_1 + \eta_2$ does not change.
So the proof works without change.
\item[(**)] In the case of \cite[Proposition 4.1]{ono} the moduli space 
${\mathcal M}_{k''+1}(B'';{\mathcal L}'';{\mathcal R}'')$ 
appeared as the second fiber product factor of the statement.
Note $k'' = k$, $B'' = B$ and  $R_{\alpha_0} = L_{\kappa_k} \cap L_{\kappa_k}$ in our case
so it is ${\mathcal M}_{k+1}(B;{\mathcal L}'';{\mathcal R}'')$.\footnote{$\mathcal R$ is the sequence of 
connected components of $L_{\kappa_k} \cap L_{\kappa_k}$. See \cite[Begining of Section 3]{ono}.}
Moreover we consider moduli space with one interior marked point, that is,
${\mathcal M}_{1;k+1}(B;{\mathcal L}'';{\mathcal R}'')$.
Since the parameter which moves the interior marked point is identified with one complex 
number, we can discuss in the same way as the proof of \cite[Proposition 4.1]{ono}
to find that
the difference between orientations of
$$
{\mathcal M}_{0,3}(0;{\mathcal L}';{\mathcal R}') {}_{{\rm ev}_1}\times_{{\rm ev}_0} {\mathcal M}_{1;k+1}(B;{\mathcal L}'';{\mathcal R}'')
$$
and of 
$$
\partial \mathcal M_{1;k+2}(B;{\mathcal L};{\mathcal R})
$$
is $\kappa$ as in  \cite[Proposition 4.1]{ono}, that is,
\begin{equation} 
\rho \equiv  (k-1)  +\sum_{i=1}^{k} \mu(R_{\alpha_i})
 + n.
\end{equation}
in our case.
(We change the symbol from $\kappa$ to $\rho$ since $\kappa$ is used for other index.)
\end{enumerate}
Now by definition, we have 
\begin{equation}\label{form2122}
\langle
\frak q^{\rm form}_B(g;\xi_1,\dots,\xi_k),\xi_{k+1}
\rangle_{\rm cyc}
= \int_{L_{\kappa_0}}\frak m^{\rm form}_{2,0}(\frak q^{\rm form}_B(g;\xi_1,\dots,\xi_k),\xi_{k+1}).
\end{equation}
We use  \cite[Lemma 4.3]{ono} to show that the integrand in \eqref{form2122} is equal to 
\begin{equation}\label{form2123}
(-1)^{\rho'} {\rm ev}_{0,!}^{(0,B)} \left({\rm ev}^{+*}g \wedge {\rm ev}_1^*\xi_1 \wedge \dots
\wedge {\rm ev}_k^*\xi_k  \wedge {\rm ev}_{k+1}^*\xi_{k+1}
\right)
\end{equation}
where
$$\aligned
\rho' 
&= \text{\rm d-}\deg g + \epsilon(\xi_1,\dots,\xi_k,\xi_{k+1}) 
+ \sum_{i=1}^k \text{\rm d-}{\rm deg~}\xi_{i} + \text{\rm d-}{\rm deg~}\xi_{k+1} - (k+1) -1 + 1 +k \\
&=\text{\rm d-}\deg g + \epsilon(\xi_1,\dots,\xi_k,\xi_{k+1})  +  \sum_{i=1}^{k+1} \text{\rm d-}{\rm deg~}\xi_{i} +1.
\endaligned$$
We have
\begin{equation}
\rho + \rho' = \text{\rm d-}\deg g + \epsilon(\xi_1,\dots,\xi_k,\xi_{k+1})  +  \sum_{i={1}}^{k+1} \text{\rm d-}{\rm deg~}\xi_{i}
+\sum_{i=1}^k \mu(R_{\alpha_i}) + k +n.
\end{equation}
We remark that 
the part of $
\partial \mathcal M_{1;k+2}(B;{\mathcal L};{\mathcal R})
$
corresponding to $
{\mathcal M}_{3}(0;{\mathcal L}';{\mathcal R}') {}_{{\rm ev}_1}\times_{{\rm ev}_0} {\mathcal M}_{1;k+1}(B;{\mathcal L}'';{\mathcal R}'')
$
appears  as the limit of objects $((D^2;z_0,\dots,z_{k+1};z^+),u)$ in $
\mathcal M_{1;k+2}(B;{\mathcal L};{\mathcal R})
$ where $z_0$ converges to $z_{k+1}$.  See Figure \ref{dualfig2}.
\begin{figure}[h]s
\centering
\includegraphics[scale=0.45]{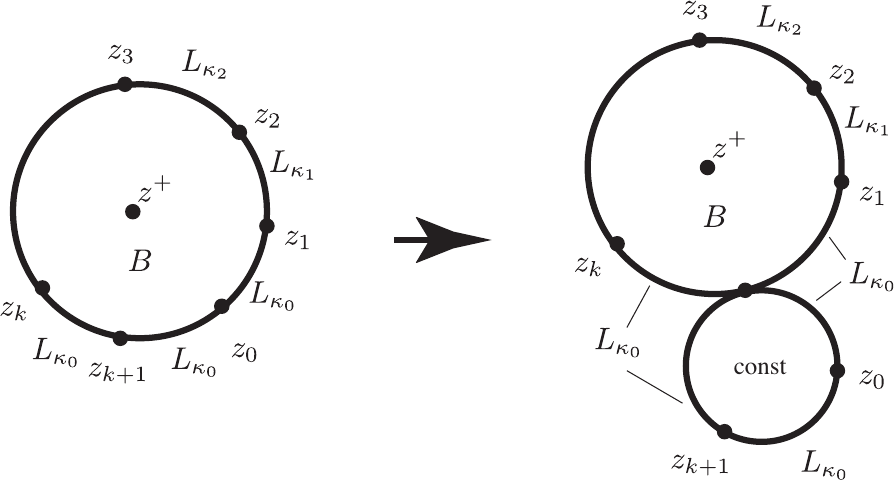}
\caption{$z_{0}$ approaches to $z_{k+1}$.}
\label{dualfig2}
\end{figure}
\par

Consider two local sections $\sigma_0$ and $\sigma_1$ of the forgetful map  $\frak{forget}_{z_0}:{\mathcal M}_{1;k+2}(B;{\mathcal L}; {\mathcal R}) \to 
{\mathcal M}_{1;k+1}(B;{\mathcal L}; {\mathcal R})$ forgetting the marked point $z_0$.  
We denote an element of ${\mathfrak x} \in {\mathcal M}_{1;k+1}(B;{\mathcal L}; {\mathcal R})$ by ${\mathfrak x}=[(u({\mathfrak x}), z_1({\mathfrak x}), \dots, z_{k+1}({\mathfrak x}))]$.  
We define $\sigma_0 ({\mathfrak x}) $ by adding a constant disk with three boundary marked points  
at $z_{k+1}({\mathfrak x})$ as in Figure \ref{dualfig1}.  
We pick $z_0(\sigma_1({\mathfrak x}))$ on the interior of the arc $\overline{z_{k+1} z_1}$ depending smoothly on ${\mathfrak x}$ in an open subset of 
${\mathcal M}_{1;k+1}(B;{\mathcal L}; {\mathcal R})$.  
We consider the smooth family of sections $\sigma_t$, $0 \leq t \leq 1$,  joining $\sigma_0$ and $\sigma_1$.  ($z_0$ moves from $z_{k+1}({\mathfrak x})$ to ${z_0(\sigma_1(\mathfrak x}))$.)   

We integrate 
\begin{equation}\label{integrand}
(-1)^{\rho'} \sigma_t^* \left({\rm ev}^{+*}g \wedge {\rm ev}_1^*\xi_1 \wedge \dots
\wedge {\rm ev}_k^*\xi_k  \wedge {\rm ev}_{k+1}^*\xi_{k+1}
\right)
\end{equation}
on ${\mathcal M}_{1;k+1}(B;{\mathcal L}; {\mathcal R})$.  
Note that ${\rm ev}_0$ on ${\mathcal M}_{1;k+2}(B;{\mathcal L}; {\mathcal R})$ does not contribute to the integrand 
and ${\rm ev}_i$ ($i=1, \dots, k+1)$ on ${\mathcal M}_{1;k+2}(B;{\mathcal L}; {\mathcal R})$ factor through  $\frak{forget}_{z_0}$.  
Hence it is equal to the integral of \eqref{integrand} over ${\mathcal M}_{1;k+1}(B;{\mathcal L}; {\mathcal R})$.  
The image of $\sigma_0$ is the boundary of the union of the images of $\sigma_t$, $0 \leq t < 1$.  
Our orientation convention is 
$$
\text{Index~}D_u\overline{\partial} \oplus {\R}_{z_0} \oplus \dots \oplus \R_{z_{k+1}} /{sl}(2,{\R}).  
$$
Therefore the boundary orientation differs from the orientation of the image of $\sigma_0$ (inherited from ${\mathcal M}_{1;k+1}(B;{\mathcal L}; \R)$) 
by 
$
(-1)^{\mu(B)+1}
$.
\begin{comment}
Therefore, by the forgetful map of $z_0$ it can be identified with 
$\mathcal M_{1;k+1}(B;{\mathcal L''};{\mathcal R}'')$ with sign correction
$
(-1)^{\mu(B)} + 1
$.
\marginpar{I want  a bit more explanation here. KF}
\end{comment}
\par 
We next recall
$$
\dim \mathcal M_{1;k+1}(B;{\mathcal L''};{\mathcal R}'')
\equiv n + \sum_{i=1}^k \mu(R_{\alpha_i}) + k. 
$$
Since $\frak q_B(g;\xi_1,\dots,\xi_k),\xi_{k+1}
\rangle_{\rm cyc}$ is nonzero only when
$$
\sum_{i=1}^{k+1}\text{\rm d-}{\rm deg~}\xi_{i} + \text{\rm d-}\deg g = \dim \mathcal M_{1;k+1}(B;{\mathcal L''};{\mathcal R}'')
$$
we have
$$
\rho + \rho' \equiv  \epsilon(\xi_1,\dots,\xi_k,\xi_{k+1}).
$$
Therefore by (\ref{form2123})
we have
$$
\aligned
&(-1)^{\epsilon(\xi_1,\dots,\xi_k,\xi_{k+1}) + \mu(B) + 1}
\langle \frak q^{\rm form}_B(g;\xi_1,\dots,\xi_k),\xi_{k+1}
\rangle_{\rm cyc}\\
&= 
\int_{\mathcal M_{1;k+1}(B;{\mathcal L''};{\mathcal R}'')}
 \left({\rm ev}^{+*}g \wedge {\rm ev}_1^*\xi_1 \wedge \dots
\wedge {\rm ev}_k^*\xi_k  \wedge {\rm ev}_{k+1}^*\xi_{k+1} 
\right) \\
&= 
\int_X g \wedge {\rm ev}_!^+  \left({\rm ev}_1^*\xi_1 \wedge \dots
\wedge {\rm ev}_k^*\xi_k  \wedge {\rm ev}_{k+1}^*\xi_{k+1} \right)
\endaligned
$$
By (\ref{form21ten2}) we have
$$
\aligned
&\frak p^{\rm form}_B\left({\rm ev}_1^*\xi_1 \wedge \dots
\wedge {\rm ev}_k^*\xi_k  \wedge {\rm ev}_{k+1}^*\xi_{k+1} \right) \\
&= (-1)^{\epsilon(\xi_1,\dots,\xi_k,\xi_{k+1})  + \text{\rm d-}\deg g +\mu(B) + 1}
 {\rm ev}_!^+  \left( {\rm ev}_1^*\xi_1 \wedge \dots
\wedge {\rm ev}_k^*\xi_k  \wedge {\rm ev}_{k+1}^*\xi_{k+1}
\right).
\endaligned
$$
Thus we conclude 
\begin{equation}
\langle
{\frak q}^{\rm form}(g;\xi_1,\dots,\xi_k),\xi_{k+1}
\rangle_{\rm cyc}
=
(-1)^{\text{\rm d-}\deg g}
\langle
g,{\frak p}^{\rm form}(\xi_1,\dots,\xi_k,\xi_{k+1})
\rangle_{\text{\rm PD}_X}.
\end{equation}
This is ${\bf b} = 0$ case of (\ref{dualityinformlevel}) (assuming that 
the $\frak p$ perturbation coincides with $\frak q$ perturbation).\footnote{In (\ref{dualityinformlevel})
we used notation $\deg g$ in place of $\text{\rm d-}\deg g$.}
The case including ${\bf b}$ is the same since bulk deformation is 
of even degree.
\qed
\subsection{Sign for the paring $\langle \rangle_{\rm cyc}$} 
\label{sec:signpair}

In this subsection we define the
sign in \eqref{form823} and (\ref{form8232}). %\marginpar{This short subsection is moved here. KF 2026 Jan.}

Let $p \in L \cap L'$ and $\widetilde{\lambda}_p$ a path of Lagrangian subspaces in $T_pX$ from $T_p L$ to $T_p L'$ equipped with 
a spin structure relative to the trivialization at the end points.  ($\lambda_p$ is the underlying path of Lagrangian subspaces.).   
We pick an orientation $o(\lambda_p)$ of $\det (\text{Index} \overline{\partial}_p)$ with the Lagrangian boundary condition given by $\lambda_p$ 
as in \eqref{operatorpL'L}.  
Denote by $\lambda_p^{\text{rev}}$ the reversal of $\lambda_p$, i.e., $\lambda_p^{\text{rev}}(t)= \lambda_p(1-t)$.  
Since the spin structure relative to 
the ends does not depend on the direction of the path, we have the reversal of $\widetilde{\lambda}^{\text{rev}}_p$.  
Pick an orientation  
$o(\lambda_p)$ (resp. $o(\lambda^{\text{rev}}_p$)) of $\text{Index} \overline{\partial}_p$ (resp. $\text{Index} \overline{\partial}_{p^{\circ}}$) with the Lagrangian boundary condition given by $\lambda_p$ (resp. $\lambda^{\text{rev}}_p$ as in \eqref{operatorpLL'} (resp. \eqref{operatorpL'L}).  
 %\marginpar{The next sentence is removed, since it is impossible to satisfy this condition for both $p$ and $p^{\circ}$. KF 2025 Dec:We may assume that $o(\lambda_p) \otimes \lambda^{\text{rev}}_p$ corresponds to the orientation of $T_p L$ in the sense mentioned just before Lemma \ref{qpropertiescat}. }
We define $c(p) \in \Z_2$ such that  $o(\lambda_p) \oplus o(\lambda^{\text{rev}}) = (-1)^{c(p)} T_pL$.
Then $c(p) + c(p^{\circ}) \equiv \mu(\lambda_p)\mu(\lambda_p^{\text{rev}})$.

For each intersection point $p \in L \cap L'$, we choose and fix $\widetilde{\lambda}_p$ and $o(\lambda_p)$.  
By abuse of notation, we write $p=(p, \widetilde{\lambda}_p, o(\lambda_p))$ and $p^{\circ}=(p, \widetilde{\lambda}^{\text{rev}}_p, o(\lambda^{\text{rev}}_p))$, 
Then we have 
$${\mathfrak m}_{2,0}(p, p^{\circ}) = (-1)^{\mu(\lambda_p)+c_p+1} \delta_{p \in L}, $$ 
where $\delta_{p \in L}$ is the distributional volume form on $L$.\footnote{Actually 
we perturb and so ${\mathfrak m}_{2,0}(p, p^{\circ})$ is a perturbation of $\pm \delta_{p \in L}$. 
The perturbation  does not affect the 
sign we discuss here.} 
Thus 
\begin{eqnarray}\label{form212888}
\langle p, p^{\circ} \rangle_{\text{cyc}} & = & \langle \frak m_{2,0}(1_L, p), p^{\circ} \rangle_{\text{cyc}} \nonumber \\
& = & (-1)^n \langle \frak m_{2,0}(p, p^{\circ}), 1_L \rangle_{\text{cyc}} \nonumber \\
& = & (-1)^{n+c_p+\mu(\lambda_p) + 1}.  
\end{eqnarray}
The sign in \eqref{form8232} is exacly the same.   

\begin{rem}
The cyclic paring should satisfy that  
\begin{equation}\label{cycpairingpp}
\langle p, p^{\circ} \rangle_{\text{cyc}} = (-1)^{\deg' p \deg' p^{\circ}  + 1} \langle p^{\circ}, p \rangle_{\text{cyc}}. 
\end{equation} 
We can check it as follows.  
%By our choice of $o(\lambda_p)$ and $o(\lambda^{\text{rev}}_p)$, we find that  
%$(-1)^{\mu(\lambda_p) \cdot \mu(\lambda^{\text{rev}}_p}) o(\lambda^{\text{rev}}_p) \otimes o(\lambda_p)$ corresponds to the orientation of $T_p L'$, hence 
In the same way as (\ref{form212888}) we have 
%$${\mathfrak m}_{2,0}(p^{\circ}, p) = (-1)^{\mu(\lambda^{\text{rev}} + 1 + \mu(\lambda_p) \cdot \mu(\lambda^{\text{rev}}_p)} \delta_{p \in L'}.$$. 
%By the same computation as above, we have 
$$\langle p^{\circ}, p \rangle_{\text{cyc}} = (-1)^{n+c(p^{\circ}) + \mu(\lambda_p^{\text{rev}}) + 1}.$$
Then \eqref{cycpairingpp} follows from $c(p) + c(p^{\circ}) \equiv \mu(\lambda_p)\mu(\lambda_p^{\text{rev}})$.  
\end{rem} 

\subsection{Sign for cyclic symmetry} 
\label{sec:signcyc}

\begin{prop}
We have: %\marginpar{This short subsection is added. KF 2026 Jan.}
\begin{enumerate}
\item
$$
\widehat{\frak p}_k^{{\bf b}}(x_0,\dots,x_k) = (-1)^{\maltese} \widehat{\frak p}^{{\bf b}}(x_k,x_0,\dots,x_{k-1})
$$
with $\maltese = \deg'x_k (\deg'x_0 + \dots + \deg' x_{k-1})$.
\item
$$
\frak m_k^{{\bf b}}
(x_0,x_1,\dots,x_K) = (-1)^\maltese
\frak m^{{\bf b}}
(x_1,\dots,x_K,x_0),
$$
where 
$\maltese = (\deg' x_0)(\deg' x_1 + \dots + \deg' x_K)$. 
\end{enumerate}
\end{prop}
In other words Theorem \ref{themp} (1) holds and $\frak m^{\bf b}$ is cyclically symmetric. 
\begin{proof}
Except the sign this is already proved in previous sections.
In the cochain level the sign in (1) can be proved by an easy calculation 
using (\ref{form21ten2}).  
Then by the duality (which we proved in Subsection \ref{sec:signduality} with sign) implies (2) holds in the cochain level.
As we discussed already the cyclic symmetries in the cohomology level follows.
Note that in case we include bulk deformation the argument for sign is the same.
\end{proof}
We can prove the sign part of the cyclic symmetries for $\frak p$ and $\frak q$ in general, in the same way.
\newpage

\part{Applications.}
\label{part4}

\section{Introduction to Part \ref{part4}.}\label{sec:intro4}

In this part except Section \ref{sec:toric} we do not use the $T$-adic topology of the Novikov field 
or of modules over it: $\Lambda$ is regarded just as  a field without any particular topology.
The filtered $A_{\infty}$ category we study is linear over $\Lambda$ and the morphism space 
is a {\it finite dimensional} vector space.
So we do not need to (and do not) take a completion of the tensor products among them.
We remark that we have constructed such a filtered $A_{\infty}$ category from Lagrangian
Floer theory in this paper through the following procedure:
\begin{enumerate}
\item
We fix a bulk class $\frak b = \frak b_0 + \frak b_2+\frak b_+$ on $X$, a background class ${\rm st} \in H^2(X, \Z_2)$, and consider a finite collection of Lagrangians $L$ equipped with primitives $\theta_L$ of the restriction of the class $\frak b_2 \in H^2(X;\F)$, as well as   ${\rm st}$-relative spin structures.
 \item Fixing an appropriate discrete monoid $G \subset \R_{\ge 0}$,
we construct, for each $n,k$,  a $G$-gapped filtered $A_{n,k}$ category with these objects.
The morphism spaces are the de Rham complexes of intersections of Lagrangian submanifolds.
(They are infinite dimensional.) %\marginpar{modified KF. 2025 Jan.}
\item
We next take the homotopy limit as $(n,k) \to \infty$, yielding a cyclic and homotopically unital
curved filtered $A_{\infty}$ category which is linear over 
$\Lambda_0$. The  morphism space is
a $T$-adic completion of the tensor product of the de Rham complex with $\Lambda_0$.
\item
We turn on  bounding cochains. 
Then we obtain a cyclic and homotopically unital weakly curvature free filtered $A_{\infty}$ category
whose objects are given by an element of the given collection of Lagrangians (with its data of primitives $\theta_L$ and st-relative Spin structures), together with 
a choice of  bounding cochain.
\item
  We  go to the canonical model.
Thus we obtain a cyclic, unital and weakly curvature free filtered $A_{\infty}$ category which is linear over 
$\Lambda_0$.  
We next change the coefficient ring to $\Lambda$.
Then the  morphism spaces are the cohomology groups of intersections of Lagrangian submanifolds 
over $\Lambda$. (In particular they are finite dimensional.)
At this stage we can forget the $T$-adic topology of $\Lambda$.
\end{enumerate}

Constructions of the open-closed map $\hat{\frak p}$ and the closed-open map $\hat{\frak q}$
are given by induction in the same way, so that the domains are  given by 
the Hochschild homology  and the Hochschild cohomology
 of the  $A_{\infty}$ category 
obtained in Step (5).
Note that the target of closed-open map $\hat{\frak q}$ is Hochschild cohomology where 
we take direct {\it product} of the $\Lambda$-vector spaces $\Hom(B_k\cL,\cL)$ with various $k$.
The domain of the open-closed map $\hat{\frak p}$  is a Hochschild homology  that is a direct {\it sum} in that sense.
They are dual each other, in the usual way that the dual to a direct product is a direct sum.
\par
Our $A_{\infty}$ category is divided into orthogonal summands according to the 
values of the potential function. The Hochschild cohomology (resp. homology) groups that we consider are given by the direct product (resp.  sum) of the groups associated to each summand. 
When we consider the $A_{\infty}$ category  
associated to a given value of the potential function we re-define $\frak m_0$ so that it vanishes (it used to be equal to the product of the unit by the potential value.) In particular, this category has vanishing curvature.
\par
However we note that we have to deal with a (curved) bi-module over objects with different 
values of potential function, that is, the morphism space among them. This plays some role in Section \ref{sec:Theorem1}.

\begin{notation}
As in the previous sections, let $\bL$ and $\bU$ be collections of Lagrangian submanifolds equipped with brane data, which respectively have potential values $\lambda_{\bL}$ and $ \lambda_{\bU} $.
Here brane data consist of a relative spin structure, $\theta_L$ and a bounding cochain.\index{brane data}

We write $\sL$, $\sU$ in\index[syindex]{Lcal@$\sL$}\index[syindex]{Ucal@$\sU$} place of  $\sL \otimes_{\Lambda_0} \Lambda$, $\sU \otimes_{\Lambda_0} \Lambda$
in Part \ref{part4}.

We sometimes omit the symbol $\frak b$ in the notation $\lambda_{\bL}$, $ \lambda_{\bU}$ and etc.
We also write $\widehat{\frak p}$ (resp. $\widehat{\frak q}$) in place of $\widehat{\frak p}^{\bf b}$ or $\widehat{\frak p}^{\frak b}$
(resp. $\widehat{\frak p}^{\bf b}$ or $\widehat{\frak q}^{\frak b}$).
\end{notation} %\marginpar{A notation added. KF. 2024 Dec.}

In Part 4 we write $CF(L,L')$ etc. for $CF^{\rm can}(L,L')$.  For example, when $L=L'$ the group
$CF(L,L)$ denotes the (co)homology group $H(L;\Lambda)$ and not the de Rham complex.\index[syindex]{CF(F@$
CF(L_{\kappa_{i-1}},L_{\kappa_{i}};\F)$}

\section{Proof of Theorem \ref{maintheorem1}.}\label{sec:Theorem1}

By the results of Subsection \ref{relative}, there is a cyclic and unital $A_{\infty}$ category, which we denote by $\Fuk$ 
and which admits fully faithful embeddings\index[syindex]{Fcal@$\Fuk$} 
\begin{equation}
  \sL \to \Fuk \gets \sU
\end{equation}
of categories whose objects are Lagrangians in $\bL$ and $\bU$.  In particular, restricting the diagonal bi-module (see Subsection \ref{sec:examples-bi-modules}) of $\Fuk$ to $\sL$ on the left and to $\sU$ on the right (or vice versa), defines and $\sU-\sL$-bi-module and a $\sL-\sU$-bidmodule which we respectively denote by ${_\sU} \Fuk_{\sL} $  \index[syindex]{UFL@${_\sU} \Fuk_{\sL} $} and $ {_\sL} \Fuk_{\sU} $.   \index[syindex]{LxFU@${_\sL} \Fuk_{\sU} $} Omitting the auxiliary data from the notation, the $\Lambda$ vector spaces underlying these bi-modules are respectively
\begin{align}
  {_\sU} \Fuk_{\sL}(U,L) & = CF^{*}(U,L)  \\
{_\sL} \Fuk_{\sU} (L,U) & = CF^{*}(L,U). 
\end{align}
\begin{rem}
We remind the reader that we use shifted gradings for all modules and bi-modules. This means that an element $x$ of $ CF^{*}(U,L)$ will be considered to have degree $\deg'(x)$ as an element of the value of the bi-module $ {_\sU} \Fuk_{\sL}$ at the pair $(U,L)$.
\end{rem}

The tensor product of these bi-modules over $\sL$ (see \eqref{tensorprodbi-modules}) is a  $\sU$-bi-module  which we denote ${_{\sU}}\Fuk_{\sL}  \displaystyle{\tensor_{\sL}} {_\sL}\Fuk_{\sU}$ with underlying modules
\begin{equation} \label{eq:L_tensor_L-U-bi-module}
\aligned
 &{_{\sU}}\Fuk_{\sL}  \tensor_{\sL} {_\sL}\Fuk_{\sU}(U,V) \\
 &= \bigoplus_{L_0,L_d \in \Ob(\sL) } CF^{*}(U,L_0) \otimes   B_{d}(\sL[1])(L_0,L_{d})\otimes   CF^{*}(L_d,V).
 \endaligned
\end{equation}
Here $d\ge 0$. %\marginpar{$d=0$ included I suppose.  KF 2025 Jan}
In this section, we reinterpret Theorem \ref{ZisPD} to yield the following result:
\begin{thm} \label{thm:commutative_diagram_p_q}
There is a  diagram
  \begin{equation} \label{eq:diagram_factor_p_q}
    \xymatrix{ HH_{*-n} ( \sL,\sL) \ar[r]^{\fphat} \ar[d] & QH^{*}_{\frak b}(X; \Lambda) \ar[d]^{\fqhat} \\
HH^{*-1}\left(\sU, {_{\sU}}\Fuk _{\sL} \displaystyle{ \tensor_{\sL}} {_\sL} \Fuk_{\sU} \right)  \ar[r] & HH^{*} ( \sU,\sU)}
  \end{equation}
which commutates up to sign depending only on the degree $*$ and the dimension $n$.
\end{thm}
\begin{rem}
Putting aside the different technical assumptions, this result should be interpreted as a generalisation and an extension of Proposition 1.3 of \cite{abouzaid:IHES} which studied the case of exact Lagrangians in Liouville manifolds.  The \emph{symplectic cohomology} of the ambient space replaces, in that situation, the quantum cohomology.  The main improvement is that the Lagrangians in \cite{abouzaid:IHES} are assumed, \emph{a priori}, to have vanishing potential value.  In addition, only the analogue of the composition
\begin{equation}
  \fq \circ \fphat  \co  HH_{*-n} ( \sL,\sL)  \to HF^{*}(U,b_{U})
\end{equation}
was studied in \cite{abouzaid:IHES} (although this map was studied in \cite{Ga} in the exact setting).  Here, the map $\fq $ from $ QH^{*}(X; \Lambda)  $  to the Lagrangian Floer cohomology of each object of $\sU$ factors through $\fqhat$.
\end{rem}
\begin{cor} \label{cor:hat_q_hat_p_equal_PO}
If $\lambda_{\bL} \neq  \lambda_{\bU}$, the composition $ \fqhat \circ \fphat $   in Equation \eqref{eq:diagram_factor_p_q} vanishes.
\end{cor}
\begin{proof}
From the commutativity of the diagram in Equation \eqref{eq:diagram_factor_p_q}, it suffices to prove that the Hochschild cohomology of the $\sU$ bi-module $  {_{\sU}}\Fuk _{\sL}  \displaystyle{\tensor_{\sL}} {_\sL}  \Fuk_{\sU}  $ vanishes.  In fact, we shall prove that for every pair of objects $U$ and $V$ in $\sU$, the complex in Equation \eqref{eq:L_tensor_L-U-bi-module} is acyclic.  To prove this, we consider the length filtration on this complex.  The $E_{1}$-page of the associated spectral sequence consists of the direct sum of the groups
\begin{equation}\label{form117}
 CF^{*}(U,L_0) \otimes   CF^{*}(L_0, L_1) \otimes \cdots \otimes   CF^{*}(L_{d-1}, L_{d}) \otimes CF^{*}(L_{d},V).
\end{equation}
and the differential on this page is the tensor product of the differentials.  As this is the tensor product of matrix factorisations with potential values $\lambda_{\bU} - \lambda_{\bL}$, 
$0$, $\lambda_{\bL} - \lambda_{\bU}$ and $\lambda_{\bL} \ne \lambda_{\bU}$, the cohomology necessarily vanishes, as we proved in Lemma \ref{lem:tensor_product_matrix_vanishes}.  We conclude that the bi-module $ {_{\sU}}\Fuk _{\sL} \displaystyle{ \tensor_{\sL}} {_\sL} \Fuk_{\sU}    $ is acyclic whenever $ \lambda_{\bU} \neq \lambda_{\bL}$.
\end{proof}

The argument for Theorem \ref{maintheorem1} can now be completed:
\begin{proof}[Proof of Theorem \ref{maintheorem1}]
Since $\fq$ is a ring homomorphism by Theorem \ref{theoremD1} and $\fq$ sends the unit of $QH^*(X;\Lambda)$ to the unit 
of $HF^*((U, b_U), (U, b_U))$, see  %\marginpar{Ref is added. KF. 2025 Jan} 
Proposition \ref{nontrivialq} and  \cite[(3.8.36.2)]{fooo09},  the non-vanishing of the self-Floer cohomology of a Lagrangian is equivalent to the non-vanishing of the image of the identity $\text{\rm id}_{X} $ under $\fq$.  
\newred{Under the assumption of Theorem  \ref{maintheorem1}, the Hochschild homology of $\sL$ 
splits into the contribution of the different potential values.
This is the consequence of the vanishing of the cohomology of (\ref{form117}) we explained above.}

Therefore there must  be some value $\lambda$ such that the composition
\begin{equation}
\fq \circ \fphat \co  HH_*(\sL_{\lambda}, \sL_{\lambda}) \to HF((U;(\frak b,b_{U})),(U;(\frak b,b_{U}));\Lambda)
\end{equation}
does not vanish.  Corollary \ref{cor:hat_q_hat_p_equal_PO} implies that this can only happen if $\lambda = \PO(U;(\frak b,b_{U}))$, which proves the desired result.
\end{proof}

\subsection{A map of $\sU$ bi-modules} \label{sec:map-su-bi-modules}
In this subsection, we construct the bottom horizontal map in Diagram \eqref{eq:diagram_factor_p_q} using the $A_{\infty}$ structure on  $\Fuk$.  It is induced by a degree $1$ morphism to the diagonal bi-module
\begin{equation}\label{formform209}
\Phi \co {_{\sU}}\Fuk _{\sL} \displaystyle{ \tensor_{\sL}} {_\sL} \Fuk_{\sU}   \to \sU
\end{equation}
with polynomial terms\index[syindex]{Phi@$\Phi$}\index[syindex]{mxfrakF@$\fm^{\Fuk}$}\footnote{$\fm^{\Fuk}$ 
is the structure operation of ${\Fuk}$.}
\begin{equation} \label{eq:map_bi-modules_2_to_1}
\Phi(\bfa \otimes   p \otimes  \bfb  \otimes  q \otimes \bfc ) = 
 \fm^{\Fuk}(\bfa \otimes   p \otimes  \bfb  \otimes  q\otimes  \bfc ).
\end{equation}
The $A_{\infty}$ relations in $\Fuk$ imply that this morphism is a chain map, and hence induces a map on Hochschild cohomology
\begin{equation}
  HH^{*}(\Phi) \co HH^{*-1}(\sU, {_{\sU}}\Fuk _{\sL} \displaystyle{ \tensor_{\sL}} {_\sL} \Fuk_{\sU}   )  \to  HH^{*}(\sU, \sU) 
\end{equation} 
which is the bottom horizontal map in Diagram \eqref{eq:diagram_factor_p_q}.

\subsection{A map of $\sL$ bi-modules}
Our goal in this section is to construct the vertical map on the left of Diagram \eqref{eq:diagram_factor_p_q}.  
The geometric intuition behind this construction is the following Figure \ref{dualfig2}.

\begin{figure}[h]
\centering
\includegraphics[scale=0.45]{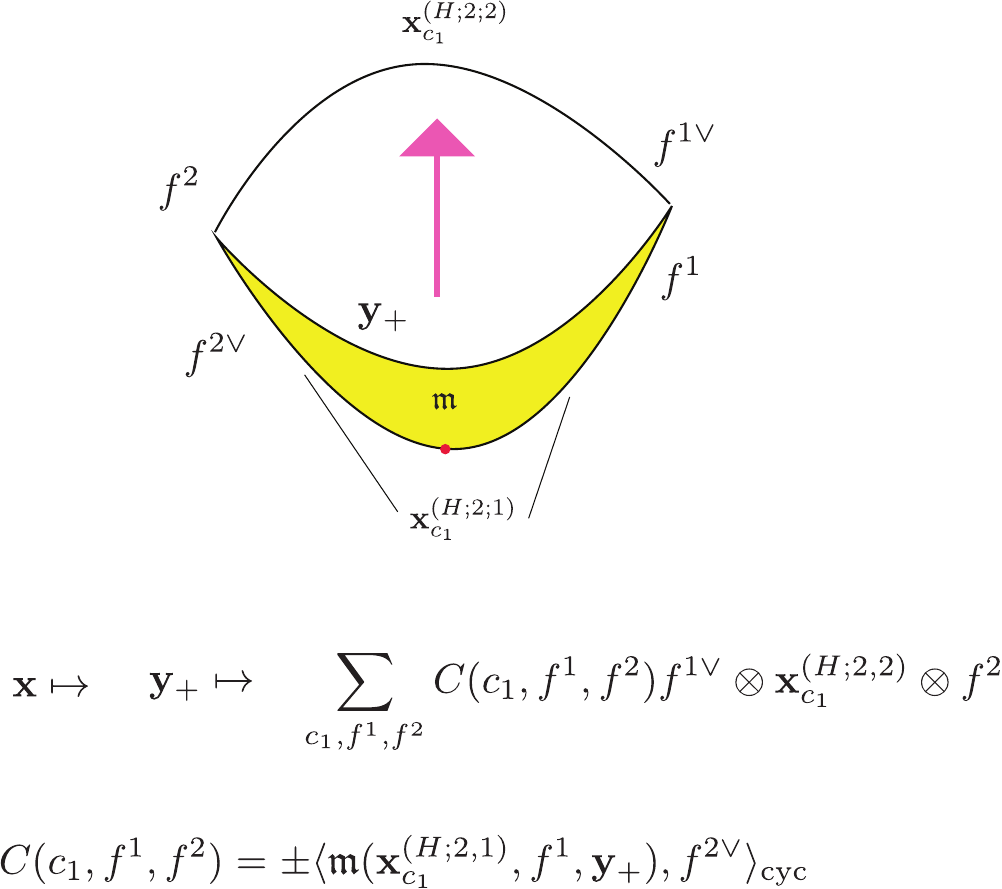}
\caption{The vertical map.}
\label{dualfig2}
\end{figure}

We start by considering the $\sL$ bi-module ${_{\sL}} \Fuk _{\sU} \displaystyle{\tensor_{\sU}} {_{\sU}}  \Fuk_{\sL} $, which also comes equipped with a degree $1$ map
\begin{equation} \label{eq:map_bi-modules_2_1-L-output}
\Phi \co {_{\sL}} \Fuk _{\sU} \tensor_{\sU} {_{\sU}}  \Fuk_{\sL} \to \sL
\end{equation}
given by Equation (\ref{eq:map_bi-modules_2_to_1}) (after swapping the r\^oles of $\sU$ and $\sL$).

The cyclic structure also gives rise to a degree $4-2n$ pairing of chain complexes with shifted gradings
\begin{equation}
CF^{*}(K,V)  \otimes CF^{*}(U,L) \otimes   CF^{*}(L',U') \otimes  CF^{*}(V',K')   \to \Lambda
\end{equation}
which vanishes whenever $U \neq U'$, $L \neq L'$, $K \neq K'$ and $V \neq V'$, or which otherwise is given by 
\begin{equation} \label{eq:pairing_two_tensor_copies_cochain_complex}
\langle  f^1 \otimes f^2 , f^{2\vee} \otimes f^{1\vee}  \rangle = (-1)^{n\deg'(f^1)} \langle f^1, f^{1\vee} \rangle \cdot \langle f^2, f^{2\vee} \rangle,
\end{equation}
where we use the fact that the degree of the pairing is $n$.\footnote{\newred{We remark that ${_{\sL}} \Fuk _{\sU} \tensor_{\sU} {_{\sU}} \Fuk_{\sL}$ 
is a direct {\it sum} over sequences of objects of $\sU$. 
$\left(  {_{\sL}} \Fuk _{\sU} \tensor_{\sU} {_{\sU}} \Fuk_{\sL}    \right)^\vee$ is 
the direct {\it sum} of the duals of those direct summands.}} %\marginpar{Footnote  added.  2025 Feb. FK.} 
\begin{rem}
Note that 
the sign in (\ref{eq:pairing_two_tensor_copies_cochain_complex}) is the Koszul sign associated to the permutation,
$$
\langle \rangle, \langle \rangle, f^1,  f^2,  f^{2\vee},  f^{1\vee}
\mapsto \langle f^1,f^{1\vee} \rangle, \langle f^2,f^{2\vee} \rangle.
$$
The variable $f^{1\vee}$ is exchanged with $f^2,  f^{2\vee}$. This process gives the sign $(-1)^{n\deg' f^{1\vee}} = (-1)^{n\deg f^1}$. 
$\langle  \rangle$ is exchanged with $f^1,  f^{1\vee}$. This process gives the sign $(-1)^{n^2}= (-1)^n$.
\end{rem}

We obtain a map
\begin{equation}
  \psi^\vee \co   \left(  {_{\sL}} \Fuk _{\sU} \tensor_{\sU} {_{\sU}} \Fuk_{\sL}    \right)^\vee \to CH^{*}(\sU,    {_{\sU}}\Fuk_{\sL}  \tensor_{\Lambda}  {_{\sL}}\Fuk_{\sU} ).
\end{equation}
which is characterised by\index[syindex]{psivee@$\psi^\vee$}
\begin{equation} \label{eq:sign_map_tensor_to_Hom}
\langle \psi^\vee (S) (\bfb), f^{2\vee} \otimes f^{1\vee} \rangle = (-1)^{ \deg'(f^{1\vee})( \deg'(\bfb) + \deg'(f^{2\vee}))} S( f^{1\vee} \otimes \bfb \otimes f^{2\vee}).    
\end{equation}
Here the inner product in the left hand side is (\ref{eq:pairing_two_tensor_copies_cochain_complex}).

Composing with the linear dual of the map in Equation (\ref{eq:map_bi-modules_2_1-L-output}), and the map from the diagonal to its dual induced by the cyclic structure, we obtain a degree $n$ morphism $\Psi^{\vee}$ of $\sL$ bi-modules as a composition\index[syindex]{Psivee@$\Psi^{\vee}$}
\begin{equation}
\psi^{\vee} \circ \Phi^{\vee} \circ  \psi \co   \sL \to \sL^{\vee} \to \left(  {_{\sL}} \Fuk _{\sU} \tensor_{\sU} {_{\sU}} \Fuk_{\sL}    \right)^\vee \to CH^{*}(\sU,    {_{\sU}}\Fuk_{\sL}  \tensor_{\Lambda}  {_{\sL}}\Fuk_{\sU} ). 
\end{equation}
% \begin{equation}
% \Psi^{\vee} \co  \cS^{-n} \sL \to   CH^{*}(\sU,    {_{\sU}}\Fuk_{\sL} \tensor_{\Lambda}  {_{\sL}}  \Fuk_{\sU} ).
% \end{equation}
\newred{Here  $\psi$ is obtained from the pairing of $\sL$,\footnote{In other words, $\psi$ is defined by $\psi(x_0)(y) = \langle x_0,y\rangle_{\rm cyc}$.} $\Phi^{\vee}$ is the dual to 
(\ref{formform209}).}
\par
Since Hochschild homology is functorial (see Section \ref{sec:hochschild-homology}), this map induces a degree $n$ map
\begin{equation} \label{eq:map-Psi-dual}
  HH_{*}(\Psi^{\vee}) \co HH_{*}( \sL,  \sL ) \to HH_{*}(\sL,  CH^{*}(\sU, {_{\sU}}\Fuk_{\sL} \tensor_{\Lambda}   {_{\sL}} \Fuk_{\sU} ))). 
\end{equation}

% Using the cyclic structure on $\Fuk$ as in Lemma \ref{cor:cyclic_structure_dual_bi-modules}, 
% %\ref{lem:pairing_induces_iso}, 
% we may identify the linear dual of $\sL$ with $\cS^{n} \sL$.  Shifting both sides by $-2n$, we obtain a map
% \begin{align}
% T \co \cS^{-n} \sL & \to \left( \cS^{2n}\sL \right)^{\vee} \\
% T (S^{-n} x)( S^{2n} a) & = (-1)^{\deg(x)}  \langle x, a \rangle. 
% \end{align}

In order to arrive at the bottom left corner in Diagram \eqref{eq:diagram_factor_p_q}, we begin by noting that we have % an isomorphism
% \begin{equation} \label{eq:shift_bi-module_shifts_homology}
%   CH_{*-n} ( \sL,\sL)  \cong CH_{*} ( \sL, \cS^{-n} \sL),
% \end{equation}
% % This seems to involve no signs!
% so we replace the shifted Hochschild homology group on the top left corner of diagram (\ref{eq:diagram_factor_p_q}) with $ HH_{*} ( \sL, \cS^{-n} \sL) $.  Moreover, at the level of complexes, we have
natural morphisms
\begin{equation} \label{eq:CH_*-almost-commutes-with-CH^*}
\aligned
  &CH_{*}\left(\sL,  CH^{*}\left(\sU, {_{\sU}}\Fuk_{\sL} \tensor_{\Lambda} {_\sL}\Fuk_{\sU} \right)\right) \\
&\to CH^{*}\left( \sU,   CH_{*}\left(\sL, {_\sL} \Fuk_{\sU}  \tensor_{\Lambda} {_{\sU}}\Fuk_{\sL}   \right)   \right) 
\to CH^{*}( \sU,  {_{\sU}}\Fuk _{\sL} \displaystyle{ \tensor_{\sL}} {_\sL} \Fuk_{\sU}  )
 \endaligned
\end{equation}
whose composition we denote by $\Iota$.\index[syindex]{iota@$\Iota$}  The formula for the first map is
\begin{equation} \label{eq:first_maps_iota}
T \otimes \bfx  \mapsto \{ \bfb \mapsto   (-1)^{\deg'(\bfb) \deg'(\bfx)}  \tau( T(\bfb)) \otimes  \bfx \} 
\end{equation}
where 
\begin{equation} \label{eq:sign_shuffle_factors}
\tau(f^1 \otimes f^2) = (-1)^{\deg'(f^1) \deg'(f^2)} f^2 \otimes f^1,
\end{equation}
and $T  \in CH^{*}\left(\sU, {_{\sU}}\Fuk_{\sL} \tensor_{\Lambda} {_\sL}\Fuk_{\sU} \right)$, ${\bf x} \in B\sL$,
${\bf b} \in  B\sU$.  %\marginpar{A line added.  KF 2025 Jan.}?
The second map in Equation \eqref{eq:CH_*-almost-commutes-with-CH^*} is obtained by composition with the map
\begin{align} \label{eq:second_maps_iota}
   CH_{*}\left(\sL, {_\sL} \Fuk_{\sU}  \tensor_{\Lambda} {_{\sU}}\Fuk_{\sL}   \right)  &  \to {_{\sU}}\Fuk _{\sL} \displaystyle{ \tensor_{\sL}} {_\sL} \Fuk_{\sU} \\
f^2 \otimes f^1 \otimes \bfx & \mapsto (-1)^{\deg'(f^2) (\deg'(f^1) + \deg'(\bfx))} f^1 \otimes \bfx \otimes  f^2.
\end{align}
\newred{Here $f^2 \in {_\sL} \Fuk_{\sU}$, $f^1 \in {_{\sU}}\Fuk_{\sL}$ and $\bfx \in B{\sL}$.}
The composition $\Iota \circ HH_{*}(\Psi^{\vee}) $  defines the left vertical map in Diagram \eqref{eq:diagram_factor_p_q}.

\begin{rem} \label{rem:ganatra-result}
Specialising to the special case $\sL = \sU$, the cyclic structure thus induces a map
\begin{equation}
      HH_{*-n} ( \sL,\sL) \to  HH^{*} ( \sL,\sL),
\end{equation}
which Theorem \ref{thm:commutative_diagram_p_q} shows factors through $\fphat$. To relate this to Ganatra's work \cite{Ga2}, %\marginpar{Whether 
%\cite{Ga2} or \cite{GTS} is a correct paper to quote. KF 2025 Sep.} 
we recall that the Calabi-Yau structure also induces a pairing 
\begin{equation}
      HH_{*-n} ( \sL,\sL) \otimes  HH^{*} ( \sL,\sL) \to \Lambda,
\end{equation}
so that we obtain a pairing
\begin{equation}
     HH_{*} ( \sL,\sL) \otimes  HH_{*} ( \sL,\sL) \to \Lambda,
\end{equation}
which  agrees with the Sklyarov pairing \cite{Sk}. For categories which are \emph{cohomologically smooth,} Sklyarov shows that this map is an isomorphism. Tracing through the above isomorphisms, we conclude from  Theorem \ref{thm:commutative_diagram_p_q} that $\fphat$ is necessarily injective. If the ambient symplectic manifold has vanishing first Chern class, the rank of $QH^0(M; \Lambda)$ is $1$, and that of $HH^0(\sL, \sL)$ is at least one, so that Ganatra automatically concludes that the existence of a smooth sub-category $\sL$ of the category of Lagrangians implies that $\sL$ split-generates. In general, one obtains a split-generation result whenever the rank of the even part of $QH^*(X; \Lambda)$ is not larger than that of the even part of $  HH^{*} ( \sL,\sL)$. 
\end{rem}

We now seek an explicit chain level formula for this map: pick bases $f^1_i$ for $ \bigoplus_{L,U} CF^{*}(U,L) $  and $f^2_{j} $ for $ \bigoplus_{K,V} CF^{*}(K,V) $, and dual bases $f_{i}^{1\vee}$ and $f_{j}^{2\vee}$ for the vector spaces $ \bigoplus_{L',U'} CF^{*}(L',U')  $ and $\bigoplus_{V',K'}  CF^{*}(V',K') $ in the sense that  
\begin{equation}
  \langle f^1_{i}, f_{i}^{1\vee} \rangle_{\rm cyc} =  %(-1)^{\deg(g_{j})} 
\langle f^2_{j}, f_{j}^{2\vee} \rangle_{\rm cyc} = 1
\end{equation}
and all other basis elements are orthogonal.%  which leads to the decompositions:

According to Equation \eqref{eq:pairing_two_tensor_copies_cochain_complex}, 
\begin{equation} \label{eq:dual_basis_tensor_prodiuct}
f^1_{i} \otimes f^2_{j} \textrm{ and } (-1)^{n \deg'(f^1_i)  }  f_{j}^{2\vee} \otimes  f_{i}^{1\vee}
\end{equation}
form dual bases with respect to the cochain pairing on the tensor product.   We thus have a cochain map
\begin{align}
  \Lambda & \to   CF^{*}(K,V)  \otimes CF^{*}(U,L) \otimes   CF^{*}(L',U')\otimes  CF^{*}(V',K') \\ \label{eq:map_to_dual_basis}
1 & \mapsto \sum_{i,j} (-1)^{\deg'(f^2_j) + (n+1) \deg'(f^1_i)  } f^1_{i} \otimes f^2_{j}\otimes  f_{j}^{2\vee} \otimes f_{i}^{1\vee}.
\end{align}
The sign difference between Equations \eqref{eq:dual_basis_tensor_prodiuct} and \eqref{eq:map_to_dual_basis} is due to Lemma \ref{lem:sign_coevaluation_pairing}.

\begin{lem}  \label{lem:compute_left_vertical_arrow}
  At the chain level, $ \Iota \circ HH_{*}(\Psi^{\vee}) $  is given by the expression
  \begin{multline}
 (\Iota \circ CH_{*}(\Psi^{\vee}))  (x_0 \otimes \bfx_+) ( \bfb )  = \\ \sum_{f^1_i, f^2_{j}}  (-1)^{~(\ref{eq:total_sign_expression_CH_Psi^vee-short})}
\langle x_{0},  \fm^{\Fuk}( \bfx^{3;1}_{\chi} \otimes  f_{i}^{1\vee} \otimes \bfb \otimes  f_{j}^{2\vee} \otimes \bfx^{3;3}_{\chi})
 \rangle_{\rm cyc}    f^1_{i} 
\otimes  \bfx^{3;2}_{\chi}   \otimes f^2_{j}.
  \end{multline}
The sign is equal to
\begin{equation} \label{eq:total_sign_expression_CH_Psi^vee-short}
\aligned
&\deg'(x_0)  +   n  \left( \deg'(f^2_j) + \deg'(\bfb) \right)  \\
&+ \deg'\left( \bfx^{3;3}_{\chi} \right) \left(1 + n + \deg'(x_0) + \deg'\left( \bfx \right)  \right)  
\\ &+ \left( \deg'(f^1_{i}) + \deg'(\bfx^{3;2}_{\chi} )  \right) \left( \deg'(\bfb) + \deg'(f^2_j) \right).
\endaligned
\end{equation}
\end{lem}
Note that ${\bf b} \in B\sU$, $x_0 \in \sL$, ${\bf x}_+ \in B\sL$
and 
$$
((\Delta \otimes 1)\circ \Delta)(\bfx_+) = \sum_{\chi} \bfx^{3;1}_{\chi} \otimes \bfx^{3;2}_{\chi} \otimes \bfx^{3;3}_{\chi}.
$$
\begin{proof}
We go through the computation step by step.  First, we use the definition of $CH_{*}( \Psi^{\vee} )$ which gives
\begin{equation}\nonumber
CH_{*}(\Psi^{\vee}  ) (x_0 \otimes \bfx_+)  = \sum (-1)^{\eqref{eq:first_sign}}(\psi^{\vee} \circ  \Phi^{\vee}) \left(  \bfx^{3;3}_{\chi} \otimes \psi( x_{0})   \otimes \bfx^{3;1}_{\chi}\right) \otimes   \bfx^{3;2}_{\chi}    \end{equation}
\begin{equation} \label{eq:first_sign}
\deg'\left( \bfx^{3;3}_{\chi} \right) \left(1+ n + \deg'(x_0) + \deg'\left( \bfx \right)  \right) \textrm{ from Equation \eqref{eq:CH_map_bi-modules}.}
\end{equation}
We remark that the roles of $p$, ${\bf a}$, ${\bf a}_{\alpha(3)}^{3;3}$ in \eqref{eq:CH_map_bi-modules}
are played by $\psi(x_0)$, ${\bf x}$, ${\bf x}_{\chi}^{3;3}$, respectively, here.
Since $\frak{deg}(\psi(x_0)) \equiv \deg' x_0 +n \mod 2$ the sign in \eqref{eq:CH_map_bi-modules} gives the above sign. 

Next, we compute that
\begin{align}
 (\psi^\vee \circ \Phi^{\vee}) \left(  \bfx^{3;3}_{\chi} \otimes  \psi( x_{0})  \otimes \bfx^{3;1}_{\chi}\right) (\bfb)   &  \nonumber\\ 
& \hspace{-1.6in} =  \sum (-1)^{\eqref{eq:1.5-sign}}\Big{\langle } (\psi^{\vee} \circ \Phi^{\vee})  \left(  \bfx^{3;3}_{\chi} \otimes \psi( x_{0})  \otimes \bfx^{3;1}_{\chi}\right) (\bfb)  ,   f_{j}^{2\vee} \otimes  f_{i}^{1\vee}    \Big{\rangle}_{\rm cyc}     f^1_{i}  \otimes f^2_{j} \nonumber\\
%Use Equation \eqref{eq:pairing_tensor_product_HH^*} 
& \hspace{-1.6in}=  \sum (-1)^{ \eqref{eq:second_sign}} \Phi^{\vee} \left(  \bfx^{3;3}_{\chi} \otimes \psi(x_{0})  \otimes \bfx^{3;1}_{\chi}\right) \left(  f_{i}^{1\vee} \otimes \bfb \otimes   f_{j}^{2\vee} \right) \cdot f^1_{i}  \otimes f^2_{j} \nonumber\\
% Formula for the dual of a morphism
& \hspace{-1.6in}  = \sum (-1)^{ \eqref{eq:third_sign}}\Big{\langle }  x_{0},    \Phi\left(\bfx^{3;1}_{\chi} \otimes f_{i}^{1\vee} \otimes \bfb  \otimes  f_{j}^{2\vee} \otimes  \bfx^{3;3}_{\chi}\right) \Big{\rangle}_{\rm cyc}  f^1_{i}  \otimes f^2_{j}.\nonumber
\end{align}
With signs
\begin{align}
\label{eq:1.5-sign} & n  \deg'(f^1_i)  \textrm{ from Equation  \eqref{eq:dual_basis_tensor_prodiuct}} \\
\label{eq:second_sign} &  \eqref{eq:1.5-sign} + \deg'(f^{1\vee}_{i})  \left( \deg'(\bfb) + \deg'(f^{2\vee}_j) \right) \textrm{ from Equation  \eqref{eq:sign_map_tensor_to_Hom}}   \\ %161121: Added n+1 
\label{eq:third_sign}&  \eqref{eq:second_sign} + n + \deg'(x_0)  \textrm{ from Equation (\ref{eq:dual_map_bi-modules})}.
\end{align}
We remark that the role of 
$T$,
${\bf b}$,
$\phi$,
${\bf a}$,
$p$ in 
 (\ref{eq:dual_map_bi-modules}) 
 is played by
 $\Phi$,
$\bfx^{3;3}_{\chi}$,
$\psi(x_{0})$,
$\bfx^{3;1}_{\chi}$,
$f_{i}^{1\vee} \otimes \bfb \otimes   f_{j}^{2\vee}$
here, respectively.
Since $\deg \Phi = 1$ the sign in  (\ref{eq:dual_map_bi-modules}) 
becomes $n + \deg'(x_0)$ here.

Finally, the contribution of the three signs in Equation \eqref{eq:first_maps_iota}, \eqref{eq:sign_shuffle_factors}, and \eqref{eq:second_maps_iota} to the map $\Iota$ is
\begin{multline}\label{eq:fifth-sign}
\deg'(\bfx^{3;2}_{\chi} ) \deg'(\bfb )  +  \deg'(f^2_{j}) \deg'(f^1_{i}) + \deg'(f^2_{j})( \deg'(f^1_{i}) + \deg'(\bfx^{3;2}_{\chi} ) )   \\
=   \deg'(\bfx^{3;2}_{\chi} )  \left(\deg'(\bfb) + \deg'(f^2_{j}) \right).
\end{multline}
Adding up the signs, we get the desired expression.
(We use the formulas $\deg f^1_i +  \deg f_i^\vee = n$, $\deg f^2_j +  \deg f_j^{2\vee} = n$ also during calculation.)
\end{proof}
\begin{rem} \label{rem:sign_ch_psivee_from_Koszul} 
Using the relation $\deg'(x_0)  = 1 + n + \deg'(\bfx^{3;1}_{\chi} ) + \deg'(f^1_i) + \deg'(\bfb) +  \deg'(f^2_j) + \deg'\left( \bfx^{3;3}_{\chi} \right) $, we can rewrite the sign in Equation \eqref{eq:total_sign_expression_CH_Psi^vee-short} as:
 %\marginpar{Formula 
%is slightly modified, since this is the form we can derive without using $\deg'(x_0)  = 1 + n + \deg'(\bfx^{3;1}_{\chi} ) + \deg'(f_i) + \deg'(\bfb) +  
%\deg'(g_j) + \deg'\left( \bfx^{3;3}_{\chi} \right) $. KF 2025 March}
\begin{equation} \label{eq:total_sign_expression_CH_Psi^vee}
\aligned
&\deg' x_0 + \deg'{\bf b} + \deg'(f^2_j) + (n+1)  \left( \deg'(f^2_j)  +  \deg'(\bfb) \right)  \\ 
&+ 
 \deg'\left( \bfx^{3;3}_{\chi} \right)  \left( \deg'(\bfx^{3;2}_{\chi})   + \deg'(f^1_{i}) \right)     \\ 
 &+\left( \deg'(f^1_{i}) + \deg'(\bfx^{3;2}_{\chi} )  + \deg'\left( \bfx^{3;3}_{\chi} \right) \right) \left( \deg'(\bfb) + \deg'(f^2_j)\right).
 \endaligned
  \end{equation}
This  can be derived by Koszul conventions by taking the product of (i) the Koszul signs for the permutation:
  \begin{equation}\label{form2339}
  \aligned
  &
 \frak m, \langle \rangle,  x_0, \bfx^{3;1}_\chi, \bfx^{3;2}_\chi, \bfx^{3,3}_\chi, \bfb, f^1_i, f^2_j,  f_j^{2\vee},  f_i^{1\vee}     \\ 
 &\mapsto
 \langle \rangle, x_0,  \frak m, \bfx^{3;1}_\chi,  f_i^{1\vee}, \bfb,  f_j^{2\vee}, \bfx^{3;3}_\chi, f^1_i, \bfx^{3;2}_\chi, f^2_j
 \endaligned
  \end{equation} %\marginpar{Formula changed, KF 2025 March}
together with the signs associated to, (ii)  the sign in Equation \eqref{eq:map_to_dual_basis}.
We can check it by a straight forward computation. 
\end{rem}

\subsection{Proof of Theorem \ref{thm:commutative_diagram_p_q}}
Note that two of the arrows in Diagram \eqref{eq:diagram_factor_p_q} are not labelled.  The following is a more precise version:
\begin{lem} \label{lem:commutative_diagram_p_q_all_arrows_labelled}
The following diagram commutes up to an overall sign depending only on the degree  $*$ 
and $n$:
\begin{equation} \label{eq:commutative_diagram_p_q_all_arrows_labelled}
    \xymatrix{ HH_{*} ( \sL, \sL) \ar[rr]^{\fphat} \ar[d]^{HH_{*}(\Psi^{\vee})} & & QH^{*}(X; \Lambda) \ar[d]^{\fqhat} \\
HH_{*}(\sL,  CH^{*}(\sU,{_{\sU}}\Fuk_\sL\displaystyle{ \tensor_{\Lambda}} {_\sL} \Fuk_{\sU}))) \ar[r]^(.55){H^{*}(\Iota)}  & HH^{*}(\sU, {_{\sU}}\Fuk_{\sL} \displaystyle{ \tensor_{\sL}} {_\sL} \Fuk_{\sU} )  \ar[r]^(.6){HH^{*}(\Phi)} & HH^{*} ( \sU,\sU).}
\end{equation}
\end{lem}

We shall derive this result from Theorem \ref{ZisPD}.  In order to do this, we use the perfect pairing
\begin{equation}
  \langle \_, \_ \rangle_{\rm HH} \co HH^{*} ( \sU,\sU) \otimes HH_{*} ( \sU, \sU) \to \Lambda
\end{equation}
defined
 %\marginpar{$\langle \_, \_ \rangle_{\rm HH}$. The notation should be unified with the notations 
%with other places.  KF 2025 Aug.},  
in Lemma \ref{lem:pairing_HH}, using the cyclic structure on $\sU$.

Our goal now is to compare the pairing $Z$ introduced in subsection \ref{inproZ}, with the pairing
\begin{align}
  CH_{*} ( \sL,\sL) & \otimes CH_{*} ( \sU,\sU) \to \Lambda \\
x_{0} \otimes \bfx_+ \otimes y_0 & \otimes \bfy_+  \mapsto \langle (CH^{*}(\Phi) \circ \Iota \circ CH_{*}(\Psi^{\vee})) \left( x_0 \otimes \bfx_+ \right),  y_0 \otimes \bfy_+ \rangle_{\rm HH}. \label{eq:pairing_homology_cohomology}
\end{align}
\begin{lem}\label{lem231010}
The right hand side of  Equation $(\ref{eq:pairing_homology_cohomology})$ is given by the expression
\begin{equation} \label{eq:pairing_expression_cyclic_chains_U_L}
\aligned
 \sum (-1)^{\maltese_1} 
&\Big{\langle}  x_{0},  \fm^{\Fuk} \left(  \bfx^{3;1}_{\chi} \otimes  f_{i}^{1\vee} \otimes \bfy^{3;2}_{\eta} \otimes  f_{j}^{2\vee} \otimes \bfx^{3;3}_{\chi} \right) \Big{ \rangle}_{\rm cyc}  \\
&\Big{  \langle}  y_{0},   \fm^{\Fuk} \left( \bfy^{3;1}_{\eta} \otimes  f^1_{i} \otimes \bfx^{3;2}_{\chi} \otimes  f^2_{j} \otimes \bfy^{3;3}_{\eta} \right)  \Big{\rangle}_{\rm cyc}  
\endaligned
\end{equation}
with the sign $\maltese_1$ given by
 %\marginpar{The way the sign  is given is changed  KF 2026 Jan.
%The rest of this subsection is rewritten. KF. 2026 Jan.}
the Koszul sign associated to the composition of permutations:
\begin{equation} \label{eq:total_sign_pairing_homological} 
\aligned
&\Phi, \iota, \Psi^{\vee}, x_0, {\bf x}_+, y_0, {\bf y}_+ \\
& \mapsto y_0, \Phi, \iota, \Psi^{\vee}, x_0, {\bf x}_+, {\bf y}_+ \\
& \mapsto y_0, \Phi, {\bf y}^{3;1}_{\eta}, \iota, \Psi^{\vee}, x_0, {\bf x}_+, {\bf y}^{3;2}_{\eta}, {\bf y}^{3;3}_{\eta}
\endaligned
\end{equation}
together with $(-1)$ and the sign $(-1)^{(\ref{eq:total_sign_expression_CH_Psi^vee-short})}$.
\end{lem}
\begin{proof}
We give a step by step computation.  First, we use Equation \eqref{eq:pairing_HH_q_first} to unwind the definition of the pairing between Hochschild homology and cohomology:
\begin{multline}
 \langle (CH^{*}(\Phi) \circ \Iota \circ CH_{*}(\Psi^{\vee})) \left(  x_{0} \otimes\bfx_+ \right), y_0 \otimes   \bfy_+ \rangle_{\rm HH}  \\ =  (-1)^{\eqref{eq:sixth_sign}}  \langle y_0,  CH^{*}(\Phi) \circ \Iota \circ CH_{*}(\Psi^{\vee}) \left(  x_{0} \otimes \bfx_+ \right)(\bfy_+)  \rangle .
\end{multline}
The sign in Equation  \eqref{eq:pairing_HH_q_first} is given by  
\begin{equation}
  \label{eq:sixth_sign}
1 + \deg'(y_0) \left(1 + n +  \deg'(\bfy_+)  \right).
\end{equation}
This is the Koszul sign associated to the permutation 
$$
\Phi, \iota, \Psi^{\vee}, x_0, {\bf x}_+, y_0, {\bf y}_+ \mapsto y_0, \Psi, \iota, \Phi^{\vee}, x_0, {\bf x}_+, {\bf y}_+
$$
times $(-1)$.  (This minus sign comes from the graded {\it anti}-commutativity of $\langle \rangle_{\rm cyc}$.)
Next, we recall that
\begin{multline}
(CH^{*}(\Phi) \circ \Iota \circ CH_{*}(\Psi^{\vee})) \left(  x_{0} \otimes \bfx_+ \right)(\bfy_+)  \\ = \sum (-1)^{~\eqref{eq:seventh_sign}} \Phi \Big{ (} \bfy^{3;1}_{\eta} \otimes  (\Iota \circ CH_{*}(\Psi^{\vee})) \left( x_{0}  \otimes \bfx_+ \right)(  \bfy^{3;2}_{\eta}  ) \otimes \bfy^{3;3}_{\eta}  \Big{)} 
\end{multline}
According to Equation \eqref{eq:induced_map_CH^*}, the sign is in the above expression is
\begin{equation} \label{eq:seventh_sign}
  \deg'(\bfy^{3;1}_{\eta} ) \left(1 + n+ \deg'(x_0) + \deg'( \bfx_+)\right) .
\end{equation}
This is the Koszul sign associated to the permutation:
$$
\iota, \Psi^{\vee}, x_0, {\bf x}_+, {\bf y}_+ \mapsto {\bf y}^{3;1}_{\eta}, \iota, \Psi^{\vee}, x_0, {\bf x}_+, {\bf y}^{3;2}_{\eta}, {\bf y}^{3;3}_{\eta}.
$$
We recall that $\bfb$ (in Lemma \ref{lem:compute_left_vertical_arrow}) is to be replaced by $ \bfy^{3;2}_{\eta} $. 
Therefore Lemma \ref{lem231010} follows from the fact that 
the assignment to associate Koszul sign to a permutation is a group homomorphism to $\{\pm 1\}$ from 
an appropriate permutation group. 
\end{proof}
We next rewrite the sign $\maltese_1$ in (\ref{eq:pairing_expression_cyclic_chains_U_L})
to a Koszul sign associated to a certain single permutation.
\par
We first remark that $\Phi$ is by definition the operator $\frak m$ which has degree $1$.
Therefore we  may rewrite $\Phi$ as $\frak m$, during the definition and calculation of 
Koszul sign.
$\iota$ is of degree $0$ and so  can be omitted during the definition and calculation of 
Koszul sign.
$\Psi^{\vee}$ is defined by $\frak m$ and $\langle \rangle_{\rm cyc}$ and 
has degree $n+1$ modulo $2$. Therefore we may replace it 
by $\frak m,\langle \rangle_{\rm cyc}$, during the definition and calculation of 
Koszul sign.
Thus the permutation (\ref{eq:total_sign_pairing_homological}) 
may be rewritten as
\begin{equation} \label{eq:total_sign_pairing_homological2} 
\aligned
&\frak m_{(1)}, \frak m_{(2)},\langle \rangle_{(2)}, x_0, {\bf x}_+, y_0, {\bf y}_+ \\
& \mapsto y_0, \frak m_{(1)}, {\bf y}^{3;1}_{\eta}, \frak m_{(2)},\langle \rangle_{(2)},  x_0, {\bf x}_+, {\bf y}^{3;2}_{\eta}, {\bf y}^{3;3}_{\eta}
\endaligned
\end{equation}
We compose it with the permutation (\ref{form2339}) that gives (\ref{eq:total_sign_expression_CH_Psi^vee-short}).
The composed permutation is:
\begin{equation}\label{form2531}
\aligned
&\langle \rangle_{(1)}, \frak m_{(1)}, \frak m_{(2)},\langle \rangle_{(2)}, x_0, {\bf x}_+, y_0, {\bf y}_+, f^1_i,f^2_j, f_j^{2\vee}, f_i^{1\vee},  
\\
&\mapsto
\langle \rangle_{(1)}, y_0, \frak m_{(1)}, {\bf y}^{3;1}_{\eta}, \frak m_{(2)},\langle \rangle_{(2)},  x_0, {\bf x}_+, {\bf y}^{3;2}_{\eta}, {\bf y}^{3;3}_{\eta}, f^1_i,f^2_j, f_j^{2\vee}, f_i^{1\vee},  
\\
&\mapsto \langle\rangle_{(1)}, y_0, \frak m_{(1)}, {\bf y}^{3;1}_{\eta}, \frak m_{(2)},\langle \rangle_{(2)},  x_0, {\bf x}_+, {\bf y}^{3;2}_{\eta}, f^1_i,f^2_j, f_j^{2\vee}, f_i^{1\vee},  {\bf y}^{3;3}_{\eta}, 
\\&
\mapsto \langle\rangle_{(1)}, y_0, \frak m_{(1)}, {\bf y}^{3;1}_{\eta}, 
\langle \rangle_{(2)}, x_0,  \frak m_{(2)}, \bfx^{3;1}_\chi,  f_i^{1\vee}, {\bf y}^{3;2}_{\eta},  f_j^{2\vee}, \bfx^{3;3}_\chi, f^1_i, \bfx^{3;2}_\chi, f^2_j
,  {\bf y}^{3;3}_{\eta}, 
\\&
\mapsto \langle\rangle_{(1)}, y_0, \frak m_{(1)}, {\bf y}^{3;1}_{\eta}, f^1_i, \bfx^{3;2}_\chi, f^2_j
,  {\bf y}^{3;3}_{\eta}, 
\langle \rangle_{(2)}, x_0,  \frak m_{(2)}, \bfx^{3;1}_\chi,  f_i^{1\vee}, {\bf y}^{3;2}_{\eta},  f_j^{2\vee}, \bfx^{3;3}_\chi, 
\endaligned
\end{equation}
Here the first permutation is (\ref{eq:total_sign_pairing_homological2}),
the second permutation does not cause sign, 
and third permutation is  (\ref{form2339}).
The fourth permutation does not cause sign since 
the shifted degree of 
$\langle \rangle_{(2)}, x_0,  \frak m_{(2)}, \bfx^{3;1}_\chi,  f_i^{1\vee}, {\bf y}^{3;2}_{\eta},  f_j^{2\vee}, \bfx^{3;3}_\chi$
is even.
Therefore the sign $\maltese_1$ is the Koszul sign 
associated to (\ref{form2531}) together the sign 
(ii)'  the sign in Equation \eqref{eq:map_to_dual_basis} and the sign (iii)'  $(-1)$.

Using cyclic symmetry and anti-graded symmetricity of 
$\langle \rangle_{\rm cyc}$ we find that (\ref{eq:pairing_expression_cyclic_chains_U_L})
is
\begin{equation}
\aligned
\sum (-1)^{\maltese_2} 
&\Big{\langle}   \fm^{\Fuk} \left(\bfx^{3;3}_{\chi}  \otimes x_{0} \otimes \bfx^{3;1}_{\chi} \otimes  f_{i}^{1\vee} \otimes \bfy^{3;2}_{\eta}     \right),
 f_{j}^{2\vee}\Big{ \rangle}_{\rm cyc}  
\\
&\Big{  \langle} \fm^{\Fuk} \left( \bfx^{3;2}_{\chi} ,  f^2_{j},  \bfy^{3;3}_{\eta} \otimes y_{0}  \otimes \bfy^{3;1}_{\eta}       \right) ,
f^1_{i}   \Big{\rangle}_{\rm cyc}  
\endaligned
\end{equation}
where $\maltese_2$ is the Koszul sign associated to the permutation:
\begin{equation}\label{2353form}
\aligned
&\langle \rangle_{(1)}, \frak m_{(1)}, \frak m_{(2)},\langle \rangle_{(2)}, x_0, {\bf x}_+, y_0, {\bf y}_+, f^1_i,f^2_j, f_j^{2\vee}, f_i^{1\vee},\\
& \mapsto 
\langle\rangle_{(1)}, \frak m_{(1)},  \bfx^{3;2}_\chi, f^2_j,  {\bf y}^{3;3}_{\eta},  y_0,  {\bf y}^{3;1}_{\eta}, f^1_i
,  
\langle \rangle_{(2)}, \frak m_{(2)}, \bfx^{3;3}_\chi, x_0,  \bfx^{3;1}_\chi,  f_i^{1\vee}, {\bf y}^{3;2}_{\eta},  f_j^{2\vee}, 
\endaligned
\end{equation}
together the sign 
(ii)'  the sign in Equation \eqref{eq:map_to_dual_basis} and the sign (iii)'  $(-1)$.

We now prove:
\begin{prop}
  \label{cor:comparing_Z_with_homological_pairing}
The pairings $Z(x_0 \otimes \bfx_+, y_0 \otimes \bfy_+)$ differs from $\langle CH^{*}(\Phi) \circ \Iota \circ CH_{*}(\Psi^{\vee}) \left( x_0 \otimes \bfx_+ \right),  y_0 \otimes \bfy_+ \rangle_{\rm HH}$ by a sign depending only on $\deg' (x_0 \otimes \bfx_+)$, 
$\deg' (y_0 \otimes \bfy_+)$  and $n$. \qed
\end{prop}
\begin{proof}
We identify 
$\bfx = x_0 \otimes \bfx_+$, $\bfy = y_0 \otimes \bfy_+$.
Then 
$$
\aligned
&\sum_{\chi}  ( \bfx^{3;3}_\chi \otimes x_0 \otimes  \bfx^{3;1}_\chi) \otimes \bfx^{3;2}_\chi
=\sum_{c_1} \text{\bf x}_{c_1}^{(H;2;1)}  \otimes
\text{\bf x}_{c_1}^{(H;2;2)}, \\
&\sum_{\eta}  ( \bfy^{3;3}_\eta \otimes y_0 \otimes  \bfy^{3;1}_\eta) \otimes \bfy^{3;2}_\eta
=\sum_{c_2} \text{\bf y}_{c_2}^{(H;2;1)}  \otimes
\text{\bf y}_{c_2}^{(H;2;2)}.
\endaligned
$$
Moreover we identify the basis $f^1$,  $f^{1\vee}$ and $f^2$,  $f^{2\vee}$ in Subsection \ref{inproZ} with
$f^{1\vee}_i$, $(-1)^{1+\deg' f^1_i(n-\deg' f^1_i)}f^1_i$ and $f^2_j$, $f^{2\vee}_j$ here.\footnote{We remark that 
$\langle f, f^{\vee}\rangle = 1$ implies $\langle f^{\vee}, f\rangle = (-1)^{1+\deg' f(n-\deg' f)}$.
Therefore $(f^{\vee})^{\vee} = (-1)^{1+\deg' f(n-\deg' f)} f$.}
Then, modulo the terms depending only on $\deg' (x_0 \otimes \bfx_+)$, 
$\deg' (y_0 \otimes \bfy_+)$  and $n$ 
the difference of Koszul sign in (\ref{form199}) and (\ref{2353form}) 
is one by the permutation 
\begin{equation}\label{2354form}
f^1_i,f^2_j, f_j^{2\vee}, f_i^{1\vee} \mapsto f^1_i, f_i^{1\vee},f^2_j, f_j^{2\vee} = f^{1\vee},  f^1, f^2, f^{2\vee} \mapsto f^1,  f^{1\vee}, f^2, f^{2\vee}
\end{equation}
This cancels with the difference between signs in Equation \eqref{eq:map_to_dual_basis}
and in Lemma \ref{lem:sign_coevaluation_pairing} together with $(-1)^{\deg' f^1(n-\deg' f^1)}$.
 %\marginpar{This last step is now checked.  KF 2026 Feb.}
In fact, the Koszul sign associated to (\ref{2354form}) is equal to
$$
(\deg' f^2_j + n\deg' f^1_i) + (\deg' f^1_i + \deg' f^2_j) + \deg' f^1_i(n-\deg' f^1_i).
$$
\end{proof}

Having computed the pairing in Equation \eqref{eq:pairing_homology_cohomology}, we can now complete the proof of Theorem \ref{thm:commutative_diagram_p_q} in its more precise version:
\begin{proof}[Proof of Lemma \ref{lem:commutative_diagram_p_q_all_arrows_labelled}]
Since the pairing  between the Hochschild homology and cohomology of $\sU$ is perfect, the commutativity of Equation \eqref{eq:commutative_diagram_p_q_all_arrows_labelled} follows from knowing that, for any element
\begin{equation}
  \bfx \otimes \bfy \in  HH_{*} ( \sL, \sL) \otimes HH_{*} ( \sU, \sU),
\end{equation}
we have an equality 
\begin{equation} \label{eq:two_pairings_HH_*}
   \langle HH^{*}(\Phi) \circ \Iota \circ HH_{*}(\Psi^{\vee}) \left( \bfx \right), \bfy \rangle_{\rm HH} =  \pm
   \langle \fqhat \left(  \fphat\left( \bfx \right) \right), \bfy \rangle_{\rm HH}.
\end{equation}
with sign depending only on $\deg' \bfx$, 
$\deg' \bfy$  and $n$. 
Since $\fphat$ and $\fqhat$ are dual by Theorem \ref{thm:duality}, the right hand side is equal to\
%marginpar{I slightly modified the sign.  KF 2025Feb}\marginpar{$\langle \rangle_{QH}$ the notation does not seem to be consistent with other place.  KF 2025 Aug.}
\begin{equation}
  (-1)^{\deg'(\bfx)}   \langle  \fphat\left( \bfx \right) ,  \fphat \left( \bfy \right) \rangle_{{\rm PD}_X}.
\end{equation}
On the other hand, applying Theorem \ref{ZisPD}, we conclude that the pairing on Quantum cohomology agrees with the pairing $Z( \bfx, \bfy)$, given in Definition  \ref{defnZ}.
\end{proof}

\section{Proof of Theorem \ref{maintheorem2}.}
\label{sec:Theorem2}
In this section we shall derive Theorem \ref{maintheorem2} from the precise version of Theorem \ref{thm:commutative_diagram_p_q} stated as Lemma \ref{lem:commutative_diagram_p_q_all_arrows_labelled}. 
We start with a category $\cL$ satisfying the  assumption of Theorem \ref{maintheorem1}, and let $\cL_{\lambda}$ denote the subcategory of objects whose potential value is $\lambda$.  As explained in Remark \ref{rem:split-generation-summand}, the proof of Theorem \ref{maintheorem2} reduces to the fact that every unobstructed Lagrangian with the potential value $\lambda$ is quasi-isomorphic to a summand of a twisted complex in $\cL_{\lambda}$.\index[syindex]{Llambda@$\cL_{\lambda}$} The relevant twisted complex has been introduced in Section \ref{sec:universal-complex}:  starting with the inclusion of categories $\cL_{\lambda} \subset \Fuk$ we defined,  for each integer $D$, and each object $X$ of $\Fuk$, a twisted complex $  \Tot_{\cL_{\lambda}}^{D}(X) $ in $ \Tw(\cL_{\lambda})  $:
\begin{prop} \label{prop:finite_generation_time}
There is an integer $D$ with the following property:  any pair $(U,b_U)$ 
of Lagrangian submanifold $U$ and its weak bounding cochain $b_U$ with potential value $ \lambda $ is a summand of $  \Tot_{\cL_{\lambda}}^{D}(U,b_U) $. 
\end{prop}

The assumption of Theorem \ref{maintheorem1} is that the unit in the Quantum cohomology ring of $X$ lies in the image of $\fphat$.  Since the Hochschild complex has an ascending filtration by length, this implies the existence of some integer $D$ such that the unit in fact lies in the image of the composition
\begin{equation} \label{eq:finite_HH_hits_identity}
  HH_{*}^{(D)}(\sL, \sL) \to  HH_{*}(\sL, \sL)  \to QH^{*}(X).
\end{equation}

By Lemma \ref{lem:commutative_diagram_p_q_all_arrows_labelled}, and using the fact that the maps in Diagram \ref{eq:commutative_diagram_p_q_all_arrows_labelled} preserve the length filtration, we conclude that the unit in the Hochschild cohomology of the category $\sU$ lies in the image of the evaluation map
\begin{equation}
  HH^{*}\left(\sU, \left( {_{\sU}}\Fuk_{\cL} \displaystyle{ \tensor_{\cL}} {_\cL} \Fuk_{\sU}  \right) ^{(D)}\right)  \underset{HH^{*}(\Phi)}{\longrightarrow} HH^{*} ( \sU,\sU).
\end{equation}

By considering the constant term of a Hochschild cochain, we obtain an evaluation map
\begin{equation}
   HH^{*}\left(\sU, \left( {_{\sU}}\Fuk_{\cL} \displaystyle{ \tensor_{\cL}} {_\cL} \Fuk_{\sU}  \right) ^{(D)}\right)  \to \blue{H^*} \left( {_{\sU}}\Fuk_{\cL} \displaystyle{ \tensor_{\cL}} {_\cL} \Fuk_{\sU}  \right) ^{(D)},
\end{equation}
 where the right hand side is the vector space underlying the bi-module.  Since the restriction from $ HH^{*} ( \sU,\sU) $   to the Floer cohomology of any object is obtained by the same procedure, we conclude that there is a commutative diagram:
 \begin{equation} \label{eq:evaluate_HH^*_underlying}
   \xymatrix{  HH^{*}\left(\sU, \left( {_{\sU}}\Fuk_{\cL} \displaystyle{ \tensor_{\cL}} {_\cL} \Fuk_{\sU}  \right) ^{(D)}\right)
   \quad   \ar[d] \ar[r]^(.6){HH^{*}(\Phi)} & HH^{*} ( \sU,\sU) \ar[d] \\ 
\blue{H^*} \left( {_{\sU}}\Fuk_{\cL} \displaystyle{ \tensor_{\cL}} {_\cL} \Fuk_{\sU}  \right) ^{(D)}   \ar[r] &  HF^{*}(U,b_U) . }
 \end{equation}

 \begin{proof}[Proof of Proposition \ref{prop:finite_generation_time}]
Let $\cY^{r}_{U}$ and $\cY^{l}_{U}$ respectively denote the right and left Yoneda modules over $\cL$ associated to $U$.    By Corollary \ref{cor:criterion_summand_tensor_bi-modules}, it suffices to prove that the unit in $ HF^{*}(U,b_U)$ lies in the image of the map 
\begin{equation}
H^{*}  \left( \cY^{r}_U[1] \otimes_{\cL} \cY^{l}_U[1]  \right)^{(D)} \to HF^{*}(U,b_U).
\end{equation}
On the other hand, $ \cY^{r}_U \otimes_{\cL} \cY^{l} _U$ is the summand of $ {_{\sU}}\Fuk_{\cL} \displaystyle{ \tensor_{\cL}} {_\cL} \Fuk_{\sU} $ corresponding to both inputs in the bi-module being $ (U,b_U) $.  In addition, the evaluation map $\Phi$ in Equation \eqref{eq:map_bi-modules_2_to_1} precisely restricts (upon shifting the bi-modules), to the map $\fm$.  We conclude from Diagram \eqref{eq:evaluate_HH^*_underlying}, and the preceding discussion that the hypothesis of  Corollary \ref{cor:criterion_summand_tensor_bi-modules} indeed holds, and hence that $ (U,b_U) $ is a summand of $ \Tot_{\cL_{\lambda}}^{D}(X) $, where $D$ is fixed, uniformly for all $(U,b_U)$, by   Equation \eqref{eq:finite_HH_hits_identity}.
\end{proof}

\section{Proof of Theorem \ref{maintheorem3}.}
\label{sec:Theorem3}

We shall prove Theorem \ref{maintheorem3} in two parts.  First, in Subsection \ref{sec:Theorem3-I}, we use the results of Section \ref{sec:ring} to show that $\hat{\frak p}  $ is surjective.  In the next two subsections, we develop further applications of the annulus argument, by elaborating on the commutative diagram in Equation \eqref{eq:diagram_factor_p_q}.  Finally, in Subsection \ref{sec:inject-hatfr-p}, we prove that  $\hat{\frak p}  $ is injective, which completes the proof of Theorem \ref{maintheorem3}.

\subsection{Surjectivity of $\hat{\frak p}  $}
\label{sec:Theorem3-I}

Recall, from Subsection \ref{sec:cup-cap-products}, that the Hochschild cohomology of $\cL$ is an algebra, over which Hochschild homology is naturally a module.  
\begin{thm}\label{pqrelation}
The map $\hat{\frak p}$ is a map of modules over  $QH^{*}_{\fb}(X;\Lambda_0)$, i.e. if $\text{\bf x} \in HH_*({\cL})$ 
and $\frak z \in QH_{\frak b}^*(X;\Lambda_0)$, then %\marginpar{$\cup$ is changed to $\cupdot$ KF 2025 Aug.}
\begin{equation}\label{form2410}
\sum_c (-1)^{\maltese_1}\hat{\frak p}(\frak m(\text{\bf x}_c^{(H;4;\blue{3})}
\otimes\frak q(\frak z;\text{\bf x}_c^{(H;4;\blue{4})})\otimes \text{\bf x}_c^{(H;4;\blue{1})})
\otimes \text{\bf x}_c^{(H;4;\blue{2})})
=
\frak z\cupdot \hat{\frak p}(\text{\bf x})
\end{equation}
where
$$
\aligned
\maltese_1 = 
&(\deg' \text{\bf x}_c^{(H;4;3)}+\deg' \text{\bf x}_c^{(H;4;4)})
(\deg' \text{\bf x}_c^{(H;4;1)}+\deg' \text{\bf x}_c^{(H;4;2)}) \\
&+ 
(\deg \frak z + 1) \deg' \text{\bf x}_c^{(H;4;3)} + 1.
\endaligned
$$
Here we use notation $(\ref{DeltaH})$.
 %\marginpar{$\cup$ is changed to $\cupdot$ in (\ref{form2410}).
%The sign in the proof should be rechecked.}
\end{thm}\

\begin{figure}[h]
\includegraphics[scale=0.5]{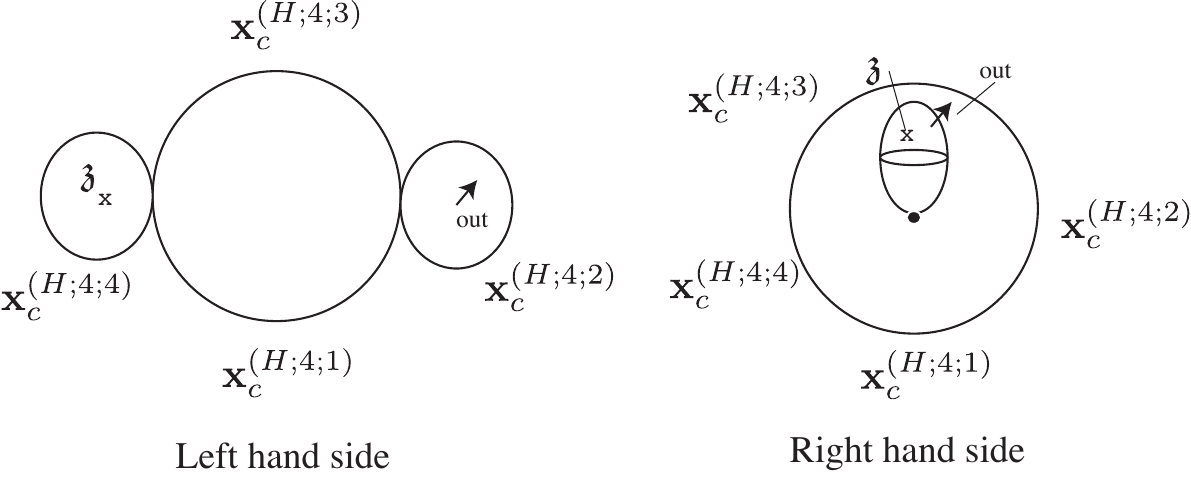}
\caption{Theorem \ref{pqrelation}}
\label{Theorem221fig}
\end{figure}

See Figure \ref{Theorem221fig}. %\marginpar{Sign put.  KF 2025 Jan.}\marginpar{Sign checked and corrected.  KF 2025 Feb.}

\begin{rem}
In the monotone case a somewhat similar result is 
obtained in Biran-Cornea \cite{bircor}.
\end{rem}
\begin{cor} \label{cor:p-surjective}
Under the assumptions of Theorem $\ref{maintheorem1}$, $\hat{\frak p}$ is surjective, and $\hat{\frak q}$ is injective.
\end{cor}
\begin{proof}
Let $\tau_{X} \in HH_{*}(\cL,\cL)$ denote the class whose image under $\hat{\frak p}$ is ${\bf e}_{X}$.  The assertion that $\hat{\frak p}$ is a module map implies that, for every class $x \in QH^{*}_{\fb}(X)$,
\begin{equation}
  \hat{\frak p} (\frak q(x) \cap \tau_{X}) = x \cup^{\frak b}  \hat{\frak p}(\tau_{X}) = x \cup^{\frak b}  {\bf e}_{X} = x,
\end{equation} from which we conclude that $\hat{\frak p}$ is surjective.  By Theorem \ref{thm:duality}, $\hat{\frak q}$ is therefore injective.
\end{proof}
\begin{rem}
  The idea of using the module structure on $HH_{*}(\cL,\cL)$ to show the surjectivity of $\hat{\frak p}$ is due to Sheel Ganatra, who implemented it in his thesis, in the setting of exact symplectic manifolds (see \cite{Ga}).
\end{rem}

Note that Corollary \ref{cor:p-surjective} completes the proof of Theorem \ref{maintheorem1}.  

\begin{proof}[Proof of  Theorem \ref{pqrelation}]
The proof is based on the duality of $\frak q$ and $\frak p$
(Theorem \ref{thm:duality}.) and the fact $\frak q$ is a ring homomorphism
(Theorem \ref{theoremD1}) and proceed as follows.
\par
Let $\frak w$ be any element of $H(X;\F)$. We then calculate
\begin{equation}\label{formulatoprove63}
\aligned
&\sum_c
(-1)^{\maltese_1}
\langle
\frak w,
\hat{\frak p}(\frak m(\text{\bf x}_c^{(H;4;3)}
\otimes\frak q(\frak z;\text{\bf x}_c^{(H;4;4)})\otimes \text{\bf x}_c^{(H;4;1)})
\otimes \text{\bf x}_c^{(H;4;2)})
\rangle_{\text{\rm PD}_X} \\
&=\sum_c
(-1)^{\maltese_2}
\langle
\frak q(\frak w;\text{\bf x}_c^{(H;4;2)}),
\frak m(\text{\bf x}_c^{(H;4;3)}
\otimes\frak q(\frak z;\text{\bf x}_c^{(H;4;4)})\otimes \text{\bf x}_c^{(H;4;1)})
\rangle_{\text{cyc}}
\\
&=
\sum_c
(-1)^{\maltese_3}
\langle
x_0,
\frak m(\text{\bf x}_c^{(H;4;1);+}
\otimes\frak q(\frak w;\text{\bf x}_c^{(H;4;2)})\otimes \text{\bf x}_c^{(H;4;3)}
\otimes \frak q(\frak z;\text{\bf x}_c^{H;4,4}) \otimes
\text{\bf x}_c^{(H;4;1);-})\rangle_{\text{\rm cyc}}
\\
&=
(-1)^{\maltese_4}
\langle
x_0,
\frak q(\frak w \cupdot \frak z;\text{\bf x}_-)
\rangle_{\text{\rm PD}_L}\\
&=
(-1)^{\maltese_5}\langle
\frak p(\text{\bf x}),
\frak w\cupdot \frak z
\rangle_{\text{\rm PD}_X}.
\endaligned
\end{equation}
Here the first equality follows from Theorem \ref{thm:duality} and Lemma \ref{lem:pairing_HH} .
\par
The symbols
$\text{\bf x}_c^{(H;4;1);\pm}$ are defined by 
$$
\text{\bf x}_c^{(H;4;1)} =
\text{\bf x}_c^{(H;4;1);-} \otimes x_0 \otimes \text{\bf x}_c^{(H;4;1);+},
$$
and the second equality is a consequence of cyclic symmetry.
(See (\ref{qismcat}).)
 %\marginpar{The position of the places of $x_0$ is different here from (\ref{qismcat}). 
%The sign is affected.  KF 2025 Jan.}
\par
We define the symbol $\text{\bf x}_-$ by
$$
\text{\bf x} = x_0 \otimes \text{\bf x}_-.
$$
We remark
\begin{equation}\label{signeq224}
\text{\bf x}_- =  (-1)^{\maltese}  \text{\bf x}_c^{(H;4;1);+} \otimes  \text{\bf x}_c^{(H;4;2)} \otimes  \text{\bf x}_c^{(H;4;3)} \otimes  \text{\bf x}_c^{(H;4;4)} \otimes \text{\bf x}_c^{(H;4;1);-},
\end{equation} %\marginpar{Formula (\ref{signeq224}) is added. KF. 2025 Feb}
with
$$
\maltese = \deg' \text{\bf x}_c^{(H;4;1);-} (\deg'x_0 + \deg' \text{\bf x}_c^{(H;4;1);+}+\deg' \text{\bf x}_c^{(H;4;2)} +\deg' \text{\bf x}_c^{(H;4;3)}
+\deg' \text{\bf x}_c^{(H;4;4)}).
$$
Then the third equality is a consequence of  Theorem  \ref{theoremD1}.
The fourth equality again follows from Theorem \ref{thm:duality}.
\par
The sign $\maltese_1$ is as in Theorem \ref{pqrelation}. Other signs are:\footnote{The 
fact that sign works is mostly obvious since all the signs are by Koszul Rule.}
$$
\aligned
\maltese_2 = &\maltese_1 + \deg \frak w + 
\deg' \text{\bf x}_c^{(H;4;2)}
(\deg' \text{\bf x}_c^{(H;4;1)}+\deg' \text{\bf x}_c^{(H;4;3)} + \deg' \text{\bf x}_c^{(H;4;4)} + \deg \frak z),
\\
\maltese_3 = &\maltese_2 + \deg'x_0(\deg' \text{\bf x}_c^{(H;4;1);-}+\deg' \text{\bf x}_c^{(H;4;2)} +\deg' \text{\bf x}_c^{(H;4;3)}\\
&\qquad\qquad\qquad\qquad\qquad\qquad + \deg' \text{\bf x}_c^{(H;4;4)} + \deg \frak z + \deg \frak w+1),
\\
&+ \deg' \text{\bf x}_c^{(H;4;1);+}(\deg' \text{\bf x}_c^{(H;4;1);-}+\deg' \text{\bf x}_c^{(H;4;2)}+\deg' \text{\bf x}_c^{(H;4;3)} \\
&\qquad\qquad\qquad\qquad\qquad\qquad\,\,\,\,\,\quad+ \deg' \text{\bf x}_c^{(H;4;4)} + \deg \frak w + \deg \frak z)\\
& +\deg' \text{\bf x}_c^{(H;4;2)} + \deg \frak w + 1,
\\
\maltese_4 = &\maltese_3 + (\deg \frak w +1)\deg' \text{\bf x}_c^{(H;4;1);+} \\
&\quad+ (\deg \frak z +1)(\deg' \text{\bf x}_c^{(H;4;1);+} + \deg' \text{\bf x}_c^{(H;4;2)} + \deg' \text{\bf x}_c^{(H;4;3)}),
\\
\maltese_5 = &\maltese_4 +  \deg'x_0 + (\deg \frak w+ \deg \frak z)\deg' {\bf x}^-= (\deg \frak w+ \deg \frak z)\deg' {\bf x}.
\endaligned
$$

\par
Theorem \ref{pqrelation} follows from (\ref{formulatoprove63}).
\end{proof}

\subsection{A consequence of the naturality of the cup product.}
\label{sec:cons-natur-cup}
We now return to the situation of Section \ref{sec:Theorem1} where $\cL$ and $\sU$ are categories associated to collections of 
Lagrangians (equipped with bounding cochains) $\bL$ and $\bU$, and admitting a fully faithful embedding into a category $\Fuk$.

Recall that the bottom horizontal map in Diagram \eqref{eq:diagram_factor_p_q} is induced, on Hochschild cohomology, by the map of $\sU$ bi-modules
\begin{equation}\label{form224}
\Phi \co {_{\sU}}\Fuk _{\cL} \displaystyle{ \tensor_{\cL}} {_\cL} \Fuk_{\sU}   \to \sU
\end{equation}
described in Section \ref{sec:map-su-bi-modules}.   In this section, we shall use Equation \eqref{eq:second_maps_iota} to send
$CH_{*}\left(\cL, {_\cL} \Fuk_{\sU}   \displaystyle{\tensor_{\Lambda}} {_{\sU}}\Fuk_{\cL}   \right)  $ to the left hand side of (\ref{form224}).
 %\marginpar{corrected please check.  KF. 2025 Jan.} 
We note that we have an isomorphism
$$
CH_{*}\left(\cL, {_\cL} \Fuk_{\sU}  \tensor_{\Lambda} {_{\sU}}\Fuk_{\cL}   \right) \cong 
{_{\sU}}\Fuk _{\cL} {\displaystyle \tensor_{\cL}} {_\cL} \Fuk_{\sU}.
$$
This induces  the $\sU$ bi-module structure of the left hand side.

This bimodule structure enables us to define the  Hochschild homology
 $$
 HH_{*}\left(\sU, CH_{*}\left(\cL, {_\cL} \Fuk_{\sU}   \displaystyle{ \tensor_{\Lambda}} {_{\sU}}\Fuk_{\cL}   \right)  \right)
 $$
 and induces a map
\begin{equation}
HH_{*}(\Phi) \co  HH_{*}\left(\sU, CH_{*}\left(\cL, {_\cL} \Fuk_{\sU}   \displaystyle{ \tensor_{\Lambda}} {_{\sU}}\Fuk_{\cL}   \right)  \right) \to HH_{*}(\sU, \sU ),
\end{equation}
which is explicitly described in Subsection \ref{sec:hochschild-homology}.
\par
In order to compare this map to the one appearing in Diagram \eqref{eq:diagram_factor_p_q}, we observe that the cap product, which we describe in Subsection \ref{sec:cup-cap-products}, defines maps
\begin{equation}
HH^{*}(\sU,\sU) \otimes HH_{*}(\sU,\sU) \to HH_{*}(\sU,\sU) 
\end{equation}
and
\begin{equation}\label{form238}
\aligned
&HH_{*}(\sU,\sU) \otimes HH^{*}\left(\sU,CH_{*}\left(\cL, {_\cL} \Fuk_{\sU} {\displaystyle \tensor_{\Lambda}} {_{\sU}}\Fuk_{\cL}   \right)  \right)   
\\& \qquad \qquad \qquad
\to HH_{*}\left(\sU,CH_{*}\left(\cL, {_\cL} \Fuk_{\sU}  \tensor_{\Lambda} {_{\sU}}\Fuk_{\cL}   \right) \right).
\endaligned
\end{equation}
We remark that we use 
$$
CH_{*}\left(\cL, {_\cL} \Fuk_{\sU}  \tensor_{\Lambda} {_{\sU}}\Fuk_{\cL}   \right) \cong 
{_{\sU}}\Fuk _{\cL} {\displaystyle \tensor_{\cL}} {_\cL} \Fuk_{\sU} 
$$
and the $\sU$ bi-module structure of the right hand side to define the cap product (\ref{form238}).

The reader may verify, via a straightforward computation the following result which asserts the compatibility of cap products with the maps induced on Hochschild homology and cohomology by morphisms of bi-modules.  The proof is a direct computation from the definition, which we omit.  
\begin{lem}
  The following diagram commutes up to sign
  \begin{equation} \label{eq:naturality_cap_product}
    \xymatrix{   HH_{*}(\sU,\sU) \otimes HH^{*}\left(\sU,CH_{*}\left(\cL, {_\cL} \Fuk_{\sU}{\displaystyle  \tensor_{\Lambda}} {_{\sU}}\Fuk_{\cL}   \right)  \right)   \ar[d]^{\id \otimes HH^{*}(\Phi)  } \ar[r]^-{\capdot} &HH_{*}\left(\sU, CH_{*}\left(\cL, {_\cL} \Fuk_{\sU}{\displaystyle  \tensor_{\Lambda}} {_{\sU}}\Fuk_{\cL}   \right)  \right) \ar[d]^{ HH_{*}(\Phi) } \\
 HH^{*}(\sU,\sU) \otimes  HH_{*}(\sU,\sU)   \ar[r]^{\capdot}  & HH_{*}(\sU,\sU).} 
  \end{equation} %\qed
\end{lem}

\subsection{A consequence of the cyclic structure on $\Fuk$}

In the previous section, we considered the Hochschild homology, over $\sU$, of the bi-module $ CH_{*}\left(\cL, {_\cL} \Fuk_{\sU}{\displaystyle  \tensor_{\Lambda}} {_{\sU}}\Fuk_{\cL}   \right)  $.  
%The cyclic chain complex which computes this group is
%\begin{equation}
% \bigoplus_{\substack{ U,V \in \Ob(\sU) \\ K,L \in \Ob(\cL)} }  CF^{*}(K,V)  \otimes  CF^{*}(U,L) \otimes B(\cL)[1](L,K) \otimes  B(\sU)[1](U,V).
%\end{equation}
%Note that, by simply rearranging the order of the terms in the above complex, we arrive at the complex computing the Hochschild homology of the $\cL$ bi-module $  CH_{*}\left(\sU, {_\sU} \Fuk_{\cL}{\displaystyle  \tensor_{\Lambda}} {_{\cL}}\Fuk_{\sU}   \right)   $:
\begin{lem}
 There is a natural isomorphism
  %\marginpar{$\cap$ is changed to $\capdot$ in several places in this 
 %subsections.  KF 2025 Aug.}
 \begin{equation}
   HH_{*}\left(\sU,  CH_{*}\left(\cL, {_\cL} \Fuk_{\sU}{\displaystyle  \tensor_{\Lambda}} {_{\sU}}\Fuk_{\cL}   \right) \right) \cong HH_{*}\left(\cL, CH_{*}\left(\sU, {_\sU} \Fuk_{\cL}{\displaystyle  \tensor_{\Lambda}} {_{\cL}}\Fuk_{\sU}   \right) \right).
 \end{equation} %\qed
\end{lem}
\begin{proof} By the functoriality of the Hochschild homology, it is enough to 
construct a natural isomorphism
$$
 HH_{*}\left(\sU,  {_\sU} \Fuk_{\cL}{\displaystyle  \tensor_{\cL}} {_{\cL}}\Fuk_{\sU}\right)
 \cong HH_{*}\left(\cL,  {_\cL} \Fuk_{\sU}{\displaystyle  \tensor_{\sU}} {_{\sU}}\Fuk_{\cL}   \right).
$$
The cyclic chain complex which computes the left hand side group is
\begin{equation}
 \bigoplus_{\substack{ U,V \in \Ob(\sU) \\ K,L \in \Ob(\cL)} }  
\big (CF^{*}(U,L) \otimes B(\cL)[1](L,K) \otimes CF^{*}(K,V) \big)  \otimes  B(\sU)[1](U,V).
\end{equation}
Note that, by simply rearranging the order of the terms in the above complex, 
we arrive at the complex computing the Hochschild homology of the $\cL$ bi-module 
$  CH_{*}\left(\sL, {_\cL} \Fuk_{\sU}{\displaystyle  \tensor_{\sU}} {_{\sU}}\Fuk_{\cL}   \right) $.
This finishes the proof.
\end{proof}

With this isomorphism at hand, we shall be able to state the main result of this section.  Its proof is a completely algebraic argument relying only on the fact that the $A_{\infty}$ structure on $\Fuk$ is cyclic.  Nonetheless, it is likely useful for the reader to keep in mind Figures \ref{Figure13-10} and \ref{Figure13-12}, 
which represent the two compositions around Diagram  \eqref{eq:commutative_diagram_tensor_HH_*}.

\begin{lem}
  The following diagram commutes up to sign.
  \begin{equation} \label{eq:commutative_diagram_tensor_HH_*}
    \xymatrix{  
  HH_{*}(\cL, \cS^{-n} \cL)  \otimes HH^{*}(\cL, {_{\cL}}\Fuk _{\sU} {\displaystyle \tensor_{\sU}} {_\sU} \Fuk_{\cL} )  \ar[r]^-{\capdot} &  HH_{*}(\cL,  {_{\cL}}\Fuk _{\sU} {\displaystyle \tensor_{\sU}} {_\sU} \Fuk_{\cL}  ) \ar[dd]^{\cong} \\
 HH_{*}(\cL, \cS^{-n} \cL) \otimes HH_{*}(\sU, \cS^{-n} \sU) \ar[u]^{\id \otimes \Iota \circ HH^*(\Phi^{\vee}) }  \ar[d]^{\Iota \circ HH^*(\Psi^{\vee}) \otimes \id}&   \\
HH^{*}(\sU, {_{\sU}}\Fuk _{\cL} {\displaystyle \tensor_{\cL}} {_\cL} \Fuk_{\sU} ) \ar[r]^-{\capdot} \otimes HH_{*}(\sU, \cS^{-n} \sU)  &   HH_{*}(\sU, {_{\sU}}\Fuk _{\cL} {\displaystyle \tensor_{\cL}} {_\cL} \Fuk_{\sU} )}
 \end{equation}
\end{lem}

\begin{proof} 
Let ${\mathbf x}=S^{-n}x_{0} \otimes \bfx'$ and ${\mathbf y}= S^{-n}y_{0} \otimes \bfy'$ denote generators of $CH_{*}(\cL, \cS^{-n} \cL)$ and $ CH_{*}(\sU, \cS^{-n} \sU) $.  
Firstly, we recall from (\ref{lem:compute_left_vertical_arrow}) that, 
for ${\bf z} \in CH_{*}(\sU,  \sU)$,
$$
\aligned
&(\Iota \circ CH(\Psi^{\vee}))({\mathbf x}) ({\mathbf z})  \\
&= \sum_{B, f^1, f^2, c}T^{ \omega([B])}(-1)^{\maltese_1} {}^+\frak m_B^{\rm f.c.u}({\mathbf x}_c^{(H;2,1)}, f^{1\vee}, {\mathbf z}, f^2) 
 f^1 \otimes {\mathbf x}_c^{(H;2,2)} \otimes f^{2\vee}.
 \endaligned
$$
Here $f^1$, $f^2$ runs over generators of $CF((U_{\kappa'_{j}}, b_{\kappa'_{j'}}),(L_{\kappa_i}, b_{\kappa_i}))$ for some $i, j$.   The operator 
${}^+\frak m_B^{\rm f.c.u}$ is defined in (\ref{form952}).  Here we omit $\vec{\kappa}$ since it is automatically determined by the input. 
The summation is taken over $B, f^1, f^2$ such that ${}^+\frak m_B^{\rm f.c.u}$  is defined.  
\par
Therefore by (\ref{eq:cap_product_chains2}) we have:
$$
\aligned
&(\Iota \circ CH(\Psi^{\vee}))({\mathbf x}) \capdot {\mathbf y} \\
&=
 \sum_{B, f^1, f^2, c_1, c_2}T^{ \omega([B+B'])}(-1)^{\maltese_2} {}^+\frak m_B^{\rm f.c.u}({\mathbf x}_{c_1}^{(H;2,1)}, f^{1\vee}, {\mathbf y}_{c_2}^{(H;3,3)}, f^2) 
 \\&\qquad\qquad\qquad\qquad\qquad\qquad
\frak n_{B'}( f^1 \otimes {\mathbf x}_{c_1}^{(H;2,2) \otimes } \otimes f^{2\vee} \otimes {\mathbf y}_{c_2}^{(H;3,1)} )  \otimes {\mathbf y}_{c_2}^{(H;3,2)}
 \endaligned
$$
Here $\frak n_{B'}$\index[syindex]{nfrakB@$\frak n_{B'}$} is a structure operation of the  bi-module structure on ${_{\sU}}\Fuk _{\cL} {\displaystyle \tensor_{\cL}} {_\cL} \Fuk_{\sU} $.
\par
Taking \eqref{eq:structure_maps_bi-module}-\eqref{eq:structure_maps_bi-module3} and $ {\mathbf y}_{c_2}^{(H;3,1)}  \ne 1$ into account, 
we obtain 
\begin{equation}
\aligned
&(\Iota \circ CH(\Psi^{\vee}))({\mathbf x}) \capdot {\mathbf y}
\\= &  \sum_{B, B', f^1, f^2, f^3, c_1, c_2} {}^+\frak m_B^{\rm f.c.u}({\mathbf x}_{c_1}^{(H;3,1)}, f^{1\vee}, {\mathbf y}_{c_2}^{(H;3,3)}, f^2)  
{}^+\frak m_{B'}^{\rm f.c.u}({\mathbf x}_{c_1}^{(H;3,3)} \otimes f^{2\vee} \otimes {\mathbf y}_{c_2}^{(H;3,1)} \otimes f^3)    \nonumber\\
&\qquad\qquad\qquad\qquad  (-1)^{\maltese_3}T^{\omega([B+B'])} 
f^1 \otimes {\mathbf x}_{c_1}^{(H;3,2)} \otimes f^{3\vee} \otimes {\mathbf y}_{c_2}^{(H;3,2)}. \nonumber
\endaligned
\end{equation} 
See Figure \ref{Figure13-10}.
\begin{figure}[h]
\centering
\includegraphics[scale=0.7]{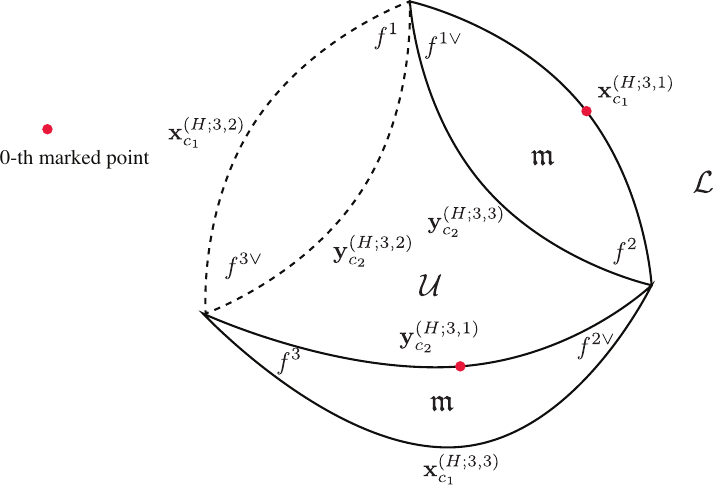}
\caption{$(\Iota \circ CH(\Psi^{\vee}))({\mathbf x}) \cap {\mathbf y}$}
\label{Figure13-10}
\end{figure}
Similarly, we have 
$$
(\Iota \circ CH(\Phi^{\vee})({\mathbf y}))({\mathbf w})
=  \sum_{B, f^1, f^2, c}T^{ \omega([B])}(-1)^{\maltese_5} {}^+\frak m_B^{\rm f.c.u}({\mathbf y}_c^{(H;2,1)}, f^1, {\mathbf w}, f^{2\vee}) 
 f^{1\vee} \otimes {\mathbf y}_c^{(H;2,2)} \otimes f^2.
$$
Therefore:
\begin{equation}
\aligned
& {\mathbf x} \capdot (\Iota \circ CH(\Phi^{\vee})({\mathbf y}))      \\
&=   \sum_{B, B', f^1, f^2, f^3, c_1, c_2} 
{}^+\frak m_B^{\rm f.c.u}({\mathbf y}_{c_2}^{(H;3,1)}, f^1,{\mathbf x}_{c_1}^{(H;3,2)}, f^{2\vee}) 
\nonumber\\&\qquad\qquad\qquad\qquad
{}^+\frak m_{B'}^{\rm f.c.u}( {\mathbf x}_{c_1}^{(H;3,1)} \otimes f^{1\vee} \otimes  {\mathbf y}_{c_2}^{(H;3,2)} \otimes f^3)
\nonumber\\
&\qquad\qquad\qquad\qquad  (-1)^{\maltese_4}T^{\omega([B+B'])} 
 {\mathbf x}_{c_1}^{(H;3,3)} \otimes f^{3\vee} \otimes {\mathbf y}_{c_2}^{(H;3,3)} \otimes f^2. \nonumber
\endaligned
\end{equation} 
\begin{figure}[h]
\centering
\includegraphics[scale=0.7]{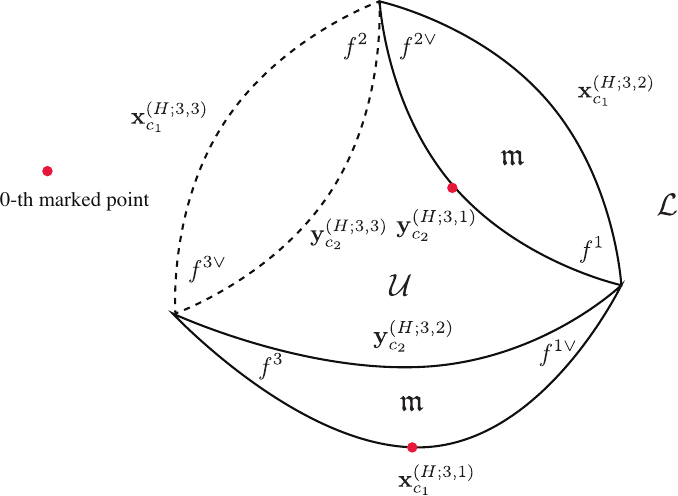}
\caption{${\mathbf x} \cap (\Iota \circ CH(\Phi^{\vee})({\mathbf y})) $}
\label{Figure13-12}
\end{figure}
See Figure \ref{Figure13-12}.
\par
Using cyclic symmetry, we can show that these two maps are chain homotopic via the chain homotopy 
$$
\aligned
&({\mathbf x},{\mathbf y}) \mapsto \\
& \sum_{B,f^1,f^2,c_1,c_2} {}^+\frak m_B^{\rm f.c.u}({\mathbf x}_{c_1}^{(H;2,1)}, f^{1\vee},{\mathbf y}_{c_2}^{(H;2,1)},f^2)T^{ \omega([B])}
\\
&\qquad\qquad\qquad\qquad\qquad(-1)^{\maltese_6}
{\mathbf x}_{c_1}^{(H;2,2)} \otimes f^1 \otimes {\mathbf y}_{c_2}^{(H;2,2)} \otimes f^{2\vee}.
\endaligned
$$ 
See Figure \ref{Figure13-13}.
\begin{figure}[h]
\centering
\includegraphics[scale=0.7]{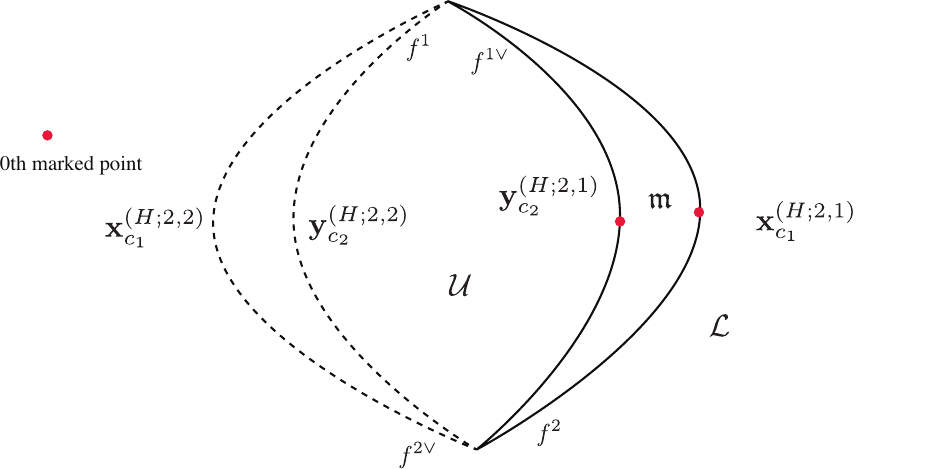}
\caption{Chain homotopy of Diagram (\ref{eq:commutative_diagram_tensor_HH_*}).}
\label{Figure13-13}
\end{figure}

Using the fact all the signs are by Koszul rule we can easily see that the sign matches.

\end{proof}

\subsection{Injectivity of $\hat{\frak p}$}
\label{sec:inject-hatfr-p}

 In this subsection, we complete the proof of Theorem \ref{maintheorem3} by proving:

\begin{prop} \label{prop:injectivity_p}
 Under the assumptions of Theorem $\ref{maintheorem1}$, the map $\hat{\frak p}$ is injective and $\hat{\frak q}$ is surjective.
\end{prop}

In order to prove that $\hat{\frak p}$ is injective, let us consider, for each pair Lagrangian  $L \in \bL$, a Hamiltonian isotopic copy, which we denote $L'$, and we now fix $\bU$ to be a collection of such Lagrangians.  We use this Hamiltonian isotopy to transport brane data from $L$ to $L'$, and obtain a category $\sU$.

The invariance of Floer homology under Hamiltonian isotopies implies that, in a category $\Fuk$ which includes both $\cL$ and $\sU$, the objects corresponding to $L$ and $L'$ are quasi-isomorphic.   In particular, we have isomorphisms
\begin{equation} \label{eq:isomorphism_cohomology_isotopic_objects}
  HF^{*}(L,K) \cong HF^{*}(L',K) \cong HF^{*}(L',K') \cong HF^{*}(L,K') 
\end{equation}
for all pairs of objects $L,K \in \Ob(\cL)$, which are compatible with the multiplication induced by $\fm^{\Fuk}_{2}$ at the cohomological level.

\begin{lem}\label{lem229}
  The morphism 
  \begin{equation}
    \Phi \co {_{\sU}}\Fuk_{\cL}    {\displaystyle\tensor_{\cL}}     {_\cL} \Fuk_{\sU}   \to \sU
  \end{equation}
is an equivalence of $\sU$ bi-modules.
\end{lem}
\begin{proof}
It suffices to prove that, for each pair of objects $L$ and $K$ of $\cL$, that we have an isomorphism on homology
\begin{equation} \label{eq:isomorphism_first_piece_bar_complex_rest}
 H^{*}\left( {_{\sU}}\Fuk_{\cL}    {\displaystyle\tensor_{\cL}}     {_\cL} \Fuk_{\sU} \right) (L',K')  \to HF^{*}(L', K').
\end{equation}
In Equation \eqref{eq:filtration_tensor_product}, we introduced a filtration of the bi-module $  {_{\sU}}\Fuk_{\cL}    {\displaystyle\tensor_{\cL}}     {_\cL} \Fuk_{\sU} $, which we also used in Section \ref{sec:Theorem2}.  We shall use the associated length spectral sequence to show that the map in Equation \eqref{eq:isomorphism_first_piece_bar_complex_rest}.  Note that the cohomology of the subquotients of the length filtration are
\begin{equation}
\aligned
  \bigoplus_{L_1, \ldots, L_d \in \Ob(\cL)} HF^{*}(L',L_1) &\otimes HF^{*}(L_1,L_2) \otimes \cdots \\
  &\otimes  HF^{*}(L_{d-1},L_d) \otimes HF^{*}(L_d,K').
\endaligned
\end{equation}
The direct sum of these groups is therefore the $E_1$ page of a spectral sequence computing the left hand side of Equation \eqref{eq:isomorphism_first_piece_bar_complex_rest}, and the next differential is induced by the product $ \fm^{\Fuk}_{2} $.  Note that the direct sum of the groups at the $E_1$ page admits a map to $HF^{*}(L', K')  $ which gives rise to the map in Equation \eqref{eq:isomorphism_first_piece_bar_complex_rest}.

At this stage, we use the isomorphisms in Equation \eqref{eq:isomorphism_cohomology_isotopic_objects} to identify the cone of the map at the $E_1$ page with
\begin{equation}
\aligned
\bigoplus_{d=0}^{\infty}  \bigoplus_{L_1, \ldots, L_d \in \Ob(\cL)} HF^{*}(L,L_1) &\otimes HF^{*}(L_1,L_2) \otimes \cdots 
\\&\otimes  HF^{*}(L_{d-1},L_d) \otimes HF^{*}(L_d,K),
\endaligned\end{equation}
equipped with the usual bar differential.  Note that the term $d=0$ corresponds to $HF^{*}(L,K)  $.  Since the category $H^{*}(\cL)$ is unital, the bar complex is acyclic, hence Equation \eqref{eq:isomorphism_first_piece_bar_complex_rest} is an isomorphism, and $\Phi$ is an equivalence of bi-modules.
\end{proof}
Applying Lemma \ref{lem:equivalence_iso_HH} to $\Phi$ and to the corresponding map in which $\cL$ and $\sU$ switch roles, we conclude:
\begin{cor} \label{cor:HH_*-of_HH_*-iso-to-HH_*}
The right vertical map in Diagram \eqref{eq:naturality_cap_product} is an isomorphism. \qed
\end{cor}
Of course, this result holds as well if we swap the r\^oles of $\sU$ and $\cL$.  

Let us now assume that every object of $\Fuk$ lies either in $\cL$ and $\sU$ (we recall that the objects of $\cL$ and $\sU$ are assumed  in this section to be Hamiltonian isotopic Lagrangians):
\begin{lem}
 The inclusion of $\cL$ and $\sU$ in $\Fuk$ induce isomorphisms on Hochschild homology which commute with $\hat{\frak p}$:
   \begin{equation} \label{eq:p-hat-commutes-with-inclusion}
     \xymatrix{  HH_{*}(\Fuk, \cS^{-n} \Fuk)  \ar[dr]^{\hat{\frak p}} & HH_{*}(\sU, \cS^{-n} \sU) \ar[l]^{\cong} \ar[d]^{\hat{\frak p}} \\
 HH_{*}(\cL, \cS^{-n} \cL) \ar[u]^{\cong} \ar[r]^{\hat{\frak p}} & QH^{*}_{\fb}(X, \Lambda).}
   \end{equation}
\end{lem}
\begin{proof}[Sketch of proof]
The fact that the inclusions induce isomorphisms on Hochschild homology follows from the more general fact that Hochschild homology is invariant under triangulated closure.  In this special case, we can filter the cyclic bar complex of $\Fuk$ by the number of morphisms belonging, say, to $\sU$.  The first step of this increasing filtration is the cyclic bar complex of $\sU$, and the acyclicity of the other subquotients follows from the acyclicity of the bar complex.

The commutativity of Diagram \eqref{eq:p-hat-commutes-with-inclusion} follows from the fact that the choice of data required to define the map $\hat{\frak p}$ in Section \ref{sec:frakq} is inductive on the objects of $\Fuk$.  In particular, having chosen data for $\sU$ and $\cL$, we can extend this data to all elements of the cyclic bar complex of $\Fuk$, and hence obtain maps which commute at the chain level.
\end{proof}

We complete this section with
\begin{proof}[Proof of Proposition \ref{prop:injectivity_p}]   
We start by summarising several previously obtained commutative diagrams in the following one:
\begin{equation} \label{eq:composition_3_diagrams}
  \xymatrix{ 
QH^{*}_{\fb}(X, \Lambda) \otimes HH_{*}(\sU, \cS^{-n} \sU)  \ar[r]^{\capdot \circ (\hat{\frak q} \otimes \id)} &  HH_{*}(\sU, \cS^{-n} \sU) \\
 & HH_{*}\left(\sU, CH_{*}\left(\cL, {_\cL} \Fuk_{\sU}{\displaystyle  \tensor_{\Lambda}} {_{\sU}}\Fuk_{\cL}   \right)  \right) \ar[u]^{ HH_{*}(\Phi)}  \ar[dd]^{\cong}  \\
HH_{*}(\cL, \cS^{-n} \cL)  \otimes   HH_{*}(\sU, \cS^{-n} \sU)   \ar[uu]^{\hat{\frak p} \otimes \id} \ar[dd]^{\id \otimes \hat{\frak p}} \ar[ur]^{*_1} \ar[dr]^{*_2} &  \\
& HH_{*}\left(\cL, CH_{*}\left(\sU, {_\sU} \Fuk_{\cL}{\displaystyle  \tensor_{\Lambda}} {_{\cL}}\Fuk_{\sU}   \right)  \right) \ar[d]^{ HH_{*}(\Phi)}   \\
HH_{*}(\cL, \cS^{-n} \cL)   \otimes  QH^{*}_{\fb}(X, \Lambda)\ar[r]^{\capdot \circ (\id \otimes \hat{q}) } &  HH_{*}(\cL, \cS^{-n} \cL) }
\end{equation}
Here
$$
*_1 =  \capdot \circ (\Iota \circ HH^*(\Psi^{\vee}) \otimes \id),  \quad
*_2 = \capdot \circ (\id \otimes \Iota \circ HH^*(\Phi^{\vee})).
$$

The middle part of the diagram, is a restatement of \eqref{eq:commutative_diagram_tensor_HH_*}, while the top and bottom commuting squares are both obtained by combining Diagrams \eqref{eq:diagram_factor_p_q} and  \eqref{eq:naturality_cap_product} (and swapping the r\^oles of $\sU$ and $\cL$).

Note that, by Corollary \ref{cor:HH_*-of_HH_*-iso-to-HH_*}, all the vertical arrows in the right most column are isomorphisms.  Let us fix an element $S^{-n} x_0 \otimes \bfx \in HH_{*}(\cL, \cS^{-n} \cL) $ whose image under $\hat{p}$ is the identity in quantum cohomology.  Since the action of quantum cohomology on $HH_{*}(\sU, \cS^{-n} \sU)  $ is unital, the  composition of the inclusion
\begin{align}
HH_{*}(\sU, \cS^{-n} \sU)  & \to  HH_{*}(\cL, \cS^{-n} \cL)  \otimes   HH_{*}(\sU, \cS^{-n} \sU)   \\
S^{-n} y_0 \otimes \bfy & \mapsto  S^{-n} x_0 \otimes \bfx  \otimes S^{-n} y_0 \otimes \bfy
\end{align}
with the top part of Diagram \eqref{eq:composition_3_diagrams} is an isomorphism.

From  the commutativity of the Diagram, we conclude that the composition of this inclusion with the bottom part is also an isomorphism.  But this implies that the map
\begin{equation}
  \hat{\frak p} \co HH_{*}(\sU, \cS^{-n} \sU)  \to QH^{*}_{\fb}(X, \Lambda)
\end{equation}
is injective.  %The commutativity of Diagram \eqref{eq:p-hat-commutes-with-inclusion}.
\end{proof}

\section{Cohomological Smoothness.}
\label{sec:smooth}

\subsection{Statement}
\label{subsec:smoothstate}

In this section we prove Theorem \ref{maintheorem3.3}.  %\marginpar{This section is new.  KF 2025 Aug 27.}
We define cohomological smoothness below.
\begin{defn}\label{defn261}
Let $\Cat$ is a curvature free $A_{\infty}$ category over a field.
We consider the DG category $\mathscr{LR}(\Cat)$\index[syindex]{LRcat@$\mathscr{LR}(\Cat)$} of left $\Cat$ right $\Cat$ bimodules.
The diagonal module $\Cat$ is its object. For any pair of objects $d_1,d_2$ 
of $\Cat$ we can define an object $\frak{Rep}(d_1,d_2)$ of $\mathscr{LR}(\Cat)$ by
$$
\frak{Rep}(d_1,d_2)(c_1,c_2) = \Cat(c_1,d_1) \otimes \Cat(d_2,c_2).
$$ 
We say $\Cat$ is {\it cohomologically smooth}\index{cohomologically smooth} if the diagonal module $\Cat$
is contained in the full subcategory of $\mathscr{LR}(\Cat)$ 
generated by objects $\frak{Rep}(d_1,d_2)$\index[syindex]{Repd1d2@$\frak{Rep}(d_1,d_2)$} via the operations:
\begin{enumerate}
\item[(i)] Degree shift.
\item[(ii)] Adding mapping cones of closed morphisms.
\item[(iii)] Idempotent closure.
\end{enumerate}
\end{defn}
We consider a finite set  $\mathscr L$ consisting of $(L,\theta_L)$
and a finite set ${\bf L}$ as in Subsection \ref{subsec:CatLag}.
Here we consider the bulk class $\frak b = \frak b_2 + \frak b_+$ as in (\ref{decomposefrakb}).
For $(L_1,\theta_{L_1}), (L_2,\theta_{L_2}) \in \mathscr L$
we consider $(L_1 \times L_2, \theta_{L_1} \otimes \theta_{L_2})$.
Here $L_1 \times L_2$ is a Lagrangian submanifold of a symplectic manifold 
$(X^2,-\pi_1^*\omega + \pi_2^*\omega)$, which we denote by $-X \times X$,\index[syindex]{Xminus@$-X$}
and $\theta_{L_1} \otimes \theta_{L_2}: = -\pi_1^*\theta_{L_1} + \pi_2^*\theta_{L_2}$\index[syindex]{theta1otimestheta2@$\theta_{L_1} \otimes \theta_{L_2}$} is a closed one form 
on $L_1 \times L_2$.
$\frak b$ induces the bulk class $\frak b^{(2)} = -\pi_1^*\frak b + \pi_2^*\frak b$
which is decomposed into  $\frak b^{(2)} = \frak b^{(2)}_2 + \frak b_+$, where 
$\frak b^{(2)}_2 = -\pi_1^*\frak b_2 + \pi_2^*\frak b_2 \in H^2(-X\times X;\F)$.
It is easy to see that $d(\theta_{L_1} \otimes \theta_{L_2}) =  \frak b^{(2)}_2$.
\par
Let ${\rm st} \in H^2(X;\Z_2)$ be the background  class 
and $V({\rm st})$ the associated background datum, a real vector bundle on the 3-skelton of $X$ such that its 2nd Stiefel-Whiteney class is ${\rm st}$. 
(See Definition \ref{def:relative_Spin}.)  We assume that, for $(L,\theta_L) \in \mathscr L$,
$L$ is given a $V({\rm st})$ relative spin structure.
As in \cite[Subsection 3.5]{fukaya:functor}, 
we regard $L$ as a Lagrangian submanifold of $-X$ 
and equip it a real oriented bundle  $TX\oplus V({\rm st})$ and its associated
relative spin structure.
We use $V^{(2)}: = \pi_1^*(TX) \oplus V({\rm st}) \oplus \pi_2^*(TX)$ as a background datum, then 
$L_1 \times L_2$ is $V^{(2)}$ relatively spin.
\par
We put
$$
\mathscr L^2 = \{(L_1 \times L_2,\theta_{L_1} \otimes \theta_{L_2}) \mid (L_i,\theta_{L_i}) \in \mathscr L\}.
$$
This is a finite set of $V^{(2)}$-relatively spin Lagrangian submanifolds of $X^2$.
\par
For any bounding cochain $b$ of $(L,\theta_L) \in \mathscr L$ with $b  = \theta_L +  b_1 + b_+$ and bulk class $\frak b$,
the class\index[syindex]{brm*@$b^*$} 
$$
b^*:=-b_1-b_+
$$  
is also a bounding cochain of $(L,-\theta_L)$ with bulk class $-\frak b$, regarded as a Lagrangian submanifold of $-X$.
\begin{rem}\label{rem261}
In case when $\frak b = b_L = 0$ this fact is given in \cite[Theorem 3.54]{fukaya:functor}
(using \cite[Theorem 1.4]{fooo:inv}).
This is based on the claim that:  
by changing $J_X$ to $J_{-X}$, $\omega_X$ to $-\omega_X$, 
$V({\rm st})$ to $(TX) \oplus V({\rm st})$ and 
applying involution $\tau$,\footnote{where $\tau(u_1(z),u_2(z)) = (u_2(\overline z),u_1(\overline z))$} the moduli space of $\beta \in \pi_2(X,L)$ holomorphic disks 
becomes the moduli space of $\tau(\beta)$ holomorphic disks. Moreover the change of orientations induces the change of signs which 
appears in the process of going from an $A_{\infty}$ category to its opposite category. 
\par
Actually there is an error in the statement of \cite[Theorem 3.54]{fukaya:functor}   related to the sign and orientation.
We will explain the error and its correction in Subsection \ref{subsec:corsign}.
 \end{rem}
% \begin{exm}
% We consider the case $\frak b = 0$, $b$ has degree $1$.
%\begin{equation}\label{forml261}
 %\frak{PO}(b) = \sum_{\beta,k} \exp(\partial \beta \cap \theta_L) \int_L\frak m_{k,\beta}(b^+,\dots,b^+).
%\end{equation}
 %is the potential function for $L \subset (X,J)$.  Here $\beta$ is a class with Maslov index $2$.
% We regard $L \subset (X,-J)$. Then $\theta'_L = - \theta_L$ is used instead of $\theta_L$.
% The potential function becomes
%\begin{equation}\label{forml262}
% \frak{PO}'(b) = \sum_{-\beta,k} \exp(\partial (-\beta) \cap \theta'_L) \int_L\frak m^{\rm op}_{k,-\beta}(b^+,\dots,b^+).
%\end{equation}
%Since Maslov index of $\beta$ is $2$ we have
%$$
%\frak m_{k,\beta}(b^+,\dots,b^+) = \frak m^{\rm op}_{k,-\beta}(b^+,\dots,b^+).
%$$
%Moreover 
%$$
%\partial \beta \cap \theta_L = \partial (-\beta) \cap \theta'_L.
%$$
%Therefore (\ref{forml261}) coincides with (\ref{forml262}).
% \end{exm}
We write 
$$
-{\bf L} = \{(L,-\theta_L,b^*) \mid (L,\theta_L,b) \in {\bf L}\}.
$$
We recall the following version of K\"unneth theorem.
(See \cite{lino} for another version of K\"unneth theorem in Lagrangian Floer theory.)
\begin{thm}{\rm (\cite{fukaya:functor}, See also \cite{lino}, \cite{LLL}.)}\label{themKuneth}
\begin{enumerate}
\item
Let $(L_1,-\theta_{L_1},b^*_1), (L_2,\theta_{L_2},b_2) \in {\bf L}$. Then  there exists 
a bounding cochain $b_1 \times b_2$  of $(L_1 \times L_2,\theta_{L_1} \otimes \theta_{L_2})$, 
which is regarded as a Lagrangian submanifold of $-X\times X$.
We have:
$$
\frak{PO}_{\frak b^{(2)}}(b_1 \times b_2) = \frak{PO}_{\frak b}(b_2) - \frak{PO}_{\frak b}(b_1).
$$
\item
We put 
$$
{\bf L}^2 = \{(L_1 \times L_2,\theta_{L_1} \otimes \theta_{L_2},b_1 \times b_2) \mid 
(L_1,-\theta_{L_1},b^*_1), (L_2,\theta_{L_2},b_2) \in {\bf L}\}.
$$
Then there exists a cyclic and unital filtered $A_{\infty}$ category, 
the set of whose objects is ${\bf L}^2$.
We denote this category as $\sL^2$.
\item
There exists a filtered $A_{\infty}$ bifunctor ${\mathscr{KU}} : \sL \times \sL \to \sL^2$.\index[syindex]{KUfrak@${\mathscr{KU}}$}
It sends a pair of objects $(L_1,-\theta_{L_1},b^*_1), (L_2,\theta_{L_2},b_2)$
to $(L_1 \times L_2,\theta_{L_1} \otimes \theta_{L_2},b_1 \times b_2)$.
\item Let ${\bf U}$ be a finite set of pairs of Lagrangians and its bounding cochains 
in $-X \times X$ such that for $(U,\theta_U,b_U) \in \bf U$ and  $L_i \in \mathscr L$, 
$U$ is transversal to $L_1 \times L_2$.
We put ${\bf F} = {\bf L}^2 \cup {\bf U}$.\index[syindex]{Fbf@${\bf F}$} Let $\Fuk$ be the cyclic and unital filtered $A_{\infty}$ category
the set of whose objects is ${\bf F}$.
Then the following two filtered $A_{\infty}$ modules over $\Fuk$ are homotopy equivalent for each 
$(L_1,-\theta_{L_1},b^*_1), (L_2,\theta_{L_2},b_2) \in {\bf L}$.
\begin{enumerate}
\item
For each $(U,\theta_U,b_U) \in {\bf F}$ it associates 
$$
CF^{\rm can}((L_1,\theta_{L_1},b_1),(U,\theta_U,b_U);(L_2,\theta_{L_2},b_2)),
$$
which is the correspodence tri-module.  (We move $(U,\theta_U,b_U)$ and then it gives a left filtered $A_{\infty}$ modules over $\Fuk$, 
for a fixed $(L_1,-\theta_{L_1},b^*_1)$ and $(L_2,\theta_{L_2},b_2)$.)
This left module over $\Fuk$ is induced from the tri-module $\mathscr{CF}(\sL,\Fuk;\sL)$ defined in 
\cite[Theorem 5.25]{fukaya:functor}
that is a left $\sL^{\rm op}$,  $\Fuk$ and right $\sL$ trimodule.
\item
For each $(U,\theta_U,b_U) \in {\bf F}$ it associates a filtered $A_{\infty}$ left  module
$$
CF^{\rm can}((U,\theta_U,b_U),(L_1 \times L_2,\theta_{L_1} \otimes \theta_{L_2},b_1 \times b_2)).
$$
This is nothing but the left Yoneda module associated to an object $(L_1 \times L_2,\theta_{L_1} \otimes \theta_{L_2},b_1 \times b_2)$ 
of $\Fuk$.
\end{enumerate}
\end{enumerate}
$\mathcal{KU}$ is called {\rm K\"unneth bifunctor} \index{K\"unneth bifunctor}\index[syindex]{KUKU@$\mathcal{KU}$}
The trimodule in item {\rm (a)} is called {\rm K\"unneth trimodule}. \index{K\"unneth trimodule} 
\end{thm}
\begin{rem}
$  $ \par
\begin{enumerate}
\item
Actually in \cite{fukaya:functor} only the case $b_i$ is a bounding cochain with $\frak{PO}(b_i) = 0$
is claimed. However as is already explained in \cite{LLL}, we can include the case of 
weak bounding cochain without changing the proof so much.
\item  We will explain some part of the construction that enters in our proof of
Theorem \ref{themKuneth} later in Subsection \ref{subsec;quilt}.
\item
The error mentioned in Remark \ref{rem261} affects the proof of this 
theorem. 
The way to modify the proof is also explained in Subsection \ref{subsec:corsign}.
\end{enumerate}
\end{rem}

The main part of the proof of Theorem \ref{maintheorem3.3} is the next theorem.
Let ${\bf L}$ be as in Theorem \ref{themKuneth}
and $\sU = \{(U,\theta_U,b_U)\}$ consist of a single embedded Lagrangian $U$ equipped with a bounding cochain.
We then consider the open-closed map:
$$
\frak p_{\sL^2} :  HH_*(\sL^2;\Lambda) \to H_*(X^2;\Lambda)
$$
and the closed-open map:
$$
\frak q_{\mathcal U} : H^*(X^2;\Lambda) \to HH^*(\mathcal U;\Lambda).
$$
\begin{thm}\label{thm253}
If $\frak p : HH_*(\sL;\Lambda) \to H_*(X;\Lambda)$ is an isomorphism then 
$\frak q_{\mathcal U} \circ \frak p_{\sL^2}$ hits the unit.
\end{thm}

\begin{proof}[Theorem \ref{thm253} $\Rightarrow$ Theorem \ref{maintheorem3.3}]
We consider the diagnal $\Delta = \{(x,x) \in -X \times X \mid x \in X\}$.
As is proved in \cite{fooo:inv}, we may take perturbations so that $0$ is a bouding cochain of $\Delta$.
We put ${\bf U} = \{(\Delta,0)\}$ and consider ${\bf F}= {\bf L}^2 \cup {\bf U}$ etc. as in Theorem \ref{maintheorem2}.

Theorem \ref{thm253} and the proof of Theorem \ref{maintheorem2} 
 imply that the diagonal  $(\Delta,0)$ is contained 
in the category which is split generated by $\sL^2$.
Since the diagonal  $(\Delta,0)$ induces an identity functor $\sL \to \sL$ 
the cohomological smoothness follows from the existence of 
filtered $A_{\infty}$ functor
$$
\Fuk \to {\mathcal{FUNC}}(\sL,\sL)
$$
to the functor category. (See \cite[Theorem 1.7]{fukaya:functor}. See also \cite{mww}, \cite{AB}.))
\end{proof}
Other than Theorem  \ref{maintheorem3.3}, 
Theorem \ref{thm253} has the following application.

\begin{cor}
In the situation of Theorem $\ref{maintheorem1}$
there exists a number $D$ such that  any
$(U,\theta_U,b_U)$  is obtained from objects of $\sL$
by taking shift, adding idempotents,  and 
taking cones iteratedly; such that one takes 
 cones  at most $D$ times.
\end{cor}

\begin{proof}
Theorems \ref{thm253} and Proposition \ref{prop:finite_generation_time}
imply that the diagonal in $-M \times M$ 
is obtained from an element of $\sL^2$ 
by the above process (i)(ii)(iii) finitely many times, $D$.
Using the fact that the diagonal induces 
the identity morphism $\sL \to \sL$ 
this implies the corollary.
\end{proof}

The rest of this section is devoted to the proof of Theoerm \ref{thm253}.

\subsection{Pairing on Hochschild homologies via $A_{\infty}$ bi-modules.}
\label{subsec;pair}

We first generalize the pairing $Z$ in Subsection \ref{inproZ} slightly.
Let $\Cat$ and $\Dat$ be unital and filtered $A_{\infty}$ categories which are curvature free 
and $\cP$ a $\Cat-\Dat$-bi-module.
It associates for each 
$c \in {\rm Ob}(\Cat)$ and $d \in {\rm Ob}(\Dat)$
a free $\Lambda_0$ module
$
\cP(c,d)
$.
For each
$c,c' \in {\rm Ob}(\mathcal C)$, $d,d' \in {\rm Ob}(\mathcal D)$ and for $r,s \ge 0$ it associates 
structure operations\index[syindex]{mPr1s@$\frak m^{\cP}_{r|1|s}$}
$$
\frak m^{\cP}_{r|1|s} : B_r\Cat[1](c',c) \otimes  \cP(c,d)  \otimes B_{s}\Dat[1](d,d')
\to  \cP (c',d')
$$ 
which is a $\Lambda_0$ module homomorphism, of degree $1$, preserves the energy filtration and satisfies the $A_{\infty}$ relation
(\ref{eq:bi-module_A_oo_equation}).
\par
We say that $\cP$ is of {\em finite type}\index{finite type} if $\cP(c,d)$ is finitely generated (free) $\Lambda_0$ module 
for each $c,d$.
 
\begin{defn}
Let $\cP$ be a $\Cat$-$\Dat$ bi-module of finite type.
We define a $\Dat$-$\Cat$   bi-module $\cP^{\vee}$ as follows.
We call it the dual bi-module.\index{dual bi-module}
\begin{enumerate}
\item $\cP^{\vee}(d,c) = \Hom_{\Lambda_0}(\cP(c,d),\Lambda_0)$.
\item
$$
\frak m^{\cP^{\vee}}_{s|1|r}({\bf y},f^{\vee},{\bf x})(f) = (-1)^{\maltese} f^{\vee}(\frak m^{\cP}_{r|1|s}({\bf x},f,{\bf y})).
$$
Here $f \in \cP(c,d)$, $f^{\vee} \in \cP^{\vee}(d',c')$  and $\frak m^{\cP^{\vee}}_{s|1|r}$ in the 
left hand side is the structure operation of $\cP^{\vee}$.
\end{enumerate}
The sign $\maltese$ is by Koszul rule. It is straightforward to check the $A_{\infty}$
relation for $\cP^{\vee}$.
\end{defn}
Suppose $\cP$ is a $\Cat$-$\Dat$ bi-module of finite type.
In the rest of this subsection we will define a bilinear map:\index[syindex]{ZP@$Z_{\cP}$}
\begin{equation}\label{formula22}
Z_{\cP} : HH_*(\Cat) \otimes HH_*(\Dat) \to \Lambda_0.
\end{equation}
We begin with certain notations.
Let 
$$
{\bf x} = x_0 \otimes \dots \otimes x_r \in CH_{r+1}(\Cat),
\quad
{\bf y} = y_0 \otimes \dots \otimes y_{s} \in CH_{s+1}(\Dat).
$$
Here $x_i \in \Cat(c_i,c_{i+1})$, $y_i \in \Dat(d_i,d_{i+1})$,
$c_i,d_i$ are objects of $\Cat$, $\Dat$ and 
$c_{k+1} = c_0$, $d_{\ell+1} = d_0$ by convention.
We take a basis $\{f_{ij;a}\mid a \in I_{ij}\}$ of $\cP(c_i,d_j)$ and its dual basis
$\{f^{\vee}_{ij;a}\mid a \in I_{ij}\}$ of $\cP^{\vee}(d_j,c_i)$.
\par
\begin{figure}[h]
\centering
\includegraphics[scale=0.6]{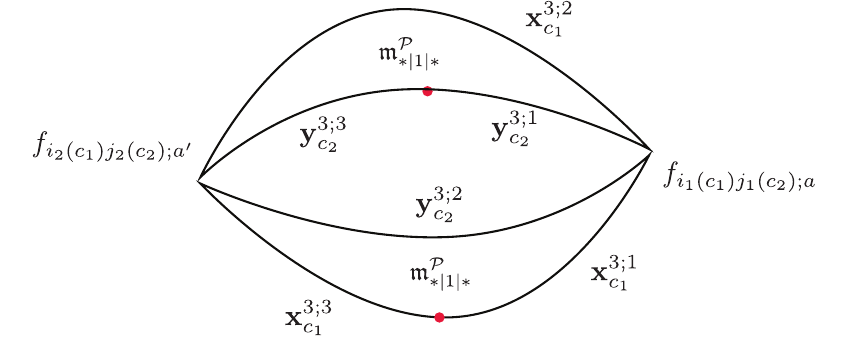}
\caption{$Z_{\cP}$}
\label{Figure13-1new}
\end{figure}

We take $i_1(c_1),i_2(c_1)$, $j_1(c_2),j_2(c_2)$ such that
$$
{\bf x}^{3;2}_{c_1} \in B\Cat[1](c_{i_1(c_1)},c_{i_2(c_1)}),
\quad
{\bf y}^{3;2}_{c_2} \in B\Dat[1](d_{j_1(c_2)},d_{j_2(c_2)}).
$$
Here we use the Sweedler's notation.
Now we define\index[syindex]{ZP@$Z_{\cP}$}
\begin{equation}\label{formula23}
\aligned
&Z_{\cP}(\text{\bf x},\text{\bf y}) \\
&= \sum_{c_1,c_2}\sum_{a,a'}
(-1)^{\maltese}
\langle 
\frak m^{\cP}_{*|1|*}(\text{\bf x}_{c_1}^{3;3} \otimes \text{\bf x}_{c_1}^{3;1},f_{i_1(c_1)j_1(c_2);a},
\text{\bf y}_{c_2}^{3;2}),f_{i_2(c_1)j_2(c_2);a'}^{\vee}
\rangle
\\
&\qquad\qquad\qquad\qquad\langle 
\frak m^{\cP}_{*|1|*}(\text{\bf x}_{c_1}^{3;2},f_{i_2(c_1)j_2(c_2);a'},
\text{\bf y}_{c_2}^{3;3}\otimes \text{\bf y}_{c_2}^{3;1}),f_{i_1(c_1)j_1(c_2);a}^{\vee}
\rangle.
\endaligned
\end{equation}
See Figure \ref{Figure13-1new}.
Here $a,a'$ run $a \in I_{i_1(c_1)j_1(c_2)}$, $a' \in I_{i_2(c_1)j_2(c_2)}$.
The sign $\maltese$ is by Koszul rule and is defined in the same way as Definition \ref{defnZ}.
It is easy to see that (\ref{formula23}) is independent of the choice of the basis $\{f_{ij;a}\}$.

\begin{lem}\label{lem23}
$$
Z_{\cP}(\delta_H \text{\bf x},\text{\bf y}) 
+ (-1)^{\deg' \text{\bf x}}Z_{\cP}(\text{\bf x},\delta_H \text{\bf y})  = 0.
$$
where $\delta_H$ is Hochschild derivative.
\end{lem}
The proof is the same as Lemma \ref{lem205}.
Actually the bilinear map $Z_{\cP}$ coincides with $Z$ in case $\cP$ is a
diagonal bi-module over $\Cat$-$\Cat$.
\par
Lemma \ref{lem23} implies that $Z_{\cP}$ induces the bilinear map (\ref{formula22}).
\par
Let $\Phi : \Dat_1 \to \Dat_2$ be a quasi-isomorphism of filtered $A_{\infty}$ categories
and let $\cP$ be a $\Cat$-$\Dat_2$ bi-module.
Then  we can pull it back by $\Phi$ to obtain a 
$\Cat$-$\Dat_1$ bi-module  $\Phi^*\cP$.
The next lemma is easy to show.

\begin{lem}\label{lem256}
For ${\bf x} \in HH_*(\Cat)$ and ${\bf y} \in HH_*(\Dat_1)$
we have
$$
Z_{\Phi^*\cP}(\text{\bf x},\text{\bf y}) = Z_{\cP}({\bf x},\Phi_*({\bf y})).
$$
\end{lem}
The next lemma is also easy to show.
\begin{lem}\label{lem257}
If $\cP$ is quasi-isomorphic to $\cP'$ then
$Z_{\cP} = Z_{\cP'}$ on $HH_*(\Cat) \otimes HH_*(\Dat)$.
\end{lem}

\subsection{Pairing and K\"unneth functor.}
\label{subsec;pairandKu}

We next apply the construction of Subsection \ref{subsec;pair} in the  situation of Subsection \ref{subsec:smoothstate}.
We first define:
\begin{defn}
Let $\Cat$ be a unital and filtered $A_{\infty}$ category.
A left $\Cat$ module $\mathcal P$ associates to each $c \in {\rm Ob}(\Cat)$ a $\Lambda_0$ module $\cP(c)$  and to each $c, c' \in {\rm Ob}(\Cat)$  
structure operations
$$
\frak m^{\cP}_{r|1|} : B_r\Cat[1](c',c) \otimes  \cP(c) \to  \cP (c'),
$$
which satisfy the $A_{\infty}$ relation. 
\par
For two  left $\Cat$ modules $\cP_1$, $\cP_2$ a left $\Cat$ module pre-homomorphism $\varphi : \cP_1 \to \cP_2$
assigns a degree $0$ map $\varphi(c) : \cP_1(c) \to \cP_2(c)$ to each $c \in {\rm Ob}(\Cat)$,
which is denoted by $y \mapsto (1_c;y)\varphi$ and 
 assigns for each ${\bf x} \in B\Cat[1](c,c')$
 $$
 \varphi({\bf x}) :  \cP_1(c) \to \cP_2(c')
 $$
 which we denote by $y \mapsto ({\bf x};y)\varphi$.
It satisfies the $A_{\infty}$ relation. We call it a $\Cat$ module homomorphism. 
\par
Here we put the symbol $\varphi$ in the right  so that the totological bi-module structure 
which we will define below looks natural.
\par 
We  define composition of homomorphisms $\varphi : \cP_1 \to \cP_2$, 
$\psi : \cP_2 \to \cP_3$ by %\marginpar{The suffix of the next formula is corrected.  KF 2026 Feb}
$$
({\bf x};y)(\varphi \circ \psi) = \sum_c ({\bf x}_c^{2;1};({\bf x}_c^{2;2};y)\varphi)\psi.
$$
\par
The DG category  $\mathbb L\Cat$ of left $\Cat$ modules\index{DG category of left modules} is defined such that
\index[syindex]{Lb1fCat@$\mathbb L\Cat$} its object is a left $\Cat$ module and morphism is a left $\Cat$ module homomorphism.
\end{defn}
There is a left $\Cat$ right $\mathbb L\Cat$ bi-module $\cQ\Cat^{\rm taut}$,  the {\em tautological bi-module}\index{tautological bi-module}
\index[syindex]{Qtot@ $\cQ\Cat^{\rm taut}$} which we define as follows.
\begin{defn}
$  $ \par
\begin{enumerate}
\item
For $c \in {\rm Ob}(\Cat)$ and $\cP \in {\rm Ob}(\mathbb L\Cat)$
we define a $\Lambda_0$ module
$$
\cQ\Cat^{\rm taut}(c,\cP) = \cP(c).
$$
\item
For ${\bf x} \in B_k\Cat[1](c,c')$, $\cP \in {\rm Ob}(\mathbb L\Cat)$, $y \in \cQ\Cat^{\rm taut}(c',\cP)$
($k=0$ is included), 
we put\index[syindex]{mxQtot@${\frak m}^{\cQ^{\rm taut}}$}
$$
{\frak m}^{\cQ\Cat^{\rm taut}}_{k|1|0}({\bf x};y) :=  \frak m_{k|1}^{\cP}({\bf x};y).
$$
Here the left hand side is a structure operation of $\cQ\Cat^{\rm taut}$ and the right hand side 
is a structure operation of $\cP$. 
\item For  ${\bf x} \in B_k\Cat[1](c,c')$, $\cP \in {\rm Ob}(\mathbb L\Cat)$, $y \in \cQ\Cat^{\rm taut}(c',\cP)$,
$\varphi \in \mathbb L\Cat(\cP,\cP')$ ($k=0$ is included), 
we put
$$
{\frak m}^{\cQ\Cat^{\rm taut}}_{k|1|1}({\bf x};y;\varphi) = ({\bf x};y)\varphi.
$$
\item For $r > 1$ we define ${\frak m}^{\cQ\Cat^{\rm taut}}_{k|1|r} = 0$.
\end{enumerate}
It is easy to check that they satisfy the axiom of a bi-module.
\par
When $\Dat$ is a full subcategory of $\mathbb L\Cat$ we restrict 
$\cQ\Cat^{\rm taut}$ to a left $\Cat$, right $\Dat$ bi-module.
We denote it by $\cQ\Cat^{\rm taut}\vert_{\Dat}$.
\end{defn}
Now we put
$
\Cat = \Fuk
$, that  is, the filtered $A_{\infty}$ category whose object set is ${\bf F} = {\bf L}^2 \cup {\bf U}$.
We also put ${\rm Ob}(\Dat_1) \cong {\rm Ob}(\Dat_2) \cong {\bf L}^2$. 
Here $\Dat_1$,$\Dat_2$  are  full subcategories 
of the DG category  $\mathbb L\Fuk$ 
of the left $\Fuk$ modules defined as follows. The object  $\in {\rm Ob}(\Dat_1)$ associated to  $(L_1,\theta_{L_1},b_1), (L_2,\theta_{L_2},b_2)$  is the left $\Fuk$ module
\index[syindex]{Da1@$\Dat_1$}
\begin{equation}\label{253form}
\aligned
&\Dat_1((L_1,\theta_{L_1},b_1), (L_2,\theta_{L_2},b_2))(U,\theta_U,b_U) \\
&: =  CF^{\rm can}((U,\theta_U,b_U),(L_1 \times L_2,\theta_{L_1} \otimes \theta_{L_2},b_1 \times b_2)).
\endaligned
\end{equation}

The left $\Fuk$ module which $\Dat_2$ associates to  $(L_1,\theta_{L_1},b_1), (L_2,\theta_{L_2},b_2)$  is\index[syindex]{Da2@$\Dat_2$}
\begin{equation}\label{254form}
\aligned
&\Dat_2((L_1,\theta_{L_1},b_1), (L_2,\theta_{L_2},b_2))(U,\theta_U,b_U) \\
&:=CF^{\rm can}((L_1,\theta_{L_1},b_1),(U,\theta_U,b_U);(L_2,\theta_{L_2},b_2)).
\endaligned
\end{equation}
The morphism part of $\Dat_1$ and $\Dat_2$ is obtained by the left $\Fuk $
module structures of the right hand sides of (\ref{253form}), (\ref{254form}), respectively.
Thus $\Dat_1$ and $\Dat_2$ are full DG subcategories of $\mathbb L\Fuk$.

Then Theorem \ref{themKuneth} (4) implies the following.
\begin{lem}\label{lem258}
There exists a quasi-equivalence  $\Phi^1 : \Dat_1 \to \Dat_2$ (which is the identity map on objects) 
such that $(\Phi^1)^*\cQ\Cat^{\rm taut}\vert_{\Dat_2}$ is quasi isomorphic to 
$\cQ\Cat^{\rm taut}\vert_{\Dat_1}$ as $\Fuk$, $\Dat_1$ bimodule.
\end{lem}
Lemmas \ref{lem256},  \ref{lem257},  \ref{lem258} imply:
\begin{lem}
For ${\bf x} \in HH_*(\Fuk)$, ${\bf y} \in HH_*(\Dat_1)$ we have
\begin{equation}\nonumber
Z_{\cQ\Cat^{\rm taut}\vert_{\Dat_1}}({\bf x},{\bf y}) = Z_{\cQ\Cat^{\rm taut}\vert_{\Dat_2}}({\bf x},\Phi^1_*({\bf y})).
\end{equation}
Here $Z_*$ is the pairing in Subsection $\ref{inproZ}$.
\end{lem}
We next consider the Yoneda embedding:
$$
\Phi_2 : \sL^2 \to \Dat_1.
$$
The left $\Fuk$ right $\sL^2$ 
 bimodule $\Phi_2^* \cQ\Cat^{\rm taut}\vert_{\Dat_1}$
is quasi isomorphic to the restriction of the diagonal module of 
$\Fuk$ to a left $\Fuk$ right $ \sL^2$ bimodule.  Therefore Lemmas \ref{lem256},  \ref{lem257} imply:
\begin{lem}
For ${\bf x} \in HH_*(\Fuk)$, ${\bf y} \in HH_*(\sL^2)$ we have
$$
Z_{\cQ\Cat^{\rm taut}\vert_{\Fuk,\Dat_1}}({\bf x},{\bf y}) = Z({\bf x},{\bf y}).
$$
Here the right hand side is the pairing in Subsection $\ref{inproZ}$ applied to $\Fuk$.
\end{lem}
We apply Theorem \ref{ZisPD} to the right hand side.
We thus conclude:
\begin{prop}
There exists a weak equivalence $\Phi : \Dat_1 \to \sL^2$ such that
for ${\bf x} \in HH_*(\Fuk)$, ${\bf y} \in HH_*(\Dat_2)$ the next equality holds.
\begin{equation}\label{form253}
Z_{\cQ\Cat^{\rm taut}\vert_{\Dat_2}}({\bf x},{\bf y}) = \langle \frak p^{\bf b}({\bf x}),\frak p^{\bf b}(\Phi_*({\bf y})) \rangle_{{\rm PD}_{X\times X}}.
\end{equation}
\end{prop}
Note that the right hand side of (\ref{form253}) is
$$
\pm \langle {\bf x},(\frak q^{\bf b}\circ \frak p^{\bf b})(\Phi_*({\bf y}))\rangle_{\rm HH},
$$
where the sign $\pm$ depends only on $\deg'{\bf x}$, $\deg'{\bf y}$ and $n$.
Thus to prove Theorem \ref{thm253} it suffices to prove the following:
\begin{prop}\label{prop2512}
Assume $\frak p^{\bf b} : HH_*(\sL;\Lambda) \to H_*(X;\Lambda)$ is an isomorphism.
Then there exists ${\bf y}_0 \in HH_*(\Dat_2)$ such that
$$
Z_{\Dat_2}({\bf x},{\bf y}_0) 
\begin{cases}
= 1   &\text{If ${\bf x} = {[{\rm vol}_U]} \in H^*(U;\Lambda)$},
\\
= 0    &\text{If ${\bf x} \in HH_*(\mathcal U;\Lambda)$ and does not contain the term $[{\rm vol}_U]$}.
\end{cases}
$$
\end{prop}
We remark that there is a canonical map from the de Rham cohomology $H^*(U;\Lambda_0)$ 
to the Hochishild homology $HH_*(\Dat_2)$.
\par
We prove Proposition \ref{prop2512} in the rest of this section. The proof is based on a variant of the 
Cardy relation.

\subsection{Pseudo-holomorphic quilt and K\"unneth trimodule: review.}
\label{subsec;quilt}

In this subsection we review the moduli space of pseudo-holomorphic quilt 
which is used to prove Theorem \ref{themKuneth}.
Let ${\bf L}$ and ${\bf F} = {\bf L}^2 \cup {\bf U}$ be as in Theorem \ref{themKuneth}.
Let ${\bf L} = \{(L_{\kappa},\theta_{L_{\kappa}},b_{\kappa}) \mid \kappa = 1,\dots,\#{\bf L}\}$,
${\bf F} = \{(U_{\rho},\theta_{U_{\rho}},b_{\rho}) \mid \rho = 1,\dots,\#{\bf F}\}$.
We take sequences $\vec{\kappa}_i = (\kappa_{i,0},\dots,\kappa_{i,K_i})$ for $i=1,2$ 
and $\vec{\rho} = (\rho_0,\dots,\rho_K)$.\footnote{We do {\it not} require $\kappa_{i,K_i+1} = \kappa_{i,0}$ this time.}
$\kappa_{i,j} \in \{1,\dots,\#{\bf L}\}$, $\rho_j \in \{1,\dots,\#{\bf F}\}$.
We also take
$p_{i,j} \in L_{\kappa_{i,{j-1}}} \cap L_{\kappa_{i,{j}}}$ if $L_{\kappa_{i,{j-1}}} \ne L_{\kappa_{i,j}}$
and $p_{i,j} = (p_{i,j}(1),p_{i,j}(2))$ is an element of $\tilde L_{\kappa_{i,j}} \times_X \tilde L_{\kappa_{i,j}} \setminus
\tilde L_{\kappa_{i,j}}$, that is, a self-intersection point if $L_{\kappa_{i,j-1}} = L_{\kappa_{i,j}}$.
We put $\vec p_{i }= (p_{i,1},\dots,p_{i,K_i})$.
$p_{12,j} \in \tilde U_{\rho_{j-1}}  \times_X \tilde U_{\rho_j}$, $\vec p_{12}$ is defined in a similar way.
\par
We also take
$\vec k_i = (k_{i,0},\dots,k_{i,K_i})$ for $i=1,2$ where $k_{i,j} \in \Z_{\ge 0}$,
and $\vec k_{12} = (k_{12,0},\dots,k_{12,K})$ with $k_{12,j} \in \Z_{\ge 0}$.
\par
We also take 
$p_{-\infty} \in \pi_0((\tilde L_{\kappa_{1,K_1}} \times \tilde L_{\kappa_{2,0}}) \times_{X\times X} U_{\rho_K})$,
$p_{+\infty} \in \pi_0((\tilde L_{\kappa_{1,0}} \times \tilde L_{\kappa_{2,K_2}}) \times_{X\times X} U_{0})$,
and let $R(p_{-\infty})$, $R(p_{+\infty})$ be the corresponding connected components.

\begin{figure}[ht]
\centering
\includegraphics[scale=0.4]{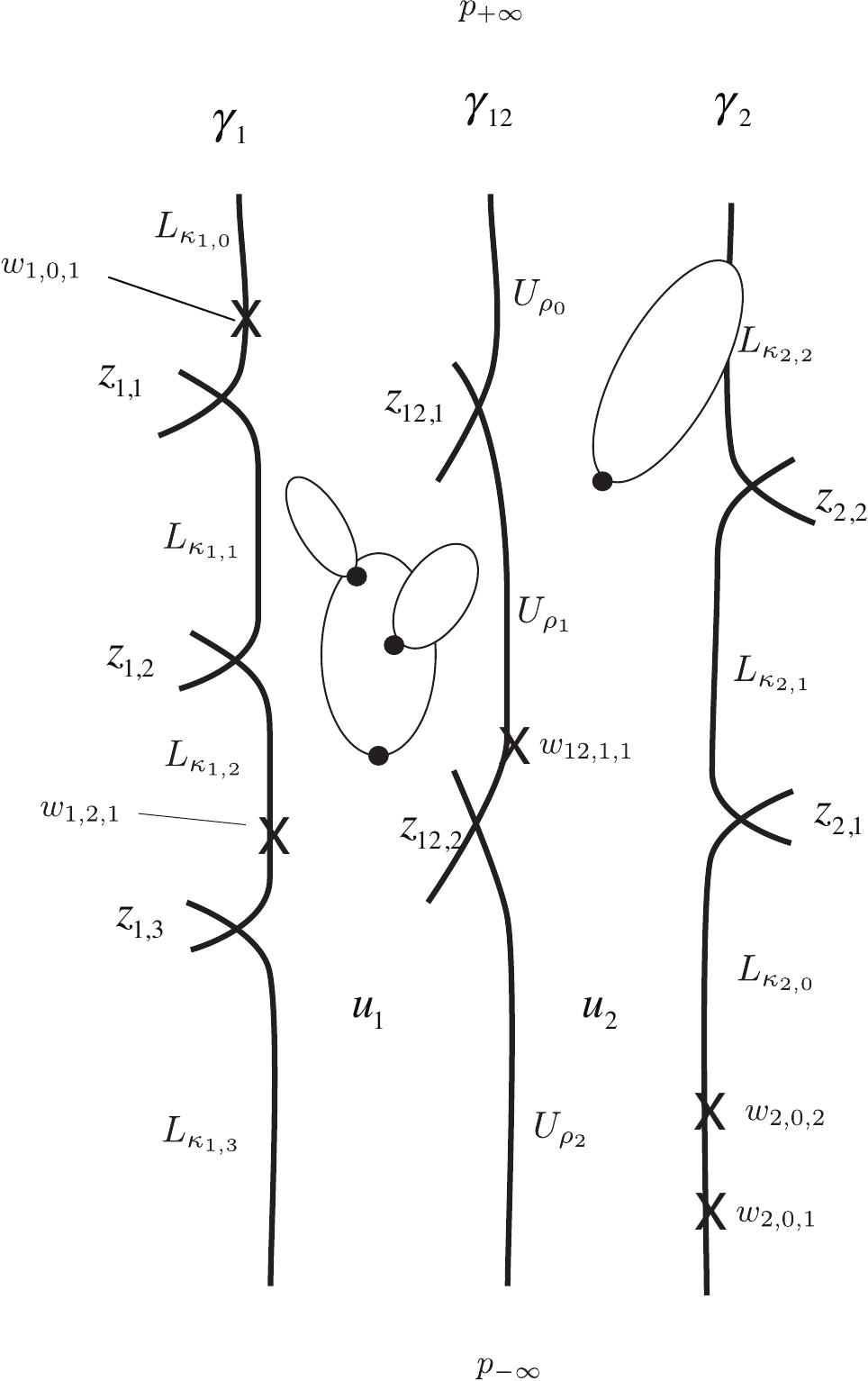}
\caption{An element of $\overset{\circ\circ}{\mathcal M}_{\rm QT}(\vec\kappa_1,\vec\kappa_2,\vec\rho;\vec p_1,\vec p_{2},\vec p_{12};\vec k_1,\vec k_2,\vec k_{12};
p_{\pm \infty};E)$.}
\label{Figure51}
\end{figure}

\begin{defn}\label{def516}{\rm (See \cite[Definition 5.27]{fukaya:functor})}
We consider objects
$$(\Sigma;\vec z^+_1,\vec z^+_2;\vec z_1,\vec z_{12},\vec z_2;\vec\vec{w}_1,\vec\vec{w}_{12},\vec\vec{w}_2;u_1,u_2;\gamma_1,\gamma_{12},\gamma_2)$$
with the following properties. (See Figure \ref{Figure51}.)
\begin{enumerate}
\item
The space 
$\Sigma$ is a bordered Riemann surface of genus $0$ with $\Sigma \supseteq ([-1,1]\times \sqrt{-1}\R)$. 
The closure of $\Sigma \setminus ([-1,1]\times \R)$ is a finite 
union of (maximal)  trees of spheres. We call  its connected component a {\it  tree of sphere components}.
\index{tree of sphere components}
We require that a tree of sphere components intersects 
with 
$[-1,1]\times \sqrt{-1}\R$ at one point, which we call its {\it root}. \index{root} 
All the roots are 
points of $((-1,0) \cup (0,1)) \times \sqrt{-1}\R$.
\footnote{In other words, we require that the root is not on $\{0,\pm 1\} \times \sqrt{-1}\R$.}
\item
Let $\Omega_1$ (resp. $\Omega_2$) be the union of $[-1,0] \times \sqrt{-1}\R$
(resp. $[0,+1] \times \R)$) and the trees of sphere components rooted on it.
\par
$u_i : \Omega_i \to X$ is a smooth map such that $u_1$ is $-J_X$ holomorhpic and 
$u_2$ is $J_X$ holomorphic.
\item
We put
$\vec z_i = (z_{i,1},\dots,z_{i,K_i})$ ($i=1,2$). Then 
$z_{1,j} \in  \{-1\} \times \sqrt{-1}\R$, $z_{2,j} \in  \{1\} \times \sqrt{-1}\R$.
$\vec z_{12} = (z_{12,1},\dots,z_{12,K})$, $z_{12,j} \in \{0\} \times \sqrt{-1}\R$.
\par
If $j_1 < j_2$ then ${\rm Im} (z_{1,j_1}) > {\rm Im} (z_{1,j_2})$, 
${\rm Im} (z_{12,j_1}) > {\rm Im} (z_{12,j_2})$
and ${\rm Im} (z_{2,j_1}) < {\rm Im} (z_{2,j_2})$.
\par
We put $\vert \vec z_i\vert = \{z_{i,1},\dots,z_{i,K_i}\}$.
$\vert \vec z_{12}\vert$ is defined in the same way.
\item
We put $\vec\vec{w}_i = (\vec w_{i,0},\dots,\vec w_{i,K_i})$, $i=1,2$.
Here $\vec w_{1,j} = (w_{1,j,1},\dots,w_{1,j,k_{1,j}})$, 
$w_{1,j,m} \in  \{-1\} \times  \sqrt{-1}({\rm Im} (z_{1,j}),{\rm Im} (z_{1,j-1}))$ ($m=1,\dots,k_{1,j}$) with decreasing ${\rm Im}(w_{1,j,m})$. And
$\vec w_{2,j} = (w_{2,j,1},\dots,w_{2,j,k_{2,j}})$, 
$w_{2,j,m} \in  \{-1\} \times  \sqrt{-1}({\rm Im} (z_{2,j}),{\rm Im} (z_{2,j+1}))$ ($m=1,\dots,k_{2,j}$) with increasing ${\rm Im}(w_{2,j,m})$.
($z_{2,K_{2}+1} = z_{1,0} = +\infty$, $z_{2,0} =z_{1,K_{1}+1}= -\infty$ by convention.)

Moreover
$\vec\vec{w}_{12}= (\vec w_{12,1},\dots,,\vec w_{12,K})$.  
Here $\vec w_{12,j} = (w_{12,j,1},\dots,w_{12,j,k_{12,j}})$, 
$w_{12,j,m} \in  \{-1\} \times  \sqrt{-1}({\rm Im} z_{12,j},{\rm Im} (z_{12,j-1}))$ ($m=1,\dots,k_{12,j}$) with decreasing ${\rm Im}(w_{12,j,m})$.
 \item
The maps
$\gamma_1 : (\{-1\}\times \R) \setminus \vert\vec z_1\vert \to \tilde L_1$, 
$\gamma_2 : (\{1\}\times \R) \setminus \vert\vec z_2\vert \to \tilde L_2$,
$\gamma_{12} : (\{0\}\times \R) \setminus \vert\vec z_{12}\vert \to \tilde L_{12}$
are smooth and satisfy
\begin{equation}
\aligned
i_{L_1}(\gamma_1(z)) = u_1(z) \qquad\qquad &\text{if $z \in (\{-1\}\times \R) \setminus \vert\vec z_1\vert$},  \\
i_{L_2}(\gamma_2(z)) = u_2(z) \qquad\qquad &\text{if $z \in (\{1\}\times \R) \setminus \vert\vec z_2\vert$},  \\
i_{L_{12}}(\gamma_{12}(z)) = (u_1(z),u_2(z)) \qquad\qquad &\text{if $z \in 
(\{0\}\times \R) \setminus \vert\vec z_{12}\vert$}.
\endaligned
\end{equation}
\item
At the points of $\vert\vec z_1\vert$, $\vert\vec z_2\vert$, $\vert\vec z_{12}\vert$, the maps $\gamma_1$, $\gamma_2$, 
$\gamma_{12}$ satisfy the switching condition, determined by $\vec p_1$, $\vec p_2$, $\vec p_{12}$
in the same way as Definition \ref{def211} (4).
\item
When $z \in [-1,1]\times \sqrt{-1}\R$, ${\rm Im}z \to \pm \infty$, 
the following asymptotic boundary condition determined by $p_{\pm \infty}$ is satisfied.
\begin{equation}
\aligned
& \lim_{\tau \to -\infty}(\gamma_1(-1-\sqrt{-1}\tau),\gamma_2(1-\sqrt{-1}\tau),\gamma_{12}(-\sqrt{-1}\tau)) \in R(p_{-\infty}), \\
& \lim_{\tau \to +\infty}(\gamma_1(-1+\sqrt{-1}\tau),\gamma_2(1+\sqrt{-1}\tau),\gamma_{12}(+\sqrt{-1}\tau)) \in R(p_{+\infty}).
\endaligned
\end{equation}
\item 
The object is stable in the sense that the group of automorphisms is a finite group.
\item
$$
-\int_{\Omega_1}u_1^*\omega + \int_{\Omega_2}u_2^*\omega = E.
$$
\end{enumerate}
\par
We  define an equivalence relation $\sim$ among those objects 
in an obvious way. We denote by
$\overset{\circ\circ}{\mathcal M}_{\rm QT}(\vec\kappa_1,\vec\kappa_2,\vec\rho;\vec p_1,\vec p_{2},\vec p_{12};\vec k_1,\vec k_2,\vec k_{12};
p_{\pm \infty};E)$
\index[syindex]{M1QTa1a12@${\mathcal M}_{\rm QT}(\vec\kappa_1,\vec\kappa_2,\vec\rho;\vec p_1,\vec p_{2},\vec p_{12};\vec k_1,\vec k_2,\vec k_{12};
p_{\pm \infty};E)$}
the set of all the equivalence classes of this equivalence relation.
We call its element a {\it pseudo-holomorphic quilt}.\index{pseudo-holomorphic quilt}
\end{defn}

We can compactify the moduli space $\overset{\circ\circ}{\mathcal M}_{\rm QT}(\vec\kappa_1,\vec\kappa_2,\vec\rho;\vec p_1,\vec p_{2},\vec p_{12};\vec k_1,\vec k_2,\vec k_{12};
p_{\pm \infty};E)$.  The compactification ${\mathcal M}_{\rm QT}(\vec\kappa_1,\vec\kappa_2,\vec\rho;\vec p_1,\vec p_{2},\vec p_{12};\vec k_1,\vec k_2,\vec k_{12};
p_{\pm \infty};E)$ carries a Kuranishi structure with corners. 
Its boundary are one of the following three types.
\begin{enumerate}
\item[(bd1)]  The disk bubble at $\{-1\} \times \R$.
\item[(bd2)]  The disk bubble at $\{+1\} \times \R$.
\item[(bd3)]  The disk bubble at $\{0\} \times \R$.
(We call $\{0\} \times \R$ a seam\index{seam}.)
\item[(bd4)]
The object splits into two pieces in the imaginary direction.
\end{enumerate}
See \cite[Theorem 5.43]{fukaya:functor}.\footnote{We need to modify the compactification 
of the moduli space of disk bubbles on the seams to obtain this Kuranishi structure.
This point is discussed in detail in \cite[Section 12]{fukaya:functor}.}
\par
We renumber the marked points as follows.
We consider $\vec z_1 \cup \vec\vec w_1$.
We then renumber them according to the descending order of the imaginary part 
and write them as $z_1(1),\dots,z_1(\mathcal K_1)$,
where $\mathcal K_1 = K_1 + \sum k_{1,j}$ is the total number of the marked points 
on $\{-1\} \times \sqrt{-1}\R$.
We also renumber $\vec z_{12} \cup \vec\vec w_{12}$
to obtain 
$z_{12}(1),\dots,z_{12}(\mathcal K)$ where $\mathcal K = K + \sum k_{12,j}$.
\par
For $\vec z_{2} \cup \vec\vec w_{2}$ we use {\it ascending} order of the imaginary part 
to renumber it and obtain $z_{2}(1),\dots,z_{2}(\mathcal K_2)$.
\par
We define $\kappa'_{i,j}$ so that $z_{i}(j) \in \tilde L_{\kappa'_{i,j-1}} \times_X  \tilde L_{\kappa'_{i,j}} $.
We define $\rho'_{j}$ so that $z_{12}(j) \in \tilde U_{\rho'_{j-1}} \times_X  \tilde U_{\rho'_{j}} $.

\begin{rem}\label{rem2616}
The process here is similar to the one in Subsection \ref{constcyclic}. 
We alert that the role of $\vec\kappa'$ and $\vec\kappa$ is opposite between 
here and Subsection \ref{constcyclic}. 
\end{rem}

We define $L^{\vec{\kappa}'_1,\vec p_1}_{\text{\rm source}}$, $L^{\vec{\kappa}'_2,\vec p_2}_{\text{\rm source}}$
and $U^{\vec{\rho}',\vec p_{12}}_{\text{\rm source}}$ in the same way as Definition \ref{defn927777}.
Namely
$$
L^{\vec{\kappa}'_1,\vec p_1}_{\text{\rm source}} \subseteq \prod_{i=1}^{\mathcal K_1} (\tilde L'_{\kappa_{1,i-1}} \times_{X} \tilde L'_{\kappa_{1,i}})
$$
etc.
\par
The evaluation maps at $z_{i}(j)$ and $z_{12}(j)$ define 
\begin{equation}
\aligned
{\rm ev} = ({\rm ev}_1,{\rm ev}_{12},{\rm ev}_2) 
: &{\mathcal M}_{\rm QT}(\vec\kappa_1,\vec\kappa_2,\vec\rho;\vec p_1,\vec p_{2},\vec p_{12};\vec k_1,\vec k_2,\vec k_{12};
p_{\pm \infty};E) \\
& \to L^{\vec{\kappa}'_1,\vec p_1}_{\text{\rm source}} \times U^{\vec{\rho}',\vec p_{12}}_{\text{\rm source}} \times
L^{\vec{\kappa}'_2,\vec p_2}_{\text{\rm source}}.
\endaligned
\end{equation}
We also have evaluation maps at $\tau \to \pm \infty$,
$$
{\rm ev}_{\pm \infty} : {\mathcal M}_{\rm QT}(\vec\kappa_1,\vec\kappa_2,\vec\rho;\vec p_1,\vec p_{2},\vec p_{12};\vec k_1,\vec k_2,\vec k_{12};
p_{\pm \infty};E)
\to R(p_{\pm \infty}).
$$
We include $\ell_1$ (resp. $\ell_2$) marked points on $\Sigma_1$ (resp. $\Sigma_2$) to define
$$
{\mathcal M}_{{\rm QT},\ell_1,\ell_2}(\vec\kappa_1,\vec\kappa_2,\vec\rho;\vec p_1,\vec p_{2},\vec p_{12};\vec k_1,\vec k_2,\vec k_{12};
p_{\pm \infty};E).
$$
It then also have evaluation maps:
$$
{\rm ev}_{\rm int} : 
{\mathcal M}_{{\rm QT},\ell_1,\ell_2}(\vec\kappa_1,\vec\kappa_2,\vec\rho;\vec p_1,\vec p_{2},\vec p_{12};\vec k_1,\vec k_2,\vec k_{12};
p_{\pm \infty};E)
\to X^{\ell_1+\ell_2}.
$$
We decompose ${\mathcal M}_{{\rm QT},\ell_1,\ell_2}(\vec\kappa_1,\vec\kappa_2,\vec\rho;\vec p_1,\vec p_{2},\vec p_{12};\vec k_1,\vec k_2,\vec k_{12};
p_{\pm \infty};E)$ by the homology class $B$ and will write 
${\mathcal M}_{{\rm QT},\ell_1,\ell_2}(\vec\kappa_1,\vec\kappa_2,\vec\rho;\vec p_1,\vec p_{2},\vec p_{12};\vec k_1,\vec k_2,\vec k_{12};
p_{\pm \infty};B)$.

We define  $BCF(\cL^{\rm form};\vec{\kappa}'_i)$ ($i=1,2$) and
$BCF(\sU^{\rm form};\vec{\rho}^{\,\prime})$ in the same way as (\ref{eq920}), (\ref{eq921}).
Namely
$$
\aligned
BCF(\cL^{\rm form};\vec{\kappa}^{\,\prime}_i) &= \bigoplus_{\vec p_i} \Omega(L^{\vec{\kappa}^{\,\prime}_i,\vec p_i}_{\text{\rm source}}) \otimes \F, \\
BCF(\sU^{\rm form};\vec{\rho}^{\,\prime})  &= \bigoplus_{\vec p_{12}} \Omega(U^{\vec{\rho}^{\,\prime},\vec p_{12}}_{\text{\rm source}})\otimes \F.
\endaligned
$$
We now define 
$$
\aligned
&\frak q^{\text{\rm form}}_{\ell_1,\ell_2;\vec{\kappa}_1,\vec{\kappa}_2,\vec{\rho},\vec p_1,\vec p_2,\vec p_{12};p_{-\infty},p_{+\infty};B} \\
&\qquad : \Omega(X)^{\otimes(\ell_1+\ell_2)}
\otimes BCF(\cL^{\rm form};\vec{\kappa}^{\,\prime}_1) \otimes BCF(\sU^{\rm form};\vec{\rho}^{\,\prime}) 
\\
&\qquad\qquad\qquad \otimes \Omega(R(p_{-\infty}))\otimes BCF(\cL^{\rm form};\vec{\kappa}^{\,\prime}_2) \to \Omega(R(p_{+\infty}))
\endaligned
$$
by (Compare Definition \ref{def827}.)
\begin{equation}\label{form269}
\aligned
&\frak q^{\rm form}_{\star}({\bf g}_1\otimes {\bf g}_2 \otimes {\bf h}_1 \otimes {\bf h}\otimes f \otimes {\bf h}_2 )\\
&=
(-1)^{{\maltese}}{\rm ev}_{+\infty !}\left( {\mathcal M}_{{\rm QT},{\star}}({\star}); {\rm ev}^*_{\rm int}({\bf g}_1\otimes {\bf g}_2)\right. \\
&\qquad\qquad\qquad\qquad \left.  {\rm ev}^*_{1}{\bf h}_1 \wedge {\rm ev}^*_{12}{\bf h} \wedge  {\rm ev}^*_{-\infty}f \wedge {\rm ev}^*_{2}{\bf h}_2   
 \right). 
\endaligned
\end{equation}
Here we use ${\star}$ in place of complicated sequence of symbols for simplicity.
\par
Note that we take an appropriate CF-perturbation on 
$$
{\mathcal M}_{{\rm QT},{\star}}({\star}) = {\mathcal M}_{{\rm QT},\ell_1,\ell_2}(\vec\kappa_1,\vec\kappa_2,\vec\rho;\vec p_1,\vec p_{2},\vec p_{12};\vec k_1,\vec k_2,\vec k_{12};
p_{\pm \infty};B)
$$
to define the integration along the fiber ${\rm ev}_{+\infty !}$.  (We require that the map ${\rm ev}_{+\infty}$ is strongly submersive with respect to this 
CF-perturbation.) The sign $\maltese$ is discussed in Subsubsection \ref{subsubsec:ref}.
\par
We next define $\rho_{\frak b,\theta}(B)$ in a similar way as Definition \ref{def936} as follows.
Let us take an object $(\Sigma;\vec z_1,\vec z_{12},\vec z_2;\vec\vec{w}_1,\vec\vec{w}_{12},\vec\vec{w}_2;u_1,u_2;\gamma_1,\gamma_{12},\gamma_2)$
belonging to the homology class $B$. We then consider
\begin{equation}\label{form2611}
\theta(B): = \sum_{i=1,2}\sum_{j=0}^{\mathcal K_i+1} \int_{z_{i,j-1}}^{z_{i,j}} \gamma_i^*\theta_{L_{\kappa_{i,j}}}  + \sum_{j=0}^{\mathcal K} 
\int_{z_{12,j-1}}^{z_{12,j}} \gamma_{12}^*\theta_{U_{\rho_{j}}} 
\end{equation}
where $z_{2,K_{2}+1} = z_{1,0} = +\infty$, $z_{2,0} =z_{1,K_{1}+1}= -\infty$ etc. by convention.\footnote{We remark that 
$\gamma_1$ and $\gamma_{12}$ is oriented downward and $\gamma_2$ is oriented upward.}
We define:
\begin{equation}\label{form2612}
\rho_{\frak b,\theta}(B) := \exp\left( -\int_{\Sigma_1} u_1^*\frak b + \int_{\Sigma_2} u_2^*\frak b + \theta(B) \right).
\end{equation}
We also put $\omega(B) = E$ where $E$ is as in Definition \ref{def516} (9).
It is easy to see that the right hand side depends only on the homology class $B$.
\par
We put
\begin{equation}\label{CFformform}
\mathscr{CF}^{\rm form}(L_{\kappa_1},U_{\rho};L_{\kappa_2};\F) =
\Omega((\tilde L_{\kappa_1} \times \tilde L_{\kappa_2}) \times_{X\times X} \tilde U_{\rho}) \otimes \F.
\end{equation}
Now we define
$$
\aligned
\frak m^{\text{\rm form},\frak b}_{\vec{\kappa}^{\,\prime}_1,\vec{\kappa}^{\,\prime}_2,\vec\rho^{\,\prime}}
: 
BCF(\cL^{\rm form};\vec{\kappa}^{\,\prime}_1) &\otimes BCF(\sU^{\rm form};\vec{\rho}^{\,\prime}) 
\otimes 
\mathscr{CF}(L_{\kappa'_{1,\mathcal K_1}},U_{\rho'_K};L_{\kappa^{\,\prime}_{2,0}};\F)
\\ &\otimes BCF(\cL^{\rm form};\vec{\kappa}^{\,\prime}_2) 
\to \mathscr{CF}(L_{\kappa'_{1,0}},U_{\rho'_0};L_{\kappa'_{2,K_2}};\F)
\otimes_{\F} \Lambda_0
\endaligned
$$
by
\begin{equation}
\aligned
&\frak m^{\text{\rm form},\frak b}_{\vec{\kappa}_1,\vec{\kappa}_2,\vec\rho}
({\bf h}_1 \otimes {\bf h}\otimes f \otimes {\bf h}_2 ) \\
&=
\sum_{\vec p_1,\vec p_2,\vec p_{12},B,\ell_1,\ell_2}T^{\omega(B) }
\frac{\rho_{\frak b,\theta}(B)
}{\ell!}\\
&\qquad \qquad \qquad
\frak q^{\rm form}_{\ell_1,\ell_2;\vec{\kappa}_1,\vec{\kappa}_2,\vec{\rho},\vec p_1,\vec p_2,\vec p_{12};\star;B}({\frak b}^{\ell_1}\otimes {\frak b}^{\ell_2} \otimes {\bf h}_1 \otimes {\bf h}\otimes f \otimes {\bf h}_2 ).
\endaligned
\end{equation}
This operation becomes a structure operation of a left $(\cL^{\rm form})^{\rm op}$, $\sU^{\rm form}$, right $\cL^{\rm form}$ filtered tri-module.
In fact the $A_{\infty}$ relation for the trimodule is the sum of
$$
\aligned
&\frak m(({\bf h}_1)^{3;1}_c\frak m(({\bf h}_1)^{3;2}_c) ({\bf h}_1)^{3;3}_c; {\bf h};f;{\bf h}_2), \\
&\frak m({\bf h}_1;{\bf h};f;({\bf h}_2)^{3;1}_c \frak m(({\bf h}_2)^{3;2}_c) ({\bf h}_2)^{3;3}_c),  \\
&\frak m({\bf h}_1; {\bf h}^{3;1}_c\frak m({\bf h}^{3;2}_c) {\bf h}^{3;3}_c;f;{\bf h}_2), \\
&\frak m(({\bf h}_1)^{2;1}_{c_1};{\bf h}^{2;1}_{c} \frak m(({\bf h}_1)^{2;2}_{c_1};{\bf h}^{2;2}_{c};f;({\bf h}_2)^{2;1}_{c_2})
({\bf h}_2)^{2;2}_{c_2}),
\endaligned
$$
with Koszul sign.

Lines 1-4 of the above formula correspond to (bd1)-(bd4) of the description of the boundary of 
the moduli space
${\mathcal M}_{\rm QT}(\vec\kappa_1,\vec\kappa_2,\vec\rho;\vec p_1,\vec p_{2},\vec p_{12};\vec k_1,\vec k_2,\vec k_{12};
p_{\pm \infty};E)$, respectively. See \cite[Proposition 5.51]{fukaya:functor} for detail.
\par
We can include forgetability data $\vec{\frak f}$ in the same way as in Section \ref{sec:unit} 
and then our trimodule becomes unital. (See \cite[Definition 5.10]{fukaya:functor} for the 
unitality of the bi-module. Unitality of the trimodule is defined in the same way.)
We then reduce the structure to the cohomologies using the $A_{\infty}$ 
functor $\frak f$ in Subsections \ref{subsec:red} and \ref{hunitcansubsec}.
Also we replace (\ref{CFformform}) by 
\begin{equation}\label{CFcancan}
\mathscr{CF}(L_{\kappa},U_{\rho};L_{\kappa'};\F) =
H((\tilde L_{\kappa} \times \tilde L_{\kappa'}) \times_{X\times X} \tilde U_{\rho};\F)
\end{equation}
in the same way as \cite[Theorem 5.4.18]{fooo09}.\footnote{\cite[Theorem 5.4.18]{fooo09} is the case of 
a bi-module over a filtered $A_{\infty}$ algebra.  It can however be easily adapted to our case of 
a trimodule over filtered $A_{\infty}$ categories.}
We thus obtain a left $\cL^{\rm op}_{\rm curv}$, $\Fuk_{\rm curv}$ right $\cL_{\rm curv}$ {\it curved} trimodule, 
which we describe below.\footnote{In \cite{fukaya:functor} reduction to the canonical model was not discussed since 
we did not use it there. Here we use bilinear map $Z_{\cP}$ for which we need finite dimensionality 
of the morphism modules.}(We omit the symbol of the bulk class $\frak b$ below since we fix it.)
\begin{notation}\label{not2020}
We denote an object $(L_{\kappa},\theta_{L_{\kappa}})$ of $\cL_{\rm curv}$ or of $\cL^{\rm op}_{\rm curv}$
by $c_{\kappa}$. We also denote an object $(U_{\rho},\theta_{U_{\rho}})$
of $\Fuk_{\rm curv}$ by $d_{\rho}$.  (Here $(U_{\rho},\theta_{U_{\rho}}) \in \mathscr U$.)
Another type of  object $(L_{\kappa_1} \times L_{\kappa_2},\theta_{L_{\kappa_1}} \otimes \theta_{L_{\kappa_2}} )$
of $\Fuk_{\rm curv}$ is denoted by $d_{(\kappa_1,\kappa_2)}$.
We remark that the set of objects of $\Fuk_{\rm curv}$ is $\mathscr L^2 \cup \mathscr U$.
We use the symbols $\sigma$ in $d_{\sigma}$ when $\sigma$ can be both  $\rho$ and $(\kappa_1,\kappa_2)$.
\end{notation}
We denote our curved trimodule by $\mathscr{CF}_{\rm curve}$.
For objects $c_{\kappa_1},c_{\kappa_2},d_{\sigma}$ it associates a $\Lambda_0$ module\index[syindex]{CFcurve@$\mathscr{CF}_{\rm curve}(c_{\kappa_1},d_{\sigma};c_{\kappa_2})$}
$$
\mathscr{CF}_{\rm curve}(c_{\kappa_1},d_{\sigma};c_{\kappa_2})
$$
in a similar way as  (\ref{CFcancan}) $\otimes \Lambda_0$.
If 
\begin{equation}
\aligned
&{\bf x} \in B_{k_1}\cL^{\rm op}[1](c_{\kappa_1},c_{\kappa'_1}), \\
&{\bf w} \in B_k\sU[1](d_{\sigma},d_{\sigma'}) \\
&{\bf y} \in B_{k_2}\cL[1](c_{\kappa_2},c_{\kappa'_2})
\endaligned
\end{equation}
then the structure operation of the trimodule is
$$
z \mapsto \frak m^{\rm curve}_{k_1,k | 1  | k_2}({\bf x},{\bf w};z;{\bf y}).
$$
This is a map:
$$
\mathscr{CF}_{\rm curve}(c_{\kappa'_1},d_{\sigma'};c_{\kappa_2}) \to \mathscr{CF}_{\rm curve}(c_{\kappa_1},d_{\sigma};c_{\kappa'_2}).
$$
We now use the next proposition to include the bounding cochain.
We remark that
$$
\mathscr{CF}_{\rm curve}(c_{\kappa_1},d_{(\kappa_1,\kappa_2)};c_{\kappa_2})
=
H((\tilde L_{\kappa_1} \times \tilde L_{\kappa_2}) \times_{X\times X} (\tilde L_{\kappa_1} \times \tilde L_{\kappa_2});\Lambda_0).
$$
It contains a fundamental class $[\tilde L_{\kappa_1} \times \tilde L_{\kappa_2}]$ of the diagonal component.
We denote it by
${\bf 1}_{(\kappa_1,\kappa_2)}$
\begin{prop}
Let $b_{\kappa_1}$, $b_{\kappa_2}$ be (weak) bounding cochains of $c_{\kappa_1}$ and $c_{\kappa_2}$, respectively.
Then there exists a unique (weak) bounding cochain $b_{\kappa_1} \times b_{\kappa_2}$ of $d_{(\kappa_1,\kappa_2)}$
with the following properties.
\begin{equation}
\sum_{k_1,k,k_2 = 0}^{\infty} \frak m^{\rm curve}_{k_1,k | 1  | k_2}(b_{\kappa_1}^{k_1},(b_{\kappa_1}\times b_{\kappa_2})^{k};{\bf 1}_{(\kappa_1,\kappa_2)};b_{\kappa_2}^{k_2}) = 0.
\end{equation}
\begin{equation}\label{formnew2515}
\frak{PO}(b_{\kappa_1} \times b_{\kappa_2}) = \frak{PO}(b_{\kappa_2}) - \frak{PO}(b_{\kappa_1}).
\end{equation}
\end{prop}
\begin{proof}
In the case when $\frak{PO}(b_{\kappa_2})=  \frak{PO}(b_{\kappa_1}) = 0$ this is \cite[Proposition 16.11]{fukaya:functor}.
In the general case of weak bounding cochains this is \cite[Proposition 2.10]{LLL}.
\end{proof}
We remark that $\frak{PO}(b_{\kappa_1})$ in  (\ref{formnew2515}) is the potential value as a
bounding cochain of $\cL$, which is minus of the potential value as a
bounding cochain of $\cL^{\rm op}$.
\par
Now using $b_{\kappa_1} \times b_{\kappa_2}$ as a bounding cochain of $d_{(\kappa_1,\kappa_2)}$
we obtain a weakly curvature free cyclic, unital and filtered $A_{\infty}$ category 
$\Fuk$.
Also we have a left $\cL^{\rm op}$, $\Fuk$ right $\cL$ trimodule $\mathscr{CF}$.
For objects $c_{\kappa_1},c_{\kappa_2},d_{\sigma}$ the trimodule $\mathscr{CF}$ associates
a $\Lambda_0$ module 
$$
\mathscr{CF}(c_{\kappa_1},d_{\sigma};c_{\kappa_2})
$$
which is the same as $\mathscr{CF}_{\rm curve}(c_{\kappa_1},d_{\sigma};c_{\kappa_2})$
as $\Lambda_0$ module.
\par
For ${\bf x} \in B_{k_1}\cL^{\rm op}[1](c_{\kappa_1},c_{\kappa'_1})$, 
${\bf w} \in B_k\sU[1](d_{\sigma},d_{\sigma'})$,
${\bf y} \in B_{k_2}\cL[1](c_{\kappa_2},c_{\kappa'_2})$ 
and $z \in \mathscr{CF}(c_{\kappa'_1},d_{\sigma'};c_{\kappa_2})$ 
we have
$$
\frak m^{\bf b}_{k_1,k | 1  | k_2}({\bf x},{\bf w};z;{\bf y}) \in \mathscr{CF}(c_{\kappa_1},d_{\sigma};c_{\kappa'_2}).
$$
The operations $\frak m^{\bf b}_{k_1,k | 1  | k_2}$ is obtained by deforming the operations\index[syindex]{mbk1k1k2@$\frak m^{\bf b}_{k_1,k | 1  | k_2}$}
$\frak m^{\rm curve}_{k_1,k | 1  | k_2}$ via bounding cochains.   
The operations $\frak m^{\bf b}_{k_1,k | 1  | k_2}$ satisfy the $A_{\infty}$ relation and are unital.
\par
As we already discussed in Subsection \ref{subsec;pairandKu}, the trimodule $\mathscr{CF}$ associates to $c_{\kappa_1} , c_{\kappa_2}$
an object   ${\rm Ob}(\Cat_2)$ of the $DG$ category $\mathbb L\Fuk$ that is the left $\Fuk$ module 
which associates $\mathscr{CF}(c_{\kappa_1},d_{\sigma};c_{\kappa_2})$ to $d_{\sigma}$.
Lemma \ref{lem258} claims that it is quasi-isomorphic to 
the left module which is a Yoneda module associated to $d_{(c_{\kappa_1},c_{\kappa_2})}$.
This is actually the way to construct the K\"unneth bi-functor 
$\mathscr{KU} : \cL \times \cL \to \Fuk$.
\par
This is the summary of the proof of Theorem \ref{themKuneth}.
We will use the moduli space ${\mathcal M}_{\rm QT}(\vec\kappa_1,\vec\kappa_2,\vec\sigma;\vec p_1,\vec p_{2},\vec p_{12};\vec k_1,\vec k_2,\vec k_{12};
p_{\pm \infty};E)$ and its variant in the next subsection also.

\subsection{A K\"unneth type homomorphism between Hochschild homologies}
\label{subsec:KuHoch}

In this section we obtain an element ${\bf y}_0$ appearing in Proposition \ref{prop2512}.
For this purpose we will construct a chain map
\begin{equation}\label{form25800}
\Phi: CH_*(\mathcal L;\Lambda_0) \otimes CH_*(\mathcal L;\Lambda_0) 
\to CH_*(\mathbb L\Fuk;\Lambda_0)
\end{equation}
which induces a homomorphism
\begin{equation}\label{form2580}
HH_*(\mathcal L;\Lambda_0) \otimes HH_*(\mathcal L;\Lambda_0) 
\to HH_*(\mathbb L\Fuk;\Lambda_0).
\end{equation}
Such a homomorphism as (\ref{form2580}) 
could be obtained as the composition
\begin{equation}\label{form258}
\aligned
HH_*(\mathcal L;\Lambda_0) \otimes HH_*(\mathcal L;\Lambda_0) 
&\to HH_*(\mathcal L \otimes \mathcal L;\Lambda_0)\\
 &\to HH_*(\Fuk;\Lambda_0) 
\to HH_*(\mathbb L\Fuk;\Lambda_0),
\endaligned
\end{equation}
where the first homomorphism is an (algebraic) K\"unneth homomorphism,
the second homomorphism is  induced from the 
$A_{\infty}$ functor 
$\mathcal L \otimes \mathcal L \to \Fuk$ constructed in \cite[Theorem 16.9]{fukaya:functor}
and the third homomorphism is induced by the Yoneda embedding.
\par
However in this paper we construct (\ref{form25800}) directly 
not as the composition (\ref{form258}).  We believe that they 
coincide.  We however do not prove it since we do not use it.
The particular construction of (\ref{form25800}) we provide in this subsection 
is used  for the proof of Proposition  \ref{prop2512}.

We first define

\begin{equation}\label{form26100}
\aligned
\varphi : B_{k_1} \sL^{\rm op}[1](c_{\kappa_{k_1}},c_{\kappa_{0}}) &\otimes B_{k_2}\sL[1](c_{\kappa'_{0}},c_{\kappa'_{k_2}}) \\
&\to
\Dat_2(\Dat_{2}(c_{\kappa_{0}},c_{\kappa'_{0}})), \Dat_{2}(c_{\kappa_{k_1}},c_{\kappa'_{k_2}}))
\endaligned
\end{equation}
by
\begin{equation}\label{form2610}
\varphi({\bf x},{\bf y})({\bf w})(z) = (-1)^{\deg'{\bf y}(\deg' {\bf w} + \deg' z)}\frak m_{k_1,k\vert 1 \vert k_2}
({\bf x},{\bf w};z;{\bf y}).
\end{equation}

Let us elaborate on (\ref{form2610}).
In case ${\bf w} = 1_{d_{\sigma}} \in B_0\Fuk[1](d_{\sigma},d_{\sigma})$
the pre-left $\Fuk$ module homomorphism $\varphi({\bf x},{\bf y})$,  gives
$$
\varphi({\bf x},{\bf y})(1_{d_{\sigma}}) \in \Hom(\mathscr{CF}(c_{\kappa_{0}},d_{\sigma};c_{\kappa'_{0}}),
\mathscr{CF}(c_{\kappa_{k_1}},d_{\sigma};c_{\kappa'_{k_2}}))
$$
that is,
$$
z \mapsto (-1)^{\deg'{\bf y}(\deg' z)}\frak m_{k_1,0\vert 1 \vert k_2}({\bf x},1_{d_{\sigma}};z;{\bf y}).
$$
When ${\bf w} \in B_k\Fuk[1](d_{\sigma},d_{\sigma'})$ ($k\ge 1$),
$\varphi({\bf x},{\bf y})$ gives
$$
\varphi({\bf x},{\bf y})({\bf w}) \in \Hom(\mathscr{CF}(c_{\kappa_{0}},d_{\sigma};c_{\kappa'_{0}}),
\mathscr{CF}(c_{\kappa_{k_1}},d_{\sigma'};c_{\kappa'_{k_2}}))
$$
that is
$$
z \mapsto (-1)^{\deg'{\bf y}(\deg' {\bf w} + \deg' z)}\frak m_{k_1,k\vert 1 \vert k_2}
({\bf x},{\bf w};z;{\bf y}).
$$
\par
We next define `Hochschild analogue' of (\ref{form26100}).
We first define a subspace of 
${\mathcal M}_{\rm QT}(\vec\kappa_1,\vec\kappa_2,\vec\sigma;\vec p_1,\vec p_{2},\vec p_{12};\vec k_1,\vec k_2,\vec k_{12};
p_{\pm \infty};E)$.
Let $k_1$ (resp. $k_2$) be the total number of the boundary marked points on $\{-1\}\times \sqrt{-1}\R$ 
(resp.   $\{+1\}\times \sqrt{-1}\R$) of  elements of 
${\mathcal M}_{\rm QT}(\vec\kappa_1,\vec\kappa_2,\vec\sigma;\vec p_1,\vec p_{2},\vec p_{12};\vec k_1,\vec k_2,\vec k_{12};
p_{\pm \infty};E)$.
Here following Notation \ref{not2020}, we use $\vec{\sigma}$ in the above notation.
Let $m_1 \in \{1,\dots, k_1\}$, $m_2 \in \{1,\dots, k_2\}$. 

\begin{defn}\label{defn2620}
We denote by 
${\mathcal M}^{H;m_1,m_2}_{\rm QT}(\vec\kappa_1,\vec\kappa_2,\vec\sigma;\vec p_1,\vec p_{2},\vec p_{12};\vec k_1,\vec k_2,\vec k_{12};
p_{\pm \infty};E)$
the subset of 
${\mathcal M}_{\rm QT}(\vec\kappa_1,\vec\kappa_2,\vec\sigma;\vec p_1,\vec p_{2},\vec p_{12};\vec k_1,\vec k_2,\vec k_{12};
p_{\pm \infty};E)$
represented by an object with an additional property that
$$
{\rm Im}(z_1(m_1)) = {\rm Im}(z_2(m_2)).
$$
Here $z_1(i)$ (resp. $z_2(i)$) is obtained by changing the enumeration 
of marked points on $\{-1\}\times \sqrt{-1}\R$ 
(resp.   $\{+1\}\times \sqrt{-1}\R$) in the way right above Remark \ref{rem2616}.
\end{defn}
${\mathcal M}^{H;m_1,m_2}_{\rm QT}(\vec\kappa_1,\vec\kappa_2,\vec\sigma;\vec p_1,\vec p_{2},\vec p_{12};\vec k_1,\vec k_2,\vec k_{12};
p_{\pm \infty};E)$
has a Kuranishi structure with corners and its (virtual) dimension is 
one less than the (virtual) dimension of 
${\mathcal M}_{\rm QT}(\vec\kappa_1,\vec\kappa_2,\vec\sigma;\vec p_1,\vec p_{2},\vec p_{12};\vec k_1,\vec k_2,\vec k_{12};
p_{\pm \infty};E)$.
\par
We use the moduli space ${\mathcal M}^{H;m_1,m_2}_{\rm QT}(\vec\kappa_1,\vec\kappa_2,\vec\sigma;\vec p_1,\vec p_{2},\vec p_{12};\vec k_1,\vec k_2,\vec k_{12};
p_{\pm \infty};E)$ in place of ${\mathcal M}_{\rm QT}(\vec\kappa_1,\vec\kappa_2,\vec\sigma;\vec p_1,\vec p_{2},\vec p_{12};\vec k_1,\vec k_2,\vec k_{12};
p_{\pm \infty};E)$ 
to obtain:
\begin{equation}
\aligned
&\frak m^{\text{\rm form}}_{(\vec{\kappa}_1^1,\vec{\kappa}_1^2),(\vec{\kappa}_2^1,\vec{\kappa}_2^2),\vec\sigma}\\
:&BCF(\cL^{\rm form};\vec{\kappa}^1_1) \otimes \cL^{\rm form}(L_{\kappa_{m_1-1}},L_{\kappa_{m_1}})
\otimes BCF(\cL^{\rm form};\vec{\kappa}^2_1)\\
&BCF(\sU^{\rm form};\vec{\sigma}) 
\otimes 
\mathscr{CF}(L_{\kappa_{1,\mathcal K_1}},U_{\rho_K};L_{\kappa_{2,0}};\F)
\\ &\otimes BCF(\cL^{\rm form};\vec{\kappa}^1_2)  \otimes \cL^{\rm form}(L_{\kappa_{m_2-1}},L_{\kappa_{m_2}})
\otimes BCF(\cL^{\rm form};\vec{\kappa}^2_2)  
\\
&\to \mathscr{CF}(L_{\kappa_{1,0}},U_{\sigma_0};L_{\kappa_{2,K_2}};\F)
\otimes_{\F} \Lambda_0.
\endaligned
\end{equation}
Here the notations are as follows.
Let
$
\vec{\kappa}_1 = (\kappa_{1,0},\dots,\kappa_{1,\mathcal K_1})
$.
$
\vec{\kappa}_2 = (\kappa_{2,0},\dots,\kappa_{2,\mathcal K_2})
$.
We then put
$$
\aligned
&\vec{\kappa}_1^1 = (\kappa_{1,0},\dots,\kappa_{1,m_1-1}),
\quad
\vec{\kappa}_1^2 = (\kappa_{1,m_1},\dots,\kappa_{1,\mathcal K_1}), \\
&\vec{\kappa}_2^1 = (\kappa_{2,0},\dots,\kappa_{2,m_2-1}),
\quad
\vec{\kappa}_2^2 = (\kappa_{2,m_2},\dots,\kappa_{2,\mathcal K_2}).
\endaligned
$$
Now as usual we proceed as follows.
\begin{enumerate}
\item[(I)]  We include interior marked points and deform 
$\frak m^{\text{\rm form}}_{(\vec{\kappa}_1^1,\vec{\kappa}_1^2),(\vec{\kappa}_2^1,\vec{\kappa}_2^2),\vec\sigma}$ 
to $\frak m^{\text{\rm form},\frak b}_{(\vec{\kappa}_1^1,\vec{\kappa}_1^2),(\vec{\kappa}_2^1,\vec{\kappa}_2^2),\vec\sigma}$.
In other words, we include the bulk deformation.
\item[(II)]
We include forgetability data and extend $\frak m^{\text{\rm form},\frak b}_{(\vec{\kappa}_1^1,\vec{\kappa}_1^2),(\vec{\kappa}_2^1,\vec{\kappa}_2^2),\vec\sigma}$ to its plus version, which admits ${\bf e}^+$ and ${\bf f}$ as inputs.
\item[(III)]
We use the $A_{\infty}$ 
functor $\frak f$ in Subsections \ref{subsec:red} and \ref{hunitcansubsec}. Moreove 
we replace $\mathscr{CF}^{\rm form}$ by its canonical model $\mathscr{CF}$ in the same way as in \cite[Theorem 5.4.18]{fooo09}.
\item[(IV)]
We eliminate the curvature using the given bounding cochains of objects in ${\bf L}$ or ${\bf F}$.
\end{enumerate}
We do not repeat the detail of those steps since in similar situations we have done it many times already in this paper.
\par
\begin{figure}[h]
\centering
\includegraphics[scale=0.35]{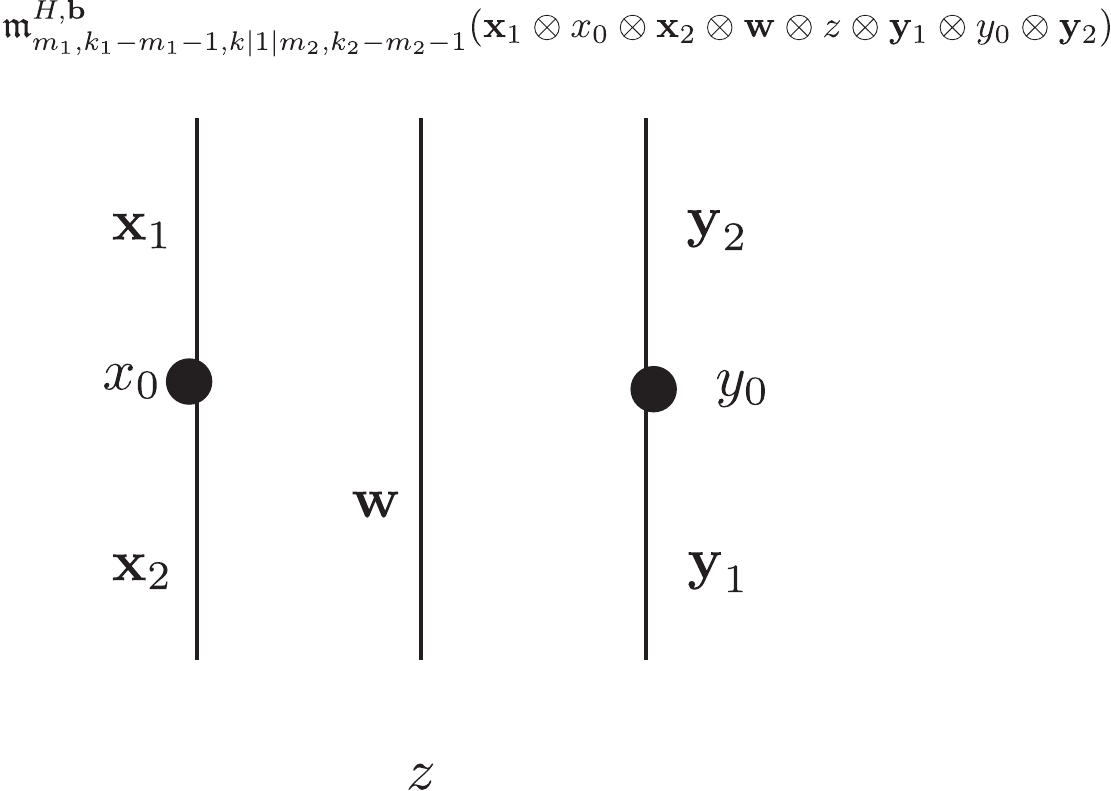}
\caption{$\frak m^{H,{\bf b}}_{m_1,k_1-m_1-1,k\vert 1 \vert m_2,k_2-m_2-1}$}
\label{FigurMHB}
\end{figure}

We thus obtain the following operation:
$$
\aligned
&\frak m^{H,{\bf b}}_{m_1,k_1-m_1-1,k\vert 1 \vert m_2,k_2-m_2-1}\\
:&
B_{m_1}\cL^{\rm op}[1](c_{\kappa^{(1)}_1},c_{\kappa^{(3)}_1}) 
\otimes \cL^{\rm op}(c_{\kappa^{(3)}_1},c_{\kappa^{(4)}_1}) \otimes
B_{k_1-m_1-1}\cL^{\rm op}[1](c_{\kappa^{(4)}_1},c_{\kappa^{(2)}_1})  \\
&\otimes  B_k\sU[1](d_{\sigma},d_{\sigma'}) \otimes \mathscr{CF}(c_{\kappa^{(2)}_{1}},d_{\sigma'};c_{\kappa^{(2)}_{k_2}}) \\
&\otimes B_{m_2}\cL[1](c_{\kappa^{(2)}_1},c_{\kappa^{(3)}_2}) 
\otimes \cL(c_{\kappa^{(3)}_2},c_{\kappa^{(4)}_2}) \otimes
B_{k_2-m_2-1}\cL^{\rm op}[1](c_{\kappa^{(4)}_2},c_{\kappa^{(1)}_2}) \\
&\to \mathscr{CF}(c_{\kappa^{(1)}_{1}},d_{\sigma};c_{\kappa^{(1)}_{2}})
\endaligned
$$
which sends 
$
{\bf x}_1 \otimes x_0 \otimes {\bf x}_2 \otimes {\bf w} \otimes z 
\otimes {\bf y}_1 \otimes y_0 \otimes {\bf y}_2
$
to
$$
\frak m^{H,{\bf b}}_{m_1,k_1-m_1-1,k\vert 1 \vert m_2,k_2-m_2-1}({\bf x}_1 \otimes x_0 \otimes {\bf x}_2 \otimes {\bf w} \otimes z 
\otimes {\bf y}_1 \otimes y_0 \otimes {\bf y}_2).
$$
See Figure \ref{FigurMHB}.

We now use it to define
$$
\aligned
\varphi^H 
& : B_{m_1}\cL^{\rm op}[1](c_{\kappa^{(1)}_1},c_{\kappa^{(3)}_1}) 
\otimes \cL^{\rm op}(c_{\kappa^{(3)}_1},c_{\kappa^{(4)}_1}) \otimes
B_{k_1-m_1-1}\cL^{\rm op}[1](c_{\kappa^{(4)}_1},c_{\kappa^{(2)}_1}) \\
& \otimes B_{m_2}\cL[1](c_{\kappa^{(2)}_2},c_{\kappa^{(3)}_2}) 
\otimes \cL(c_{\kappa^{(3)}_2},c_{\kappa^{(4)}_2})\otimes
B_{k_2-m_2-1}\cL^{\rm op}[1](c_{\kappa^{(4)}_2},c_{\kappa^{(1)}_2}) \\
&\to
\Dat_2(\Dat_{2}((c_{\kappa_1^{(2)}},c_{\kappa_2^{(2)}}),d_{\sigma'}), 
\Dat_{2}((c_{\kappa_1^{(1)}},c_{\kappa_2^{(1)}}),d_{\sigma}))
\endaligned
$$
by
\begin{equation}\label{form2990}
\aligned
&\varphi^H(
{\bf x}_1 \otimes x_0 \otimes {\bf x}_2  
\otimes {\bf y}_1 \otimes y_0 \otimes {\bf y}_2)({\bf w})(z) \\
& = (-1)^{\maltese}\frak m^H_{m_1,k_1-m_1-1,k\vert 1 \vert m_2,k_2-m_2-1} \\
&\qquad\qquad\qquad\qquad({\bf x}_1 \otimes x_0 \otimes {\bf x}_2 \otimes {\bf w} \otimes z 
\otimes {\bf y}_1 \otimes y_0 \otimes {\bf y}_2),
\endaligned
\end{equation}
with 
$\maltese = (\deg'{\bf w} + \deg z)(\deg'{\bf y}_1+\deg y_0 + \deg'{\bf y}_2)$.

We next use $\varphi$, $\varphi^H$ to define $\Phi$ in (\ref{form25800}) as follows.
\par
We first introduce a few notations. We consider 
${\bf x}= x_0 \otimes {\bf x}_+,  {\bf y} = y_0 \otimes {\bf y}_+
\in CH_*(\mathcal L;\Lambda_0)$.
We then put
\begin{equation}\label{DeltaH2}
\Delta_H^{k-1}(\text{\bf x})
=
\sum_{c_1}(-1)^{\maltese} \text{\bf x}_{c_1}^{(H;k;1)} \otimes \cdots \otimes 
\text{\bf x}_{c_1}^{(H;k;k)}
\end{equation}
as in (\ref{DeltaH}). (The sign $(-1)^{\maltese} $ is as in (\ref{DeltaH}).) 
$\text{\bf x}_{c_1}^{(H;k;1)}$ contains $x_0$ and we write
$$
\text{\bf x}_{c_1}^{(H;k;1)} =  \text{\bf x}_{c_1}^{(H;k;1,1)} \otimes x_0 \otimes \text{\bf x}_{c_1}^{(H;k;1,2)}.
$$
We define $\text{\bf y}_{c_1}^{(H;k;1,1)}$, $\text{\bf y}_{c_1}^{(H;k;1,2)}$ and $\text{\bf y}_c^{(H;k;i)}$
in the same way.
\begin{defn}
\begin{equation}
\aligned
\Phi({\bf x},{\bf y})
= \sum_{k=0}^{\infty} \sum_{c_1,c_2}\!{}'
&\varphi^H(\text{\bf x}_{c_1}^{(H;k;1,1)} \otimes x_0 \otimes \text{\bf x}_{c_1}^{(H;k;1,2)},
\text{\bf y}_{c_2}^{(H;k;1,1)}  \otimes y_0 \otimes \text{\bf y}_{c_2}^{(H;k;1,2)}) \\
&\otimes \varphi(\text{\bf x}_{c_1}^{(H;k;2)},\text{\bf y}_{c_2}^{(H;k;2)})  \otimes \dots \otimes \varphi(\text{\bf x}_{c_1}^{(H;k;k)},\text{\bf y}_{c_2}^{(H;k;k)})\\
&\!\!\!\!\!\!\!\!\!\!\!\!\!\!\!\!\!\!\!\! \in CH_*(\mathbb L\Fuk;\Lambda_0).
\endaligned
\end{equation}
Here $\sum\!{}'$ is a sum over all pairs $(c_1,c_2)$ such that there is no
$j =2,\dots,k$ with $\text{\bf x}_{c_1}^{(H;k;1,j)} = \text{\bf y}_{c_2}^{(H;k;1,j)} =1 \in B_0\sL$.
\end{defn}

\begin{prop}\label{prop1621}
$\Phi$ is a chain map.
\end{prop}
\begin{proof}
We first study $\varphi^H$.
As usual we  study the codimension $1$ boundary of the moduli space 
${\mathcal M}^{H;m_1,m_2}_{\rm QT}(\vec\kappa_1,\vec\kappa_2,\vec\sigma;\vec p_1,\vec p_{2},\vec p_{12};\vec k_1,\vec k_2,\vec k_{12};
p_{\pm \infty};E)$.
It is classified by the 9 types shown in the next Figures 
\ref{Figureabc}, \ref{Figuredef}, \ref{Figureghi}.

\begin{figure}[h]
\centering
\includegraphics[scale=0.35]{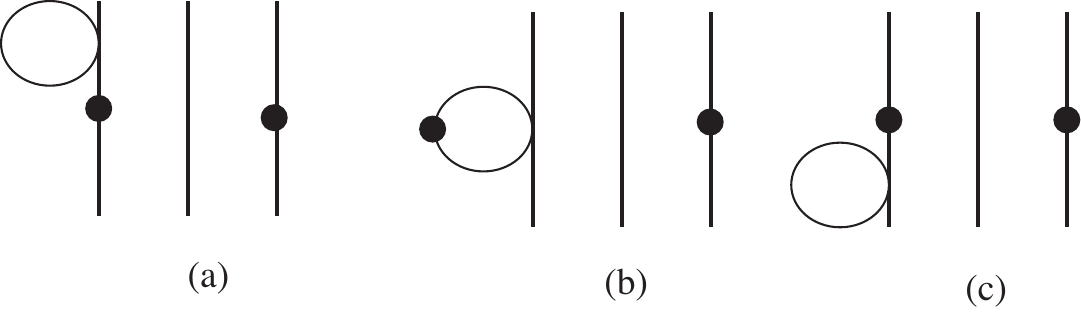}
\caption{Bubble at $\{-1\} \times \sqrt{-1}\R$}
\label{Figureabc}
\end{figure}
\begin{figure}[h]
\centering
\includegraphics[scale=0.35]{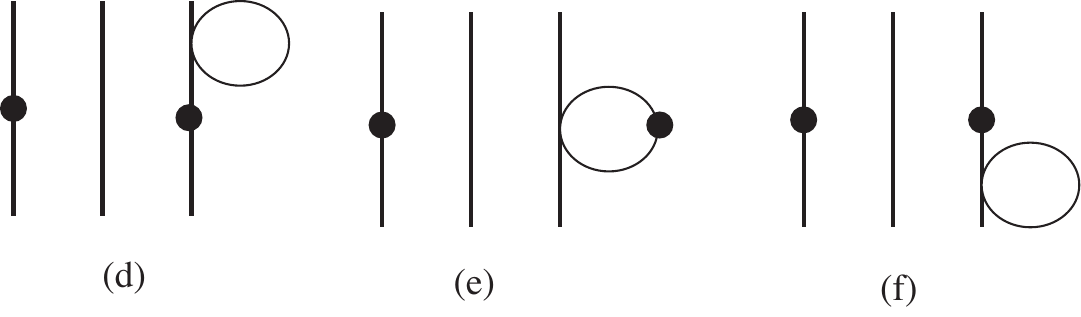}
\caption{Bubble at $\{+1\} \times \sqrt{-1}\R$}
\label{Figuredef}
\end{figure}
\begin{figure}[h]
\centering
\includegraphics[scale=0.35]{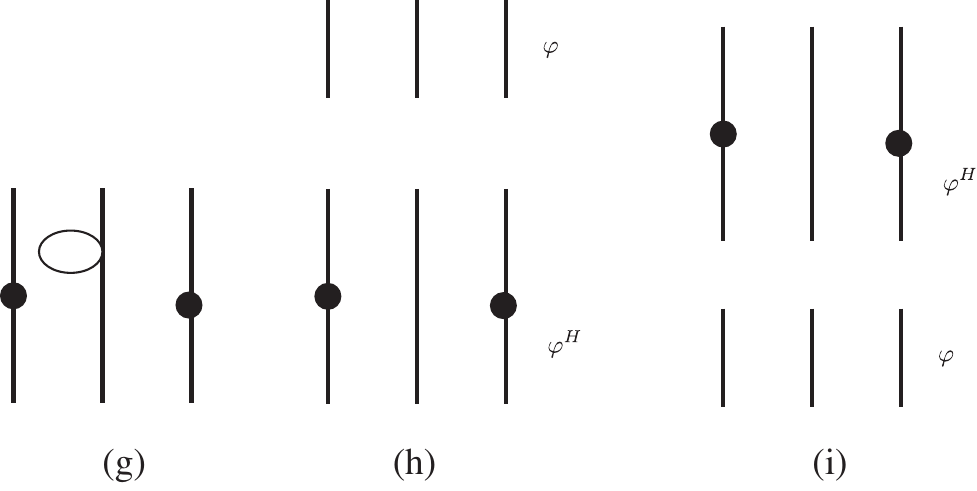}
\caption{Other types of boundaries.}
\label{Figureghi}
\end{figure}
For $\varphi$, the boundary of the moduli space, we need to study is (bd1)-(bd4)
appearing right after Definition \ref{def516}.
\par
The Hochschild boundary $\delta_H\circ \Phi$.
It is sum of 
$$
\aligned
&\left(\varphi^H(\cdots) 
\circ \varphi(\text{\bf x}_{c_1}^{(H;k;2)},\text{\bf y}_{c_2}^{(H;k;2)})\right) \\
&\qquad\otimes \varphi(\text{\bf x}_{c_1}^{(H;k;3)},\text{\bf y}_{c_2}^{(H;k;3)})  \otimes \dots \otimes \varphi(\text{\bf x}_{c_1}^{(H;k;k)},\text{\bf y}_{c_2}^{(H;k;k)}),
\\
&\left(\varphi(\text{\bf x}_{c_1}^{(H;k;k)},\text{\bf y}_{c_2}^{(H;k;k)}) \circ \varphi^H(\cdots)  \right) \\
&\qquad\otimes \varphi(\text{\bf x}_{c_1}^{(H;k;2)},\text{\bf y}_{c_2}^{(H;k;2)})  \otimes \dots \otimes \varphi(\text{\bf x}_{c_1}^{(H;k;k-1)},\text{\bf y}_{c_2}^{(H;k;k)}),
\\
&\varphi^H(\cdots)  \otimes
\dots  \otimes 
\left(\varphi(\text{\bf x}_{c_1}^{(H;k;i)},\text{\bf y}_{c_2}^{(H;k;i)}) \circ  \varphi(\text{\bf x}_{c_1}^{(H;k;i+1)},\text{\bf y}_{c_2}^{(H;k;i+1)})\right) \\
&\qquad\otimes   \dots \otimes \varphi(\text{\bf x}_{c_1}^{(H;k;k)},\text{\bf y}_{c_2}^{(H;k;k)})
\endaligned
$$
(here we omit the variables in $\varphi^H$ and write $\varphi^H(\cdots)$ simplicity)
and 
terms obtained by replacing $\varphi^H$ or various $\varphi$ in the definition of $\Phi$ by their commutators with 
$\frak m_{0,k'\vert 1\vert 0}$.
The commutator with $\varphi^H$ is obtained from Figure \ref{Figureghi} (h)(i) but in the case there is no
marked points on $\{\pm 1\} \times \sqrt{-1}\R$ in the part written $\varphi$.
The commutator with $\varphi$ is obtained from (bd4) in a similar way.
\par
Now we can see that those terms cancel each other to give $0$.
The proof of Proposition \ref{prop1621} is complete.
\end{proof}
Our next task is to calculate
$
Z_{{\cQ\Cat^{\rm taut}\vert_{\Dat_2}}}(\Phi({\bf x},{\bf y}),{\bf w})
$
for ${\bf x},  {\bf y} 
\in CH_*(\mathcal L;\Lambda_0)$
and ${\bf w} 
\in CH_*(\Fuk;\Lambda_0)$
in terms of the structure operations $\frak m^{\bf b}_{k_1,k | 1  | k_2}$ of trimodule. 
%and  
%the structure operations $\frak m^{H,{\bf b}}_{m_1,k_1-m_1-1,k\vert 1 \vert m_2,k_2-m_2-1}$ of quatromodule

\begin{figure}[h]
\centering
\includegraphics[scale=0.5]{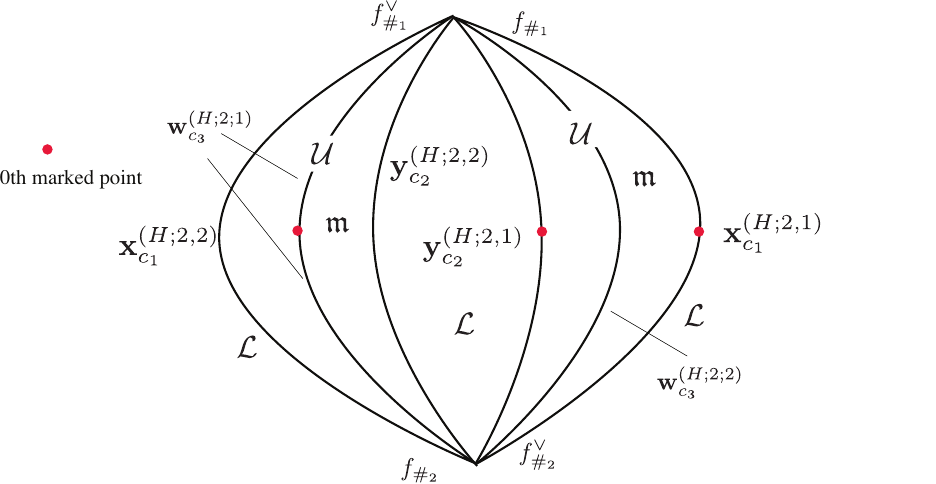}
\caption{$Z_{\rm QT}$}
\label{Figure67}
\end{figure}

We will define\index[syindex]{ZQT@$Z_{\rm QT}$} 
$$
Z_{\rm QT}: CH_*(\mathcal L;\Lambda_0)  \otimes CH_*(\Fuk;\Lambda_0) \otimes CH_*(\mathcal L;\Lambda_0)
\to \Lambda_0
$$ 
for this purpose.
We put
$$
\aligned
&\Delta_H(\text{\bf x})
=
\sum_{c_1}(-1)^{\maltese}  \text{\bf x}_{c_1}^{(H;2;1)} \otimes 
\text{\bf x}_{c_1}^{(H;2;2)} \\
&\text{\bf x}_{c_1}^{(H;2;1)} = \text{\bf x}_{c_1}^{(H;2;1,1)} \otimes x_0 \otimes \text{\bf x}_{c_1}^{(H;2;1,2)}
\endaligned
$$
and define $\text{\bf y}_{c_2}^{(H;2;1)}$, $\text{\bf y}_{c_2}^{(H;2;2)}$,
$\text{\bf w}_{c_3}^{(H;2;1)}$, $\text{\bf w}_{c_3}^{(H;2;2)}$
in the same way.
Now we put
\begin{equation}\label{formula23999}
\aligned
&Z_{\rm QT}(\text{\bf x},\text{\bf w},\text{\bf y}) \\
&= \sum_{c_1,c_2,c_3}\sum_{a,a'}
(-1)^{\maltese}
\langle 
\frak m^{{\bf b}}_{*|1|*}(\text{\bf x}_{c_1}^{(H;2;1)},{\bf w}^{(H;2;2)}_{c_3};f_{\#_1},
\text{\bf y}_{c_2}^{(H;2;1)}),f_{\#_2}^{\vee}
\rangle
\\
&\qquad\qquad\qquad\qquad\langle 
\frak m^{{\bf b}}_{*|1|*}(\text{\bf x}_{c_1}^{(H;2;2)},{\bf w}^{(H;2;1)}_{c_3} ;f_{\#_2},
\text{\bf y}_{c_2}^{(H;2;2)}),f_{\#_1}^{\vee}
\rangle.
\endaligned
\end{equation}
See Figure \ref{Figure67}.
We define the sign $\maltese$ later in Subsubsection \ref{subsubsec:ref}.
Here $f_{\#_1}$ (resp. $f_{\#_2}$) is a basis of 
$\mathscr{CF}(c_{\kappa_{1}},d_{\sigma_1};c_{\kappa'_{2}})$
(resp. $\mathscr{CF}(c_{\kappa'_{1}},d_{\sigma'_1};c_{\kappa_{2}})$),
where 
$\text{\bf x}_{c_1}^{(H;2;1)} \in B\sL(c_{\kappa_{1}},c_{\kappa'_{1}})$,
$\text{\bf y}_{c_2}^{(H;2;1)} \in B\sL(c_{\kappa'_{2}},c_{\kappa_{2}})$,
${\bf w}^{(H;2;1)}_{c_3} \in B\sU(d_{\rho_{1}},d_{\rho'_{1}})$.
$f^{\vee}_{\#_1}$ and $f^{\vee}_{\#_2}$  are their dual basis.

\begin{lem}\label{lem2622}
We have
$$
Z_{{\cQ\Cat^{\rm taut}\vert_{\Dat_2}}}(\Phi({\bf x},{\bf y}),{\bf w})
= Z_{\rm QT}(\text{\bf x},\text{\bf w};\text{\bf y}).
$$
\end{lem}
In fact this is immediate from the definition.
(We use the fact that $\Dat_2$ is a DG category.
Namely all the higher compositions $\frak m_k$, $k > 2$ vanish.) 

\begin{prop}\label{prop2623}
Suppose $\frak p^{\bf b}({\bf x}) = [{\rm 1}_X]$
and 
${\bf w} \in HH_*(\mathcal U;\Lambda)$. Then 
$$
Z_{\rm QT}({\bf x},{\bf w};{\bf x}) 
\begin{cases}
= 1   &\text{If ${\bf w} =  [{\rm vol}_{U}]\in H^*(U;\Lambda)$},
\\
= 0    &\text{If ${\bf w} \in HH_*(\mathcal U;\Lambda)$ and does not contain the term $[{\rm vol}_{U}]$}.
\end{cases}
$$
\end{prop}
Proposition \ref{prop2623} and Lemma \ref{lem2622} immediately imply Proposition \ref{prop2512}.\footnote{
We remark that the symbols ${\bf x}$, ${\bf w}$ in  Proposition \ref{prop2623} 
correspond to ${\bf y}_0$ and ${\bf x}$, respectively, in Proposition  \ref{prop2512}.}
We will prove Proposition \ref{prop2623}
in the next subsection.

\subsection{The Cardy relation for Lagrangian correspondences}
\label{subsec;cardy}

\subsubsection{Quilted annuli}
\label{subsubsec:qanu}

The proof of Proposition \ref{prop2623} is based on certain moduli spaces of pseudo-holomorphic curves,
which we call  {\it quilted annuli}.\index{quilted annuli}
We start with defining them.  

For $j=1,2$, let $(\vec{\kappa}^{(j)},\vec k^{(j)}) \in \text{\rm seq}'_{K(j)}(\cL)$.  Here  $\text{\rm seq}'_{K(j)}(\cL)$
is defined at the beginning of Subsection \ref{sec:moduli-space-annuli}.
Let $(\vec{{\upsilon}},\vec r) \in \text{\rm seq}'_{R}(\sU)$ where $\text{\rm seq}'_{R}(\sU)$
is  defined in the same way.

We fix an intersection point $L_{\kappa^{(j)}_{\blue{i-1}}} \cap L_{\kappa^{(j)}_{\blue{i}}}$ (or self-intersection point of 
$L_{\kappa^{(j)}_{\blue{i-1}}} = L_{\kappa^{(j)}_{\blue{i}}}$),
which we write $p^{(j)}_i$.
Similarly %$y_i$ is identified with
we fix $q_i \in U_{\rho_{\blue{i-1}}} \cap U_{\rho_{\blue{i}}}$ (or self-intersection point). %whenever these Lagrangians are different.
We denote by $\vec p^{\,(j)}$, (resp. $\vec q$) the totality of $p^{(j)}$'s, (resp $q$'s).
\par
We put $k_{(j)}-1 = \vert\vec k^{(j)}\vert  = \sum_{i=0}^{K(j)} k^{(j)}_i$, 
$r -1 = \vert\vec r\vert  = \sum_{i=0}^R r_i$.

Let  $T \in (0,\infty)$ and
$\Sigma(T)$ be  $[-T,T] \times S^1$ together (possibly) with trees of sphere bubbles 
rooted on $((-T,T) \setminus \{0\}) \times S^1$.
$\partial\Sigma(T)$ is a union of two circles $\partial_{\pm}\Sigma$ where 
$\partial_-\Sigma(T) = \{-T\} \times S^1$, $\partial_+\Sigma(T) = \{+T\} \times S^1$.  The circle $\{0\} \times S^1$ is 
called a seam.\index{seam}
\par
We denote by $\Sigma_-(T)$ (resp.  $\Sigma_+(T)$) the union of $[-T,0] \times S^1$ 
(resp. $[0,T] \times S^1$) together with trees of spheres 
rooted on it.

\begin{defn}\label{2624dsef}
 $\overset{\circ}{\mathcal M}_{k_{(1)}+1,k_{(2)}+1,r+1;T}$
 is the set of $(\Sigma;\vec z_1,\vec z_2,\vec z)$ such that:
 \begin{enumerate}
 \item $\Sigma = \Sigma(T)$ is as above.
 \item  $\vec z_j = (z_{j}(0),\dots,z_{j}(k_{(j)}))$.  $z_{1}(i) \in \partial_-\Sigma$, $z_{2}(i) \in \partial_+\Sigma$.
 The points $z_{j}(0),\dots,z_{j}(k_{(j)})$ are mutually distinct and  respect  the cyclic order determined by the boundary orientation of 
 $\partial_{\pm}\Sigma$.   In other words for $j=1$ the orientation is the negative direction and for $j=2$ the orientation 
 is positive direction.
 \item $\vec z = (z(0),\dots,z(r)$.  $z(i) \in \{0\} \times S^1$.  The points $z(0),\dots,z(r)$ are mutually distinct and respect the 
 downward cyclic order of 
 $S^1$.  (The orientation is of negative direction.)
 \item  $z_{1}(0) = (-T,0)$, $z_{2}(0) = (T,0)$, $z(0) = (0,1/2)$. 
 Here we identify $S^1 \cong \R/\Z$.
 \end{enumerate}

Let
\begin{equation}\label{formnew21644}
\aligned
\vec z^{\,(1)} &= (z^{(1)}_0,w^{(1)}_{{0},1},\dots,w^{(1)}_{{0} ,k^{(1)}_1},
z^{(1)}_1,,\dots,z^{(1)}_K,
w^{(1)}_{{K}(1) ,1},\dots,w^{(1)}_{{K}(1),k_{K(1)}}), \\
\vec z^{\,(2)} &= (z^{(2)}_0,w^{(1)}_{{0},2},\dots,w^{(1)}_{{0} ,k^{(2)}_1},
z^{(2)}_1,,\dots,z^{(2)}_K,
w^{(2)}_{{K}(2) ,1},\dots,w^{(2)}_{{K}(2),k_{K(2)}})\\
\vec z^{\,(12)} &= (z^{(12)}_0,w^{(12)}_{{0},1},\dots,w^{(12)}_{{0} ,k^{(12)}_1},
z^{(12)}_1,,\dots,z^{(12)}_K,
w^{(12)}_{{K} },\dots,w^{(12)}_{{K},k_{K}})
\endaligned
\end{equation}
be points on $\partial_-\Sigma(T)$, $\partial_+\Sigma(T)$, $\{0\} \times S^1$, respectively.
We require that they respect the cyclic order of $\partial_-\Sigma(T)$, $\partial_+\Sigma(T)$,
and of $\{0\} \times S^1$, respectively.
\par 
Let  $m_j \in \{0,\dots,k^{(j)}_0\}$ ($j=1,2$),  $m \in \{0,\dots,r_0\}$.
We renumber $\vec z^{\,(j)}$ to $\vec z_j =  (z_{j}(0),\dots,z_{j}(k_{(j)}))$
and $\vec z^{\,(12)}$ to $\vec z =(z(0),\dots,z(r))$ such that:
\begin{enumerate}
\item[(I)] $\vec z_j$ and $\vec z$ are as in (2)(3)(4) above.
\item [(II)]
$z_j(0) = z^{\,(j)}_0$ if $m_j =0$. If $m_j > 0$ then
$z_j(0) = w^{\,(j)}_{0,m_j}$.
\item[(III)]
$z(0) = z_0$ if $m =0$. If $m > 0$ then
$z(0) = w_{0,m}$.
\end{enumerate}
\end{defn}

\begin{figure}[h]
\centering
\includegraphics[scale=0.4]{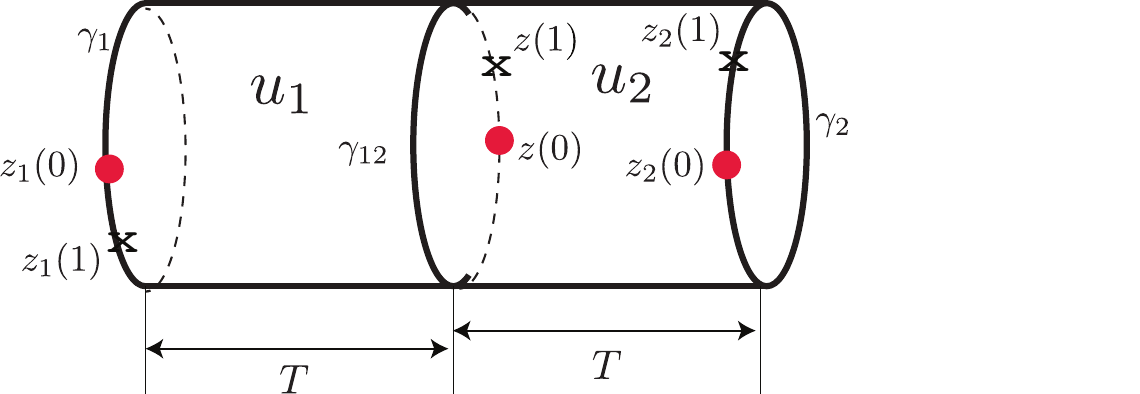}
\caption{${\mathcal M}_{(\vec{\kappa}^{(1)},\vec k^{(1)},m_1), (\vec{\kappa}^{(2)},\vec k^{(2)},m_2),(\vec{{\upsilon}},\vec r,m)}(\sL;\sU;\sL;T;E)$}
\label{Figure11}
\end{figure}

In the next definition we assume $U_{\rho}$ is transversal to $L_{\kappa} \times L_{\kappa'}$.
For the proof of Proposition \ref{prop2623}, we use the case $U =$ diagonal. In such a case we perturb the diagonal a bit 
so that it is transversal to $L_{\kappa} \times L_{\kappa'}$.\index[syindex]{MsLsUsL@${\mathcal M}_{(\vec{\kappa}^{(1)},\vec k^{(1)},m_1), (\vec{\kappa}^{(2)},\vec k^{(2)},m_2),(\vec{{\upsilon}},\vec r,m)}(\sL;\sU;\sL;R;E)$}

\begin{defn}\label{def1625}
We define $\overset{\circ}{\mathcal M}_{(\vec{\kappa}^{(1)},\vec k^{(1)},m_1), (\vec{\kappa}^{(2)},\vec k^{(2)},m_2),(\vec{{\upsilon}},\vec r,m)}(\sL;\sU;\sL;T;E)$ to be the set of equivalence classes of 
the objects 
$((\Sigma;\vec z_1,\vec z_2,\vec z_{12});u_1,u_2;\gamma_1,\gamma_2,\gamma_{12})$ satisfying the conditions we describe below.
See Figure \ref{Figure11}.
\begin{enumerate}
\item $(\Sigma;\vec z_1,\vec z_2,\vec z)$ represents an element of $\overset{\circ}{\mathcal M}_{k_{(1)}+1,k_{(2)}+1,r+1;T}$.
We change the enumeration of $\vec z_1,\vec z_2,\vec z$ and obtain $\vec z^{\,(1)}$, $\vec z^{\,(2)}$, $\vec z^{\,(12)}$, $m_1$, $m_2$, $m$ 
such that $\vec z_1,\vec z_2,\vec z$ are obtained from them as in the second half of Definition \ref{2624dsef}.
\item
$u_1 : \Sigma_-(T) \setminus (\vec z^{\,(1)}\cup \vec z^{\,(12)}) \to X$ is a $-J_{X}$-holomorphic map 
and $u_2 : \Sigma_+(T) \setminus (\vec z^{\,(2)}\cup \vec z^{\,(12)}) \to X$ is a $J_{X}$-holomorphic map.
\item 
$\gamma_1 : (\{-T\} \times S^1) \setminus \vec z^{\,(1)} \to \tilde L$, 
$\gamma_2 : (\{T\} \times S^1) \setminus \vec z^{\,(2)} \to \tilde L$,  
$\gamma_{12} : (\{0\} \times S^1) \setminus \vec z^{\,(12)} \to \tilde U$ are smooth maps.
Here $\tilde L$ and  $\tilde U$ are disjoint union of $\tilde L_{\kappa}$, $\tilde U_{\rho}$
respectively.
\item 
$$
\gamma_j\left(\overline{z^{(j)}_{i-1},z^{(j)}_{i}}\right) \subset \tilde L_{\kappa^{(j)}_{i-1}},
\qquad 
\gamma_{12}\left(\overline{z^{(12)}_{i-1},z^{(12)}_{i}}\right) \subset \tilde U_{\rho_{i-1}}.
$$
Here $\kappa^{(j)}_{K(j)+1} = \kappa^{(j)}_{0}$, $\rho_{K+1} = \rho_{0}$ by convention.
\item  For $j=1,2$ we have an equality $u_j = i_{L} \circ \gamma_j$ when both sides are defined.
Here $i_L = i_{L_{\kappa}}$ on $\tilde L_{\kappa}$.
We also have $(u_1,u_2) = i_{U} \circ \gamma_{12}$ when both sides are defined.
Here $i_U = i_{U_{\rho}}$ on $\tilde U_{\rho}$.
\item We require $u_j(z^{(j)}_i) = p^{(j)}_i$.  $(u_1(z^{(12)}_i), u_2(z^{(12)}_i)) =  q_i$.
(Here we take various choices of $\vec p^{\,(j)}$ and $\vec q$ and take the union.
When we fix a homology class $B$ we fix the choice.)
\item
We require a similar switching condition as Definition \ref{def211} (4)
at $z^{(j)}_i$ and $z^{(12)}_i$.
\item 
We require the stability condition defined below.
\item
$$
-\int u_1^*\omega_{X} + \int u_2^*\omega_{X} = E.
$$
\end{enumerate}

We say an object $((\Sigma;\vec z_1,\vec z_2,\vec z_{12});u_1,u_2;\gamma_1,\gamma_2,\gamma_{12})$ is equivalent to 
another object $((\Sigma';\vec z^{\,\prime}_1,\vec z^{\,\prime}_2,\vec z^{\,\prime}_{12});u'_1,u'_2;\gamma'_1,\gamma'_2,\gamma'_{12})$ 
if there exists a biholomorphic map $v : \Sigma \to \Sigma'$ such that
$v(\vec z_1) = \vec z^{\,\prime}_1$, $v(\vec z_2) = \vec z^{\,\prime}_2$, $v(\vec z_{12}) = \vec z^{\,\prime}_{12}$, $u'_1 \circ v =  u_1$, 
$u'_2 \circ v =  u_2$, $u'_{12} \circ v =  u_{12}$, $\gamma'_1 \circ v =  \gamma_1$, 
$\gamma'_2 \circ v =  \gamma_2$, $\gamma'_{12} \circ v =  \gamma_{12}$.
\par
We say  $((\Sigma;\vec z_1,\vec z_2,\vec z_{12});u_1,u_2;\gamma_1,\gamma_2,\gamma_{12})$
is stable
if the group of automorphisms in the above sense is a finite group.
\par
We decompose 
$\overset{\circ}{\mathcal M}_{(\vec{\kappa}^{(1)},\vec k^{(1)},m_1), (\vec{\kappa}^{(2)},\vec k^{(2)},m_2),(\vec{{\upsilon}},\vec r,m_)}(\sL;\sU;\sL;T;E)$ 
according to the homology classes $B$ and use the symbol 
$\overset{\circ}{\mathcal M}_{(\vec{\kappa}^{(1)},\vec k^{(1)},m_1), (\vec{\kappa}^{(2)},\vec k^{(2)},m_2),(\vec{{\upsilon}},\vec r,m_)}(\sL;\sU;\sL;T;B)$ to denote the components. \index[syindex]{MsLsUsLB@${\mathcal M}_{(\vec{\kappa}^{(1)},\vec k^{(1)},m_1), (\vec{\kappa}^{(2)},\vec k^{(2)},m_2),(\vec{{\upsilon}},\vec r,m_)}(\sL;\sU;\sL;T;B)$}
\end{defn} 

We can compactify $\overset{\circ}{\mathcal M}_{(\vec{\kappa}^{(1)},\vec k^{(1)},m_1), (\vec{\kappa}^{(2)},\vec k^{(2)},m_2),(\vec{{\upsilon}},\vec r,m_)}(\sL;\sU;\sL;T;E)$ 
to obtain its compactification ${\mathcal M}_{(\vec{\kappa}^{(1)},\vec k^{(1)},m_1), (\vec{\kappa}^{(2)},\vec k^{(2)},m_2),(\vec{{\upsilon}},\vec r,m_)}(\sL;\sU;\sL;T;E)$ 
in the same way as \cite[Subsection 5.2 and Section 12]{fukaya:functor}.
They carry Kuranishi structures with a boundary and corners.
Its boundary is classified into the following three cases.
 \begin{enumerate}
\item[(bds1)]  The disk bubble at $\{-T\} \times S^1$.
\item[(bds2)]  The disk bubble at $\{+T\} \times S^1$.
\item[(bds3)]  The disk bubble at $\{0\} \times S^1$.
\end{enumerate}
The evaluation maps at boundaries $\partial_-\Sigma$, $\partial_+\Sigma$ and the seam 
define maps ${\rm ev}^{\partial}_-$, ${\rm ev}^{\partial}_+$, ${\rm ev}^{\rm sm}$, 
respectively. They are maps
\begin{equation}
\aligned
&{\rm ev}^{\partial}_-: 
{\mathcal M}_{(\vec{\kappa}^{(1)},\vec k^{(1)},m_1), (\vec{\kappa}^{(2)},\vec k^{(2)},m_2),(\vec{{\upsilon}},\vec r,m_)}(\sL;\sU;\sL;T;E) \to 
\bigcup_{\vec p^{(1)}}L^{\vec{\kappa}^{(1)},\vec p^{(1)}}_{\text{\rm source}}, \\
&{\rm ev}^{\partial}_+: 
{\mathcal M}_{(\vec{\kappa}^{(1)},\vec k^{(1)},m_1), (\vec{\kappa}^{(2)},\vec k^{(2)},m_2),(\vec{{\upsilon}},\vec r,m_)}(\sL;\sU;\sL;T;E) \to \bigcup_{\vec p^{(2)}}L^{\vec{\kappa}^{(2)},\vec p^{(2)}}_{\text{\rm source}}, \\
&{\rm ev}^{\rm sm} : {\mathcal M}_{(\vec{\kappa}^{(1)},\vec k^{(1)},m_1), (\vec{\kappa}^{(2)},\vec k^{(2)},m_2),(\vec{{\upsilon}},\vec r,m)}(\sL;\sU;\sL;T;E) \to \bigcup_{\vec q}U^{\vec{\kappa},\vec q}_{\text{\rm source}}.
\endaligned
\end{equation}
Here targets are defined in the same way as Definition \ref{defn927777}.
Let ${\bf h}^{(j)} \in \Omega(L^{\vec{\kappa}^{(j)},\vec p^{(j)}}_{\text{\rm source}})$ 
and 
${\bf h} \in \Omega(U^{\vec{\kappa},\vec q}_{\text{\rm source}})$.
We define:\footnote{Recall as usual ${\rm f.c.}$ stands for from, curved, respectively.}
$$
Z_{T,B}^{\rm f.c.}({\bf h}^{(1)},{\bf h};{\bf h}^{(2)})
= 
\int_{\star} ({\rm ev}^{\partial}_-)^*({\bf h}^{(1)}) \wedge ({\rm ev}^{\rm sm})^*{\bf h} \wedge({\rm ev}^{\partial}_+)^*({\bf h}^{(2)})
$$
where $\star ={\mathcal M}_{(\vec{\kappa}^{(1)},\vec k^{(1)},m_1), (\vec{\kappa}^{(2)},\vec k^{(2)},m_2),(\vec{{\upsilon}},\vec r,m_)}(\sL;\sU;\sL;T;B)$.
We take and fix an appropriate CF-perturbation to define the integration.
We now consider the sum
$$
Z_{T}^{\rm f.c.} = \sum_B e^{\rho_{\frak b,\theta}(B)} T^{B(\omega)} Z_{T,B}^{\rm f.c.}.
$$
Here we regard $\vec p^{\,(j)}$ and $\vec{q}$ as 
parts of the input as in the case of the definition of $\frak q$ etc. and
$\rho_{\frak b,\theta}(B)$ is defined in the same way as (\ref{form2612}).
\par
Now as usual we proceed:
\begin{enumerate}
\item[(i)]  We include interior marked points and deform 
$Z_{T}^{\rm f.c.}$ 
to $Z_{T}^{\rm f.c.\frak b}$.
In other words, we include bulk deformations.
\item[(ii)]
We include forgetability data and extend $Z_{T}^{\rm f.c.\frak b}$ to its plus version, which admits ${\bf e}^+$ and ${\bf f}$ as inputs.
\item[(iii)]
We use the $A_{\infty}$ 
functor $\frak f$ in Subsections \ref{subsec:red} and \ref{hunitcansubsec}. 
\item[(iv)]
We eliminate the curvature using the given bounding cochains of objects in ${\bf L}$ or ${\bf F} = {\bf L}
\cup {\bf U}$.
\end{enumerate}
We thus obtain\index[syindex]{ZRb@$Z_{R}^{\bf b}$} 
\begin{equation}\label{form2630}
Z_{T}^{\bf b} :  CH(\sL,\sL;\Lambda_0) \otimes CH(\sU,\sU;\Lambda_0)  \otimes CH(\sL,\sL;\Lambda_0) \to \Lambda_0.
\end{equation}
\begin{lem}
$(\ref{form2630})$ induces a map
$$
Z_{T}^{\bf b} :  HH_*(\sL,\sL;\Lambda_0) \otimes HH_*(\sU,\sU;\Lambda_0)  \otimes HH_*(\sL,\sL;\Lambda_0) \to \Lambda_0.
$$
\end{lem}
\begin{proof}
This is a consequence of the fact that the contributions of (bds1), (bds2), (bds3) correspond
to the Hochishild boundary operators on the first, third and second tensor factors, respectively.
\end{proof}

\subsubsection{A variant of the Cardy relation}
\label{subsubsec:varcard}

We now modify the argument of Section \ref{sec:annuli} to our situation to prove Proposition \ref{prop2623}.

\begin{lem}\label{lem2628}
The map $Z_{T}^{\bf b}$ on the Hochishild homologies is independent of $T \in (0,\infty)$.
\end{lem}
We can prove it in the same way as Lemma \ref{lem414}.

\begin{lem}\label{lem2629}
We have
$$
\lim_{T \to 0}Z_{T}^{\bf b} = Z_{\rm QT}.
$$
\end{lem}
We can prove it in the same way as Lemma \ref{lem412} together with the $c=2$ case of 
Lemma \ref{annlimit}, by comparing the limit $T\to 0$ of  
Figure \ref{Figure11} to  Figure \ref{Figure67}.
\par
Now to complete the proof of Proposition \ref{prop2623} 
it remains to study the limit:
$\lim_{T \to \infty}Z_{T}^{\bf b}$.
For this purpose  we next study the behavior of our moduli space 
${\mathcal M}_{(\vec{\kappa}^{(1)},\vec k^{\,(1)}), (\vec{\kappa}^{\,(2)},\vec k^{(2)}),(\vec{{\upsilon}},\vec r)}(\sL;\sU;\sL;\vec p^{\,(1)},\vec p^{\,(2)},\vec q;T;B)$ when $T$ goes to infinity.
We observe that an object of this limit can be depicted by the next Figure \ref{Figure69}.

\begin{figure}[h]
\centering
\includegraphics[scale=0.3]{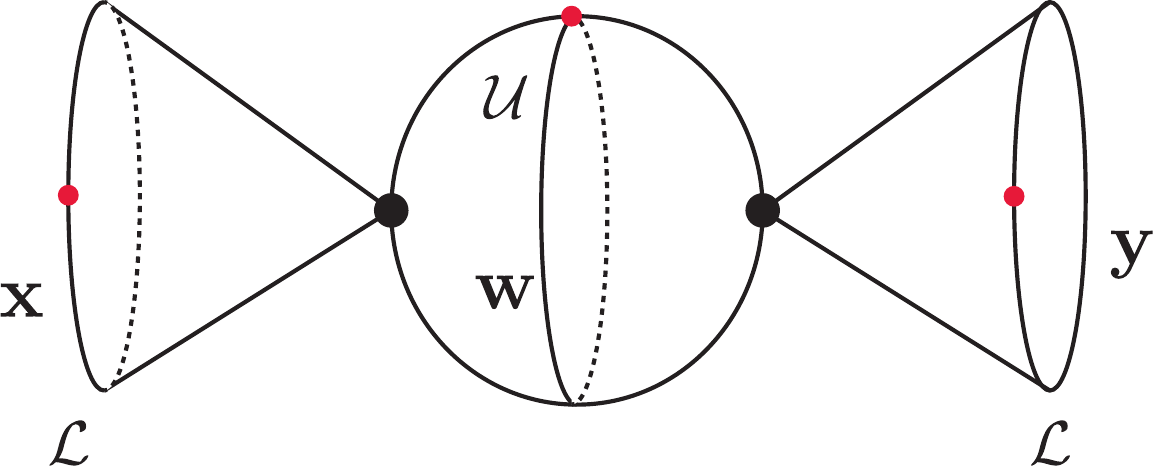}
\caption{$\lim_{T \to \infty}Z_{T}^{\bf b}({\bf x},{\bf w};{\bf y})$}
\label{Figure69}
\end{figure}

Let us define precisely the moduli space of objects which are depicted by Figure \ref{Figure69}.
We begin with defining the moduli space of  holomorphic spheres appearing in the 
middle of Figure \ref{Figure69}.
We consider $S^2$ and its equator $S^1_{\rm eq}$, the north pole $\frak{np} \in S^2$, 
and the south  pole $\frak{sp} \in S^2$.
We decompose $S^2 \setminus S^1_{\rm eq}$ into the northern hemisphere $S^2_{\frak{nh}}$
and the southern hemisphere $S^2_{\frak{sh}}$.
We also consider $\vec{{\upsilon}},\vec r,m$ and $\vec q$ as in Definition \ref{def1625}.

\begin{defn}
We consider $((\Sigma^s,\vec z,z^+_{\frak n},z^+_{\frak s}),u_1^s,u_2^s,\gamma_{12})$ with the following 
properties.
\begin{enumerate}
\item
$\Sigma^s$ is a 2-sphere (possibly) together with trees of sphere bubbles
rooted at $S^2 \setminus S^1_{\rm eq}$.
We denote by $\Sigma^s_{\frak{nh}}$ (resp. $\Sigma^s_{\frak{sh}}$) the 
union of $S^2_{\frak{nh}}$ (resp. $S^2_{\frak{sh}}$) and the trees of sphere bubbles
rooted at it.
\item
$z^+_{\frak n} \in \Sigma^s_{\frak{nh}} \setminus S^1_{\rm eq}$ and
$z^+_{\frak s} \in \Sigma^s_{\frak{nh}} \setminus S^1_{\rm eq}$ are north and south poles of $S^2$ respectively.
\item $\vec z = (z(0),\dots,z(r))$ are points on $S^1_{\rm eq}$.
They are disjoint and respects the cyclic order defined by the 
orientation of $S^1_{\rm eq} = \partial \Sigma^s_{\frak{nh}}$.
We change the enumeration of $\vec z$ and define $\vec z^{\,(12)}$ as in Definition \ref{2624dsef}.
\item
$u_1^s  : \Sigma^s_{\frak{sh}} \setminus \{z^{(12)}_0,\dots,z^{(12)}_K\}
\to X$ is $-J_X$ holomorphic 
and 
$u_2^s  : \Sigma^s_{\frak{nh}} \setminus \{z^{(12)}_0,\dots,z^{(12)}_K\}
\to X$ is $J_X$ holomorphic.  Here $z^{(12)}_i$ is as in the third 
line of (\ref{formnew21644}).
\item 
$\gamma_{12} : S^1_{\rm eq} \setminus \{z^{(12)}_0,\dots,z^{(12)}_K\}
\to \tilde U$ is a smooth map. Here $\tilde U$ is the disjoint union of $\tilde U_{\rho}$
($\rho =0,\dots,K$) and is the equator minus finitely many points.
\item
$\gamma_{12}\left(\overline{z^{(12)}_{i-1},z^{(12)}_{i}}\right) \subset \tilde U_{\rho_i}$.
Here $\rho_{K} = \rho_{-1}$ by convention.
\item
$(u_1,u_2) = i_{U} \circ \gamma_{12}$ when both sides are defined.
Here $i_U = i_{U_{\rho}}$ on $\tilde U_{\rho}$.
\item We require $(u_1(z^{(12)}_i), u_2(z^{(12)}_i)) =  q_i$.
\item
We require a similar switching condition as Definition \ref{def211} (4)
at $z^{(12)}_i$.
\item It is stable in the sense that the group of automorphisms (which we can define in an 
obvious way) is a finite group.
\item
$$
-\int_{\Sigma_{\frak{sh}}^s}(u^s_1)^*\omega + \int_{\Sigma_{\frak{nh}}^s}(u^s_2)^*\omega = E_s.
$$
\end{enumerate}
The isomorphism class of such an object is called a {\it quilted sphere}\index{quilted sphere} and 
the set of all such objects are denoted by
$\overset{\circ}{\mathcal M}^{\rm QS}_{(\vec{{\upsilon}},\vec r,m)}(\sU;E_s)$.
(Here we take the union of various choices of $\vec q$.)
Then we decompose it by homology classes $B_s$ and define
$\overset{\circ}{\mathcal M}^{\rm QS}_{(\vec{{\upsilon}},\vec r,m)}(\sU;B_s)$.
(Here the choice of $\vec q$ is fixed and is included in $B$.)
\end{defn}
We remark that an element of $\overset{\circ}{\mathcal M}^{\rm QS}_{(\vec{{\upsilon}},\vec r,m)}(\sU;E_s)$\index[syindex]{MQSupsilon
@${\mathcal M}^{\rm QS}_{(\vec{{\upsilon}},\vec r,m)}(\sU;E_s)$}
has a representative such that $z(0)$ is the given point on the equator.
\par
We can define their compactifications and denote them as
${\mathcal M}^{\rm QS}_{(\vec{{\upsilon}},\vec r,m)}(\sU;E_s)$
and ${\mathcal M}^{\rm QS}_{(\vec{{\upsilon}},\vec r,m)}(\sU;B_s)$.
Their boundaries correspond to disk bubbles on the seam $= S^1_{\rm eq}$.
\par
We have  evaluation maps:
$$
\aligned
{\rm ev}^{\rm sm} : &{\mathcal M}^{\rm QS}_{(\vec{{\upsilon}},\vec r,m)}(\sU;E_s) \to \bigcup_{\vec q}U^{\vec{\kappa},\vec q}_{\text{\rm source}},
\\
({\rm ev}^{\frak s},{\rm ev}^{\frak n}) : &{\mathcal M}^{\rm QS}_{(\vec{{\upsilon}},\vec r,m)}(\sU;E_s) \to X^2.
\endaligned
$$
Here the map in the first line is obtained by the evaluation at the marked points $\vec z$,  
and the map in the second line is obtained by the evaluation at the marked points $z^+_{\frak s}, z^+_{\frak n}$.
Similar evaluation maps are defined on ${\mathcal M}^{\rm QS}_{(\vec{{\upsilon}},\vec r,m)}(\sU;B_s)$. 
\par
We next discuss moduli spaces which describe the disk components of Figure \ref{Figure69}.
Actually they are moduli spaces of disks appeared in the definition of $\frak p$.
We consider $(\vec{\kappa}^{(j)},\vec k^{(j)},m_j)$, 
$\vec p^{\,(j)}$, as in  Definition \ref{def1625}.  (Here $j=1,2$.)
We then consider the moduli spaces
$
{\mathcal M}_{1,\vec k^{(j)}}((\vec{\kappa}^{\,(j)},\vec p^{\,(j)},m_j);B_j)
$
as in (\ref{formdef915}). ($1 = \ell$ is the number of interior marked points.)
It comes with evaluation maps
\begin{equation}
\aligned
&{\rm ev}^{\partial} : {\mathcal M}_{1,\vec k^{(j)}}((\vec{\kappa}^{\,(j)},\vec p^{\,(j)},m_j);B_j),
\to L^{\vec{\kappa}^{(j)},\vec p^{\,(j)}}_{\text{\rm source}} \\
&{\rm ev}^{+} : {\mathcal M}_{1,\vec k^{(j)}}((\vec{\kappa}^{\,(j)},\vec p^{\,(j)},m_j);B_j) 
\to X.
\endaligned
\end{equation}
We now consider the next fiber product over $X$:
\begin{equation}\label{form2632}
\aligned
& {\mathcal M}_{1,\vec k^{(1)}}((\vec{\kappa}^{(1)},\vec p^{\,(1)},m_1);B_1) \\
& {}_{{\rm ev}^{+}}\times_{{\rm ev}^{\frak s}} {\mathcal M}^{\rm QS}_{(\vec{{\upsilon}},\vec r,m)}(\sU;B_s)
\\
&{}_{{\rm ev}^{\frak n}}\times_{{\rm ev}^{+}}
{\mathcal M}_{1,\vec k^{(2)}}((\vec{\kappa}^{(2)},\vec p^{\,(1)},m_2);B_2)
\endaligned
\end{equation}
and denote it by
$$
{\mathcal M}^{\rm DSD}_{(\vec{\kappa}^{(1)},\vec k^{(1)},m_1), (\vec{\kappa}^{(2)},\vec k^{(2)},m_2),(\vec{{\upsilon}},\vec r,m)}(\sL;\sU;\sL;B_1,B_s,B_2).
$$
The is the moduli space of objects depicted by Figure \ref{Figure69}.
It comes with evaluation maps:
\begin{equation}
\aligned
&{\rm ev}^{\partial}_-: 
{\mathcal M}^{\rm DSD}_{(\vec{\kappa}^{(1)},\vec k^{(1)},m_1), (\vec{\kappa}^{(2)},\vec k^{(2)},m_2),(\vec{{\upsilon}},\vec r,m_)}(\sL;\sU;\sL;B_1,B_s,B_2) \to 
L^{\vec{\kappa}^{(1)},\vec p^{\,(1)}}_{\text{\rm source}}, \\
&{\rm ev}^{\partial}_+: 
{\mathcal M}^{\rm DSD}_{(\vec{\kappa}^{(1)},\vec k^{(1)},m_1), (\vec{\kappa}^{(2)},\vec k^{(2)},m_2),(\vec{{\upsilon}},\vec r,m_)}(\sL;\sU;\sL;B_1,B_s,B_2) \to L^{\vec{\kappa}^{(2)},\vec p^{\,(2)}}_{\text{\rm source}}, \\
&{\rm ev}^{\rm sm} : {\mathcal M}^{\rm DSD}_{(\vec{\kappa}^{(1)},\vec k^{(1)},m_1), (\vec{\kappa}^{(2)},\vec k^{(2)},m_2),(\vec{{\upsilon}},\vec r,m)}(\sL;\sU;\sL;B_1,B_s,B_2) \to U^{\vec{\kappa},\vec q}_{\text{\rm source}}.
\endaligned
\end{equation}
For ${\bf h}^{(j)} \in \Omega(L^{\vec{\kappa}^{(j)},\vec p^{(j)}}_{\text{\rm source}})$ 
and 
${\bf h} \in \Omega(U^{\vec{\kappa},\vec q}_{\text{\rm source}})$
we define:
$$
Z_{B_1,B_s,B_2}^{\rm f.c.}({\bf h}^{(1)},{\bf h};{\bf h}^{(2)})
= 
\int_{\star} ({\rm ev}^{\partial}_-)^*({\bf h}^{(1)}) \wedge ({\rm ev}^{\rm sm})^*{\bf h} \wedge({\rm ev}^{\partial}_+)^*({\bf h}^{(2)})
$$
where $\star ={\mathcal M}^{\rm DSD}_{(\vec{\kappa}^{(1)},\vec k^{(1)},m_1), (\vec{\kappa}^{(2)},\vec k^{(2)},m_2),(\vec{{\upsilon}},\vec r,m_)}(\sL;\sU;\sL;B_1,B_s,B_2)$.
We then put
$$
Z_{\rm DSD}^{\rm f.c.} = \sum_{B_1,B_2,B_s} e^{\rho(\theta,B)} T^{B(\omega)} Z_{B_1,B_s,B_2}^{\rm f.c.}.
$$
Here $B$ is the homology class obtained by concatenating $B_1,B_2,B_s$.
$\rho(\theta,B)$ is defined in the same way as (\ref{form2612}).
We then take the steps (i)(ii)(iii)(iv) as before to obtain \index[syindex]{ZDSD@$Z_{\rm DSD}^{\bf b}$}
\begin{equation}\label{form26302}
Z_{\rm DSD}^{\bf b} :  CH(\sL,\sL;\Lambda_0) \otimes CH(\sU,\sU;\Lambda_0)  \otimes CH(\sL,\sL;\Lambda_0) \to \Lambda_0
\end{equation}
in a similar way as (\ref{form2630}).
\begin{lem}\label{lem2631}
We have:
$$
\lim_{T \to \infty}Z_{T}^{\bf b} = Z^{\bf b}_{\rm DSD}.
$$
\end{lem}
We can prove the lemma in the same way as the $c=1$ case of Lemma \ref{annlimit} 
by comparing Figures \ref{Figure69} and \ref{Figure11}.

Now Proposition \ref{prop2623} follows from the next lemma.

\begin{figure}[h]
\centering
\includegraphics[scale=0.3]{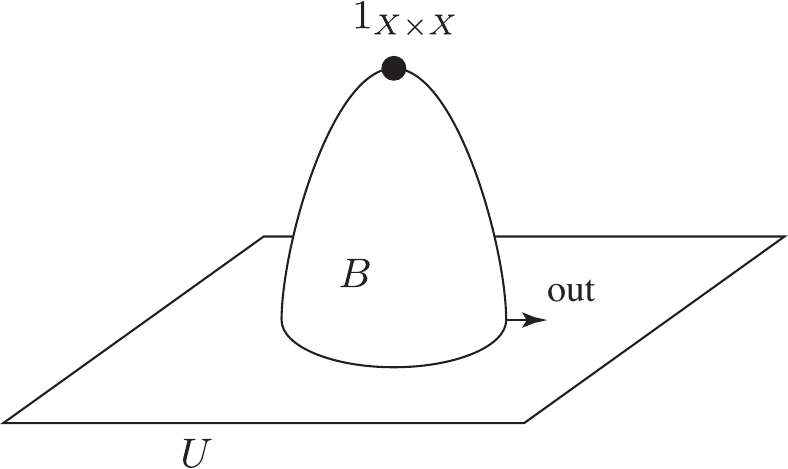}
\caption{$
\langle \frak q^{\bf b}(1_{X\times X}),{\bf w}\rangle_{\rm HH}
$}
\label{Figure70}
\end{figure}

\begin{lem}
Let ${\bf x} \in HH_*(\sL,\sL;\Lambda)$ be an element such that $\frak p^{\bf b}({\bf x}) = {\bf e}_X$.
Then we have:
$$
Z^{\bf b}_{\rm DSD}({\bf x},{\bf w};{\bf x}) = 0
$$
unless ${\bf w}$ contains a term  $c[{\rm vol}_{U_{\rho}}]$ for some $\rho$.
Moreover
\begin{equation}\label{form2635}
Z^{\bf b}_{\rm DSD}({\bf x},c[{\rm vol}_{U_{\rho}}];{\bf x}) = \pm c.
\end{equation}
\end{lem}
\begin{proof}
We remark that ${\bf x}$, ${\bf w}$ here corresponds ${\bf y}_0$, $\bf z$ in 
Proposition \ref{prop2512}.
\par
For simplicity of notation we consider the case $\frak b = b_{U_{\rho}} = 0$ and $\theta_{U_{\rho}} = 0$.
We can modify the proof so that it works in the general case in a straightforward way.
(However then the notation becomes complicated.)
\par
We first observe that the map $Z_{\rm DSD}$ in the homology level is independent of the perturbation.
Among fiber product factors in (\ref{form2632}) 
we take $\frak p$-perturbation (See Proposition \ref{existmultipolu2}) for the 1st and the third factors.
Then using the fact that $\frak p^{\bf b}({\bf x}) = 1_X$
we obtain the next formula. Suppose ${\bf w}$ is represented by ${\bf h} 
\in \Omega(U^{\vec{\kappa},\vec q}_{\text{\rm source}})$.
$$
Z_{\rm DSD}({\bf x},{\bf w};{\bf x})
= \sum_{B_s} T^{\omega(B_s)} \int_{{\mathcal M}^{\rm QS}_{(\vec{{\upsilon}},\vec r,m)}(\sU;B_s)}
 ({\rm ev}^{\rm sm})^*{\bf h}.
$$
By taking appropriate perturbation the right hand side is 
equal to 
$
\langle \frak q^{\bf b}(1_X),{\bf w}\rangle_{\rm HH}
$
up to sign. 
(Here $\frak q^{\bf b}$ is the closed open map of $X \times X$.)  
See Figure \ref{Figure70}, 
whose double is the configuration appearing 
in Figure \ref{Figure69} when ${\bf x} ={\bf y}$ with $\frak p^{\bf b}({\bf x}) = 1_X$. Therefore by a similar cobordism argument 
as  in Section \ref{nontrivialsec} we obtain the lemma.
\end{proof}

Thus the proof of Theorem \ref{thm253} is completed except orientation and sign issue, 
which we will discuss in the next subsubsection.

\subsubsection{Reflexion principle and orientation}
\label{subsubsec:ref}

We discuss orientation and sign in this subsubsection.
The strategy is to use the reflexion principle to 
reduce the orientation of various moduli spaces involved 
to those we studied in Section \ref{sec:signCardy}.
\par
We first consider the moduli space
  ${\mathcal M}_{\rm QT}(\vec\kappa_1,\vec\kappa_2,\vec\rho;\vec p_1,\vec p_{2},\vec p_{12};\vec k_1,\vec k_2,\vec k_{12};
p_{\pm \infty};E)$ in Definition \ref{def516}.
Their elements are  equivalence classes of
$$(\Sigma;\vec z^+_1,\vec z^+_2;\vec z_1,\vec z_{12},\vec z_2;\vec\vec{w}_1,\vec\vec{w}_{12},\vec\vec{w}_2;u_1,u_2;\gamma_1,\gamma_{12},\gamma_2).$$
Since our concern is only on orientations it suffices to consider only 
the case $\Sigma = [-1,1] \times \R = [-1,1] \times \sqrt{-1}\R \subset \C$.
\par
We renumber the marked points on $\{\pm 1\} \times \R$ to obtain $\vec z_1 = (z_1(0),\dots,z_1(k_1))$ and 
$\vec z_2 = (z_2(0),\dots,z_2(k_2))$. Note that $\vec z_1$ is enumerated downward and $\vec z_2$ is enumerated
upward with respect to the $\sqrt{-1}\R$ coordinate.
We revert the enumeration of $\vec z_1$ to obtain $\vec z^{\,\prime}_1$ such that
$z'_1(i) = z_1(k_1-i)$.
\par
We consider a map $z \mapsto \tilde z = -\overline z$,\index[syindex]{ztilde@$\tilde z$} where  $\overline z$ is the complex conjugate.
It sends  $\{-1\} \times \sqrt{-1}\R$ to $\{+1\} \times  \sqrt{-1}\R$.
We further restrict ourselves to the subset such that:
\begin{enumerate}
\item[(*)]$\tilde z'_1(i) \ne z_2(j)$ for any $i,j$. 
\end{enumerate}
We renumber the disjoint union 
$$
\{ \tilde z'_1(i) \mid i \} \cup \{z_2(j) \mid j\} = \{z_3(i) \mid i= 0,\dots,k_1+k_2+1\}
$$
such that the imaginary part of $z_3(i)$ is increasing.
Note the way $\vec z^{\,\prime}_1$ and $\vec z_2$ determine $\vec z_3$ is 
parametrized by shuffles:\index{shuffle}
$\sigma = (\sigma_1,\sigma_2) : \{0,\dots,k_1\} \sqcup \{0,\dots,k_2\} \to \{0,\dots,k_1+k_2+1\}$,
that is, the bijection from the disjoint union such that each of $\sigma_i$ is order preserving.
\par
We next define a $J_{X\times X}$ holomorphic map $u_3 : [0,1] \times \sqrt{-1}\R \to X^2$ by
\begin{equation}\label{form2636}
u_3(z) = (u_1(\tilde z),u_2(z)),
\end{equation}
and $\gamma_3 : \{1\} \times \sqrt{-1}\R \setminus \vec z_3 \to \tilde L \times \tilde L$
by
\begin{equation}\label{form2637}
\gamma_3(z) = (\gamma_1(\tilde z),\gamma_2(z)).
\end{equation}
For each shuffle $\sigma$ we obtain $\kappa^{\sigma}_i = (\kappa^{\sigma}_{1,i},\kappa^{\sigma}_{2,i})$
such that 
$$
u_3(z_3(i)) \in (L_{\kappa^{\sigma}_{1,i-1}} \times L_{\kappa^{\sigma}_{2,i-1}}) \cap
(L_{\kappa^{\sigma}_{1,i}} \times L_{\kappa^{\sigma}_{2,i}}).
$$
Thus $(([0,1]\times \R,\vec z_3);u_3;\gamma_{12},\gamma_3)$ is an element in the moduli space 
which defines the bi-module structure 
\begin{equation}\label{form2638}
\aligned
B\sU(U_{\rho_0},U_{\rho_r})&\otimes CF(U_{\rho_r},(L_{\kappa_{1,0}}\times L_{\kappa_{2,0}})) \\
&\otimes B\sL^2((L_{\kappa_{1,0}}\times L_{\kappa_{2,0}}),(L_{\kappa_{1,k_1+k_2-1}}\times L_{\kappa_{2,k_1+k_2-1}}))\\
&\longrightarrow
CF(U_{\rho_r},(L_{\kappa_{1,k_1+k_2-1}}\times L_{\kappa_{2,k_1+k_2-1}})).
\endaligned
\end{equation}
The orientation of this moduli space is defined in the previous sections.\footnote{Actually the intersection between
$L_{\kappa^{\sigma}_{1,i-1}} \times L_{\kappa^{\sigma}_{2,i-1}}$ and
 $L_{\kappa^{\sigma}_{1,i}} \times L_{\kappa^{\sigma}_{2,i}}$ may not be 
 transversal.  They may not  coincide either.  On the other hand, the intersection is always of Morse-Bott type.
 The study of orientation in \cite{ono} includs the Morse-Bott case. The argument of Section \ref{sec:signCardy}
 can be applied to such a case also. Note that the intersection of $U_{\rho}$ and 
 $L_{\kappa} \times L_{\kappa'}$ is transversal.}
 We can use the identification to define the orientation of ${\mathcal M}_{\rm QT}(\vec\kappa_1,\vec\kappa_2,\vec\rho;\vec p_1,\vec p_{2},\vec p_{12};\vec k_1,\vec k_2,\vec k_{12};
p_{\pm \infty};E)$ by this isomorphism.
\par
Note that we revert the enumeration of the marked points on $\{-1\} \times \sqrt{-1}\R$.  This corresponds to 
changing from a left $\sL$ module to a right $\sL$ module.  We also use the 
Koszul conventions.  So the identification by the reflexion principle converts $A_{\infty}$ formulas with sign to 
$A_{\infty}$ formulas with sign.  (See \cite[Example 17.3]{fukaya:functor}.)
\par
We remark that we need to check that the induced orientation on 
the moduli space ${\mathcal M}_{\rm QT}(\vec\kappa_1,\vec\kappa_2,\vec\sigma;\vec p_1,\vec p_{2},\vec p_{12};\vec k_1,\vec k_2,\vec k_{12};
p_{\pm \infty};E)$
is independent of the choice of shuffles $\sigma$.  We can prove it by using \cite[Lemma 8.4.3]{fooo092}.

We next discuss the orientation of
${\mathcal M}^{H;m_1,m_2}_{\rm QT}(\vec\kappa_1,\vec\kappa_2,\vec\sigma;\vec p_1,\vec p_{2},\vec p_{12};\vec k_1,\vec k_2,\vec k_{12};
p_{\pm \infty};E)$
appearing in Definition \ref{defn2620}.
We enumerate the marked points on $\{-1\} \times \sqrt{-1}\R$  (resp. $\{+1\} \times \sqrt{-1}\R$) 
by downward (resp. upward) order  and denote it by $\vec z_1 = (z_1(0),\dots,z_1(k_1))$ and 
$\vec z_2 = (z_2(0),\dots,z_2(k_2))$.
We remark that
$\tilde z_1(m_1) =  z_2(m_2)$. We restrict ourselves to the subset such that:
\begin{enumerate}
\item[(**)]$\tilde z_1(i) \ne  z_2(j)$ for any $i \ne m_1, j \ne m_2$. 
\end{enumerate}
We put $z'_1(j) = z_1(k_1 - j)$.
We renumber the disjoint union 
$$
\aligned
&\{ \tilde z'_1(i) \mid i=0,\dots,k_1-m_1-1 \} \cup \{z_2(j) \mid j=0,\dots, m_2-1\}  \\
&= \{z_{3-}(i) \mid i= 1,\dots,k_1-m_1+m_2\}
\endaligned
$$
such that the imaginary part of $z_{3-}(i)$ is increasing.
We renumber the disjoint union 
$$
\aligned
&\{ \tilde z'_1(i) \mid i=k_1-m_1+1,\dots,k_1 \} \cup \{z_2(j) \mid j=m_2+1,\dots, k_2\}  \\
&= \{z_{3+}(i) \mid i= 1,\dots,k_2-m_2+m_1\}
\endaligned
$$
such that the imaginary part of $z_{3+}(i)$ is increasing.
We put
$$
z_3(0) = \tilde z_1(m_1) = z_2(m_2)
$$
and
$$
\vec z_3 = (z_{3-}(1),\dots,z_{3-}(k_1-m_1+m_2),z_3(0),z_{3+}(1),\dots,z_{3+}(k_2-m_2+m_1)).
$$
Here the right hand side is a $k_1+k_2+1$-tuple of points on $\{+1\} \times \sqrt{-1}\R$
with increasing imaginary part.
\par
The combinatorial type of such $\vec z_3$ is classified by  
pairs of shuffles
$$
\aligned
\sigma_- = (\sigma_{1-},\sigma_{2-}) : \{0,\dots,k_1-m_1-1\} 
&\sqcup \{0,\dots,m_2-1\}\\ &\to \{1,\dots,k_1-m_1+m_2\}, \\
\sigma_+ = (\sigma_{1+},\sigma_{2+}) : \{k_1-m_1+1,\dots,k_1\} &\sqcup \{m_2+1,\dots,k_2\} \\
&\to \{k_1-m_1+m_2+2,\dots,k_1+k_2+1\}
\endaligned
$$
in an obvious way.

We define $u_3 : [0,1] \times \sqrt{-1}\R \to X^2$ and $\gamma_3 : \{1\} \times \sqrt{-1}\R \setminus \vec z_3 \to \tilde L \times \tilde L$
by (\ref{form2636}), (\ref{form2637}).

For each $\vec\sigma = (\sigma_-,\sigma_+)$ we obtain $\kappa^{\vec\sigma}_i = (\kappa^{\vec\sigma}_{1,i},\kappa^{\vec\sigma}_{2,i})$
such that 
$$
u_3(z_3(i)) \in (L_{\kappa^{\vec\sigma}_{1,i-1}} \times L_{\kappa^{\vec\sigma}_{2,i-1}}) \cap
(L_{\kappa^{\vec\sigma}_{1,i}} \times L_{\kappa^{\vec\sigma}_{2,i}}).
$$
Thus $(([0,1]\times \R,\vec z_3);u_3;\gamma_{12},\gamma_3)$ is an element of the moduli spaces 
which define the bi-module structure 
(\ref{form2638}).
\par
Therefore the orientation of ${\mathcal M}^{H;m_1,m_2}_{\rm QT}(\vec\kappa_1,\vec\kappa_2,\vec\rho;\vec p_1,\vec p_{2},\vec p_{12};\vec k_1,\vec k_2,\vec k_{12};
p_{\pm \infty};E)$ is defined by identifying it with a moduli space appearing in the previous sections.
\par
We have defined the sign and the orientation appearing in the operations $\frak m^H_{*,*\vert 1\vert *}$ and
$\frak m_{*,*\vert 1\vert *}$.
\par
We now consider the map $Z_{\rm QT}$.
The moduli space which is used to define it appears in Figure \ref{Figure67}, 
where the disk in the right  (together with marked points and maps) is an element of
${\mathcal M}^{H;m_1,m_2}_{\rm QT}(\star)$ for some $\star$, $m_1$, $m_2$
and the disk in the left  (together with marked points and maps) is 
an element of ${\mathcal M}_{\rm QT}(\star)$ for some $\star$.
Restricting to a subset where there is no bubble and Conditions (*),(**) are satisfied,
we can use the reflexion principle to each of those disks to obtain an object 
depicted in Figure \ref{Figure67revrev} below.

\begin{figure}[h]
\centering
\includegraphics[scale=0.5]{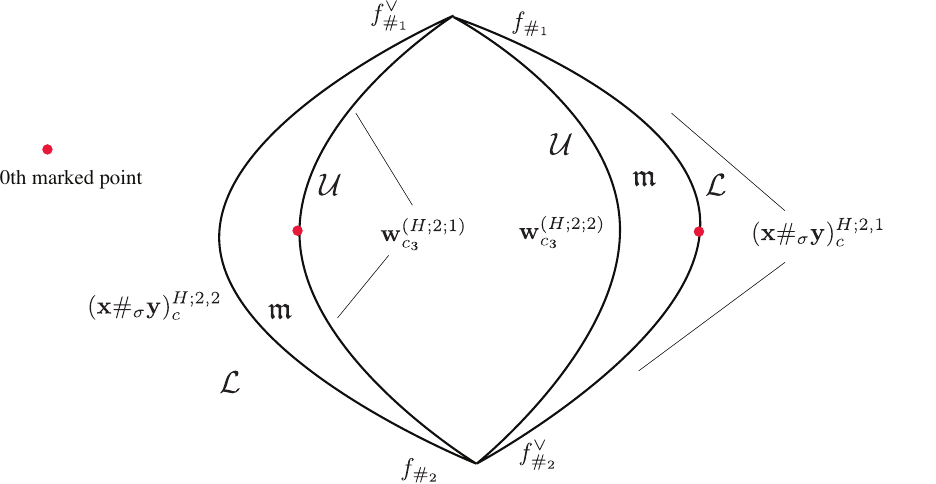}
\caption{Reflexion principle applied to $Z_{\rm QT}$}
\label{Figure67revrev}
\end{figure}

Here
$
{\bf x}\#_{\sigma}{\bf y}
$
is a shuffle of ${\bf x}$ and ${\bf y}$ which we define as follows.
Let ${\bf x} = x_0 \otimes {\bf x}_+$, ${\bf y} = y_0 \otimes {\bf y}_+$,
where ${\bf x}_+ \in B_{k_1}\sL[1]$, ${\bf y}_+ \in B_{k_2}\sL[1]$.
We consider a shuffle
\begin{equation}\label{form2639}
\sigma = (\sigma_1,\sigma_2) : \{1,\dots,k_1\} \sqcup \{1,\dots,k_2\} \to \{1,\dots,k_1+k_2\}.
\end{equation}
We now put: $
{\bf x}\#_{\sigma}{\bf y}
= (x_0 \otimes y_0) \otimes ({\bf x}^{\rm op}_+\#_{\sigma}{\bf y}_+)$
where
$$
({\bf x}^{\rm op}_+\#_{\sigma}{\bf y}_+)_i =
\begin{cases}
x_{+j} \otimes {\bf e}    & \text{if $\sigma_1(k_1-j+1) = i$,} \\
{\bf e} \otimes  y_{+j}  & \text{if $\sigma_2(j) = i$.}
\end{cases}
$$
Here ${\bf e}$ is the volume form on $\tilde L_{\kappa}$.
($\kappa$ is determined automatically from the condition that ${\bf x}\#_{\sigma}{\bf y}$ 
becomes a Hochishild chain.)
\par
$({\bf x}\#_{\sigma}{\bf y})^{(H;2,1)}$, $({\bf x}\#_{\sigma}{\bf y})^{(H;2,2)}$
then are defined by (\ref{DeltaH}).
\par
Now the reflexion principle {\it could} imply the next formula:
\begin{equation}\label{form2640}
Z_{\rm QT}(\text{\bf x},\text{\bf w};\text{\bf y})
= \sum_{\sigma} Z(\text{\bf w};{\bf x}\#_{\sigma}{\bf y}).
\end{equation}
Here $Z$ is as in Definition \ref{eq:formula_Z-pairing}.
The sum is taken over all shuffles (\ref{form2639}).
\par
We however immediately remark the following.
Our discussion on the reflexion principle is {\it not} enough to show 
Formula (\ref{form2640}).  In fact the reflexion principle is applied only 
on a subset where there is no bubble\footnote{See however Remark \ref{rem2633}.} and 
conditions (*), (**) are satisfied.
To prove (\ref{form2640}) we need to extend it to the compactifications 
and find the way to obtain Kuranishi structures and their CF-perturbations 
(and their cobordisms) of the compactification so that it is compatible with reflexion.\footnote{This kinds of argument 
is worked out in \cite[Subsection 6.2]{fooo:inv} in a different but somewhat similar situation.}
We believe one can do it but to work out its detail is much cumbersome and lengthy. We therefore take a short cut in this paper 
to avoid proving (\ref{form2640}) itself.

We remark that here we only concern with orientations and sign issue.
So, for our purpose, it suffices to compare the local contributions of the left and the right hand sides of 
(\ref{form2640}) at the points where the reflexion principle is applicable.
Note that the maps $Z_{\rm QT}$ and $Z$ are defined by the integration on the moduli spaces. 
Therefore  the notion of  local contributions makes sense.  Moreover in the part 
where the reflexion principle is applicable we can compare local 
contributions.

We recall that the sign for the map $Z$ is given in Definition \ref{defnZ} and 
discussion thereafter.  We also remark that the orientation 
of the moduli space of the left hand side is defined by using the reflexion principle.
Therefore it is obvious that we can find the sign ${\maltese}$ in (\ref{formula23999})
uniquely so that (\ref{form2640}) holds with sign in the above sense.  (Namely 
local contributions at the point where the reflexion principle applies 
coincide together with sign.)

We have thus fixed the sign  ${\maltese}$ in (\ref{formula23999}).
(We do not provide an explicit formula since we do not need it.)

\begin{rem}\label{rem2633}
To prove various properties of $Z_{\rm QT}$ and $\varphi$, $\varphi^H$ 
we use the reflexion principle and reduce the proof to the corresponding statements 
on $Z$ or $\frak m$.  For this purpose we need to study certain codimension $1$ 
components, such as the one with disk bubbles.  We can easily observe 
that to the codimension $1$ components appearing in those discussions, 
the reflexion principle extends.
\end{rem}

We now proceed to the discussion on orientations and signs appearing in the Cardy relation 
(Subsubsections \ref{subsubsec:qanu} and \ref{subsubsec:varcard}).
We observe that to the moduli spaces we use there (which are depicted in 
Figures \ref{Figure11} and \ref{Figure69}) we can apply the reflexion principle in a 
similar way using appropriate shuffles.
So we can prove a similar formula as (\ref{form2640}) for
$Z_{R}^{\bf b}$ and $Z_{\rm DSD}$.
More precisely we can identify local contributions to 
$Z_{R}^{\bf b}$ and $Z_{\rm DSD}$ to ones of corresponding 
operators obtained by the moduli spaces used to prove the Cardy relation in  Section \ref{sec:annuli}.
Then the sign and orientation part of the proof of Lemmas \ref{lem2628},  \ref{lem2629}, 
\ref{lem2631} follow from the corresponding arguments in Section \ref{sec:signCardy}.
We have now completed the orientation and sign part of the proof of Theorem \ref{thm253}.
\qed

\subsection{Correction of a sign in \cite{fukaya:functor}}
\label{subsec:corsign}
\subsubsection{$J \mapsto -J$, orientation and opposite category.}
\label{subsubsec:Jto-J}
Here we explain an error mentioned in Remark \ref{rem261}. %\marginpar{Subsubsection is added. 2026 Feb. KF}
 \cite[Theorem 3.54]{fukaya:functor} claims the following (*).
 Let $(X,\omega)$ be a symplectic manifold
 and ${\bf L}$   a finite set of $V$-relatively spin Lagrangian 
 submanifolds.
 We regard $L \in {\bf L}$ also as a Lagrangian submanifold
 of $(X,-\omega)$.
It is then $TX \oplus V$-relatively spin.
 We write ${\bf L}'$ the collection of $L$'s regarded as $TX \oplus V$-relatively spin 
 Lagrangian submanifolds of $(X,-\omega)$.
 \par
 Let $\sL$ (resp. $\sL'$) be the filtered $A_{\infty}$ category associated to 
 $(X,\omega,{\bf L})$ (resp. $(X,-\omega,{\bf L}')$.
 \begin{enumerate}
 \item[(*)]
 $\sL'$ is isomorphic to the opposite category $\sL^{\rm op}$ of $\sL$.
 \end{enumerate}
In \cite{fukaya:functor} there is an error in the proof of this claim.
 Our correction is the next Proposition \ref{2636prop}.
\begin{defn} 
For a filtered $A_{\infty}$ category $\Cat$ we define 
$-\Cat$ as follows.
\begin{enumerate}
\item The set of objects of $-\Cat$ is the same as $\Cat$.
\item  The structure operations $\frak m^{\Cat}$, $\frak m^{-\Cat}$ 
are related by
$$
\frak m^{-\Cat} = -\frak m^{\Cat}.
$$
\end{enumerate}
\end{defn}
It is easy to see that if $\frak m^{\Cat}$ satisfies the $A_{\infty}$ relation then 
$-\frak m^{\Cat}$ satisfies the $A_{\infty}$ relation.

\begin{prop}\label{2636prop}
$\sL'$ is isomorphic to $-\sL^{\rm op}$.
\end{prop}

\begin{proof}
The error of the proof of (*) (= \cite[Theorem 3.54]{fukaya:functor}) lies 
in the proof of \cite[Lemma 3.57]{fukaya:functor}.
Following \cite{fukaya:functor} we use the trick to replace a finite of mutually
transversal Lagrangians by a single immersed Lagrangian. 
Let $L$ be an immersed Lagrangian. We take $\vec a = (a_0,\dots,a_k)$ where 
$a_i$ is either a self-intersection point or the symbol ${\rm d}$ which 
means diagonal intersection.
We define 
$$
\mathcal M((X,J),\vec a;E;d)
$$
to be the isomorphism classes of stable pseudo-holomorphic disks,
with boundary condition specified by $\vec a$, energy $E$
and in the components of virtual dimension $d$.
Let $h_i$ be $1$ if $a_i$ is a self intersection point and 
$h_i$ be a differential form of $\tilde L$ if $a_i = {\rm d}$.
We assume $\sum \deg h_i = d$.
We put
$$
\vec a^{\rm op} = (a_0,a_{k},\dots,a_1).
$$
The corrected version of \cite[Lemma 3.57]{fukaya:functor} is the following.
\begin{lem}\label{lem2637}
\begin{equation}
\aligned
&(-1)^{\epsilon(h_k,\dots,h_1)}\int_{{{\mathcal M}}((X,-J_X);L;\vec a^{\rm op};E;d)}
{\rm ev}^*(h_0 \times h_k \times \dots \times h_1) \\
&= (-1)^{\maltese} \int_{{{\mathcal M}}((X,J_X);L;\vec a;E;d)}
{\rm ev}^*(h_0 \times h_1 \times \dots \times h_k)
\endaligned
\end{equation}
where
\begin{equation}\label{form2642}
\maltese = 1 + k + \sum_{1\le i<j\le k}\deg' h_i\deg' h_j + \varepsilon + \epsilon(h_1,\dots,h_k).
\end{equation}
Here $\varepsilon =0$ if and only if
$
d - (k-2)
$
is divisible by $4$. Otherwise, ${\varepsilon} = 1$. $\epsilon(h_1,\dots,h_k)$
is as in $(\ref{212plus1})$.
\end{lem}
In \cite[Lemma 3.57]{fukaya:functor} it was written 
$
\maltese = 1 + \sum_{1\le i<j\le k}\deg' h_i\deg' h_j + \varepsilon.
$\footnote{$\epsilon(h_1,\dots,h_k)$, $\epsilon(h_k,\dots,h_1)$ were not explicitly written 
in \cite[Lemma 3.57]{fukaya:functor}. However they were there.  In fact the sign 
of \cite[Lemma 3.57]{fukaya:functor} is based on \cite{ono}.}
\begin{proof}
The sign $\maltese$ coincides with one in \cite[Theorem 4.1]{fooo:inv}.
We identify  ${{\mathcal M}}((X,J_X);L;\vec a;E;d)$
and ${{\mathcal M}}((X,-J_X);L;\vec a^{\rm op};E;d)$ 
by the map
$$
(D^2;\vec z,u) \to (D^2;\vec z',u')
$$
where $\vec z = (z_0,\dots,z_k)$ is the boundary marked points 
and 
$$
\vec z' = (\overline z_0, \overline z_k,\dots,\overline z_1).
$$
$u : D^2 \to X$ is a $J$-holomorphic map satisfying the boundary 
condition specified by $L$ and $\vec a$.
Then
$$
u'(z) = z(\overline z).
$$
The sign comes from the following:
\begin{enumerate}
\item 
The sign coming from $u \mapsto u'$.
\item
The sign coming from $\vec z \mapsto \vec z'$.
\item
The sign coming from $h_0 \times h_1 \times \dots \times h_k \mapsto h_0 \times h_k \times \dots \times h_1$.
\end{enumerate}
Using the fact  $u \mapsto u'$ induces the complex conjugation on the deformation space 
and the index calculation, (1) gives $\varepsilon$.
\par
(2) gives $k+1 + \#\{(i,j) \mid 0 \le i < j \le k\}$.
\par
$\#\{(i,j) \mid 0 \le i < j \le k\}$ (3) together with difference of correction terms 
$\epsilon(h_1,\dots,h_k) - \epsilon(h_k,\dots,h_1)$ gives 
$\sum_{1\le i<j\le k}\deg' h_i\deg' h_j$.
Hence the proof.
\end{proof}
We remark that in \cite[Lemma 3.57]{fukaya:functor} the term $k$ which 
appears by the change of the directions of $z_1,\dots,z_k$ was overlooked.

As is explained in the proof of \cite[Lemma 3.57]{fukaya:functor} the term
$\varepsilon$ in (\ref{form2642}) cancels with the effect of changing $V$ relative spin structures of $L$ 
to $V \oplus TX$ relative spin structure. Therefore we have
$$
\frak m^{\sL'}(h_k,\dots,h_1) = (-1)^{k+1 + \sum_{1\le i<j\le k}\deg' h_i\deg' h_j} \frak m^{\sL}(h_1,\dots,h_k)
$$
On the other hand, by definition
$$
\frak m^{\sL^{\rm op}}(h_k,\dots,h_1) = (-1)^{1 + \sum_{1\le i<j\le k}\deg' h_i\deg' h_j} \frak m^{\sL}(h_1,\dots,h_k).
$$
We define a map $\frak{minus} : \Cat'(c,c') \to \Cat'(c,c')$ by 
$$
\frak{minus}(x) = -x.
$$
Then we have
$$
\frak m^{\sL'} \circ (\frak{minus} \otimes \dots \otimes \frak{minus}) =  - \frak{minus} \circ \frak m^{\sL^{\rm op}}.
$$
Proposition \ref{2636prop} follows.
\end{proof}
We have the next corollary from the proof of Lemma \ref{lem2637}.

\begin{cor}\label{cor2638}
If $b \in H(\tilde L \otimes_X\tilde L;\Lambda_0)$ is a (weak) bounding cochain of 
the object $L$ of  $\sL$ then $-b$ is a (weak) bounding cochain of the object 
$L'$ of $\sL'$.
\end{cor}

\subsubsection{A correction to the construction of $A_{\infty}$ trimodule.}
\label{subsubsec:cortrimodule}

The correction of the previous subsubsection affects the 
construction of the $A_{\infty}$ trimodule
in \cite[Section 5.2]{fukaya:functor}, which we use as Theorem \ref{themKuneth}.
In this subsection we explain how we modify the construction and adapt the correction 
in the previous subsubsection.
We use the following simple algebraic lemma.
Let $\Cat_i$, $i=1,2,3$ be filtered $A_{\infty}$ categories 
and $\Dat$ a left $\Cat_1$, $\Cat_2$ and right $\Cat_3$ trimodule.
For $c_i \in {\rm OB}(\Cat_i)$ we have a free $\Lambda_0$ graded module
$$
\Dat(c_1,c_2;c_3)
$$
and also structure operations:
$$
\aligned
\frak m_{k_1,k_2\vert 1\vert k_3} : B_{k_1}\Cat_1[1](c_1,c'_1) \otimes B_{k_2}\Cat_2[1](c_2,c'_2) &\otimes \Dat(c'_1,c'_2;c'_3)
 \otimes B_{k_3}\Cat_2[1](c'_3,c_3) \\ &\to \Dat(c_1,c_2;c_3).
 \endaligned
$$
We put
$$
\Dat[1]^k(c_1,c_2;c_3) = \Dat^{k+1}(c_1,c_2;c_3)
$$
and
\begin{equation}
\frak m'_{k_1,k_2\vert 1 \vert k_3}({\bf x},{\bf w};v[1];{\bf y}) = - (-1)^{\deg'{\bf w}} \frak m_{k_1,k_2\vert 1\vert k_3}({\bf x},{\bf w};v;{\bf y})[1].
\end{equation}
\begin{lem}\label{lem1639}
$(\Dat[1],\frak m')$ is a left $-\Cat_1$, $\Cat_2$ right $\Cat_3$ trimodule.
\end{lem}
The proof is immediate from the definition.
\begin{proof}[Proof of Theorem \ref{themKuneth}]
For simplicity we consider the case when $\theta_L = 0$, $\frak b =0$.
(The case $\theta_L  \ne 0$ can be proved in the same way.)
We use the moduli space 
${\mathcal M}_{\rm QT}(\vec\kappa_1,\vec\kappa_2,\vec\rho;\vec p_1,\vec p_{2},\vec p_{12};\vec k_1,\vec k_2,\vec k_{12};
p_{\pm \infty};E)$ in Definition \ref{def516}
and define $\frak q^{\rm form}_{\star}$ as in (\ref{form269}).
In the case ${\bf g}_1 =  {\bf g}_2 = 0$ it gives a 
structure operation of left $\sL'$, $\Fuk$ right $\sL$ (curved) trimodule structure.  
Here $\sL'$ is a filtered $A_{\infty}$ category associated to $(X,-J)$, $V \oplus TX$ and ${\bf L}'$.
(Here we use $(V \oplus TX) \times V$ as a background class of relative spin structure on $X\times X$
to define orientation of ${\mathcal M}_{\rm QT}(\vec\kappa_1,\vec\kappa_2,\vec\rho;\vec p_1,\vec p_{2},\vec p_{12};\vec k_1,\vec k_2,\vec k_{12};
p_{\pm \infty};E)$ .)
Now we apply Proposition \ref{2636prop} to obtain a left  $-\sL^{\rm op}$, $\Fuk$ right $\sL$ (curved) trimodule.
We then apply Lemma \ref{lem1639} to obtain a left $\sL^{\rm op}$, $\Fuk$ right $\sL$ (curved) trimodule.
We remark $-b$ is a bounding cochain of $\sL^{\rm op}$ via the isomorphism in Lemma \ref{lem1639}.
(This is Corollary \ref{cor2638}.)
The rest of the proof of Theorem \ref{themKuneth} is the same as written in \cite{fukaya:functor}.
\end{proof}
We put a similar sign correction to define $Z_{T,B}^{\rm f.c.}$, $Z_{T}^{\rm f.c.}$.
We can thus 
incorporate
 the sign correction of Subsection \ref{subsec:corsign} to the discussion of the proof of 
 Theorem \ref{thm253}.
 
 \begin{rem}
 In this section we explain the sign correction to the paper \cite{fukaya:cyc} until its Section 7.
 We need similar modification to the other part of \cite{fukaya:cyc}.
 The explanation of such modification is postponed to the documents which will appear  near future.
 Here we only present the part we need for the purpose of this paper, that is, the proof of Theorem \ref{thm253}.
 \end{rem}

\section{Degeneration Hodge-to-de Rham.}
\label{sec:Hodge-to-de Rham}

In this section we prove Theorem \ref{maintheorem3.5}.  %\marginpar{This section is new. KF 2025 July}
In this section we study a filtered $A_{\infty}$ category $\cL$ which is linear over $\Lambda_0$.
We first recall a few definitions.
Let $\cL$ be a unital and cyclic $A_{\infty}$ category linear over $\Lambda_0$ which is curvature free.
Its Hochishild chain complex is 
$$
CH_*(\sL,\sL;\Lambda_0) := \bigoplus_{c,c' \in {\rm Ob}(\sL)} \cL(c,c') \otimes_{\Lambda_0} B^{\rm red}\sL(c',c)
$$
together with the Hochishild differential $\delta_H$.
\begin{defn}
The (negative) cyclic chain complex\index{cyclic complex} $CC(\sL,\sL;\Lambda_0)$ is\index[syindex]{CCStar()@$CC(\sL,\sL;\Lambda_0)$}
$$
CC_*(\sL,\sL;\Lambda_0) = CH_*(\sL,\sL;\Lambda_0) \otimes_{\Lambda_0} \Lambda_0[\![u]\!]
$$
together with differential $\delta_H + uB$, where $u$ is a formal parameter of degree $2$ and $B : CH_*(\sL,\sL;\Lambda_0) \to CH_{*-1}(\sL,\sL;\Lambda_0)$
is defined by\index[syindex]{B@$B$}
\begin{equation}\label{form241}
B(x_0 \otimes {\bf x}) = \sum_{i=0}^k {\bf e}_{c_{i}} \otimes \mathcal S_i([x_0] \otimes {\bf x}).
\end{equation}
Here\index[syindex]{Scali@$\mathcal S_i$} 
$$
\mathcal S_i(x_0 \otimes {\bf x}) = (-1)^{{\maltese}}x_i \otimes x_{i+1} \otimes \dots \otimes x_k \otimes x_0  \otimes \dots \otimes x_{i-1}
$$
with ${\bf x} = x_1 \otimes \dots \otimes x_k$ and $x_i \in \sL(c_{i},c_{i+1})$,
${\maltese} = (\deg'x_i + \dots \deg' x_k)(\deg'x_0 + \dots \deg' x_{i-1})$. 
It is easy to see that $(\delta_H + uB)^2 = 0$. The homology of $\delta_H + uB$
is called the negative cyclic homology\index{cyclic homology}\index{negative cyclic homology} $HC_*(\sL,\sL;\Lambda_0)$ of $\sL$.\index[syindex]{HC@$HC_*(\cdot,\cdot)$}
\par
We put $CC(\sL,\sL) = CC(\sL,\sL;\Lambda_0) \otimes_{\Lambda_0} \Lambda$
\end{defn}
Obviously there exists an exact sequence of chain complexes
\begin{equation}
0 \longrightarrow CC_*(\sL,\sL) \overset{\times u}{\longrightarrow} CC_*(\sL,\sL) \longrightarrow CH_*(\sL,\sL) \longrightarrow 0.
\end{equation}
It induces a long exact sequence:
\begin{equation}
\longrightarrow HC_{*-2}(\sL,\sL) {\longrightarrow} HC_*(\sL,\sL) \longrightarrow HH_*(\sL,\sL) \longrightarrow HC_{*-1}(\sL,\sL) \longrightarrow
\end{equation}
 (\cite{Co,LQ}.)
We also consider the filtration $F^nCC(\sL,\sL) = u^nCC(\sL,\sL)$
to obtain a spectral sequence $E$ with 
$$
E^2_k \cong HH(\sL,\sL),  \quad  E^{\infty}_k =  u^kHC(\sL,\sL)/u^{k+1}HC(\sL,\sL).
$$
which is called Hodge-to-de Rham spectral sequence.\index{Hodge-to-de Rham spectral sequence}
Theorem \ref{maintheorem3.5} is its degeneracy at $E_2$ level.
\par
We consider the case of $\sL = \oplus_{\lambda}\sL_{\lambda}$.
\begin{lem}\label{lem242}
We have $\widehat{\frak p}\circ B = 0$ where $B$ is $(\ref{form241})$.
\end{lem}
\begin{proof}
$
\widehat{\frak p}({\bf e}_{c_i} \otimes \mathcal S_i([x_0] \otimes {\bf x})) =  0,
$
by Theorem \ref{themp} (2)(3).
\end{proof}
We define $\frak p_c : CC_*(\sL,\sL;\Lambda_0) \to \Omega^*(X) \widehat{\otimes}_{\F} \Lambda_0[\![u]\!]$
by $\frak p_c = \frak p \otimes {\rm id}$.
We define a filtration of the codomain by $F^n(\Omega^*(X) \widehat{\otimes}_{\F} \Lambda[\![u]\!])
= \Omega^*(X) \widehat{\otimes}_{\F} u^n\Lambda[\![u]\!]$. \index[syindex]{pfrac@${\frak p}_c$}

\begin{lem}\label{lem243}
${\frak p}_c : CC_*(\sL,\sL) \to \Omega^*(X) \widehat{\otimes}_{\F} \Lambda[\![u]\!]$ is a 
filtration preserving chain map which induces an 
isomorphism on homologies. Here the differential of $\Omega^*(X) \widehat{\otimes}_{\F} \Lambda$ 
is $d  \otimes {\rm id}$.
\end{lem}
\begin{proof}
It is obvious from definition that ${\frak p}_c$ preserves the filtrations.  
Lemma \ref{lem242} and the fact that $\widehat{\frak p}$ is a chain map imply that ${\frak p}_c$
is a chain map.  
Therefore it induces a morphism of the spectral sequences associated to  the filtrations.
To $E_2$ terms, it induces an isomorphism (Theorem \ref{maintheorem3}),
$
{\frak p}_c : HH_*(\cL,\cL;\Lambda) \to H_*(X;\Lambda).
$
Therefore ${\frak p}_c$ induces an isomorphism on homologies.
\end{proof}
Theorem \ref{maintheorem3.5} is an immediate consequence of 
Lemma \ref{lem243}.
\qed

\begin{rem}
It is easy to see that the above proof also shows  Remark \ref{rem119}.  This version can be used to 
study a component $\sL_{\lambda}$ of $\sL$.
\end{rem}
 
We can also prove the next property of $B$.

\begin{lem}
$$
Z(B \text{\bf x},\text{\bf y}) 
+ (-1)^{\deg' \text{\bf x}}Z(\text{\bf x},B \text{\bf y})  = 0.
$$
\end{lem}
\begin{proof}
We have
$$
\aligned
Z(B \text{\bf x},\text{\bf y})  &= \sum_i (-1)^{\maltese} Z({\bf e} \otimes S_i({\bf x}),\text{\bf y}) \\
&= 
\sum_{f,g,i,j} (-1)^{\maltese} \langle\frak m_2(f,{\bf e}),g^{\vee} \rangle_{\rm cyc } \langle g, 
\frak m(\mathcal S_i({\bf x}),f^{\vee},\mathcal S_j({\bf y}))\rangle_{\rm cyc}\\
&= \sum_{f,i,j} (-1)^{\maltese}  \langle f, \frak m(\mathcal S_i({\bf x}),f^{\vee},\mathcal S_j({\bf y}))\rangle_{\rm cyc}.
\endaligned
$$
We calculate $Z(\text{\bf x},B \text{\bf y})$ in a similar way and prove the lemma.
\end{proof}
Therefore $Z$ defines an inner product of $\delta_H + uB$ homology,\index[syindex]{Zhat@$\hat Z$}
that is,
$$
\hat Z : HC_*(\sL;\Lambda_0) \otimes HC_*(\sL;\Lambda_0) \to \Lambda_0[\![u]\!],
$$
which extends $Z$ and satisfies
$
\hat Z(u {\bf x},{\bf y}) = - \hat Z({\bf x},u{\bf y}) = u \hat Z({\bf x},{\bf y}).
$
This pairing is expected to correspond to K. Saito's higher residue pairing (\cite{Sa}) in certain cases.

\section{The toric case.}\label{sec:toric}
\subsection{The strategy of the proof}

In this section, we prove Theorem \ref{maintheorem4}.
We first recall some notations and results from \cite{toric3}.
Let $(X,\omega)$ be a compact toric manifold 
and $\mu : X \to \R^n$ its moment map whose image 
is a polytope $P$, with interior denoted by $\overset{\circ}P $.
We take the coordinates of $\R^n$ so that the composition 
of $\mu$ and the coordinate functions are moment maps of 
$S^1$ actions which are factors of $T^n$.
This choice determines a factorization 
of $T^n$ into the direct product $(S^1)^n$ 
and we denote by $e_1,\dots,e_n$  the basis of 
$H_1(T^n;\Z)$ corresponding to this 
factorization.
It determines a basis of $H_1(L(u);\Z)$
where $L(u) = \mu^{-1}(u)$ is the Lagrangian fibers.
We take $\frak b \in H^{\rm even}(X;\Lambda_0^{\C})$ and 
use it as the bulk class of our 
filtered $A_{\infty}$ category.
Then there exist 
finitely many pairs $(L(u_{\kappa}),b_{\kappa})$
where $u_{\kappa} \in {\overset{\circ}P }$ and 
$b_{\kappa} \in H^1(L(u_{\kappa});\Lambda_0^{\C})/2\pi \sqrt{-1}H^1(L(u_{\kappa});\Z)$
which are regarded as  weak bounding cochains,
such that
\begin{equation}\label{HFisbig}
HF((L(u_{\kappa}),b_{\kappa},\frak b),(L(u_{\kappa}),b_{\kappa},\frak b);\Lambda^{\C}) \cong H(L(u_{\kappa});\Lambda^{\C}).
\end{equation}
The number of such pairs (counted with appropriate multiplicity) is equal to the Betti number of $X$.
(\cite[Theorem 1.1.1]{toric3}.)
We consider the set  
${\bf L} = \{(L(u_{\kappa}),b_{\kappa}) \mid \kappa =1,\dots,K\}$
 with ($K \le {\operatorname{rank}}\, H(X)$) consisting of all such pairs
and will show that it satisfies the conclusion of Theorem \ref{maintheorem4}.
From now on in this section we write $\Lambda_0$, $\Lambda$ etc. 
in place of $\Lambda^{\C}_0$, $\Lambda^{\C}$.
We take $\frak b$ so that (other than $H^0$ part) it is supported in a small 
neighborhood of the toric divisor.  So $i_{L_{u_{\kappa}}}^*(\frak b) = 0$ and 
we take $\theta_{u_{\kappa}} \equiv b_{\kappa} \mod \Lambda_{+}$, which is a closed one form on $L(u_{\kappa})$.
\par
For the proof we study the relationship between Hochschild cohomology and the Jacobian ring.
We take the origin $0 = (0,\dots,0) \in P$ and 
use the above mentioned basis to identify $H^1(L(0);\Lambda_0) \cong (\Lambda_0)^n$. 
We denote by $x_1,\dots,x_n$ the coordinate functions of  $H^1(L(0);\Lambda_0)$ corresponding to this basis.
For $u \in P$ we denote the (corresponding) coordinates of $H^1(L(u);\Lambda_0)$
by $x_1(u),\dots,x_n(u)$. Each element
$b(u) = \sum x_i(u)e_i \in H^1(L(u);\Lambda_0)$ is a weak bounding 
cochain as proved in \cite[Proposition 3.1]{fooo:bulk}, and\index[syindex]{POfra@$\frak{PO}_{\frak b}$} we define the 
{\it potential function}\index{potential function} by:
$$
\frak{PO}_{\frak b}(L(u),b(u)){\bf e} = \frak m^{b(u),\frak b}(1)
= \sum_{k,\ell} \frak q(\frak b^{\otimes\ell};b(u)^{\otimes k})
$$
here $\frak m$, $\frak q$ are the $A_{\infty}$ operation and the closed-open map defined for 
$L(u)$, ${\bf e}$ is the unit in $H^0(L(u);\Lambda_0)$.
\par
It is proved in \cite[Lemma 1.3.7]{toric3} that there exists a function
$$
\frak{PO}_{\frak b} \in \Lambda\langle\!\langle y,y^{-1} \rangle\!\rangle_0^{\overset{\circ}P}
$$
such that
\begin{equation}\label{form141}
\frak{PO}_{\frak b}(L(u),b(u)) = \frak{PO}_{\frak b}(y).
\end{equation}
The explanation of the notations in the above formula is in order.
\par
We define the valuation $v_T(a) = \lambda_1$ when $a = \sum_{i=1}^{\infty} a_i T^{\lambda_i}$
with $a_i \in \C\setminus \{0\}$, $\lambda_i \in \R_{\ge 0}$, $\lambda_i < \lambda_{i+1}$, $\lim_{i\to \infty} \lambda_i = \infty$.
The ring of {\it strictly convergent power series}\index{strictly convergent power series ring}  $\Lambda\langle\!\langle y,y^{-1} \rangle\!\rangle_0^{\overset{\circ}P}$
is the set of formal sums
$$
\sum_{I \in \Z^n} a_I y_1^{i_1}\dots y_n^{i_n}
$$
where $a_I \in \Lambda_0$ such that 
$$
\lim_{\max \{\vert i_j\vert \mid j=1,\dots,n\} \to \infty} \left(v_T(a_I) + \sum_{j=1}^n i_ju_{j}\right)= + \infty,
$$
for any $u = (u_1,\dots,u_n)\in \overset{\circ}P$. See \cite[Definition 1.3.1]{toric3}. In other words 
for any $\overline y_i \in \Lambda$ with $(v_T(\overline y_1),\dots,v_T(\overline y_n)) \in \overset{\circ}P $
we require for the formal sum
$
\sum_{I \in \Z^n} a_I \overline y_1^{i_1}\dots \overline y_n^{i_n}
$
to converge in the $T$ adic topology.
Given $u$, $b(u) = \sum x_i(u)e_i$ 
we put
$
y_i = T^{u_i}e^{x_i(u)}.
$
Then
$$
\frak{PO}_{\frak b}(y) = \frak{PO}_{\frak b}(y_1,\dots,y_n)
$$
is defined and becomes an element of $\Lambda\langle\!\langle y,y^{-1} \rangle\!\rangle_0^{\overset{\circ}P}$.
It is important that the function $\frak{PO}_{\frak b}$ in the right hand side of (\ref{form141}) 
is independent of $u \in \overset{\circ}P$.
\par
The finitely many pairs $(L(u_{\kappa}),b_{\kappa})$ correspond to the 
critical points $\frak y_{\kappa} = (y_{\kappa,1},\dots,y_{\kappa,n})$
of $\frak{PO}_{\frak b}$ by the following formula:
$$
y_{\kappa,i} = T^{u_{\kappa,i}} e^{x(u_{\kappa,i})}.
$$
Using the fact that $\frak y_{\kappa}$ is a critical point we can show 
(\ref{HFisbig}).
\par
We define the {\it Jacobian ring}\index{Jacobian ring} ${\rm Jac}(\frak{PO}_{\frak b})$ of $\frak{PO}_{\frak b}$
as the quotient\index[syindex]{Jac@${\rm Jac}(\frak{PO}_{\frak b})$} of
$\Lambda\langle\!\langle y,y^{-1} \rangle\!\rangle_0^{\overset{\circ}P}$
by the closure of the ideal  generated by $\partial \frak{PO}_{\frak b}/\partial y_i$, $i=1,\dots n$.
(\cite[Definition 1.3.10]{toric3}.)
\par
The {\it Kodaira-Spencer map}\index{Kodaira-Spencer map} \index[syindex]{ksfra@$\frak{ks}_{\frak b}$} 
$$
\frak{ks}_{\frak b} : H^*(X;\Lambda_0) \to {\rm Jac}(\frak{PO}_{\frak b})
$$
defined in \cite[Section 2.4]{toric3} is obtained by taking the derivative of $\frak{PO}_{\frak b}$.
This is a ring homomorphism where the ring structure is the quantum cup product deformed by $\frak b$.
It is proved in \cite[Theorem 1.1.1]{toric3} that $\frak{ks}_{\frak b}$ is an isomorphism.
\par
We next define the formal power series version of the Jacobian ring.
We consider one of the critical points 
$\frak y_{\kappa} = (y_{\kappa,1},\dots,y_{\kappa,n})$.
We consider coordinates $w^{\kappa}_i$ such that $y_i = y_{\kappa,i} e^{w^{\kappa}_i}$
and formal power series ring:
$$
\Lambda_0 [\![ w^{\kappa}_1,\dots,w^{\kappa}_n ]\!],
$$
whose element is a formal sum
\begin{equation}\label{formula1433}
\sum_{I \in \Z_{\ge 0}^n} a_I (w^{\kappa}_1)^{i_1}\dots (w^{\kappa}_n)^{i_n}
\end{equation}
where $a_I \in \Lambda_0$.
The ring
$$
\Lambda_0 [\![ w^{\kappa}_1,\dots,w^{\kappa}_n ]\!] \otimes_{\Lambda_0} \Lambda
$$
is the set of formal sum (\ref{formula1433}) such that $a_I \in \Lambda$ with 
$$
\inf \{ v_T(a_I) \mid I \in \Z_{\ge 0}^n\}  \ne -\infty.
$$
(This ring is smaller than the
formal power series ring $\Lambda [\![ w^{\kappa}_1,\dots,w^{\kappa}_n ]\!]$ over $\Lambda$.)
\par
The strictly convergent power series $\frak{PO}_{\frak b}$ 
defines an element of 
$
\Lambda_0 [\![ w^{\kappa}_1,\dots,w^{\kappa}_n ]\!] \otimes_{\Lambda_0} \Lambda
$ by putting $y_i = y_{\kappa,i} e^{w^{\kappa}_i}$.\index[syindex]{Jacform@${\rm Jac}^{\rm fm}_{\frak y_{\kappa}}(\frak{PO}_{\frak b})$}
 %\marginpar{{Topology is removed. KF}}
\begin{defn}
The {\it formal Jacobian ring},\index{formal Jacobian ring} ${\rm Jac}^{\rm fm}_{\frak y_{\kappa}}(\frak{PO}_{\frak b})$, 
is the quotient of the ring 
$
\Lambda_0 [\![ w^{\kappa}_1,\dots,w^{\kappa}_n ]\!] \otimes_{\Lambda_0} \Lambda
$
by the  ideal generated by $\partial \frak{PO}_{\frak b}/\partial w^{\kappa}_i$, $i=1,\dots n$.
\end{defn}
We remark that here  the above definition does not use the topology of $\Lambda$.
\par
We next define a map
\begin{equation}\label{form144}
{\rm Jac}(\frak{PO}_{\frak b}) \otimes_{\Lambda_0} \Lambda \to 
\prod_{\kappa} {\rm Jac}^{\rm fm}_{\frak y_{\kappa}}(\frak{PO}_{\frak b}).
\end{equation}
We take $P' \subset \overset{\circ}P $, a polytope which contains the 
valuation vectors %$(\in \overset{\circ}P)$
of all the critical points $\frak y_{\kappa}$.
There is an obvious map
$$
\Lambda\langle\!\langle y,y^{-1} \rangle\!\rangle_0^{\overset{\circ}P}
\to
\Lambda\langle\!\langle y,y^{-1} \rangle\!\rangle_0^{P'}.
$$
Note that the target is the ring of (rigid) analytic functions 
which converge on $P'$. We emphasise that here we require the convergence on $P'$, rather than only in its interior.
\par
By \cite[Lemma 2.9.8 and Remark 2.9.9]{toric3}  %\marginpar{The number 2.9.8 to be checked KF 2025 March} 
the ideal generated by $\partial \frak{PO}_{\frak b}/\partial y_i$, $i=1,\dots n$
is closed in $\Lambda\langle\!\langle y,y^{-1} \rangle\!\rangle_0^{P'}$.
Therefore we obtain
$$
{\rm Jac}(\frak{PO}_{\frak b}) \to
\frac{\Lambda\langle\!\langle y,y^{-1} \rangle\!\rangle_0^{P'}}
{\left(\partial \frak{PO}_{\frak b}/\partial y_i:
i=1,\dots n
\right)}.
$$
By associating a formal power series to an element of 
$\Lambda\langle\!\langle y,y^{-1} \rangle\!\rangle_0^{P'}$
we obtain 
$$
\frac{\Lambda\langle\!\langle y,y^{-1} \rangle\!\rangle_0^{P'}}
{\left(\partial \frak{PO}_{\frak b}/\partial y_i:
i=1,\dots n
\right)}
\to 
\prod_{\kappa} {\rm Jac}^{\rm fm}_{\frak y_{\kappa}}(\frak{PO}_{\frak b}).
$$
The composition is the map (\ref{form144}).
 %\marginpar{An argument is added. KF}

In Subsection \ref{subsec:HHtoJacobi} we will construct the next commutative diagram.

\begin{equation}\label{formalizationsquare}
\xymatrix{ 
H(X;\Lambda_0) \ar[rr]^{\widehat{\frak q}}\ar[d]^{\frak{ks}_{\frak b}} &  &   \ar[d]^{}  HH^{*}_{\rm red} ( \cL,\cL) \\
{\rm Jac}(\frak{PO}_{\frak b}) \otimes_{\Lambda_0} \Lambda\ar[rr]^{}  & &  \prod_{\kappa} {\rm Jac}^{\rm fm}_{\frak y_{\kappa}}(\frak{PO}_{\frak b}).
}
\end{equation}
Here the lower horizontal arrow is the map (\ref{form144}). %\marginpar{$HH$ is changed to $HH_{\rm red}$.  KF 2025 March}
\par
In Subsection \ref{subsec:Jringcomp} we prove the next result.

\begin{prop}\label{jaccompthm}
The map $(\ref{form144})$ is an isomorphism.
 %\marginpar{injective is changed to isomorphism.
%We need injective only for the proof of main theorem. But iso looks to be a better result.  KF 2025 July.}
\end{prop}

\begin{proof}[Proposition \ref{jaccompthm} $\Rightarrow$ Theorem \ref{maintheorem4}]
The map 
$\frak{ks}_{\frak b} : H(X;\Lambda_0) \to {\rm Jac}(\frak{PO}_{\frak b})$ is 
an isomorphism. Therefore ${\rm Jac}(\frak{PO}_{\frak b})$ is a free $\Lambda_0$ module.
Hence the left virtical arrow of Diagram (\ref{formalizationsquare}) is injective.
The lower horizontal arrow is injective by Proposition \ref{jaccompthm}.
Therefore $\widehat{\frak q}$ is injective.
The map (\ref{mapphat(1.3)}) is its dual and so is surjective, as required.
\end{proof}

\begin{rem}
In \cite{fooo:bulk, toric3}, the bulk deformation is defined by using equivariant cycles on the 
toric divisor. In this paper we use differential forms to define it.
However they give the same result as is proved in \cite[Proposition 20.7]{spectre}.
See Subsection \ref{perturbcomp}.
\par
Actually even though we use both de Rham model here and in \cite{fooo:bulk, toric3} 
the sign conventions may not  coincide among them.  In \cite{fooo:bulk, toric3}  we used 
the sign in \cite{fooo092} (and transform  it to differential form). 
However it is easy to see that they exactly coincide in case inputs are one form (for forms on $L$)
and two forms (for forms on $X$).  Therefore the calculation of the leading order term of 
the potential function (which we borrow from \cite{cho-oh}) coincides in two convention.
Then the rest of the argument of \cite{fooo:bulk, toric3} uses only the general property of 
$\frak q$ which holds for both sign convention.  Thus we can prove the results 
of \cite{fooo:bulk, toric3}  in our sign convention by the same proof. 
\end{rem}

\begin{rem}
Theorem \ref{maintheorem4} is stated and proved for a toric manifold $X$.
However the same argument can be used to derive the same conclusion 
in various other cases, where we know the number of critical points 
of the potential function coincides with the Betti-number, as implemented e.g. in
\cite{NNU,CKO,YX}. 
\end{rem}

\subsection{Hochschild cohomology and Jacobian ring}\label{subsec:HHtoJacobi}

The purpose of this subsection is to construct  Diagram 
(\ref{formalizationsquare}).
\par
The construction of the right vertical arrow in (\ref{formalizationsquare})
is given in \cite[Lemma 4.7.5]{toric3}. We first repeat it here for completeness.
Let $(u_{\kappa},b_{\kappa})$ be a pair corresponding to a critical point $\frak y_{\kappa}$ 
of $\frak{PO}_{\frak b}$.
We take an element $\vec{\varphi}$ of the Hochschild cochain complex.
It induces a map
$$
{\varphi}_k : H(L(u_{\kappa});\Lambda)^{\otimes k} \to H(L(u_{\kappa});\Lambda).
$$
We  assume that $\vec{\varphi}$ is a {\it reduced} Hochschild cochain, 
that is
$$
{\varphi}_k(x_1,\dots,x_{i-1},{\bf e},x_{i+1},\dots,x_k) = 0. 
$$
For $w = (w_1,\dots,w_n) \in H^1(L(u_{\kappa});\C)$ we put
$$
\mathcal K_{\kappa}(\vec{\varphi})(w) =  \sum_{k=1}^{\infty} 
\int_{L(u_{\kappa})}{\varphi}_k(w,\dots,w).
$$
We remark that the right hand side is well defined as a formal power series.
(It does not necessarily converge in the $T$ adic topology 
since $v_T(w_i)$ is not positive.) 
\par
This defines a map
$$
\mathcal K_{\kappa} : CH^{*}_{\rm red} ( \cL,\cL) \to {\rm Jac}^{\rm fm}_{\frak y_{\kappa}}(\frak{PO}_{\frak b}).
$$
Let $\vec{\psi}$ be a reduced 
cochain, and set $\vec{\varphi} = \delta^H(\vec{\psi})$. 
We will prove that $\mathcal K_{\kappa}(\vec{\varphi})$ is in the 
Jacobian ideal.
By definition we have
\begin{equation}
\aligned
\varphi_k(x_1,\dots,x_k) 
= 
&\sum (-1)^{\maltese_1}\frak m^{\frak b,b_{\kappa}}_{k-j+1}(x_1,\dots,\psi_{j-i+1}(x_i,\dots,x_j),\dots,x_k)
\\
&+ \sum (-1)^{\maltese_2}\psi_{k-j+1}(x_1,\dots,\frak m^{\frak b,b_{\kappa}}_{j-i+1}(x_i,\dots,x_j),\dots,x_k)
\endaligned
\end{equation}
with $\maltese_1 = (\deg' x_1 + \dots + \deg' x_{i-1})\deg \psi$
and $\maltese_2 = \deg' x_1 + \dots + \deg' x_{i-1}$.
Therefore
$$
\aligned
\mathcal K_{\kappa}(\vec{\varphi})(w) = 
&\sum (-1)^{\maltese_1} \int \frak m^{\frak b,b_{\kappa}}_{k-j+1}(w,\dots,\psi_{j-i+1}(w,\dots,w),\dots,w)
\\
&+ \sum (-1)^{\maltese_2}\int  \psi_{k-j+1}(w,\dots,\frak m^{\frak b,b_{\kappa}}_{j-i+1}(w,\dots,w),\dots,w).
\endaligned
$$
The first sum becomes
$$
\int \frak m_1^{\frak b,b_{\kappa} + w}(\vec\psi(e^w)).
$$
As is proved in \cite[Proposition 2.4.16]{toric3} 
this is contained in the formal Jacobian ideal.
\par
Since an element of $H^1(L(u_{\kappa}))$ gives a  
bounding cochain it follows that
$$
\sum_k \frak m_k^{\frak b,b_{\kappa}}(w,\dots,w)
$$
is proportional to the unit.
Therefore, since $\vec{\psi}$ is reduced, the second term vanishes.

We have thus constructed the map, the right vertical arrow in (\ref{formalizationsquare}).
We next prove the commutativity of the diagram.
Let $\Delta_{\frak b} \in H^*(X;\Lambda_0)$.
The upper horizontal arrow (closed-open map) in Diagram (\ref{formalizationsquare})
maps $\Delta_{\frak b}$ to
$\hat{\varphi}^{\Delta_{\frak b}}$ such that
$$
\varphi^{\Delta_{\frak b}}_k(x_1,\dots,x_k) = 
\left.\frac{d}{d\epsilon} \frak q^{\frak b,b_{\kappa}}(e^{\epsilon \Delta_{\frak b} };x_1,\dots,x_k)\right\vert_{\epsilon=0}.
$$
Therefore
$$
\mathcal K_{\kappa}(\vec{\varphi}^{\Delta_{\frak b}})(w) =
\sum_k
\int
\left.\frac{d}{d\epsilon} \frak q^{\frak b,b_{\kappa}}(e^{\epsilon \Delta_{\frak b} };
\underbrace{w,\dots,w}_{k})\right\vert_{\epsilon=0}.
$$
Using \cite[Remark 2.3.19]{toric3} (which is proved in  \cite[Subsection 4.1]{toric3}) we can 
easily show the first equality of: %\marginpar{modified here. KF 2024 Dec}
\begin{equation}\label{eq23777}
\sum_k
\frak q^{\frak b,b_{\kappa}}(e^{\epsilon \Delta_{\frak b} };\underbrace{w,\dots,w}_{k})
=
\frak q^{\frak b_0}(e^{\frak b_+ + \epsilon \Delta_{\frak b}}; e^{b_{\kappa}+ w})
=
\frak q(e^{\frak b + \epsilon \Delta_{\frak b}}; e^{b_{\kappa}+ w}).
\end{equation}
For the second equality, there is an issue mentioned in Subsection \ref{divisor}.  However
this equality can be proved in this situation as we explain in Subsection \ref{perturbcomp}.

Therefore by definition
$$
\mathcal K_{\kappa}(\vec{\varphi}^{\Delta_{\frak b}})(w) 
= \left.\frac{d}{d\epsilon}\frak{PO}(\frak b + \epsilon 
\Delta_{\frak b};\frak y_{\kappa,1}e^{w_1},\dots,\frak y_{\kappa,n}e^{w_n})\right\vert_{\epsilon=0}.
$$

By definition
$$
\frak{ks}_{\frak b}(y) 
= \left.\frac{d}{d\epsilon}\frak{PO}(\frak b + \epsilon 
\Delta_{\frak b};y)\right\vert_{\epsilon=0}.
$$
Therefore the lower horizontal arrow sends it to 
$\mathcal K_{\kappa}(\vec{\varphi}^{\Delta_{\frak b}})$.
We have thus proved the commutativity of Diagram (\ref{formalizationsquare}).

\subsection{Comparison between two Jacobian rings 1}\label{subsec:Jringcomp}

In this subsection we prove Proposition \ref{jaccompthm}. This is a 
``rigid analytic analogue'' of a standard result, which we first explain.
\par
Let $M$ be a Stein manifold and $W$  a holomorphic function on it.
We assume $W$ has a finitely many isolated critical points, 
which we denote by $\{ z_{\kappa} \}$.
Let $\frak O$ be the sheaf of holomorphic functions on $M$.
We consider the sheaf $\frak J_W$ of ideals generated by 
the first derivatives of $W$.
We put
$$
{\rm Jac}(W) = \frac{\frak O(M)}{\frak J_W(M)}.
$$
Next we consider for each $\kappa$ the ring of formal power series 
at $z_{\kappa}$ and consider its quotient by the ideal generated by the 
first derivatives of $W$ at $z_{\kappa}$.
We denote this quotient by ${\rm Jac}^{\rm fm}_{z_{\kappa}}(W)$.
We consider the obvious map
\begin{equation}\label{map1466}
{\rm Jac}(W) \to \prod_{\kappa} {\rm Jac}^{\rm fm}_{z_{\kappa}}(W).
\end{equation}
This map is an analogue of (\ref{form144}) and is an isomorphism.
We can prove this fact (which is standard) as follows.
We first consider the quotient of the convergent power series ring 
at $z_{\kappa}$ by the ideal generated by the 
first derivatives of $W$. We denote it by ${\rm Jac}^{\rm conv}_{z_{\kappa}}(W)$.
Using the fact that $z_{\kappa}$ is an isolated critical point, 
we can easily show that 
${\rm Jac}^{\rm fm}_{z_{\kappa}}(W)$ is isomorphic to ${\rm Jac}^{\rm conv}_{z_{\kappa}}(W)$.
(We will prove   it later as Lemma \ref{lem14444}.)
\par
We next consider the 
exact sequence of coherent sheaves
$$
0 \longrightarrow \frak J_W \longrightarrow \frak O \longrightarrow \frac{\frak O}{\frak J_W}
\longrightarrow 0.
$$
Then the cohomology exact sequence becomes
$$
H^0(\frak J_W) \to H^0(\frak O) \to H^0\left(\frac{\frak O}{\frak J_W}\right) \to H^1(\frak J_W) \to
$$
Note $H^0(\frak J_W) = \frak J_W(M)$, $H^0(\frak O) = \frak O(M)$
and $H^0(\frac{\frak O}{\frak J_W}) = \prod_{\kappa} {\rm Jac}^{\rm conv}_{z_{\kappa}}(W)$.
On the other hand $H^1(\frak J_W) = 0$ since $M$ is Stein. Thus (\ref{map1466}) is an isomorphism.
 %\marginpar{
%Injective is changed to isomorphism here in some other places in this subsection.  to be checked.  KF 2025 July.}
\par
It is very likely that we can prove Proposition \ref{jaccompthm} 
in a similar way using sheaf theory in rigid 
analytic geometry. 
Below we provide a different proof without using sheaf theory over rigid analytic space.
\par\medskip
First we recall the decomposition:
\begin{equation}
{\rm Jac}(\frak{PO}_{\frak b}) \otimes_{\Lambda_0} \Lambda
\cong \prod_{\kappa}  {\rm Jac}(\frak{PO}_{\frak b};\frak y_{\kappa})
\end{equation}
where
$$
\begin{aligned}
&{\rm Jac}(\frak{PO}_{\frak b};\frak y_{\kappa}) \\
&= 
\{[\mathcal P] \in {\rm Jac}(\frak{PO}_{\frak b}) \otimes_{\Lambda_0} \Lambda \mid 
(y_i - y_{\kappa,i})^N [\mathcal P] = 0 \,\,\,\, \text{for large $N$ and any $i=1,\dots,n$.}\}
\end{aligned}
$$
(See \cite[Proposition 1.3.16]{toric3}.) and
$\frak y_{\kappa} = (y_{\kappa,1},\dots,y_{\kappa,n})$.
Therefore it suffices to show that
$$
F_{\kappa} : {\rm Jac}(\frak{PO}_{\frak b};\frak y_{\kappa}) \to {\rm Jac}^{\rm fm}_{\frak y_{\kappa}}(\frak{PO}_{\frak b})
$$
is an isomorphism.
We take $\kappa$ and fix it.
By changing $P$ to $P -  u_{\kappa}$ and $y_i$ to $T^{-u_{\kappa,i}}y_i$ 
we may assume $u_{\kappa}  = 0$. (In other words $v_T(y_{\kappa,i}) = 0$.)
We define $y'_i$ by 
$y'_i = y_i - y_{\kappa,i}$.
Then $y_i^{-1} = (y_{\kappa,i} + y'_i)^{-1} \in \Lambda_0[\![y'_i ]\!]$.
Therefore we may regard 
$\frak{PO}_{\frak b} \in \Lambda_0[\![y'_1,\dots,y'_n ]\!] \otimes_{\Lambda_0} \Lambda$.
We put
$$
{\rm Jac}^{\rm fm,\prime}_{\frak y_{\kappa}}(\frak{PO}_{\frak b})
= 
\frac{\Lambda_0[\![y'_1,\dots,y'_n ]\!] \otimes_{\Lambda_0} \Lambda}{\frak I}
$$
where $\frak I$ is an ideal generated by the first derivatives of $\frak{PO}_{\frak b}$.
Since $y'_i = e^{w_i} - 1$ induces an isomorphism of formal power
series rings, we have 
$
{\rm Jac}^{\rm fm}_{\frak y_{\kappa}}(\frak{PO}_{\frak b})
\cong
{\rm Jac}^{\rm fm,\prime}_{\frak y_{\kappa}}(\frak{PO}_{\frak b}).
$
Thus it suffices to prove that
$$
F'_{\kappa} : {\rm Jac}(\frak{PO}_{\frak b};y_{\kappa}) \to {\rm Jac}^{\rm fm,\prime}_{\frak y_{\kappa}}(\frak{PO}_{\frak b})
$$
is an isomorphism.
\par
As is proved in \cite[Theorem 2.9.2]{toric3} there is a coordinate change (which 
converges in $\Lambda\langle\!\langle y,y^{-1} \rangle\!\rangle_0^{\overset{\circ}P}$)
which transforms $\frak{PO}_{\frak b}$ to a Laurent polynomial $W$. 
It suffices to study
$$
F'_{\kappa} : {\rm Jac}(W;\frak y_{\kappa}) \to {\rm Jac}^{\rm fm,\prime}_{\frak y_{\kappa}}(W).
$$
We consider the map
$$
\frac{\Lambda[y_1,\dots,y_n,y^{-1}_1,\dots,y_n^{-1}]}{(\partial W/\partial y_1,\dots,\partial W/\partial y_n)}
\to 
{\rm Jac}(W;\frak y_{\kappa}),
$$
where the domain is the quotient of Laurent {\it polynomial} ring.
%Note $\frak y_{\kappa}$ is a zero of $\partial W/\partial y_i$ $i=1,\dots,n$.
This map is surjective.
Therefore there exists an ideal $\alpha$ of  $\Lambda[y_1,\dots,y_n,y^{-1}_1,\dots,y_n^{-1}]$
containing the ideal $(\partial W/\partial y_1,\dots,\partial W/\partial y_n)$
such that
$$
{\rm Jac}(W;\frak y_{\kappa})
\cong
\frac{\Lambda[y_1,\dots,y_n,y^{-1}_1,\dots,y_n^{-1}]}
{\alpha}.
$$
Note that $\frak y_{\kappa}$ is the unique geometric point of
$Spec(\frac{\Lambda[y_1,\dots,y_n,y^{-1}_1,\dots,y_n^{-1}]}
{\alpha})$.
\par
Thus it suffices to show that
$$
\frac{\Lambda[y_1,\dots,y_n,y^{-1}_1,\dots,y_n^{-1}]}{\alpha}
\to 
{\rm Jac}^{\rm fm,\prime}_{\frak y_{\kappa}}(W)
$$
is an isomorphism.
The fact that this map (which sends $y_i$ to $y_i$) is 
well-defined implies that the ideal $\alpha$ is contained in the ideal 
$\frak I$ generated by the first derivatives of $\frak{PO}_{\frak b}$
in the formal power series ring $\Lambda_0[\![y'_1,\dots,y'_n ]\!] \otimes_{\Lambda_0} \Lambda$.
\par
Note the $T$-adic topology no longer plays a role here.
Therefore we may replace $\Lambda$ by a field $K \subset \Lambda$ 
which is the algebraic closure of the extension of $\Q$ 
by all the elements appearing as the coefficients of 
$W$.
$K$ is an algebraically closed field of characteristic $0$ and of 
countable transcendental degree.
Therefore there exists an embedding $K \to \C$.
Thus it suffices to show that
$$
\frac{\C[y_1,\dots,y_n,y^{-1}_1,\dots,y_n^{-1}]}{\alpha}
\to 
{\rm Jac}^{\rm fm,\prime}_{\frak y_{\kappa}}(W)
$$
is an isomorphism.
Here $W$ is a Laurent polynomial with complex coefficients 
and $\frak y_{\kappa}$ is a critical point in $\C^n$.
\par
Let $O$ be the ring of holomorphic functions on $(\C^*)^n$ and 
$\frak A$ its ideal generated by $\alpha$.
By GAGA we have
$$
\frac{\C[y_1,\dots,y_n,y^{-1}_1,\dots,y_n^{-1}]}{\alpha}
\cong
\frac{O}{\frak A}.
$$
On the other hand, applying the argument of the beginning of this 
subsection to $M = (\C^*)^n$ we find that
$$
\frac{O}{\frak A} \to {\rm Jac}^{\rm fm,\prime}_{\frak y_{\kappa}}(W)
$$
is an isomorphism. 
Namely we regard $\frak A$ as coherent analytic sheaf and 
use the cohomology exact sequence associated to the 
sheaf exact sequence 
$$
0 \longrightarrow \frak A \longrightarrow \frak O \longrightarrow \frac{\frak O}{\frak A}
\longrightarrow 0.
$$

Thus to complete the proof of Proposition \ref{jaccompthm}
it remains to prove the next lemma.

\begin{lem}\label{lem14444}
Let $W$ be a holomorphic function defined on a neighborhood of $0$ in $\C^n$, with an isolated critical point at $0$.
Let $\C\{z_1,\dots,z_n\}$ be the ring of convergent power series 
and $\C[\![z_1,\dots,z_n]\!]$ the ring of formal power series.
Let $\frak I$ be the ideal of $\C[\![z_1,\dots,z_n]\!]$ generated by $\partial W/\partial z_i$, 
$i=1,\dots,n$
and 
 $\frak A$ an ideal of $\C\{z_1,\dots,z_n\}$ 
 containing $\partial W/\partial z_i$, $i=1,\dots,n$.
 We assume $\frak A \subseteq \frak I$.
Then
$$
\frac{\C\{z_1,\dots,z_n\}}{\frak A} \to \frac{\C[\![z_1,\dots,z_n]\!]}{\frak I}
$$
is an isomorphism.
\end{lem}
\begin{proof}
We provide a proof of this standard fact for completeness. 
Since $0$ is an isolated critical point, the quotient
$\frac{\C\{z_1,\dots,z_n\}}{\frak A}$ is finite dimensional and $z_i$ is 
contained in its unique maximal ideal.
Therefore there exists $N$ such that $z_i^N \in \frak A$ for all $i=1,\dots,n$.
The surjectivity follows. Let us prove the injectivity.
Suppose $f \in \C\{z_1,\dots,z_n\}$ such that $[f]$ goes to zero.
Then there exist formal power series $Q_i$ such that 
$$
f = \sum_{i=1}^{n} Q_i \frac{\partial W}{\partial z_i}.
$$
We write $Q_i = Q_i^0 + Q_i'$ such that 
$Q_i^0$ is a polynomial and all the terms of $Q_i'$ contain 
$z_i^m$ with $m\ge N$ for some $i$.
We may replace $f$ by 
$f - \sum_{i=1}^{n} Q^0_i \frac{\partial W}{\partial z_i}$.
Thus $f$ is of the form
$$
f = \sum_I a_I z^I
$$
such that the sum is taken over all $I =(i_1,\dots,i_n)$ 
with $\max \{i_1,\dots,i_n\} \ge N$.
It suffices to show that all such convergent power series 
$f$ are contained in $\frak A$.
We prove the next statement by induction on $m$.
\par
\begin{enumerate}
\item[(*)]
Suppose in addition $f$ does not contain $z_m, \dots, z_n$
then $f \in \frak A$
\end{enumerate}
\par
When $m=1$ this is obvious.
Suppose we proved (*) for $m$.
Let $f$ be a convergent power series over all $I =(i_1,\dots,i_n)$ 
with $\max \{i_1,\dots,i_n\} \ge N$ such that $f$ does not contain $z_{m+1}, \dots, z_n$.
We obtain
$$
f = f_N z_m^N + f_{N-1} z_m^{N-1} + \dots + f_1 z_m + f_0,
$$
where $f_i$ are all convergent power series and 
$f_0,\dots,f_{N-1}$ do not contain $z_m,\dots,z_{n}$.
In fact 
$$
f_i = \frac{1}{i!}\left.\frac{\partial^i f}{\partial z_m^i}\right\vert_{z_m=0}
$$
for $i=0,\dots,N-1$ and 
$$
f_N = \frac{f - f_{N-1} z_m^{N-1} - \dots - f_1 z_m - f_0}{z_m^N}.
$$
Since $f_N$ is a convergent power series  $f_N z_m^N \in \frak A$.
On the other hand, by assumption, each term of $f_i$, $i<N$ must 
contain $z_j^k$ with $j < m$, $k \ge N$.
Therefore by induction hypothesis $f_i \in \frak A$ for $i<N$.
Thus $f \in \frak A$ as required.
\end{proof}

\subsection{Comparison between two Jacobian rings 2}\label{subsec:Jringcomp2}

In this subsection, we give an alternative proof of 
Proposition \ref{jaccompthm}.
Let $\frak y_k$ be a critical point of $\frak{PO}_{\frak b}$.
We first remark that there exists 
a finitely generated ideal $\alpha$ of  ${\rm Jac}(\frak{PO}_{\frak b}) \otimes_{\Lambda_0} \Lambda$
such that
$$
{\rm Jac}(\frak{PO}_{\frak b};\frak y_{\kappa})
=
\frac{{\rm Jac}(\frak{PO}_{\frak b}) \otimes_{\Lambda_0} \Lambda}{\alpha}.
$$
There exists $N$ such that $(y_i - y_{\kappa,i})^N = 0$
in ${\rm Jac}(\frak{PO}_{\frak b};\frak y_{\kappa})$.
On the other hand for each $\kappa'\ne \kappa$ there exists 
$i \in \{1,\dots,n\}$ such that  $y_i - y_{\kappa,i}$
is invertible in ${\rm Jac}(\frak{PO}_{\frak b};\frak y_{\kappa'})$.
Thus $\alpha$ contains the first derivatives of $\frak{PO}_{\frak b}$
and $(y_i - y_{\kappa,i})^N$, $i=1,\dots,n$.
\par
By a coordinate change we may assume that $v_T(y_{\kappa,i}) = 0$, and moreover that $y_{\kappa,i} \in \C$.
We put $y'_i = y_i - y_{\kappa,i}$.

\begin{lem}\label{lem145}
Let $f \in \Lambda\langle\!\langle y,y^{-1} \rangle\!\rangle_0^{\overset{\circ}P} \otimes_{\Lambda_0} \Lambda$
and assume that there exists $g \in \Lambda[\![ y'_1,\dots,y'_n]\!]$ with 
$$
f = y'_i g
$$
in $\Lambda[\![ y'_1,\dots,y'_n]\!]$. Then 
$g \in  \Lambda\langle\!\langle y,y^{-1} \rangle\!\rangle_0^{\overset{\circ}P} \otimes_{\Lambda_0} \Lambda$.
\end{lem}
\begin{proof}
We write
$$
f = \sum T^{\lambda_k} P_k(y_1,\dots,y_n).
$$
Here $P_k \in \C[y_1,\dots,y_n,y^{-1}_1,\dots,y^{-1}_n]$ are 
Laurent {\it polynomials} with complex coefficients and $\lambda_k \in \R$, 
is an increasing sequence such that $\lim_{k\to \infty} \lambda_k = + \infty$.
Put:
\begin{equation}
g = \sum T^{\lambda_k} \frac{P_k(y_1,\dots,y_n)}{y_i - y_{\kappa,i}}.
\end{equation}
We first claim
$$
 \frac{P_k(y_1,\dots,y_n)}{y_i - y_{\kappa,i}} \in \C[y_1,\dots,y_n,y^{-1}_1,\dots,y^{-1}_n].
$$
In fact the left hand side is a meromorphic function on $\C P^n$ which is holomorphic 
on $(\C^*)^n \setminus \{(y_1,\dots,y_n) \mid y_i = y_{\kappa,i} \}$.
Since the left hand side can be identified with a formal power series of $y'_1,\dots,y'_n$ by assumption, it is holomorphic 
at $(y'_1,\dots,y'_n) = (0,\dots,0)$. Therefore it is a holomorphic function of $y_i$ on $(\C^*)^n$. Moreover 
it is a meromorphic function on $\C P^n$. Therefore 
it is a Laurent polynomial.
\par
We define $Q_k = \frac{P_k(y_1,\dots,y_n)}{y_i - y_{\kappa,i}}$ a Laurent polynomial 
of $y_1,\dots,y_n$.
To complete the proof it suffices to show that 
\begin{equation}\label{form14999}
\sum_{k=1}^{\infty} T^{\lambda_k} Q_k(\frak y_1,\dots,\frak y_n)
\end{equation}
converges in the $T$-adic topology for any $\frak y_1,\dots,\frak y_n$ 
with $(v_T(\frak y_1),\dots,v_T(\frak y_n)) \in \overset{\circ}P $. 
\par
By assumption
$$
\lim_{k\to \infty} (\lambda_k + P_k(\frak y_1,\dots,\frak y_n))
= 
+\infty
$$
for $(v_T(\frak y_1),\dots,v_T(\frak y_n)) \in \overset{\circ}P $.
On the other hand
$$
v_T(Q_k(\frak y_1,\dots,\frak y_n))
+ v_T(\frak y_i-y_{\kappa,i}) = v_T(P_k(\frak y_1,\dots,\frak y_n)).
$$
If $v_T(\frak y_i) \ne 0$ then 
$$
v_T(\frak y_i-y_{\kappa,i})
=
\min \{ v_T(\frak y_i), 0\}
$$ 
and is uniformly bounded in absolute value.
Therefore (\ref{form14999}) converges.
\par
We can use the fact that the convex hull of 
$$
\overset{\circ}P \setminus \{(u_1,\dots,u_n) \mid u_i = 0\}
$$
is $\overset{\circ}P $ and prove that (\ref{form14999}) converges 
for $(v_T(\frak y_1),\dots,v_T(\frak y_n)) \in \overset{\circ}P $.
\end{proof}
Now we prove that
$$
\frak F : \frac{{\rm Jac}(\frak{PO}_{\frak b}) \otimes_{\Lambda_0} \Lambda}{\alpha}
\to {\rm Jac}^{\rm fm}_{\frak y_{\kappa}}(\frak{PO}_{\frak b})
$$
is an isomorphism.
\par
Let $f \in {\rm Jac}(\frak{PO}_{\frak b}) \otimes_{\Lambda_0} \Lambda$
be such that $\frak F(f) =0$ in 
${\rm Jac}^{\rm fm,\prime}_{\frak y_{\kappa}}(\frak{PO}_{\frak b})$.
We will prove that $f \in \alpha$.
We have
$$
f = \sum_i Q_i \frac{\partial \frak{PO}_{\frak b}}{\partial y'_i}
$$
in $\Lambda[\![ y'_1,\dots,y'_n]\!]$.
There exists $N$ such that $(y'_i)^N \in \alpha$ for all $i=1,\dots,n$.
\par
We write $Q_i = Q_i^0 + Q_i'$ such that 
$Q_i^0$ is a polynomial of $y'_j$, $j=1,\dots,n$ and all the terms of $Q_i'$ contain some
$(y'_j)^m$ with $m\ge N$ and $j=1,\dots,n$.
We may replace $f$ by 
$f - \sum_{i=1}^{n} Q^0_i \frac{\frak{PO}_{\frak b}}{\partial y'_i}$.
Thus $f$ is of the form
$$
f = \sum_I a_I (y')^I
$$
such that the sum is taken over all $I =(i_1,\dots,i_n)$ 
with $\max \{i_1,\dots,i_n\} \ge N$.
It suffices to show that all such $f$ are contained in $\alpha$.
We prove the next statement by induction on $m$.
\par
\begin{enumerate}
\item[(*)]
Suppose in addition that $f$ does not contain $y'_m, \dots, y'_n$
then $f \in \alpha$
\end{enumerate}
\par
When $m=1$ this is obvious.
Suppose we proved (*) for $m$.
Let $f$ be a convergent power series over all $I =(i_1,\dots,i_n)$ 
with $\max \{i_1,\dots,i_n\} \ge N$ such that $f$ does not contain $y'_{m+1}, \dots, y'_n$.
We obtain
$$
f = f_N (y'_m)^N + f_{N-1} (y'_m)^{N-1} + \dots + f_1 y'_m + f_0,
$$
where $f_i \in {\rm Jac}(\frak{PO}_{\frak b}) \otimes_{\Lambda_0} \Lambda$ and 
$f_0,\dots,f_{N-1}$ do not contain $y'_m,\dots,y'_n$.
In fact 
$$
f_i = \frac{1}{i!}\left.\frac{\partial^i f}{\partial (y'_m)^i}\right\vert_{y'_m=0}
$$
for $i=0,\dots,N-1$ and 
$$
f_N = \frac{f - f_{N-1} (y'_m)^{N-1} - \dots - f_1 (y'_m) - f_0}{(y'_m)^N}.
$$
Here we use Lemma \ref{lem145}.
Since $f_N \in {\rm Jac}(\frak{PO}_{\frak b}) \otimes_{\Lambda_0} \Lambda$,  $f_N (y'_m)^N \in \alpha$.
On the other hand, by assumption, each term of $f_i$, $i<N$ must 
contain $(y'_i)^m$ with $i\ge m+1$.
Therefore by induction hypothesis $f_i \in \alpha$ for $i<N$.
Thus $f \in \alpha$ as required.
We thus proved injectivity.
The proof of surjectivity is easier and so is omitted. %\marginpar{added to be checked KF 2025 July}
\qed

\subsection{Comparison of the perturbations of this paper and of \cite{toric3}.}
\label{perturbcomp}

In this paper we construct $A_{\infty}$ categories, closed-open maps, and open-closed maps 
using Kuranishi structures and CF-perturbations of the moduli spaces. %\marginpar{Subsection added.  KF 2024 Dec.}
In this section we use the result of \cite[Theorem 1.1.1]{toric3}, that is, the Kodaira-Spencer map from the 
quantum cohomology of the ambient space to the Jacobian ring 
is a ring isomorphism.   We use also Kuranishi structures and CF-perturbations in \cite{toric3}.
However there are small differences between two constructions, which we discuss in this subsection.

First we use the divisor axiom (\ref{divbundary}) in \cite{toric3}.
We can safely do so since we need it only in the case of finitely many disjoint Lagrangian tori 
which are orbits of $T^n$ action.
To apply our main theorems Theorems \ref{maintheorem1} and  \ref{maintheorem2} we need 
to consider the case when an `unknown' Lagrangian submanifold $U$ is involved. 
So we need to study the case when intersection of $U$ with an orbit is non-trivial.
However the construction of the perturbations 
for the finite set of Lagrangians, that is, the union of finitely many orbits and $U$, 
is constructed by extending the existing perturbations for the given finitely many $T^n$ orbits. 
Therefore, 
the divisor axiom (\ref{divbundary}) is satisfied as far as finitely many $T^n$ orbits we start with concerns.

Another point concerns the way to how to involve bulk deformations.
In this paper we choose differential forms which represent bulk classes.
In \cite{toric3} we take singular cycles which are union of irreducible components of 
toric divisors or their intersections, and cut the moduli space by requiring 
the evaluation maps to hit such cycles.  (See  \cite[Section 2.3]{toric3}.)
We explain that we can still apply the results (\cite[Theorem 1.1.1]{toric3}) 
in our situation to prove (\ref{eq23777}).
Here we  compare two potential functions 
(associated to finitely many $T^n$ orbits). 
Namely:
\begin{equation}\label{form2311}
\frak{PO}_{\frak b}(b) = \langle \frak m^{\frak b,b}_0(1), {\rm vol}_L \rangle_{\rm cyc}  \in \Lambda_0
\end{equation}
for the operations $\frak m_0^{\frak b,b}$ obtained by the two (different) perturbations.
Note that (\ref{form2311}) is an element of $\Lambda_0$ which is well 
defined in the homology level (once a system of perturbations is given).
Therefore it suffices to compare them modulo $T^E$ for arbitrary but fixed $E$.
To be  a bit more detailed,
we first observe that the contribution of an element  $\frak b_2$ of $H^2(X;\R)$
is the same for the two perturbations.
Note there is a non trivial relation between cohomology classes 
represented by irreducible components of 
toric divisors. We fix finitely many irreducible components which 
generate the
subspace isomorphic to $H^2(X;\R)$. (The proof of the independence of such 
choices is a part of \cite[Sections 2.4 and 2.5]{toric3}.)
Then we take smooth closed 2 forms representing a basis of $H^2(X;\R)$
so that they are smoothing of the currents corresponding those 
irreducible components.  
In particular their supports are disjoint from the Lagrangian torus involved.
By this choice 
it is easy to see the contribution is by the cap product $\beta \cap \frak b_2$
in both cases.

Now it suffices to consider  finitely many moduli spaces to compare
(\ref{form2311}) for two different choices of perturbations.
The argument leading to a proof of such coincidence is 
the same as \cite[Section 28.3]{spectre}.  
\par
The third point we need to discuss is that the curved $A_{\infty}$ categories 
(bulk deformed by the same class $\frak b$) are homotopy equivalent 
when we use two different perturbations.
Namely the part which does not contain $T$ and of degree $0$ contribute 
in the same way.
Then we can use \cite[Section 28.3]{spectre} to compare the bulk deformation 
using irreducible components of toric divisors or their intersections and one by using
their smoothing (degree two forms). 

We remark that using the perturbations of  \cite{spectre} 
it is proved that any cohomology class  $b_1 \in H^1(L;\Lambda_0)$
gives rise to a  bounding cochain
\cite[Proposition 4.3]{fooo:toric1}.
If we use the bulk deformation by differential forms (the ones we use here),
this may not be the case.\footnote{In the Fano case it is still true.}
In fact the  curved $A_{\infty}$ algebras   $(H^*(L;\Lambda_0),\frak m^{\frak b})$
obtained by these two different types of bulk deformations  are homotopy equivalent.
However the homotopy equivalence may not preserve the $\Z$ grading.
Therefore the set of bounding cochains may not coincide with  
the cohomology group $H^1(L;\Lambda_0)$ 
under the bulk deformation we use in this paper, 
although it is identified with $H^1(L;\Lambda_0)$ for the toric case.
Neverthless there is still an embedding 
$H^1(L;\Lambda_0) \hookrightarrow \mathcal M_{\rm weak}(L)$ since $\mathcal M_{\rm weak}(L)$ 
is invariant under the homotopy equivalence.
We can therefore use the result of \cite{fooo:toric1}.

\section{Quantum cup product between $\frak p$ images.}
\label{sec:Theorem2+1}
Let $\bL$ and $\bU$ be collections of Lagrangian submanifolds
 equipped with brane data, which respectively have potential values $\lambda_{\bL}$ and $ \lambda_{\bU} $.
(We assume Assumption \ref{LUtrans}.)\index[syindex]{lambd0aLbf@$\lambda_{\bL}$}
Let $\text{\bf x} \in HH_*({\cL})$, 
$\text{\bf y} \in HH_*({\sU})$.
Theorem \ref{maintheorem1} which we proved in Section  \ref{sec:Theorem1} 
 implies the equality
\begin{equation}\label{pnormal}
\langle \widehat{\frak p}(\text{\bf x}),\hat{\frak p}(\text{\bf y}) \rangle_{\text{\rm PD}_X}
=0
\end{equation}
if $\lambda_{\bL} \ne \lambda_{\bU} $.
\par
Here and hereafter in this and the  next sections we write 
$\hat{\frak p}$, $\hat{\frak q}$, ${\frak p}$ in place of 
$\hat{\frak p}^{\text{\bf b}}$, $\hat{\frak q}^{\text{\bf b}}$, ${\frak p}^{\text{\bf b}}$
for simplicity.
\par
We also have the following:
\begin{prop}
Let $\lambda_{\bL} = \lambda_{\bU}$ and 
assume that
\begin{equation}\label{noFLoer}
HF((L_{\kappa},b_{\kappa}),(U_{\upsilon},b_{{\upsilon}});\Lambda) 
= 0
\end{equation}
for any $\kappa$ and ${\upsilon}$.
Then we have
\begin{equation}\label{pnormal2}
\langle \hat{\frak p}(\text{\bf x}),\hat{\frak p}(\text{\bf y}) \rangle_{\text{\rm PD}_X}
=0
\end{equation}
for all $\text{\bf x} \in HH_*({\cL})$, 
$\text{\bf y} \in HH_*({\sU})$.
\end{prop}
\begin{proof}
Note
Theorem \ref{thm:commutative_diagram_p_q}
holds in our situation also.
By assumption, ${_{\sU}}\Fuk _{\cL} \displaystyle{ \tensor_{\cL}} {_\cL} \Fuk_{\sU}$ 
is acyclic. 
\end{proof}
In this section we improve this proposition as follows.
\begin{thm}\label{thmpcup0}
Assume that either {\rm (i)} $\lambda_{\bL} \ne \lambda_{\bU} $ 
or {\rm (ii)} $\lambda_{\bL} = \lambda_{\bU} $ and $(\ref{noFLoer})$.
Then we have:
\begin{equation}\label{pcup0}
\hat{\frak p}(\text{\bf x}) \cup^{\frak b} \hat{\frak p}(\text{\bf y})=0
\end{equation}
for all $\text{\bf x} \in HH_*({\cL})$, 
$\text{\bf y} \in HH_*({\sU})$.
\end{thm}
% the following:
%\begin{thm}\label{pqrelation}
%Let $\text{\bf x} \in HH_*({\cL})$ 
%and $\frak z \in QH_{\frak b}^*(X;\Lambda)$. Then
%\begin{equation}
%\sum_c (-1)^{\maltese}\hat{\frak p}(\frak m(\text{\bf x}_c^{(H;4;1)}
%\otimes\frak q(\frak z;\text{\bf x}_c^{(H;4;2)})\otimes \text{\bf x}_c^{(H;4;3)})
%\otimes \text{\bf x}_c^{(H;4;4)})
%=
%\frak z\cup^{\frak b} \hat{\frak p}(\text{\bf x}) .
%\end{equation}
%where
%$
%\maltese = ????$.
%Here we use notation $(\ref{DeltaH})$.
%\end{thm}
%\begin{rem}
%In the monotone case somewhat similar result is 
%obtained in Biran-Cornea \cite{bircor}.
%\end{rem}
\begin{proof}[Proof of Theorem
% \ref{pqrelation} + (\ref{pnormal}),(\ref{pnormal2}) $\Rightarrow$ Theorem 
 \ref{thmpcup0}]
Let $\frak z \in H(X;\Lambda)$. Then, by Theorem \ref{pqrelation}, we have
$$
\aligned
&\langle 
\frak z, \hat{\frak p}(\text{\bf x}) \cup^{\frak b} \hat{\frak p}(\text{\bf y})
\rangle_{\text{\rm PD}_X} 
=
\langle \frak z \cup^{\frak b} \hat{\frak p}(\text{\bf x}), \hat{\frak p}(\text{\bf y})
\rangle_{\text{\rm PD}_X} \\
&= 
\sum_c (-1)^{\maltese} \langle \hat{\frak p}(\frak m(\text{\bf x}_c^{(H;4;1)}
\otimes\frak q(\frak z;\text{\bf x}_c^{(H;4;2)})\otimes \text{\bf x}_c^{(H;4;3)})
\otimes \text{\bf x}_c^{(H;4;4)}),   \hat{\frak p}(\text{\bf y})
\rangle_{\text{\rm PD}_X} \\
\endaligned$$
where $\maltese$ is the Koszule sign.
This is zero by (\ref{pnormal}) or (\ref{pnormal2}).
\end{proof}
\begin{rem}
Equations (\ref{pnormal}) and (\ref{pnormal2}) follow from (\ref{pcup0}). In fact
$$
\langle \hat{\frak p}(\text{\bf x}),\hat{\frak p}(\text{\bf y}) \rangle_{\text{\rm PD}_X}
=
%(-1)^{\maltese}
\int_X \hat{\frak p}(\text{\bf x}) \cup^{\frak b} \hat{\frak p}(\text{\bf y}).
$$  
%$\maltese=???.$
\end{rem}
% Theorem \ref{thmpcup0} follows from (\ref{pnormal}), 
% (\ref{pnormal2}) and Theorem \ref{pqrelation}

\begin{rem}
The deformed cup product $\cup^{\frak b}$ is 
defined by using the moduli space of pseudo-holomorphic spheres.
In a similar way, we can define the $Comm_{\infty}$ structure 
$$
\frak{GW}^{\frak b}_{\ell} : H^*(X;\Lambda_0)^{\otimes \ell} \to H^*(X;\Lambda_0)
$$
by
$$
\aligned
&\langle
\frak{GW}^{\frak b}_{\ell}(\frak z_1,\dots,\frak z_{\ell}),
\frak z_0
\rangle_{X}\\
&=
\sum_{\ell'=0}^{\infty}\sum_{\alpha \in \pi_2(X)}
\frac{T^{\omega(\alpha)}}{\ell'!}
\int_{\mathcal M_{\ell+\ell'+1}^{\text{\rm cl}}(\alpha)}
(\text{\rm ev}^+)^*
(\frak z_1 \times \dots \times \frak z_{\ell}  \times \frak z_0 
\times \frak b^{\ell'}).
\endaligned
$$
See Propositions \ref{Kuraproduct}, \ref{multiproduct}.
\par
We remark that Theorem \ref{thmpcup0}
{\it does not} generalize to $\frak{GW}^{\frak b}_{\ell}$ for $\ell \ge 3$.
Namely let $\text{\bf x}_i \in HH_*(\cL_{\lambda_i})$ 
such that not all $\lambda_i$ coincide with each other. 
However
\begin{equation}\label{prodzeroll}
\frak{GW}^{\frak b}_{\ell}(\hat{\frak p}(\text{\bf x}_1),\dots,\hat{\frak p}(\text{\bf x}_{\ell})) = 0
\end{equation}
does not hold in general.  
In fact let us consider the case $X = \C P^2$ and $\frak b = 0$.
It is well known that
$$
QH(\C P^2;\Lambda) = \Lambda^3.
$$
The Clifford torus $T^2 \subset \C P^2$ has three elements 
$b_i \in H^2(T^2;\Lambda_0)$ for which Floer cohomology 
$HF((T^2,b_i),(T^2,b_i);\Lambda)$ is nonzero.
The values $\{\frak{PO}(b_i)\}$ are the product of $3T^{1/3}$ with distinct cubic roots of unity.
So we have $\cL_{\lambda_i}$ $i=1,2,3$. Since $HF((T^2,b_i),(T^2,b_i);\Lambda)$ 
is a Clifford algebra as a ring (Cho \cite{Cho05II}) the Hochschild homology 
of $\cL_{\lambda_i}$ is $\Lambda$.
(See \cite{kassel}.) It is easy to see 
(and follows from Theorem \ref{maintheorem3}) that 
$$
\hat{\frak p} : \prod_{i=1}^3 HH_*(\cL_{\lambda_i}) \to H(\C P^2;\Lambda)
$$
is an isomorphism.
Since {\it big} quantum cohomology of $ H(\C P^2;\Lambda)$ does not 
split into a direct product of three factors, 
it follows that (\ref{prodzeroll}) can not hold for \blue{some} $\ell$.
\par
This fact seems to imply that Fukaya category of $X$ does {\it not} 
determine $Comm_{\infty}$ structure of $X$.
\end{rem}

\section{Decomposition of quantum cohomology and of $A_{\infty}$ category.}
\label{sec: decomp}

In this section we apply our main results 
to study the relationship between the 
decompositions of the category $\cL$ 
and of quantum cohomology.
\par
Let $X$ be a symplectic manifold and 
$\frak b \in H^{\rm even}(X;\Lambda_0)$ be 
a bulk class. 
We consider the set $\text{\rm Crit}(X,\text{\rm st},\frak b)$ 
($\subset \Lambda$)\index[syindex]{CritX@$\text{\rm Crit}(X,\text{\rm st},\frak b)$}
of all values $\frak{PO}_{\frak b}(b)$ such that 
there exists a $\text{\rm st}$-relatively spin Lagrangian submanifold  
$L$ and $b \in \mathcal M_{\text{\rm weak}} (L;\frak b,\Lambda_0)$
satisfying %\marginpar{$(b,\frak b)$ or $(\frak b,b)$ ?}
$$
HF((L,b),(L,b);\Lambda) \ne 0.
$$
\begin{conj}\label{betti}
The following holds.
\begin{equation}\label{betti}
\#\text{\rm Crit}(X,\text{\rm st},\frak b)
\le
\sum \text{\rm rank}_{\Q} H^k(X;\Q).
\end{equation}
If the equality holds then $QH_{\frak b}(X;\Lambda)$ 
is semi-simple.
\end{conj}
We next consider the $\frak b$-deformed\index[syindex]{QHb@$QH_{\frak b}(X;\Lambda)$} 
quantum cohomology ring
$QH_{\frak b}(X;\Lambda)$.\index[syindex]{QHba@$QH_{\frak b}(X;\Lambda)_a$}
Let
\begin{equation}\label{decompQH}
QH_{\frak b}(X;\Lambda)
= 
\prod_{a=1}^A QH_{\frak b}(X;\Lambda)_a
\end{equation}
be its direct product decomposition  
(as a ring) into indecomposable $QH_{\frak b}(X;\Lambda)_a$. 
Let $e_a$ be the unit of the factor
$QH_{\frak b}(X;\Lambda)_a$.
It is easy to see that the decomposition 
(\ref{decompQH}) is unique.
\par
We take a finite set $\text{\bf L}$ of pairs 
$(L,b)$ of Lagrangian submanifolds $L$ and 
weak bounding cochains $b$.
Let $\cL$ be the cyclic $A_{\infty}$
category the set of whose object is
$\text{\bf L}$.
We then define the associated cyclic $A_{\infty}$ category $D^{\pi}({\cL})$
as in the introduction.
We say an object
$c$ of $D^{\pi}({\cL})$ is {\it indecomposable}\index{indecomposable} if the idempotent of $HF(c,c)$ is the unit or $0$.

\begin{rem}
For a pair $(L,b)$ of Lagrangian submanifold and its weak bounding cochain, 
 %\marginpar{Following Ono's suggestion, definition of indecomposability is 
%changed and remark is added. 2025 June KF}
it frequently happens that it decomposes into finitely many mutually isomorphic objects.
For example sometimes $HF((L,b),(L,b);\Lambda^{\C})$ is a matrix algebra $M(N)$.  (In the case of a $T^n$ orbit of a toric 
manifold, for example, its self-Floer homology is a Clifford algebra.  Sometimes it is a matrix algebra.) 
In such a case $(L,b)$ is isomorphic to the direct sum of $N$ copies of indecomposable objects.
Proposition \ref{factorunique} below applies to such $(L,b)$ since direct summands are mutually isomorphic.
\end{rem}

It is easy to see that any object of $D^{\pi}(\cL)$ is a 
direct sum of  indecomposable objects.
If $c$ is indecomposable, then there exists a unique $\lambda \in \Lambda_0$ such that 
$c \in D^{\pi}(\cL_{\lambda})$.
We write this $\lambda$ as $\lambda(c)$.
\par
We remark that
$1_c \in B_0(HF(c,c)[1])$ and is different from $\text{\bf e}_c \in B_1(HF(c,c)[1])$ that is the unit.
\begin{prop}\label{factorunique}
If $c$ is indecomposable then  there exists a  
unique $a(c) \in A$ such that $
{\frak q}_c(e_{a(c)};1_c) \in HF(c,c)
$
is nonzero. 
\end{prop}
\begin{proof}
By assumption $\frak q(e_a;1_c) \in HF(c,c)$ is either $\text{\bf e}_c$ or $0$.
Since $e_a \cup^{\frak b} e_{a'} = 0$ for $a\ne a'$,
${\frak q}(e_a;1_c)=\text{\bf e}_c$ for at most one $a$. 
Since ${\frak q}(\sum e_a;1_c) =\text{\bf e}_c$ there exists 
exactly one $a$ such that ${\frak q}(e_a;1_c)=\text{\bf e}_c$.
\end{proof}
Proposition \ref{factorunique} implies that
$D^{\pi}(\cL)$ is decomposed into \index[syindex]{DpiLa@$D^{\pi}(\cL;a)$}
direct sum
\begin{equation}\label{catedecompose}
D^{\pi}(\cL) = \prod_{a\in A} D^{\pi}(\cL;a)
\end{equation}
where an object of $D^{\pi}(\cL;a)$ is a direct sum of the objects $c$ with 
$a(c) = a$.
\par
We remark
\begin{equation}
HH_*(\cL) 
\cong HH_*(D^{\pi}(\cL))
\cong
\prod_{a\in A} HH_*(D^{\pi}(\cL;a)).
\end{equation}

\begin{lem}\label{ppreservesa}
$$
\hat{\frak p} \left( HH_*(D^{\pi}(\cL;a) \right)
\subset QH_{\frak b}(X;\Lambda)_a.
$$
\end{lem}
\begin{proof}
Let $[\text{\bf x}] \in HH_*(D^{\pi}(\cL;a))$.
We need to show 
\begin{equation}\label{115toshow}
\langle\hat{\frak p}(\text{\bf x}),e_{a'}\rangle_{\text{\rm PD}_X} = 0,
\end{equation}
for $a' \ne a$. 
In fact then for any element $y$ of $QH_{\frak b}(X;\Lambda)_{a'}$ we will have
$$
\langle \widehat{\frak p}({\bf x}),y\rangle_{\text{\rm PD}_X}
= \pm \langle \widehat{\frak p}({\bf x})\cup^{\frak b}y, e_{a'}\rangle_{\text{\rm PD}_X}
= \pm \langle\widehat{\frak p}({\bf x} \cup \widehat{\frak q}(y)),e_{a'}\rangle_{\text{\rm PD}_X} = 0.
$$

To prove (\ref{115toshow}), we first consider the case 
$\text{\bf x} = x 
\in B_1(HF(c,c)[1])$.
Using the fact that
${\frak q}(e_{a'};1_c) = 0$ for any object $c$ in $D^{\pi}(\cL;a)$,
we calculate:
$$
\langle\hat{\frak p}(x),e_{a'}\rangle_{\text{\rm PD}_X}
= 
\langle x,{\frak q}(e_{a'};1_c)\rangle_{\rm HH} = 0.
$$

\par
We next consider the length filtration:
\begin{equation}\label{numberfiltHH}
\prod_{\lambda}F^k CH_*(D^{\pi}(\cL_{\lambda};a))
=
\prod_{\lambda}\bigoplus_{k'\ge k}
CH_{k'}(D^{\pi}(\cL_{\lambda};a)[1])
\end{equation}
of $\prod_{\lambda}CH_*(D^{\pi}(\cL_{\lambda};a))$.
Since we put $\frak m_0 = 0$ on each of 
$D^{\pi}(\cL_{\lambda};a)$
(See Subsection \ref{cvnoncomon}.)
the Hochschild boundary operator preserves the length 
filtration.
\par
We next study 
\begin{equation}\label{qtobezero}
{\frak q}(e_{a'};\text{\bf x})
\end{equation}
for $\text{\bf x} \in (\ref{numberfiltHH})$ and show that 
$(\ref{qtobezero})$ is zero by induction on $k$.
\par
By Theorem \ref{theoremD1} and $e_{a'}\cup^{\frak b} e_{a'} 
= e_{a'}$ we have
\begin{equation}\label{113calcu}
\aligned
&{\frak q}(e_{a'};\text{\bf x})
=
{\frak q}(e_{a'}\cup^{\frak b} e_{a'};\text{\bf x})\\
&=
\sum_{c}(-1)^{\maltese_c}
\frak m(\text{\bf x}^{(5;1)}_c \otimes
{\frak q}(e_{a'};\text{\bf x}^{(5;2)}_c) \otimes
\text{\bf x}^{(5;3)}_c \otimes
\frak q(e_{a'};\text{\bf x}^{(5;4)}_c) \otimes
\text{\bf x}^{(5;5)}_c). \\
\endaligned
\end{equation}
 The right hand side 
vanishes by induction hypothesis.
Thus we have proved that $(\ref{qtobezero})$ is zero.
The duality then implies (\ref{115toshow}).
\end{proof}
\par
We conjecture that each of the 
factors $D^{\pi}(\cL;a)$ contains at most one category $\cL_{\lambda}$. % is contained in one of the 
%factors $D^{\pi}(\cL;a)$. Namely:
\begin{conj}\label{lambdacoinc}
If $a(c) = a(c')$ then 
$\lambda(c) = \lambda(c')$.
\end{conj}
\begin{prop}
Conjecture  $\ref{lambdacoinc}$ implies Conjecture $\ref{betti}$.
In fact Conjecture  $\ref{lambdacoinc}$ implies
\begin{equation}\label{Criandbetti}
\#\text{\rm Crit}(X,\text{\rm st},\frak b)
\le A.
\end{equation}
\par
If  the equality holds in Conjecture $\ref{betti}$ in addition
then the map
$\widehat{\frak p} : HH_*(\cL,\cL) \to H(X;\Lambda)$ is an isomorphism.
\end{prop}
\begin{proof}
It is easy to see that Conjecture  $\ref{lambdacoinc}$ implies the inequality (\ref{Criandbetti}).
Suppose that equality holds in (\ref{betti}).
Then $A$ must be the Betti number. Therefore each factor of 
(\ref{decompQH}) is rank one, which implies semi-simplicity.
\par
Since $\widehat{\frak p} : HH_*(D^{\pi}(\cL;a)) \to H(X;\Lambda)$
is nonzero, Lemma \ref{ppreservesa} implies that 
$\widehat{\frak p} : HH_*(\cL) \to H(X;\Lambda)$ is an isomorphism.
\end{proof}
\begin{prop}
Suppose that $\frak b$ is concentrated in degree $2$  
and that  $c$ (resp. $c'$) is a direct factor of  $(L,b)$ (resp. $(L',b')$) with $b,b'$ concentrated in degree $1$.
Then
Conjecture $\ref{lambdacoinc}$ holds. 
\end{prop}
\begin{proof}
This is an easy consequence of Theorem \ref{c1crit}.
In fact Theorem \ref{c1crit} implies
$$
c_1(X) \cup^{\frak b} e_a = \lambda(c) e_a
=  \lambda(c') e_a.
$$
\end{proof}
\begin{thm}
If 
$$
\text{\bf e}_X \in \text{\rm Im}(\hat{\frak p} : HH_*(\cL) \to H_*(X;\Lambda))
$$
then Conjecture $\ref{lambdacoinc}$ holds.
\end{thm}
\begin{proof}
Let $a = a(c) = a(c')$.
We take $\lambda_i$ and $\text{\bf x}_i \in HH_*(\cL_{\lambda_i})$
such that $e_a = \sum_{i=1}^k\hat{\frak p}(\text{\bf x}_i)$, and 
$\lambda_i \ne \lambda_j$ for $i\ne j$.
Then
$$
e_a = e_a \cup^{\frak b} e_a
= 
\sum_i  \hat{\frak p}(\text{\bf x}_i)\cup^{\frak b} \hat{\frak p}(\text{\bf x}_i),
$$
since $\langle\hat{\frak p}(\text{\bf x}_i),\hat{\frak p}(\text{\bf x}_j)\rangle_{\text{\rm PD}_X}
= 0$ for $i\ne j$. (We use Theorem \ref{thmpcup0}  here.)
Thus we have
$$
\hat{\frak p}(\text{\bf x}_i)\cup^{\frak b} \hat{\frak p}(\text{\bf x}_i)
=
\hat{\frak p}(\text{\bf x}_i).
$$
The indecomposability of $QH_{\frak b}(X;\Lambda)_a$ now implies
$e_a = \hat{\frak p}(\text{\bf x}_i)$ for some $i$.
\par
Then
$
{\frak q}(\widehat{\frak p}(\text{\bf x}_i);1_c) \ne 0
$
implies $\lambda(c) = \lambda_i$.
We also have $\lambda(c') = \lambda_i$ in the same way. 
\end{proof}
\begin{thm}\label{semsimpleOK}
Conjecture $\ref{lambdacoinc}$ holds if $QH_{\frak b}(X;\Lambda)$ is semi-simple. 
%Here $a=a(c) = a(c')$.
\end{thm}
\begin{proof}
Since ${\frak q}(e_a;1_c) \ne 0$, ${\frak q}(e_a;1_{c'}) \ne 0$ there exists 
$x_c \in HF(c,c)$ and $x_{c'} \in HF(c',c')$ such that
$$
\langle \hat{\frak p}(x_c),e_a\rangle_{\text{\rm PD}_X} \ne 0
\ne \langle \hat{\frak p}(x_{c'}),e_a\rangle_{\text{\rm PD}_X}.
$$
We use Theorem \ref{thm:duality} here.
By assumption $QH_{\frak b}(X;\Lambda)_a \cong \Lambda$. 
\par
By Proposition \ref{factorunique}, 
$$
\langle \hat{\frak p}(x_c),e_{a'}\rangle_{\text{\rm PD}_X} = 0
= \langle \hat{\frak p}(x_{c'}),e_{a'}\rangle_{\text{\rm PD}_X}.
$$
for any $a' \ne a$. Therefore
$$
\langle \hat{\frak p}(x_c),\hat{\frak p}(x_{c'})\rangle_{\text{\rm PD}_X} = 0.
$$
Theorem \ref{semsimpleOK} now follows from Proposition \ref{factorunique}.
\end{proof}

\begin{thm}\label{ssfinite}
We assume that $QH_{\frak b}(X;\Lambda)$ is semi-simple.
Let $b_i \in \mathcal M_{\text{\rm weak}}(L_i,\frak b)$ 
be such that, for $i=1,\dots,m$,
$$
HF((L_i, b_i),(L_i,b_i);\Lambda)\ne 0.
$$ 
We also assume that one of the following holds for each $i\ne j$:
\begin{enumerate}
\item 
$\frak{PO}_{\frak b}(b_i) \ne \frak{PO}_{\frak b}(b_j)$.
\item 
$\frak{PO}_{\frak b}(b_i) = \frak{PO}_{\frak b}(b_j)$ and
$$
HF((L_i,b_i),(L_j,b_j);\Lambda) = 0.
$$
\end{enumerate}
\par
Then 
\begin{equation}\label{misestirank}
m \le \sum_k\text{\rm rank}_{\Q} H^k(X;\Q).
\end{equation}
\par
If equality holds in $(\ref{misestirank})$ then 
the map
$\widehat{\frak p} : HH_*(\cL) \to H(X;\Lambda)$ is an isomorphism  for $\text{\bf L} =\{(L_i,b_i) \mid i=1,\dots,m\}$.
\end{thm}
\begin{proof}
Since $QH_{\frak b}(X;\Lambda)$ is semi-simple we have
$QH_{\frak b}(X;\Lambda)_a \cong \Lambda$ for each $a 
= 1, \dots, A$ and $A$ is the Betti number of $X$.
For each $i$ let us take 
an indecomposable component  $c_i$ of $(L_i,(\frak b,b_i))$.
We put $a_i \in A$ to be the semi-simple factor mapping non-trivially to the Floer homology of $c_i$. 
The assumption and Theorem \ref{thmpcup0} imply that 
$$
\langle \widehat{\frak p}(x_i),   \widehat{\frak p}(x_j)\rangle_{{\rm PD}_X} = 0
$$
for $x_i \in HF(c_i,c_i)$, $x_j \in HF(c_j,c_j)$.
Using $QH_{\frak b}(X;\Lambda)_a \cong \Lambda$ 
and Proposition \ref{nontriialp}, this implies that 
$a_i \ne a_j$.
Therefore 
$$
m \le \# A \le \sum_k\text{\rm rank}_{\Q} H^k(X;\Q).
$$
Suppose the equality holds. By semi-simplicity, 
we can find $x_i \in HF(c_i,c_i)$ such that
$\frak p(x_i) = w_ie_{a(c_i)}$ for $w_i \in \Lambda \setminus \{0\}$.
Then
$$
\text{\bf e}_X = \sum_{i} \widehat{\frak p}\left(\frac{x_i}{w_i}\right).
$$
\end{proof}
\begin{rem}
$  $ \par
\begin{enumerate}
\item
In the case $X = T^2$, the conclusion of Theorem \ref{ssfinite} does not hold.
Similarly it does not hold in many of Calabi-Yau manifolds.
In fact in many such cases we find infinitely many Lagrangian tori
which are mutually disjoint and have nontrivial Floer homology.
\par
For example let us consider the case $X=T^{2n} = \C^n/\Z^n$ and 
$L_i = T^n = \R^n/\Z^n + a_i\sqrt{-1}$. ($a_i \in \R^n)$. 
In this case $\hat{\frak p} : HF(L_i,L_i;\Lambda) \to H(X;\Lambda)$
coincides with Gysin map $H^k(T^n;\Lambda) \to H^{k+n}(T^{2n};\Lambda)$.
Therefore obviously 
$$
\langle \widehat{\frak p}(x_i), \widehat{\frak p}(x_j)\rangle_{{\rm PD}_X}  = 0,
$$
for $x_i\in HF(L_i,L_i;\Lambda)$, $x_j\in HF(L_j,L_j;\Lambda)$.
Moreover the open closed map $\widehat{\frak p}$ 
(defined on the Hochschild homology of the $A_{\infty}$ category whose 
objects are Lagrangian $T^n$ as above) never hits $\text{\bf e}_X\in H^0(X)$.
\item
In \cite[Corollary 24.4]{spectre} we proved that the conclusion of Theorem \ref{ssfinite} does not hold 
for the cubic surface and $\frak b =0$. In fact the quantum cohomology of 
cubic surface is not semi-simple.
\item In case $X$ is $S^2 \times S^2$ or certain $k\ge 2$ points toric
blow up of $\C P^2$, we find infinitely many 
Lagrangian tori that are mutually disjoint and have 
nontrivial Floer homology (\cite{fooo:bulk,fooo:S2S2}).
However their quantum cohomology are semi-simple.
In those examples, we varied the bulk class $\frak b$ 
to obtain nontrivial Floer homology.
So Theorem \ref{ssfinite} does not apply to that case. 
%\item
%We remark that the last part of Theorem \ref{ssfinite} together with
%Theorem \ref{maintheorem2} in the introduction can be used to prove 
%homological Mirror symmetry in certain cases.
%For example in many of the toric cases we can check the assumption 
%directly. 
%(See \cite{afooo} the general toric case etc..)
\end{enumerate}
\end{rem}
Let $L \subset X$ be a Lagrangian submanifold.
We take its Weinstein neighborhood and 
take  a family $L(u)$ of Lagrangian submanifolds
parametrized by $u$ in a small neighborhood of 
$0$ in $H^1(L;\R)$. Namely $L(u)$ is the graph 
of a closed one form representing $u$.
(Here we identify the Weinstein neighborhood with a neighborhood of 
zero section in $T^*L$.)
The Lagrangian $L(u)$ is unique up to Hamiltonian isotopy
supported in the Weinstein neighborhood.
\begin{conj}\label{finitenessLag}
Suppose $QH_{\frak b}(X;\Lambda)$ is semi-simple.
Then there exist only a finite number of $u$ that enjoy the following property:
\par
There exists $b \in \mathcal M_{\text{\rm weak}}(L(u);\frak b)$ 
such that
$$
HF((L(u),b),(L(u),b);\Lambda) \ne 0.
$$ 
\end{conj}
Theorem \ref{ssfinite} implies Conjecture \ref{finitenessLag} 
in case when 
there exists a fiber bundle $\pi : L\to T^b$ 
such that $b = \text{\rm rank} \,H^1(L;\R)$ and that
the first Betti number of the fiber is zero.
In fact, in such case, we may take $L(u)$ so that $L(u) \cap L(u') = \emptyset$ for 
$u\ne u'$.
The general case seems to be related to the 
Lagrangian version of the Flux conjecture. (See \cite{On}.)
\par
We next mention a relationship of Theorem \ref{maintheorem2} to the 
(super) heaviness of Lagrangian submanifolds.
We review its definition here.
Let $e_a$ be an idempotent in $QH_{\frak b}(X;\Lambda)$ that is 
the unit of $QH_{\frak b}(X;\Lambda)_a$. 
By \cite{EP:morphism,usher:talk,spectre} we obtain 
$$
\zeta_{e_a} : \widetilde{\text{\rm Ham}}(X,\omega) \to \R
$$
where $\widetilde{\text{\rm Ham}}(X,\omega)$ is the universal cover of the group 
of Hamiltonian diffeomorphisms of $X$.\index[syindex]{HamXomega@${\text{\rm Ham}}(X,\omega)$}
Let $\widetilde\varphi_H \in \widetilde{\text{\rm Ham}}(X,\omega)$ be the 
element corresponding to the time dependent Hamiltonian $H : [0,1] \times X \to \R$.
Namely 
$$
\frac{d}{dt}\varphi_H^t = X_{H_t} \circ \varphi_H^t, 
\quad \varphi_H^0 = identity,
$$
$t \mapsto  \widetilde\varphi_H^t$ 
is its lift to $\widetilde{\text{\rm Ham}}(X,\omega)$ 
and $\widetilde\varphi_H = \widetilde\varphi_H^1$.
\par
We call a subset $Y$ is $\zeta_{e_a}$-{\it heavy}\index{heavy} if 
\begin{equation}
\zeta(\widetilde\varphi_H)  \leq \sup \{H(t,p) \mid t\in [0,1], p \in Y\}.
\end{equation}
We call $Y$ is $\zeta_{e_a}$-{\it superheavy}\index{superheavy} if 
\begin{equation}
\zeta(\widetilde\varphi_H)  \geq \inf \{H(t,p) \mid t\in [0,1], p \in Y\}.
\end{equation}
\begin{rem}
Here we use the convention of \cite{spectre} that is equivalent to but is 
slightly different from Entov-Polterovich's.  
Namely, our $\zeta^\text{\rm FOOO}$ and Entov-Polterovich's $\zeta^\text{\rm EP}$ 
are related by $\zeta^\text{\rm FOOO}(F)=-\zeta^\text{\rm EP}(-F)$.   
See \cite[\newred{Subsection 1.3}]{spectre}.  

\end{rem}
Let $\{(L_i,b_i)\mid i=1,\dots,m\}  = \bL$ be as in the introduction. 
We define the filtered $A_{\infty}$ category $\cL$  as in Section 
\ref{sec:cyclicfil}.
\begin{thm}\label{spectrethm}
$  $ \par
\begin{enumerate}
\item
If $\hat{\frak q}(e_a) \in HH^*(\cL,\cL)$ is nonzero,
then $\bigcup_{i=1}^m L_i$ is $\zeta_{e_a}$-heavy.
\item
If 
\begin{equation}\label{phaseaimage}
PD(e_a) \in \text{\rm Im}(\hat{\frak p} : HH_*(\cL) \to H_*(X;\Lambda)),
\end{equation}
then  $\bigcup_{i=1}^m L_i$ is $\zeta_{e_a}$-superheavy.
\end{enumerate}
\end{thm}

\begin{proof}
The proof is based on the following slight enhancement of \cite[Proposition 25.4]{spectre} 
and is similar to its \cite[Subsection 5.3]{fooo09}. 
\begin{prop}\label{spectreestimate}
For any $a \in H^*(X;\Lambda)$
we have
\begin{equation}\label{estimatespectre}
\rho^{\frak b}(H;a)
\ge
\inf \{- H(t,p)\mid (t,p)\in S^1\times \bigcup_{i=1}^m L_i \}
- \frak v_T(\hat{\frak q}(a)).
\end{equation}
\end{prop}
Here $\rho^{\frak b}$ is the spectral invariant with bulk 
that is defined in \cite[Definition 7.6]{spectre}.
The proof of Proposition \ref{spectreestimate} is 
the same as that in \cite[Proposition 25.4]{spectre}.
In fact, the only difference is that we use several Lagrangian submanifolds,  
while in \cite[Proposition 25.4]{spectre} we considered only one 
Lagrangian submanifold. However the proof works without any change 
and so can be safely omitted.
\par
The proof of (1) is the same as \cite[Theorem 18.8 (1)]{spectre}  
using Proposition \ref{spectreestimate} in place of \cite[Proposition 18.9]{spectre}.
So we omit it.
\par
We next prove (2).
We consider the set $\bL_a$ of all indecomposable components $c_j$ of  
members $(L_i,b_i)$ of $\bL$ such that $a(c_j) = a$. Let 
$\cL(\bL_a)$ be the full subcategory whose object set is $\bL_a$.
\par
Our assumption (\ref{phaseaimage}) 
implies 
$$
e_a \in \text{\rm Im}(\hat{\frak p} : HH_*(\cL(\bL_a)) \to H(X;\Lambda)).
$$
Therefore an analogue of Theorem \ref{maintheorem3} implies that

$$
\hat{\frak q} : QH(X;\Lambda)_a \to HH^*(\cL(\bL_a),\cL(\bL_a))
$$
is an isomorphism.
Therefore
$$
\hat{\frak q} : QH(X;\Lambda)_a \to HH^*(\cL,\cL)
$$
is injective.
\par
We next use  \cite[Theorem 15.1]{spectre} (that is 
a minor modification of \cite{EP:morphism}) to show
\begin{equation}\label{25100}
\rho^{\frak b}(\widetilde \psi_H^n;e_{a})
=
- \inf \{
\rho^{\frak b}(\widetilde \psi_H^{-n};b)
\mid \Pi(e_{a},b)\ne 0
\}.
\end{equation}
Here 
$\Pi : QH(X;\Lambda) \otimes QH(X;\Lambda) \to \C$
is the coefficient of $T^0$ of the Poincar\'e duality pairing.
\par
Since
the image of the map
$\hat{\frak q}$
is a finite dimensional vector space over $\Lambda$, we can apply the
argument of \cite[Subsection 8.1]{spectre}
to find a standard basis 
$\hat{\frak q}(e_1),\dots,\hat{\frak q}(e_k)$ of the image of
$\hat{\frak q}$. (See  \cite[Definition 8.3]{spectre} for the definition of standard basis.)
Then we have 
\begin{equation}\label{259formu}
\aligned
\frak v_T\left(\hat{\frak q}\left(\sum_{i=1}^k x_i e_i \right)\right)
& = \min \{\frak v_T(x_i \hat{\frak q} (e_i)) \mid i=1,\dots,k \} \\
& \le \min \{\frak v_T(x_i ) \mid i=1,\dots,k \}+ C_1
\endaligned
\end{equation}
where $C_1$ is independent of $x_i$.
\par
Let
$\widetilde \psi_H \in \widetilde{\text{\rm Ham}}(M,\omega)$.
We estimate the right hand side of (\ref{25100}).
Suppose $\Pi(e_{a},b) \ne 0$.
We put
$$
e_{a}\cup^{\frak b}b = \sum_{i=1}^k x_i e_i,
\qquad x_i \in \Lambda.
$$
Using $\Pi(e_{a},b) \ne 0$,  \newred{it is easy to see}
$
\frak v_T(e_{a}\cup^{\frak b} b)\le 0.
$
Therefore
\begin{equation}\label{2550}
\min\{\frak v_T(x_i) \mid i = 1,\dots, k\}
\le C_2,
\end{equation}
where
$C_2 = -\min\{\frak v_T(e_i) \mid i = 1,\dots, k\}$.
\par
By the triangle inequality (\cite[Theorem 7.8 (5)]{spectre}),
$$
\rho^{\frak b}(\widetilde \psi_H^{-n};b)
\ge
\rho^{\frak b}(\widetilde \psi_H^{-n};\sum_{i=1}^k x_ie_i)
- \rho^{\frak b}(\underline 0;e_{a}).
$$
Using Proposition \ref{spectreestimate} the right hand side is
not larger than
$$
n \inf \{- \widetilde H(t,p) \mid (t,p) \in S^1\times M\}
-\frak v_T\left(\hat{\frak q} \left(\sum_{i=1}^k x_i e_i\right)\right)
-  \rho^{\frak b}(\underline 0;e_{\frak y}).
$$
By (\ref{259formu}), this is
not smaller than 
$$
n \inf \{H(t,p) \mid (t,p) \in S^1\times M\}
-\min\{ \frak v_T(x_i) \mid i=1,\dots, k\}
- C_1 - \rho^{\frak b}(\underline 0;e_{\frak y}).
$$
Using (\ref{2550}) we have
$$
\rho^{\frak b}(\widetilde \psi_H^{-n};b)
\ge
n\inf\{H(t,p) \mid (t,p) \in S^1\times M\}
- C_2 - C_1 - \rho^{\frak b}(\underline 0;e).
$$
Therefore by \eqref{25100}, we obtain
$$
-\frac{\rho^{\frak b}(\widetilde \psi_H^n;e_{\frak y})}{n}
\ge \inf \{H(t,p) \mid (t,p) \in S^1\times M\}
- \frac{C_3}{n},
$$
where $C_3$ is independent of $n$. Therefore we obtain
$$
\aligned
\frac{\rho^{\frak b}(\widetilde \psi_H^n;e_{\frak y})}{n}
& \le  -\inf \{H(t,p) \mid (t,p) \in S^1\times M\}
+ \frac{C_3}{n}\\
& =  \sup \{- H(t,p) \mid (t,p) \in S^1\times M\}
+ \frac{C_3}{n}.
\endaligned
$$
By letting $n \to \infty$, we have finished the proof of Theorem \ref{spectrethm}.
(Note
$$
\zeta_{e_a}(\widetilde \psi_H)
= -\lim_{n\to \infty}\frac{\rho^{\frak b}(\widetilde \psi_H^n;e_{\frak y})}{n}.
$$
See \cite[(14.4)]{spectre}.)
\end{proof}
Theorem  \ref{spectrethm} and \cite[Theorem 1.4]{EP:rigid} implies the following 
intersection result.
Let $\text{\bf U} = \{(U_i,b^U_i) \mid i=1,\dots,k\}$ be a set of 
$\text{\rm st}$-relatively spin Lagrangian submanifolds and 
weak bounding chains. We define the cyclic filtered $A_{\infty}$
category $\sU$ from it.
\begin{cor}\label{intersectionUL}
Suppose
$\hat{\frak q}(e_a) \in HH^*(\sU,\sU)$ is nonzero and 
$e_a \in \text{\rm Im}(\hat{\frak p} : HH_*(\cL, \cL) \to H(X;\Lambda))$.
Then there exists $i,j$ such that
$$
U_j \cap L_i \ne \emptyset.
$$
\end{cor}
We remark that Corollary \ref{intersectionUL} also follows 
from Theorem \ref{thm:commutative_diagram_p_q}.  In fact the assumption implies that 
$$
\hat{\frak q}\circ \hat{\frak p} :
HH_*(\cL, \cL) \to HH^*(\sU,\sU)
$$
is nonzero.
\par
Corollary \ref{intersectionUL} also follows from Theorem \ref{maintheorem2}. In fact the second 
assumption implies that $\bL$ split generates $D^{\pi}(\cL;a)$.
The first assumption implies that one of $(U_j,b^U_j)$ 
contains a direct summand that is an element of $D^{\pi}(\cL;a)$.
\par\smallskip
The results of this section are mostly concerned with the case when
quantum cohomology is close to be semi-simple.
In the Calabi-Yau case, namely in case $c^1(X) = 0$, the situation is different.
\begin{conj}\label{conj2717}
If  $c^1(X) = 0$ then $\text{\rm Crit}(X,\text{\rm st},\frak b) =\{0\}$.
\end{conj}
Let $L$ be a Lagrangian submanifold of Calabi-Yau manifold $X$.
Suppose that
\begin{enumerate}
\item[(i)]
$b \in 
\mathcal M_{\text{\rm weak}}(L;\frak b;\Lambda_0)$.
\item[(ii)] 
The Maslov index $\pi_2(X;L) \to \Z$ vanishes. 
\end{enumerate}
Then by dimension counting of the moduli space, it is easy to show
\begin{equation}\label{PO0}
\frak{PO}_{\frak b}(b) = 0.
\end{equation}
By Theorem \ref{c1crit},  (\ref{PO0}) also follows from the following assumption.
\begin{enumerate}
\item [(a)]
$b \in 
H^1(L;\Lambda_0) \cap 
\mathcal M_{\text{\rm weak}}(L;\frak b;\Lambda_0)$.
\item [(b)]
$\frak b \in H^2(X;\Lambda_0)$.
\item[(c)]
$HF((L,(\frak b,b)),(L,(\frak b,b));\Lambda) \ne 0$.
\end{enumerate}
\par
We point out that there is an example of Lagrangian submanifold in 
K3 surface where (a)(b)(c) above hold but (ii) does not hold.
This suggests that Conjecture \ref{conj2717} is not so obvious.
\par
To construct such an example we study Kummer surface below.
We consider 
$
X_0 = {T^4}/{\{\pm 1\}}  
$
where $-1$ acts by $[x_1,\dots,x_4] \mapsto [-x_1,\dots,-x_4]$.
We put any symlectic form $\omega$ on it that is induced by \blue{an invariant} symplectic 
form on $T^4$.
$X_0$ have $16$ singular points of type $A_1$ (that is locally $\R^4/\{\pm 1\}$).
We replace small neighborhoods of those 16 points by the Milnor fiber
that is symplectomorphic to a neighborhood of zero section 
in $T^*S^2$. 
We thus obtain a smooth symplectic manifold $X$, that is actually 
a Kummer surface.
(It contains 16 Lagrangian spheres.)
\par
$T^*S^2$ contains a family of Lagrangian submanifolds $T(u)$, that we describe below.
For $u=1$ we take
$
L(1) \subset T^*S^2 
$
as the union of unit tangent vector of closed geodesics of $S^2$ which contains 
both north and south poles. (See \cite[appendix 2]{fooo:S2S2}.)
(Here we identify $TS^2$ with $T^*S^2$ using the standard metric of $S^2$.)
\par
For general $u$ we put $L(u) = \{uv \mid v \in L(1)\}$.
\par
We may identify a small neighborhood of zero section in $T^*S^2$
to an open set in $X$.
Therefore $L(u)$ may be regarded as a Lagrangian submanifold of $X$ for 
sufficiently small $u$, that we write 
$T(u)$.  (Actually there are 16 different choices of $T(u)$ for given $u$.) 
\par
It is easy to see that the minimal Maslov number of $T(u)$ is $2$ and is not $0$.
On the other hand, we have the following:
\begin{prop}\label{propKummer}
If $u$ is sufficiently small then there exists $\frak b(u) \in H^2(X;\Lambda_0)$ 
and $b(u) \in H^1(T(u);\Lambda_0)$ such that
$$
HF((T(u),(\frak b(u),b(u))),(T(u),(\frak b(u),b(u)));\Lambda_0)
= H(T^2;\Lambda_0).
$$
\end{prop}
\begin{rem}
$  $ \par
\begin{enumerate}
\item
Since $\dim X = 4$ we have $$H^1(T(u);\Lambda_0) \subset \mathcal M_{\text{\rm weak}}(L(u);\frak b;\Lambda_0),$$
by \cite[Theorem A.1]{fooo:S2S2}.)
\item
(\ref{PO0}) holds in our case, since (a), (b), (c) above are satisfied.
\end{enumerate}\end{rem}
\begin{proof}[Proof of Proposition \ref{propKummer}.]
We take compatible almost complex structure $J_U$ on $U$ 
so that $J_U$ is exactly the same as one used in \cite[Section 6]{fooo:S2S2}.
We choose an almost complex structure $J$ on $X$ so that 
there is no $J$-holomorphic disk bounding $T(u)$ with 
Maslov index $\le 0$ and $J$ coincides with $J_U$ on $U$. (We can prove it by using starndard 
dimension counting argument.)
\par
Let $\beta \in \Pi_2(X;T(u))$ be an element of Maslov index $2$.
We define the number $n_{\beta}$ as the degree of the 
evaluation map 
$$
\mathcal M_{1}(T(u);\beta) \to T(u).
$$
The degree is well-defined since there is no disk with Maslov number $\le 0$.
(We remark that $n_{\beta}$ may depend on the choice of almost complex structure $J$, 
however.)
\par
By \cite[Section 6]{fooo:S2S2} there are four  classes $\beta_i \in \Pi_2(X;T(u))$, 
$i=1,\dots,4$ such that $n_{\beta_i} = 1$ and 
$\beta_i$ is realized by a holomorphic disk contained in $U$.
Moreover $\beta_i \cap \omega = u$.
There is a positive number independent of $u$ such that 
if $\beta\ne \beta_i$ and if $n_{\beta} \ne 0$ then 
$\beta \cap \omega  > \lambda$.
(This is a consequence of the fact that a $J$-disk realizing the class 
$\beta$ cannot be contained in $U$.)
\par
We may choose the basis of $H^1(T(u);\Z)$ such that
$$
\partial\beta_1 = (1,0),
\quad 
\partial\beta_2 = \partial\beta_3 = (0,-1),
\quad
\partial\beta_4 = (-1,-2).
$$
Therefore we have
$$
\frak{PO}(x_1,x_2) \equiv 
T^u(y_1 + 2y_2^{-1} + y_1^{-1} y_2^{-2})
\equiv 
T^u(y_1 + y_1^{-1} y_2^{-1})^2
\mod T^{\lambda}.
$$
Here $(x_1,x_2) \in H^1(T(u);\Lambda_0)$ and we use the 
above mentioned basis here. 
The variable $y_i$ is defined by
$y_i = e^{x_i}$.
\par
We next introduce a bulk deformation.
We choose cycles $D_1,D_2$ of $X$ representing classes in $H_2(X;\Q)$ as follows.
The zero section $=S^2$ is contained in $T^*(S^2)$.
We may regard it as a submanifold of $X$ and denote it by $D_1$.
We choose $D_2$ so that 
$$
D_2 \cap U = U \cap T_{\text{\rm sp}}^*S^{2}
$$ 
where $T_{\text{\rm sp}}^*S^{2}$ is the fiber of south pole.
(We can find such $D_2$ because the unit sphere bundle of $T^*S^2$ is 
a rational homology sphere.)
Then we have
$$
\aligned
&\beta_1\cap D_1 = 0 = \beta_4\cap D_1, 
\quad \beta_2\cap D_1 =1, \quad \beta_3\cap D_1 =- 1 \\
&\beta_1\cap D_2 = 1, \quad \beta_2\cap D_2=\beta_3\cap D_2=\beta_4\cap D_2=0. 
\endaligned
$$
We put
$$
\frak b(w_1,w_2) = w_1 PD([D_1]) + w_2 PD([D_2]) \in H^2(X;\Lambda_+).
$$
Therefore, in the same way as \cite[Section 7]{fooo:S2S2}, we have
$$
\frak{PO}_{\frak b(w_1,w_2)}(x_1,x_2) \equiv 
T^u(e^{w_2}y_1 + (e^{w_1}+e^{-w_1})y_2^{-1} + y_1^{-1} y_2^{-2})
\mod T^{\lambda}.
$$
We put $b = (0,\pi \sqrt{-1})$. Then we have
\begin{equation}\label{criteq}
\aligned
\frac{\partial \frak{PO}_{\frak b(w_1,w_2)}}{\partial y_1}(b)
&= T^u (e^{w_2} - 1) + T^{\lambda}P_1(w_1,w_2), \\
\frac{\partial \frak{PO}_{\frak b(w_1,w_2)}}{\partial y_2}(b)
&= T^u (2 - e^{w_1} - e^{-w_1}) + T^{\lambda}P_2(w_1,w_2), 
\endaligned
\end{equation} %\marginpar{{Formula corrected. KF}}
where $P_1,P_2$ are power series of $w_1,w_2$ with coefficient in $\Lambda_0$.
It is easy to see that if $u<\lambda$ then there exists $w_1,w_2 \in \Lambda_+$ 
such that the right hand side of (\ref{criteq}) vanishes.
Proposition \ref{propKummer} now follows from \cite[Theorem 2.3]{fooo:S2S2}.
\end{proof}
\begin{rem}
In the above proof we chose $J$ and fixed it. The existence of $b$ 
satisfying the conclusion of Proposition \ref{propKummer} is independent of 
$J$. However $b$ itself may depend on $J$.
\end{rem}

\section{$c_1$ and the value of potential function.}
\label{sec: value}
In this section we consider a bulk class $\frak b \in H^2(X;\Lambda_0)$, with a chosen representative differential form for  $\frak b_2 \in H^2(X;\F)$ which vanishes on a Lagrangian $L$. We also consider a class $b \in H^1(L;\Lambda_0)$, where the leading order term $b_1 \in H^1(L; \F)$ determines the class of the primitive $\theta_L$ for the restriction of $\frak b_2$ to $L$, as in Equation \eqref{fixthetaL}.
\begin{thm}\label{c1crit}
%We decompose $b= b_0 + b_+$. (See $(\ref{bdecomp})$.)
Suppose the pair  $(\frak b,b)$ satisfies weak Maurer-Cartan equation, namely
$$
\sum_{k=0}^{\infty}\sum_{\beta} \exp(\partial\beta \cap b_0)
\frak m^{\frak b}_{\beta,k}(b_+^k) =
\frak{PO}_{\frak b}(b) \cdot \text{\bf e}_L.
$$
Suppose also that $HF((L,b),(L,b);\Lambda) \ne 0$.
Then there exists a non zero $v \in H(X;\Lambda_0)$ such that
\begin{equation}\label{formulaeigen}
c_1(X) \cup^{\frak b} v = \frak{PO}_{\frak b}(b) \cdot v.
\end{equation}
\end{thm}
\begin{rem}
Theorem \ref{c1crit} was mentioned by M. Kontsevich as a conjecture.
It was proved in certain cases by Auroux \cite{auroux} and in \cite{fooo:bulk},  \cite[Theorem 23.12]{spectre}. The result can be generalised in a straightforward way to drop the condition that $\frak b_0$ vanishes on $L$; in that case, we replace the choice of $b_0$ by that of the primitive $\theta_L$.
\end{rem}
\begin{rem}
We remark that the condition $b \in H^1(L;\Lambda_0)$ 
depends on various choices in general.
Note that we have
$
\widehat{\mathcal M}_{\text{\rm weak}}(L;\Lambda_0) 
\subset H^{\rm odd}(L;\Lambda_0),
$
and that the moduli space $\widehat{\mathcal M}_{\text{\rm weak}}(L;\Lambda_0)$
of  bounding cochains is independent of the choices.
Precisely speaking, this independence means the following.
Let $J_1$, $J_2$ be two different choices of compatible almost 
complex structures (we also need to fix Kuranishi structure and CF-perturbation also, but we do not include them in the 
notation). % but this set may depend on them.) 
We then define   $\widehat{\mathcal M}_{\text{\rm weak}}(L;\Lambda_0;J_i)$
as a subset of the cohomology groups of odd degrees. 
Now there exists a map 
$$
\Psi : H^{\rm odd}(L;\Lambda_0) 
\to H^{\rm odd}(L;\Lambda_0) 
$$
that is non-linear and does not preserve degree in general 
such that
$$
\widehat{\mathcal M}_{\text{\rm weak}}(L;\Lambda_0;J_2)
= \Psi(\widehat{\mathcal M}_{\text{\rm weak}}(L;\Lambda_0;J_1)).
$$
Thus the condition $b \in H^1(L;\Lambda_0)$ depends 
on the choices and is not an intrinsic condition.
In the following cases this condition becomes more reasonable:
\begin{enumerate}
\item 
If $\dim_{\C}X = 2$ then 
$H^{\rm odd}(L;\Lambda_0) = H^{1}(L;\Lambda_0)$.
Hence the condition is automatically satisfied.
\item
Let $L$ be a fiber of toric manifold. We proved in 
\cite{fooo:toric1, fooo:bulk} that $H^{1}(L;\Lambda_0)
\subset \widehat{\mathcal M}_{\text{\rm weak}}(L;\Lambda_0)$
in case we take the almost complex structure, the Kuranishi structure, and the multisection to be $T^n$ equivariant.
So the condition holds in this case, if we reformulate the story using singular homology rather 
than de Rham cohomology.
\item
If $L$ is the fixed locus of an anti-symplectic involution, then we can take an involution anti-invariant almost complex structure, and appropriately equivariant Kuranishi structures and multisections.
In this case $b=0$ is a good candidate for the 
choice of $\widehat{\mathcal M}_{\text{\rm weak}}(L;\Lambda_0)$.
\item
If $L$ is monotone again $b=0$ is an element of 
$\widehat{\mathcal M}_{\text{\rm weak}}(L;\Lambda_0)$.
\end{enumerate}
\end{rem}
We deduce Theorem \ref{c1crit} from the following:
\begin{prop}\label{wherec1goes}
If $\frak b\in H^2(X;\Lambda_0)$ and 
$b \in H^1(L;\Lambda_0)$ then 
$$
{\frak q}(c_1(X);1_{(L,b)}) = \frak{PO}_{\frak b}(b) \cdot \text{\bf e}_L.
$$ 
\end{prop}
\begin{proof}[Proof of Proposition \ref{wherec1goes}
$\Rightarrow$ Theorem \ref{c1crit}]
We decompose 
$$
H^*(X;\Lambda) = \bigoplus_{\lambda} V_{\lambda}
$$
where
$$
V_{\lambda} = \{ x \mid 
\exists N, \,\,(c_1(X) - \lambda)^N x = 0 \}.
$$
(Here we use the product structure $\cup^{\frak b}$.)
If $x \in V_{\lambda}$ and $\frak{PO}_{\frak b}(b) \ne \lambda$ then 
by Proposition \ref{wherec1goes} and Theorem \ref{theoremD1} we have
$$
(\frak{PO}_{\frak b}(b) - \lambda)^N {\frak q}(x;1_{(L,b)}) = 0.
$$
Therefore ${\frak q}(x;1_{(L,b)}) = 0$.
Theorem  \ref{c1crit} now follows from Proposition \ref{nontrivialq}.
(The non-triviality of $\hat{\frak q}$.)
\end{proof}
\begin{cor}\label{c1eigen}
Let $\bL =\{(L_i,b_i)\}$ be the set of pairs $(L,b)$ satisfying the conditions 
of Theorem {\rm \ref{c1crit}}. We assume that $\lambda = \frak{PO}_{\frak b}(b_i)$ 
for all $i$.  
Let 
$\text{\bf x} \in  HH_*(\cL)$. Then
\begin{equation}
\hat{\frak p}^{\text{\bf b}}(\text{\bf x}) \cup^{\frak b}
\underbrace
{( c_1(X) - \lambda)
 \cup^{\frak b} \cdots  \cup^{\frak b} ( c_1(X) - \lambda)}_{N} = 0
\end{equation}
for sufficiently large $N$.
\end{cor}
\begin{proof}
We first consider the case:
$\text{\bf x} = x_0 \in HF_*((L,b),(L,b);\Lambda)
\subset HH_*(\cL)$
for some object $(L,b)$ 
and prove
\begin{equation}\label{E1simplecase}
\hat{\frak p}^{\text{\bf b}}(x_0) \cup^{\frak b} c_1(X) = \lambda \hat{\frak p}^{\text{\bf b}}(x_0).
\end{equation}
Let $\frak z \in QH^*(X;\Lambda)$. Then by 
Theorem \ref{pqrelation} we have
$$
\aligned
& \langle 
c_1(X) \cup^{\frak b} \hat{\frak p}^{\text{\bf b}}(x_0),\frak z
\rangle_{\text{\rm PD}_X}
= \langle 
c_1(X), \hat{\frak p}^{\text{\bf b}}(x_0) \cup^{\frak b}\frak z
\rangle_{\text{\rm PD}_X}\\
&=
\langle 
c_1(X), 
\hat{\frak p}^{\text{\bf b}}(\frak m^{\text{\bf b}}(x_0
\otimes 
\frak q^{\text{\bf b}}(\frak z;1_{(L,b)})))
\rangle_{\text{\rm PD}_X}
\\
&=
\langle 
\frak q^{\text{\bf b}}(c_1(X);1_{(L,b)}), \frak m^{\text{\bf b}}(x_0
\otimes 
\frak q^{\text{\bf b}}(\frak z;1_{(L,b)}))
\rangle_{\text{\rm cyc}}
\\
&= \lambda
\langle
\text{\bf e}_L,
\frak m^{\text{\bf b}}(x_0
\otimes 
\frak q^{\text{\bf b}}(\frak z;1_{(L,b)}))
\rangle_{\text{\rm cyc}}
\\
&= \lambda
\langle
x_0,
\frak m^{\text{\bf b}}(\frak q^{\text{\bf b}}(\frak z;1_{(L,b)}),
\text{\bf e}_L)
\rangle_{\text{\rm cyc}}
\\
&=
\lambda
\langle
x_0,
\frak q^{\text{\bf b}}(\frak z;1_{(L,b)})
\rangle_{\rm cyc}
=
\lambda
\langle
\hat{\frak p}^{\text{\bf b}}(x_0),\frak z
\rangle_{\text{\rm PD}_X}.
\endaligned
$$
Note the equality of the second line is a consequence of Theorem \ref{pqrelation}.
The equality of the third line is a consequence of Theorem \ref{thm:duality}.
The equality of the fourth line is a consequence of Proposition \ref{wherec1goes}.
The equality of the fifth line is a consequence of cyclic symmetry.
We have thus proved (\ref{E1simplecase}).
\par
To study the general case we consider the decomposition 
(\ref{decompQH}).
Since $c^1(X)$ has even degree it is in the center 
of $QH_{\frak b}(X;\Lambda)$.
Therefore each space
$$
QH_{\frak b}(X;\Lambda)_{\rho} = 
\{ x \in QH_{\frak b}(X;\Lambda) \mid \exists N, \,\, (c_1(X) - \rho)^N x = 0\}
$$
is a direct product factor of $QH_{\frak b}(X;\Lambda)$.
Therefore each direct factor $QH_{\frak b}(X;\Lambda)_a$ 
in the decomposition (\ref{decompQH}) is contained in one of 
the generalized eigenspaces $QH_{\frak b}(X;\Lambda)_{\rho}$.
\par
We decompose $D^{\pi}(\cL)$ into the product of 
$D^{\pi}(\cL;a)$ as in (\ref{catedecompose}).
Then Lemma \ref{ppreservesa} implies that
$$
\hat{\frak p}^{\text{\bf b}}(D^{\pi}(\cL;a))
\subset  
QH_{\frak b}(X;\Lambda)_{a}.
$$
Therefore 
$$
\hat{\frak p}^{\text{\bf b}}(D^{\pi}(\cL;a))
\subset   
QH_{\frak b}(X;\Lambda)_{\rho},
$$
for some $\rho$. By (\ref{E1simplecase})
 $\rho = \lambda$ if $D^{\pi}(\cL;a) 
\ne 0$.
The proof of Corollary \ref{c1eigen} is complete.
\end{proof}
We remark that Corollary \ref{c1eigen} 
implies Theorem \ref{maintheorem1} in case 
the assumption of Theorem {\rm \ref{c1crit}} is satisfied.
\begin{proof}[Proof of Proposition \ref{wherec1goes}]
We first remark that there exists a cycle $Q$ of $X \setminus L$ 
such that 
\begin{equation}\label{condE3}
\aligned
&\beta \cap Q = \mu(\beta), \qquad \beta \in \pi_2(X,L), \\
&PD([Q]) = 2c_1(X) \in H^2(X;\Z).
\endaligned
\end{equation}
See \cite[Lemma 23.11]{spectre}.
We regard $Q$ as a current and approximate it by a smooth differential $2$ form $g$
whose support is disjoint from $L$.
\par
Let $g_1,\dots,g_{B_2}$ be smooth forms which represent a basis of 
$H^2(X;\R)$.
Then we have
$$
\frak b = \frak b_2 + \sum_{i=1}^{B_2} w_i g_i
$$
where $\frak b_2 \in H^2(X;\C)$ and $w_i \in \Lambda_+$.
\par
Consider the moduli space 
${\mathcal M}_{\ell;k+1}(L;\beta)$ which we used to define the operator $\frak q$ and  the forgetful map
\begin{equation}\label{forget1}
\frak{forget} : {\mathcal M}_{\ell+1;k+1}(L;\beta) \to {\mathcal M}_{\ell;k+1}(L;\beta)
\end{equation}
forgetting the first interior marked point $z_0^+$.
In Section \ref{sec:cyclicfil} we introduced a system of Kuranishi structures 
and \blue{CF-perturbations} on ${\mathcal M}_{\ell;k+1}(L;\beta)$.
We consider another system of
Kuranishi structures 
and \blue{CF-perturbations} on 
${\mathcal M}_{\ell+1;k+1}(L;\beta)$ that is induced from the first 
one by the forgetful map (\ref{forget1}).
This second system of Kuranishi structures and CF-perturbations are {\it not} 
invariant under the permutation of interior marked points.
(It is invariant under the permutations of 2nd-$(\ell+1)$st interior marked points.)
When we use this 2nd system of Kuranishi structures we write 
${\mathcal M}_{\ell+1;k+1}(L;\beta)_{\text{\rm II}}$ 
in place of ${\mathcal M}_{\ell+1;k+1}(L;\beta)$.
\par
Let  
$
b = \theta + b_+
$, where
$\theta $ is a complex valued one form with $d\theta = \frak b_2$ and 
$b_+$ is a representative of  $H^1(L;\Lambda_+)$.
Let $\text{\bf h} = (h_1,\dots,h_k)$, $h_i \in \Omega(L)$. 
We put
$$
\aligned
{\frak q}^{\rm form,(\frak b,b_0)}_{\beta}(g;\text{\bf h}) 
=
&\sum_{i_1,\dots,i_B=0}^{\infty} \sum_k
(-1)^{\maltese}
\frac{\rho_{\frak b,\theta}(\beta)}%\rho_{b}(\beta)}
{i_1!\cdots i_{B}!}
w_1^{i_1}\cdots w_B^{i_B}\\
&\text{\rm Corr}({\mathcal M}_{\ell+1;k+1}(L;\beta);{\rm ev}^+,({\rm ev}_1,\dots,{\rm ev}_k);{\rm ev}_0)
(g_1^{i_1}\cdots g_B^{i_B};\text{\bf h}).
\endaligned
$$
(where $\maltese$ is the same as (\ref{form21ten1}).)
We define 
${\frak q}^{\rm form,(\frak b,b_0)}_{\beta;\text{\rm II}}(g;\text{\bf h}) $
to be the one obtained by using ${\mathcal M}_{\ell+1;k+1}(L;\beta)_{\text{\rm II}}$  in place of 
${\mathcal M}_{\ell+1;k+1}(L;\beta)$  in the above formula.
\par
We also define
$$
\aligned
\frak q^{\rm form,(\frak b,b_0)}_{\beta}(1_X;\text{\bf h}) 
=
&\sum_{i_1,\dots,i_B=0}^{\infty} \sum_k(-1)^{\maltese}
\frac{\rho_{\frak b,\theta}(\beta)}%\rho_{b}(\beta)}
{i_1!\cdots i_{B}!}
w_1^{i_1}\cdots w_B^{i_B}\\
&\text{\rm Corr}({\mathcal M}_{\ell+1;k+1}(L;\beta);{\rm ev}^+,({\rm ev}_1,\dots,{\rm ev}_k);{\rm ev}_0)
(g_1^{i_1}\cdots g_B^{i_B};\text{\bf h}).
\endaligned
$$
Here $1_X \in E_0(\Omega(X)[2])$. (It is different from $\text{\bf e}_X$.)
Using the operators
$\frak q^{\rm form,(\frak b,b_0)}_{\beta_i}(g;\text{\bf h})$,
$\frak q^{\rm form,(\frak b,b_0)}_{\beta_i}(1_X;\text{\bf h})$
and  the Green kernel $G_L$ in the way explained in Subsection \ref{frakqoncoho},
we obtain
$\hat{\frak q}^{\text{\rm can},\text{\bf b}}_{\beta}(g) $.
(Here the sum of $\beta_i$ associated to each of the vertex in the tree 
appearing in the construction of subsection \ref{frakqoncoho} is $\beta$ 
and we put $b_+$ on each of the exterior vertices.)
\par
We use 
$\frak q^{\rm form,(\frak b,b_0)}_{\beta_i;\text{\rm II}}(g;\text{\bf h}) $
and 
$\frak q^{\rm form,(\frak b,b_0)}_{\beta_i}(1_X;\text{\bf h}) $
instead to define 
$\hat{\frak q}^{\text{\bf b}}_{\beta;\text{\rm II}}(g) $.
\par
Since $b$ is a  bounding cochain (in the canonical model) we have
$$
\sum_{\beta}T^{\omega(\beta)} \frak q^{\text{\bf b}}_{\beta}(1_X)
= \frak{PO}_{\frak b}(b) \cdot \text{\bf e}_L.
$$
\begin{rem}
Note that we assumed this formula only for the canonical model.
\end{rem}
Now we have
\begin{lem}\label{gIIand1}
$$
\sum_{\beta}T^{\omega(\beta)} \hat{\frak q}^{\text{\bf b}}_{\beta;\text{\rm II}}(g)
= \frak{PO}_{\frak b}(b)\cdot \text{\bf e}_L.
$$
\end{lem}
\begin{proof}
We consider a decorated 
ribbon tree $\Gamma = (\mathcal T,\beta(\cdot))$ 
used to define $\frak q^{\text{\rm can},\text{\bf b}}_{\beta}$.
Note that we have only the data $\mathcal T$ and $\beta(\cdot)$ since we have only 
one Lagrangian submanifold.
We use $\frak q^{\rm form,(\frak b,b_0)}_{\beta(v)}(1_X;\cdot)$ at
each of the vertices to define ${\frak q}^{\text{\rm can},(\frak b,b_0)}_{\frak b}$. We replace it by 
${\frak q}^{\rm form,(\frak b,b_0)}_{\beta(v);\text{\rm II}}(g;\cdot)$ at 
one of the vertices. We sum them up for various $\Gamma$ and their vertices. We thus obtain 
${\frak q}^{\text{\bf b}}_{\beta;\text{\rm II}}(g)$.
Then since $\sum \mu_L(\beta(v)) = 2$ (here $\mu_L$ is the Maslov index), 
the second formula of (\ref{condE3}) together with the choice of 
our perturbation to define $ \frak q^{\rm form,\text{\bf b}}_{\beta;\text{\rm II}}$ 
implies 
$$
\hat{\frak q}^{\text{\bf b}}_{\beta;\text{\rm II}}(g)
 =  
\hat{\frak q}^{\text{\bf b}}_{\beta}(1_X).
$$
Note that $\sum \mu(\beta(v)) = 2$ follows from degree counting.
Namely since the output is supported in degree $0$,
$b$ has degree $1$ and $\frak b$ has degree $2$, the Maslov 
index should be $2$.
In other words we use the fact that $b$ is a  bounding cochain 
and all the outputs other than those of degree $0$ must cancel each other.
We remark that this cancellation continues to hold for 
$ \hat{\frak q}^{\text{\rm can},\text{\bf b}}_{\beta;\text{\rm II}}(g)$.
In fact, the ratio of 
$\hat{\frak q}^{\text{\rm can},\text{\bf b}}_{\beta}(1_X)$ and 
$\hat{\frak q}^{\text{\rm can},\text{\bf b}}_{\beta;\text{\rm II}}(g)$, 
for any class $\beta$, is proportional to the Maslov index of $\beta$.
So since the cancellation occurs in each degree for 
$\hat{\frak q}^{\text{\bf b}}_{\beta}(1_X)$, 
it occurs also for $\hat{\frak q}^{\text{\bf b}}_{\beta;\text{\rm II}}(g)$.
\par
The proof of Lemma \ref{gIIand1} is complete.
\end{proof}
\begin{lem}
$$
\sum_{\beta}T^{\omega(\beta)} \hat{\frak q}^{\text{\bf b}}_{\beta;\text{\rm II}}(g)
=
\sum_{\beta}T^{\omega(\beta)}\hat{\frak q}^{\text{\bf b}}_{\beta}(g).
$$
\end{lem}
\begin{proof}
The argument is similar to, for example, the proof of Proposition \ref{dualityinformlevel}.
The only difference between 
$\hat{\frak q}^{\text{\bf b}}_{\beta;\text{\rm II}}$ and 
$\hat{\frak q}^{\text{\bf b}}_{\beta}$ is that we use 
different perturbation at one of the disks appearing in the definition.
We can find a homotopy between two perturbations.
(The perturbation we put on the `bubbles' are the same.)
Using it we can show that the left side is homologous to the right side.
Since we are working in the canonical model, the output is always a harmonic form.
So homologous means equal.
\end{proof}
The proof of Proposition \ref{wherec1goes} and Theorem \ref{c1crit} is now complete.
\end{proof}

\section{Quantitative version and its application.}
\label{sec: quatitative}

\subsection{Mod $T^E$ version of 
Theorem \ref{maintheorem1}}
\label{subsecmodTEstatement}

Let $\bL$ be a finite collection of 
$(L_{\kappa},\theta_{\kappa},b_{\kappa})$ where
$L_{\kappa}$ is a $\text{\rm st}$-relatively 
spin Lagrangian submanifolds and 
$b_{\kappa} \in H^{\rm odd}(L_{\kappa};\Lambda_0)$
such that
\begin{equation}
\sum_{k=0}^{\infty} 
\frak m_k^{(\frak b,0)}(b_{\kappa}^{k}) \equiv 
\lambda_{\kappa} \cdot \text{\bf e}_{L_{\kappa}}
\mod T^{E_0}. 
\end{equation}
We define
$\frak{PO}(b_{\kappa}) = \lambda_{\kappa} \in \Lambda_0/T^{E_0}\Lambda_0$

From such a collection we defined in Section \ref{sec:cyclicfil} a unital and cyclic filtered $A_{\infty}$
category mod $T^{E_0}$ for each $E_0 \ge 0$, that we write $\cL^{E_0}$.
We also defined in Subsection \ref{frakponcoho}
a map\index[syindex]{pfrakhatE0@$\hat{\frak p}^{E_0}$}
\begin{equation}\label{pmodTE}
\hat{\frak p}^{E_0}: 
HH_*(\cL^{E_0};\Lambda_0/T^{E_0}\Lambda_0) \to H(X;\Lambda_0/T^{E_0}\Lambda_0).
\end{equation}

\begin{thm}\label{TEversion}
Assume that the image of $(\ref{pmodTE})$ contains $T^{\rho}\text{\bf e}_{X}$.
Let $(U,\theta_U,b_U)$ be another pair consisting of a $\text{\rm st}$-relatively 
spin Lagrangian submanifold $U$ and a class $b_{U} \in H^{\rm odd}(U;\Lambda_0)$
such that
\begin{equation}\label{condbU1}
\sum_{k=0}^{\infty} 
\frak m_k^{(\frak b,0)}(b^{k}_U) \equiv 
\lambda^{U} \cdot \text{\bf e}_U
\mod T^{E_0}. 
\end{equation}
Assume in addition that there exists $E_1 < E_0$ such that
$$
T^{E_1}HF((U,b_{U}),(U,b_{U});\Lambda_{0}/T^{E_0}\Lambda_0)
\ne 0.
$$
Then there exists $\kappa$ such that
$$
\lambda^U \equiv \frak{PO}(b_{\kappa}) \mod T^{E_1-\rho}.
$$
\end{thm}

\begin{proof}
We consider the following variant of Diagram \ref{eq:diagram_factor_p_q}.
  \begin{equation} \label{eq:diagram_factor_p_q2}
    \xymatrix{ HH_{*-n} ( \cL,\cL) \otimes \Lambda_0/T^{E_0} \ar[r]^{\fphat} \ar[d] & QH^{*}(X; \Lambda_0/T^E)  \ar[d]^{\fqhat} \\
HH^{*-1}\left(\sU, {_{\sU}}\Fuk _{\cL} \displaystyle{ \tensor_{\cL}} {_\cL} \Fuk_{\sU} \right)\otimes \Lambda_0/T^{E_0}  \ar[r] & HH^{*} ( \sU,\sU)\otimes \Lambda_0/T^{E_0}.}
  \end{equation}
Here $\mathcal U$ is the filtered $A_{\infty}$ category with one object $(U,b_{U})$.
By assumption there must be some value $\lambda$ such that
\begin{equation}
T^{E_1 - \rho}\fq \circ \fphat \co  HH_*(\cL_{\lambda}, \cL_{\lambda}) \otimes \Lambda_0/T^{E_0} \to HF((U;(\frak b,b_{U})),(U;(\frak b,b_{U}));\Lambda_0/T^{E_0})
\end{equation}
does not vanish.
\par
We put $r = \lambda^U - \frak{PO}_{\frak b}(b_{\kappa})$. Then 
$$
r \cdot HH^{*-1}\left(\sU, {_{\sU}}\Fuk _{\cL} \displaystyle{ \tensor_{\cL}} {_\cL} \Fuk_{\sU} \right)\otimes \Lambda_0/T^{E_0}
= 0.
$$
Therefore by the commutativity of (\ref{eq:diagram_factor_p_q2}) 
$r \equiv 0 \mod T^{E_1-\rho}$ as required.
\end{proof} %\marginpar{Proof added. KF}

\subsection{An application}
\label{subsecmodTEappli}

In this subsection we demonstrate how we can apply Theorem \ref{TEversion} in an example.
We consider the case $X = \C P^n$ in this subsection. We can 
apply the same argument for any compact toric $X$.
We normalize the symplectic form $\omega$ so that 
$\int_{\C P^1}\omega = 1$.

\begin{prop}
Let $U$ be a relatively spin Lagrangian submanifold in $\C P^n$.
Then there exists a non-constant holomorphic map 
$
u : (D^2,\partial D^2) \to (\C P^n,U)
$
such that 
\begin{equation}\label{bcond1}
\mu_U([u]) \le n+1, \qquad 
\int_{D^2} u^* \omega \le 1.
\end{equation}
Here $\mu_U : H_2(\C P^n,U) \to \Z$ is the Maslov index.
\end{prop}

\begin{proof}
Let $L_0 \subset \C P^n$ be the Clifford torus. We take $\frak b =0$.
We put $\chi_k = \exp(2\pi \sqrt{-1}k/(n+1))$, $k=0,\dots,n$
$$
b_k = (\log \chi_k,\dots,\log \chi_k)\in H^1(L_0;\C).
$$
As was shown by Cho (see \cite[Section 1.4]{toric3}), the class
$b_k$ lies in $\mathcal M_{\text{\rm weak}}(\C P^n;\Lambda_0)$
and the corresponding self-Floer homology group is isomorphic to ordinary cohomology
$$
HF((L_0,b_k),(L_0,b_k);\Lambda_0)
\cong H(T^n;\Lambda_0).
$$
Moreover $\frak{PO}(b_k) = T^{1/(n+1)}(n+1)\chi_k$.
\par
The quantum cohomology $QH(\C P^n;\Lambda)$
splits to $n+1$ copies of $\Lambda$.
We write $e_k$ ($k=0,\dots,n$) for the unit of $k$-th factor.
We may choose them so that
$$
{\frak q}(e_k;1_{L_0}) \ne 0 \in HF((L_0,b_k),(L_0,b_k);\Lambda).
$$
(This is well-known and follows for example from Theorem \ref{ssfinite}.)
Let $\text{\rm vol}_{(L_0,b_k)}$ be the 
generator of $H^n(L_0;\Z)$ regarded  as an element of 
$HF((L_0,b_k),(L_0,b_k);\Lambda)$.
By Theorem \ref{thm:duality} we have
$\hat{\frak p}(\text{\rm vol}_{(L_0,b_k)}) = c_k e_k$.
Then
$$
\aligned
&1 = \langle {\frak q}(e_k;1_{L_0}),\text{\rm vol}_{(L_0,b_k)}\rangle
= \langle e_k,\hat{\frak p}(\text{\rm vol}_{(L_0,b_k)})\rangle_{\text{\rm PD}_X}
= c_k\langle e_k,e_k\rangle_{\text{\rm PD}_X}, \\
&
c_k =c_k^2\langle e_k,e_k \rangle_{\text{\rm PD}_X}= \langle \hat{\frak p}(\text{\rm vol}_{(L_0,b_k)}),
\hat{\frak p}(\text{\rm vol}_{(L_0,b_k)}) \rangle_{\text{\rm PD}_X}.
\endaligned$$
By \cite[Theorems 1.3.25 and 3.4.1]{toric3} we have
\begin{equation}\label{hessianapper}
\langle \hat{\frak p}(\text{\rm vol}_{(L_0,b_k)}),
\hat{\frak p}(\text{\rm vol}_{(L_0,b_k)}) \rangle_{\text{\rm PD}_X}
= \text{\rm det}\left(y_iy_j \frac{\partial^2\frak{PO}}{\partial y_i\partial y_j}\right)
(\frak y_k)
\end{equation}
where
$
\frak{PO} = y_1+\dots+y_k + \frac{T}{y_1\cdots y_k}
$
and 
$
\frak y_k = T^{1/(n+1)}(\chi_k,\dots,\chi_k).
$
We have $(\ref{hessianapper}) = (n+1)T^{n/(n+1)}\chi_k$.
(see \cite[Section 1.4]{toric3}).
Therefore
$$
T^{n/(n+1)}\text{\bf e}_{X} \in \text{\rm Im}(\hat{\frak p}
: HH_*(\cL;\Lambda_0) \to H(\C P^n;\Lambda_0))
$$
here $\bL =\{(L_0,b_k)\mid k=0,\dots,n\}$.
\par
Now suppose that there does not exist a map $u$ satisfying (\ref{bcond1}).
Then $b_U = 0$  satisfies (\ref{condbU1}) for $E_0 > 1$ and $\lambda^U =0$.
(We can show that any holomorphic disk with Maslov index $> 2$ does not 
contribute to $\frak m_0$ by a dimension counting argument. See \cite[Subsection 3.6.3]{fooo09}
for a related argument.)
\par
Again by dimension counting the holomorphic disk with Maslov index $> n+1$ does 
not contribute to the differential of the spectral sequence of \cite[Theorem B]{fooo09}.
Therefore 
$$
HF((U,b^{U}),(U,b^{U});\Lambda_{0}/T^{E_0}\Lambda_0)
= H(U;\Lambda_{0}/T^{E_0}\Lambda_0).
$$
It follows that there exists $b_k$ such that
$$
\frak{PO}(b_k) \equiv \frak{PO}(b^{U}) \mod T^{E_0 - n/(n+1)}
$$
Since $\frak{PO}(b_k)$ is not congruent to $0$ modulo $T^{1/(n+1)}$, 
this is a contradiction.
\end{proof}

%\appendix

\printindex[syindex]
\printindex
%\section{Sketch of the proof of Proposition
%\ref{spectreestimate}
%\label{sec:spectreestimate}

%\input{reference.tex}
\bibliographystyle{amsalpha}

\end{document}